\documentclass[reqno,11pt]{amsart} 
\usepackage{amsmath,amsfonts,amsthm,amscd,amssymb,graphicx}
\usepackage{subfigure}

\usepackage{soul} 

\usepackage{kotex}
\usepackage{csquotes}

  \usepackage[margin=
1
 in]{geometry}

\usepackage{comment}

\usepackage[makeroom]{cancel}

\usepackage[english]{babel}

\addto\captionsenglish{
	\renewcommand{\contentsname}%
	{Contents}%
}

\newtheorem*{theorem*}{Theorem}

\usepackage{CJKutf8}

\usepackage{stmaryrd}

\usepackage{cancel,accents,todonotes}

\numberwithin{equation}{section}

\usepackage{hyperref}
\usepackage{verbatim}
\usepackage{color}

\usepackage{stackengine,scalerel}

\newcommand\obullet[1]{\ThisStyle{\ensurestackMath{%
			\stackon[1pt]{\SavedStyle#1}{\SavedStyle\kern.6\LMpt\bullet}}}}
\newcommand\ocirc[1]{\ThisStyle{\ensurestackMath{%
			\stackon[1pt]{\SavedStyle#1}{\SavedStyle\kern.6\LMpt\circ}}}}

\newtheorem{theorem}{Theorem}[section]

\newtheorem{proposition}[theorem]{Proposition}

\newtheorem{corollary}[theorem]{Corollary}

\newtheorem{assumption}[theorem]{Assumption}
\newtheorem{definition}[theorem]{Definition}

\newtheorem{lemma}{Lemma}[section]

\theoremstyle{remark}
\newtheorem{remark}{Remark}[section]

\theoremstyle{plain}

\let\hide\iffalse
\let\unhide\fi

\newcommand{\bk}{\mbox{\boldmath $k$}}
\newcommand{\bl}{\mbox{\boldmath $l$}}

\def\eps{\varepsilon }

\renewcommand{\div}{{\rm div}}

\newcommand{\pt}{\partial}

\newcommand{\p}{\partial}

\newcommand{\II}{{\mathbb I}}

\newcommand{\e}{{\varepsilon}}

\newcommand{\Be}{\begin{equation}}
	\newcommand{\Ee}{\end{equation}}

\newcommand{\AC}{\accentset{\circ}}

\newcommand{\TbT}{\mathcal T \! \!\!\!\ b\! \!\mathcal T}

\def\beq{\begin{equation}}
\def\eeq{\end{equation}}
\def\bb1{{1\!\!1}}

\def\bW{{\bf W}}

\def\bPhi{{\bf \Phi}}
\def\bw{\boldsymbol{\omega}}

\def\R{\mbox{Re }}

\def\w{{\omega}}

\def\p{\partial}

\def\les{\lesssim}

\newcommand{\T}{{\mathbb T}}

\def\diag{\mbox{diag}}

\def\eps{\varepsilon}

\def\triangle{\Delta}
\def\bega{\begin{aligned}}
\def\enda{\end{aligned}}
\def\R{\mathbb{R}^2}

\def\lw{\left}
\def\rw{\right}

\def\R{\mathbb{R}}

\def\la{\langle}
\def\ra{\rangle}

\def\Z{\mathbb{Z}}

\def\ta{\theta}

\def\bcase{\begin{cases}}
\def\ecase{\end{cases}}
\def\al{\alpha}
\def\bmx{\begin{bmatrix}}
\def\emx{\end{bmatrix}}

\def\P{\mathbf{P}}
\def\II{\mathbf{I}}

\makeatletter
\newcommand{\StartNoTOC}{%
  \begingroup
  \let\addcontentsline\@gobblethree
  \let\addtocontents\@gobbletwo
}
\newcommand{\StopNoTOC}{\endgroup}
\makeatother

\begin{document}

\title[Quantitative Closure analysis]{Quantitative Closure Analysis toward Ideal Fluids  \\[1.5em]{\it{\small{In Honor of Professor Yan Guo’s 60th Birthday}}}
 }


\author{Gi-Chan Bae}
\address{Department of Mathematical Sciences, Seoul National University, Seoul 08826, Korea.}
\email{gcbae02@snu.ac.kr}

\author{Chanwoo Kim}
\address{Department of Mathematics, University of Wisconsin-Madison, Madison, WI, USA, 53717}
\email{chanwookim.math@gmail.com}

 \thanks{Authors’ comment: This arXiv preprint is the detailed companion to the compact journal-submission version of the same title and is posted on arXiv.org for readers' convenience. This version is not intended for separate journal publication.}


\maketitle

\begin{abstract}
We establish the incompressible low--Mach/high--Reynolds limit for the Boltzmann equation for a broad class of initial data, without recourse to any asymptotic expansion. Exploiting the local Maxwellian manifold and the macro--micro decomposition in a new quasi-linear analysis, we derive quantitative estimates for the purely microscopic fluctuation, as well as bounds for the kinetic vorticity and the entropic fluctuation in terms of the initial data. As a consequence, in two space dimensions, the rescaled velocity and temperature converge to a global solution of the incompressible Euler equations coupled to a transported temperature, within the frameworks of DiPerna--Lions--Majda and Delort.



  

\hide

Understanding the specifics of dimension 2 in the hydrodynamic limit toward the incompressible Euler equation, without any prior knowledge of the limiting solution, has been an important open question in the theory of Boltzmann equations. In this paper, we establish that solutions of the scaled Boltzmann equations converge to renormalized solutions of the incompressible Euler equations as the Knudsen number tends to zero, with the convergence holding globally in time. Our results apply under the conditions that the initial Boltzmann vorticity is uniformly in $L^p(\mathbb{T}^2)$ and the entropy remains uniformly bounded. Additionally, we show that the limiting Euler solution conserves energy.

A key component of our proof is a novel weak stability result for solutions of the Boltzmann equations, which remains valid even when the gradients of momentum and temperature are unbounded. Specifically, we analyze the vorticity at the microscopic level through the dissipative momentum equation, successfully decoupling the nonlinear macro-macro interactions in the process.



Understanding the specificity of dimension 2 in the hydrodynamic limit toward the incompressible Euler equation \textit{without} any prior information about the limiting solution has been an important open question in the theory of Boltzmann equations. In this paper, we prove that the solutions of the scaled Boltzmann equations converge to renormalized solutions of the incompressible Euler equations as the Knudsen number vanishes, globally in time. Our convergence holds when the initial Boltzmann vorticity is uniformly in $L^p(\mathbb{T}^2)$ and the entropy is uniformly bounded. Additionally, we establish that this limiting Euler solution conserves energy.

A crucial part of our proof involves a new weak stability result for the solutions to the Boltzmann equations, even when the gradients of momentum and temperature can be unbounded. In particular, we examine the vorticity at the microscopic level via the dissipative momentum equation, for which we are able to construct an energy-dissipation pair to effectively ``decouple'' the nonlinear macro-macro interaction.


It has been an important open question in the theory of Boltzmann equations to understand the specificity of dimension 2 in the incompressible inviscid regime without using any expansions (such as those by Hilbert or Chapman-Enskog). In this paper, we prove that the solutions of the scaled Boltzmann equations converge to a renormalized solution of the incompressible Euler equations as the Knudsen number $\e$ vanishes, when the initial Boltzmann vorticity belongs to $L^p(\mathbb{T}^2)$ uniformly-in-$\e$. We also establish that our limiting Euler solution conserves energy.


As a key of the proof, we establish a new weak stability result for the solutions to the Boltzmann equations, even when the gradients of momentum and temperature may be unbounded. In particular, we understand the vorticity at the microscopic level throughout the dissipative momentum equation, and we are able to construct an effective energy-dissipation pair for the whole system. 
 .


\unhide

\hide
In this paper, we rigorously prove that, without employing any expansions (such as those by Hilbert or Chapman-Enskog), the solutions of the scaled Boltzmann equations converge to a renormalized solution (\`a la Diperna-Lions) of the incompressible Euler equations provided the vorticity is contained in $L^p(\mathbb{T}^2)$ as the Knudsen number $\e \downarrow 0$. We also prove that this limiting solution conserves the energy: $\|u(t)\|_{L^2_x} = \|u(0)\|_{L^2_x}$. Our proof is based one a new weak stability result of the solutions to the Boltzmann equations even when the gradient of the momentum and temperature can blow up in the limit, e.g. $ \| \nabla_x  u^\e \|_{L^1_tL^\infty_x} \rightarrow \infty$ as $\e \downarrow 0$. In particular, we realize the vorticity in microscopic level throughout the dissipative momentum equation, and therefore can establish full dissipation for both macroscopics and macroscopics.

In this paper, we rigorously demonstrate that, without using any expansions (such as those by Hilbert or Chapman-Enskog), the solutions of the scaled Boltzmann equations converge to a renormalized solution (à la Diperna-Lions) of the incompressible Euler equations, provided that the vorticity is contained in \( L^p(\mathbb{T}^2) \) as the Knudsen number \( \varepsilon \) approaches 0.

We also establish that this limiting solution conserves energy, meaning \( \|u(t)\|_{L^2_x} = \|u(0)\|_{L^2_x} \). Our proof relies on a new weak stability result for the solutions to the Boltzmann equations, even when the gradients of momentum and temperature may blow up in the limit, for example, when \( \| \nabla_x  u^\varepsilon \|_{L^1_t L^\infty_x} \) approaches infinity as \( \varepsilon \) tends to 0. 

In particular, we address the vorticity at the microscopic level throughout the dissipative momentum equation, which allows us to establish full dissipation for both macroscopic and microscopic scales.

In particular, we address the vorticity at the microscopic level throughout the dissipative momentum equation, which allows us to establish full dissipation for both macroscopic and microscopic scales.

The specificity of dimension 2 in the incompressible inviscid regime has been an important open question in the theory of Boltzmann equations. In this paper, we rigorously demonstrate that, without using any expansions (such as those by Hilbert or Chapman-Enskog), the solutions of the scaled Boltzmann equations converge to a renormalized solution (à la Diperna-Lions) of the incompressible Euler equations, provided that the vorticity is contained in \( L^p(\mathbb{T}^2) \) as the Knudsen number \( \varepsilon \) approaches 0. 

We also establish that this limiting solution conserves energy, meaning $\|u(t)\|_{L^2_x} = \|u(0)\|_{L^2_x}$.

\unhide


\end{abstract}
 

\tableofcontents

 \section*{Introduction}
 
 Kinetic equations can be viewed as an \textit{infinite moment hierarchy}: multiplying by test functions in the
velocity variable and integrating produces an unclosed cascade of macroscopic relations, whereas fluid
descriptions close at the level of finitely many macroscopic fields.\footnote{\small\itshape
This viewpoint is classical; see, for example, Grad's \emph{method of moments} of \cite{Gr49}, which truncates the infinite
hierarchy to finitely many moments and seeks constitutive relations (closure laws) justified by near local
equilibrium.}
In this work we focus on the low--Mach/high--Reynolds regime and study a \emph{singularly perturbed family}
$F^\varepsilon=F^\e(t,x,v)\ge 0$ governed by the scaled Boltzmann equation for hard spheres:
\begin{equation}\label{BE}
\varepsilon\,\partial_t F^\varepsilon + v\cdot\nabla_x F^\varepsilon
=
\frac{1}{\kappa \varepsilon}\,\mathcal{N}(F^\varepsilon,F^\varepsilon),
\qquad \kappa=\kappa(\varepsilon)\to 0 \ \text{as } \varepsilon\to 0.
\end{equation}
Here the Mach number $\varepsilon$ measures the smallness of the macroscopic flow speed relative to the thermal
speed, while the collision frequency is of order $1/(\kappa\varepsilon)$ (equivalently, the Knudsen number is
$\mathrm{Kn}=\kappa\varepsilon$).

A fundamental obstruction in the kinetic-to-fluid passage is \emph{closure}: the five moments determine
$(\mathrm P^\varepsilon,\mathrm U^\varepsilon,\Theta^\varepsilon)$, but the stress and heat flux involve higher
velocity moments and hence genuinely microscopic information. We define the macroscopic fields by the five moments
\begin{equation}\label{eq:macro_fields_intro}
\int_{\mathbb{R}^3} F^\varepsilon
\begin{pmatrix}
1\\ v\\ |v|^2/2
\end{pmatrix}\,dv
=
\begin{pmatrix}
\mathrm P^\varepsilon\\ \mathrm P^\varepsilon \mathrm U^\varepsilon\\
\tfrac{3k_B}{2}\mathrm P^\varepsilon\Theta^\varepsilon + \tfrac{1}{2}\mathrm P^\varepsilon|\mathrm U^\varepsilon|^2
\end{pmatrix}.
\end{equation}In the low--Mach regime, it is convenient to parametrize these fields by logarithmic fluctuations,
\begin{equation}\label{expform}
\mathrm{P}^{\e} = e^{\e \rho^{\e}}, \qquad
\mathrm{U}^{\e}= \e u^{\e}, \qquad
\mathrm{\Theta}^{\e} = e^{\e \ta^{\e}}.
\end{equation}

The collision operator preserves the collision invariants:
\begin{equation}\label{eq:collision_invariants_intro}
\int_{\mathbb{R}^3}\mathcal{N}(F^\e,F^\e)(v)\,(1,\ v,\ |v|^2/2)^T\,dv =(0,0,0)^T.
\end{equation}
These collision invariants yield exact balance laws:
\begin{align}\label{loccon}
\bega
 \e\p_t \mathrm{P}^{\e} + \nabla_x\cdot(\mathrm{P}^{\e}\mathrm{U}^{\e})=0, \cr
 \e\p_t(\mathrm{P}^{\e}\mathrm{U}^{\e})
+\nabla_x\cdot\big(\mathrm{P}^{\e}\,\mathrm{U}^{\e}\otimes\mathrm{U}^{\e}+k_B \mathrm{P}^{\e} \Theta^{\e}I+\mathbf{r}^{\e}\big)=0,
\cr
 \e\p_t \mathrm{E}^{\e}
+\nabla_x\cdot\big(\mathrm{U}^{\e}\,\mathrm{E}^{\e}+\big(k_B \mathrm{P}^{\e} \Theta^{\e}I+\mathbf{r}^{\e}\big)\mathrm{U}^{\e}+ \mathfrak{q}^{\e}\big)=0,
\enda
\end{align}
where $\mathrm{E}^{\varepsilon}:=\frac{3}{2}k_B\,\mathrm{P}^{\varepsilon}\mathrm{\Theta}^{\varepsilon}
+\frac{1}{2}\mathrm{P}^{\varepsilon}|\mathrm{U}^{\varepsilon}|^2$. The system \eqref{loccon} is not closed: the
missing fluxes are encoded in the Burnett functionals $\mathbf{r}^{\varepsilon}=\langle F^{\e},  (v-\mathrm{U}^{\e})\otimes(v-\mathrm{U}^{\e})-\frac{|v-\mathrm{U}^{\e} |^2}{3} I \rangle$ (stress deviator) and
$\mathfrak{q}^{\varepsilon}=\langle F^{\e}, (v-\mathrm{U}^{\e})\frac{(|v-\mathrm{U}^{\e}|^2-5k_B\mathrm{\Theta}^{\e})}{2}   \rangle$ (heat flux).

Collisions provide the mechanism that can enable an effective closure: they single out the \emph{local Maxwellian
manifold} as the set of collision equilibria. Given $(\mathrm P^\varepsilon,\mathrm U^\varepsilon,\Theta^\varepsilon)$,
we denote by $\mathbf P F^\varepsilon$ or $M^\varepsilon$ the corresponding local Maxwellian,
\begin{equation}\label{M-def}
\mathbf P F^\e =
M^\varepsilon 
=
\frac{\mathrm P^\varepsilon(t,x)}{\big(2\pi k_B \Theta^\varepsilon(t,x)\big)^{3/2}}
\exp\!\Big(-\frac{|v-\mathrm U^\varepsilon(t,x)|^2}{2k_B\Theta^\varepsilon(t,x)}\Big),
\qquad
\mathcal{N}(\mathbf P F^\e,\mathbf P F^\e)=0.
\end{equation}
This motivates the macro--micro decomposition
\begin{equation}\label{eq:macro_micro_intro}
F^\varepsilon = \mathbf P F^\e + \AC{\mathbf P} F^\varepsilon,
\qquad
\int_{\mathbb{R}^3} \AC{\mathbf P} F^\varepsilon \,(1,\ v,\ |v|^2/2)^T\,dv  =(0,0,0)^T,
\end{equation}
which decomposes $F^\varepsilon$ into its component on the local Maxwellian manifold and a transverse microscopic fluctuation.
The resulting non-closure (stress and heat flux) is carried by the microscopic fluctuation $\AC{\mathbf P}F^\varepsilon$,
so quantitative control of $\AC{\mathbf P}F^\varepsilon$ through collisional relaxation provides the natural
closure mechanism.

\hide
This motivates the macro--micro decomposition
\begin{equation}\label{eq:macro_micro_intro}
F^\varepsilon = \mathbf P F^\e   
+  \AC{\mathbf P} F^\varepsilon,
\qquad
\int_{\mathbb{R}^3} \AC{\mathbf P} F^\varepsilon \,(1,\ v,\ |v|^2/2)^T\,dv  =(0,0,0)^T,
\end{equation}
which decomposes $F^\varepsilon$ into its component on the local Maxwellian manifold and a transverse
fluctuation. All non-closure is concentrated in the microscopic fluctuation $\AC{\mathbf P} F^\varepsilon$. In particular, the stress
and heat flux are explicit functionals of $\AC{\mathbf P} F^\varepsilon$ (measuring the departure from the local Maxwellian
manifold), and quantitative control of $\AC{\mathbf P} F^\varepsilon$ via collisional relaxation is therefore the natural route
to a closed macroscopic description.
\unhide
 
\StartNoTOC

\subsection*{A. Main Results} We prove a moment--closure route to the incompressible Euler limit directly from the Boltzmann equation, without asymptotic expansions, for rough two-dimensional vorticity classes including nonunique regimes. Our results are organized in three layers. First, we prove an
$\varepsilon$-uniform quantitative closure estimate for the microscopic fluctuation,
which is the core mechanism of the analysis. Second, in two space dimensions, this
yields global convergence to weak notions of Euler solutions in rough vorticity classes,
including nonunique regimes. Third, under additional regularity assumptions, the
limiting dynamics is better posed and the convergence strengthens to quantitative
rates and full-family convergence.

\hide The key difficulty is \textbf{$\boldsymbol{\varepsilon}$-uniform closure}: controlling the microscopic fluctuation
$\AC{\mathbf P}F^\varepsilon$ strongly enough to close the macroscopic fluxes despite fast acoustic time scales
and large transport. Our main results show that this closure can be achieved for a broad class of initial data
in two space dimensions, yielding \emph{quantitative decay of $\AC{\mathbf P}F^\varepsilon$}. A key ingredient is the analysis of the \emph{kinetic vorticity} $\omega^\varepsilon := \nabla_x^\perp\cdot u^\varepsilon$ and optimal dissipation control. \unhide

Formally, the limiting incompressible dynamics is given by the Euler equations coupled to an advected entropic fluctuation 
\begin{align}\label{incompE}
\begin{aligned}
 \partial_t u^E + u^E\cdot\nabla_x u^E + \nabla_x p^E  = 0,\\
 \partial_t \Big(\frac{3}{2} \theta^E - \rho^E \Big) + u^E\cdot\nabla_x \Big(\frac{3}{2} \theta^E - \rho^E \Big)  = 0,\\
 \nabla_x\cdot u^E  = 0, \ \  \ \nabla_x(\rho^E+\theta^E)=0,
\end{aligned}
\end{align}
where the last constraint reflects the absence of an acoustic component. In particular, $u^E$ satisfies the vorticity formulation in 2D
\begin{equation}\label{weqnE2D}
\partial_t \omega^E + u^E\cdot\nabla_x \omega^E = 0,
\qquad
\omega^E = \nabla_x^\perp\cdot u^E.
\end{equation}

The major difficulty is \textbf{$\boldsymbol{\varepsilon}$-uniform closure}: controlling the
microscopic fluctuation $\mathring{\mathbf P}F^\varepsilon$ strongly enough to close the
macroscopic fluxes despite fast acoustic time scales and large transport. Our main
results show that this closure can be achieved for a broad class of initial data in
two space dimensions, yielding \emph{quantitative decay of
$\mathring{\mathbf P}F^\varepsilon$}. A key ingredient is the analysis of the
\emph{kinetic vorticity} $\omega^\varepsilon:=\nabla_x^\perp\cdot u^\varepsilon$ and
optimal dissipation control. The following quantitative microscopic
dissipation estimate provides the corresponding $\varepsilon$-uniform closure
mechanism and is the cornerstone of our global-in-time two-dimensional analysis,
underpinning the convergence results below.

\newtheorem*{theoremG}{Theorem G}
\begin{theoremG}[Informal statement of Theorem \ref{T.2D.global}]\label{theoremG}
Let $\Omega=\mathbb{R}^2$. Assume that the initial data $\{F^\varepsilon_0\}$ satisfy the modulated entropy bound 
(see \eqref{L2unif}), while allowing blow-up in
strong topologies (
see \eqref{ABC1}). We prove that the microscopic fluctuation $\AC{\P}F^\varepsilon$ is controlled in a weighted higher-order topology by the initial energy, up to the growth factor $\exp\!\big(\exp(\exp(\|\omega_0^\varepsilon\|_{\infty}\, t))\big)$ dictated by the initial kinetic vorticity scale.
In particular, $\e^{-2}\AC{\P}F^\varepsilon$ vanishes as $\varepsilon\to 0$ in this weighted topology for any fixed time.\end{theoremG}


We now state the principal hydrodynamic consequences in two space dimensions. The closure estimate above yields global convergence in rough weak-solution classes, including renormalized Euler solutions in the sense of DiPerna--Lions--Majda for $L^p$ vorticity (\cite{DiLi,DiMa}) and Delort-type solutions for nonnegative Radon-measure vorticity (\cite{Delort}). Precise definitions are recalled in Sections \ref{sec:DLM} and \ref{sec:RM}.


\hide

\begin{definition}[DiPerna--Lions--Majda \cite{DiLi,DiMa}]\label{D.soluw}
A pair $(u^E,\omega^E)$ is a \emph{renormalized solution} of \eqref{weqnE2D} on $[0,T]$ with initial vorticity
$\omega^E|_{t=0}=\omega^E_0$ if:
\begin{itemize}
\item $u^E={\bf K}\ast \omega^E$ and $\nabla_x\cdot u^E = 0$ in the sense of distributions, where
${\bf K}(x):= \frac{1}{2\pi}\frac{x^{\perp}}{|x|^2}$;
\item for every $\psi \in C_c^{\infty}([0,T)\times \mathbb{R}^2)$,
\[
\int_0^T \int_{\mathbb{R}^2} \boldsymbol{\beta}(\omega^E)\,(\partial_t\psi + u^E \cdot \nabla_x \psi)\,dxdt
+ \int_{\mathbb{R}^2}\boldsymbol{\beta}(\omega^E_0)\,\psi(0,x)\,dx = 0,
\]
for all $\boldsymbol{\beta} \in C^1(\mathbb{R}) \cap L^\infty(\mathbb{R})$ that vanish in a neighborhood of $0$.
\end{itemize}
\end{definition}

 \unhide

\newtheorem*{theoremP}{Theorem P}

\begin{theoremP}[Informal statement of Theorem \ref{T.C.Lp}]\label{theoremP}
Assume the hypotheses of Theorem~G. Assume moreover that the initial kinetic vorticity
$\omega^\varepsilon_0:=\nabla^\perp\!\cdot u^\varepsilon_0$ and entropic fluctuation
$\mathfrak s^\varepsilon_0:=\frac{3}{2}\theta^\varepsilon_0-\rho^\varepsilon_0$ converge strongly to
$(\omega_0,\mathfrak s_0)\in (L^p\cap L^1) \times L^2$ for some $p\in[1,\infty)$. 
Then, globally in time, $(u^\varepsilon,\omega^\varepsilon,\mathfrak s^\varepsilon)$ converges strongly (up to a
subsequence) and the microscopic fluctuation $\e^{-2}\AC{\P}F^\varepsilon$ decays as $\varepsilon\to 0$.
Any such limit $(u,\omega, \mathfrak s)$ is a global renormalized solution in the sense recalled in Section \ref{sec:DLM}. When $\omega_0  \in L^p \cap L^1$ for $p>1$, this solution conserves the energy $\|u(t)\|_{L_x^2}$ (Proposition~\ref{P.Econs}).
\end{theoremP}


\newtheorem*{theoremR}{Theorem R} 
\begin{theoremR}[Informal statement of Theorem \ref{T.Radon}]\label{theoremR}
Assume the hypotheses of Theorem~G, and suppose in addition that the initial vorticity
$\omega^\varepsilon_0\geq0$ converges weakly$^\ast$ in $\mathcal{M}(\mathbb{R}^2)$ to a nonnegative Radon measure
$\omega_0\ge 0$, that $\|\bw_0^{\e}\|_{L^1\cap H^{-1}_{loc}}<\infty$ and $\mathfrak{s}^{\e}_0 \to \mathfrak{s}_0$ strongly in $L^2(\R^2)$. Then $(u^\varepsilon,\mathfrak s^\e)$ converges up to a
subsequence and the microscopic fluctuation
$\e^{-2}\AC{\P}F^\varepsilon$ decays as $\varepsilon\to 0$, globally in time. Any such limit $(u,\mathfrak s)$ is a
global weak solution in primitive variables in the sense recalled in Section \ref{sec:RM}, equivalently a Delort solution \cite{Delort}.
\end{theoremR}

Additional regularity leads to stronger conclusions. In this regime, we identify the limit in a stronger solution class, establish quantitative convergence rates, and upgrade subsequential convergence to convergence of the full family. Our results include
\begin{itemize}
  \item
  In two space dimensions, if the initial vorticity is bounded, then the full family converges to the unique Yudovich solution. In particular, the velocity convergence is quantitative.

  \item
  If, in addition, the initial vorticity has additional Besov-type (Triebel--Lizorkin) regularity, then the vorticity convergence is quantitative as well.

  \item
  If the initial data are smooth, then in three space dimensions, as well as in two space dimensions, we obtain a quantitative convergence rate on a short time interval.
\end{itemize}
The full statements are given in Theorem~\ref{T.C.Linf}, Theorem~\ref{T.C.TL}, and Theorem~\ref{T.C.Hk}, respectively.
\hide
Additional regularity leads to stronger conclusions. We construct stronger limiting solutions, quantify convergence rates, and upgrade subsequential convergence to convergence of the full family. Our results include

\begin{itemize}[leftmargin=5.5em]
  \item
  In two space dimensions, if the initial vorticity is bounded, then the sequence converges
to the unique Yudovich solution. In particular, the velocity convergence becomes quantitative.
  \item
  In addition, if the initial vorticity has Besov-type regularity then the vorticity convergence also becomes quantitative.

  \item
  In three space dimensions (and also in two space dimensions), if the initial data are smooth,
then we obtain a sharp quantitative convergence rate on a short time interval.
\end{itemize}

The corresponding full statements can be found in Theorem~\ref{T.C.Linf}, Theorem~\ref{T.C.TL}, and Theorem~\ref{T.C.Hk}, respectively. 
\unhide
\hide
\hide 
\newtheorem*{theoremY}{Theorem Y}
\begin{theoremY}[Informal statement of Theorem~\ref{T.C.Linf}]
In two space dimensions, if the initial vorticity is bounded, then the sequence converges
to the unique Yudovich solution. Moreover, the velocity convergence is quantitative: $\| (\mathbb P u^\varepsilon - u^E)(t) \|_{L^2_x}
\lesssim_t
 (
\| \mathbb P u^\varepsilon_0 - u^E_0 \|_{L^2_x}
+ \kappa^{\frac12-}
 )^{e^{-Ct}}$.
\end{theoremY}

\newtheorem*{theoremB}{Theorem B}
\begin{theoremB}[Informal statement of Theorem~\ref{T.C.TL}]
In addition, if the initial vorticity has Besov-type regularity (Triebel--Lizorkin regularity $\mathcal F^s_2(\mathbb R^2)$ for some $s\in(0,1)$), then the vorticity convergence is quantitative: $\| (\omega^\varepsilon - \omega^E)(t) \|_{L^q_x}
\lesssim_t
 (
\| \mathbb P u^\varepsilon_0 - u^E_0 \|_{L^2_x}
+ \kappa^{\frac12-}
 )^{\widetilde C_{s,q} e^{-Ct}}$ when $q\in[2,\infty)$.
\end{theoremB}

\newtheorem*{theoremS}{Theorem S}
\begin{theoremS}[Informal statement of Theorem~\ref{T.C.Hk}]
In three space dimensions (and also in two space dimensions), if the initial data are smooth,
then we obtain a sharp quantitative convergence rate on a short time interval.
\end{theoremS}
\unhide

Additional regularity leads to stronger conclusions. The limiting fluid dynamics is better posed, and the convergence can be strengthened: we construct stronger limiting solutions, quantify convergence rates, and upgrade subsequential convergence to convergence of the full family. Our results include

\newtheorem*{theoremY}{Theorem Y} 

\begin{theoremY}[Informal statement of Theorem~\ref{T.C.Linf}]
In two space dimensions, if the initial vorticity is bounded, then the sequence converges
to the unique Yudovich solution. Moreover, the velocity convergence is quantitative:
  $\| (\mathbb P u^\e - u^E)(t) \|_{L^2_x}
    \lesssim_{t} \left(
    \|  \mathbb P u^\e_0 - u^E_0  \|_{L^2_x} + \kappa^{\frac{1}{2}-}
    \right)^{e^{-C t}}$.
\end{theoremY}

\newtheorem*{theoremB}{Theorem B}

    \begin{theoremB}[Informal statement of Theorem~\ref{T.C.TL}]In addition, if the initial vorticity has a Besov-type regularity (e.g. Triebel--Lizorkin regularity $\mathcal F^{s}_2(\R^2)$ for any $s\in(0,1)$),
    then the convergence upgrades to a \emph{quantified vorticity convergence}: $\| (  \omega^\e - \omega^E)(t) \|_{L^q_x}
    \lesssim_{t} \left(
    \|  \mathbb P u^\e_0 - u^E_0  \|_{L^2_x} + \kappa^{\frac{1}{2}-}
    \right)^{\tilde{C}_{s, q}e^{-C t}}$ for $q \in [2, \infty)$.
    \end{theoremB}

\newtheorem*{theoremS}{Theorem S} 
\begin{theoremS}[Informal statement of Theorem~\ref{T.C.Hk}]
 In 3D (and 2D) with smooth initial data, we obtain a sharp quantitative convergence rate on a short time interval.
\end{theoremS}

 These stronger regimes include:
(Y) the unique Yudovich regime in 2D, with a quantitative velocity estimate;
(B) sharper vorticity convergence under additional Besov-type regularity;
(S) a short-time sharp convergence theory for smooth data in 3D (and in 2D).


\begin{itemize}
    \item[\textbf{(Theorem Y)}] In 2D with bounded initial vorticity, convergence to the unique Yudovich solution with a \emph{quantified velocity convergence} (Theorem~\ref{T.C.Linf}):
    \[
    \| (\mathbb P u^\e - u^E)(t) \|_{L^2_x}
    \lesssim_{t} \left(
    \|  \mathbb P u^\e_0 - u^E_0  \|_{L^2_x} + \kappa^{\frac{1}{2}-}
    \right)^{e^{-C t}} .
    \]

    \item[\textbf{(Theorem S)}] In 3D (and 2D) with smooth initial data, we obtain a sharp quantitative convergence rate on a short time interval (Theorem~\ref{T.C.Hk}).
\end{itemize}
 
\unhide

\hide
\begin{definition}\label{D.weaksol}
A function $u(t,x)$ is called a \emph{weak solution} of the  Euler equation in primitive variable form if the following conditions hold:
\begin{itemize}
\item $u\in L^1([0,T]\times B_R)$ and $u\otimes u \in L^1([0,T]\times B_R)$. 
\item $u\in Lip([0,T];H^{-s}_{loc}(\R^2))$ and $u_0 \in H^{-s}_{loc}(\R^2)$ for some $s>0$.
\item $\div(u) = 0$ in the sense of distributions.
\item For every test function $\Lambda \in C_c^{\infty}(\R^+\times \R^2)$ with $\nabla_x \cdot \Lambda = 0$, 
\begin{align*}
\int_0^T \int_{\R^2} (\p_t\Lambda \cdot u + \nabla_x \Lambda : u\otimes u ) dxdt = 0,
\end{align*}
where $X:Y$ denotes $\sum_{i,j}X_{ij}Y_{ij}$. 
\end{itemize}
\end{definition}

\begin{theorem}\label{T.Info.Radon}[Informal statement of Theorem \ref{T.Radon}]
For $\Omega=\R^2$, consider
\begin{align*}
&\bw_0 \in \mathcal{M}(\R^2), \qquad \bw_0 \geq 0 \quad \mbox{or} \quad \bw_0 \leq 0, \qquad \mbox{supp}~ \bw_0 \subset \{x~|~|x|<R\}, \quad u_0 \in L^2_{\mathrm{loc}}(\R^2).
\end{align*}
Suppose that the mollified initial data satisfy the assumptions of Theorem~\ref{T.Info.global}. 
Then the following convergence results hold:
\begin{align*}
\bega
m^{\e}_A \quad &\rightarrow \quad m^{\#} \qquad \mbox{strongly in}, \qquad L^r([0,T]\times B_R(0)), \quad \mbox{for} \quad 1\leq r<2, \cr 
\bw_A^{\e} \quad &\overset{\ast}{\rightharpoonup} \quad d\bw^{\#} \qquad \mbox{weakly* in} \qquad L^\infty(0,T;\mathcal{M}(\R^2)).
\enda
\end{align*}
Here $m^{\#}$ is a \emph{weak solution of the Euler equations in primitive variable form} in the sense of Definition~\ref{D.weaksol}. Moreover,
$m^{\#} \in L^\infty(0,T;L^2_{\mathrm{loc}}(\R^2))$ and 
$\bw^{\#} \in L^\infty(0,T;\mathcal{M}(\R^2))$. In addition, The forcing part vanishes in the limit:
\begin{align*}
&m^{\e}_B, \ \bw^{\e}_B \quad \rightarrow \quad 0 \ \ \text{strongly in} \quad L^\infty(0,T; L^p(\R^2)), \quad \mbox{for} \quad  2\leq p \leq \infty \quad \mbox{with rate.} 
\end{align*}

\end{theorem}

\unhide

\subsection*{B. Relation to Previous Work}
We restrict ourselves to previous works directly related to our results. For general background on kinetic
theory and fluid dynamics, we refer to
\cite{CIP,MaBe,Villani02,VillaniBourbaki02}.

\subsubsection*{B.1. Moment--closure proofs.} 
A \emph{moment--closure proof} derives the fluid limit by closing the macroscopic moment system through uniform control of the microscopic fluctuation, without prescribing an \emph{a priori} limiting flow and without relying on asymptotic expansions. Consequently, the limit procedure provides a \textit{kinetic existence theory} for the limiting fluid dynamics.

Notably, the program initiated by
Bardos--Golse--Levermore (building on DiPerna--Lions renormalized solutions \cite{DL89}) develops hydrodynamic
limits from \emph{uniform} physical bounds via moment and compactness arguments \cite{BGL1,BGL2}. In particular,
the incompressible Navier--Stokes limit admits a well-developed theory\footnote{Nishida's compressible Euler limit \cite{Nishida} belongs to the moment-closure proofs relying on a
Cauchy--Kowalevskaya framework (hence not designed for shocks or low regularity) without an asymptotic expansion.}, with major
advances by Lions--Masmoudi and Golse--Saint-Raymond in the dynamical setting \cite{LM01a,LM01b,Golse-Saint}. At this viscous scaling, the entropy/energy inequality provides (i) a uniform $L^2$-type control of the fluctuation and (ii) a coercive \emph{collision dissipation} that strongly damps the microscopic component. This dissipation identifies the leading non-closed fluxes (stress/heat flux) at the correct scale and yields a genuinely \emph{parabolic} compactness mechanism at the macroscopic level, which is what ultimately allows one to pass to the nonlinear limit. The corresponding \emph{steady} problem is challenging, since the time-dependent entropy inequality is no longer available as a Lyapunov mechanism. In this direction, Esposito--Guo--Kim--Marra develop a different approach to
establish steady hydrodynamic limits via a micro--macro decomposition and stationary coercivity estimates for
the collision operator \cite{EGKM18}.

\subsubsection*{B.2. The incompressible Euler limit and the entropy method.} 
The incompressible Euler (low--Mach, high--Reynolds) scaling is substantially more singular than viscous limits. Even the \emph{existence} of a basic uniform $L^2$-type control of
fluctuations around a \emph{global} Maxwellian is not available in general, so the compactness mechanism
underlying moment methods for Navier--Stokes does not directly apply. This is one of the reasons why a \emph{moment--closure proof} of the incompressible Euler limit has long remained out of reach. A different line of attack \cite{Go, LM01b,Saint} is the \emph{entropy method of Yau.}\footnote{
Yau's entropy method was introduced for stochastic particle systems \cite{Y1,OVY,Va} and later adapted
to the kinetic framework in \cite{Go}, where it is observed that the method is compatible with Lions' concept of
\emph{dissipative solutions}, tailored to weak--strong stability.
} 
However, this weak-strong stability type result does not by itself provide a
moment-closure/compactness mechanism that would produce weak Euler limits without
inserting the fluid solution into the comparison functional.

Our results can be viewed as filling precisely this gap: we establish a \textit{moment--closure route} to the
incompressible Euler limit \emph{directly from the Boltzmann dynamics}, without prescribing an \emph{a priori}
Euler flow and without using asymptotic expansions. Moreover, the limiting Euler dynamics we obtain is broad: in
two space dimensions it covers the classical well-posed Yudovich regime, as well as genuinely non-unique
settings described by the DiPerna--Lions--Majda renormalized vorticity formulation and by the Delort class. As a
byproduct, the limit procedure yields a \emph{kinetic existence theory} for incompressible Euler solutions in
these weak vorticity classes.

\hide

In the incompressible Euler (low--Mach/high--Reynolds) regime, this parabolic compactness disappears, and one no longer has an $L^2$-based uniform control strong enough around the global Maxwellian to close the moment system by compactness alone.

This is precisely why a moment--closure proof of the incompressible Euler limit has been considered out of reach in the renormalized-solution framework and why alternative strategies (notably relative-entropy methods relying on a prescribed smooth Euler flow) have played a central role \cite{VillaniBourbaki}.

Grad's \emph{method of moments} starts from the (infinite) hierarchy of moment equations and asks for a \emph{closure} justified by near local equilibrium. In this spirit, the ``moment method'' program initiated by Bardos--Golse--Levermore and developed further with Lions--Masmoudi, Golse--Saint-Raymond, and others has yielded robust hydrodynamic limits in regimes where a uniform energy/entropy framework provides enough compactness and effective closure at the macroscopic level. In particular, the incompressible Navier--Stokes limit admits a well-developed moment/compactness theory for \emph{dynamic} solutions, while related \emph{steady} limits have also been treated (e.g.\ by Esposito--Guo--Kim--Marra). For compressible Euler limits, results in the moment/closure direction typically require strong regularity assumptions on the limiting flow (for instance, analytic frameworks as in Nishida), and therefore do not address regimes involving shocks or low regularity.

\paragraph{(ii) The incompressible Euler limit as a genuinely singular closure problem.} The inviscid incompressible Euler scaling (low Mach and high Reynolds) is markedly more singular. In this regime, even a basic uniform control of fluctuations in an $L^2$-based topology around a global Maxwellian is not available in general, and the standard compactness route underlying moment methods for viscous limits does not directly apply. This is one of the reasons why a \emph{moment--closure proof} of the incompressible Euler limit has long remained out of reach. A different line of attack is the \emph{relative entropy} strategy (Yau and its adaptations), which compares the kinetic solution to a \emph{prescribed} smooth solution of the target fluid system. This approach yields a weak--strong type stability principle and quantitative convergence as long as the limiting Euler flow is smooth, but it is not designed to produce rough Euler limits (nor to close the moment hierarchy without inserting the fluid solution by hand). Our results can be viewed as filling precisely this gap: we establish a moment--closure route to the incompressible Euler limit \emph{directly from the Boltzmann dynamics}, without prescribing an \emph{a priori} Euler flow and without using asymptotic expansions. Moreover, the limiting Euler dynamics we reach is broad: it covers, in two space dimensions, the classical well-posed Yudovich regime, as well as genuinely non-unique settings in which one must work with the DiPerna--Lions--Majda renormalized vorticity formulation and with the Delort class.

In contrast, for the \emph{inviscid} incompressible Euler limit, Villani emphasizes that the available approach
has been the relative-entropy method initiated by Yau and adapted to the kinetic setting by Golse: it proceeds by
introducing an explicit \emph{limiting} hydrodynamic solution and controlling the relative entropy with a Gr\"onwall
argument. :contentReference[oaicite:7]{index=7} :contentReference[oaicite:8]{index=8}
The price is that the limit is naturally formulated in very weak solution frameworks (e.g.\ Lions' dissipative
solutions), with a weak--strong principle when a smooth Euler flow exists. :contentReference[oaicite:9]{index=9}

\medskip
\noindent
\textbf{Our contribution (Euler by moment closure).}
The present work targets precisely this gap: we provide a \emph{moment--closure proof} of the incompressible Euler
limit \emph{directly from the Boltzmann dynamics}, without hinging on a prescribed limiting Euler flow and without
any Hilbert/Chapman--Enskog expansion. At the level of the limit equation, we accommodate a broad range of 2D
vorticity data---from classes where uniqueness holds (e.g.\ Yudovich) to genuinely nonunique regimes, treated in
the DiPerna--Lions--Majda setting and, for measure-valued vorticity, in the Delort class.

\subsubsection*{C.3. Hilbert expansions: validity and beyond.}
Hilbert and Chapman--Enskog expansions are primarily \emph{validity} results rather than mere convergence arguments:  they yield \emph{high-order} asymptotics and
thereby validate refined hydrodynamic descriptions (e.g. boundary layers \cite{JangKim, NKim1, NKim2}, ghost effects \cite{ghost}, shock profiles \cite{TYang}). As a byproduct of revealing the fine structure of the kinetic solutions, a convergence proof follows as a corollary, which is usually quantitative in a strong topology. 강한 결과와 스트럭쳘에 대한 fine 정보는 줌에도 불구하고, as a convergence argument, the method has clear limitation: typically require strong structural assumptions on the initial data and rely on the regularity of the
limiting fluid solution.   
\unhide

\subsubsection*{B.3. Quantitative hydrodynamic theory beyond Hilbert expansions.}

Hilbert and Chapman--Enskog expansions occupy a central place in hydrodynamic limit theory, as they provide much more than mere convergence: they deliver \emph{high-order} asymptotic descriptions and thereby justify refined hydrodynamic models, including phenomena beyond the leading-order closure \cite{dMEL,Guo-NS,KimLa}, such as boundary layers
\cite{JangKim,KimN1,KimN2}, ghost effects \cite{ghost}, and shock profiles. Precisely because they capture the fine structure of kinetic solutions, they often yield convergence as a byproduct, frequently with \emph{quantitative} bounds in strong topologies. Their strength, however, is naturally accompanied by a correspondingly structured regime of applicability, typically requiring well-prepared initial data and sufficient regularity of the limiting fluid solution over the time interval under consideration.

In particular, a \emph{quantitative} control of microscopic fluctuations away from the local Maxwellian manifold provides a concrete basis for closure in asymptotic hydrodynamic limits.\footnote{\label{Villani}\small\itshape 
«\,Le sujet est cependant bien loin d'\^etre \'epuis\'e, et l'on pourrait sugg\'erer maintes am\'eliorations de grande ampleur,
faisant na\^\i tre des difficult\'es colossales.
Tout d'abord, seul un th\'eor\`eme quantitatif pourrait donner une base physique incontestable \`a ces th\'eor\`emes
asymptotiques, par exemple : \`a quelles conditions sur la donn\'ee initiale et sur le nombre de Knudsen peut-on assurer
que les \'equations hydrodynamiques sont satisfaites avec une erreur relative n'exc\'edant pas, disons, 1\% ?
M\^eme si l'on parvenait \`a rendre quantitatifs les arguments des preuves actuelles, on serait amen\'e \`a des majorations en $O(1/\sqrt{\log |\log \e|})$, donc \`a des nombres de Knudsen d\'eraisonnables ($10^{10^{1000}} \cdots$).\,»
(from Section~6, \emph{Conclusions et probl\`emes ouverts}, in Villani's Bourbaki expos\'e \cite{VillaniBourbaki02}).}
Our analysis provides a new direct quantitative control of microscopic fluctuations at the \emph{optimal dissipation scale}, without relying on Hilbert or Chapman--Enskog expansions, or more generally on any high-order asymptotic construction around a smooth hydrodynamic profile, even when the limiting Euler solution belongs only to very weak regularity classes. In more regular regimes (for instance, bounded vorticity or slightly enhanced Besov regularity), we further obtain quantitative convergence rates for the velocity and, under stronger regularity assumptions on the vorticity, for the vorticity itself. In this sense, our results provide a new quantitative perspective on the hydrodynamic limit, complementary to the classical expansion method.

\hide Hilbert and Chapman--Enskog expansions are primarily \emph{validity} results rather than mere convergence
arguments: they yield \emph{high-order} asymptotics and thereby justify refined hydrodynamic descriptions (and,
in particular, moment information beyond the leading-order closure), such as boundary layers
\cite{JangKim,KimN1,KimN2}, ghost effects \cite{ghost}, and shock profiles. By revealing the fine
structure of kinetic solutions, these expansions typically yield convergence as a corollary, often with
\emph{quantitative} bounds in strong topologies. Even when they provide quantitative bounds, these methods have
clear limitations as general convergence arguments: they usually require stringent structural assumptions
(well-preparedness) on the initial data and rely on the regularity of the limiting fluid solution over the time
interval of validity. 
Hilbert and Chapman--Enskog expansions occupy a central place in hydrodynamic limit theory, as they provide much more than mere convergence: they deliver \emph{high-order} asymptotic descriptions and thereby justify refined hydrodynamic models, including phenomena beyond the leading-order closure, such as boundary layers
\cite{JangKim,KimN1,KimN2}, ghost effects \cite{ghost}, and shock profiles. Precisely because they capture the fine structure of kinetic solutions, they often yield convergence as a byproduct, frequently with \emph{quantitative} bounds in strong topologies. Their strength, however, is naturally accompanied by a correspondingly structured regime of applicability, typically requiring well-prepared initial data and sufficient regularity of the limiting fluid solution over the time interval under consideration.

 In particular, a \emph{quantitative} control of microscopic fluctuations away from the local Maxwellian
manifold provides a concrete quantitative basis for closure in asymptotic hydrodynamic limits.\footnote{\label{Villani}\small\itshape 
«\,Le sujet est cependant bien loin d'\^etre \'epuis\'e, et l'on pourrait sugg\'erer maintes am\'eliorations de grande ampleur,
faisant na\^\i tre des difficult\'es colossales.
Tout d'abord, seul un th\'eor\`eme quantitatif pourrait donner une base physique incontestable \`a ces th\'eor\`emes
asymptotiques, par exemple : \`a quelles conditions sur la donn\'ee initiale et sur le nombre de Knudsen peut-on assurer
que les \'equations hydrodynamiques sont satisfaites avec une erreur relative n'exc\'edant pas, disons, 1\% ?
M\^eme si l'on parvenait \`a rendre quantitatifs les arguments des preuves actuelles, on serait amen\'e \`a des majorations en $O(1/\sqrt{\log |\log \e|})$, donc \`a des nombres de Knudsen d\'eraisonnables ($10^{10^{1000}} \cdots$).\,»
(from Section~6, \emph{Conclusions et probl\`emes ouverts}, in Villani's Bourbaki expos\'e \cite{VillaniBourbaki02}).} Our analysis provides a direct quantitative control of the microscopic fluctuation at the \emph{optimal dissipation scale}, without recourse to Hilbert or Chapman--Enskog expansions, nor to any high-order asymptotic construction around a smooth hydrodynamic profile, even when the limiting Euler solution belongs only to very weak regularity classes. In more regular regimes (for instance, bounded vorticity or slightly enhanced Besov regularity), we further obtain quantitative convergence rates for the velocity and, under stronger regularity assumptions on the vorticity, for the vorticity itself. In other words, our results can be viewed as providing a new quantitative perspective on the hydrodynamic limit,
complementary to the classical expansion method.

\unhide

\hide
In particular, a \emph{quantitative} control of microscopic fluctuations away from the local Maxwellian
manifold provides a concrete quantitative basis for closure in asymptotic hydrodynamic limits.\footnote{\label{Villani}\small\itshape 
«\,Le sujet est cependant bien loin d'\^etre \'epuis\'e, et l'on pourrait sugg\'erer maintes am\'eliorations de grande ampleur,
faisant na\^\i tre des difficult\'es colossales.
Tout d'abord, seul un th\'eor\`eme quantitatif pourrait donner une base physique incontestable \`a ces th\'eor\`emes
asymptotiques, par exemple : \`a quelles conditions sur la donn\'ee initiale et sur le nombre de Knudsen peut-on assurer
que les \'equations hydrodynamiques sont satisfaites avec une erreur relative n'exc\'edant pas, disons, 1\% ?
M\^eme si l'on parvenait \`a rendre quantitatifs les arguments des preuves actuelles, on serait amen\'e \`a des majorations en $O(1/\sqrt{\log |\log \e|})$, donc \`a des nombres de Knudsen d\'eraisonnables ($10^{10^{1000}} \cdots$).\,»
(from Section~6, \emph{Conclusions et probl\`emes ouverts}, in Villani's Bourbaki expos\'e \cite{VillaniBourbaki02}).} In this quantitative perspective, our analysis yields quantitative control of the microscopic fluctuation at the \emph{optimal dissipation scale}, even when the limiting Euler solution is sought in
very weak classes. Moreover, in better-posed regimes (e.g.\ bounded vorticity, or slight additional regularity
in Besov spaces) we obtain quantified convergence rates for the velocity and, under higher vorticity regularity,
for the vorticity itself.
 \unhide

 \hide

 {Hydrodynamic limits and Previous Frameworks}
The mathematical study of the hydrodynamic limit primarily addresses {\it validating the closure of fluid models}, while the main challenge is to prove that {\it pure microscopic part vanishes} in the scaling limit. 

For fixed hydrodynamic fields   and 
T
T, the distribution   called the local Maxwellian, maximizes the entropy  
H(f) (Gibbs' lemma); this is why it is referred to as thermodynamic equilibrium.

Of course, it makes little physical sense to introduce a family of solutions depending on a parameter (
Kn
Kn) that tends to 0; this is a qualitative way of mathematically formalizing the smallness of 
Kn
Kn. It would be more satisfying to establish a quantitative estimate of the error made when replacing the kinetic equation with the hydrodynamic equation.

 it is natural to think that the distribution  is close to a local equilibrium in a physical regime where collisions are very frequent, i.e., with a low Knudsen number. The study of the hydrodynamic limits of the Boltzmann equation consists precisely of:

Proving a form of local thermodynamic equilibrium in the asymptotic regime.
Deriving limiting equations for the hydrodynamic fields associated with 
The lecture notes [Go], written for mathematician readers, provide a very accessible introduction to these problems.

From the Boltzmann equation \eqref{BE}, we formally have the local conservation laws: 
\Be\label{LCL}
\textbf{St} \p_t  \int_{\R^3}F(t,x,v) \begin{bmatrix}
1 \\ v \\ \frac{|v|^2}{2}
\end{bmatrix}
 dv + \nabla_x  \cdot   \int_{\R^3}F(t,x,v) \begin{bmatrix}
 v \\ v \otimes v  \\ v \frac{|v|^2}{2}
\end{bmatrix}
 dv
 = \begin{bmatrix}
0 \\ 0  \\ 0
\end{bmatrix}.
\Ee

Obviously the flux term contains microscopic parts.

The hydrodynamic limit of the Boltzmann equation has become a central focus of mathematical kinetic theory, as it systematically leads to the derivation of fluid equations such as the incompressible Navier-Stokes equations and the compressible Euler equations, among others. Notable contributions to this field have been made by researchers such as Bardos, Caflisch, Golse, Guo, Levermore, Nishida, Saint-Raymond, and many others. Their work addresses various aspects of both diffusive and compressible limits and has provided critical insights into the interplay between microscopic and macroscopic scales, offering rigorous frameworks for understanding convergence, stability, and compatibility conditions. Variuos mathematical frameworks for hydrodynamic limits can be categorized into three primary methodologies: the moment method, the asymptotic expansion method, and the relative entropy method. Each of these approaches provides distinct perspectives and tools for tackling the challenges of convergence, stability, and regularity within the context of the Boltzmann equation.

\hide

 The moment method involves deriving moment equations from the Boltzmann equation and demonstrating their convergence to the target fluid equations as the Knudsen number approaches zero. By examining these moment equations, one can establish the existence of solutions to the fluid equations through the Boltzmann framework. Key results in this area include the following. Bardos, Golse, and Levermore (BGL) present a rigorous proof of the fluid dynamic limits of the Boltzmann equation. Their work provides a foundational framework for the hydrodynamic limit, highlighting the convergence of solutions and systematically addressing transitions to the incompressible and compressible Euler and Navier-Stokes equations. The BGL approach relies on compactness methods and uniform estimates, which serve as a cornerstone for further developments in this area. Golse and Saint-Raymond, in their study published in Inventiones Mathematicae, build upon this foundational work, refining the analysis of fluid dynamic limits by delving into higher-order corrections and interactions between kinetic and fluid scales. Their study extends the understanding of hydrodynamic limits by incorporating detailed asymptotic expansions and addressing second-order effects in the Knudsen number. Guo extends the diffusive limit of the Boltzmann equation beyond the classical Navier-Stokes approximation. This work provides uniform estimates and compactness arguments, which are crucial for demonstrating convergence, especially in regimes where the standard fluid dynamic approximations fail. Esposito, Guo, Kim, and Marra investigate stationary solutions to the Boltzmann equation within the hydrodynamic limit. Their research emphasizes the impact of boundary effects on the solutions and their convergence properties, shedding light on the interaction between boundaries and fluid behavior in kinetic theory. Nishida explores the compressible Euler limit of the nonlinear Boltzmann equation, particularly addressing the regime up to the formation of shocks. This work bridges the gap between kinetic theory and compressible fluid dynamics, offering insights into how shocks and nonlinear phenomena emerge in the hydrodynamic limit.

\unhide

\subsubsection*{C.1. Moment method}

\subsubsection*{C.2. Yau's method}

\subsubsection*{C.3. Hilbert Expansion}  

The moment method involves deriving moment equations from the Boltzmann equation and demonstrating their convergence to the target fluid equations as the Knudsen number approaches zero. By examining these moment equations, one can establish the existence of solutions to the fluid equations through the Boltzmann framework. Key results in this area include the following.

Bardos, Golse, and Levermore (BGL) present a rigorous proof of the fluid dynamic limits of the Boltzmann equation. Their work provides a foundational framework for the hydrodynamic limit, highlighting the convergence of solutions and systematically addressing transitions to the incompressible and compressible Euler and Navier-Stokes equations. The BGL approach relies on compactness methods and uniform estimates, which serve as a cornerstone for further developments in this area.

Building upon this foundational framework, Lions and Masmoudi contribute significantly to the understanding of hydrodynamic limits by extending the analysis to low-regularity solutions and nonlinear interactions in the Boltzmann equation. Their work introduces new techniques that go beyond the compactness-based arguments used in BGL, employing functional analytic methods and renormalization techniques. These innovations enable the study of weak convergence in settings where traditional approaches face difficulties, particularly in capturing nonlinear effects and broader solution classes.

The advancements made by Lions and Masmoudi provide a strong foundation for the subsequent work of Golse and Saint-Raymond, who refine the understanding of hydrodynamic limits further. In their study published in Inventiones Mathematicae, Golse and Saint-Raymond build upon the results of BGL and Lions-Masmoudi by focusing on higher-order corrections and detailed asymptotic expansions. They delve into second-order effects in the Knudsen number and explore the interactions between kinetic and fluid scales in more complex regimes. This refinement integrates the foundational results of BGL with the extended scope introduced by Lions and Masmoudi, offering a comprehensive view of the hydrodynamic limit under various physical and mathematical conditions. Together, these contributions form a coherent narrative, with each step building on the previous to advance the field systematically.

Esposito, Guo, Kim, and Marra make a significant contribution by investigating stationary solutions of the Boltzmann equation within the hydrodynamic limit, particularly focusing on the impact of boundary effects. Their work complements and extends earlier results by analyzing the interplay between the Boltzmann equation and macroscopic fluid equations in bounded domains. By introducing refined techniques for handling boundary layers and demonstrating uniform stability of solutions, this study addresses challenges not fully resolved by earlier works such as BGL or Golse and Saint-Raymond. Their research highlights how boundary interactions influence the convergence of kinetic solutions to incompressible Navier-Stokes equations, providing a critical perspective on the moment method. Notably, this work builds upon the foundational analysis provided by Guo [2006], whose exploration of the diffusive limit introduced essential uniform estimates and compactness arguments. While Guo's 2006 results provided the groundwork for understanding diffusive scaling, the study by Esposito et al. extends these ideas to stationary states, emphasizing the role of physical boundary constraints in the hydrodynamic limit.

Nishida explores the compressible Euler limit of the nonlinear Boltzmann equation, particularly addressing the regime up to the formation of shocks. This work bridges the gap between kinetic theory and compressible fluid dynamics, offering insights into how shocks and nonlinear phenomena emerge in the hydrodynamic limit. Nishida's study extends the ideas introduced by BGL and complements the detailed asymptotic analysis of Golse-Saint-Raymond by tackling highly nonlinear regimes where compressibility and singularity formation are significant.

Taken together, these works illustrate the development and interconnectedness of the moment method, with each contribution building on and extending the foundational ideas laid out by Bardos, Golse, and Levermore. The progression from foundational convergence proofs to detailed asymptotic expansions and the inclusion of boundary, low-regularity, and nonlinear effects highlights the versatility and depth of the moment method in addressing the hydrodynamic limits of the Boltzmann equation.

\subsubsection{Asymptotic Expansion method}
The Asymptotic Expansion Method involves expanding the Boltzmann equation and its solutions as formal power series in terms of the Knudsen number, systematically deriving high-order fluid equations as the Knudsen number approaches zero. Unlike the moment method, which focuses on the convergence of solutions, the asymptotic expansion method emphasizes the ``convergence of equations," providing a hierarchical description of macroscopic dynamics derived from kinetic theory. This method is instrumental in constructing fluid equations and exploring the transition from kinetic to fluid regimes in greater detail.

Caflisch's foundational work provides a rigorous framework for the fluid dynamic limit of the nonlinear Boltzmann equation using asymptotic expansions. His study introduces formal expansions in the Knudsen number and proves the validity of these expansions for deriving fluid equations. By analyzing higher-order corrections, Caflisch establishes a systematic procedure for linking the Boltzmann equation to classical hydrodynamic models such as the Euler and Navier-Stokes equations. His approach underscores the importance of formal asymptotics in understanding the interplay between microscopic and macroscopic scales.

Building on these ideas, the work of de Masi, Esposito, and Lebowitz delves deeper into the incompressible Navier-Stokes and Euler limits of the Boltzmann equation. Their research focuses on specific scaling regimes and demonstrates how asymptotic expansions can rigorously connect the Boltzmann framework to the incompressible fluid equations. They explore the stability and structure of the derived equations, addressing challenges such as boundary conditions and the role of initial data in the convergence process.

Guo (2006) extends these results by addressing the Boltzmann diffusive limit beyond the classical Navier-Stokes approximation. His work introduces uniform estimates and compactness arguments to construct solutions that remain valid in more complex regimes, such as those involving high-order fluid equations. Guo's analysis also provides crucial insights into the robustness of asymptotic expansions under varying physical conditions, including higher-order corrections in diffusive scaling.

Recent studies on vorticity convergence further refine the application of the asymptotic expansion method. These works explore the transition from the Boltzmann equation to two-dimensional incompressible Euler equations below the Yudovich class, providing critical results for initial data with analytic regularity. Additionally, research on the incompressible Euler limit from the Boltzmann equation with diffuse boundary conditions investigates how boundary interactions influence the formal expansions and their convergence. These studies extend the asymptotic framework to include boundary layers, demonstrating how asymptotics can accommodate complex physical phenomena at kinetic-fluid interfaces.

Together, these works illustrate the development and refinement of the asymptotic expansion method, highlighting its role in connecting kinetic theory with fluid dynamics. Caflisch provides the foundational framework, de Masi, Esposito, and Lebowitz extend its scope to incompressible limits, Guo incorporates high-order corrections, and the vorticity convergence studies address boundary effects and analytic regularity. These contributions collectively deepen our understanding of the hydrodynamic limits of the Boltzmann equation, emphasizing the power of formal expansions in bridging scales.

\subsubsection{Relative Entropy method}
The relative entropy method offers a fundamentally different perspective by establishing the weak stability of Boltzmann solutions relative to target fluid solutions. This approach relies on the construction of a Lyapunov functional, which quantifies the deviation between kinetic and fluid solutions. A comprehensive treatment of this method is provided in Saint-Raymond’s monograph, Hydrodynamic Limits of the Boltzmann Equation. This work systematically applies the relative entropy method to various fluid regimes, offering a unifying framework for analyzing stability and convergence.

\newpage

\hide

Hilbert's sixth problem poses a profound challenge to both the mathematical and physical
communities: to establish a rigorous framework connecting microscopic dynamics to macroscopic
phenomena through the principles of mechanics and probability.
Within this program, a central objective is the derivation of fluid equations from the kinetic
theory of gases, most prominently from the Boltzmann equation.
In a dimensionless formulation, the Boltzmann equation takes the multiscale form
\begin{equation}\label{sBoltzmann}
\textbf{St}\,\partial_t F + v\cdot \nabla_x F \;=\; \frac{1}{\textbf{Kn}} \,\mathcal{N}(F,F),
\end{equation}
where $F=F(t,x,v)\ge 0$ is the one-particle distribution on phase space, the Strouhal number
$\textbf{St}$ characterizes the macroscopic temporal scale, and the Knudsen number $\textbf{Kn}$
measures the mean free path relative to the macroscopic length scale.
Monatomic intermolecular collisions are encoded by the bilinear collision operator
\begin{align}\label{Qdef}
\mathcal{N}(F,G)
:= \frac12 \int_{\mathbb{R}^3}\int_{\mathbb{S}^2_+} B(v-v_*,w)
\Big[&\,F(v_*')G(v') + F(v')G(v_*') \\
&\,- F(v_*)G(v) - F(v)G(v_*) \Big]\,dw\,dv_*,
\nonumber
\end{align}
where the post-collisional velocities are given by
\[
v^\prime = v- \big((v-v_*)\cdot w\big)\, w,
\qquad
v_*^\prime = v_*+ \big((v-v_*)\cdot w\big)\, w,
\]
and $\mathbb{S}^2_+$ is the half-sphere preventing double counting of collisions.
A structural cornerstone is the \emph{collision invariance}:
\begin{equation}\label{eq:collision_invariance_intro}
\Big\langle \mathcal{N}(F,G)(v),
\begin{pmatrix}
1\\ v\\ |v|^2/2
\end{pmatrix}\Big\rangle
=
\begin{pmatrix}
0\\ 0\\ 0
\end{pmatrix},
\end{equation}
expressing conservation of mass, momentum, and energy through binary collisions.
Throughout, $\langle \cdot,\cdot\rangle$ denotes the standard $L^2_v$ inner product.

\medskip
\noindent
\textbf{Phase space and the $\infty$-many moments system.}
The macroscopic fields arise as velocity moments of $F$ against the collision invariants:
\begin{equation}\label{eq:macro_fields_intro}
\rho(t,x):=\langle F,1\rangle,\qquad
\rho u(t,x):=\langle F,v\rangle,\qquad
\frac{3}{2}\rho\theta(t,x)+\frac12\rho|u|^2 := \Big\langle F,\frac{|v|^2}{2}\Big\rangle.
\end{equation}
More generally, for any polynomial (or suitable test function) $\phi(v)$, the moment
$m_\phi(t,x):=\langle F,\phi\rangle$ satisfies an evolution law obtained by multiplying
\eqref{sBoltzmann} by $\phi$ and integrating in $v$.
This yields an \emph{infinite hierarchy} of coupled PDEs---an $\infty$-many moments system---whose
intrinsic lack of closure reflects the multiscale character of \eqref{sBoltzmann}.
At the level of macroscopic conservation laws, the hierarchy closes only up to higher-order
moments (stress and heat flux), whose control is precisely the analytical bottleneck in
hydrodynamic limits.

\medskip
\noindent
\textbf{Collisional relaxation and Maxwellian equilibria.}
The collision operator drives $F$ toward local thermodynamic equilibrium.
The equilibrium manifold is given by the Maxwellians
\begin{equation}\label{eq:maxwellian_intro}
M[\rho,u,\theta](v)
=
\frac{\rho}{(2\pi\theta)^{3/2}}
\exp\!\Big(-\frac{|v-u|^2}{2\theta}\Big),
\end{equation}
characterized by $\mathcal{N}(M,M)=0$ and uniquely determined by the moments
\eqref{eq:macro_fields_intro}.
Thus, the mechanism by which the moment hierarchy becomes effectively closed is the
\emph{relaxation by collision}: for regimes in which collisions are sufficiently frequent,
the solution remains close to the Maxwellian manifold, and the only persistent degrees of
freedom are the macroscopic fields $(\rho,u,\theta)$.

\medskip
\noindent
\textbf{Macro--micro decomposition.}
A convenient framework to quantify this closure is the macro--micro splitting
\begin{equation}\label{eq:macro_micro_intro}
F = M + G,\qquad M:=M[\rho,u,\theta],\qquad
\Big\langle G,
\begin{pmatrix}
1\\ v\\ |v|^2/2
\end{pmatrix}\Big\rangle
=
\begin{pmatrix}
0\\ 0\\ 0
\end{pmatrix},
\end{equation}
which separates fluid variables (the Maxwellian part) from kinetic fluctuations (the micro part).
Projecting \eqref{sBoltzmann} onto the collision invariants yields the compressible balance laws
for $(\rho,u,\theta)$, but with fluxes containing the microscopic stress and heat flux, which are
linear functionals of $G$.
On the other hand, the equation for $G$ involves the linearized collision operator around $M$;
its coercivity (modulo the collision invariants) is the quantitative expression of relaxation and
is the source of dissipation that controls the micro component in norms compatible with the
macroscopic nonlinearities.

\medskip
\noindent
\textbf{Incompressible regime (low Mach number).}
A particularly delicate and fundamental asymptotic limit is the incompressible (low Mach number)
regime, in which the macroscopic velocity is small compared to the thermal speed, while acoustic
modes oscillate on fast time scales.
In dimensionless variables, this corresponds to a coupled scaling of the form
\begin{equation}\label{eq:lowMach_scaling_intro}
\textbf{St}=\varepsilon,\qquad \textbf{Kn}=\kappa\,\varepsilon,
\qquad \varepsilon\to 0,\quad \kappa=\kappa(\varepsilon)\to 0,
\end{equation}
so that collisions enforce near-equilibrium on the kinetic level, yet viscous effects vanish at
the macroscopic level.
In this scaling, the macro--micro decomposition \eqref{eq:macro_micro_intro} becomes the natural
tool to disentangle the incompressible component from acoustic oscillations and to identify the
effective limiting dynamics.
Formally, after appropriate filtering of the fast acoustic time scale, one expects the limiting
velocity field to satisfy the incompressible Euler equations, while the fluctuation $G$ is
controlled by collisional relaxation.
The rigorous implementation of this program requires a robust closure mechanism:
one must propagate uniform bounds on the macroscopic fields and exploit coercivity of the
linearized collision operator to control $G$ strongly enough to pass to the limit in the
nonlinear fluxes of the projected conservation laws.

Hilbert's sixth problem poses a profound challenge to both the mathematical and physical communities: to establish a rigorous framework that connects microscopic dynamics to macroscopic phenomena using the principles of mechanics and probability. A central aspect of this problem involves deriving macroscopic fluid equations from the foundational kinetic theory of gases, particularly through the Boltzmann equations. In its dimensionless form, the Boltzmann equation is expressed as a multiscale model:
\Be\label{sBoltzmann}
\textbf{St} \p_t F + v\cdot \nabla_x F = \frac{1}{\textbf{Kn}} \mathcal{N}(F,F),
\Ee
where the Strouhal number ($\textbf{St}$) characterizes the temporal scale, while the Knudsen number ($\textbf{Kn}$) represents the mean free path. Monatomic intermolecular collisions are accounted for by the collision operator:
\begin{align} 
\mathcal{N}(F,G) 
:=   \frac{1}{2}\int_{\mathbb{R}^3}\int_{\mathbb{S}^2_+}B(v-v_*,w)\big[ F(v_*')G(v')+F(v')G(v_*') 
-F(v_*)G(v)-F(v)G(v_*)\big]dwdv_*, \label{Qdef}
\end{align} 
where $v^\prime = v- ((v-v_*) \cdot w) w$ and $v_*^\prime = v_*+ ((v-v_*) \cdot w) w$. This operator enjoys a crucial {\it collision invariance} $ \langle \mathcal{N}(F,G) (v), [1, v , \frac{|v|^2}{2}]^T \rangle= [0,0,0]^T$. In this paper, we denote a standard $L^2_v$-inner product by $\langle   \ ,  \   \rangle $.

\subsection{Phase space, moment hierarchies, and the closure problem}

A collisional gas is described at the kinetic level by a nonnegative distribution
\[
F=F(t,x,\xi)\ge 0, \qquad (t,x,\xi)\in [0,T]\times \Omega \times \mathbb{R}^3,
\]
where $x\in\Omega\subset \mathbb{R}^d$ ($d=2,3$) and $\xi\in\mathbb{R}^3$ is the microscopic velocity.
Macroscopic observables are encoded by velocity moments: for any test function $\phi(\xi)$,
\begin{equation}\label{eq:moment_def}
m_\phi(t,x) := \int_{\mathbb{R}^3} F(t,x,\xi)\,\phi(\xi)\,d\xi .
\end{equation}
The collection $\{m_\phi\}$ may be viewed as an \emph{infinite moment system} in phase space.
For instance, expanding $F$ in an orthogonal basis in $\xi$ (e.g.\ Hermite polynomials) leads formally to
a countable set of coupled PDEs for the coefficients.
The fundamental difficulty is that the moment equations are not closed: lower-order moments
are forced by higher-order moments, reflecting the intrinsically multiscale nature of kinetic dynamics.
A quantitative understanding of \emph{how} and \emph{when} this hierarchy effectively closes
in singular asymptotic regimes is a central theme of hydrodynamic limits.

Among all moments, the collision invariants $(1,\xi,|\xi|^2/2)$ play a distinguished role.
They define the macroscopic density, bulk velocity, and temperature through
\begin{equation}\label{eq:macro_fields}
\int_{\mathbb{R}^3} F
\begin{pmatrix}
1\\ \xi\\ |\xi|^2/2
\end{pmatrix}\,d\xi
=
\begin{pmatrix}
P_F\\ P_F U_F\\ \tfrac{3k_B}{2}P_F\Theta_F + \tfrac{1}{2}P_F|U_F|^2
\end{pmatrix},
\end{equation}
where $k_B$ is Boltzmann's constant.
The triple $(P_F,U_F,\Theta_F)$ is the natural macroscopic state associated to $F$.

\subsection{Collisional kinetic equations and Maxwellian equilibria}

The evolution of $F$ is governed by a transport--collision balance of the form
\begin{equation}\label{eq:kinetic}
\partial_t F + \xi\cdot\nabla_x F = \mathcal{N}(F,F),
\end{equation}
where $\mathcal{N}$ is a bilinear collision operator encoding particle--particle interactions.
For concreteness, one may keep in mind the hard-sphere Boltzmann operator,
\begin{equation}\label{eq:boltzmann_hardsphere}
\mathcal{N}(F,H)(\xi)
=
\int_{\mathbb{R}^3}\int_{\mathbb{S}^2_{+}}
|(\xi-\xi_*)\cdot \omega|\,
\Big( F(\xi')H(\xi_*') - F(\xi)H(\xi_*)\Big)\,d\omega\,d\xi_*,
\end{equation}
where $(\xi',\xi_*')$ are the post-collisional velocities produced by specular reflection along $\omega$,
so that momentum and kinetic energy are conserved.
The operator $\mathcal{N}$ is irreversible in the sense that it creates correlations among outgoing pairs,
even when incoming particles are approximately independent.

The collision invariants imply that $\mathcal{N}$ conserves mass, momentum, and energy:
\begin{equation}\label{eq:collision_invariants}
\int_{\mathbb{R}^3}
\begin{pmatrix}
1\\ \xi\\ |\xi|^2/2
\end{pmatrix}
\mathcal{N}(\cdot,\cdot)(\xi)\,d\xi
=
\begin{pmatrix}
0\\ 0\\ 0
\end{pmatrix}.
\end{equation}
Equilibria are local Maxwellians, characterized by $\mathcal{N}(M,M)=0$:
\begin{equation}\label{eq:maxwellian}
M[\rho,u,\theta](\xi)
=
\frac{\rho}{(2\pi k_B \theta)^{3/2}}
\exp\!\Big(-\frac{|\xi-u|^2}{2k_B\theta}\Big),
\qquad (\rho,u,\theta)\in \mathbb{R}_+\times \mathbb{R}^3\times \mathbb{R}_+.
\end{equation}
The Maxwellian manifold is the natural candidate for the macroscopic closure of \eqref{eq:kinetic}.

\subsection{Macro--micro decomposition and quantitative relaxation}

A robust way to separate fluid variables from kinetic fluctuations is the macro--micro decomposition
\begin{equation}\label{eq:macro_micro}
F = M + G,
\qquad
M := M[P_F,U_F,\Theta_F],
\qquad
\int_{\mathbb{R}^3} G
\begin{pmatrix}
1\\ \xi\\ |\xi|^2/2
\end{pmatrix}\,d\xi = 0.
\end{equation}
Projecting \eqref{eq:kinetic} onto the collision invariants yields the local conservation laws
\begin{align}
\partial_t P_F + \nabla_x\cdot(P_F U_F) &= 0, \label{eq:mass}\\
\partial_t(P_F U_F)
+ \nabla_x\cdot\Big(P_F U_F\otimes U_F + k_B P_F\Theta_F I + r_F\Big) &= 0,\label{eq:momentum}\\
\partial_t\Big(\tfrac{3k_B}{2}P_F\Theta_F + \tfrac12 P_F|U_F|^2\Big)
+ \nabla_x\cdot\Big(\tfrac{5k_B}{2}P_F\Theta_F U_F + \tfrac12 P_F|U_F|^2U_F + q_F + r_F U_F\Big) &= 0,
\label{eq:energy}
\end{align}
where the stress deviator $r_F$ and heat flux $q_F$ are \emph{purely microscopic}:
\begin{equation}\label{eq:stress_heat}
r_{F,ij} = \int_{\mathbb{R}^3} R_{ij}(\xi)\,G\,d\xi,
\qquad
q_{F,j} = \int_{\mathbb{R}^3} Q_j(\xi)\,G\,d\xi,
\end{equation}
with
\[
R(\xi):=(\xi-U_F)\otimes(\xi-U_F)-\frac{|\xi-U_F|^2}{3}I,
\qquad
Q(\xi):=(\xi-U_F)\Big(\frac{|\xi-U_F|^2-5k_B\Theta_F}{2}\Big).
\]
Thus, the moment hierarchy closes \emph{exactly} at the Euler level if and only if $G\equiv 0$;
more generally, the deviations from ideal fluid dynamics are entirely encoded by the micro part $G$.

The dynamics of $G$ is governed by a transport equation with a linearized collision operator
around the evolving Maxwellian $M$:
\begin{equation}\label{eq:micro_eq}
\partial_t G + \xi\cdot\nabla_x G - \big(\mathcal{N}(M,G)+\mathcal{N}(G,M)\big)
= \mathcal{N}(G,G) + \big(\partial_t M + \xi\cdot\nabla_x M\big).
\end{equation}
A key mechanism is the coercivity of the linearized collision operator,
which can be read as a quantitative relaxation toward local equilibrium.
At the level of entropy production one expects an estimate of the form
\begin{equation}\label{eq:coercivity}
-\int_{\mathbb{R}^3}\big(\mathcal{N}(M,G)+\mathcal{N}(G,M)\big)\,\frac{G}{M}\,d\xi
\ \gtrsim\ 
\Big\|\frac{G}{\sqrt{M}}\Big\|_{L^2_\xi}^2,
\end{equation}
modulo the null space generated by the collision invariants.
In singular scaling limits, the central analytical problem is to convert
\eqref{eq:coercivity} into \emph{effective closure}: uniform control of the macro fields
together with a quantitative decay (or smallness) of $G$ in norms compatible with
the macroscopic nonlinearities in \eqref{eq:mass}--\eqref{eq:energy}.

\subsection{Low-Mach incompressible regime and hydrodynamic scaling}

We now outline the incompressible (low Mach number) regime.
Fix reference scales $(T_0,L_0,\Theta_0)$ and the associated thermal speed $c\sim \sqrt{\Theta_0}$.
The Mach number is the ratio of macroscopic and thermal velocities,
\[
\mathrm{Ma}:=\frac{L_0/T_0}{c},
\]
and the Knudsen number $\mathrm{Kn}$ measures the mean free path relative to $L_0$.
After nondimensionalization one arrives at the scaled kinetic equation
\begin{equation}\label{eq:scaled_boltzmann}
\mathrm{Ma}\,\partial_t F + \xi\cdot\nabla_x F = \frac{1}{\mathrm{Kn}}\,\mathcal{N}(F,F).
\end{equation}
In the low-Mach regime the macroscopic state is close to a global Maxwellian,
and one may parametrize fluctuations by
\begin{equation}\label{eq:low_mach_param}
(P,U,\Theta) = (1,0,1) + O(\mathrm{Ma}).
\end{equation}
To access the incompressible Euler dynamics, one considers a coupled limit in which
\begin{equation}\label{eq:euler_scaling}
\mathrm{Ma}=\varepsilon,
\qquad
\mathrm{Kn}=\kappa\,\varepsilon,
\qquad
\kappa\to 0 \ \text{as}\ \varepsilon\to 0,
\end{equation}
so that collisions are fast enough to enforce near-equilibrium
while viscosity effects vanish at the macroscopic level.
In this scaling, it is convenient to represent the macroscopic fields as
\begin{equation}\label{eq:macro_eps}
P^\varepsilon = e^{\varepsilon \rho^\varepsilon},
\qquad
U^\varepsilon = \varepsilon u^\varepsilon,
\qquad
\Theta^\varepsilon = e^{\varepsilon \theta^\varepsilon},
\end{equation}
which isolates the acoustic component $(\rho^\varepsilon+\theta^\varepsilon)$
and the incompressible component of $u^\varepsilon$.
Formally, the leading-order balance in \eqref{eq:scaled_boltzmann} drives $F^\varepsilon$
toward a Maxwellian and yields, after filtering acoustic oscillations,
the incompressible Euler equations for the limiting velocity field.
The rigorous justification of this scenario requires a quantitative closure mechanism:
one must propagate sufficiently strong bounds on the macroscopic fields while exploiting
collisional coercivity to control $\AC{\P}F^{\e}$ in a manner compatible with
the low-Mach singular structure (fast time scales and large transport).

1. Hilbert 6

 키네틱 에서 플루이드로

2. Convergence in equation:  admissible class 를 특정하는 것이 일반적 

NS 리밋

Open question 1: 인컴프레서블 오일러 리밋 오픈

솔루션의 컨스트럭션

3. Convergence in solution: stability+uniqueness

\unhide

Open question 2

목적: ``The works presented here bear witness to considerable progress in the treatment of hydrodynamic limits. The subject is, however, very far from being exhausted, and one could suggest many large-scale improvements (불어 원어), giving rise to colossal difficulties (불어 원어). First of all, only a quantitative theorem could give an indisputable physical basis to these asymptotic theorems (불어 원어), for example: under what conditions on the initial data and on the Knudsen number can one guarantee that the hydrodynamic equations are satisfied with a relative error not exceeding, say, $1\%$? Even if one succeeded in making the arguments of the current proofs quantitative,
this would lead to estimates in $O(1/\sqrt{\log |\log \e|})$, and therefore to unreasonably large Knudsen numbers $10^{{10}^{1000}}$.''

왜 오픈이고, 다른 결과들과 어떠한 관계에 있는가?

$\nabla u  \in L^\infty$ control 이 핵심

4. Validity of asymptotic expansion: key difficulty raises when singularity of the target solution such as Shock layer, Prandtl layer, 

매끈한(또는 안정한) 유체 해 주변에서의 정확한 해의 local manifold를 구성하는 방법

expansion

우리 결과 설명

한 문단으로 컨샙셜 노블티 + 태크니컬 노블 +

대표적인 방식이 모멘텀 메쏘드이다. 이는 어프라이오리 유니폼 바운드와 컴팩트니스를 기반으로 약 수렴성을 공부, 볼츠만 솔루션의 리밋으로서 유체 해의 존재성을 준다. 또한 싱굴러 리밋 문제의 스테빌리티 어널리스로 볼수 있는데, 근방을 최대한 크게 (더 큰 싸이즈, 더 약한 토폴로지) 잡아주는 것이 미덕이다. 
\footnote{The second author first learned this perspective from Professor Yan Guo’s lectures for the course “Topics in Kinetic Theory” at Brown University (2010; CK’s personal notes).}
그림

Golse--Saint-Raymond

Guo 2006

스테디 문제는 EGKM 의 결과

\subsection{Two Major Open Questions}

\subsubsection{a quantitative theorem}

유체극한의 정량분석에서 무엇이 알려져 있는가? 힐버트 익스펜젼 등의 방식으로 유체 솔루션의 다단계을 알게되면 고차 유체솔루션과 가까운 볼츠만 방정식을 구성, 따라서 이는 넌리니어 유체와 고차 유체의 벨리디티를 준다. 이는 볼츠만 방정식의 해 구성 방법으로 볼수도 있으면, 특정 유체 프로파일과 어심토틱컬리 가까운 해를 구성하는데 사용할 수 있다. 따라서 근방을 키우는 것이 미덕이 아니라 유체의 해에 가깝게 (구조적으로, 토폴로지적으로) 근사하는 것이 메쏘드가 가지는 본연의 강점이다. 물론 G의 입장에선, 두번째 방식엔 초기 조건에 제약이 따른다. 

Well-prepared 를 그림으로 설명

\subsubsection{Large Reynolds number limit}

Saint-Raymond

 그중에 인컴프레써블 나비어스톡스는 이 둘이 잘 작동하는 스케일링으로 수학적으로도 가장 잘 발전되어 있다.

 Main difficulties: 
1) 따라서 스케일 자체는 인컴프레써블 오일러 리밋이 가장 어렵다고 할 수 있다. 포말하게 보면 이는 꽤 쉽게 보이는데, 
\[
\text{ 에스티메잇: 에너지+ 디씨페이션 $\lesssim  \frac{1}{\sqrt{\kappa}}  \sqrt{\text{에너지}}$ 디씨페이션 }\]
물론 나비어스톡 (카파가 오더 원)의 경우 수학적으로 트리비얼 하지 않은 에버리지 렘마와 $L^6$ bound를 사용해서 리궈러스하게 닫을 수 있다.  

 반면 나비어스톡스의 하이 레이놀드 수렴인 인컴프레서블 오일러는 스케일이 좋지 않다. 위의 식에서 $1/\sqrt{\kappa}$ 만큼 싱굴러한 것이 포말레벨에서 보인다.

2) 또한 플루이드 입장에서는 페널라이제이션 텀이 싱굴러하다. $\frac{1}{\e} \nabla \cdot pressure$. 시간미분에 대한 조건을 강하게 주면 문제가 쉬워지나, 이럴 경우 이니셜 데이터에 강한 제약을 줘야한다. 

따라서 모멘텀 메쏘드가 작동하지 않는다고 믿었다. 

유니폼 바운드가 좋지 않아 모멘텀 메쏘드로는 안되었고, 디시파티브 솔루션으로의 존재성이 알려졌는데 이는 렐라티브 엔트로피+위크스트롱 알규먼트라고 할 수 있다. (굴스의 써베이에 많이 인용해서 이 점을 명확히 하자) 따라서 립쉬츠 혹은 유도비치 클래스가 아닌 오일러 솔류션의 경우, 수렴성의 의미가 직관적이지 않다. 더욱이 그 수렴하는 솔루션의 구조적인 특징을 알수가 없다. 수렴의 정량화 등등

\subsection{Main Results}

Motivation: 1) resolve these two major open questions introducing a new quasilinear method; 

Substantial goals

2) microscopic fluctuation control when $\omega$ control without well-preparedness: 2D tubulence theorem, $\omega \notin $ measure, $L^p$, $L^\infty$, Besov 

3) Construction of Euler solutions from Boltzmann solutions as a limit without well-preparedness: Energy conserved-solution

Full Statement

\subsection{New Methods and Ideas to overcome Difficulties}

Major difficulty: lack of low frequency dissipation control

Acoustic wave control

\subsection{Various Expansions for the Well-Prepared data}
$\bullet$ 하이드로 다이나믹의 두 가지 다른 접근법들: 해의 수렴, 방정식의 수렴, 해의 수렴은 유체 존재성을 줄 수 있고, 방정식의 수렴은 어심토틱 어널리시스의 장점이 있다.
\begin{itemize}
    \item

    모멘텀 메쏘드로는 오일러 방정식을 공부하지 못했다.

    \item Our results: large perturbation (not well-prepared), momentum method approach, 
\end{itemize}

\subsection{Open Questions}

$\bullet$ 오픈된 싱굴러 스케일 문제들: 싱굴러 스케일 문제는 스케일링에 따라 문제가 더 어려워지며, 기대되는 리밋 솔루션의 유니크니스/레귤러리티에 따라 또한 문제가 어려워진다. 물론 어느 방식이냐에 따라서도 난이도에 썹틀한 차이가 있다.
\begin{itemize}
\item 컴프레써블 스케일링은 볼츠만에서 오는 자연스러운 스케일이나 기대되는 리밋 솔루션의 정칙성때문에 문제가 어렵다.

\item 터뷸런스 관련 2차원에선 볼티시티를 공부해야한다. 

키네틱 볼티시티 에스티멧은 $\nabla^\perp \cdot \div \hat \alpha$ 을 공붛해야하는데, 이는 하드 스피어에선, 하이오더 소볼래브 에스티멧이 필수적

2차원 터불런스 theory enthropy

\item 정량화: 빌라니 부르바키 논문을 인용하여 정량화된 하이드로 리밋의 중요성을 강조하자.

어려운점: 만약 $1/\sqrt \kappa$ 블로업이 있으면 이니셜 레귤라리제이션을 써도 빌로우 유도비치에서 문제를 풀수 없다. $\sqrt{\kappa}e^{- \frac{t}{\sqrt \kappa}}$. 그리고 $1/\e$ 곱하기 프레셔의 페넌트레이션 텀이 존재하기 때문에, 이니셜데이터가 웰프리페어 되어 있지 않는 경우, 이는 또다른 어려움이다.

\end{itemize}

$\bullet$ 이 논문의 모티베이션과 우리의 결과: 이 논문의 모티베이션은 이 지점에서 출발한다. 이 논문의 목적은 다음과 같다.

\subsection{Validity of Closure}
 {\color{red} “벨리티이 오브 클로져”: 퓨어 마이크로 스코픽의 정량적 감쇄를 보인다. 즉 유체 매니폴드에서의 초기 디비에이션이 스케일링과 시간에 따라 변하는 것을 정량적으로 규명한다. 즉, 꽤 일반적인 퍼터베이션에서 유체 매니폴드가 마이크로 퍼터베이션에 대해 로컬 스테이블 매니폴드임을 정량적으로 보인다. }

 \textbf{*Informal Statement of $G$ estimate} 백그라운드 매크로스코픽의 이니셜 데이터 가정하에, 정량적 수렴을 보일수 있다. 

리마크들....

(Bootstrap 을 쓰면)
너무 테크니컬 해보임

{\color{red}

\subsection{Existence proof of the Fluids}
 유체 해의 존재성을 어프라이오리 가정하지 않고, 순전히 볼츠만에서 스케일 리밋으로 유체의 해를 구성한다.

\textbf{*Informal Statement: 2D 라돈 메져 볼티시티 결과}

\subsection{Invariant Property}
  리밋 해가 유니크하지 않는 상황에서도, 스케일링 리밋에서 인바리언트한 프라퍼티를 규명한다. 구체적으로 립쉬츠 이니셜 데이터가 아닌 볼츠만 솔루션의 스케일 리밋이 리노말라이즈드 해로 수렴함을 보인다. 즉 수렴해는 라그랑지안이어야한다 (리지디티: 그로모브 리지디티, 디펄나-리옹). 

  \textbf{*Informal Statement: 2D $\omega_0 \in L^p$ for $p \in [1, \infty]$, 리노말라이즈 솔루션 즉 라그랑지안이 인바이언트함 }

\subsection{Quantitative Analysis}
 모멘텀 시스템의 정확한 정량적 구조 공부한다. 예를 들어 유체에서 잘 알려진 더블 익스포넨셜 그로스 등등을 관찰할 수 있다. 또한 레귤러리티가 좋은 경우, 컨벌젼스 레잇을 공부할수 있다. 나레이티브하게 서술하자. 

  \textbf{*Informal Statement: 3D 로컬인 타임 수렴+ $H^s$}

    \textbf{*Informal Statement: 더블 익스포넨셜 그로스........}

 }

\begin{theorem}\label{T.Info.global}[Informal statement of Theorem \ref{T.2D.global}]

For $\kappa=\e^q$ with $0<q<2$, there exists $\e_0>0$ such that, if the initial data satisfy
\begin{align*}
&\sup_{0<\e<\e_0}
\Big(
\|(\rho^{\e}_0,u^{\e}_0-\bar{u},\ta^{\e}_0)\|_{L^2(\R^2)}
+\|\w^{\e}_0\|_{L^{p_1}(\R^2)}
\Big)
<\infty, \quad 
\sup_{0<\e<\e_0}
\bigg\|
\frac{1}{\e}
\frac{(F_0-M_{[\rho_0,u_0,\ta_0]})}{\sqrt{\tilde{\mu}}}
\bigg\|_{L^2_x(\R^2 ; L^2_v(\R^3))}
<\infty,
\\[0.5em]
&\sup_{0<\e<\e_0}
\bigg\|
\frac{\e^{1-}}{\e}
\frac{(F|_{t=0}-\mu)}{\sqrt{\tilde{\mu}}}
\bigg\|_{L^\infty_{x,v}(\R^2\times\R^3)}
\les 1,
\qquad \hspace{16mm}
\sup_{0<\e<\e_0}
\bigg\|
\frac{\e^{1-}}{\e}
\frac{(F|_{t=0}-M_{[\rho_0,u_0,\ta_0]})}{\sqrt{\tilde{\mu}}}
\bigg\|_{L^\infty_{x,v}(\R^2\times\R^3)}
\les 1,
\end{align*}
where $\tilde{\mu}=M_{[1,0,1-c_0]}$ for any $0<c_0\ll 1$.
Then the family of Boltzmann solutions $\{F^{\e}\}_{0<\e<\e_0}$ satisfies

\begin{align*}
&\frac{1}{\e^2}
\frac{\nabla_x^n \AC{\P}F^{\e}}{\sqrt{M^\e}}
\to 0
\quad
\text{in $L^2(0,T; L^2(\R^2 \times \R^3))$},
\\[0.3em]
&\sup_{0\leq t\leq T}
\|\w^{\e}_A(t)\|_{L^{p_1}(\R^2)}
\leq C\big(\|\w^{\e}|_{t=0}\|_{L^{p_1}(\R^2)},T\big),
\qquad
\w^{\e}_B \to 0
\ \ \text{in} \quad
L^\infty(0,T; L^{p_2}(\R^2)),
\quad \mbox{with rate},
\\[0.3em]
&\sup_{0\leq t\leq T}
\|(\rho^{\e}-3/2\ta^{\e})(t)\|_{L^2(\R^2)}
\leq
C\big(|(\rho^{\e}-3/2\ta^{\e})|_{t=0}\|_{L^2(\R^2)},T\big),
\\[0.3em]
&\mathbb{P}^{\perp}u^{\e},
\ \rho^{\e}+\ta^{\e}
\to 0
\quad
\text{in $ L^r (0,T; \dot{B}_{p_3,1}^{s+2(\frac{1}{p_3}-\frac{1}{2})+\frac{1}{r}}(\R^2))$ with rate}.
\end{align*}
This holds for any $0\leq n\leq 3$, $1\leq p_1\leq \infty$, $2\leq p_2\leq \infty$, $2<p_3<\infty$, 
$\frac{1}{r}\leq 2(\frac{1}{2}-\frac{1}{p})$, and for any $s\in(0,3)$.

\hide
\begin{align*}
&\frac{1}{\e^2} \frac{\nabla_x^n \AC{\P}F^{\e} }{\sqrt{M^\e}}\to   0   
  \ \ \text{in $L^2 (0,T; L^2 (\R^2 \times \R^3))$}, \cr 
&\sup_{0\leq t\leq T}\|\w^{\e}_A(t)\|_{L^{p_1}(\R^2)} \leq C(\|\w^{\e}|_{t=0}\|_{L^{p_1}(\R^2)},T) , \qquad \w^{\e}_B \to 0 \ \ \text{in} \quad L^\infty(0,T; L^{p_2}(\R^2)),\quad \mbox{with rate.} 
\cr 
&\sup_{0\leq t\leq T}\|(\rho^{\e}-3/2\ta^{\e})(t)\|_{L^2(\R^2)} \leq C\big(|(\rho^{\e}-\frac{3}{2}\ta^{\e})|_{t=0}\|_{L^2(\R^2)},T\big) \cr 
&\mathbb{P}^{\perp}u^{\e}, \ \rho^{\e}+\ta^{\e}    \to  0  \ \ \text{in $ L^r (0,T; \dot{B}_{p_3,1}^{s+2(\frac{1}{p_3}-\frac{1}{2})+\frac{1}{r}}(\R^2))$ with rate}
\end{align*}
for any $0\leq n\leq 3$, $1\leq p_1\leq \infty$, $2\leq p_2\leq \infty$, $2<p_3<\infty$ and $\frac{1}{r}\leq 2(\frac{1}{2}-\frac{1}{p})$ and for any $s\in(0,3)$.
\begin{equation}
\mbox{If} \quad \sup_{0<\e<\e_0}\frac{1}{\e} \bigg\|\frac{\nabla_x^n \AC{\P}F^{\e}}{\sqrt{M^\e}}\bigg|_{t=0}\bigg\|_{L^2_x(\R^2;L^2_v(\R^3))}<\infty , \quad \mbox{then} \quad 
\frac{1}{\e^2} \frac{\nabla_x^n \AC{\P}F^{\e} }{\sqrt{M^\e}}\to   0   
  \ \ \text{in $L^2 (0,T; L^2 (\R^2 \times \R^3))$}.
\end{equation}
for any $n\in\mathbb{N}\cup\{0\}$
\begin{equation}
\mbox{If} \quad \sup_{0<\e<\e_0}\|\w^{\e}|_{t=0}\|_{X_1(\R^2)}<\infty , \quad \mbox{then} \quad 
\begin{cases}
&\sup_{0\leq t\leq T}\|\w^{\e}_A(t)\|_{X_1(\R^2)} \leq C(\|\w^{\e}|_{t=0}\|_{X_1(\R^2)},T) , \\
&\w^{\e}_B \to 0 \ \ \text{in} \quad L^\infty(0,T; L^p(\R^2)), \quad \mbox{for} \quad  2\leq p<\infty \quad \mbox{with rate.} 
\end{cases}
\end{equation}
where $X_1= L^p$ for $1\leq p \leq \infty$ or $X_1=H^k$ for any $k\geq0$.

\begin{equation}
\mbox{If} \quad \sup_{0<\e<\e_0}\bigg\|\Big(\rho^{\e}-\frac{3}{2}\ta^{\e}\Big)\bigg|_{t=0}\bigg\|_{X_2(\R^2)}<\infty , \quad \mbox{then} \quad \sup_{0\leq t\leq T}\bigg\|\Big(\rho^{\e}-\frac{3}{2}\ta^{\e}\Big)(t)\bigg\|_{X_2(\R^2)} \leq C\bigg(\bigg\|\Big(\rho^{\e}-\frac{3}{2}\ta^{\e}\Big)\bigg|_{t=0}\bigg\|_{X_2(\R^2)},T\bigg)
\end{equation}
where $X_2= L^p$ for $1\leq p \leq \infty$ or $X_2=H^k$ for any $k\geq0$.

\begin{align}
\mathbb{P}^{\perp}u^{\e}, \ \rho^{\e}+\ta^{\e}    \to  0  \ \ &\text{in $ L^r (0,T; \dot{B}_{p,1}^{s+d(\frac{1}{p}-\frac{1}{2})+\frac{1}{r}}(\R^d))$ with rate} ,\label{conv_irrorhota}\\
 \frac{3}{2} \theta^\e - \rho^\e  \to \frac{3}{2} \theta^E - \rho^E  \ \ &\text{in  $L^\infty (0, T; {H}^{k}(\R^d))$ },  \label{conv_entr}
\end{align}
for $2<p<\infty$ and $\frac{1}{r}\leq d(\frac{1}{2}-\frac{1}{p})$ and for any $s\in(0,3)$.
\unhide
\end{theorem}

(아래는 Sec 9에 다 적어놓은 Remark들 입니다. 뭘 가져올지 취사선택을 여기서 해야할것 같습니다)
\begin{remark}
가정 둘째줄은 훨씬 약하다.
볼티시티 $L^p_1$ 는 $H^k$로 대체되도 성립한다.
$\rho-3/2\ta$의 $L^2$ 도 $H^k$로 대체되도 성립한다.
$n$ 과 $s$도 사실 arbitrary 큰 수에 대해서도 성립한다. 
\end{remark}

\begin{remark}
G가 0으로 가는 sense가 moments closoure를 가능하게 해준다 
그런데, rho,u,ta 각각 모멘트 이퀘이션에서는 안되고,
rho-3/2ta 혹은 $\w$의 이퀘이션에서 가능하게 해준다.
\end{remark}

\begin{remark}
3d에서 각 초기조건을 rho,u,ta,G는 $H^3$에서, $\w$ 는 $H^2$에서 주면, 각 초기조건에 해당하는 norm을 propagate시킬수 있고, G도 똑같이 0으로간다 유한한 시간동안.
$p_3$와 $r$ relation은 바뀜.
$\sup_{0\leq t\leq T_*}\|\nabla_x(\rho^{\e},u^{\e},\ta^{\e})\|<\infty$ 가 되는 시간 $T_*$동안 Hydrodynamic limit을 보일 수 있다. 
\end{remark}

\begin{remark}
볼티시티가 $H^2$에서 더블익스포넨셜 그로우 하는것을 생각해보면,
rho,ta가 살아있으면 볼츠만 에너지 디시페이션는 트리플 익스포넨셜 그로우,
rho,ta가 죽어있으면 볼츠만 에너지 디시페이션도 더블 익스포넨셜 그로우 한다.
\end{remark}

\begin{theorem}\label{T.Info.Hk}[Informal statement of Theorem \ref{T.C.Hk}]
For $\Omega=\R^2$ or $\Omega=\R^3$, consider
\begin{align*}
\w_0\in H^2(\R^d), 
\quad (u_0-\bar{u}) \in L^2(\R^d), 
\quad \rho_0, \ta_0 \in H^3(\R^d), 
\quad \frac{1}{\e}\frac{\AC{\P}F_0}{\sqrt{\tilde{\mu}}}
\in H^{3}_x(\R^d ; L^2_v(\R^3)).
\end{align*}
Suppose that the mollified initial data satisfy the assumptions of Theorem~\ref{T.Info.global}. 
Then the following convergence results hold:
\begin{align*}
&\mathbb{P}u^{\e} - u^E \to 0  
\quad \text{in $L^\infty (0,T; H^3(\R^d))$},
\\[0.3em]
&\frac{3}{2} \theta^\e - \rho^\e  
\to 
\frac{3}{2} \theta^E - \rho^E  
\quad \text{in $L^\infty (0,T; H^{3}(\R^d))$},
\end{align*}
where the limit point 
$(\rho^E,u^E,\ta^E)$ is the unique classical solution of the incompressible Euler equations~\eqref{incompE} and the convergence holds with rate for $H^2$.
\end{theorem}

Yudovich convergence
\begin{definition}\label{D.weaksolYudo}
A functions $(u,\w)$ is called a \emph{weak solution} to the vorticity-stream formulation of the Euler equation corresponding to the initial data $\w_0 \in L^1(\R^2)\cap L^\infty(\R^2)$ if the following conditions hold:
\begin{itemize}
\item $\w\in L^\infty(0,T;L^1(\R^2)\cap L^\infty(\R^2))$.
\item $u={\bf K}\ast \w$ and $\div(u) = 0$ in the sense of distributions.
\item For every test function $\Lambda \in C^1(0,T;C^1_0(\R^2))$, 
\begin{align*}
\int_{\R^2}\Lambda(T,x)\w(T,x)dx-\int_{\R^2}\Lambda(0,x)\w(0,x)dx = \int_0^T\int_{\R^2} (\p_t\Lambda  + u\cdot \nabla_x \Lambda) \w dxdt.
\end{align*}
\end{itemize}
\end{definition}

\begin{theorem}\label{T.Info.Yudo}[Informal statement of Theorem \ref{T.C.Linf}]
For $\Omega=\R^2$, consider
\begin{align*}
\w_0\in L^\infty(\R^2) \cap L^1(\R^2).
\end{align*}
Suppose that the mollified initial data satisfy the assumptions of Theorem~\ref{T.Info.global}. 
Then the following convergence results hold:
\begin{align*}
&\w^{\e}_A \quad \rightarrow \quad \w^{E} \quad \mbox{in} \quad  L^\infty(0,T;L^p(\R^2)), \quad \mbox{for any} \quad 1\leq p < \infty \\
&\w^{\e}_B \quad \rightarrow \quad 0 \quad \mbox{in} \quad  L^\infty (0, T; L^p(\R^2)) \mbox{with rate}, 
\end{align*}
where the limit point $(u^E,\w^E)$ with $\w^E:=\nabla_x^{\perp}\cdot u^E$ is the \emph{unique weak solution to the vorticity--stream formulation} of the Euler equation~\eqref{incompE} with initial data $\w^E|_{t=0}=\w_0$. 
If we further assume that $\w_0 \in \mathcal{F}^s_2(\R^2)$ for some $0<s<1$, where $\mathcal{F}^s_2$ is defined in Definition~\ref{D.Fsdef}, then for any $2\leq p<\infty$, we have the following convergence rate.
\end{theorem}

DiPerna-Majda Convergence:

\begin{theorem}\label{T.Info.Lp}[Informal statement of Theorem \ref{T.C.Lp}]
For $\Omega=\R^2$, consider
\begin{align*}
\w_0\in L^p_c(\R^2) \cap H^{-1}_{loc}(\R^2).
\end{align*}
Suppose that the mollified initial data satisfy the assumptions of Theorem~\ref{T.Info.global}. 
Then the following convergence results hold:
\begin{align*}
&u^{\e}_A \quad \rightarrow \quad u^{E} \quad \mbox{strongly in} \quad  \begin{cases}L^2(0,T;L^2_{loc}(\R^2)), \quad &\mbox{for} \quad p>1, \\ 
L^r(0,T;L^r_{loc}(\R^2)), \quad \mbox{for} \quad 1\leq r<2 , \quad &\mbox{for} \quad p=1,
\end{cases}  \\ 
&\w^{\e}_A \quad \rightarrow \quad \w^{E} \quad \mbox{strongly in} \quad C(0,T;L^p(\R^2)), \quad \mbox{for} \quad 1\leq p<\infty.  
\end{align*}
The pair $(u^{E},\w^{E})$ satisfies the incompressible Euler equations in the renormalized sense as in Definition \ref{D.soluw}. In addition, The forcing part vanishes in the limit:
\begin{align*}
&\w^{\e}_B \quad \rightarrow \quad 0 \ \ \text{strongly in} \quad L^\infty(0,T; L^p(\R^2)), \quad \mbox{for} \quad  2\leq p<\infty \quad \mbox{with rate,} 
\end{align*}
and 
\begin{align*}
3/2\theta^\e - \rho^\e  \to 3/2\theta^E - \rho^E  \ \ &\text{in  $L^\infty (0, T; L^2(\R^2))$ }.
\end{align*}
the pair $(u^{E},3/2\ta^{E}-\rho^{E})$ satisfies the incompressible Euler equations in the renormalized sense as in Definition \ref{D.soluw}.
\end{theorem}

\hide
\newpage

\section{A New Quasi-Linear Method: Difficulties and Key Ideas}


\subsection{Global Validity: Takeaway 3--Macro-Micro decoupling}

Altogether, the energy–dissipation inequality \eqref{EDintro} shows that the growth of the solution is governed by 
$\|\nabla_x(\rho^{\e},u^{\e},\ta^{\e})\|_{L^\infty_x}$.
This quantity is controlled by decomposing the dynamics into the de-penalized macroscopic variables $\w^{\e}$ and $\mathfrak{s}^{\e}$ together with the acoustic variables. 
We emphasize that the macroscopic variables $\w^{\e}$ and $\mathfrak{s}^{\e}$ are almost decoupled from both the microscopic dynamics and the acoustic components, thanks to the recovery of the full dissipation scale from the microscopic estimates in \eqref{Geqn0}, combined with the favorable scaling provided by the Strichartz estimates for the acoustic system.
As a consequence, the growth of 
$\|\nabla_x\mathbb{P}u^{\e}(t)\|_{L^\infty_x}$ and $\|\nabla_x \mathfrak{s}^{\e}(t)\|_{L^\infty_x}$
can be controlled solely in terms of the time interval and the initial data through a suitable bootstrap argument, under which the small parameters $\e$ and $\kappa$ absorb the forcing contributions appearing in the energy and dissipation estimates.
The remaining components $\|\nabla_x\mathbb{P}^{\perp}u^{\e}(t)\|_{L^\infty_x}$ and $\|\nabla_x p^{\e}(t)\|_{L^\infty_x}$ 
are controlled through dispersive effects.


Consequently, estimate \eqref{EDintro} implies that, whenever
$\mathcal{E}_{tot}^{\mathrm N}(F_0^{\e})<\infty$, the Boltzmann energy may grow at most triple-exponentially in time.
Under the additional assumption $\sup_{\e>0}\|\nabla_x \mathfrak{s}^{\e}_0\|_{L^\infty_x} \leq \kappa$, this growth rate improves to a double-exponential bound.


Furthermore, the structure of \eqref{EDintro} allows the analysis to extend below the classical Yudovich regularity threshold.
By imposing an appropriate admissible blow-up condition on the initial data, the right-hand side of \eqref{EDintro} can be controlled by a sufficiently small fraction of $\kappa^{-1}$, which in turn enables the propagation of vorticity with merely $L^p$ initial data over arbitrarily long time intervals.

(찬우가할일: 렐라티브엔트로피메쏘드와 차이써주기)


\subsection{Quantitative convergence term by term}



Regarding convergence, we emphasize that our results provide not only qualitative convergence toward the incompressible limit, but also explicit convergence rates. 
A distinctive feature of the present analysis is that different components of the solution—namely the de-penalized variables, the acoustic variables, and the microscopic part—converge in distinct functional topologies, each governed by a different underlying mechanism.

For the de-penalized variables $\w^{\e}$ and $\mathfrak{s}^{\e}$, convergence toward solutions of the incompressible Euler equations is obtained directly at the level of the evolution equations. 
Indeed, the forcing terms appearing in \eqref{wrtaeqn} vanish in $L_T^2L_x^\infty$, so that the limiting dynamics are governed purely by the Euler structure. 
As a consequence, once uniform bounds on the initial data sequence are available, the convergence can be quantified.

In particular, if the initial data satisfy
$\mathbb{P}u^{\e}_0-u^E_0 \to 0$ in $H^3(\R^d)$ and the macroscopic and microscopic components remain uniformly bounded in the corresponding functional spaces, then the solenoidal velocity admits the convergence rate
\begin{align*}
\bega
\sup_{0 \leq t \leq T}&\|(\mathbb{P}u^{\e}-u^E)(t)\|_{H^2(\R^d)} \leq C_T\|(\mathbb{P}u^{\e}-u^E)(0)\|_{H^2(\R^d)}  + C\Big(\e^{\frac{(d-1)-}{4}-} +\kappa^{\frac{1}{2}-}\Big),
\enda
\end{align*}
where $u^E$ denotes the unique classical Euler solution associated with $u^E_0$.
Here the rate $\e^{\frac{(d-1)-}{4}-}$ originates from dispersive Strichartz estimates for the acoustic variables, while the contribution $\kappa^{\frac12-}$ reflects the microscopic dissipation $\kappa^{\frac{1}{2}}\mathcal{D}^{\frac{1}{2}}$. 
This convergence rate becomes possible due to the refined dissipation scale recovered from the microscopic equation for $\AC{\P}F^{\e}$ in \eqref{Geqn0}.

In the two-dimensional Yudovich regime, where the initial vorticity satisfies $\sup_{\e>0}\|\w^{\e}_0\|_{L^\infty\cap L^1(\R^2)}<\infty $ and $\w^{\e}_0 \to \w^E_0$ in $L^2(\R^2)$, we obtain convergence toward the unique Yudovich solution in
$L^\infty(0,T;L^p(\R^2))$ for all $2\leq p<\infty$.
Moreover, assuming additional fractional regularity
$\sup_{\e>0}[\w^{\e}_0]_{\mathcal{F}_2^{s}(\R^2)} <+\infty$ for $0<s<1$, where $\mathcal{F}^s_2$ is defined in Definition~\ref{D.Fsdef},
the convergence becomes quantitative and yields the rate
\begin{align*}
\sup_{0 \leq t \leq T}&\|(\w^{\e}-\w^E)(t)\|_{L^p(\R^2)} \les \bigg( \bigg(\|\mathbb{P}u^{\e}_0-u^{E}_0\|_{L^2_x}+C \Big(\e^{\frac{1}{4}-} +\kappa^{\frac{1}{2}-}\Big) \bigg)^{e^{-Ct}}\bigg)^{\frac{2}{p}\frac{se^{-CT}}{1+se^{-CT}}-},
\end{align*}
for all $2\leq p<\infty$. Here, $\w^{E}$ denotes the unique Yudovich solution associated with the initial data $\w_0^E$.

The acoustic variables, on the other hand, converge to zero through dispersive effects. As a consequence of the Strichartz estimates, they decay in suitable Besov spaces $L_T^{r}\dot{B}_{p,1}^{s}$ (see \eqref{Besovdef}), even when the initial acoustic component is nonvanishing.

Meanwhile, the microscopic component converges to zero due to strong dissipation induced by coercivity. More precisely, the scaled microscopic term $\frac{1}{\e}\AC{\P}F^{\e}(t)$ vanishes in the $|M^{\e}|^{-1}$-weighted $L_T^{2}H_x^{\mathrm N}L_v^{2}$ norm at rate $\e\kappa^{\frac12}$, even when the initial microscopic fluctuation is of Mach-number size, namely $\AC{\P}F_0^{\e}=O(\e)$.
Such convergence is achievable because the present framework permits controlled growth of high-order initial norms through an admissible blow-up condition.



Regarding the convergence rate, we observe a noteworthy phenomenon originating from the microscopic equation \eqref{Geqn0}. 
Up to the $H_x^{\mathrm N}$ level, the microscopic forcing term is controlled by 
$\kappa^{\frac12}\mathcal{D}^{\frac12}$, and the dissipation remains scale-independent provided that the initial macroscopic and microscopic data are uniformly bounded in $H_x^{\mathrm N}$. 
Otherwise, under admissible blow-up initial data, the dissipation becomes singular through small fractions of $\kappa^{-1}$, which explains the appearance of the reduced rate $\kappa^{\frac12-}$ in the general setting.

A crucial observation is that, at low frequencies, the microscopic component does not exhibit temporal growth but only propagates the initial data while preserving strong dissipation. 
As a consequence, the microscopic contribution generates the full convergence rate $\kappa^{\frac{1}{2}}\int_0^t\mathcal{D}^{\frac{1}{2}}(s)ds \leq C \kappa^{\frac{1}{2}} T^{\frac{1}{2}}$ for the Burnett functionals
$\frac{1}{\e^2}\nabla_x^{\ell}\mathbf r_{ij}^{\e}$ and
$\frac{1}{\e^2}\nabla_x^{\ell}\mathfrak q_j^{\e}$
appearing in \eqref{wrtaeqn}, provided that
$\sup_{\e>0}\frac{1}{\e}\|(\nabla_x^{\ell}\AC{\P}F^{\e}) |M^{\e}|^{-\frac{1}{2}}|_{t=0}\|_{L^2_{x,v}} < \infty$. 
This rate is consistent with the inviscid-limit scaling observed in the Navier–Stokes framework and shows that the convergence rate $\kappa^{\frac12}T^{\frac12}$ is ultimately governed by the regularity of the initial microscopic data.

\subsection{Difficulties/Ideas in the classes of DiPerna-Majda} 





A major difficulty for solutions in the DiPerna–Majda class (namely, when the vorticity belongs to $L^p(\R^2)$ for $1\leq p<\infty$) arises from the fact that the limiting Euler solution is no longer unique.
Despite this lack of uniqueness, we prove that the velocity and vorticity limits obtained from the Boltzmann equation satisfy the incompressible Euler equations in the renormalized sense introduced by DiPerna–Lions \cite{DiLi}.

When the initial vorticity sequence is uniformly bounded in $L^p_c(\R^2)$ for $1<p<\infty$, the associated velocity fields admit uniform $W^{1,p}(\R^2)$ control.
In this regime, we exploit the DiPerna–Lions theory of renormalized solutions for the continuity equation developed in \cite{DiLi}, which provides a notion of uniqueness determined by the underlying vector field.
Although the limiting Euler velocity itself may not be unique, one can extract, up to subsequences, a limit velocity field.
For this fixed limiting vector field, the continuity equation admits a unique renormalized solution.
This allows us to identify the limit of the vorticity sequence as the renormalized solution associated with the limit velocity field.
Moreover, by considering the continuity equation associated with this fixed limit velocity field, we establish strong convergence of the vorticity.

A key additional difficulty in our setting is that the vorticity generated from the Boltzmann equation is not uniformly controlled in $L^p$ for $1\le p<2$.
To overcome this issue and to send both the acoustic variables and microscopic contributions to zero in a suitable topology, we decompose the vorticity into an initial-data part and a forcing part.
The initial-data component converges to a renormalized solution, while the forcing part vanishes with an explicit convergence rate.
This refined analysis becomes possible due to the improved dissipation scale obtained from the microscopic equation \eqref{Geqn0}, which removes the loss originating from macro–macro interactions.
As a further consequence, we prove that the limiting velocity $u^{E}$ satisfies energy conservation, $\|u^{\e}\|_{L^2_x} = \|u_0\|_{L^2_x}$.

When the initial vorticity sequence is only uniformly bounded in
$L^1_c \cap H^{-1}_{\mathrm{loc}}(\R^2)$, the velocity field becomes a singular integral of an $L^1$ function.
In this case, we employ the theory of regular Lagrangian flows developed by Bouchut–Crippa \cite{B-Crippa}.
Their framework ensures uniqueness of the associated Lagrangian solution, which in turn coincides with the renormalized solution of the continuity equation, allowing us again to identify the Euler limit.

The above arguments crucially rely on the structure of the velocity equation derived from the Boltzmann dynamics.
By removing the penalized term in the vorticity equation and exploiting the evolution equation for the microscopic component $\AC{\P}F^{\e}$, we recover dissipation at the full hydrodynamic scale, which ultimately enables the DiPerna–Majda–type limiting procedure described above.





\subsection{Difficulties/Ideas in the classes Delort}

For even rougher initial data, we consider the case where the initial vorticity is a Radon measure with distinguished sign and the approximating sequence is uniformly bounded in $L^1(\R^2)$. In this regime, substantially new difficulties arise compared to the classical Euler theory.

A first obstruction is that the distinguished-sign structure, which plays a crucial role in Delort’s compactness framework \cite{Delort}, is no longer preserved along the Boltzmann evolution due to microscopic and acoustic effects.
A second difficulty concerns spatial decay. Delort’s argument relies on precise decay properties of the vorticity maximal function, whereas the Boltzmann vorticity equation contains microscopic contributions weighted by $1/\mathrm{P}^{\e}$, preventing direct control of spatial moments such as $\int |x|^2\w^{\e}dx$.

To overcome these difficulties, we introduce the Boltzmann momentum
$m^{\e}:= \frac{1}{\e}\mathrm{P}^{\e}\mathrm{U}^{\e}$
and its associated vorticity
$\bw^{\e} := \nabla_x^{\perp}\cdot m^{\e}$. 
At the level of this momentum vorticity, the microscopic contribution appears in divergence form,
$\nabla_x^{\perp}\nabla_x\cdot\mathbf{r}^{\e}$, allowing the weight $|x|^2$ to be removed by integration by parts. As a consequence, the moment $\int |x|^2\w^{\e}dx$ is conserved, restoring the key compactness mechanism underlying Delort’s theory.

To compensate for the loss of sign preservation, we further decompose the momentum vorticity into initial-data and forcing components, $\bw^{\e}=\bw^{\e}_A+\bw^{\e}_B$, so that the initial-data component retains the distinguished-sign structure.
By estimating the vorticity maximal function associated with $\bw^{\e}_A$ and exploiting the conservation property satisfied by the full momentum vorticity $\bw^{\e}$, we show that the initial-data part $\bw^{\e}_A$ satisfies the same sharp decay estimate for the vorticity maximal function, up to an additional forcing contribution.

Since this forcing term vanishes in a suitable topology, we obtain convergence of the momentum $m^{\e}$, whose limit generates a weak solution of the Euler equations.

\unhide


\newpage
.
\newpage

\hide
In this paper, we define the \emph{well-prepared initial data} if the quantity $|F^{\e}-M_{[1,0,1]}|$ goes to 0 faster than $\e$. For example, if $\frac{1}{\e}\AC{\P}F^{\e}|_{t=0}\to0$ in some topology, then this is well-prepared initial data. 
If we use the wave equation for the acoustic variables, then we need a control of the time derivative of initial data $\e\p_t \mathbb{P}^{\perp}u^{\e}_0$. 
However, the initial data corresponding to the scaled time derivative $\e\p_t$ is too much well-prepared.
(Equation)
This aspect is different from the low Mach number limit problems in [..,..] where $\e\p_t u^{\e}_0$ is not well-prepared at least in the sense of the scale. 

well-prepared를 완전히 피할수 있다. (여러가지 센스로)
$1/\e \AC{\P}F^{\e}_0$ 도 커도 됨
$\mathbb{P}^{\perp}u^{\e}_0$ 도 0이 아님.
시간미분을 안써도 됨
\unhide



\hide
The moment closure methods is constructing solutions to the Boltzmann equation by identifying a topology in which the singular components vanish, without relying on any a priori information on the limiting Euler solution.
Near a global Maxwellian $F^{\e}=\mu+\e\sqrt{\mu}f^{\e}$, one may expect that a suitable dissipation scale is enough to suppress the microscopic part on the moments equations \eqref{locconNew} to zero. 
\begin{align}\label{globalexpand}
\frac{d}{dt}\|f^{\e}\|_{L^2}^2 + \frac{1}{\kappa \e^2}\|\sqrt{\nu}(\II-\P_{\mu})f^{\e}\|_{L^2}^2 &\les \frac{1}{\kappa \e} \la \Gamma(f^{\e},f^{\e}), (\II-\P_{\mu})f^{\e}\ra.
\end{align}
where $\P_{\mu}$ is $5$-dimensional projection with respect to the basis $(1,v,|v|^2)\sqrt{\mu}$. 
(See the moments equations \eqref{locconNew} that the dissipation suppresses the microscopic part since $\|\frac{1}{\e^2}\nabla_x\cdot \mathbf{r}^{\e}\|\sim\frac{1}{\e}\|(\II-\P)f^{\e}\|\sim\kappa^{\frac{1}{2}}$.)
If such dissipation is available, then after canceling the penalized terms $\frac{1}{\e}\nabla_x\cdot u^{\e}$ and $\frac{1}{\e}\nabla_x(\rho^{\e}+\ta^{\e})$, we expect that the resulting equations for suitable macroscopic quantities such as $\omega^{\e}$ and $\rho^{\e}-\frac32\theta^{\e}$ are no longer singular.
However, despite its apparent simplicity, this closure has not been achieved at the incompressible Euler scaling. The main obstruction is the highly singular collision operator of order $1/(\e\kappa)$, which is more singular than that appearing in the incompressible Navier–Stokes limit. Obtaining sufficiently strong dissipation to overcome this difficulty is one of the central contributions of this paper and requires exploiting a structure that has not been used previously. For a more detailed discussion, we describe two main difficulties. \\
\textbf{Difficulty 1:}
At the level of a global Maxwellian in \eqref{globalexpand}, the estimates can not be closed at all, because the macro-macro interaction inside the nonlinear term has scaling factor $1/\sqrt{\kappa}$. 
Alternatively, following the approach of [Liu–Yu–Yang] and working near a local Maxwellian, the term $\mathcal{N}(M^{\e},M^{\e})$ is canceled; however, starting from second-order derivatives, one encounters a loss of $\kappa$ again:
\begin{align*}
\frac{d}{dt}\frac{1}{\eps^2}\int|\p^{\alpha}F^{\e}|^2|M^{\e}|^{-1} &+ \frac{1}{\eps^4\kappa}\int \nu |\pt^\al \AC{\P}F^{\e}|^2|M^{\e}|^{-1} \leq \frac{1}{\eps^2}\bigg|\int|\p^{\alpha}F^{\e}|^2|M^{\e}|^{-2}\left(\pt_t M^{\e} + \frac{1}{\eps}v\cdot\nabla_x M^{\e}\right)\bigg| \\
& +\frac{1}{\sqrt{\kappa}}\sum_{0< \beta< \alpha}\int \mathcal{N}\Big(\frac{1}{\e}\mathbf{P}(\p^{\beta}M^{\e}),\frac{1}{\e}\mathbf{P}(\p^{\alpha-\beta}M^{\e})\Big)\frac{1}{\eps^2\sqrt{\kappa}}\p^{\alpha}\AC{\P}F^{\e}|M^{\e}|^{-1}  + \cdots
\end{align*}
Both approaches—writing $F^{\e}=M^{\e}+\AC{\P}F^{\e}$ and then differentiating, or differentiating $F^{\e}$ first and subsequently decomposing it as $F^{\e}=M^{\e}+\AC{\P}F^{\e}$—lead to the same difficulty. 
Expanding near a local Maxwellian improves marginally the situation: a $1/\sqrt{\kappa}$ loss appears after the $H^2$ level. However, this indicates that the dissipation scale deteriorates as the regularity increases, and moment closure in terms of $(\rho,u,\theta)$ becomes impossible already from $H^1$, due to the relation $\partial_t(\rho^{\e},u^{\e},\theta^{\e})\sim \frac{1}{\e^2}\nabla_x\cdot \mathbf{r}^{\e}$. Which means that the microscopic part of the vorticity equation $\w^{\e}$ is singular.\\

\textbf{Difficulty 2:}
To overcome these two difficulties, we identify a new cancellation mechanism intrinsic to the equation, which we call macro–micro cancellation. In the moment equations, penalized terms are eliminated by combining energy estimates at carefully chosen scales. Throughout the analysis, we work with perturbations around a local Maxwellian, $F^{\e}=M^{\e}+\AC{\P}F^{\e}$

and focus exclusively on the equation satisfied by the microscopic component $\AC{\P}F^{\e}$.
\unhide

\subsection{Symmetric Hyperbolic System of Moments}

*타이핑 리우 논문 하고 비교 완전히 다르다. 이사람들은 아무것도 안했다. 
$G$ 의 방정식이 다르다. 컨서베이션 로를 썼다!!

나비어스톡스로 가지 않았다.

\subsection{Microscopic estimate and Key Cancellation}

$G$의 에스티멧이랑 SH 에스티멧 합치면 캔슬 (버넷 캔슬)

캔슬이 없으면, 스케일 잃는다. 나비어 스톡스를 안써도 된다. 
타이핑 (에러처리), 코랙션(라준현)

[Rmk: 볼츠만 방정식만 가지고, 에너지에스티메잍을 돌리면, $H^2_x$ 부터 에너지와 디시페이션에서 $\kappa$를 1개씩 잃어야한다. - 이럴경우 vorticty equation \& symmetric hyperbolic system 둘다 오일러 방정식 형태의 forcing term이 singular하다. 그래서 NS를 쓸수밖에 없다. 
$G$의 에너지 에스티메잍 ($1/M^{\e}$ 로 잘 나눈) + symmetric hyperbolic system을 해서 캔슬을 해야만, 최대 스케일의 에너지,디시페이션을 얻을 수 있다. ($\e^4\kappa$ dissipation). 그러면 Euler만 쓰고 닫을 수 있다.]

\subsection{Kinetic Vorticity formulation}
캔슬이 있었기때문에 나비어스톡스를 안써도 된다.

\subsection{Acoustic system and Well-Preparedness} 이는 큰 어려움인 페넌트레이션 텀을 이니셜 데이터의 웰프리페어드니스 없이 해결하는데 결정적인 과정이다. 

Explanation about the well-preparedness and time derivative. (In the low Mach number limit $\e\p_t\rho \sim O(1)$. However, our case, BE $\p_tF_0 \sim \e\p_t(\rho_0,u_0,\ta_0) \sim O(1)$ is well-prepared because of microscopic part) ($\frac{1}{\e}G \sim 0 $ )
Once we take a norm to the scaled Boltzmann equation \eqref{BE}, 
\Be \notag
\eps \|\pt_t F^\e_0 \| +\|v\cdot\nabla_x F^\e_0 \|=\frac{1}{\kappa\eps}\|\mathcal{N}(F^\e_0,F^\e_0)\|.
\Ee
\begin{align*}
\nabla_x \cdot u &\approx \e \p_t \rho, \quad
\nabla_x(\rho+\ta) \approx \e \p_t u.
\end{align*}

\subsection{High Moments Control}
acoustic wave + STR estimate, 다이벌젼스 꽤 크다. 
시간 미분 필요 없다. 

[Rmk: incompressibility를 보일때: Dachin method의 장점. Wave equation을 쓰면 시간미분을 필히 써야, $\div(u)$ 와 $\nabla_x(\rho+\ta)$ 의 convergence를 보일 수 있다. Danchin은 Transport type \eqref{varsigeqn}을 썼기때문에, 시간미분을 안써도 된다.]

Note that $\div(u^{\e})$ and $(\rho^{\e}+\ta^{\e})$ satisfy the linear wave equation with forcing term:
\begin{align}
\bega
\begin{cases} \p_t^2 u - \Delta_x u =0 , \cr 
(u,\p_tu)|_{t=0} = u_0, u_1, \end{cases}
\enda
\end{align}
If we use the linear wave equation, we will have dispersion term $e^{\pm it|\xi|}$ with the time derivative of the initial data $\hat{\p_tu_0}$:
\begin{align}
\bega
u(t,x) &= \mathcal{F}^{-1}(e^{it|\xi|}\varphi) \ast \mathcal{F}^{-1}\bigg(\frac{1}{2}\bigg(\hat{u}_0(\xi)+\frac{1}{i|\xi|}\hat{\p_tu_0}(\xi)\bigg)\varphi \bigg) \cr 
&+ \mathcal{F}^{-1}(e^{-it|\xi|}\varphi) \ast \mathcal{F}^{-1}\bigg(\frac{1}{2}\bigg(\hat{u}_0(\xi)-\frac{1}{i|\xi|}\hat{\p_tu_0}(\xi)\bigg)\varphi \bigg)
\enda
\end{align}
On the other hand, if we use the transport equation \eqref{varsigeqn} (\cite{Danchin} type), we have the dispersion $e^{it|\xi|}$ without time derivative of the initial data:
\begin{align*}
\bega
\p_t\left( \begin{array}{c} \hat{\varsigma}_1 \\ \hat{\varsigma}_2 \end{array} \right) = i |\xi|  \left( \begin{array}{cc} 0 & -1 \\ 1 & 0  \end{array} \right) \left( \begin{array}{c} \hat{\varsigma}_1 \\ \hat{\varsigma}_2 \end{array} \right) 
, \qquad J:= \left( \begin{array}{cc} 0 & -1 \\ 1 & 0  \end{array} \right).
\enda
\end{align*}
Let $\hat{U}(t,\xi):= (\hat{\varsigma}_1(t,\xi),\hat{\varsigma}_2(t,\xi))^T$. Then $\hat{U}(t,\xi)$ satisfies 
\begin{align*}
\bega
\hat{U}(t,\xi) = e^{it|\xi|J}\hat{U}(0,\xi) = \left( \begin{array}{cc} \cos(t|\xi|) & -\sin(t|\xi|) \\ \sin(t|\xi|) & \cos(t|\xi|) \end{array} \right) \hat{U}(0,\xi).
\enda
\end{align*}

[Rmk: purely spatial dervatives를 사용할때, 3D local-in-time uniform boundedness를 보이기 위해서는, $\div(u^{\e})$ 와 $\rho^{\e}+\ta^{\e}$가 0으로 간다는것을 사용하지 않아도 된다. 먼저 uniform boundedness를 보이고, $\div(u^{\e})$ 와 $\rho^{\e}+\ta^{\e}$가 0으로 가는것을 보일수 있다. 이때, 에너지에서 필요한 최소한의 미분갯수는 (\mathrm{N},\mathrm{N}+1) = (3,4) 이다. 2D global-in-time uniform boundedness를 보이기 위해서는, $\|\nabla_x\div(u^{\e})\|_{L^\infty_x}$와 $\|\nabla_x^2(\rho^{\e}+\ta^{\e})\|_{L^\infty_x}$가 0으로 간다는 것을 사용해서 $\|\w^{\e}\|_{H^2_x}$를 컨트롤 해야한다. 그래서 Strichartz estimate때문에 최소한의 필요한 미분갯수는 (\mathrm{N},\mathrm{N}+1) = (4,5) 이다. (space-time derivatives를 고려할경우 위설명과는 조금 달라집니다.)]

\subsection{Discussion} {\color{red} 아래빨간부분이 추가로 포함되어야함}
$\bullet$ 메쏘드에 대한 논의: 이 메쏘드의 해석학적 특징과 장점은 명확하다. 
\begin{itemize}
\item[1.] 쿼지 리니어 문제 공부: 결정적인 어려움, 로우프리퀀시 인스테빌리티 차단: 떨모다이나믹 인터액션에서 매크로파트가 완전히 탈락한다 $\Gamma(b\cdot v, b\cdot v)$가 폴싱하는 인스테빌리티를 구조적으로 캔슬시김.
\item[3.] {\color{red}최대 장점은 에너지 에스테잇이 가능하다는 점이다. 어떤 의미에선 레라티브 엔트로피와 에너지에스티멧의 장점을 동시에 누리는 방식이다 (단 구조적 장점을 최대한 이용하기 위해 쿼지문제가 되었다는 점에서 수학적 난이도는 올라갔다고 볼 수 있다. ) }
\item[4.] 모멘텀 시스템의 시메트릭 하이퍼볼릭 구조가 명확히 나타난다. (볼티시티 구조를 완전히 이용, 어쿠스틱 시스템의 에너지 에스티메잇 가능)

\item[5.] 디씨파티브 성질을 완전히 이용할 수 있다
{\color{red}[Rmk: low regularity $0\leq |\al|\leq \mathrm{N}-2$의 microscopic part는 time growth 하지 않다는 것 까지도 보일 수 있다. (initial data propagte 한다.)]}

$F^{\e}=\mu+\e\sqrt{\mu}f$ 로 decompose 하면, $\mathcal{V}^{\al}_{\ell}= \mathcal{E}_{tot}+\e^{3/2}\kappa(something)$ 이 된다. ($\eqref{L.Vdecomp}_1$).
$F^{\e}=M^{\e}+ \e\sqrt{\bar{\mu}}g$ 로 decompose하면 $\mathcal{E}_{tot}$ 이 날라간다.(스케일 좋아짐) 대신 $\p_tM^{\e}+\frac{v}{\e}\cdot\nabla_xM^{\e}$이 나와서 $0\leq |\al|\leq \mathrm{N}-2$ 까지밖에 못쓴다.

\item[6.] 기존의 에너지 에스티멧과 다른점
\begin{itemize}
    \item {\color{red} 이니셜 데이터 클래스를 구조적으로 제약을 줄였다. 인컴프레서블 조건은 스트릭할츠 적용, 초기치 시간 미분에 대약 제약을 극적으로 줄임 }
\item 이니셜 데이터 유니폼 바운드 조건이 굉장히 낮을때도 정량적 어널리시스가 가능하다.  {\color{red} 나비어스톡스를 쓰게되면 초기치의 시간미분 가정이 들어간다.  }

\item {\color{red} 파라볼릭 에스티메잇 없어 오로지 하이퍼볼릭 에스티메잇만으로 문제를 해결, 파라볼릭 에스티멧을 쓰면 (나비어스톡스 중간과정) 필연적으로 시간미분을 사용해야하고, 그 결과 이니셜 데이터의 웰프리페어드니스를 가정해야한다. }
\end{itemize}
\end{itemize}

\subsection{Problems of interest and their Challenges}

In this paper, we are interested in the hydrodynamic limit toward the incompressible Euler equations. The diffusive limit typically involves a rescaling of time and space to capture the slow evolution of macroscopic quantities in high Reynolds number regime:  
  \Be\begin{split}
  \textbf{St} = \e = \textbf{Ma}, \ \ & \textbf{Kn} = \kappa \e, \ \ \text{for $\kappa = \kappa (\e) \to 0$ as $\e \to 0$.}
\end{split}  \Ee

more singular scaling (vs NS)

No priori information about fluid

target equation does not have uniqueness

Lack of smoothness

\hide
\subsection{
Key Ideas, New Framework and Results}

no use of approximation of local Maxwellian.

Liu Yang Yu: not the top order,

Let us first discuss \textit{validating the closure of fluid models} for finitely many macroscopic quantities derived from the Boltzmann equation, which inherently leads to an unclosed infinite hierarchy of moment equations. Macroscopic quantities naturally emerge from the microscopic description of states such as mass moment $ \mathrm{P}^{\e} 
$ (an uppercase `rho'), velocity moments $\mathrm{U}^{\e}
$, and temperature moment $\mathrm{\Theta}^{\e}$ via

Here, monatomic ideal gases holds for the internal energy law $\mathrm{E}^{\e} = \frac{3}{2}k_B\mathrm{\Theta}^{\e}$ with the Boltzmann constant $k_B$.

This system is unclosed as its forcing terms is component of microscopics: $\mathbf{r}_{ij}^{\e} = \langle \mathfrak{R}^{\e}_{ij},
F  \rangle,  
\mathfrak{q}_j^{\e} = \langle \mathcal{Q}^{\e}_j , 
F \rangle$. Here $\mathfrak{R}^{\e}_{ij}=(v_i-\mathrm{U}^{\e}_i)(v_j-\mathrm{U}^{\e}_j)-\frac{|v-\mathrm{U}^{\e} |^2}{3}\delta_{ij},$ and $ \mathcal{Q}^{\e}_i = (v_i-\mathrm{U}^{\e}_i)\frac{(|v-\mathrm{U}^{\e}|^2-5k_B\mathrm{\Theta}^{\e})}{2}.$ Specifically, this entails deriving macroscopic fluid equations for $(\mathrm{P}, \mathrm{U}, \mathrm{\Theta})$ from the mesoscopic Boltzmann equation as the Knudsen number vanishes under appropriate parameter scalings.

Besides the Knudsen number and Strouhal number, an additional dimensionless parameter, Mach number ($\mathbf{Ma}$), plays an important role in this context, which quantifies the velocity magnitude relative to a reference Maxwellian: $\mathbf{Ma} \sim  F - M_{[1,0,1]}$. In fact, the ratio $\mathbf{Ma}/\mathbf{St}$ determine the size of convection in the system \eqref{loccon}. Moreover, the Mach number and Knudsen number together determine the Reynolds number, which characterizes macroscopic viscosity, via the von K\'arm\'an relation ($\text{Reynolds number} \sim  {\textbf{Ma}} / {\textbf{Kn}}$).

Second, the hydrodynamic limit addresses the convergence toward the associated local Maxwellians $M_{[ \mathrm{P},\mathrm{U},\mathrm{\Theta}]}$ with the distribution function satisfying $F= M_{[ \mathrm{P},\mathrm{U},\mathrm{\Theta}]}  + o(\textbf{Ma})$ under appropriate parameter scalings as $\textbf{Kn}\rightarrow 0$. In other words, the major main in this context is to prove \textit{asymptotic stability in the scaling limit} of an associate local Maxwellin. 
 
The connection between these macroscopic quantities and the microscopic dynamics is established through the local Maxwellian (a.k.a. $F-$Maxwellian):

$\sigma_L$ plays an important role in quantifying the convergence of the solution toward a Maxwellian as time grows at a fixed scale. However, in the context of our problem—where the main focus is on the $\e$-dependence of the size of $\AC{\P}F^{\e}$—it becomes a relatively secondary issue, essentially corresponding to determining the constant in the Knudsen number in the $\e$-scale limit. Therefore, we do not address the problem of quantifying $\sigma_L$ here.

\hide
\begin{align}\label{Gform_intro}
\bega
\AC{\P}^{\e}F&=-\textbf{Kn}  \bigg\{ \sum_{i,j}\frac{\pt_i \mathrm{U}^{\e}_j}{k_B\mathrm{\Theta}^{\e}} \mathcal{L}^{-1} (\mathfrak{R}^{\e}_{ij}M^{\e})+\sum_j \frac{\pt_j\mathrm{\Theta}^{\e}}{k_B|\mathrm{\Theta}^{\e}|^2}\mathcal{L}^{-1} (\mathcal{Q}^{\e}_j M^{\e}) \bigg\}\cr 
&\quad -(\mathcal{L}^{\e} )^{-1}\lw\{\textbf{Kn}  \lw(\textbf{St} \pt_t \AC{\P}^{\e}F+\AC{\P}^{\e}(v\cdot\nabla_x \AC{\P}^{\e}F)\rw)-\mathcal{N}(\AC{\P}^{\e}F,\AC{\P}^{\e}F)\rw\} .
\enda
\end{align}
\unhide

Moreover, this analysis ensures that the derived macroscopic equations are consistent with the microscopic dynamics, reinforcing their physical relevance and predictive accuracy in the hydrodynamic regime.

validity of hydrodynamic limit

The mathematical study of hydrodynamic limit focuses on two issues mainly. First, validification of closure in fluids for finitely many macroscopic quantities derived from Boltzmann equations which leads to an infinite hierarchy of moment equations where each moment depends on the next higher moment, creating an unclosed system. In our case, this corresponds to deriving macroscopic fluid equations for $(\mathrm{P}, \mathrm{U}, \mathrm{\Theta})$ from the mesoscopic Boltzmann equation as the Knudsen number vanishes ($\textbf{Kn}\rightarrow 0$) under appropriate parameter scalings. An important dimensionless parameter in this context is the Mach number ($\mathbf{Ma}$), which quantifies the velocity magnitude relative to a reference Maxwellian: $\mathbf{Ma} \sim  F - M_{[1,0,1]}$.

The closure problem for kinetic equations constitutes one of the fundamental challenges in statistical physics. Based on the kinetic theory of Maxwell and Boltzmann, Hilbert [1] asked a humble question: Is there a self-consistent way to derive the evolution of macroscopic variables directly from particle-based models, most prominently, the Boltzmann equation?
Classically, fluid models such as the Euler equations or the Navier–Stokes equations are derived from kinetic equations through expansions in Knudsen number—the famous Chapman–Enskog expansion [2]—and were shown to be consistent with well-established local fluid models in the high-collision limit [3]. Higher-order terms in Knudsen number, such as the Burnett equation, however, exhibit nonphysical behavior, such as instabilities [4], and cannot be extended to the medium- or low-collisional regime [5]. Furthermore, there exist thermodynamic effects, such as thermally induced creep or Knudsen diffusion, which cannot be explained within the Navier–Stokes equations [6]. Nevertheless, a reduced description in terms of finitely many macroscopic quantities is highly desirable both in light of moments closures, as well as for efficient numerical methods

Motivation

convergence of equations

convergence of solutions

existence proof

\unhide


\hide
This relation provides a crucial connection between macroscopic fluid dynamics and the underlying kinetic theory, encapsulating the interplay of dimensionless parameters that govern the transition between scales. The process of deriving macroscopic equations from kinetic theory, known as the hydrodynamic limit, bridges the statistical description of particle interactions with the continuum description in fluid mechanics. 

In recent decades, substantial progress has been made in understanding this limit, driven by both theoretical interest and practical applications in physics and engineering.

\hide

As the Boltzmann equation suggested, $\textbf{St}$ is a temporal scale and $\textbf{Kn}$ is the mean-free path. The Mach number refer to a size of velocity such that
\Be
 \textbf{Ma}  \sim |U|  \sim  |F- M_{[1,0,1]}|.
 \Ee
 The Reynolds number is determined by the von K\'arm\'an relation:
 \Be\label{vKarman}
 \text{Reynold number} = \frac{\textbf{Ma}}{\textbf{Kn}} 
 \Ee

To connect the Boltzmann equation to the fluid models, we rescale the equation. For the macroscopic length scale $L$ and time scale $T$, we introduce $
\hat{t} = \frac{t}{T}, \hat{x} = \frac{x}{L}, \hat{v} = \frac{v}{c}$ where $c$ is the speed of sound.
We define 
\begin{align*}
\hat{F}(\hat{t},\hat{x},\hat{v}) = \frac{L^3c^2}{\mathcal{N}}F(t,x,v)
\end{align*}
where $\mathcal{N}$ is the number of molecules in the volume $L^3$. 
\begin{align*}
\mathcal{Q}^{\e}(\hat{z},w) = \frac{1}{c\pi r^2} B(z,w), \quad \hat{z}= \frac{z}{c}
\end{align*}
Dimensionless form of the Boltzmann equation: 
\begin{align*}
\frac{L}{cT}\pt_t \hat{F}+\hat{v}\cdot\nabla_{\hat{x}} \hat{F}=\frac{\mathcal{N}\pi r^2}{L^2}\int_{\mathbb{R}^3}\int_{\mathbb{S}^2_+}\mathcal{Q}^{\e}(\hat{v}-\hat{v}_*,w)\big[\hat{F}(\hat{v}_*')\hat{F}(\hat{v}')+\hat{F}(\hat{v}')\hat{F}(\hat{v}_*')\big]dwd\hat{v}_*.
\end{align*}\unhide

\hide

 \subsection{Main Theorem}

The moment method involves deriving moment equations from the Boltzmann equation and demonstrating their convergence to the target fluid equations as the Knudsen number approaches zero. By examining these moment equations, one can establish the existence of solutions to the fluid equations through the Boltzmann framework. Key results in this area include the following. Bardos, Golse, and Levermore present a rigorous proof of the fluid dynamic limits of the Boltzmann equation. Their work highlights the convergence of solutions and lays a strong theoretical foundation for this approach. Later, Golse and Saint-Raymond refine the analysis of fluid dynamic limits. They focus on higher-order effects and provide insights into the interaction between kinetic and fluid scales. Guo extends the diffusive limit of the Boltzmann equation beyond the classical Navier-Stokes approximation. This work provides uniform estimates and compactness arguments, which are crucial for demonstrating convergence. Esposito, Guo, Kim, and Marra investigate stationary solutions to the Boltzmann equation within the hydrodynamic limit. Their research emphasizes the impact of boundary effects on the solutions and their convergence properties. Nishida explores the compressible Euler limit of the nonlinear Boltzmann equation, particularly addressing the regime up to the formation of shocks. This work bridges the gap between kinetic theory and compressible fluid dynamics.

\unhide


\hide
2) Previous methods

2-1) Moments method

Proving that the scale of solution guarantees the microscopic part $\nabla_x\cdot \mathbf{r}^{\e}$ and $\nabla_x\cdot \mathfrak{q}^{\e}$ goes to $0$ in some sense. 

By using such way, the hydrodynamic limit is proved in incompressible Euler equation [Saint Raymond], incompressible navier-Stokes [Golse, Yan], [Bardos,Golse,Levermore], compressible Euler before shock [Nishida].

2-2) Asymptotic expansion 
$F= \mu + \sqrt{\mu}(\e f_1 + \e^2f_2 + \cdots \e^kf_k)$

3) Our strategy and what is different

4) key difficulty and idea

Two main terms. a) Coming from the commutator of $1/M^{\e}$. which produce $\p_t M^{\e} + \frac{v}{\e}\cdot\nabla_xM^{\e}$. Forcing $\|\nabla_x u\|_{L^\infty_x}$, $\|\nabla_x \ta\|_{L^\infty_x}$ be order constant and necessity of high order moments estimate $|v|^2F$ and $|v|^3F$.

High order Boltzmann energy estimate is singular.
 
Local conservation laws are singular. We should use $\w$ and $\rho-\frac{3}{2}\ta$.

Microscopic part in the conservation laws is singular. We should expand it getting one $\e$ and lossing one regularity with respect to $t$ and $x$.
\unhide


\hide
\begin{lemma}[Local conservation laws] We have 
\end{lemma}
\begin{proof}
From \eqref{ABdef}, we can compute that 
\Be\begin{split}\notag
\int_{\R^3} v_i v_j  F  dv & =\int_{\R^3} \left( \mathfrak{R}^{\e}_{ij}   - \mathrm{U}^{\e}_i \mathrm{U}^{\e}_j  
 +  \mathrm{U}^{\e}_i v_j   +  \mathrm{U}^{\e}_j v_i   
+\delta_{ij}  \frac{|v-\mathrm{U}^{\e} |^2}{3}  \right)  F dv \\
& = \mathbf{r}_{ij}^{\e}  - \mathrm{P}^{\e}\mathrm{U}^{\e}_i \mathrm{U}^{\e}_j
+2\mathrm{P}^{\e}\mathrm{U}^{\e}_i \mathrm{U}^{\e}_j+ \frac{2\delta_{ij}}{3} \Big(\mathrm{P}^{\e} \mathrm{E}^{\e} + \frac{1}{2} \mathrm{P}^{\e} |\mathrm{U}^{\e}|^2\Big)- \frac{2 \delta_{ij}}{3} \mathrm{P}^{\e}  |\mathrm{U}^{\e}|^2 + \frac{\delta_{ij}}{3} \mathrm{P}^{\e} |\mathrm{U}^{\e}|^2 \\
& =  \mathbf{r}_{ij}^{\e}  +  \mathrm{P}^{\e}\mathrm{U}^{\e}_i \mathrm{U}^{\e}_j + \delta_{ij} \frac{2}{3} \mathrm{P}^{\e} \mathrm{E}^{\e},
\end{split}
\Ee
and 
\Be\begin{split}\notag
\int_{\R^3} v_i  \frac{|v|^2}{2} F  dv & =\int_{\R^3} \left( \mathcal{Q}^{\e}_{i}  +  \frac{ |\mathrm{U}^{\e}|^2- 5R \mathrm{\Theta}^{\e}  }{2} \mathrm{U}^{\e}_i  -   \frac{  |\mathrm{U}^{\e}|^2 - 5R \mathrm{\Theta}^{\e}  }{2}  v_i 
- \mathrm{U}^{\e}_i \mathrm{U}^{\e} \cdot v + \frac{\mathrm{U}^{\e}_i}{2}|v|^2 + \mathrm{U}^{\e} \cdot v v_i \right) F d v \\
&= \mathfrak{q}^{\e}_i
- \mathrm{U}^{\e}_i \mathrm{P}^{\e} |\mathrm{U}^{\e}|^2 + \mathrm{U}^{\e}_i \left(\mathrm{P}^{\e} \mathrm{E}^{\e} + \frac{1}{2} \mathrm{P}^{\e} |\mathrm{U}^{\e}|^2\right) + \mathrm{U}^{\e} \cdot \int v v_i F dv \\
& = \mathfrak{q}^{\e}_i + \mathrm{P}^{\e} \mathrm{E}^{\e} \mathrm{U}^{\e}_i - \frac{1}{2} \mathrm{P}^{\e} |\mathrm{U}^{\e}|^2 \mathrm{U}^{\e}_i + \sum_j  \mathrm{U}^{\e}_j ( \mathbf{r}_{ij}^{\e}  +  \mathrm{P}^{\e}\mathrm{U}^{\e}_i \mathrm{U}^{\e}_j + \delta_{ij} \frac{2}{3} \mathrm{P}^{\e} \mathrm{E}^{\e} )\\
& = \mathfrak{q}^{\e}_i + \frac{5}{3} \mathrm{P}^{\e} \mathrm{E}^{\e} \mathrm{U}^{\e}_i
+ \frac{1}{2} \mathrm{P}^{\e} |\mathrm{U}^{\e}|^2 \mathrm{U}^{\e}_i + \sum_j \mathbf{r}_{ij}^{\e} \mathrm{U}^{\e}_j .
\end{split}
\Ee
Combining with \eqref{LCL}, these lead \eqref{loccon}. 

\hide
 {\color{red} Do we need the whole following stuff? \bigg[\bigg[

Using \eqref{M-moment}, we have 
\beq\label{random2}
\bega
&\int \eps\pt_t M^{\e} \bmx 1\\
v_i-\mathrm{U}^{\e}_i\\
|v-\mathrm{U}^{\e}|^2
\emx +\int v_j \pt_jM^{\e} \bmx 1\\v_i-\mathrm{U}^{\e}_i\\|v-\mathrm{U}^{\e}|^2
\emx +\int v_j \pt_j \AC{\P}^{\e}F\bmx 1\\v_i-\mathrm{U}^{\e}_i\\|v-\mathrm{U}^{\e}|^2\emx =0
\enda
\eeq
For just one derivative ($\p=\p_t$ or $\p_{x_i}$), we have 
\Be 
\bega
\p M^{\e}&=\lw(\frac{\p \mathrm{P}^{\e}}{\mathrm{P}^{\e}}-\frac{3}{2} \frac{\p\mathrm{\Theta}^{\e}}{\mathrm{\Theta}^{\e}}+\frac{\p \mathrm{U}^{\e} \cdot (v-\mathrm{U}^{\e})}{k_B\mathrm{\Theta}^{\e}}+\frac{\p\mathrm{\Theta}^{\e}}{\mathrm{\Theta}^{\e}}\cdot \frac{|v-\mathrm{U}^{\e}|^2}{2k_B\mathrm{\Theta}^{\e}} \rw)M
\enda
\Ee
Hence, we get
\[\bega
\int \eps\pt_t M^{\e} \bmx 1\\
v_i-\mathrm{U}^{\e}_i\\
|v-\mathrm{U}^{\e}|^2
\emx 
&=\eps\bmx \pt_t \mathrm{P}^{\e} \\ \mathrm{P}^{\e} \pt_t \mathrm{U}^{\e}_i\\3k_B\p_t(\mathrm{P}^{\e}\mathrm{\Theta}^{\e})
\emx
\enda
\]
Next we substitute \eqref{M1deri} for $\p_jM$ to have 
\beq\label{random1}
\bega
&\int v_j \pt_j M^{\e} \bmx 1\\v_i-\mathrm{U}^{\e}_i\\|v-\mathrm{U}^{\e}|^2\emx =\int (v_j-\mathrm{U}^{\e}_j)\pt_j M^{\e}\bmx 1\\v_i-\mathrm{U}^{\e}_i\\|v-\mathrm{U}^{\e}|^2\emx+\int \mathrm{U}^{\e}_j \pt_j M^{\e} \bmx 1\\v_i-\mathrm{U}^{\e}_i\\|v-\mathrm{U}^{\e}|^2\emx\\
&= \left[ \left(\begin{array}{c} \mathrm{P}^{\e} \nabla_x\cdot\mathrm{U}^{\e}  \\ \nabla_x\mathrm{P}^{\e} k_B\mathrm{\Theta}^{\e} +\nabla_x\mathrm{\Theta}^{\e} \mathrm{P}^{\e} k_B \\ \nabla_x \cdot \mathrm{U}^{\e} (5\mathrm{P}^{\e} k_B\mathrm{\Theta}^{\e} ) \end{array} \right) + \left(\begin{array}{c} \mathrm{U}^{\e} \cdot \nabla_x\mathrm{P}^{\e}   \\ \mathrm{P}^{\e}  \mathrm{U}^{\e} \cdot\nabla_x \mathrm{U}^{\e}  \\ \mathrm{U}^{\e} \cdot\nabla_x\mathrm{P}^{\e} (3k_B\mathrm{\Theta}^{\e} )+\mathrm{U}^{\e} \cdot \nabla_x\mathrm{\Theta}^{\e} (3\mathrm{P}^{\e} k_B) \end{array} \right) \right]
\enda
\eeq
Next we compute the last term in \eqref{random2}. We have
\[\bega
\int v_j \pt_j \AC{\P}^{\e}F\bmx 1\\v_i-\mathrm{U}^{\e}_i\\|v-\mathrm{U}^{\e}|^2\emx &=\int (v_j-\mathrm{U}^{\e}_j)\pt_j \AC{\P}^{\e}F\bmx 1\\v_i-\mathrm{U}^{\e}_i\\|v-\mathrm{U}^{\e}|^2\emx \cr 
&= \bmx 0\\ \pt_{x_j}\la \AC{\P}^{\e}F,\mathfrak{R}^{\e}_{ij}(v-\mathrm{U}^{\e})\ra_{L^2_v}\\ 2\pt_{x_j}\la \AC{\P}^{\e}F,\mathcal{Q}^{\e}_j(v-\mathrm{U}^{\e})\ra_{L^2_v}+2 \pt_j \mathrm{U}^{\e}_i\la \AC{\P}^{\e}F,\mathfrak{R}^{\e}_{ij}(v-\mathrm{U}^{\e})\ra_{L^2_v}
\emx
\enda
\]
Here we note that
\[\bega
&\int (v_j-\mathrm{U}^{\e}_j)\pt_j \AC{\P}^{\e}F|v-\mathrm{U}^{\e}|^2dv=\int \pt_j\lw((v_j-\mathrm{U}^{\e}_j)|v-\mathrm{U}^{\e}|^2\AC{\P}^{\e}F
\rw)dv-\int(v_j-\mathrm{U}^{\e}_j)\AC{\P}^{\e}F\cdot 2(v_i-\mathrm{U}^{\e}_i)(-\pt_j \mathrm{U}^{\e}_i)\\
&=\pt_j\lw(\int (v_j-\mathrm{U}^{\e}_j)|v-\mathrm{U}^{\e}|^2 \AC{\P}^{\e}F
\rw)+2 \pt_j \mathrm{U}^{\e}_i\int (v_j-\mathrm{U}^{\e}_j)(v_i-\mathrm{U}^{\e}_i)\AC{\P}^{\e}F
\enda
\]
In summary, we obtain
\begin{align}\label{consloc}
\bega
&\eps\p_t\mathrm{P}^{\e} +\nabla_x \cdot (\mathrm{P}^{\e}\mathrm{U}^{\e})=0 \cr
&\eps\mathrm{P}^{\e} \p_t\mathrm{U}^{\e} +\mathrm{P}^{\e}  \mathrm{U}^{\e}
\cdot\nabla_x \mathrm{U}^{\e} +k_B\nabla_x(\mathrm{P}^{\e} \mathrm{\Theta}^{\e})
+\nabla_x \cdot \mathbf{r}^{\e} =0 \cr  
&\eps3k_B\p_t(\mathrm{P}^{\e}\mathrm{\Theta}^{\e})
+3k_B \mathrm{U}^{\e}
\cdot\nabla_x(\mathrm{P}^{\e} \mathrm{\Theta}^{\e})+ 5k_B\nabla_x \cdot \mathrm{U}^{\e}
(\mathrm{P}^{\e} \mathrm{\Theta}^{\e}) \cr 
&\quad +2\nabla_x\cdot \mathfrak{q}^{\e} + 2\nabla_x \mathrm{U}^{\e} : \mathbf{r}^{\e}=0
\enda
\end{align}
Combining the first two equation, we can have the equation for $\p_t(\mathrm{P}^{\e} \mathrm{U}^{\e})$:
\begin{align}\label{Wa1}
&\eps \p_t(\mathrm{P}^{\e} \mathrm{U}^{\e}) +\mathrm{U}^{\e}\nabla_x \cdot (\mathrm{P}^{\e}\mathrm{U}^{\e})+\mathrm{P}^{\e}  \mathrm{U}^{\e}
\cdot\nabla_x \mathrm{U}^{\e} +k_B\nabla_x(\mathrm{P}^{\e} \mathrm{\Theta}^{\e})
+\nabla_x \cdot \mathbf{r}^{\e} =0 
\end{align}
\bigg]\bigg]
}
\unhide
\end{proof}
\hide
We need local conservation laws for $(\rho,u,\ta)$.
\begin{lemma} We have
\begin{align}\label{locconS}
\bega
&\e\p_t \rho + \e u \cdot \nabla_x \rho + \mathrm{P}^{\e} (\nabla_x \cdot u) = 0 \cr 
&\e\p_t u + \e u\cdot\nabla_x u +k_B\mathrm{\Theta}^{\e} \lw(\frac{\nabla_x\rho}{\mathrm{P}^{\e}} + \frac{\nabla_x \ta}{\mathrm{\Theta}^{\e}}\rw)+\frac{1}{\e \mathrm{P}^{\e}}\nabla_x \cdot \mathbf{r}^{\e} =0  \cr 
&\e\p_t \ta + \e u\cdot \nabla_x \ta 
+ \frac{2}{3}\mathrm{\Theta}^{\e} \nabla_x \cdot u  +\frac{2}{\eps3k_B\mathrm{P}^{\e}}\nabla_x\cdot \mathfrak{q}^{\e} + \frac{2}{\eps3k_B\mathrm{P}^{\e}}\nabla_x \mathrm{U}^{\e} : \mathbf{r}^{\e}=0
\enda
\end{align}
\end{lemma}
\begin{proof}
(1) Dividing $\eqref{consloc}_1$ by $\e$ gives 
\begin{align}
\e\p_t \rho + \e u \cdot \nabla_x \rho + \mathrm{P}^{\e} (\nabla_x \cdot u) = 0
\end{align}
(2) See momentum conservation law :
\begin{align}\label{}
\bega
&\eps\mathrm{P}^{\e} \p_t\mathrm{U}^{\e} +\mathrm{P}^{\e}  \mathrm{U}^{\e}
\cdot\nabla_x \mathrm{U}^{\e} +k_B\nabla_x(\mathrm{P}^{\e} \mathrm{\Theta}^{\e})
+\nabla_x \cdot \mathbf{r}^{\e} =0
\enda
\end{align}
We divide each side of $\eqref{consloc}_2$ by $\e\mathrm{P}^{\e}$: 
\begin{align}\label{}
\bega
\lw(\e\p_t u + \e u
\cdot\nabla_x u\rw) +\frac{k_B}{\mathrm{P}^{\e}}(\nabla_x\rho \mathrm{\Theta}^{\e}+ \mathrm{P}^{\e} \nabla_x \ta)
+\frac{1}{\e \mathrm{P}^{\e}}\nabla_x \cdot \mathbf{r}^{\e} =0 
\enda
\end{align}
We can write it 
\begin{align}\label{}
\bega
\lw(\e\p_t u + \e u\cdot\nabla_x u\rw) +k_B\mathrm{\Theta}^{\e} \lw(\frac{\nabla_x\rho}{\mathrm{P}^{\e}} + \frac{\nabla_x \ta}{\mathrm{\Theta}^{\e}}\rw)+\frac{1}{\e \mathrm{P}^{\e}}\nabla_x \cdot \mathbf{r}^{\e} =0 
\enda
\end{align}
(3) Consider $\eqref{loccon}_3$:
\begin{align}\label{}
\bega
&\eps3k_B\p_t(\mathrm{P}^{\e}\mathrm{\Theta}^{\e})
+3k_B \mathrm{U}^{\e}
\cdot\nabla_x(\mathrm{P}^{\e} \mathrm{\Theta}^{\e})+ 5k_B\nabla_x \cdot \mathrm{U}^{\e}
(\mathrm{P}^{\e} \mathrm{\Theta}^{\e}) \cr 
&\quad + 2\nabla_x\cdot \mathfrak{q}^{\e} + 2\nabla_x \mathrm{U}^{\e} : \mathbf{r}^{\e}=0
\enda
\end{align}
Make cancel on $\p_t\mathrm{P}^{\e}$ part: 
\begin{align}\label{}
\bega
&\eps3k_B\mathrm{P}^{\e}\p_t \mathrm{\Theta}^{\e}
+3k_B\mathrm{P}^{\e} \mathrm{U}^{\e}
\cdot \nabla_x \mathrm{\Theta}^{\e}+ 2k_B\nabla_x \cdot \mathrm{U}^{\e}
(\mathrm{P}^{\e} \mathrm{\Theta}^{\e}) \cr 
&\quad + 2\nabla_x\cdot \mathfrak{q}^{\e} + 2\nabla_x \mathrm{U}^{\e} : \mathbf{r}^{\e}=0
\enda
\end{align}
Divide each side by $\eps3k_B\mathrm{P}^{\e}$: 
\begin{align}\label{}
\bega
&\e\p_t \ta + \e u\cdot \nabla_x \ta 
+ \frac{2}{3}\mathrm{\Theta}^{\e} \nabla_x \cdot u  +\frac{1}{\eps3k_B\mathrm{P}^{\e}} 2\nabla_x\cdot \mathfrak{q}^{\e} + \frac{2}{\eps3k_B\mathrm{P}^{\e}}\nabla_x \mathrm{U}^{\e} : \mathbf{r}^{\e}=0
\enda
\end{align}
\end{proof}
\unhide
\unhide


\hide
\begin{theorem}\label{Info1}[Informal statement of Theorem \ref{Thm1}: Below Yudovich solution] Assume the initial relative entropy is bounded $\frac{1}{\e^2}\{ \mathcal{H}(F^{\e}_0) - \mathcal{H}(\mu)\}\leq C$, $\|(\rho^{\e},\ta^{\e})\|_{L^2_x}\leq$ and $\|\w^{\e}_0\|_{L^p_x} \leq C$ and let the other higher order quantities have certain rate of singular $\kappa$ until $5$ derivatives (with maximally one time derivative). Furthermore, if the initial data is bounded in high-moments 
\begin{align*}
&\sum_{0\leq |\al_x|\leq 5}\kappa^{\frac{(|\al_x|-1)_+}{2}}\left\|  w \p^{\al_x}  \left(\frac{F^\e_0 - \mu}{ \sqrt \mu }
\right) \right\|_{L^\infty_{x,v}} \leq \frac{1}{C},
\end{align*}
for a positive constant $C$. Then for each initial datum $F^\e_0$, and for arbitraly time $T_*$, we can choose $(\e,\kappa)$ which guarantees global validity \( T_* \to \infty \) as \( \e, \kappa \to 0 \) such that there exists a unique non-negative solution $F^\e (t,x,v)$ to the Boltzmann equation \eqref{BE} for $0\leq t\leq T_*$ satisfying the following uniform estimate: 
\begin{align*}
\sup_{0\leq t \leq T}\|(\rho^{\e},u^{\e},\ta^{\e})(t)\|_{L^2_x}^2 \leq C, \qquad \sup_{0\leq t \leq T}\|\w^{\e}(t)\|_{L^p_x} \leq C.
\end{align*}
\end{theorem}

For the two dimensional case, see Theorem \ref{T.2D.global}. 
For the three dimensional case, see Theorem \ref{T.3D.unif}.

\begin{theorem}\label{Info2}[Informal statement 2] 
Let $\Omega = \R^2$. We assume that the initial data share the same radial part, namely,
\begin{align*}
u^{\e}_0 = \widetilde{u}^{\e}_0+ \bar{u} + \mathbb{P}^{\perp}u^{\e}_0, \qquad
u^E_0 = \widetilde{u}^E_0+ \bar{u}.
\end{align*}
Let $F^{\e}$ be the solution to the Boltzmann equation constructed in Theorem \ref{Info1}, defined on the time interval $t\in[0,T]$, with initial data satisfying 
\[
(u^{\e}_0-\bar{u}) \in L^2(\R^2), 
\qquad (\rho^{\e}_0,\ta^{\e}_0) \in H^{-\mathfrak{j}}(\R^2) \quad \text{for any } \mathfrak{j}\in\R.
\]
Then, for any $T>0$, the following convergence holds:
\begin{itemize}
\item If $\w_0,\w^E_0 \in H^2(\R^2)$ and $\w^{\e}_0 \rightarrow \w^E_0$ strongly in $H^2(\R^2)$, then 
\begin{align*}
\sup_{0\leq t\leq T}\|(\w^{\e}-\w^{E})(t)&\|_{H^2_x} \to 0 ,\quad \mbox{as} \quad \e \to 0.
\end{align*}

\item If $\w_0,\w^E_0 \in L^2\cap L^\infty(\R^2)$ and $\w^{\e}_0 \rightarrow \w_0$ strongly in $L^2(\R^2)$, then 
\begin{align*}
\sup_{0\leq t\leq T}\|(\w^{\e}-\w^{E})(t)&\|_{L^p_x} \to 0,\quad \mbox{as} \quad \e \to 0 ,\quad \mbox{for all} \quad 1\leq p <\infty,
\end{align*}
where $\w^E$ is the unique weak solution of the Euler equation with respect to the initial data $\w^E_0$.

\item If $\w_0\in L^p_c(\R^2) \cap H^{-1}_{loc}(\R^2)$, and $\w^{\e}_0 \rightarrow \w_0$ strongly in $L^p(\R^2)$ then (up to subsequence)
\begin{align*}
&\w^{\e}_A \quad \rightarrow \quad \w^{E} \quad \mbox{strongly in} \quad C(0,T;L^p(\R^2)), \quad \mbox{for} \quad 1\leq p<\infty,
\end{align*}
In addition, $\nabla_x \cdot u^{E}(t,x) = 0$, and $u^{E}(t,x) = {\bf K} \ast \w^{E}(t,x)$ for ${\bf K}(x):= \frac{1}{2\pi}\frac{x^{\perp}}{|x|^2}$ and the pair $(u^{E}, \w^{E})$ satisfies the incompressible Euler equations in the renormalized sense.

\item If $\bw_0 \in \mathcal{M}(\R^2)$ with $\bw_0\geq0$ or $\bw_0\leq 0$ and $\mbox{supp}~ \bw_0 \subset \{x~|~|x|<R\}$, then (up to subsequence)
\begin{align*}
m^{\e}_A \quad &\rightarrow \quad m^{\#} \qquad \mbox{strongly in}, \qquad L^r([0,T]\times B_R(0)), \quad \mbox{for} \quad 1\leq r<2,
\end{align*}
and $m^{\#}$ is weak solution of the Euler equation.

\end{itemize}
\end{theorem}
\unhide

\hide
\subsection{Difficulties and Main Ideas}

If we start from the near-global Maxwellian regime, the most singular term will appear when the macroscopic and microscopic components collide:
\begin{align}
\frac{d}{dt}\|f\|_{L^2}^2 + \frac{1}{\kappa \e^2}\|\sqrt{\nu}(\II-\P_{\mu})f\|_{L^2}^2 &\leq \frac{1}{\kappa \e} \la \Gamma(f,f), (\II-\P_{\mu})f\ra \cr 
&\leq\frac{1}{\sqrt{\kappa}} \|\Gamma(\P_{\mu}f,\P_{\mu}f)\|_{L^2}\bigg\|\frac{1}{\sqrt{\kappa} \e}\sqrt{\nu}(\II-\P_{\mu})f\bigg\|_{L^2} + l.o.t.
\end{align}
To close the system, we require \(\|u_i u_j\|_{L^2_tL^2_x} + \|u\ta\|_{L^2_tL^2_x} + \|\ta^2\|_{L^2_tL^2_x} \approx \kappa^{\frac{1}{2}}\), which results in \(u = \ta = 0\). This indicates that closure is not feasible under this approach.
Once we follow the idea of \cite{LYY}, by splitting \(F\) near a local Maxwellian, then the collisions among the macroscopic components can be removed:
\begin{align*}
\p_t (M^{\e}+\AC{\P}F^{\e})+\frac{v}{\e}\cdot\nabla_x (M^{\e}+\AC{\P}F^{\e})=\frac{2}{\kappa\eps^2}\mathcal{N}(M^{\e},\AC{\P}F^{\e}) + \frac{1}{\kappa\eps^2}\mathcal{N}(\AC{\P}F^{\e},\AC{\P}F^{\e}).
\end{align*}
This approach is valid for the low-regularity part.
For the high-regularity regime, where \(|\al| \geq 1\), we obtain
\begin{align}\label{localexpand}
\bega
\frac{1}{2\e^2}&\frac{d}{dt} \int_{\Omega \times \R^3}|\p^{\al}F|^2|M^{\e}|^{-1}dvdx + \frac{1}{\kappa\e^4}\int_{\Omega \times \R^3}\mathcal{L}(\p^{\al}F)(\p^{\al}F)|M^{\e}|^{-1} dvdx \cr
&= -\frac{1}{2\e^2}\int_{\Omega \times \R^3} |\p^{\al}F|^2|M^{\e}|^{-2}\bigg(\p_t M^{\e} + \frac{v}{\eps}\cdot\nabla_x M^{\e}\bigg)dvdx \cr
&\quad+\frac{1}{\kappa\e^4}\sum_{0< \beta< \al} \binom{\al}{\beta} \int_{\Omega \times \R^3}\mathcal{N}(\P(\p^{\beta}M^{\e}),\P(\p^{\al-\beta}M^{\e}))\AC{\P}(\p^{\al}F)|M^{\e}|^{-1}dvdx + l.o.t.
\enda
\end{align}
(For the detailed proof, see Lemma \ref{Eest}).
We note that the macroscopic-macroscopic collision term now appears only after the second derivative. Since \(\p^{\al}M^{\e} \approx \e\p^{\al}(\rho,u,\ta)(1+|v|)^nM^{\e}\) for some power \(n>0\), we encounter the same problem as in \eqref{globalexpand} after the second derivative. Specifically, we still require \(\|\p^1u_i\p^1u_j\|_{L^2_tL^2_x}+\|\p^1u\p^1\ta\|_{L^2_tL^2_x}+\|\p^1\ta\p^1\ta\|_{L^2_tL^2_x} \approx \kappa^{\frac{1}{2}}\). However, the situation is more favorable than in \eqref{globalexpand} since \eqref{localexpand} becomes singular only after the \(H^2_x\) level.
Our key idea is to lose the \(\kappa\)-scale after the second derivative, thereby enabling us to hope for closing the estimate. Another issue arises in the second line of \eqref{localexpand}. When the transport operator acts on \(M^{\e}\), it produces terms like \(|v-\mathrm{U}^{\e}|^2M^{\e}\) and \(|v-\mathrm{U}^{\e}|^3M^{\e}\), which cannot be controlled by dissipation.

By multiplying \(\kappa\) to \eqref{localexpand}, we derive the following estimate:
\begin{align}\label{Eintro}
\frac{d}{dt}\mathcal{E}_{top} +\mathcal{D}_{top}  \leq \lw(\|\nabla_xu\|_{L^{\infty}_x}+\|\nabla_x\ta\|_{L^{\infty}_x}+\kappa\|\nabla_x(u,\ta)\|_{H^3_x}\rw)(\mathcal{E}_{top}+\mathcal{E}_{top}^{\frac{1}{2}}\mathcal{D}_{top}^{\frac{1}{2}}) + \|\nabla_xu\|_{L^\infty_x}\mathcal{V}_2 + \|\nabla_x\ta\|_{L^\infty_x}\mathcal{V}_3 + l.o.t,
\end{align}
where \(\mathcal{V}_{\ell}\) is a dissipation-like term with \((1+|v|^{\ell})\)-growth (the precise definition is given in \eqref{largev}). 
Since we lose a scale after \(H^2_x\), estimating \(\|\nabla_xu\|_{L^{\infty}_x}\) and \(\|\nabla_x\ta\|_{L^{\infty}_x}\) using embedding makes the inequality still singular. (The energy scale in \eqref{localexpand} yields \(\|\nabla_xu\|_{H^2_x}+\|\nabla_x\ta\|_{H^2_x} \les \kappa^{-\frac{1}{2}}\mathcal{E}_{top}^{\frac{1}{2}}\).)
Moreover, the inequality becomes singular in time even for a fixed \(\kappa\): \(\frac{d}{dt}\mathcal{E}_{top} \leq \kappa^{-\frac{1}{2}}\mathcal{E}_{top}^{\frac{3}{2}}\). Our key argument is that the macroscopic part can be estimated almost independently from the microscopic part. If \(\|\nabla_xu\|_{L^{\infty}_x}\) and \(\|\nabla_x\ta\|_{L^{\infty}_x}\) can be estimated in a manner similar to the incompressible Euler equation—meaning their growth depends only on time and initial data—then we can close the estimate in the Yudovich class. 
However, using the equations for \((\rho, u, \ta)\) directly is not feasible because they are too singular with respect to \(\nabla_x \cdot u\) and \(\nabla_x(\rho + \ta)\):
\Be\begin{split}\label{locconS}
\bega
&\p_t \rho + u \cdot \nabla_x \rho + \frac{1}{\e}\mathrm{P}^{\e} (\nabla_x \cdot u) = 0, \cr 
&\p_t u + u\cdot\nabla_x u +\frac{1}{\e}k_B\mathrm{\Theta}^{\e} \lw(\frac{\nabla_x\rho}{\mathrm{P}^{\e}} + \frac{\nabla_x \ta}{\mathrm{\Theta}^{\e}}\rw)+\frac{1}{\e^2 \mathrm{P}^{\e}}\nabla_x \cdot \mathbf{r}^{\e} =0  ,\cr 
&\p_t \ta + u\cdot \nabla_x \ta 
+ \frac{1}{\e}\frac{2}{3}\mathrm{\Theta}^{\e} \nabla_x \cdot u  +\frac{2}{\e^2 3k_B\mathrm{P}^{\e}}\nabla_x\cdot \mathfrak{q}^{\e} + \frac{2}{\e^2 3k_B\mathrm{P}^{\e}}\nabla_x \mathrm{U}^{\e} : \mathbf{r}^{\e}=0.
\enda
\end{split}
\Ee
To prove \(\frac{1}{\e}(\nabla_x \cdot u) \rightarrow 0\) and \(\frac{1}{\e}\nabla_x(\rho + \ta) \rightarrow 0\) as \(\e \rightarrow 0\), a higher level of well-preparedness than \(\e\) is required. Our goal is to close the system with the minimal well-preparedness. To achieve this, we choose the scaled time derivative \(\p_{\tilde{t}} = \e \varpi^{-1}(\e,\kappa)\p_t\), which introduces \(\varpi(\e,\kappa)\) as the level of well-preparedness. For the Yudovich class and below the Yudovich class, we set
\begin{align*}
&\varpi(\e,\kappa) = |\log \e|^{-1}, \quad \text{Case B: Below Yudovich class,} \cr
&\varpi(\e,\kappa) = \kappa^{\frac{3}{2}}, \quad \quad \text{Case C: Yudovich class.}
\end{align*}
This choice guarantees that \(\|\nabla_x \cdot u\|_{H^s_x} \rightarrow 0\) and \(\|\nabla_x(\rho + \ta)\|_{H^s_x} \rightarrow 0\) as \(\e \rightarrow 0\) for a certain regularity \(s\), while ensuring that \(\frac{1}{\e} \nabla_x \cdot u\) and \(\frac{1}{\e} \nabla_x(\rho + \ta)\) diverge to infinity.
Our main idea is to utilize the vorticity \(\w\) and the  specific entropy fluctuation \(\rho - \frac{3}{2}\ta\) to bypass the singularity. The equations for \(\w\) and \(\rho - \frac{3}{2}\ta\) are given by
\begin{align*}
&\p_t \w +u\cdot\nabla_x \w -\eta_0 k_B^{\frac{1}{2}} \kappa  \Delta_x \w = g_{\w}, \cr
&\p_t \Big(\rho-\frac{3}{2}\ta\Big) + u \cdot \nabla_x \Big(\rho-\frac{3}{2}\ta\Big)  -\eta_1 k_B^{\frac{1}{2}} \kappa  \Delta_x\Big(\rho-\frac{3}{2}\ta\Big) = g_{\mathfrak{s}}.
\end{align*}
(For the exact forms of these equations, please refer to Proposition \ref{meqnlem}.)
Then, we can estimate the macroscopic quantities. 
Here, we note that the microscopic parts \(\nabla_x\cdot \mathbf{r}^{\e}\) and \(\nabla_x\cdot \mathfrak{q}^{\e}\) need to be expanded. Since a scale is lost after \(H^2_x\), the microscopic part becomes singular in any Sobolev space \(H^s_x\) for \(s > 0\). By expanding the microscopic part, we can replace one derivative with one scale \(\e\). Consequently, we can obtain only an \(H^4_x\) estimate for the macroscopic part, relying on an \(H^5_x\) and \(H^1_tH^4_x\) estimate for the Boltzmann equation.
Notably, the forcing terms \(g_{\w}\) and \(g_{\mathfrak{s}}\) always include at least one component among \(\e\), \(\nabla_x \cdot u\), or \(\nabla_x(\rho + \ta)\). Therefore, if we choose a sufficient scale for \(\nabla_x \cdot u\) and \(\nabla_x(\rho + \ta)\) such that all forcing terms are absorbed into the scale, we can estimate the macroscopic quantities for \(g_{\w}\) and \(g_{\mathfrak{s}}\), nearly decoupled from the forcing terms. By applying the maximum principle and parabolic estimates, we can obtain estimates for \(u\) similar to those obtained for the incompressible Euler equation:
\begin{align}\label{ugrowintro}
\bega
\|\nabla_x u(t)\|_{L^\infty_x} 
&\les e^{(1+\|\w_0\|_{L^\infty_x})t}\lw(1+\ln^+\|\w_0\|_{H^2_x}\rw)\lw(1+\|\w_0\|_{L^\infty_x}\rw), \cr
\|\nabla_x u(t)\|_{H^3_x} 
&\les \exp\lw[e^{(1+\|\w_0\|_{L^\infty_x})t}\rw]  \Big(\|\w_0\|_{H^3_x}^2\Big)^{e^{(1+\|\w_0\|_{L^\infty_x})t}}.
\enda
\end{align}

Assuming \(H^3_x\) initial data for \((\rho, u, \ta)\), we can estimate the first part of \eqref{Eintro} in the Yudovich class. However, if we consider data below the Yudovich class, the inequality \eqref{Eintro} becomes implausible because \(\nabla_x u \notin L^\infty_x\). To address this, we select initial data below the Yudovich class with \(\|\w_0\|_{L^p} \leq C\) and define \(\w_0^{\e} := \mathcal{F}^{-1}(1_{|\xi| \leq 1/\sqrt{\kappa}} \mathcal{F}\w_0)\), where \(\hat{f} := \int_{\R^2} e^{-2\pi i \xi \cdot x} f(x) \, dx\). This sequence of initial data is singular with respect to \(\kappa\) in higher regularity spaces, satisfying \(\|\w^{\e}_0\|_{W^{s,p}} \leq C\kappa^{-s}\) for \(0 \leq s \leq 5\) and \(1 < p < 2\), or \(0 \leq s \leq 4\) and \(2 \leq p\).

We leverage the fact that the forcing terms \(g_{\w}\) and \(g_{\mathfrak{s}}\) always include at least one component among \(\e\), \(\nabla_x \cdot u\), or \(\nabla_x(\rho + \ta)\). Even though more well-preparedness is required than in the Yudovich class, choosing small \(\e\) well-preparedness, such as \(\varpi(\e,\kappa) = \e^{0.0001}\), \(|\log \e|^{-1}\), or \((\log (|\log \e| + 1))^{-1}\), can absorb the forcing terms. In this case, we can still achieve the estimate \eqref{ugrowintro}, even if the upper bound remains singular in \(\kappa\). Consequently, we can specify the growth of \(\|\nabla_x u\|_{L^{\infty}_x} + \|\nabla_x \ta\|_{L^{\infty}_x} + \kappa\|\nabla_x(u, \ta)\|_{H^3_x}\) in terms of time, initial data, and singular \(\kappa\). 
Thus, the energy and dissipation are bounded by singular \(\kappa\), providing no control in the limit sense. Nevertheless, we can achieve uniform boundedness of \(\|(\rho, u, \ta)(t)\|_{L^2_x}\) and \(\|\w(t)\|_{L^p_x}\) using the entropy and the vorticity equation. This is feasible because the forcing term in the vorticity equation contains \(\e\), or we select small \(\e\) well-preparedness that absorbs the forcing term arising from \(\nabla_x \cdot u\).

For the large velocity part, we decompose it near a global Maxwellian,

We also mention two technical difficulties. 
1) In the usual energy estimate of the Boltzmann equation, when microscopic and microscpic collide, we will have $\int_0^t \mathcal{D}_{top}^{\frac{3}{2}}(s)ds$. When we estimate it, if we use embedding for time $\|f\|_{L^\infty_t} \leq C_T\|f\|_{H^1_t}$, then we loss a scale $\e$ and get small well-preparedness $\varpi(\e,\kappa)$, since we are using scaled time derivative $\p_{\tilde{t}}=\e\varpi^{-1}(\e,\kappa)\p_t $. We should carefully estimate the scale when we use embedding for time. 
2) The second difficulty arise when we estimate the forcing terms $(g_{\w},g_{\mathfrak{s}})$ of the macroscopic quantities. We note that for the linear operator $\mathcal{L}H=-2\mathcal{N}(M^{\e},H)$, the inverse $\mathcal{L}^{-1}$ can be defined for each fixed $x$, since it only depends on $(\mathrm{P}^{\e},\mathrm{U}^{\e},\mathrm{\Theta}^{\e})$ (See Lemma \ref{Llem1}). So that, when we take a derivative to $\mathcal{L}^{-1}$, we should consider it inductively. We prove it in Proposition \ref{P.Xi} and Lemma \ref{L-1G2}.

Hydrodynamic..
Add references: Guo 2006, Golse, Saint-Raymond..

Hilbert expansion..

Difference from [Golse, Saint-Raymond].
Regularity, limit on vorticity equation.

One derivative can be changed by one $\e$. -This is not enough to close the system.
As the derivative increases, $\kappa$ is inevitably lost.

$5$-derivative is the least derivative that can gurantee the $L^\infty_x$ estimate of the vorticity from the maximum principle.

Filtering method: present three conservation laws.  more singular part.

How can we get $H^3_x$ estimate of $\ta$?
Energy conservation plus Boussinesq.


\unhide

This paper is organized as follows: \\
Section 2: Conservation laws, vorticity equation, temperature equation.\\
Section 3: Energy estimate of BE\\
Section 4: Wave equation - Boussinesq + div free \\
Section 5: $H^3$ estimate, $L^\infty$ estimate of macroscopic fiedls  \\
Section 6: Proof of the main theorem

\unhide

\hide
Curl for 2D: For 
\begin{equation}
\nabla_x \times F= \nabla^\perp \cdot F =  - \p_2 F_1  + \p_1 F_2, \qquad \nabla^\perp= \bigg( \begin{array}{c} -\p_2 \\ \p_1 \end{array} \bigg).
\end{equation}
Curl for 3D:
\begin{equation}
\nabla_x \times F= 
\end{equation}
\unhide 

\hide
$\bullet$ Sobolev norm and homogeneous sobolev norm
\begin{align*}
\bega
\|f\|_{H^s_x} := \bigg( \sum_{0\leq |\al_x| \leq s}\|\p^{\al_x}f\|_{L^2_x}^2 \bigg)^{\frac{1}{2}}, \qquad  \|f\|_{\dot{H}^s_x} := \bigg( \sum_{|\al_x| = s}\|\p^{\al_x}f\|_{L^2_x}^2 \bigg)^{\frac{1}{2}}
\enda
\end{align*}
$\bullet$ $L^\infty$ and $L^2$ norm in time
\begin{align*}
\bega
\|h(t)\|_{L^\infty_t} := \sup_{0\leq s\leq t}|h(s)|, \qquad \|h\|_{L^2_t} := \lw(\int_0^t |h(s)|^2 ds \rw)^{\frac{1}{2}}.
\enda
\end{align*}
\unhide

\unhide

\subsection*{C. New Quasi-Linear Method: Difficulties and Key Ideas}

Our contribution can be viewed as a \emph{moment--closure proof} of the incompressible Euler limit directly from
the Boltzmann dynamics: we rely on neither a prescribed \emph{a priori} Euler flow nor asymptotic expansions.
The main obstruction is the intrinsic singularity of the low--Mach/high--Reynolds scaling, on both the microscopic and
macroscopic levels. Any linearization around a \emph{fixed} Maxwellian that does not fully incorporate the exact
macroscopic fields $(\rho^{\varepsilon},u^{\varepsilon},\theta^{\varepsilon})$ leaves an uncancelled
macro--macro collision term at a \emph{supercritical} $\kappa$-scale. At the same time, the macroscopic moment
system carries genuinely singular penalization terms that obstruct a direct energy closure:
\begin{align}\label{local_intro}
\bega
 \p_t \rho^{\e} +u^{\e}\cdot\nabla_x\rho^{\e} + \frac{1}{\e}\nabla_x\cdot u^{\e} &=0 ,\cr
 \p_tu^{\e} + u^{\e}\cdot \nabla_x u^{\e} +\frac{1}{\e}\nabla_x(\rho^{\e}+\theta^{\e})
+\frac{1}{\e^2}\nabla_x \cdot \mathbf{r}^{\e} &= \text{\emph{l.o.t.}}, \cr
 \frac{3}{2}\big(\p_t \theta^{\e} +u^{\e}\cdot\nabla_x\theta^{\e}\big) + \frac{1}{\e}\nabla_x\cdot u^{\e}
+\frac{1}{\e^2}\nabla_x \cdot  \mathfrak{q}^{\e} &= \text{\emph{l.o.t.}}.
\enda
\end{align}
Classical approaches typically bypass such singularities by hinging on a prescribed limiting fluid solution via
Hilbert-type expansions and by imposing strict structural well-preparedness at both the micro and macro levels.
Achieving a theory beyond these restrictions is one of the main outcomes of the present work.

\subsubsection*{C.1. Micro--macro decomposition and quasi-linear framework}

To overcome the \emph{supercritical} macro--macro interaction, we linearize along the exact \emph{evolving local
Maxwellian manifold} determined by the Boltzmann macroscopic fields, rather than around a prescribed or frozen
equilibrium. We then choose the micro--macro decomposition \eqref{eq:macro_micro_intro} so that the
perturbation is precisely the geometric \emph{transverse fluctuation} to the manifold, motivated by \cite{DesvillettesVillani, Guo-NS}. In this formulation the
collision operator takes a quasi-linear form,
\[
\mathcal{L}(\AC{\P}F^{\varepsilon}) := -2\,\mathcal{N}\big(M^{\varepsilon},\AC{\P}F^{\varepsilon}\big).
\]

While this approach is effective at the \emph{low-frequency} level, the resulting structure does not
automatically propagate to higher regularity, since the underlying local-Maxwellian manifold varies in
space--time. Indeed, a \emph{projection--differentiation commutator obstruction} forces Burnett-type
contributions into the energy identity, destroying the dissipation--testing alignment. Each differentiation
incurs a $\kappa$-loss (at a supercritical scale); if accumulated, this loss breaks the energy closure and
precludes global (even local) convergence. 

To overcome this obstruction, we employ a coupled macro--micro energy method. We treat the macroscopic moment
system as a symmetric hyperbolic system with singular penalization and design an energy test that exploits its
internal cancellations. On the microscopic side, the quasi-linear equation yields dissipation at the
\emph{optimal scale} $\varepsilon^{2}\kappa^{1/2}$ \emph{without any $\kappa$-loss}, even at higher derivatives.
The main remaining difficulty is derivative loss, which reappears through Burnett-type terms. The key
observation is a structural cancellation: although the Burnett functionals couple the macroscopic and
microscopic dynamics at the level of the \emph{individual} equations, their leading contributions cancel once
the two estimates are \emph{combined} in a single coupled energy balance. This mechanism eliminates the final
$\kappa$-loss at principal order and allows us to close the estimates at the optimal scale. An additional difficulty is that the convection of the local Maxwellian generates both the cubic velocity weight
$|v|^3$ and gradients of the macroscopic fields. To handle this, we use a Caflisch-type local-Maxwellian
argument to obtain the low-frequency weighted $L^\infty_x$ control and exploit the $L^2_x$-integrability of the
macroscopic field gradients. However, the Caflisch argument cannot be used at the high-order level
because of derivative loss; this particular difficulty is absent in the Hilbert expansion framework. At the
high-order level we instead work with a global-Maxwellian $\mu$ formulation, and the resulting error
$M^\varepsilon-\mu$ (which is not $L^2_x$ for the infinite velocity energy case) is controlled in the Lagrangian $L^\infty$-estimate by exploiting the strong singular
dissipation. \hide

To overcome this obstruction, we employ a coupled macro--micro energy method. We treat the macroscopic moment
system as a symmetric hyperbolic system with singular penalization and design an energy test that exploits its
internal cancellations. On the microscopic side, the quasi-linear equation yields dissipation at the
\emph{optimal scale} $\varepsilon^{2}\kappa^{1/2}$ \emph{without any $\kappa$-loss}, even at higher derivatives.
An additional difficulty is that the convection of the local Maxwellian generates both the cubic velocity weight
$|v|^3$ and gradients of the macroscopic fields. To handle this, we use a Caflisch-type local-Maxwellian
argument to obtain the low-frequency weighted $L^\infty_x$ control and exploit the $L^2_x$-integrability of the
macroscopic field gradients. However, this local-Maxwellian argument cannot be used at the high-order level because of derivative loss; this particular difficulty is absent in the Hilbert expansion framework. At the high-order level we instead work with a global-Maxwellian formulation; the
resulting error $M^\varepsilon-\mu$ is then controlled in the Lagrangian $L^\infty$ estimate by exploiting the
strong singular dissipation of order $(\kappa\varepsilon^{2})^{-1}$. A further difficulty is derivative loss,
which reappears through Burnett-type terms. The key observation is a structural cancellation: although the
Burnett functionals couple the macroscopic and microscopic dynamics at the level of the \emph{individual}
equations, their leading contributions cancel once the two estimates are \emph{combined} in a single coupled
energy balance. This mechanism eliminates the final $\kappa$-loss at principal order and allows us to close the
estimates at the optimal scale.

To overcome this obstruction, we employ a coupled macro--micro energy method. We treat the macroscopic moment
system as a symmetric hyperbolic system with singular penalization and design an energy test that exploits its
internal cancellations. On the microscopic side, the quasi-linear equation yields dissipation at the
\emph{optimal scale} $\varepsilon^{2}\kappa^{1/2}$ \emph{without any $\kappa$-loss}, even at higher derivatives.
An additional difficulty is that the convection of the local Maxwellian generates both the cubic velocity weight
$|v|^3$ and the macroscopic field gradients. To handle this, we use a Caflisch-type
local-Maxwellian argument to obtain the low-frequency weighted $L^\infty_x$ control and use the integrability of the macroscopic field gradients. But this local-Maxwellian
argument is suffered by derivative loss. at the
high-order level we instead work with a global-Maxwellian formulation; the resulting error $M^\varepsilon-\mu$
is then controlled in the Lagrangian $L^\infty$ estimate by exploiting the strong singular dissipation of order
$(\kappa\varepsilon^{2})^{-1}$. A further difficulty is derivative loss, which reappears through Burnett-type
terms. The key observation is a structural cancellation: although the Burnett functionals couple the macroscopic
and microscopic dynamics at the level of the \emph{individual} equations, their leading contributions cancel
once the two estimates are \emph{combined} in a single coupled energy balance. This mechanism eliminates the
final $\kappa$-loss at principal order and allows us to close the estimates at the optimal scale.

To overcome this obstruction, we employ a coupled macro--micro energy method. We treat the macroscopic moment
system as a symmetric hyperbolic system with singular penalization and design an energy test that exploits its
internal cancellations. On the microscopic side, the quasi-linear equation yields dissipation at the
\emph{optimal scale} $\varepsilon^{2}\kappa^{1/2}$ \emph{without any $\kappa$-loss}, even at higher derivatives.
The remaining difficulty is derivative loss, which reappears through Burnett-type terms. The key observation is
a structural cancellation: although the Burnett functionals couple the macroscopic and microscopic dynamics at
the level of the \emph{individual} equations, their leading contributions cancel once the two estimates are
\emph{combined} in a single coupled energy balance. This mechanism eliminates the final $\kappa$-loss at
principal order and allows us to close the estimates at the optimal scale.

At the differentiated level, the transport of the local Maxwellian produces both the cubic velocity weight $|v|^3$ and differentiated macroscopic fields. While the integrability of the latter is a key ingredient in a local-Maxwellian perturbative description, our top-order energy estimates are instead carried out in the global-Maxwellian formulation, and a Caflisch-type local-Maxwellian argument is used only for the non-differentiated low-frequency $L^\infty_x$ control. 

\unhide
At the end of the macro--micro analysis one arrives at a Lyapunov inequality with a Gr\"onwall-type growth
mechanism of the schematic form
\begin{equation}\label{EDintro}
\text{(Energy)}+\int_0^t \text{(Dissipation)}
\;\lesssim\;
\int_0^t \|\nabla_x(\rho^\varepsilon,u^\varepsilon,\theta^\varepsilon) \|_{L^\infty_x}\,\text{(Energy)}    
\;+\;\textit{l.o.t.}.
\end{equation}
Thus the Beale–Kato–Majda-type quantity
\begin{equation}\label{Lipschitz}
\|\nabla_x(\rho^\varepsilon,u^\varepsilon,\theta^\varepsilon)\|_{L^1(0,T_*;L^\infty_x)}
\end{equation}
governs the lifespan of the energy method. In three space dimensions, for smooth data and sufficiently high
Sobolev regularity, Sobolev embedding controls this threshold by the energy itself, yielding local-in-time
construction; in two dimensions the same mechanism yields local-in-time control for smooth data, while
global-in-time bounds require additional structure, discussed next.

\subsubsection*{C.2. Global convergence at low regularity}

In two space dimensions, our goal is to obtain a global-in-time hydrodynamic limit for vorticity data beyond the
Yudovich class. At this level of regularity, two fundamental difficulties arise: obtaining a priori control of
the BKM-type quantity in terms of the initial data and the small-scale parameters is delicate, and uniqueness
of the limiting Euler dynamics is in general unavailable.

The key idea is to obtain the BKM-type control from the local conservation laws \eqref{local_intro}. The moment
system for $(\rho^\varepsilon,u^\varepsilon,\theta^\varepsilon)$ contains singular penalization
$\varepsilon^{-1}\nabla_x\!\cdot u^\varepsilon$ and $\varepsilon^{-1}\nabla_x(\rho^\varepsilon+\theta^\varepsilon)$,
as well as derivatives of microscopic forcing terms. To bypass the penalization, we work with variables that
eliminate the singular gradients: the kinetic vorticity $\omega^\varepsilon:=\nabla_x^\perp\!\cdot u^\varepsilon$
and the entropic fluctuation $\mathfrak{s}^\varepsilon:=\frac{3}{2}\theta^\varepsilon-\rho^\varepsilon$, together
with the acoustic variables $(\mathbb{P}^\perp u^\varepsilon,\rho^\varepsilon+\theta^\varepsilon)$. Assuming
dissipation control at the optimal scale, the BKM-type control for $u^\varepsilon$ follows from a
maximum-principle estimate for the kinetic vorticity equation,
\begin{equation}\notag
 \p_t\omega^{\e} + u^{\e}\cdot \nabla_x \omega^{\e} + (\nabla_x\cdot u^{\e})\,\omega^{\e}
 +\frac{1}{\e^2}\nabla_x^{\perp}\cdot(\nabla_x \cdot  \mathbf{r}^{\e})
 = \text{\emph{l.o.t.}}.
\end{equation}
The acoustic variables $(\mathbb{P}^\perp u^\varepsilon,\rho^\varepsilon+\theta^\varepsilon)$ satisfy a forced
half-wave system (see Lemma~\ref{L.locconP}). Strichartz estimates yield dispersive decay and quantitative
smallness in Besov norms, so the acoustic component converges to zero even when the initial acoustic modes do
not vanish.

Once the BKM-type quantity is controlled in terms of the initial data, we finally exploit the Lyapunov
inequality \eqref{EDintro} to recover dissipation at the optimal scale, at the expense of a super-exponential
growth factor in time stemming from the BKM-type control.

It is worth emphasizing that we do not invoke any parabolic estimates at the fluid level (e.g.\ via a
Navier--Stokes expansion). Such an approach typically requires controlling time derivatives of microscopic
quantities and therefore imposes additional well-preparedness conditions on the initial data. By contrast, our
argument relies solely on the underlying hyperbolic structure (Euler-type conservation laws) and avoids any
assumptions on time derivatives at $t=0$. 
\hide

\subsection*{B. New Quasi-Linear Method: Difficulties and Key Ideas} 
Our contribution can be viewed as a \emph{moment--closure proof} of the incompressible Euler limit directly from
the Boltzmann dynamics: we rely on neither a prescribed \emph{a priori} Euler flow nor asymptotic expansions.
The main obstruction is the intrinsic singularity of the low--Mach/high--Reynolds scaling, on both the microscopic and
macroscopic levels. Any linearization around a \emph{fixed} Maxwellian that does not fully incorporate the exact
macroscopic fields $(\rho^{\varepsilon},u^{\varepsilon},\theta^{\varepsilon})$ leaves an uncancelled
macro--macro collision term at the \emph{supercritical} $\kappa^{-1/4}$ scale. At the same time, the macroscopic
moment system carries genuinely singular penalization terms that obstruct a direct energy closure:
\begin{align}\label{local_intro}
\bega
 \p_t \rho^{\e} +u^{\e}\cdot\nabla_x\rho^{\e} + \frac{1}{\e}\nabla_x\cdot u^{\e} &=0 ,\cr
 \p_tu^{\e} + u^{\e}\cdot \nabla_x u^{\e} +\frac{1}{\e}\nabla_x(\rho^{\e}+\ta^{\e})
+\frac{1}{\e^2}\nabla_x \cdot \mathbf{r}^{\e} &= \text{\emph{l.o.t.}}, \cr
 \frac{3}{2}\big(\p_t \ta^{\e} +u^{\e}\cdot\nabla_x\ta^{\e}\big) + \frac{1}{\e}\nabla_x\cdot u^{\e}
+\frac{1}{\e^2}\nabla_x \cdot  \mathfrak{q}^{\e} &= \text{\emph{l.o.t.}}.
\enda
\end{align}

Classical approaches typically bypass such singularities by hinging on a prescribed limiting fluid solution via
Hilbert-type expansions and by imposing strict structural well-preparedness at both the micro and macro levels.
Achieving a theory beyond these restrictions is one of the main outcomes of the present work.

 \subsubsection*{B.1. Micro-Macro decomposition and Quasi-linear framework} 
To overcome the \emph{supercritical} macro--macro interaction, we linearize along the exact \emph{evolving local
Maxwellian manifold} determined by the Boltzmann macroscopic fields, rather than around a prescribed or frozen
equilibrium. We then choose the micro-macro decomposition \eqref{eq:macro_micro_intro} so that the perturbation is precisely the geometric \emph{transverse fluctuation} to the manifold. In this formulation the collision operator takes a quasi-linear form,
$\mathcal{L}(\AC{\P}F^{\varepsilon}) := -2\,\mathcal{N}\big(M^{\varepsilon},\AC{\P}F^{\varepsilon}\big)$, with
$M^\varepsilon$ in \eqref{M-def}.

While this approach is effective at the \emph{low-frequency} level, the resulting structure does not
automatically propagate to higher regularity, since the underlying local-Maxwellian manifold varies in
space--time. Indeed, a \emph{projection--differentiation commutator obstruction} forces Burnett-type contributions into the
energy identity, destroying the dissipation--testing alignment. Each differentiation costs a fatal $\kappa^{-1/2}$
factor; accumulating this loss breaks the energy estimates and precludes global (even local) convergence.

To overcome this obstruction, we employ a coupled macro--micro energy method. We treat the macroscopic moment system as a symmetric hyperbolic system with singular penalization,
and we design an energy test that exploits its internal cancellations. On the microscopic side, the quasi-linear
equation yields dissipation at the optimal $\varepsilon^{2}\kappa^{1/2}$ scale \emph{without any scale loss}, even at
higher derivatives. The remaining difficulty is derivative loss, which reintroduces a $\kappa^{-1/4}$ defect through Burnett-type
terms. The key observation is a structural cancellation: although the Burnett functionals couple the macroscopic
and microscopic dynamics at the level of the \emph{individual} equations, their leading contributions cancel
once the two estimates are \emph{combined} in a single coupled energy balance. This mechanism eliminates the
final $\kappa^{-1/4}$ loss (and, in effect, the Burnett contribution at principal order) and allows us to close
the estimates without a loss of scale.

At the end of the macro--micro analysis one arrives at a Lyapunov inequality with a Gr\"onwall-type growth
mechanism of the schematic form
\begin{equation}\label{EDintro}
\text{(Energy)}+\int_0^t \text{(Dissipation)}
\;\lesssim\;
\int_0^t \|\nabla_x(\rho^\varepsilon,u^\varepsilon,\theta^\varepsilon)(s)\|_{L^\infty_x}\,\text{(Energy)}(s)\,ds
\;+\;\textit{l.o.t.}.
\end{equation}
Thus the BKM-type quantity
\begin{equation}\label{Lipschitz}
\|\nabla_x(\rho^\varepsilon,u^\varepsilon,\theta^\varepsilon)\|_{L^1(0,T_*;L^\infty_x)}
\end{equation}
governs the lifespan of the energy method. In three space dimensions, for smooth data and sufficiently high
Sobolev regularity, Sobolev embedding controls this threshold by the energy itself, yielding local-in-time
construction in 2D and 3D.

\subsubsection*{B.2. Global convergence in low regularity}
In two space dimensions, our goal is to obtain a global-in-time hydrodynamic limit for vorticity data beyond the
Yudovich class. At this level of regularity, two fundamental difficulties arise: obtaining a priori control of
the BKM-type quantity in terms of the initial data and the small-scale parameters is delicate, and uniqueness
of the limiting Euler dynamics is in general unavailable.

The key idea is to obtain the BKM-type control from the local conservation laws \eqref{local_intro}. The moment
system for $(\rho^\varepsilon,u^\varepsilon,\theta^\varepsilon)$, however, contains singular penalization
$\varepsilon^{-1}\nabla_x\!\cdot u^\varepsilon$ and $\varepsilon^{-1}\nabla_x(\rho^\varepsilon+\theta^\varepsilon)$,
as well as derivatives of microscopic forcing terms. To bypass the penalization, we work with
variables that eliminate the singular gradients: the kinetic vorticity
$\omega^\varepsilon:=\nabla_x^\perp\!\cdot u^\varepsilon$ and the entropic fluctuation
$\mathfrak{s}^\varepsilon=\frac{3}{2}\theta^\varepsilon-\rho^\varepsilon$, together with the acoustic variables
$(\mathbb{P}^\perp u^\varepsilon,\rho^\varepsilon+\theta^\varepsilon)$. Assuming dissipation control at the
optimal scale, the BKM-type control for $u^\varepsilon$ follows from a maximum-principle estimate for the
kinetic vorticity equation,
\begin{equation}\notag
 \p_t\omega^{\e} + u^{\e}\cdot \nabla_x \omega^{\e} + (\nabla_x\cdot u^{\e})\,\omega^{\e}
 +\frac{1}{\e^2}\nabla_x^{\perp}\cdot(\nabla_x \cdot  \mathbf{r}^{\e})
 = \text{\emph{l.o.t.}}.
\end{equation}
The acoustic variables $(\mathbb{P}^\perp u^\varepsilon,\rho^\varepsilon+\theta^\varepsilon)$ satisfy a forced
half-wave system (see Lemma~\ref{L.locconP}). Strichartz estimates yield dispersive decay and quantitative
smallness in Besov norms, so the acoustic component converges to zero even when the initial acoustic modes do
not vanish.

Once the BKM-type quantity is controlled in terms of the initial data, we finally exploit the Lyapunov
inequality \eqref{EDintro} to recover dissipation at the optimal scale, at the expense of a super-exponential
growth factor in time stemming from the BKM-type control.

It is worth emphasizing that we do not invoke any parabolic regularization at the fluid level (e.g.\ via a
Navier--Stokes expansion). Such an approach typically requires controlling time derivatives of microscopic
quantities and therefore imposes additional well-preparedness conditions on the initial data. By contrast, our
argument relies solely on the underlying hyperbolic structure (Euler-type conservation laws) and avoids any
assumptions on time derivatives at $t=0$.

\hide
, which yields
quantitative control of $\AC{\P}F^\varepsilon$, provided the time-integrated Lipschitz threshold is controlled.
In turn, this microscopic control bounds the Burnett-type forcing terms in the macroscopic system and allows the
coupled macro--micro estimates to close.

through the \emph{kinetic vorticity} and the \emph{maximum principle}, complemented by potential-theoretic estimates, while the \emph{acoustic variables} are shown to decay by dispersive effects.

Controlling this forcing at the correct scale is therefore a central ingredient, but it must be combined with a
robust mechanism to handle the singular low--Mach/high--Reynolds structure (penalization and derivative loss).

The required $L^\infty_x$ control of $\nabla_x u^\varepsilon$ is then
obtained by combining a maximum-principle estimate for the kinetic vorticity equation with potential-theoretic
bounds, while the acoustic variables decay by dispersive effects.

A key idea is to exploit the Lyapunov inequality \eqref{EDintro} to recover dissipation at the optimal scale, thereby obtaining quantitative control of the microscopic fluctuation $\AC{\P}F^\varepsilon$ (once the time-integrated Lipschitz threshold is controlled). This microscopic control, in turn, bounds the Burnett-type forcing terms that enter the macroscopic equations as a forcing term.

 \medskip In the global analysis the main object is the macroscopic system. However, as seen from \eqref{locconNew}, the evolution of $(\rho^\varepsilon,u^\varepsilon,\theta^\varepsilon)$ contains \emph{singularly penalized terms}, notably $\varepsilon^{-1}\nabla_x\!\cdot u^\varepsilon$ and $\varepsilon^{-1}\nabla_x(\rho^\varepsilon+\theta^\varepsilon)$. To bypass this obstruction, we work with macroscopic variables free of the leading penalization: the kinetic vorticity $\omega^\varepsilon:=\nabla_x^\perp\!\cdot u^\varepsilon$ and the entropic fluctuation $\mathfrak{s}^\varepsilon=\frac{3}{2}\theta^\varepsilon-\rho^\varepsilon$. Taking curl eliminates the penalized gradient, yielding an equation of the form \eqref{wrtaeqn} in which the only remaining forcing enters through Burnett functionals. For $(\omega^\varepsilon,\mathfrak{s}^\varepsilon)$ we combine a maximum-principle estimate with potential theory to control $\|\nabla_x\mathbb{P}u^\varepsilon(t)\|_{L^\infty_x}$, and the optimal microscopic dissipation is used precisely to bound the Burnett forcing.

\newpage

 A central mechanism is the recovery of the full dissipation scale from the microscopic equation \eqref{Geqn0}: since $\AC{\P}F^\varepsilon$ is controlled at the optimal $\varepsilon^4\kappa$ scale, the microscopic forcing in the vorticity/entropy equations remains compatible with the macro energy method and no cumulative $\kappa$-loss occurs.

provided the crucial bootstrap assumption and control the forcing term (and Gronwall growth factor) of the system of kinetic vorticity, a gradient of entropic fluctuation, and acoustic variables. For kinetic vorticity, a gradient of entropic fluctuation, we use maximum principle and acoustic variables the dispersion estimate (Strichartz).

In two dimensions, the global-in-time argument hinges on controlling the growth factor
$\|\nabla_x(\rho^\varepsilon,u^\varepsilon,\theta^\varepsilon)\|_{L^\infty_x}$ in the Gronwall-type inequality
governing the total energy--dissipation (cf.\ \eqref{EDintro}). The penalized macroscopic system obstructs a
direct bound, so we exploit variables that are free from the leading penalization:
the vorticity $\omega^\varepsilon=\nabla_x^\perp\!\cdot u^\varepsilon$ and the entropic fluctuation
$\mathfrak s^\varepsilon=\frac32\theta^\varepsilon-\rho^\varepsilon$.

\smallskip
\noindent\emph{(i) Acoustic variables.}
The acoustic pair $(\mathbb P^\perp u^\varepsilon,\rho^\varepsilon+\theta^\varepsilon)$ satisfies a half-wave
system with forcing (see \eqref{locconP}). Strichartz estimates yield dispersive decay in Besov-type norms,
providing quantitative smallness of the acoustic component even when the initial acoustic modes do not vanish.

\newpage

The low--Mach/high--Reynolds scaling is intrinsically singular: large transport competes with fast relaxation, and the
macro--micro interaction must be treated at the correct scale. A simple linearization around a \emph{fixed}
Maxwellian $\mu$ fails to capture this structure: the macro--macro part of the collision term is not cancelled
and one encounters a $\kappa^{-1/2}$-singular energy imbalance:  
\[
\text{(Energy)}+\text{(Dissipation)}
\;\lesssim\;
\kappa^{-\frac12}\,
\text{(macro--macro interaction)}\cdot \text{(Dissipation)}^{\frac12},
\]
Classical approaches
avoid this either by prescribing the macroscopic fields through a given fluid solution (Hilbert expansion) or by
imposing strict structural well-preparedness to suppress the singular interactions. In contrast, we assume neither; instead we
control the non-closed fluxes---stress and heat flux---directly
within our quasi-linear framework.

The low--Mach/high--Reynolds hydrodynamic scaling couples large transport with fast relaxation:
\[
\varepsilon\,\partial_t F^\varepsilon + v\cdot\nabla_x F^\varepsilon
=
\frac{1}{\kappa\varepsilon}\,\mathcal N(F^\varepsilon,F^\varepsilon),
\qquad \kappa=\kappa(\varepsilon)\to 0,
\]
and hence produces an intrinsic singular structure in the macro--micro interaction. In particular, if one
linearizes around a \emph{fixed} Maxwellian $\mu$ (or any Maxwellian not incorporating the evolving macroscopic
fields), writing $F^\varepsilon=\mu+\varepsilon f^\varepsilon$, then macro--macro interactions are not
damped at the correct scale. At the level of energy estimates this manifests itself through a schematic
imbalance of the form
\[
\text{(Energy)}+\text{(Dissipation)}
\;\lesssim\;
\kappa^{-\frac12}\,
\text{(macro--macro interaction)}\cdot \text{(Dissipation)}^{\frac12},
\]
which is incompatible with closing an $\varepsilon$--uniform theory when the macroscopic component is not
prescribed.

Classical Hilbert-type frameworks avoid this obstruction by expanding around a Maxwellian whose parameters are
taken from a \emph{given} fluid solution; the singular interactions are then cancelled order-by-order by the
expansion. Another standard route is to enforce well-preparedness so that the singular acoustic interactions are
suppressed from the initial time. In the present work we avoid both assumptions: the macroscopic fields are not
fixed \emph{a priori}, and the data are not assumed to be well-prepared. Thus the macroscopic variables
generated by the Boltzmann dynamics interact nonlinearly with the microscopic component through higher velocity
moments (Burnett-type functionals), and one must control these interactions directly.

\subsubsection*{B.1. Singular scale limits and the failure of naive linearization}

The low--Mach/high--Reynolds hydrodynamic scaling couples large transport with fast relaxation:
\[
\varepsilon\,\partial_t F^\varepsilon + v\cdot\nabla_x F^\varepsilon
=
\frac{1}{\kappa\varepsilon}\,\mathcal N(F^\varepsilon,F^\varepsilon),
\qquad \kappa=\kappa(\varepsilon)\to 0,
\]
and hence produces an intrinsic singular structure in the macro--micro interaction. In particular, if one
linearizes around a \emph{fixed} Maxwellian $\mu$ (or any Maxwellian not incorporating the evolving macroscopic
fields), writing $F^\varepsilon=\mu+\varepsilon f^\varepsilon$, then macro--macro interactions are not
damped at the correct scale. At the level of energy estimates this manifests itself through a schematic
imbalance of the form
\[
\text{(Energy)}+\text{(Dissipation)}
\;\lesssim\;
\kappa^{-\frac12}\,
\text{(macro--macro interaction)}\cdot \text{(Dissipation)}^{\frac12},
\]
which is incompatible with closing an $\varepsilon$--uniform theory when the macroscopic component is not
prescribed.

Classical Hilbert-type frameworks avoid this obstruction by expanding around a Maxwellian whose parameters are
taken from a \emph{given} fluid solution; the singular interactions are then cancelled order-by-order by the
expansion. Another standard route is to enforce well-preparedness so that the singular acoustic interactions are
suppressed from the initial time. In the present work we avoid both assumptions: the macroscopic fields are not
fixed \emph{a priori}, and the data are not assumed to be well-prepared. Thus the macroscopic variables
generated by the Boltzmann dynamics interact nonlinearly with the microscopic component through higher velocity
moments (Burnett-type functionals), and one must control these interactions directly.

\subsubsection*{B.2. Local Maxwellian micro--macro decomposition and quasi-linearization}

The starting point is the exact cancellation $\mathcal N(M^\varepsilon,M^\varepsilon)=0$ for the \emph{local}
Maxwellian $M^\varepsilon=M_{[\mathrm P^\varepsilon,\mathrm U^\varepsilon,\Theta^\varepsilon]}$. We introduce the
macroscopic projection $\mathbf P$ associated with $M^\varepsilon$ (see \eqref{Pdef}) and decompose
\[
F^\varepsilon=\mathbf P F^\varepsilon + \AC{\mathbf P}F^\varepsilon,
\qquad \mathbf P F^\varepsilon =: M^\varepsilon,
\qquad
\int_{\mathbb R^3}\AC{\mathbf P}F^\varepsilon(1,v,|v|^2/2)^T\,dv=0.
\]
Because the reference Maxwellian already carries the evolving macroscopic fields, the singular macro--macro
interaction in the collision term is cancelled at the \emph{basic} energy level, and the linearized collision
operator becomes \emph{quasi-linear}:
\[
\mathcal L_{M^\varepsilon}(\AC{\mathbf P}F^\varepsilon)
=
-2\,\mathcal N(M^\varepsilon,\AC{\mathbf P}F^\varepsilon).
\]
This quasi-linearization is the first essential step toward an $\varepsilon$--uniform closure.

However, a new difficulty appears at higher regularity: differentiation does not commute with the
quasi-linear structure. When spatial derivatives fall on $M^\varepsilon$, one encounters terms of the form
$\AC{\mathbf P}(\partial^\alpha M^\varepsilon)$ inside the weighted inner product, which reintroduce a mismatch
between dissipation and testing. In particular, while dissipation naturally acts on
$\partial^\alpha\AC{\mathbf P}F^\varepsilon$, the energy identity requires testing against
$\AC{\mathbf P}(\partial^\alpha F^\varepsilon)$, producing contributions involving $\partial^\alpha M^\varepsilon$
and hence higher velocity moments (Burnett-type quantities). This dissipation--testing mismatch is precisely the
mechanism by which singular scale-loss can reemerge beyond the lowest derivative levels (see
Proposition~\ref{P.F.Energy}).

\subsubsection*{B.3. Quasi-linear framework: separating macroscopic penalization and microscopic dissipation}

To overcome the derivative-level mismatch, we separate the analysis into two coupled parts:

\smallskip
\noindent\emph{(i) Macroscopic penalized system.}
Taking collision-invariant moments yields a symmetric hyperbolic system for
$(\rho^\varepsilon,u^\varepsilon,\theta^\varepsilon)$ with singular penalization terms of the form
$\varepsilon^{-1}\nabla_x\cdot u^\varepsilon$ and $\varepsilon^{-1}\nabla_x(\rho^\varepsilon+\theta^\varepsilon)$
(see \eqref{locconNew}). Testing this system against the macroscopic variables themselves produces a crucial
cancellation of the leading penalized contributions, which motivates the nonlinear parametrization
$(\mathrm P^\varepsilon,\mathrm U^\varepsilon,\Theta^\varepsilon)=(e^{\varepsilon\rho^\varepsilon},\varepsilon u^\varepsilon,
e^{\varepsilon\theta^\varepsilon})$ introduced in \eqref{expform}.

\smallskip
\noindent\emph{(ii) Microscopic dissipative equation.}
The microscopic equation for $\AC{\mathbf P}F^\varepsilon$ (see \eqref{Geqn0}) provides strong dissipation at the
optimal scale $\varepsilon^{4}\kappa$ without incurring any additional $\kappa$-loss, even at higher derivative
levels. The key point is that the microscopic energy estimate tests directly against
$\partial^\alpha\AC{\mathbf P}F^\varepsilon/M^\varepsilon$ and therefore avoids the problematic appearance of
$\partial^\alpha M^\varepsilon$ that plagues the macroscopic testing.

\smallskip
By combining these two estimates at each derivative level, we recover the full dissipation scale and eliminate
the singular macro--macro interactions that would otherwise be generated by the Burnett functionals. At the
top-order level, where both the macroscopic and microscopic systems suffer an unavoidable derivative loss due to
the splitting of transport between them, we close the argument by estimating the full solution $F^\varepsilon$
as a whole and allowing a controlled loss only at the very top derivative level (cf.\ the $\dot H^{\mathrm N+1}$
estimate in Section~\ref{sec:toporder}).

\subsubsection*{B.4. Global 2D strategy: dispersive acoustics and vorticity/entropy variables}

In two dimensions, the global-in-time argument hinges on controlling the growth factor
$\|\nabla_x(\rho^\varepsilon,u^\varepsilon,\theta^\varepsilon)\|_{L^\infty_x}$ in the Gronwall-type inequality
governing the total energy--dissipation (cf.\ \eqref{EDintro}). The penalized macroscopic system obstructs a
direct bound, so we exploit variables that are free from the leading penalization:
the vorticity $\omega^\varepsilon=\nabla_x^\perp\!\cdot u^\varepsilon$ and the entropic fluctuation
$\mathfrak s^\varepsilon=\frac32\theta^\varepsilon-\rho^\varepsilon$.

\smallskip
\noindent\emph{(i) Acoustic variables.}
The acoustic pair $(\mathbb P^\perp u^\varepsilon,\rho^\varepsilon+\theta^\varepsilon)$ satisfies a half-wave
system with forcing (see \eqref{locconP}). Strichartz estimates yield dispersive decay in Besov-type norms,
providing quantitative smallness of the acoustic component even when the initial acoustic modes do not vanish.

\smallskip
\noindent\emph{(ii) Vorticity and entropic fluctuation.}
Taking curl removes the gradient penalization, and the $\mathfrak s^\varepsilon$-equation eliminates the
divergence penalization (see \eqref{wrtaeqn}). This reduction yields transport-type equations with microscopic
forcing through Burnett functionals. The microscopic dissipation estimate (Theorem~\ref{theoremG}) ensures that
these forcing terms are compatible with the energy--dissipation scale, and hence the vorticity equation yields
a closed control of $\|\omega^\varepsilon(t)\|_{L^\infty_x}$ (via a maximum principle) together with the
potential-theory bound for $\|\nabla_x \mathbb P u^\varepsilon(t)\|_{L^\infty_x}$. Once
$\|\nabla_x u^\varepsilon(t)\|_{L^\infty_x}$ is controlled, the transport structure further controls
$\|\nabla_x\mathfrak s^\varepsilon(t)\|_{L^\infty_x}$.

\subsubsection*{B.5. Weighted moment control and high-velocity regimes}

Finally, the $M^\varepsilon$-weighted energy method generates commutators involving the streaming operator
acting on the Maxwellian, $\partial_t M^\varepsilon + v\cdot\nabla_x M^\varepsilon$, which carry high velocity
moments. To control these contributions we employ an $L^2$--$L^\infty$ strategy in the spirit of the
energy/large-velocity decomposition initiated by Caflisch and developed in the modern framework of Guo: the
low-velocity region is controlled by energy--dissipation, while the large-velocity region recovers favorable
scaling by distributing integrability between $\nabla_x(\rho^\varepsilon,u^\varepsilon,\theta^\varepsilon)$ and
$\AC{\mathbf P}F^\varepsilon$ and invoking weighted $L^\infty_{x,v}$ bounds.

\subsubsection*{B.6. Local-in-time validity in 3D}

In three dimensions, the same quasi-linear framework yields local-in-time construction and quantitative
convergence under smooth initial data. Here the energy controls
$\|\nabla_x(\rho^\varepsilon,u^\varepsilon,\theta^\varepsilon)\|_{L^\infty_x}$ for $\mathrm N>\frac d2+1$,
so the Gronwall mechanism closes on a short time interval, leading to the local 3D results stated in
Theorem~\ref{T.C.Hk}.

\newpage

.
\newpage

\subsubsection*{B.1. Singular scale limits}
 A rigorous construction of solutions to the Boltzmann equation, together with the justification of the hydrodynamic limit without any \emph{a priori} information on the corresponding fluid solution, remains largely unexplored and constitutes one of the main open problems related to Hilbert’s sixth problem.
In contrast to the classical framework based on the Hilbert expansion—where macroscopic fields are prescribed by a given fluid solution—the macroscopic variables arising directly from the Boltzmann dynamics are not fixed \emph{a priori}. Instead, they evolve through nonlinear interactions with the microscopic component via the Burnett functionals.

As a consequence, the coupled macro–micro system inherently carries singular scaling structures, making the analysis fundamentally different from perturbative approaches around predetermined fluid states. This intrinsic interaction prevents a straightforward closure of the system and constitutes a major analytical obstacle in deriving hydrodynamic limits directly from the Boltzmann equation itself.

It is readily observed that, if one linearizes the Boltzmann equation around a Maxwellian of the form $F^{\e}= M_{[*,*,*]} + \e g^{\e}$ which does not fully incorporate the evolving macroscopic fields $(\mathrm{P}^{\e}, \mathrm{U}^{\e}, \mathrm{\Theta}^{\e})$, then the resulting microscopic dissipation is no longer sufficient to control the macro–macro interactions arising in the nonlinear collision term.
At the level of energy estimates, one encounters a singular structure of the form 
\[
\text{   (Energy) + (Dissipation) $\lesssim  \frac{1}{\sqrt{\kappa}}   \text{에너지의 매크로 파트} \text{(Dissipation)}^{1/2}$ }\]

\begin{align}\label{globalexpand}
\frac{d}{dt}\|f^{\e}\|_{L^2}^2 + \frac{1}{\kappa \e^2}\|\sqrt{\nu}(\II-\P_{\mu})f^{\e}\|_{L^2}^2 &\les \frac{1}{\sqrt \kappa} \la \Gamma(\P_{\mu} f^{\e},\P_{\mu} f^{\e}),  \frac{1}{\sqrt \kappa \e} (\II-\P_{\mu})f^{\e}\ra.
\end{align} 

A classical way to avoid such a singular behavior is to employ the Hilbert expansion, where the macroscopic fields are prescribed by a limiting fluid solution. However, this approach necessarily relies on \emph{a priori} information on the corresponding fluid dynamics. Another commonly adopted strategy is to impose well-prepared initial data so that the singular interactions are suppressed from the outset.

In contrast, the present work avoids both assumptions. We neither assume prior knowledge of the limiting fluid solution nor restrict ourselves to well-prepared initial configurations. Instead, we investigate how solutions to the Boltzmann equation evolve from general initial data and dynamically generate the hydrodynamic behavior.

Furthermore, without any knowledge of the limiting fluid solution, the macroscopic moment equations contain singular penalized terms
$\frac{1}{\e}\nabla_x\cdot u^{\e}$ and
$\frac{1}{\e}\nabla_x(\rho^{\e}+\ta^{\e})$
(see \eqref{locconNew}).
If one imposes stronger regularity assumptions on the time derivatives—for instance,
$\|\partial_t(\rho^{\e},u^{\e},\ta^{\e})\|<\infty$—the analysis becomes substantially simpler, since the moment equations then directly provide control of the acoustic variables.
However, such bounds in fact correspond to assuming \emph{well-prepared initial data}.
Indeed, at the initial time one formally obtains
\begin{align*}
\|\partial_t(\rho^{\e},u^{\e},\ta^{\e}) |_{t=0}\|
&\sim  
\frac{1}{\e}
\Big\|(\partial_tF^{\e})|M^{\e}|^{-\frac{1}{2}}|_{t=0}\Big\| 
\sim
\frac{1}{\e^2\kappa}
\Big\|\frac{1}{\e}\mathcal{L}(\AC{\P}F^{\e})|M^{\e}|^{-\frac{1}{2}}|_{t=0} \Big\|,
\end{align*}
which remains finite only if the microscopic component converges faster than the Mach-number scale, namely $\frac{1}{\e}\AC{\P}F^{\e}_0 \to 0$.
In this sense, strong time-derivative bounds implicitly enforce well-preparedness.
Notably, the weaker condition
$\|\e \partial_t(\rho^{\e},u^{\e},\ta^{\e})|_{t=0}\|<\infty$, which is typically regarded as not well-prepared in the low-Mach-number limit of compressible fluid equations, becomes well-prepared in the hydrodynamic limit.
Therefore, the initial data are not well-prepared unless one either avoids time derivatives altogether or works with at least the scaled derivative $\e^2\kappa\p_t$.

\hide
It is worth emphasizing that our approach requires \emph{no time derivatives of the initial data}. 
This leads to three fundamental novelties:
(1) We do not rely on a Navier–Stokes expansion for the macroscopic equations, 
which would otherwise be necessary to control terms such as 
$\kappa \partial_t \mathcal{L}^{-1}\AC{\P}F^{\e}$.
(2) Control of the acoustic variables is achieved without assuming bounds on 
their initial time derivatives. Instead, we exploit the half-wave structure 
satisfied by $\mathbb{P}^{\perp}u^{\e}$ and $\rho^{\e}+\ta^{\e}$.
(3) In estimating $F^{\e}$, particularly in the interaction of three 
microscopic components $\AC{\P}F^{\e}$, we avoid time-variable embeddings 
by simultaneously exploiting both the energy and the dissipation 
of the microscopic part.
\unhide

\subsubsection*{B.2. Micro-Macro decomposition}

To overcome the difficulty arising from \eqref{globalexpand}, we exploit the fundamental cancellation property of the collision operator, $\mathcal{N}(M^{\e},M^{\e}) = 0$. Motivated by the macro–micro decomposition framework developed in \cite{GuoInvent} and \cite{LYY}, we introduce a macroscopic projection associated with the local Maxwellian $M^{\e}$, (\eqref{Pdef}) where $\{e_i^{\e}\}_{i=0}^4$ forms an orthonormal basis of the space spanned by $\{1,v,|v|^2\}$ in velocity space with respect to the weighted inner product $\langle\cdot,\frac{*}{M^{\e}}\rangle$.
The basis functions correspond to the five collision invariants associated with the local Maxwellian $M^{\e}$.

We then decompose the solution as $F^{\e} = M^{\e}+ \AC{\P}F^{\e}$ where $M^{\e}=M_{[\mathrm{P}^{\e},\mathrm{U}^{\e},\mathrm{\Theta}^{\e}]}$ 
so that the macroscopic component is incorporated directly into the reference Maxwellian.
This decomposition, previously introduced in \cite{LYY}, cancels the singular macro–macro interactions in the compressible Euler scaling. In the present setting, the same structure yields an exact cancellation of the singular factor $\kappa^{-1/2}$ at the level of the basic energy estimate, effectively transforming the linearized collision operator into a quasi-linear one.

However, this cancellation mechanism does not persist under differentiation. Once spatial derivatives act on the quasi-linear operator $\mathcal{L}(\AC{\P}F^{\e})=-2\mathcal{N}(M^{\e},\AC{\P}F^{\e})$, a structural mismatch arises between dissipation and testing.
While the dissipation naturally acts on $\partial^{\alpha}\AC{\P}F^{\e}$, the corresponding energy estimate requires testing against $\AC{\P}(\partial^{\alpha}F^{\e})$. As a consequence, derivatives falling on the Maxwellian generate contributions involving $\AC{\P}(\partial^{\alpha}M^{\e})$ in the inner product.
The leading-order contribution of $\AC{\P}(\partial^{\alpha}M^{\e})$ is of order $O(\e^2)$, and it contains fourth-order velocity moments corresponding to the Burnett functionals.
\begin{align*}
\frac{1}{\e^4 \kappa }
\bigg\langle
\mathcal L(\partial^\alpha \AC{\P}F^\e),
\frac{\AC{\P}\partial^\alpha M^{\e}}{M^{\e}}
\bigg\rangle
\sim
\frac{1}{\sqrt{\kappa}}
(\text{macro--macro interactions})
\text{(Dissipation)}^{1/2}
\end{align*}  
This dissipation–testing mismatch reproduces the singular $\kappa^{-1/2}$ scaling. Consequently, already from the second-order derivative level, singular macro–macro interactions reemerge, constituting one of the main analytical difficulties of the present problem (see Proposition~\ref{P.F.Energy}).

If one continues the estimates so that one loses an additional power of $\kappa^{\frac{1}{2}}$ per derivative beginning at the second-order level, the vorticity equation derived from the macroscopic Boltzmann system becomes ineffective. 
Consequently, even local-in-time convergence appears out of reach, 
as the dissipation degenerates:
$\frac{d}{dt}\|\w^{\e}\|_{H^k} \sim \frac{1}{\e^2}\|\nabla_x^2\mathbf{r}^{\e}\|_{H^k} 
\les \kappa^{-k}\mathcal{D}^{\frac{1}{2}}$. 
One may use the macroscopic equations in Navier--Stokes form by expanding the microscopic part. However, in this case, the initial data must be prepared so as to control the time derivative of the microscopic part $\kappa \partial_t \mathcal{L}^{-1}\AC{\P}F^{\e}$.

\subsubsection*{B.3. Quasi-linear Framework}
To overcome this difficulties, we separate the macroscopic and microscopic components in the equation. The macroscopic system derived from the Boltzmann equation then becomes a symmetric hyperbolic system with penalization terms. (\eqref{locconNew})
Here we used $\mathrm{P}^{\e} = e^{\e \rho^{\e}}$, $\mathrm{U}^{\e}= \e u^{\e}$, and $\mathrm{\Theta}^{\e} = e^{\e \ta^{\e}}$. 
The quantities $\mathbf{r}^{\e}$ and $\mathfrak{q}^{\e}$ denote the Burnett functionals defined in \eqref{albe-def}.
Meanwhile, the microscopic equation provides strong dissipation without any loss of $\kappa$, even at higher derivative levels.
(\eqref{Geqn0})
where $\mathfrak{R}^{\e}_{ij}$ and $\mathcal{Q}^{\e}_i$ are defined in \eqref{ABdef}.

We stress that, although both estimates for $\p^{\al}F^{\e}$ and $\p^{\al}\AC{\P}F^{\e}$ contain commutators associated with the quasi-linear operator $\mathcal{L}\AC{\P}F^{\e} = -2\mathcal{N}(M^{\e},\AC{\P}F^{\e})$, the underlying inner-product structures differ essentially. Testing the macroscopic equation against $\p^{\al}F^{\e}$ introduces higher-order derivatives of the Maxwellian, namely $\p^{\al}M^{\e}$, which are responsible for the singular $\kappa$-loss. In contrast, the microscopic estimate obtained from \eqref{Geqn0} avoids this mechanism and therefore does not produce any singular dependence on $\kappa$.
As a consequence, by exploiting the microscopic equation \eqref{Geqn0} for $\AC{\P}F^{\e}$, one can eliminate the singular $\kappa$-scaling without sacrificing dissipation, thereby recovering the maximal scale $\e^{4}\kappa$.

추가: Liu-Yu-Yang \cite{LYY}, Tong Yang 등등 과 뭐가 다른가? 우리는 플루이드를 풀고 들어오지 않았다. 따라서 L 도 넌리니어 오퍼레이터임. 또한 나비어스톡스를 쓰지 않음. 우리의 어프로치는 얀/리우의 뷰를 따르나 로컬막스웰리안과 컨써베이션로우를 사용하고 켄설레이션을 사용했다는 점에서 데빌레-빌라니의 스피릿과도 매우 가깝다다.    )

Although the dissipation estimate is recovered at the optimal scaling, the decomposition $F^{\e} = M^{\e}+\AC{\P}F^{\e}$ renders the macroscopic equations more singular by generating additional penalized terms $\frac{1}{\e}\nabla_x\cdot u^{\e}$ and $\frac{1}{\e}\nabla_x(\rho^{\e}+\ta^{\e})$.
Moreover, since the transport structure associated with $\frac{v}{\e}\cdot\nabla_x F^{\e}$ is split between the macroscopic and microscopic dynamics, both systems \eqref{locconNew} and \eqref{Geqn0} inevitably suffer from derivative loss. As a consequence, at the top-order derivative level there is no available mechanism to control the resulting macro–macro interactions.

To handle the singular penalized terms $\frac{1}{\e}\nabla_x\cdot u^{\e}$ and $\frac{1}{\e}\nabla_x(\rho^{\e}+\ta^{\e})$, we exploit the symmetric hyperbolic structure of the macroscopic system by testing each equation in \eqref{locconNew} against the variables $(\rho^{\e},u^{\e},\ta^{\e})$, which produces a crucial cancellation.
This mechanism motivates the choice of the nonlinear formulation $(\mathrm{P}^{\e},\mathrm{U}^{\e},\mathrm{\Theta}^{\e}) = (e^{\e \rho^{\e}},\e u^{\e},e^{\e \ta^{\e}})$, under which the cancellation becomes transparent at the energy level.

*To address the derivative loss present in both systems \eqref{locconNew} and \eqref{Geqn0}, we found that the interactions mediated by the Burnett functionals cancel at leading order once the macroscopic and microscopic estimates are combined. 
This reveals a notable structural mechanism where the Burnett functionals link the macroscopic and microscopic dynamics at the level of individual equations, while their contributions disappear in the coupled energy estimate, effectively eliminating the Burnett effect from the final energy balance.

*Since this cancellation occurs only at the principal level, the remaining derivative-loss terms must still be controlled. After establishing $H_x^{\mathrm N}$ estimates for both \eqref{locconNew} and \eqref{Geqn0}, we allow a loss of $\kappa^{\frac12}$ only at the top-order level $\dot{H}_x^{\mathrm N+1}$ by estimating the full solution $F^{\e}$ as a whole.
This behavior is reminiscent of the Navier–Stokes framework, in which $\|\w^{\e}\|_{H^{\mathrm{N}}_x}$ is controlled, while the corresponding dissipation is available only at the weaker $\kappa^{\frac12}$ scale through $\kappa^{\frac{1}{2}}\|\nabla_x\w^{\e}\|_{H^{\mathrm{N}}_x}$.

\subsubsection*{B.4. High-order moment control in local in time validity}  
Finally, to control the high-order moments generated by the streaming term $\pt_t M^{\e} + \frac{v}{\eps}\cdot\nabla_x M^{\e}$ arising from the commutator associated with the weight $1/M^{\e}$, we follow the approach initiated by Caflisch \cite{Caflisch} and later developed into an $L^2$–$L^\infty$ framework in \cite{GuoJJCPAM}.
Expanding the solution around a global Maxwellian and decomposing the phase space into small- and large-velocity regions, the low-velocity part is handled by energy–dissipation estimates, while favorable scaling is recovered in the large-velocity regime by distributing integrability between $\nabla_x u^{\e}$, $\nabla_x\theta^{\e}$ and $\AC{\P}F^{\e}$ and extracting the $L^\infty_{x,v}$ bound of the perturbation.

\subsubsection*{B.5. Local-in-time Validity in 3d}

In the energy estimates for $F^{\e}$, the macro–macro interaction terms initially produce a singular loss at the $\kappa$-scale.
By decomposing the solution into macroscopic and microscopic components and exploiting the precise structure of the microscopic equation, this singular behavior is removed through an exact cancellation that restores the full dissipation scale.
The derivative loss arising from this decomposition is subsequently compensated by additional cancellations between the macroscopic and microscopic systems.

Altogether, combining the $H^{\mathrm N}$ estimates for \eqref{locconNew} and \eqref{Geqn0} with the $\dot H^{\mathrm N+1}$ estimate for $F^{\e}$ allows us to characterize the growth of the total energy–dissipation structure, which is governed by $\|\nabla_x(\rho^{\e},u^{\e},\ta^{\e})\|_{L^\infty_x}$:
\begin{align}\label{EDintro}
\bega
\mathcal{E}_{tot}(t)+\int_0^t\mathcal{D}_{tot}(s)ds &\leq C\int_0^t\big(1+\|\nabla_x(\rho^{\e},u^{\e},\ta^{\e})(s)\|_{L^\infty_x}^2\big) \mathcal{E}_{tot}(s)ds + \e t.
\enda
\end{align}
Here $\mathcal{E}_{tot}^{\mathrm{N}}(F^\e(t))$ and $\mathcal{D}_{tot}^{\mathrm{N}}(F^\e(t))$ are defined in \eqref{EDtotdef}. 
In the 3-dimensional case, since the energy functional already controls $\|\nabla_x(\rho^{\e},u^{\e},\ta^{\e})\|_{L^\infty_x}\leq C\mathcal{E}_M^{\frac{1}{2}}$ for $\mathrm{N}>\frac{d}{2}+1$, this yields the construction of local-in-time solutions.

To the best of our knowledge, this provides the first result on the incompressible hydrodynamic limit in which solutions are constructed solely from the Boltzmann equation, without assuming any \emph{a priori} information on the limiting fluid variables. 
This approach is in the same spirit as the work of Nishida~\cite{Nishida}, where the hydrodynamic limit toward the compressible Euler equations was established.
We now aim to extend this result to a global-in-time hydrodynamic limit, particularly in the two-dimensional case, while overcoming the above regularity constraint. 
The main difficulties are summarized as follows.

\subsubsection*{B.6. New Difficulity: Penalization terms}

To extend the analysis to the global-in-time two-dimensional setting, obtaining a uniform bound for $\|\nabla_x(\rho^{\e},u^{\e},\ta^{\e})\|_{L^\infty_x}$ constitutes a critical threshold.
However, as seen from the macroscopic system \eqref{locconNew}, the evolution equations for $(\rho^{\e},u^{\e},\ta^{\e})$ contain singularly penalized terms, $\frac{1}{\e}\nabla_x\cdot u^{\e}$ and $\frac{1}{\e}\nabla_x(\rho^{\e}+\ta^{\e})$ which obstruct a direct control of this quantity.

To overcome this difficulty, we exploit macroscopic variables free from the singular penalization, namely the vorticity $\w^{\e}:= \nabla_x^\perp \cdot u^{\e}$ and the thermodynamic variable $\mathfrak{s}^{\e}=\frac{3}{2}\ta^{\e}-\rho^{\e}$.

To estimate $\|\nabla_xu^{\e}\|_{L^\infty_x}$, we decompose the velocity field into vorticity part $\w^{\e}$ and irrotational part $\nabla_x \mathbb{P}^{\perp}u^{\e}$.
For the estimate of $\|\nabla_x\rho^{\e}\|_{L^\infty_x}$ and $\|\nabla_x\ta^{\e}\|_{L^\infty_x}$, we decompose the estimate by thermodynamic variables $\|\nabla_xp^{\e}\|_{L^\infty_x}$ and $\|\nabla_x\mathfrak{s}^{\e}\|_{L^\infty_x}$.

\subsubsection*{B.7. Global Validity: Takeaway 1-- Quantitative control of the acoustic variables}

For the acoustic variables, the local conservation laws 
\eqref{locconNew} yield the following half-wave system with forcing terms: (\eqref{locconP})
Precise definitions of $\bPhi^{\e}_{\rho+\ta}$ and $\bPhi^{\e}_{\mathbb{P}^{\perp}u}$ are given in Lemma~\ref{L.locconP}.
By exploiting the associated Strichartz estimates, the acoustic variables converge to zero due to dispersive effects, even when the initial acoustic components are not vanishing. 
Moreover, this approach provides the explicit convergence rate 
$\e^{\frac{1}{r}}$
in some Besov space ($L^r_T \dot{B}^{s}_{p,1}$ where $\dot{B}_{p,q}^s$ is defined in \eqref{Besovdef}) for $d=2,3$.

\subsubsection*{B.8. Global Validity: Takeaway 2--Control of the vorticity and the entropic fluctuation $\mathfrak{s}^{\e} = \frac{3}{2} \theta^\e - \rho^\e$ 
}

The vorticity equation eliminates the singular penalized gradient term $\frac{1}{\e}\nabla_x(\rho^{\e}+\ta^{\e})$, while the evolution equation for the thermodynamic variable $\mathfrak{s}^{\e}$ removes the divergence penalization $\frac{1}{\e}\nabla_x\cdot u^{\e}$:
\begin{align}\label{wrtaeqn}
\bega
&\p_t\w^{\e} + u^{\e}\cdot \nabla_x \w^{\e} + (\nabla_x\cdot u^{\e})\w^{\e} +k_B \mathrm{\Theta}^{\e} (\nabla_x^{\perp}\ta^{\e})\cdot\nabla_x(\rho^{\e}+\ta^{\e}) +\frac{1}{\e^2}\nabla_x^{\perp}\cdot\bigg(\frac{1}{\mathrm{P}^{\e} }\nabla_x \cdot \mathbf{r}^{\e}\bigg) = 0, \cr
&\p_t \mathfrak{s}^{\e} + u^{\e}\cdot \nabla_x \mathfrak{s}^{\e} +\frac{1}{\e^2}\frac{1}{k_B\mathrm{P}^{\e}\mathrm{\Theta}^{\e} } \big(\nabla_x\cdot \mathfrak{q}^{\e} + \nabla_x \mathrm{U}^{\e} : \mathbf{r}^{\e} \big)=0.
\enda
\end{align}
This reduction allows the macroscopic dynamics to be analyzed through transport-type equations without the singular penalization.

For the vorticity, we employ the maximum principle to control $\|\w^{\e}(t)\|_{L^\infty_x}$ while potential-theory estimates provide bounds for $\|\w^{\e}(t)\|_{H^2_x}$ and $\|\nabla_x \mathbb{P}u^{\e}(t)\|_{L_x^\infty}$.

A central point of the present work concerns the scaling of the microscopic forcing appearing in the vorticity equation.
If, as discussed in Section~3.2, the dissipation were to lose a factor of $\kappa^{1/2}$ at each derivative level beyond $H^2$, even the basic $L^2$ estimate for the vorticity would fail. 
In that case, one would obtain $\frac{d}{dt}\|\w^{\e}\|_{L^2_x} \sim \mathcal{D}^{\frac{1}{2}}$ and inserting the resulting bound for $\|\nabla_x\mathbb{P}u^{\e}\|_{L^\infty_x}$ into the growth factor in \eqref{EDintro} would destroy the global-in-time control of the solution.
At higher regularity levels, the dissipation would become increasingly singular, preventing both $H^2$ and $L^\infty$ estimates for the vorticity since $\frac{d}{dt}\|\w^{\e}\|_{H^k} \sim \frac{1}{\e^2}\|\nabla_x^2\mathbf{r}^{\e}\|_{H^k} 
\les \kappa^{-k}\mathcal{D}^{\frac{1}{2}}$.

The key mechanism preventing this loss is the recovery of the full dissipation scale from the microscopic equation \eqref{Geqn0}.
Since the microscopic component $\AC{\P}F^{\e}$ is controlled at the optimal $\e^4\kappa$ scale by dissipation, the microscopic forcing appearing in the vorticity equation remains compatible with the energy–dissipation framework.
In contrast, one may derive a viscous term such as $-\kappa\Delta_x\w^{\e}$ through a Navier–Stokes expansion. However, this approach requires control of time derivatives of microscopic quantities, namely
$\kappa \nabla_x \mathcal{L}^{-1} (\partial_t \AC{\P}F^{\e})$,
and therefore implicitly assumes well-prepared initial data.
Our approach avoids this requirement and relies solely on the Euler-type moment equations, thereby allowing global-in-time validity without preparation assumptions.

A further difficulty arises from the fact that, unlike in the incompressible Euler equation, the velocity field transporting the vorticity remains compressible.
Although dispersive estimates yield favorable scaling for the acoustic variables in Besov spaces, such control is available only below the top-order derivative level due to forcing terms involving $\nabla_x\cdot \mathbf{r}^{\e}$ and $\nabla_x\cdot \mathfrak{q}^{\e}
$ in \eqref{locconP}.
Consequently, the estimates must be arranged so that the compressible component carries as few derivatives as possible.

Combining the refined dissipation scale with the dispersive gain obtained from the acoustic dynamics, the vorticity equation yields control of the growth of
$\|\w^{\e}(t)\|_{L^\infty_x}$ and
$\|\nabla_x\mathbb{P}u^{\e}(t)\|_{L^\infty_x}$
in terms only of the initial data, the time interval, and energy–dissipation quantities endowed with favorable scaling.
Once $\|\nabla_x u^{\e}(t)\|_{L^\infty_x}$ is controlled, the maximum principle further determines the growth of $\|\nabla_x\mathfrak{s}^{\e}(t)\|_{L^\infty_x}$.

 \newtheorem*{theoremG}{Theorem G}

\begin{theoremG}[Informal statement of Theorem \ref{T.2D.global}]\label{theoremG}
Let $\Omega=\mathbb{R}^2$. Assume that the initial data $\{F^\varepsilon_0\}$ satisfy a \emph{uniform-in-$\varepsilon$ scaled modulated entropy bound}, while they are allowed to blow up in strong topologies in the sense of the admissible blow-up condition (Definition~\ref{ABC}). Then
\begin{align}\label{eq:informal_micro_dissipation}
 \sum_{|\alpha| \leq N+1}  \frac{\kappa^{(|\alpha|-\mathrm{N})_+}}{\eps^4\kappa}
 \int_0^{T_*}\!\!\int_{\Omega \times \R^3} \nu \,
 \frac{\big|\p^{\alpha}\AC{\mathbf P} F^\varepsilon (t)\big|^2}{M^{\e}}
 \,dvdxdt
 \;\lesssim\;
\exp\!\Big(\exp\!\big(\exp(\|\omega_0^\varepsilon\|_\infty\, T_*)\big)\Big)\,\mathcal{E}_{\mathrm{tot}}(0)
 + \textit{l.o.t.}
\end{align}
\end{theoremG}

\hide

\subsection*{B. New Quasi-Linear Method: Difficulties and Key Ideas}
\subsubsection*{B.1. Singular scale limits}
A rigorous construction of solutions to the Boltzmann equation, together with the justification of the hydrodynamic limit without any \emph{a priori} information on the corresponding fluid solution, remains largely unexplored and constitutes one of the main open problems related to Hilbert’s sixth problem.
In contrast to the classical framework based on the Hilbert expansion—where macroscopic fields are prescribed by a given fluid solution—the macroscopic variables arising directly from the Boltzmann dynamics are not fixed \emph{a priori}. Instead, they evolve through nonlinear interactions with the microscopic component via the Burnett functionals.

As a consequence, the coupled macro–micro system inherently carries singular scaling structures, making the analysis fundamentally different from perturbative approaches around predetermined fluid states. This intrinsic interaction prevents a straightforward closure of the system and constitutes a major analytical obstacle in deriving hydrodynamic limits directly from the Boltzmann equation itself.

It is readily observed that, if one linearizes the Boltzmann equation around a Maxwellian of the form $F^{\e}= M_{[*,*,*]} + \e g^{\e}$ which does not fully incorporate the evolving macroscopic fields $(\mathrm{P}^{\e}, \mathrm{U}^{\e}, \mathrm{\Theta}^{\e})$, then the resulting microscopic dissipation is no longer sufficient to control the macro–macro interactions arising in the nonlinear collision term.
At the level of energy estimates, one encounters a singular structure of the form
\[
\text{   (Energy) + (Dissipation) $\lesssim  \frac{1}{\sqrt{\kappa}}   \text{에너지의 매크로 파트} \text{(Dissipation)}^{1/2}$ }\]

\begin{align}\label{globalexpand}
\frac{d}{dt}\|f^{\e}\|_{L^2}^2 + \frac{1}{\kappa \e^2}\|\sqrt{\nu}(\II-\P_{\mu})f^{\e}\|_{L^2}^2 &\les \frac{1}{\sqrt \kappa} \la \Gamma(\P_{\mu} f^{\e},\P_{\mu} f^{\e}),  \frac{1}{\sqrt \kappa \e} (\II-\P_{\mu})f^{\e}\ra.
\end{align}

A classical way to avoid such a singular behavior is to employ the Hilbert expansion, where the macroscopic fields are prescribed by a limiting fluid solution. However, this approach necessarily relies on \emph{a priori} information on the corresponding fluid dynamics. Another commonly adopted strategy is to impose well-prepared initial data so that the singular interactions are suppressed from the outset.

In contrast, the present work avoids both assumptions. We neither assume prior knowledge of the limiting fluid solution nor restrict ourselves to well-prepared initial configurations. Instead, we investigate how solutions to the Boltzmann equation evolve from general initial data and dynamically generate the hydrodynamic behavior.

Furthermore, without any knowledge of the limiting fluid solution, the macroscopic moment equations contain singular penalized terms
$\frac{1}{\e}\nabla_x\cdot u^{\e}$ and
$\frac{1}{\e}\nabla_x(\rho^{\e}+\ta^{\e})$
(see \eqref{locconNew}).
If one imposes stronger regularity assumptions on the time derivatives—for instance,
$\|\partial_t(\rho^{\e},u^{\e},\ta^{\e})\|<\infty$—the analysis becomes substantially simpler, since the moment equations then directly provide control of the acoustic variables.
However, such bounds in fact correspond to assuming \emph{well-prepared initial data}.
Indeed, at the initial time one formally obtains
\begin{align*}
\|\partial_t(\rho^{\e},u^{\e},\ta^{\e}) |_{t=0}\|
&\sim  
\frac{1}{\e}
\Big\|(\partial_tF^{\e})|M^{\e}|^{-\frac{1}{2}}|_{t=0}\Big\| 
\sim
\frac{1}{\e^2\kappa}
\Big\|\frac{1}{\e}\mathcal{L}(\AC{\P}F^{\e})|M^{\e}|^{-\frac{1}{2}}|_{t=0} \Big\|,
\end{align*}
which remains finite only if the microscopic component converges faster than the Mach-number scale, namely $\frac{1}{\e}\AC{\P}F^{\e}_0 \to 0$.
In this sense, strong time-derivative bounds implicitly enforce well-preparedness.
Notably, the weaker condition
$\|\e \partial_t(\rho^{\e},u^{\e},\ta^{\e})|_{t=0}\|<\infty$, which is typically regarded as not well-prepared in the low-Mach-number limit of compressible fluid equations, becomes well-prepared in the hydrodynamic limit.
Therefore, the initial data are not well-prepared unless one either avoids time derivatives altogether or works with at least the scaled derivative $\e^2\kappa\p_t$.

\hide
It is worth emphasizing that our approach requires \emph{no time derivatives of the initial data}. 
This leads to three fundamental novelties:
(1) We do not rely on a Navier–Stokes expansion for the macroscopic equations, 
which would otherwise be necessary to control terms such as 
$\kappa \partial_t \mathcal{L}^{-1}\AC{\P}F^{\e}$.
(2) Control of the acoustic variables is achieved without assuming bounds on 
their initial time derivatives. Instead, we exploit the half-wave structure 
satisfied by $\mathbb{P}^{\perp}u^{\e}$ and $\rho^{\e}+\ta^{\e}$.
(3) In estimating $F^{\e}$, particularly in the interaction of three 
microscopic components $\AC{\P}F^{\e}$, we avoid time-variable embeddings 
by simultaneously exploiting both the energy and the dissipation 
of the microscopic part.
\unhide



\subsubsection*{B.2. Micro-Macro decomposition}


To overcome the difficulty arising from \eqref{globalexpand}, we exploit the fundamental cancellation property of the collision operator, $\mathcal{N}(M^{\e},M^{\e}) = 0$. Motivated by the macro–micro decomposition framework developed in \cite{GuoInvent} and \cite{LYY}, we introduce a macroscopic projection associated with the local Maxwellian $M^{\e}$, (\eqref{Pdef})
where $\{e_i^{\e}\}_{i=0}^4$ forms an orthonormal basis of the space spanned by $\{1,v,|v|^2\}$ in velocity space with respect to the weighted inner product $\langle\cdot,\frac{*}{M^{\e}}\rangle$.
The basis functions correspond to the five collision invariants associated with the local Maxwellian $M^{\e}$.

We then decompose the solution as $F^{\e} = M^{\e}+ \AC{\P}F^{\e}$ where $M^{\e}=M_{[\mathrm{P}^{\e},\mathrm{U}^{\e},\mathrm{\Theta}^{\e}]}$ 
so that the macroscopic component is incorporated directly into the reference Maxwellian.
This decomposition, previously introduced in \cite{LYY}, cancels the singular macro–macro interactions in the compressible Euler scaling. In the present setting, the same structure yields an exact cancellation of the singular factor $\kappa^{-1/2}$ at the level of the basic energy estimate, effectively transforming the linearized collision operator into a quasi-linear one.


However, this cancellation mechanism does not persist under differentiation. Once spatial derivatives act on the quasi-linear operator $\mathcal{L}(\AC{\P}F^{\e})=-2\mathcal{N}(M^{\e},\AC{\P}F^{\e})$, a structural mismatch arises between dissipation and testing.
While the dissipation naturally acts on $\partial^{\alpha}\AC{\P}F^{\e}$, the corresponding energy estimate requires testing against $\AC{\P}(\partial^{\alpha}F^{\e})$. As a consequence, derivatives falling on the Maxwellian generate contributions involving $\AC{\P}(\partial^{\alpha}M^{\e})$ in the inner product.
The leading-order contribution of $\AC{\P}(\partial^{\alpha}M^{\e})$ is of order $O(\e^2)$, and it contains fourth-order velocity moments corresponding to the Burnett functionals.
\begin{align*}
\frac{1}{\e^4 \kappa }
\bigg\langle
\mathcal L(\partial^\alpha \AC{\P}F^\e),
\frac{\AC{\P}\partial^\alpha M^{\e}}{M^{\e}}
\bigg\rangle
\sim
\frac{1}{\sqrt{\kappa}}
(\text{macro--macro interactions})
\text{(Dissipation)}^{1/2}
\end{align*}
\hide
\begin{align*}
\frac{\AC{\P} (\p^{\al} M^{\e}) }{M^{\e}} \sim \sum_{0<\beta<\al} \bigg( \p^{\alpha- \beta}\mathrm{U}^\e \otimes \p^\beta \mathrm{U}^\e  :\mathfrak R ^\e +  \p^{\alpha - \beta } \mathrm{U}^\e \p^\beta \mathrm{\Theta}^\e  \cdot \mathcal Q ^\e + \p^{\al-\beta}\mathrm{\Theta}^\e \p^{\beta}\mathrm{\Theta}^\e(1+|v|^4) \bigg)
\end{align*}
\begin{align*}
\frac{1}{\e^4 \kappa } \left\langle \mathcal L ( \p^\alpha \AC{\P} F^\e) , \frac{\AC{\P} \p^{\al}M^{\e}}{M^{\e}} \right\rangle &\sim  \frac{1}{\sqrt{\kappa}}  \sum_{0<\beta<\al} \frac{1}{\e^2}\bigg( |\p^{\alpha-\beta}\mathrm{U}^{\e}\p^\beta \mathrm{U}^{\e}|+  |\p^{\alpha - \beta } \mathrm{U}^\e \p^\beta \mathrm{\Theta}^\e| + \p^{\al-\beta}\mathrm{\Theta}^\e \p^{\beta}\mathrm{\Theta}^\e| \bigg) \cr 
&\quad \times \bigg(\frac{1}{\e^4\kappa} \int_{\Omega\times\R^3} \frac{|\p^{\al}\AC{\P}F^{\e}|^2}{M^{\e}} dvdx \bigg)^{\frac{1}{2}}.
\end{align*}
\unhide
\hide
\[
\frac{1}{\e^4 \kappa }\langle \mathcal L (  \p^\alpha \AC{\P} F^\e) , \AC{\P} \p^{\al}M^{\e}\rangle  \sim \sum_{0 < \beta}  \p^{\alpha- \beta}U^\e \otimes \p^\beta  U^\e:  
\frac{1}{\e^4 \kappa }\langle \mathcal L (  \p^\alpha \AC{\P} F^\e) , \mathfrak R^\e M^\e\rangle 
+  \p^{\alpha - \beta }  \Theta^\e \p^\beta U^\e \cdot \frac{1}{\e^4 \kappa }\langle \mathcal L (  \p^\alpha \AC{\P} F^\e) , \mathcal Q^\e M^\e\rangle 
\]
\unhide
This dissipation–testing mismatch reproduces the singular $\kappa^{-1/2}$ scaling. Consequently, already from the second-order derivative level, singular macro–macro interactions reemerge, constituting one of the main analytical difficulties of the present problem (see Proposition~\ref{P.F.Energy}).

\hide
\begin{align}\label{FE00}
\bega
\frac{d}{dt}\frac{1}{\eps^2}\int \frac{|\p^{\alpha}F^{\e}|^2}{M^{\e}} &+ \frac{1}{\eps^4\kappa}\int \nu \frac{|\pt^\al \AC{\P}F^{\e}|^2}{M^{\e}} \leq \frac{1}{\eps^2}\bigg|\int \frac{|\p^{\al}F^{\e}|^2}{M^{\e}}\left(\pt_t M^{\e} + \frac{v}{\eps}\cdot\nabla_x M^{\e}\right)|M^{\e}|^{-2} \bigg| \cr
& +\frac{1}{\sqrt{\kappa}}\sum_{0< \beta< \alpha}\int \mathcal{N}\Big(\frac{1}{\e}\mathbf{P}(\p^{\beta}M^{\e}),\frac{1}{\e}\mathbf{P}(\p^{\alpha-\beta}M^{\e})\Big)\frac{1}{\eps^2\sqrt{\kappa}}\p^{\alpha}\AC{\P}F^{\e}|M^{\e}|^{-1}  + \cdots
\enda
\end{align}
 \unhide

If one continues the estimates so that one loses an additional power of $\kappa^{\frac{1}{2}}$ per derivative beginning at the second-order level, the vorticity equation derived from the macroscopic Boltzmann system becomes ineffective. 
Consequently, even local-in-time convergence appears out of reach, 
as the dissipation degenerates:
$\frac{d}{dt}\|\w^{\e}\|_{H^k} \sim \frac{1}{\e^2}\|\nabla_x^2\mathbf{r}^{\e}\|_{H^k} 
\les \kappa^{-k}\mathcal{D}^{\frac{1}{2}}$. 
One may use the macroscopic equations in Navier--Stokes form by expanding the microscopic part. However, in this case, the initial data must be prepared so as to control the time derivative of the microscopic part $\kappa \partial_t \mathcal{L}^{-1}\AC{\P}F^{\e}$.

\subsubsection*{B.3. Quasi-linear Framework}
To overcome this difficulties, we separate the macroscopic and microscopic components in the equation. The macroscopic system derived from the Boltzmann equation then becomes a symmetric hyperbolic system with penalization terms. (\eqref{locconNew})
Here we used $\mathrm{P}^{\e} = e^{\e \rho^{\e}}$, $\mathrm{U}^{\e}= \e u^{\e}$, and $\mathrm{\Theta}^{\e} = e^{\e \ta^{\e}}$. 
The quantities $\mathbf{r}_{ij}^{\e}$ and $\mathfrak{q}_j^{\e}$ denote the Burnett functionals defined in \eqref{albe-def}.
Meanwhile, the microscopic equation provides strong dissipation without any loss of $\kappa$, even at higher derivative levels.
(\eqref{Geqn0})
where $\mathfrak{R}^{\e}_{ij}$ and $\mathcal{Q}^{\e}_i$ are defined in \eqref{ABdef}.

We stress that, although both estimates for $\p^{\al}F^{\e}$ and $\p^{\al}\AC{\P}F^{\e}$ contain commutators associated with the quasi-linear operator $\mathcal{L}\AC{\P}F^{\e} = -2\mathcal{N}(M^{\e},\AC{\P}F^{\e})$, the underlying inner-product structures differ essentially. Testing the macroscopic equation against $\p^{\al}F^{\e}$ introduces higher-order derivatives of the Maxwellian, namely $\p^{\al}M^{\e}$, which are responsible for the singular $\kappa$-loss. In contrast, the microscopic estimate obtained from \eqref{Geqn0} avoids this mechanism and therefore does not produce any singular dependence on $\kappa$.
As a consequence, by exploiting the microscopic equation \eqref{Geqn0} for $\AC{\P}F^{\e}$, one can eliminate the singular $\kappa$-scaling without sacrificing dissipation, thereby recovering the maximal scale $\e^{4}\kappa$.


추가: Liu-Yu-Yang \cite{LYY}, Tong Yang 등등 과 뭐가 다른가? 우리는 플루이드를 풀고 들어오지 않았다. 따라서 L 도 넌리니어 오퍼레이터임. 또한 나비어스톡스를 쓰지 않음. 우리의 어프로치는 얀/리우의 뷰를 따르나 로컬막스웰리안과 컨써베이션로우를 사용하고 켄설레이션을 사용했다는 점에서 데빌레-빌라니의 스피릿과도 매우 가깝다다.    )



Although the dissipation estimate is recovered at the optimal scaling, the decomposition $F^{\e} = M^{\e}+\AC{\P}F^{\e}$ renders the macroscopic equations more singular by generating additional penalized terms $\frac{1}{\e}\nabla_x\cdot u^{\e}$ and $\frac{1}{\e}\nabla_x(\rho^{\e}+\ta^{\e})$.
Moreover, since the transport structure associated with $\frac{v}{\e}\cdot\nabla_x F^{\e}$ is split between the macroscopic and microscopic dynamics, both systems \eqref{locconNew} and \eqref{Geqn0} inevitably suffer from derivative loss. As a consequence, at the top-order derivative level there is no available mechanism to control the resulting macro–macro interactions.


To handle the singular penalized terms $\frac{1}{\e}\nabla_x\cdot u^{\e}$ and $\frac{1}{\e}\nabla_x(\rho^{\e}+\ta^{\e})$, we exploit the symmetric hyperbolic structure of the macroscopic system by testing each equation in \eqref{locconNew} against the variables $(\rho^{\e},u^{\e},\ta^{\e})$, which produces a crucial cancellation.
This mechanism motivates the choice of the nonlinear formulation $(\mathrm{P}^{\e},\mathrm{U}^{\e},\mathrm{\Theta}^{\e}) = (e^{\e \rho^{\e}},\e u^{\e},e^{\e \ta^{\e}})$, under which the cancellation becomes transparent at the energy level.

To address the derivative loss present in both systems \eqref{locconNew} and \eqref{Geqn0}, we found that the interactions mediated by the Burnett functionals cancel at leading order once the macroscopic and microscopic estimates are combined. 
This reveals a notable structural mechanism where the Burnett functionals link the macroscopic and microscopic dynamics at the level of individual equations, while their contributions disappear in the coupled energy estimate, effectively eliminating the Burnett effect from the final energy balance.

Since this cancellation occurs only at the principal level, the remaining derivative-loss terms must still be controlled. After establishing $H_x^{\mathrm N}$ estimates for both \eqref{locconNew} and \eqref{Geqn0}, we allow a loss of $\kappa^{\frac12}$ only at the top-order level $\dot{H}_x^{\mathrm N+1}$ by estimating the full solution $F^{\e}$ as a whole.
This behavior is reminiscent of the Navier–Stokes framework, in which $\|\w^{\e}\|_{H^{\mathrm{N}}_x}$ is controlled, while the corresponding dissipation is available only at the weaker $\kappa^{\frac12}$ scale through $\kappa^{\frac{1}{2}}\|\nabla_x\w^{\e}\|_{H^{\mathrm{N}}_x}$.

\subsubsection*{B.4. High-order moment control in local in time validity} 
Finally, to control the high-order moments generated by the streaming term $\pt_t M^{\e} + \frac{v}{\eps}\cdot\nabla_x M^{\e}$ arising from the commutator associated with the weight $1/M^{\e}$, we follow the approach initiated by Caflisch \cite{Caflisch} and later developed into an $L^2$–$L^\infty$ framework in \cite{GuoJJCPAM}.
Expanding the solution around a global Maxwellian and decomposing the phase space into small- and large-velocity regions, the low-velocity part is handled by energy–dissipation estimates, while favorable scaling is recovered in the large-velocity regime by distributing integrability between $\nabla_x u^{\e}$, $\nabla_x\theta^{\e}$ and $\AC{\P}F^{\e}$ and extracting the $L^\infty_{x,v}$ bound of the perturbation.

\subsubsection*{B.5. Local-in-time Validity in 3d}



In the energy estimates for $F^{\e}$, the macro–macro interaction terms initially produce a singular loss at the $\kappa$-scale.
By decomposing the solution into macroscopic and microscopic components and exploiting the precise structure of the microscopic equation, this singular behavior is removed through an exact cancellation that restores the full dissipation scale.
The derivative loss arising from this decomposition is subsequently compensated by additional cancellations between the macroscopic and microscopic systems.

Altogether, combining the $H^{\mathrm N}$ estimates for \eqref{locconNew} and \eqref{Geqn0} with the $\dot H^{\mathrm N+1}$ estimate for $F^{\e}$ allows us to characterize the growth of the total energy–dissipation structure, which is governed by $\|\nabla_x(\rho^{\e},u^{\e},\ta^{\e})\|_{L^\infty_x}$:
\begin{align}\label{EDintro}
\bega
\mathcal{E}_{tot}(t)+\int_0^t\mathcal{D}_{tot}(s)ds &\leq C\int_0^t\big(1+\|\nabla_x(\rho^{\e},u^{\e},\ta^{\e})(s)\|_{L^\infty_x}^2\big) \mathcal{E}_{tot}(s)ds + \e t.
\enda
\end{align}
Here $\mathcal{E}_{tot}^{\mathrm{N}}(F^\e(t))$ and $\mathcal{D}_{tot}^{\mathrm{N}}(F^\e(t))$ are defined in \eqref{EDtotdef}. 
In the 3-dimensional case, since the energy functional already controls $\|\nabla_x(\rho^{\e},u^{\e},\ta^{\e})\|_{L^\infty_x}\leq C\mathcal{E}_M^{\frac{1}{2}}$ for $\mathrm{N}>\frac{d}{2}+1$, this yields the construction of local-in-time solutions.

To the best of our knowledge, this provides the first result on the incompressible hydrodynamic limit in which solutions are constructed solely from the Boltzmann equation, without assuming any \emph{a priori} information on the limiting fluid variables. 
This approach is in the same spirit as the work of Nishida~\cite{Nishida}, where the hydrodynamic limit toward the compressible Euler equations was established.
We now aim to extend this result to a global-in-time hydrodynamic limit, particularly in the two-dimensional case, while overcoming the above regularity constraint. 
The main difficulties are summarized as follows.


\subsubsection*{B.6. New Difficulity: Penalization terms}


To extend the analysis to the global-in-time two-dimensional setting, obtaining a uniform bound for $\|\nabla_x(\rho^{\e},u^{\e},\ta^{\e})\|_{L^\infty_x}$ constitutes a critical threshold.
However, as seen from the macroscopic system \eqref{locconNew}, the evolution equations for $(\rho^{\e},u^{\e},\ta^{\e})$ contain singularly penalized terms, $\frac{1}{\e}\nabla_x\cdot u^{\e}$ and $\frac{1}{\e}\nabla_x(\rho^{\e}+\ta^{\e})$ which obstruct a direct control of this quantity.

To overcome this difficulty, we exploit macroscopic variables free from the singular penalization, namely the vorticity $\w^{\e}:= \nabla_x^\perp \cdot u^{\e}$ and the thermodynamic variable $\mathfrak{s}^{\e}=\frac{3}{2}\ta^{\e}-\rho^{\e}$. 

To estimate $\|\nabla_xu^{\e}\|_{L^\infty_x}$, we decompose the velocity field into vorticity part $\w^{\e}$ and irrotational part $\nabla_x \mathbb{P}^{\perp}u^{\e}$.
For the estimate of $\|\nabla_x\rho^{\e}\|_{L^\infty_x}$ and $\|\nabla_x\ta^{\e}\|_{L^\infty_x}$, we decompose the estimate by thermodynamic variables $\|\nabla_xp^{\e}\|_{L^\infty_x}$ and $\|\nabla_x\mathfrak{s}^{\e}\|_{L^\infty_x}$.

\hide
As a consequence, all singular penalized effects become confined to the acoustic variables
$\mathbb{P}^{\perp}u^{\e}$ and
$p^{\e}=\rho^{\e}+\ta^{\e}$.
This leads to a structural decomposition of the main obstacle quantity,
\begin{align}\label{rutdecomp-intro}
\bega
|\nabla_xu^{\e}(t)| \sim  \underbrace{| (\w^{\e}(t) , \nabla_x\mathbb{P}^{\perp}u^{\e}(t) ) | }_{ \textit{de-penalized variable} } + \underbrace{ | (\nabla_x\mathbb{P}^{\perp}u^{\e}(t) , \nabla_xp^{\e}(t) ) | }_{ \textit{acoustic (penalized) dynamics} }.
\enda
\end{align}
\begin{align}\label{rutdecomp-intro}
\bega
|\nabla_x(\rho^{\e},\ta^{\e})(t)| \sim  \underbrace{| \nabla_x\mathfrak{s}^{\e}(t) | }_{ \textit{de-penalized variable} } + \underbrace{ | (\nabla_x\mathbb{P}^{\perp}u^{\e}(t) , \nabla_xp^{\e}(t) ) | }_{ \textit{acoustic (penalized) dynamics} }.
\enda
\end{align}
This decomposition separates the dynamics into penalization-free and acoustic components.

However, such a separation alone is not sufficient: if the resulting estimates still depend on the energy–dissipation quantities without favorable scaling, the inequality \eqref{EDintro} cannot provide global-in-time control of the solution.
\unhide

\subsubsection*{B.7. Global Validity: Takeaway 1-- Quantitative control of the acoustic variables}
\hide
{\color{blue} 강하게 축약: It is well known that the acoustic variables $\mathbb{P}^{\perp}u^{\e}$ and $\rho^{\e}+\ta^{\e}$ satisfy a wave-type equation. However, deriving Strichartz estimates from the wave equation requires boundedness of the scaled time derivative of the initial data, namely $\|\e\p_t \mathbb{P}^{\perp}u^{\e}_0\| < \infty$. 
Such a condition is typically not well-prepared in the low-Mach-number limit of compressible fluid equations, since 
$\|\e\p_t(\rho^{\e},u^{\e},\ta^{\e})|_{t=0}\|
\sim \big\|\big(\nabla_x(\rho^{\e}+\ta^{\e}), \nabla_x\cdot u^{\e}\big)|_{t=0} \big\|$. 
On the other hand, in the hydrodynamic limit derived from the Boltzmann equation, this requirement becomes excessively restrictive due to the microscopic contribution:
\begin{align*}
\|\e\p_t(\rho^{\e},u^{\e},\ta^{\e})|_{t=0}\| &\sim  
\frac{1}{\e}\Big\|(\e\p_tF^{\e})|M^{\e}|^{-\frac{1}{2}}|_{t=0}\Big\| \sim
\frac{1}{\e\kappa}\Big\|\frac{1}{\e}\mathcal{L}(\AC{\P}F^{\e})|M^{\e}|^{-\frac{1}{2}}|_{t=0}\Big\| <\infty.
\end{align*}
Therefore, the initial data are not well-prepared unless one either avoids time derivatives altogether or works with at least the scaled derivative $\e^2\kappa\p_t$. 
For this reason, in hydrodynamic limit problems it is advantageous to employ the Strichartz estimates for the half-wave equation developed in \cite{Danchin}.} 
\unhide

For the acoustic variables, the local conservation laws 
\eqref{locconNew} yield the following half-wave system with forcing terms: (\eqref{locconP})
Precise definitions of $\bPhi^{\e}_{\rho+\ta}$ and $\bPhi^{\e}_{\mathbb{P}^{\perp}u}$ are given in Lemma~\ref{L.locconP}.
By exploiting the associated Strichartz estimates, the acoustic variables converge to zero due to dispersive effects, even when the initial acoustic components are not vanishing. 
Moreover, this approach provides the explicit convergence rate 
$\e^{\frac{1}{r}}$
in some Besov space ($L^r_T \dot{B}^{s}_{p,1}$ where $\dot{B}_{p,q}^s$ is defined in \eqref{Besovdef}) for $d=2,3$. 


\subsubsection*{B.8. Global Validity: Takeaway 2--Control of the vorticity and the entropic fluctuation $\mathfrak{s}^{\e} = \frac{3}{2} \theta^\e - \rho^\e$
}



The vorticity equation eliminates the singular penalized gradient term $\frac{1}{\e}\nabla_x(\rho^{\e}+\ta^{\e})$, while the evolution equation for the thermodynamic variable $\mathfrak{s}^{\e}$ removes the divergence penalization $\frac{1}{\e}\nabla_x\cdot u^{\e}$:
\begin{align}\label{wrtaeqn}
\bega
&\p_t\w^{\e} + u^{\e}\cdot \nabla_x \w^{\e} + (\nabla_x\cdot u^{\e})\w^{\e} +k_B \mathrm{\Theta}^{\e} (\nabla_x^{\perp}\ta^{\e})\cdot\nabla_x(\rho^{\e}+\ta^{\e}) +\frac{1}{\e^2}\nabla_x^{\perp}\cdot\bigg(\frac{1}{\mathrm{P}^{\e} }\nabla_x \cdot \mathbf{r}^{\e}\bigg) = 0, \cr
&\p_t \mathfrak{s}^{\e} + u^{\e}\cdot \nabla_x \mathfrak{s}^{\e} +\frac{1}{\e^2}\frac{1}{k_B\mathrm{P}^{\e}\mathrm{\Theta}^{\e} } \big(\nabla_x\cdot \mathfrak{q}^{\e} + \nabla_x \mathrm{U}^{\e} : \mathbf{r}^{\e} \big)=0.
\enda
\end{align}
This reduction allows the macroscopic dynamics to be analyzed through transport-type equations without the singular penalization.

For the vorticity, we employ the maximum principle to control $\|\w^{\e}(t)\|_{L^\infty_x}$ while potential-theory estimates provide bounds for $\|\w^{\e}(t)\|_{H^2_x}$ and $\|\nabla_x \mathbb{P}u^{\e}(t)\|_{L_x^\infty}$.

A central point of the present work concerns the scaling of the microscopic forcing appearing in the vorticity equation.
If, as discussed in Section~3.2, the dissipation were to lose a factor of $\kappa^{1/2}$ at each derivative level beyond $H^2$, even the basic $L^2$ estimate for the vorticity would fail. 
In that case, one would obtain $\frac{d}{dt}\|\w^{\e}\|_{L^2_x} \sim \mathcal{D}^{\frac{1}{2}}$ and inserting the resulting bound for $\|\nabla_x\mathbb{P}u^{\e}\|_{L^\infty_x}$ into the growth factor in \eqref{EDintro} would destroy the global-in-time control of the solution.
At higher regularity levels, the dissipation would become increasingly singular, preventing both $H^2$ and $L^\infty$ estimates for the vorticity since $\frac{d}{dt}\|\w^{\e}\|_{H^k} \sim \frac{1}{\e^2}\|\nabla_x^2\mathbf{r}^{\e}\|_{H^k} 
\les \kappa^{-k}\mathcal{D}^{\frac{1}{2}}$.

The key mechanism preventing this loss is the recovery of the full dissipation scale from the microscopic equation \eqref{Geqn0}.
Since the microscopic component $\AC{\P}F^{\e}$ is controlled at the optimal $\e^4\kappa$ scale by dissipation, the microscopic forcing appearing in the vorticity equation remains compatible with the energy–dissipation framework.
In contrast, one may derive a viscous term such as $-\kappa\Delta_x\w^{\e}$ through a Navier–Stokes expansion. However, this approach requires control of time derivatives of microscopic quantities, namely
$\kappa \nabla_x \mathcal{L}^{-1} (\partial_t \AC{\P}F^{\e})$,
and therefore implicitly assumes well-prepared initial data.
Our approach avoids this requirement and relies solely on the Euler-type moment equations, thereby allowing global-in-time validity without preparation assumptions.

A further difficulty arises from the fact that, unlike in the incompressible Euler equation, the velocity field transporting the vorticity remains compressible.
Although dispersive estimates yield favorable scaling for the acoustic variables in Besov spaces, such control is available only below the top-order derivative level due to forcing terms involving $\nabla_x\cdot \mathbf{r}^{\e}$ and $\nabla_x\cdot \mathfrak{q}^{\e}
$ in \eqref{locconP}.
Consequently, the estimates must be arranged so that the compressible component carries as few derivatives as possible.

Combining the refined dissipation scale with the dispersive gain obtained from the acoustic dynamics, the vorticity equation yields control of the growth of
$\|\w^{\e}(t)\|_{L^\infty_x}$ and
$\|\nabla_x\mathbb{P}u^{\e}(t)\|_{L^\infty_x}$
in terms only of the initial data, the time interval, and energy–dissipation quantities endowed with favorable scaling.
Once $\|\nabla_x u^{\e}(t)\|_{L^\infty_x}$ is controlled, the maximum principle further determines the growth of $\|\nabla_x\mathfrak{s}^{\e}(t)\|_{L^\infty_x}$.


 \newtheorem*{theoremG}{Theorem G}

\begin{theoremG}[Informal statement of Theorem \ref{T.2D.global}]\label{theoremG}
Let $\Omega=\mathbb{R}^2$. Assume that the initial data $\{F^\varepsilon_0\}$ satisfy a \emph{uniform-in-$\varepsilon$ scaled modulated entropy bound}, while they are allowed to blow up in strong topologies in the sense of the admissible blow-up condition (Definition~\ref{ABC}). Then
\begin{align}\label{eq:informal_micro_dissipation}
 \sum_{|\alpha| \leq N+1}  \frac{\kappa^{(|\alpha|-\mathrm{N})_+}}{\eps^4\kappa}
 \int_0^{T_*}\!\!\int_{\Omega \times \R^3} \nu \,
 \frac{\big|\p^{\alpha}\AC{\mathbf P} F^\varepsilon (t)\big|^2}{M^{\e}}
 \,dvdxdt
 \;\lesssim\;
\exp\!\Big(\exp\!\big(\exp(\|\omega_0^\varepsilon\|_\infty\, T_*)\big)\Big)\,\mathcal{E}_{\mathrm{tot}}(0)
 + \textit{l.o.t.}
\end{align}
\end{theoremG}

\unhide
\unhide

\unhide

\StopNoTOC

\section{Preliminaries}

In this section, we introduce several notations and equations for the macroscopic fields (local conservation laws), as well as the equation for the microscopic component 
$\AC{\P}F^{\e}$. We also define the energy and dissipation.


\begin{definition}
(1) We define a $5$-orthonormal basis with respect to the inner product $\langle \cdot ,  \frac{*}{M^{\e}} \rangle$, where $\langle   \ ,  \   \rangle $ denotes the standard $L^2_v$ inner product. For $i=1,2,3$,
\begin{align*}
\bega
e_0^{\e} = \frac{1}{\sqrt{\mathrm{P}^{\e}}}M^{\e}, \qquad e_{i}^{\e} = \frac{v_i-\mathrm{U}^{\e}_i}{\sqrt{k_B\mathrm{P}^{\e} \mathrm{\Theta}^{\e}}}M^{\e},  \qquad e_{4}^{\e} = \frac{1}{\sqrt{6 \mathrm{P}^{\e}}}\lw(\frac{|v-\mathrm{U}^{\e}|^2}{k_B\mathrm{\Theta}^{\e}}-3\rw)M^{\e}.
\enda
\end{align*}
We define the macroscopic and microscopic projections as follows:
\begin{align}\label{Pdef}
\bega
\P F&=\sum_{i=0}^4 \lw\la F, \frac{e_{i}^{\e}}{M^{\e}} \rw\ra e_{i}^{\e}, \qquad \AC{\P}F= (\mathbf{I}-\P)F.
\enda
\end{align}
Note that $\langle \AC{\P}F, (1, v, |v|^2) \rangle = (0,0,0)$.  

(2) We write the local Maxwellian in \eqref{M-def} depending on the macroscopic fields $(\mathrm{P}^\e,\mathrm{U}^{\e} ,\mathrm{\Theta}^{\e})$ by $M^{\e}=M_{[\mathrm{P}^\e,\mathrm{U}^{\e} ,\mathrm{\Theta}^{\e}]}$ and we denote the global Maxwellian $\mu(v):=M_{[1,0,1]}(v)$.

(3) Monatomic intermolecular collisions are described by the collision operator
\begin{align}\label{Qdef}
\mathcal{N}(F,G) 
:=   \frac{1}{2}\int_{\mathbb{R}^3}\int_{\mathbb{S}^2_+}B(v-v_*,w)\big[ F(v_*')G(v')+F(v')G(v_*') 
-F(v_*)G(v)-F(v)G(v_*)\big]dwdv_*, 
\end{align} 
where $v^\prime = v-((v-v_*) \cdot w) w$ and $v_*^\prime = v_*+((v-v_*) \cdot w) w$ denote the post-collisional velocities. 
Here, we consider the hard-sphere collision kernel $B(v-v_*,w):= |(v-v_*)\cdot w|$.

(4) We define the collision frequency
\begin{align}\label{nudef}
\nu(v)= \int_{\mathbb{R}^3}\int_{\mathbb{S}^2_+}B(v-v_*,w)\mu(v_*)dwdv_* , \qquad \nu^{\e}(t,x,v):= \nu\bigg(\frac{v-\mathrm{U}^{\e}(t,x)}{\sqrt{\mathrm{\Theta}^{\e}(t,x)}}\bigg).
\end{align}
\end{definition}


\begin{definition}\label{D.ABdef} Suppose $F^{\e}$ is given and the macroscopic fields \eqref{eq:macro_fields_intro} are well-defined. We define 
\begin{align}\label{albe-def}
\mathbf{r}_{ij}^{\e} := \langle F^{\e},  \mathfrak{R}^{\e}_{ij}   
  \rangle, \qquad
\mathfrak{q}_j^{\e} := \langle F^{\e}, \mathcal{Q}^{\e}_j   \rangle, 
\end{align}
\beq\label{ABdef}
\mathfrak{R}^{\e}_{ij} :=(v_i-\mathrm{U}^{\e}_i)(v_j-\mathrm{U}^{\e}_j)-\frac{|v-\mathrm{U}^{\e} |^2}{3}\delta_{ij},\qquad 
\mathcal{Q}^{\e}_i := (v_i-\mathrm{U}^{\e}_i)\frac{(|v-\mathrm{U}^{\e}|^2-5k_B\mathrm{\Theta}^{\e})}{2}.
\eeq
\end{definition}

\begin{definition} 
For a vector field $u:\R^d\rightarrow \R^d$ for $d\geq2$, we define the Hodge decomposition: 
\begin{align}\label{Hodge}
\bega
u(x)&= \mathbb{P}u(x) +\mathbb{P}^{\perp}u(x),
\enda
\end{align}
for the Leray projection and orthogonal Leray projection as follows: 
\begin{align}\label{LerayPdef}
\mathbb{P}u = u - \nabla_x\Delta_x^{-1}(\nabla_x\cdot u), \qquad \mathbb{P}^{\perp}u = \nabla_x\Delta_x^{-1}(\nabla_x\cdot u).
\end{align}
In particular, 
$\|u\|_{H^k}^2 = \|\mathbb{P}u\|_{H^k}^2 + \|\mathbb{P}^{\perp}u\|_{H^k}^2$.
\end{definition}

Recall \eqref{eq:macro_fields_intro} and \eqref{expform}. 
We refer to $(\rho^{\e},u^{\e},\ta^{\e})$ as the \emph{primitive variables}. In contrast, the macroscopic convergence established in this work is most naturally expressed in terms of the \emph{thermodynamic variables}
\begin{equation}\label{thermoV}
    p^{\e}:=\rho^{\e}+\ta^{\e},
\qquad
\mathfrak{s}^{\e}:=\frac{3}{2}\ta^{\e}-\rho^{\e},
\end{equation}
which represent, to first order in $\e$, the fluctuations of pressure and (dimensionless) specific entropy.
Indeed, for an ideal gas in physical units $\text{(full pressure)}=k_B\,\mathrm{P}^{\e}\mathrm{\Theta}^{\e}$, hence
\[
\text{full pressure}=k_B\,\mathrm{P}^{\e}\mathrm{\Theta}^{\e}
= k_B e^{\e(\rho^{\e}+\ta^{\e})}
= k_B\bigl(1+\e(\rho^{\e}+\ta^{\e})\bigr)+O(\e^2),
\]
so the linear pressure fluctuation (normalized by $k_B$) is $p^{\e}$, i.e.\ $\delta (\text{full pressure}) = k_B p^{\e}$.
For a monatomic ideal gas one may take (up to an additive constant) the specific entropy $
S^{\e} = k_B\Bigl(\frac{3}{2}\log \mathrm{\Theta}^{\e} - \log \mathrm{P}^{\e} \Bigr)$, and therefore
$
S^{\e}
= k_B\Bigl(\frac{3}{2}\log(e^{\e\ta^{\e}})-\log(e^{\e\rho^{\e}})\Bigr)
= \e k_B\Bigl(\frac{3}{2}\ta^{\e}-\rho^{\e}\Bigr),
$
so the linear entropy fluctuation is $\delta S^{\e} = k_B \mathfrak{s}^{\e}$.
The change of variables $(\rho^{\e},\ta^{\e})\leftrightarrow(p^{\e},\mathfrak{s}^{\e})$ is invertible: $\theta^{\e}=\frac{2}{5}\bigl(p^{\e}+\mathfrak{s}^{\e}\bigr)$, $\rho^{\e}=\frac{3}{5}p^{\e}-\frac{2}{5}\mathfrak{s}^{\e}.$


\hide
Due to the collision invariance, these moments formally solve the following system:
\begin{align*}
&\textbf{St}\p_t \mathrm{P}^{\e} + \nabla_x\cdot(\mathrm{P}^{\e}\mathrm{U}^{\e})=0, \cr
&\textbf{St}\p_t(\mathrm{P}^{\e}\mathrm{U}^{\e})
+\nabla_x\cdot\Big(\mathrm{P}^{\e}\,\mathrm{U}^{\e}\otimes\mathrm{U}^{\e}+p^{\e}I+\mathbf{r}^{\e}\Big)=0,
\cr
&\textbf{St}\p_t \mathrm{E}^{\e}
+\nabla_x\cdot\bigg(\mathrm{U}^{\e}\,\mathrm{E}^{\e}+\Big(p^{\e}I+\mathbf{r}^{\e}\Big)\mathrm{U}^{\e}+ \mathfrak{q}^{\e}\bigg)=0.
\end{align*}
where $\mathrm{E}^{\e}:= \frac{3}{2}k_B\mathrm{P}^{\e}\mathrm{\Theta}^{\e}+\frac{1}{2}\mathrm{P}^{\e}|\mathrm{U}^{\e}|^2 $ and $p^{\e}:= k_B \mathrm{P}^{\e} \Theta^{\e}$.
\textcolor{red}{($E$를 $1/2|v|^2F$에 해당하는 폼으로 썼습니다. 위에 처럼 쓰려면 $p^{\e}$이 중복 notation이 됨)}

$\p_t(\mathrm{P}^{\e},\mathrm{U}^{\e},\mathrm{\Theta}^{\e})$ 꼴은 \eqref{locconRUT}에 있습니다
\hide
\begin{CJK}{UTF8}{mj}\color{red}[CK6/5:  
추가해주세여: We can recast the second equation in the Euler–Reynolds type form:
\[
\textbf{St} \p_t\Big(\mathrm{P}^{\e} \mathrm{U}^{\e} \Big) +   \nabla \cdot \Big(\mathrm{P}^{\e}  \mathrm{U}^{\e}  \otimes  \mathrm{U}^{\e}
+ k_B \mathrm P ^{\e} \Theta^{\e} I + \mathbf{r}^{\e}
\Big)=0  \]
\unhide

마지막 방정식도 컨서베이션 폼으로 쓴다음에 마지막 식을 유도해주세여.
\[
\p_t  \Big( \frac{3 k_B}{2} \mathrm P^{\e} \Theta^{\e} + \frac{1}{2} \mathrm P ^{\e}| \mathrm U^{\e}|^2 \Big)
 + \nabla \cdot \Big(
 \frac{5}{2} k_B \mathrm P^{\e} \Theta^{\e} \mathrm U^{\e} + \frac{1}{2} \mathrm P^{\e} |\mathrm U^{\e}|^2 \mathrm U^{\e}
 + \mathbf q^{\e} + \mathbf{r}^{\e} \mathrm U^{\e}\Big)=0
\]

이 방정식도 써주시고 
\[
\p_t \Big( \frac{1}{2} \mathrm P^{\e} |\mathrm U^{\e}|^2\Big)
+ \nabla \cdot \Big(
\frac{1}{2} \mathrm P^{\e} |\mathrm U^{\e}|^2 \mathrm U^{\e} + k_B \mathrm P^{\e} \Theta^{\e} \mathrm U^{\e}
\Big) = k_B \mathrm P^{\e} \Theta ^{\e} \nabla \cdot \mathrm U^{\e} - \nabla \cdot \mathbf r^{\e} \cdot \mathrm U^{\e}
\]

그렇게 해서 직관적으로 유도되는 아래 폼이 더 좋은 듯 하네요. (결국 다 같은것이지만요)
\[
\p_t \Big( \frac{3 k_B}{2} \mathrm P^{\e} \Theta^{\e}\Big)
+ \nabla \cdot \Big(  \frac{3 k_B}{2} \mathrm P^{\e} \Theta ^{\e} \mathrm U^{\e} \Big)
+ k_B \mathrm P^{\e} \Theta^{\e} \nabla \cdot \mathrm U^{\e} = - \nabla \cdot \mathbf q^{\e}- \sum_{i,j} \p_{x_i} \mathrm U^{\e}_j \mathbf r_{ij}^{\e}
\]
여기서 $ \p_{x_i} \mathrm U^{\e}_j \mathbf r_{ij}^{\e} = \nabla \mathrm U^{\e} \! :  \mathbf r^{\e}$ 로 쓰면 더 깔끔하겠죠

그리고 
\[
\Big[\p_t + \mathrm U^{\e} \cdot \nabla_x  + \frac{2}{3} \nabla \cdot \mathrm U^{\e}  \Big]\Big( \frac{3 k_B}{2} \Theta^{\e}\Big)= -\frac{ \nabla \cdot \mathbf q^{\e}}{ \mathrm P^{\e}}  - \frac{\nabla \mathrm U^{\e}\! : \mathbf r^{\e} }{ \mathrm P^{\e}}
\]

요것도 어디 써져 있으면 좋을듯

\[
\p_t \mathrm U^{\e} + \mathrm U^{\e} \cdot \nabla \mathrm U^{\e}
+ k_B \frac{\Theta^{\e} \nabla \mathrm P^{\e}}{\mathrm P^{\e}} + k_B \nabla \Theta^{\e} = -  \frac{ \nabla \cdot \mathbf r ^{\e}}{\mathrm P^{\e}} 
\]

Compute 
\begin{align*}
    &\frac{[ \p_t +  \xi \cdot \nabla_x ] M}{M} \\
   = & \frac{ \p_t P  + U \cdot \nabla_x P + (\xi - U ) \cdot \nabla P }{P}\\
    & +   \Big(
    \frac{|\xi - U|^2}{ 2 k_B \Theta }
    - \frac{3}{2}\Big) \frac{ \p_t \Theta  + U \cdot \nabla_x \Theta + (\xi - U ) \cdot \nabla \Theta }{\Theta}\\
    &  + \frac{ (\p_t U  + U \cdot \nabla_x  U) \cdot (\xi - U) }{k_B \Theta}\\
    &+ \frac{ (\xi - U ) \cdot \nabla U \cdot (\xi -  U)}{k_B \Theta}\\
  =  &  \frac{ \p_t P + U \cdot \nabla P + P \nabla \cdot U  }{P}\\
  &+ (\xi - U) \cdot \Big(
  \frac{\nabla P}{P} + \frac{\p_t U + U \cdot \nabla U}{k_B \Theta} + \frac{\nabla \Theta}{\Theta}
  \Big)\\
  &+ \left(\frac{|\xi - U|^2 - 3 k_B \Theta }{2 k_B \Theta}\right) \frac{ \p_t \Theta + U \cdot \nabla \Theta
  + \frac{2}{3} \Theta \nabla \cdot U
  }{\Theta}\\
  &+ \frac{|\xi -U|^2 - 5 k_B \Theta}{2k_B \Theta} (\xi - U) \cdot \frac{\nabla \Theta}{\Theta}\\
  & + \frac{\nabla U}{k_B \Theta}: \left(
  (\xi - U) \otimes (\xi - U) - \frac{| \xi - U|^2}{3}I
  \right) 
\end{align*}

\begin{align*}
     \frac{[ \p_t +  \xi \cdot \nabla_x ] M}{M}  
    =&  -\frac{ \nabla \cdot \mathbf r^{\e}}{k_B \mathrm P^{\e} \Theta^{\e} }\cdot(\xi -  \mathrm U^{\e})  - \frac{2 (  {\nabla \cdot q^{\e}}  +   {\nabla  \mathrm U^{\e}: \mathbf r^{\e}} )   }{ 3 k_B \mathrm  P^{\e} \Theta^{\e}  } \frac{|\xi -  \mathrm U^{\e}|^2 - 3 k_B \Theta^{\e}}{2 k_B \Theta^{\e}} \\
    & + \frac{\nabla \Theta^{\e}}{ k_B (\Theta^{\e})^2} \cdot  \mathcal Q^{\e}
   + \frac{\nabla  \mathrm U^{\e}}{k_B \Theta^{\e}}:
    \mathfrak{R}^{\e}
    \end{align*}

] \end{CJK}
\unhide

\begin{lemma}\label{L.loccon}
Recall that the macroscopic fields $(\mathrm{P}^{\e},\mathrm{U}^{\e},\mathrm{\Theta}^{\e})$ associated with the Boltzmann solution $F^{\e}$ formally solve \eqref{loccon}, \hide
\begin{align}\label{loccon}
\bega
&\e\p_t \mathrm{P}^{\e} + \nabla_x\cdot(\mathrm{P}^{\e}\mathrm{U}^{\e})=0, \cr
&\e\p_t(\mathrm{P}^{\e}\mathrm{U}^{\e})
+\nabla_x\cdot\Big(\mathrm{P}^{\e}\,\mathrm{U}^{\e}\otimes\mathrm{U}^{\e}+k_B \mathrm{P}^{\e} \Theta^{\e}I+\mathbf{r}^{\e}\Big)=0,
\cr
&\e\p_t \mathrm{E}^{\e}
+\nabla_x\cdot\bigg(\mathrm{U}^{\e}\,\mathrm{E}^{\e}+\Big(k_B \mathrm{P}^{\e} \Theta^{\e}I+\mathbf{r}^{\e}\Big)\mathrm{U}^{\e}+ \mathfrak{q}^{\e}\bigg)=0,
\enda
\end{align}
where $\mathrm{E}^{\e}:= \frac{3}{2}k_B\mathrm{P}^{\e}\mathrm{\Theta}^{\e}+\frac{1}{2}\mathrm{P}^{\e}|\mathrm{U}^{\e}|^2$, 
$I$ denotes the identity matrix, and $\mathbf{r}^{\e}$ and $\mathfrak{q}^{\e}$ are defined in \eqref{albe-def}.\unhide
\hide
\begin{align}
\bega
&\e\p_t\mathrm{P}^{\e} +\nabla_x \cdot (\mathrm{P}^{\e}\mathrm{U}^{\e})=0, \cr
&\e\p_t(\mathrm{P}^{\e} \mathrm{U}^{\e}) + \nabla_x\cdot\big(\mathrm{P}^{\e}\,\mathrm{U}^{\e}\otimes\mathrm{U}^{\e}\big) +k_B\nabla_x(\mathrm{P}^{\e} \mathrm{\Theta}^{\e}) + \nabla_x\cdot \mathbf{r}^{\e} =0, \cr  
&\e\p_t\bigg(\frac{3k_B}{2}\mathrm{P}^{\e}\mathrm{\Theta}^{\e}\bigg)
+  \nabla_x\cdot\bigg(\frac{3k_B}{2}\mathrm{P}^{\e}\mathrm{\Theta}^{\e} \mathrm{U}^{\e}\bigg) + k_B
\mathrm{P}^{\e} \mathrm{\Theta}^{\e}(\nabla_x \cdot \mathrm{U}^{\e})  + \nabla_x\cdot \mathfrak{q}^{\e} +\nabla \mathrm U^{\e} \! : \mathbf{r}^{\e} =0,
\enda
\end{align}
\unhide
 which are equivalent to the following system if we express the variables in the form \eqref{expform} for $(\rho^\e, u^\e, \theta^\e)$:
\begin{align}\label{locconNew}
\bega
&\p_t \rho^{\e} +u^{\e}\cdot\nabla_x\rho^{\e} + \frac{1}{\e}\nabla_x\cdot u^{\e} =0 ,\cr
&\frac{1}{k_B\mathrm{\Theta}^{\e}}\bigg[\p_tu^{\e} + u^{\e}\cdot \nabla_x u^{\e} \bigg] +\frac{1}{\e}\nabla_x(\rho^{\e}+\ta^{\e})
+\frac{1}{\e^2}\frac{1}{k_B\mathrm{P}^{\e} \mathrm{\Theta}^{\e}}\nabla_x \cdot \mathbf{r}^{\e} =0, \cr  
&\frac{3}{2}\bigg[\p_t \ta^{\e} +u^{\e}\cdot\nabla_x\ta^{\e} \bigg] + \frac{1}{\e}\nabla_x\cdot u^{\e}  +\frac{1}{\e^2}\frac{1}{k_B\mathrm{P}^{\e}\mathrm{\Theta}^{\e} } \big(\nabla_x\cdot \mathfrak{q}^{\e} + \nabla_x \mathrm{U}^{\e} : \mathbf{r}^{\e} \big)=0.
\enda
\end{align}
Here, $\mathbf{r}^{\e}$ and $\mathfrak{q}^{\e}$ are defined in \eqref{albe-def}.
\end{lemma}

\hide 기본 계산이므로 유도 생략
\begin{proof}
From \eqref{loccon}, we can rewrite the equation as 
\begin{align*}
\bega
&\eps\p_t\mathrm{P}^{\e} +\nabla_x \cdot (\mathrm{P}^{\e}\mathrm{U}^{\e})=0, \cr
&\eps\mathrm{P}^{\e} \p_t\mathrm{U}^{\e} +\mathrm{P}^{\e}  \mathrm{U}^{\e}
\cdot\nabla_x \mathrm{U}^{\e} +k_B\nabla_x(\mathrm{P}^{\e} \mathrm{\Theta}^{\e})
+\sum_{j} \p_{x_j} \mathbf{r}_{ij}^{\e} =0, \cr  
&\eps\mathrm{P}^{\e}\p_t\mathrm{\Theta}^{\e}
+\mathrm{P}^{\e} \mathrm{U}^{\e}
\cdot\nabla_x \mathrm{\Theta}^{\e}+ \frac{2}{3}\nabla_x \cdot \mathrm{U}^{\e}
(\mathrm{P}^{\e} \mathrm{\Theta}^{\e})  +\frac{2}{3k_B}\sum_{j} \p_{x_j} \mathfrak{q}_j^{\e} + \frac{2}{3k_B}\sum_{i,j} \p_{x_i}\mathrm{U}^{\e}_j \mathbf{r}_{ij}^{\e}=0.
\enda
\end{align*}
We divide the second and third equation by $k_B\mathrm{P}^{\e} \mathrm{\Theta}^{\e}$ and $\mathrm{P}^{\e}$, respectively. Then, applying \eqref{expform}, we have 
\begin{align*}
\bega
&\mathrm{P}^{\e}\eps^2 \bigg[ \p_t \rho^{\e} +u^{\e}\cdot\nabla_x\rho^{\e} + \frac{1}{\e}\nabla_x\cdot u^{\e} \bigg] =0 ,\cr
&\frac{1}{k_B\mathrm{\Theta}^{\e}} \e^2 \bigg[ \p_tu^{\e} + u^{\e}\cdot \nabla_x u^{\e} \bigg]  +\frac{1}{\mathrm{P}^{\e} \mathrm{\Theta}^{\e}}\nabla_x(\mathrm{P}^{\e} \mathrm{\Theta}^{\e}) 
+\frac{1}{k_B\mathrm{P}^{\e} \mathrm{\Theta}^{\e}}\sum_{j} \p_{x_j} \mathbf{r}_{ij}^{\e} =0 ,\cr  
&\mathrm{\Theta}^{\e} \e^2 \bigg[ \p_t \ta^{\e} +u^{\e}\cdot\nabla_x\ta^{\e} + \frac{2}{3}\frac{1}{\e}\nabla_x\cdot u^{\e} \bigg]  +\frac{2}{3k_B\mathrm{P}^{\e}}\sum_{j} \p_{x_j} \mathfrak{q}_j^{\e} + \frac{2}{3k_B\mathrm{P}^{\e}}\sum_{i,j} \p_{x_i}\mathrm{U}^{\e}_j \mathbf{r}_{ij}^{\e}=0.
\enda
\end{align*}
Using \eqref{expform}, the pressure in the second equation becomes
\begin{align*}
\bega
\frac{1}{\mathrm{P}^{\e} \mathrm{\Theta}^{\e}}\nabla_x(\mathrm{P}^{\e} \mathrm{\Theta}^{\e}) = \nabla_x \log (\mathrm{P}^{\e} \mathrm{\Theta}^{\e}) = \e\nabla_x(\rho^{\e}+\ta^{\e}).
\enda
\end{align*}
Thus we have \eqref{locconNew}.
\end{proof}\unhide

\begin{definition}\label{D.bWdef}
For $x\in \R^d$, we define the $(d+2)$-dimensional macroscopic fields as
\begin{align}\label{bWdef}
\bW^{\e}(t,x) =[\rho^{\e}(t,x),u^{\e}(t,x),\ta^{\e}(t,x)]^T.
\end{align}
We define the coefficients matrix of \eqref{locconNew} and the linear differential operator as
\begin{align}\label{Sdef}
\mathcal{S}^{\e}(t,x) = \lw[ \begin{array}{ccc} 1 & 0 & 0 \\ 0 & \frac{1}{k_B\mathrm{\Theta}^{\e}(t,x)}I_{d\times d} & 0 \\  0 & 0 & \frac{3}{2} \end{array} \rw], \qquad L(\p_x) = \lw[ \begin{array}{ccc} 0 & \div & 0 \\ 
\nabla_x &0_{d\times d} & \nabla_x \\  
0 & \div & 0 \end{array} \rw].
\end{align}
For brevity, we also write $\mathcal{S}^{\e}$ and $(\mathcal{S}^{\e})^{-1}$ as follows: 
\begin{align*}
\mathcal{S}^{\e}(t,x) = \diag\bigg(1,\frac{1}{k_B\mathrm{\Theta}^{\e}(t,x)}\vec{1}_d,\frac{3}{2}\bigg), \qquad (\mathcal{S}^{\e})^{-1}(t,x) = \diag\bigg(1,k_B\mathrm{\Theta}^{\e}(t,x)\vec{1}_d,\frac{2}{3}\bigg),
\end{align*}
where $\vec{1}_d=(1,\cdots,1)$ is a $d$-dimensional vector whose all components are 1.
We denote the macroscopic forcing term of \eqref{locconNew} by 
\begin{align}\label{gWdef}
\bPhi^{\e}_{\bW}(t,x) = \lw[ \begin{array}{ccc} 0 \\ -\frac{1}{\e^2}\frac{1}{k_B\mathrm{P}^{\e}(t,x) \mathrm{\Theta}^{\e}(t,x)}\sum_{j} \p_{x_j} \mathbf{r}_{ij}^{\e}(t,x)
\\  -\frac{1}{\e^2}\frac{1}{k_B\mathrm{P}^{\e}(t,x)\mathrm{\Theta}^{\e}(t,x) }\Big(\sum_{j} \p_{x_j} \mathfrak{q}_j^{\e}(t,x) + \sum_{i,j} \p_{x_i}\mathrm{U}^{\e}_j(t,x) \mathbf{r}_{ij}^{\e}(t,x)\Big) \end{array} \rw].
\end{align}

Then a system of the local conservation laws \eqref{locconNew} is equivalent to 
\begin{align}\label{bWeqn}
\mathcal{S}^{\e} (\p_t \bW^{\e} + u^{\e} \cdot \nabla_x \bW^{\e}) + \frac{1}{\e} L(\p_x) \bW^{\e} = \bPhi^{\e}_{\bW}.
\end{align}
\end{definition}

\begin{lemma}
For macroscopic fields satisfying the local conservation laws \eqref{loccon}, the local Maxwellian defined in \eqref{M-def} satisfies
\begin{align}\label{Phitx}
\bega
\p_t M^{\e} + \frac{v}{\eps}\cdot\nabla_x M^{\e} &=  \frac{1}{k_B\mathrm{\Theta}^{\e}}\sum_{i,j}\p_{x_i}u^{\e}_j \mathfrak{R}^{\e}_{ij}M^{\e} +\frac{1}{k_B\mathrm{\Theta}^{\e}}\sum_i \p_{x_i}\ta^{\e} \mathcal{Q}^{\e}_iM^{\e}  \cr 
&-\frac{1}{\e k_B\mathrm{P}^{\e}\mathrm{\Theta}^{\e}}\sum_{i,j} \p_{x_j} \mathbf{r}_{ij}^{\e} (v_i-\mathrm{U}^{\e}_i) M^{\e}\cr 
&-\frac{2}{\eps3k_B \mathrm{P}^{\e}\mathrm{\Theta}^{\e}} \bigg(\sum_{j} \p_{x_j} \mathfrak{q}_j^{\e} + \sum_{i,j} \p_{x_i}\mathrm{U}^{\e}_j \mathbf{r}_{ij}^{\e} \bigg) \lw(\frac{|v-\mathrm{U}^{\e}|^2-3k_B \mathrm{\Theta}^{\e}}{2k_B \mathrm{\Theta}^{\e}}\rw)M^{\e}.
\enda
\end{align}
The microscopic component $\AC{\P}F^{\e}$ satisfies
\begin{align}\label{Geqn0}
\bega
&\p_t \AC{\P}F^{\e}+\frac{v}{\e}\cdot\nabla_x \AC{\P}F^{\e}  + \frac{1}{\kappa\e^2}\mathcal{L}(\AC{\P}F^{\e})
- \frac{1}{\kappa\e^2}\mathcal{N}(\AC{\P}F^{\e},\AC{\P}F^{\e}) 
\\
& = - \frac{1}{k_B\mathrm{\Theta}^{\e}}\sum_{i,j}\p_{x_i}u^{\e}_j \mathfrak{R}^{\e}_{ij}M^{\e}-\frac{1}{k_B\mathrm{\Theta}^{\e}}\sum_i \p_{x_i}\ta^{\e} \mathcal{Q}^{\e}_iM^{\e} + \frac{1}{\e k_B\mathrm{P}^{\e}\mathrm{\Theta}^{\e}}\sum_{i,j} \p_{x_j} \mathbf{r}_{ij}^{\e} (v_i-\mathrm{U}^{\e}_i)M^{\e} \cr 
& +\frac{2}{\e 3k_B\mathrm{P}^{\e}\mathrm{\Theta}^{\e}} \bigg(\sum_{j} \p_{x_j} \mathfrak{q}_j^{\e} + \sum_{i,j} \p_{x_i}\mathrm{U}^{\e}_j \mathbf{r}_{ij}^{\e} \bigg) \lw( \frac{|v-\mathrm{U}^{\e}|^2-3k_B\mathrm{\Theta}^{\e}}{2k_B\mathrm{\Theta}^{\e}}\rw)M^{\e},
\enda
\end{align}
where $\mathbf{r}^{\e}_{ij}$, $\mathfrak{q}^{\e}_i$, $\mathfrak{R}^{\e}_{ij}$, and $\mathcal{Q}^{\e}_i$ are defined in Definition \ref{D.ABdef}.
\end{lemma}
\begin{proof}
Once \eqref{Phitx} holds, it is direct to see \eqref{Geqn0} from the decomposition \eqref{Pdef} and equation \eqref{BE}. 
Now we prove \eqref{Phitx}. We rewrite the local conservation laws \eqref{loccon} as follows:
\begin{align}\label{locconRUT}
\bega
&\eps\p_t\mathrm{P}^{\e} +\nabla_x \cdot (\mathrm{P}^{\e}\mathrm{U}^{\e})=0, \cr
&\eps \p_t\mathrm{U}^{\e} + \mathrm{U}^{\e}
\cdot\nabla_x \mathrm{U}^{\e} +\frac{1}{\mathrm{P}^{\e}}k_B\nabla_x(\mathrm{P}^{\e} \mathrm{\Theta}^{\e})
+\frac{1}{\mathrm{P}^{\e}}\sum_{j} \p_{x_j} \mathbf{r}_{ij}^{\e} =0, \cr  
&\eps\p_t \mathrm{\Theta}^{\e}
+\mathrm{U}^{\e}
\cdot \nabla_x \mathrm{\Theta}^{\e}+ \frac{2}{3} \mathrm{\Theta}^{\e} \nabla_x \cdot \mathrm{U}^{\e} +\frac{2}{3k_B\mathrm{P}^{\e}}\bigg(\sum_{j} \p_{x_j} \mathfrak{q}_j^{\e} + \sum_{i,j} \p_{x_i}\mathrm{U}^{\e}_j \mathbf{r}_{ij}^{\e}\bigg)=0.
\enda
\end{align}
By a direct computation, we have
\begin{align*}
\bega
\bigg(\frac{\p M^{\e}}{\p\mathrm{P}^{\e}},\frac{\p M^{\e}}{\p\mathrm{U}^{\e}},\frac{\p M^{\e}}{\p\mathrm{\Theta}^{\e}}\bigg) = \bigg(\frac{1}{\mathrm{P}^{\e}}, \frac{(v-\mathrm{U}^{\e})}{k_B \mathrm{\Theta}^{\e}},\frac{|v -\mathrm{U}^{\e}|^2-3k_B \mathrm{\Theta}^{\e}}{2k_B |\mathrm{\Theta}^{\e}|^2}\bigg) M^{\e},
\enda
\end{align*}
and
\begin{align*}
&\Big(\p_t + \frac{v}{\e}\cdot\nabla_x\Big)M^{\e}  = \bigg(\p_t + \frac{\mathrm{U}^{\e}}{\e}\cdot\nabla_x + \frac{v-\mathrm{U}^{\e}}{\e}\cdot\nabla_x\bigg)M^{\e} \cr
&=  \frac{1}{\e}\bigg[\frac{(\e\p_t\mathrm{P}^{\e}+\mathrm{U}^{\e}\cdot\nabla_x\mathrm{P}^{\e})}{\mathrm{P}^{\e}}+\frac{(\e\p_t\mathrm{U}^{\e}+\mathrm{U}^{\e}\cdot\nabla_x\mathrm{U}^{\e})\cdot (v-\mathrm{U}^{\e})}{k_B\mathrm{\Theta}^{\e}} \cr 
&+\frac{(\e\p_t \mathrm{\Theta}^{\e}+\mathrm{U}^{\e}\cdot\nabla_x \mathrm{\Theta}^{\e})}{\mathrm{\Theta}^{\e}} \lw( \frac{|v-\mathrm{U}^{\e}|^2-3k_B\mathrm{\Theta}^{\e}}{2k_B\mathrm{\Theta}^{\e}} \rw)\bigg]M^{\e} \cr 
&+ \frac{1}{\e}\bigg[\frac{\nabla_x\mathrm{P}^{\e}}{\mathrm{P}^{\e}}\cdot (v-\mathrm{U}^{\e}) + \sum_{i,j}\frac{\p_i\mathrm{U}^{\e}_j(v_i-\mathrm{U}^{\e}_i)(v_j-\mathrm{U}^{\e}_j)}{k_B\mathrm{\Theta}^{\e}}+\frac{\nabla_x \mathrm{\Theta}^{\e}}{\mathrm{\Theta}^{\e}}\cdot(v-\mathrm{U}^{\e})\frac{|v-\mathrm{U}^{\e}|^2-3k_B\mathrm{\Theta}^{\e}}{2k_B\mathrm{\Theta}^{\e}}\bigg]M^{\e}.
\end{align*}
Applying the local conservation laws \eqref{locconRUT}, together with 
$\frac{\nabla_x \mathrm{P}^{\e}}{\mathrm{P}^{\e}} = \e\nabla_x\rho^{\e}$, 
$\frac{\nabla_x \mathrm{\Theta}^{\e}}{\mathrm{\Theta}^{\e}} = \e\nabla_x\ta^{\e}$, 
and the definitions of $\mathfrak{R}^{\e}_{ij}$ and $\mathcal{Q}^{\e}_i$ in \eqref{ABdef}, 
we obtain the result \eqref{Phitx}.

\hide
where we used $\mathrm{P}^{\e}=e^{\e\rho^{\e}}$, $\mathrm{\Theta}^{\e}=e^{\e\ta^{\e}}$, so that $\p\mathrm{P}^{\e}=\e(\p\rho^{\e}) e^{\e\rho^{\e}}$, $\p\mathrm{\Theta}^{\e}=\e(\p\ta^{\e}) e^{\e\ta^{\e}}$. Using .., we have
\begin{align}\notag
\bega
&\left(\p_tM^{\e} + \frac{v}{\e}\cdot\nabla_x M^{\e}\right) \cr
& = \bigg[\e(\p_t\rho^{\e}+u^{\e}\cdot\nabla_x\rho^{\e})+\e(\p_tu^{\e}+u^{\e}\cdot\nabla_xu^{\e})\cdot\frac{(v-\mathrm{U}^{\e})}{k_B\mathrm{\Theta}^{\e}}+\e(\p_t\ta^{\e}+u^{\e}\cdot\nabla_x\ta^{\e}) \lw(\frac{|v-\mathrm{U}^{\e}|^2-3k_B \mathrm{\Theta}^{\e}}{2k_B \mathrm{\Theta}^{\e}}\rw)\bigg] M^{\e} \cr 
&+ \bigg[\nabla_x\rho^{\e}\cdot (v-\mathrm{U}^{\e})+\sum_{i,j}\frac{\p_{x_i}u^{\e}_j(v_i-\mathrm{U}^{\e}_i)(v_j-\mathrm{U}^{\e}_j)}{k_B\mathrm{\Theta}^{\e}}+\nabla_x\ta^{\e}\cdot(v-\mathrm{U}^{\e}) \lw(\frac{|v-\mathrm{U}^{\e}|^2-3k_B \mathrm{\Theta}^{\e}}{2k_B \mathrm{\Theta}^{\e}}\rw)\bigg]M^{\e} .
\enda
\end{align}
Applying \eqref{locconNew} and using the definition of $\mathfrak{R}^{\e}_{ij}$ and $\mathcal{Q}^{\e}_i$ in \eqref{ABdef}, we have \eqref{Phitx}. 
\unhide
\hide \\(Proof of \eqref{Geqn0})
We substitute $M^{\e}+\AC{\P}F^{\e}$ for $F^{\e}$ to \eqref{BE} to have
\beq \notag
\p_t \AC{\P}F^{\e}+\frac{v}{\e}\cdot\nabla_x \AC{\P}F^{\e}= -\bigg(\p_t M^{\e}+\frac{v}{\e}\cdot\nabla_x M^{\e}\bigg) + \frac{2}{\kappa\eps^2}\mathcal{N}(\AC{\P}F^{\e},M^{\e})+\frac 1 {\kappa\eps^2}\mathcal{N}(\AC{\P}F^{\e},\AC{\P}F^{\e}).
\eeq
Applying \eqref{Phitx}, we derive \eqref{Geqn0}.\unhide
\end{proof}

\subsection{Energy and Dissipation}

Before defining the energy and dissipation, we distinguish two cases depending on whether the velocity field belongs to $L^2$. This distinction arises from the structure of the initial vorticity. In particular, in the 2-dimensional case the velocity field generated by a nonzero total vorticity may fail to lie in $L^2$. Accordingly, we introduce the finite velocity energy case and the infinite velocity energy case.

\begin{definition}
We define the ``\emph{finite velocity energy case}" by
\begin{align}\label{caseEC}
\bega
\int_{\R^2} |u|^2 dx <\infty, \quad \mbox{for} \quad d=2,3.
\enda
\end{align}
We define the ``\emph{infinite velocity energy case}" for $d=2$ by
\begin{align}\label{caseECX}
\int_{\R^2} |u|^2 dx =\infty, \quad \int_{K} |u|^2 dx \leq C , \quad \mbox{for any compact set} \quad K\subset \R^2.
\end{align}
\end{definition}

\hide
\begin{remark}\label{Rmk.infH}
When the velocity energy is infinite, the initial relative entropy is infinite. This can be seen from the following bound:
\begin{align}\label{F-mu,split2}
\bega
\|\rho^{\e}_0\|_{L^2(\R^d)}^2 +\|u^{\e}_0\|_{L^2(\R^d)}^2 + \|\ta^{\e}_0\|_{L^2(\R^d)}^2 &\leq \frac{C}{\e^2}\{ \mathcal{H}(F^{\e}_0) - \mathcal{H}(\mu)\} + \frac{C}{\e^2}\int_{\R^d \times\R^3}\frac{| \AC{\P}F^{\e}(0)|^2}{M^{\e}(0)} dvdx \cr 
&+ \e^2\Big(\|\rho^{\e}_0\|_{L^4(\R^d)}^4 +\|u^{\e}_0\|_{L^4(\R^d)}^4 + \|\ta^{\e}_0\|_{L^4(\R^d)}^4\Big).
\enda
\end{align}
\end{remark}
\unhide

\hide
$\bullet$ Note: When the global energy is infinite \eqref{caseECX}, the initial relative entropy is infinite. 

\begin{lemma} {\color{red} 정말 필요? 위의 리마크로 충분하지 않음?} Let $\rho_0^{\e}, \ta_0^{\e} \in C_c^{\infty}(\R^d)$ and $\AC{\P}F^{\e} \in C_c^{\infty}(\R^d\times\R^d)$. For the initial velocity $u_0$ having infinite velocity energy  \eqref{caseECX}, the mollified initial data satisfies  
\begin{align}\label{F-mu,split2}
\bega
\|\rho^{\e}_0\|_{L^2(\R^d)}^2 +\|u^{\e}_0\|_{L^2(\R^d)}^2 + \|\ta^{\e}_0\|_{L^2(\R^d)}^2 &\leq \frac{C}{\e^2}\{ \mathcal{H}(F^{\e}_0) - \mathcal{H}(\mu)\} + \frac{C}{\e^2}\int_{\R^d \times\R^3}\frac{| \AC{\P}F^{\e}(0)|^2}{M^{\e}(0)} dvdx \cr 
&+ \e^2\Big(\|\rho^{\e}_0\|_{L^4(\R^d)}^4 +\|u^{\e}_0\|_{L^4(\R^d)}^4 + \|\ta^{\e}_0\|_{L^4(\R^d)}^4\Big).
\enda
\end{align}
Which means 
\begin{align}
\bega
\infty &= \frac{C}{\e^2}\{ \mathcal{H}(F^{\e}_0) - \mathcal{H}(\mu)\} 
\enda
\end{align}
for each $\e$ ensuring the existence of the solution.
\end{lemma} \unhide

\hide
{\color{blue}
\begin{itemize}
\item 본문에 쓴 노테이션: $\widetilde{u}^{\e}$, $\bar{u}$.
\item 이니셜데이터에 대해 Radial energy decomp를 정의해주면, 위 노테이션들이 정의가 안됨.
\item full velocity $u^{\e}$ 으로부터 정의를 하려면, $\mathbb{P}u^{\e}$ 에서 정의해야됨. 
\item 제너럴 $u$에 대해 radial-energy decomposition을 정의해주고, 우리의 $u^{\e}$이 나뉘어지는데, $\bar{u}$가 $\e$, $t$에 디펜드 안한다를 보여주는게 나을거 같습니다. 
\end{itemize}
}
\unhide

\begin{definition}\label{D.Ra-E}\cite{MaBe} For the infinite velocity energy case defined in \eqref{caseECX}, and for a divergence free vector field $u:\R^2\to\R^2$ satisfying $\w(x)= \nabla_x^{\perp}\cdot u(x)$ with $\w \in L^1(\R^2)$, there exist its radial-energy decomposition as
\begin{align*}
\bega
u(x)&= \widetilde{u}(x) + \bar{u}(x),
\enda
\end{align*}
where $\bar{u}$ (called the ``radial eddy'') and $\widetilde{u}$ satisfy
\begin{align*}
\bega
\int_{\R^2} |\widetilde{u}(x)|^2 dx <\infty, \qquad \div(\widetilde{u})=0, \qquad \bar{u}(x) = \frac{x^{\perp}}{|x|^2} \int_0^{|x|}r\bar{\w}(r)dr, \qquad \div(\bar{u})=0.
\enda
\end{align*}
Here, $\bar{\w}(\cdot)$ is a smooth radially symmetric vorticity satisfying
\begin{align}\label{rw-def}
\bega
&\int_{\R^2} \w(y) dy = \int_{\R^2} \bar{\w}(|y|) dy.
\enda
\end{align}
In particular, we can choose $\bar{\w}(|x|)\sim |x|^{\mathrm{N}}$ near $|x|=0$, $\| |x|^{-k} |\nabla_x^{\mathrm{N}-k}\bar{\w}(|x|)|\|_{L^2_x}\leq C$, for $0\leq k\leq\mathrm{N}$.
\end{definition}

\begin{lemma}
Let $\{u^{\e}\}_{\e>0}$ satisfy the infinite velocity energy case \eqref{caseECX}, and assume that the initial total vorticity $\int_{\R^2} \w^{\e}_0 dx$ is independent of $\e$. Then the velocity $u^{\e}$ admits the following radial-energy decomposition:
\begin{align}\label{uexpand}
\bega
u^{\e}(t,x)= \widetilde{u}^{\e}(t,x)+\bar{u}(x)+\mathbb{P}^{\perp}u^{\e}(t,x).
\enda
\end{align}
Here, the radial eddy $\bar{u}$ is independent of $\e$ and $t$, and is generated by an $\e$-independent radial vorticity $\bar{\w}$:
\begin{align}\label{barudef}
\bar{u}(x) := \begin{cases} 
0, \quad &\text{for the finite velocity energy case} \quad \eqref{caseEC}, \\
\frac{x^{\perp}}{|x|^2} \int_0^{|x|} r \bar{\w}(r) \, dr, \quad &\text{for the infinite velocity energy case} \quad \eqref{caseECX}.
\end{cases}
\end{align}
\end{lemma}
\begin{proof}
For each $\e>0$ and $t\geq 0$, applying the Hodge decomposition \eqref{Hodge} and the radial energy decomposition in Definition~\ref{D.Ra-E}, we obtain $u^{\e}(t,x) = \widetilde{u}^{\e}(t,x)+\bar{u}^{\e}(t,x) + \mathbb{P}^{\perp}u^{\e}(t,x)$. 
We note that the Boltzmann vorticity $\w^{\e}(t,x):=\nabla_x^\perp \cdot u^{\e}(t,x)$ preserves its circulation $\int_{\R^2}\w^{\e}(t,x)dx$ (See Lemma \ref{L.util}).
Moreover, by assumption, the initial circulation is independent of $\e>0$. Hence, there exists a radially symmetric vorticity $\bar{\w}$, independent of $\e$, such that $\int_{\R^2}\w^{\e}(t)dx=\int_{\R^2}\w^{\e}_0dx = \int_{\R^2} \bar{\w}(|y|)dy$.
By the construction in \eqref{rw-def}, it follows that $\bar{u}^{\e}(t,x)$ is independent of $\e$ and $t$.
\end{proof}

\begin{remark}
We can approximate initial data to $\int_{\R^2} \w^{\e}_0 dx$ is independent of $\e$. See Appendix \ref{A.molli}.
\end{remark}

We say that $u_0^{\e}$ satisfies the finite (or infinite) velocity energy case if the family $\{u_0^{\e}\}_{\e>0}$ satisfies the corresponding condition.

\begin{lemma}\label{L.barubdd}
In the infinite velocity energy case \eqref{caseEC}, the radial eddy $\bar{u}$ defined in Definition \ref{D.Ra-E} satisfies
\begin{align*}
\bega
\|\bar{u}\|_{L^\infty_x} \leq C, \quad \|\nabla_x \bar{u}\|_{H^{\mathrm{N}}_x} \leq C, \quad \|\nabla_x \bar{u}\|_{L^\infty_x} \leq C, \quad \|\bar{u}\|_{L^p_x} \leq C\Big(\frac{1}{p-2}\Big)^{\frac{1}{p}}, \quad p>2,
\enda
\end{align*}
for some positive constant $C>0$.
\end{lemma}
\hide
\begin{proof}
By the definition of $\bar{u}$ in Definition~\ref{D.Ra-E}, we have
\begin{align*}
\bega
|\nabla_x^{\mathrm{N}+1}\bar{u}| &\leq \frac{1}{|x|^{\mathrm{N}+2}}\bigg|\int_0^{|x|}r\bar{\w}(r)dr\bigg| + \sum_{0\leq k\leq\mathrm{N}}\frac{1}{|x|^k}|\nabla_x^{\mathrm{N}-k}\bar{\w}(|x|)| .
\enda
\end{align*}
From \eqref{rw-def}, we have
\begin{align*}
\bega
\|\nabla_x^{\mathrm{N}+1}\bar{u}\|_{L^2_x} &\leq \bigg\|\mathbf{1}_{|x|\leq1}\frac{1}{|x|^{\mathrm{N}+2}}\bigg|\int_0^{|x|}r\bar{\w}(r)dr\bigg\|_{L^2_x} dx +\bigg\|\mathbf{1}_{|x|>1}\frac{1}{|x|^{\mathrm{N}+2}}\bigg|\int_0^{|x|}r\bar{\w}(r)dr\bigg\|_{L^2_x} dx \cr 
&+ \sum_{0\leq k\leq\mathrm{N}}\bigg\|\frac{1}{|x|^k}|\nabla_x^{\mathrm{N}-k}\bar{\w}(|x|)|\bigg\|_{L^2_x} \cr 
&\leq C,
\enda
\end{align*}
where we used $\bar{\w}(|x|)\sim |x|^{\mathrm{N}}$ near $|x|=0$, $\int_{\R^2} \bar{\w}(|y|) dy \leq C$, and \eqref{rw-def}.
\end{proof}
\unhide

We now introduce the definitions of the energy and dissipation.

\hide 다시쓰기For the vortex patch and vortex sheet solution, we need to consider the infinite velocity energy case \eqref{caseECX}. For that we define some definitions for this case: 
\begin{align*}
\p_t F^{\e} + \frac{v}{\e}\cdot \nabla_x F^{\e} = \frac{1}{\e^2\kappa} \mathcal{N}(F^{\e},F^{\e}), \qquad F^{\e}|_{t=0}:= F^{\e}_0(x,v)
\end{align*}
We can decompose the initial data 
\begin{align*}
F^{\e}_0(x,v) = M^{\e}_0(x,v) + \AC{\P}F^{\e}_0(x,v), \qquad M^{\e}_0(x,v):= M_{[e^{\e \rho^{\e}_0(x)},\e u^{\e}_0(x),e^{\e \ta^{\e}_0(x)}]},
\end{align*}
where $M_{[\cdot,\cdot,\cdot]}$ is defined in \eqref{M-def}. \unhide

\hide
We can decompose $u_0 (x)$ via the Hodge decomposition \eqref{Hodge} together with the radial–energy decomposition in Definition~\ref{D.Ra-E} as 
\begin{align}\label{iniRa-E}
\bega
u_0(x)&= \widetilde{u}_0(x) + \bar{u}_0(x) +\mathbb{P}^{\perp}u_0(x), \ \   \int_{\R^2} \w_0(x) dx = \int_{\R^2} \nabla_x^{\perp}\cdot \bar{u}_0(x)dx,\ \ \int_{\R^2} \nabla_x^{\perp}\cdot \widetilde{u}_0(x) dx =0.
\enda
\end{align}
\unhide

\hide
By using the definition \ref{D.tiludef}, to unify the notation for both cases \eqref{caseEC} and \eqref{caseECX}, we can write 
\begin{align}\label{uexpand}
\bega
u^{\e}(t,x)= \widetilde{u}^{\e}(t,x)+\bar{u}(x)+\mathbb{P}^{\perp}u^{\e}(t,x).
\enda
\end{align}
We note that $\bar{u}(x)$ is non-zero only for the case \eqref{caseECX}: 
\begin{align}\label{barudef}
\bar{u}(x) := \begin{cases} 
0, \quad &\text{for the finite velocity energy case,} \quad \eqref{caseEC} \\
\frac{x^{\perp}}{|x|^2} \int_0^{|x|} r \bar{\w}(r) \, dr, \quad &\text{for the infinite velocity energy case,} \quad \eqref{caseECX}
\end{cases}
\end{align}
where $\bar{\w}(r)$ is defined in Definition \ref{D.Ra-E}.\unhide

\begin{definition}\label{D.al}
We define spatial multi-indices and their corresponding derivatives as
\begin{align}\label{mindexx}
\al_x = (\al_1, \cdots, \al_d), \qquad \p^{\al_x} = \p_{x_1}^{\al_1}\cdots\p_{x_d}^{\al_d}.
\end{align}
For scaled space-time derivatives, we define 
\begin{align}\label{mindex}
\al = (\al_0, \al_1, \cdots,\al_d), \qquad \p^{\al} = 
\p_{\tilde{t}}^{\al_0}\p_{x_1}^{\al_1}\cdots\p_{x_d}^{\al_d}, \quad \text{where} \quad \p_{\tilde{t}} = \e^{\mathfrak{n}}\p_t,
\end{align}
for arbitrary $\mathfrak{n} \geq 0$.
In this paper, when space–time derivatives are used, we restrict the temporal derivatives to at most first order, i.e., $\alpha_0 = 0,1$:
\begin{align}
&\textit{Purely spatial derivatives} \quad 
\p^{\al} = \p_{x_1}^{\al_1}\cdots\p_{x_d}^{\al_d}, \qquad \qquad d=2,3,  \label{caseA}
\\
&\textit{Space-time derivatives} 
\quad \hspace{5mm} \p^{\al} = (\e^{\mathfrak{n}}\p_t)^{\al_0}\p_{x_1}^{\al_1}\cdots\p_{x_d}^{\al_d}, \hspace{3mm} d=2,3. \label{caseB}
\end{align}

Let $\Omega=\R^d$ with $d=2$ or $3$. For $\mathrm{N}\geq0$, we define the lower-order energy and dissipation as follows:
\begin{align}\label{N-EDdef}
\bega
\mathcal{E}_{M}^{\mathrm{N}}(F^\e(t))=
\mathcal{E}_{M}^{\mathrm{N}}(t) &:= \bigg(\|\rho^{\e}(t)\|_{L^2_x}^2+\|u^{\e}(t)-\bar{u}\|_{L^2_x}^2 +\|\ta^{\e}(t)\|_{L^2_x}^2\bigg) \\
&+ \sum_{1\leq|\al|\leq\mathrm{N}} \bigg(\|\p^{\al}\rho^{\e}(t)\|_{L^2_x}^2+\|\p^{\al}u^{\e}(t)\|_{L^2_x}^2+\|\p^{\al}\ta^{\e}(t)\|_{L^2_x}^2\bigg), \cr
\mathcal{E}_G^{\mathrm{N}}(F^\e (t)) =\mathcal{E}_G^{\mathrm{N}}(t)&:= \sum_{0\leq|\al|\leq\mathrm{N}}\frac{1}{\e^2} \int_{\Omega\times\R^3}|\p^{\al}\AC{\P}F^{\e}|^2|M^{\e}|^{-1}dvdx, \cr 
\mathcal{D}_G^{\mathrm{N}}(F^\e(t))=
\mathcal{D}_G^{\mathrm{N}}(t)&:= \sum_{0\leq|\al|\leq\mathrm{N}} \frac{1}{\kappa\eps^4}\int_{\Omega \times \R^3} \nu\lw(\frac{v-\mathrm{U}^{\e}}{\sqrt{\mathrm{\Theta}^{\e}}}\rw) |\p^{\al}\AC{\P}F^{\e}|^2|M^{\e}|^{-1} dvdx .
\enda
\end{align}
\hide
For later convenience, we also define the energy and dissipation of the microscopic part at each derivative level, for $0\leq|\al|\leq \mathrm{N}$, as follows:
\begin{align}\label{N-EDdef3}
\bega
\mathcal{E}^{\al}_G(t)&:= \frac{1}{\e^2} \int_{\Omega\times\R^3}|\p^{\al}\AC{\P}F^{\e}|^2|M^{\e}|^{-1}dvdx, \cr 
\mathcal{D}^{\al}_G(t)&:=  \frac{1}{\kappa\eps^4}\int_{\Omega \times \R^3} \nu\lw(\frac{v-\mathrm{U}^{\e}}{\sqrt{\mathrm{\Theta}^{\e}}}\rw) |\p^{\al}\AC{\P}F^{\e}|^2|M^{\e}|^{-1} dvdx.
\enda
\end{align}
\unhide
For the top-order derivatives, we define
\begin{align}\label{N-EDdef2}
\bega
\mathcal{E}_{top}^{\mathrm{N}}(F^\e(t))=
\mathcal{E}_{top}^{\mathrm{N}}(t)&:=
\sum_{|\al|=\mathrm{N}+1}\kappa\lw(\|\p^{\al}\rho^{\e}\|_{L^2_x}^2+\|\p^{\al} u^{\e}\|_{L^2_x}^2+\|\p^{\al}\ta^{\e}\|_{L^2_x}^2\rw)  \cr 
&\quad + \sum_{|\al|=\mathrm{N}+1}\frac{\kappa}{\e^2} \int_{\Omega \times \R^3} | \AC{\P} \p^{\al} F^{\e} |^2 |M^{\e}|^{-1} dv dx, \cr 
\mathcal{D}_{top}^{\mathrm{N}}(F^\e (t)) = \mathcal{D}_{top}^{\mathrm{N}}(t)&:= \sum_{|\al|=\mathrm{N}+1} \frac{1}{\eps^4}\int_{\Omega \times \R^3} \nu\left( \frac{v-\mathrm{U}^{\e}}{\sqrt{\mathrm{\Theta}^{\e}}} \right) | \AC{\P} \p^{\al} F^{\e} |^2 |M^{\e}|^{-1} dv dx.
\enda
\end{align}
The total energy and total dissipation are defined by
\begin{align}\label{EDtotdef}
\bega
\mathcal{E}_{tot}^{\mathrm{N}}(F^\e(t)) =\mathcal{E}_{tot}^{\mathrm{N}}(t) &:= \mathcal{E}_M^{\mathrm{N}}(t) + \mathcal{E}_G^{\mathrm{N}}(t) + \mathcal{E}_{top}^{\mathrm{N}}(t), \\
\mathcal{D}_{tot}^{\mathrm{N}}(F^\e(t)) = \mathcal{D}_{tot}^{\mathrm{N}}(t) &:= \mathcal{D}_G^{\mathrm{N}}(t) + \mathcal{D}_{top}^{\mathrm{N}}(t).
\enda
\end{align}
\end{definition}

\subsection{A modulated entropy bound and the $\mathrm{N}$-admissible blow-up condition}

We close the energy--dissipation estimate under the following assumptions on the initial data.

\begin{definition}[Modulated entropy bound]
We say that the initial data satisfy the modulated entropy bound if
\begin{equation}\label{L2unif}
\sup_{\e>0}
\bigg(
\|(\rho_0^{\e}, u_0^{\e}-\bar{u}, \theta_0^{\e})\|_{L^2(\R^2)}
+ \left\|\frac{1}{\e}\frac{\AC{\P}F^{\e}_0}{\sqrt{\tilde{\mu}}}\right\|_{L^2_x(\R^2;L_v^2(\R^3))}
\bigg)
<+\infty.
\end{equation}
Here $\tilde{\mu}=M_{[1,0,1-c_0]}$ for some $0<c_0\ll1$, and $\bar{u}$
denotes the radial eddy defined in Definition~\ref{D.Ra-E} when $d=2$, while we set $\bar{u}=0$ for $d=3$.
\end{definition}

\begin{definition}[$\mathrm{N}$-admissible blow-up condition (ABC)]\label{ABC}
We say that the initial data satisfy the $\mathrm{N}$-admissible blow-up condition if the following bounds hold\footnote{
Villani stresses the importance of quantitative hydrodynamic limits and cautions that quantitative versions are likely to yield only extremely poor iterated-logarithmic bounds; see Section~6 of
\cite{VillaniBourbaki02} and footnote~\ref{Villani}. In particular, the restriction in the first line of
\eqref{ABC1} is consistent with the iterated-logarithmic barrier highlighted in \cite{VillaniBourbaki02}.
}:
\begin{align}\label{ABC1}
\bega
&\sum_{0 \leq |\alpha| \leq \mathrm{N}+1}\!\!\! \left(
\kappa^{(|\al|-\mathrm{N})_+}\| \p^\alpha (\rho^\e_0, u^\e_0, \theta^\e_0) \|_{L^2_x}^2
 + \frac{\kappa^{(|\al|-\mathrm{N})_+}}{\e^2}\left\| \frac{\p^\alpha \AC{\P}F^{\e}_0}{\sqrt{\tilde{\mu}}} \right\|_{L^2_{x,v}}^2
 \right)
\lesssim \log\!\Big(\log\!\big(\log(1/\e)\big)\Big),
\cr
&\sum_{0 \leq |\alpha| \leq \mathrm{N}+1}\!
\left\|
\p^\alpha \left( \frac{F_0^\e - \mu}{\sqrt{\tilde{\mu}}} \right)
 \right\|_{L^\infty_{x,v}} < +\infty,
\enda
\end{align}
where $\tilde{\mu}=M_{[1,0,1-c_0]}$ for some $0<c_0\ll 1$.
\end{definition}

\hide
\subsection{A modulated entropy bound and $\mathrm{N}$-admissible blow-up condition}
We close the energy--dissipation estimate under the following assumptions on the initial data. 

\begin{definition}[A modulated entropy bound]
We say that the initial data satisfy 
the modulated entropy bound if the following inequality holds: 
\begin{equation}\label{L2unif}
\sup_{\e>0}
\bigg(  \|( \rho_0^{\e}, u_0^{\e}-\bar{u}, \ta_0^{\e}
)\|_{L^2(\R^2)}
+ \left\|\frac{1}{\e}\frac{\AC{\P}F^{\e}_0}{\sqrt{\tilde{\mu}}}\right\|_{L^2_x(\R^2;L_v^2(\R^3))}
\bigg)
<+\infty.
\end{equation}
Here $\tilde{\mu}=M_{[1,0,1-c_0]}$ for some $0<c_0\ll1$, and $\bar{u}$
denotes the radial eddy defined in Definition~\ref{D.Ra-E} when $d=2$, while we set $\bar{u}=0$ for $d=3$.
\end{definition}

\begin{definition}[$\mathrm{N}$-Admissible Blow-up Condition (ABC)]\label{ABC}
We say that the initial data satisfy the \emph{Admissible Blow-up Condition (ABC)} if the following inequalities hold\footnote{
Villani stresses the importance of quantitative hydrodynamic limits and cautions that quantitative versions would likely yield only extremely poor iterated-logarithmic bounds; see Section~6 of
\cite{VillaniBourbaki02} and footnote~\ref{Villani}. In particular, the restriction in the first line of
\eqref{ABC1} is consistent with the iterated-logarithmic barrier highlighted in \cite{VillaniBourbaki02}.
}:
\begin{align}\label{ABC1}
\bega
&\sum_{0 \leq |\alpha| \leq \mathrm{N}+1}\!\!\! \left(
\kappa^{(|\al|-\mathrm{N})_+}\| \p^\alpha (\rho^\e_0, u^\e_0, \theta^\e_0) \|_{L^2_x}^2
 + \frac{\kappa^{(|\al|-\mathrm{N})_+}}{\e^2}\left\| \frac{\p^\alpha \AC{\P}F^{\e}_0}{\sqrt{\tilde{\mu}}} \right\|_{L^2_{x,v}}^2
 \right)
\lesssim \log\!\Big(\log\!\big(\log(1/\e)\big)\Big),
\cr
&\sum_{0 \leq |\alpha| \leq \mathrm{N}+1}\!
\left\|
\p^\alpha \left( \frac{F_0^\e - \mu}{\sqrt{\tilde{\mu}}} \right)
 \right\|_{L^\infty_{x,v}} <+ \infty
 , 
 \enda
\end{align}
where $\tilde{\mu}=M_{[1,0,1-c_0]}$ for some $0<c_0\ll 1$. 
\end{definition}

\unhide

\hide
\begin{remark}From the definition, 
    \begin{align}\label{Crutadef}
\bega
\mathcal{E}_{tot}^{\mathrm{N}}(0)&\les  \|(\rho^{\e}_0,u^{\e}_0-\bar{u},\ta^{\e}_0)\|_{L^2_x}^2 + \sum_{1\leq|\al|\leq \mathrm{N}+1}\kappa^{(|\al|-\mathrm{N})_+}\|\p^{\al}(\rho^{\e}_0,u^{\e}_0,\ta^{\e}_0)\|_{L^2_x}^2 \cr 
&+\sum_{0\leq|\al|\leq\mathrm{N}+1}\frac{\kappa^{(|\al|-\mathrm{N})_+}}{\e^2} \int_{\Omega\times\R^3}\frac{|\p^{\al}\AC{\P}F^{\e}(0)|^2}{|M^{\e}(0)|}dvdx. 
\enda
\end{align}
\end{remark}
\unhide

\section{Microscopic Energy estimates}\label{Sec.micro}

In this section, we present two main propositions.
The first concerns the energy estimate for the microscopic part
$\p^{\al}\AC{\P}F^{\e}$ for $0 \leq |\al|\leq \mathrm{N}$
(see Proposition~\ref{P.G.Energy} in Subsection~\ref{Sec.3.G}).
This estimate cannot be closed on its own, since the bound
involves top-order derivatives of order $\mathrm{N}+1$.
We denote this contribution by $\mathfrak{S}_G^{\al}(t)$,
which cancels with the corresponding leading-order term in the
macroscopic energy estimate.
The second proposition establishes the estimate for $\p^{\al}F^{\e}$
at the top-order derivative level $|\al|=\mathrm{N}+1$
with a loss of $\kappa$
(see Proposition~\ref{P.F.Energy} in Subsection~\ref{Sec.3.Top}).


%

\subsection{Purely microscopic estimate}\label{Sec.3.G}
We define two quantities to state the main proposition:
\begin{definition}
We define the high-order moments control:
\begin{align}\label{largev}
\bega
\mathcal{V}_{\ell}(F^{\e}(t,x)) &:= \sum_{0\leq|\al|\leq\mathrm{N}+1}\frac{\kappa^{(|\al|-\mathrm{N})_+}}{\eps^2} \int_{\R^3} \la v \ra^{\ell} | \p^{\al} \AC{\P}F^{\e} |^2 |M^{\e}|^{-1} dv  , 
\enda
\end{align}
where $\la v \ra^{\ell} := (1+|v|^2)^{\frac{\ell}{2}}$ and $x_+$ denotes $x_+=x$ if $x\geq0$ and $x_+=0$ if $x<0$.

We define a control of the turbulent transport as
\begin{align}\label{RSdef}
\bega
\TbT(t,x) := \sum_{i, j} |\p_{x_j} \mathbf{r}_{ij}^{\e}(t,x)| + \sum_{j} |\p_{x_j} \mathfrak{q}_j^{\e}(t,x)| + \sum_{i, j} |\p_{x_i} \mathrm{U}^{\e}_j \mathbf{r}_{ij}^{\e}(t,x)|, 
\enda
\end{align}
where $\mathbf{r}_{ij}^{\e}$ and $\mathfrak{q}_j^{\e}$ are defined in \eqref{albe-def}.
\end{definition}

To prove main propositions and theorems, we employ the following bootstrap assumption.
\begin{definition}[Weak Bootstrap Assumption]
We assume that the following bootstrap assumptions hold for some $T>0$ and for some sufficiently small constant $0<c_0 \ll 1 $:
\begin{align}\label{condition}
\bega
(\mathcal{B}_1):& \quad \e \kappa^{-\frac{1}{2}} \sup_{t \in [0,T]} (\mathcal{E}_{tot}^{\mathrm{N}}(t))^{\frac{1}{2}} \leq \frac{1}{4}, \cr 
(\mathcal{B}_2):& \quad \sup_{t \in [0,T]} \Big(|\mathrm{P}^{\e}(t,x) -1|,~|\mathrm{U}^{\e}(t,x)|,~|\mathrm{\Theta}^{\e}(t,x) -1|\Big) < \frac{c_0}{2}.
\enda
\end{align} 
\end{definition}

We present the main proposition of this subsection.

\begin{proposition}\label{P.G.Energy} 
Let $\Omega=\R^d$ with $d=2$ or $3$ and $\mathrm{N}>d/2+1$.
Assume that the initial data $u_0^{\e}$ satisfies either the finite velocity energy condition \eqref{caseEC} or the infinite velocity energy condition \eqref{caseECX}.
Under the bootstrap assumption \eqref{condition}, the following estimate holds in both the purely spatial derivative case \eqref{caseA} and the space-time derivative case \eqref{caseB} with arbitrary $\mathfrak{n}>0$.
\begin{align}\label{totalGt}
\bega
\frac{d}{dt}&\mathcal{E}_G^{\mathrm{N}}(F^{\e}(t)) +\frac{\sigma_L}{C} \mathcal{D}_G^{\mathrm{N}}(F^{\e}(t)) \leq  \sum_{0\leq|\al|\leq\mathrm{N}}\mathfrak{S}_G^{\al}(t) \cr 
&+ C_{}\bigg[ \int_{\Omega}|\nabla_xu^{\e}(t)| \mathcal{V}_2(F^{\e}(t)) dx + \int_{\Omega}|\nabla_x\ta^{\e}(t)| \mathcal{V}_3(F^{\e}(t)) dx + \frac{1}{\eps} \int_{\Omega} \TbT(t) \mathcal{V}_2(F^{\e}(t)) dx \bigg] \cr 
&+ C_{}\bigg[\e \kappa^{\frac{1}{2}} \mathcal{E}_{M}^{\mathrm{N}}(F^{\e}(t)) (\mathcal{D}_G^{\mathrm{N}}(F^{\e}(t)))^{\frac{1}{2}} +\e (\mathcal{E}_{M}^{\mathrm{N}}(F^{\e}(t)))^{\frac{1}{2}} \mathcal{D}_G^{\mathrm{N}}(F^{\e}(t)) \bigg] \cr
&+ C_{}\frac{1}{\kappa^{\frac{1}{2}}\e} \sum_{0\leq \al \leq \lfloor \mathrm{N}/2 \rfloor}\Big\|\|\la v \ra^{\frac{1}{2}} \p^{\al}\AC{\P}F^{\e} |M^{\e}|^{-\frac{1}{2}}(t)\|_{L^2_v}\Big\|_{L^\infty_x} (\mathcal{E}_{M}^{\mathrm{N}}(F^{\e}(t)))^{\frac{1}{2}} (\mathcal{D}_G^{\mathrm{N}}(F^{\e}(t)))^{\frac{1}{2}}   \cr 
&+ C_{}\sum_{0\leq|\al|\leq \lfloor\mathrm{N}/2\rfloor}\Big\|\|\la v \ra^{\frac{1}{2}} \p^{\al}\AC{\P}F^{\e} |M^{\e}|^{-\frac{1}{2}}(t)\|_{L^2_v}\Big\|_{L^\infty_x} \mathcal{D}_G^{\mathrm{N}}(F^{\e}(t)),
\enda
\end{align}
where $\sigma_L>0$ is the constant satisfying the coercivity estimate in Lemma~\ref{L.coer}, and $C>0$ is a universal constant.
Here, $\lfloor \cdot \rfloor$ denotes the greatest integer function.
The momentum-flux alignment contribution
$\mathfrak{S}_G^{\al}(t)$ is defined by
\begin{align}\label{Amidef}
\bega
\mathfrak{S}_G^{\al}(t) &:=-\frac{1}{\e^2}  \int_{\Omega} \frac{1}{k_B\mathrm{\Theta}^{\e}}\sum_{i,j}\p^{\al}\p_{x_i}u^{\e}_j \p^{\al} \mathbf{r}_{ij}^{\e} dx \cr 
&-\frac{1}{\e^2}  \int_{\Omega} \frac{1}{k_B\mathrm{\Theta}^{\e}}\sum_i \p^{\al}\p_{x_i}\ta^{\e} \bigg(\p^{\al}\mathfrak{q}_i^{\e}+\sum_{1\leq\beta\leq \al}\binom{\al}{\beta} \sum_{k} \p^{\beta}\mathrm{U}^{\e}_k \p^{\al-\beta}\mathbf{r}_{ik}^{\e}\bigg) dx.
\enda
\end{align}
\end{proposition}

In addition to the proposition, we will make essential use of the following lemma when closing the estimate. Lemma \ref{L.ABG} specifically addresses the contribution of $\TbT(t)$ appearing in \eqref{totalGt}.

\begin{lemma}\label{L.ABG} 
Under the bootstrap assumption \eqref{condition}, with $\mathrm{N}>d/2+1$, the following estimates hold for both the purely spatial derivatives \eqref{caseA} and the space-time derivatives \eqref{caseB}.
\begin{align}\label{ABGscale}
\bega
&(1)~ \|\p^{\al}\mathbf{r}_{ij}^{\e}(t)\|_{L^2_x}+\|\p^{\al}\mathfrak{q}_j^{\e}(t)\|_{L^2_x} + \sum_{0\leq \beta<\al} \|\p^{\al-\beta}\mathrm{U}^{\e}_j \p^{\beta}\mathbf{r}_{ij}^{\e}(t)\|_{L^2_x} \cr 
& \hspace{3cm} \leq \begin{cases} C \e^2\kappa^{\frac{1}{2}}(\mathcal{D}_G^{\mathrm{N}}(F^{\e}(t)))^{\frac{1}{2}}, \quad &\mbox{when} \quad 0\leq|\al|\leq\mathrm{N}, \\ 
C\e^2\big((\mathcal{D}_{tot}^{\mathrm{N}}(F^{\e}(t)))^{\frac{1}{2}}+\mathcal{E}_M^{\mathrm{N}}(F^{\e}(t))(t)\big), \quad &\mbox{when} \quad |\al|=\mathrm{N}+1.
\end{cases} \cr 
&(2) ~~ \|\TbT(t)\|_{L^\infty_x} 
\leq C\e^2\kappa^{\frac{1}{2}}(\mathcal{D}_G^{\mathrm{N}}(F^{\e}(t)))^{\frac{1}{2}}.
\enda
\end{align} 
Here, $\mathbf{r}_{ij}^{\e}$, $\mathfrak{q}_j^{\e}$, and $\TbT(t)$ are defined in \eqref{albe-def} and \eqref{RSdef}.
\end{lemma}
\begin{proof}
See the proof in Section \ref{Sec.Gembed}.
\end{proof}

\subsubsection{Linear coercivity}

\begin{definition}
For a given $M^{\e}$, we define the quasi-linear operator $\mathcal{L}$ by
\Be \label{Ldef}
\mathcal{L}(G)(v)=-2 \mathcal{N}(M^{\e},G)(v) . 
\Ee
\end{definition}
We state some standard properties of this operator. 
For the reader’s convenience, the detailed proofs are presented in Appendix \ref{A.Coll}.

\begin{lemma}\label{HH}
If $\int_{\R^3} (1,v,|v|^2)G dv =0$, then for both types of derivatives $\p^{\al}$ defined in \eqref{caseA} and \eqref{caseB}, we have $\int_{\R^3} (1,v,|v|^2) \p^{\al}G dv =0$.
\end{lemma}
\begin{proof}
It is a direct observation from 
$\int_{\R^3} (1,v,|v|^2) \p^{\al} G dv =  \p^{\al} \left(\int_{\R^3} (1,v,|v|^2) G dv\right)=0.$
\end{proof}

\begin{lemma}\label{Llem1}  
The operator $\mathcal{L}$ defined in \eqref{Ldef} satisfies the following properties:
\begin{enumerate}
\item $\mathcal{L}$ is symmetric: $\la \mathcal{L}(F), G  \frac{1}{M^{\e}} \ra  = \la F, \mathcal{L}(G)\frac{1}{M^{\e}}\ra $. 
\item On the functional subspace satisfying $
\int_{\R^3} \nu^{\e}(v)|F(v)|^2 \frac{1}{M^{\e}} dv < \infty$, where $\nu^{\e}$ is defined in \eqref{nudef}, the null space and its orthogonal complement are given by 
\begin{align*}
\bega
&\mathrm{Ker}(\mathcal{L})  :=  \{\mathcal{L}F=0  \}= span \left\{ M^{\e},vM^{\e},|v|^2M^{\e} \right\} , \cr 
&(\mathrm{Ker}(\mathcal{L}))^\perp := \left\{F~\bigg|~ \int_{\R^3} F G \frac{1}{M^{\e}} dv =0, \quad \mbox{for} \quad \forall G\in \mathrm{Ker}(\mathcal{L}) \right\} = \left\{F ~|~  \P F= 0 \right\},
\enda
\end{align*}
where projections $\P$ and $\AC{\P}$ are defined in \eqref{Pdef}. Moreover, the operator satisfies $\P \mathcal{L}(F) = 0$.

\item 
For any $G \in (\mathrm{Ker}(\mathcal{L}))^\perp$, there exists a unique $F \in (\mathrm{Ker}(\mathcal{L}))^\perp$ such that $\mathcal{L} F = G$. We denote this inverse by $F = \mathcal{L}^{-1} G$, and interpret the inverse operator as $\mathcal{L}^{-1}: (\mathrm{Ker}(\mathcal{L}))^\perp \to (\mathrm{Ker}(\mathcal{L}))^\perp$.

\item 
If $|\mathrm{P}^{\e} - 1|, |\mathrm{U}^{\e}|, |\mathrm{\Theta}^{\e} - 1| < 1/2$, then for any $n \geq 0$, the operators $\mathcal{L}$ and $\mathcal{L}^{-1}$ satisfy the following estimates:
\begin{align}\label{L-1L2}
\bega
&\int_{\R^3} (1+|v|^{n}) |\mathcal{L}(F)|^2 |M^{\e}|^{-1} dv \lesssim 
\int_{\R^3}(1+|v|^{n+2})|F|^2 |M^{\e}|^{-1} dv, \cr 
&\int_{\R^3}(1+|v|^{n+2})|\mathcal{L}^{-1}(F)|^2 |M^{\e}|^{-1} dv \lesssim \int_{\R^3} (1+|v|^{n}) |F|^2 |M^{\e}|^{-1} dv.
\enda
\end{align}
\end{enumerate}
\end{lemma}
\begin{proof}
    See the proof in Appendix \ref{A.Coll}.
\end{proof}

\begin{lemma}\label{L.coer}
For any given functions $G, H, K$ such that the right-hand side is bounded, the quasi-linear operator $\mathcal{L}$, defined in \eqref{Ldef}, satisfies
\begin{align}\label{coer1}
\la \mathcal{L}G, G|M^{\e}|^{-1} \ra_{L^2_v} 
&\geq \sigma_L \mathrm{P}^{\e}\sqrt{\mathrm{\Theta}^{\e}} \int_{\R^3} \nu\lw(\frac{v-\mathrm{U}^{\e}}{\sqrt{\mathrm{\Theta}^{\e}}}\rw) |\AC{\P} G|^2 |M^{\e}|^{-1} \, dv,
\end{align}
for a positive constant $\sigma_L > 0$. Here, $\nu(\cdot)$ is defined in \eqref{nudef}. \\
The collision operator $\mathcal N(\cdot,\cdot)$, defined in \eqref{Qdef}, satisfies
\begin{align}\label{nonlin}
\bega
\int_{\R^3}\mathcal N(G,H)K |M^{\e}|^{-1} dv &\les  \sqrt{\mathrm{P}^{\e}\mathrm{\Theta}^{\e}} \bigg(\left\|\sqrt{\nu^{\e}}G|M^{\e}|^{-\frac{1}{2}}\right\|_{L^2_v}\left\|H|M^{\e}|^{-\frac{1}{2}}\right\|_{L^2_v} \cr 
& +\left\|G|M^{\e}|^{-\frac{1}{2}}\right\|_{L^2_v}\left\|\sqrt{\nu^{\e}}H|M^{\e}|^{-\frac{1}{2}}\right\|_{L^2_v}\bigg)\left\|\sqrt{\nu^{\e}}K|M^{\e}|^{-\frac{1}{2}}\right\|_{L^2_v}.
\enda
\end{align}
If we further assume $\sup_{t \in [0,T]} (|\mathrm{P}^{\e}(t,x) -1|,|\mathrm{U}^{\e}(t,x) |,|\mathrm{\Theta}^{\e}(t,x) -1|) <1$, then we have 
\begin{align}\label{nonlin2}
\bega
\int_{\R^3}\la v \ra^n |\mathcal N(G,H)|^2|M^{\e}|^{-1} dv &\les  \int_{\R^3}(1+|v|^{n+2})\frac{|G(v)|^2}{M^{\e}(v)}dv\int_{\R^3}\frac{|H(v)|^2}{M^{\e}(v)}dv \cr 
&+ \int_{\R^3}\frac{|G(v)|^2}{M^{\e}(v)}dv\int_{\R^3}(1+|v|^{n+2})\frac{|H(v)|^2}{M^{\e}(v)}dv,
\enda
\end{align}
for any $n \geq 0$. 
\end{lemma}
\begin{proof}
    See the proof in Appendix \ref{A.Coll}.
\end{proof}

\subsubsection{Proof of Proposition \ref{P.G.Energy}}

In the proof, for brevity, we slightly abuse notation by writing $\mathcal{E}(t)$, $\mathcal{D}(t)$, and $\mathcal{V}_{\ell}(t)$ for $\mathcal{E}^{\mathrm{N}}(F^{\e}(t))$, $\mathcal{D}^{\mathrm{N}}(F^{\e}(t))$, and $\mathcal{V}_{\ell}(F^{\e}(t))$, respectively.

\begin{lemma}\label{Gcomb} 
For the microscopic part $\AC{\P}F^{\e}$ satisfying the equation \eqref{Geqn0}, and for both types of derivatives $\p^{\al}$ defined in \eqref{caseA} and \eqref{caseB}, we have the following energy estimate:
\begin{equation}\label{ACE-1} 
\begin{split}
\frac{d}{dt}  \bigg(\sum_{0\leq|\al|\leq\mathrm{N}} \frac{1}{2\e^2}\int_{\Omega \times \R^3} \frac{|\p^{\al}\AC{\P}F^{\e}|^2}{M^{\e}}&dvdx\bigg) + \sum_{0\leq|\al|\leq\mathrm{N}}\frac{1}{\kappa\e^4}\int_{\Omega \times \R^3}\mathcal{L}(\p^{\al}\AC{\P}F^{\e})\frac{\p^{\al}\AC{\P}F^{\e}}{M^{\e}} dvdx  \cr
&= \AC{I}_1(t)+\AC{I}_2(t)+\AC{I}_3(t)+\AC{I}_4(t)+\AC{I}_5(t)+\AC{I}_6(t)+\AC{I}_7(t),
\end{split}
\end{equation}
where
\begin{align}\label{ACEdef}
\bega
\AC{I}_1(t)&= -\sum_{0\leq|\al|\leq\mathrm{N}}\frac{1}{2\e^2}\int_{\Omega \times \R^3} \bigg(\p_t M^{\e} + \frac{v}{\eps}\cdot\nabla_x M^{\e}\bigg) \frac{|\p^{\al}\AC{\P}F^{\e}|^2}{|M^{\e}|^2} dvdx,  \cr 
\AC{I}_2(t)&=-\sum_{0\leq|\al|\leq\mathrm{N}}\frac{1}{\e^2} \int_{\Omega \times \R^3}\p^{\al}\bigg[\frac{1}{k_B\mathrm{\Theta}^{\e}}\sum_{i,j}\p_{x_i}u^{\e}_j \mathfrak{R}^{\e}_{ij}M^{\e}\bigg] \frac{\p^{\al}\AC{\P}F^{\e}}{M^{\e}} dvdx, \cr 
\AC{I}_3(t)&=-\sum_{0\leq|\al|\leq\mathrm{N}}\frac{1}{\e^2} \int_{\Omega \times \R^3}\p^{\al}\bigg[\frac{1}{k_B\mathrm{\Theta}^{\e}}\sum_i \p_{x_i}\ta^{\e} \mathcal{Q}^{\e}_iM^{\e}\bigg]\frac{\p^{\al}\AC{\P}F^{\e}}{M^{\e}} dvdx, \cr 
\AC{I}_4(t)&=  \sum_{0\leq|\al|\leq\mathrm{N}}\frac{1}{\kappa\e^4}\int_{\Omega \times \R^3}\sum_{1\leq \beta \leq \al}\binom{\al}{\beta}\mathcal N(\p^{\beta}M^{\e},\p^{\al-\beta}\AC{\P}F^{\e}) \frac{\p^{\al}\AC{\P}F^{\e}}{M^{\e}} dvdx,  \cr
\AC{I}_5(t)&=  \sum_{0\leq|\al|\leq\mathrm{N}}\frac{1}{\kappa\e^4} \int_{\Omega \times \R^3} \sum_{0\leq\beta\leq\al}\binom{\al}{\beta}\mathcal N(\p^{\beta}\AC{\P}F^{\e},\p^{\al-\beta}\AC{\P}F^{\e}) \frac{\p^{\al}\AC{\P}F^{\e}}{M^{\e}} dvdx, \cr 
\AC{I}_6(t)&= \sum_{0\leq|\al|\leq\mathrm{N}}\frac{1}{\e^3} \int_{\Omega \times \R^3} \p^{\al}\bigg[\frac{1}{ k_B\mathrm{P}^{\e}\mathrm{\Theta}^{\e}}\sum_{i,j} \p_{x_j} \mathbf{r}_{ij}^{\e} (v_i-\mathrm{U}^{\e}_i)M^{\e}\bigg] \frac{\p^{\al}\AC{\P}F^{\e}}{M^{\e}} dvdx,  \cr 
\AC{I}_7(t)&= \!\!\!\!\!\sum_{0\leq|\al|\leq\mathrm{N}}\! \frac{1}{\e^3} \int_{\Omega \times \R^3} \!\! \p^{\al}\bigg[\frac{2}{3k_B\mathrm{P}^{\e}\mathrm{\Theta}^{\e}} \bigg(\!\sum_{j} \p_{x_j} \mathfrak{q}_j^{\e} + \sum_{i,j} \p_{x_i}\mathrm{U}^{\e}_j \mathbf{r}_{ij}^{\e} \bigg) \!\lw( \frac{|v-\mathrm{U}^{\e}|^2}{2k_B\mathrm{\Theta}^{\e}}-\frac{3}{2}\! \rw)\!M^{\e}\bigg]   \frac{\p^{\al}\AC{\P}F^{\e}}{M^{\e}} dvdx.
\enda
\end{align}
\end{lemma}
\begin{proof}  
We apply $\p^{\al}$  to the equation  \eqref{Geqn0} and distribute the derivatives over the identity $\mathcal{L}(\AC{\P}F^{\e}) = -2\mathcal{N}(M^{\e}, \AC{\P}F^{\e})$ to obtain
\begin{align*}
\bega
\p_t &\p^{\al}\AC{\P}F^{\e}+\frac{v}{\e}\cdot\nabla_x \p^{\al}\AC{\P}F^{\e} + \frac{1}{\kappa\e^2}\mathcal{L}(\p^{\al}\AC{\P}F^{\e}) = \frac{2}{\kappa\e^2}\sum_{1\leq \beta \leq \al}\binom{\al}{\beta}\mathcal{N}(\p^{\beta}M^{\e},\p^{\al-\beta}\AC{\P}F^{\e}) \cr 
&- \p^{\al}\bigg(\frac{1}{k_B\mathrm{\Theta}^{\e}}\sum_{i,j}\p_{x_i}u^{\e}_j \mathfrak{R}^{\e}_{ij}M^{\e}\bigg)-\p^{\al}\bigg(\frac{1}{k_B\mathrm{\Theta}^{\e}}\sum_i \p_{x_i}\ta^{\e} \mathcal{Q}^{\e}_iM^{\e}\bigg) \cr 
&+\frac{1}{\kappa\e^2}\sum_{0\leq\beta\leq\al}\binom{\al}{\beta}\mathcal{N}(\p^{\beta}\AC{\P}F^{\e},\p^{\al-\beta}\AC{\P}F^{\e}) + \frac{1}{\e}\p^{\al}\bigg(\frac{1}{k_B\mathrm{P}^{\e}\mathrm{\Theta}^{\e}}\sum_{i,j} \p_{x_j} \mathbf{r}_{ij}^{\e} (v_i-\mathrm{U}^{\e}_i)M^{\e}\bigg) \cr 
&+\frac{1}{\e}\p^{\al}\bigg(\frac{2}{3k_B\mathrm{P}^{\e}\mathrm{\Theta}^{\e}} \bigg(\sum_{j} \p_{x_j} \mathfrak{q}_j^{\e} + \sum_{i,j} \p_{x_i}\mathrm{U}^{\e}_j \mathbf{r}_{ij}^{\e} \bigg) \lw( \frac{|v-\mathrm{U}^{\e}|^2}{2k_B\mathrm{\Theta}^{\e}}-\frac{3}{2} \rw)M^{\e}\bigg).
\enda
\end{align*}
We then test the above equation against  $\frac{1}{\e^2}\p^{\al}\AC{\P}F^{\e}|M^{\e}|^{-1}$ and sum over all multi-indices with $0 \leq |\al| \leq \mathrm{N}$ to obtain the desired result.
\end{proof}

To proceed with the energy estimate in Lemma \ref{Gcomb}, we expand $\p^{\al}M^{\e}$ using an order expansion.
\begin{lemma}\label{Mal}
For the two types of derivatives $\p^{\al}$ defined in \eqref{caseA} and \eqref{caseB}, and for $|\al| \geq 1$, 
\begin{align}
\p^{\al} M^{\e} = \e\Phi_{\al}^1M^{\e} +\sum_{2\leq i\leq |\al|}\eps^i \Phi_{\al}^i M^{\e}.
\label{RpM-def}
\end{align}
In particular, $\Phi_{\al}^1$ and $\Phi_{\al}^2$ have the following explicit forms:
\beq\label{Phi1def}
\bega
\Phi_\al^1 = \p^{\al} \rho^{\e}+\frac{\p^{\al} u^{\e} \cdot (v-\mathrm{U}^{\e})}{k_B\mathrm{\Theta}^{\e}}+ \p^{\al} \ta^{\e} \lw( \frac{|v-\mathrm{U}^{\e}|^2}{2k_B\mathrm{\Theta}^{\e}}-\frac{3}{2} \rw),
\enda
\eeq
and
\begin{align}\label{Phi2def}
\bega
\Phi_{\al}^2&= \sum_{0<\beta<\al}\binom{\al}{\beta}\Phi_{\beta}^1\Phi_{\al-\beta}^1 -\sum_{0<\beta<\al}\binom{\al}{\beta}\bigg(\frac{\p^{\al-\beta}u^{\e}\cdot\p^{\beta}u^{\e}}{k_B \mathrm{\Theta}^{\e}} + \p^{\al-\beta}\ta^{\e}\frac{\p^{\beta}u^{\e}\cdot(v - \mathrm{U}^{\e})}{k_B \mathrm{\Theta}^{\e}} \cr 
& \quad + \p^{\beta}\ta^{\e}\frac{\p^{\al-\beta}u^{\e}\cdot(v - \mathrm{U}^{\e})}{k_B \mathrm{\Theta}^{\e}} \bigg) + \p^{\al-\beta}\ta^{\e}\p^{\beta}\ta^{\e}\frac{|v - \mathrm{U}^{\e}|^2}{2k_B \mathrm{\Theta}^{\e}} .
\enda
\end{align}
For $t \in [0,T]$ satisfying the bootstrap assumption \eqref{condition}, and for $1 \leq i \leq |\al|$, we have
\begin{align}
\e^i|\Phi_{\al}^i|M^{\e} &\les \e^i(1+|v|^{2i})  \sum_{\substack{\al_1+\cdots+\al_i=\al \\ \al_i>0}}|\p^{\al_1}(\rho^{\e},u^{\e},\ta^{\e})|\times \cdots \times |\p^{\al_i}(\rho^{\e},u^{\e},\ta^{\e})| M^{\e}, \label{Phiinf0} \\ 
|\p^{\al}M^{\e}| &\les \sum_{1\leq i\leq |\al|}\eps^i (1+|v|^{2i})\sum_{\substack{\al_1+\cdots+\al_i=\al \\ \al_i>0}}|\p^{\al_1}(\rho^{\e},u^{\e},\ta^{\e})|\times \cdots \times |\p^{\al_i}(\rho^{\e},u^{\e},\ta^{\e})| M^{\e}. \label{al-M}
\end{align}
For any $n \geq 0$, $0\leq n_0 < 1$, and $1 \leq p \leq \infty$, the following estimates hold:
\begin{align}\label{Phiscale}
\bega
\|\la v\ra^n |\p^{\al}M^{\e}||M^{\e}|^{-n_0}(t)\|_{L^2_xL^p_v} &\les \e\mathcal{E}_{M}^{\frac{1}{2}}(t), \quad &\mbox{for} \quad 1\leq|\al|\leq\mathrm{N}, \cr
\|\la v\ra^n |\p^{\al}M^{\e}||M^{\e}|^{-n_0}(t)\|_{L^2_xL^p_v} &\les \e\kappa^{-\frac{1}{2}}\mathcal{E}_{tot}^{\frac{1}{2}}(t), \quad &\mbox{for} \quad |\al|=\mathrm{N}+1, 
\enda
\end{align}
and
\begin{align}\label{Phiinfscale}
\bega
\|\la v\ra^n |\p^{\al}M^{\e}||M^{\e}|^{-n_0}(t)\|_{L^\infty_xL^p_v} &\les \e\mathcal{E}_{M}^{\frac{1}{2}}(t), \qquad \hspace{3mm} \mbox{for} \quad 1\leq|\al|\leq\mathrm{N}-2, \cr
\|\la v\ra^n |\p^{\al}M^{\e}||M^{\e}|^{-n_0}(t)\|_{L^\infty_xL^p_v} &\les \e\kappa^{-\frac{1}{2}}\mathcal{E}_{tot}^{\frac{1}{2}}(t), \quad \mbox{for} \quad |\al|=\mathrm{N}-1.
\enda
\end{align}
The estimates \eqref{Phiscale} and \eqref{Phiinfscale} also hold when replacing $\p^{\al} M^{\e}$ with $\e \Phi_{\al}^1 M^{\e}$. 

Furthermore, the remainder $\sum_{2\leq i\leq |\al|}\eps^i \Phi_{\al}^i M^{\e}$ in \eqref{RpM-def} satisfies
\begin{align}
&\Big\|\la v\ra^n   \sum_{2\leq i\leq |\al|}\eps^i |\Phi_{\al}^i M^{\e} ||M^{\e}|^{-n_0}(t)\Big\|_{L^2_xL^p_v} \les \e^2\mathcal{E}_M(t), \quad \mbox{for} \quad  |\al|\leq\mathrm{N}+1, \label{Rscale} \\
&\Big\|\la v\ra^n   \sum_{2\leq i\leq |\al|}\eps^i |\Phi_{\al}^i M^{\e} ||M^{\e}|^{-n_0}(t) \Big\|_{L^\infty_xL^p_v} \les \e^2\mathcal{E}_M(t), \quad \mbox{for} \quad |\al|\leq\mathrm{N}-1. \label{Rinfscale}
\end{align}

\hide
For $\Phi_\al^i$, $i=3,4,5$, we denote the form of the terms  
\begin{align}\label{Phi3def}
\Phi_\al^i = \sum_{\substack{\al_1+\cdots+\al_i=\al\\0<|\al_i|<|\al| \\ x_1,\cdots,x_i\in\{\rho^{\e},u^{\e},\ta^{\e}\}}} \frac{\mathcal{P}_{\al_1,\cdots,\al_i}^i(\mathrm{P}^{\e},\mathrm{\Theta}^{\e},v-\mathrm{U}^{\e})}{\mathcal{R}_{\al_1,\cdots,\al_i}^i(\mathrm{P}^{\e},\mathrm{\Theta}^{\e})}  |\p^{\al_1}x_1| \times \cdots \times |\p^{\al_i}x_i|
\end{align}
Here the polynomial $\mathcal{P}_{\al_1,\cdots,\al_k}^i$ and monomial $\mathcal{R}_{\al_1,\cdots,\al_k}^i$ satisfy the following condition:
\begin{itemize}
\item Once we write  $(\mathrm{P}^{\e},\mathrm{\Theta}^{\e},v-\mathrm{U}^{\e})=(x_1,\cdots,x_5)$, then the generic polynomial $\mathcal{P}_{\al_1,\cdots,\al_i}^i$ satisfies  $\mathcal{P}_{\al_1,\cdots,\al_i}^i(x_1,\cdots,x_5)= \sum_{n}C_n x_1^{n_1}\cdots x_5^{n_5}$. And the maximum velocity growth is $|v-\mathrm{U}^{\e}|^{2i}$ with $n_1,\cdots,n_5 \in \mathbb{Z}_{\geq0}$. 
\item $\mathcal{R}_{\al_1,\cdots,\al_i}^i$ is generic monomial with $\mathcal{R}_{\al_1,\cdots,\al_i}^i(\mathrm{P}^{\e},\mathrm{\Theta}^{\e})=C_{\al_1,\cdots,\al_i}|\mathrm{P}^{\e}|^{n_1}|\mathrm{\Theta}^{\e}|^{n_2}$ such that $n_1,n_2 \in \mathbb{Z}_{\geq0}$ depending on $(\al_1,\cdots,\al_i)$ with $C_{\al_1,\cdots,\al_i}\in\mathbb{R}/0$.
\end{itemize}
\unhide
\end{lemma}

\hide
\begin{proof}
To decompose $\p^{\al}M^{\e}$ to $\e, \cdots,\e^5$ terms, we note that if one derivative acts on $(\mathrm{P}^{\e},\mathrm{U}^{\e},\mathrm{\Theta}^{\e})$, then it makes $\e\p^1(\rho^{\e},u^{\e},\ta^{\e})$. Otherwise, if one derivative acts on $(\rho^{\e},u^{\e},\ta^{\e})$, then just derivative increases $\p^1(\rho^{\e},u^{\e},\ta^{\e})$ without additional scale. \\
(i) Computation of $\e\Phi_{\al}^1$:
The $\e\Phi_{\al}^1$ part of $\p^{\al}M^{\e}$ is the term in which the derivative acts $(\mathrm{P}^{\e},\mathrm{U}^{\e},\mathrm{\Theta}^{\e})$ only once, and the remaining derivatives act $(\rho^{\e},u^{\e},\ta^{\e})$. Once we take $\p^{\beta}$ on $M^{\e}$, for $|\beta|=1$, then we get 
\beq\label{M1deri}
\bega
\p^{\beta}M^{\e} = \eps \lw(\frac{\p^{\beta} \rho^{\e}}{\mathrm{P}^{\e}}+\frac{\p^{\beta} u^{\e} \cdot (v-\mathrm{U}^{\e})}{k_B\mathrm{\Theta}^{\e}}+\frac{\p^{\beta}\ta^{\e}}{\mathrm{\Theta}^{\e}} \lw( \frac{|v-\mathrm{U}^{\e}|^2}{2k_B\mathrm{\Theta}^{\e}}-\frac{3}{2} \rw)\rw) M^{\e} = \eps \Phi_{\beta}^1 M^{\e}
\enda
\eeq
Taking the remaining derivatives $\p^{\al-\beta}$ on \eqref{M1deri}, applying the Leibniz rule gives
\beq\label{M2deri0}
\bega
\p^{\al} M^{\e} &=  \sum_{|\beta|=1,\beta\leq\al}\p^{\al-\beta}(\e\Phi_{\beta}^1 M^{\e} ) = \eps \sum_{|\beta|=1,\beta\leq\al}\sum_{0\leq \gamma \leq \al-\beta}\binom{\al}{\beta}\binom{\al-\beta}{\gamma}(\p^{\gamma}\Phi_{\beta}^1 \p^{\al-\beta-\gamma} M^{\e} ) 
\enda
\eeq
From \eqref{M2deri0}, the $O(\e)$ part of $\p^{\al}M^{\e}$ arise when the every derivative acts on $(\rho^{\e},u^{\e},\ta^{\e})$ on $\Phi_{\beta}^1$ with $\gamma=\al-\beta$, that is 
\beq\label{M5de}
\bega
O(\e) ~\mbox{of}~ \p^{\al} M^{\e} &= \eps \Phi_{\al}^1 M^{\e}  
\enda
\eeq
Similarly, $O(\e^2)$ part of $\p^{\al}M^{\e}$ is the term in which the derivative acts $(\mathrm{P}^{\e},\mathrm{U}^{\e},\mathrm{\Theta}^{\e})$ twice, and the remaining derivatives act $(\rho^{\e},u^{\e},\ta^{\e})$. That is, 
\beq\label{M5de2}
\bega
O(\e^2) ~\mbox{of}~ \p^{\al} M^{\e} &= \eps \sum_{|\beta|=1,\beta\leq\al}\sum_{0\leq \gamma \leq \al-\beta}\binom{\al}{\beta}\binom{\al-\beta}{\gamma} \bigg[ [ O(1) ~\mbox{of}~ \p^{\gamma}\Phi_{\beta}^1 ] \times [O(\e) ~\mbox{of}~ \p^{\al-\beta-\gamma}M^{\e}] \cr 
&\quad + [O(1) ~\mbox{of}~ \p^{\al-\beta-\gamma}M^{\e}] \times [ O(\e) ~\mbox{of}~ \p^{\gamma}\Phi_{\beta}^1 ] \bigg] 
\enda
\eeq
By the same way as in \eqref{M2deri0} and \eqref{M5de}, we have
\beq \label{O1e}
\bega
& [ O(1) ~\mbox{of}~ \p^{\gamma}\Phi_{\beta}^1 ] = \Phi_{\beta+\gamma}^1 \qquad \hspace{19mm} [O(\e) ~\mbox{of}~ \p^{\al-\beta-\gamma}M^{\e}] = \e\Phi_{\al-\beta-\gamma}^1M^{\e} \cr 
& [O(1) ~\mbox{of}~ \p^{\al-\beta-\gamma}M^{\e}] = \mathbf{1}_{\gamma=\al-\beta}M^{\e} \qquad [ O(\e) ~\mbox{of}~ \p^{\al-\beta}\Phi_{\beta}^1 ] = \sum_{\substack{|\gamma|=1\\\gamma \leq \al-\beta}} \binom{\al-\beta}{\gamma}\p^{\gamma}(\mathrm{P}^{\e},\mathrm{U}^{\e},\mathrm{\Theta}^{\e})\cdot\nabla_{(\mathrm{P}^{\e},\mathrm{U}^{\e},\mathrm{\Theta}^{\e})}\Phi_{\al-\gamma}^1
\enda
\eeq
Substituting \eqref{O1e} for \eqref{M5de2} gives 
\beq\label{M5de22}
\bega
O(\e^2) ~\mbox{of}~ \p^{\al} M^{\e} &= \eps^2 \sum_{|\beta|=1,\beta\leq\al}\sum_{0\leq \gamma \leq \al-\beta}\binom{\al}{\beta}\binom{\al-\beta}{\gamma}  \Phi_{\beta+\gamma}^1 \Phi_{\al-\beta-\gamma}^1 M^{\e} + \eps^2 \sum_{|\gamma|=1} \binom{\al}{\beta} \binom{\al-\beta}{\gamma} \p^{\gamma}(\rho^{\e},u^{\e},\ta^{\e})\cdot\nabla_{(\mathrm{P}^{\e},\mathrm{U}^{\e},\mathrm{\Theta}^{\e})}\Phi_{\al-\gamma}^1 M^{\e} 
\enda
\eeq
We can compute the last term of \eqref{M5de22} as follows: 
\beq \label{Phideri}
\bega
\sum_{|\gamma|=1} \p^{\gamma}(\rho^{\e},u^{\e},\ta^{\e})\cdot & \nabla_{(\mathrm{P}^{\e},\mathrm{U}^{\e},\mathrm{\Theta}^{\e})}\Phi_{\al-\gamma}^1 = \sum_{|\gamma|=1} \bigg(-\frac{\p^{\gamma}\rho^{\e}\p^{\al-\gamma}\rho^{\e}}{|\mathrm{P}^{\e}|^2}-\frac{\p^{\gamma}u^{\e}\cdot\p^{\al-\gamma}u^{\e}}{k_B\mathrm{\Theta}^{\e}}-\frac{\p^{\gamma}u^{\e}\p^{\al-\gamma}\ta^{\e} \cdot (v-\mathrm{U}^{\e})}{k_B|\mathrm{\Theta}^{\e}|^2} \cr
&\quad -\frac{\p^{\gamma}\ta^{\e}\p^{\al-\gamma}\ta^{\e}}{|\mathrm{\Theta}^{\e}|^2}\lw( \frac{|v-\mathrm{U}^{\e}|^2}{2k_B\mathrm{\Theta}^{\e}}-\frac{3}{2} \rw) - \frac{\p^{\gamma}\ta^{\e}}{\mathrm{\Theta}^{\e}}\lw(\frac{2(v-\mathrm{U}^{\e})\cdot\p^{\al-\gamma} u^{\e}}{2k_B\mathrm{\Theta}^{\e}} +\frac{\p^{\al-\gamma}\ta^{\e} |v-\mathrm{U}^{\e}|^2}{2k_B|\mathrm{\Theta}^{\e}|^2}\rw) \bigg)
\enda
\eeq
Combining \eqref{M5de22} and \eqref{Phideri}, we obtain 
\beq\label{M5derie2'}
\bega
O(\e^2) ~\mbox{of}~ \p^{\al} M^{\e} &=\e^2 \sum_{0<\beta<\al}\Phi_{\beta}^1\Phi_{\al-\beta}^1M^{\e} -\e^2\sum_{|\beta|=1,\beta\leq\al}\bigg[\frac{\p^{\beta}\rho^{\e}\p^{\al-\beta}\rho^{\e}}{|\mathrm{P}^{\e}|^2}+\frac{\p^{\beta}u^{\e}\cdot\p^{\al-\beta}u^{\e}}{k_B\mathrm{\Theta}^{\e}}+\frac{\p^{\beta}u^{\e}\p^{\al-\beta}\ta^{\e}\cdot(v-\mathrm{U}^{\e})}{k_B|\mathrm{\Theta}^{\e}|^2} \cr
& +\frac{\p^{\beta}\ta^{\e}\p^{\al-\beta}\ta^{\e}}{|\mathrm{\Theta}^{\e}|^2}\lw( \frac{|v-\mathrm{U}^{\e}|^2}{2k_B\mathrm{\Theta}^{\e}}-\frac{3}{2} \rw) + \frac{\p^{\beta}\ta^{\e}}{\mathrm{\Theta}^{\e}}\lw(\frac{2(v-\mathrm{U}^{\e})\cdot\p^{\al-\beta} u^{\e}}{2k_B\mathrm{\Theta}^{\e}} +\frac{\p^{\al-\beta}\ta^{\e} |v-\mathrm{U}^{\e}|^2}{2k_B|\mathrm{\Theta}^{\e}|^2}\rw) \bigg] M^{\e}
\enda
\eeq
This corresponds to the definition of $\Phi_{\al}^2$ in \eqref{Phi2def}. Considering that the scale $\e$ occurs as many times as the differential acts on $(\mathrm{P}^{\e},\mathrm{U}^{\e},\mathrm{\Theta}^{\e})$, by the same way as in \eqref{M2deri0}, we can have the form of $\Phi_{\al}^i$ for $i=3,4,5$ inductively. The term $\Phi_{\al}^i$ has the following property: The term $\Phi_{\al}^i$ has the form of multiplying $\p(\rho^{\e},u^{\e},\ta^{\e})$ $i$-times. And the maximum velocity growth of $\Phi_{\al}^i$ is $|v-\mathrm{U}^{\e}|^{2i}$ and which appear in the maximum number of product of $\p\ta^{\e}$ term.
\end{proof}
\unhide

\begin{proof}
We prove this in Appendix \ref{A.B}.
\end{proof}

\begin{lemma}\label{lemscale}
For the two types of derivatives $\p^{\al}$ defined in \eqref{caseA} and \eqref{caseB}, and for $d = 2, 3$, we have
\begin{align}
&(1)~\|\p^{\al}(\rho^{\e},u^{\e},\ta^{\e})(t)\|_{L^2_x} \leq \begin{cases} \mathcal{E}_{M}^{\frac{1}{2}}(t), \quad &\mbox{for} \quad 1\leq |\al| \leq \mathrm{N}, \\ \kappa^{-\frac{1}{2}}\mathcal{E}_{top}^{\frac{1}{2}}(t), \quad &\mbox{for} \quad |\al| =\mathrm{N}+1, \end{cases} 
\label{rutscale}\\
&(2)~ \|\p^{\al}(\rho^{\e},u^{\e},\ta^{\e})(t)\|_{L^\infty_x}) \leq \begin{cases} C \mathcal{E}_{M}^{\frac{1}{2}}(t),  \quad  &\mbox{for}  \quad  1\leq |\al| \leq \mathrm{N}-2, \\ 
C \kappa^{-\frac{1}{4}}\mathcal{E}_{tot}^{\frac{1}{2}}(t),  \quad  &\mbox{for}  \quad  |\al| = \mathrm{N}-1,
\end{cases} \label{rutinf}
\end{align}
for a constant $C > 0$ determined by the Agmon inequality \eqref{Agmon}. Here, $\mathcal{E}_M$, $\mathcal{E}_{top}$, and $\mathcal{E}_{tot}$ are defined in Definition \ref{D.al}.
\end{lemma}
\begin{proof}
(1) This estimate follows directly from the definition of the energy $\mathcal{E}_M$ given in \eqref{N-EDdef}. \\
(2) For any $\zeta \in \{\rho^{\e}, u^{\e}, \ta^{\e}\}$, we apply Agmon's inequality \eqref{Agmon} for $d = 3$ and $|\al| = \mathrm{N} - 1$ to obtain
\begin{align*}
\|\p^{\al}\zeta\|_{L^\infty_x} \leq C\|\p^{\al}\zeta\|_{H^1_x}^{\frac{1}{2}}\|\p^{\al}\zeta\|_{H^2_x}^{\frac{1}{2}} \leq C\Big(\mathcal{E}_{M}^{\frac{1}{2}}\Big)^{\frac{1}{2}} \Big(\kappa^{-\frac{1}{2}}\mathcal{E}_{top}^{\frac{1}{2}}\Big)^{\frac{1}{2}} \leq C\kappa^{-\frac{1}{4}}\mathcal{E}_{M}^{\frac{1}{4}}\mathcal{E}_{top}^{\frac{1}{4}}.
\end{align*}
For $d=2$, the result follows from the same argument.
\end{proof}

\begin{lemma}\label{L.G1} 
For $t\in[0,T]$ satisfying the bootstrap assumption \eqref{condition}, and for both the purely spatial derivatives \eqref{caseA} and the space-time derivatives \eqref{caseB}, the term $\AC{I}_1$ defined in \eqref{ACEdef} satisfies the following upper bound:
\begin{align*}
\AC{I}_1(t) &\les \int_{\Omega}|\nabla_xu^{\e}(t)| \mathcal{V}_2(t) dx + \int_{\Omega}|\nabla_x\ta^{\e}(t)| \mathcal{V}_3(t) dx + \frac{1}{\eps} \int_{\Omega} \TbT(t) \mathcal{V}_2(t) dx,
\end{align*}
where \( \mathcal{V}_\ell(t,x) \) and the control of the turbulence transport $\TbT(t,x)$ are defined in \eqref{largev} and \eqref{RSdef}, respectively.
\end{lemma}
\begin{proof}
Applying \eqref{Phitx} to $\AC{I}_1$, we get
\begin{align*}
\bega
\AC{I}_1 &=- \sum_{0\leq|\al|\leq\mathrm{N}}\frac{1}{2\eps^2}\int_{\Omega}\frac{1}{k_B\mathrm{\Theta}^{\e}}\sum_{i,j}\p_{x_i}u^{\e}_j \int_{\R^3} \frac{|\p^{\al}\AC{\P}F^{\e}|^2}{M^{\e}}\mathfrak{R}^{\e}_{ij}dvdx \cr
&-\sum_{0\leq|\al|\leq\mathrm{N}}\frac{1}{2\eps^2}\int_{\Omega}\frac{1}{k_B\mathrm{\Theta}^{\e}}\sum_i \p_{x_i}\ta^{\e} \int_{\R^3} \frac{|\p^{\al}\AC{\P}F^{\e}|^2}{M^{\e}}\mathcal{Q}^{\e}_idvdx \cr 
&+\sum_{0\leq|\al|\leq\mathrm{N}}\frac{1}{2\eps^2}\int_{\Omega}\frac{2}{\eps3k_B\mathrm{P}^{\e}\mathrm{\Theta}^{\e}} \bigg(\sum_{j} \p_{x_j} \mathfrak{q}_j^{\e} + \sum_{i,j} \p_{x_i}\mathrm{U}^{\e}_j \mathbf{r}_{ij}^{\e} \bigg) \int_{\R^3} \frac{|\p^{\al}\AC{\P}F^{\e}|^2}{M^{\e}}\lw( \frac{|v-\mathrm{U}^{\e}|^2}{2k_B\mathrm{\Theta}^{\e}}-\frac{3}{2} \rw)dvdx \cr 
&+\sum_{0\leq|\al|\leq\mathrm{N}}\frac{1}{2\eps^2}\int_{\Omega} \frac{1}{\e k_B\mathrm{P}^{\e}\mathrm{\Theta}^{\e}}\sum_{i,j} \p_{x_j} \mathbf{r}_{ij}^{\e} \int_{\R^3} \frac{|\p^{\al}\AC{\P}F^{\e}|^2}{M^{\e}}(v_i-\mathrm{U}^{\e}_i)dvdx.
\enda
\end{align*}
For the first two lines, we use the pointwise bounds 
$|\mathfrak{R}^{\e}_{ij}|\leq C(1+|v|^2)$ and $\mathcal{Q}^{\e}_i\leq C(1+|v|^3)$ from the bootstrap assumption $(\mathcal{B}_2)$ in \eqref{condition}. 
Then, by the definition of $\mathcal{V}_2$ and $\mathcal{V}_3$ in \eqref{largev}, we obtain the desired bounds for the first two lines. 
For the last two lines, we apply $L^\infty_x$ bounds to the terms 
$\sum_{j} \p_{x_j} \mathfrak{q}_j^{\e} + \sum_{i,j} \p_{x_i}\mathrm{U}^{\e}_j \mathbf{r}_{ij}^{\e}$ 
and $\p_{x_j} \mathbf{r}_{ij}^{\e}$. 
Then, by the definition of $\TbT(t)$ in \eqref{RSdef}, the last two lines are bounded by 
$\frac{1}{\eps} \int_{\Omega} \TbT(t) \mathcal{V}_2(t) dx$.
\end{proof}

\begin{lemma}\label{L.G23}
For $t\in[0,T]$ satisfying the bootstrap assumption \eqref{condition}, and for both the purely spatial derivatives \eqref{caseA} and the space-time derivatives \eqref{caseB}, the terms $\AC{I}_2$ and $\AC{I}_3$ defined in \eqref{ACEdef} satisfy the following estimates:
\begin{align*}
&\bigg|\AC{I}_2(t)  -\frac{-1}{\e^2} \!\!\!\sum_{0\leq|\al|\leq\mathrm{N}} \int_{\Omega} \frac{1}{k_B\mathrm{\Theta}^{\e}}\sum_{i,j}\p^{\al}\p_{x_i}u^{\e}_j(t) \p^{\al} \mathbf{r}_{ij}^{\e}(t)  \bigg| \lesssim 
\e \kappa^{\frac{1}{2}} \mathcal{E}_{M}(t) \mathcal{D}_G^{\frac{1}{2}}(t), \cr 
& \bigg|\AC{I}_3(t)- \frac{-1}{\e^2} \!\!\!\sum_{0\leq|\al|\leq\mathrm{N}} \int_{\Omega} \frac{1}{k_B\mathrm{\Theta}^{\e}}\sum_i \p^{\al}\p_{x_i}\ta^{\e}(t) \bigg(\p^{\al}\mathfrak{q}_i^{\e}(t)+\sum_{1\leq\beta\leq \al}\binom{\al}{\beta} \sum_{k} \p^{\beta}\mathrm{U}^{\e}_k(t) \p^{\al-\beta}\mathbf{r}_{ik}^{\e}(t)\bigg) \bigg| \cr
&\quad \les \e \kappa^{\frac{1}{2}} \mathcal{E}_{M}(t) \mathcal{D}_G^{\frac{1}{2}}(t).
\end{align*} 
\end{lemma}
\begin{proof}
(Estimate of $\AC{I}_2$) 
We decompose $\AC{I}_2$ according to whether all derivatives act on $\p_i u^{\e}_j$ or not, writing $\AC{I}_2 = \AC{I}_{2,1} + \AC{I}_{2,2}$:
\begin{align*}
\AC{I}_{2,1}&:=-\frac{1}{\e^2} \sum_{0\leq|\al|\leq\mathrm{N}} \int_{\Omega} \frac{1}{k_B\mathrm{\Theta}^{\e}}\sum_{i,j}\p^{\al}\p_{x_i}u^{\e}_j \int_{\R^3}\mathfrak{R}^{\e}_{ij} \p^{\al}\AC{\P}F^{\e} dv dx, \cr 
\AC{I}_{2,2}&:=-\frac{1}{\e^2} \sum_{0\leq|\al|\leq\mathrm{N}}\sum_{0<\beta\leq\al}\binom{\al}{\beta} \int_{\Omega} \sum_{i,j}\p^{\al-\beta}\p_{x_i}u^{\e}_j \int_{\R^3}\p^{\beta}\bigg(\frac{1}{k_B\mathrm{\Theta}^{\e}}\mathfrak{R}^{\e}_{ij}M^{\e} \bigg) \p^{\al}\AC{\P}F^{\e} |M^{\e}|^{-1} dv dx.  
\end{align*}
To estimate $\AC{I}_{2,1}$, we observe that, since $\p \mathfrak{R}^{\e}_{ij}$ is spanned by $(1,v)$ for any spatial or temporal derivative $\p^{\al}$, Lemma~\ref{HH} yields
\begin{align}\label{alhata}
\p^{\al} \mathbf{r}_{ij}^{\e} = \p^{\al} \la \mathfrak{R}^{\e}_{ij} , \AC{\P}F^{\e} \ra = \sum_{0\leq\beta\leq\al}\binom{\al}{\beta}\la \p^{\beta}\mathfrak{R}^{\e}_{ij},\p^{\al-\beta}\AC{\P}F^{\e} \ra =  \la \mathfrak{R}^{\e}_{ij} , \p^{\al}\AC{\P}F^{\e} \ra.
\end{align}
This gives 
\begin{align}\label{claimal}
\int_{\R^3}\mathfrak{R}^{\e}_{ij} \p^{\al}\AC{\P}F^{\e} dv = \p^{\al} \mathbf{r}_{ij}^{\e}.
\end{align}

For $\AC{I}_{2,2}$, applying H\"{o}lder's inequality in $x$ and $v$, we obtain
\begin{align*}
\AC{I}_{2,2}&\les \frac{1}{\e^2} \sum_{0\leq|\al|\leq\mathrm{N}}\sum_{1\leq\beta\leq\al} \bigg(\int_{\Omega}|\p^{\al-\beta}\nabla_xu^{\e}|^2 \int_{\R^3}\bigg|\p^{\beta}\bigg(\frac{1}{k_B\mathrm{\Theta}^{\e}}\mathfrak{R}^{\e}_{ij}M^{\e} \bigg)\bigg|^2 |M^{\e}|^{-1}dv dx\bigg)^{\frac{1}{2}} \cr 
&\times \bigg(\int_{\Omega \times \R^3} |\p^{\al}\AC{\P}F^{\e}|^2 |M^{\e}|^{-1} dvdx\bigg)^{\frac{1}{2}}.
\end{align*}
Using the definition of the dissipation $\mathcal{D}_G$ in \eqref{N-EDdef}, we get
\begin{align*}
\AC{I}_{2,2}&\les \kappa^{\frac{1}{2}} \sum_{0\leq|\al|\leq\mathrm{N}}\sum_{1\leq\beta\leq\al} \bigg(\int_{\Omega}|\p^{\al-\beta}\nabla_xu^{\e}|^2 \int_{\R^3}\bigg|\p^{\beta}\bigg(\frac{1}{k_B\mathrm{\Theta}^{\e}}\mathfrak{R}^{\e}_{ij}M^{\e} \bigg)\bigg|^2 |M^{\e}|^{-1}dv dx\bigg)^{\frac{1}{2}} \mathcal{D}_G^{\frac{1}{2}}.
\end{align*}
We apply $L^2_x$ or $L^\infty_x$ to $\p^{\al - \beta} \nabla_x u^{\e}$ depending on $1\leq |\beta| \leq \lfloor \mathrm{N}/2 \rfloor+1$ or $\lfloor \mathrm{N}/2 \rfloor+1< |\beta| \leq \mathrm{N}$. We then expand $\p^{\beta}(\frac{1}{k_B \mathrm{\Theta}^{\e}} \mathfrak{R}^{\e}_{ij} M^{\e})$ using \eqref{al-vM} and apply the estimates \eqref{Phiscale} and \eqref{Phiinfscale}. This yields
\begin{align*}
\AC{I}_{2,2}&
\les \kappa^{\frac{1}{2}} \mathcal{E}_{M}^{\frac{1}{2}} (\e\mathcal{E}_{M}^{\frac{1}{2}})\mathcal{D}_G^{\frac{1}{2}}  \les \e \kappa^{\frac{1}{2}} \mathcal{E}_{M} \mathcal{D}_G^{\frac{1}{2}}.
\end{align*}
Combining this with \eqref{claimal} for $\AC{I}_{2,1}$ gives the desired estimate for $\AC{I}_2$. \\
(Estimate of $\AC{I}_3$) 
In the same way, we split $\AC{I}_3$ as $\AC{I}_3 = \AC{I}_{3,1} + \AC{I}_{3,2}$, where
\begin{align*}
\bega
\AC{I}_{3,1}&=-\frac{1}{\e^2}\sum_{0\leq|\al|\leq\mathrm{N}} \int_{\Omega}\frac{1}{k_B\mathrm{\Theta}^{\e}}\sum_i \p^{\al}\p_{x_i}\ta^{\e} \int_{\R^3} \mathcal{Q}^{\e}_i\p^{\al}\AC{\P}F^{\e} dvdx, \cr 
\AC{I}_{3,2}&=-\frac{1}{\e^2}\sum_{0\leq|\al|\leq\mathrm{N}}\sum_{1\leq\beta\leq\al}\binom{\al}{\beta} \int_{\Omega} \sum_i \p^{\al-\beta}\p_{x_i}\ta^{\e} \int_{\R^3}\p^{\beta}\bigg(\frac{1}{k_B\mathrm{\Theta}^{\e}} \mathcal{Q}^{\e}_iM^{\e}\bigg)\p^{\al}\AC{\P}F^{\e}|M^{\e}|^{-1} dvdx.
\enda
\end{align*}
For the estimate of $\AC{I}_{3,1}$, we claim that
\begin{align}\label{claimbeta}
\int_{\R^3} \mathcal{Q}^{\e}_i\p^{\al}\AC{\P}F^{\e} dv = \p^{\al} \mathfrak{q}_i^{\e} + \sum_{
0<\beta\leq \al}\binom{\al}{\beta} \sum_{k} \p^{\beta}\mathrm{U}^{\e}_k \p^{\al-\beta}\mathbf{r}_{ik}^{\e}.
\end{align}
By direct computation, for $|\al| \geq 1$, we have
\begin{align*}
\bega
\p^{\al}\mathcal{Q}^{\e}_j &= -\frac{1}{2} \p^{\al}\mathrm{U}^{\e}_j(|v-\mathrm{U}^{\e}|^2-5k_B\mathrm{\Theta}^{\e}) + \frac{1}{2}\sum_{0<\beta<\al}\binom{\al}{\beta} \p^{\beta}(v_j-\mathrm{U}^{\e}_j)\p^{\al-\beta}(|v-\mathrm{U}^{\e}|^2-5k_B\mathrm{\Theta}^{\e})\cr 
&+ \frac{1}{2}(v_j-\mathrm{U}^{\e}_j) \bigg(-2 \p^{\al}\mathrm{U}^{\e} \cdot (v-\mathrm{U}^{\e})  + \sum_{0<\beta<\al}\binom{\al}{\beta}\p^{\beta}\mathrm{U}^{\e}\cdot \p^{\al-\beta}\mathrm{U}^{\e}\bigg).
\enda
\end{align*}
We observe that the right-hand side of $\p^{\al} \mathcal{Q}^{\e}_j$ is spanned by $(1, v, |v|^2)$, except for the term  $\p^{\al}\mathrm{U}^{\e}_k(v_j-\mathrm{U}^{\e}_j)(v_k-\mathrm{U}^{\e}_k)$. Therefore, using Lemma \ref{HH}, we get
\begin{align}\label{alhatb}
\bega
\p^{\al}\mathfrak{q}_i^{\e} = \p^{\al} \la \mathcal{Q}^{\e}_i , \AC{\P}F^{\e} \ra &= \la \mathcal{Q}^{\e}_i , \p^{\al}\AC{\P}F^{\e} \ra + \sum_{1\leq\beta\leq \al}\binom{\al}{\beta}\langle \p^{\beta}\mathcal{Q}^{\e}_i, \p^{\al-\beta}\AC{\P}F^{\e} \rangle \cr 
&= \la \mathcal{Q}^{\e}_i , \p^{\al}\AC{\P}F^{\e} \ra - \sum_{1\leq\beta\leq \al}\binom{\al}{\beta} \sum_{k} \p^{\beta}\mathrm{U}^{\e}_k \langle \mathfrak{R}^{\e}_{ik} , \p^{\al-\beta}\AC{\P}F^{\e} \rangle.
\enda
\end{align}
Applying \eqref{claimal} to the last term in \eqref{alhatb} yields the desired identity \eqref{claimbeta}.  
The estimate for $\AC{I}_{3,2}$ can then be obtained in the same way as the estimate for $\AC{I}_{2,2}$.
\end{proof}

\begin{lemma}\label{L.G45}
For $t\in[0,T]$ satisfying the bootstrap assumption \eqref{condition}, and for both the purely spatial derivatives \eqref{caseA} and the space--time derivatives \eqref{caseB}, the terms $\AC{I}_4$ and $\AC{I}_5$, defined in \eqref{ACEdef}, satisfy the following upper bound:
\begin{align*}
\AC{I}_4(t)&\les \e \mathcal{E}_{M}^{\frac{1}{2}}(t) \mathcal{D}_G(t)+ \frac{1}{\kappa^{\frac{1}{2}}\e} \sum_{0\leq \al \leq \lfloor \mathrm{N}/2 \rfloor}\Big\|\|\la v \ra^{\frac{1}{2}} \p^{\al}\AC{\P}F^{\e} |M^{\e}|^{-\frac{1}{2}}(t)\|_{L^2_v}\Big\|_{L^\infty_x} \mathcal{E}_{M}^{\frac{1}{2}}(t)\mathcal{D}_G^{\frac{1}{2}}(t), \cr 
\AC{I}_5(t)&\les \sum_{0\leq|\al|\leq \lfloor\mathrm{N}/2\rfloor}\Big\|\|\la v \ra^{\frac{1}{2}} \p^{\al}\AC{\P}F^{\e} |M^{\e}|^{-\frac{1}{2}}(t)\|_{L^2_v}\Big\|_{L^\infty_x} \mathcal{D}_G(t),
\end{align*}
where $\lfloor \cdot \rfloor$ denotes the greatest integer function.
\end{lemma}
\begin{proof}
(Estimate of $\AC{I}_4(t)$)
We first note that, as long as $t \in [0, T]$ satisfies the bootstrap assumption \eqref{condition}, the collision frequency defined in \eqref{nudef} satisfies the following pointwise equivalence for some constant $C>0$:
\begin{align}\label{nusimilar}
\frac{1}{C}\sqrt{1+|v|^2} \leq 
\nu^{\e}(t,x,v)
\leq C\sqrt{1+|v|^2 }.
\end{align}
Applying \eqref{nonlin} to $\AC{I}_4$, we obtain the following upper bound for $\AC{I}_4$:
\begin{align*}
\bega 
\sum_{0\leq|\al|\leq\mathrm{N}}\frac{C}{\kappa\e^4} \int_{\Omega} \sum_{1\leq \beta \leq \al} \Big\|\la v \ra^{\frac{1}{2}} \p^{\al-\beta}\AC{\P}F^{\e}|M^{\e}|^{-\frac{1}{2}}\Big\|_{L^2_v}\Big\|\la v \ra^{\frac{1}{2}} \p^{\beta}M^{\e}|M^{\e}|^{-\frac{1}{2}}\Big\|_{L^2_v} \Big\|\la v \ra^{\frac{1}{2}} \p^{\al}\AC{\P}F^{\e}|M^{\e}|^{-\frac{1}{2}}\Big\|_{L^2_v}dx .
\enda
\end{align*}
We apply the H\"{o}lder inequality and take the $L^\infty_x$ norm on the factor with lower regularity between $\p^{\beta}$ and $\p^{\al - \beta}$. Then,
\begin{align*}
\bega 
\AC{I}_4&\les  \sum_{0\leq|\al|\leq\mathrm{N}}\frac{1}{\kappa^{\frac{1}{2}}\e^2} \sum_{1\leq \beta \leq \lfloor |\al|/2  \rfloor} \Big\|\la v \ra^{\frac{1}{2}} \p^{\al-\beta}\AC{\P}F^{\e}|M^{\e}|^{-\frac{1}{2}}\Big\|_{L^2_{x,v}}\Big\|\la v \ra^{\frac{1}{2}} \p^{\beta}M^{\e}|M^{\e}|^{-\frac{1}{2}}\Big\|_{L^\infty_xL^2_v} \mathcal{D}_G^{\frac{1}{2}} \cr 
&+\sum_{0\leq|\al|\leq\mathrm{N}}\frac{1}{\kappa^{\frac{1}{2}}\e^2} \sum_{\lfloor |\al|/2  \rfloor< \beta \leq |\al|} \Big\|\la v \ra^{\frac{1}{2}} \p^{\al-\beta}\AC{\P}F^{\e}|M^{\e}|^{-\frac{1}{2}}\Big\|_{L^\infty_xL^2_v}\Big\|\la v \ra^{\frac{1}{2}} \p^{\beta}M^{\e}|M^{\e}|^{-\frac{1}{2}}\Big\|_{L^2_{x,v}} \mathcal{D}_G^{\frac{1}{2}},
\enda
\end{align*}
where we used the dissipation $\mathcal{D}_G$ defined in \eqref{N-EDdef}. 
Applying $\mathcal{D}_G$ again, along with the estimates \eqref{Phiscale} and \eqref{Phiinfscale}, we obtain
\begin{align*}
\bega 
\AC{I}_4&\les \frac{1}{\kappa^{\frac{1}{2}}\e^2} (\e^2\kappa^{\frac{1}{2}}\mathcal{D}_G^{\frac{1}{2}}) (\e\mathcal{E}_{M}^{\frac{1}{2}}) \mathcal{D}_G^{\frac{1}{2}} + \frac{1}{\kappa^{\frac{1}{2}}\e^2} \sum_{0\leq \al \leq \lfloor \mathrm{N}/2  \rfloor} \Big\|\la v \ra^{\frac{1}{2}} \p^{\al}\AC{\P}F^{\e}|M^{\e}|^{-\frac{1}{2}}\Big\|_{L^\infty_xL^2_v} (\e\mathcal{E}_{M}^{\frac{1}{2}})\mathcal{D}_G^{\frac{1}{2}}.
\enda
\end{align*}
This completes the estimate. \\
(Estimate of $\AC{I}_5(t)$) 
Applying \eqref{nonlin} to $\AC{I}_5$ gives 
\begin{align*}
\bega 
\sum_{0\leq|\al|\leq\mathrm{N}}\frac{1}{\kappa\e^4} \int_{\Omega} \sum_{0\leq\beta\leq\al} \Big\|\la v \ra^{\frac{1}{2}} \p^{\al-\beta}\AC{\P}F^{\e}|M^{\e}|^{-\frac{1}{2}}\Big\|_{L^2_v}\Big\|\la v \ra^{\frac{1}{2}} \p^{\beta}\AC{\P}F^{\e} |M^{\e}|^{-\frac{1}{2}}\Big\|_{L^2_v} \Big\|\la v \ra^{\frac{1}{2}} \p^{\al}\AC{\P}F^{\e}|M^{\e}|^{-\frac{1}{2}}\Big\|_{L^2_v}dx .
\enda
\end{align*}
Since the upper bound is symmetric with respect to the exchange $(\al - \beta, \beta) \leftrightarrow (\beta, \al - \beta)$, we may assume without loss of generality that $0 \leq |\beta| \leq \lfloor \mathrm{N}/2 \rfloor$. Applying the H\"{o}lder inequality and taking the $L^\infty_x$ norm on the term involving $\p^{\beta} \AC{\P}F^{\e}$, we obtain
\begin{align*}
\bega 
\AC{I}_5
&\les \sum_{0\leq|\beta|\leq \lfloor\mathrm{N}/2\rfloor}\bigg\|\Big\|\la v \ra^{\frac{1}{2}} \p^{\beta}\AC{\P}F^{\e} |M^{\e}|^{-\frac{1}{2}}\Big\|_{L^2_v}\bigg\|_{L^\infty_x} \mathcal{D}_G,
\enda
\end{align*}
where we used the dissipation $\mathcal{D}_G$ defined in \eqref{N-EDdef}. This completes the proof.
\end{proof}

\begin{lemma}\label{L.G67}
For $t\in[0,T]$ satisfying the bootstrap assumption \eqref{condition}, and for both the purely spatial derivatives \eqref{caseA} and the space-time derivatives \eqref{caseB}, the terms $\AC{I}_6$ and $\AC{I}_7$ defined in \eqref{ACEdef} satisfy the following upper bound:
\begin{align*}
\AC{I}_6(t)&\leq C\e^2 \kappa^{\frac{1}{2}} \bigg(\sum_{0\leq|\al|\leq\mathrm{N}-1}\frac{1}{\e^2}\bigg\|\sum_{i,j}\p^{\al}\p_{x_j} \mathbf{r}_{ij}^{\e}(t)\bigg\|_{L^2_x}\bigg) \mathcal{E}_{M}^{\frac{1}{2}}(t) \mathcal{D}_G^{\frac{1}{2}}(t), \cr 
\AC{I}_7(t)&\leq C\e^2 \kappa^{\frac{1}{2}} \bigg(\sum_{0\leq|\al|\leq\mathrm{N}-1}\frac{1}{\e^2}\bigg\|\sum_{j} \p^{\al}\p_{x_j} \mathfrak{q}_j^{\e}(t) + \sum_{i,j} \p^{\al}\big(\p_{x_i}\mathrm{U}^{\e}_j(t) \mathbf{r}_{ij}^{\e}(t)\big)\bigg\|_{L^2_x}\bigg) \mathcal{E}_{M}^{\frac{1}{2}}(t) \mathcal{D}_G^{\frac{1}{2}}(t).
\end{align*}
\end{lemma}
\begin{proof}
We note that in the terms $\AC{I}_6$ and $\AC{I}_7$, if all derivatives fall on $\mathbf{r}_{ij}^{\e}$ or $\mathfrak{q}_j^{\e}$, then the corresponding integrals vanish due to Lemma \ref{HH}, which gives
\begin{align*}
\AC{I}_6&= \frac{1}{\e^3}\sum_{0\leq|\al|\leq\mathrm{N}} \sum_{1\leq\beta\leq\al}\int_{\Omega \times \R^3}   \binom{\al}{\beta} \p^{\al-\beta}\bigg(\sum_{i,j}\p_{x_j} \mathbf{r}_{ij}^{\e}\bigg) \p^{\beta}\bigg(\frac{1}{ k_B\mathrm{P}^{\e}\mathrm{\Theta}^{\e}} (v_i-\mathrm{U}^{\e}_i)M^{\e}\bigg) |M^{\e}|^{-1} \p^{\al}\AC{\P}F^{\e} dvdx, \cr 
\AC{I}_7&=\frac{1}{\e^3}\sum_{0\leq|\al|\leq\mathrm{N}} \sum_{1\leq\beta\leq\al} \int_{\Omega \times \R^3} \binom{\al}{\beta} \p^{\al-\beta}\bigg(\sum_{j} \p_{x_j} \mathfrak{q}_j^{\e} + \sum_{i,j} \p_{x_i}\mathrm{U}^{\e}_j \mathbf{r}_{ij}^{\e} \bigg) \cr 
&\quad \times \p^{\beta} \bigg( \frac{2}{ 3k_B\mathrm{P}^{\e}\mathrm{\Theta}^{\e}} \lw( \frac{|v-\mathrm{U}^{\e}|^2}{2k_B\mathrm{\Theta}^{\e}}-\frac{3}{2} \rw)M^{\e} \bigg) |M^{\e}|^{-1} \p^{\al}\AC{\P}F^{\e} dvdx.
\end{align*}
Since the estimates for both terms are analogous, we provide the detailed estimate only for $\AC{I}_6$.
Applying H\"{o}lder's inequality in $x$ and $v$, we have
\begin{align*}
\AC{I}_6& \les \frac{1}{\e^3} \sum_{0\leq|\al|\leq\mathrm{N}}\sum_{1\leq\beta\leq\al} \bigg(\int_{\Omega}\bigg|\p^{\al-\beta}\bigg(\sum_{i,j}\p_{x_j} \mathbf{r}_{ij}^{\e}\bigg)\bigg|^2 \int_{\R^3}\bigg| \p^{\beta}\bigg(\frac{1}{ k_B\mathrm{P}^{\e}\mathrm{\Theta}^{\e}} (v_i-\mathrm{U}^{\e}_i)M^{\e}\bigg)\bigg|^2 |M^{\e}|^{-1}dv dx\bigg)^{\frac{1}{2}} \cr 
&\quad \times \bigg(\int_{\Omega \times \R^3} |\p^{\al}\AC{\P}F^{\e}|^2 |M^{\e}|^{-1} dvdx\bigg)^{\frac{1}{2}}. 
\end{align*}
Following the same approach used for estimating $\AC{I}_{2,2}$ in Lemma \ref{L.G23} — namely, invoking the dissipation norm $\mathcal{D}_G$ from \eqref{N-EDdef}, and using the estimates \eqref{al-vM}, \eqref{Phiscale}, and \eqref{Phiinfscale} — we obtain
\begin{align*}
\AC{I}_6&\les \e\kappa^{\frac{1}{2}} \bigg(\sum_{0\leq|\al|\leq\mathrm{N}-1}\frac{1}{\e^2}\bigg\|\sum_{i,j}\p^{\al}\p_{x_j} \mathbf{r}_{ij}^{\e}\bigg\|_{L^2_x}\bigg) (\e\mathcal{E}_{M}^{\frac{1}{2}}) \mathcal{D}_G^{\frac{1}{2}},
\end{align*}
where we used the Agmon inequality \eqref{Agmon} to estimate $\|\p^{\al-\beta}\p_{x_j} \mathbf{r}_{ij}^{\e}\|_{L^\infty_x} \les \|\p^{\al-\beta}\p_{x_j} \mathbf{r}_{ij}^{\e}\|_{H^2_x}$ when $\lfloor \mathrm{N}/2 \rfloor +1 \leq |\beta| \leq \mathrm{N}$.
\end{proof}

\begin{proof}[\textbf{Proof of Proposition \ref{P.G.Energy}}] 
We combine the estimates of $\AC{I}_i$ for $i = 1, \dots, 7$ from Lemmas \ref{L.G1}, \ref{L.G23}, \ref{L.G45}, and \ref{L.G67}, and insert them into \eqref{ACE-1}:
\begin{align*}
\bega
\sum_{0\leq|\al|\leq\mathrm{N}}& \frac{1}{2\e^2}\frac{d}{dt} \int_{\Omega \times \R^3}|\p^{\al}\AC{\P}F^{\e}|^2|M^{\e}|^{-1}dvdx + \sum_{0\leq|\al|\leq\mathrm{N}}\frac{1}{\kappa\e^4}\int_{\Omega \times \R^3}\mathcal{L}(\p^{\al}\AC{\P}F^{\e})(\p^{\al}\AC{\P}F^{\e})|M^{\e}|^{-1} dvdx  \cr
&\leq C\bigg( \int_{\Omega}|\nabla_xu^{\e}| \mathcal{V}_2 dx + \int_{\Omega}|\nabla_x\ta^{\e}| \mathcal{V}_3 dx + \frac{1}{\eps} \int_{\Omega} \TbT \mathcal{V}_2 dx \Bigg) + \sum_{0\leq|\al|\leq\mathrm{N}}\mathfrak{S}_G^{\al} \cr 
&+ C\e \kappa^{\frac{1}{2}} \mathcal{E}_{M} \mathcal{D}_G^{\frac{1}{2}} + C\e \mathcal{E}_{M}^{\frac{1}{2}} \mathcal{D}_G+ \frac{C}{\kappa^{\frac{1}{2}}\e}\sum_{0\leq \al \leq \lfloor \mathrm{N}/2 \rfloor}\Big\|\|\la v \ra^{\frac{1}{2}} \p^{\al}\AC{\P}F^{\e} |M^{\e}|^{-\frac{1}{2}}\|_{L^2_v}\Big\|_{L^\infty_x} \mathcal{E}_{M}^{\frac{1}{2}}\mathcal{D}_G^{\frac{1}{2}} \cr 
&+ C\sum_{0\leq|\al|\leq \lfloor\mathrm{N}/2\rfloor}\Big\|\|\la v \ra^{\frac{1}{2}} \p^{\al}\AC{\P}F^{\e} |M^{\e}|^{-\frac{1}{2}}\|_{L^2_v}\Big\|_{L^\infty_x} \mathcal{D}_G \cr 
&+ C\e^2 \kappa^{\frac{1}{2}} \sum_{0\leq|\al|\leq\mathrm{N}-1}\frac{1}{\e^2}\bigg(\bigg\|\sum_{i,j}\p^{\al}\p_{x_j} \mathbf{r}_{ij}^{\e}\bigg\|_{L^2_x}+\bigg\|\sum_{j} \p_{x_j} \mathfrak{q}_j^{\e} + \sum_{i,j} \p_{x_i}\mathrm{U}^{\e}_j \mathbf{r}_{ij}^{\e}\bigg\|_{L^2_x}\bigg) \mathcal{E}_{M}^{\frac{1}{2}} \mathcal{D}_G^{\frac{1}{2}}.
\enda
\end{align*}
For the second term on the first line, we apply the coercivity estimate in Lemma \ref{L.coer}, then use the bootstrap assumption $(\mathcal{B}_2)$ in \eqref{condition}. By the definitions of $\mathcal{E}_G$ and $\mathcal{D}_G$ in \eqref{N-EDdef}, we derive a lower bound of the first line as $\frac{d}{dt} \mathcal{E}_G(t) + \frac{\sigma_L}{C} \mathcal{D}_G(t)$.
Finally, we apply $\eqref{ABGscale}_1$ to the last line, which completes the proof of Proposition \ref{P.G.Energy}.
\end{proof}

\subsection{Top-order microscopic estimate}\label{Sec.3.Top}

In this subsection, we derive the top-order energy estimate for
$\p^{\al}F^{\e}$ with $|\al|=\mathrm{N}+1$,
which incurs a loss of a factor of $\kappa$ compared to the
lower-order estimates.
Unlike the estimate for $\AC{\P}F^{\e}$ in
Subsection~\ref{Sec.3.G}, whose right-hand side of~\eqref{Geqn0}
does not involve macro--macro interactions,
we must here control the full nonlinear contribution.
The leading terms arise from the macro--macro interaction and from the commutator
generated by the weight $1/M^{\e}$.

\begin{proposition}\label{P.F.Energy}
Let $\Omega=\R^d$ with $d=2$ or $3$ and $\mathrm{N}>d/2+1$.
Assume that the initial data $u_0^{\e}$ satisfies either the finite velocity energy condition \eqref{caseEC} or the infinite velocity energy condition \eqref{caseECX}.
Under the bootstrap assumption \eqref{condition}, the following estimate holds in both the purely spatial derivative case \eqref{caseA} and the space-time derivative case \eqref{caseB} with arbitrary $\mathfrak{n} \geq 1$.
\begin{align}\label{totalEt}
\bega
&\sum_{|\al|=\mathrm{N}+1}\frac{\kappa}{2\e^2}\frac{d}{dt} \int_{\Omega \times \R^3}|\p^{\al}F^{\e}|^2 |M^{\e}|^{-1}(t)dvdx  +\frac{\sigma_L}{C} \mathcal{D}_{top}^{\mathrm{N}}(F^{\e}(t)) \cr 
&\les \big(\|\nabla_xu^{\e}(t)\|_{L^\infty_x}+\|\nabla_x\ta^{\e}(t)\|_{L^\infty_x}\big)\big(\mathcal{E}_{tot}^{\mathrm{N}}(F^{\e}(t))+(\mathcal{E}_M^{\mathrm{N}}(F^{\e}(t)))^{\frac{1}{2}} (\mathcal{D}_{top}^{\mathrm{N}}(F^{\e}(t)))^{\frac{1}{2}}\big)  \cr
&+ \int_{\Omega}|\nabla_xu^{\e}(t)| \mathcal{V}_2(F^{\e}(t)) dx + \int_{\Omega}|\nabla_x\ta^{\e}(t)| \mathcal{V}_3(F^{\e}(t)) dx + \frac{1}{\eps} \int_{\Omega} \TbT(t) \mathcal{V}_2(F^{\e}(t)) dx \cr 
&+ \e\kappa^{\frac{1}{2}} \Big( \mathcal{E}_{tot}^{\mathrm{N}}(F^{\e}(t))\mathcal{D}_{tot}^{\mathrm{N}}(F^{\e}(t)))^{\frac{1}{2}}+(\mathcal{E}_{tot}^{\mathrm{N}}(F^{\e}(t)))^{\frac{1}{2}} \mathcal{D}_{tot}^{\mathrm{N}}(F^{\e}(t)) \Big) \cr 
&+\e(\mathcal{E}_{tot}^{\mathrm{N}}(F^{\e}(t)))^{\frac{3}{2}} (\mathcal{D}_{top}^{\mathrm{N}}(F^{\e}(t)))^{\frac{1}{2}} +\mathcal{Z}_{top}^{time}(t)\cr 
&+ \frac{1}{\e}\Big\|\|\la v \ra^{\frac{1}{2}} \AC{\P}F^{\e} |M^{\e}|^{-\frac{1}{2}}(t)\|_{L^2_v}\Big\|_{L^\infty_x} \Big(\kappa^{-\frac{1}{2}} (\mathcal{E}_{tot}^{\mathrm{N}}(F^{\e}(t)))^{\frac{1}{2}} + \e(\mathcal{D}_{top}^{\mathrm{N}}(F^{\e}(t)))^{\frac{1}{2}}\Big) (\mathcal{D}_{top}^{\mathrm{N}}(F^{\e}(t)))^{\frac{1}{2}} \cr 
& + \frac{1}{\e}\sum_{1\leq |\al| \leq\lfloor(\mathrm{N}+1)/2\rfloor} \Big\|\|\la v \ra^{\frac{1}{2}} \p^{\al}\AC{\P}F^{\e} |M^{\e}|^{-\frac{1}{2}}(t)\|_{L^2_v}\Big\|_{L^\infty_x} (\mathcal{E}_{M}^{\mathrm{N}}(F^{\e}(t)))^{\frac{1}{2}} (\mathcal{D}_{top}^{\mathrm{N}}(F^{\e}(t)))^{\frac{1}{2}} \cr 
&+ \kappa^{\frac{1}{2}} \sum_{0\leq|\al|\leq \lfloor(\mathrm{N}+1)/2\rfloor}\Big\|\|\la v \ra^{\frac{1}{2}} \p^{\al}\AC{\P}F^{\e} |M^{\e}|^{-\frac{1}{2}}(t)\|_{L^2_v}\Big\|_{L^\infty_x} (\mathcal{D}_G^{\mathrm{N}}(F^{\e}(t)))^{\frac{1}{2}} (\mathcal{D}_{top}^{\mathrm{N}}(F^{\e}(t)))^{\frac{1}{2}},
\enda
\end{align}
where $\sigma_L>0$ is the constant given in Lemma~\ref{L.coer}, and $C>0$ is a universal constant.

The contribution from time derivatives, denoted by 
$\mathcal{Z}_{top}^{\mathrm{time}}(t)$, is defined by 
$\mathcal{Z}_{top}^{\mathrm{time}}(t)=0$ in the purely spatial case \eqref{caseA}, and by
\begin{align}\label{Bmidef}
\begin{aligned}
\mathcal{Z}_{top}^{time}(t) &:=\e^{\mathfrak{n}-1}  \Big[\big\|\big((\nabla_x\cdot u^{\e}),\nabla_x(\rho^{\e}+\ta^{\e})\big)(t)\big\|_{L^\infty_x} + \e \|u^{\e}(t)\|_{L^\infty_x}\, (\mathcal{E}_{M}^{\mathrm{N}}(F^{\e}(t)))^{\tfrac{1}{2}} 
\cr 
&\quad + \e\kappa^{\frac{1}{2}} (\mathcal{D}_G^{\mathrm{N}}(F^{\e}(t)))^{\frac{1}{2}} \Big]  \times
(\mathcal{E}_{M}^{\mathrm{N}}(F^{\e}(t)))^{\tfrac{1}{2}} (\mathcal{D}_{top}^{\mathrm{N}}(F^{\e}(t)))^{\tfrac{1}{2}},
\end{aligned}
\end{align}
in the space--time derivative case \eqref{caseB}.

\hide 
We can choose different $\mathcal{Z}_{top}^{time}(t)$ depending on whether derivatives are purely spatial or not. For the purely spatial derivative case \eqref{caseA}, $\mathcal{Z}_{top}^{time}(t)=0$. On ther other hand, for the space-time derivative case \eqref{caseB}, we have that 
\begin{align}\label{Bmidef}
\bega
\mathcal{Z}_{top}^{time}(t) &:= \e^{\mathfrak{n}-1} \Big( \big\|\big((\nabla_x\cdot u^{\e}),\nabla_x(\rho^{\e}+\ta^{\e})\big)(t)\big\|_{L^\infty_x} +\e\|u^{\e}\|_{L^\infty_x}\mathcal{E}_{M}^{\frac{1}{2}} + \e \|\bPhi^{\e}_{\bW}\|_{L^\infty_x} \Big)\mathcal{E}_{M}^{\frac{1}{2}}\mathcal{D}_{top}^{\frac{1}{2}},
\enda
\end{align}\unhide  

In addition, we have the following equivalence relation for the top-order energy:
\begin{align}\label{E0eq}
\bega
\frac{1}{C} \Big(\mathcal{E}_{top}^{\mathrm{N}}(F^{\e}(t)) -\e\kappa^{\frac{1}{2}} (\mathcal{E}_{tot}^{\mathrm{N}}(F^{\e}(t)))^{\frac{3}{2}}(t)\Big) &\leq \bigg|\sum_{|\al|=\mathrm{N}+1} \frac{\kappa}{2\e^2} \int_{\Omega \times \R^3} |\p^{\al}F^{\e}|^2 |M^{\e}|^{-1}(t) \, dv dx\bigg|\cr 
&\leq C \Big(\mathcal{E}_{top}^{\mathrm{N}}(F^{\e}(t)) +\e\kappa^{\frac{1}{2}} (\mathcal{E}_{tot}^{\mathrm{N}}(F^{\e}(t)))^{\frac{3}{2}}(t)\Big) , 
\enda
\end{align}
for some constant $C > 0$.
\end{proposition}

We present the proof of Proposition \ref{P.F.Energy} at the end of this section. The proof involves the analysis of three leading-order terms:
\begin{itemize}
\item The equivalence of the left-hand side of \eqref{totalEt} and $\mathcal{E}_{top}$. (See Lemma \ref{L.E0}.)
\item Control of the nonlinear macro-macro interaction term. (See Lemma \ref{L.E1}.)
\item  Control of the convection against the local Maxwellian. (See Lemma \ref{L.E2}.)
\end{itemize}
We first present an identity for the energy estimate in Lemma \ref{L.EFtop}, and then estimate each term individually.
Before proceeding, we begin by decomposing $\p^{\al}F^{\e}$ into its macroscopic and microscopic components.

\begin{lemma}\label{PalF} 
For both the purely spatial derivative case \eqref{caseA} and the space-time derivative case \eqref{caseB} for arbitrary $\mathfrak{n} \geq 0$, the projection operators $\P$ and $\AC{\P}$ defined in \eqref{Pdef} satisfy the following properties:

(1) Whenever defined, $\P\big(\p^{\al} (\AC{\P} F^{\e})\big) = 0$.

(2) We have
\begin{equation}\label{palF}
\p^{\al}F^{\e} = \P(\p^{\al} M^{\e}) + \AC{\P}(\p^{\al}F^{\e}).
\end{equation}
\end{lemma}
\begin{proof}
Both statements are trivial when $\al = 0$. Assume now that $\al \neq 0$. From the definition of $\P$ in \eqref{Pdef}, we note that $\P(\p^{\al}(\AC{\P}F^{\e}))$ is a linear combination of moments of the form
$\int_{\R^3} \p^{\al}(\AC{\P}F^{\e})(1,v,|v|^2) \, dv$. 
Since $\AC{\P}F^{\e}$ is orthogonal to the null space of the linear operator, it satisfies
\begin{align*}
\int_{\R^3} \p^{\al}(\AC{\P}F^{\e})(1,v,|v|^2) \, dv = \p^{\al} \left( \int_{\R^3} \AC{\P}F^{\e} (1,v,|v|^2) \, dv \right) = 0,
\end{align*}
which proves (1). 
To prove (2), we recall the decomposition $F^{\e} = M^{\e} + \AC{\P}F^{\e}$, which gives $
\P(\p^{\al}F^{\e}) = \P(\p^{\al} M^{\e}) + \P\big(\p^{\al} (\AC{\P}F^{\e})\big)$.
Thus, (1) of this lemma implies $\P\big(\p^{\al} (\AC{\P}F^{\e})\big) = 0$, and we obtain \eqref{palF}.
\end{proof}

In the proof of Proposition~\ref{P.F.Energy}, for brevity, we slightly abuse notation by writing $\mathcal{E}(t)$, $\mathcal{D}(t)$, and $\mathcal{V}_{\ell}(t)$ for $\mathcal{E}^{\mathrm{N}}(F^{\e}(t))$,
$\mathcal{D}^{\mathrm{N}}(F^{\e}(t))$, and $\mathcal{V}_{\ell}(F^{\e}(t))$, respectively.

\begin{lemma}\label{L.EFtop}
For the distribution function $F^{\e}$ satisfying the scaled Boltzmann equation \eqref{BE}, and for both types of derivatives $\p^{\al}$ defined in \eqref{caseA} and \eqref{caseB}, we have the following energy estimate:
\begin{align}\label{E'0}
\bega
\sum_{|\al|=\mathrm{N}+1}\frac{\kappa}{2\e^2}\frac{d}{dt} \int_{\Omega \times \R^3}\frac{|\p^{\al}F^{\e}|^2}{M^{\e}} dvdx  &+ \sum_{|\al|=\mathrm{N}+1}\frac{1}{\e^4}\int_{\Omega \times \R^3}\mathcal{L}(\p^{\al}F^{\e})\frac{\p^{\al}F^{\e}}{M^{\e}} dvdx \cr 
&= I^F_1(t)+I^F_2(t)+I^F_3(t)+I^F_4(t)+I^F_5(t),
\enda
\end{align}
where
\begin{align}\label{E123def}
\bega
I^F_1(t)&=\sum_{|\al|=\mathrm{N}+1}\frac{1}{\e^4} \int_{\Omega \times \R^3} \mathcal{L}(\p^{\al}M^{\e}) \frac{\AC{\P}(\p^{\al}F^{\e})}{M^{\e}} dvdx, \cr 
I^F_2(t)&= -\sum_{|\al|=\mathrm{N}+1}\frac{\kappa}{2\e^2}\int_{\Omega \times \R^3} \left(\p_t M^{\e} + \frac{v}{\eps}\cdot\nabla_x M^{\e}\right) \frac{|\p^{\al}F^{\e}|^2}{|M^{\e}|^2} dvdx, \cr 
I^F_3(t)&= \sum_{|\al|=\mathrm{N}+1}\frac{2}{\e^4}\int_{\Omega \times \R^3} \mathcal{N}(\p^{\al}F^{\e},\AC{\P}F^{\e}) \frac{\AC{\P}(\p^{\al}F^{\e})}{M^{\e}} dvdx,  \cr 
I^F_4(t)&=\sum_{|\al|=\mathrm{N}+1}\frac{2}{\e^4}\sum_{0< \beta< \al} \binom{\al}{\beta} \int_{\Omega \times \R^3}\mathcal{N}(\p^{\beta}\AC{\P}F^{\e},\p^{\al-\beta}M^{\e}) \frac{\AC{\P}(\p^{\al}F^{\e})}{M^{\e}} dvdx, \cr 
I^F_5(t)&=\sum_{|\al|=\mathrm{N}+1}\frac{1}{\e^4}\sum_{0< \beta< \al} \binom{\al}{\beta} \int_{\Omega \times \R^3}\mathcal{N}(\p^{\beta}\AC{\P}F^{\e},\p^{\al-\beta}\AC{\P}F^{\e}) \frac{\AC{\P}(\p^{\al}F^{\e})}{M^{\e}} dvdx.
\enda
\end{align}
\end{lemma}
\begin{proof}
We apply $\p^{\al}$ to \eqref{BE} and decompose the collision operator into two parts: one where all derivatives act on a single factor, and the other where the derivatives are distributed between both factors.
\[ 
\p_t \p^{\al}F^{\e}+\frac{v}{\e}\cdot\nabla_x \p^{\al}F^{\e} = \frac{2}{\kappa\e^2}\mathcal{N}(\p^{\al}F^{\e},F^{\e}) + \frac{1}{\kappa\e^2}\sum_{0< \beta< \al}\binom{\al}{\beta}\mathcal{N}(\p^{\beta}F^{\e},\p^{\al-\beta}F^{\e}).
\]
We further decompose the first term on the right-hand side by using the relation $F = M^{\e} + \AC{\P}F^{\e}$, and apply the definition of the linear operator $\mathcal{L}$ from \eqref{Ldef} to get
\begin{align*}
\p_t \p^{\al}F^{\e}+\frac{v}{\e}\cdot\nabla_x \p^{\al}F^{\e} + \frac{1}{\kappa\e^2}\mathcal{L}(\p^{\al}F^{\e}) = \frac{2}{\kappa\e^2}\mathcal{N}(\p^{\al}F^{\e},\AC{\P}F^{\e}) + \frac{1}{\kappa\e^2}\sum_{0< \beta< \al}\binom{\al}{\beta}\mathcal{N}(\p^{\beta}F^{\e},\p^{\al-\beta}F^{\e}).
\end{align*}
Then, testing the above equation with $\frac{\kappa}{\e^2}\p^{\al}F^{\e}|M^{\e}|^{-1}$ yields
\begin{align}\label{FEesti}
\bega
\frac{\kappa}{2\e^2}&\frac{d}{dt} \int_{\Omega \times \R^3}|\p^{\al}F^{\e}|^2 |M^{\e}|^{-1}dvdx  + \frac{1}{\e^4}\int_{\Omega \times \R^3}\mathcal{L}(\p^{\al}F^{\e})\p^{\al}F^{\e}|M^{\e}|^{-1}  dvdx \cr
&= -\frac{\kappa}{2\e^2}\int_{\Omega \times \R^3}  \left(\p_t M^{\e} + \frac{v}{\eps}\cdot\nabla_x M^{\e}\right) |\p^{\al}F^{\e}|^2 |M^{\e}|^{-2} dvdx \cr 
&\quad + \frac{2}{\e^4}\int_{\Omega \times \R^3} \mathcal{N}(\p^{\al}F^{\e},\AC{\P}F^{\e})\AC{\P}(\p^{\al}F^{\e}) |M^{\e}|^{-1}  dvdx  \cr 
&\quad +\frac{1}{\e^4}\sum_{0< \beta< \al} \binom{\al}{\beta} \int_{\Omega \times \R^3}\mathcal{N}(\p^{\beta}F^{\e},\p^{\al-\beta}F^{\e})\AC{\P}(\p^{\al}F^{\e}) |M^{\e}|^{-1}  dvdx,
\enda
\end{align}
where we have used the identity $\int_{\R^3} \mathcal{N}(\cdot,\cdot)(1,v,|v|^2) \, dv = 0$ to the last two terms. Summing over all multi-indices with $|\al| = \mathrm{N}+1$, we find that the second and third lines of \eqref{FEesti} correspond to the terms $I^F_2$ and $I^F_3$, respectively.
For the last line in \eqref{FEesti}, we further decompose the collision operator using the decomposition $F = M^{\e} + \AC{\P}F^{\e}$ as follows: 
\begin{align*}
\bega
\sum_{|\al|=\mathrm{N}+1}\frac{1}{\e^4}\sum_{0< \beta< \al} \binom{\al}{\beta} \int_{\Omega \times \R^3}\mathcal{N}\Big(\p^{\beta}(M^{\e}+\AC{\P}F^{\e}),\p^{\al-\beta}(M^{\e}+\AC{\P}F^{\e})\Big)\AC{\P}(\p^{\al}F^{\e}) |M^{\e}|^{-1} dvdx.
\enda
\end{align*}
When the collision operator includes interactions between $\AC{\P}F^{\e}$ and $M^{\e}$, the resulting contribution is assigned to the term $I^F_4$. When the interaction is purely between $\AC{\P}F^{\e}$ terms, the contribution corresponds to $I^F_5$. When $M^{\e}$ interacts with $M^{\e}$ inside the collision operator, the corresponding term is expressed as
\begin{align*}
\bega
\sum_{|\al|=\mathrm{N}+1}\frac{1}{\e^4}\sum_{0< \beta< \al} \binom{\al}{\beta} \int_{\Omega \times \R^3}\mathcal{N}\Big(\p^{\beta}M^{\e},\p^{\al-\beta}M^{\e}\Big)\AC{\P}(\p^{\al}F^{\e}) |M^{\e}|^{-1} dvdx.
\enda
\end{align*}
Since $\mathcal{N}(M^{\e},M^{\e}) = 0$, we have
\begin{align*}
\sum_{0<\beta<\alpha}\binom{\al}{\beta} \mathcal{N}\Big(\p^{\beta}M^{\e},\p^{\al-\beta}M^{\e}\Big)  =  \p^{\al} \mathcal{N}(M^{\e},M^{\e}) - \mathcal{N}(\p^{\al} M^{\e}, M^{\e} ) - \mathcal{N}( M^{\e}, \p^{\al} M^{\e} ) 
 = \mathcal{L}(\p^{\al} M^{\e}),
\end{align*}
where we have used the definition of \(\mathcal{L}\) given in \eqref{Ldef}. This corresponds to the term \( I^F_1 \).
\end{proof}

\begin{remark}
In contrast to the proof of Lemma \ref{L.EFtop}, one may first decompose
$F^{\e}=M^{\e}+\AC{\P}F^{\e}$ and then apply $\p^{\al}$ to the equation:
\begin{align*}
\p_t F^{\e}+\frac{v}{\e}\cdot\nabla_x F^{\e} + \frac{1}{\kappa\e^2}\mathcal{L}(\AC{\P}F^{\e}) =  \frac{1}{\kappa\e^2}\mathcal{N}(\AC{\P}F^{\e},\AC{\P}F^{\e}).
\end{align*}
In this case, a singular macro--macro interaction appears through the term
$\mathcal{N}(M^{\e},\p^{\al}\AC{\P}F^{\e}) \AC{\P}(\p^{\al}M^{\e})$,
which has the same structure and scaling as the term $I^F_1$.
However, the resulting dissipation $|\p^{\al}\AC{\P}F^{\e}|^2|M^{\e}|^{-1}$ is weaker than the dissipation
$|\AC{\P}\p^{\al}F^{\e}|^2|M^{\e}|^{-1}$ obtained in Lemma~\ref{L.EFtop}.
\end{remark}

Now, we present the estimate showing that the first term of \eqref{E'0} is equivalent to \(\mathcal{E}_{\mathrm{top}}(t)\), as stated in \eqref{E0eq}.

\begin{lemma}\label{L.E0}
For $t\in[0,T]$ satisfying the bootstrap assumption \eqref{condition}, the inequality \eqref{E0eq} holds for both the purely spatial derivatives \eqref{caseA} and the space-time derivatives \eqref{caseB}.
\end{lemma}
\begin{proof}
We define $I^F_L(t)$ as follows:
\begin{align*}
I^F_L(t) := \sum_{|\al|=\mathrm{N}+1} \frac{\kappa}{2\e^2} \int_{\Omega \times \R^3} |\p^{\al}F^{\e}|^2 |M^{\e}|^{-1} \, dv \, dx.
\end{align*}
Using \eqref{palF} in Lemma \ref{PalF}, we have 
\begin{align}\label{E0a}
\bega
I^F_L(t) &=  \sum_{|\al|=\mathrm{N}+1} \frac{\kappa}{\eps^2} \int_{\Omega \times \R^3} \Big(|\P(\p^{\al}M^{\e})|^2 +2\P(\p^{\al}M^{\e}) \AC{\P}(\p^{\al}F^{\e}) +  |\AC{\P}(\p^{\al}F^{\e})|^2 \Big) |M^{\e}|^{-1} \, dv \, dx.
\enda
\end{align}

Since \( \P \) and \( \AC{\P} \) are orthogonal in the sense of \eqref{Pproj}, the middle term in \eqref{E0a} vanishes. We now decompose the first term in \eqref{E0a} using Lemma \ref{Mal}:
\begin{align}\label{Psqdecomp}
\bega
|\P(\p^{\al}M^{\e})|^2 = |\e\Phi_{\al}^1 M^{\e}|^2 + \Big( |\P(\p^{\al}M^{\e})|^2 - |\e\Phi_{\al}^1 M^{\e}|^2 \Big).
\enda
\end{align}
For the first term in \eqref{Psqdecomp}, we have
\begin{align}\label{Philower}
\bega
\int_{\R^3} |\Phi_{\al}^1|^2M^{\e}dv 
&= \int_{\R^3} \frac{|\p^{\al} \rho^{\e}|^2}{|\mathrm{P}^{\e}|^2}M^{\e}+\sum_i\frac{|\p^{\al} u^{\e}_i|^2  |v_i-\mathrm{U}^{\e}_i|^2}{|k_B\mathrm{\Theta}^{\e}|^2}M^{\e}+\frac{|\p^{\al} \ta^{\e}|^2}{|\mathrm{\Theta}^{\e}|^2} \lw( \frac{|v-\mathrm{U}^{\e}|^2}{2k_B\mathrm{\Theta}^{\e}}-\frac{3}{2} \rw)^2M^{\e}dv \cr
&= \lw(\frac{|\p^{\al} \rho^{\e}|^2}{\mathrm{P}^{\e}} + \frac{\mathrm{P}^{\e}|\p^{\al} u^{\e}|^2}{k_B\mathrm{\Theta}^{\e}} + \frac{3}{2} \mathrm{P}^{\e} \frac{|\p^{\al} \ta^{\e}|^2}{|\mathrm{\Theta}^{\e}|^2}\rw).
\enda
\end{align}
For the last term in \eqref{Psqdecomp}, using the identity $a^2 - b^2 = (a + b)(a - b)$, the expansion \eqref{RpM-def}, and Lemma \ref{Pprop}, we obtain
\begin{align}\label{Psqremain}
\bega
&\int_{\Omega \times \R^3} \Big(|\P(\p^{\al}M^{\e})|^2-|\e\Phi_{\al}^1M^{\e}|^2\Big)|M^{\e}|^{-1}dvdx  \cr 
&\les \int_{\Omega \times \R^3} \bigg(2\e\Phi_{\al}^1M^{\e} + \P\bigg( \sum_{2\leq i\leq |\al|}\eps^i \Phi_{\al}^i M^{\e}\bigg)\bigg) \P\bigg( \sum_{2\leq i\leq |\al|}\eps^i \Phi_{\al}^i M^{\e}\bigg) |M^{\e}|^{-1}dvdx  \cr 
&\les \bigg(\int_{\Omega \times \R^3} \bigg|2\e\Phi_{\al}^1M^{\e} + \P\bigg( \sum_{2\leq i\leq |\al|}\eps^i \Phi_{\al}^i M^{\e}\bigg)\bigg|^2|M^{\e}|^{-1}dvdx\bigg)^{\frac{1}{2}}  \cr 
&\quad \times \bigg(\int_{\Omega \times \R^3} \bigg|\P\bigg(\sum_{2\leq i\leq |\al|}\eps^i \Phi_{\al}^i M^{\e}\bigg)\bigg|^2|M^{\e}|^{-1}dvdx\bigg)^{\frac{1}{2}}  \cr
&\les (\e\kappa^{-\frac{1}{2}}\mathcal{E}_{top}^{\frac{1}{2}}+\e^2 \mathcal{E}_M) (\e^2 \mathcal{E}_M) \les \e^3 \kappa^{-\frac{1}{2}}\mathcal{E}_M\mathcal{E}_{tot}^{\frac{1}{2}},
\enda
\end{align}
where we used \eqref{Phiscale}, \eqref{Rscale}, bootstrap assumption \( (\mathcal{B}_1) \) in \eqref{condition}, and the identity $\P (\Phi_{\al}^1 M^{\e}) = \Phi_{\al}^1 M^{\e}$.
Applying \eqref{Philower} and \eqref{Psqremain} to \eqref{E0a}, we conclude
\begin{align*}
\bega
\bigg|I^F_L(t) &-\sum_{|\al|=\mathrm{N}+1}\frac{\kappa}{\eps^2}\int_{\Omega \times \R^3}|\AC{\P}(\p^{\al}F^{\e})|^2|M^{\e}|^{-1}dvdx \cr 
&- \sum_{|\al|=\mathrm{N}+1}\kappa \int_{\Omega} \lw(\frac{|\p^{\al} \rho^{\e}|^2}{\mathrm{P}^{\e}} - \frac{\mathrm{P}^{\e}|\p^{\al} u^{\e}|^2}{k_B\mathrm{\Theta}^{\e}} + \frac{3}{2} \mathrm{P}^{\e} \frac{|\p^{\al} \ta^{\e}|^2}{|\mathrm{\Theta}^{\e}|^2}\rw) dx \bigg| \les \e\kappa^{\frac{1}{2}} \mathcal{E}_{tot}^{\frac{3}{2}}(t).
\enda
\end{align*}
Thus, under the bootstrap assumption $\sup_{t \in [0,T]} \left( | \mathrm{P}^{\e} - 1|, |\mathrm{U}^{\e}|, |\mathrm{\Theta}^{\e} - 1| \right) \ll 1$ from \( (\mathcal{B}_2) \) in \eqref{condition}, we conclude that \( I^F_L(t) \) is equivalent to the energy \( \mathcal{E}_{top}(t) \) defined in \eqref{N-EDdef}, up to a difference of order $\e\kappa^{\frac{1}{2}} \mathcal{E}_{tot}^{\frac{3}{2}}(t)$:
\begin{align*}
\bega
\frac{1}{C} \Big(\mathcal{E}_{top}(t) -\e\kappa^{\frac{1}{2}} \mathcal{E}_{tot}^{\frac{3}{2}}(t)\Big) \leq \big|I^F_L(t)\big| \leq C \Big(\mathcal{E}_{top}(t) +\e\kappa^{\frac{1}{2}} \mathcal{E}_{tot}^{\frac{3}{2}}(t)\Big),
\enda
\end{align*}
for some constant \( C > 0 \). This gives the result \eqref{E0eq}.
\end{proof}

\subsubsection{Nonlinear Macro-Macro Interaction}

Before estimating the term $I^F_1$ defined in \eqref{E123def}, we first state a lemma that establishes the relations between $\e \p_t u^{\e}$ and $\nabla_x (\rho^{\e} + \ta^{\e})$; and between $\e \p_t \ta^{\e}$ and $\nabla_x \cdot u^{\e}$ in suitable norms.

\begin{lemma}\label{L.ept}
For $t\in[0,T]$ satisfying the bootstrap assumption \eqref{condition}, and for $\mathrm{N} > d/2 + 1$, we have 
\begin{align}\label{eptineq}
\bega
\|\e\p_t (u^{\e},\ta^{\e}) \|_{L^\infty_x} &\les \e \|u^{\e}\|_{L^\infty_x} \mathcal{E}_M^{\frac{1}{2}} + \|(\nabla_x\cdot u^{\e},\nabla_x(\rho^{\e}+\ta^{\e}))\|_{L^\infty_x} + \e \|\bPhi^{\e}_{\bW}\|_{L^\infty_x}, \cr
\|\e\p_t (u^{\e},\ta^{\e}) \|_{H^k_x} &\les \e \big(\|u^{\e}\|_{L^\infty_x} \mathcal{E}_M^{\frac{1}{2}} + \mathcal{E}_M\big) + \|(\nabla_x\cdot u^{\e},\nabla_x(\rho^{\e}+\ta^{\e}))\|_{H^k_x} + \e \|\bPhi^{\e}_{\bW}\|_{H^k_x}, 
\enda
\end{align}
for $0\leq k \leq \mathrm{N}-1$. 
\end{lemma}
\begin{proof}
Using the local conservation laws \eqref{locconNew}$_2$ and \eqref{locconNew}$_3$, we get
\begin{align*}
\|\e\p_tu^{\e}\| &\les \|\e u^{\e}\cdot \nabla_x u^{\e}\| +C\|\nabla_x(\rho^{\e}+\ta^{\e})\|
+ \frac{1}{\e} \sum_{j} \| \p_{x_j} \mathbf{r}_{ij}^{\e}\|, \cr 
\|\e\p_t\ta^{\e}\| &\les \|\e u^{\e}\cdot\nabla_x\ta^{\e}\|+ \frac{2}{3}\|\nabla_x\cdot u^{\e}\| + \frac{1}{\e}\sum_{j} \| \p_{x_j} \mathfrak{q}_j^{\e}\| + \frac{1}{\e}\sum_{i,j} \|\p_{x_i}\mathrm{U}^{\e}_j \mathbf{r}_{ij}^{\e}\|,
\end{align*}
for either \( \|\cdot\|_{L^\infty_x} \) or \( \|\cdot\|_{H^k_x} \), where we used  \eqref{pTaleq} and the bootstrap assumption \eqref{condition}.
For the convection term, applying \eqref{uvHk} yields
\begin{align*}
\|u^{\e}\cdot \nabla_x \bW^{\e}\|_{H^{\mathrm{N}-1}_x} &\leq \|u^{\e}\|_{L^\infty_x} \|\nabla_x \bW^{\e}\|_{H^{\mathrm{N}-1}_x} + \|\nabla_x u^{\e}\|_{L^\infty_x}\|\nabla_x \bW^{\e}\|_{H^{\mathrm{N}-2}_x} + \|\nabla_x \bW^{\e}\|_{L^\infty_x}\|\nabla_x u^{\e}\|_{H^{\mathrm{N}-2}_x} \cr 
&\leq \|u^{\e}\|_{L^\infty_x} \mathcal{E}_M^{\frac{1}{2}} + \mathcal{E}_M,
\end{align*}
where we used \( \bW^{\e} = [\rho^{\e}, u^{\e}, \ta^{\e}]^T \).
Moreover,
\begin{align*}
\|u^{\e}\cdot \nabla_x \bW^{\e}\|_{L^\infty_x} &\leq \|u^{\e}\|_{L^\infty_x} \|\nabla_x \bW^{\e}\|_{L^\infty_x} \leq \|u^{\e}\|_{L^\infty_x} \mathcal{E}_M^{\frac{1}{2}}.
\end{align*}
By the definition of $\bPhi^{\e}_{\bW}$ in \eqref{gWdef}, we obtain the desired result.
\end{proof}

\begin{lemma}\label{L.E1}
For $t\in[0,T]$ satisfying the bootstrap assumption \eqref{condition}, depending on the type of derivative—either the purely spatial derivative case \eqref{caseA} or the space--time derivative case \eqref{caseB} with $\mathfrak{n} \geq 1$—the term $I^F_1$ defined in \eqref{E123def} satisfies the following upper bounds:
\begin{align*}
\bega
&I^F_1(t) \leq \|\nabla_x(u^{\e},\ta^{\e})(t)\|_{L^\infty_x} \mathcal{E}_M^{\frac{1}{2}}(t)\mathcal{D}_{top}^{\frac{1}{2}}(t) + \e\mathcal{E}_{M}^{\frac{3}{2}}(t)\mathcal{D}_{top}^{\frac{1}{2}}(t) +\mathcal{Z}_{top}^{time}(t),
\enda
\end{align*}
where $\mathcal{Z}_{top}^{time}(t)$ is defined in \eqref{Bmidef}.
\hide
\begin{align*}
\bega
&I^F_1(t) \leq \|\nabla_x(u^{\e},\ta^{\e})(t)\|_{L^\infty_x} \mathcal{E}_M^{\frac{1}{2}}(t)\mathcal{D}_{top}^{\frac{1}{2}}(t) + \e\mathcal{E}_{M}^{\frac{3}{2}}(t)\mathcal{D}_{top}^{\frac{1}{2}}(t) \cr 
&+\begin{cases} 0,  &\mbox{for} \quad \eqref{caseA}, \\ 
\e^{\mathfrak{n}-1} \Big( \|((\nabla_x\cdot u^{\e}),\nabla_x(\rho^{\e}+\ta^{\e}))(t)\|_{L^\infty_x} +\e\|u^{\e}(t)\|_{L^\infty_x}\mathcal{E}_{M}^{\frac{1}{2}}(t) + \e \kappa^{\frac{1}{2}}\mathcal{D}_G^{\frac{1}{2}}(t) \Big)\mathcal{E}_{M}^{\frac{1}{2}}(t)\mathcal{D}_{top}^{\frac{1}{2}}(t) ,  &\mbox{for} \quad \eqref{caseB}.
\end{cases} 
\enda
\end{align*}
\unhide
\end{lemma}

\begin{remark}\label{Rmk.scale}
The result of Lemma~\ref{L.E1} depends on the temporal scaling. 
If the scaled time derivative $\e^{\mathfrak{n}} \p_t$ is included in $\p^{\al}$ for some $\mathfrak{n} \geq 1$, then applying the local conservation laws in \eqref{locconNew} yields an additional scaling factor of $\e^{\mathfrak{n}-1}$. 
On the other hand, if the scaling $\e^{\mathfrak{n}} \p_t$ is used with $0 \leq \mathfrak{n} < 1$, the local conservation laws do not produce any improvement in scaling. 
In this case, we must estimate the time derivative directly using the energy bound in \eqref{zzL2}, namely, $\|\e^{\mathfrak{n}} \p_t u^{\e}\|_{L^\infty_x} \leq \mathcal{E}_M^{1/2}$. 
The remaining term to be controlled is then $\|\nabla_x(u^{\e}, \ta^{\e})\|_{H^{\mathrm{N}-1}_x}$, which leads to a different type of upper bound for the case $0 \leq \mathfrak{n} < 1$.
\end{remark}
\begin{proof}[Proof of Lemma \ref{L.E1}] Applying the H\"{o}lder inequality and \eqref{L-1L2} to $I^F_1$ yields
\begin{align*}
\bega
I^F_{3}
&\leq \sum_{|\al|=\mathrm{N}+1}\frac{1}{\e^4} \bigg(\int_{\Omega \times \R^3} \Big(\mathcal{L}(\p^{\al}M^{\e})\Big)^2|M^{\e}|^{-1} dvdx\bigg)^{\frac{1}{2}} \bigg(\int_{\Omega \times \R^3}|\AC{\P}(\p^{\al}F^{\e})|^2 |M^{\e}|^{-1} dvdx\bigg)^{\frac{1}{2}} \cr 
&\leq \sum_{|\al|=\mathrm{N}+1}\frac{1}{\e^2} \bigg(\int_{\Omega \times \R^3} (1+|v|^2)\Big(\AC{\P}(\p^{\al}M^{\e})\Big)^2|M^{\e}|^{-1} dvdx\bigg)^{\frac{1}{2}} \mathcal{D}_{top}^{\frac{1}{2}}
,
\enda
\end{align*}
where we used the dissipation $\mathcal{D}_{top}$ defined in \eqref{N-EDdef2}. 
Recall the expansion $\p^{\al} M^{\e} = \sum_{1\leq i\leq |\al|}\eps^i \Phi_{\al}^i M^{\e}$ from \eqref{RpM-def}. Note that $\AC{\P}(\eps \Phi_{\al}^1 M^{\e}) = 0$, and the $\rho^{\e}$-dependent part of $\Phi_{\al}^2 M^{\e}$ lies in the span of $(1, v, |v|^2)$, so the $\rho^{\e}$-contribution in $\AC{\P}(\Phi_{\al}^2 M^{\e})$ vanishes due to definitions \eqref{Phi1def} and \eqref{Phi2def}. Thus, we have
\begin{align*}
\bega
\big|\AC{\P}(\p^{\al}M^{\e}) \big| &\les \e^2\sum_{0<\beta<\al} (|\p^{\beta}u^{\e} \p^{\al-\beta}u^{\e}| +|\p^{\beta}u^{\e} \p^{\al-\beta}\ta^{\e}| + |\p^{\beta}\ta^{\e} \p^{\al-\beta}\ta^{\e}|)(1+|v|^4)M^{\e} \cr 
&+\sum_{i=3}^{|\al|}\e^i  \sum_{\substack{\al_1+\cdots+\al_i=\al \\ \al_i>0}}|\p^{\al_1}(\rho^{\e},u^{\e},\ta^{\e})|\times \cdots \times |\p^{\al_i}(\rho^{\e},u^{\e},\ta^{\e})| (1+|v|^{2i}) M^{\e},
\enda
\end{align*}
where we also used \eqref{Phiinf0} in the second line. 
For the terms involving three or more factors of 
$\p^{\al_i}(\rho^{\e},u^{\e},\ta^{\e})$, 
the maximal derivative index is bounded by $|\al_i|\leq \mathrm{N}-1$. 
Using the bound  $\|\p^{\al_i}(\rho^{\e},u^{\e},\ta^{\e})\|_{L^2_x} \leq \mathcal{E}_{M}^{\frac{1}{2}}$ for $|\al_i|=\mathrm{N}-1$ from \eqref{Phiinfscale}, and arguing as in the proof of \eqref{Phiscale} in Appendix \ref{A.B}, we obtain
\begin{align}\label{I3Fbdd}
\bega
I^F_{3}&\les \sum_{|\al|=\mathrm{N}+1} \bigg(\int_{\Omega} \sum_{0<\beta<\al}\Big| |\p^\beta u^{\e}| |\p^{\al-\beta} u^{\e}| + |\p^\beta u^{\e}| |\p^{\al-\beta} \ta^{\e}| + |\p^\beta \ta^{\e}| |\p^{\al-\beta} \ta^{\e}| \Big|^2 dx \bigg)^{\frac{1}{2}} \mathcal{D}_{top}^{\frac{1}{2}} \cr 
&+ \frac{1}{\e^2}\sum_{i=3}^{|\al|} \e^i \mathcal{E}_M^{\frac{i}{2}} \mathcal{D}_{top}^{\frac{1}{2}}.
\enda
\end{align}
Let $\zeta_1, \zeta_2 \in \{u^{\e}, \ta^{\e} \}$, and consider estimating $\|\p^{\beta}\zeta_1 \p^{\al-\beta} \zeta_2\|_{L^2_x}$. Let $\p^{\bullet}$ denote either $\p_{\tilde{t}} = \e^{\mathfrak{n}} \p_t$ or $\nabla_x$. Since $|\beta| \geq 1$ and $|\al - \beta| \geq 1$, we have  
\begin{align}\label{zzL2}
\bega
\sum_{|\al|=\mathrm{N}+1}  \sum_{0<\beta<\al} \left\| |\partial^\beta \zeta_1| |\partial^{\al-\beta} \zeta_2| \right\|_{L^2_x} 
&\les \|\nabla_x\zeta_1 \p^{\bullet}\zeta_2\|_{H^{\mathrm{N}-1}_x} \cr 
&\les \|\nabla_x\zeta_1\|_{L^\infty_x} \|\p^{\bullet}\zeta_2\|_{H^{\mathrm{N}-1}_x} + \|\nabla_x\zeta_1\|_{H^{\mathrm{N}-1}_x} \|\p^{\bullet}\zeta_2\|_{L^\infty_x} \cr 
&\les \|\nabla_x(u^{\e},\ta^{\e})\|_{L^\infty_x} \mathcal{E}_M^{\frac{1}{2}} +  \|\p^{\bullet}\zeta_2\|_{L^\infty_x}\mathcal{E}_M^{\frac{1}{2}} ,
\enda
\end{align}
where we used \eqref{uvHk} in the second line and $\|\p^{\bullet}\zeta\|_{H^{\mathrm{N}-1}_x} \leq \mathcal{E}_M^{\frac{1}{2}}$ in the last line.\\
(Case A): When there are no time derivatives, i.e., $\p^{\bullet} = \nabla_x$, it is clear that the second term in \eqref{zzL2} is also bounded by $\|\nabla_x(u^{\e}, \ta^{\e})\|_{L^\infty_x} \mathcal{E}_M^{\frac{1}{2}}$. \\ 
(Case B): In the case when $\p^{\bullet} = \p_{\tilde{t}} = \e^{\mathfrak{n}} \p_t$, we claim that
\begin{align}\label{MMAB}
\bega
\|\p_{\tilde{t}}\zeta_2\|_{L^\infty_x} \leq 
\e^{\mathfrak{n}-1} \Big( (\|\nabla_x\cdot u^{\e}\|_{L^\infty_x}+\|\nabla_x(\rho^{\e}+\ta^{\e})\|_{L^\infty_x}) +\e\|u^{\e}\|_{L^\infty_x}\mathcal{E}_{M}^{\frac{1}{2}} + \e \kappa^{\frac{1}{2}}\mathcal{D}_G^{\frac{1}{2}} \Big).
\enda
\end{align}
Using the temporal scaling $\e^{\mathfrak{n}}$, we apply Lemma \ref{L.ept} to estimate $\e \p_t(u^{\e}, \ta^{\e})$ in the $L^\infty_x$ norm:
\begin{align*}
\bega
\|\e\p_t u^{\e} \|_{L^\infty_x} &\les \e \|u^{\e}\|_{L^\infty_x} \mathcal{E}_M^{\frac{1}{2}} + \|\nabla_x(\rho^{\e}+\ta^{\e})\|_{L^\infty_x} + \e \|\bPhi^{\e}_{\bW}\|_{L^\infty_x}, \cr
\|\e\p_t \ta^{\e} \|_{L^\infty_x} &\les \e \|u^{\e}\|_{L^\infty_x} \mathcal{E}_M^{\frac{1}{2}} + \|\nabla_x\cdot u^{\e}\|_{L^\infty_x} + \e \|\bPhi^{\e}_{\bW}\|_{L^\infty_x}.
\enda
\end{align*}
Using $\|\bPhi^{\e}_{\bW}\|_{L^\infty_x} \leq C\kappa^{\frac{1}{2}}\mathcal{D}_G^{\frac{1}{2}}$ from \eqref{ABGscale}$_2$, we obtain the bound \eqref{MMAB} for any $\zeta_2 \in \{u^{\e}, \ta^{\e}\}$. Combining \eqref{zzL2} and \eqref{MMAB} with \eqref{I3Fbdd} completes the estimate.
\end{proof}

\subsubsection{Contribution of Convection Against the Local Maxwellian} 

In this section, we aim to prove the estimate for $I^F_2$ defined in \eqref{E123def}, which arises from the commutator $\llbracket\p_t +\frac{v}{\eps} \cdot \nabla_x, \frac{1}{M^\e} \rrbracket$. 
Before estimating $I^F_2$, we first analyze the commutator between \( \p^{\al} \) and \( \AC{\P} \) in order to control the dissipation arising from this commutator.

\begin{lemma}
For both types of derivatives \( \p^{\al} \) defined in \eqref{caseA} and \eqref{caseB}, and for \( F^{\e} = M^{\e} + \AC{\P}F^{\e} \), the following commutation relations hold:
\begin{align}\label{pP}
\bega
\p^{\al} \P F^{\e} - \P (\p^{\al}F^{\e}) &= \begin{cases}
0, \quad &\text{for} \quad |\al| = 0,1, \\
\AC{\P} \left(  \sum_{2\leq i\leq |\al|}\eps^i \Phi_{\al}^i M^{\e} \right), \quad &\text{for} \quad |\al| \geq 2,
\end{cases}
\enda
\end{align}
and
\begin{align}\label{pACP}
\p^{\al} \AC{\P}F^{\e} - \AC{\P} (\p^{\al}F^{\e}) &= \begin{cases}
0, \quad &\text{for} \quad |\al| = 0,1, \\
-\AC{\P} \lw( \sum_{2\leq i\leq |\al|}\eps^i \Phi_{\al}^i M^{\e} \rw), \quad &\text{for} \quad |\al| \geq 2.
\end{cases}
\end{align}
\end{lemma}
\begin{proof}
To prove \eqref{pP}, we use the identity \( \P F^{\e} = M^{\e} \) and apply Lemma \ref{PalF} to get
\begin{align*}
\bega
\p^{\al}\P F^{\e}-\P(\p^{\al}F^{\e})= \p^{\al}M^{\e}-\P(\p^{\al}M^{\e}) &= \AC{\P}(\p^{\al}M^{\e}) = \AC{\P}\sum_{2\leq i\leq |\al|}\eps^i \Phi_{\al}^i M^{\e},
\enda
\end{align*}
where we used \( \AC{\P}(\Phi_{\al}^1 M^{\e}) = 0 \) and the definition of \( \sum_{2\leq i\leq |\al|}\eps^i \Phi_{\al}^i M^{\e} \) in \eqref{RpM-def}. This completes the proof of \eqref{pP}.
Next, using the definition of \( \AC{\P} \) in \eqref{Pdef}, we compute
\begin{align*}
\p^{\al}\AC{\P}F^{\e} = \p^{\al}(\mathbf{I}-\P)F^{\e}
= \p^{\al}F^{\e} - \p^{\al}\P F^{\e} 
= \AC{\P}(\p^{\al}F^{\e}) +\lw(\P(\p^{\al}F^{\e})- \p^{\al}\P F^{\e} \rw) .
\end{align*}
Applying \eqref{pP} to the last term yields \eqref{pACP}.
\end{proof}

\begin{lemma}\label{L.E2}
For $t\in[0,T]$ satisfying the bootstrap assumption \eqref{condition}, and for both the purely spatial derivatives \eqref{caseA} and the space-time derivatives \eqref{caseB}, the term $I^F_2$ defined in \eqref{E123def} satisfies the following upper bound:
\begin{align}\label{E1est}
\bega
I^F_2(t) &\les \big(\|\nabla_xu^{\e}(t)\|_{L^\infty_x}+\|\nabla_x\ta^{\e}(t)\|_{L^\infty_x}\big) \mathcal{E}_{tot}(t) \cr 
&+ \int_{\Omega}|\nabla_xu^{\e}(t)| \mathcal{V}_2(t) dx + \int_{\Omega}|\nabla_x\ta^{\e}(t)| \mathcal{V}_3(t) dx + \frac{1}{\eps} \int_{\Omega} \TbT(t) \mathcal{V}_2(t) dx \cr 
&+ \e\kappa^{\frac{1}{2}}\mathcal{D}_G^{\frac{1}{2}}(t) \Big(\mathcal{E}_{tot}(t) +\e\kappa^{\frac{1}{2}}\mathcal{E}_{top}^{\frac{1}{2}}(t)\mathcal{D}_{top}^{\frac{1}{2}}(t) \Big)  + \e\kappa^{\frac{1}{2}}  \mathcal{E}_{tot}(t)\mathcal{D}_{top}^{\frac{1}{2}}(t) ,
\enda
\end{align}
where \( \mathcal{V}_\ell(t,x) \) and the control of the turbulence transport $\TbT(t,x)$ are defined in \eqref{largev} and \eqref{RSdef}, respectively.
\end{lemma}

\begin{proof}
We prove Lemma \ref{L.E2} in three steps.
In the first step, we decompose $I^F_2$ into the macroscopic contribution and the microscopic projection. 
In the second step, we estimate $I^F_{2,1}$. 
In the third step, we estimate $I^F_{2,2}$, $I^F_{2,3}$, and $I^F_{2,4}$. \\
\noindent (Step 1) We first claim that $I^F_2$ can be decomposed as follows $I^F_2 = I^F_{2,1} + I^F_{2,2} + I^F_{2,3} + I^F_{2,4}$:
\begin{align}\label{E1idef}
\bega
I^F_{2,1}(t)=& -\sum_{|\al|=\mathrm{N}+1}\frac{\kappa}{2}\int_{\Omega \times \R^3} |\Phi_{\al}^1|^2  \left(\p_t M^{\e} + \frac{v}{\eps}\cdot\nabla_x M^{\e}\right) dvdx,  \cr
I^F_{2,2}(t)= &-\sum_{|\al|=\mathrm{N}+1}\frac{\kappa}{2\eps^2}\int_{\Omega \times \R^3} |\AC{\P}(\p^{\al}F^{\e})|^2|M^{\e}|^{-2}   \left(\p_t M^{\e} + \frac{v}{\eps}\cdot\nabla_x M^{\e}\right) dvdx,  \cr
I^F_{2,3}(t)=& \sum_{|\al|=\mathrm{N}+1}\frac{\kappa}{2\eps^2}\int_{\Omega \times \R^3} \Big(|\P(\p^{\al}M^{\e})|^2-|\e\Phi_{\al}^1M^{\e}|^2\Big) |M^{\e}|^{-2}   \left(\p_t M^{\e} + \frac{v}{\eps}\cdot\nabla_x M^{\e}\right) dvdx,  \cr 
I^F_{2,4}(t)=&-\sum_{|\al|=\mathrm{N}+1}\frac{\kappa}{\eps^2}\int_{\Omega \times \R^3} \P(\p^{\al}M^{\e})\AC{\P}(\p^{\al}F^{\e})    |M^{\e}|^{-2}\left(\p_t M^{\e} + \frac{v}{\eps}\cdot\nabla_x M^{\e}\right)dvdx. 
\enda
\end{align}
Here, recall that $\left( \p_t M^{\e} + \frac{v}{\eps} \cdot \nabla_x M^{\e} \right)$ is given by \eqref{Phitx}. The term $I^F_{2,1}$ satisfies the following form:
\begin{align}\label{E11equal}
\bega
&I^F_{2,1}(t)=
\sum_{|\al|=\mathrm{N}+1}\frac{\kappa}{2} \int_{\Omega} \frac{2}{\e k_B\mathrm{\Theta}^{\e}} \lw( \frac{\p^{\al} \rho^{\e}}{\mathrm{P}^{\e}}  + \frac{\p^{\al} \ta^{\e}}{\mathrm{\Theta}^{\e}}\rw) \sum_{i,j} \p_{x_j} \mathbf{r}_{ij}^{\e} \p^{\al} u^{\e}_i dx \cr 
&+\sum_{|\al|=\mathrm{N}+1}\frac{\kappa}{2} \int_{\Omega}\frac{2}{\eps3k_B\mathrm{\Theta}^{\e}} \lw(\frac{1}{k_B\mathrm{\Theta}^{\e}} |\p^{\al} u^{\e}|^2 + 3\lw(\frac{\p^{\al} \rho^{\e} }{\mathrm{P}^{\e}}\frac{\p^{\al} \ta^{\e}}{\mathrm{\Theta}^{\e}}+\frac{|\p^{\al} \ta^{\e}|^2}{|\mathrm{\Theta}^{\e}|^2}\rw) \rw) \bigg(\sum_{j} \p_{x_j} \mathfrak{q}_j^{\e} + \sum_{i,j} \p_{x_i}\mathrm{U}^{\e}_j \mathbf{r}_{ij}^{\e} \bigg) dx  \cr 
&-\sum_{|\al|=\mathrm{N}+1}\frac{\kappa}{2} \int_{\Omega}\bigg[ \frac{\mathrm{P}^{\e}}{k_B\mathrm{\Theta}^{\e}}\bigg(-\frac{2}{3}(\nabla_x\cdot u^{\e})|\p^{\al} u^{\e}|^2 + 2\sum_{i,j}\p_{x_i}u^{\e}_j \p^{\al} u^{\e}_i\p^{\al} u^{\e}_j \bigg)  +\frac{5\mathrm{P}^{\e}}{|\mathrm{\Theta}^{\e}|^2}\nabla_x\ta^{\e} \cdot \p^{\al} u^{\e} \p^{\al} \ta^{\e} \bigg] dx.
\enda
\end{align}
(Proof of \eqref{E1idef})
Using \( \p^{\al}F^{\e} = \P(\p^{\al} M^{\e}) + \AC{\P}(\p^{\al}F^{\e}) \) from Lemma \ref{PalF}, we can express \( I^F_2 \) as
\begin{align}\label{E1sp}
\bega
I^F_2 &= -\sum_{|\al|=\mathrm{N}+1}\frac{\kappa}{2\eps^2}\int_{\Omega \times \R^3} |\P(\p^{\al} M^{\e})|^2 |M^{\e}|^{-2} \left(\p_t M^{\e} + \frac{v}{\eps}\cdot\nabla_x M^{\e}\right) dvdx \cr
&\quad - \sum_{|\al|=\mathrm{N}+1}\frac{\kappa}{\eps^2}\int_{\Omega \times \R^3} \P(\p^{\al} M^{\e}) \AC{\P}(\p^{\al}F^{\e}) |M^{\e}|^{-2} \left(\p_t M^{\e} + \frac{v}{\eps}\cdot\nabla_x M^{\e}\right) dvdx \cr
&\quad - \sum_{|\al|=\mathrm{N}+1}\frac{\kappa}{2\eps^2}\int_{\Omega \times \R^3} |\AC{\P}(\p^{\al}F^{\e})|^2 |M^{\e}|^{-2} \left(\p_t M^{\e} + \frac{v}{\eps}\cdot\nabla_x M^{\e}\right) dvdx.
\enda
\end{align}
The second and third lines of \eqref{E1sp} correspond to $I^F_{2,4}$ and $I^F_{2,2}$, respectively.
In the first line of \eqref{E1sp}, applying the decomposition of \( |\P(\p^{\al} M^{\e})|^2 \) from \eqref{Psqdecomp} yields the terms $I^F_{2,1}$ and $I^F_{2,3}$. \\
(Proof of \eqref{E11equal}) By the definition of \( \Phi_{\al}^1 \) in \eqref{Phi1def}, we expand \( |\Phi_\al^1|^2 \) explicitly as
\begin{align*}
\bega
|\Phi_\al^1|^2&= \frac{|\p^{\al} \rho^{\e}|^2}{|\mathrm{P}^{\e}|^2}+\sum_{k,l}\frac{\p^{\al} u^{\e}_k\p^{\al} u^{\e}_l}{k_B^2|\mathrm{\Theta}^{\e}|^2}(v_k-\mathrm{U}^{\e}_k)(v_l-\mathrm{U}^{\e}_l) +\frac{|\p^{\al} \ta^{\e}|^2}{4k_B^2|\mathrm{\Theta}^{\e}|^4} \lw(|v-\mathrm{U}^{\e}|^2-3k_B\mathrm{\Theta}^{\e}\rw)^2  \cr 
&+2\frac{\p^{\al} \rho^{\e} \p^{\al} u^{\e} }{\mathrm{P}^{\e}k_B\mathrm{\Theta}^{\e}}\cdot (v-\mathrm{U}^{\e}) + \frac{\p^{\al} \rho^{\e} \p^{\al} \ta^{\e}}{\mathrm{P}^{\e}k_B|\mathrm{\Theta}^{\e}|^2}\lw(|v-\mathrm{U}^{\e}|^2-3k_B\mathrm{\Theta}^{\e}\rw) \cr 
&+ \frac{\p^{\al} u^{\e} \p^{\al} \ta^{\e}}{k_B^2|\mathrm{\Theta}^{\e}|^3}\cdot(v-\mathrm{U}^{\e}) \lw(|v-\mathrm{U}^{\e}|^2-3k_B\mathrm{\Theta}^{\e}\rw).
\enda
\end{align*}
Since the term \( \p_t M^{\e} + \frac{v}{\eps} \cdot \nabla_x M^{\e} \) has the explicit form given in \eqref{Phitx}, we only need to compute the following moments of \( |\Phi_{\al}^1|^2 M^{\e} \):
\begin{align}\label{Phi-moment}
\bega
&\int_{\R^3} (v_i-\mathrm{U}^{\e}_i)|\Phi_{\al}^1|^2M^{\e}dv = 2\mathrm{P}^{\e}\p^{\al} u^{\e}_i \lw( \frac{\p^{\al} \rho^{\e}}{\mathrm{P}^{\e}}  + \frac{\p^{\al} \ta^{\e}}{\mathrm{\Theta}^{\e}}\rw),
\cr
&\int_{\R^3} \lw(\frac{|v-\mathrm{U}^{\e}|^2}{2k_B\mathrm{\Theta}^{\e}}-\frac{3}{2}\rw) |\Phi_{\al}^1|^2M^{\e}dv = \frac{\mathrm{P}^{\e}}{k_B\mathrm{\Theta}^{\e}} |\p^{\al} u^{\e}|^2 + 3\mathrm{P}^{\e}\lw(\frac{\p^{\al} \rho^{\e} }{\mathrm{P}^{\e}}\frac{\p^{\al} \ta^{\e}}{\mathrm{\Theta}^{\e}}+\frac{|\p^{\al} \ta^{\e}|^2}{|\mathrm{\Theta}^{\e}|^2}\rw), \cr 
&\int_{\R^3} \mathfrak{R}^{\e}_{ij} |\Phi_{\al}^1|^2M^{\e}dv = -\frac{2}{3}\delta_{ij}\lw(\mathrm{P}^{\e} |\p^{\al} u^{\e}|^2 \rw) + 2\mathrm{P}^{\e} \p^{\al} u^{\e}_i\p^{\al} u^{\e}_j, \cr 
&\int_{\R^3} \mathcal{Q}^{\e}_i |\Phi_{\al}^1|^2M^{\e}dv = 5k_B\mathrm{P}^{\e} \p^{\al} u^{\e}_i \p^{\al} \ta^{\e}.
\enda
\end{align}
In particular, for the third line of \eqref{Phi-moment}, we used 
\begin{align*}
&\int_{\R^3} (v_i-\mathrm{U}^{\e}_i)(v_j-\mathrm{U}^{\e}_j) |\Phi_{\al}^1|^2M^{\e}dv = \int_{\R^3} (v_i-\mathrm{U}^{\e}_i)(v_j-\mathrm{U}^{\e}_j)M^{\e} \cr 
&\times \bigg[ \frac{|\p^{\al} \rho^{\e}|^2}{|\mathrm{P}^{\e}|^2}\delta_{ij} +\frac{|\p^{\al} \ta^{\e}|^2}{4k_B^2|\mathrm{\Theta}^{\e}|^4} \lw(|v-\mathrm{U}^{\e}|^2-3k_B\mathrm{\Theta}^{\e}\rw)^2\delta_{ij} + \frac{\p^{\al} \rho^{\e} \p^{\al} \ta^{\e}}{\mathrm{P}^{\e}k_B|\mathrm{\Theta}^{\e}|^2}\lw(|v-\mathrm{U}^{\e}|^2-3k_B\mathrm{\Theta}^{\e}\rw)\delta_{ij} \cr 
&+\sum_{k,l}\frac{\p^{\al} u^{\e}_k\p^{\al} u^{\e}_l}{k_B^2|\mathrm{\Theta}^{\e}|^2}(v_k-\mathrm{U}^{\e}_k)(v_l-\mathrm{U}^{\e}_l)(\delta_{i=j \neq k=l}+\delta_{i=k\neq j=l}+\delta_{i=l\neq j=k}+\delta_{i=j=k=l})\bigg] dv,
\end{align*}
so that
\begin{align*}
&\int_{\R^3} (v_i-\mathrm{U}^{\e}_i)(v_j-\mathrm{U}^{\e}_j) |\Phi_{\al}^1|^2M^{\e}dv =  2\mathrm{P}^{\e} (1-\delta_{ij}) \p^{\al} u^{\e}_i\p^{\al} u^{\e}_j \cr &+\delta_{ij}\lw(\frac{k_B\mathrm{\Theta}^{\e}}{\mathrm{P}^{\e}} |\p^{\al} \rho^{\e}|^2 + \sum_{i\neq k}\mathrm{P}^{\e} |\p^{\al} u^{\e}_k|^2 + 3\mathrm{P}^{\e} |\p^{\al} u^{\e}_i|^2 +\frac{7}{2} \frac{k_B\mathrm{P}^{\e}}{\mathrm{\Theta}^{\e}} |\p^{\al} \ta^{\e}|^2  + 2k_B \p^{\al} \rho^{\e} \p^{\al}\ta^{\e} \rw) \cr 
&= \delta_{ij}\lw(k_B \mathrm{P}^{\e} \mathrm{\Theta}^{\e} \lw( \frac{|\p^{\al} \rho^{\e}|^2}{|\mathrm{P}^{\e}|^2} +2\frac{\p^{\al} \rho^{\e}}{\mathrm{P}^{\e}}\frac{\p^{\al}\ta^{\e}}{\mathrm{\Theta}^{\e}} +\frac{7}{2} \frac{|\p^{\al} \ta^{\e}|^2}{|\mathrm{\Theta}^{\e}|^2}\rw) + \mathrm{P}^{\e} |\p^{\al} u^{\e}|^2 \rw) + 2\mathrm{P}^{\e} \p^{\al} u^{\e}_i\p^{\al} u^{\e}_j.
\end{align*}
\hide
First line of \eqref{vPhi2}
\begin{align}
\bega
&\int_{\R^3} (v_i-\mathrm{U}^{\e}_i)|\Phi_{\al}^1|^2M^{\e}dv \cr 
&=  \int_{\R^3}(v_i-\mathrm{U}^{\e}_i) \bigg[2\frac{\p^{\al} \rho^{\e} \p^{\al} u^{\e}_i }{\mathrm{P}^{\e}k_B\mathrm{\Theta}^{\e}} (v_i-\mathrm{U}^{\e}_i) + \frac{\p^{\al} u^{\e}_i \p^{\al} \ta^{\e}}{k_B^2|\mathrm{\Theta}^{\e}|^3}(v_i-\mathrm{U}^{\e}_i) \lw(|v-\mathrm{U}^{\e}|^2-3k_B\mathrm{\Theta}^{\e}\rw) \bigg] M^{\e} dv \cr 
&= \bigg[2\frac{\p^{\al} \rho^{\e} \p^{\al} u^{\e}_i }{\mathrm{P}^{\e}k_B\mathrm{\Theta}^{\e}}k_B \mathrm{P}^{\e}\mathrm{\Theta}^{\e} + \frac{\p^{\al} u^{\e}_i \p^{\al} \ta^{\e}}{k_B^2|\mathrm{\Theta}^{\e}|^3}\lw(5k_B^2 \mathrm{P}^{\e}|\mathrm{\Theta}^{\e}|^2-3k_B\mathrm{\Theta}^{\e}(k_B \mathrm{P}^{\e}\mathrm{\Theta}^{\e})\rw) \bigg] \cr 
&= 2\p^{\al} \rho^{\e} \p^{\al} u^{\e}_i + \frac{2\mathrm{P}^{\e}}{\mathrm{\Theta}^{\e}}\p^{\al} u^{\e}_i \p^{\al} \ta^{\e} \cr 
&= 2\mathrm{P}^{\e}\p^{\al} u^{\e}_i \lw( \frac{\p^{\al} \rho^{\e}}{\mathrm{P}^{\e}}  + \frac{\p^{\al} \ta^{\e}}{\mathrm{\Theta}^{\e}}\rw)
\enda
\end{align}
Second line of \eqref{vPhi2}
\begin{align}
\bega
&\int_{\R^3} \lw(\frac{|v-\mathrm{U}^{\e}|^2}{2k_B\mathrm{\Theta}^{\e}}-\frac{3}{2}\rw) |\Phi_{\al}^1|^2M^{\e}dv \cr 
&=  \int_{\R^3}\lw(\frac{|v-\mathrm{U}^{\e}|^2}{2k_B\mathrm{\Theta}^{\e}}-\frac{3}{2}\rw) \bigg[\frac{|\p^{\al} \rho^{\e}|^2}{|\mathrm{P}^{\e}|^2}+\sum_{k}\frac{|\p^{\al} u^{\e}_k|^2}{k_B^2|\mathrm{\Theta}^{\e}|^2}(v_k-\mathrm{U}^{\e}_k)^2 +\frac{|\p^{\al} \ta^{\e}|^2}{4k_B^2|\mathrm{\Theta}^{\e}|^4} \lw(|v-\mathrm{U}^{\e}|^2-3k_B\mathrm{\Theta}^{\e}\rw)^2  \cr 
&+ \frac{\p^{\al} \rho^{\e} \p^{\al} \ta^{\e}}{\mathrm{P}^{\e}k_B|\mathrm{\Theta}^{\e}|^2}\lw(|v-\mathrm{U}^{\e}|^2-3k_B\mathrm{\Theta}^{\e}\rw) \bigg] M^{\e} dv \cr 
&=\sum_{k}\frac{|\p^{\al} u^{\e}_k|^2}{k_B^2|\mathrm{\Theta}^{\e}|^2}\lw(\frac{5k_B^2 \mathrm{P}^{\e}|\mathrm{\Theta}^{\e}|^2}{2k_B\mathrm{\Theta}^{\e}}-\frac{3}{2}k_B \mathrm{P}^{\e}\mathrm{\Theta}^{\e}\rw) \cr
&+\frac{|\p^{\al} \ta^{\e}|^2}{4k_B^2|\mathrm{\Theta}^{\e}|^4} \bigg(\lw(\frac{105 k_B^3  \mathrm{P}^{\e} |\mathrm{\Theta}^{\e}|^3-6k_B\mathrm{\Theta}^{\e}(15k_B^2 \mathrm{P}^{\e}|\mathrm{\Theta}^{\e}|^2) + 9k_B^2|\mathrm{\Theta}^{\e}|^2(3k_B \mathrm{P}^{\e}\mathrm{\Theta}^{\e})}{2k_B\mathrm{\Theta}^{\e}}\rw)\cr 
&\quad -\frac{3}{2}\lw(15k_B^2 \mathrm{P}^{\e}|\mathrm{\Theta}^{\e}|^2-6k_B\mathrm{\Theta}^{\e}(3k_B \mathrm{P}^{\e}\mathrm{\Theta}^{\e}) + 9k_B^2|\mathrm{\Theta}^{\e}|^2(\mathrm{P}^{\e})\rw) \bigg) \cr 
&+ \frac{\p^{\al} \rho^{\e} \p^{\al} \ta^{\e}}{\mathrm{P}^{\e}k_B|\mathrm{\Theta}^{\e}|^2}\lw(\frac{15k_B^2 \mathrm{P}^{\e}|\mathrm{\Theta}^{\e}|^2}{2k_B\mathrm{\Theta}^{\e}}-\frac{3}{2}3k_B \mathrm{P}^{\e}\mathrm{\Theta}^{\e} \rw) \cr 
&=\frac{\mathrm{P}^{\e}}{k_B\mathrm{\Theta}^{\e}} \sum_{k}|\p^{\al} u^{\e}_k|^2 + \frac{3\mathrm{P}^{\e}}{|\mathrm{\Theta}^{\e}|^2}|\p^{\al} \ta^{\e}|^2 + \frac{3}{\mathrm{\Theta}^{\e}} \p^{\al} \rho^{\e} \p^{\al} \ta^{\e} \cr 
&= \frac{\mathrm{P}^{\e}}{k_B\mathrm{\Theta}^{\e}} |\p^{\al} u^{\e}|^2 + 3\mathrm{P}^{\e}\lw(\frac{\p^{\al} \rho^{\e} }{\mathrm{P}^{\e}}\frac{\p^{\al} \ta^{\e}}{\mathrm{\Theta}^{\e}}+\frac{|\p^{\al} \ta^{\e}|^2}{|\mathrm{\Theta}^{\e}|^2}\rw)
\enda
\end{align}
Third line of \eqref{vPhi2} - first term:
\begin{align}
\bega
&\int_{\R^3} (v_i-\mathrm{U}^{\e}_i)(v_j-\mathrm{U}^{\e}_j) |\Phi_{\al}^1|^2M^{\e}dv \cr 
&=  \int_{\R^3}(v_i-\mathrm{U}^{\e}_i)(v_j-\mathrm{U}^{\e}_j) \bigg[ \frac{|\p^{\al} \rho^{\e}|^2}{|\mathrm{P}^{\e}|^2}+\sum_{k,l}\frac{\p^{\al} u^{\e}_k\p^{\al} u^{\e}_l}{k_B^2|\mathrm{\Theta}^{\e}|^2}(v_k-\mathrm{U}^{\e}_k)(v_l-\mathrm{U}^{\e}_l) +\frac{|\p^{\al} \ta^{\e}|^2}{4k_B^2|\mathrm{\Theta}^{\e}|^4} \lw(|v-\mathrm{U}^{\e}|^2-3k_B\mathrm{\Theta}^{\e}\rw)^2  \cr 
& + \frac{\p^{\al} \rho^{\e} \p^{\al} \ta^{\e}}{\mathrm{P}^{\e}k_B|\mathrm{\Theta}^{\e}|^2}\lw(|v-\mathrm{U}^{\e}|^2-3k_B\mathrm{\Theta}^{\e}\rw) \bigg] M^{\e} dv \cr 
\enda
\end{align}

\begin{align}
&\int_{\R^3} (v_i-\mathrm{U}^{\e}_i)(v_j-\mathrm{U}^{\e}_j) |\Phi_{\al}^1|^2M^{\e}dv \cr 
&= \int_{\R^3} (v_i-\mathrm{U}^{\e}_i)(v_j-\mathrm{U}^{\e}_j) \bigg[ \frac{|\p^{\al} \rho^{\e}|^2}{|\mathrm{P}^{\e}|^2}\delta_{ij}+\sum_{k,l}\frac{\p^{\al} u^{\e}_k\p^{\al} u^{\e}_l}{k_B^2|\mathrm{\Theta}^{\e}|^2}(v_k-\mathrm{U}^{\e}_k)(v_l-\mathrm{U}^{\e}_l)(\delta_{i=j \neq k=l}+\delta_{i=k\neq j=l}+\delta_{i=l\neq j=k}+\delta_{i=j=k=l}) \cr 
&+\frac{|\p^{\al} \ta^{\e}|^2}{4k_B^2|\mathrm{\Theta}^{\e}|^4} \lw(|v-\mathrm{U}^{\e}|^2-3k_B\mathrm{\Theta}^{\e}\rw)^2\delta_{ij}  + \frac{\p^{\al} \rho^{\e} \p^{\al} \ta^{\e}}{\mathrm{P}^{\e}k_B|\mathrm{\Theta}^{\e}|^2}\lw(|v-\mathrm{U}^{\e}|^2-3k_B\mathrm{\Theta}^{\e}\rw)\delta_{ij} \bigg] M^{\e} dv \cr 
&= \frac{|\p^{\al} \rho^{\e}|^2}{|\mathrm{P}^{\e}|^2}(k_B \mathrm{P}^{\e}\mathrm{\Theta}^{\e})\delta_{ij}\cr 
&+\frac{1}{k_B^2|\mathrm{\Theta}^{\e}|^2}\sum_{i\neq k}\int_{\R^3}(v_i-\mathrm{U}^{\e}_i)^2(v_k-\mathrm{U}^{\e}_k)^2M^{\e} dv  |\p^{\al} u^{\e}_k|^2\delta_{ij}\cr 
&+\frac{1}{k_B^2|\mathrm{\Theta}^{\e}|^2}(1-\delta_{ij})\int_{\R^3}(v_i-\mathrm{U}^{\e}_i)^2(v_j-\mathrm{U}^{\e}_j)^2M^{\e} dv \p^{\al} u^{\e}_i\p^{\al} u^{\e}_j \cr 
&+\frac{1}{k_B^2|\mathrm{\Theta}^{\e}|^2}(1-\delta_{ij})\int_{\R^3}(v_i-\mathrm{U}^{\e}_i)^2(v_j-\mathrm{U}^{\e}_j)^2M^{\e} dv \p^{\al} u^{\e}_j\p^{\al} u^{\e}_i \cr 
&+\frac{1}{k_B^2|\mathrm{\Theta}^{\e}|^2}\int_{\R^3}(v_i-\mathrm{U}^{\e}_i)^4 M^{\e} dv |\p^{\al} u^{\e}_i|^2 \delta_{ij}  \cr 
&+\frac{|\p^{\al} \ta^{\e}|^2}{4k_B^2|\mathrm{\Theta}^{\e}|^4} \frac{1}{3}\lw(105 k_B^3  \mathrm{P}^{\e} |\mathrm{\Theta}^{\e}|^3-6k_B\mathrm{\Theta}^{\e}(15k_B^2 \mathrm{P}^{\e}|\mathrm{\Theta}^{\e}|^2) + 9k_B^2|\mathrm{\Theta}^{\e}|^2(3k_B \mathrm{P}^{\e}\mathrm{\Theta}^{\e})\rw)\delta_{ij} \cr
& + \frac{\p^{\al} \rho^{\e} \p^{\al} \ta^{\e}}{\mathrm{P}^{\e}k_B|\mathrm{\Theta}^{\e}|^2}\frac{1}{3}\lw(15k_B^2 \mathrm{P}^{\e}|\mathrm{\Theta}^{\e}|^2-3k_B\mathrm{\Theta}^{\e}(3k_B \mathrm{P}^{\e}\mathrm{\Theta}^{\e})\rw)\delta_{ij} 
\end{align}
($\delta_{ij}$ yields $\nabla_x \cdot u^{\e}$)
Arrange it: 
\begin{align}
&\int_{\R^3} (v_i-\mathrm{U}^{\e}_i)(v_j-\mathrm{U}^{\e}_j) |\Phi_{\al}^1|^2M^{\e}dv \cr 
&=\frac{k_B\mathrm{\Theta}^{\e}}{\mathrm{P}^{\e}} |\p^{\al} \rho^{\e}|^2 \delta_{ij} \cr 
&+\frac{1}{k_B^2|\mathrm{\Theta}^{\e}|^2}\sum_{i\neq k}k_B^2\mathrm{P}^{\e} |\mathrm{\Theta}^{\e}|^2 |\p^{\al} u^{\e}_k|^2\delta_{ij}\cr 
&+\frac{1}{k_B^2|\mathrm{\Theta}^{\e}|^2}k_B^2\mathrm{P}^{\e} |\mathrm{\Theta}^{\e}|^2 (1-\delta_{ij})\p^{\al} u^{\e}_i\p^{\al} u^{\e}_j \cr 
&+\frac{1}{k_B^2|\mathrm{\Theta}^{\e}|^2}k_B^2\mathrm{P}^{\e} |\mathrm{\Theta}^{\e}|^2 (1-\delta_{ij})\p^{\al} u^{\e}_j\p^{\al} u^{\e}_i \cr 
&+\frac{1}{k_B^2|\mathrm{\Theta}^{\e}|^2}3k_B^2\mathrm{P}^{\e} |\mathrm{\Theta}^{\e}|^2 |\p^{\al} u^{\e}_i|^2 \delta_{ij} \cr
&+\frac{7}{2} \frac{k_B\mathrm{P}^{\e}}{\mathrm{\Theta}^{\e}} |\p^{\al} \ta^{\e}|^2 \delta_{ij}  + 2k_B \p^{\al} \rho^{\e} \p^{\al}\ta^{\e} \delta_{ij} \cr 
&= \delta_{ij}\lw(\frac{k_B\mathrm{\Theta}^{\e}}{\mathrm{P}^{\e}} |\p^{\al} \rho^{\e}|^2 + \sum_{i\neq k}\mathrm{P}^{\e} |\p^{\al} u^{\e}_k|^2 + 3\mathrm{P}^{\e} |\p^{\al} u^{\e}_i|^2 +\frac{7}{2} \frac{k_B\mathrm{P}^{\e}}{\mathrm{\Theta}^{\e}} |\p^{\al} \ta^{\e}|^2  + 2k_B \p^{\al} \rho^{\e} \p^{\al}\ta^{\e} \rw) + 2\mathrm{P}^{\e} (1-\delta_{ij}) \p^{\al} u^{\e}_i\p^{\al} u^{\e}_j \cr 
&= \delta_{ij}\lw(k_B \mathrm{P}^{\e} \mathrm{\Theta}^{\e} \lw( \frac{|\p^{\al} \rho^{\e}|^2}{|\mathrm{P}^{\e}|^2} +2\frac{\p^{\al} \rho^{\e}}{\mathrm{P}^{\e}}\frac{\p^{\al}\ta^{\e}}{\mathrm{\Theta}^{\e}} +\frac{7}{2} \frac{|\p^{\al} \ta^{\e}|^2}{|\mathrm{\Theta}^{\e}|^2}\rw) + \sum_{i\neq k}\mathrm{P}^{\e} |\p^{\al} u^{\e}_k|^2 + \mathrm{P}^{\e} |\p^{\al} u^{\e}_i|^2 \rw) + 2\mathrm{P}^{\e} \p^{\al} u^{\e}_i\p^{\al} u^{\e}_j \cr 
&= \delta_{ij}\lw(k_B \mathrm{P}^{\e} \mathrm{\Theta}^{\e} \lw( \frac{|\p^{\al} \rho^{\e}|^2}{|\mathrm{P}^{\e}|^2} +2\frac{\p^{\al} \rho^{\e}}{\mathrm{P}^{\e}}\frac{\p^{\al}\ta^{\e}}{\mathrm{\Theta}^{\e}} +\frac{7}{2} \frac{|\p^{\al} \ta^{\e}|^2}{|\mathrm{\Theta}^{\e}|^2}\rw) + \mathrm{P}^{\e} |\p^{\al} u^{\e}|^2 \rw) + 2\mathrm{P}^{\e} \p^{\al} u^{\e}_i\p^{\al} u^{\e}_j
\end{align}

Third line of \eqref{vPhi2} - second term:
\begin{align}
\bega
&\int_{\R^3} |v-\mathrm{U}^{\e} |^2 |\Phi_{\al}^1|^2M^{\e}dv \cr 
&=  \int_{\R^3}|v-\mathrm{U}^{\e} |^2 \bigg[ \frac{|\p^{\al} \rho^{\e}|^2}{|\mathrm{P}^{\e}|^2}+\sum_{k,l}\frac{\p^{\al} u^{\e}_k\p^{\al} u^{\e}_l}{k_B^2|\mathrm{\Theta}^{\e}|^2}(v_k-\mathrm{U}^{\e}_k)(v_l-\mathrm{U}^{\e}_l) +\frac{|\p^{\al} \ta^{\e}|^2}{4k_B^2|\mathrm{\Theta}^{\e}|^4} \lw(|v-\mathrm{U}^{\e}|^2-3k_B\mathrm{\Theta}^{\e}\rw)^2  \cr 
&+2\frac{\p^{\al} \rho^{\e} \p^{\al} u^{\e} }{\mathrm{P}^{\e}k_B\mathrm{\Theta}^{\e}}\cdot (v-\mathrm{U}^{\e}) + \frac{\p^{\al} \rho^{\e} \p^{\al} \ta^{\e}}{\mathrm{P}^{\e}k_B|\mathrm{\Theta}^{\e}|^2}\lw(|v-\mathrm{U}^{\e}|^2-3k_B\mathrm{\Theta}^{\e}\rw) + \frac{\p^{\al} u^{\e} \p^{\al} \ta^{\e}}{k_B^2|\mathrm{\Theta}^{\e}|^3}\cdot(v-\mathrm{U}^{\e}) \lw(|v-\mathrm{U}^{\e}|^2-3k_B\mathrm{\Theta}^{\e}\rw) \bigg] M^{\e} dv \cr 
&=  \int_{\R^3} |v-\mathrm{U}^{\e}|^2 \bigg[\frac{|\p^{\al} \rho^{\e}|^2}{|\mathrm{P}^{\e}|^2}+\sum_{k}\frac{|\p^{\al} u^{\e}_k|^2}{k_B^2|\mathrm{\Theta}^{\e}|^2}(v_k-\mathrm{U}^{\e}_k)^2 +\frac{|\p^{\al} \ta^{\e}|^2}{4k_B^2|\mathrm{\Theta}^{\e}|^4} \lw(|v-\mathrm{U}^{\e}|^2-3k_B\mathrm{\Theta}^{\e}\rw)^2  \cr 
&+ \frac{\p^{\al} \rho^{\e} \p^{\al} \ta^{\e}}{\mathrm{P}^{\e}k_B|\mathrm{\Theta}^{\e}|^2}\lw(|v-\mathrm{U}^{\e}|^2-3k_B\mathrm{\Theta}^{\e}\rw) \bigg] M^{\e} dv \cr 
&= \frac{|\p^{\al} \rho^{\e}|^2}{|\mathrm{P}^{\e}|^2}(3k_B \mathrm{P}^{\e}\mathrm{\Theta}^{\e}) + \sum_{k}\frac{|\p^{\al} u^{\e}_k|^2}{k_B^2|\mathrm{\Theta}^{\e}|^2}(5k_B^2 \mathrm{P}^{\e}|\mathrm{\Theta}^{\e}|^2)\cr
&+\frac{|\p^{\al} \ta^{\e}|^2}{4k_B^2|\mathrm{\Theta}^{\e}|^4} (105 k_B^3  \mathrm{P}^{\e} |\mathrm{\Theta}^{\e}|^3-6k_B\mathrm{\Theta}^{\e}(15k_B^2 \mathrm{P}^{\e}|\mathrm{\Theta}^{\e}|^2) + 9k_B^2|\mathrm{\Theta}^{\e}|^2(3k_B \mathrm{P}^{\e}\mathrm{\Theta}^{\e}))\cr 
&+ \frac{\p^{\al} \rho^{\e} \p^{\al} \ta^{\e}}{\mathrm{P}^{\e}k_B|\mathrm{\Theta}^{\e}|^2}\lw(15k_B^2 \mathrm{P}^{\e}|\mathrm{\Theta}^{\e}|^2-3k_B\mathrm{\Theta}^{\e}(3k_B \mathrm{P}^{\e}\mathrm{\Theta}^{\e}) \rw) \cr 
&=3k_B \mathrm{P}^{\e}\mathrm{\Theta}^{\e} \frac{|\p^{\al} \rho^{\e}|^2}{|\mathrm{P}^{\e}|^2} + 5\mathrm{P}^{\e} |\p^{\al} u^{\e}|^2 +\frac{21}{2} \frac{k_B\mathrm{P}^{\e}}{\mathrm{\Theta}^{\e}}|\p^{\al} \ta^{\e}|^2 + 6k_B \p^{\al} \rho^{\e} \p^{\al} \ta^{\e} \cr 
&= 3k_B \mathrm{P}^{\e}\mathrm{\Theta}^{\e} \lw(\frac{|\p^{\al} \rho^{\e}|^2}{|\mathrm{P}^{\e}|^2} + 2\frac{\p^{\al} \rho^{\e}}{\mathrm{P}^{\e}} \frac{\p^{\al} \ta^{\e}}{\mathrm{\Theta}^{\e}} +\frac{7}{2} \frac{|\p^{\al} \ta^{\e}|^2}{|\mathrm{\Theta}^{\e}|^2} \rw) + 5\mathrm{P}^{\e} |\p^{\al} u^{\e}|^2
\enda
\end{align}
Combining two part of the third line of \eqref{vPhi2}:
\begin{align}
\bega
&\int_{\R^3} \lw((v_i-\mathrm{U}^{\e}_i)(v_j-\mathrm{U}^{\e}_j)-\frac{|v-\mathrm{U}^{\e} |^2}{3}\delta_{ij}\rw) |\Phi_{\al}^1|^2M^{\e}dv \cr 
&= \delta_{ij}\lw(k_B \mathrm{P}^{\e} \mathrm{\Theta}^{\e} \lw( \frac{|\p^{\al} \rho^{\e}|^2}{|\mathrm{P}^{\e}|^2} +2\frac{\p^{\al} \rho^{\e}}{\mathrm{P}^{\e}}\frac{\p^{\al}\ta^{\e}}{\mathrm{\Theta}^{\e}} +\frac{7}{2} \frac{|\p^{\al} \ta^{\e}|^2}{|\mathrm{\Theta}^{\e}|^2}\rw) + \mathrm{P}^{\e} |\p^{\al} u^{\e}|^2 \rw) + 2\mathrm{P}^{\e} \p^{\al} u^{\e}_i\p^{\al} u^{\e}_j \cr 
&-\frac{1}{3}\delta_{ij}\lw(3k_B \mathrm{P}^{\e}\mathrm{\Theta}^{\e} \lw(\frac{|\p^{\al} \rho^{\e}|^2}{|\mathrm{P}^{\e}|^2} + 2\frac{\p^{\al} \rho^{\e}}{\mathrm{P}^{\e}} \frac{\p^{\al} \ta^{\e}}{\mathrm{\Theta}^{\e}} +\frac{7}{2} \frac{|\p^{\al} \ta^{\e}|^2}{|\mathrm{\Theta}^{\e}|^2} \rw) + 5\mathrm{P}^{\e} |\p^{\al} u^{\e}|^2\rw) \cr 
&= -\frac{2}{3}\delta_{ij}\lw(\mathrm{P}^{\e} |\p^{\al} u^{\e}|^2 \rw) + 2\mathrm{P}^{\e} \p^{\al} u^{\e}_i\p^{\al} u^{\e}_j
\enda
\end{align}
Fourth line of \eqref{vPhi2}
\begin{align}
\bega
&\int_{\R^3} \lw((v_i-\mathrm{U}^{\e}_i)\frac{(|v-\mathrm{U}^{\e}|^2-5k_B\mathrm{\Theta}^{\e})}{2}\rw) |\Phi_{\al}^1|^2M^{\e}dv \cr 
&=  \int_{\R^3}\lw((v_i-\mathrm{U}^{\e}_i)\frac{(|v-\mathrm{U}^{\e}|^2-5k_B\mathrm{\Theta}^{\e})}{2}\rw) \bigg[ 2\frac{\p^{\al} \rho^{\e} \p^{\al} u^{\e}_i }{\mathrm{P}^{\e}k_B\mathrm{\Theta}^{\e}} (v_i-\mathrm{U}^{\e}_i) + \frac{\p^{\al} u^{\e}_i \p^{\al} \ta^{\e}}{k_B^2|\mathrm{\Theta}^{\e}|^3}(v_i-\mathrm{U}^{\e}_i) \lw(|v-\mathrm{U}^{\e}|^2-3k_B\mathrm{\Theta}^{\e}\rw) \bigg] M^{\e} dv \cr 
&= \bigg[\frac{\p^{\al} \rho^{\e} \p^{\al} u^{\e}_i }{\mathrm{P}^{\e}k_B\mathrm{\Theta}^{\e}} \lw( 5k_B^2 \mathrm{P}^{\e}|\mathrm{\Theta}^{\e}|^2-5k_B\mathrm{\Theta}^{\e}(k_B \mathrm{P}^{\e}\mathrm{\Theta}^{\e})\rw) + \frac{\p^{\al} u^{\e}_i \p^{\al} \ta^{\e}}{k_B^2|\mathrm{\Theta}^{\e}|^3}\int_{\R^3}\frac{1}{2}(v_i-\mathrm{U}^{\e}_i)^2(|v-\mathrm{U}^{\e}|^4-8k_B\mathrm{\Theta}^{\e}|v-\mathrm{U}^{\e}|^2+15k_B^2|\mathrm{\Theta}^{\e}|^2)M^{\e} dv\bigg]  \cr 
&= \frac{\p^{\al} u^{\e}_i \p^{\al} \ta^{\e}}{k_B^2|\mathrm{\Theta}^{\e}|^3}\int_{\R^3}\frac{1}{2}(v_i-\mathrm{U}^{\e}_i)^2(|v-\mathrm{U}^{\e}|^4-8k_B\mathrm{\Theta}^{\e}|v-\mathrm{U}^{\e}|^2+15k_B^2|\mathrm{\Theta}^{\e}|^2) M^{\e} dv \cr 
&= \frac{\p^{\al} u^{\e}_i \p^{\al} \ta^{\e}}{k_B^2|\mathrm{\Theta}^{\e}|^3}\frac{1}{2}\frac{1}{3}(105 k_B^3  \mathrm{P}^{\e} |\mathrm{\Theta}^{\e}|^3-8k_B\mathrm{\Theta}^{\e}(15k_B^2 \mathrm{P}^{\e}|\mathrm{\Theta}^{\e}|^2)+15k_B^2|\mathrm{\Theta}^{\e}|^2(3k_B \mathrm{P}^{\e}\mathrm{\Theta}^{\e})) \cr 
&= \frac{\p^{\al} u^{\e}_i \p^{\al} \ta^{\e}}{k_B^2|\mathrm{\Theta}^{\e}|^3}\frac{30}{6}k_B^3 \mathrm{P}^{\e} |\mathrm{\Theta}^{\e}|^3 \cr 
&= 5k_B\mathrm{P}^{\e} \p^{\al} u^{\e}_i \p^{\al} \ta^{\e}
\enda
\end{align}
\unhide
We substitute \eqref{Phitx} for \( \p_t M^{\e} + \frac{v}{\eps}\cdot\nabla_x M^{\e} \) into the quantity \( \int_{\R^3} |\Phi_{\al}^1|^2  \left(\p_t M^{\e} + \frac{v}{\eps}\cdot\nabla_x M^{\e}\right) dv \), and then, using \eqref{Phi-moment}, we derive the desired result \eqref{E11equal}.
\\ (Step 2) We claim the following estimate. 
\begin{align}\label{I2,1est}
\bega
I^F_{2,1}(t)   & \les \left( \|\nabla_xu^{\e}(t)\|_{L^\infty_x} + \| \nabla_x \ta^{\e}(t)\|_{L^\infty_x} + \frac{1}{\e}\|\TbT(t)\|_{L^\infty_x} \right) \mathcal{E}_{top}(t) \,  .
\enda
\end{align}
In particular, when summing over \( |\al_x| = 1 \) for spatial derivatives only (i.e., excluding time derivatives) in the 2-dimensional case, we have
\begin{align}\label{E11,al1}
\bega
&\sum_{|\al_x|=1}\int_{\Omega \times \R^3} |\Phi_{\al_x}^1|^2  \left(\p_t M^{\e} + \frac{v}{\eps}\cdot\nabla_x M^{\e}\right) dvdx  \cr 
&\les \left( \|\nabla_x\cdot u^{\e}(t)\|_{L^\infty_x}+\| \nabla_x\ta^{\e}(t)\|_{L^\infty_x} +\frac{1}{\e}\|\TbT(t)\|_{L^\infty_x}\right)\mathcal{E}_{M}(t)+ \|\nabla_x \mathbb{P}u^{\e}(t)\|_{L^2_x}^2\|\nabla_x\mathbb{P}^{\perp}u^{\e}(t)\|_{L^\infty_x} \cr 
& +\|\nabla_x\mathbb{P}u^{\e}(t)\|_{L^2_x}\|\nabla_x\mathbb{P}^{\perp}u^{\e}(t)\|_{L^4_x}^2 + \|\nabla_x \mathbb{P}^{\perp}u^{\e}(t)\|_{L^3_x}^3 ,
\enda
\end{align}
where $\mathbb{P}$ is defined in \eqref{LerayPdef}. \\
(Proof of \eqref{I2,1est})
From the expression for \( I^F_{2,1} \) in \eqref{E11equal}, we apply $(\mathcal{B}_2)$ from \eqref{condition} and use Hölder's inequality to obtain
\begin{align}\label{E11,1}
\bega
&I^F_{2,1} \les 
\sum_{|\al|=\mathrm{N}+1}\kappa (\|\p^{\al}\rho^{\e}\|_{L^2_x}+\|\p^{\al}\ta^{\e}\|_{L^2_x})\|\p^{\al}u^{\e}\|_{L^2_x} \frac{1}{\e} \sum_{i,j} \|\p_{x_j} \mathbf{r}_{ij}^{\e} \|_{L^\infty_x} \cr 
&+\sum_{|\al|=\mathrm{N}+1}\kappa \lw(\|\p^{\al}\rho^{\e}\|_{L^2_x}^2+\|\p^{\al}u^{\e}\|_{L^2_x}^2+\|\p^{\al}\ta^{\e}\|_{L^2_x}^2\rw) \frac{1}{\e}\bigg(\sum_{j} \|\p_{x_j} \mathfrak{q}_j^{\e}\|_{L^\infty_x} + \sum_{i,j} \|\p_{x_i}\mathrm{U}^{\e}_j \mathbf{r}_{ij}^{\e}\|_{L^\infty_x}\bigg) \cr 
&+\sum_{|\al|=\mathrm{N}+1}\kappa  \bigg(\|\nabla_x\cdot u^{\e}\|_{L^\infty_x}\|\p^{\al} u^{\e}\|_{L^2_x}^2 +\| \nabla_x\ta^{\e}\|_{L^\infty_x} \|\p^{\al} u^{\e}\|_{L^2_x} \|\p^{\al} \ta^{\e} \|_{L^2_x} + \int_{\Omega}\Big|\sum_{i,j}\p_{x_i}u^{\e}_j \p^{\al} u^{\e}_i\p^{\al} u^{\e}_j\Big| dx\bigg).
\enda
\end{align}
Using \eqref{rutscale} and the definition of \( \TbT \) in \eqref{RSdef}, we have
\begin{align*}
\bega
I^F_{2,1} &\les 
\Big(\|\nabla_xu^{\e}\|_{L^\infty_x}+\| \nabla_x\ta^{\e}\|_{L^\infty_x}+\frac{1}{\e}\|\TbT\|_{L^\infty_x}\Big)\mathcal{E}_{top} .
\enda
\end{align*}
(Proof of \eqref{E11,al1}) 
We focus on the term \( \sum_{i,j} \p_i u^{\e}_j \p^{\al} u^{\e}_i \p^{\al} u^{\e}_j \) in the third line of \eqref{E11,1}. After summing over \( |\al_x| = 1 \) in dimension \( d=2 \), we obtain
\begin{align}\label{uuubdd}
\bega
\sum_{|\al_x|=1}\sum_{i,j}&\p_{x_i}u^{\e}_j \p^{\al_x} u^{\e}_i\p^{\al_x} u^{\e}_j  = \sum_{1\leq i,j,k \leq 2}\p_{x_i}u^{\e}_j \p_k u^{\e}_i\p_k u^{\e}_j 
\cr 
&= (\p_1u^{\e}_1)^3+(\p_2u^{\e}_2)^3 +  (\nabla_x \cdot u^{\e})\big((\p_1u^{\e}_2)^2+\p_1u^{\e}_2\p_2u^{\e}_1+(\p_2u^{\e}_1)^2\big).
\enda
\end{align}
Using the Hodge decomposition \( u^{\e} = \mathbb{P}u^{\e} + \mathbb{P}^{\perp}u^{\e} \) from \eqref{Hodge}, we get
\begin{align}\label{uuubdd2}
\bega
(\p_1u^{\e}_1)^3+(\p_2u^{\e}_2)^3 &= (\p_1\mathbb{P}u^{\e}_1 + \p_1\mathbb{P}^{\perp}u^{\e}_1)^3+(\p_2\mathbb{P}u^{\e}_2 + \p_2\mathbb{P}^{\perp}u^{\e}_2)^3 \cr 
&= 3(\p_1\mathbb{P}u^{\e}_1)^2(\nabla_x\cdot \mathbb{P}^{\perp}u^{\e}) +3(\p_1\mathbb{P}u^{\e}_1)((\p_1\mathbb{P}^{\perp}u^{\e}_1)^2-(\p_2\mathbb{P}^{\perp}u^{\e}_2)^2) \cr 
&+ (\p_1\mathbb{P}^{\perp}u^{\e}_1)^3+(\p_2\mathbb{P}^{\perp}u^{\e}_2)^3,
\enda
\end{align}
where we used $\nabla_x\cdot \mathbb{P}u^{\e}=0$.
Applying \eqref{uuubdd} and \eqref{uuubdd2}, we can estimate the last term of \eqref{E11,1} as follows:
\begin{align*}
\bega
\sum_{|\al_x|=1} \int_{\Omega}\Big|\sum_{i,j}\p_{x_i}u^{\e}_j \p^{\al} u^{\e}_i\p^{\al} u^{\e}_j\Big| dx &\les \|\nabla_x\cdot u^{\e}\|_{L^\infty_x}\|\nabla_x u^{\e}\|_{L^2_x}^2 + \|\nabla_x \mathbb{P}u^{\e}\|_{L^2_x}^2\|\nabla_x\mathbb{P}^{\perp}u^{\e}\|_{L^\infty_x} \cr 
& +\|\nabla_x\mathbb{P}u^{\e}\|_{L^2_x}\|\nabla_x\mathbb{P}^{\perp}u^{\e}\|_{L^4_x}^2 + \|\nabla_x \mathbb{P}^{\perp}u^{\e}\|_{L^3_x}^3.
\enda
\end{align*}
Since the remaining terms are already estimated in \eqref{E11,1}, we conclude the proof of \eqref{E11,al1}.
\\ (Step 3) 
We claim the following estimates.
\begin{align}\label{I2,234}
\bega
I^F_{2,2}(t) 
&\les  \int_{\Omega}\sum_{i,j}|\p_{x_i}u^{\e}_j(t)| \mathcal{V}_{2}(t) dx + \int_{\Omega}\sum_i |\p_{x_i}\ta^{\e}(t)| \mathcal{V}_{3}(t) dx + \frac{1}{\eps} \int_{\Omega} \TbT(t) \mathcal{V}_{2}(t) dx \cr 
& \quad + \e^2\kappa \big(\|\nabla_xu^{\e}\|_{L^\infty_x}+\|\nabla_x\ta^{\e}\|_{L^\infty_x} \big)\mathcal{E}_M^2(t) + \e\kappa \|\TbT(t)\|_{L^\infty_x} \mathcal{E}_M^2(t), \cr
I^F_{2,3}(t)  &\les \e \kappa^{\frac{1}{2}}  \lw(\|\nabla_xu^{\e}(t)\|_{L^\infty_x}+\|\nabla_x\ta^{\e}(t)\|_{L^\infty_x} +\frac{1}{\e}\|\TbT(t)\|_{L^\infty_x}\rw) \mathcal{E}_M(t)\mathcal{E}_{tot}^{\frac{1}{2}}(t) ,  \cr
I^F_{2,4}(t)  &\les \e\kappa^{\frac{1}{2}}  \lw(\|\nabla_xu^{\e}(t)\|_{L^\infty_x}+\|\nabla_x\ta^{\e}(t)\|_{L^\infty_x} +\frac{1}{\e}\|\TbT(t)\|_{L^\infty_x}\rw) \mathcal{E}_{tot}^{\frac{1}{2}}(t)\mathcal{D}_{top}^{\frac{1}{2}}(t).
\enda
\end{align}
(Estimate of $I^F_{2,2}$)
Applying \eqref{Phitx} to \( I^F_{2,2} \), we obtain
\begin{align*}
\bega
&I^F_{2,2}=- \sum_{|\al|=\mathrm{N}+1}\frac{\kappa}{2\eps^2}\int_{\Omega}\frac{1}{k_B\mathrm{\Theta}^{\e}}\sum_{i,j}\p_{x_i}u^{\e}_j \int_{\R^3} \frac{|\AC{\P}(\p^{\al}F^{\e})|^2}{M^{\e}}\mathfrak{R}^{\e}_{ij}dvdx \cr
&-\sum_{|\al|=\mathrm{N}+1}\frac{\kappa}{2\eps^2}\int_{\Omega}\frac{1}{k_B|\mathrm{\Theta}^{\e}|^2}\sum_i \p_{x_i}\ta^{\e} \int_{\R^3} \frac{|\AC{\P}(\p^{\al}F^{\e})|^2}{M^{\e}}\mathcal{Q}^{\e}_idvdx \cr 
&+\sum_{|\al|=\mathrm{N}+1}\frac{\kappa}{2\eps^2}\int_{\Omega}\frac{2}{\eps3k_B\mathrm{P}^{\e}\mathrm{\Theta}^{\e}} \bigg(\sum_{j} \p_{x_j} \mathfrak{q}_j^{\e} + \sum_{i,j} \p_{x_i}\mathrm{U}^{\e}_j \mathbf{r}_{ij}^{\e} \bigg) \int_{\R^3} \frac{|\AC{\P}(\p^{\al}F^{\e})|^2}{M^{\e}}\lw( \frac{|v-\mathrm{U}^{\e}|^2}{2k_B\mathrm{\Theta}^{\e}}-\frac{3}{2} \rw)dvdx \cr 
&+\sum_{|\al|=\mathrm{N}+1}\frac{\kappa}{2\eps^2}\int_{\Omega} \frac{1}{\e k_B\mathrm{P}^{\e}\mathrm{\Theta}^{\e}}\sum_{i,j} \p_{x_j} \mathbf{r}_{ij}^{\e} \int_{\R^3} \frac{|\AC{\P}(\p^{\al}F^{\e})|^2}{M^{\e}}(v_i-\mathrm{U}^{\e}_i)dvdx.
\enda
\end{align*}
Using condition \( (\mathcal{B}_2) \) in \eqref{condition}, the macroscopic fields are bounded \( \sup_{t \in [0,T]} \big(|\mathrm{P}^{\e} - 1|, |\mathrm{U}^{\e}|, |\mathrm{\Theta}^{\e} - 1|\big) \ll 1 \). 
In addition, we have the pointwise bounds \( |\mathfrak{R}^{\e}_{ij}| \leq C(1+|v|^2) \) and \( |\mathcal{Q}^{\e}_i| \leq C(1+|v|^3) \). 
We first take the \( L^\infty_x \) norm for \( \sum_{j} \p_{x_j} \mathfrak{q}_j^{\e} + \sum_{i,j} \p_{x_i} \mathrm{U}^{\e}_j \mathbf{r}_{ij}^{\e} \), and \( \sum_{i,j} \p_{x_j} \mathbf{r}_{ij}^{\e} \). 
Then, using \eqref{pACP}, we apply the inequality $|\AC{\P}(\p^{\al}F^{\e})|^2 \les |\p^{\al}\AC{\P}F^{\e}|^2 + |\AC{\P} ( \sum_{2\leq i\leq |\al|}\eps^i \Phi_{\al}^i M^{\e} )|^2 $ 
to each of the four lines above. Finally, by applying \eqref{Rscale} along with the definitions of \( \mathcal{V}_{\ell} \) and \( \TbT \) from \eqref{largev} and \eqref{RSdef}, we conclude the desired estimate.
\\ (Estimate of $I^F_{2,3}$)
We first observe that, by taking the \( L^\infty_{x,v} \) norm of \eqref{Phitx} and using the bootstrap assumption 
\( \sup_{t \in [0,T]} \big(|\mathrm{P}^{\e} - 1|, |\mathrm{U}^{\e}|, |\mathrm{\Theta}^{\e} - 1|\big) \ll 1 \) from \eqref{condition}, together with the bound 
\( \| \la v \ra^n |M^{\e}|^{1-n_0} \|_{L^\infty_{x,v}} \leq C \) for any \( n \geq 0 \) and \( 0 \leq n_0 < 1 \), we obtain
\begin{align}\label{Phitxinf}
\bega
\bigg\|\left(\p_t M^{\e} + \frac{v}{\eps}\cdot\nabla_x M^{\e}\right) &|M^{\e}|^{-n_0}(t)\bigg\|_{L^\infty_{x,v}} \les \|\nabla_xu^{\e}(t)\|_{L^\infty_x} + \|\nabla_x\ta^{\e}(t)\|_{L^\infty_x} + \frac{1}{\e}\|\TbT(t)\|_{L^\infty_x}.
\enda
\end{align}
By following the same approach as in \eqref{Psqremain}, and using \eqref{Phitxinf}, we obtain
\begin{align*}
I^F_{2,3}&\leq \sum_{|\al|=\mathrm{N}+1}\frac{\kappa}{2\eps^2}\int_{\Omega \times \R^3} \Big(|\P(\p^{\al}M^{\e})|^2-|\e\Phi_{\al}^1M^{\e}|^2\Big) |M^{\e}|^{-\frac{3}{2}} dvdx \bigg\|\frac{\left(\p_t M^{\e} + \frac{v}{\eps}\cdot\nabla_x M^{\e}\right)}{|M^{\e}|^{\frac{1}{2}}} \bigg\|_{L^\infty_{x,v}}\cr 
&\les \sum_{|\al|=\mathrm{N}+1}\frac{\kappa}{2\eps^2} \big(\e^3 \kappa^{-\frac{1}{2}}\mathcal{E}_M\mathcal{E}_{tot}^{\frac{1}{2}}\big)\bigg[\|\nabla_xu^{\e}\|_{L^\infty_x}+\|\nabla_x\ta^{\e}\|_{L^\infty_x} +\frac{1}{\e}\|\TbT\|_{L^\infty_x}\bigg],
\end{align*}
(Estimate of $I^F_{2,4}$) For the term \( I^F_{2,4} \), using the Hölder inequality, \eqref{Phiscale}, and \eqref{Phitxinf}, we get
\begin{align*}
I^F_{2,4}&\les \sum_{|\al|=\mathrm{N}+1}\frac{\kappa}{\eps^2} \|\P(\p^{\al}M^{\e}) |M^{\e}|^{-\frac{3}{4}}\|_{L^2_{x,v}} \| \AC{\P}(\p^{\al}F^{\e})|M^{\e}|^{-\frac{1}{2}}\|_{L^2_{x,v}}\bigg\|\frac{\left(\p_t M^{\e} + \frac{v}{\eps}\cdot\nabla_x M^{\e}\right)}{|M^{\e}|^{\frac{3}{4}}} \bigg\|_{L^\infty_{x,v}} \cr 
&\les \sum_{|\al|=\mathrm{N}+1}\frac{\kappa}{\eps^2} (\e\kappa^{-\frac{1}{2}}\mathcal{E}_{tot}^{\frac{1}{2}}) (\e^2 \mathcal{D}_{top}^{\frac{1}{2}}) \bigg[\|\nabla_xu^{\e}\|_{L^\infty_x}+\|\nabla_x\ta^{\e}\|_{L^\infty_x} +\frac{1}{\e}\|\TbT\|_{L^\infty_x}\bigg]. 
\end{align*}

Combining \eqref{I2,1est} and \eqref{I2,234} with the decomposition
\( I^F_2 = I^F_{2,1}+I^F_{2,2}+I^F_{2,3}+I^F_{2,4} \) in \eqref{I2,234}, we obtain
\begin{align*}
\bega
I^F_2(t) &\les \big(\|\nabla_xu^{\e}(t)\|_{L^\infty_x}+\|\nabla_x\ta^{\e}(t)\|_{L^\infty_x}\big) \Big(\mathcal{E}_{top}(t)+\e\kappa^{\frac{1}{2}}\mathcal{E}_M(t)\mathcal{E}_{tot}^{\frac{1}{2}}(t)+\e^2\kappa\mathcal{E}_M^2(t)\Big) \cr 
&+ \int_{\Omega}\sum_{i,j}|\p_{x_i}u^{\e}_j(t)| \mathcal{V}_{2}(t) dx + \int_{\Omega}\sum_i |\p_{x_i}\ta^{\e}(t)| \mathcal{V}_{3}(t) dx + \frac{1}{\eps} \int_{\Omega} \TbT(t) \mathcal{V}_{2}(t) dx, \cr 
&+ \frac{1}{\e} \|\TbT(t)\|_{L^\infty_x} \Big(\mathcal{E}_{top}(t)+\e\kappa^{\frac{1}{2}}\mathcal{E}_M(t)\mathcal{E}_{tot}^{\frac{1}{2}}(t) +\e\kappa^{\frac{1}{2}}\mathcal{E}_{top}^{\frac{1}{2}}(t)\mathcal{D}_{top}^{\frac{1}{2}}(t) +\e^2\kappa\mathcal{E}_M^2(t) \Big)  \cr 
&+ \e\kappa^{\frac{1}{2}}  \mathcal{E}_M^{\frac{1}{2}}(t) \mathcal{E}_{top}^{\frac{1}{2}}(t)\mathcal{D}_{top}^{\frac{1}{2}}(t),
\enda
\end{align*}
where we used \( \|\nabla_xu^{\e}\|_{L^\infty_x}+\|\nabla_x\ta^{\e}\|_{L^\infty_x} \leq C\mathcal{E}_M^{\frac{1}{2}} \) from \eqref{rutinf} in the last line to estimate \( I^F_{2,4} \). 
We also note that, by using \( \e\mathcal{E}_M^{\frac{1}{2}} \leq C \) from \( (\mathcal{B}_1) \) in \eqref{condition}, we have
\begin{align*}
\bega
\mathcal{E}_{top}(t)+\e\kappa^{\frac{1}{2}}\mathcal{E}_M(t)\mathcal{E}_{tot}^{\frac{1}{2}}(t)+\e^2\kappa\mathcal{E}_M^2(t) \leq C\mathcal{E}_{tot}(t).
\enda
\end{align*}
Finally, applying \( \|\TbT(t)\|_{L^\infty_x} \leq C\e^2\kappa^{\frac{1}{2}}\mathcal{D}_G^{\frac{1}{2}}(t) \) from \eqref{ABGscale}$_2$, we complete the proof of the estimate \eqref{E1est}.
\end{proof}

\subsubsection{Estimate of other terms}
In this part, we estimate \( I^F_3 \), \( I^F_4 \), and \( I^F_5 \), as defined in~\eqref{E123def}. 
When taking the \( L^\infty_x \)-norm of the microscopic part, for instance $\big\|\langle v \rangle^{\tfrac{1}{2}} \AC{\P}F^{\e}\, |M^{\e}|^{-\tfrac{1}{2}}(t)\big\|_{L^\infty_xL^2_v},$ and applying the Sobolev embedding, the derivative may fall either on \( \AC{\P}F^{\e} \) or on \( |M^{\e}|^{-\tfrac{1}{2}} \). 
In the latter case, it generates high moments \( \mathcal{V}_\ell \), which will later be controlled by means of a pointwise estimate derived from the characteristic approach in Section \ref{Sec.V}.


\begin{lemma}\label{L.E345}
For $t\in[0,T]$ satisfying the bootstrap assumption \eqref{condition}, and for both the purely spatial derivatives \eqref{caseA} and the space-time derivatives \eqref{caseB}, the terms $I^F_3$, $I^F_4$ and $I^F_5$, defined in \eqref{E123def} satisfies the following upper bound:
\begin{align*}
\bega 
I^F_3(t) &\les \frac{1}{\e}\Big\|\|\la v \ra^{\frac{1}{2}} \AC{\P}F^{\e} |M^{\e}|^{-\frac{1}{2}}(t)\|_{L^2_v}\Big\|_{L^\infty_x} \Big(\kappa^{-\frac{1}{2}}\mathcal{E}_{tot}^{\frac{1}{2}}(t) + \e\mathcal{D}_{top}^{\frac{1}{2}}(t)\Big)\mathcal{D}_{top}^{\frac{1}{2}}(t), \cr 
I^F_4(t) &\les  \e\kappa^{\frac{1}{2}}\mathcal{E}_{M}^{\frac{1}{2}}(t)\mathcal{D}_G^{\frac{1}{2}}(t)\mathcal{D}_{top}^{\frac{1}{2}}(t) + \frac{1}{\e}\sum_{1\leq |\al| \leq\lfloor(\mathrm{N}+1)/2\rfloor} \Big\|\|\la v \ra^{\frac{1}{2}} \p^{\al}\AC{\P}F^{\e} |M^{\e}(t)|^{-\frac{1}{2}}\|_{L^2_v}\Big\|_{L^\infty_x} \mathcal{E}_{M}^{\frac{1}{2}}(t)\mathcal{D}_{top}^{\frac{1}{2}}(t), \cr 
I^F_5(t) & \les \kappa^{\frac{1}{2}} \sum_{0\leq|\al|\leq \lfloor(\mathrm{N}+1)/2\rfloor}\Big\|\|\la v \ra^{\frac{1}{2}} \p^{\al}\AC{\P}F^{\e} |M^{\e}|^{-\frac{1}{2}}(t)\|_{L^2_v}\Big\|_{L^\infty_x} \mathcal{D}_G^{\frac{1}{2}}(t)\mathcal{D}_{top}^{\frac{1}{2}}(t).
\enda
\end{align*}
\end{lemma}
\begin{proof}
(Estimate of \( I^F_3 \))
Applying \eqref{nonlin} to \( I^F_3 \), we have 
\begin{align*}
\bega 
I^F_3&\les  \sum_{|\al|=\mathrm{N}+1}\frac{2}{\e^4} \int_{\Omega} \Big\|\la v \ra^{\frac{1}{2}} \p^{\al}F^{\e}|M^{\e}|^{-\frac{1}{2}}\Big\|_{L^2_v}\Big\|\la v \ra^{\frac{1}{2}} \AC{\P}F^{\e} |M^{\e}|^{-\frac{1}{2}}\Big\|_{L^2_v} \Big\|\la v \ra^{\frac{1}{2}} \p^{\al}\AC{\P}F^{\e}|M^{\e}|^{-\frac{1}{2}}\Big\|_{L^2_v}dx.
\enda
\end{align*}
Then, we take the \( L^\infty_x \) norm of the middle term \( \|\la v \ra^{\frac{1}{2}} \AC{\P}F^{\e} |M^{\e}|^{-\frac{1}{2}} \|_{L^2_v} \) and apply Hölder's inequality in \( x \), yielding
\begin{align}\label{IF2leq}
\bega 
I^F_3&\les  \sum_{|\al|=\mathrm{N}+1}\frac{2}{\e^4} \Big\|\la v \ra^{\frac{1}{2}} \p^{\al}F^{\e}|M^{\e}|^{-\frac{1}{2}}\Big\|_{L^2_{x,v}} \Big\|\|\la v \ra^{\frac{1}{2}} \AC{\P}F^{\e} |M^{\e}|^{-\frac{1}{2}}\|_{L^2_v}\Big\|_{L^\infty_x} (\e^2\mathcal{D}_{top}^{\frac{1}{2}}),
\enda
\end{align}
where we used the dissipation \( \mathcal{D}_{top} \) defined in \eqref{N-EDdef2}.
Using the decomposition \( \p^{\al}F^{\e} = \P(\p^{\al}M^{\e}) + \AC{\P}(\p^{\al}F^{\e}) \) from \eqref{PalF}, and applying \eqref{Phiscale} and \eqref{N-EDdef2}, the first term in \eqref{IF2leq} can be estimated by
\begin{align}\label{IF2leq2}
\bega 
\Big\|\la v \ra^{\frac{1}{2}} \p^{\al}F^{\e}|M^{\e}|^{-\frac{1}{2}}\Big\|_{L^2_{x,v}} &\leq \Big\|\la v \ra^{\frac{1}{2}} \Big(\P(\p^{\al}M^{\e})+\AC{\P}(\p^{\al}F^{\e})\Big)|M^{\e}|^{-\frac{1}{2}}\Big\|_{L^2_{x,v}}  \les \e\kappa^{-\frac{1}{2}}\mathcal{E}_{tot}^{\frac{1}{2}} + \e^2\mathcal{D}_{top}^{\frac{1}{2}}.
\enda
\end{align}
Combining \eqref{IF2leq} and \eqref{IF2leq2} yields the desired result.
\\ (Estimate of \( I^F_4 \)) By following the same approach as in Lemma \ref{L.G45}, we have 
\begin{align*}
\bega 
I^F_{4}&\les  \sum_{|\al|=\mathrm{N}+1}\frac{1}{\e^2} \sum_{1\leq \beta \leq \lfloor(\mathrm{N}+1)/2\rfloor} \Big\|\la v \ra^{\frac{1}{2}} \p^{\beta}\AC{\P}F^{\e}|M^{\e}|^{-\frac{1}{2}}\Big\|_{L^\infty_xL^2_v}\Big\|\la v \ra^{\frac{1}{2}} \p^{\al-\beta}M^{\e}|M^{\e}|^{-\frac{1}{2}}\Big\|_{L^2_{x,v}} \mathcal{D}_{top}^{\frac{1}{2}} \cr 
&+\sum_{|\al|=\mathrm{N}+1}\frac{1}{\e^2} \sum_{\lfloor(\mathrm{N}+1)/2\rfloor< \beta <\al} \Big\|\la v \ra^{\frac{1}{2}} \p^{\beta}\AC{\P}F^{\e}|M^{\e}|^{-\frac{1}{2}}\Big\|_{L^2_{x,v}}\Big\|\la v \ra^{\frac{1}{2}} \p^{\al-\beta}M^{\e}|M^{\e}|^{-\frac{1}{2}}\Big\|_{L^\infty_xL^2_v} \mathcal{D}_{top}^{\frac{1}{2}} \cr 
&\leq \frac{1}{\e^2}\sum_{1\leq \beta \leq\lfloor(\mathrm{N}+1)/2\rfloor} \Big\|\la v \ra^{\frac{1}{2}} \p^{\al-\beta}\AC{\P}F^{\e}|M^{\e}|^{-\frac{1}{2}}\Big\|_{L^\infty_xL^2_v} (\e\mathcal{E}_{M}^{\frac{1}{2}})\mathcal{D}_{top}^{\frac{1}{2}} + \frac{1}{\e^2}(\e^2\kappa^{\frac{1}{2}}\mathcal{D}_G^{\frac{1}{2}}) (\e\mathcal{E}_{M}^{\frac{1}{2}})\mathcal{D}_{top}^{\frac{1}{2}} .
\enda
\end{align*}
(Estimate of \( I^F_5 \)) 
We note that depending on the number of derivatives, we apply different scales of dissipation, namely \( \mathcal{D}_G \) and \( \mathcal{D}_{top} \), defined in \eqref{N-EDdef} and \eqref{N-EDdef2}, respectively. Without loss of generality, we assume \( 0 \leq |\beta| \leq \lfloor(\mathrm{N}+1)/2\rfloor \). Then, using \( \mathcal{D}_{top} \) and \( \mathcal{D}_G \) to control \( \AC{\P}(\p^{\al}F^{\e}) \) and \( \p^{\al - \beta} \AC{\P}F^{\e} \), respectively, we obtain
\begin{align*}
\bega 
I^F_{5} &\les \kappa^{\frac{1}{2}} \sum_{0\leq|\al|\leq \lfloor(\mathrm{N}+1)/2\rfloor}\Big\|\|\la v \ra^{\frac{1}{2}} \p^{\al}\AC{\P}F^{\e} |M^{\e}|^{-\frac{1}{2}}\|_{L^2_v}\Big\|_{L^\infty_x} \mathcal{D}_G^{\frac{1}{2}} \mathcal{D}_{top}^{\frac{1}{2}} .
\enda
\end{align*}
\end{proof}

\subsubsection{Proof of Lemma \ref{L.ABG}}\label{Sec.Gembed}

In this part, we prove Lemma \ref{L.ABG}.


\begin{proof}[\textbf{Proof of Lemma \ref{L.ABG}}]
We apply the \( L^2_x \) norm to the identities in \eqref{alhata} and \eqref{alhatb}, and use the H\"{o}lder inequality to get
\begin{align*}
\bega
\|\p^{\al} \mathbf{r}_{ij}^{\e}\|_{L^2_x} &\leq  \Big\|\|\mathfrak{R}^{\e}_{ij}|M^{\e}|^{\frac{1}{2}}\|_{L^2_v}\|\p^{\al}\AC{\P}F^{\e}|M^{\e}|^{-\frac{1}{2}}\|_{L^2_v} \Big\|_{L^2_x}, \cr 
\|\p^{\al}\mathfrak{q}_i^{\e}\|_{L^2_x} &\leq \Big\|\|\mathcal{Q}^{\e}_i|M^{\e}|^{\frac{1}{2}}\|_{L^2_v}\|\p^{\al}\AC{\P}F^{\e}|M^{\e}|^{-\frac{1}{2}}\|_{L^2_v} \Big\|_{L^2_x}  + \sum_{1\leq\beta\leq \al}\binom{\al}{\beta} \sum_{k} \Big\|\p^{\beta}\mathrm{U}^{\e}_k \p^{\al-\beta} \mathbf{r}_{ik}^{\e} \Big\|_{L^2_x}.
\enda
\end{align*}
Since the macroscopic fields are bounded as \( (|\mathrm{P}^{\e} -1|, |\mathrm{U}^{\e}|, |\mathrm{\Theta}^{\e}-1|) \ll 1 \) by \eqref{condition}, it follows that 
$\|\mathfrak{R}^{\e}_{ij}e^{-C|v|^2}\|_{L^2_v}^2 \leq C$ and $\|\mathcal{Q}^{\e}_ie^{-C|v|^2}\|_{L^2_v}^2 \leq C$. Thus, we obtain the estimates
\begin{align}\label{albe.1}
\bega
\|\p^{\al} \mathbf{r}_{ij}^{\e}\|_{L^2_x} &\les \|\p^{\al}\AC{\P}F^{\e}|M^{\e}|^{-\frac{1}{2}}\|_{L^2_{x,v}}, \cr 
\|\p^{\al}\mathfrak{q}_i^{\e}\|_{L^2_x} &\les \|\p^{\al}\AC{\P}F^{\e}|M^{\e}|^{-\frac{1}{2}}\|_{L^2_{x,v}} + \sum_{1\leq\beta\leq \al}\binom{\al}{\beta} \sum_{k} \Big\|\p^{\beta}\mathrm{U}^{\e}_k \p^{\al-\beta} \mathbf{r}_{ik}^{\e} \Big\|_{L^2_x}.
\enda
\end{align}
For the estimate of the first term in \eqref{albe.1}, we note that different dissipation norms are used depending on the number of the derivative, namely \( |\al| = \mathrm{N} \) and \( |\al| = \mathrm{N}+1 \). We claim that
\begin{align}\label{albe.claim}
\bega
\|\p^{\al}\AC{\P}F^{\e}|M^{\e}|^{-\frac{1}{2}}\|_{L^2_{x,v}} \leq \begin{cases} \e^2\kappa^{\frac{1}{2}}\mathcal{D}_G^{\frac{1}{2}}, \quad &\mbox{when} \quad 0\leq|\al|\leq\mathrm{N}, \\ 
\e^2(\mathcal{D}_{top}^{\frac{1}{2}}(t)+\mathbf{1}_{|\al|\geq2}\mathcal{E}_M(t)), \quad &\mbox{when} \quad |\al|=\mathrm{N}+1.
\end{cases}
\enda
\end{align}
When \( 0\leq|\al|\leq\mathrm{N} \), the definition of \( \mathcal{D}_G \) in \eqref{N-EDdef} directly yields the first estimate in \eqref{albe.claim}. For the case \( |\al|=\mathrm{N}+1 \), we apply the commutator identity between \( \p^{\al} \) and \( \AC{\P} \) given in \eqref{pACP}. This gives
\begin{align}\label{al5decomp}
\bega
\int_{\Omega}\int_{\R^3}|\p^{\al}\AC{\P}F^{\e}|^2|M^{\e}|^{-1} dvdx &\leq \int_{\Omega}\int_{\R^3}|\AC{\P} \p^{\al}F^{\e}|^2|M^{\e}|^{-1} dvdx \cr 
&+ \int_{\Omega}\int_{\R^3}\bigg|\AC{\P} \bigg( \sum_{2\leq i\leq |\al|}\eps^i \Phi_{\al}^i M^{\e} \bigg)\bigg|^2|M^{\e}|^{-1} dvdx \cr 
&\les \e^4\mathcal{D}_{top} + \e^4 \mathbf{1}_{|\al|\geq2}\mathcal{E}_M^2,
\enda
\end{align}
where the first term is bounded using \( \mathcal{D}_{top} \) from \eqref{N-EDdef2}, and the second term is estimated using \eqref{ACPX} and \eqref{Rscale}. We also used the fact that \( \sum_{2\leq i\leq |\al|}\eps^i \Phi_{\al}^i M^{\e} \) is nonzero only when \( |\al|\geq2 \), as defined in \eqref{RpM-def}. This completes the proof of \eqref{albe.claim}.
For the second line of \eqref{albe.1}, we also claim that
\begin{align}\label{albe.claim2}
\bega
\sum_{1\leq\beta\leq \al} \Big\|\p^{\beta}\mathrm{U}^{\e}_k \p^{\al-\beta} \mathbf{r}_{ik}^{\e} \Big\|_{L^2_x} &\les \begin{cases} \e^3\kappa^{\frac{1}{2}}\mathcal{E}_M^{\frac{1}{2}}\mathcal{D}_G^{\frac{1}{2}}, \quad &\mbox{when} \quad 0\leq|\al|\leq\mathrm{N}, \\ 
\e^3\mathcal{E}_{tot}^{\frac{1}{2}}(\mathcal{D}_{tot}^{\frac{1}{2}}+\mathbf{1}_{|\al|\geq2}\mathcal{E}_M(t)), \quad &\mbox{when} \quad |\al|=\mathrm{N}+1.
\end{cases}
\enda
\end{align}
To estimate this term, we apply Agmon's inequality \eqref{Agmon}, yielding
\begin{align*}
\bega
\sum_{1\leq\beta\leq \al} \Big\|\p^{\beta}\mathrm{U}^{\e}_k \p^{\al-\beta} \mathbf{r}_{ik}^{\e} \Big\|_{L^2_x} &\les \sum_{\lfloor |\al|/2 \rfloor < |\beta| \leq |\al|} \e \|\p^{\beta}u^{\e}_k\|_{L^2_x} \|\p^{\al-\beta} \mathbf{r}_{ik}^{\e}\|_{H^2_x} \cr 
&+ \sum_{1\leq |\beta| \leq [|\al|/2]} \e  \|\p^{\beta}u^{\e}_k\|_{L^\infty_x} \|\p^{\al-\beta} \mathbf{r}_{ik}^{\e}\|_{L^2_x}.
\enda
\end{align*}
Using the bounds from \eqref{rutscale}, \eqref{rutinf}, and the estimate \eqref{albe.claim}, we get
\begin{align*}
\bega
&\sum_{1\leq\beta\leq \al} \Big\|\p^{\beta}\mathrm{U}^{\e}_k \p^{\al-\beta} \mathbf{r}_{ik}^{\e} \Big\|_{L^2_x} \cr 
&\les \begin{cases}
\e\mathcal{E}_M^{\frac{1}{2}}\e^2\kappa^{\frac{1}{2}}\mathcal{D}_G^{\frac{1}{2}},  &\mbox{when} \quad 0\leq|\al|\leq\mathrm{N}, \\ 
\e (\kappa^{-\frac{1}{2}}\mathcal{E}_{tot}^{\frac{1}{2}}) (\e^2\kappa^{\frac{1}{2}}\mathcal{D}_G) + \e\mathcal{E}_M^{\frac{1}{2}}\Big(\e^2\kappa^{\frac{1}{2}}\mathcal{D}_G+\e^2\mathcal{D}_{top}^{\frac{1}{2}}(t)+\mathbf{1}_{|\al|\geq2}\e^2\mathcal{E}_M(t)\Big) ,  &\mbox{when} \quad |\al|=\mathrm{N}+1. \end{cases}
\enda
\end{align*}
Using the definitions of \( \mathcal{E}_{tot} \) and \( \mathcal{D}_{tot} \) in \eqref{EDtotdef}, this proves the claim \eqref{albe.claim2}.
Combining the estimates \eqref{albe.claim} and \eqref{albe.claim2} with \eqref{albe.1}, and using the bound \( \e\mathcal{E}_{tot}^{\frac{1}{2}} \leq 1 \) from \eqref{condition}, we conclude the proof of the estimate \eqref{ABGscale}$_1$. \\
(2) Applying Agmon's inequality \eqref{Agmon} and using $\eqref{ABGscale}_1$, we obtain the desired result in $\eqref{ABGscale}_2$. 
\end{proof}

\begin{proof}[\textbf{Proof of Proposition \ref{P.F.Energy}}]
Combining the estimates for $I^F_1$, $I^F_2$, $I^F_3$, $I^F_4$, and $I^F_5$ obtained in Lemma~\ref{L.E1}, Lemma~\ref{L.E2}, and Lemma~\ref{L.E345}, we apply them to \eqref{E'0}.
\begin{align*}
\bega
\sum_{|\al|=\mathrm{N}+1}&\frac{\kappa}{2\e^2}\frac{d}{dt} \int_{\Omega \times \R^3}|\p^{\al}F^{\e}|^2 |M^{\e}|^{-1}dvdx  + \sum_{|\al|=\mathrm{N}+1}\frac{1}{\e^4}\int_{\Omega \times \R^3}\mathcal{L}(\p^{\al}F^{\e})\p^{\al}F^{\e}|M^{\e}|^{-1}  dvdx \cr 
&\les  \big(\|\nabla_xu^{\e}\|_{L^\infty_x} +\|\nabla_x\ta^{\e}\|_{L^\infty_x}\big)  \big(\mathcal{E}_{tot}+\mathcal{E}_M^{\frac{1}{2}}\mathcal{D}_{top}^{\frac{1}{2}}\big)\cr 
&+ \int_{\Omega}|\nabla_xu^{\e}| \mathcal{V}_2 dx + \int_{\Omega}|\nabla_x\ta^{\e}| \mathcal{V}_3 dx + \frac{1}{\eps} \int_{\Omega} \TbT \mathcal{V}_2 dx \cr 
&+ \e\kappa^{\frac{1}{2}}\mathcal{D}_G^{\frac{1}{2}} \Big(\mathcal{E}_{tot}+\e\kappa^{\frac{1}{2}}\mathcal{E}_{top}^{\frac{1}{2}}\mathcal{D}_{top}^{\frac{1}{2}} \Big)  + \e\kappa^{\frac{1}{2}} \mathcal{E}_{tot}\mathcal{D}_{top}^{\frac{1}{2}} +\e\mathcal{E}_{M}^{\frac{3}{2}}\mathcal{D}_{top}^{\frac{1}{2}} +\mathcal{Z}_{top}^{time}  \cr 
&+ \frac{1}{\e}\Big\|\|\la v \ra^{\frac{1}{2}} \AC{\P}F^{\e} |M^{\e}|^{-\frac{1}{2}}\|_{L^2_v}\Big\|_{L^\infty_x} \Big(\kappa^{-\frac{1}{2}}\mathcal{E}_{tot}^{\frac{1}{2}} + \e\mathcal{D}_{top}^{\frac{1}{2}}\Big)\mathcal{D}_{top}^{\frac{1}{2}} \cr 
&+ \e\kappa^{\frac{1}{2}}\mathcal{E}_{M}^{\frac{1}{2}}\mathcal{D}_G^{\frac{1}{2}}\mathcal{D}_{top}^{\frac{1}{2}} + \frac{1}{\e}\sum_{1\leq |\al| \leq\lfloor(\mathrm{N}+1)/2\rfloor} \Big\|\|\la v \ra^{\frac{1}{2}} \p^{\al}\AC{\P}F^{\e} |M^{\e}|^{-\frac{1}{2}}\|_{L^2_v}\Big\|_{L^\infty_x} \mathcal{E}_{M}^{\frac{1}{2}}\mathcal{D}_{top}^{\frac{1}{2}} \cr 
&+ \kappa^{\frac{1}{2}} \sum_{0\leq|\al|\leq \lfloor(\mathrm{N}+1)/2\rfloor}\Big\|\|\la v \ra^{\frac{1}{2}} \p^{\al}\AC{\P}F^{\e} |M^{\e}|^{-\frac{1}{2}}\|_{L^2_v}\Big\|_{L^\infty_x} \mathcal{D}_G^{\frac{1}{2}} \mathcal{D}_{top}^{\frac{1}{2}}.
\enda
\end{align*}
Applying the coercivity estimate from Lemma~\ref{L.coer} to the second term on the first line, we obtain the desired estimate \eqref{totalEt} in Proposition~\ref{P.F.Energy}.
\end{proof}

\section{Macroscopic estimate}\label{Sec.macro}

In this section, we estimate the macroscopic variables $(\rho^{\e}, u^{\e}, \ta^{\e})$ for multi-indices $0 \leq |\al| \leq \mathrm{N}$, considering both purely spatial derivatives~\eqref{caseA} and space--time derivatives~\eqref{caseB}. The key element in the analysis is the momentum--flux alignment term $\mathfrak{S}_W^{\al}$ from~\eqref{Amadef}, which contains the highest--order derivative of order $\mathrm{N}+1$. 
At leading order, we discover that this contribution is crucially canceled by $\mathfrak{S}_G^{\al}$ from~\eqref{Amidef}. Another important point is that, by exploiting the symmetric hyperbolic structure of equation \eqref{bWeqn}, we were able to cancel the contributions of the penalized terms $\frac{1}{\e}\nabla_x \cdot u^{\e}$ and $\frac{1}{\e}\nabla_x(\rho^{\e}+\ta^{\e})$ in the energy estimate.

The estimates in this section apply to both the finite velocity energy case in \eqref{caseEC} and the infinite velocity energy case in \eqref{caseECX}.
For the infinite velocity energy case, since $u^{\e} \notin L^2_x$, we apply the radial-energy decomposition to $u^{\e}$ and then use the cancellation between $\frac{1}{\e} \nabla_x \cdot u^{\e}$ and $\frac{1}{\e} \nabla_x (\rho^{\e} + \ta^{\e})$ through the energy estimate for $(\rho^{\e}, u^{\e}-\bar{u}, \ta^{\e})$.
In Section \ref{Sec.macro.Stri}, we estimate the acoustic variables $\mathbb{P}^{\perp}u^{\e}$ and $(\rho^{\e} + \ta^{\e})$ using the Strichartz estimates.

\subsection{Symmetric hyperbolic estimate}\label{Sec.macro.1}
In this part, we estimate $(\rho^{\e}, u^{\e}, \ta^{\e})$ using the local conservation laws \eqref{locconNew}.
For the infinite velocity energy case \eqref{caseECX}, we decompose $u^{\e}$ using both the Hodge decomposition and the radial-energy decomposition, as in \eqref{uexpand}, since $u^{\e} \notin L^2$.
We also recall the definition of $\bar{u}(x)$ in \eqref{barudef}, which depends on the case considered in \eqref{caseEC} and \eqref{caseECX}.

%
%

\begin{proposition}\label{P.EW} 
Let $\Omega=\R^d$ with $d=2$ or $3$ and $\mathrm{N}>d/2+1$.
Assume that the initial data $u_0^{\e}$ satisfies either the finite velocity energy condition \eqref{caseEC} or the infinite velocity energy condition \eqref{caseECX}.
Under the bootstrap assumption \eqref{condition}, the following estimate holds in both the purely spatial derivative case \eqref{caseA} and the space-time derivative case \eqref{caseB} with arbitrary $\mathfrak{n}\geq1$.
\begin{align}\label{EWresult}
\bega
\frac{1}{2}&\bigg(\|\rho^{\e}(t)\|_{L^2_x}^2+\frac{1}{k_B}\|u^{\e}(t)-\bar{u}\|_{L^2_x}^2+\frac{3}{2}\|\ta^{\e}(t)\|_{L^2_x}^2\bigg) +\frac{1}{C}\sum_{1\leq|\al|\leq\mathrm{N}}\|\p^{\al}(\rho^{\e},u^{\e},\ta^{\e})(t)\|_{L^2_x}^2 \cr 
&\leq \frac{1}{2}\bigg(\|\rho^{\e}_0\|_{L^2_x}^2+\frac{1}{k_B}\|u^{\e}_0-\bar{u}\|_{L^2_x}^2+\frac{3}{2}\|\ta^{\e}_0\|_{L^2_x}^2\bigg) +C\sum_{1\leq|\al|\leq\mathrm{N}}\|\p^{\al}(\rho^{\e}_0,u^{\e}_0,\ta^{\e}_0)\|_{L^2_x}^2 \cr 
&+C \int_0^t\Big(1+\|\nabla_x\bar{u}\|_{L^\infty_x} +\|\nabla_x(\rho^{\e},u^{\e},\ta^{\e}))(s)\|_{L^\infty_x} \Big) \mathcal{E}_M^{\mathrm{N}}(F^{\e}(s))ds +\sum_{0\leq|\al|\leq\mathrm{N}}\int_0^t\mathfrak{S}_W^{\al}(s)ds \cr 
&+ \int_0^t\bigg(\e(\mathcal{E}_M^{\mathrm{N}}(F^{\e}(s)))^{\frac{3}{2}}+\e(\mathcal{E}_M^{\mathrm{N}}(F^{\e}(s)))^2 + \e \kappa^{\frac{1}{2}}\mathcal{E}_M^{\mathrm{N}}(F^{\e}(s))(\mathcal{D}_G^{\mathrm{N}}(F^{\e}(s)))^{\frac{1}{2}} + \mathcal{Z}_W^{time}(s)\bigg)ds.
\enda
\end{align}
Here, $\bar{u}$, $\mathcal{E}_M^{\mathrm{N}}(F^{\e}(t))$, and 
$\mathcal{D}_G^{\mathrm{N}}(F^{\e}(t))$ are defined in
\eqref{barudef} and \eqref{N-EDdef}.
Here, the momentum-flux alignment contribution of $\mathfrak{S}_W^{\al}(t)$ is given by
\begin{align}\label{Amadef}
\bega
\mathfrak{S}_W^{\al}(t) &:= -\frac{1}{\e^2}\sum_{i,j} \int_{\Omega} \p^{\al}\bigg[\frac{1}{k_B\mathrm{P}^{\e} \mathrm{\Theta}^{\e}} \p_{x_j} \mathbf{r}_{ij}^{\e}\bigg] \p^{\al}u^{\e}_i dx \cr 
&-\frac{1}{\e^2}\int_{\Omega} \p^{\al}\bigg[\frac{1}{k_B\mathrm{P}^{\e}\mathrm{\Theta}^{\e}} \Big(\sum_{j} \p_{x_j} \mathfrak{q}_j^{\e} +\sum_{i,j} \p_{x_i}\mathrm{U}^{\e}_j \mathbf{r}_{ij}^{\e}\Big)\bigg] \p^{\al}\ta^{\e} dx,
\enda
\end{align}
for $1 \leq |\al| \leq \mathrm{N}$, and
\begin{align}\label{Amadef0}
\bega
\mathfrak{S}_W^{0}(t)&: = -\frac{1}{\e^2}\sum_{i,j} \int_{\Omega}\frac{1}{k_B\mathrm{P}^{\e}} \p_{x_j} \mathbf{r}_{ij}^{\e}  (u^{\e}-\bar{u}) dx \cr 
&- \frac{1}{\e^2}\int_{\Omega}\frac{1}{k_B\mathrm{P}^{\e}\mathrm{\Theta}^{\e} }\bigg(\sum_{j} \p_{x_j} \mathfrak{q}_j^{\e} + \sum_{i,j} \p_{x_i}\mathrm{U}^{\e}_j \mathbf{r}_{ij}^{\e}\bigg) \ta^{\e} dx,
\enda
\end{align}
for $|\al| = 0$. 
The contribution from the time derivatives, denoted by 
$\mathcal{Z}_W^{\mathrm{time}}(t)$, is defined by
\begin{align}\label{Bmadef}
\bega
\mathcal{Z}_W^{time}(t) &:= \begin{cases}
0, \quad &\text{for $\eqref{caseA}$}, \\
\e^{\mathfrak{n}-1}\Big(\|\nabla_x\nabla_x\cdot u^{\e}(t)\|_{L^\infty_x}\mathcal{E}_M^{\mathrm{N}}(F^{\e}(t))  +\e^2\kappa(\mathcal{E}_{M}^{\mathrm{N}}(F^{\e}(t)))^\frac{1}{2}\mathcal{D}_G^{\mathrm{N}}(F^{\e}(t)) \Big), \quad &\text{for $\eqref{caseB}$}.
\end{cases}
\enda
\end{align}

\end{proposition}

We provide the proof of Proposition \ref{P.EW} at the end of this section, after presenting several lemmas toward it.

\begin{lemma}\label{L.EW} 
For the purely spatial derivatives \eqref{caseA} and space-time derivatives \eqref{caseB}, and for $0 \leq |\al| \leq \mathrm{N}$, we have
\begin{align}\label{EWineq}
\frac{d}{dt}\int_{\Omega} &(\p^{\al}\bW^{\e})^T \mathcal{S}^{\e} (\p^{\al}\bW^{\e}) dx \leq \mathfrak{S}_W^{\al}(t)
 + C\bigg[I_{W}^{\al}(t) + II_{W}^{\al}(t) + III_{W}^{\al}(t) \bigg],
\end{align}
where $\bW^{\e}$, $\mathcal{S}^{\e}$, and $\mathfrak{S}_W^{\al}$
are defined in~\eqref{bWdef}, \eqref{Sdef}, and \eqref{Amadef},
respectively, and
\begin{align}\label{EW12def}
\bega
I_{W}^{\al}(t) &=\Big(\| \nabla_x \cdot (tr(\mathcal{S}^{\e})u^{\e})\|_{L^\infty_x} + \|\p_t\mathcal{S}^{\e}\|_{L^\infty_x}\Big) \|\p^{\al}\bW^{\e}\|_{L^2_x}^2, \cr
II_{W}^{\al}(t) &= \sum_{1\leq \beta \leq \al} \sum_{0\leq \gamma\leq\beta}\|\p^{\gamma}\mathcal{S}^{\e} \p^{\beta-\gamma}u^{\e} \cdot \nabla_x \p^{\al-\beta}\bW^{\e}\|_{L^2_x}\|\p^{\al}\bW^{\e}\|_{L^2_x}, \cr
III_{W}^{\al}(t) &= \sum_{1\leq\beta\leq\al}\|\p^{\beta}\mathcal{S}^{\e} \p_t\p^{\al-\beta}\bW^{\e}\|_{L^2_x}\|\p^{\al}\bW^{\e}\|_{L^2_x}.
\enda
\end{align}
\end{lemma}

\begin{proof}
We apply the operator $\p^{\al}$ to the local conservation laws \eqref{bWeqn}:
\begin{align}\label{altilw}
\mathcal{S}^{\e} (\p_t \p^{\al}\bW^{\e} + u^{\e} \cdot \nabla_x \p^{\al}\bW^{\e}) + \frac{1}{\e} L(\p_x) \p^{\al}\bW^{\e} = \p^{\al} \bPhi^{\e}_{\bW} - \llbracket\p^{\al},\mathcal{S}^{\e} (\p_t + u^{\e} \cdot \nabla_x )\rrbracket \bW^{\e},
\end{align}
where the commutator is defined as $\llbracket A, B \rrbracket := AB - BA$.  
Since $\int_{\Omega}\frac{1}{\e} (L(\p_x)\p^{\al}\bW^{\e}) \p^{\al}\bW^{\e} dx =0$, multiplying \eqref{altilw} by $(\p^{\al} \bW^{\e})^T$ and integrating yields
\begin{align}\label{EW.1}
\bega
\frac{1}{2}\frac{d}{dt} \int_{\Omega} (\p^{\al}\bW^{\e})^T \mathcal{S}^{\e} (\p^{\al}\bW^{\e})& dx = \frac{1}{2}\int_{\Omega} \nabla_x \cdot (tr(\mathcal{S}^{\e})u^{\e})|\p^{\al}\bW^{\e}|^2 dx + \frac{1}{2}\int_{\Omega} (\p^{\al}\bW^{\e})^T \p_t\mathcal{S}^{\e}  (\p^{\al}\bW^{\e}) dx \cr 
&+ \int_{\Omega} (\p^{\al}\bW^{\e})^T \llbracket \p^{\al},\mathcal{S}^{\e} (\p_t + u^{\e} \cdot \nabla_x )  \rrbracket \bW^{\e} dx + \int_{\Omega} (\p^{\al}\bW^{\e})^T \p^{\al}\bPhi^{\e}_{\bW} dx.
\enda
\end{align}
For the first line of \eqref{EW.1}, taking $L^\infty_x$ norms of $\nabla_x \cdot (\mathrm{tr}(\mathcal{S}^{\e}) u^{\e})$ and $\p_t \mathcal{S}^{\e}$ gives rise to the term $I_W^{\al}(t)$.
For the first term in the second line of \eqref{EW.1}, the commutator expands as
\begin{align*}
\llbracket \p^{\al},\mathcal{S}^{\e} (\p_t + u^{\e} \cdot \nabla_x ) \rrbracket \bW^{\e} &= \sum_{1\leq\beta\leq\al} \binom{\al}{\beta} \p^{\beta}\mathcal{S}^{\e} \p_t\p^{\al-\beta}\bW^{\e} \cr 
&+ \sum_{1\leq \beta \leq \al} \sum_{0\leq \gamma\leq\beta} \binom{\al}{\beta} \binom{\beta}{\gamma} \p^{\gamma}\mathcal{S}^{\e} \p^{\beta-\gamma}u^{\e} \cdot \nabla_x \p^{\al-\beta}\bW^{\e} .
\end{align*}
Then, by applying the H\"{o}lder inequality, we obtain the terms $II_W^{\al}(t)$ and $III_W^{\al}(t)$, respectively.  
For the last term in \eqref{EW.1}, we substitute the definition of $\bPhi^{\e}_{\bW}$ from \eqref{gWdef} to obtain
\begin{align*}
\int_{\Omega} (\p^{\al}\bW^{\e})^T \p^{\al}\bPhi^{\e}_{\bW} dx &= -\frac{1}{\e^2}\sum_i \int_{\Omega} \p^{\al}\bigg[\frac{1}{k_B\mathrm{P}^{\e} \mathrm{\Theta}^{\e}}\sum_{j} \p_{x_j} \mathbf{r}_{ij}^{\e}\bigg] \p^{\al}u^{\e}_i dx \cr 
&\quad -\frac{1}{\e^2}\int_{\Omega} \p^{\al}\bigg[\frac{1}{k_B\mathrm{P}^{\e}\mathrm{\Theta}^{\e}} \Big(\sum_{j} \p_{x_j} \mathfrak{q}_j^{\e} +\sum_{i,j} \p_{x_i}\mathrm{U}^{\e}_j \mathbf{r}_{ij}^{\e}\Big)\bigg] \p^{\al}\ta^{\e} dx .
\end{align*}
This matches the definition of $\mathfrak{S}_W^{\al}$ given in \eqref{Amadef}.  
Therefore, we have obtained the desired result.
\end{proof}



In contrast to the finite-velocity energy case \eqref{caseEC}, in the infinite-velocity energy case \eqref{caseECX} we are only able to control $u^{\e}-\bar{u}$ in the low-frequency $L^2$.


\begin{lemma}\label{L.rutL2w} 
Let $\Omega=\R^d$ with $d=2$ or $3$. For $(\rho^{\e}, u^{\e}, \ta^{\e})$ satisfying the local conservation laws \eqref{locconNew}, we have
\begin{align*}
\bega
&\frac{1}{2}\frac{d}{dt}\bigg(\|\rho^{\e}\|_{L^2_x}^2+\frac{1}{k_B}\|u^{\e}-\bar{u}\|_{L^2_x}^2+\frac{3}{2}\|\ta^{\e}\|_{L^2_x}^2\bigg) \cr 
&\leq C \Big(1+\|\nabla_x\cdot u^{\e}\|_{L^\infty_x}+\|\nabla_x \bar{u}\|_{L^\infty_x}\Big) \bigg(\|\rho^{\e}\|_{L^2_x}^2+\frac{1}{k_B}\|u^{\e}-\bar{u}\|_{L^2_x}^2+\frac{3}{2}\|\ta^{\e}\|_{L^2_x}^2\bigg) \cr 
&+\|\nabla_x(\rho^{\e}+\ta^{\e})\|_{L^\infty_x} \frac{\|\mathrm{\Theta}^{\e}-1\|_{L^2_x}}{\e}  \|u^{\e}-\bar{u}\|_{L^2_x} +\mathfrak{S}_W^{0}(t),
\enda
\end{align*}
where $\bar{u}(x)$, and $\mathfrak{S}_W^{0}(t)$ are defined in \eqref{barudef}, and \eqref{Amadef0}, respectively.
\end{lemma}
\begin{proof}
Recall the decomposition of $u^{\e}$ given in \eqref{uexpand}. For brevity, we introduce the following notation, used only within this lemma:
\begin{align*}
\bega
\mathfrak{u}^{\e}(t,x):=u^{\e}(t,x)-\bar{u}(x) = \widetilde{u}^{\e}(t,x)+\mathbb{P}^{\perp}u^{\e}(t,x).
\enda
\end{align*}
Substituting \eqref{uexpand} into the velocity equation $\eqref{locconNew}_2$ yields
\begin{align*}
\bega
\p_t\mathfrak{u}^{\e} + \mathfrak{u}^{\e}\cdot \nabla_x \mathfrak{u}^{\e} + \mathfrak{u}^{\e}\cdot \nabla_x \bar{u} + \bar{u}\cdot\nabla_x \mathfrak{u}^{\e} +  \bar{u}\cdot\nabla_x \bar{u} +\frac{1}{\e}k_B\mathrm{\Theta}^{\e}\nabla_x(\rho^{\e}+\ta^{\e})
+\frac{1}{\e^2}\frac{1}{\mathrm{P}^{\e}}\sum_{j} \p_{x_j} \mathbf{r}_{ij}^{\e} =0.
\enda
\end{align*}
We now derive the energy estimate with respect to $\mathfrak{u}^{\e}$:
\begin{align}\label{muL2}
\bega
\frac{1}{2}\frac{d}{dt}\|\mathfrak{u}^{\e}\|_{L^2_x}^2 &\leq C \Big(1+\|\nabla_x\cdot\mathfrak{u}^{\e}\|_{L^\infty_x}+\|\nabla_x \bar{u}\|_{L^\infty_x}\Big)\|\mathfrak{u}^{\e}\|_{L^2_x}^2 \cr 
&-\frac{1}{\e}k_B \int_{\Omega}\mathrm{\Theta}^{\e}\nabla_x(\rho^{\e}+\ta^{\e}) \mathfrak{u}^{\e}dx 
-\frac{1}{\e^2}\int_{\Omega}\frac{1}{\mathrm{P}^{\e}}\sum_{j} \p_{x_j} \mathbf{r}_{ij}^{\e}  \mathfrak{u}^{\e}dx.
\enda
\end{align}
Here, we used
\begin{align*}
\bega
&\int_{\Omega} (\bar{u}\cdot\nabla_x \mathfrak{u}^{\e}) \mathfrak{u}^{\e}dx = -\frac{1}{2}\int_{\Omega} (\nabla_x\cdot\bar{u})|\mathfrak{u}^{\e}|^2dx =0, \cr 
&\int_{\Omega}(\bar{u}\cdot\nabla_x \bar{u}) \mathfrak{u}^{\e}dx = \int_{\Omega}\frac{x^{\perp}}{|x|^4}\bigg(\int_0^{|x|}r\bar{\w}(r)dr\bigg)^2 \mathfrak{u}^{\e}dx  \leq  \|\mathfrak{u}^{\e}\|_{L^2_x}^2,
\enda
\end{align*}
since $\frac{1}{|x|^3}\big(\int_0^{|x|}r\bar{\w}(r)dr\big)^2 \in L^2_x$ due to the choice of the radial vorticity $\bar{\w}(|x|) \sim |x|^{\mathrm{N}}$ near $|x| = 0$, as specified in \eqref{rw-def}.
We also derive the energy estimates for $\rho^{\e}$ and $\ta^{\e}$ from \eqref{locconNew}, as follows:
\begin{align}\label{rtaL2-E}
\bega
\frac{1}{2}\frac{d}{dt}\|\rho^{\e}\|_{L^2_x}^2 &= \frac{1}{2}\int_{\Omega} (\nabla_x \cdot u^{\e}) |\rho^{\e}|^2 dx -\frac{1}{\e} \int_{\Omega} (\nabla_x \cdot u^{\e}) \rho^{\e} dx , \cr
\frac{1}{2}\frac{d}{dt}\frac{3}{2}\|\ta^{\e}\|_{L^2_x}^2 &= \frac{1}{2}\int_{\Omega} (\nabla_x \cdot u^{\e}) \frac{3}{2}|\ta^{\e}|^2 dx -\frac{1}{\e} \int_{\Omega} (\nabla_x \cdot u^{\e}) \ta^{\e} dx \cr 
&\quad - \frac{1}{\e^2}\int_{\Omega}\bigg(\frac{1}{k_B\mathrm{P}^{\e}\mathrm{\Theta}^{\e} }\sum_{j} \p_{x_j} \mathfrak{q}_j^{\e} + \frac{1}{k_B\mathrm{P}^{\e}\mathrm{\Theta}^{\e} }\sum_{i,j} \p_{x_i}\mathrm{U}^{\e}_j \mathbf{r}_{ij}^{\e}\bigg) \ta^{\e} dx.
\enda
\end{align}
Combining the energy estimates $\frac{1}{k_B} \eqref{muL2}$ and \eqref{rtaL2-E} yields
\begin{align*}
\bega
\frac{1}{2}\frac{d}{dt}&\bigg(\|\rho^{\e}\|_{L^2_x}^2+\frac{1}{k_B}\|\mathfrak{u}^{\e}\|_{L^2_x}^2+\frac{3}{2}\|\ta^{\e}\|_{L^2_x}^2\bigg) \leq C\|\nabla_x\cdot u^{\e}\|_{L^\infty_x} \Big(\|\rho^{\e}\|_{L^2_x}^2+\frac{3}{2}\|\ta^{\e}\|_{L^2_x}^2\Big) \cr 
&+ C \Big(1+\|\nabla_x\cdot\mathfrak{u}^{\e}\|_{L^\infty_x}+\|\nabla_x \bar{u}\|_{L^\infty_x}\Big)\frac{1}{k_B}\|\mathfrak{u}^{\e}\|_{L^2_x}^2 \cr 
&-\frac{1}{\e} \int_{\Omega}(\mathrm{\Theta}^{\e}-1)\nabla_x(\rho^{\e}+\ta^{\e}) \mathfrak{u}^{\e}dx 
-\frac{1}{\e^2}\int_{\Omega}\frac{1}{k_B\mathrm{P}^{\e}}\sum_{j} \p_{x_j} \mathbf{r}_{ij}^{\e}  \mathfrak{u}^{\e}dx \cr 
&- \frac{1}{\e^2}\int_{\Omega}\bigg(\frac{1}{k_B\mathrm{P}^{\e}\mathrm{\Theta}^{\e} }\sum_{j} \p_{x_j} \mathfrak{q}_j^{\e} + \frac{1}{k_B\mathrm{P}^{\e}\mathrm{\Theta}^{\e} }\sum_{i,j} \p_{x_i}\mathrm{U}^{\e}_j \mathbf{r}_{ij}^{\e}\bigg) \ta^{\e} dx.
\enda
\end{align*}
Here, we decomposed $\mathrm{\Theta}^{\e} = 1 + (\mathrm{\Theta}^{\e} - 1)$ in the first term on the second line of \eqref{muL2}, and then used the following cancellation, relying on the identities $\nabla_x \cdot \widetilde{u}^{\e} = 0$ and $\nabla_x \cdot u^{\e} = \nabla_x \cdot \mathbb{P}^{\perp} u^{\e}$:
\begin{align*}
\bega
&\frac{1}{\e}\int_{\Omega}(\nabla_x\cdot u^{\e}) \rho^{\e}dx + \frac{1}{\e}\int_{\Omega}\nabla_x(\rho^{\e}+\ta^{\e})\mathfrak{u}^{\e} dx + \frac{1}{\e}\int_{\Omega}(\nabla_x\cdot u^{\e}) \ta^{\e}dx \cr 
&=\frac{1}{\e}\int_{\Omega}(\nabla_x\cdot u^{\e}) \rho^{\e}dx - \frac{1}{\e}\int_{\Omega}(\rho^{\e}+\ta^{\e})\nabla_x\cdot(\widetilde{u}^{\e}+\mathbb{P}^{\perp}u^{\e}) dx + \frac{1}{\e}\int_{\Omega}(\nabla_x\cdot u^{\e}) \ta^{\e}dx = 0.
\enda
\end{align*}
Since $\nabla_x \cdot \mathfrak{u}^{\e} = \nabla_x \cdot \mathbb{P}^{\perp} u^{\e} = \nabla_x \cdot u^{\e}$, we apply the H\"{o}lder inequality to obtain
\begin{align*}
\bega
\frac{1}{2}\frac{d}{dt}&\bigg(\|\rho^{\e}\|_{L^2_x}^2+\frac{1}{k_B}\|\mathfrak{u}^{\e}\|_{L^2_x}^2+\frac{3}{2}\|\ta^{\e}\|_{L^2_x}^2\bigg) \cr 
&\leq C \Big(1+\|\nabla_x\cdot u^{\e}\|_{L^\infty_x}+\|\nabla_x \bar{u}\|_{L^\infty_x}\Big) \bigg(\|\rho^{\e}\|_{L^2_x}^2+\frac{1}{k_B}\|\mathfrak{u}^{\e}\|_{L^2_x}^2+\frac{3}{2}\|\ta^{\e}\|_{L^2_x}^2\bigg) \cr 
&+\frac{\|\mathrm{\Theta}^{\e}-1\|_{L^2_x}}{\e} \|\nabla_x(\rho^{\e}+\ta^{\e})\|_{L^\infty_x} \|\mathfrak{u}^{\e}\|_{L^2_x} +\mathfrak{S}_W^{0}.
\enda
\end{align*}
This yields the desired result.
\hide
\begin{remark}
For the finite velocity energy case, we have 
\begin{align}
\bega
&\p_tu^{\e} + u^{\e}\cdot \nabla_x u^{\e} +\frac{k_B\mathrm{\Theta}^{\e}}{\e}\nabla_x(\rho^{\e}+\ta^{\e})
+\frac{1}{\e^2}\frac{1}{\mathrm{P}^{\e}}\sum_{j} \p_{x_j} \mathbf{r}_{ij}^{\e} =0, 
\enda
\end{align}
\begin{align*}
\bega
\frac{1}{2}\frac{d}{dt}\|u^{\e}\|_{L^2_x}^2 &\leq \frac{1}{2}\|\nabla_x\cdot u^{\e}\|_{L^\infty_x}\|u^{\e}\|_{L^2_x}^2-\frac{1}{\e}k_B \int_{\Omega}\nabla_x(\rho^{\e}+\ta^{\e}) u^{\e}dx \cr 
&
-\frac{k_B}{\e} \int_{\Omega}(\mathrm{\Theta}^{\e}-1)\nabla_x(\rho^{\e}+\ta^{\e}) \mathfrak{u}^{\e}dx 
-\frac{1}{\e^2}\int_{\Omega}\frac{1}{\mathrm{P}^{\e}}\sum_{j} \p_{x_j} \mathbf{r}_{ij}^{\e}  \mathfrak{u}^{\e}dx.
\enda
\end{align*}
\begin{align*}
\bega
\frac{1}{2}\frac{d}{dt}\bigg(\|\rho^{\e}\|_{L^2_x}^2+\frac{1}{k_B}\|u^{\e}\|_{L^2_x}^2+\frac{3}{2}\|\ta^{\e}\|_{L^2_x}^2\bigg) &\leq \frac{1}{2}\|\nabla_x\cdot u^{\e}\|_{L^\infty_x} \bigg(\|\rho^{\e}\|_{L^2_x}^2+\frac{1}{k_B}\|u^{\e}\|_{L^2_x}^2+\frac{3}{2}\|\ta^{\e}\|_{L^2_x}^2\bigg) \cr 
&+\|\nabla_x(\rho^{\e}+\ta^{\e})\|_{L^\infty_x} \frac{\|\mathrm{\Theta}^{\e}-1\|_{L^2_x}}{\e}  \|u^{\e}\|_{L^2_x} +\mathfrak{S}_W^{0}(t),
\enda
\end{align*}
where $\mathfrak{S}_W^{0}(t)$ is defined in \eqref{Amadef0}.
\end{remark}
\unhide
\end{proof}

\subsubsection{Key cancellation within the momentum-flux alignment}

Next, we present the crucial cancellation at the leading order between the momentum-flux alignment contributions $\mathfrak{S}_G^{\al}(t)$ and $\mathfrak{S}_W^{\al}(t)$, which arise from $\AC{\P}F^{\e}$ and $\bW^{\e}$, respectively, and are defined in \eqref{Amidef} and \eqref{Amadef}.
In the proof, for brevity, we slightly abuse notation by writing $\mathcal{E}(t)$, and $\mathcal{D}(t)$ for $\mathcal{E}^{\mathrm{N}}(F^{\e}(t))$, and $\mathcal{D}^{\mathrm{N}}(F^{\e}(t))$, respectively.

\begin{lemma}\label{L.cancel} For both the purely spatial derivative case \eqref{caseA} and the space-time derivative case \eqref{caseB}, we have the following leading-order cancellation between the momentum-flux alignment terms $\mathfrak{S}_G^{\al}(t)$ and $\mathfrak{S}_W^{\al}(t)$ defined in \eqref{Amidef} and \eqref{Amadef}, respectively:
\begin{enumerate}
\item For the finite velocity energy case \eqref{caseEC}, we have
\begin{align*}
\bega
|\mathfrak{S}_G^{\al}(t) + \mathfrak{S}_W^{\al}(t)| \leq C\e \mathcal{E}_{tot}^{\mathrm{N}}(F^{\e}(t)) (\mathcal{D}_{tot}^{\mathrm{N}}(F^{\e}(t)))^{\frac{1}{2}}, \quad \mbox{for} \quad  0\leq|\al|\leq \mathrm{N},
\enda
\end{align*}
for some constant $C > 0$.
\item For the infinite velocity energy case \eqref{caseECX}, we have
\begin{align*}
\bega
&|\mathfrak{S}_G^{\al}(t) + \mathfrak{S}_W^{\al}(t)| \leq C\e \mathcal{E}_{tot}^{\mathrm{N}}(F^{\e}(t)) (\mathcal{D}_{tot}^{\mathrm{N}}(F^{\e}(t)))^{\frac{1}{2}},  & \hspace{-5mm} 1\leq|\al|\leq \mathrm{N}, \cr 
&|\mathfrak{S}_G^{0}(t) + \mathfrak{S}_W^0(t)| \leq C\kappa^{\frac{1}{2}}\|\nabla_x\bar{u}(t)\|_{L^2_x}(\mathcal{D}_G^{\mathrm{N}}(F^{\e}(t)))^{\frac{1}{2}}+ C\e \mathcal{E}_{tot}^{\mathrm{N}}(F^{\e}(t)) (\mathcal{D}_{tot}^{\mathrm{N}}(F^{\e}(t)))^{\frac{1}{2}},  &\hspace{-5mm} |\al|=0,
\enda
\end{align*}
for some constant $C > 0$.
\end{enumerate}
\end{lemma}
\begin{proof}
(i) For the finite velocity energy case \eqref{caseEC}, we note that $\bar{u} = 0$. We decompose $\mathfrak{S}_W^{\al}(t)$, as defined in \eqref{Amadef}, into the leading-order term $\mathfrak{S}_{W,1}^{\al}(t)$ and the remainder term $\mathfrak{S}_{W,2}^{\al}(t)$ for $0 \leq |\al| \leq \mathrm{N}$, as follows:
\begin{align*}
\bega
\mathfrak{S}_{W,1}^{\al}(t) &:= -\frac{1}{\e^2} \int_{\Omega} \frac{1}{k_B\mathrm{P}^{\e} \mathrm{\Theta}^{\e}}\sum_{i,j} \p^{\al}\p_{x_j} \mathbf{r}_{ij}^{\e} \p^{\al}u_i dx -\frac{1}{\e^2}\int_{\Omega} \frac{1}{k_B\mathrm{P}^{\e}\mathrm{\Theta}^{\e}} \sum_{j} \p^{\al}\p_{x_j} \mathfrak{q}_j^{\e} \p^{\al}\ta dx, \cr 
\mathfrak{S}_{W,2}^{\al}(t) &:= -\frac{1}{\e^2} \int_{\Omega} \sum_{1\leq\beta\leq\al}\binom{\al}{\beta}\p^{\beta}\bigg(\frac{1}{k_B\mathrm{P}^{\e} \mathrm{\Theta}^{\e}}\bigg)\sum_{i,j} \p^{\al-\beta}\p_{x_j} \mathbf{r}_{ij}^{\e} \p^{\al}u_i dx \cr 
&\quad -\frac{1}{\e^2}\int_{\Omega}\sum_{1\leq\beta\leq\al}\binom{\al}{\beta}\p^{\beta}\bigg(\frac{1}{k_B\mathrm{P}^{\e}\mathrm{\Theta}^{\e}}\bigg) \sum_{j} \p^{\al-\beta}\p_{x_j} \mathfrak{q}_j^{\e} \p^{\al}\ta dx \cr 
&\quad  -\frac{1}{\e^2}\int_{\Omega} \p^{\al}\bigg[\frac{1}{k_B\mathrm{P}^{\e}\mathrm{\Theta}^{\e}} \Big(\sum_{i,j} \p_{x_i}\mathrm{U}^{\e}_j \mathbf{r}_{ij}^{\e}\Big)\bigg] \p^{\al}\ta dx.
\enda
\end{align*}
For the estimate of the remainder $\mathfrak{S}_{W,2}^{\al}$, we observe that the maximum number of derivatives acting on the terms $\p^{\al-\beta} \p_{x_j} \mathbf{r}_{ij}^{\e}$ and $\p^{\al-\beta} \p_{x_j} \mathfrak{q}_j^{\e}$ is at most $\mathrm{N}$, since $1 \leq |\beta|$.  
Therefore, applying \eqref{pTaleq} and \eqref{ABGscale}, we obtain
\begin{align}\label{RW2est}
\bega
|\mathfrak{S}_{W,2}^{\al}(t)| &\les \e\mathcal{E}_{M}^{\frac{1}{2}} \frac{1}{\e^2} \sum_{0\leq|\al|\leq \mathrm{N}}\sum_{i,j}\|\p^{\al}\mathbf{r}_{ij}^{\e}\|_{L^2_x} \|\p^{\al}u^{\e}_i\|_{L^2_x} +\e\mathcal{E}_{M}^{\frac{1}{2}} \frac{1}{\e^2} \sum_{0\leq|\al|\leq \mathrm{N}}\sum_{j} \|\p^{\al} \mathfrak{q}_j^{\e}\|_{L^2_x} \|\p^{\al}\ta^{\e}\|_{L^2_x} \cr 
&+\frac{1}{\e^2}\sum_{0\leq|\al|\leq \mathrm{N}}\bigg\|\p^{\al}\bigg[\frac{1}{k_B\mathrm{P}^{\e}\mathrm{\Theta}^{\e}} \Big(\sum_{i,j} \p_{x_i}\mathrm{U}^{\e}_j \mathbf{r}_{ij}^{\e}\Big)\bigg]\bigg\|_{L^2_x} \|\p^{\al}\ta^{\e}\|_{L^2_x} \cr 
&\leq \e\kappa^{\frac{1}{2}} \mathcal{E}_{M}\mathcal{D}_G^{\frac{1}{2}} + \frac{1}{\e^2}\e \Big(\mathcal{E}_M^{\frac{1}{2}}+\kappa^{-\frac{1}{2}}\mathcal{E}_{top}^{\frac{1}{2}}\Big)\Big(\e^2\kappa^{\frac{1}{2}}\mathcal{D}_G^{\frac{1}{2}}\Big)\mathcal{E}_M^{\frac{1}{2}}.
\enda
\end{align}
Here, we used \eqref{rutscale} for $0 \leq |\al| \leq \mathrm{N} + 1$ in the third line.  
Similarly, for the term $\mathfrak{S}_G^{\al}(t)$ defined in \eqref{Amidef}, we decompose it into the leading-order term $\mathfrak{S}_{G,1}^{\al}(t)$ and the remainder term $\mathfrak{S}_{G,2}^{\al}(t)$ for $0 \leq |\al| \leq \mathrm{N}$, after applying integration by parts to the derivative $\p_i$ in the terms $\p_i u^{\e}_j$ and $\p_i \ta^{\e}$:
\begin{align*}
\bega
\mathfrak{S}_{G,1}^{\al}(t)&:=\frac{1}{\e^2}  \int_{\Omega} \frac{1}{k_B\mathrm{\Theta}^{\e}}\sum_{i,j} \p^{\al} \p_i\mathbf{r}_{ij}^{\e} \p^{\al}u^{\e}_j dx +\frac{1}{\e^2}  \int_{\Omega} \frac{1}{k_B\mathrm{\Theta}^{\e}}\sum_i \p^{\al}\p_i\mathfrak{q}_i^{\e} \p^{\al}\ta^{\e} dx , \cr 
\mathfrak{S}_{G,2}^{\al}(t)&:=\frac{1}{\e^2}  \int_{\Omega} \sum_{i,j} \p_i\frac{1}{k_B\mathrm{\Theta}^{\e}} \p^{\al} \mathbf{r}_{ij}^{\e} \p^{\al}u^{\e}_j dx +\frac{1}{\e^2}  \int_{\Omega} \sum_i \p_i\frac{1}{k_B\mathrm{\Theta}^{\e}} \p^{\al}\mathfrak{q}_i^{\e} \p^{\al}\ta^{\e} dx \cr 
&-\frac{1}{\e^2}  \int_{\Omega} \frac{1}{k_B\mathrm{\Theta}^{\e}}\sum_{i,k} \bigg(\sum_{1\leq\beta\leq \al}\binom{\al}{\beta} \p^{\beta}\mathrm{U}^{\e}_k \p^{\al-\beta}\mathbf{r}_{ik}^{\e}\bigg) \p^{\al}\p_{x_i}\ta^{\e} dx.
\enda
\end{align*}
Proceeding in the same manner as in \eqref{RW2est}, we obtain the estimate
\begin{align}\label{RG2est}
\bega
|\mathfrak{S}_{G,2}^{\al}(t)|&\les \e\kappa^{\frac{1}{2}} \mathcal{E}_{M}\mathcal{D}_G^{\frac{1}{2}} + \frac{1}{\e^2}\Big(\e\mathcal{E}_M^{\frac{1}{2}}\Big)\Big(\e^2\kappa^{\frac{1}{2}}\mathcal{D}_G^{\frac{1}{2}}\Big)\Big(\mathcal{E}_M^{\frac{1}{2}}+\kappa^{-\frac{1}{2}}\mathcal{E}_{top}^{\frac{1}{2}}\Big).
\enda
\end{align}
Next, we consider the combination of the leading-order terms $\mathfrak{S}_{W,1}^{\al}(t)$ and $\mathfrak{S}_{G,1}^{\al}(t)$:
\begin{align}\label{R+R}
\bega
|\mathfrak{S}_{G,1}^{\al}(t)+\mathfrak{S}_{W,1}^{\al}(t)|&=\frac{1}{\e^2}  \int_{\Omega} \bigg(1-\frac{1}{\mathrm{P}^{\e}}\bigg)\frac{1}{k_B\mathrm{\Theta}^{\e}}\bigg(\sum_{i,j} \p^{\al} \p_i\mathbf{r}_{ij}^{\e} \p^{\al}u^{\e}_j +\sum_i  \p^{\al}\p_i\mathfrak{q}_i^{\e}\p^{\al}\ta^{\e}\bigg) dx \cr 
&\leq \frac{1}{\e^2}\Big(\e\mathcal{E}_M^{\frac{1}{2}}\Big)\Big(\e^2\kappa^{\frac{1}{2}}\mathcal{D}_G^{\frac{1}{2}} + \e^2\mathcal{D}_{top}^{\frac{1}{2}}\Big)\mathcal{E}_M^{\frac{1}{2}}.
\enda
\end{align}
Here, we used \eqref{ABGscale} for $0 \leq |\al| \leq \mathrm{N} + 1$, as well as \eqref{rutscale}, and the bound $|1 - \tfrac{1}{\mathrm{P}^{\e}}| \leq \e \mathcal{E}_M^{\frac{1}{2}}$.  
Combining the estimates \eqref{RW2est}, \eqref{RG2est}, and \eqref{R+R}, and applying the notations $\mathcal{E}_{tot}$ and $\mathcal{D}_{tot}$ defined in \eqref{EDtotdef}, we obtain the desired result. \\
(ii) We observe that for $1 \leq |\al| \leq \mathrm{N}$, the terms $\mathfrak{S}_G^{\al}(t)$ and $\mathfrak{S}_W^{\al}(t)$ satisfy the same estimate as established in the proof of part (i) of this lemma.
Furthermore, in the case $|\al| = 0$, the second lines of both $\mathfrak{S}_G^{0}(t)$ and $\mathfrak{S}_W^{0}(t)$ exhibit the same cancellation structure.  
Therefore, it remains only to estimate the first lines of $\mathfrak{S}_G^{0}(t)$ and $\mathfrak{S}_W^{0}(t)$.
Following the same approach as in the proof of part (i) of this lemma, we decompose the leading-order term and apply integration by parts to obtain
\begin{align*}
\bega
|\mathfrak{S}_G^{0}(t) + \mathfrak{S}_W^0(t)| &\les \frac{1}{\e^2}\int_{\Omega}\frac{1}{k_B\mathrm{P}^{\e}}\sum_{i,j} \bigg(\mathbf{r}_{ij}^{\e}  \p_{x_j}(u^{\e}-\bar{u}) - \mathbf{r}_{ij}^{\e} \p_ju_i^{\e}\bigg) dx  + \e \mathcal{E}_{tot} \mathcal{D}_{tot}^{\frac{1}{2}} \cr 
&\les \frac{1}{\e^2}\|\nabla_x\bar{u}\|_{L^2_x}\|\mathbf{r}_{ij}^{\e}\|_{L^2_x}+ \e \mathcal{E}_{tot} \mathcal{D}_{tot}^{\frac{1}{2}} \cr 
&\les \kappa^{\frac{1}{2}}\|\nabla_x\bar{u}\|_{L^2_x}\mathcal{D}_G^{\frac{1}{2}}+ \e \mathcal{E}_{tot} \mathcal{D}_{tot}^{\frac{1}{2}},
\enda
\end{align*}
where we used \eqref{ABGscale}. 
\end{proof}

We now estimate $I_W^{\al}$, $II_W^{\al}$, and $III_W^{\al}$,
defined in~\eqref{EW12def}. To this end, we recall from
Lemma~\ref{L.ept} the relations $\e \p_t u^{\e} \sim \nabla_x (\rho^{\e} + \ta^{\e})$ and $\e \p_t \ta^{\e} \sim \nabla_x \cdot u^{\e}$.


\begin{lemma}\label{L.EW12} 
Let the initial data $u_0^{\e}$ satisfy either the finite velocity energy condition \eqref{caseEC} or the infinite velocity energy condition \eqref{caseECX}.
For $t\in[0,T]$ satisfying the bootstrap assumption \eqref{condition}, and for both the purely spatial derivatives \eqref{caseA} and the space-time derivatives \eqref{caseB}, the terms $I_W^{\al}(t)$, $II_W^{\al}(t)$ and $III_{W}^{\al}(t)$ defined in \eqref{EW12def} satisfy the following upper bounds:
\begin{align*}
\bega
\sum_{1\leq|\al|\leq\mathrm{N}}&I_{W}^{\al}(t) \les \|\nabla_x\cdot u(t)\|_{L^\infty_x}\mathcal{E}_M(t) + \e\big(\mathcal{E}_M^{\frac{3}{2}}(t)+\mathcal{E}_M^2(t)\big) + \e \|\bPhi^{\e}_{\bW}\|_{L^\infty_x}\mathcal{E}_M(t), \cr 
\sum_{1\leq|\al|\leq\mathrm{N}}&II_{W}^{\al}(t) \les \|\nabla_x\bW^{\e}(t)\|_{L^\infty_x} \mathcal{E}_M(t) + \e\big(\mathcal{E}_M^{\frac{3}{2}}(t)+\mathcal{E}_M^2(t)\big) \cr 
&+ \begin{cases} 0 , \quad & \mbox{for} \quad \eqref{caseA}, \\ 
\e^{\mathfrak{n}}  \|\bPhi^{\e}_{\bW}(t)\|_{L^\infty_x}\mathcal{E}_M(t), \quad & \mbox{for} \quad \eqref{caseB}, 
\end{cases} 
\cr 
\sum_{1\leq|\al|\leq\mathrm{N}}&III_{W}^{\al}(t) \les \Big(\|\nabla_x\ta^{\e}(t)\|_{L^\infty_x} + \|\nabla_x(\rho^{\e}+\ta^{\e})(t)\|_{L^\infty_x}\Big) \mathcal{E}_M(t) \cr 
&+ \e\big(\mathcal{E}_M^{\frac{3}{2}}(t)+\mathcal{E}_M^2(t)\big) + \e \mathcal{E}_M(t) \sum_{0\leq|\al|\leq\mathrm{N}-1}\|\p^{\al}\bPhi^{\e}_{\bW}(t)\|_{L^2_x} \cr 
&+ \begin{cases}
0 , \quad & \mbox{for} \quad \eqref{caseA}, \\ 
\e^{\mathfrak{n}-1}\Big(\|\nabla_x\nabla_x\cdot u^{\e}(t)\|_{L^\infty_x}\mathcal{E}_M(t)  +\e^2\mathcal{E}_{M}^\frac{1}{2}(t)\|\bPhi^{\e}_{\bW}(t)\|_{L^\infty_x}\|\bPhi^{\e}_{\bW}(t)\|_{H^{\mathrm{N}-1}_x} \Big), \quad & \mbox{for} \quad \eqref{caseB}, \end{cases} 
\enda
\end{align*}
for $\mathfrak{n} \geq 1$, where $\mathcal{E}_M(t)$ is defined in \eqref{N-EDdef}.
\end{lemma}
\begin{proof}
We note that for the finite velocity energy case \eqref{caseEC}, the energy includes the $L^2_x$ norm of $u^{\e}$ since $\bar{u} = 0$.  
On the other hand, in the infinite velocity energy case \eqref{caseECX}, we cannot bound $\|u^{\e}\|_{L^2_x}$, as $\|\bar{u}\|_{L^2_x} = \infty$.  
Therefore, in proving this lemma, we use $\|u^{\e}\|_{L^\infty_x}$ to avoid the divergence of $\|u^{\e}\|_{L^2_x}$ in the infinite velocity energy case.
Depending on whether the finite velocity energy case \eqref{caseEC} or the infinite velocity energy case \eqref{caseECX} applies, we estimate $\|u^{\e}\|_{L^\infty_x}$ as follows:
\begin{align}\label{uLinfEC}
\bega
\|u^{\e}(t)\|_{L^\infty_x} \leq 
\begin{cases}
C\mathcal{E}_M^{\frac{1}{2}}(t), \quad &\mbox{for finite velocity energy case } 
\eqref{caseEC}, \\ 
C\big(\mathcal{E}_{M}^{\frac{1}{2}}(t) +1\big), \quad &\mbox{for infinite velocity energy case } \eqref{caseECX}, \end{cases}
\enda
\end{align}
for some constant $C>0$.
Here, we used Agmon's inequality \eqref{Agmon} for the finite velocity energy case \eqref{caseEC}:
\begin{align*}
\bega
\|u^{\e}\|_{L^\infty_x} \leq \|u^{\e}\|_{H^2_x} \leq \mathcal{E}_M^{\frac{1}{2}}, \quad d=2,3,
\enda
\end{align*}
and also Agmon's inequality \eqref{Agmon} for the infinite velocity energy case \eqref{caseECX}:
\begin{align*}
\bega
\|u^{\e}\|_{L^\infty_x} \leq \|u^{\e}-\bar{u}\|_{L^\infty_x} +\|\bar{u}\|_{L^\infty_x} \leq  \|u^{\e}-\bar{u}\|_{H^2_x} +\|\bar{u}\|_{L^\infty_x} \leq \|\widetilde{u}^{\e}+\mathbb{P}^{\perp}u^{\e}\|_{H^2_x} +\|\bar{u}\|_{L^\infty_x} \leq \mathcal{E}_{M}^{\frac{1}{2}} +\|\bar{u}\|_{L^\infty_x},
\enda
\end{align*}
where we also used Lemma \ref{L.barubdd}.
\\(Proof of the estimate for $I_W^{\al}$)  
To prove the estimate for $I_W^{\al}$, we claim that
\begin{align}\label{claimE0W}
\bega
\| \nabla_x \cdot (tr(\mathcal{S}^{\e})u^{\e})\|_{L^\infty_x} &\les \|\nabla_x\cdot u^{\e}\|_{L^\infty_x} + \e \|u^{\e}\|_{L^\infty_x} \mathcal{E}_M^{\frac{1}{2}}, \cr 
\|\p_t\mathcal{S}^{\e}\|_{L^\infty_x} &\les  \|\nabla_x\cdot u^{\e}\|_{L^\infty_x} +\e \big( \|u^{\e}\|_{L^\infty_x} \mathcal{E}_M^{\frac{1}{2}} + \mathcal{E}_M\big) + \e \|\bPhi^{\e}_{\bW}\|_{L^\infty_x}.
\enda
\end{align}
If the claim \eqref{claimE0W} holds, then the result follows directly from the bound $\|\p^{\al} \bW^{\e}\|_{L^2_x}^2 \leq \mathcal{E}_M$ for $1 \leq |\al| \leq \mathrm{N}$ together with \eqref{uLinfEC}.
From the definition of $\mathcal{S}^{\e}$ in \eqref{Sdef}, we observe that
$tr(\mathcal{S}^{\e}) = \frac{5}{2} + \frac{d}{k_B\mathrm{\Theta}^{\e}}$. Since $\p_i (\mathrm{\Theta}^{\e})^{-1} = -\e \p_i \ta^{\e} (\mathrm{\Theta}^{\e})^{-1}$, we have $\|\nabla_x(tr(\mathcal{S}^{\e}))\|_{L^\infty_x}\leq \|\e\nabla_x\ta^{\e}\|_{L^\infty_x}\|\frac{d}{k_B\mathrm{\Theta}^{\e}}\|_{L^\infty_x} \leq C\e\|\nabla_x\ta^{\e}\|_{L^\infty_x}$ and  $\|tr(\mathcal{S}^{\e})\|_{L^\infty_x}\leq C$ by \eqref{pTaleq} and \eqref{condition}. Thus, we get
\begin{align*}
\| \nabla_x \cdot (tr(\mathcal{S}^{\e})u^{\e})\|_{L^\infty_x} &\les \|tr(\mathcal{S}^{\e})\|_{L^\infty_x} \|\nabla_x\cdot u^{\e}\|_{L^\infty_x} + \|\e\nabla_x\ta^{\e}\|_{L^\infty_x}\|u^{\e}\|_{L^\infty_x}  \cr 
&\les \|\nabla_x\cdot u^{\e}\|_{L^\infty_x}+ \e\|u^{\e}\|_{L^\infty_x}\mathcal{E}_{M}^{\frac{1}{2}},
\end{align*}
where we used $\|\nabla_x\ta^{\e}\|_{L^\infty_x} \leq \mathcal{E}_M^{\frac{1}{2}}$ by \eqref{rutscale}. This establishes the first inequality in \eqref{claimE0W}.
To prove the second inequality in \eqref{claimE0W}, we differentiate $\mathcal{S}^{\e}$ and its inverse using the definition \eqref{Sdef}, which gives
\begin{align}\label{pS}
\p^{\al}\mathcal{S}^{\e} = \p^{\al}\bigg(\frac{1}{k_B\mathrm{\Theta}^{\e}}\bigg) \diag\big(0,1,1,1,0\big), \qquad \p^{\al}(\mathcal{S}^{\e})^{-1} = \p^{\al}(k_B\mathrm{\Theta}^{\e}) \diag\big(0,1,1,1,0\big) ,
\end{align}
by definition of $\mathcal{S}^{\e}$ in \eqref{Sdef}. Applying \eqref{condition} and the estimate \eqref{eptineq}, we obtain
\begin{align*}
\bega
\|\p_t\mathcal{S}^{\e}\|_{L^\infty_x} = \bigg\|(\e\p_t\ta)\frac{1}{k_B\mathrm{\Theta}^{\e}}\bigg\|_{L^\infty_x} &\leq \|\e\p_t\ta\|_{L^\infty_x} \bigg\|\frac{1}{k_B\mathrm{\Theta}^{\e}}\bigg\|_{L^\infty_x} \cr 
&\leq \e \big( \|u^{\e}\|_{L^\infty_x} \mathcal{E}_M^{\frac{1}{2}} + \mathcal{E}_M\big) + \|\nabla_x\cdot u^{\e}\|_{L^\infty_x} + \e \|\bPhi^{\e}_{\bW}\|_{L^\infty_x}.
\enda
\end{align*}
This gives the proof of the second inequality in \eqref{claimE0W}. \\
(Proof of the estimate of $II_{W}^{\al}$)
We decompose $II_{W}^{\al}$ into the following two parts:
\begin{align}\label{EW2decomp}
II_{W}^{\al} = \sum_{1\leq \beta \leq \al} \sum_{0\leq \gamma\leq\beta} \Big\|\p^{\gamma}\mathcal{S}^{\e} \p^{\beta-\gamma}u^{\e} \cdot \nabla_x \p^{\al-\beta}\bW^{\e} \Big\|_{L^2_x} \|\p^{\al}\bW^{\e}\|_{L^2_x} = 
(R_{2,1}+R_{2,2}) \|\p^{\al}\bW^{\e}\|_{L^2_x},
\end{align}
where
\begin{align*}
R_{2,1} &= \sum_{1\leq \beta \leq \al}  \Big\| \mathcal{S}^{\e} \p^{\beta}u^{\e} \cdot \nabla_x \p^{\al-\beta}\bW^{\e} \Big\|_{L^2_x}, \qquad 
R_{2,2} = \sum_{1\leq \beta \leq \al} \sum_{1\leq \gamma\leq\beta} \Big\| \p^{\gamma}\mathcal{S}^{\e} \p^{\beta-\gamma}u^{\e} \cdot \nabla_x \p^{\al-\beta}\bW^{\e} \Big\|_{L^2_x}.
\end{align*}
We first estimate the term $R_{2,1}$. Using the bound $\|\mathcal{S}^{\e}\|_{L^\infty_x} \leq C$ from \eqref{condition}, we have
\begin{align*}
R_{2,1} &\les \sum_{1\leq \beta \leq \al} \|\mathcal{S}^{\e}\|_{L^\infty_x}\Big\|\p^{\beta}u^{\e} \cdot \nabla_x \p^{\al-\beta}\bW^{\e} \Big\|_{L^2_x} \les \sum_{1\leq \beta \leq \al} \Big\|\p^{\beta}u^{\e} \cdot \nabla_x \p^{\al-\beta}\bW^{\e} \Big\|_{L^2_x}.
\end{align*}
(Case A) When $\al_0 = 0$, applying \eqref{uvHk}, we get 
\begin{align*}
R_{2,1} &\les \|\nabla_xu^{\e}\|_{L^\infty_x} \|\nabla_x\bW^{\e}\|_{H^{\mathrm{N}-1}_x} + \|\nabla_x\bW^{\e}\|_{L^\infty_x} \|\nabla_xu^{\e}\|_{H^{\mathrm{N}-1}_x} \leq \|\nabla_x\bW^{\e}\|_{L^\infty_x} \mathcal{E}_M^{\frac{1}{2}}.
\end{align*}
(Case B) The scaled-time derivative $\p_{\tilde{t}}$ can act on either $u^{\e}$ or $\bW^{\e}$. Let $\p^{\bullet}$ denote either $\p_{\tilde{t}}$ or $\nabla_x$. Then we estimate $R_{2,1}$ as follows:
\begin{align*}
R_{2,1} \les \sum_{0\leq|\beta_x+\gamma_x|\leq \mathrm{N}-1}\|(\p^{\beta_x}\p^{\bullet}u^{\e}) (\p^{\gamma_x}\nabla_x\bW^{\e}) \|_{L^2_x} + \sum_{0\leq|\beta_x+\gamma_x|\leq \mathrm{N}-1} \| (\p^{\beta_x}\nabla_xu^{\e}) (\p^{\gamma_x}\p^{\bullet}\bW^{\e}) \|_{L^2_x} .
\end{align*}
Since the total number of derivatives is $\mathrm{N} + 1$, the upper bound for $R_{2,1}$ corresponds to the term \eqref{zzL2} in Lemma \ref{L.E1}. Therefore, we obtain the same upper bounds as in \eqref{zzL2} and \eqref{MMAB}:
\begin{align}\label{R21est}
\bega
R_{2,1} \leq \begin{cases}
\|\nabla_x\bW^{\e}\|_{L^\infty_x} \mathcal{E}_M^{\frac{1}{2}}, \quad & \mbox{for} \quad \eqref{caseA},  \\
\|\nabla_x\bW^{\e}\|_{L^\infty_x} \mathcal{E}_M^{\frac{1}{2}} + 
\e^{\mathfrak{n}-1}  \Big(\e \big( \|u^{\e}\|_{L^\infty_x} \mathcal{E}_M^{\frac{1}{2}} + \mathcal{E}_M\big) + \e \|\bPhi^{\e}_{\bW}\|_{L^\infty_x}\mathcal{E}_M^{\frac{1}{2}} \Big), \quad & \mbox{for} \quad \eqref{caseB},
\end{cases}
\enda
\end{align}
for any $\mathfrak{n} \geq 1$, where we used the fact that $\|\nabla_x \cdot u^{\e}\|_{L^\infty_x}$ and $\|\nabla_x (\rho^{\e} + \ta^{\e})\|_{L^\infty_x}$ can be absorbed into $\|\nabla_x \bW^{\e}\|_{L^\infty_x}$.
For the estimate of $R_{2,2}$, we decompose the terms according to whether at least one derivative falls on $u^{\e}$ or not:
\begin{align*}
R_{2,2} = \sum_{1\leq \beta \leq \al} \sum_{1\leq \gamma<\beta} \Big\| \p^{\gamma}\mathcal{S}^{\e} \p^{\beta-\gamma}u^{\e} \cdot \nabla_x \p^{\al-\beta}\bW^{\e} \Big\|_{L^2_x} + \sum_{1\leq \beta \leq \al} \Big\| \p^{\beta}\mathcal{S}^{\e} u^{\e} \cdot \nabla_x \p^{\al-\beta}\bW^{\e} \Big\|_{L^2_x}.
\end{align*}
Using $|\p^{\al}\mathcal{S}^{\e}|\leq|\p^{\al}(k_B\mathrm{\Theta}^{\e})^{-1}|$, together with the estimates \eqref{pTa}, \eqref{pTaleq}, and Agmon's inequality \eqref{Agmon}, we obtain
\begin{align}\label{R22est}
R_{2,2} \les \e\mathcal{E}_M^{\frac{3}{2}} + \e \|u^{\e}\|_{L^\infty_x}\mathcal{E}_M.
\end{align}
Combining the estimates \eqref{R21est} and \eqref{R22est} into the decomposition \eqref{EW2decomp}, and using $\|\p^{\al} \bW^{\e}\|_{L^2_x} \leq \mathcal{E}_M^{\frac{1}{2}}$ together with \eqref{uLinfEC}, we conclude the proof of the desired estimate. 
\\(Proof of the estimate for $III_W^{\al}$)  
We define the leading $L^2_x$-norm of $III_{W}^{\al}(t)$ by $R_3 $, as follows:
\begin{align*}
III_{W}^{\al}(t) &= R_3\|\p^{\al}\bW^{\e}\|_{L^2_x}, \qquad R_3:=\sum_{1\leq\beta\leq\al}\Big\|\p^{\beta}\mathcal{S}^{\e} \p_t\p^{\al-\beta}\bW^{\e}\Big\|_{L^2_x}.
\end{align*}
By the definition of \( \mathcal{S}^{\e} \) in \eqref{Sdef} and the expression \eqref{pS}, and since \( 1 \leq |\beta| \), we have
\begin{align*}
R_{3} &= \sum_{1\leq\beta\leq\al}\bigg\|\frac{1}{\e}\p^{\beta}\frac{1}{k_B\mathrm{\Theta}^{\e}} \p^{\al-\beta}(\e\p_tu^{\e})\bigg\|_{L^2_x}. 
\end{align*}
(Case A) Using \eqref{condition} and the estimate \eqref{uvHk}, we get
\begin{align*}
R_{3}&\les \bigg(\Big\|\frac{1}{\e}\nabla_x\frac{1}{k_B\mathrm{\Theta}^{\e}}\Big\|_{L^\infty_x} \big\|\e\p_t u^{\e}\big\|_{H^{\mathrm{N}-1}_x}+\Big\|\frac{1}{\e}\nabla_x\frac{1}{k_B\mathrm{\Theta}^{\e}}\Big\|_{H^{\mathrm{N}-1}_x} \big\|\e\p_t u^{\e}\big\|_{L^\infty_x}\bigg) .
\end{align*}
To estimate \( \frac{1}{\e} \nabla_x \frac{1}{k_B \mathrm{\Theta}^{\e}} \), we apply \eqref{pTa} to have
\begin{align}\label{pxta}
\bigg\|\frac{1}{\e}\nabla_x\frac{1}{k_B\mathrm{\Theta}^{\e}}\bigg\|_{H^{\mathrm{N}-1}_x} \les \mathcal{E}_{M}^{\frac{1}{2}}, \qquad \bigg\|\frac{1}{\e}\nabla_x\frac{1}{k_B\mathrm{\Theta}^{\e}}\bigg\|_{L^\infty_x} \les \|\nabla_x\ta^{\e}\|_{L^\infty_x},
\end{align}
where we also used the bound from \eqref{condition}. Therefore, by Lemma \ref{L.ept} and the estimate \eqref{pxta}, we conclude
\begin{align}\label{R3-A}
R_{3} &\les \Big(\|\nabla_x\ta^{\e}\|_{L^\infty_x} + \|\nabla_x(\rho^{\e}+\ta^{\e})\|_{L^\infty_x}\Big) \mathcal{E}_M^{\frac{1}{2}} + \e\|u^{\e}\|_{L^\infty_x}\mathcal{E}_M + \e \mathcal{E}_M^{\frac{1}{2}} \Big(\|\bPhi^{\e}_{\bW}\|_{L^\infty_x}+\|\bPhi^{\e}_{\bW}\|_{H^{\mathrm{N}-1}_x}\Big),
\end{align}
where we have used $\|\nabla_x(\rho^{\e}+\ta^{\e})\|_{H^{\mathrm{N}-1}_x}\leq \mathcal{E}_M^{\frac{1}{2}}$.\\
(Case B) We only need to consider the case when $\p^{\al}$ includes the scaled time derivative $\p_{\tilde{t}} = \e^{\mathfrak{n}}\p_t$ defined in \eqref{mindex}. There are two subcases to analyze: when $\p_{\tilde{t}}$ acts on $\frac{1}{k_B\mathrm{\Theta}^{\e}}$, and when it acts on $\e\p_t u^{\e}$.
\begin{align*}
R_{3}&\les \Big\|\frac{1}{\e}\p_{\tilde{t}}\frac{1}{k_B\mathrm{\Theta}^{\e}} (\e\p_t u^{\e})\Big\|_{H^{\mathrm{N}-1}_x}+\Big\|\frac{1}{\e}\nabla_x\frac{1}{k_B\mathrm{\Theta}^{\e}} \p_{\tilde{t}}(\e\p_t u^{\e})\Big\|_{H^{\mathrm{N}-2}_x} := R_{3,1} + R_{3,2} .
\end{align*}
For the estimate of $R_{3,1}$, we claim that
\begin{align}\label{R31-BC}
\bega
R_{3,1} 
&\les \e^{\mathfrak{n}-1}\bigg[\big(\|\nabla_x\cdot u^{\e}\|_{L^\infty_x}+\|\nabla_x(\rho^{\e}+\ta^{\e})\|_{L^\infty_x}\big)\mathcal{E}_M^{\frac{1}{2}}\bigg] \cr 
& + \e^{\mathfrak{n}-1} \bigg[ \e\|u^{\e}\|_{L^\infty_x}\mathcal{E}_M + \e^2\|u^{\e}\|_{L^\infty_x}^2\mathcal{E}_M +\e\|\bPhi^{\e}_{\bW}\|_{H^{\mathrm{N}-1}_x} \mathcal{E}_M^{\frac{1}{2}} + \e^2 \|\bPhi^{\e}_{\bW}\|_{L^\infty_x}\|\bPhi^{\e}_{\bW}\|_{H^{\mathrm{N}-1}_x}\bigg] .
\enda
\end{align}
Using \eqref{uvHk} and $\p\frac{1}{k_B\mathrm{\Theta}^{\e}} = -\e\p\ta^{\e} \frac{1}{k_B\mathrm{\Theta}^{\e}}$, we get 
\begin{align*}
R_{3,1} &\les \Big\|\frac{1}{k_B\mathrm{\Theta}^{\e}}\p_{\tilde{t}}\ta^{\e}\Big\|_{L^\infty_x} \|(\e\p_t u^{\e})\|_{H^{\mathrm{N}-1}_x} + \Big\|\frac{1}{k_B\mathrm{\Theta}^{\e}}\p_{\tilde{t}}\ta^{\e}\Big\|_{H^{\mathrm{N}-1}_x} \|(\e\p_t u^{\e})\|_{L^\infty_x}.
\end{align*}
Applying \eqref{uvHk} and \eqref{pTaleq}, we have
\begin{align*}
\Big\|\frac{1}{k_B\mathrm{\Theta}^{\e}}\p_{\tilde{t}}\ta^{\e}\Big\|_{H^{\mathrm{N}-1}_x} &\les \Big\|\frac{1}{k_B\mathrm{\Theta}^{\e}}\Big\|_{H^{\mathrm{N}-1}_x} \|\p_{\tilde{t}}\ta^{\e}\|_{L^\infty_x} + \Big\|\frac{1}{k_B\mathrm{\Theta}^{\e}}\Big\|_{L^\infty_x} \|\p_{\tilde{t}}\ta^{\e}\|_{H^{\mathrm{N}-1}_x} \les \|\p_{\tilde{t}}\ta^{\e}\|_{H^{\mathrm{N}-1}_x},
\end{align*}
where we used Agmon's inequality \eqref{Agmon}
$\|\cdot\|_{L^\infty_x} \leq \|\cdot\|_{H^{2}_x}$ and $\mathrm{N} > d/2 + 1$. Therefore, we obtain 
\begin{align*}
R_{3,1} &\les \|\p_{\tilde{t}}\ta^{\e}\|_{L^\infty_x} \|(\e\p_t u^{\e})\|_{H^{\mathrm{N}-1}_x} + \|\p_{\tilde{t}}\ta^{\e}\|_{H^{\mathrm{N}-1}_x} \|(\e\p_t u^{\e})\|_{L^\infty_x}.
\end{align*}
Then, applying Lemma \ref{L.ept} and using $\|\nabla_x(\rho^{\e}+\ta^{\e})\|_{H^{\mathrm{N}-1}_x},\|\nabla_x\cdot u^{\e}\|_{H^{\mathrm{N}-1}_x} \les \mathcal{E}_M^{\frac{1}{2}}$, the bound $\e \mathcal{E}_M^{\frac{1}{2}}\leq C$ from \eqref{condition}, and the Sobolev embedding $\|\cdot\|_{L^\infty_x} \leq \|\cdot\|_{H^{\mathrm{N}-1}_x}$, we obtain the desired estimate \eqref{R31-BC}. \\
For the estimate of $R_{3,2}$, we claim that 
\begin{align}\label{R32-BC}
R_{3,2}&\les \e^{\mathfrak{n}-1}\big(\|\nabla_x\ta^{\e}\|_{L^\infty_x}+\|\nabla_x\nabla_x\cdot u^{\e}\|_{L^\infty_x}\big)\mathcal{E}_M^{\frac{1}{2}} \cr 
&+ \bigg(\e^{\mathfrak{n}-1}\Big(\e (\|u^{\e}\|_{L^\infty_x}+\mathcal{E}_M^{\frac{1}{2}})\mathcal{E}_M^{\frac{1}{2}}\Big) + \e \sum_{0\leq|\al|\leq\mathrm{N}-1}\|\p^{\al}\bPhi^{\e}_{\bW}\|_{L^2_x}\bigg)\mathcal{E}_M^{\frac{1}{2}}.
\end{align}
Using \eqref{uvHk}, we also get 
\begin{align}\label{R32leq}
R_{3,2} &\leq \Big\|\frac{1}{\e}\nabla_x\frac{1}{k_B\mathrm{\Theta}^{\e}}\Big\|_{L^\infty_x} \|\p_{\tilde{t}}(\e\p_t u^{\e})\|_{H^{\mathrm{N}-2}_x} + \Big\|\frac{1}{\e}\nabla_x\frac{1}{k_B\mathrm{\Theta}^{\e}}\Big\|_{H^{\mathrm{N}-2}_x} \|\p_{\tilde{t}}(\e\p_t u^{\e})\|_{L^\infty_x}.
\end{align}
Following the same approach as in Lemma \ref{L.ept}, and using the local conservation law $\eqref{locconNew}_2$, we find
\begin{align}\label{ttu}
\bega
\|\p_{\tilde{t}}(\e\p_t u^{\e})\| &= \bigg\| \p_{\tilde{t}} \bigg(\e u^{\e}\cdot \nabla_x u^{\e} +k_B\mathrm{\Theta}^{\e}\nabla_x(\rho^{\e}+\ta^{\e})
+\frac{1}{\e}\frac{1}{\mathrm{P}^{\e}}\sum_{j} \p_{x_j} \mathbf{r}_{ij}^{\e} \bigg) \bigg\| \cr 
&\les \e (\|u^{\e}\|_{L^\infty_x}+\mathcal{E}_M^{\frac{1}{2}})\mathcal{E}_M^{\frac{1}{2}} + \|\p_{\tilde{t}}\nabla_x(\rho^{\e}+\ta^{\e})\| + \e \|\p_{\tilde{t}}\bPhi^{\e}_{\bW}\|,
\enda
\end{align}
for either the $\|\cdot\|_{L^\infty_x}$ or $\|\cdot\|_{H^{\mathrm{N}-2}_x}$ norms.
We use the local conservation laws $\eqref{locconNew}_1$ and $\eqref{locconNew}_3$ to estimate the term $\nabla_x(\rho^{\e} + \ta^{\e})$:
\begin{align}\label{pt.rta}
\bega
\|\p_{\tilde{t}}\nabla_x(\rho^{\e}+\ta^{\e})\| \les \e^{\mathfrak{n}-1}\Big(\e (\|u^{\e}\|_{L^\infty_x}+\mathcal{E}_M^{\frac{1}{2}})\mathcal{E}_M^{\frac{1}{2}} + \|\nabla_x\nabla_x\cdot u^{\e}\| + \e \|\nabla_x\bPhi^{\e}_{\bW}\|\Big), 
\enda
\end{align}
again for $\|\cdot\|_{L^\infty_x}$ or $\|\cdot\|_{H^{\mathrm{N}-2}_x}$.
Substituting \eqref{pt.rta} into \eqref{ttu}, we obtain
\begin{align}\label{ttu2}
\bega
\|\p_{\tilde{t}}(\e\p_t u^{\e})\| &\les  \e^{\mathfrak{n}-1}\Big(\e (\|u^{\e}\|_{L^\infty_x}+\mathcal{E}_M^{\frac{1}{2}})\mathcal{E}_M^{\frac{1}{2}} + \|\nabla_x\nabla_x\cdot u^{\e}\| + \e \|\nabla_x\bPhi^{\e}_{\bW}\|\Big) \cr 
& + \e (\|u^{\e}\|_{L^\infty_x}+\mathcal{E}_M^{\frac{1}{2}})\mathcal{E}_M^{\frac{1}{2}}  + \e \|\p_{\tilde{t}}\bPhi^{\e}_{\bW}\|.
\enda
\end{align}
By applying \eqref{ttu2} and \eqref{pxta} to \eqref{R32leq}, we obtain the desired bound \eqref{R32-BC}.
Combining \eqref{R3-A}, \eqref{R31-BC}, and \eqref{R32-BC}, and then applying \eqref{condition}, \eqref{uLinfEC}, together with the bound $\|\p^{\al}\bW^{\e}\|_{L^2_x} \leq \mathcal{E}_M^{\frac{1}{2}}$ in the definition of $III_{W}^{\al}$, yields the desired result.
\end{proof}

\begin{proof}[\textbf{Proof of Proposition \ref{P.EW}}] 
We first go back to the case $|\al|=0$ in Lemma \ref{L.rutL2w}. 
Since 
\begin{align*}
\bega
\bigg(\|\rho^{\e}(t)\|_{L^2_x}^2+\frac{1}{k_B}\|u^{\e}-\bar{u}\|_{L^2_x}^2+\frac{3}{2}\|\ta^{\e}(t)\|_{L^2_x}^2\bigg) \leq C\mathcal{E}_M(t),
\enda
\end{align*}
we can write Lemma \ref{L.rutL2w} as
\begin{align}\label{PEW2-a}
\bega
&\frac{1}{2}\frac{d}{dt}\bigg(\|\rho^{\e}\|_{L^2_x}^2+\frac{1}{k_B}\|u^{\e}-\bar{u}\|_{L^2_x}^2+\frac{3}{2}\|\ta^{\e}\|_{L^2_x}^2\bigg) \cr 
&\leq C \Big(1+\|\nabla_x\cdot u^{\e}\|_{L^\infty_x} +\|\nabla_x(\rho^{\e}+\ta^{\e})\|_{L^\infty_x} + \|\nabla_x\bar{u}\|_{L^\infty_x} \Big) \mathcal{E}_M(t) +\mathfrak{S}_W^{0}(t),
\enda
\end{align}
where we used $\frac{\|\mathrm{\Theta}^{\e}-1\|_{L^2_x}}{\e}\leq \mathcal{E}_M^{\frac{1}{2}}$ and $\mathfrak{S}_W^{0}(t)$ defined in \eqref{Amadef0}.
For $1\leq|\al|\leq\mathrm{N}$, from Lemma \ref{L.EW}, taking summation $\sum_{1\leq|\al|\leq\mathrm{N}}$ to \eqref{EWineq} gives 
\begin{align*}
\frac{d}{dt}\sum_{1\leq|\al|\leq\mathrm{N}}\int_{\Omega} &(\p^{\al}\bW^{\e})^T \mathcal{S}^{\e} (\p^{\al}\bW^{\e}) dx \leq \sum_{1\leq|\al|\leq\mathrm{N}}\mathfrak{S}_W^{\al}(t) + C \sum_{1\leq|\al|\leq\mathrm{N}} \bigg[I_{W}^{\al} + II_{W}^{\al} + III_{W}^{\al} \bigg],
\end{align*}
where $I_{W}^{\al}$, $II_{W}^{\al}$, and $III_{W}^{\al}$ are defined in \eqref{EW12def}.
We combine the estimates of $I_{W}^{\al}$, $II_{W}^{\al}$, and $III_{W}^{\al}$ in Lemma \ref{L.EW12} to have
\begin{align}\label{PEW2-b}
\bega
&\frac{d}{dt}\sum_{1\leq|\al|\leq\mathrm{N}}\int_{\Omega} (\p^{\al}\bW^{\e})^T \mathcal{S}^{\e} (\p^{\al}\bW^{\e}) dx \leq \sum_{1\leq|\al|\leq\mathrm{N}}\mathfrak{S}_W^{\al}(t) +\|\nabla_x\bW^{\e}(t)\|_{L^\infty_x} \mathcal{E}_M(t) \cr
& + \e\big(\mathcal{E}_M^{\frac{3}{2}}(t)+\mathcal{E}_M^2(t)\big) + \e \mathcal{E}_M(t) \sum_{0\leq|\al|\leq\mathrm{N}-1}\|\p^{\al}\bPhi^{\e}_{\bW}(t)\|_{L^2_x} + \mathcal{Z}_W^{time}(t),
\enda
\end{align}
where $\mathcal{Z}_W^{time}(t)$ is defined in \eqref{Bmadef}. 
For the term $\p^{\al}\bPhi^{\e}_{\bW}$, combining \eqref{ABGscale}$_2$ and \eqref{ABGscale}$_1$ with \eqref{condition} yields  
\begin{align}\label{Phiest}
\|\bPhi^{\e}_{\bW}(t)\|_{L^\infty_x} \leq C\kappa^{\frac{1}{2}}\mathcal{D}_G^{\frac{1}{2}}, \qquad \sum_{0\leq|\al|\leq\mathrm{N}-1}\|\p^{\al}\bPhi^{\e}_{\bW}(t)\|_{L^2_x} \leq C\kappa^{\frac{1}{2}}\mathcal{D}_G^{\frac{1}{2}}.
\end{align}
For the left-hand side of \eqref{PEW2-b}, by the definition of $\mathcal{S}^{\e}$ in \eqref{Sdef} and the bound $\sup_{t \in [0,T]}|\mathrm{\Theta}^{\e}(t,x) -1| \ll 1$ from \eqref{condition}, we get 
\begin{align}\label{Ssimil}
\bega
\frac{1}{C}\|\p^{\al}(\rho^{\e},u^{\e},\ta^{\e})(t)\|_{L^2_x}^2 \leq \int_{\Omega} (\p^{\al}\bW^{\e})^T \mathcal{S}^{\e} (\p^{\al}\bW^{\e}) dx \leq C\|\p^{\al}(\rho^{\e},u^{\e},\ta^{\e})(t)\|_{L^2_x}^2,
\enda
\end{align}
for $1\leq|\al|\leq\mathrm{N}$.
Combining \eqref{PEW2-a} and \eqref{PEW2-b}, and then integrating in time, we apply \eqref{Phiest} and \eqref{Ssimil} to obtain the desired result \eqref{EWresult}.

\end{proof}

\subsection{Acoustic dispersion}\label{Sec.macro.Stri}

In this section, we estimate $\mathbb{P}^{\perp}u^{\e}$ and $\rho^{\e} + \ta^{\e}$, which solve a dispersion system:
\begin{lemma}\label{L.locconP}
For $(\rho^{\e},u^{\e},\ta^{\e})$ satisfying the local conservation laws \eqref{locconNew},
the quantities $\rho^{\e}+\ta^{\e}$ and $\mathbb{P}^{\perp}u^{\e}$ satisfy the following equations:
\begin{align}\label{locconP}
\bega
&\p_t \big(\rho^{\e}+\ta^{\e}\big) + \frac{5}{3}\frac{1}{\e}\nabla_x\cdot u^{\e} =\bPhi^{\e}_{\rho+\ta}, \qquad \p_t \mathbb{P}^{\perp}u^{\e}  +\frac{k_B}{\e}\nabla_x(\rho^{\e}+\ta^{\e})=\bPhi^{\e}_{\mathbb{P}^{\perp}u},
\enda
\end{align}
where
\begin{align}\label{bPhilocconP}
\bega
\bPhi^{\e}_{\rho+\ta}(t,x)&:=-u^{\e}\cdot\nabla_x\big(\rho^{\e}+\ta^{\e}\big)-\frac{2}{3}\frac{1}{\e^2}\frac{1}{k_B\mathrm{P}^{\e}\mathrm{\Theta}^{\e} }\sum_{j} \p_{x_j} \mathfrak{q}_j^{\e} - \frac{2}{3}\frac{1}{\e^2}\frac{1}{k_B\mathrm{P}^{\e}\mathrm{\Theta}^{\e} }\sum_{i,j} \p_{x_i}\mathrm{U}^{\e}_j \mathbf{r}_{ij}^{\e}, \cr 
\bPhi^{\e}_{\mathbb{P}^{\perp}u}(t,x)&:=\mathbb{P}^{\perp}\bigg[-u^{\e}\cdot \nabla_x u^{\e} + \frac{k_B(1-\mathrm{\Theta}^{\e})}{\e}\nabla_x(\rho^{\e}+\ta^{\e})- \frac{1}{\e^2}\frac{1}{\mathrm{P}^{\e}}\sum_{j} \p_{x_j} \mathbf{r}_{ij}^{\e}\bigg].
\enda
\end{align}
Here, the Leray projector $\mathbb{P}$ is defined in \eqref{LerayPdef}.
\end{lemma}
\begin{proof}
We obtain the first equation of \eqref{locconP} by adding $\eqref{locconNew}_1$ and $\frac{2}{3}\eqref{locconNew}_3$. To derive the second equation of \eqref{locconP}, we multiply equation $\eqref{locconNew}_2$ by $k_B\mathrm{\Theta}^{\e}$:
\begin{align*}
\bega
&\p_tu^{\e} + u^{\e}\cdot \nabla_x u^{\e} +\frac{k_B\mathrm{\Theta}^{\e}}{\e}\nabla_x(\rho^{\e}+\ta^{\e})
+\frac{1}{\e^2}\frac{1}{\mathrm{P}^{\e}}\sum_{j} \p_{x_j} \mathbf{r}_{ij}^{\e} =0.
\enda
\end{align*}
Next, we decompose $\mathrm{\Theta}^{\e} = 1 + (\mathrm{\Theta}^{\e} - 1)$ in front of $\nabla_x(\rho^{\e}+\ta^{\e})$, and then apply the Leray projector $\mathbb{P}^{\perp}$ to obtain the second equation of \eqref{locconP}.
\end{proof}


\hide

잘 요약해서 인트로로
Different from the low Mach number limit, where the condition $\|\e\p_t u^{\e}_0\|_{L^2_x} \leq C$ does not necessarily require a well-prepared initial data, the scaled Boltzmann equation requires the initial data $\|\e\p_t F^{\e}_0 |M_0|^{-\frac{1}{2}}\|_{L^2_{x,v}} \leq C$ to be well-prepared due to the singular scaling of the Knudsen number $\textbf{Kn} = \e\kappa$. Therefore, in the context of the incompressible Euler limit from the Boltzmann equation, it is preferable to avoid using time derivatives. From this perspective, the dispersion relation derived from the linearized equation \eqref{varsigeqn} in \cite{Danchin} is more suitable than the wave equation approach in \cite{JMR}, since the Strichartz estimate for \eqref{varsigeqn} does not require time derivatives of the initial data.

The author in \cite{Danchin} proved that if $\varsigma_1$ and $\varsigma_2$ satisfy the following equation,
\begin{align}\label{varsigeqn}
\bega
&\p_t \varsigma_1 + \e^{-1}|D|\varsigma_2 = \bPhi_{\varsigma_1}, \cr 
&\p_t \varsigma_2 - \e^{-1}|D|\varsigma_1 = \bPhi_{\varsigma_2}, \cr 
&(\varsigma_1, \varsigma_2)|_{t=0} = (\varsigma_1(0), \varsigma_2(0)),
\enda
\end{align}
where $|D|$ is the Fourier multiplier defined in \eqref{fourierDdef}, 
then $\varsigma_1$ and $\varsigma_2$ can be bounded in terms of the scale $\e$ via Strichartz estimates.  
The main idea is that the dispersion relation of \eqref{varsigeqn} matches that of the wave equation associated with the oscillations $e^{\pm it|\xi|}$.

\unhide


\begin{proposition}[Proposition 2.2 of \cite{Danchin1}] \label{P.Book}
Let $(\varsigma_1,\varsigma_2)$ be the solution of 
\begin{align}\label{varsigeqn}
\bega
&\p_t \varsigma_1 + \e^{-1}|D|\varsigma_2 = \bPhi_{\varsigma_1}, \qquad \varsigma_1|_{t=0} = \varsigma_1(0),  \cr 
&\p_t \varsigma_2 - \e^{-1}|D|\varsigma_1 = \bPhi_{\varsigma_2}, \qquad \varsigma_2|_{t=0} = \varsigma_2(0),
\enda
\end{align}
where $|D|$ is the Fourier multiplier defined in \eqref{fourierDdef}.

Then, for any $s \in \mathbb{R}$, we have 
\begin{align*}
\bega
\|(\varsigma_1,\varsigma_2)\|_{L^r_T\dot{B}_{p,1}^{s+d(\frac{1}{p}-\frac{1}{2})+\frac{1}{r}}} \leq C \e^{\frac{1}{r}}\|(\varsigma_1(0),\varsigma_2(0))\|_{\dot{B}_{2,1}^s} + C\e^{\frac{1}{r}}\|(\bPhi_{\varsigma_1},\bPhi_{\varsigma_2})\|_{L^1_T\dot{B}_{2,1}^{s}},
\enda
\end{align*}
where
\begin{align}\label{rpd-def}
\bega
2 \leq p \leq \infty, \qquad \frac{2}{r} \leq \min \left\{ 1, (d-1)\bigg(\frac{1}{2}-\frac{1}{p}\bigg)\right\}, \qquad (r,p,d) \neq (2,\infty,3).
\enda
\end{align}
Here, the Besov norm is defined in \eqref{Besovdef}.
\end{proposition}

The proof can be found in \cite{Danchin,Danchin1}; in particular, a detailed proof is given in Proposition 10.30 of \cite{Danchin}. For the reader's convenience, we sketch the argument in Appendix \ref{A.Stri}.

We now apply this Strichartz estimate to the system \eqref{locconP} to control the lower frequencies.
Recall that $\mathrm{N}+1$ denotes the highest order of derivatives appearing in our energy–dissipation framework \eqref{N-EDdef2}.
\begin{proposition}\label{P.div.ineq}
For any real number $s \in (0, \mathrm{N})$, and for $d = 2, 3$, there exists a positive constant $C > 0$ such that 
\begin{align}\label{div.ineq}
\bega
\Big\|&\Big((\rho^{\e}+\ta^{\e}),\mathbb{P}^{\perp}u^{\e}\Big)\Big\|_{L^{r}_T\dot{B}_{p,1}^{s+d(\frac{1}{p}-\frac{1}{2})+\frac{1}{r}}} \leq C \e^{\frac{1}{r}}\Big\|\Big((\rho^{\e}_0+\ta^{\e}_0),(\mathbb{P}^{\perp}u^{\e}_0)\Big)\Big\|_{\dot{B}_{2,1}^s} \cr 
&+ \begin{cases}
C\e^{\frac{1}{r}}\int_0^T\Big(\mathcal{E}_{tot}^{\mathrm{N}}(
F^\e(\tau)
) +\kappa^{\frac{1}{2}}\Big(\mathcal{D}_{tot}^{\mathrm{N}}(F^\e(\tau))\Big)^{\frac{1}{2}}\Big) d\tau, \quad &0\leq s < \mathrm{N}-1, \\
C\e^{\frac{1}{r}}\kappa^{-\frac{\mathfrak{\eta}}{2}}\int_0^T\Big(\mathcal{E}_{tot}^{\mathrm{N}}(F^\e(\tau)) +\kappa^{\frac{1}{2}}\Big(\mathcal{D}_{tot}^{\mathrm{N}}(F^\e(\tau))\Big)^{\frac{1}{2}}\Big) d\tau, \quad & s = \mathrm{N}-1, \\
C\e^{\frac{1}{r}}\kappa^{-\frac{s+1-\mathrm{N}}{2}}\int_0^T\Big(\mathcal{E}_{tot}^{\mathrm{N}}(F^\e(\tau)) +\kappa^{\frac{1}{2}}\Big(\mathcal{D}_{tot}^{\mathrm{N}}(F^\e(\tau))\Big)^{\frac{1}{2}}\Big) d\tau, \quad &\mathrm{N}-1 < s < \mathrm{N},
\end{cases} 
\enda
\end{align}
where $(r,p,d)$ satisfies \eqref{rpd-def}, and for any $0<\eta\ll1$.
Here, the Fourier multiplier $|D|$, the Leray projector $\mathbb{P}$, the Besov norm, and the energy and dissipation terms $\mathcal{E}_{tot}^N$ and $\mathcal{D}_{tot}^N$ are defined in \eqref{fourierDdef}, \eqref{LerayPdef}, \eqref{Besovdef}, and \eqref{EDtotdef}, respectively.
\end{proposition}

\begin{remark}
The maximal scaling exponents are obtained as $\e^{\frac{1}{4}}$ 
by choosing $p=\infty$ for $d=2$, and as $\e^{\frac{1}{2}-}$ 
by letting $p\to\infty$ for $d=3$, respectively.
\end{remark}

\hide
\begin{lemma}\label{L.divinf} 
Let $d = 2,3$, and let $\mathrm{N}$ be given.
Then the following estimate holds for all $t \in (0,T)$:
\begin{align*}
\bega
\Big\|\nabla_x^{\ell}\Big((\rho^{\e}+\ta^{\e}),&\mathbb{P}^{\perp}u^{\e}\Big)\Big\|_{L^{r}_TL^\infty_x} \leq C \e^{\frac{1}{r}}\Big\|\Big((\rho^{\e}_0+\ta^{\e}_0),(\mathbb{P}^{\perp}u^{\e}_0)\Big)\Big\|_{\dot{B}_{2,1}^s} \cr 
&+ \begin{cases}
C\e^{\frac{1}{r}}\int_0^T\Big(\mathcal{E}_{tot}^{\mathrm{N}}(\tau) +\kappa^{\frac{1}{2}}\Big(\mathcal{D}_{tot}^{\mathrm{N}}(\tau)\Big)^{\frac{1}{2}}\Big) d\tau, \qquad &0\leq s < \mathrm{N}-1, \\
C\e^{\frac{1}{r}}\kappa^{-\frac{\mathfrak{\eta}}{2}}\int_0^T\Big(\mathcal{E}_{tot}^{\mathrm{N}}(\tau) +\kappa^{\frac{1}{2}}\Big(\mathcal{D}_{tot}^{\mathrm{N}}(\tau)\Big)^{\frac{1}{2}}\Big) d\tau, \quad \mbox{for any $0<\eta\ll1$} \qquad & s = \mathrm{N}-1, \\
C\e^{\frac{1}{r}}\kappa^{-\frac{s+1-\mathrm{N}}{2}}\int_0^T\Big(\mathcal{E}_{tot}^{\mathrm{N}}(\tau) +\kappa^{\frac{1}{2}}\Big(\mathcal{D}_{tot}^{\mathrm{N}}(\tau)\Big)^{\frac{1}{2}}\Big) d\tau, \qquad &\mathrm{N}-1 < s < \mathrm{N},
\end{cases} 
\enda
\end{align*}
where
\begin{align*}
\bega
s= \ell + \frac{d}{2} - \frac{1}{r} < \mathrm{N}.
\enda
\end{align*}
\end{lemma}
\begin{proof}
Applying the embedding of Besov spaces in \eqref{Besov-Lp} with $q=\infty$, we obtain
\begin{align*}
\bega
\Big\|\nabla_x^{\ell}\Big((\rho^{\e}+\ta^{\e}),\mathbb{P}^{\perp}u^{\e}\Big)\Big\|_{L^{r}_TL^\infty_x} \leq C \Big\|\nabla_x^{\ell}\Big((\rho^{\e}+\ta^{\e}),\mathbb{P}^{\perp}u^{\e}\Big)\Big\|_{L^{r}_T\dot{B}_{p,1}^{\frac{d}{p}}}.
\enda
\end{align*}
Choosing $s := \ell + \frac{d}{2} - \frac{1}{r}$ and applying Proposition~\ref{P.div.ineq} yield the desired estimate.
\end{proof}
\unhide

\begin{remark}\label{Rmk.divuLinf}
Applying the embedding of Besov spaces in \eqref{Besov-Lp} with $q=\infty$, for $\ell \in (-d/2+1/r, \mathrm{N}-d/2+1/r)$ we obtain
\begin{align*}
\bega
\Big\|\nabla_x^{\ell}&\Big((\rho^{\e}+\ta^{\e}),\mathbb{P}^{\perp}u^{\e}\Big)\Big\|_{L^{r}_TL^\infty_x} \leq C \e^{\frac{1}{r}}\Big\|\Big((\rho^{\e}_0+\ta^{\e}_0),(\mathbb{P}^{\perp}u^{\e}_0)\Big)\Big\|_{\dot{B}_{2,1}^s} \cr 
&+ \begin{cases}
C\e^{\frac{1}{r}}\int_0^T\Big(\mathcal{E}_{tot}^{\mathrm{N}}(F^\e(\tau)) +\kappa^{\frac{1}{2}}\Big(\mathcal{D}_{tot}^{\mathrm{N}}(F^\e(\tau))\Big)^{\frac{1}{2}}\Big) d\tau, \quad &0\leq s < \mathrm{N}-1, \\
C\e^{\frac{1}{r}}\kappa^{-\frac{\mathfrak{\eta}}{2}}\int_0^T\Big(\mathcal{E}_{tot}^{\mathrm{N}}(F^\e(\tau)) +\kappa^{\frac{1}{2}}\Big(\mathcal{D}_{tot}^{\mathrm{N}}(F^\e(\tau))\Big)^{\frac{1}{2}}\Big) d\tau,  \quad & s = \mathrm{N}-1, \\
C\e^{\frac{1}{r}}\kappa^{-\frac{s+1-\mathrm{N}}{2}}\int_0^T\Big(\mathcal{E}_{tot}^{\mathrm{N}}(F^\e(\tau)) +\kappa^{\frac{1}{2}}\Big(\mathcal{D}_{tot}^{\mathrm{N}}(F^\e(\tau))\Big)^{\frac{1}{2}}\Big) d\tau, \quad &\mathrm{N}-1 < s < \mathrm{N},
\end{cases} 
\enda
\end{align*}
for any $0<\eta\ll1$ when $s=\mathrm{N}-1$, where $s= (\ell + \frac{d}{2} - \frac{1}{r}) \in (0,\mathrm{N})$. In particular, we apply this Strichartz inequality to prove Theorem \ref{T.2D.global} in the case $(\mathrm{N}, r, d) = (4,4,2)$ and $0\leq \ell \leq 2$. 
\end{remark}

\begin{remark}
Depending on the spatial dimension and the top-order derivative $\mathrm{N}+1$, the singular $\kappa$-scaling differs.
Heuristically, examining the $\e$ and $\kappa$ scaling in Remark \ref{Rmk.divuLinf}, and choosing $\ell=2$ and $p\to\infty$, we obtain
\begin{align*}
\bega
\begin{cases} 
\|\nabla_x^2 \mathbb{P}^{\perp} u^{\e}\|_{L^{4}_TL^\infty_x} \les 
\e^{\frac{1}{4}}\kappa^{-\frac{3}{8}}, \quad &d=2,~\mathrm{N}=3, \\
\|\nabla_x^2 \mathbb{P}^{\perp} u^{\e}\|_{L^{4}_TL^\infty_x} \les \e^{\frac{1}{4}}, \quad &d=2,~\mathrm{N}=4, \\
\|\nabla_x^2 \mathbb{P}^{\perp} u^{\e}\|_{L^{2+}_TL^\infty_x} \les \e^{\frac{1}{2}-}\kappa^{-}, \quad &d=3,~\mathrm{N}=4.
\end{cases} 
\enda
\end{align*}
\end{remark}

\hide
\begin{remark} (If we use the Besov space $\dot{B}_{s+\frac{1}{\delta},2})$)
\begin{align*}
&\big\|\mathbb{P}^{\perp}u^{\e}\big\|_{L^{4+}_T\dot{B}_{2+\frac{1}{\delta},2}^{3-\frac{3}{4}}} \les C \e^{\frac{1}{4}-}\kappa^{-\frac{1}{2}}\int_0^T\Big(\mathcal{E}_{tot}(\tau) +\kappa^{\frac{1}{2}}\mathcal{D}_{tot}^{\frac{1}{2}}(\tau)\Big) d\tau , \quad \mbox{when} \quad s=3 \cr 
&\big\|\mathbb{P}^{\perp}u^{\e}\big\|_{L^{4+}_T\dot{B}_{2+\frac{1}{\delta},2}^{2-\frac{3}{4}}} \les C \e^{\frac{1}{4}-} \int_0^T\Big(\mathcal{E}_{tot}(\tau) +\kappa^{\frac{1}{2}}\mathcal{D}_{tot}^{\frac{1}{2}}(\tau)\Big) d\tau, \quad \mbox{when} \quad s=2
\end{align*}
\end{remark}
\unhide

In the proof, for brevity, we slightly abuse notation by writing $\mathcal{E}(t)$, and $\mathcal{D}(t)$ for $\mathcal{E}^{\mathrm{N}}(F^{\e}(t))$, and $\mathcal{D}^{\mathrm{N}}(F^{\e}(t))$, respectively.

\begin{proof}[\textbf{Proof of Proposition \ref{P.div.ineq}}]
Applying the operator $|D|^{-1}\div$ to the second equation of \eqref{locconP}, we get
\begin{align}\label{locconP2}
\bega
&\p_t \big(\rho^{\e}+\ta^{\e}\big) + \frac{5}{3}\frac{1}{\e}\nabla_x\cdot u^{\e} =\bPhi^{\e}_{\rho+\ta}, \quad \p_t |D|^{-1}\div\mathbb{P}^{\perp}u^{\e}  -\frac{k_B}{\e}|D|(\rho^{\e}+\ta^{\e})= |D|^{-1}\div \bPhi^{\e}_{\mathbb{P}^{\perp}u}.
\enda
\end{align}
Since $\nabla_x\cdot u^{\e} = |D|\big(|D|^{-1}\div \mathbb{P}^{\perp}u^{\e}\big)$, the system \eqref{locconP2} has the same structure as \eqref{varsigeqn}. 
Therefore, we can apply Proposition \ref{P.Book} to \eqref{locconP2}, identifying $\varsigma_1 = \rho^{\e}+\ta^{\e}$ and $\varsigma_2 = |D|^{-1}\div(\mathbb{P}^{\perp}u^{\e})$:
\begin{align}\label{bvineq}
\bega
\Big\|\Big((\rho^{\e}+\ta^{\e}),|D|^{-1}\div(\mathbb{P}^{\perp}u^{\e})\Big)\Big\|_{L^r_T\dot{B}_{p,1}^{s+d(\frac{1}{p}-\frac{1}{2})+\frac{1}{r}}} &\leq C \e^{\frac{1}{r}}\Big\|\Big((\rho^{\e}_0+\ta^{\e}_0),|D|^{-1}\div(\mathbb{P}^{\perp}u^{\e}_0)\Big)\Big\|_{\dot{B}_{2,1}^s} \cr 
&+ C\e^{\frac{1}{r}}\Big\|(\bPhi^{\e}_{\rho+\ta},|D|^{-1}\div \bPhi^{\e}_{\mathbb{P}^{\perp}u})\Big\|_{L^1_T\dot{B}_{2,1}^{s}},
\enda
\end{align}
for any $s\in\mathbb{R}$, where $\bPhi^{\e}_{\rho+\ta}$ and $\bPhi^{\e}_{\mathbb{P}^{\perp}u}$ are defined in \eqref{bPhilocconP}.
Since $|D|^{-1}\div$ is a homogeneous multiplier of degree $0$ and $\mathbb{P}^{\perp}u^{\e}$ is curl free, it suffices to estimates $\big\|(\bPhi^{\e}_{\rho+\ta},\bPhi^{\e}_{\mathbb{P}^{\perp}u})\big\|_{L^{1}_T\dot{B}_{2,1}^{s}}$.
We note that $\bPhi^{\e}_{\rho+\ta} \sim \p_{x_j} \mathfrak{q}_j^{\e}$, and we estimated $\|\p^{\al}\mathfrak{q}_j^{\e}\|_{L^2_x} \les \e^2\kappa^{\frac{1}{2}}\mathcal{D}_G^{\frac{1}{2}}$ for $0\leq|\al|\leq \mathrm{N}$ and $\|\p^{\al}\mathfrak{q}_j^{\e}\|_{L^2_x} \les \e^2\mathcal{D}_{top}^{\frac{1}{2}}$ for $|\al|=\mathrm{N}+1$ in \eqref{ABGscale}. Thus, we interpolate the Besov norm $\|\cdot\|_{\dot{B}_{2,1}^{s}}$ depending on the regularity parameter $0< s < \mathrm{N}$ and $\mathrm{N} \leq s < \mathrm{N}+1$.
For simplicity, we write $\bPhi^{\e}=\bPhi^{\e}_{\rho+\ta} $ or $|D|^{-1}\div \bPhi^{\e}_{\mathbb{P}^{\perp}u}$ only in this proof. 
Applying the interpolation inequality in \eqref{Besov-inter}, we have
\begin{align*}
\bega
\|\bPhi^{\e}\|_{\dot{B}_{2,1}^{s}} &\leq \begin{cases} 
\frac{C}{\mathrm{N}-1}\bigg(\frac{\mathrm{N}-1}{\mathrm{N}-s-1}+\frac{\mathrm{N}-1}{s}\bigg)\|\bPhi^{\e}\|_{\dot{B}_{2,\infty}^{0}}^{\frac{\mathrm{N}-s-1}{\mathrm{N}-1}} \|\bPhi^{\e}\|_{\dot{B}_{2,\infty}^{\mathrm{N}-1}}^{\frac{s}{\mathrm{N}-1}}, \qquad &0< s < \mathrm{N}-1, \\ 
\frac{C}{\eta}\|\bPhi^{\e}\|_{\dot{B}_{2,\infty}^{\mathrm{N}-1-\mathfrak{\eta}}}^{\frac{1}{2}} \|\bPhi^{\e}\|_{\dot{B}_{2,\infty}^{\mathrm{N}-1+\mathfrak{\eta}}}^{\frac{1}{2}}, \qquad &s= \mathrm{N}-1, \\ 
C\bigg(\frac{1}{\mathrm{N}-s}+\frac{1}{s+1-\mathrm{N}}\bigg)\|\bPhi^{\e}\|_{\dot{B}_{2,\infty}^{\mathrm{N}-1}}^{\mathrm{N}-s} \|\bPhi^{\e}\|_{\dot{B}_{2,\infty}^{\mathrm{N}}}^{s+1-\mathrm{N}}, \qquad &\mathrm{N}-1 < s < \mathrm{N},
\end{cases}
\enda
\end{align*}
where we choose $(s_1,s_2,\tau)= (0,\mathrm{N}-1,\frac{\mathrm{N}-s-1}{\mathrm{N}-1})$ for the first case and $(s_1,s_2,\tau)= (\mathrm{N}-1-\mathfrak{\eta},\mathrm{N}-1+\mathfrak{\eta},\frac{1}{2})$ with $0<\eta\ll1$ for the second case and $(s_1,s_2,\tau)= (\mathrm{N}-1,\mathrm{N},\mathrm{N}-s)$ for the third case.  
Since $\|\cdot\|_{\dot{B}^s_{2,\infty}} \leq \|\cdot\|_{\dot{B}^s_{2,2}} \les \|\cdot\|_{H^s_x}$, it follows from \eqref{Besov-H} that
\begin{align}\label{GGineq}
\bega
\|\bPhi^{\e}\|_{\dot{B}_{2,1}^{s}} &\leq \begin{cases} C\bigg(\frac{1}{\mathrm{N}-s-1}+\frac{1}{s}\bigg)\|\bPhi^{\e}\|_{L^2_x}^{\frac{\mathrm{N}-s-1}{\mathrm{N}-1}} \|\bPhi^{\e}\|_{H^{\mathrm{N}-1}_x}^{\frac{s}{\mathrm{N}-1}}, \qquad &0< s < \mathrm{N}-1, \\
\frac{C}{\eta}\|\bPhi^{\e}\|_{H_x^{\mathrm{N}-1-\mathfrak{\eta}}}^{\frac{1}{2}} \|\bPhi^{\e}\|_{H_x^{\mathrm{N}-1+\mathfrak{\eta}}}^{\frac{1}{2}}, \qquad &s= \mathrm{N}-1, \\ 
C\bigg(\frac{1}{\mathrm{N}-s}+\frac{1}{s+1-\mathrm{N}}\bigg)\|\bPhi^{\e}\|_{H_x^{\mathrm{N}-1}}^{\mathrm{N}-s} \|\bPhi^{\e}\|_{H_x^{\mathrm{N}}}^{s+1-\mathrm{N}}, \qquad &\mathrm{N}-1 < s < \mathrm{N}.
\end{cases}
\enda
\end{align}
Then, we estimate $\bPhi^{\e}_{\rho+\ta}$ in the $H^k$ norm for $0 \leq k \leq \mathrm{N}$. Using \eqref{uvHk}, we get
\begin{align*}
\bega
\big\|u^{\e}\cdot\nabla_x\big(\rho^{\e}+\ta^{\e}\big)\big\|_{H^k_x} &\leq  \|\nabla_xu^{\e}\|_{L^\infty_x}\|\nabla_x(\rho^{\e}+\ta^{\e})\|_{H^{k-1}_x} +\|\nabla_x(\rho^{\e}+\ta^{\e})\|_{L^\infty_x}\|\nabla_xu^{\e}\|_{H^{k-1}_x}  \cr 
&\quad + \|u^{\e}\|_{L^\infty_x}\|(\rho^{\e}+\ta^{\e})\|_{H^{k+1}_x}\cr 
&\leq \begin{cases}  C\mathcal{E}_{tot}, & \quad 0\leq k\leq\mathrm{N}-1, \\ 
C\kappa^{-\frac{1}{2}}\mathcal{E}_{tot},  & \quad k=\mathrm{N},
\end{cases}
\enda
\end{align*}
and
\begin{align*}
\bega
&\bigg\|\frac{1}{\e^2}\frac{1}{k_B\mathrm{P}^{\e}\mathrm{\Theta}^{\e} }\sum_{j} \p_{x_j} \mathfrak{q}_j^{\e} +\frac{1}{\e^2}\frac{1}{k_B\mathrm{P}^{\e}\mathrm{\Theta}^{\e} }\sum_{i,j} \p_{x_i}\mathrm{U}^{\e}_j \mathbf{r}_{ij}^{\e} \bigg\|_{H^k_x}  \leq \begin{cases} C \kappa^{\frac{1}{2}}\mathcal{D}_G^{\frac{1}{2}}, & \quad 0\leq k\leq\mathrm{N}-1, \\ 
C(\mathcal{D}_{tot}^{\frac{1}{2}}+\mathcal{E}_M),  & \quad k=\mathrm{N},
\end{cases}
\enda
\end{align*}
where we used $\eqref{ABGscale}_1$. Therefore, we obtain
\begin{align}\label{G12est}
\bega
\|\bPhi^{\e}_{\rho+\ta}\|_{H^k_x} &\les \begin{cases} C (\mathcal{E}_{tot} +\kappa^{\frac{1}{2}}\mathcal{D}_{tot}^{\frac{1}{2}}), & \quad 0\leq k\leq\mathrm{N}-1, \\ 
C \kappa^{-\frac{1}{2}}(\mathcal{E}_{tot} +\kappa^{\frac{1}{2}}\mathcal{D}_{tot}^{\frac{1}{2}}),  & \quad k=\mathrm{N}.
\end{cases}
\enda
\end{align}
For the estimate of $|D|^{-1}\div \bPhi^{\e}_{\mathbb{P}^{\perp}u}$, we observe that $\||D|^{-1}\div\cdot\|_{H^s_x} \leq \|\cdot\|_{H^s_x}$ and $\|\mathbb{P}^{\perp}\cdot\|_{H^s_x} \leq \|\cdot\|_{H^s_x}$, since $\widehat{|D|^{-1}\div g} = \frac{\xi}{|\xi|}\cdot\hat{g}$ and $\widehat{\mathbb{P}^{\perp}g} = \sum_j\frac{\xi_i\xi_j}{|\xi|^2}\hat{g}_j$. 
Then, applying $\|\frac{k_B(1-\mathrm{\Theta}^{\e})}{\e}\|_{H^k_x} \leq \mathcal{E}_{tot}^{\frac{1}{2}}$ from Lemma \ref{L.1/P}, we can readily conclude that $|D|^{-1}\div \bPhi^{\e}_{\mathbb{P}^{\perp}u}$ satisfies the same estimate as in \eqref{G12est}, following the same argument used for $\bPhi^{\e}_{\rho+\ta}$. 
Applying \eqref{GGineq} and \eqref{G12est} to \eqref{bvineq}, we obtain the desired result.
\end{proof}

\hide
\begin{lemma} Let $u$ is purely gradient part i.e. $u = \mathbb{P}^{\perp}u$. Then we have 
\begin{align*}
\bega
\frac{1}{C}\|\mathbb{P}^{\perp}u\|_{\dot{B}_{p,q}^{s}} \leq \||D|^{-1}\div \mathbb{P}^{\perp}u\|_{\dot{B}_{p,q}^{s}} &\leq C\|\mathbb{P}^{\perp}u\|_{\dot{B}_{p,q}^{s}}
\enda
\end{align*}
\end{lemma}
\begin{proof}
Since $u$ is purely gradient part, defining $g:= |D|^{-1}(\div u)$, and $R_j:=i \xi_j/|\xi|$ we can write 
\begin{align*}
\bega
\mathbb{P}^{\perp} u = \nabla_x |D|^{-2}\div u, \quad \widehat{\mathbb{P}^{\perp}u} = \frac{\xi\otimes\xi}{|\xi|^2}\hat{u} = \frac{\xi_i}{|\xi|} \frac{\xi_j}{|\xi|}\hat{u}_j = -i \frac{\xi_i}{|\xi|} \hat{g}(\xi)
\enda
\end{align*}
Hence, we can write $u = -i R g$. We claim the following inequality
\begin{align}\label{claimRf}
\bega
\|Rf\|_{\dot{B}_{p,q}^{s}} \leq C(p,d)\|f\|_{\dot{B}_{p,q}^{s}}
\enda
\end{align}
If this claim hold, then we can have the result since 
\begin{align*}
\bega
\|Ru\|_{\dot{B}_{p,q}^{s}} \leq C(p,d)\|u\|_{\dot{B}_{p,q}^{s}}, \quad \|u\|_{\dot{B}_{p,q}^{s}} =\|R g\|_{\dot{B}_{p,q}^{s}}  \leq C(p,d)\|g\|_{\dot{B}_{p,q}^{s}}
\enda
\end{align*}
Now we prove the claim \eqref{claimRf}. Let $\tilde{\varphi}$ be a smooth cutoff function $\tilde{\varphi}=1$ on $supp \varphi$. Then we have 
\begin{align*}
\bega
\widehat{\dot{\Delta}_k (Rf)} =  \frac{\xi}{|\xi|} \varphi(2^{-k}\xi) \hat{f}(\xi) = \frac{\xi}{|\xi|} \varphi(2^{-k}\xi) \tilde{\varphi}(2^{-k}\xi) \hat{f}(\xi)
\enda
\end{align*}
Applying Young's convolution inequality gives 
\begin{align*}
\bega
\|\dot{\Delta}_k (Rf)\|_{L^p} = \Big\| \mathcal{F}^{-1}(\frac{\xi}{|\xi|} \varphi(2^{-k}\xi))\Big\|_{L^1} \Big\| \mathcal{F}^{-1}(\tilde{\varphi}(2^{-k}\xi) \hat{f}(\xi)) \Big\|_{L^p} 
\enda
\end{align*}
Using the change of variable $2^k\xi = \eta$, we get 
\begin{align*}
\bega
\mathcal{F}^{-1}(\frac{\xi}{|\xi|} \varphi(2^{-k}\xi)) &= \frac{1}{(2\pi)^d}\int_{\R^d} e^{ix\cdot\xi} \frac{\xi}{|\xi|} \varphi(2^{-k}\xi) d\xi \cr 
&= 2^{kd} \frac{1}{(2\pi)^d}\int_{\R^d} e^{ix\cdot 2^k\eta} \frac{\eta}{|\eta|} \varphi(\eta) d\eta \cr 
&= 2^{kd} K(2^kx)
\enda
\end{align*}
Hence, taking $L^1_x$ norm gives 
\begin{align*}
\bega
\Big\| \mathcal{F}^{-1}(\frac{\xi}{|\xi|} \varphi(2^{-k}\xi))\Big\|_{L^1} \leq 2^{kd} \int_{\R^d} |K(2^kx)| dx = \int_{\R^d} |K(y)| dy = C
\enda
\end{align*}
where we used the change of variable $2^kx=y$ and 
\begin{align*}
\bega
K(y) = \frac{1}{(2\pi)^d}\int_{\R^d} e^{iy\cdot \eta} \frac{\eta}{|\eta|} \varphi(\eta) d\eta 
\enda
\end{align*}
Since $\tilde{\varphi}(2^{-k\xi})$ is supported in $\varphi(2^{-k}\xi)$, $\varphi(2^{-k+1}\xi)$ and $\varphi(2^{-k-1}\xi)$, we get 
\begin{align*}
\bega
\|\dot{\Delta}_k (Rf)\|_{L^p} \leq C \Big\| \mathcal{F}^{-1}(\tilde{\varphi}(2^{-k}\xi) \hat{f}(\xi)) \Big\|_{L^p} \leq \sum_{i=+1,-1,0}\Big\| \mathcal{F}^{-1}(\tilde{\varphi}(2^{-k-i}\xi) \hat{f}(\xi)) \Big\|_{L^p}
\enda
\end{align*}
Taking $\sum_{k\in\mathbb{Z}}2^{ks}$ then taking $l^q$, we proved the claim \eqref{claimRf}. 
\end{proof}
\unhide

\begin{remark}
We compare the end-point Strichartz estimate from \cite{JMR} and Proposition \ref{P.Book} from \cite{Danchin} for $d=3$.
In the context of estimating $\|\nabla_x\cdot u^{\e}\|_{L^\infty_x}$, we examine both the achievable scale and the required initial data.
\begin{itemize}
\item From~\cite{JMR}, the function $\mathrm{u}^{\e}$ satisfying the wave equation fulfills the estimate
\begin{align*}
\| \mathrm{u}^{\e}  \|_{L^2_TL_x^\infty}  &\leq  \sqrt{\e} \left(\sqrt{ \ln (1+ 2   T/\e ) }    \left(
\|O(1)\|_{L^1 (0,T; L^2_x)} + \|\mathrm{u}^{\e}_0 \|_{\dot{H}^1}  + \| \e\p_t \mathrm{u}^{\e}_0 \|_{L^2_x}
\right) + \| \p^2_x \mathrm{u}^{\e} \|_{L^2_TL^2_x} \right).
\end{align*}
The optimal scaling 
$\e^{\frac{1}{2}}\sqrt{ \ln (1+ 2   T/\e ) }$ is attainable only when the scaled time derivative $\e\p_t$ and the corresponding initial data are used.
For estimating $\|\nabla_x\cdot u^{\e}\|_{L^\infty}$, the required initial data includes $\|\nabla_x\cdot u^{\e}_0\|_{\dot{H}^1_x}$ and $\|\e\p_t\nabla_x\cdot u^{\e}_0\|_{L^2_x}$.
\item From~\cite{Danchin}, we employ the embedding
\begin{align*}
\bega
\|\nabla_x\cdot u^{\e}\|_{L^\infty} \leq C\|\nabla_x\cdot u^{\e}\|_{\dot{B}_{p,1}^{\frac{3}{p}}}, \quad \mbox{for} \quad d=3.
\enda
\end{align*}
To control the right-hand side of $\|\nabla_x\cdot u^{\e}\|_{L^\infty_x}$, we require the bound
\begin{align*}
\bega
\e^{\frac{1}{r}}\Big\|\Big(\nabla_x(\rho^{\e}_0+\ta^{\e}_0),(\nabla_x\mathbb{P}^{\perp}u^{\e}_0)\Big)\Big\|_{\dot{B}_{2,1}^{\frac{3}{2}-\frac{1}{r}}}. 
\enda
\end{align*}
As $p \to \infty$, this estimate yields the scaling $\e^{\frac{1}{2}-}$,
and the required initial data becomes
\begin{align*}
\bega
\Big\|\div(\mathbb{P}^{\perp}u^{\e}_0)\Big\|_{\dot{B}_{2,1}^{1+}} \sim \|\nabla_x\cdot u^{\e}_0\|_{H^{1+}_x}.
\enda
\end{align*}
\item Thus, the result in \cite{JMR} provides a slightly better scale with less restrictive initial data, but only when using the scaled time derivative $\e^{\mathfrak{n}}\p_t$ with $\mathfrak{n}=1$. 
\end{itemize}
\end{remark}

\section{High-order moment estimate in $L^\infty$}\label{Sec.V}

\StartNoTOC

\subsection{High-Order Moment Estimate Including the Top Order Derivatives}\label{Sec.V.Top}

In this section, we estimate the large velocity component $\mathcal{V}_{\ell}$, which arises in the Boltzmann energy estimates given in Proposition \ref{P.G.Energy} and Proposition \ref{P.F.Energy}. Recall that in the terms $\AC{I}_1$ and $I^F_2$ from Lemma \ref{L.G1} and Lemma \ref{L.E2}, respectively, the growth rates $|v|^2$ and $|v|^3$ appear due to the term \(\p_t M^{\e}+ \frac{v}{\e} \cdot \nabla_x M^{\e}\). This growth exceeds the allowable velocity growth \(\nu(v)\), which can be controlled through dissipation. The aim of this section is to provide control over \(\mathcal{V}_{\ell}\).
For the detailed estimate, we decompose $\mathcal{V}_{\ell}(F^{\e}(t))$, defined in \eqref{largev}, according to the order of derivatives as follows:  
\begin{align}\label{Vldef}
\bega
\mathcal{V}_{\ell}^{\al}(t,x) &:= \frac{\kappa^{(|\al|-\mathrm{N})_+}}{\eps^2} \int_{\R^3} (1 + |v|)^{\ell} | \p^{\al} \AC{\P}F^{\e} |^2 |M^{\e}|^{-1} dv , \quad \mbox{for} \quad 0\leq|\al|\leq \mathrm{N}+1,
\enda
\end{align}
so that $\mathcal{V}_{\ell}(F^{\e}(t,x)) = \sum_{0\leq|\al|\leq\mathrm{N}+1}\mathcal{V}_{\ell}^{\al}$.
Here, $x_+$ denotes $x_+ = x$ if $x \geq 0$ and $x_+ = 0$ if $x < 0$.
We decompose \( F^{\e} = \mu + \e\sqrt{\mu}f^{\e} \) near the global Maxwellian $\mu(v)=M_{[1,0,1]}$, and estimate \(\mathcal{V}_{\ell}^{\al}\) in terms of \(f^{\e}\).
We define an exponentially weighted function \(h^{\e}\) as follows:
\begin{align}\label{hdef}
\bega
h^{\e}(t,x,v):=w(v)f^{\e}(t,x,v), \quad \text{where} \quad w(v):= e^{c_1|v|^2}, \quad \frac{c_0}{4k_B(1-c_0/2)}< c_1 <\frac{1}{4k_B},
\enda
\end{align}
where \(c_0 \ll 1\) is a small constant from \eqref{condition}, satisfying \(|\mathrm{P}^{\e} -1|, |\mathrm{U}^{\e}|, |\mathrm{\Theta}^{\e} -1| < c_0/2 \).
Multiplying the equation of $f^{\e}$ by $w(v)$ gives
\Be \label{heqn}
\bigg[\p_t + \frac{v}{\e} \cdot \nabla_x + \frac{\nu}{\kappa \e^2} \bigg] h^{\e}
= \frac{1}{\kappa \e^2} \mathbf{K}_w(h^{\e}) + \frac{w}{\kappa \e} \Gamma\bigg(\frac{h^{\e}}{w},\frac{h^{\e}}{w}\bigg).
\Ee
Here the linear operator is defined by
\begin{align}\label{kdef}
\bega
\mathbf{K}_w h(v) &:= \int_{\R^3}\mathbf{k}_w(v,v_*)h(v_*)dv_*,
\qquad
\mathbf{k}_w(v,v_*) := \mathbf{k}(v,v_*)\frac{w(v)}{w(v_*)}, \cr
\mathbf{k}(v,v*)&:= c_1|v-v_*|e^{-\frac{|v|^2+|v_*|^2}{2}} - \frac{c_2}{|v-v_*|}e^{-\frac{1}{4}|v-v_*|^2-\frac{1}{4}\frac{(|v|^2-|v_*|^2)^2}{|v-v_*|^2}}, 
\enda
\end{align}
for some positive constants $c_1$ and $c_2$.
And the nonlinear operator is given by $\Gamma(f,f) = \frac{1}{\sqrt{\mu}}\mathcal{N}(\sqrt{\mu}f,\sqrt{\mu}f)$.
Applying the derivative \(\p^{\al}\), we have
\Be \label{wBE}
\bigg[\p_t + \frac{v}{\e} \cdot \nabla_x + \frac{\nu}{\kappa \e^2} \bigg] \p^{\al}h^{\e}
= \frac{1}{\kappa \e^2} \mathbf{K}_w \p^{\al}h^{\e} + \frac{w}{\kappa \e } \sum_{0\leq\beta\leq\al}\binom{\al}{\beta}\Gamma\bigg(\frac{\p^{\beta}h^{\e}}{w},\frac{\p^{\al-\beta}h^{\e}}{w}\bigg).
\Ee


We note that the energy defined in \eqref{N-EDdef} differs depending on whether the initial macroscopic velocity corresponds to the case of finite velocity energy \eqref{caseEC} or infinite velocity energy \eqref{caseECX}. In particular, 
\begin{align*}
\bega
\|u^{\e}(t)\|_{L^2_x} &\leq \mathcal{E}_M^{\frac{1}{2}}(t), \quad &\mbox{for the finite velocity energy case \eqref{caseEC}}, \cr 
\|u^{\e}(t)\|_{L^2_x} &= \infty, \quad &\mbox{for the infinite velocity energy case \eqref{caseECX}}.
\enda
\end{align*}
However, when estimating $\mathcal{V}^{\al}_{\ell}$, a difficulty arises. 
Indeed, the quantity $\mathcal{V}^{\al}_{\ell}$ contains the term
$\frac{1}{\eps^2}\int_{\Omega\times\R^3}|(\mu-M^{\e})|^2|M^{\e}|^{-1} dvdx$,
which forces $u^{\e}$ to lie in $L^2$. Hence, this decomposition cannot be applied in the infinite velocity energy case \eqref{caseECX} at the level without derivatives.
For this reason, we present the main propositions of this section depending on whether \eqref{caseEC} or \eqref{caseECX} applies.

\begin{proposition}\label{P.hLinf}
Let $\Omega=\R^d$ with $d=2$ or $3$.
For both the purely spatial derivative case \eqref{caseA} and the space-time derivative case \eqref{caseB}, for any $\mathfrak{n} \geq 0$ and $\mathrm{N} > d/2 + 1$, and for $t\in[0,T]$ satisfying the bootstrap assumption \eqref{condition}, the following estimates hold:
\begin{align}\label{L.Vdecomp}
\bega
&\int_{\Omega} \bigg(|\nabla_xu^{\e}(t)| \mathcal{V}_{2}^{\al}(t) + |\nabla_x\ta^{\e}(t)| \mathcal{V}_{3}^{\al}(t)\bigg) dx \les \|\nabla_x (u^{\e},\ta^{\e})(t)\|_{L^\infty_x}\Big(\mathcal{E}_{tot}^{\mathrm{N}}(F^{\e}(t)) +\e^{\frac{3}{2}}\kappa \mathcal{D}_{tot}^{\mathrm{N}}(F^{\e}(t))\Big) \cr 
&\qquad + e^{-c_1\e^{-{\frac{1}{3}}}}\frac{1}{\e} \mathcal{E}_{tot}^{\mathrm{N}}(F^{\e}(t)) \kappa^{(|\al|-\mathrm{N})_+} \|\p^{\al}h^{\e}(t)\|_{L^\infty_{x,v}},  \cr
&\frac{1}{\eps} \int_{\Omega} \TbT(t) \mathcal{V}_{2}^{\al}(t) dx \les \e^{\frac{1}{2}}\kappa^{\frac{1}{2}}(\mathcal{D}_{G}^{\mathrm{N}}(F^{\e}(t)))^{\frac{1}{2}} \mathcal{E}_{tot}^{\mathrm{N}}(F^{\e}(t)) \cr 
&\qquad + \e\kappa^{\frac{1}{2}} e^{-c_1\e^{-{\frac{1}{2}}}}  (\mathcal{D}_{G}^{\mathrm{N}}(F^{\e}(t)))^{\frac{1}{2}}(\mathcal{E}_{tot}^{\mathrm{N}}(F^{\e}(t)))^{\frac{1}{2}} \kappa^{(|\al|-\mathrm{N})_+}\|\p^{\al}h^{\e}(t)\|_{L^\infty_{x,v}},
\enda
\end{align}
\hide
\begin{align}
\bega
\mathcal{V}^{\al}_{\ell}(t) &\leq C\bigg[\mathcal{E}_{tot}^{\mathrm{N}}(F^{\e}(t)) +\e^{\frac{3}{2}}\kappa \Big(\mathcal{D}_{tot}^{\mathrm{N}}(F^{\e}(t))+\mathbf{1}_{|\al|\geq2}(\mathcal{E}_M^{\mathrm{N}}(F^{\e}(t)))^2\Big) + e^{-c_1\e^{-{\frac{1}{\ell}}}}\kappa^{(|\al|-\mathrm{N})_+}  \|\p^{\al}h^{\e}(t)\|_{L^\infty_{x,v}}^2\bigg], \cr 
\mathcal{V}^{\al}_{\ell}(t) &\leq C\bigg[\e^{-\frac{1}{2}}\mathcal{E}_{tot}^{\mathrm{N}}(F^{\e}(t)) + e^{-c_1\e^{-{\frac{1}{\ell}}}}\kappa^{(|\al|-\mathrm{N})_+}  \|\p^{\al}h^{\e}(t)\|_{L^\infty_{x,v}}^2\bigg],
\enda
\end{align}
\unhide
for $0 \leq |\al| \leq \mathrm{N}+1$ in the finite velocity energy case \eqref{caseEC}, 
and for $1 \leq |\al| \leq \mathrm{N}+1$ in the infinite velocity energy case \eqref{caseECX}.
\end{proposition}

\begin{remark}\label{Rmk.VLt}
In the first line of \eqref{L.Vdecomp}, we decomposed $\mathcal{V}_{\ell}^{\al}$ using the dissipation in $L^2_T$, which yields a favorable scaling. 
However, this decomposition cannot be used for the second inequality of \eqref{L.Vdecomp}, since $\|\TbT(t)\|_{L^\infty_x}$ is controlled through the dissipation term $\e^2 \kappa^{\frac{1}{2}} \mathcal{D}_G^{\frac{1}{2}}$. 
Therefore, we decompose $\mathcal{V}_{\ell}^{\al}$ using only the energy and estimate it in $L^\infty_T$, which leads to a weaker scaling than in the first line. 
Nevertheless, this is still sufficient to close the argument.
\hide
Recall the estimate \eqref{totalGt} from Proposition \ref{P.G.Energy}. When estimating the term $\|\nabla_x(u^{\e}, \ta^{\e})\|_{L^\infty_x} \mathcal{V}_{\ell}$, we can use the first line of \eqref{L.Vdecomp}, which includes the dissipation term $\mathcal{D}_{tot}(t)$. However, for estimating the term $\TbT(t)\mathcal{V}^{\al}_2(t)$, the first line of \eqref{L.Vdecomp} is not applicable, since $\TbT$ is controlled by the dissipation term $\e^2 \kappa^{\frac{1}{2}} \mathcal{D}_G^{\frac{1}{2}}$, as shown in (2) of Lemma \ref{L.ABG}.
Therefore, even though the second line of \eqref{L.Vdecomp} has a weaker scaling compared to the first, it must be used in this context, and it suffices to close the argument.
(We could alternatively use time embeddings if we employ space-time derivatives as in \eqref{caseB}. However, this approach leads to worse scaling for $\mathfrak{n} \geq 1$, as it requires losing a factor of $\e^{\mathfrak{n}}$.)
\unhide
\end{remark}

Different from the assumption in Proposition \ref{P.hLinf}, the following proposition holds for any number of derivatives $0 \leq |\al| \leq \mathrm{N}+1$ in both the finite velocity energy case \eqref{caseEC} and the infinite velocity energy case \eqref{caseECX}.

\begin{proposition}\label{P.hLinf2}
Let $\Omega = \R^2$ or $\R^3$. Suppose the initial data $u_0^{\e}$ satisfies either the finite velocity energy condition \eqref{caseEC} or the infinite velocity energy condition \eqref{caseECX}. Depending on the spatial dimension and on whether one considers
the purely spatial derivative case~\eqref{caseA} or the space--time derivative case~\eqref{caseB}, we assume that
\begin{align}\label{hassume}
\bega
\e\sum_{0\leq |\al_x|\leq \mathrm{N}}\|\p^{\al_x}h^{\e}
_0\|_{L^\infty_{x,v}} + \mathcal{C}_{(\e,\kappa)} \sup_{t \in [0,T]} (\mathcal{E}_{tot}^{\mathrm{N}}(F^{\e}(t)))^{\frac{1}{2}} +\e \sup_{t \in [0,T]}\big(1+\mathcal{E}_{tot}^{\mathrm{N}}(F^{\e}(t))\big)^{\frac{1}{2}} < \frac{1}{4C},
\enda
\end{align}
where $C>0$ is a constant to be chosen, and $\alpha_x$ denotes the purely spatial multi-index defined in~\eqref{mindexx}.
Here, $\mathcal{C}_{(\e,\kappa)}$ is defined by
\begin{align*}
\bega
\mathcal{C}_{(\e,\kappa)} = 
&\begin{cases}
\e(\e\kappa)^{-\frac{2}{p}}, \quad \mbox{for any} \quad 1\leq p<\infty, \quad &d=2 \quad \mbox{with} \quad \eqref{caseA}, \\
\e\kappa^{-\frac{1}{2}}(\e\kappa)^{-\frac{2}{p}}, \quad \mbox{for any} \quad 1\leq p<\infty, \quad &d=2 \quad \mbox{with} \quad \eqref{caseB}, \\
\e\kappa^{-\frac{1}{4}}, \quad &d=3 \quad  \mbox{with} \quad \eqref{caseA}, \\
\e^{\frac{1}{4}}\kappa^{-\frac{9}{8}}, \quad &d=3 \quad  \mbox{with} \quad \eqref{caseB}.
\end{cases}
\enda
\end{align*}
Then the function \(h^{\e}\), defined in \eqref{hdef}, satisfies
\begin{align}\label{claimx4}
\bega
\begin{cases} 
\displaystyle \e\sum_{0\leq |\al_x|\leq \lfloor (\mathrm{N}+1)/2 \rfloor}\|\p^{\al_x}h^{\e}(t)\|_{L^\infty_tL^\infty_{x,v}} \leq \frac{1}{2C},  \quad &\mbox{for the purely spatial derivative }~\eqref{caseA}, \\
\displaystyle \e\sum_{0\leq |\al_x|\leq \mathrm{N}}\|\p^{\al_x}h^{\e}(t)\|_{L^\infty_tL^\infty_{x,v}} \leq \frac{1}{2C},  \quad &\mbox{for the space-time derivative }~\eqref{caseB},
\end{cases}
\enda
\end{align}
and
\begin{align}
&\sum_{0\leq|\al|\leq\mathrm{N}+1}\!\!\!\kappa^{\frac{(|\al|-\mathrm{N})_+}{2}}\|\p^{\al}h^{\e}\|_{L^2_tL^\infty_{x,v}}  \les \e\kappa^{\frac{1}{2}} \!\!\!\!\!\!\sum_{0\leq|\al|\leq\mathrm{N}+1}\!\!\!\!\kappa^{\frac{(|\al|-\mathrm{N})_+}{2}}\|\p^{\al}h^{\e}_0\|_{L^\infty_{x,v}}+ \frac{1}{(\e\kappa)^{\frac{d}{2}}}\bigg(\int_0^t  \big(1+\mathcal{E}_{tot}^{\mathrm{N}}(F^{\e}(s))\big) ds\bigg)^{\frac{1}{2}}, \label{hL2}   \\
&\sum_{0\leq|\al|\leq\mathrm{N}+1}\kappa^{\frac{(|\al|-\mathrm{N})_+}{2}}\|\p^{\al}h^{\e}\|_{L^\infty_tL^\infty_{x,v}}  \les\! \!\!\!\sum_{0\leq|\al|\leq\mathrm{N}+1}\!\!\kappa^{\frac{(|\al|-\mathrm{N})_+}{2}}\|\p^{\al}h^{\e}_0\|_{L^\infty_{x,v}} +\frac{1}{(\e\kappa)^{\frac{d}{2}}}\|\big(1+\mathcal{E}_{tot}^{\mathrm{N}}(F^{\e}(t))\big)\|_{L^\infty_t}^{\frac{1}{2}} ,\label{hLinf} 
\end{align}
for both the purely spatial derivative case \eqref{caseA} and the space-time derivative case \eqref{caseB}, for any $\mathfrak{n} \geq 0$ and $\mathrm{N} > d/2 + 1$.
\end{proposition}

In the proof, for brevity, we slightly abuse notation by writing 
$h$, $f$, $\mathcal{E}(t)$, and $\mathcal{D}(t)$ for 
$h^{\e}$, $f^{\e}$, $\mathcal{E}^{\mathrm{N}}(F^{\e}(t))$, 
and $\mathcal{D}^{\mathrm{N}}(F^{\e}(t))$, respectively.





\begin{proof}[\textbf{Proof of Proposition \ref{P.hLinf}}]
We divide the velocity integration domain of $\mathcal{V}^{\al}_{\ell}$ into the regions $\{ |v| \leq \e^{-\frac{1}{2\ell}} \}$ and $\{ |v| > \e^{-\frac{1}{2\ell}} \}$ as follows:
\begin{align*}
\bega
\mathcal{V}^{\al}_{\ell}(t,x)&=\frac{\kappa^{(|\al|-\mathrm{N})_+}}{\eps^2}\bigg[\int_{|v| \leq \e^{-\frac{1}{2\ell}}} + \int_{|v| > \e^{-\frac{1}{2\ell}}}\bigg]\la v \ra^{\ell}|\p^{\al}\AC{\P} F|^2|M^{\e}|^{-1} dv = \mathcal{V}^{\al}_{\ell,s}(t,x) + \mathcal{V}^{\al}_{\ell,L}(t,x).
\enda
\end{align*}
Using the decomposition $F^{\e} = M^{\e} + \AC{\P}F^{\e} = \mu + \e \sqrt{\mu} f^{\e}$, we split $\mathcal{V}^{\al}_{\ell,L} \leq \mathcal{V}^{\al}_{\ell,L,1}+\mathcal{V}^{\al}_{\ell,L,2}$, where
\begin{align}\label{VlLdef}
\bega
\mathcal{V}^{\al}_{\ell,L,1}(t,x) &:= \frac{\kappa^{(|\al|-\mathrm{N})_+}}{\eps^2}\int_{\R^3}(1+|v|)^{\ell}|\p^{\al}(\mu-M^{\e})|^2|M^{\e}|^{-1} dv, \cr
\mathcal{V}^{\al}_{\ell,L,2}(t,x) &:=\frac{\kappa^{(|\al|-\mathrm{N})_+}}{\eps^2}\int_{|v| \geq \e^{-\frac{1}{2\ell}}}(1+|v|)^{\ell}|\e\sqrt{\mu}\p^{\al}f^{\e}|^2|M^{\e}|^{-1} dv.
\enda
\end{align}
As explained in Remark \ref{Rmk.VLt}, depending on the available time-space norms $L^2_T$ or $L^\infty_T$, we require two types of estimates. For each part, we claim the following estimates:
\begin{align}\label{smallv}
\bega
\int_{\Omega}\mathcal{V}^{\al}_{\ell,s}(t)dx &\les \begin{cases} \e^{\frac{3}{2}}\kappa \Big(\mathcal{D}_{tot}(t)+\mathcal{E}_M^2(t)\Big), \\ \e^{-\frac{1}{2}}\mathcal{E}_{tot}(t),
\end{cases}
\enda
\end{align}
and
\begin{align}\label{Vl1b}
\bega
\int_{\Omega} \mathcal{V}^{\al}_{\ell,L,1} dx 
&\les \mathcal{E}_{tot}(t), \quad \begin{cases} \mbox{for} \quad 0 \leq |\al| \leq \mathrm{N}+1, \quad \mbox{for} \quad \eqref{caseEC}, \\ 
\mbox{for} \quad 1 \leq |\al| \leq \mathrm{N}+1, \quad \mbox{for} \quad \eqref{caseECX}.
\end{cases}
\enda
\end{align}
And, for $\mathfrak{X}(t,x) \in \{\p_iu_j^{\e}(t,x), \p_i \ta^{\e}(t,x) , \TbT(t,x) \} $, we claim 
\begin{align}\label{Vl2a}
\bega
&\int_{\Omega} |\mathfrak{X}(t)| \mathcal{V}^{\al}_{\ell,L,2}(t) dx \les \kappa^{(|\al|-\mathrm{N})_+} \frac{1}{\e}e^{-c_1\e^{-{\frac{1}{\ell}}}} \|\mathfrak{X}(t)\|_{L^2_x} \mathcal{E}_{tot}^{\frac{1}{2}}(t) \|\p^{\al}h^{\e}(t)\|_{L^\infty_{x,v}}.
\enda
\end{align}
If the above three claims \eqref{smallv}, \eqref{Vl1b} and \eqref{Vl2a} hold, then applying $\|\nabla_x(u^{\e},\ta^{\e})\|_{L^2_x}  \leq \mathcal{E}_M^{\frac{1}{2}}$ and $\|\TbT\|_{L^2_x} \leq \e^2\kappa^{\frac{1}{2}}\mathcal{D}_G^{\frac{1}{2}}$ from \eqref{ABGscale}, we directly obtain the result of Proposition \ref{P.hLinf}. \\
(Proof of the claim \eqref{smallv}) For the small velocity region, we estimate it by noting the loss of scaling due to $\la v \ra^{\ell} \les 1 + \e^{-\frac{1}{2}} \les \e^{-\frac{1}{2}}$:
\begin{align}\label{smallv0}
\bega
\mathcal{V}^{\al}_{\ell,s}&\les  \e^{-\frac{1}{2}}\frac{\kappa^{(|\al|-\mathrm{N})_+}}{\eps^2}\int_{|v| \leq \e^{-\frac{1}{2\ell}}}|\p^{\al}\AC{\P} F|^2|M^{\e}|^{-1} dv.
\enda
\end{align}
To derive the first line of \eqref{smallv}, we use the dissipation $\mathcal{D}_G$ for $0 \leq |\al| \leq \mathrm{N}$ and apply estimate \eqref{al5decomp} for $|\al| = \mathrm{N}+1$:
\begin{align*}
\bega
\int_{\Omega}\mathcal{V}^{\al}_{\ell,s}(t)dx &\les \e^{-\frac{1}{2}}\bigg(\frac{1}{\e^2}\Big(\e^4\kappa\mathcal{D}_G\Big) + \frac{\kappa}{\e^2}\Big(\e^4\mathcal{D}_{top} + \e^4 \mathbf{1}_{|\al|\geq2}\mathcal{E}_M^2\Big) \bigg) \cr 
&\les \e^{\frac{3}{2}}\kappa \Big(\mathcal{D}_G+\mathcal{D}_{top}+\mathbf{1}_{|\al|\geq2}\mathcal{E}_M^2\Big).
\enda
\end{align*}
Then, the first line of \eqref{smallv} follows from the definition of $\mathcal{D}_{tot}$ in \eqref{EDtotdef}.
To derive the second line of \eqref{smallv}, particularly for the case $|\al| = \mathrm{N}+1$, we apply estimates \eqref{pACP}, \eqref{ACPX}, and \eqref{Rscale} to have
\begin{align}\label{Eal=5}
\bega
&\kappa^{(|\al|-\mathrm{N})_+}\int_{\Omega \times \R^3} | \p^{\al} \AC{\P} F^{\e} |^2 |M^{\e}|^{-1} dv dx\\
&\les \kappa^{(|\al|-\mathrm{N})_+}\int_{\Omega \times \R^3} | \AC{\P} \p^{\al} F^{\e} |^2 |M^{\e}|^{-1} dv dx + \kappa^{(|\al|-\mathrm{N})_+}\int_{\Omega \times \R^3} \bigg|\AC{\P} \bigg( \sum_{2\leq i\leq |\al|}\eps^i \Phi_{\al}^i M^{\e} \bigg)\bigg|^2 |M^{\e}|^{-1} dv dx \cr 
&\les \e^2\mathcal{E}_{top} + \e^4\kappa \mathcal{E}_M^2.
\enda
\end{align}
For $0 \leq |\al| \leq \mathrm{N}$, we use the definition of 
$\mathcal{E}_G(t)$ in \eqref{N-EDdef}. 
Combining this with \eqref{Eal=5} in \eqref{smallv0}, we obtain
\begin{align*}
\bega 
\int_{\Omega}\mathcal{V}^{\al}_{\ell,s}(t) dx 
&\les \e^{-\frac{1}{2}}\bigg(\frac{1}{\e^2}\Big(\e^2\mathcal{E}_G\Big) + \frac{1}{\e^2}\Big(\e^2\mathcal{E}_{top} + \e^4\kappa \mathcal{E}_M^2\Big) \bigg) \les \e^{-\frac{1}{2}}\mathcal{E}_{tot},
\enda
\end{align*}
where we used $(\mathcal{B}_1)$ in \eqref{condition}. This proves the claim of the second line of \eqref{smallv}. \\
(Proof of the claim \eqref{Vl1b}) To estimate $\mathcal{V}^{\al}_{\ell,L,1}$, we use the expansion of $M^{\e}$ from Lemma \ref{L.MFME}:
\begin{align*}
\bega
\p^{\al}(M^{\e}-\mu) &=  \e\p^{\al}\left(\rho^{\e} + u^{\e} \cdot \frac{v}{k_B} + \ta^{\e} \frac{|v|^2 - 3k_B}{2k_B}\right) \mu \\
&+ \e^2 \p^{\al}\bigg((\mathrm{P}^{\vartheta}\rho^{\e}, u^{\e}, \mathrm{\Theta}^{\vartheta}\ta^{\e})^T \int_0^1 \left\{D^2_{(\mathrm{P}^{\vartheta}, \mathrm{U}^{\vartheta}, \mathrm{\Theta}^{\vartheta})}M(\vartheta)\right\} (1 - \vartheta) d\vartheta (\mathrm{P}^{\vartheta}\rho^{\e}, u^{\e}, \mathrm{\Theta}^{\vartheta}\ta^{\e})^T \bigg).
\enda
\end{align*}
We note that $\|u^{\e}\|_{L^2_x}$ is bounded by $\mathcal{E}_M^{\frac{1}{2}}$, defined in \eqref{N-EDdef}, in the case of finite velocity energy \eqref{caseEC}, but it is unbounded in the case of infinite velocity energy \eqref{caseECX}. 
Hence, we consider $0 \leq |\al| \leq \mathrm{N}+1$ in the finite velocity energy case, whereas $1 \leq |\al| \leq \mathrm{N}+1$ in the infinite velocity energy case.
Using the definitions of $\mathcal{E}_M$ and $\mathcal{E}_{top}$ from \eqref{N-EDdef} and \eqref{N-EDdef2}, and applying \eqref{Phiscale} and condition $(\mathcal{B}_1)$ in \eqref{condition}, we obtain the result in \eqref{Vl1b}. \\ 
(Proof of the claim \eqref{Vl2a}) 
We apply H\"{o}lder inequality twice by utilizing the integrability of spatial variable:
\begin{align}\label{lL2-0}
\bega
\int_{\Omega} |\mathfrak{X}| \mathcal{V}^{\al}_{\ell,L,2}dx &\les \frac{\kappa^{(|\al|-\mathrm{N})_+}}{\eps^2} \|\mathfrak{X}\|_{L^2_x} \bigg(\int_{\Omega}\int_{|v| \geq \e^{-\frac{1}{2\ell}}}\frac{|\e\sqrt{\mu}\p^{\al}f^{\e}|^2}{M^{\e}} dvdx\bigg)^{\frac{1}{2}} \cr 
&\quad \times \bigg\|\bigg(\int_{|v| \geq \e^{-\frac{1}{2\ell}}}(1+|v|)^{2\ell} \frac{|\e\sqrt{\mu}\p^{\al}f^{\e}|^2}{M^{\e}} dv\bigg)^{\frac{1}{2}} \bigg\|_{L^\infty_x} .
\enda
\end{align}
For the first velocity integral, we again use the decomposition $F^{\e} = M^{\e} + \AC{\P}F^{\e} = \mu + \e \sqrt{\mu} f^{\e}$ to have 
\begin{align}\label{lL2-1}
\bega
\bigg(\int_{\Omega}\int_{|v| \geq \e^{-\frac{1}{2\ell}}}\frac{|\e\sqrt{\mu}\p^{\al}f^{\e}|^2}{M^{\e}} dvdx\bigg)^{\frac{1}{2}} &\les \bigg(\int_{\Omega}\int_{|v| \geq \e^{-\frac{1}{2\ell}}}\frac{|\p^{\al}(M^{\e}-\mu)|^2}{M^{\e}} + \frac{|\p^{\al}\AC{\P}F^{\e}|^2}{M^{\e}} dvdx\bigg)^{\frac{1}{2}}  \cr
&\les \e^2 \mathcal{E}_M+\e^2 \mathcal{E}_G + \e^4\mathcal{E}_M^2.
\enda
\end{align}
For the last velociy integral, we insert a factor of $w^2(v)$ by multiplying and dividing inside the integral:
\begin{align}\label{lL2-2}
\bega
\bigg\|\bigg(\int_{|v| \geq \e^{-\frac{1}{2\ell}}}(1+|v|)^{2\ell} \frac{|\e\sqrt{\mu}\p^{\al}f^{\e}|^2}{M^{\e}} &dv\bigg)^{\frac{1}{2}} \bigg\|_{L^\infty_x} \les \int_{|v| \geq \e^{-\frac{1}{2\ell}}}(1+|v|)^{2\ell}|\p^{\al}h^{\e}|^2 \frac{1}{w^2}\mu|M^{\e}|^{-1} dv \cr 
&\hspace{-2mm} \les e^{-c_1\e^{-{\frac{1}{\ell}}}} \|\p^{\al}h^{\e}\|_{L^\infty_{x,v}} \int_{|v| \geq \e^{-\frac{1}{2\ell}}}(1+|v|)^{2\ell}\frac{1}{w}\mu|M^{\e}|^{-1} dv ,
\enda
\end{align}
where we used $\frac{1}{w(v)} \leq e^{-c_1 \e^{-\frac{1}{\ell}}}$. For $c_1 > \frac{c_0}{4k_B(1 - c_0/2)}$ and using $(\mathcal{B}_2)$ in \eqref{condition}, the integral is bounded as follows:
\begin{align*}
\int_{\R^3}(1+|v|)^{2\ell}\frac{1}{w}\mu|M^{\e}|^{-1} dv &\les \int_{\R^3} (1+|v|)^{2\ell}e^{-c_1|v|^2}e^{-\frac{|v|^2}{2k_B}}e^{\frac{|v|^2 +|\mathrm{U}^{\e}|^2}{2k_B(1-c_0/2)}}dv \leq C.
\end{align*}
Applying \eqref{lL2-1} and \eqref{lL2-2} to \eqref{lL2-0} and using $(\mathcal{B}_1)$ in \eqref{condition} proves the claim \eqref{Vl2a}. 
\end{proof}

To prove Proposition \ref{P.hLinf2}, we first introduce the following two auxiliary lemmas.

\begin{lemma}\label{L.htt} 
Let \( h^{\e}(t,x,v) \) solve the equation \eqref{wBE}. For both types of derivatives \( \p^{\al} \) defined in \eqref{caseA} and \eqref{caseB}, for any $\mathfrak{n} \geq 0$, and for \( (x,v) \in \R^d \times \R^3 \) with \( d = 2,3 \), we have
\begin{align}
\|\p^{\al}h^{\e}\|_{L^2_tL^\infty_{x,v}} &\leq C\e\kappa^{\frac{1}{2}}\|\p^{\al}h^{\e}_0\|_{L^\infty_{x,v}} + \frac{C_\mathfrak{N}}{(\e\kappa)^{\frac{d}{p}}} \|\mathbf{1}_{|v|\leq2\mathfrak{N}}\p^{\al}f^{\e}\|_{L^2_tL^p_xL^2_v} \label{hL2t} \\
& + C \e \sum_{0\leq\beta\leq\al} \big\|\|\p^{\beta}h^{\e}(t)\|_{L^\infty_{x,v}} \|\p^{\al-\beta}h^{\e}(t)\|_{L^\infty_{x,v}}\big\|_{L^2_t},  \cr
\|\p^{\al}h^{\e}\|_{L^\infty_tL^\infty_{x,v}} &\leq C\|\p^{\al}h^{\e}_0\|_{L^\infty_{x,v}} + \frac{C_\mathfrak{N}}{(\e\kappa)^{\frac{d}{p}}} \|\mathbf{1}_{|v|\leq2\mathfrak{N}}\p^{\al}f^{\e}\|_{L^\infty_tL^p_xL^2_v} \label{hLinft} \\
&+ C \e \sum_{0\leq\beta\leq\al} \|\p^{\beta}h^{\e}\|_{L^\infty_tL^\infty_{x,v}} \|\p^{\al-\beta}h^{\e}\|_{L^\infty_tL^\infty_{x,v}}, \notag 
\end{align}
for any \( p \in [1, \infty] \) and a constant \( C > 0 \), where \( C \) is the constant from \eqref{wgammav}. The constant \( C_\mathfrak{N} > 0 \) depends on \( \mathfrak{N} \), which will be determined in the proof.
\end{lemma}

\begin{proof}
Along the characteristic curve, we can rewrite \eqref{wBE} as follows: 
\begin{align}\label{1duhamel}
\bega
\p^{\al}h(t,x,v) &= e^{ \frac{-\nu(v)t}{\e^2 \kappa}}\p^{\al}h_0\bigg(x-\frac{v}{\e}t,v\bigg) \cr 
&+ \frac{1}{\e^2 \kappa}\int_0^t e^{ \frac{-\nu(v)(t-s)}{\e^2 \kappa}} \int_{\R^3}\mathbf{k}_w(v,v_*) \p^{\al}h \bigg(s, x-\frac{v}{\e}(t-s),v_*\bigg) dv_* ds \cr 
&+ \frac{1}{\e \kappa}\sum_{0\leq\beta\leq\al}\binom{\al}{\beta}\int_0^t e^{ \frac{-\nu(v)(t-s)}{\e^2 \kappa}} w\Gamma\bigg(\frac{\p^{\beta}h}{w},\frac{\p^{\al-\beta}h}{w}\bigg)\bigg(s, x-\frac{v}{\e}(t-s),v\bigg)ds.
\enda
\end{align}
We apply Duhamel’s principle once more to the term $\p^{\al}h$ appearing in the second line of \eqref{1duhamel}, obtaining
\begin{align}\label{alh}
\p^{\al}h(t,x,v) &= I^h_1 + I^h_2 + I^h_3 + I^h_4 + I^h_5,
\end{align}
where
\begin{align*}
I^h_1(t,x,v) &= e^{ \frac{-\nu(v)t}{\e^2 \kappa}}\p^{\al}h_0\bigg(x-\frac{v}{\e}t,v\bigg) ,\cr 
I^h_2(t,x,v)&= \frac{1}{\e^2 \kappa}\int_0^t e^{ \frac{-\nu(v)(t-s)}{\e^2 \kappa}} \int_{\R^3}\mathbf{k}_w(v,v_*) e^{ \frac{-\nu(v_*)s}{\e^2 \kappa}}\p^{\al}h_0\bigg(x-\frac{v}{\e}(t-s)-\frac{v_*}{\e}s,v_*\bigg) dv_* ds ,\cr 
I^h_3(t,x,v)&= \frac{1}{\e^2 \kappa}\int_0^t e^{ \frac{-\nu(v)(t-s)}{\e^2 \kappa}} \int_{\R^3}\mathbf{k}_w(v,v_*) \frac{1}{\e^2 \kappa} \int_0^s e^{\frac{-\nu(v_*)(s-\tau)}{\e^2 \kappa}}  \cr 
&\quad \int_{\R^3}\mathbf{k}_w(v_*,v_{**}) \p^{\al}h\bigg(\tau, x-\frac{v}{\e}(t-s)-\frac{v_*}{\e}(s-\tau),v_{**}\bigg)dv_{**} d\tau dv_* ds ,\cr 
I^h_4(t,x,v)&= \frac{1}{\e^2 \kappa}\int_0^t e^{ \frac{-\nu(v)(t-s)}{\e^2 \kappa}} \int_{\R^3}\mathbf{k}_w(v,v_*) \frac{1}{\e \kappa} \sum_{0\leq\beta\leq\al}\binom{\al}{\beta} \int_0^s e^{\frac{-\nu(v_*)(s-\tau)}{\e^2 \kappa}} \cr 
&\quad \times  w\Gamma\bigg(\frac{\p^{\beta}h}{w},\frac{\p^{\al-\beta}h}{w}\bigg)\bigg(\tau, x-\frac{v}{\e}(t-s)-\frac{v_*}{\e}(s-\tau),v_*\bigg)d\tau dv_* ds ,\cr 
I^h_5(t,x,v)&= \frac{1}{\e \kappa}\sum_{0\leq\beta\leq\al}\binom{\al}{\beta}\int_0^t e^{ \frac{-\nu(v)(t-s)}{\e^2 \kappa}} w\Gamma\bigg(\frac{\p^{\beta}h}{w},\frac{\p^{\al-\beta}h}{w}\bigg)\bigg(s, x-\frac{v}{\e}(t-s),v\bigg)ds .
\end{align*}
We estimate the terms \( I^h_1, \dots, I^h_5 \) one by one. We present the proof for the \( L^2_t \) case \eqref{hL2t}; the \( L^\infty_t \) case \eqref{hLinft} follows by a similar argument. 
Taking the \( L^2_t \) norm of \( I^h_1 \), we have
\begin{align}\label{hI1}
\bega
\|I^h_1\|_{L^2_t} &\les \lw(\int_0^t \Big|e^{ \frac{-\nu(v)t}{\e^2 \kappa}}\Big|^2 dt\rw)^{\frac{1}{2}} \|\p^{\al}h_0\|_{L^\infty_{x,v}} \les \e\kappa^{\frac{1}{2}}\|\p^{\al}h_0\|_{L^\infty_{x,v}}.
\enda
\end{align}
Applying Young’s convolution inequality to \( I^h_2 \), we get
\begin{align}\label{hI2}
\bega
\|I^h_2\|_{L^2_t} 
&\les \bigg\|\frac{1}{\e^2 \kappa}e^{ \frac{-\nu(v)t}{\e^2 \kappa}}\bigg\|_{L^1_t} \bigg\|\int_{\R^3}\mathbf{k}_w(v,v_*) e^{ \frac{-\nu(v_*)t}{\e^2 \kappa}} dv_* \bigg\|_{L^2_t} \|\p^{\al}h_0\|_{L^\infty_{x,v}} 
\les \e\kappa^{\frac{1}{2}}\|\p^{\al}h_0\|_{L^\infty_{x,v}},
\enda
\end{align}
where we used property $\eqref{kw}_2$ of the kernel.
Similarly, applying Young’s inequality and the estimate \eqref{wgammav} to \( I^h_5 \), we obtain
\begin{align}\label{hI5}
\bega
\|I^h_5\|_{L^2_t} &\leq C \e \bigg\|\frac{1}{\e^2\kappa}e^{\frac{-\nu(v)t}{\e^2 \kappa}}\bigg\|_{L^1_t} \sum_{0\leq\beta\leq\al}  \nu(v)\Big\|\|\p^{\beta}h(s)\|_{L^\infty_{x,v}}\|\p^{\al-\beta}h(s)\|_{L^\infty_{x,v}}\Big\|_{L^2_s} \cr 
&\leq C \e \sum_{0\leq\beta\leq\al}  \Big\|\|\p^{\beta}h(s)\|_{L^\infty_{x,v}}\|\p^{\al-\beta}h(s)\|_{L^\infty_{x,v}}\Big\|_{L^2_s}.
\enda
\end{align}
To estimate \( I^h_4 \), we use Young’s inequality, Minkowski’s integral inequality, the bound \eqref{wgammav}, and the kernel property $\eqref{kw}_2$:
\begin{align}\label{hI4}
\bega
\|I^h_4\|_{L^2_t}
&\leq C\bigg\|\frac{1}{\e^2\kappa}e^{\frac{-\nu(v)t}{\e^2 \kappa}}\bigg\|_{L^1_t} \sum_{0\leq\beta\leq\al}\int_{\R^3}|\mathbf{k}_w(v,v_*)|  \e\bigg\|\frac{1}{\e^2\kappa}e^{ \frac{-\nu(v_*)s}{\e^2 \kappa}}\bigg\|_{L^1_s} \cr 
&\quad \times \nu(v_*)\Big\|\|\p^{\beta}h(\tau)\|_{L^\infty_{x,v}}\|\p^{\al-\beta}h(\tau)\|_{L^\infty_{x,v}}\Big\|_{L^2_\tau} dv_* \cr 
&\leq C\e\sum_{0\leq\beta\leq\al} \Big\|\|\p^{\beta}h(\tau)\|_{L^\infty_{x,v}}\|\p^{\al-\beta}h(\tau)\|_{L^\infty_{x,v}}\Big\|_{L^2_\tau}.
\enda
\end{align}
For \( I^h_3 \), we split it into two parts, \( I^h_{31} \) and \( I^h_{32} \), based on the integration region in \( d\tau \): namely, \( \tau \in [s - \bar{\delta} \e^2 \kappa, s] \) and \( \tau \in [0, s - \bar{\delta} \e^2 \kappa] \), respectively.
Using Young's convolution inequality, Minkowski’s integral inequality, and Hölder's inequality, we estimate \( I^h_{31} \) as
\begin{align}\label{hI31}
\bega
\|I^h_{31}\|_{L^2_t} &\les \bigg\|\frac{1}{\e^2 \kappa}e^{\frac{-\nu(v)t}{\e^2 \kappa}} \bigg\|_{L^1_t } \bigg\|\int_{\R^3}\mathbf{k}_w(v,v_*) \frac{1}{\e^2 \kappa} \int_{s-\bar{\delta}\e^2\kappa}^s e^{\frac{-\nu(v_*)(s-\tau)}{\e^2 \kappa}} \frac{1}{1+|v_*|} \|\p^{\al}h(\tau)\|_{L^\infty_{x,v}} d\tau dv_* \bigg\|_{L^2_s} \cr 
&\les \bigg\|\int_{\R^3}\mathbf{k}_w(v,v_*) \frac{1}{\e^2 \kappa} \lw(\int_{s-\bar{\delta}\e^2\kappa}^s  e^{\frac{-2\nu(v_*)(s-\tau)}{\e^2 \kappa}} d\tau\rw)^{\frac{1}{2}} \lw(\int_{s-\bar{\delta}\e^2\kappa}^s\|\p^{\al}h(\tau)\|_{L^\infty_{x,v}}^2 d\tau\rw)^{\frac{1}{2}} dv_* \bigg\|_{L^2_s} \cr 
&\les \sqrt{\bar{\delta}}\|\p^{\al}h\|_{L^2_tL^\infty_{x,v}},
\enda
\end{align}
where we used
\begin{align*}
\int_0^t\int_{s-\bar{\delta}\e^2\kappa}^s\|\p^{\al}h(\tau)\|_{L^\infty_{x,v}}^2 d\tau ds = \int_0^t \int_{\tau}^{\tau+\bar{\delta}\e^2\kappa}\|\p^{\al}h(\tau)\|_{L^\infty_{x,v}}^2 ds d\tau = \bar{\delta}\e^2\kappa \|\p^{\al}h\|_{L^2_tL^\infty_{x,v}}^2.
\end{align*}
For \( I^h_{32} \), we decompose the kernel \( \mathbf{k}_w(v,v_*) \) as \(\mathbf{k}_w = \mathbf{k}_{w,\mathfrak{N}} + (\mathbf{k}_w - \mathbf{k}_{w,\mathfrak{N}})\), where
\begin{align*}
\mathbf{k}_{w,\mathfrak{N}}(v,v_*) &:= 1_{\frac{1}{\mathfrak{N}}<|v-v_*|<\mathfrak{N}, |v|<\mathfrak{N}} \mathbf{k}_w(v,v_*), \\
\mathbf{k}_{w,R}(v,v_*) &:= \mathbf{k}_w(v,v_*) - \mathbf{k}_{w,\mathfrak{N}}(v,v_*).
\end{align*}
To simplify the notation, we define the following term:
\begin{align*}
I^h_{32}(\mathbf{k}_w^1,\mathbf{k}_w^2) &:= \frac{1}{\e^2 \kappa} \int_0^t e^{ \frac{-\nu(v)(t-s)}{\e^2 \kappa}} \int_{\R^3} \mathbf{k}_w^1(v,v_*) \frac{1}{\e^2 \kappa} \int_0^{s-\bar{\delta}\e^2\kappa} e^{\frac{-\nu(v_*)(s-\tau)}{\e^2 \kappa}} \\
&\quad \int_{\R^3} \mathbf{k}_w^2(v_*,v_{**}) \p^{\al}h\bigg(\tau, x-\frac{v}{\e}(t-s)-\frac{v_*}{\e}(s-\tau),v_{**}\bigg) dv_{**} d\tau dv_* ds,
\end{align*}
for \( \{\mathbf{k}_w^1, \mathbf{k}_w^2\} \in \{ \mathbf{k}_{w,\mathfrak{N}}, \mathbf{k}_{w,R} \} \).
Then \( I^h_{32}(\mathbf{k}_w, \mathbf{k}_w) \) can be written as
\begin{align*}
I^h_{32}(\mathbf{k}_w,\mathbf{k}_w) = I^h_{32}(\mathbf{k}_{w,\mathfrak{N}},\mathbf{k}_{w,\mathfrak{N}}) + I^h_{32}(\mathbf{k}_{w,R},\mathbf{k}_{w,\mathfrak{N}}) + I^h_{32}(\mathbf{k}_{w,\mathfrak{N}},\mathbf{k}_{w,R}) + I^h_{32}(\mathbf{k}_{w,R},\mathbf{k}_{w,R}).
\end{align*}
For the last three terms, each involving at least one \( \mathbf{k}_{w,R} \), we use estimate \eqref{kw} to have
\begin{align}\label{hI32a}
\int_{\R^3} \mathbf{k}_{w,R}(v,v_*) dv_* &\les \lw(\frac{1}{\mathfrak{N}} + e^{-\mathfrak{N}^2}\rw).
\end{align}
Applying Young’s inequality twice in time and using Minkowski’s integral inequality, the quantity \( I^h_{32}(\mathbf{k}_w^1, \mathbf{k}_w^2) \) satisfies
\begin{align}\label{hI32b}
\bega
\|I^h_{32}(\mathbf{k}_w^1,\mathbf{k}_w^2)\|_{L^2_t}&\les \bigg\| \frac{1}{\e^2 \kappa}e^{\frac{-\nu(v)t}{\e^2 \kappa}}\bigg\|_{L^1_t} \int_{\R^3}\mathbf{k}_w^1(v,v_*) \bigg\|\frac{1}{\e^2 \kappa} e^{\frac{-\nu(v_*)s}{\e^2 \kappa}}\bigg\|_{L^1_s} \cr 
&\quad \times \int_{\R^3}\mathbf{k}_w^2(v_*,v_{**}) \|\p^{\al}h\|_{L^2_tL^\infty_{x,v}}dv_{**} d\tau dv_* \cr 
&\les \int_{\R^3}\mathbf{k}_w^1(v,v_*) \int_{\R^3}\mathbf{k}_w^2(v_*,v_{**}) dv_{**} d\tau dv_* \|\p^{\al}h\|_{L^2_tL^\infty_{x,v}}.
\enda
\end{align}
Combining \eqref{hI32a} and \eqref{hI32b}, we obtain
\begin{align}\label{hI32R}
\|I^h_{32}(\mathbf{k}_{w,R},\mathbf{k}_{w,\mathfrak{N}})\|_{L^2_t} + \|I^h_{32}(\mathbf{k}_{w,\mathfrak{N}},\mathbf{k}_{w,R})\|_{L^2_t} + \|I^h_{32}(\mathbf{k}_{w,R},\mathbf{k}_{w,R})\|_{L^2_t} &\les \lw(\frac{1}{\mathfrak{N}} + e^{-\mathfrak{N}^2}\rw) \|\p^{\al}h\|_{L^2_tL^\infty_{x,v}}.
\end{align}
We note that \( I^h_{32}(\mathbf{k}_{w,\mathfrak{N}}, \mathbf{k}_{w,\mathfrak{N}}) \) is the main contribution. In the definition of \( \mathbf{k}_{w,\mathfrak{N}} \), the conditions \( |v_* - v_{**}| \leq \mathfrak{N} \) and \( |v_*| \leq \mathfrak{N} \) imply \( |v_{**}| \leq 2\mathfrak{N} \) and hence \( w(v_{**}) \leq e^{4c_1 \mathfrak{N}^2} \).
Moreover, since \( |v - v_*| \geq \frac{1}{\mathfrak{N}} \), it follows from $\eqref{kw}_1$ that
\begin{align*}
\mathbf{k}_{w,\mathfrak{N}}(v,v_*) \leq 1_{|v-v_*| \geq \frac{1}{\mathfrak{N}}}\frac{1}{|v-v_*|}e^{-C|v-v_*|^2} \leq \mathfrak{N}e^{-\frac{C}{\mathfrak{N}^2}}.
\end{align*}
Thus, we obtain
\begin{align*}
\bega
I^h_{32}&(\mathbf{k}_{w,\mathfrak{N}},\mathbf{k}_{w,\mathfrak{N}})\leq C_\mathfrak{N} \frac{1}{\e^2 \kappa}\int_0^t e^{ \frac{-\nu(v)(t-s)}{\e^2 \kappa}} \frac{1}{\e^2 \kappa} \int_0^{s-\bar{\delta}\e^2\kappa} e^{ \frac{-(s-\tau)}{\e^2 \kappa}} \cr 
&\quad \times \int_{|v_*|\leq \mathfrak{N}} \int_{|v_{**}|\leq 2\mathfrak{N}}\bigg|\p^{\al}f\bigg(\tau, x-\frac{v}{\e}(t-s)-\frac{v_*}{\e}(s-\tau),v_{**}\bigg)\bigg| dv_{**}dv_* d\tau  ds .
\enda
\end{align*}
We apply Hölder’s inequality in the \( v_{**} \)- and \( v_* \)-integrals, and perform the change of variables \( x_i - \frac{v_i}{\e}(t - s) - \frac{v_{*i}}{\e}(s - \tau) = y_i \) for \( i = 1, \ldots, d \). Then \( \left( \frac{s - \tau}{\e} \right)^2 dv_{*1}dv_{*2} = dy \) for $d=2$ and \( \left( \frac{s - \tau}{\e} \right)^3 dv_* = dy \) for $d=3$. Since $\int_{|v_*|\leq \mathfrak{N}}1 dv_{*3}\leq C_\mathfrak{N}$ for $d=2$, we have 
\begin{align*}
I^h_{32}(\mathbf{k}_{w,\mathfrak{N}},\mathbf{k}_{w,\mathfrak{N}}) &\leq C_\mathfrak{N} \frac{1}{\e^2 \kappa}\int_0^t e^{ \frac{-\nu(v)(t-s)}{\e^2 \kappa}} \frac{1}{\e^2 \kappa} \int_0^{s-\bar{\delta}\e^2\kappa} e^{ \frac{-(s-\tau)}{\e^2 \kappa}} \lw(\int_{|v_*|\leq \mathfrak{N}} \lw(\int_{|v_{**}|\leq 2\mathfrak{N}}1^2 dv_{**}\rw)^{\frac{q}{2}}dv_*\rw)^{\frac{1}{q}}\cr 
& \times \lw((1+\mathfrak{N}) \int_{\R^d} \lw(\frac{\e}{s-\tau}\rw)^d\lw(\int_{|v_{**}|\leq 2\mathfrak{N}}|\p^{\al}f(\tau, y,v_{**})|^2 dv_{**}\rw)^{\frac{p}{2}}dy\rw)^{\frac{1}{p}}  d\tau  ds ,
\end{align*}
where \( (p, q) \in [1, \infty]^2 \) satisfies \( \frac{1}{p} + \frac{1}{q} = 1 \). 
Since \( \tau \leq s - \bar{\delta} \e^2 \kappa \), we estimate
\begin{align*}
\lw(\frac{\e}{s-\tau}\rw)^d \les \lw(\frac{\e}{\bar{\delta}\e^2\kappa}\rw)^d \les \frac{1}{(\e\kappa)^d}.
\end{align*}
Thus, we derive
\begin{align*}
I^h_{32}&(\mathbf{k}_{w,\mathfrak{N}},\mathbf{k}_{w,\mathfrak{N}})\leq C_\mathfrak{N} \frac{1}{\e^2 \kappa}\int_0^t e^{ \frac{-\nu(v)(t-s)}{\e^2 \kappa}} \frac{1}{\e^2 \kappa} \int_0^{s-\bar{\delta}\e^2\kappa} e^{ \frac{-(s-\tau)}{\e^2 \kappa}} \frac{1}{(\e\kappa)^{\frac{d}{p}}} \|\mathbf{1}_{|v|\leq2\mathfrak{N}}\p^{\al}f(\tau)\|_{L^p_xL^2_v}  d\tau  ds.
\end{align*}
Applying Young’s inequality twice, we obtain
\begin{align}\label{I32r}
\big\|I^h_{32}(\mathbf{k}_{w,\mathfrak{N}},\mathbf{k}_{w,\mathfrak{N}})\big\|_{L^2_t} 
&\les C_\mathfrak{N} \frac{1}{(\e\kappa)^{\frac{d}{p}}} \big\|\mathbf{1}_{|v|\leq2\mathfrak{N}}\p^{\al}f\big\|_{L^2_tL^p_xL^2_v}.
\end{align}
Combining the estimates \eqref{hI1}, \eqref{hI2}, \eqref{hI5}, \eqref{hI4}, \eqref{hI31}, \eqref{hI32R}, and \eqref{I32r}, we conclude
\begin{align*}
\bega
\lw(1-C\sqrt{\bar{\delta}}-\frac{C}{\mathfrak{N}}-Ce^{-\mathfrak{N}^2}\rw)\|\p^{\al}h\|_{L^2_tL^\infty_{x,v}}&\les \e\kappa^{\frac{1}{2}}\|\p^{\al}h_0\|_{L^\infty_{x,v}} +C_\mathfrak{N} \frac{1}{(\e\kappa)^{\frac{d}{p}}} \|\mathbf{1}_{|v|\leq2\mathfrak{N}}\p^{\al}f\|_{L^2_tL^p_{x,v}} \cr 
&+\e \sum_{0\leq\beta\leq\al} \big\|\|\p^{\beta}h(t)\|_{L^\infty_{x,v}} \|\p^{\al-\beta}h(t)\|_{L^\infty_{x,v}}\big\|_{L^2_t}.
\enda
\end{align*}
By choosing \( \bar{\delta} > 0 \) sufficiently small and \( \mathfrak{N} > 0 \) sufficiently large such that \(\lw(C\sqrt{\bar{\delta}}+\frac{C}{\mathfrak{N}}+Ce^{-\mathfrak{N}^2}\rw)\leq 1/2\), we conclude the desired estimate \eqref{hL2t}.

\end{proof}

Now, we estimate the linear term in \eqref{hL2t}, namely \( \|\mathbf{1}_{|v| \leq 2\mathfrak{N}} \p^{\al}f^{\e}\|_{L^p_x L^2_v} \).
In the infinite velocity energy case \eqref{caseECX}, attempting to obtain an $L^2_x$ estimate for $f^{\e}$ again produces contributions involving $M^{\e}-\mu$, which in turn require control of $\|u^{\e}\|_{L^2_x}$.
To avoid imposing an $L^2_x$ condition on $u^{\e}$, we instead derive an $L^\infty$ estimate for $f^{\e}$ in the case without derivatives.
This is possible because, in Lemma~\ref{L.htt}, we obtained estimates in $L^p_x$ for any $p\in[1,\infty]$, and we only need to control the velocity cut-off part of $f^{\e}$ in the $L^\infty$ norm.

\begin{lemma}\label{L.ftoE}
For $t\in[0,T]$ satisfying the bootstrap assumption \eqref{condition}, and for both types of derivatives \( \p^{\al} \) defined in \eqref{caseA} and \eqref{caseB}, and for the constant \( \mathfrak{N} > 0 \) chosen in Lemma \ref{L.htt}, we have 
\begin{align}\label{ftoE}
\begin{split}
\kappa^{(|\al|-\mathrm{N})_+}\big\|\mathbf{1}_{|v|\leq2\mathfrak{N}}\p^{\al}f^{\e}(t)\big\|_{L^2_{x,v}}^2 &\les C_{\mathfrak{N}} \mathcal{E}_{tot}^{\mathrm{N}}(F^{\e}(t)),
\end{split}
\end{align}
for $0 \leq |\al| \leq \mathrm{N}+1$ in the finite velocity energy case \eqref{caseEC}, 
and for $1 \leq |\al| \leq \mathrm{N}+1$ in the infinite velocity energy case \eqref{caseECX}. In addition, when $|\al|=0$ for the infinite velocity energy case \eqref{caseECX}, we have 
\begin{align}\label{ftoE2}
\begin{split}
\big\|\mathbf{1}_{|v|\leq2\mathfrak{N}}f^{\e}(t)\big\|_{L^\infty_xL^2_v}^2 &\les C_{\mathfrak{N}} \big(1+\mathcal{E}_{tot}^{\mathrm{N}}(F^{\e}(t))\big), \quad \mbox{for} \quad \eqref{caseECX}.
\end{split}
\end{align}

\end{lemma}
\begin{proof}
We first prove \eqref{ftoE}. Using the decomposition \( F = M^{\e} + \AC{\P}F^{\e} = \mu + \e \sqrt{\mu} f^{\e} \), we estimate
\begin{align}\label{u<N}
\begin{split}
\kappa^{(|\al|-\mathrm{N})_+}\big\|\mathbf{1}_{|v|\leq2\mathfrak{N}}\p^{\al}f^{\e} \big\|_{L^2_{x,v}}^2 &\leq \kappa^{(|\al|-\mathrm{N})_+}\lw\|\mathbf{1}_{|v|\leq2\mathfrak{N}}\frac{\p^{\al}(M^{\e}-\mu)}{\e\sqrt{\mu}}\rw\|_{L^2_{x,v}}^2 \cr 
&+ \kappa^{(|\al|-\mathrm{N})_+}\lw\|\mathbf{1}_{|v|\leq2\mathfrak{N}}\frac{\p^{\al}\AC{\P}F^{\e}}{\e\sqrt{\mu}}\rw\|_{L^2_{x,v}}^2.
\end{split}
\end{align}
For the first term in \eqref{u<N}, as in the estimate for \eqref{Vl1b}, we have
\begin{align*}
\kappa^{(|\al|-\mathrm{N})_+}\lw\|\mathbf{1}_{|v|\leq2\mathfrak{N}}\frac{\p^{\al}(M^{\e}-\mu)}{\e\sqrt{\mu}}\rw\|_{L^2_{x,v}}^2 &= \kappa^{(|\al|-\mathrm{N})_+}\frac{1}{\eps^2}\int_{\Omega}\int_{|v|\leq2\mathfrak{N}}|\p^{\al}(M^{\e}-\mu)|^2|\mu|^{-1} dvdx  \cr 
&\les \kappa^{(|\al|-\mathrm{N})_+}\|\p^{\al}(\rho^{\e},u^{\e},\ta^{\e})\|_{L^2_x}^2 \les \mathcal{E}_{tot},
\end{align*}
for $0 \leq |\al| \leq \mathrm{N}+1$ in the finite velocity energy case \eqref{caseEC}, 
and for $1 \leq |\al| \leq \mathrm{N}+1$ in the infinite velocity energy case \eqref{caseECX}.
For the second term in \eqref{u<N}, we use the definition of \( \mathcal{E}_G \) in \eqref{N-EDdef} for \( 0 \leq |\al| \leq \mathrm{N} \), and apply \eqref{Eal=5} along with the definition of \( \mathcal{E}_{\mathrm{top}} \) in \eqref{N-EDdef2} for \( |\al| = \mathrm{N} + 1 \), to obtain
\begin{align}\label{Gv<N}
\bega
\kappa^{(|\al|-\mathrm{N})_+}\lw\|\mathbf{1}_{|v|\leq2\mathfrak{N}}\frac{\p^{\al}\AC{\P}F^{\e}}{\e\sqrt{\mu}}\rw\|_{L^2_{x,v}}^2 &\leq \kappa^{(|\al|-\mathrm{N})_+}\lw\|\mathbf{1}_{|v|\leq2\mathfrak{N}} \frac{\sqrt{M^{\e}}}{\sqrt{\mu}}\rw\|_{L^\infty_{x,v}}^2  \lw\|\frac{\p^{\al}\AC{\P}F^{\e}}{\e \sqrt{M^{\e}}} \rw\|_{L^2_{x,v}}^2 \cr 
&\les \begin{cases} C_{\mathfrak{N}} \mathcal{E}_G, \quad &\mbox{when} \quad 0\leq|\al|\leq\mathrm{N}, \\ 
C_{\mathfrak{N}}(\mathcal{E}_{top}+\e^2\mathcal{E}_M^2), \quad &\mbox{when} \quad |\al|=\mathrm{N}+1, \end{cases}
\enda
\end{align}
where we used \( \lw\|\mathbf{1}_{|v|\leq2\mathfrak{N}} \frac{\sqrt{M^{\e}}}{\sqrt{\mu}}\rw\|_{L^\infty_{x,v}} \leq C \), which follows from condition \( (\mathcal{B}_2) \) in \eqref{condition}. \\ 
For the proof of \eqref{ftoE2}, we take the $L^\infty_x$ norm instead of $L^2_x$ in \eqref{u<N}, which gives 
\begin{align}\label{u<N2}
\bega
\big\|\mathbf{1}_{|v|\leq2\mathfrak{N}}f\big\|_{L^\infty_xL^2_v}^2 &\leq \lw\|\mathbf{1}_{|v|\leq2\mathfrak{N}}\frac{(M^{\e}-\mu)}{\e\sqrt{\mu}}\rw\|_{L^\infty_xL^2_v}^2 + \lw\|\mathbf{1}_{|v|\leq2\mathfrak{N}}\frac{\AC{\P}F^{\e}}{\e\sqrt{\mu}}\rw\|_{L^\infty_xL^2_v}^2 .
\enda
\end{align}
For the first term in \eqref{u<N2}, we use \eqref{rutscale} and \eqref{uLinfEC} to obtain 
\begin{align*}
\lw\|\mathbf{1}_{|v|\leq2\mathfrak{N}}\frac{(M^{\e}-\mu)}{\e\sqrt{\mu}}\rw\|_{L^\infty_xL^2_v}^2 \leq \|(\rho^{\e},u^{\e},\ta^{\e})\|_{L^\infty_x}^2 \les (1+\mathcal{E}_{tot}).
\end{align*}
For the second term in \eqref{u<N2}, we apply Agmon's inequality \eqref{Agmon} to the microscopic part, which yields
\begin{align*}
\lw\|\mathbf{1}_{|v|\leq2\mathfrak{N}}\frac{\AC{\P}F^{\e}}{\e\sqrt{\mu}}\rw\|_{L^\infty_xL^2_v}^2 
\leq \lw\|\mathbf{1}_{|v|\leq2\mathfrak{N}} \frac{\sqrt{M^{\e}}}{\sqrt{\mu}}\rw\|_{L^\infty_{x,v}}^2  \sum_{0\leq|\al_x|\leq2}\lw\|\frac{\p^{\al_x}\AC{\P}F^{\e}}{\e \sqrt{M^{\e}}} \rw\|_{L^2_{x,v}}^2 \leq \mathcal{E}_{G}.
\end{align*}
Thus, the desired result \eqref{ftoE2} follows.
\end{proof}

\begin{proof}[\textbf{Proof of Proposition \ref{P.hLinf2}}]
We first establish \eqref{hL2} and \eqref{hLinf} under the assumption that \eqref{claimx4} holds, and subsequently prove \eqref{claimx4}. Since the proofs of \eqref{hL2} and \eqref{hLinf} are similar, we begin with the proof of \eqref{hL2} and then briefly explain the modifications required for \eqref{hLinf}.
Depending on whether the finite velocity energy condition \eqref{caseEC} or the infinite velocity energy condition \eqref{caseECX} applies, we employ different estimates for the linear term in \eqref{ftoE} and \eqref{ftoE2}, respectively. In the finite velocity energy case \eqref{caseEC}, we apply  \(\sum_{0\leq|\al|\leq\mathrm{N}+1}\kappa^{\frac{(|\al|-\mathrm{N})_+}{2}}\) to \eqref{hL2t} with \(p = 2\) in Lemma \ref{L.htt}, which yields
\begin{align}\label{h=I+II}
\sum_{0\leq|\al|\leq\mathrm{N}+1}\kappa^{\frac{(|\al|-\mathrm{N})_+}{2}}\|\p^{\al}h\|_{L^2_tL^\infty_{x,v}}&\leq C \e\kappa^{\frac{1}{2}}\sum_{0\leq|\al|\leq\mathrm{N}+1}\kappa^{\frac{(|\al|-\mathrm{N})_+}{2}}\|\p^{\al}h_0\|_{L^\infty_{x,v}} + I_h + II_h,
\end{align}
where
\begin{align*}
I_h&= \frac{C}{(\e\kappa)^{\frac{d}{2}}}\sum_{0\leq|\al|\leq\mathrm{N}+1}\kappa^{\frac{(|\al|-\mathrm{N})_+}{2}} \|\mathbf{1}_{|v|\leq2\mathfrak{N}}\p^{\al}f\|_{L^2_tL^2_{x,v}}, \cr
II_h&= C \e \sum_{0\leq|\al|\leq\mathrm{N}+1}\sum_{\beta+\gamma=\al} \Big\|\Big(\kappa^{\frac{(|\beta|-N)_+}{2}}\|\p^{\beta}h(t)\|_{L^\infty_{x,v}}\Big) \Big(\kappa^{\frac{(|\gamma|-N)_+}{2}}\|\p^{\gamma}h(t)\|_{L^\infty_{x,v}}\Big) \Big\|_{L^2_t}.
\end{align*}
For the infinite velocity energy case \eqref{caseECX}, we apply 
\(\sum_{1\leq|\al|\leq\mathrm{N}+1}\kappa^{\frac{(|\al|-\mathrm{N})_+}{2}}\) 
to \eqref{hL2t} with \(p = 2\) in Lemma \ref{L.htt}, and additionally combine it with \eqref{hL2t} for \(|\al|=0\) with \(p = \infty\). 
Then the inequality \eqref{h=I+II} holds if we replace \(I_h\) in \eqref{h=I+II} with the following modified term \(I_h'\):
\begin{align*}
I_h'&= \frac{C}{(\e\kappa)^{\frac{d}{2}}}\sum_{1\leq|\al|\leq\mathrm{N}+1}\kappa^{\frac{(|\al|-\mathrm{N})_+}{2}} \|\mathbf{1}_{|v|\leq2\mathfrak{N}}\p^{\al}f\|_{L^2_tL^2_{x,v}} + \|\mathbf{1}_{|v|\leq2\mathfrak{N}}f\|_{L^2_tL^\infty_xL^2_v}.
\end{align*}
For the term $II_h$, the weight 
\(\kappa^{\frac{(|\al| - \mathrm{N})_+}{2}}\) 
is distributed to the nonlinear term in \(II_h\) according to the following rule:
\begin{align*}
\kappa^{\frac{(|\al|-\mathrm{N})_+}{2}}\p^{\al}(XY) &= \begin{cases} 
\p^{\beta}X \p^{\al-\beta}Y, \quad &\text{if} \quad 0 \leq |\al|, |\al-\beta|\leq \mathrm{N}, \\ 
(\kappa^{\frac{1}{2}}\p^{\al}X) Y + X(\kappa^{\frac{1}{2}}\p^{\al}Y), \quad &\text{otherwise}.  
\end{cases} 
\end{align*}
We apply \eqref{ftoE} in Lemma \ref{L.ftoE} to $I_h$, and apply both \eqref{ftoE} and \eqref{ftoE2} to $I_h'$, to obtain
\begin{align}\label{Ih}
\bega
I_h &\leq \frac{C}{(\e\kappa)^{\frac{d}{2}}}\bigg(\int_0^t  \mathcal{E}_{tot}(s) ds\bigg)^{\frac{1}{2}}, & \mbox{for the case} \quad \eqref{caseEC}, \cr 
I_h' &\leq \frac{C}{(\e\kappa)^{\frac{d}{2}}}\bigg(\int_0^t  \mathcal{E}_{tot}(s) ds\bigg)^{\frac{1}{2}} +\bigg(\int_0^t  \big(1+\mathcal{E}_{tot}(s)\big) ds\bigg)^{\frac{1}{2}} , & \mbox{for the case} \quad \eqref{caseECX}.
\enda
\end{align}

To estimate \( II_h \), we divide the case when $\p^{\al}$ is purely spatial derivative case \eqref{caseA} or the space-time derivative case \eqref{caseB}. We distribute the derivatives in $II_h$ as follows:
\begin{align}\label{IIh-00}
II_h &\leq C\e \sum_{0\leq|\al|\leq\mathrm{N}+1}\sum_{\beta+\gamma=\al} \Big(\kappa^{\frac{(|\beta|-\mathrm{N})_+}{2}}\|\p^{\beta}h(t)\|_{L^2_tL^\infty_{x,v}}\Big) \Big(\kappa^{\frac{(|\gamma|-\mathrm{N})_+}{2}}\|\p^{\gamma}h(t)\|_{L^\infty_tL^\infty_{x,v}}\Big).
\end{align}
(1) When $\p^{\al}$ is purely spatial derivative case \eqref{caseA}, the lower number of derivative should be less that $\lfloor (\mathrm{N}+1)/2 \rfloor$. Hence, applying the estimates for \( I_h \) (or \( I_h' \)) and \( II_h \) from \eqref{Ih} and \eqref{IIh-00} to \eqref{h=I+II}, and then move the nonlinear term to the left-hand side yields
\begin{align}\label{h-hh-(1)}
\bega
\Big(1-&C\e\sum_{0\leq |\gamma_x|\leq \lfloor (\mathrm{N}+1)/2 \rfloor}\|\p^{\gamma_x}h(t)\|_{L^\infty_tL^\infty_{x,v}}\Big)\sum_{0\leq|\al_x|\leq\mathrm{N}+1}\kappa^{\frac{(|\al_x|-\mathrm{N})_+}{2}}\|\p^{\al}h\|_{L^2_tL^\infty_{x,v}} \cr 
&\les C \e\kappa^{\frac{1}{2}}\sum_{0\leq|\al_x|\leq\mathrm{N}+1}\kappa^{\frac{(|\al|-\mathrm{N})_+}{2}}\|\p^{\al_x}h_0\|_{L^\infty_{x,v}}+\frac{C}{(\e\kappa)^{\frac{d}{2}}}\bigg(\int_0^t  \big(1+\mathcal{E}_{tot}(s)\big) ds\bigg)^{\frac{1}{2}}.
\enda
\end{align}
(2) When $\p^{\al}$ is space-time derivative case \eqref{caseB}, recall the notation $\p_{\tilde{t}}=\e^{\mathfrak{n}}\p_t$. Without loss of generality, we assume that $\p^\beta$ in \eqref{IIh-00} include the derivative $\p_{\tilde{t}}$. Then by the same way, we have 
\begin{align}\label{h-hh}
\bega
\Big(1-&C\e\sum_{0\leq |\gamma_x|\leq \mathrm{N}}\|\p^{\gamma_x}h(t)\|_{L^\infty_tL^\infty_{x,v}}\Big)\sum_{0\leq|\al|\leq\mathrm{N}+1}\kappa^{\frac{(|\al|-\mathrm{N})_+}{2}}\|\p^{\al}h\|_{L^2_tL^\infty_{x,v}} \cr 
&\les C \e\kappa^{\frac{1}{2}}\sum_{0\leq|\al|\leq\mathrm{N}+1}\kappa^{\frac{(|\al|-\mathrm{N})_+}{2}}\|\p^{\al}h_0\|_{L^\infty_{x,v}}+\frac{C}{(\e\kappa)^{\frac{d}{2}}}\bigg(\int_0^t  \big(1+\mathcal{E}_{tot}(s)\big) ds\bigg)^{\frac{1}{2}}.
\enda
\end{align}

If the inequality \eqref{claimx4} holds, then \eqref{h-hh-(1)} or \eqref{h-hh} implies the $L^2_t$ result in \eqref{hL2}.
To prove $L^\infty_t$ estimate in \eqref{hLinf}, we follow the same approach. From \eqref{hLinft}, we write the corresponding part of \eqref{h-hh-(1)} or \eqref{h-hh} in \( L^\infty_t \) as
\begin{align*}
\bega
\Big(1-&C\e\sum_{0\leq |\gamma_x|\leq \mathrm{N}_*}\|\p^{\gamma_x}h(t)\|_{L^\infty_tL^\infty_{x,v}}\Big)\sum_{0\leq|\al|\leq\mathrm{N}+1}\kappa^{\frac{(|\al|-\mathrm{N})_+}{2}}\|\p^{\al}h\|_{L^\infty_tL^\infty_{x,v}} \cr 
&\les C \sum_{0\leq|\al|\leq\mathrm{N}+1}\kappa^{\frac{(|\al|-\mathrm{N})_+}{2}}\|\p^{\al}h_0\|_{L^\infty_{x,v}}+\frac{C}{(\e\kappa)^{\frac{d}{2}}}\|(1+\mathcal{E}_{tot})\|_{L^\infty_t}^{\frac{1}{2}},
\enda
\end{align*}
here $\mathrm{N}_*=\lfloor (\mathrm{N}+1)/2 \rfloor$ for the purely spatial derivative case \eqref{caseA} or $\mathrm{N}_*=\mathrm{N}$ for the space-time derivative case \eqref{caseB}. By the same way, if \eqref{claimx4} holds, we conclude the $L^\infty_t$ result in \eqref{hLinf}. \newline

{\bf Proof of \eqref{claimx4}} We now prove two inequalities in \eqref{claimx4} depending on the corresponding two cases. \\ 
(1) When $\p^{\al}$ is purely spatial derivative case \eqref{caseA}, we apply \( \sum_{1\leq |\al_x|\leq \lfloor (\mathrm{N}+1)/2 \rfloor } \) to \eqref{hLinft} in Lemma \ref{L.htt} and add \eqref{hLinft} for \( |\al_x| = 0 \) and \( p = \infty \) to have 
\begin{align}\label{10.20}
\bega
\bigg(&1-C\e\sum_{0\leq |\al_x|\leq \lfloor (\mathrm{N}+1)/2 \rfloor}\|\p^{\al_x}h\|_{L^\infty_tL^\infty_{x,v}}\bigg)\sum_{0\leq |\al_x|\leq \lfloor (\mathrm{N}+1)/2 \rfloor}\|\p^{\al_x}h\|_{L^\infty_tL^\infty_{x,v}} \cr 
&\les \sum_{0\leq |\al_x|\leq \lfloor (\mathrm{N}+1)/2 \rfloor}\|\p^{\al_x}h_0\|_{L^\infty_{x,v}}+\|(1+\mathcal{E}_{tot})\|_{L^\infty_t}^{\frac{1}{2}}+\frac{1}{(\e\kappa)^{\frac{d}{p}}}\sum_{1\leq |\al_x|\leq \lfloor (\mathrm{N}+1)/2 \rfloor} \|\mathbf{1}_{|v|\leq2\mathfrak{N}}\p^{\al_x}f\|_{L^\infty_tL^p_xL^2_v},
\enda
\end{align}
for $p \in [1,\infty]$.
Depending on the dimension, we choose different number $p$ for the linear term. We choose $p \gg 1 $ and apply interpolation inequality~\eqref{Ga-Ni} for $d=2$ and choose $p=\infty$ and use Agmon inequality \eqref{Agmon} for $d=3$ to have 
\begin{align*}
\|\mathbf{1}_{|v|\leq2\mathfrak{N}}\p^{\al_x}f\|_{L^p_xL^2_v} &\leq C\|\mathbf{1}_{|v|\leq2\mathfrak{N}}\nabla_x\p^{\al_x}f\|_{L^2_{x,v}}^{1-\frac{2}{p}} \|\mathbf{1}_{|v|\leq2\mathfrak{N}}\p^{\al_x}f\|_{L^2_{x,v}}^{\frac{2}{p}} \leq C_{\mathfrak{N}} \mathcal{E}_{tot}(t), &\mbox{for} \quad d=2, \cr 
\|\mathbf{1}_{|v|\leq2\mathfrak{N}}\p^{\al_x}f\|_{L^\infty_xL^2_v} &\leq C\|\mathbf{1}_{|v|\leq2\mathfrak{N}}\nabla_x\p^{\al_x}f\|_{H^1_xL^2_v}^{\frac{1}{2}} \|\mathbf{1}_{|v|\leq2\mathfrak{N}}\p^{\al_x}f\|_{H^2_xL^2_v}^{\frac{1}{2}} \leq C_{\mathfrak{N}} \kappa^{-\frac{1}{4}}\mathcal{E}_{tot}(t) , &\mbox{for} \quad d=3,
\end{align*}
where we used Lemma \ref{L.ftoE} for $|\al_x|\leq \lfloor (\mathrm{N}+1)/2 \rfloor$. 
Multiplying \eqref{10.20} by \( \e \), we obtain
\begin{align*}
\bigg(1-&C\e\sum_{0\leq |\al_x|\leq \lfloor (\mathrm{N}+1)/2 \rfloor}\|\p^{\al_x}h\|_{L^\infty_tL^\infty_{x,v}}\bigg)\bigg(\e\sum_{0\leq |\al_x|\leq \lfloor (\mathrm{N}+1)/2 \rfloor}\|\p^{\al_x}h\|_{L^\infty_tL^\infty_{x,v}}\bigg) \cr 
&\les \e\sum_{0\leq |\al_x|\leq \lfloor (\mathrm{N}+1)/2 \rfloor}\|\p^{\al_x}h_0\|_{L^\infty_{x,v}} +\e \|(1+\mathcal{E}_{tot})\|_{L^\infty_t}^{\frac{1}{2}} +Y_1,
\end{align*}
where
\begin{align*}
Y_1= 
\begin{cases} 
\e(\e\kappa)^{-\frac{d}{p}}\|\mathcal{E}_{tot}\|_{L^\infty_t}^{\frac{1}{2}}, \quad \mbox{for} \quad d=2, \\ 
\e\kappa^{-\frac{1}{4}}\|\mathcal{E}_{tot}\|_{L^\infty_t}^{\frac{1}{2}}, \quad \hspace{5mm} \mbox{for}  \quad d=3.
\end{cases}
\end{align*}
Now, define \( Z := \e\sum_{0\leq |\al_x|\leq \mathrm{N}}\|\p^{\al_x}h\|_{L^\infty_tL^\infty_{x,v}} \). Under the Bootstrap assumption in Assumption \ref{Boot}, for sufficiently large $p$ for $d=2$ this leads to the quadratic inequality\( (1 - CZ)Z < \frac{1}{4C} \), which yields \eqref{claimx4}.

(2) When $\p^{\al}$ is space-time derivative case \eqref{caseB},
We apply \( \sum_{1\leq |\al_x|\leq \mathrm{N}} \) to \eqref{hLinft} in Lemma \ref{L.htt}, 
and combine it with \eqref{hLinft} for \( |\al_x| = 0 \) and \( p = \infty \) as follows:
\begin{align}\label{10.20-t}
\bega
\sum_{0\leq |\al_x|\leq \mathrm{N}}&\|\p^{\al_x}h\|_{L^\infty_tL^\infty_{x,v}}\les \sum_{0\leq |\al_x|\leq \mathrm{N}}\|\p^{\al_x}h_0\|_{L^\infty_{x,v}}+\frac{1}{(\e\kappa)^{\frac{d}{p}}}\sum_{1\leq |\al_x|\leq \mathrm{N}} \|\mathbf{1}_{|v|\leq2\mathfrak{N}}\p^{\al_x}f\|_{L^\infty_tL^p_xL^2_v} \cr 
&+\|(1+\mathcal{E}_{tot})\|_{L^\infty_t}^{\frac{1}{2}} + \e  \sum_{0\leq |\al_x|\leq \mathrm{N}}\sum_{0\leq\beta_x\leq\al_x} \|\p^{\beta_x}h\|_{L^\infty_tL^\infty_{x,v}}\|\p^{\al_x-\beta_x}h\|_{L^\infty_tL^\infty_{x,v}}.
\enda
\end{align}
To minimize the singular scale of $\kappa$, we again choose different $p$ depending on the dimension. We choose $p\gg1$ for $d=2$ and apply interpolation inequality \eqref{Ga-Ni}, and choose $p=4$ for $d=3$ and apply Ladyzhenskaya's inequality \eqref{Lady}:
\begin{align*}
\|\mathbf{1}_{|v|\leq2\mathfrak{N}}\p^{\al_x}f\|_{L^p_xL^2_v} &\leq C\|\mathbf{1}_{|v|\leq2\mathfrak{N}}\nabla_x\p^{\al_x}f\|_{L^2_{x,v}}^{1-\frac{2}{p}} \|\mathbf{1}_{|v|\leq2\mathfrak{N}}\p^{\al_x}f\|_{L^2_{x,v}}^{\frac{2}{p}} \leq C_{\mathfrak{N}}\kappa^{-\frac{2}{p}}\mathcal{E}_{tot}(t),  &\mbox{for} \quad d=2, \cr 
\|\mathbf{1}_{|v|\leq2\mathfrak{N}}\p^{\al_x}f\|_{L^4_xL^2_v} &\leq C\|\mathbf{1}_{|v|\leq2\mathfrak{N}}\nabla_x\p^{\al_x}f\|_{L^2_{x,v}}^{\frac{3}{4}} \|\mathbf{1}_{|v|\leq2\mathfrak{N}}\p^{\al_x}f\|_{L^2_{x,v}}^{\frac{1}{4}} \leq C_{\mathfrak{N}} \kappa^{-\frac{3}{8}}\mathcal{E}_{tot}(t) ,  &\mbox{for} \quad d=3,
\end{align*}
where we used Lemma \ref{L.ftoE} for $|\al_x|\leq \mathrm{N}$. Hence, multiplying \eqref{10.20-t} by \( \e \), we obtain
\begin{align*}
\bigg(1-&C\e\sum_{0\leq |\al_x|\leq \mathrm{N}}\|\p^{\al_x}h\|_{L^\infty_tL^\infty_{x,v}}\bigg)\bigg(\e\sum_{0\leq |\al_x|\leq \mathrm{N}}\|\p^{\al_x}h\|_{L^\infty_tL^\infty_{x,v}}\bigg) \cr 
&\les \e\sum_{0\leq |\al_x|\leq \mathrm{N}}\|\p^{\al_x}h_0\|_{L^\infty_{x,v}} +\e \|(1+\mathcal{E}_{tot})\|_{L^\infty_t}^{\frac{1}{2}} + Y_2,
\end{align*}
where
\begin{align*}
Y_2= 
\begin{cases} 
\e\kappa^{-\frac{1}{2}}(\e\kappa)^{-\frac{d}{p}}\|\mathcal{E}_{tot}\|_{L^\infty_t}^{\frac{1}{2}}, \quad \mbox{for} \quad d=2, \\ 
\e^{\frac{1}{4}}\kappa^{-\frac{9}{8}}\|\mathcal{E}_{tot}\|_{L^\infty_t}^{\frac{1}{2}}, \quad \hspace{5mm} \hspace{5mm} \mbox{for}  \quad d=3.
\end{cases}
\end{align*}
Now, define \( Z := \e\sum_{0\leq |\al_x|\leq \mathrm{N}}\|\p^{\al_x}h\|_{L^\infty_tL^\infty_{x,v}} \). Under the Bootstrap assumption in Assumption \ref{Boot}, for sufficiently large $p$ for $d=2$ this leads to the quadratic inequality\( (1 - CZ)Z < \frac{1}{4C} \), which yields \eqref{claimx4}.
\end{proof}

\refstepcounter{subsection}
\subsection*{\thesubsection\quad High-Order Moments Estimate for Low-Order Derivatives}
\label{Sec.V.L}

In this subsection, we estimate the high-order moments $\mathcal{V}^{\al}_{\ell}$ by a method different from that in Subsection~\ref{Sec.V.Top}. 
More precisely, we employ the decomposition $F^{\e} = M^{\e} + \e \sqrt{\bar{\mu}}g^{\e}$ where the global Maxwellian $\bar{\mu}$ corresponds to the parameters $(\mathrm{P},\mathrm{U},\mathrm{\Theta}) = (1,0,1-c_0/2)$ in \eqref{M-def}. In particular,
\begin{align}\label{muMdef}
\bar{\mu}(v) := M_{[1,0,1-c_0/2]} =  \frac{1}{(2\pi k_B (1-c_0/2))^{\frac{3}{2}}}\exp\left(-\frac{|v|^2}{2k_B(1-c_0/2)}\right).
\end{align}
Here, \(c_0 \ll 1\) is the small constant introduced in \eqref{condition}, which ensures that
\(|\mathrm{P}^{\e} -1|, |\mathrm{U}^{\e}|, |\mathrm{\Theta}^{\e} -1| < c_0/2 \).
There are two main reasons for adopting this alternative approach:
\begin{enumerate}
\item In the infinite velocity energy case \eqref{caseECX}, the norm $\|u^{\e}\|_{L^2_x}$ cannot be controlled by the energy. Consequently, the decomposition of $\mathcal{V}^{\al}_{\ell}$ provided in Proposition \ref{P.hLinf} is not applicable for $|\al|=0$. By employing the decomposition $F^{\e} = M^{\e} + \e \sqrt{\bar{\mu}}g^{\e}$, we can avoid the term $|M^{\e}-\mu|/\sqrt{\mu}$ that arises in Proposition \ref{P.hLinf2}.
\item The high-order moments $\mathcal{V}^{\al}_{\ell}$ can be estimated with a better scaling if we apply the decomposition $F^{\e} = M^{\e} + \e \sqrt{\bar{\mu}}g^{\e}$ only to the low-order derivatives $0 \leq |\al| \leq \mathrm{N}-2$. This leads to sharper bounds for the microscopic part associated with low-order derivatives, as derived from \eqref{totalGt} in Proposition \ref{P.G.Energy}. (See \eqref{Grefine.2D} in Theorem~\ref{T.2D.global}.)
\end{enumerate}
With the same exponential weight $w(v)$ defined in \eqref{hdef}, we introduce
\begin{align}\label{hdef2}
\bega
\mathfrak{h}^{\e}(t,x,v):=w(v)g^{\e}(t,x,v).
\enda
\end{align}
Unlike \eqref{wBE}, substituting $F^{\e}=M^{\e}+\e\sqrt{\bar{\mu}}g^{\e}$ into \eqref{BE} and then multiplying by $\tfrac{w(v)}{\e\sqrt{\bar{\mu}(v)}}$ yields
\begin{align}\label{h2eqn}
\bega
\bigg(\p_t \mathfrak{h}^{\e} + \frac{v}{\e}\cdot \nabla_x\mathfrak{h}^{\e} \bigg) 
&= -\frac{w}{\e\sqrt{\bar{\mu}}}\bigg(\p_tM^{\e}+\frac{v}{\e}\cdot\nabla_x M^{\e}\bigg) \cr 
&+ \frac{2}{\kappa \e^2} \frac{w}{\sqrt{\bar{\mu}}}\mathcal{N}(M^{\e},\sqrt{\bar{\mu}}g^{\e}) + \frac{1}{\kappa \e} \frac{w}{\sqrt{\bar{\mu}}}\mathcal{N}(\sqrt{\bar{\mu}}g^{\e},\sqrt{\bar{\mu}}g^{\e}).
\enda
\end{align}
Here, in contrast to the equation for $h^{\e}$ in \eqref{heqn}, the transport contribution of the local Maxwellian appears explicitly. Here, we define 
\begin{align}\label{LkMdef}
\bega
&\bar{L}g:= -\frac{2}{\sqrt{\bar{\mu}}}\mathcal{N}(M^{\e},\sqrt{\bar{\mu}}g) = \bar{\nu} g - \bar{\mathbf{K}} g,  \qquad \bar{\mathbf{K}} g(v):= \int_{\R^3}\bar{\mathbf{k}}(v,v_*) g(v_*)dv_*, \cr 
&\bar{\mathbf{k}}_{w}(v,v_*):= \bar{\mathbf{k}}(v,v_*)\frac{w(v)}{w(v_*)} ,  \qquad \hspace{1.7cm} \bar{\mathbf{K}}_{w} h(v) := \int_{\R^3}\bar{\mathbf{k}}_{w}(v,v_*) h(v_*)dv_*, \cr 
&\bar{\Gamma} (g,g) := \frac{1}{\sqrt{\bar{\mu}}}\mathcal{N}(\sqrt{\bar{\mu}}g,\sqrt{\bar{\mu}}g).
\enda
\end{align}
Thus, the equation can be rewritten as
\begin{align}\label{wBE2}
\bigg[\p_t + \frac{v}{\e} \cdot \nabla_x + \frac{\bar{\nu}}{\kappa \e^2} \bigg] \p^{\al}\mathfrak{h}^{\e}
&= \frac{1}{\kappa \e^2} \bar{\mathbf{K}}_{w} \p^{\al}\mathfrak{h}^{\e} + \frac{w}{\kappa \e } \sum_{0\leq\beta\leq\al}\binom{\al}{\beta}\bar{\Gamma}\bigg(\frac{\p^{\beta}\mathfrak{h}^{\e}}{w},\frac{\p^{\al-\beta}\mathfrak{h}^{\e}}{w}\bigg) -\frac{1}{\e}\mathrm{M}_{\text{trans}}^{\al}(t) ,
\end{align}
where the transport term of the local Maxwellian is written as
\begin{align}\label{TMdef}
\bega
\mathrm{M}_{\text{trans}}^{\al}(t) &:= \frac{w}{\sqrt{\bar{\mu}}}\p^{\al}\bigg(\p_tM^{\e}+\frac{v}{\e}\cdot\nabla_x M^{\e}\bigg).
\enda
\end{align}
When estimating the $L^\infty_x$ norm of $\mathrm{M}_{\text{trans}}^{\al}(t)$, we encounter terms of the form 
$\|\p^{\al}\nabla_x u^{\e}\|_{L^\infty_x}$ and $\|\p^{\al}\nabla_x \ta^{\e}\|_{L^\infty_x}$. 
This explains why the decomposition $F^{\e} = M^{\e} + \e \sqrt{\bar{\mu}}\, g^{\e}$ can be used only for derivatives satisfying $0 \leq |\al| \leq \mathrm{N}-2$, in contrast to \eqref{wBE}.

\begin{lemma}\label{L.muM}
For $t \in [0,T]$, where $T>0$ is defined in~\eqref{condition}, and on any bounded velocity domain $\{\,|v|\leq C\,\}$ with $C>0$, there exist constants $C_1>0$ and $C_2>0$ such that
\begin{align*}
C_1M^{\e}(t,x,v) \leq \bar{\mu}(v) \leq C_2M^{\e}(t,x,v),
\end{align*}
where $\bar{\mu}(v)$ is defined in \eqref{muMdef}.
\end{lemma}
\begin{proof}
By $(\mathcal{B}_2)$ in \eqref{condition}, we have 
\begin{align*}
1-c_0/2<\mathrm{\Theta}^{\e}(t,x) < 1+c_0/2 ,\quad \mbox{for} \quad (t,x)\in[0,T]\times \Omega.
\end{align*}
This immediately implies the desired result.
\end{proof}

We remark that the $L^\infty_{x,v}$ estimate of $\mathfrak{h}^{\e}$ is available only for $0 \leq |\al| \leq \mathrm{N}-2$, whereas the decomposition of $\mathcal{V}^{\al}_{\ell}(t)$ is valid for $0 \leq |\al| \leq \mathrm{N}+1$. 
Accordingly, we state Proposition \ref{P.hLinf3} for $0 \leq |\al| \leq \mathrm{N}+1$ and Proposition \ref{P.hLinf4} for $0 \leq |\al| \leq \mathrm{N}-2$.

\begin{proposition}\label{P.hLinf3}
Let $\Omega=\R^d$ with $d=2$ or $3$. Suppose the initial data $u_0^{\e}$ satisfies either the finite velocity energy condition \eqref{caseEC} or the infinite velocity energy condition \eqref{caseECX}. Then, for both the purely spatial derivative case \eqref{caseA} and the space-time derivative case \eqref{caseB}, for any $\mathfrak{n} \geq 0$ and $0 \leq |\al| \leq \mathrm{N}+1$ with $\mathrm{N} > d/2 + 1$, the following estimates hold uniformly for all $t \in [0,T]$, where $T > 0$ is the time defined by the bootstrap assumption in \eqref{condition}:
\begin{align}\label{L.Vdecomp3}
\bega
&\int_{\Omega} \bigg(|\nabla_xu^{\e}(t)| \mathcal{V}_{2}^{\al}(t) + |\nabla_x\ta^{\e}(t)| \mathcal{V}_{3}^{\al}(t)\bigg) dx \les \|\nabla_x (u^{\e},\ta^{\e})(t)\|_{L^\infty_x} \e^{\frac{3}{2}}\kappa \cr 
&\qquad \times \Big(\mathcal{D}_{tot}^{\mathrm{N}}(F^{\e}(t))+(\mathcal{E}_M^{\mathrm{N}}(F^{\e}(t)))^2\Big)  + e^{-c_1\e^{-{\frac{1}{3}}}}\frac{1}{\e} \mathcal{E}_{tot}^{\mathrm{N}}(F^{\e}(t)) \kappa^{(|\al|-\mathrm{N})_+} \|\p^{\al}\mathfrak{h}^{\e}(t)\|_{L^\infty_{x,v}},  \cr 
&\frac{1}{\eps} \int_{\Omega} \TbT(t) \mathcal{V}_{2}^{\al}(t) dx \les \e^{\frac{1}{2}}\kappa^{\frac{1}{2}}(\mathcal{D}_{G}^{\mathrm{N}}(F^{\e}(t)))^{\frac{1}{2}} \mathcal{E}_{tot}^{\mathrm{N}}(F^{\e}(t)) \cr 
&\qquad + \e\kappa^{\frac{1}{2}} e^{-c_1\e^{-{\frac{1}{2}}}}  (\mathcal{D}_{G}^{\mathrm{N}}(F^{\e}(t)))^{\frac{1}{2}}(\mathcal{E}_{tot}^{\mathrm{N}}(F^{\e}(t)))^{\frac{1}{2}} \kappa^{(|\al|-\mathrm{N})_+}\|\p^{\al}\mathfrak{h}^{\e}(t)\|_{L^\infty_{x,v}}.
\enda
\end{align}
\hide
\begin{align}
\bega
\mathcal{V}^{\al}_{\ell}(t) &\leq C\bigg[\e^{\frac{3}{2}}\kappa \Big(\mathcal{D}_{tot}^{\mathrm{N}}(F^{\e}(t))+\mathbf{1}_{|\al|\geq2}(\mathcal{E}_M^{\mathrm{N}}(F^{\e}(t)))^2\Big) + e^{-c_1\e^{-{\frac{1}{\ell}}}}\kappa^{(|\al|-\mathrm{N})_+}  \|\p^{\al}\mathfrak{h}^{\e}(t)\|_{L^\infty_{x,v}}^2\bigg], \cr 
\mathcal{V}^{\al}_{\ell}(t) &\leq C\bigg[\e^{-\frac{1}{2}}\mathcal{E}_{tot}^{\mathrm{N}}(F^{\e}(t)) + e^{-c_1\e^{-{\frac{1}{\ell}}}}\kappa^{(|\al|-\mathrm{N})_+}  \|\p^{\al}\mathfrak{h}^{\e}(t)\|_{L^\infty_{x,v}}^2\bigg],
\enda
\end{align}
\unhide
\end{proposition}
\begin{proof}
Since the argument is almost identical to that of Proposition \ref{P.hLinf}, we only provide a sketch of the proof. As in Proposition \ref{P.hLinf}, we split the velocity integration domain of $\mathcal{V}^{\al}_{\ell}$ into $\{ |v| \leq \e^{-\frac{1}{2\ell}} \}$ and $\{ |v| > \e^{-\frac{1}{2\ell}} \}$. 
The small velocity contribution $\mathcal{V}^{\al}_{\ell,s}$ satisfies the same estimate as in \eqref{smallv}. For the large velocity contribution, using the decomposition $F^{\e} = M^{\e} + \AC{\P}F^{\e} = M^{\e} + \e \sqrt{\bar{\mu}} g^{\e}$ implies that $\mathcal{V}^{\al}_{\ell,L,1}=0$ in \eqref{VlLdef}. Since $\mathcal{V}^{\al}_{\ell,L,2}$ satisfies the same estimate as in \eqref{Vl2a}, with $h^{\e}$ replaced by $\mathfrak{h}^{\e}$, the desired result follows.
\end{proof}

\begin{proposition}\label{P.hLinf4}
Assume the same hypotheses as in Proposition \ref{P.hLinf3}. In addition, suppose that 
\begin{align}\label{h2assume}
\bega
\e\sum_{0\leq |\al_x|\leq \mathrm{N}-3}\|\p^{\al_x}\mathfrak{h}^{\e}_0\|_{L^\infty_{x,v}} + \e \sup_{t \in [0,T]} (\mathcal{E}_{tot}^{\mathrm{N}}(F^{\e}(t)))^{\frac{1}{2}}  < \frac{1}{4C},
\enda
\end{align}
for $d=2,3$, where \(C\) is defined in \eqref{wgammav}, and \(C > 0\) is chosen in Proposition \ref{P.hLinf2}. Then the function \(\mathfrak{h}^{\e}\), defined in \eqref{hdef2}, satisfies
\begin{align}\label{claimx42}
\e\sum_{0\leq |\al_x|\leq \mathrm{N}-3}\|\p^{\al_x}\mathfrak{h}^{\e}(t)\|_{L^\infty_tL^\infty_{x,v}} \leq \frac{1}{2C},
\end{align}
and
\begin{align}
\sum_{0\leq|\al|\leq\mathrm{N}-2}\|\p^{\al}\mathfrak{h}^{\e}\|_{L^2_tL^\infty_{x,v}} &\les \e\kappa^{\frac{1}{2}}\sum_{0\leq|\al|\leq\mathrm{N}-2}\|\p^{\al}\mathfrak{h}^{\e}_0\|_{L^\infty_{x,v}} + \frac{1}{(\e\kappa)^{\frac{d}{2}}}\bigg(\int_0^t \mathcal{E}_{tot}^{\mathrm{N}}(F^{\e}(s)) ds\bigg)^{\frac{1}{2}}, \label{h2L2} \\ 
\sum_{0\leq|\al|\leq\mathrm{N}-2}\|\p^{\al}\mathfrak{h}^{\e}\|_{L^\infty_tL^\infty_{x,v}} &\les \sum_{0\leq|\al|\leq\mathrm{N}-2}\|\p^{\al}\mathfrak{h}^{\e}_0\|_{L^\infty_{x,v}} +\frac{1}{(\e\kappa)^{\frac{d}{2}}}\|\mathcal{E}_{tot}^{\mathrm{N}}(F^{\e}(t))\|_{L^\infty_t}^{\frac{1}{2}}. \label{h2Linf}
\end{align}
\end{proposition}

\begin{remark}\label{Rmk.h0}
We note that Proposition \ref{P.hLinf4} also holds in the case without derivatives.  
Suppose the same hypotheses as in Proposition \ref{P.hLinf3} hold, and in addition
\begin{align}\label{hassume0}
\bega
\e\|\mathfrak{h}^{\e}_0\|_{L^\infty_{x,v}} + \e \sup_{t \in [0,T]} (\mathcal{E}_{tot}^{\mathrm{N}}(F^{\e}(t)))^{\frac{1}{2}}  < \frac{1}{4C},
\enda
\end{align}
for $d=2,3$. Then the function \(\mathfrak{h}\) satisfies the following estimates:
\begin{align*}
\bega
\e\|\mathfrak{h}^{\e}\|_{L^\infty_tL^\infty_{x,v}} &\leq \frac{1}{2C}, \cr 
\|\mathfrak{h}^{\e}\|_{L^2_tL^\infty_{x,v}} &\les \e\kappa^{\frac{1}{2}}\|\mathfrak{h}^{\e}_0\|_{L^\infty_{x,v}} + \frac{1}{(\e\kappa)^{\frac{d}{2}}}\bigg(\int_0^t \mathcal{E}_{tot}^{\mathrm{N}}(F^{\e}(s)) ds\bigg)^{\frac{1}{2}}, \cr
\|\mathfrak{h}^{\e}\|_{L^\infty_tL^\infty_{x,v}} &\les \|\mathfrak{h}^{\e}_0\|_{L^\infty_{x,v}} +\frac{1}{(\e\kappa)^{\frac{d}{2}}}\|\mathcal{E}_{tot}^{\mathrm{N}}(F^{\e}(t))\|_{L^\infty_t}^{\frac{1}{2}}. 
\enda
\end{align*}
\end{remark}

\begin{remark}\label{Rmk.h0-2}
Note that we cannot estimate $\p^{\al}\mathfrak{h}^{\e}$ in \eqref{L.Vdecomp3} for high-order derivatives $|\al|\geq \mathrm{N}-1$. Moreover, the estimate in Proposition \ref{P.hLinf} is not applicable for $|\al|=0$ in the infinite velocity energy case \eqref{caseECX}. 
Therefore, when estimating $\mathcal{V}_{\ell}$, we proceed as follows:
\begin{itemize}
\item For the finite velocity energy case \eqref{caseEC}, we use $\eqref{L.Vdecomp}_2$ for $0 \leq |\al|\leq \mathrm{N}+1$.  
\item For the infinite velocity energy case \eqref{caseECX}, we use $\eqref{L.Vdecomp}_2$ for $1 \leq |\al|\leq \mathrm{N}+1$ and $\eqref{L.Vdecomp3}_2$ for $|\al|=0$.  
\end{itemize}
This yields
\begin{align}\label{L.Vdecompsum}
\bega
&\sum_{0\leq|\al|\leq\mathrm{N}+1} \int_{\Omega} \bigg(|\nabla_xu^{\e}(t)| \mathcal{V}_{2}^{\al}(t) + |\nabla_x\ta^{\e}(t)| \mathcal{V}_{3}^{\al}(t)\bigg) dx \cr 
&\qquad \les \|\nabla_x (u^{\e},\ta^{\e})(t)\|_{L^\infty_x}\Big(\mathcal{E}_{tot}^{\mathrm{N}}(F^{\e}(t)) +\e^{\frac{3}{2}}\kappa \mathcal{D}_{tot}^{\mathrm{N}}(F^{\e}(t))\Big)  + e^{-c_1\e^{-{\frac{1}{3}}}}\frac{1}{\e} \mathcal{E}_{tot}^{\mathrm{N}}(F^{\e}(t)) \mathrm{Y}(t),  \cr 
&\sum_{0\leq|\al|\leq\mathrm{N}+1}\frac{1}{\eps} \int_{\Omega} \TbT(t) \mathcal{V}_{2}^{\al}(t) dx \cr 
&\qquad \les \e^{\frac{1}{2}}\kappa^{\frac{1}{2}}(\mathcal{D}_{G}^{\mathrm{N}}(F^{\e}(t)))^{\frac{1}{2}} \mathcal{E}_{tot}^{\mathrm{N}}(F^{\e}(t)) + \e\kappa^{\frac{1}{2}} e^{-c_1\e^{-{\frac{1}{2}}}}  (\mathcal{D}_{G}^{\mathrm{N}}(F^{\e}(t)))^{\frac{1}{2}}(\mathcal{E}_{tot}^{\mathrm{N}}(F^{\e}(t)))^{\frac{1}{2}} \mathrm{Y}(t),
\enda
\end{align}
\hide
\begin{align}
\bega
\mathcal{V}_{\ell}(F^{\e}(t)) &\leq C\bigg[\mathcal{E}_{tot}^{\mathrm{N}}(t) +\e^{\frac{3}{2}}\kappa \Big(\mathcal{D}_{tot}^{\mathrm{N}}(F^{\e}(t))+(\mathcal{E}_M^{\mathrm{N}}(F^{\e}(t)))^2\Big) + e^{-c_1\e^{-{\frac{1}{\ell}}}}\mathrm{Y}(t)\bigg], \cr 
\mathcal{V}_{\ell}(F^{\e}(t)) &\leq C\bigg[\e^{-\frac{1}{2}}\mathcal{E}_{tot}^{\mathrm{N}}(F^{\e}(t)) + e^{-c_1\e^{-{\frac{1}{\ell}}}}\mathrm{Y}(t)\bigg],
\enda
\end{align}
\unhide
where the total $h^{\e}$--part $\mathrm{Y}(t)$ is defined by
\begin{align}\label{Thdef}
\bega
\mathrm{Y}(t) := \begin{cases}
\displaystyle \sum_{0\leq|\al|\leq\mathrm{N}+1}\kappa^{(|\al|-\mathrm{N})_+}  \|\p^{\al}h^{\e}(t)\|_{L^\infty_{x,v}}, & \mbox{for} \quad \eqref{caseEC}, \\
\displaystyle \|\mathfrak{h}^{\e}(t)\|_{L^\infty_{x,v}}+\sum_{1\leq|\al|\leq\mathrm{N}+1}\kappa^{(|\al|-\mathrm{N})_+}  \|\p^{\al}h^{\e}(t)\|_{L^\infty_{x,v}}, & \mbox{for} \quad \eqref{caseECX}.
\end{cases}
\enda
\end{align}

\end{remark}

In addition to the above proposition, we will use the following lemma to estimate the last two terms in \eqref{totalGt} and \eqref{totalEt}.
In the lemma below, the factor $\kappa^{-\frac{1}{8}(4-\mathrm{N})_+}$ reflects the following dichotomy: when $\mathrm{N}=3$, a singular scaling $\kappa^{-\frac{1}{8}}$ appears, whereas for $\mathrm{N}\geq 4$, no singular scaling occurs.

\begin{lemma}\label{L.Ginft}
Assume the same hypotheses as in Proposition \ref{P.hLinf2}. In addition, for the infinite velocity energy case \eqref{caseECX}, suppose the additional assumption \eqref{hassume0}. Then the following estimates hold:
\begin{align}\label{GL2}
\bega
&\sum_{0\leq|\al|\leq \lfloor(\mathrm{N}+1)/2\rfloor}\Big\|\|\la v \ra^{\frac{1}{2}} \p^{\al}\AC{\P}F^{\e} |M^{\e}|^{-\frac{1}{2}}(t)\|_{L^2_v}\Big\|_{L^\infty_x} \cr 
&\les \e^2(\mathcal{D}_{tot}^{\mathrm{N}}(F^{\e}(t)))^{\frac{1}{2}} + \e^2\kappa^{-\frac{1}{8}(4-\mathrm{N})_+}\Big(\mathcal{E}_{tot}^{\mathrm{N}}(F^{\e}(t)) + \e^{\frac{3}{2}}\kappa(\mathcal{D}_{tot}^{\mathrm{N}}(F^{\e}(t)))^{\frac{1}{2}}\Big)  \cr 
&+ \e^{\frac{3}{2}} (\mathcal{E}_{tot}^{\mathrm{N}}(F^{\e}(t)))^{\frac{1}{2}}e^{-\frac{c_1}{2}\e^{-{\frac{1}{9}}}}\bigg[ \mathrm{Y}^{\frac{1}{2}}(0) +\frac{1}{(\e\kappa)^{\frac{d}{2}}}\|\mathcal{E}_{tot}^{\mathrm{N}}(F^{\e}(t))\|_{L^\infty_t}^{\frac{1}{2}}\bigg],
\enda
\end{align}
and
\begin{align}\label{GLinfty}
\bega
&\sum_{0\leq|\al|\leq \lfloor(\mathrm{N}+1)/2\rfloor}\sup_{0\leq t\leq T}\Big\|\|\la v \ra^{\frac{1}{2}} \p^{\al}\AC{\P}F^{\e} |M^{\e}|^{-\frac{1}{2}}\|_{L^2_v}\Big\|_{L^\infty_x} \cr 
&\leq \e^{\frac{3}{4}}
 \sqrt{\mathcal{E}_{tot}^{\mathrm{N}}(F^{\e}(t))} 
+ \e^{\frac{1}{2}} e^{-\frac{c_1}{2}\e^{-{\frac{1}{9}}}}\bigg[\mathrm{Y}^{\frac{1}{2}}(0) +\frac{1}{(\e\kappa)^{\frac{d}{2}}}\|\mathcal{E}_{tot}^{\mathrm{N}} (F^{\e}(t))\|_{L^\infty_t}^{\frac{1}{2}}\bigg],
\enda
\end{align}
where $\mathrm{Y}(t)$ is defined in \eqref{Thdef}. 
\end{lemma}
\begin{proof}
The proof of this lemma is provided at the end of this section.
\end{proof}

\begin{remark}
(1) As shown in the estimate \eqref{totalGt} of Proposition~\ref{P.G.Energy}, 
we require the $L^2_t$ estimate \eqref{GL2} and the $L^\infty_t$ estimate \eqref{GLinfty} for the fourth and fifth lines of \eqref{totalGt}, respectively. 
Therefore, although the scaling in \eqref{GLinfty} is weaker, the estimate can still be closed.
(2) In the purely spatial derivative case~\eqref{caseA}, time embedding cannot be applied. However, if we adopt the scaled space--time derivatives $\p^{\al}  = (\e^{\mathfrak{n}}\p_t)^{\al_0}\p_{x_1}^{\al_1}\cdots\p_{x_d}^{\al_d}$, then time embedding can be used at the cost of losing the scale $\e^{\mathfrak{n}}$.
\hide 
following the same approach as in \eqref{Gembed}, we can prove
\begin{align*}
&\sup_{0 \leq t \leq T}\Big\| \|\la v \ra^{\frac{1}{2}} \p^{\al_x}\AC{\P}F^{\e}|M^{\e}|^{-1/2}(t) \|_{L^2_v} \Big\|_{L^\infty_x} \les \frac{1}{\sqrt{T}}\bigg( \e^2\kappa^{\frac{1}{2}}\|\mathcal{D}_G^{\frac{1}{2}}\|_{L^2_T}+\sum_{i=1}^2\e\kappa^{-\frac{1}{4}}\Big\|(\e^i\mathcal{E}_M^{\frac{i}{2}})\mathcal{V}_{4i+1}^{\frac{1}{2}}\Big\|_{L^2_T} \bigg) \cr 
&+ \begin{cases} \frac{\sqrt{T}}{\e^{\mathfrak{n}}}\Big\|\e^2\kappa^{\frac{1}{2}}\mathcal{D}_G^{\frac{1}{2}}+\sum_{i=1}^3\e(\e^i\mathcal{E}_M^{\frac{i}{2}})\mathcal{V}_{4i+1}^{\frac{1}{2}} \Big\|_{L^2_T},  &\mbox{when} \quad 0\leq |\al|\leq \mathrm{N}-3, \\ 
\frac{\sqrt{T}}{\e^{\mathfrak{n}}}\Big\|\e^2\mathcal{D}_{top}^{\frac{1}{2}} +\e^2\kappa^{\frac{1}{2}}\mathcal{D}_G^{\frac{1}{2}} + \e^2\mathcal{E}_M+\sum_{i=1}^3\e\kappa^{-\frac{1}{4}}(\e^i\mathcal{E}_M^{\frac{i}{2}})\mathcal{V}_{4i+1}^{\frac{1}{2}} \Big\|_{L^2_T},  &\mbox{when} \quad |\al|= \mathrm{N}-2. \end{cases}
\end{align*}
For simplicity, however, we will only use the estimate \eqref{GLinfty}.
\unhide
\end{remark}

We will prove Proposition \ref{P.hLinf4} at the end of this subsection, after establishing several auxiliary lemmas.
In the proof, for brevity, we slightly abuse notation by writing 
$\mathfrak{h}$, $g$, $\mathcal{E}(t)$, and $\mathcal{D}(t)$ for 
$\mathfrak{h}^{\e}$, $g^{\e}$, $\mathcal{E}^{\mathrm{N}}(F^{\e}(t))$, 
and $\mathcal{D}^{\mathrm{N}}(F^{\e}(t))$, respectively.

\begin{lemma}\label{L.htt2}
Let \( \mathfrak{h}^{\e}(t,x,v) \) be a solution of \eqref{wBE2}.  
For both types of derivatives \( \p^{\al} \) defined in \eqref{caseA} and \eqref{caseB}, for any $\mathfrak{n} \geq 0$, and for \( (x,v) \in \R^d \times \R^3 \) with \( d = 2,3 \), the following estimates hold:
\begin{align}
\|\p^{\al}\mathfrak{h}^{\e}\|_{L^2_tL^\infty_{x,v}} &\leq C\e\kappa^{\frac{1}{2}}\|\p^{\al}\mathfrak{h}^{\e}_0\|_{L^\infty_{x,v}} + \frac{C_\mathfrak{N}}{(\e\kappa)^{\frac{d}{p}}} \|\mathbf{1}_{|v|\leq2\mathfrak{N}}\p^{\al}g^{\e}\|_{L^2_tL^p_xL^2_v} \label{h2L2t} \\
& + C \e \sum_{0\leq\beta\leq\al} \big\|\|\p^{\beta}\mathfrak{h}^{\e}(t)\|_{L^\infty_{x,v}} \|\p^{\al-\beta}\mathfrak{h}^{\e}(t)\|_{L^\infty_{x,v}}\big\|_{L^2_t} + \e\kappa \Big\|\|\mathrm{M}_{\text{trans}}^{\al}\|_{L^\infty_{x,v}}\Big\|_{L^2_t},  \cr
\|\p^{\al}\mathfrak{h}^{\e}\|_{L^\infty_tL^\infty_{x,v}} &\leq C\|\p^{\al}\mathfrak{h}^{\e}_0\|_{L^\infty_{x,v}} + \frac{C_\mathfrak{N}}{(\e\kappa)^{\frac{d}{p}}} \|\mathbf{1}_{|v|\leq2\mathfrak{N}}\p^{\al}g^{\e}\|_{L^\infty_tL^p_xL^2_v} \label{h2Linft} \\
&+ C \e \sum_{0\leq\beta\leq\al} \|\p^{\beta}\mathfrak{h}^{\e}\|_{L^\infty_tL^\infty_{x,v}} \|\p^{\al-\beta}\mathfrak{h}^{\e}\|_{L^\infty_tL^\infty_{x,v}} + \e\kappa \Big\|\|\mathrm{M}_{\text{trans}}^{\al}\|_{L^\infty_{x,v}}\Big\|_{L^\infty_t}, \notag 
\end{align}
for any \( p \in [1, \infty] \), where \( C > 0 \) is a constant, \( C \) is the constant from \eqref{wgammav}, and \( C_\mathfrak{N} > 0 \) depends on \( \mathfrak{N} \).
\end{lemma}
\begin{proof}
In the same manner as in Lemma \ref{L.htt}, we apply Duhamel's principle twice to \eqref{wBE2}, which yields
\begin{align*}
\p^{\al}\mathfrak{h}(t,x,v) &= I^{\mathfrak{h}}_1 + I^{\mathfrak{h}}_2 + I^{\mathfrak{h}}_3 + I^{\mathfrak{h}}_4 + I^{\mathfrak{h}}_5 + I^{\mathfrak{h}}_6 + I^{\mathfrak{h}}_7,
\end{align*}
where
\begin{align*}
I^{\mathfrak{h}}_1(t,x,v)&:= e^{ \frac{-\nu(v)t}{\e^2 \kappa}}\p^{\al}\mathfrak{h}_0\bigg(x-\frac{v}{\e}t,v\bigg) ,\cr 
I^{\mathfrak{h}}_2(t,x,v)&:= \frac{1}{\e^2 \kappa}\int_0^t e^{ \frac{-\nu(v)(t-s)}{\e^2 \kappa}} \int_{\R^3}\bar{\mathbf{k}}_{w}(v,v_*) e^{ \frac{-\nu(v_*)s}{\e^2 \kappa}}\p^{\al}\mathfrak{h}_0\bigg(x-\frac{v}{\e}(t-s)-\frac{v_*}{\e}s,v_*\bigg) dv_* ds ,\cr 
I^{\mathfrak{h}}_3(t,x,v)&:=\frac{1}{\e^2 \kappa}\int_0^t e^{ \frac{-\nu(v)(t-s)}{\e^2 \kappa}} \int_{\R^3}\bar{\mathbf{k}}_{w}(v,v_*) \frac{1}{\e^2 \kappa} \int_0^s e^{\frac{-\nu(v_*)(s-\tau)}{\e^2 \kappa}}  \cr 
&\quad \int_{\R^3}\bar{\mathbf{k}}_{w}(v_*,v_{**}) \p^{\al}\mathfrak{h}\bigg(\tau, x-\frac{v}{\e}(t-s)-\frac{v_*}{\e}(s-\tau),v_{**}\bigg)dv_{**} d\tau dv_* ds ,\cr 
I^{\mathfrak{h}}_4(t,x,v)&:= \frac{1}{\e^2 \kappa}\int_0^t e^{ \frac{-\nu(v)(t-s)}{\e^2 \kappa}} \int_{\R^3}\bar{\mathbf{k}}_{w}(v,v_*) \frac{1}{\e \kappa} \sum_{0\leq\beta\leq\al}\binom{\al}{\beta} \int_0^s e^{\frac{-\nu(v_*)(s-\tau)}{\e^2 \kappa}} \cr 
&\quad \times  w\bar{\Gamma}\bigg(\frac{\p^{\beta}\mathfrak{h}}{w},\frac{\p^{\al-\beta}\mathfrak{h}}{w}\bigg)\bigg(\tau, x-\frac{v}{\e}(t-s)-\frac{v_*}{\e}(s-\tau),v_*\bigg)d\tau dv_* ds ,\cr 
I^{\mathfrak{h}}_5(t,x,v)&:= \frac{1}{\e \kappa}\sum_{0\leq\beta\leq\al}\binom{\al}{\beta}\int_0^t e^{ \frac{-\nu(v)(t-s)}{\e^2 \kappa}} w\bar{\Gamma}\bigg(\frac{\p^{\beta}\mathfrak{h}}{w},\frac{\p^{\al-\beta}\mathfrak{h}}{w}\bigg)\bigg(s, x-\frac{v}{\e}(t-s),v\bigg)ds, \cr 
I^{\mathfrak{h}}_{6}(t,x,v)&:= -\int_0^t e^{ \frac{-\nu(v)(t-s)}{\e^2 \kappa}}\frac{1}{\e}\mathrm{M}_{\text{trans}}^{\al}\bigg(s, x-\frac{v}{\e}(t-s),v_*\bigg) ds ,\cr 
I^{\mathfrak{h}}_{7}(t,x,v)&:= -\frac{1}{\e^2 \kappa}\int_0^t e^{ \frac{-\nu(v)(t-s)}{\e^2 \kappa}} \int_{\R^3}\bar{\mathbf{k}}_{w}(v,v_*)  \int_0^s e^{\frac{-\nu(v_*)(s-\tau)}{\e^2 \kappa}}\cr 
&\quad \times  \frac{1}{\e}\mathrm{M}_{\text{trans}}^{\al}\bigg(\tau, x-\frac{v}{\e}(t-s)-\frac{v_*}{\e}(s-\tau),v_{**}\bigg)dv_{**} d\tau dv_* ds .
\end{align*}
The estimates for $I^{\mathfrak{h}}_1, \ldots, I^{\mathfrak{h}}_5$ are the same as in Lemma \ref{L.htt}, since $\bar{\mathbf{k}}_{w}$ and $\bar{\Gamma}$ satisfy the similar properties as $\mathbf{k}_{w}$ and $\Gamma$ in \eqref{kw2} and \eqref{wgammav2}. 
Therefore, it remains to estimate $I^{\mathfrak{h}}_{6}$ and $I^{\mathfrak{h}}_{7}$. \\
Applying Young’s convolution inequality to \( I^{\mathfrak{h}}_{6} \), we obtain
\begin{align*}
\bega
\|I^{\mathfrak{h}}_6\|_{L^2_t} 
&\les \Big\|e^{ \frac{-\nu(v)t}{\e^2 \kappa}}\Big\|_{L^1_t} \bigg\|\Big\|\frac{1}{\e}\mathrm{M}_{\text{trans}}^{\al}\Big\|_{L^\infty_{x,v}} \bigg\|_{L^2_t} \les \e\kappa \bigg\| \Big\|\mathrm{M}_{\text{trans}}^{\al}\Big\|_{L^\infty_{x,v}} \bigg\|_{L^2_t}.
\enda
\end{align*}
To estimate \( I^{\mathfrak{h}}_7 \), we use Young’s inequality, Minkowski’s integral inequality, the bound \eqref{wgammav2}, and the kernel property $\eqref{kw2}_2$:
\begin{align*}
\bega
\|I^{\mathfrak{h}}_7\|_{L^2_t}&\les \bigg\|\frac{1}{\e^2 \kappa}e^{\frac{-\nu(v)t}{\e^2 \kappa}}\bigg\|_{L^1_t} \bigg\|\int_{\R^3}\bar{\mathbf{k}}_{w}(v,v_*)  \int_0^s e^{\frac{-\nu(v_*)(s-\tau)}{\e^2 \kappa}}\frac{1}{\e}\|\mathrm{M}_{\text{trans}}^{\al}(\tau)\|_{L^\infty_{x,v}}d\tau dv_* \bigg\|_{L^2_s} \cr 
&\les \bigg\|\frac{1}{\e^2\kappa}e^{\frac{-\nu(v)t}{\e^2 \kappa}}\bigg\|_{L^1_t} \int_{\R^3}|\bar{\mathbf{k}}_{w}(v,v_*)|  \e^2\kappa\bigg\|\frac{1}{\e^2\kappa}e^{ \frac{-\nu(v_*)s}{\e^2 \kappa}}\bigg\|_{L^1_s} \frac{1}{\e}\Big\|\|\mathrm{M}_{\text{trans}}^{\al}(\tau)\|_{L^\infty_{x,v}}\Big\|_{L^2_\tau} dv_* \cr 
&\les \e\kappa \Big\|\|\mathrm{M}_{\text{trans}}^{\al}(\tau)\|_{L^\infty_{x,v}}\Big\|_{L^2_\tau}.
\enda
\end{align*}
This completes the proof of the desired estimates.
\end{proof}

\begin{lemma}\label{L.ftoE2}
For $t\in[0,T]$ satisfying the bootstrap assumption \eqref{condition}, for both types of derivatives \( \p^{\al} \) defined in \eqref{caseA} and \eqref{caseB}, and for the constant \( \mathfrak{N} > 0 \) chosen in Lemma \ref{L.htt2}, we have 
\begin{align}\label{ftoE3}
\begin{split}
\kappa^{(|\al|-\mathrm{N})_+}\big\|\mathbf{1}_{|v|\leq2\mathfrak{N}}\p^{\al}g^{\e}(t)\big\|_{L^2_{x,v}}^2 &\les C_{\mathfrak{N}} \mathcal{E}_{tot}^{\mathrm{N}}(F^{\e}(t)),
\end{split}
\end{align}
for $0 \leq |\al| \leq \mathrm{N}+1$, in both the finite velocity energy case \eqref{caseEC} 
and the infinite velocity energy case \eqref{caseECX}.
\end{lemma}
\begin{proof}
Different from Lemma \ref{L.ftoE}, using the decomposition \( F = M^{\e} + \AC{\P}F^{\e} = M^{\e} + \e \sqrt{\mu} g^{\e} \), we have only the microscopic contribution:
\begin{align*}
\begin{split}
\kappa^{(|\al|-\mathrm{N})_+}\big\|\mathbf{1}_{|v|\leq2\mathfrak{N}}\p^{\al}g^{\e}\big\|_{L^2_{x,v}}^2 &= \kappa^{(|\al|-\mathrm{N})_+}\lw\|\mathbf{1}_{|v|\leq2\mathfrak{N}}\frac{\p^{\al}\AC{\P}F^{\e}}{\e\sqrt{\mu}}\rw\|_{L^2_{x,v}}^2.
\end{split}
\end{align*}
Proceeding as in \eqref{Gv<N}, and using the definitions of the energy $\mathcal{E}_{G}$ and $\mathcal{E}_{top}$ given in \eqref{N-EDdef} and \eqref{N-EDdef2}, we directly obtain the desired estimate.
\end{proof}

\begin{lemma}\label{L.TMbdd} 
For $t\in[0,T]$ satisfying the bootstrap assumption \eqref{condition}, and for both types of derivatives \( \p^{\al} \) defined in \eqref{caseA} and \eqref{caseB}, and for \( |\al| \leq \mathrm{N}-2 \), the following estimate holds:
\begin{align*}
\bega
\|\mathrm{M}_{\text{trans}}^{\al}(t)\|_{L^\infty_{x,v}} \les \begin{cases} (\mathcal{E}_{tot}^{\mathrm{N}}(F^{\e}(t)))^{\frac{1}{2}},  &\mbox{for} \quad 0\leq |\al|\leq \mathrm{N}-3, \\
\kappa^{-\frac{1}{4}}(\mathcal{E}_{tot}^{\mathrm{N}}(F^{\e}(t)))^{\frac{1}{2}},  &\mbox{for} \quad |\al|=\mathrm{N}-2,
\end{cases}
\enda
\end{align*}
where \( \mathrm{M}_{\text{trans}}^{\al}(t) \) is defined in \eqref{TMdef}.
\end{lemma}
\begin{proof}
We first claim the following estimate:
\begin{align}\label{TMbd1}
\bega
\|\mathrm{M}_{\text{trans}}^{\al}&(t)\|_{L^\infty_{x,v}} \les \sum_{0\leq \beta \leq \al}\bigg[\|\p^{\beta}\nabla_xu^{\e}(t)\|_{L^\infty_x} + \|\p^{\beta}\nabla_x\ta^{\e}(t)\|_{L^\infty_x} \cr 
&+ \sum_{i,j}\|\p^{\beta} \p_{x_j} \mathbf{r}_{ij}^{\e}(t)\|_{L^\infty_x} + \sum_{j}\| \p^{\beta}\p_{x_j} \mathfrak{q}_j^{\e}(t)\|_{L^\infty_x} + \sum_{i,j} \|\p^{\beta}\big(\p_{x_i}\mathrm{U}^{\e}_j \mathbf{r}_{ij}^{\e}\big)(t)\|_{L^\infty_x} \bigg].
\enda
\end{align}
Recall that we have an explicit expression for $\p_t M^{\e} + \frac{v}{\e}\cdot\nabla_x M^{\e}$ in \eqref{Phitx}. 
Taking $\frac{w}{\sqrt{\bar{\mu}}}\p^{\al}$ and then the 
\(L^\infty_{x,v}\)-norm of the first term of \eqref{Phitx} yields
\begin{align*}
\bigg\|\frac{w}{\sqrt{\bar{\mu}}}\p^{\al}\bigg(\frac{1}{k_B\mathrm{\Theta}^{\e}}\sum_{i,j}\p_{x_i}u^{\e}_j \mathfrak{R}^{\e}_{ij}M^{\e}\bigg)\bigg\|_{L^\infty_{x,v}} &\les \sum_{0\leq\beta\leq\al}\binom{\al}{\beta}\bigg\|\frac{w}{\sqrt{\bar{\mu}}}\p^{\al-\beta}\bigg(\frac{1}{k_B\mathrm{\Theta}^{\e}}\sum_{i,j} \mathfrak{R}^{\e}_{ij}M^{\e}\bigg) \p^{\beta}\p_{x_i}u^{\e}_j\bigg\|_{L^\infty_{x,v}} \cr 
&\les \sum_{0\leq\beta\leq\al} (1+\e\mathcal{E}_M^{\frac{1}{2}}+\cdots+\e^{|\al-\beta|}\mathcal{E}_M^{\frac{|\al-\beta|}{2}})\|\p^{\beta}\p_{x_i}u^{\e}_j\|_{L^\infty_x} \cr 
&\les \sum_{0\leq\beta\leq\al}\|\p^{\beta}\p_{x_i}u^{\e}_j\|_{L^\infty_x},
\end{align*}
where we have used \(\e \mathcal{E}_M^{1/2} \leq 1\) in \((\mathcal{B}_2)\) from \eqref{condition} and the estimate \eqref{al-vM}.  
By applying the same argument to the remaining terms of \eqref{Phitx}, we obtain the claim \eqref{TMbd1}. \\
For the estimate of the first line of \eqref{TMbd1}, we use \eqref{rutinf} to have 
\begin{align}\label{utaE}
\bega
\sum_{0\leq \beta \leq \al}\bigg[\|\p^{\beta}\nabla_xu^{\e}(t)\|_{L^\infty_x} + \|\p^{\beta}\nabla_x\ta^{\e}(t)\|_{L^\infty_x} \bigg] \les \begin{cases} \mathcal{E}_{M}^{\frac{1}{2}},  &\mbox{for} \quad 0\leq |\al|\leq \mathrm{N}-3, \\ 
\kappa^{-\frac{1}{4}}\mathcal{E}_{tot}^{\frac{1}{2}},  &\mbox{for} \quad |\al|=\mathrm{N}-2.
\end{cases}
\enda
\end{align}
For the second line of \eqref{TMbd1}, unlike Lemma \ref{L.ABG} where we used dissipation estimates, 
here we use the energy bound since we need an \(L^\infty_t\) estimate. 
Recall the inequality in \eqref{albe.1}.  
By the same argument as in \eqref{albe.claim}, we obtain the following estimate:
\begin{align*}
\bega
\|\p^{\al}\AC{\P}F^{\e}|M^{\e}|^{-\frac{1}{2}}\|_{L^2_{x,v}} \les \begin{cases} \e\mathcal{E}_G^{\frac{1}{2}}, \quad &\mbox{when} \quad 0\leq|\al|\leq\mathrm{N}, \\ 
\e\kappa^{-\frac{1}{2}}\mathcal{E}_{top}^{\frac{1}{2}}(t) +\e^2\mathbf{1}_{|\al|\geq2}\mathcal{E}_M(t), \quad &\mbox{when} \quad |\al|=\mathrm{N}+1.
\end{cases}
\enda
\end{align*}
Consequently, 
\begin{align*}
\bega
\|\p^{\al} \mathbf{r}_{ij}^{\e}\|_{L^2_x} , \|\p^{\al}\mathfrak{q}_i^{\e}\|_{L^2_x} &\les \begin{cases} \e\mathcal{E}_G^{\frac{1}{2}}, \quad &\mbox{when} \quad 0\leq|\al|\leq\mathrm{N}, \\ 
\e\kappa^{-\frac{1}{2}}\mathcal{E}_{top}^{\frac{1}{2}}(t) +\e^2\mathbf{1}_{|\al|\geq2}\mathcal{E}_M(t), \quad &\mbox{when} \quad |\al|=\mathrm{N}+1.
\end{cases}
\enda
\end{align*}
Applying Agmon’s inequality \eqref{Agmon} to the second line of \eqref{TMbd1} and using \eqref{albe.E}, 
we conclude, in the same way as in the proof of \eqref{rutinf}, that
\begin{align}\label{albe.E}
\bega
\sum_{0\leq \beta \leq \al}&\bigg[\sum_{i,j}\|\p^{\beta} \p_{x_j} \mathbf{r}_{ij}^{\e}(t)\|_{L^\infty_x} + \sum_{j}\| \p^{\beta}\p_{x_j} \mathfrak{q}_j^{\e}(t)\|_{L^\infty_x} + \sum_{i,j} \|\p^{\beta}\big(\p_{x_i}\mathrm{U}^{\e}_j \mathbf{r}_{ij}^{\e}\big)(t)\|_{L^\infty_x} \bigg] \cr 
&\les \begin{cases} \e\mathcal{E}_{G}^{\frac{1}{2}},  &\mbox{for} \quad 0\leq |\al|\leq \mathrm{N}-3, \\ 
\e\kappa^{-\frac{1}{4}}\mathcal{E}_{tot}^{\frac{1}{2}},  &\mbox{for} \quad |\al|=\mathrm{N}-2,
\end{cases}
\enda
\end{align}
where we have again used \(\e\mathcal{E}_M^{1/2}\leq 1\) in $(\mathcal{B}_2)$ from \eqref{condition}.  
Finally, combining \eqref{utaE} and \eqref{albe.E} into \eqref{TMbd1}, we obtain the desired result.

\end{proof}

\begin{proof}[\textbf{Proof of Proposition \ref{P.hLinf4}}]
The proof is essentially the same as that of Proposition \ref{P.hLinf2}. 
Here we only establish \eqref{h2L2} and \eqref{claimx42}.  
We apply $\sum_{0\leq|\al|\leq\mathrm{N}-2}$ to \eqref{h2L2t} in Lemma \ref{L.htt2} with $p=2$, having
\begin{align*}
\sum_{0\leq|\al|\leq\mathrm{N}-2}\|\p^{\al}\mathfrak{h}\|_{L^2_tL^\infty_{x,v}}&\leq C \e\kappa^{\frac{1}{2}}\sum_{0\leq|\al|\leq\mathrm{N}-2}\|\p^{\al}\mathfrak{h}_0\|_{L^\infty_{x,v}} + I_{\mathfrak{h}} + II_{\mathfrak{h}},
\end{align*}
where
\begin{align*}
I_{\mathfrak{h}}&= \frac{C}{(\e\kappa)^{\frac{d}{2}}}\sum_{0\leq|\al|\leq\mathrm{N}-2} \|\mathbf{1}_{|v|\leq2\mathfrak{N}}\p^{\al}g\|_{L^2_tL^2_{x,v}} + \e\kappa \sum_{0\leq|\al|\leq\mathrm{N}-2}\Big\|\|\mathrm{M}_{\text{trans}}^{\al}\|_{L^\infty_{x,v}}\Big\|_{L^2_t}, \cr
II_{\mathfrak{h}}&= C \e \sum_{0\leq|\al|\leq\mathrm{N}-2}\sum_{\beta+\gamma=\al} \Big\|\|\p^{\beta}\mathfrak{h}(t)\|_{L^\infty_{x,v}} \|\p^{\gamma}\mathfrak{h}(t)\|_{L^\infty_{x,v}} \Big\|_{L^2_t}.
\end{align*}
We apply \eqref{ftoE3} together with Lemma \ref{L.TMbdd} to the term $I_{\mathfrak{h}}$ and obtain 
\begin{align*}
I_{\mathfrak{h}}&\les  \frac{C}{(\e\kappa)^{\frac{d}{2}}}\bigg(\int_0^t \mathcal{E}_{tot}(s)ds\bigg)^{\frac{1}{2}} + \e\kappa^{\frac{3}{4}}\bigg(\int_0^t \mathcal{E}_{tot}(s)ds\bigg)^{\frac{1}{2}} \les \frac{C}{(\e\kappa)^{\frac{d}{2}}}\bigg(\int_0^t \mathcal{E}_{tot}(s)ds\bigg)^{\frac{1}{2}}.
\end{align*}
Then, in the same manner as in the proof of Proposition \ref{P.hLinf2}, we absorb the nonlinear term into the left-hand side.
If the inequality \eqref{claimx42} holds, then we get 
\begin{align*}
\frac{1}{2}\sum_{0\leq|\al|\leq\mathrm{N}-2}\|\p^{\al}\mathfrak{h}\|_{L^2_tL^\infty_{x,v}} 
&\les \e\kappa^{\frac{1}{2}}\sum_{0\leq|\al|\leq\mathrm{N}-2}\|\p^{\al}\mathfrak{h}_0\|_{L^\infty_{x,v}} +\frac{C}{(\e\kappa)^{\frac{d}{2}}}\bigg(\int_0^t \mathcal{E}_{tot}(s)ds\bigg)^{\frac{1}{2}}.
\end{align*}
This directly yields the desired estimate \eqref{hL2}. \\
We now prove \eqref{claimx42}. 
Applying \( \sum_{0\leq |\al_x|\leq \mathrm{N}-3} \) to \eqref{h2Linft} in Lemma \ref{L.htt}, we get
\begin{align}\label{10.20-2}
\bega
&\sum_{0\leq |\al_x|\leq \mathrm{N}-3}\|\p^{\al_x}\mathfrak{h}\|_{L^\infty_tL^\infty_{x,v}}\les \sum_{0\leq |\al_x|\leq \mathrm{N}-3}\|\p^{\al_x}\mathfrak{h}_0\|_{L^\infty_{x,v}}+\frac{1}{(\e\kappa)^{\frac{d}{p}}}\sum_{0\leq |\al_x|\leq \mathrm{N}-3} \|\mathbf{1}_{|v|\leq2\mathfrak{N}}\p^{\al_x}g\|_{L^\infty_tL^p_xL^2_v} \cr 
&+ \e  \sum_{0\leq |\al_x|\leq \mathrm{N}-3}\sum_{0\leq\beta_x\leq\al_x} \|\p^{\beta_x}\mathfrak{h}\|_{L^\infty_tL^\infty_{x,v}}\|\p^{\al_x-\beta_x}\mathfrak{h}\|_{L^\infty_tL^\infty_{x,v}} + \e\kappa \sum_{0\leq |\al_x|\leq \mathrm{N}-3}\Big\|\|\mathrm{M}_{\text{trans}}^{\al}\|_{L^\infty_{x,v}}\Big\|_{L^\infty_t},
\enda
\end{align}
for any $p \in [1,\infty]$.
For the linear term involving $\p^{\al_x} g$, we choose $p=\infty$ and apply Agmon's inequality~\eqref{Agmon} together with~\eqref{ftoE3}, which yields
\begin{align*}
\sum_{0\leq |\al_x|\leq \mathrm{N}-3} \|\mathbf{1}_{|v|\leq2\mathfrak{N}}\p^{\al_x}g\|_{L^\infty_tL^\infty_xL^2_v} \les \|\mathcal{E}_{tot}\|_{L^\infty_t}^{\frac{1}{2}}.
\end{align*}
We also employ Lemma \ref{L.TMbdd} to estimate the term $\mathrm{M}_{\text{trans}}^{\al}$.  
Multiplying \eqref{10.20-2} by \( \e \) and then moving the nonlinear contribution to the left-hand side, we have
\begin{align*}
\bigg(1-&C\e\sum_{0\leq |\al_x|\leq \mathrm{N}-3}\|\p^{\al_x}\mathfrak{h}\|_{L^\infty_tL^\infty_{x,v}}\bigg)\bigg(\e\sum_{0\leq |\al_x|\leq \mathrm{N}-3}\|\p^{\al_x}\mathfrak{h}\|_{L^\infty_tL^\infty_{x,v}}\bigg) \cr 
&\les \e\sum_{0\leq |\al_x|\leq \mathrm{N}-3}\|\p^{\al_x}\mathfrak{h}_0\|_{L^\infty_{x,v}}+C\e\|\mathcal{E}_{tot}\|_{L^\infty_t}^{\frac{1}{2}},
\end{align*}
for $d=2,3$. Under the assumption \eqref{h2assume}, this reduces to the quadratic inequality \( (1 - CZ)Z < \frac{1}{4C} \) for \( Z := \e\sum_{0\leq |\al_x|\leq \mathrm{N}-3}\|\p^{\al_x}\mathfrak{h}\|_{L^\infty_tL^\infty_{x,v}} \). This proves the desired estimate \eqref{claimx42}.
\end{proof}

\begin{proof}[\textbf{Proof of Lemma \ref{L.Ginft}}] 
We first claim the following inequality:
\begin{align}\label{Gembed}
\bega
\Big\|& \|\la v \ra^{\frac{1}{2}} \p^{\al}\AC{\P}F^{\e}|M^{\e}|^{-1/2} (t)\|_{L^2_v} \Big\|_{L^\infty_x} \cr 
&\les \begin{cases} \e^2\kappa^{\frac{1}{2}}(\mathcal{D}_G^{\mathrm{N}}(F^{\e}(t)))^{\frac{1}{2}}+\mathfrak{A}(t) , \quad &\mbox{when} \quad 0\leq |\al|\leq \mathrm{N}-2, \\ 
\e^2(\mathcal{D}_{tot}^{\mathrm{N}}(F^{\e}(t)))^{\frac{1}{2}} + \e^2\mathcal{E}_M^{\mathrm{N}}(F^{\e}(t)) +\mathfrak{A}(t), \quad &\mbox{when} \quad |\al|=\mathrm{N}-1,
\end{cases}
\enda
\end{align}
where
\begin{align*}
\bega
\mathfrak{A}(t):= \e^2\Big(\int_{\Omega} \mathcal{V}_5 \big(|\nabla_x(\rho^{\e},u^{\e},\ta^{\e})|+|\nabla_x^2(\rho^{\e},u^{\e},\ta^{\e})|\big)dx\Big)^{\frac{1}{2}} + \e^3 \Big(\int_{\Omega} \mathcal{V}_9 |\nabla_x(\rho^{\e},u^{\e},\ta^{\e})|^2dx\Big)^{\frac{1}{2}}.
\enda
\end{align*}

We control the term via Agmon's inequality \eqref{Agmon}, and estimate the growth of the velocity weight arising from $\p |M^{\e}|^{-1/2}$.  
We first apply Minkowski’s integral inequality 
, and then apply Agmon's inequality \eqref{Agmon}. Distributing the two spatial derivatives to $\p^{\al}\AC{\P}F^{\e}$ and $|M^{\e}|^{-1/2}$, we obtain
\begin{align*}
\Big\| \| \la v \ra^{\frac{1}{2}} \p^{\al}\AC{\P}F^{\e}|M^{\e}|^{-1/2}(t)\|_{L^2_v} \Big\|_{L^\infty_x} &\les \Big\| \| \la v \ra^{\frac{1}{2}} \p^{\al}\AC{\P}F^{\e}|M^{\e}|^{-1/2}(t) \|_{H^2_x} \Big\|_{L^2_v} = I_{p}^{\al}(t)+ II_{p}^{\al}(t),
\end{align*}
where
\begin{align}\label{P12def}
\bega
I_{p}^{\al}(t) &:= \sum_{0\leq|\beta_x|\leq2} \left\| \la v \ra^{\frac{1}{2}}\p^{\al+\beta_x}\AC{\P}F^{\e} (|M^{\e}|^{-1/2})(t)\right\|_{L^2_{x,v}}, \cr 
II_{p}^{\al}(t) &:=  \sum_{1\leq|\beta_x|\leq2}\sum_{\substack{\beta_x^1+\beta_x^2=\beta_x\\ \beta_x^2>0}} \left\| \la v \ra^{\frac{1}{2}}\p^{\al+\beta_x^1}\AC{\P}F^{\e} \p^{\beta_x^2}(|M^{\e}|^{-1/2})(t)\right\|_{L^2_{x,v}}.
\enda
\end{align}
To estimate $I_{p}^{\al}$, we note that different dissipation norms are used depending on the number of derivatives, specifically whether $|\al| = \mathrm{N}$ or $|\al| = \mathrm{N}+1$. 
When $|\al| \leq \mathrm{N} - 2$, we directly apply the definition of $\mathcal{D}_G$ in \eqref{N-EDdef} to obtain
\begin{align*}
I_{p}^{\al} &\les \e^2\kappa^{\frac{1}{2}}\mathcal{D}_G^{\frac{1}{2}}, \quad \mbox{when} \quad 0\leq |\al|\leq \mathrm{N}-2.
\end{align*}
When $|\al| = \mathrm{N} - 1$, we use the commutator between $\p^{\al}$ and $\AC{\P}$ given in \eqref{pACP}, which yields
\begin{align*}
I_{p}^{\al} &\les \sum_{0\leq|\beta_x|\leq2} \left\| \la v \ra^{\frac{1}{2}}\AC{\P}(\p^{\al+\beta_x}F) (|M^{\e}|^{-1/2})\right\|_{L^2_{x,v}} \\
& \ + \sum_{0\leq|\beta_x|\leq2} \bigg\| \la v \ra^{\frac{1}{2}}\AC{\P} \bigg(\sum_{2\leq i\leq |\al|}\eps^i \Phi_{\al}^i M^{\e} \bigg) (|M^{\e}|^{-1/2})\bigg\|_{L^2_{x,v}} \cr 
&\les \e^2\mathcal{D}_{top}^{\frac{1}{2}} +\e^2\kappa^{\frac{1}{2}}\mathcal{D}_G^{\frac{1}{2}} + \e^2\mathcal{E}_M, \quad \mbox{when} \quad |\al|=\mathrm{N}-1,
\end{align*}
where we used \eqref{N-EDdef2} to estimate the top-order dissipation terms involving $\mathrm{N} + 1$ derivatives, and applied \eqref{Rscale} to control the remainder term involving $\sum_{2\leq i\leq |\al|}\eps^i \Phi_{\al}^i M^{\e}$.
To estimate $II_{p}^{\al}$, we focus on the velocity growth of the term $\p^{\beta_x^2}(|M^{\e}|^{-1/2})$. Depending on the number of derivatives acting on $|M^{\e}|^{-1/2}$, applying \eqref{al-1/M} gives 
\begin{align*}
\bega
II_{p}^{\al}(t) &\les  \e \sum_{0\leq \beta_x^1 \leq1}\left\| \la v \ra^{2+\frac{1}{2}}\p^{\al+\beta_x^1}\AC{\P}F^{\e} \big(|\nabla_x(\rho^{\e},u^{\e},\ta^{\e})|+|\nabla_x^2(\rho^{\e},u^{\e},\ta^{\e})|\big)|M^{\e}|^{-1/2}(t)\right\|_{L^2_{x,v}} \cr 
&+ \e^2\left\| \la v \ra^{4+\frac{1}{2}}\p^{\al}\AC{\P}F^{\e} |\nabla_x(\rho^{\e},u^{\e},\ta^{\e})|^2|M^{\e}|^{-1/2}(t)\right\|_{L^2_{x,v}} \cr 
&\les \e^2\Big(\int_{\Omega} \mathcal{V}_5 \big(|\nabla_x(\rho^{\e},u^{\e},\ta^{\e})|+|\nabla_x^2(\rho^{\e},u^{\e},\ta^{\e})|\big)dx\Big)^{\frac{1}{2}} + \e^3 \Big(\int_{\Omega} \mathcal{V}_9 |\nabla_x(\rho^{\e},u^{\e},\ta^{\e})|^2dx\Big)^{\frac{1}{2}}.
\enda
\end{align*}
Combining the estimate of $I_{p}^{\al}(t)$ and $II_{p}^{\al}(t)$ proves the claim \eqref{Gembed}. \\

(Proof of \eqref{GL2})
By the same way to the proof of $\eqref{L.Vdecomp}_1$ in Proposition \ref{P.hLinf}, ($\mathrm{Y}(t)$ is derived by the same way to \eqref{L.Vdecompsum}), we have 
\begin{align*}
\bega
\mathfrak{A}(t) &\les \e^2\big(\|\nabla_x (\rho^{\e},u^{\e},\ta^{\e})\|_{L^\infty_x}+\|\nabla_x^2 (\rho^{\e},u^{\e},\ta^{\e})\|_{L^\infty_x} \big)^{\frac{1}{2}}\Big(\mathcal{E}_{tot} +\e^{\frac{3}{2}}\kappa \mathcal{D}_{tot}\Big)^{\frac{1}{2}} \cr 
&+ \e^3 \|\nabla_x (\rho^{\e},u^{\e},\ta^{\e})\|_{L^\infty_x}\Big(\mathcal{E}_{tot} +\e^{\frac{3}{2}}\kappa \mathcal{D}_{tot}\Big)^{\frac{1}{2}} + \e^2 \bigg(e^{-c_1\e^{-{\frac{1}{9}}}}\frac{1}{\e} \mathcal{E}_{tot} \mathrm{Y}(t) \bigg)^{\frac{1}{2}},
\enda
\end{align*}
where the total $h$--part $\mathrm{Y}(t)$ is defined in \eqref{Thdef}.
Then, we use $\|\nabla_x (\rho^{\e},u^{\e},\ta^{\e})\|_{L^\infty_x}\leq \mathcal{E}_M^{\frac{1}{2}}$ and
\begin{align*}
\bega
\|\nabla_x^2 (\rho^{\e},u^{\e},\ta^{\e})\|_{L^\infty_x} \les \begin{cases} C \mathcal{E}_{M}^{\frac{1}{2}}(t),  \quad  &\mbox{for}  \quad  \mathrm{N}\geq 4, \\ 
C \kappa^{-\frac{1}{4}}\mathcal{E}_{tot}^{\frac{1}{2}}(t),  \quad  &\mbox{for}  \quad  \mathrm{N}=3,
\end{cases}
\enda
\end{align*}
because of Agmon's inequality \eqref{Agmon} and \eqref{rutinf}. 
To estimate the term $\mathrm{Y}(t)$, we use \eqref{hLinf} for $0 \leq |\al| \leq \mathrm{N}+1$ in the finite velocity energy case \eqref{caseEC}, and use \eqref{hLinf} for $1 \leq |\al| \leq \mathrm{N}+1$ together with Remark \ref{Rmk.h0} for $|\al|=0$ in the infinite velocity energy case \eqref{caseECX}. This yields 
\begin{align}\label{Thbdd}
\bega
\mathrm{Y}(t) &\les  \begin{cases}
\displaystyle \sum_{0\leq|\al|\leq\mathrm{N}+1}\kappa^{(|\al|-\mathrm{N})_+}  \|\p^{\al}h_0\|_{L^\infty_{x,v}}^2 +\frac{1}{(\e\kappa)^{d}}\|(1+\mathcal{E}_{tot})\|_{L^\infty_t}, & \mbox{for} \quad \eqref{caseEC}, \\
\displaystyle \|\mathfrak{h}_0\|_{L^\infty_{x,v}}^2+ \sum_{1\leq|\al|\leq\mathrm{N}+1}\kappa^{(|\al|-\mathrm{N})_+}  \|\p^{\al}h_0\|_{L^\infty_{x,v}}^2 +\frac{1}{(\e\kappa)^{d}}\|(1+\mathcal{E}_{tot})\|_{L^\infty_t}, & \mbox{for} \quad \eqref{caseECX},
\end{cases} \cr 
&\les \mathrm{Y}(0) + \frac{1}{(\e\kappa)^{d}}\|(1+\mathcal{E}_{tot})\|_{L^\infty_t}.
\enda
\end{align}
Applying \eqref{Thbdd}, we obtain the result \eqref{GL2}.
\\
(Proof of \eqref{GLinfty})  
When we prove \eqref{GL2}, we estimated $I_{p}^{\al}(t)$ in \eqref{P12def} by using the dissipation.  
Here, instead, we estimate it by using the energy. 
Writting \eqref{P12def} again, we have  
\begin{align*}
\bega 
\Big\| \| \la v \ra^{\frac{1}{2}} \p^{\al}\AC{\P}F^{\e}&|M^{\e}|^{-1/2}(t)\|_{L^2_v} \Big\|_{L^\infty_x} \les \sum_{0\leq\beta_1+\beta_2\leq2} \sum_{i=0}^{|\beta_2|} \Big\| \la v \ra^{2i+\frac{1}{2}} \p^{\al+\beta_1}\AC{\P}F^{\e} |M^{\e}|^{-1/2} \cr 
& \times \sum_{\substack{\beta_2^1+\cdots+\beta_2^i=\beta_2 \\ \beta_2^i>0}}|\p^{\beta_2^1}(\rho^{\e},u^{\e},\ta^{\e})|\times \cdots \times |\p^{\beta_2^i}(\rho^{\e},u^{\e},\ta^{\e})| \Big\|_{L^2_{x,v}}.
\enda
\end{align*}
For the finite velocity energy case \eqref{caseEC}, we apply the decomposition $F^{\e} = \mu + \e\sqrt{\mu}f^{\e}$
 for $0\leq |\al|\leq \mathrm{N}+1$.  
On the other hand, for the infinite velocity energy case \eqref{caseECX}, we apply $F^{\e} = \mu + \e\sqrt{\mu}f^{\e}$
 for $1\leq |\al|\leq \mathrm{N}+1$ and $F^{\e} = M^{\e} + \e \sqrt{\bar{\mu}}g^{\e}$ for $|\al|=0$.
Then, by the same way to the proof of $\eqref{L.Vdecomp}_2$, we have
\begin{align*}
\bega 
\Big\| \| \la v \ra^{\frac{1}{2}} \p^{\al}\AC{\P}F^{\e}|M^{\e}|^{-1/2}(t)\|_{L^2_v} \Big\|_{L^\infty_x} &\les \e \bigg(\e^{-\frac{1}{2}}\mathcal{E}_{tot} + e^{-c_1\e^{-{\frac{1}{9}}}}\frac{1}{\e} \mathcal{E}_{tot} \mathrm{Y} \bigg)^{\frac{1}{2}} \sum_{i=0}^{2} \Big(\e^i\mathcal{E}_{tot}^{\frac{i}{2}}\Big) \cr 
&\les \e^{\frac{3}{4}}\mathcal{E}_{tot}^{\frac{1}{2}} + \e^{\frac{1}{2}} e^{-\frac{c_1}{2}\e^{-{\frac{1}{9}}}}\mathrm{Y}^{\frac{1}{2}},
\enda
\end{align*}
where we used $\e^i\mathcal{E}_{tot}^{\frac{i}{2}} \leq C$ from $(\mathcal{B}_1)$ in \eqref{condition}, and the fact that $\lfloor(\mathrm{N}+1)/2\rfloor + 2 \leq \mathrm{N}+1$ for $\mathrm{N} > \tfrac{d}{2}+1$.
Finally, taking $\sup_{0 \leq t \leq T}$ and applying \eqref{Thbdd} completes the proof of \eqref{GLinfty}.
\end{proof}

\StopNoTOC

\section{
Kinetic vorticity  
}\label{Sec.macroglobal}


The objective of this section is to establish control of $\|\nabla_x (\rho^{\e},u^{\e},\ta^{\e})\|_{L^\infty_x}$ in the two-dimensional case.
The key idea is that, by considering the vorticity $\w^{\e}$ and specific entropy fluctuation $\mathfrak{s}^{\e}=\tfrac{3}{2}\ta^{\e}-\rho^{\e}$, we can eliminate singular penalized terms $\tfrac{1}{\e}\nabla_x\cdot u^{\e}$ and $\tfrac{1}{\e}\nabla_x(\rho^{\e}+\ta^{\e})$ in the local conservation laws \eqref{locconNew}. 
We estimate the velocity $u^{\e}$ using the Hodge decomposition \eqref{Hodge}, 
and estimate $\rho^{\e}$ and $\ta^{\e}$ through the acoustic variable $p^{\e}=\rho^{\e}+\ta^{\e}$ and the specific entropy fluctuation 
$\mathfrak{s}^{\e}=\tfrac{3}{2}\ta^{\e}-\rho^{\e}$:
\begin{align}\label{rutadecomp}
\bega
\|\nabla_xu^{\e}\|_{L^\infty_x} &\leq \|\nabla_x\mathbb{P}u^{\e}\|_{L^\infty_x} +\|\nabla_x\mathbb{P}^{\perp}u^{\e}\|_{L^\infty_x}, \cr 
\|\nabla_x\rho^{\e}\|_{L^\infty_x}+\|\nabla_x\ta^{\e}\|_{L^\infty_x} &\leq \|\nabla_xp^{\e}\|_{L^\infty_x} + \|\nabla_x\mathfrak{s}^{\e}\|_{L^\infty_x}.
\enda
\end{align}
We have already estimated $\|\nabla_x\mathbb{P}^{\perp}u^{\e}\|_{L^\infty_x}$ and $\|\nabla_xp^{\e}\|_{L^\infty_x}$ using the Strichartz estimate \eqref{div.ineq} in Proposition \ref{P.div.ineq}.  
For the estimates of $\|\nabla_x\mathbb{P}u^{\e}\|_{L^\infty_x}$ and $\|\nabla_x\mathfrak{s}^{\e}\|_{L^\infty_x}$, we employ the maximum principle and the energy estimate for $\mathbb{P}u^{\e}$ via the vorticity equation for $\w^{\e}$.
As a result, we obtain double exponential growth of $\|\w^{\e}\|_{H^2_x}$ and exponential growth of $\|\nabla_x\mathbb{P}u^{\e}\|_{L^\infty_x}$, both with additional forcing terms.


\begin{definition}\label{def_vorticity}
The (kinetic) vorticity $\w^{\e}$ is defined by
\begin{equation}\label{wdef}
\w^{\e}(t,x) := \nabla_x^\perp \cdot u^{\e}(t,x) = -\p_2 u^{\e}_1(t,x) + \p_1 u^{\e}_2(t,x).
\end{equation}
\end{definition}

\begin{lemma}\label{L.wrtaeqn}
The vorticity $\w^{\e}$ and the specific entropy fluctuation
$\rho^{\e} - \tfrac{3}{2}\ta^{\e}$ satisfy the following equations:
\begin{align}
&\p_t\w^{\e} + u^{\e}\cdot \nabla_x \w^{\e} + (\nabla_x\cdot u^{\e})\w^{\e} =\varPi_{\w}^{\e},  \label{weqnnew} \\ 
&\p_t \mathfrak{s}^{\e} + u^{\e}\cdot \nabla_x \mathfrak{s}^{\e} = \varPi_{\mathfrak{s}}^{\e}, \label{rtaeqnnew}
\end{align}
where the forcing terms are given by
\begin{align}\label{gwrta}
\bega
\varPi_{\w}^{\e}(t,x) &:= -k_B\mathrm{\Theta}^{\e}\nabla_x^{\perp}\ta^{\e}\cdot\nabla_x(\rho^{\e}+\ta^{\e})
-\frac{1}{\e^2}\nabla_x^{\perp}\cdot\bigg(\frac{1}{\mathrm{P}^{\e} }\sum_{j} \p_{x_j} \mathbf{r}_{ij}^{\e}\bigg), \cr 
\varPi_{\mathfrak{s}}^{\e}(t,x) &:= -\frac{1}{\e^2}\frac{1}{k_B\mathrm{P}^{\e}\mathrm{\Theta}^{\e} }\bigg(\sum_{j} \p_{x_j} \mathfrak{q}_j^{\e} + \sum_{i,j} \p_{x_i}\mathrm{U}^{\e}_j \mathbf{r}_{ij}^{\e} \bigg).
\enda
\end{align}
\end{lemma}
\begin{proof}
We multiply $k_B\mathrm{\Theta}^{\e}$ to $\eqref{locconNew}_2$ and then apply $\nabla_x^{\perp}\cdot$.  
Using the identity $\nabla_x \times (u^{\e}\cdot\nabla_x u^{\e}) = u^{\e}\cdot \nabla_x \w^{\e} + (\nabla_x \cdot u^{\e})\w^{\e}$, we get
\begin{align}\label{weqnnew-1}
\bega
&\p_t\w^{\e} + u^{\e}\cdot \nabla_x \w^{\e} + (\nabla_x\cdot u^{\e})\w^{\e}  \cr  &+\nabla_x^{\perp}\cdot\bigg(\frac{k_B}{\e}\mathrm{\Theta}^{\e}\nabla_x(\rho^{\e}+\ta^{\e})\bigg)
+\frac{1}{\e^2}\nabla_x^{\perp}\cdot\bigg(\frac{1}{\mathrm{P}^{\e} }\sum_{j} \p_{x_j} \mathbf{r}_{ij}^{\e}\bigg) =0.
\enda
\end{align}
Since $\nabla_x^{\perp}\cdot\frac{k_B}{\e}\nabla_x(\rho^{\e}+\ta^{\e})=0$, this reduces to the vorticity equation \eqref{weqnnew}.  
Similarly, subtracting $\eqref{locconNew}_3$ from $\eqref{locconNew}_1$ yields the equation \eqref{rtaeqnnew}.
\end{proof}

\hide
\begin{remark} Recall that from the vorticity equation \eqref{weqnnew} in Lemma \ref{L.util}, we proved that the mean vorticity $ \int_{\R^d}\w^{\e}(t,x) dx $ is conserved, since the forcing term is of the stream form. If we directly apply the curl operator to the second equation of \eqref{locconNew}, then the corresponding mean vorticity is not conserved. 
\hide
and then multiply by $k_B\mathrm{\Theta}^{\e}$, we obtain the following alternative vorticity equation:
\begin{align}\label{weqnnew2}
\bega
\p_t\w^{\e} + u^{\e}\cdot \nabla_x \w^{\e} + (\nabla_x\cdot u^{\e})\w^{\e} &=- k_B\mathrm{\Theta}^{\e} \nabla_x^{\perp}\ta^{\e} \cdot \bigg[\nabla_x(\rho^{\e}+\ta^{\e}) 
+\frac{1}{\e}\frac{1}{k_B\mathrm{P}^{\e} \mathrm{\Theta}^{\e}}\sum_{j} \p_{x_j} \mathbf{r}_{ij}^{\e}\bigg] \cr 
&
-\frac{1}{\e^2} k_B\mathrm{\Theta}^{\e} \nabla_x^{\perp}\bigg[\frac{1}{k_B\mathrm{P}^{\e} \mathrm{\Theta}^{\e}}\sum_{j} \p_{x_j} \mathbf{r}_{ij}^{\e}\bigg].
\enda
\end{align}
In contrast, for the equation \eqref{weqnnew2}, this cancellation is no longer immediate, and it is not trivial to deduce $\frac{d}{dt}\int_{\R^d}\w^{\e} dx =0$.
\unhide\end{remark}
\unhide

\hide
\begin{lemma}\label{weqnlem} We have 
\Be\bega\label{weqn}
\p_t \w +u\cdot\nabla_x \w -\eta_0 k_B^{\frac{1}{2}} \kappa \Delta_x \w = - (\nabla_x \cdot u)\w + \nabla_x \times V_m  
\enda
\Ee
where $V_m$ is defined in Lemma \ref{meqnlem}.
\end{lemma}
\begin{proof}
Taking curl on the equation \eqref{meqn} yields 
\begin{align*}
\p_t \w +\nabla_x \times (u\cdot\nabla_x m) -\eta_0 k_B^{\frac{1}{2}} \kappa \Delta_x \w =\nabla_x \times V_m  
\end{align*}
Since 

we have the result. 
\end{proof}

\begin{remark}
We can write vorticity equation without expanding $\AC{\P}F^{\e}$ as follows: 
\begin{align}\label{weqntransport}
\bega
\p_t \w +u\cdot\nabla_x \w  = - (\nabla_x \cdot u)\w -u_i \nabla_x \cdot u + \e m
\cdot\nabla_x (\rho u_i)  -\frac 1 {\eps^2} \nabla_x \times \sum_j \p_j \la \AC{\P}F^{\e},\mathfrak{R}^{\e}_{ij} \ra_{L^2_v}
\enda
\end{align}
\end{remark}

The equation for $\ta$ has penalized terms such as $\frac{1}{\e}\div(u)$. We should carefully consider it. We use some cancellation and Boussinesq relation to obtain the $L^\infty$ estimate of $\nabla_x \ta$.

\begin{lemma}\label{taeqnlem} We have 
\begin{align}\label{taeqn}
\bega
&\p_t\ta+ \frac{5\mathrm{P}^{\e}}{3\mathrm{P}^{\e}+2\mathrm{\Theta}^{\e}} u\cdot\nabla_x\ta -\eta_1k_B^{1/2} \kappa \frac{5\sqrt{\mathrm{\Theta}^{\e}}}{3\mathrm{P}^{\e}+2\mathrm{\Theta}^{\e}} \Delta_x\ta =V_{\ta}
\enda
\end{align}
where
\begin{align*}
\bega
V_{\ta}&=\frac{2\mathrm{\Theta}^{\e}}{3\mathrm{P}^{\e}+2\mathrm{\Theta}^{\e}} \lw(\p_t\rho+\p_t\ta\rw) +\frac{1}{\e}\frac{2}{3\mathrm{P}^{\e}+2\mathrm{\Theta}^{\e}}
 u
\cdot\nabla_x(\mathrm{P}^{\e} \mathrm{\Theta}^{\e}) +\frac{5}{3\mathrm{P}^{\e}+2\mathrm{\Theta}^{\e}}
\eta_1k_B^{1/2} \kappa  \frac{\e|\nabla_x\ta|^2}{2\sqrt{\mathrm{\Theta}^{\e}}}\cr 
& -\frac{1}{k_B\e^2}\frac{2}{3\mathrm{P}^{\e}+2\mathrm{\Theta}^{\e}}
\sum_{i,j} \p_{x_i}\mathrm{U}^{\e}_j \mathbf{r}_{ij}^{\e} \cr 
&+\frac{1}{k_B\e^2}\frac{2}{3\mathrm{P}^{\e}+2\mathrm{\Theta}^{\e}}
\sum_j \p_j\lw\{\lw\la \mathcal{L}^{-1}\lw\{\kappa\eps \lw(\eps\p_t \AC{\P}F^{\e}+(\II-\P)(v\cdot\nabla_x \AC{\P}F^{\e})\rw)-\mathcal{N}(\AC{\P}F^{\e},\AC{\P}F^{\e})\rw\},\mathcal{Q}^{\e}_j\rw\ra\rw\}
\enda
\end{align*}
\end{lemma}
\begin{proof}
By the same way with Lemma \ref{meqnlem}, applying \eqref{Gform} and Lemma \ref{ABcomp} gives 
\begin{align}\bega \label{GBform}
\la \AC{\P}F^{\e}, \mathcal{Q}^{\e}_j\ra 
&=-\kappa\eps \sum_{m}\frac{\p_m\mathrm{\Theta}^{\e}}{k_B|\mathrm{\Theta}^{\e}|^2} (k_B\mathrm{\Theta}^{\e})^{5/2} \frac{5}{2}\eta_1 \delta_{mj}  \\
&\quad-\lw\la \mathcal{L}^{-1}\lw\{\kappa\eps \lw(\eps\p_t \AC{\P}F^{\e}+(\II-\P)(v\cdot\nabla_x \AC{\P}F^{\e})
\rw)-\mathcal{N}(\AC{\P}F^{\e},\AC{\P}F^{\e})
\rw\},\mathcal{Q}^{\e}_j\rw\ra
\enda
\end{align}
We also denote the first line of \eqref{GBform} as $GB_1$. Taking $\sum_j \p_j$ on $GB_1$ gives 
\[\bega
\sum_j \p_j(\AC{\P}F^{\e}B_1)
&=\eta_1 \kappa\eps \frac{5}{2} k_B^{\frac{3}{2}} \sum_j\lw(\p_j^2\mathrm{\Theta}^{\e} \sqrt{\mathrm{\Theta}^{\e}} + \frac{(\p_j\mathrm{\Theta}^{\e})^2}{2\sqrt{\mathrm{\Theta}^{\e}}}\rw)
\enda\]
This combining with \eqref{GBform} and 
$1/(3k_B\e)\times \eqref{loccon}_3$ gives desired result.
\begin{align}\label{ta-La}
\bega
&\p_t(\mathrm{P}^{\e}\mathrm{\Theta}^{\e})
+\frac{1}{\e} \lw(\mathrm{U}^{\e}
\cdot\nabla_x(\mathrm{P}^{\e} \mathrm{\Theta}^{\e})+ \frac{5}{3}\nabla_x \cdot \mathrm{U}^{\e}
(\mathrm{P}^{\e} \mathrm{\Theta}^{\e}) \rw)  \cr 
&\quad + \frac{2}{3k_B\e}\sum_{i,j} \p_{x_i}\mathrm{U}^{\e}_j \mathbf{r}_{ij}^{\e} - \frac{5}{3}\eta_1k_B^{1/2} \kappa \lw(\Delta_x\mathrm{\Theta}^{\e} \sqrt{\mathrm{\Theta}^{\e}} + \frac{|\nabla_x\mathrm{\Theta}^{\e}|^2}{2\sqrt{\mathrm{\Theta}^{\e}}}\rw) \\
&\quad-\frac{2}{3k_B\e}\sum_j \p_j\lw\{\lw\la \mathcal{L}^{-1}\lw\{\kappa\eps \lw(\eps\p_t \AC{\P}F^{\e}+(\II-\P)(v\cdot\nabla_x \AC{\P}F^{\e})
\rw)-\mathcal{N}(\AC{\P}F^{\e},\AC{\P}F^{\e})
\rw\},\mathcal{Q}^{\e}_j\rw\ra\rw\}=0
\enda
\end{align}
From the local energy conservation law $\eqref{loccon}_3$, dividing each side by $\e$, we have 
\begin{align}\label{}
\bega
&\mathrm{\Theta}^{\e}\p_t\rho +\mathrm{P}^{\e}\p_t\ta
+\frac{1}{\e^2} \lw(\mathrm{U}^{\e}
\cdot\nabla_x(\mathrm{P}^{\e} \mathrm{\Theta}^{\e})+ \frac{5}{3}\nabla_x \cdot \mathrm{U}^{\e}
(\mathrm{P}^{\e} \mathrm{\Theta}^{\e}) \rw)  \cr 
&\quad + \frac{2}{3k_B\e^2}\sum_{i,j} \p_{x_i}\mathrm{U}^{\e}_j \mathbf{r}_{ij}^{\e} - \frac{5}{3}\eta_1k_B^{1/2} \kappa \lw(\Delta_x\ta \sqrt{\mathrm{\Theta}^{\e}} + \frac{\e|\nabla_x\ta|^2}{2\sqrt{\mathrm{\Theta}^{\e}}}\rw) \\
&\quad-\frac{2}{3k_B\e^2}\sum_j \p_j\lw\{\lw\la \mathcal{L}^{-1}\lw\{\kappa\eps \lw(\eps\p_t \AC{\P}F^{\e}+(\II-\P)(v\cdot\nabla_x \AC{\P}F^{\e})
\rw)-\mathcal{N}(\AC{\P}F^{\e},\AC{\P}F^{\e})
\rw\},\mathcal{Q}^{\e}_j\rw\ra\rw\}=0
\enda
\end{align}
If we try to use the maximum principle with repect to $\p_t\ta$ or $\p_t\nabla_x \ta$, then the terms $\mathrm{\Theta}^{\e}\p_t\rho$ and $\frac{1}{\e^2}\nabla_x \cdot \mathrm{U}^{\e}(\mathrm{P}^{\e} \mathrm{\Theta}^{\e})$ blow up in $L^\infty_x$. So, we use the cancellation by using the density conservation law $\eqref{loccon}_1$, 
\begin{align*}
\bega
&\p_t\rho = -\frac{1}{\e^2}\nabla_x \cdot (\mathrm{P}^{\e}\mathrm{U}^{\e})
\enda
\end{align*}
to extract $-\p_t \rho$ on the left-side: 
\begin{align}\label{}
\bega
&-\frac{2}{3}\mathrm{\Theta}^{\e}\p_t\rho+\mathrm{P}^{\e}\p_t\ta
+\frac{1}{\e^2} \lw(-\frac{2}{3}\mathrm{U}^{\e}
\cdot\nabla_x(\mathrm{P}^{\e} \mathrm{\Theta}^{\e})
+\frac{5}{3}(\mathrm{P}^{\e}\mathrm{U}^{\e})\cdot\nabla_x\mathrm{\Theta}^{\e} \rw)  \cr 
&\quad + \frac{2}{3k_B\e^2}\sum_{i,j} \p_{x_i}\mathrm{U}^{\e}_j \mathbf{r}_{ij}^{\e} - \frac{5}{3}\eta_1k_B^{1/2} \kappa \lw(\Delta_x\ta \sqrt{\mathrm{\Theta}^{\e}} + \frac{\e|\nabla_x\ta|^2}{2\sqrt{\mathrm{\Theta}^{\e}}}\rw) \\
&\quad-\frac{2}{3k_B\e^2}\sum_j \p_j\lw\{\lw\la \mathcal{L}^{-1}\lw\{\kappa\eps \lw(\eps\p_t \AC{\P}F^{\e}+(\II-\P)(v\cdot\nabla_x \AC{\P}F^{\e})
\rw)-\mathcal{N}(\AC{\P}F^{\e},\AC{\P}F^{\e})
\rw\},\mathcal{Q}^{\e}_j\rw\ra\rw\}=0
\enda
\end{align}
We add $\frac{2}{3}\mathrm{\Theta}^{\e}\lw(\p_t\rho+\p_t\ta\rw)$ on both sides, then divide each side by $\frac{3\mathrm{P}^{\e}+2\mathrm{\Theta}^{\e}}{3}$ to have 
\begin{align*}
\bega
&\p_t\ta+\frac{5\mathrm{P}^{\e}}{3\mathrm{P}^{\e}+2\mathrm{\Theta}^{\e}} u\cdot\nabla_x\ta -\frac{5\sqrt{\mathrm{\Theta}^{\e}}}{3\mathrm{P}^{\e}+2\mathrm{\Theta}^{\e}}\eta_1k_B^{1/2} \kappa \Delta_x\ta   \cr 
&=\frac{2\mathrm{\Theta}^{\e}}{3\mathrm{P}^{\e}+2\mathrm{\Theta}^{\e}} \lw(\p_t\rho+\p_t\ta\rw) +\frac{1}{\e}\frac{2}{3\mathrm{P}^{\e}+2\mathrm{\Theta}^{\e}}
 u
\cdot\nabla_x(\mathrm{P}^{\e} \mathrm{\Theta}^{\e}) +\frac{5}{3\mathrm{P}^{\e}+2\mathrm{\Theta}^{\e}}
\eta_1k_B^{1/2} \kappa  \frac{\e|\nabla_x\ta|^2}{2\sqrt{\mathrm{\Theta}^{\e}}}\cr 
& -\frac{1}{k_B\e^2}\frac{2}{3\mathrm{P}^{\e}+2\mathrm{\Theta}^{\e}}
\sum_{i,j} \p_{x_i}\mathrm{U}^{\e}_j \mathbf{r}_{ij}^{\e} \cr 
&+\frac{1}{k_B\e^2}\frac{2}{3\mathrm{P}^{\e}+2\mathrm{\Theta}^{\e}} \sum_j \p_j\lw\{\lw\la \mathcal{L}^{-1}\lw\{\kappa\eps \lw(\eps\p_t \AC{\P}F^{\e}+(\II-\P)(v\cdot\nabla_x \AC{\P}F^{\e})\rw)-\mathcal{N}(\AC{\P}F^{\e},\AC{\P}F^{\e})\rw\},\mathcal{Q}^{\e}_j\rw\ra\rw\}
\enda
\end{align*}
Adding and subtracting some coefficients on the first line, we have 

\end{proof}

\unhide

\begin{lemma}\label{L.util}
Let $\w^{\e}(t,x)$ be the solution to the equation \eqref{weqnnew}. Then we have
\begin{align*}
\bega
\frac{d}{dt}\int_{\R^2} \w^{\e}(t,x) dx =0, \qquad \int_{\R^2} \w^{\e}_0(x) dx = \int_{\R^2} \nabla_x^{\perp}\cdot\bar{u}(x) dx.
\enda
\end{align*}
\end{lemma}
\begin{proof}
Integrating \eqref{weqnnew} over $\R^2$ and applying integration by parts, we obtain the first identity, since
\begin{align*}
\bega
\int_{\R^2} \Big(u^{\e}\cdot \nabla_x \w^{\e} + (\nabla_x\cdot u^{\e})\w^{\e}\Big) dx = - \int_{\R^2} (\nabla_x\cdot u^{\e})\w^{\e}dx + \int_{\R^2}(\nabla_x\cdot u^{\e})\w^{\e} dx = 0.
\enda
\end{align*}
By the definition of the radial-energy decomposition in Definition~\ref{D.Ra-E}, we have
\begin{align*}
\bega
\int_{\R^2} \w^{\e}(t,x) dx = \int_{\R^2} \w^{\e}_0(x) dx = \int_{\R^2} \nabla_x^{\perp}\cdot\bar{u}(x) dx.
\enda
\end{align*}
\end{proof}

\subsection{Maximum principle}

In this subsection, we estimate \(\|\w^{\e}\|_{L^\infty_x}\) and \(\|\nabla_x \mathfrak{s}^{\e}\|_{L^\infty_x}\) by applying the maximum principle.

\begin{lemma}\label{maximum} 
Let \(W(t, \cdot) \in C^2(\R^2)\) be a solution to the convection–diffusion equation
\beq\label{heat-max}
\p_t W +U(t,x)\cdot\nabla_x W- \mathrm{J}(t,x)\triangle W=\mathrm{K}(t,x)W+ \mathrm{H},
\eeq
where \[
\mathrm{J}(t,x)\ge 0, \qquad \|\mathrm{K}(t,x)\|_{L^\infty_x} < \infty.
\]
Then the following estimate holds:
\[
\|W(t)\|_{C(\R^2)}\le e^{\int_0^t \|\mathrm{K}(s)\|_{C(\R^2)}ds}\|W (0)\|_{C(\R^2)}+\int_0^t  e^{\int_s^t \|\mathrm{K}(\tau)\|_{C(\R^2)}d\tau} \|\mathrm{H}(s)\|_{C(\R^2)} ds,
\]
where $
\|W\|_{C(\R^2)}=\sup_{x\in \R^2}|W(x)|$.
\end{lemma}
\begin{proof}
Although the proof is classical, we include it here for completeness.  
When \(t = 0\), the inequality is immediate. Hence we assume \(t > 0\).  
Let \(|W(t,x)|\) attain its maximum value at the point \(x = x_\star \in \R^2\), and denote $W_{max}(t) := \|W(t)\|_{C(\R^2)} = |W(t,x_\star)|$.
We first consider the case \(W(t,x_\star) = W_{max}(t) \ge 0\) (the case of a minimum can be treated by considering \(-W\)).  
At \(x = x_\star\), we have the standard maximum principle properties
\[
\nabla_x W(t, x)|_{x=x_\star} = 0, \qquad \mathrm{J}(t, x_\star) \triangle W(t, x)|_{x=x_\star} \leq 0.
\]
Evaluating \eqref{heat-max} at \((t,x_\star)\), we obtain
\begin{align*}
\frac{d}{dt}W(t,x_\star)
& = (-U\cdot\nabla_x W + \mathrm{J}\triangle W + \mathrm{K}W + \mathrm{H})|_{(t,x_\star)}\\
& \le \mathrm{K}(t,x_\star)W(t,x_\star) + \mathrm{H}(t,x_\star) \cr 
&\le \|\mathrm{K}(t)\|_{C(\R^2)}\|W(t)\|_{C(\R^2)} + \|\mathrm{H}(t)\|_{C(\R^2)}.
\end{align*}
Thus, \(W_{max}(t)\) satisfies the differential inequality
\[
\frac{d}{dt} W_{max}(t) \le \|\mathrm{K}(t)\|_{C(\R^2)}\, W_{max}(t) + \|\mathrm{H}(t)\|_{C(\R^2)}.
\]
Applying Grönwall’s inequality yields the desired estimate.
\end{proof}

\begin{proposition}\label{P.max}
Let $\Omega=\R^2$, and let \(\w^{\e}\) and \(\rho^{\e} - \frac{3}{2}\ta^{\e}\) be solutions of \eqref{weqnnew} and \eqref{rtaeqnnew}, respectively. Suppose that the bootstrap assumption \eqref{condition} holds on the interval $[0,T]$. Then, for all $t \in [0,T]$, the following estimates hold:
\begin{align}\label{wmax}
\bega
\|&\w^{\e}(t)\|_{L^\infty_x} \leq e^{\int_0^t \|(\nabla_x \cdot u^{\e})(s)\|_{L^\infty_x}ds}\|\w^{\e}_0\|_{L^\infty_x} \cr
&+C \int_0^t e^{\int_s^t \|(\nabla_x \cdot u^{\e})(\tau)\|_{L^\infty_x}d\tau} \Big((\mathcal{E}_M^{\mathrm{N}}(F^{\e}(s)))^{\frac{1}{2}} \|\nabla_x(\rho^{\e}+\ta^{\e})(s)\|_{L^\infty_x} + \kappa^{\frac{1}{2}}(\mathcal{D}_G^{\mathrm{N}}(F^{\e}(s)))^{\frac{1}{2}}\Big) ds, 
\enda
\end{align}
and
\begin{align}\label{rtamax}
\bega
\|\nabla_x\mathfrak{s}^{\e}(t)\|_{L^\infty_x}  \leq e^{\int_0^t \|\nabla_xu^{\e}(s)\|_{L^\infty_x}ds }\|\nabla_x\mathfrak{s}^{\e}_0\|_{L^\infty_x}   
 +C\kappa^{\frac{1}{2}} \int_0^t  e^{\int_s^t \|\nabla_xu^{\e}(\tau)\|_{L^\infty_x}d\tau }  (\mathcal{D}_G^{\mathrm{N}}(F^{\e}(s)))^{\frac{1}{2}}ds.
\enda
\end{align}
Here, $\mathcal{E}_M^{\mathrm{N}}(F^{\e}(t))$ and $\mathcal{D}_G^{\mathrm{N}}(F^{\e}(t))$ are defined in \eqref{N-EDdef} for $\mathrm{N}\geq 4$.
\end{proposition}
\begin{proof}
In the proof, for brevity, we slightly abuse notation by writing $\mathcal{E}(t)$, and $\mathcal{D}(t)$ for $\mathcal{E}^{\mathrm{N}}(F^{\e}(t))$, and $\mathcal{D}^{\mathrm{N}}(F^{\e}(t))$, respectively.
For the estimate of \(\varPi_{\w}^{\e}\) and \(\varPi_{\mathfrak{s}}^{\e}\) defined in \eqref{gwrta}, we claim that
\begin{align}\label{Gesti}
\bega
&\|\varPi_{\w}^{\e}(t) \|_{L^\infty_x} \leq C\Big(\mathcal{E}_M^{\frac{1}{2}} \|\nabla_x(\rho^{\e}+\ta^{\e})(t)\|_{L^\infty_x} + \kappa^{\frac{1}{2}}\mathcal{D}_G^{\frac{1}{2}}(t)\Big), \cr 
&\|\nabla_x\varPi_{\mathfrak{s}}^{\e}(t)\|_{L^\infty_x} \leq C\kappa^{\frac{1}{2}}\mathcal{D}_G^{\frac{1}{2}}(t).
\enda
\end{align}
Taking the $L^\infty_x$-norm of each term in $\varPi_{\w}^{\e}(t,x)$, we first obtain
\begin{align*}
\bega
\frac{1}{\e}\big\|\nabla_x^{\perp}\mathrm{\Theta}^{\e}\cdot\nabla_x(\rho^{\e}+\ta^{\e})\big\|_{L^\infty_x} \leq \mathcal{E}_M^{\frac{1}{2}} \|\nabla_x(\rho^{\e}+\ta^{\e})\|_{L^\infty_x},
\enda
\end{align*}
and
\begin{align}\label{ahatesti}
\bega
\bigg\|\nabla_x^{\perp}\cdot\bigg(\frac{1}{\mathrm{P}^{\e}}\sum_{j} \p_{x_j} \mathbf{r}_{ij}^{\e}\bigg)\bigg\|_{L^\infty_x} &\leq \bigg\|\nabla_x\frac{1}{\mathrm{P}^{\e}}\bigg\|_{L^\infty_x}\|\nabla_x\mathbf{r}_{ij}^{\e}\|_{L^\infty_x} + \bigg\|\frac{1}{\mathrm{P}^{\e}}\bigg\|_{L^\infty_x}\|\nabla_x^2\mathbf{r}_{ij}^{\e}\|_{L^\infty_x} + \kappa^{\frac{1}{2}}\mathcal{D}_G^{\frac{1}{2}}(t) \cr 
&\les \e\mathcal{E}_M^{\frac{1}{2}}\big(\e^2\kappa^{\frac{1}{2}}\mathcal{D}_G^{\frac{1}{2}}(t)\big) + \e^2\kappa^{\frac{1}{2}}\mathcal{D}_G^{\frac{1}{2}}(t) \les \e^2\kappa^{\frac{1}{2}}\mathcal{D}_G^{\frac{1}{2}}(t),
\enda
\end{align}
where we used \eqref{pTaleq}, $\eqref{ABGscale}_1$, and $(\mathcal{B}_2)$ in \eqref{condition}.  
This proves the first inequality in \eqref{Gesti}.  
The second inequality of \eqref{Gesti} follows from the same reasoning.
Next, applying Lemma \ref{maximum} to equation \eqref{weqnnew} with the choices \(W = \w^{\e}\), \(U = u^{\e}\), \(\mathrm{J} = 0\), \(\mathrm{K} = -\nabla_x \cdot u^{\e}\), and \(\mathrm{H} = \varPi_{\w}^{\e}\), together with $\eqref{Gesti}_1$, yields the bound \eqref{wmax}.
For \eqref{rtamax}, we first apply \(\nabla_x\) to equation \eqref{rtaeqnnew}, obtaining
\begin{align*}
&\p_t \nabla_x\mathfrak{s}^{\e} + u^{\e} \cdot \nabla_x \big(\nabla_x\mathfrak{s}^{\e}\big) = - \nabla_xu^{\e} \cdot \nabla_x \mathfrak{s}^{\e} + \nabla_x\varPi_{\mathfrak{s}}^{\e}. 
\end{align*}
Applying Lemma \ref{maximum} with \(W = \nabla_x\mathfrak{s}^{\e}\), \(U = u^{\e}\), \(\mathrm{J} = 0\), \(\mathrm{K} = -\nabla_x u^{\e}\), and \(\mathrm{H} = \nabla_x \varPi_{\mathfrak{s}}^{\e}\), together with $\eqref{Gesti}_2$, gives the desired result \eqref{rtamax}.
\end{proof}

\hide
\begin{remark} Once we directly compute $\eqref{locconS}_1-\frac{3}{2}\eqref{locconS}_3$ without expanding $\mathfrak{q}_j^{\e}$, then we can have 
\begin{align*}
\p_t (\rho-\frac{3}{2}\ta) + u \cdot \nabla_x (\rho-\frac{3}{2}\ta) + (\rho-\ta) (\nabla_x \cdot u) = \e\rho u \cdot \nabla_x (\rho-\frac{3}{2}\ta) +\frac{1}{\e^2k_B\mathrm{P}^{\e}}\sum_{j} \p_{x_j} \mathfrak{q}_j^{\e} + \frac{1}{\e^2k_B\mathrm{P}^{\e}}\sum_{i,j} \p_{x_i}\mathrm{U}^{\e}_j \mathbf{r}_{ij}^{\e}.
\end{align*}
By the same reason in Remark \ref{R.wtransp}, the quantity $\nabla_x\p_j\mathfrak{q}_j^{\e}$ cannot be bounded in $L^\infty_x$:
\begin{align*}
\bega
\Big\|\frac 1 {\eps^2} \nabla_x\p_j \la \AC{\P}F^{\e},\mathcal{Q}^{\e}_{j} \ra_{L^2_v}\Big\|_{L^\infty_x} &\les  \Big\|\frac 1 {\eps^2} \nabla_x\p_j \la \AC{\P}F^{\e},\mathcal{Q}^{\e}_{j} \ra_{L^2_v}\Big\|_{L^2_x}^{\frac{1}{2}} \Big\|\frac 1 {\eps^2} \nabla_x\p_j \la \AC{\P}F^{\e},\mathcal{Q}^{\e}_{j} \ra_{L^2_v}\Big\|_{H^2_x}^{\frac{1}{2}} \les \kappa^{-\frac{1}{2}}\mathcal{D}_{top}^{\frac{1}{2}} \rightarrow \infty \quad \mbox{as} \quad \kappa \downarrow 0.
\enda
\end{align*}
\end{remark}
\unhide

\subsection{Vorticity energy}
In this subsection, we derive an explicit exponential growth estimate by combining the bounds for $\|\w^{\e}\|_{H^2_x}$ and $\|\nabla_x \mathbb{P}u^{\e}\|_{L^\infty_x}$.
In two dimensions, the potential estimate is crucial, providing logarithmic control of the vorticity energy. However, it alone is not sufficient for global-in-time bounds, due to the growth of the Boltzmann energy–dissipation through $\|\nabla_x(\rho^{\e},\ta^{\e})\|_{L^\infty_x}$. Hence, together with strong dissipation from macro–micro cancellation, it enables global control by effectively decoupling the macroscopic and microscopic components.
Since the acoustic variables are controlled only below order $\mathrm{N}+1$, we minimize the derivatives acting on them.

\begin{lemma}
For the vorticity $\w^{\e}$ defined in \eqref{wdef} and the Leray projector defined in \eqref{LerayPdef}, we have 
\begin{align}\label{pmw}
\bega
\|\nabla_x \mathbb{P}u^{\e}\|_{H^k_x} \leq C\|\w^{\e}\|_{H^k_x}, \qquad 
\|\w^{\e}\|_{H^k_x} \leq \|\nabla_x \mathbb{P}u^{\e}\|_{H^k_x}, \quad \text{for} \quad k \in \mathbb{N} \cup \{0\},
\enda
\end{align}
for some positive constant $C$. Moreover,
\begin{align}\label{Be-Ma,m}
\bega
\|\nabla_x \mathbb{P}u^{\e}\|_{L^\infty_x} \les \big(1 + \ln^+\|\w^{\e}\|_{H^2_x}\big)\big(1 + \|\w^{\e}\|_{L^\infty_x}\big).
\enda
\end{align}
\end{lemma}

\begin{proof}
By the definition of the Leray projector \eqref{LerayPdef}, the vector field \(\mathbb{P}u^{\e}\) is divergence-free, and moreover \(\nabla_x^{\perp}\cdot\mathbb{P}u^{\e} = \w^{\e}\) since $\nabla_x^{\perp}\cdot\mathbb{P}^{\perp}u^{\e}=0$. Therefore, applying \eqref{puw} with \(u_* = \mathbb{P}u^{\e}\), we obtain
\begin{align*}
\bega
\|\w^{\e}\|_{H^k_x} \les \|\nabla_x\mathbb{P}u^{\e}\|_{H^k_x} \les \|\w^{\e}\|_{H^k_x},
\enda
\end{align*}
which proves \eqref{pmw}.  
Since $\nabla_x^{\perp}\cdot \mathbb{P}u^{\e} = \w^{\e}$, applying the inequality \eqref{Be-Ma} yields \eqref{Be-Ma,m}.
\end{proof}

\hide
\begin{align*}
\bega
\varPi_{\w}^{\e} &:=\eta_0 k_B^{\frac{1}{2}} \kappa(|\mathrm{\Theta}^{\e}|^{\frac 1 2}-1)  \Delta_x \w^{\e} - (\nabla_x \cdot u^{\e})\w^{\e} +\eta_0 \kappa 
 \nabla^\perp (k_B^{\frac 1 2}|\mathrm{\Theta}^{\e}|^{\frac 1 2}) \cdot \Big(\Delta_x m  +  \nabla_x\nabla_x\cdot u^{\e}
\Big) \cr 
& + \nabla_x^\perp \cdot \bigg[-u_i \nabla_x \cdot u^{\e} + \e u\cdot\nabla_x (\rho u_i) -\eta_0 k_B^{\frac 1 2} \kappa\eps|\mathrm{\Theta}^{\e}|^{\frac 1 2}\lw(\frac{1}{3}\p_i\nabla_x\cdot (\rho u_i) + \Delta_x(\rho u_i)\rw) \cr 
&+\eta_0 \kappa \sum_j \p_j(k_B^{\frac 1 2}|\mathrm{\Theta}^{\e}|^{\frac 1 2})\lw(\p_i u_j+\p_j u_i-\frac 2 3 \delta_{ij}(\nabla_x\cdot u^{\e})\rw) \bigg] \cr 
&+ \nabla_x^\perp \cdot \frac{1}{\e^2}\sum_j \p_j\lw\la \mathcal{L}^{-1}\lw\{\kappa\eps \lw(\eps\p_t \AC{\P}F^{\e}+\AC{\P}(v\cdot\nabla_x \AC{\P}F^{\e})\rw)-\mathcal{N}(\AC{\P}F^{\e},\AC{\P}F^{\e})\rw\},\mathfrak{R}^{\e}_{ij}\rw\ra_{L^2_v}
\enda
\end{align*}
\begin{align*}
\bega 
\varPi_{\mathfrak{s}} &:= \eta_1 k_B^{\frac{1}{2}} \kappa \Big(\frac{|\mathrm{\Theta}^{\e}|^{\frac 1 2}}{\mathrm{P}^{\e}}-1\Big)  \Delta_x(\rho-\frac{3}{2}\ta) -(\rho-\ta) (\nabla_x \cdot u^{\e}) -\eta_1\frac{k_B^{\frac 1 2}|\mathrm{\Theta}^{\e}|^{\frac 1 2}}{\mathrm{P}^{\e}} \kappa \Delta_x(\rho+\ta) +\e\rho u \cdot \nabla_x \rho \cr 
&- \frac{3}{2}\bigg[\e(\rho u)\cdot \nabla_x \ta + \frac{5}{6}\eta_1 \e\kappa   \frac{k_B^{1/2}}{\mathrm{P}^{\e}\sqrt{\mathrm{\Theta}^{\e}}}|\nabla_x\ta|^2-\frac{2}{\eps^23k_B\mathrm{P}^{\e}}\sum_{i,j} \p_{x_i}\mathrm{U}^{\e}_j \mathbf{r}_{ij}^{\e} \bigg] \cr 
&- \frac{1}{\e^2 k_B \mathrm{P}^{\e}} \sum_j \p_j \lw\la \mathcal{L}^{-1}\lw\{\kappa\eps \lw(\eps\p_t \AC{\P}F^{\e}+\AC{\P}(v\cdot\nabla_x \AC{\P}F^{\e})\rw)-\mathcal{N}(\AC{\P}F^{\e},\AC{\P}F^{\e})\rw\},\mathcal{Q}^{\e}_j\rw\ra_{L^2_v}
\enda
\end{align*}
\unhide

\begin{lemma}\label{L.wrtaHs}
Let $\Omega=\R^2$, and let \(\w^{\e}\) be the solution to \eqref{weqnnew}. Then $\w^{\e}$ satisfies the following estimate:
\begin{align}\label{wHs}
\bega
\frac{d}{dt}\|\w^{\e}\|_{H^2_x}^2 
&\les \Big(\|\nabla_x \mathbb{P}u^{\e}\|_{L^\infty_x}+\|\nabla_x \mathbb{P}^{\perp}u^{\e}\|_{L^\infty_x}+\|\w^{\e}\|_{L^\infty_x}+\|\nabla_x^2 \mathbb{P}^{\perp}u^{\e}\|_{L^\infty_x}\Big) \|\w^{\e}\|_{H^2_x}^2  \cr 
& +\bigg(\|\nabla_x\nabla_x\cdot u^{\e}\|_{L^\infty_x}\|\w^{\e}\|_{L^2_x}+ \bigg\| k_B\mathrm{\Theta}^{\e}\nabla_x^{\perp}\ta^{\e}\cdot\nabla_x(\rho^{\e}+\ta^{\e})\bigg\|_{H^1_x}\bigg)\|\w^{\e}\|_{H^3_x}  \cr 
&+ \bigg\|\frac{1}{\e^2}\nabla_x^{\perp}\cdot\bigg(\frac{1}{\mathrm{P}^{\e} }\sum_{j} \p_{x_j} \mathbf{r}_{ij}^{\e}\bigg)\bigg\|_{H^2_x}\|\w^{\e}\|_{H^2_x}.
\enda
\end{align}
\end{lemma}
\begin{proof}
We perform an energy estimate in \(H^2_x\) for the vorticity equation \eqref{weqnnew} as follows:
\begin{align}\label{wE1}
\bega
\frac{1}{2}\frac{d}{dt}\|\w^{\e}\|_{H^2_x}^2 
&\les I_{\w} + \sum_{0 \leq |\al_x| \leq 2}\int_{\Omega}\p^{\al_x}\varPi_{\w}^{\e} \p^{\al_x}\w^{\e} dx,
\enda
\end{align}
where
\begin{align*}
\bega
I_{\w}(t)&:= -\sum_{0 \leq |\al_x| \leq 2} \int_{\Omega} \p^{\al_x}(u^{\e} \cdot \nabla_x \w^{\e}) \p^{\al_x} \w^{\e} \, dx -\sum_{0 \leq |\al_x| \leq 2} \int_{\Omega} \p^{\al_x}\big((\nabla_x\cdot u^{\e})\w^{\e}\big) \p^{\al_x} \w^{\e} \, dx .
\enda
\end{align*}
For the term $I_{\w}(t)$, an integration by parts yields
\begin{align}\label{wE2}
\bega
I_{\w}(t) &\leq \frac{1}{2}\|\nabla_x \cdot u^{\e} \|_{L^\infty_x} \|\w^{\e}\|_{H^2_x}^2 + I_{\w}^{1}(t) + I_{\w}^{2}(t) ,
\enda
\end{align}
where we separated the contributions involving $\nabla_x\cdot u^{\e}$:
\begin{align*}
\bega
I_{\w}^{1}(t)&:= \sum_{0 \leq |\al_x| \leq 2}\int_{\Omega} \bigg(\big(\p^{\al_x}(u^{\e} \cdot \nabla_x \w^{\e}) - u^{\e} \cdot \nabla_x \p^{\al_x} \w^{\e}\big)\p^{\al_x} \w^{\e} \bigg) dx, \cr 
I_{\w}^{2}(t)&:= \sum_{0 \leq |\al_x| \leq 2}\int_{\Omega} \bigg(\big(\p^{\al_x}\big((\nabla_x\cdot u^{\e})\w^{\e}\big) - (\nabla_x\cdot u^{\e})\p^{\al_x}\w^{\e} \big)\p^{\al_x} \w^{\e}\bigg) dx .
\enda
\end{align*}
For $I_{\w}^{1}(t)$, the velocity $u^{\e}$ necessarily carries at least one derivative, so that
\begin{align*}
\bega
I_{\w}^{1}(t)&\les \|\nabla_x u^{\e}\|_{L^\infty_x} \|\w^{\e}\|_{H^2_x}^2 + \|\nabla_x^2 u^{\e} \nabla_x \w^{\e}\|_{L^2_x} \|\w^{\e}\|_{H^2_x}.
\enda
\end{align*}
We decompose the last term using the Hodge decomposition \eqref{Hodge}:
\begin{align*}
\bega
\|\nabla_x^2 u^{\e} \nabla_x \w^{\e}\|_{L^2_x} &\les \|\nabla_x^2 \mathbb{P}u^{\e} \nabla_x \w^{\e}\|_{L^2_x} + \|\nabla_x^2 \mathbb{P}^{\perp}u^{\e} \nabla_x \w^{\e}\|_{L^2_x}.
\enda
\end{align*}
Using \eqref{uvHk} with \(u_* = \nabla_x u^{\e}\), \(v_* = \w^{\e}\) for the $\mathbb{P}u^{\e}$ part, and estimating in $L^\infty_x$ for the $\mathbb{P}^{\perp}u^{\e}$ part, we obtain
\begin{align}\label{wE3}
\bega
\|\nabla_x^2 \mathbb{P}u^{\e} \nabla_x \w^{\e}\|_{L^2_x} &\les \|\nabla_x\mathbb{P}u^{\e}\|_{H^2_x}\|\w^{\e}\|_{L^\infty_x}+\|\nabla_x\mathbb{P}u^{\e}\|_{L^\infty_x}\|\w^{\e}\|_{H^2_x}, \cr 
\|\nabla_x^2 \mathbb{P}^{\perp}u^{\e} \nabla_x \w^{\e}\|_{L^2_x} &\les  \|\nabla_x^2 \mathbb{P}^{\perp}u^{\e}\|_{L^\infty_x} \|\nabla_x \w^{\e}\|_{L^2_x}.
\enda
\end{align}
Since $\|\nabla_x \mathbb{P}u^{\e}\|_{H^2_x} \les \|\w^{\e}\|_{H^2_x}$ by \eqref{pmw}, this yields
\begin{align*}
\bega
I_{\w}^{1}(t)&\les \|\nabla_x u^{\e}\|_{L^\infty_x} \|\w^{\e}\|_{H^2_x}^2 + \Big(\|\w^{\e}\|_{L^\infty_x}+\|\nabla_x\mathbb{P}u^{\e}\|_{L^\infty_x}+\|\nabla_x^2 \mathbb{P}^{\perp}u^{\e}\|_{L^\infty_x}\Big) \|\w^{\e}\|_{H^2_x}^2.
\enda
\end{align*}
For $I_{\w}^{2}(t)$, the term $\nabla_x\cdot u^{\e}$ also carries at least one derivative, giving
\begin{align*}
\bega
I_{\w}^{2}(t)&\les \bigg|\int_{\Omega} (\nabla_x^2\nabla_x\cdot u^{\e})
\w^{\e} \nabla_x^2\w^{\e} dx\bigg| + \|\nabla_x\nabla_x\cdot u^{\e}\|_{L^\infty_x} \|\w\|_{H^2_x}^2.
\enda
\end{align*}
Integrating by parts, we further have
\begin{align*}
\bega
\bigg|\int_{\Omega} (\nabla_x^2\nabla_x\cdot u^{\e})
\w^{\e} \nabla_x^2\w^{\e} dx\bigg| &\leq \bigg|\int_{\Omega} (\nabla_x\nabla_x\cdot u^{\e})
\nabla_x \w^{\e} \nabla_x^2\w^{\e} dx\bigg| + \bigg|\int_{\Omega} (\nabla_x\nabla_x\cdot u^{\e})
\w^{\e} \nabla_x^3\w^{\e} dx\bigg| \cr 
&\les \|\nabla_x\nabla_x\cdot u^{\e}\|_{L^\infty_x}\Big(\|\w^{\e}\|_{\dot{H}^1_x}\|\w^{\e}\|_{\dot{H}^2_x}+\|\w^{\e}\|_{L^2_x}\|\w^{\e}\|_{\dot{H}^3_x} \Big).
\enda
\end{align*}
Hence, we get
\begin{align}\label{wE4}
\bega
I_{\w}^{2}(t)&\les \|\nabla_x\nabla_x\cdot u^{\e}\|_{L^\infty_x} \|\w\|_{H^2_x}^2 +\|\nabla_x\nabla_x\cdot u^{\e}\|_{L^\infty_x}\|\w^{\e}\|_{L^2_x}\|\w^{\e}\|_{H^3_x} .
\enda
\end{align}
Finally, for the last term in \eqref{wE1}, we integrate by parts only in the first contribution of $\varPi_{\w}^{\e}$ to obtain
\begin{align}\label{wE5}
\bega
\sum_{0 \leq |\al_x| \leq 2}\int_{\Omega}\p^{\al_x}\varPi_{\w}^{\e} \p^{\al_x}\w^{\e} dx &\leq \bigg\| k_B\mathrm{\Theta}^{\e}\nabla_x^{\perp}\ta^{\e}\cdot\nabla_x(\rho^{\e}+\ta^{\e})\bigg\|_{H^1_x}\|\w^{\e}\|_{H^3_x} \cr 
&+ \bigg\|\frac{1}{\e^2}\nabla_x^{\perp}\cdot\bigg(\frac{1}{\mathrm{P}^{\e} }\sum_{j} \p_{x_j} \mathbf{r}_{ij}^{\e}\bigg)\bigg\|_{H^2_x}\|\w^{\e}\|_{H^2_x}.
\enda
\end{align}
Combining \eqref{wE2}, \eqref{wE3}, \eqref{wE4}, and \eqref{wE5} with \eqref{wE1} gives the desired estimate \eqref{wHs}.
\end{proof}

\begin{proposition}\label{P.macro.u}
Let $\Omega=\R^2$, and let \(\w^{\e}\) be the solution of \eqref{weqnnew}. Suppose that the bootstrap assumption \eqref{condition} holds on the interval $[0,T]$. Then, for every $t \in [0,T]$, the following estimates hold: 
\begin{align}\label{pmLinfe}
\bega
\lw\|\nabla_x \mathbb{P}u^{\e}(t)\rw\|_{L^\infty_x} &\leq C \bigg[e^{C\int_0^t \lw(1+\|\w^{\e}(s)\|_{L^\infty_x}\rw) ds}\lw(1+\ln^+\|\w^{\e}_0\|_{H^2_x}^2\rw) \cr 
&+ \int_0^t e^{C\int_s^t \lw(1+\|\w^{\e}(\tau)\|_{L^\infty_x}\rw) d\tau}\mathfrak{O}(s) ds\bigg]\lw(1+\|\w^{\e}(t)\|_{L^\infty_x}\rw), 
\enda
\end{align}
and
\begin{align}\label{wH2e}
\bega
\|\w^{\e}(t)\|_{H^2_x}^2 &\leq C \exp\lw[e^{C\int_0^t \lw(1+\|\w^{\e}(s)\|_{L^\infty_x}\rw) ds}\rw] \times \Big(\|\w^{\e}_0\|_{H^2_x}^2\Big)^{e^{C\int_0^t \lw(1+\|\w^{\e}(s)\|_{L^\infty_x}\rw) ds}} \cr 
&\times \exp \lw[\int_0^t e^{C\int_s^t \lw(1+\|\w^{\e}(\tau)\|_{L^\infty_x}\rw) d\tau}\mathfrak{O}(s) ds \rw],
\enda
\end{align}
where 
\begin{align}\label{GH2esti}
\mathfrak{O}(t) &\les \sum_{0\leq k \leq 2}\Big\|\nabla_x^k\Big((\rho^{\e}+\ta^{\e}),\mathbb{P}^{\perp}u^{\e}\Big)(t)\Big\|_{L^\infty_x}(1+\mathcal{E}_{tot}^{\mathrm{N}}(F^{\e}(t)))+ \kappa^{\frac{1}{2}}(\mathcal{D}_G^{\mathrm{N}}(F^{\e}(t)))^{\frac{1}{2}}.
\end{align}
Here, $\mathcal{E}_{tot}^{\mathrm{N}}(F^{\e}(t))$ and $\mathcal{D}_G^{\mathrm{N}}(F^{\e}(t))$ are defined in Definition \ref{D.al} for $\mathrm{N}\geq 4$.
\end{proposition}

\begin{proof}
In the proof, for brevity, we slightly abuse notation by writing $\mathcal{E}(t)$, and $\mathcal{D}(t)$ for $\mathcal{E}^{\mathrm{N}}(F^{\e}(t))$, and $\mathcal{D}^{\mathrm{N}}(F^{\e}(t))$, respectively.
We first prove \eqref{wH2e}.  
In the energy estimate of the vorticity equation \eqref{wHs}, we apply \eqref{Be-Ma,m} to handle the term \(\|\nabla_x\mathbb{P}u^{\e}\|_{L^\infty_x}\), which yields
\begin{align}\label{omegaH2}
\bega
\frac{d}{dt}\|\w^{\e}\|_{H^2_x}^2
&\les \Big(\big(1 + \ln^+ \|\w^{\e}\|_{H^2_x}\big) \big(1 + \|\w^{\e}\|_{L^\infty_x}\big)+\|\nabla_x \mathbb{P}^{\perp}u^{\e}\|_{L^\infty_x}+\|\nabla_x^2 \mathbb{P}^{\perp}u^{\e}\|_{L^\infty_x}\Big) \|\w^{\e}\|_{H^2_x}^2 \cr 
& +\bigg(\|\nabla_x\nabla_x\cdot u^{\e}\|_{L^\infty_x}\|\w^{\e}\|_{L^2_x}+ \big\| k_B\mathrm{\Theta}^{\e}\nabla_x^{\perp}\ta^{\e}\cdot\nabla_x(\rho^{\e}+\ta^{\e})\big\|_{H^1_x}\bigg)\|\w^{\e}\|_{H^3_x} \cr 
&+ \bigg\|\frac{1}{\e^2}\nabla_x^{\perp}\cdot\bigg(\frac{1}{\mathrm{P}^{\e} }\sum_{j} \p_{x_j} \mathbf{r}_{ij}^{\e}\bigg)\bigg\|_{H^2_x}\|\w^{\e}\|_{H^2_x}.
\enda
\end{align}
While \(\|\w^{\e}(t)\|_{H^2_x}^2 > 1\), dividing both sides by \(\|\w^{\e}(t)\|_{H^2_x}^2\) gives
\begin{align*}
\frac{d}{dt}\lw(1+\ln^+\|\w^{\e}\|_{H^2_x}^2\rw) &\les \lw(1+\|\w^{\e}\|_{L^\infty_x}\rw)\lw(1+\ln^+\|\w^{\e}\|_{H^2_x}^2\rw) +\mathfrak{O}(t),
\end{align*}
where
\begin{align}\label{GH2def}
\bega
\mathfrak{O}(t) &:=\|\nabla_x \mathbb{P}^{\perp}u^{\e}\|_{L^\infty_x}+\|\nabla_x^2 \mathbb{P}^{\perp}u^{\e}\|_{L^\infty_x} +\|\nabla_x\nabla_x\cdot u^{\e}\|_{L^\infty_x}\|\w^{\e}\|_{H^3_x} \cr 
&+ \bigg\| k_B\mathrm{\Theta}^{\e}\nabla_x^{\perp}\ta^{\e}\cdot\nabla_x(\rho^{\e}+\ta^{\e})\bigg\|_{H^1_x}\|\w^{\e}\|_{H^3_x} + \bigg\|\frac{1}{\e^2}\nabla_x^{\perp}\cdot\bigg(\frac{1}{\mathrm{P}^{\e} }\sum_{j} \p_{x_j} \mathbf{r}_{ij}^{\e}\bigg)\bigg\|_{H^2_x}.
\enda
\end{align}
Applying Grönwall’s inequality, we obtain
\begin{align}\label{lnwHs}
\bega
\big(1+\ln^+\|\w^{\e}(t)\|_{H^2_x}^2\big) &\leq e^{C\int_0^t \lw(1+\|\w^{\e}(s)\|_{L^\infty_x}\rw) ds}\lw(1+\ln^+\|\w^{\e}_0\|_{H^2_x}^2\rw) \cr 
&+ C\int_0^t e^{C\int_s^t \lw(1+\|\w^{\e}(\tau)\|_{L^\infty_x}\rw) d\tau}\mathfrak{O}(s) ds.
\enda
\end{align}
This proves \eqref{wH2e}.  
For \eqref{pmLinfe}, substituting \eqref{lnwHs} into \eqref{Be-Ma,m} yields the desired estimate.
Next, we estimate $\mathfrak{O}(t)$ in \eqref{GH2def}.  
Using \eqref{pmw} we observe that
\begin{align*}
\|\w^{\e}\|_{H^3_x} \leq \|\nabla_x\mathbb{P}u^{\e}\|_{H^3_x} \leq \mathcal{E}_{M}^{\frac{1}{2}}, \quad &\mbox{for} \quad \mathrm{N}\geq4,
\end{align*}
where we used that the Hodge decomposition is orthogonal in $H^k$, which implies $\|\nabla_x\mathbb{P}u^{\e}\|_{H^k_x} +\|\nabla_x\mathbb{P}^{\perp}u^{\e}\|_{H^k_x} = \|\nabla_xu^{\e}\|_{H^k_x}$.
For the remaining terms, by \eqref{pTaleq} we estimate
\begin{align*}
\bega
\bigg\| k_B\mathrm{\Theta}^{\e}\nabla_x^{\perp}\ta^{\e}\cdot\nabla_x(\rho^{\e}+\ta^{\e})\bigg\|_{H^1_x} &\leq \bigg\| k_B\mathrm{\Theta}^{\e}\nabla_x^{\perp}\ta^{\e}\bigg\|_{H^1_x}  \sum_{0\leq|\al_x|\leq1}\|\p^{\al_x}\nabla_x(\rho^{\e}+\ta^{\e})\|_{L^\infty_x} \cr &\leq \mathcal{E}_M^{\frac{1}{2}}\sum_{0\leq|\al_x|\leq1}\|\p^{\al_x}\nabla_x(\rho^{\e}+\ta^{\e})\|_{L^\infty_x}.
\enda
\end{align*}
Finally, by the same reasoning as in \eqref{ahatesti}, we also have
\begin{align*}
\bega
\bigg\|\frac{1}{\e^2}\nabla_x^{\perp}\cdot\bigg(\frac{1}{\mathrm{P}^{\e} }\sum_{j} \p_{x_j} \mathbf{r}_{ij}^{\e}\bigg)\bigg\|_{H^2_x} \les \kappa^{\frac{1}{2}}\mathcal{D}_G^{\frac{1}{2}}(t), \quad &\mbox{for} \quad \mathrm{N}\geq4.
\enda
\end{align*}
Combining these estimates in \eqref{GH2def} yields the bound \eqref{GH2esti}.
\end{proof}

\begin{remark}  
(1) The singular dependence on $\kappa$ in the estimate~\eqref{GH2esti}
varies with the number of derivatives $\mathrm{N}$ due to the two terms
$\|\nabla_x^2 \mathbb{P}^{\perp}u^{\e}\|_{L^\infty_x}$
and $\|\w^{\e}\|_{H^3_x}$.
We note that, when $\mathrm{N}=3$, the term
$\|\nabla_x^2 \mathbb{P}^{\perp}u^{\e}\|_{L^\infty_x}$
can be controlled with a singular factor in $\kappa$,
whereas for $\mathrm{N}\ge4$ it can be bounded without such a singular
dependence, as shown in Remark~\ref{Rmk.divuLinf}.
For the term $\|\w^{\e}\|_{H^3_x}$ in~\eqref{omegaH2},
if $\mathrm{N}=3$, the definition of the top-order energy
$\mathcal{E}_{top}$ implies
\begin{align*}
\|\w^{\e}\|_{H^3_x} \leq \|\nabla_x\mathbb{P}u^{\e}\|_{H^3_x} \leq \|\nabla_xu^{\e}\|_{H^3_x} \leq  \kappa^{-\frac{1}{2}}(\mathcal{E}_{top}^{\mathrm{N}}(F^{\e}(t)))^{\frac{1}{2}}.
\end{align*}
If $\mathrm{N}\ge4$, this estimate is no longer singular in $\kappa$; $\|\w^{\e}\|_{H^3_x} \leq (\mathcal{E}_M^{\mathrm{N}}(F^{\e}(t)))^{\frac{1}{2}}$. \\
(2) Since the $H^2_x$-norm of the vorticity grows only in terms of $\|\nabla_x \mathbb{P}u^{\e}\|_{L^\infty_x}$, it can be nearly closed on its own by applying the Biot–Savart law \eqref{Be-Ma,m}, in contrast to Proposition \ref{P.EW}, where the growth depends on $\|\nabla_x (\rho^{\e},u^{\e},\ta^{\e})\|_{L^\infty_x}$.
\end{remark}

\subsubsection{Basic estimates for the vorticity and the specific entropy fluctuation}

In this section, we derive estimates for the vorticity $\w^{\e}$ in $L^p_x$ and for the specific entropy fluctuation $\mathfrak{s}^{\e} = \tfrac{3}{2}\ta^{\e} - \rho^{\e}$ in $L^2_x$.
Note that there is no control of $\varPi_{\w}^{\e}$ in $L^p$ for $1 \leq p < 2$ due to the presence of the microscopic term $\p_{x_j} \mathbf{r}_{ij}^{\e}$.
Hence, to estimate $\w^{\e}$ for $1 \leq p < 2$, we decompose the solution into the contribution from the initial data and a forcing part.
By writing $u^{\e}= \mathbb{P}u^{\e}+\mathbb{P}^{\perp}u^{\e}$ in the term $u^{\e}\cdot\nabla_x \w^{\e}$, the equation \eqref{weqnnew} can be rewritten as
\begin{align}\label{weqnP}
\bega
&\p_t\w^{\e} + \mathbb{P}u^{\e}\cdot \nabla_x \w^{\e} = \bar{\varPi}_{\w}^{\e},
\enda
\end{align}
where 
\hide
\begin{align}\label{barvarPiwdef}
\bega
\bar{\varPi}_{\w}^{\e}(t,x)&:= -\mathbb{P}^{\perp}u^{\e}\cdot \nabla_x \w^{\e} -(\nabla_x\cdot u^{\e})\w^{\e} \cr 
&\quad -k_B\mathrm{\Theta}^{\e}\nabla_x^{\perp}\ta^{\e}\cdot\nabla_x(\rho^{\e}+\ta^{\e})
-\frac{1}{\e^2}\nabla_x^{\perp}\cdot\bigg(\frac{1}{\mathrm{P}^{\e} }\sum_{j} \p_{x_j} \mathbf{r}_{ij}^{\e}\bigg), 
\enda
\end{align}
\unhide
\begin{align*}
\bega
\bar{\varPi}_{\w}^{\e}(t,x)&:=  -\nabla_x\cdot (\mathbb{P}^{\perp}u^{\e}\w^{\e}) -k_B\mathrm{\Theta}^{\e}\nabla_x^{\perp}\ta^{\e}\cdot\nabla_x(\rho^{\e}+\ta^{\e})
-\frac{1}{\e^2}\nabla_x^{\perp}\cdot\big(\frac{1}{\mathrm{P}^{\e} } \nabla_x \cdot  \mathbf{r} ^{\e}\big).
\enda
\end{align*}

\begin{definition}\label{D.wAwB}
For each $\e>0$, we define $\w^{\e}_A$ and $\w^{\e}_B$ as the solutions of the following equations:
\begin{align}
&\p_t\w^{\e}_A + \mathbb{P}u^{\e}\cdot \nabla_x \w^{\e}_A = 0 , \qquad \hspace{3mm} \w^{\e}_A(t,x)|_{t=0} = \w^{\e}_0(x), \label{weqnA} \\ 
&\p_t\w^{\e}_B + \mathbb{P}u^{\e}\cdot \nabla_x \w^{\e}_B = \bar{\varPi}_{\w}^{\e} , \qquad \w^{\e}_B(t,x)|_{t=0} = 0 . \label{weqnB}
\end{align}
We define the associated divergence-free velocity fields via the Biot--Savart law by
\begin{align*}
\bega
u^{\e}_A(t,x) := {\bf K} \ast \w^{\e}_A(t,x), \qquad u^{\e}_B(t,x) := {\bf K} \ast \w^{\e}_B(t,x),
\enda
\end{align*}
where $\ast$ denotes the standard convolution, $(f\ast g)(x)=\int_{\R^d} f(x-y)g(y)dy$ and the kernel ${\bf K}(x):= \frac{1}{2\pi}\frac{x^{\perp}}{|x|^2}$.
\end{definition}

By definition, it follows immediately that 
$\w^{\e}=\w^{\e}_A+\w^{\e}_B$ is a solution of~\eqref{weqnP}, and that 
$\mathbb{P}u^{\e}= u^{\e}_A+u^{\e}_B$.

\begin{lemma}\label{L.w-Lp}
Let $\Omega=\R^2$, and let \(\w^{\e}\) be the solution of \eqref{weqnP}. Suppose that the bootstrap assumption \eqref{condition} holds on the interval $[0,T]$. Then, for every $t \in [0,T]$, the following estimates hold:
\begin{align}
&\|\w^{\e}_A(t)\|_{L^p_x} = \|\w^{\e}_0\|_{L^p_x},  \hspace{2cm} \hspace{2mm} \mbox{for} \quad 1\leq p\leq \infty, \label{w-Lp} \\
&\|\w^{\e}_B(t)\|_{L^p_x}\leq \int_0^t \big\|\bar{\varPi}_{\w}^{\e}(s,\cdot)\big\|_{L^p_x} ds,  \quad \mbox{for} \quad 2\leq p \leq \infty. \label{wABLp}
\end{align}
Moreover, for $\mathrm{N}\geq4$, we have 
\begin{align}\label{barGest}
\bega
\big\|\bar{\varPi}_{\w}^{\e}(t)\big\|_{L^p_x} &\les \Big(
\|\mathbb{P}^{\perp}u^{\e}\|_{L^\infty_x} + \|\nabla_x\cdot u^{\e}\|_{L^\infty_x}+\|\nabla_x(\rho^{\e}+\ta^{\e})\|_{L^\infty_x}\Big)(\mathcal{E}_{tot}^{\mathrm{N}}(F^{\e}(t)))^{\frac{1}{2}}+ \kappa^{\frac{1}{2}}(\mathcal{D}_G^{\mathrm{N}}(F^{\e}(t)))^{\frac{1}{2}}.
\enda
\end{align}
\hide
where
\begin{align*}
\mathcal{C}_{\kappa} := 
\begin{cases}
1,\quad &2\leq p< \infty, \\
\kappa^{-\frac{(4-\mathrm{N})_+}{4}}, \quad &p= \infty.
\end{cases}
\end{align*}
\unhide
\end{lemma}
\begin{proof}
In the proof, for brevity, we slightly abuse notation by writing $\mathcal{E}(t)$, and $\mathcal{D}(t)$ for $\mathcal{E}^{\mathrm{N}}(F^{\e}(t))$, and $\mathcal{D}^{\mathrm{N}}(F^{\e}(t))$, respectively.

Since $\mathbb{P}u^{\e}$ is divergence-free, the estimate~\eqref{w-Lp} follows from standard arguments for the continuity equation.
\hide
We multiply \(p \w |\w|^{p-2}\) (for \(p > 1\)) to \eqref{weqnnew}, and then integrate over \(\Omega\) to derive \eqref{w-Lp}. Similarly, applying the same method to equation \eqref{prtaeqn} yields \eqref{rtaLp}.
Applying the $L^p$ estimate, we have 
\begin{align}\label{NS-Lp}
&\|\rho(t)\|_{L^p}^p \leq \|\rho
(0)\|_{L^p}^p + \int_0^t \|\nabla_x\cdot b(s)\|_{L^\infty_x} \|\rho(s)\|_{L^p}^pds +p \int_0^t \|\varPi_{\w}^{\e}(s)\|_{L^p}\|\rho(s)\|_{L^p}^{p-1} ds.
\end{align}
\begin{align}\label{dw-Lp}
&\frac{d}{dt}\int_{\Omega}|\rho^{\nu}|^pdx -\int_{\Omega}(\nabla_x\cdot b^{\nu}) |\rho^{\nu}|^p dx  +\nu p(p-1)\int_{\Omega} |\nabla_x\rho^{\nu}|^2 |\rho^{\nu}|^{p-2}dx = p\int_{\Omega}\varPi_{\w}^{\e}\rho^{\nu}|\rho^{\nu}|^{p-2} dx.
\end{align}
\begin{align}
&\frac{d}{dt}\int_{\Omega}|\rho^{\nu}|^pdx +\nu p(p-1)\int_{\Omega} |\nabla_x\rho^{\nu}|^2 |\rho^{\nu}|^{p-2}dx \leq \int_{\Omega}(\nabla_x\cdot b^{\nu}) |\rho^{\nu}|^p dx +p\|\varPi_{\w}^{\e}\|_{L^p_x}\|\rho^{\nu}\|_{L^p_x}^{p-1}.
\end{align}
Dividing each side by $p\|\rho^{\nu}\|_{L^p_x}^{p-1}$ gives 
\begin{align}
&\frac{d}{dt}\|\rho^{\nu}\|_{L^p} \leq \frac{1}{p}\|\nabla_x\cdot b^{\nu}\|_{L^\infty_x}\|\rho^{\nu}\|_{L^p} +\|\varPi_{\w}^{\e}\|_{L^p_x}.
\end{align}
\begin{align}
\|\rho^{\nu}(t)\|_{L^p}&\leq e^{\int_0^t \|\nabla_x\cdot b^{\nu}(s)\|_{L^\infty_x}ds}\|\rho^{\nu}_0\|_{L^p} + \int_0^t e^{\int_s^t \|\nabla_x\cdot b^{\nu}(\tau)\|_{L^\infty_x}d\tau} \|\varPi_{\w}^{\e}(s)\|_{L^p_x}ds
\end{align}
We use the form \eqref{weqnnew}: 
\begin{align*}
\bega
&\p_t\w^{\e} + u^{\e}\cdot \nabla_x \w^{\e} + (\nabla_x\cdot u^{\e})\w^{\e} =\varPi_{\w}^{\e},
\enda
\end{align*}
We denote the second line by $\varPi_{\w}$. 
We multiply \(p \w^{\e} |\w^{\e}|^{p-2}\) (for \(p > 1\)) to \eqref{weqnnew}, and then integrate over \(\Omega\) to derive
\begin{align*}
&\frac{d}{dt}\int_{\Omega}|\w^{\e}|^pdx = -(p-1)\int_{\Omega}(\nabla_x\cdot u^{\e}) |\w^{\e}|^p dx - p\int_{\Omega}\varPi_{\w}^{\e}\w^{\e}|\w^{\e}|^{p-2} dx.
\end{align*}
Thus we have
\begin{align*}
&\frac{d}{dt}\|\w^{\e}\|_{L^p_x}^p \leq (p-1)\|\nabla_x\cdot u^{\e}\|_{L^\infty_x}\|\w^{\e}\|_{L^p_x}^p +  Cp\|\varPi_{\w}^{\e}\|_{L^p}\|\w^{\e}\|_{L^p_x}^{p-1}
\end{align*}
(Thus uniform $L^p_x$ estimate hold for infinite energy case.)
\unhide
\hide we first observe that $\w^{\e}_A$ satisfies the transport equation
\begin{align}\label{wAeqnu}
\bega
\p_t\w^{\e}_A + \mathbb{P}u^{\e}\cdot \nabla_x \w^{\e}_A = 0 , \qquad \w^{\e}_A(t,x)|_{t=0} = \w^{\e}_0(x).
\enda
\end{align}
(Note that $\nabla_x^{\perp}\cdot\mathbb{P}u^{\e} \neq \w^{\e}_A$.)  
Since $\mathbb{P}u^{\e}$ is divergence-free, the $L^p_x$–norm of $\w^{\e}_A$ is conserved. This proves the first inequality in \eqref{wABLp}. \unhide  
For the proof of \eqref{wABLp}, taking the $L^p_x$–norm of $\w^{\e}_B$ and applying Minkowski’s integral inequality yields the inequality in \eqref{wABLp}.
\\ 
Next, we estimate $\bar{\varPi}_{\w}^{\e}(t)$ in $L^p_x$ for $2 \leq p \leq \infty$. 
From the definition~\eqref{barvarPidef}, we obtain
\begin{align*}
\bega
\big\|\bar{\varPi}_{\w}^{\e}(t)\big\|_{L^p_x} &\leq \|\mathbb{P}^{\perp}u^{\e}\|_{L^\infty_x}\|\nabla_x \w^{\e}\|_{L^p_x} + \|\nabla_x\cdot u^{\e}\|_{L^\infty_x}\|\w^{\e}\|_{L^p_x} + \mathcal{E}_M^{\frac{1}{2}}\|\nabla_x(\rho^{\e}+\ta^{\e})\|_{L^\infty_x} + \kappa^{\frac{1}{2}}\mathcal{D}_G^{\frac{1}{2}}(t),
\enda
\end{align*}
where we have used 
$\nabla_x^{\perp}\!\cdot\nabla_x(\rho^{\e}+\theta^{\e}) = 0$. 
For the term $\|\nabla_x \w^{\e}\|_{L^p_x}$, we apply the interpolation inequality~\eqref{Ga-Ni} for $2\le p<\infty$ and Agmon’s inequality~\eqref{Agmon} for $p=\infty$ to have
\begin{align*}
\bega
\|\nabla_x \w^{\e}\|_{L^p_x} \leq \begin{cases} \|\nabla_x^2 \w^{\e}\|_{L^2_x}^{1-\frac{2}{p}}\|\nabla_x \w^{\e}\|_{L^2_x}^{\frac{2}{p}}, \quad &2\leq p<\infty, \\ 
\|\nabla_x \w^{\e}\|_{L^2_x}^{\frac{1}{2}}\|\nabla_x \w^{\e}\|_{H^2_x}^{\frac{1}{2}}, \quad &p=\infty.\end{cases}
\enda
\end{align*}
Moreover, using~\eqref{pmw} together with the definitions of $\mathcal{E}_M$ and $\mathcal{E}_{top}$, we have
\begin{align*}
\bega
\|\nabla_x \w^{\e}\|_{H^2_x} \les \|\nabla_x \mathbb{P}u^{\e}\|_{H^3_x} \les
\mathcal{E}_M^{1/2}, \quad &\mbox{for} \quad \mathrm{N}\geq 4.
\enda
\end{align*}
This yields the desired estimate~\eqref{barGest}.
\end{proof}

\begin{lemma}\label{L.rho-ta}
Suppose that the bootstrap assumption \eqref{condition} holds up to time $t \in [0,T]$. 
Let $\rho^{\e}-\frac{3}{2}\ta^{\e}$ solve equation \eqref{rtaeqnnew}. Then, for $t \in [0,T]$, we have
\begin{align*}
\bega
\|\mathfrak{s}^{\e}(t)\|_{L^2_x} &\leq e^{\frac{1}{2}\int_0^t \|(\nabla_x\cdot u^{\e})(s)\|_{L^\infty_x}ds}\|\mathfrak{s}^{\e}_0\|_{L^2_x}  + \int_0^t e^{\frac{1}{2}\int_s^t \|(\nabla_x\cdot u^{\e})(\tau)\|_{L^\infty_x}d\tau} \kappa^{\frac{1}{2}} (\mathcal{D}_{G}^{\mathrm{N}}(F^{\e}(s)))^{\frac{1}{2}} ds.
\enda
\end{align*}
\end{lemma}
\begin{proof}
Taking the $L^2$ inner product of \eqref{rtaeqnnew} with $\mathfrak{s}^{\e}$ yields the energy inequality
\begin{align*}
\bega
&\frac{1}{2}\frac{d}{dt}\|\mathfrak{s}^{\e}\|_{L^2_x}^2 \leq \frac{1}{2}\|\nabla_x\cdot u^{\e}\|_{L^\infty_x}\|\mathfrak{s}^{\e}\|_{L^2_x}^2 + \big\|\varPi_{\mathfrak{s}}^{\e}\big\|_{L^2_x}\|\mathfrak{s}^{\e}\|_{L^2_x}.
\enda
\end{align*}
Dividing both sides by $\|\mathfrak{s}^{\e}\|_{L^2_x}$, and using $(\mathcal{B}_2)$ in \eqref{condition} together with $\eqref{ABGscale}_1$, we have $\|\varPi_{\mathfrak{s}}^{\e}\|_{L^2_x} \les \kappa^{\frac{1}{2}}\mathcal{D}_G^{\frac{1}{2}} $. Therefore, an application of Grönwall’s inequality yields the desired estimate.
\end{proof}

\hide
\subsection{Vorticity equation in 3D}

\begin{lemma}\label{} In 3D, $\w$ satisfies 
\Be\bega\label{weqn3D}
&\p_t \w +u\cdot\nabla_x \w -\eta_0 k_B^{\frac 1 2} \kappa|\mathrm{\Theta}^{\e}|^{\frac 1 2}  \Delta_x \w = \w\cdot\nabla_x u - \w(\nabla_x\cdot u) +V^{\w}.
\enda
\Ee
\end{lemma}
\begin{proof}
Since the following identity holds $
\nabla_x \times (u\cdot\nabla_x u) 
 = u\cdot \nabla_x \w + \w(\nabla_x\cdot u)-\w\cdot\nabla_x u$ ,
we derive \eqref{weqn3D}.
\end{proof}

When we try $H^m_x$ estimate on the 3D vorticity equation, the only problem is: for $m>d/2+1$
\begin{align}
\bega
\int_{\Omega} \p^{\al}(\w\cdot\nabla_x m)  \p^{\al}\w dx &\les (\|\w\|_{L^\infty_x}\|\nabla_xu\|_{H^2_x}+\|\nabla_xu\|_{L^\infty_x}\|\w\|_{H^2_x})\|\w\|_{H^m_x} \les \|\w\|_{H^m_x}^3 \cr 
\int_{\Omega} \p^{\al}(\w(\nabla_x\cdot u)) \p^{\al}\w dx &\les \|\nabla_x\cdot u\|_{H^m_x}\|\w\|_{H^m_x}^2
\enda
\end{align}
where we used \eqref{uvHk} and \eqref{puw}. 
\unhide

\section{Uniform global estimate in 2D}\label{Sec.2D}

\hide
Since the theorem presented here remains valid even when $(\rho_0,u_0,\ta_0,\AC{\P}F_0)\in H^{-\mathfrak{j}}$, we introduce the setting of mollifiers in order to clarify the assumptions on the initial data.
We mollify the initial data with respect to the the scale $\mathfrak{d}(\e)$, which is a function depending on $\e$. 
\begin{align*}
\bega
\rho^{\e}_0 := \varphi^{\mathfrak{d}(\e)} \ast \rho_0, \qquad \ta^{\e}_0 := \varphi^{\mathfrak{d}(\e)} \ast \ta_0, \qquad u^{\e}_0(x) := \varphi^{\mathfrak{d}(\e)} \ast \widetilde{u}_0(x) + \bar{u}_0(x) + \varphi^{\mathfrak{d}(\e)} \ast \mathbb{P}^{\perp}u_0,
\enda
\end{align*}
where $\varphi^{\mathfrak{d}(\e)}$ is a standard mollifier satisfying  
$\int_{\R^d} \varphi(x) dx=1, \quad \varphi(x) \geq0$, $\varphi(x) \in C_c^{\infty}(\R^d),$ $\varphi^{\mathfrak{d}(\e)}(x) = \frac{1}{\mathfrak{d}(\e)^2}\varphi(\frac{x}{\mathfrak{d}(\e)})$. Here, $\bar{u}(x)$ is non-zero only for the infinite velocity energy case as defined in \eqref{barudef}.
For instance, when $(\rho_0,u_0,\ta_0)$ belong to $H^{-\mathfrak{j}}_x$, their mollification becomes singular at a rate depending on $\mathfrak{d}$: $\|(\rho^{\e}_0,u^{\e}_0,\ta^{\e}_0)\|_{H^{\mathrm{N}}_x} \leq C\mathfrak{d}^{-\mathrm{N}-\mathfrak{j}}$. 
For the general initial data, we assume  
\begin{align*}
\mathcal{E}_{tot}(F_0^{\e}) = \mathcal{E}_G(F_0^{\e})+ \mathcal{E}_{M}(F_0^{\e})+ \mathcal{E}_{top}(F_0^{\e}) \les \mathfrak{d}^{-\mathfrak{p}},\quad \mbox{for any} \quad \mathfrak{p}\geq0. 
\end{align*}
With a suitable choice of the parameters $(\e,\kappa(\e),\mathfrak{d}(\e))$, we show that even starting from $H^{-\mathfrak{j}}_x$ initial data, the microscopic component $F^{\e}-M^{\e}$ converges to zero by the following sense:
\begin{align*}
\bega 
\sum_{0\leq|\al|\leq\mathrm{N}}& \frac{1}{\eps^2}\int_0^t\int_{\Omega \times \R^3} \nu\lw(\frac{v-\mathrm{U}^{\e}}{\sqrt{\mathrm{\Theta}^{\e}}}\rw) \frac{|\p^{\al}(F^{\e}-M^{\e})|^2}{M^{\e}} dvdxds \leq C\e^2 \kappa^{1-}\to 0, \quad \mbox{as} \quad \e\to0.
\enda
\end{align*}
(For the top order $|\al|=\mathrm{N}+1$, we have one $\kappa$ worse scale.)


The two theorems presented in this section holds in both settings, depending on whether we use purely spatial derivatives \eqref{caseA} or space-time derivatives \eqref{caseB}, and whether we consider the finite velocity energy case \eqref{caseEC} or the infinite velocity energy case \eqref{caseECX}. In the case of space–time derivatives $\e^{\mathfrak{n}}\p_t$ with $\mathfrak{n}<1$, the proof requires a different approach; this will be discussed later in the discussion part.
\unhide

\hide
----------------------------------
We take an initial data satisfying 
\begin{align}
\bega
\rho_0 \in L^2(\R^2) , \qquad  u_0-\bar{u} \in L^2(\R^2), \qquad \ta_0 \in L^2(\R^2), \qquad \frac{1}{\e^2} \int_{\Omega\times\R^3}|\AC{\P}F(0)|^2|M(0)|^{-1}dvdx \leq C, 
\enda
\end{align}
Our theorem hold for various kind of initial data satisfying  
\begin{itemize}
\item Sobolev space: $\w_0\in H^k_x$ for any $k\in \mathbb{Z}^+\cup\{0\}$.
\item Yudovich class: $\w_0\in L^{\infty}_x$. 
\item Below Yudovich class: $\w_0\in L^p_x$, for any $1\leq p <\infty$.
\item Space of the Radon measures on $\R^d$: $\w_0 \in \mathcal{M}(\R^2)$.
\end{itemize}
\unhide


In this section, we present the main theorem.
Before stating the main theorem, we introduce a growth function that describes the growth of the total energy and dissipation.


\begin{definition}\label{D.Psi-def}
Choose a sufficiently large constant $C>0$, depending (in principle) on the universal constant appearing in Proposition~\ref{P.macro.u}. For the given initial data $\w^{\e}_0$ and $\mathfrak{s}^{\e}_0=\frac{3}{2}\ta^{\e}_0-\rho^{\e}_0$, we define the associated growth function:
\begin{align}\label{Psi-def}
\bega
&\Psi[\w^{\e}_0,\mathfrak{s}^{\e}_0](t):= 1 +  t \big(1+\|\nabla_x\mathfrak{s}^{\e}_0\|_{L^\infty_x}\big)^2 
\lw(1+ \|\w^{\e}_0\|_{H^2_x}\rw) 
\exp\left(e^{ C(1+\|\w^{\e}_0\|_{L^\infty_x})t}\right).
\enda
\end{align}

For brevity, we often denote $\Psi(t) := \Psi[\w^{\e}_0,\mathfrak{s}^{\e}_0](t) $.
\end{definition}


The theorem presented below holds for the purely spatial derivative case~\eqref{caseA} as well as for the space--time derivative case~\eqref{caseB} with arbitrary $\mathfrak{n}>1$. 
The slight differences arising in the case $\mathfrak{n}\leq1$ will be discussed in a remark. We first state the following two bootstrap assumptions.

\begin{assumption}[$\mathrm{N}$-Bootstrap Assumption]\label{Boot}
We introduce the following bootstrap assumption:
\begin{enumerate}
\item[(1)] For $0<\mathfrak{y}\ll1$,
\begin{align}\label{Boot1}
\sup_{0\leq t \leq T_*}
\big(\e^{\mathfrak{y}}+\e\kappa^{-\frac{1}{2}}\big)\big(1+t\big)\bigg(\mathcal{E}^N_{tot}(F^\e(t))+\int_0^t\mathcal{D}^N_{tot}(F^\e(s))ds\bigg)  < \infty.
\end{align}
The energy and dissipation functionals are defined in~\eqref{EDtotdef}.

\item[(2)] For $\mu=M_{[1,0,1]}$ and $\bar{\mu}=M_{[1,0,1-c_0/2]}$ with $0<c_0 \ll 1$,
\begin{align}\label{Boot2}
&\sum_{0\leq|\al|\leq\mathrm{N}+1}\kappa^{\frac{(|\al|-\mathrm{N})_+}{2}}\bigg\|w\p^{\al}\bigg(\frac{F^{\e}-\mu}{\sqrt{\mu}}\bigg)\bigg\|_{L^\infty_tL^\infty_{x,v}} + \sum_{0\leq|\al|\leq\mathrm{N}-2}\bigg\|w\p^{\al}\bigg(\frac{F^{\e}-M^{\e}}{\sqrt{\bar{\mu}}}\bigg)\bigg\|_{L^\infty_tL^\infty_{x,v}} \les \frac{1}{\e\kappa^2},
\end{align}
where $w(v)=e^{c_1|v|^2}$ and the constant $c_1$ is defined in \eqref{hdef}.

\hide
We expect that $\|\p^{\al}h^{\e}\|_{L^2_tL^\infty_{x,v}}$ will blow up as follows: (See \eqrefhLinf})
\begin{align}
\sum_{0\leq|\al|\leq\mathrm{N}+1}\kappa^{\frac{(|\al|-\mathrm{N})_+}{2}}\|\p^{\al}h^{\e}\|_{L^\infty_tL^\infty_{x,v}}  \les\! \!\!\!\sum_{0\leq|\al|\leq\mathrm{N}+1}\!\!\kappa^{\frac{(|\al|-\mathrm{N})_+}{2}}\|\p^{\al}h^{\e}_0\|_{L^\infty_{x,v}} +\frac{1}{(\e\kappa)^{\frac{d}{2}}}\|\big(1+\mathcal{E}_{tot}^{\mathrm{N}}(F^{\e}(t))\big)\|_{L^\infty_t}^{\frac{1}{2}} ,
\end{align}
Expect $\frac{1}{(\e\kappa)^{\frac{d}{2}}}\kappa^{-\frac{1}{2}}$ blow-up
\unhide
\end{enumerate}
\end{assumption}

\hide
(2-Boots) (Bootstrap은 $F^{\e}(t)$에 대해 얘기해주기 때문에, 미분한 값에 대해서밖에 얘기가 안됩니다.)
For any $T_*>1$ the following assumption hold
\begin{align}\label{Boot2}
\sum_{0\leq |\al_x|\leq \mathrm{N}} \sup_{0\leq t \leq T_*}\bigg\|\frac{\p^{\al_x}(F^{\e}(t)-\mu)}{\sqrt{\tilde{\mu}}}\bigg\|_{L^\infty_tL^\infty_{x,v}} \ll 1, \quad \sum_{0\leq |\al_x|\leq \mathrm{N}-3} \sup_{0\leq t \leq T_*} \bigg\|\frac{\p^{\al_x}(F^{\e}-M^{\e})(t)}{\sqrt{\tilde{\mu}}}\bigg\|_{L^\infty_tL^\infty_{x,v}} \ll 1 
\end{align}
where $\tilde{\mu}=M_{[1,0,1-c_0]}$ for any $0<c_0\ll 1$. \\
(2-initial) 이니셜 데이터에 대해 준다면: 
We assume that the initial velocity tail is controlled in a weak $L^\infty_{x,v}$ sense:
\begin{align}\label{ini2D.tail}
\bega
&\bigg\|\frac{\e^{1-}}{\e}\frac{(F|_{t=0}-\mu)}{\sqrt{\tilde{\mu}}}\bigg\|_{L^\infty_{x,v}}\les 1, \quad \mbox{and} \quad \bigg\|\frac{\e^{1-}}{\e}\frac{(F|_{t=0}-M_{[\rho_0,u_0,\ta_0]})}{\sqrt{\tilde{\mu}}}\bigg \|_{L^\infty_{x,v}}\les 1,
\enda
\end{align}
where $\tilde{\mu}=M_{[1,0,1-c_0]}$ for any $0<c_0\ll 1$.
\unhide

\hide
{\color{red}[다른 것은 바꿀 필요가 없고, 아래만 수정요:

\begin{itemize} 
    \item Since you choose difference $N$ for other cases, let's indicate N in the energy and the dissipation ONLY in the Theorem/Proposition/Lemma/Remark. In the proof, we keep the current form by abusing the notation (omit $N,\e$, etc)
\[
\mathcal{E}^N_{tot}(F^\e(t)),  \mathcal{D}^N_{tot}(F^\e(t)) ,
\mathcal{E}^N_{G}(F^\e(t)),  \mathcal{D}^N_{G}(F^\e(t))  
\text{등등}
\]


\begin{cases}
\mbox{If } \quad \mathrm{N}=3, \quad q = \frac{7}{2}- , \\[2mm]
\mbox{If } \quad \mathrm{N}\geq 4, \quad q = \frac{1}{2}- .
\end{cases}
\end{align*}
\end{itemize}

]}
\unhide

\begin{theorem}\label{T.2D.global}
Let $\Omega=\R^2$. Suppose that the $\mathrm{N}$-Bootstrap Assumptions~\eqref{Boot1} and~\eqref{Boot2} hold on $t \in [0,T_*]$ for some $\mathrm{N}\geq4$. 
Then the solution $F^{\e}$ of the Boltzmann equation~\eqref{BE} with $\kappa=\e^q$ for some $0<q<2$ satisfies the following uniform-in-$\e$ estimates for all $0 \leq t \leq T_*$:
\begin{itemize}
\item 
The total energy and dissipation defined in~\eqref{EDtotdef} satisfy
\begin{align}\label{EDbdd2D}
\bega
\mathcal{E}_{tot}^{\mathrm{N}}(F^\e (t))+\int_0^t\mathcal{D}_{tot}^{\mathrm{N}}(F^\e (s))ds&\leq  C\exp\Big(C\Psi[\w^{\e}_0,\mathfrak{s}^{\e}_0](t) \Big)(\mathcal{E}_{tot}^{\mathrm{N}}(F^\e (0) )+\e t).
\enda
\end{align}
where $\Psi[\w^{\e}_0,\mathfrak{s}^{\e}_0](t)$ is defined in Definition~\ref{D.Psi-def}.   

\hide \begin{align}\label{Crutadef}
\bega
\mathcal{E}_{tot}(0)&\les  \|(\rho^{\e}_0,u^{\e}_0-\bar{u},\ta^{\e}_0)\|_{L^2_x}^2 + \sum_{1\leq|\al|\leq \mathrm{N}+1}\kappa^{(|\al|-\mathrm{N})_+}\|\p^{\al}(\rho^{\e}_0,u^{\e}_0,\ta^{\e}_0)\|_{L^2_x}^2 \cr 
&+\sum_{0\leq|\al|\leq\mathrm{N}+1}\frac{\kappa^{(|\al|-\mathrm{N})_+}}{\e^2} \int_{\Omega\times\R^3}\frac{|\p^{\al}\AC{\P}F^{\e}(0)|^2}{|M^{\e}(0)|}dvdx. 
\enda
\end{align}
\unhide

\item  
The acoustic/compressible components are uniformly controlled uniformly-in-$\e$ for any $s\in [0,\mathrm{N}-1)$:
\hide
\begin{align}
\bega
\Big\|\Big((\rho^{\e}+\ta^{\e}),\mathbb{P}^{\perp}u^{\e}\Big)\Big\|_{L^{4+8\delta}_T\dot{B}_{2+\frac{1}{\delta},1}^{s-\frac{3}{4+8\delta}}} 
&\les \e^{\frac{1}{4}-}\kappa^{-\frac{1}{2}}  \int_0^{T}\Big(\mathcal{E}_{tot}(s) +\kappa^{\frac{1}{2}}\mathcal{D}_{tot}^{\frac{1}{2}}(s)\Big) ds \to 0, 
\enda
\end{align}
\unhide
\begin{align}\label{incomp.2D}
\bega
\Big\|\Big((\rho^{\e}+\ta^{\e}),\mathbb{P}^{\perp}u^{\e}\Big)&\Big\|_{L^r_T\dot{B}_{p,1}^{s+2(\frac{1}{p}-\frac{1}{2})+\frac{1}{r}}} \leq C \e^{\frac{1}{r}}\Big\|\Big((\rho^{\e}_0+\ta^{\e}_0),\mathbb{P}^{\perp}u^{\e}_0\Big)\Big\|_{\dot{B}_{2,1}^s} \cr 
&+ C\e^{\frac{1}{r}}\int_0^{T}\Big(\mathcal{E}_{tot}^{\mathrm{N}}(F^\e (s)) +\kappa^{\frac{1}{2}} \sqrt{\mathcal{D}_{tot}^{\mathrm{N}}(F^\e (s))}\Big) ds.
\enda
\end{align}
Here $2\leq p\leq \infty$ and $\frac{2}{r} \leq \frac{1}{2}-\frac{1}{p}$, while the Leray projector $\mathbb{P}$ and the Besov norms are defined in~\eqref{LerayPdef} and~\eqref{Besovdef}, respectively.

\item The microscopic energy at lower regularity is controlled by
\begin{align}\label{Grefine.2D}
\bega
\mathcal{E}^{\mathrm{N}-2}_G( F^\e (t) ) +\frac{\sigma_L}{C}\int_0^t\mathcal{D}^{\mathrm{N}-2}_G(F^\e (s) ) ds &\leq \mathcal{E}^{\mathrm{N}-2}_G(F^\e (0)) \cr 
&+ \kappa^{\frac{1}{2}}\mathcal P (\mathcal{E}_{tot}^{\mathrm{N}}(F^\e (t)),\mathcal{D}_{tot}^{\mathrm{N}}(F^\e (t)) ),
\enda
\end{align}
where $\mathcal P (a,b)$ is a polynomial in $a$ and $b$.

\hide
\begin{align}\label{k1/2-}
\bega
\int_0^T \bigg(\Big\|\Big((\rho^{\e}+\ta^{\e})(t),\mathbb{P}^{\perp}u^{\e}(t)\Big)\Big\|_{L^\infty_x}\mathcal{E}_{tot}^{\frac{1}{2}}(t) \bigg)dt + \kappa^{\frac{1}{2}}\Big(\int_0^T\mathcal{D}_G(t)dt\Big)^{\frac{1}{2}} \leq  C \kappa^{\frac{1}{2}-}.
\enda
\end{align}
\unhide

\end{itemize}
\end{theorem}

\begin{remark}
(1) By the definition of $\mathcal{D}_{tot}$ in \eqref{N-EDdef}, the estimate \eqref{EDbdd2D} implies that the microscopic component is controlled in the limit $\e \to 0$ for all $0 \leq |\al| \leq \mathrm{N}+1$:
\begin{align}\label{Gto0}
\bega 
\frac{\kappa^{(|\al|-\mathrm{N})_+}}{\eps^4\kappa}&\int_0^T\int_{\Omega \times \R^3} \nu \frac{|\p^{\al}\AC{\P}F^{\e}(t)|^2}{M^{\e}}  dvdxdt &\leq C\exp\bigg(C\Psi[\w^{\e}_0,\mathfrak{s}^{\e}_0](t) \bigg)(\mathcal{E}_{tot}(0)+\e t).
\enda
\end{align}
By the bootstrap assumption \eqref{Boot1}, for $0 \leq |\al| \leq \mathrm{N}$, we obtain
\begin{align}\label{Gto02}
\bega 
\frac{1}{\e^4}&\int_0^T\int_{\Omega \times \R^3} \nu \frac{|\p^{\al}\AC{\P}F^{\e}(t)|^2}{M^{\e}} dvdxdt &\leq C_T \kappa^{\frac{1}{2}-}, \qquad \Big\|\frac{1}{\e^2}\p^{\al} \mathbf{r}_{ij}^{\e}\Big\|_{L^2_TL^2_x} \leq C_T \kappa^{\frac{1}{2}-}.
\enda
\end{align}
\\
(2) In perspective of scale, the bootstrap assumption \eqref{Boot2} is easily verified because of the Mach number scale $|F^{\e}-M^{\e}|\sim O(\e)$. \\
\hide (3) The estimate \eqref{EDbdd2D} does not provide boundedness in a limiting process when $(\rho_0,u_0,\ta_0,\AC{\P}F_0)$ is not uniformly bounded in $H^3_x$. However, the uniform estimates in Lemma \ref{L.unif} remain valid in the hydrodynamic limit.
\\\unhide
(3) For the case $\mathcal{E}_{tot}(0) < \infty$, inequality \eqref{EDbdd2D} ensures that both the energy and the dissipation remain uniformly bounded in $\e$. Moreover, if $\|\nabla_x(\rho^{\e}_0-\tfrac{3}{2}\ta^{\e}_0)\|_{L^\infty_x} \geq C$, then the right-hand side of \eqref{EDbdd2D} exhibits triple-exponential growth. On the other hand, if $\|\nabla_x(\rho^{\e}_0-\tfrac{3}{2}\ta^{\e}_0)\|_{L^\infty_x} \leq \kappa$, then the right-hand side of \eqref{EDbdd2D} grows at most double-exponentially, which is known to be the optimal growth rate \cite{Denisov}.  
\end{remark}

\begin{remark}\label{Rmk.Grefine}
(1) For $0 \leq |\al| \leq \mathrm{N}-2$, the inequality \eqref{Grefine.2D} implies
\begin{align*}
\sup_{\e>0}\bigg\|\frac{1}{\e} \frac{\p^{\al}\AC{\P}F^{\e}}{\sqrt{M^{\e}}}\bigg|_{t=0}\bigg\|_{L^2_{x,v}} <\infty \quad \Rightarrow  \quad \frac{1}{\e^2}\|\p^{\al} \mathbf{r}_{ij}^{\e}\|_{L^2_TL^2_x} \leq \kappa^{\frac{1}{2}}\Big(\int_0^T\mathcal{D}_G^{\al}(t)dt\Big)^{\frac{1}{2}} 
\leq C_T \kappa^{\frac{1}{2}}.
\end{align*}
Comparing this with \eqref{Gto02}, we see that this estimate yields exactly the same convergence rate as that arising from the inviscid limit of the incompressible Navier--Stokes equations. \\
(2) When deriving the energy estimate for $\p^{\al} \AC{\P}F^{\e}$, in contrast to the energy estimate for $\p^{\al} F^{\e}$, the nonlinear term appears in the form $\int_{\Omega }\p_{x_i}u^{\e}_j\mathcal{V}_{\ell}dx$ rather than $\|\nabla_xu^{\e}\|_{L^\infty_x}\mathcal{E}+\int_{\Omega }\p_{x_i}u^{\e}_j\mathcal{V}_{\ell}dx$
see Proposition~\ref{P.G.Energy} and Proposition~\ref{P.F.Energy}.
Moreover, the quantity $\AC{\P}F^{\e}$ contained in $\mathcal{V}_{\ell}$, after substituting the decomposition $F^{\e} = M^{\e} + \e \sqrt{\bar{\mu}}g^{\e}$ produces transport terms involving $M^{\e}$ on the right-hand side. (See \eqref{h2eqn} in Section \ref{Sec.V.L}.)
Since this requires an $L^\infty_{x,v}$ estimate, the energy estimate for $\p^{\al} \AC{\P}F^{\e}$ can be further refined up to when
$\|\p^{\al}\nabla_x u^{\e}\|_{L^\infty_x}$ remains controllable.
This explains why the maximal number of derivatives in~\eqref{Grefine.2D} is restricted to $\mathrm{N}-2$.
\end{remark}

\begin{remark}\label{Rmk.prepared}
We emphasize that Theorem~\ref{T.2D.global} holds even without using any time derivatives. 
Here we briefly discuss the relation between well-preparedness and the scaling in front of the time derivative $\p_{\tilde t}=\e^{\mathfrak n}\p_t$ when the space--time derivative case~\eqref{caseB} is considered with arbitrary $\mathfrak n>0$.
In the low Mach number limit for the compressible Euler equations, 
the following initial data are not regarded as well-prepared with respect to the scaling ($\p_{\tilde{t}}=\e\p_t$):
\begin{align*}
\|\e\p_t(\rho^{\e},u^{\e},\ta^{\e})|_{t=0}\| <\infty \quad \sim \quad \big\|\big(\nabla_x(\rho^{\e}+\ta^{\e}), \nabla_x\cdot u^{\e}\big)|_{t=0} \big\| <\infty.
\end{align*}
In contrast, for the incompressible Euler limit from the Boltzmann equation~\eqref{BE}, 
the initial data are well-prepared since the microscopic component vanishes faster than the Mach number. 
More precisely,
\begin{align*}
\frac{1}{\e}\Big\|(\e\p_t)F^{\e}|M^{\e}|^{-\frac{1}{2}}|_{t=0}\Big\| <\infty \quad &\sim \quad 
\frac{1}{\e\kappa}\Big\|\frac{1}{\e}\mathcal{L}(\AC{\P}F^{\e})|M^{\e}|^{-\frac{1}{2}}|_{t=0}\Big\| <\infty.
\end{align*}
\hide
Decomposing $F^{\e}=M^{\e}+\AC{\P}F^{\e}$ in the collision operator yield
\begin{align*}
\p_t F^{\e} + \frac{v}{\e} \cdot \nabla_x F^{\e} + \frac{1}{\e^2\kappa}\mathcal{L}(\AC{\P}F^{\e})= \frac{1}{\e^2\kappa} \mathcal{N}(\AC{\P}F^{\e},\AC{\P}F^{\e}).
\end{align*}
\unhide
While the scaled time derivative $\e \p_t$ is compatible with $\e\kappa$--well-preparedness, 
the use of the scaled time derivative $\e^2\kappa \p_t$ implies that, 
at least from the viewpoint of scaling, the data are not well-prepared.
\end{remark}


\begin{remark}
(1) If we use the space--time derivative case $\p_{\tilde{t}}=\e^{\mathfrak{n}}\p_t$ in \eqref{caseB}, then for $\mathfrak{n}=1$ the proof of Theorem~\ref{T.2D.global} holds when $\mathrm{N}\geq 4$.. 
This is because the argument requires an estimate of the form 
$\|\nabla_x^2\mathbb{P}^{\perp}u^{\e}\|_{L^\infty_x}$ appearing in \eqref{Bmadef}, which cannot be obtained from the Strichartz estimate when $\mathrm{N}=3$. \\
(2) For $0<\mathfrak{n}<1$, the growth rate in \eqref{EDbdd2D} changes, since the use of the moments equation \eqref{locconNew} does not provide an improved scaling. 
The main reason is that the estimates of $\mathcal{Z}_{top}^{time}$ and $\mathcal{Z}_W^{time}$ are modified, and consequently the growth of the energy and dissipation in \eqref{EBE0} can no longer be controlled solely by $\|\nabla_x(\rho,u,\ta)\|_{L^\infty_x}$. \\
(3) The maximum regularity at which we can establish the incompressibility condition is the space $L^r_T\dot{B}_{p,1}^{s+2(\frac{1}{p}-\frac{1}{2})+\frac{1}{r}}$. For $p \to \infty$ and $s<\mathrm{N}$, the maximal regularity is close to $\mathrm{N}-\frac{d+1}{4}$, which is relatively low compared to the top order regularity $H^{\mathrm{N}+1}_x$. Therefore, it is essential to reduce the regularity requirements on $(\rho^{\e}+\ta^{\e})$ and $\mathbb{P}^{\perp}u^{\e}$ in the terms that need to be controlled. 

\hide
(1) Depending on the order of the highest derivative, the admissible range of $q$ in the relation $\kappa = \e^q$ changes as follows: \textcolor{red}{Need to revise}
\begin{align*}
\begin{cases}
\mbox{If } \quad \mathrm{N}=3, \quad q = \frac{2}{7}- , \\[2mm]
\mbox{If } \quad \mathrm{N}\geq 4, \quad q = \frac{1}{2}- .
\end{cases}
\end{align*}
The main reason is the term $\|\nabla_x\nabla_x\cdot u^{\e}\|_{L^\infty_x}\|\w^{\e}\|_{L^2_x}\|\w^{\e}\|_{H^3_x}$ in \eqref{wHs} of Lemma~\ref{L.wrtaHs}.  
Indeed, when $\mathrm{N}=3$, the $H^3_x$-norm of the vorticity can only be estimated as $\|\w^{\e}\|_{H^3_x} \leq \kappa^{-\frac{1}{2}}\mathcal{E}_{top}$. In contrast, if $\mathrm{N}=4$, we have the improved bound $\|\w^{\e}\|_{H^3_x} \leq \mathcal{E}_{M}$. 
\unhide
\end{remark}

\hide
\begin{remark}
(2) If $(\rho_0,u_0,\ta_0)$ and $\AC{\P}F$ are regular up to $\mathrm{N}$th order derivatives, that is, 
\begin{align*}
&\bigg(\|\rho^{\e}_0\|_{L^2_x}^2+\|u^{\e}_0-\bar{u}\|_{L^2_x}^2 +\|\ta^{\e}_0\|_{L^2_x}^2\bigg) 
+ \sum_{1\leq|\al|\leq\mathrm{N}} \bigg(\|\p^{\al}\rho^{\e}_0\|_{L^2_x}^2+\|\p^{\al}u^{\e}_0\|_{L^2_x}^2+\|\p^{\al}\ta^{\e}_0\|_{L^2_x}^2\bigg) < \mathfrak{m}, \\
&\sum_{0\leq|\al|\leq \mathrm{N}}\frac{1}{\e^2} \int_{\Omega\times\R^3}\frac{|\p^{\al}\AC{\P}F|^2}{M^{\e}}\,dv\,dx < \mathfrak{m},
\end{align*}
then its mollification satisfies $\mathcal{E}_{tot}(F_0^{\e}) \leq \mathfrak{m}$. This corresponds to $\mathfrak{p}=0$ in \eqref{iniE2D}. Otherwise, the initial data depends on the singular scale $\mathfrak{d}^{-1}$. 
\end{remark}
\unhide

\subsection{Microscopic fluctuation (Theorem \ref{T.2D.global})}

\subsubsection{Step 1: Bootstrap assumption}

In the proof, we choose any $\mathrm{N}$ for $\mathrm{N}\geq4$ and abuse the notations by writing 
$\mathcal{E}(t)$ and $\mathcal{D}(t)$ for 
$\mathcal{E}^{\mathrm{N}}(F^{\e}(t))$ and $\mathcal{D}^{\mathrm{N}}(F^{\e}(t))$, respectively.
We first note that the Bootstrap assumption \eqref{Boot1} and \eqref{Boot2} imply the weak Bootstrap assumption \eqref{condition}.

\hide
The $\e$, $\kappa$, $\mathfrak{d}$ relation should be intersection of the three estimates. 
\begin{itemize}
\item Proposition \ref{P.Csmall}
Recall that \eqref{ek-rel:P}:
\begin{align*}
\bega
&\e \kappa^{-\frac{3}{4}}\sup_{t \in [0,T]}\mathcal{E}_{tot}(t) \leq C_{min}', 
\enda
\end{align*}
\item \eqref{hassume} for $d=2$
\begin{align*}
\bega
\e\sum_{0\leq |\al_x|\leq \mathrm{N}}\|\p^{\al_x}h_0\|_{L^\infty_{x,v}} + \e^{1-\frac{d}{4}}\kappa^{-\frac{3d}{8}} \sup_{t \in [0,T]} \mathcal{E}_{tot}^{\frac{1}{2}}(t) +\e \sup_{t \in [0,T]}(1+\mathcal{E}_{tot}(t))^{\frac{1}{2}} < \frac{1}{4C},
\enda
\end{align*}
We need for $d=2$ that 
\begin{align*}
\bega
\e^{\frac{1}{2}}\kappa^{-\frac{3}{4}} \sup_{t \in [0,T]} \mathcal{E}_{tot}^{\frac{1}{2}}(t) < \frac{1}{16C}, \qquad \e\kappa^{-\frac{3}{2}} \sup_{t \in [0,T]} \mathcal{E}_{tot}(t) < \frac{1}{(16C)^2}
\enda
\end{align*}

\item Strichartz estimate for $d=2$ guaranting the estimate of $\|\nabla_x\nabla_x\cdot u^{\e}\|_{L^\infty_x}$. Recall Remark \ref{Rmk.divLinf}:
\begin{align*}
\bega
\|\nabla_x\nabla_x\cdot u^{\e}\|_{L^\infty} &\leq C\|\nabla_x\nabla_x\cdot u^{\e}\|_{\dot{B}_{2+\frac{1}{\delta},1}^{\frac{d\delta}{1+2\delta}}}, \quad d=2,3.
\enda
\end{align*}
We need $s=2+\frac{d+1}{2}\frac{1}{2+4\delta}+\frac{2\delta}{1+2\delta}$. 
\begin{align*}
\bega
s=2+\frac{d+1}{2}\frac{1}{2+4\delta}+\frac{2\delta}{1+2\delta} = 2+ \frac{3}{2}\frac{1}{2+4\delta}+\frac{2\delta}{1+2\delta} = 2+ \frac{3+8\delta}{4+8\delta}<3
\enda
\end{align*}
In addition for the estimate of the forcing term in Proposition \ref{P.div.ineq} we need $s<\mathrm{N}-1$. So that, we need $4^-\leq \mathrm{N}$. 
\begin{align*}
\bega
\Big\|\Big((\rho^{\e}+\ta^{\e}),&\mathbb{P}^{\perp}u^{\e}\Big)\Big\|_{L^{\frac{4+8\delta}{d-1}}_T\dot{B}_{2+\frac{1}{\delta},1}^{s-\frac{d+1}{2}\frac{1}{2+4\delta}}} \leq C \e^{\frac{1}{r}}\Big\|\Big((\rho^{\e}_0+\ta^{\e}_0),\mathbb{P}^{\perp}u^{\e}_0\Big)\Big\|_{\dot{B}_{2,1}^s} \cr 
&+ \begin{cases}
C\e^{\frac{1}{r}}\int_0^T\Big(\mathcal{E}_{tot}(s) +\kappa^{\frac{1}{2}}\mathcal{D}_{tot}^{\frac{1}{2}}(s)\Big) ds, \qquad &0\leq s+1 < \mathrm{N}-1, \\
C\e^{\frac{1}{r}}\kappa^{-\frac{s+1}{2\mathrm{N}}}\int_0^T\Big(\mathcal{E}_{tot}(s) +\kappa^{\frac{1}{2}}\mathcal{D}_{tot}^{\frac{1}{2}}(s)\Big) ds, \qquad &\mathrm{N}-1 \leq s+1 < \mathrm{N},
\end{cases} 
\enda
\end{align*}

\end{itemize}

The bootstrap assumption should be satisfy the following three inequality:
\begin{align*}
T_{max}  = \sup\bigg\{&t\geq 0: 
\underbrace{\frac{\e}{\kappa^{ {3}/{2}} }\big(1+t\big)\bigg(\mathcal{E}_{tot}(t)+\int_0^t\mathcal{D}_{tot}(s)ds\bigg)^2 \leq C_{min}^{2d}}_{(\mathcal{A}_1)}, \cr
&\mbox{and}, \quad \underbrace{\frac{\e^{\frac{1}{r}}}{\kappa^{ {1}/{2}} }\big(1+t\big)\bigg(\mathcal{E}_{tot}(t)+\int_0^t\mathcal{D}_{tot}(s)ds\bigg)^2 \leq C_{min}^{2d}}_{(\mathcal{A}_2)}, \cr 
&\mbox{and}, \quad \underbrace{\kappa \big(1+t\big)\int_0^t\mathcal{D}_{tot}(s)ds \leq C_{min}^{2d}}_{(\mathcal{A}_3)} \bigg\}.
\end{align*}
for any $0<\delta<1$. 
We will prove that $T_{max}\geq T_*$ under the assumptions of Theorem \ref{T.2D.global}. 

Let $\e = \kappa^{2+\mathfrak{s}}$ for any positive number $\mathfrak{s}>0$. If $\kappa$ satisfies the following relation, then we can have the result...
We choose $4\delta = \frac{\mathfrak{s}}{2}$. I want to make a statement that 
\begin{align}
&\underbrace{\kappa\big(1+t\big)\bigg(\mathcal{E}_{tot}(t)+\int_0^t\mathcal{D}_{tot}(s)ds\bigg)^2 \leq C_{min}^{2d}}_{(\mathcal{A}_1)}, \cr
&\underbrace{\kappa^{\frac{\mathfrak{s}}{10}} \big(1+t\big)\bigg(\mathcal{E}_{tot}(t)+\int_0^t\mathcal{D}_{tot}(s)ds\bigg)^2 \leq C_{min}^{2d}}_{(\mathcal{A}_2)},
\end{align}
or $\kappa^{\frac{\mathfrak{s}}{24}}$. Then we can have the result. 
If we choose $\delta:= \min\{\frac{1}{8},\frac{\mathfrak{s}}{8}\}$, then we will have $\kappa^{\frac{2+\mathfrak{s}}{5}-\frac{1}{2}}$ or $\kappa^{\frac{\mathfrak{s}}{2}\frac{1}{4+8\delta}}$. If $\kappa^{\frac{\mathfrak{s}}{10}}$ absorb some energy, then it is fine.
We can write it by one line 
\begin{align}
&\underbrace{\kappa^{\min\{1,\frac{\mathfrak{s}}{10}\}}\big(1+t\big)\bigg(\mathcal{E}_{tot}(t)+\int_0^t\mathcal{D}_{tot}(s)ds\bigg)^2 \leq C_{min}^{2d}}_{(\mathcal{A}_1)},
\end{align}
\unhide

{\bf Claim: Assembly of the Three Energy Estimates.}
In the first step, we combine the three energy estimates in Proposition \ref{P.G.Energy}, Proposition \ref{P.F.Energy}, and Proposition \ref{P.EW}. 
Let $T_*$ be the maximal time such that the Bootstrap assumption \eqref{Boot1} and \eqref{Boot2} hold. We claim the following bound:
\begin{align}\label{EBE0}
\bega   
\mathcal{E}_{tot}(t)&+\int_0^t\mathcal{D}_{tot}(s)ds \leq Ce^{C\int_0^t(1+\|\nabla_x(\rho^{\e},u^{\e},\ta^{\e})(s)\|_{L^\infty_x}^2)ds} (\mathcal{E}_{tot}(0)+\e t),
\enda
\end{align}
for $0\leq t< T_*$.
By combining Proposition \ref{P.G.Energy}, Proposition \ref{P.F.Energy}, and Proposition \ref{P.EW}, we obtain
\begin{align}\label{E.sum-1}
\bega
&\frac{1}{C}\mathcal{E}_{tot}(s)+\frac{1}{C}\int_0^t\mathcal{D}_{tot}(s)ds \leq \sum_{0\leq|\al|\leq\mathrm{N}}\int_0^t\Big(\mathfrak{S}_G^{\al}(s)+\mathfrak{S}_W^{\al}(s)\Big)ds \cr 
&+\int_0^t\big(\|\nabla_xu^{\e}(s)\|_{L^\infty_x}+\|\nabla_x\ta^{\e}(s)\|_{L^\infty_x}\big)\big(\mathcal{E}_{tot}(s)+\mathcal{E}_M^{\frac{1}{2}}(s)\mathcal{D}_{top}^{\frac{1}{2}}(s)\big)ds \cr 
&+ \int_0^t\Big(1+\|\nabla_x\bar{u}\|_{L^\infty_x} +\|\nabla_x(\rho^{\e},u^{\e},\ta^{\e})(s)\|_{L^\infty_x} \Big) \mathcal{E}_M(s) + \e\big(\mathcal{E}_M^{\frac{3}{2}}(s)+\mathcal{E}_M^2(s)\big)ds  \cr 
&+ \int_{\Omega}|\nabla_xu^{\e}(t)| \mathcal{V}_{2}(t) dx + \int_{\Omega}|\nabla_x\ta^{\e}(t)| \mathcal{V}_{3}(t) dx + \frac{1}{\eps} \int_{\Omega} \TbT(t) \mathcal{V}_{2}(t) dx \cr
&+ \int_0^t \e\kappa^{\frac{1}{2}} \mathcal{E}_{tot}(s)\mathcal{D}_{tot}^{\frac{1}{2}}(s) + \e\mathcal{E}_{tot}^{\frac{1}{2}}(s)\mathcal{D}_{tot}(s) +\e\mathcal{E}_{tot}^{\frac{3}{2}}(s)\mathcal{D}_{top}^{\frac{1}{2}}(s) ds \cr 
&+ \int_0^t\frac{1}{\kappa^{\frac{1}{2}}\e}\sum_{0\leq |\al| \leq\lfloor(\mathrm{N}+1)/2\rfloor} \Big\|\|\la v \ra^{\frac{1}{2}} \p^{\al}\AC{\P}F^{\e} |M^{\e}|^{-\frac{1}{2}}(s)\|_{L^2_v}\Big\|_{L^\infty_x} \mathcal{E}_{tot}^{\frac{1}{2}}(s)\mathcal{D}_{tot}^{\frac{1}{2}}(s) ds\cr 
&+  \int_0^t\sum_{0\leq|\al|\leq \lfloor(\mathrm{N}+1)/2\rfloor}\Big\|\|\la v \ra^{\frac{1}{2}} \p^{\al}\AC{\P}F^{\e} |M^{\e}|^{-\frac{1}{2}}(s)\|_{L^2_v}\Big\|_{L^\infty_x} \mathcal{D}_{tot}(s) ds \cr
& +\int_0^t\big(\mathcal{Z}_{top}^{time}(s) +\mathcal{Z}_W^{time}(s)\big)ds,
\enda
\end{align}
where $\mathfrak{S}_G^{\al}(t)$ and $\mathfrak{S}_W^{\al}(t)$ are defined in \eqref{Amidef} and \eqref{Amadef}. Here, $\mathcal{Z}_{top}^{time}(t)$ and $\mathcal{Z}_W^{time}(t)$ are defined in \eqref{Bmidef} and \eqref{Bmadef}, and they appear only when we use space–time derivatives $\p_{\tilde{t}}=\e^{\mathfrak{n}}\p_t$. We focus only on the main contributions. 

(1) By applying the cancellation between the leading-order momentum-flux alignment terms $\mathfrak{S}_G^{\al}$ and $\mathfrak{S}_W^{\al}$ from Lemma \ref{L.cancel}, we get
\begin{align*}
\sum_{0\leq|\al|\leq\mathrm{N}}\bigg|\Big(\mathfrak{S}_G^{\al}(t)+\mathfrak{S}_W^{\al}(t)\Big)\bigg| 
\les \begin{cases} \e \mathcal{E}_{tot}(t) \mathcal{D}_{tot}^{\frac{1}{2}}(t), & \hspace{-4mm} \mbox{ finite velocity energy case}, \\
\kappa^{\frac{1}{2}}\|\nabla_x\bar{u}(t)\|_{L^2_x}\mathcal{D}_G^{\frac{1}{2}}(t)+ \e \mathcal{E}_{tot}(t) \mathcal{D}_{tot}^{\frac{1}{2}}(t), & \hspace{-4mm} \mbox{ infinite velocity energy case.} \end{cases}
\end{align*}
By Lemma \ref{L.barubdd}, we have $\|\nabla_x\bar{u}(t)\|_{L^2_x}\leq C$. 

(2) Estimate of $\mathcal{V}_{3}$: 
Recall the definitions of $h^{\e}$ and $\mathfrak{h}^{\e}$ in \eqref{hdef} and \eqref{hdef2}, respectively. The bootstrap assumption \eqref{Boot2} is equivalent to
\begin{align}\label{Boot2-h}
\sum_{0\leq|\al|\leq\mathrm{N}+1}&\kappa^{\frac{(|\al|-\mathrm{N})_+}{2}}\big\|\p^{\al}h^{\e}\big\|_{L^\infty_tL^\infty_{x,v}} + \sum_{0\leq|\al|\leq\mathrm{N}-2}\big\|\p^{\al}\mathfrak{h}^{\e}\big\|_{L^\infty_tL^\infty_{x,v}} \les \frac{1}{\e^2\kappa^2}.
\end{align}
Applying \eqref{Boot2-h} to \eqref{L.Vdecompsum}, we obtain
\begin{align}\label{B1}
\bega
\int_{\Omega}&|\nabla_xu^{\e}(t)| \mathcal{V}_{2}(t) dx + \int_{\Omega}|\nabla_x\ta^{\e}(t)| \mathcal{V}_{3}(t) dx + \frac{1}{\eps} \int_{\Omega} \TbT(t) \mathcal{V}_{2}(t) dx \cr 
&\les \|\nabla_x (u^{\e},\ta^{\e})\|_{L^\infty_x}\Big(\mathcal{E}_{tot} +\e^{\frac{3}{2}}\kappa \mathcal{D}_{tot}\Big)  + e^{-c_1\e^{-{\frac{1}{3}}}}\frac{1}{\e} \mathcal{E}_{tot} \frac{1}{\e^2\kappa^2} \cr 
&+\e^{\frac{1}{2}}\kappa^{\frac{1}{2}}(\mathcal{D}_{G})^{\frac{1}{2}} \mathcal{E}_{tot} + \e\kappa^{\frac{1}{2}} e^{-c_1\e^{-{\frac{1}{2}}}}  (\mathcal{D}_{G})^{\frac{1}{2}}(\mathcal{E}_{tot})^{\frac{1}{2}} \frac{1}{\e^2\kappa^2}.
\enda
\end{align}
Since the exponential factor $e^{-c_1\e^{-{\frac{1}{\ell}}}}$ 
absorbs any polynomially singular power of $\e$, including those 
arising from $\kappa=\e^q$ for any $q>0$, the right-hand side of 
\eqref{B1} can be bounded by 
$\|\nabla_x(u^{\e},\ta^{\e})\|_{L^\infty_x}\mathcal{E}_{tot}$ 
and $\e^{\frac{3}{2}}\mathcal{D}_{tot}$.

(3) For the sixth and seventh lines of \eqref{E.sum-1}, we apply
$L^2_T$ estimate \eqref{GL2} and $L^\infty_T$ estimate \eqref{GLinfty} to the term $\|\|\la v \ra^{\frac{1}{2}} \p^{\al}\AC{\P}F^{\e} |M^{\e}|^{-\frac{1}{2}}\|_{L^2_v}\|_{L^\infty_x}$, respectively.

For brevity, we first consider the purely spatial derivative case \eqref{caseA}, in which the last line of \eqref{E.sum-1} vanishes.  
Applying Young’s inequality $\|\nabla_x(u^{\e},\ta^{\e})\|_{L^\infty_x}\mathcal{E}_M^{\frac{1}{2}}\mathcal{D}_{top}^{\frac{1}{2}}
\leq C \|\nabla_x(u^{\e},\ta^{\e})\|_{L^\infty_x}^2 \mathcal{E}_M + \frac{1}{C}\mathcal{D}_{top}$ to the second line of \eqref{E.sum-1}, then by the Bootstrap assumption, we obtain
\begin{align}\label{EDbddG}
\bega
\mathcal{E}_{tot}(t)+\int_0^t\mathcal{D}_{tot}(s)ds &\leq C\int_0^t\big(1+\|\nabla_x(\rho^{\e},u^{\e},\ta^{\e})(s)\|_{L^\infty_x}^2\big) \mathcal{E}_{tot}(s)ds + \e t.
\enda
\end{align}
Applying Gr\"onwall’s inequality then yields the result \eqref{EBE0}.

\subsubsection{Step 2: Determining the time growth from the macroscopic estimate}

In this step, we claim that for any $t \in [0,T_*)$, the growth rate of the Boltzmann energy obtained in \eqref{EBE0} satisfies 
\begin{align}\label{growth}
\int_0^t \Big(1+\|\nabla_x(\rho,u,\ta)(s)\|_{L^\infty_x}^2\Big) ds \leq \Psi[\w^{\e}_0,\mathfrak{s}^{\e}_0](t),
\end{align}
where $\Psi[\w^{\e}_0,\mathfrak{s}^{\e}_0](t)$ is defined in Definition \ref{D.Psi-def}. If the claim \eqref{growth} holds, then we obtain the result in \eqref{EDbdd2D}:
\begin{align}\label{EDbdd-T}
\bega
\mathcal{E}_{tot}(t)&+\int_0^t\mathcal{D}_{tot}(s)ds \leq Ce^{C\Psi[\w^{\e}_0,\mathfrak{s}^{\e}_0](t)} (\mathcal{E}_{tot}(0)+\e t).
\enda
\end{align}
To prove \eqref{growth}, we decompose $(\rho^{\e},u^{\e},\ta^{\e})$ via the Hodge decomposition \eqref{Hodge} and thermodynamic variables:
\begin{align*}
\int_0^t \|\nabla_x(\rho,u,\ta)(s)\|_{L^\infty_x}^2 ds \leq \mathfrak{G}_1(t) +\mathfrak{G}_2(t) ,
\end{align*}
where
\begin{align*}
\mathfrak{G}_1(t) &:=  2\int_0^t \Big( \|\nabla_x\mathbb{P}^{\perp}u^{\e} \|_{L^\infty_x}^2 +\|\nabla_xp^{\e} \|_{L^\infty_x}^2 \Big)ds, \cr
\mathfrak{G}_2(t) &:= 2\int_0^t \Big(\|\nabla_x\mathbb{P}u^{\e} \|_{L^\infty_x}^2 + \|\nabla_x\mathfrak{s}^{\e}\|_{L^\infty_x}^2\Big)ds .
\end{align*}

(1) Estimate of the compressible parts: 
To estimate $\mathfrak{G}_1(t)$, it suffices to control derivatives of 
$\mathbb{P}^{\perp}u^{\e}$ and $\rho^{\e}+\ta^{\e}$ up to first order.  
However, estimating $\mathfrak{G}_2(t)$ requires derivatives up to second order.  
Therefore, for second-order derivatives $0\leq \ell\leq 2$, by substituting 
$d=2$, $r=4$, and $\mathrm{N}\geq4$ into Remark~\ref{Rmk.divuLinf}, we obtain
\begin{align}\label{divbouclaim}
\bega
\int_0^t&\Big\|\nabla_x^{\ell}\Big((\rho^{\e}+\ta^{\e})(\tau),(\mathbb{P}^{\perp}u^{\e})(\tau)\Big)\Big\|_{L^\infty_x} d\tau \cr 
&\leq C\e^{\frac{1}{4}}(1+T_*)\bigg( \mathcal{E}_{M}(0)+ \int_0^{T_*}\Big(\mathcal{E}_{tot}(s) +\mathcal{D}_{tot}^{\frac{1}{2}}(s)\Big) ds\bigg) \leq  \e^{\frac{1}{4}-\mathfrak{y}},
\enda
\end{align}
for $0\leq t< T_*$ and $0<\mathfrak{y}\ll1$, where we used the bootstrap assumption~\eqref{Boot1}.  
For the initial data, we used the interpolation inequality in \eqref{Besov-inter} together with \eqref{Besov-H}:
\begin{align*}
\big\|\big((\rho^{\e}_0+\ta^{\e}_0),(\mathbb{P}^{\perp}u^{\e}_0)\big)\big\|_{\dot{B}_{2,1}^\ell} \leq \big\|\big((\rho^{\e}_0+\ta^{\e}_0),(\mathbb{P}^{\perp}u^{\e}_0)\big)\big\|_{\dot{B}_{2,2}^0}^{\tau_0} \big\|\big((\rho^{\e}_0+\ta^{\e}_0),(\mathbb{P}^{\perp}u^{\e}_0)\big)\big\|_{\dot{B}_{2,2}^\mathrm{N}}^{1-\tau_0} \leq \mathcal{E}_{M}(0),
\end{align*}
for some $0<\tau_0<1$.

(2) Estimate of the solenoidal parts and $\rho^{\e}-3/2\ta^{\e}$: 
We claim the following three estimates:
\begin{align}
\|\w^{\e}(t)\|_{L^\infty_x} &\leq C\big(\|\w^{\e}_0\|_{L^\infty_x}+1\big), \label{w-T} \\
\lw\|\nabla_x \mathbb{P}u^{\e}(t)\rw\|_{L^\infty_x} &\leq C e^{C(1+\|\w^{\e}_0\|_{L^\infty_x})t}\lw(1+\ln^+\|\w^{\e}_0\|_{H^2_x}\rw)\lw(1+\|\w^{\e}_0\|_{L^\infty_x}\rw), \label{nablau-T} \\
\Big\|\nabla_x\Big(\rho^{\e}-\frac{3}{2}\ta^{\e}\Big)(t)\Big\|_{L^\infty_x} &\leq C\exp\bigg(Ct e^{C(1+\|\w^{\e}_0\|_{L^\infty_x})t}\lw(1+\ln^+\|\w^{\e}_0\|_{H^2_x}\rw)\lw(1+\|\w^{\e}_0\|_{L^\infty_x}\rw)\bigg) \label{drta-T}  \\
&\times \bigg(1+\Big\|\nabla_x\Big(\rho^{\e}_0-\frac{3}{2}\ta^{\e}_0\Big)\Big\|_{L^\infty_x}\bigg), \notag
\end{align}
for $0\leq t < T_*$. 
To prove \eqref{w-T}, we apply \eqref{wmax} in Proposition \ref{P.max}. 
Then, using \eqref{divbouclaim} and
\begin{align*}
&\int_0^t \Big(\mathcal{E}_M^{\frac{1}{2}}(s) \|\nabla_x(\rho^{\e}+\ta^{\e})(s)\|_{L^\infty_x} + \kappa^{\frac{1}{2}}\mathcal{D}_G^{\frac{1}{2}}(s)\Big)ds \leq \frac{1}{C},
\end{align*}
together with \eqref{Boot1}, we obtain \eqref{w-T}. 
To prove \eqref{nablau-T}, we apply \eqref{w-T} to \eqref{pmLinfe} for $\mathrm{N}\geq4$ in Proposition \ref{P.macro.u}. 
Then, using
\begin{align}\label{GH2est-T}
\bega
\int_0^t \mathfrak{O}(s)ds &\leq \int_0^t \bigg(\sum_{0\leq k \leq 2}\Big\|\nabla_x^k\Big((\rho^{\e}+\ta^{\e}),\mathbb{P}^{\perp}u^{\e}\Big)(s)\Big\|_{L^\infty_x}(1+\mathcal{E}_{tot}^{\mathrm{N}}(F^{\e}(s)))+ \kappa^{\frac{1}{2}}\mathcal{D}_G^{\frac{1}{2}}(s)\bigg) ds \cr 
&\leq 1,
\enda
\end{align}
where $\mathfrak{O}(t)$ is defined in \eqref{GH2esti}, we obtain \eqref{nablau-T}. 
To prove \eqref{drta-T}, we apply \eqref{rtamax} in Proposition \ref{P.max}. 
For the exponential factor $\int_0^t \|\nabla_xu^{\e}(s)\|_{L^\infty_x}ds$, we decompose $\nabla_xu^{\e}$ into its Leray-projected part and its orthogonal complement. 
Applying \eqref{nablau-T} to $\|\nabla_x\mathbb{P}u^{\e}\|_{L^\infty_x}$ and \eqref{divbouclaim} to $\|\nabla_x\mathbb{P}^{\perp}u^{\e}\|_{L^\infty_x}$ yields \eqref{drta-T}.

Combining \eqref{divbouclaim}, \eqref{nablau-T}, and \eqref{drta-T} with Definition~\ref{D.Psi-def} completes the proof of the claim \eqref{growth}.

\medskip
When considering the space--time derivative case \eqref{caseB} with $\mathfrak{n}\geq1$, we must also take into account the additional contribution $\mathcal{Z}_{top}^{time}(t)+\mathcal{Z}_W^{time}(t)$ in \eqref{E.sum-1}. 
Using \eqref{divbouclaim} together with the bootstrap assumption \eqref{Boot1}, we obtain
\begin{align*}
\bega
\mathcal{Z}_{top}^{time}(t) +\mathcal{Z}_W^{time}(t) 
&\leq \e^{\mathfrak{n}-1}\Big(C\mathcal{E}_{tot}(t) + \frac{1}{C}\mathcal{D}_{tot}\Big).
\enda
\end{align*}
Therefore, we recover the same inequality as in \eqref{EBE0}. 
Proceeding as in the purely spatial case, we conclude that the estimate \eqref{EDbdd-T} also holds (possibly with different constants).

\hide
\subsubsection{Verification of the Bootstrap assumption} 
To close the estimate, we must show that the growth of $\mathcal{E}_{tot}$ up to time $T_*$ still lies within the space \eqref{space2D}. Recalling \eqref{EDbdd-T}, the definitions of $\kappa$ and $\mathfrak{d}$ in \eqref{ek-rel:2D} ensure that
\begin{align}\label{ekphiT}
\bega
&\e^{\frac{\mathfrak{s}}{20}}\big(1+T_*\big)^2 \bigg(Ce^{C\Psi[\rho_0^\e,\ta^\e_0,\omega^\e_0](T_*)} (\mathcal{E}_{tot}(0)+\e T_*)\bigg)^2\leq C\e^{\frac{\mathfrak{s}}{60}}< \frac{1}{1+\mathfrak{m}}.
\enda
\end{align}
This implies $T_*\leq T_*$, where $T_*$ is defined in \eqref{space2D}. In other words, up to the chosen time $T_*$, the solution satisfies the inequality in \eqref{space2D}. Consequently, the estimate \eqref{EDbdd-T} holds for all $0\leq t \leq T_*$, for any chosen $T_*>1$. This completes the proof of \eqref{EDbdd2D}.

\unhide

\subsubsection{Step 3: Proof of \eqref{incomp.2D} and \eqref{Grefine.2D}}

(Proof of \eqref{incomp.2D})
Applying Proposition \ref{P.div.ineq} with $d=2$ and $s \in [0, \mathrm{N}-1)$ yields the result \eqref{incomp.2D}.  \\
(Proof of \eqref{Grefine.2D})
Using the same argument as in the microscopic energy estimate of Proposition \ref{P.G.Energy}, we sum the energy estimates for 
$\frac{1}{\e^2} \int_{\Omega\times\R^3}|\p^{\al}\AC{\P}F^{\e}|^2|M^{\e}|^{-1}dvdx$ over $0\leq|\al|\leq\mathrm{N}-2$ to obtain
\begin{align}\label{G.low}
\bega
&\mathcal{E}^{\mathrm{N}-2}_G(t) + \int_0^t\mathcal{D}^{\mathrm{N}-2}_G(s)ds
\leq \mathcal{E}^{\mathrm{N}-2}_G(0) +  \sum_{0\leq|\al|\leq\mathrm{N}-2}\int_0^t\mathfrak{S}_G^{\al}(s)ds \cr 
&\quad + C_{}\sum_{0\leq|\al|\leq\mathrm{N}-2}\int_0^t \bigg(\int_{\Omega}|\nabla_xu^{\e}(s)| \mathcal{V}_{2}^{\al}(s) dx + \int_{\Omega}|\nabla_x\ta^{\e}(s)| \mathcal{V}_{3}^{\al}(s) dx\bigg) ds +C\e t. 
\enda
\end{align}
For the second line of \eqref{G.low}, we apply $\eqref{L.Vdecomp3}_1$ together with the bound 
$\|\nabla_x(u^{\e},\ta^{\e})\|_{L^\infty_x}\leq \mathcal{E}_M^{\frac{1}{2}}$ from \eqref{rutinf}. Then
\begin{align*}
\bega
&\int_{\Omega}|\nabla_xu^{\e}| \mathcal{V}_{2}^{\al} dx + \int_{\Omega}|\nabla_x\ta^{\e}| \mathcal{V}_{3}^{\al} dx \cr 
&\les \mathcal{E}_M^{\frac{1}{2}}\bigg[\e^{\frac{3}{2}}\kappa \Big(\mathcal{D}_{tot}^{\mathrm{N}-2}(t)+(\mathcal{E}_M^{\mathrm{N}-2}(t))^2\Big) + e^{-c_1\e^{-{\frac{1}{3}}}}\frac{1}{\e} \mathcal{E}_{tot}^{\mathrm{N}-2}(F^{\e}(t))  \sum_{0\leq|\al|\leq\mathrm{N}-2} \|\p^{\al}\mathfrak{h}^{\e}(t)\|_{L^\infty_{x,v}}\bigg] \to 0,
\enda
\end{align*}
as $\e\to0$. Here we used the bootstrap assumptions \eqref{Boot1} and \eqref{Boot2}. 
For the last term in the first line of \eqref{G.low}, recall the definition of 
$\mathfrak{S}_G^{\al}(t)$ in \eqref{Amadef}.  
Using \eqref{rutscale} and $\eqref{ABGscale}_1$, we obtain
\begin{align*}
\bega
\sum_{0\leq|\al|\leq\mathrm{N}-2}\int_0^t\mathfrak{S}_G^{\al}(s)ds &\les \frac{1}{\e^2}\int_0^t\big(\mathcal{E}_M^{\frac{1}{2}}(s)\big) \big(\e^2\kappa^{\frac{1}{2}}\mathcal{D}_G^{\frac{1}{2}}(s)\big)ds \les \kappa^{\frac{1}{2}}\int_0^t\mathcal{E}_M^{\frac{1}{2}}(s)\mathcal{D}_G^{\frac{1}{2}}(s)ds \to 0, \quad \mbox{as} \quad \e \to 0,
\enda
\end{align*}
where we again used \eqref{Boot1}. 
Since the right-hand side of \eqref{G.low} tends to zero except for the initial data term, we obtain the desired result~\eqref{Grefine.2D}.



\subsection{Propagation of Uniform bounds}\label{S.unif}

In this part, we show that for a class of initial data satisfying a specific blow-up condition, the bootstrap assumptions \eqref{Boot1} and \eqref{Boot2} are ensured (see Lemma~\ref{L.ini}). 
For the solution constructed in Theorem~\ref{T.2D.global}, a further assumption of uniform-in-$\e$ boundedness of the initial data yields the propagation of this uniform boundedness (see Lemma~\ref{L.unif}).

\begin{lemma}\label{L.ini}
Assume a family of initial data $\{F_0^{\e}\}_{\e>0}$ satisfies a modulated entropy bound \eqref{L2unif}, 
\hide
\begin{equation}\label{L2unif}
\sup_{\e>0}
\bigg(  \|( \rho_0^{\e}, u_0^{\e}-\bar{u}, \ta_0^{\e}
)\|_{L^2(\R^2)}
+ \left\|\frac{1}{\e}\frac{\AC{\P}F^{\e}_0}{\sqrt{\tilde{\mu}}}\right\|_{L^2_x(\R^2;L_v^2(\R^3))}
\bigg)
<+\infty.
\end{equation}
Here $\tilde{\mu}=M_{[1,0,1-c_0]}$ for some $0<c_0\ll1$, and $\bar{u}$
denotes the radial eddy defined in Definition~\ref{D.Ra-E}when $d=2$, while we set $\bar{u}=0$ for $d=3$.

\begin{equation*}
\sup_{\e>0}
\bigg(
\|(\rho_0^{\e},u_0^{\e}-\bar{u},\ta_0^{\e})\|_{L^2(\R^d)}
+ \left\|\frac{1}{\e}\frac{\AC{\P}F^{\e}_0}{\sqrt{M^{\e}_0}}\right\|_{L^2(\R^d_x  \times \R^3_v)} 
\bigg)
<+\infty.
\end{equation*}
Here, $\bar{u}$ denotes the radial eddy defined in Definition~\ref{D.Ra-E} when $d=2$, while we set $\bar{u}=0$ for $d=3$.\unhide
and the $N$-Admissible Blow-up Condition in~\eqref{ABC1} of Definition~\ref{ABC}. Then, the family of strong Boltzmann solutions $\{F^{\e}\}_{\e>0}$ to \eqref{BE} constructed in Theorem \ref{T.2D.global} satisfies the Bootstrap asumptions \eqref{Boot1} and \eqref{Boot2}.
\end{lemma}
\begin{proof}
(Proof of \eqref{Boot1})
Applying the initial data \eqref{ABC1} to the estimate \eqref{EDbdd2D}, we deduce that there exists a small constant $0<\mathfrak{y}_0\ll1$ such that
\begin{align*}
\bega
\mathcal{E}_{tot}^{\mathrm{N}}(F^\e (t))+\int_0^t\mathcal{D}_{tot}^{\mathrm{N}}(F^\e (s))ds&\leq \e^{-\mathfrak{y}_0}\Big(\log\!\Big(\log\!\big(\log(1/\e)\big)\Big)+\e t\Big).
\enda
\end{align*}
For any $T>0$ chosen in Theorem \ref{T.2D.global}, multiplying $\e^{\mathfrak{y}}$ for $\mathfrak{y}>\mathfrak{y}_0$ gives 
\begin{align*}
\e^{\mathfrak{y}}\big(1+t\big)\bigg(\mathcal{E}^N_{tot}(F^\e(t))+\int_0^t\mathcal{D}^N_{tot}(F^\e(s))ds\bigg)  < \e^{\mathfrak{y}-\mathfrak{y}_0}\Big(\log\!\Big(\log\!\big(\log(1/\e)\big)\Big)+\e T\Big) < \infty.
\end{align*}
Since there exists $\mathfrak{y}>0$ such that $\e\kappa^{-\frac{1}{2}} = \e^{\mathfrak{y}}$, we can have the same estimate for the scale $\e\kappa^{-\frac{1}{2}}$ in \eqref{Boot1}. This implies the inequality in \eqref{Boot1}. \\
(Proof of \eqref{Boot2}) Recall the definition of $w(v)$, $h^{\e}(t,x,v)$, $\mathfrak{h}^{\e}(t,x,v)$ in \eqref{hdef} and \eqref{hdef2}.
The second inequality of the Admissible Blow-up Condition \eqref{ABC1} implies
\begin{align*}
\bega
&\sum_{0 \leq |\alpha| \leq \mathrm{N}+1}\!
\left\|
\p^\alpha h^{\e}_0 \right\|_{L^\infty_{x,v}} = \sum_{0 \leq |\alpha| \leq \mathrm{N}+1}\!
\left\|w\p^\alpha \left( \frac{F_0^\e - \mu}{\sqrt{\mu}} \right) \right\|_{L^\infty_{x,v}} \ll 1, \cr 
&\sum_{0 \leq |\alpha| \leq \mathrm{N}+1}\!
\left\|
\p^\alpha \mathfrak{h}^{\e}_0 \right\|_{L^\infty_{x,v}} = \sum_{0 \leq |\alpha| \leq \mathrm{N}+1}\!
\left\|w\p^\alpha \left( \frac{F^{\e}_0 - M^{\e}_0}{\sqrt{\bar{\mu}}} \right) \right\|_{L^\infty_{x,v}} \ll 1,
\enda
\end{align*}
where we used $\frac{w}{\sqrt{\mu}} \leq C\frac{w}{\sqrt{\bar{\mu}}} \leq C \frac{1}{\sqrt{\tilde{\mu}}}$. 
Here $\mu=M_{[1,0,1]}$, $\bar{\mu}=M_{[1,0,1-c_0/2]}$, $\tilde{\mu}=M_{[1,0,1-c_0]}$.
Together with the bootstrap assumption~\eqref{Boot1}, the hypotheses of 
Proposition~\ref{P.hLinf2} in~\eqref{hassume} and of 
Proposition~\ref{P.hLinf4} in~\eqref{h2assume} are satisfied. 
Hence, the conclusions \eqref{hLinf} in Propositions~\ref{P.hLinf2} and \eqref{h2Linf} in Proposition \ref{P.hLinf4} follow, which implies the second bootstrap assumption \eqref{Boot2}.
\hide 
\begin{align}
&\sum_{0\leq|\al|\leq\mathrm{N}+1}\kappa^{\frac{(|\al|-\mathrm{N})_+}{2}}\bigg\|w\p^{\al}\bigg(\frac{F^{\e}-\mu}{\sqrt{\mu}}\bigg)\bigg\|_{L^\infty_tL^\infty_{x,v}}  \les \frac{1}{(\e\kappa)^{\frac{d}{2}}} \kappa^{-\frac{1}{2}}
\end{align}
\begin{align}
\sum_{0\leq|\al|\leq\mathrm{N}-2}\bigg\|w\p^{\al}\bigg(\frac{F^{\e}-M^{\e}}{\sqrt{\bar{\mu}}}\bigg)\bigg\|_{L^\infty_tL^\infty_{x,v}} &\les \frac{1}{(\e\kappa)^{\frac{d}{2}}} \kappa^{-\frac{1}{2}} 
\end{align}
\unhide
\end{proof}

\begin{remark}
Assume a family of initial data $\{F_0^{\e}\}_{\e>0}$ satisfies a modulated entropy bound \eqref{L2unif}, and the $N$-Admissible Blow-up Condition in~\eqref{ABC1} of Definition~\ref{ABC}.
In Theorem \ref{T.2D.global}, by choosing $T_*\sim\Big(\log\!\Big(\log\!\big(\log(1/\e)\big)\Big)^{-1}$, we can ensure that the two bootstrap assumptions \eqref{Boot1} and \eqref{Boot2} hold globally. Consequently, the solution is globally valid, and \( T_* \to \infty \) as \( \e  \to 0 \).
\end{remark}

\begin{lemma}
Assume that the family of initial data $\{F_0^{\e}\}_{\e>0}$ satisfies the $\mathrm{N}$-Admissible Blow-up Condition~\eqref{ABC1} in Definition~\ref{ABC}.
Then the compressible parts of the family of solution $\{F^{\e}\}_{\e>0}$ constructed in Theorem \ref{T.2D.global} satisfies  
\begin{align}\label{C*def}
\bega
\sup_{0\leq t< T_*}\mathcal{E}_{tot}&(t) \int_0^{T_*} \|\p^{\al_x}\big((\rho^{\e}+\ta^{\e}),(\mathbb{P}^{\perp}u^{\e})\big)(s)\|_{L^\infty_x}ds + \kappa^{\frac{1}{2}}\Big(\int_0^{T_*} \mathcal{D}_G(s)ds\Big)^{\frac{1}{2}} \cr 
&\leq \begin{cases}
\mathcal{C}_* \Big(\e^{\frac{1}{4}-} +\kappa^{\frac{1}{2}-}\Big) \quad &\mbox{for} \quad |\al_x| \leq \mathrm{N}-2,\\
\mathcal{C}_* \Big(\e^{\frac{1}{4}-}\kappa^{-\frac{1}{2}}+\kappa^{\frac{1}{2}-}\Big) \quad &\mbox{for} \quad |\al_x| = \mathrm{N}-1.
\end{cases}
\enda
\end{align}
\end{lemma}

\begin{lemma}\label{L.unif}
Assume that the family of initial data $\{F_0^{\e}\}_{\e>0}$ satisfies the $\mathrm{N}$-Admissible Blow-up Condition~\eqref{ABC1} in Definition~\ref{ABC} for $\mathrm{N}\geq4$.
If, in addition, we assume that the following bounds hold uniformly in $\e>0$ in each case below, then the corresponding uniform bounds propagate in time as follows:

\hide
Under assumptions \eqref{iniH-s} and \eqref{ini2D.tail}, for mollified initial data associated with initial conditions satisfying the respective independent assumptions in the statement below, the Boltzmann solutions $F^{\e}$ constructed through Theorem \ref{T.2D.global} enjoy the following uniform in $(\e,\kappa(\e)$ estimates, independently of each assumption.  

\begin{itemize}
\item[(1)] If $\w_0 \in H^2(\R^2)$, then
\begin{align*}
\bega
\sup_{0<\e\ll1}\sup_{0\leq t \leq T}\|\w^{\e}(t)\|_{H^2_x} \leq C(\|\w_0\|_{H^2},T),
\enda
\end{align*}
where the constant $C(\|\w_0\|_{H^2},T)$ depends only on $\|\w_0\|_{H^2_x}$ and $T$.


\item[(2)] If $\w_0 \in L^p(\R^2)$ for any $1\leq p \leq \infty$, then
\begin{align*}
\bega
\sup_{0<\e\ll1}\sup_{0\leq t \leq T}\|\w^{\e}(t)\|_{L^p_x} \leq C(\|\w_0\|_{L^p_x},T), \quad \mbox{for} \quad 1\leq p \leq \infty,
\enda
\end{align*}
where the constant $C(\|\w_0\|_{L^p_x},T)$ depends only on $\|\w_0\|_{L^p_x}$ and $T$.

\item[(3)] If $\w_0 \in \mathcal{M}(\R^2)$, then
\begin{align*}
\bega
\sup_{0<\e\ll1}\sup_{0\leq t \leq T}\|\w^{\e}(t)\|_{L^1_x} \leq C(\|\w_0\|_{\mathcal{M}},T).
\enda
\end{align*}
where the constant $C(\|\w_0\|_{\mathcal{M}},T)$ depends only on $\|\w_0\|_{\mathcal{M}}$ and $T$.

\item[(4)] If we further assume $\rho_0, \ta_0 \in H^3(\R^2)$, then the corresponding solution also satisfies
\begin{align*}
\bega
\sup_{0<\e\ll1}\sup_{0\leq t \leq T}\Big\|\Big(\rho^{\e}-\frac{3}{2}\ta^{\e}\Big)(t)\Big\|_{H^3_x} \leq C\bigg(\Big\|\Big(\rho_0-\frac{3}{2}\ta_0\Big)\Big\|_{H^3_x},\|\w_0\|_{H^2_x},T\bigg).
\enda
\end{align*}
Here the constant
$C\!\left(\| \rho_0-\frac{3}{2}\ta_0\|_{H^3_x},\|\w_0\|_{H^2_x},T\right)$
depends only on $\|\rho_0-\frac{3}{2}\ta_0\|_{H^3_x}$, $\|\w_0\|_{H^2_x}$, and $T$.
\end{itemize}
\unhide

\begin{enumerate}
\item (a) For the finite velocity energy case \eqref{caseEC}, if $\sup_{\e>0}\|(\rho^{\e}_0, u^{\e}_0, \ta^{\e}_0)\|_{L^2_x}<\infty$, then $\rho^{\e}(t)$, $u^{\e}(t)$, and $\ta^{\e}(t)$ are uniformly bounded in $L^2_x$: 
\begin{equation}\label{U.L2}
\|(\rho^\e, u^\e, \theta^\e)(t)\|_{L^2_x} \leq C\|(\rho^{\e}_0, u^{\e}_0, \ta^{\e}_0)\|_{L^2_x}+ \mathcal{C}_* \Big(\e^{\frac{1}{4}-} +\kappa^{\frac{1}{2}-}\Big).
\end{equation}
(b) For the infinite velocity energy case \eqref{caseECX}, if
$\sup_{\e>0}\|(\rho^{\e}_0,u^{\e}_0-\bar{u},\ta^{\e}_0)\|_{L^2_x}<\infty$, then $\rho^{\e}(t)$, $u^{\e}(t)-\bar{u}$, and $\ta^{\e}(t)$ are uniformly bounded in $L^2_x$:
\begin{align}\label{U.L2inf}
\bega
\|(\rho^{\e},u^{\e}-\bar{u},\ta^{\e})(t)\|_{L^2_x}
&\leq Ce^{Ct}\|(\rho^{\e}_0,u^{\e}_0-\bar{u},\ta^{\e}_0)\|_{L^2_x}+ \mathcal{C}_* \Big(\e^{\frac{1}{4}-} +\kappa^{\frac{1}{2}-}\Big).
\enda
\end{align}

\item For the finite velocity energy case \eqref{caseEC}, if $\sup_{\e>0}\|\mathbb{P}u^{\e}_0\|_{L^2_x}<\infty$, then $\mathbb{P}u^{\e}(t)$ is uniformly bounded in $L^2_x$:
\begin{align}\label{U.PuL2}
\bega
\|\mathbb{P}u^{\e}(t)\|_{L^2_x} &\leq C\|\mathbb{P}u^{\e}_0\|_{L^2_x} +\mathcal{C}_* \Big(\e^{\frac{1}{4}-} +\kappa^{\frac{1}{2}-}\Big).
\enda
\end{align}

\item If $\sup_{\e>0}\|\mathfrak{s}^{\e}_0\|_{L^2_x}<\infty$, then $\mathfrak{s}^{\e}(t)=\frac{3}{2}\ta^{\e}(t)-\rho^{\e}(t)$ is uniformly bounded in $L^2_x$:
\begin{align}\label{U.rtaL2}
\bega
\|\mathfrak{s}^{\e}(t)\|_{L^2_x} \leq C\|\mathfrak{s}^{\e}_0\|_{L^2_x} +\mathcal{C}_* \Big(\e^{\frac{1}{4}-} +\kappa^{\frac{1}{2}-}\Big).
\enda
\end{align}

\item If $\sup_{\e>0}\|\w^{\e}_0\|_{L^p_x}<\infty$ for $1 \leq p \leq \infty$, then $\w^{\e}_A(t)$ and $\w^{\e}_B(t)$ defined in Definition \ref{D.wAwB} satisfy 
\begin{align}\label{U.wLp<2}
\bega
&\|\w^{\e}_A(t)\|_{L^p_x} = \|\w^{\e}_0\|_{L^p_x}, \quad \mbox{for} \quad 1\leq p \leq \infty, \cr 
&\sup_{t\in[0,T_*]}\|\w^{\e}_B(t)\|_{L^p_x} \leq \mathcal{C}_* \Big(\e^{\frac{1}{4}-} +\kappa^{\frac{1}{2}-}\Big) \to 0 , \quad \mbox{as} \quad \e\to0, \quad &\mbox{for} \quad 2\leq p \leq \infty.
\enda
\end{align}

\item If $\sup_{\e>0}\|\w^{\e}_0\|_{L^\infty_x \cap H^2_x}<\infty$, then
\begin{align}\label{U.uwgrow}
\bega
\lw\|\nabla_x \mathbb{P}u^{\e}(t)\rw\|_{L^\infty_x} &\leq Ce^{C \lw(1+\|\w^{\e}_0\|_{L^\infty_x}\rw)t}\lw(1+\ln^+\|\w^{\e}_0\|_{H^2_x}^2\rw)\lw(1+\|\w^{\e}_0\|_{L^\infty_x}\rw),\cr 
\|\w^{\e}(t)\|_{H^2_x} &\leq e^{e^{C (1+\|\w^{\e}_0\|_{L^\infty_x})t}} \times \Big(\|\w^{\e}_0\|_{H^2_x}\Big)^{e^{C (1+\|\w^{\e}_0\|_{L^\infty_x})t}}.
\enda
\end{align}

\item If $\sup_{\e>0}\|\w^{\e}_0\|_{L^\infty_x \cap H^2_x}<\infty$ and $\sup_{\e>0}\|\nabla_x\mathfrak{s}^{\e}_0\|_{L^\infty_x}<\infty$, then \eqref{U.uwgrow} holds and
\begin{align}\label{U.rtainf}
\bega
\|\nabla_x\mathfrak{s}^{\e}(t)\|_{L^\infty_x} \leq C\exp\bigg(Ct e^{C(1+\|\w^{\e}_0\|_{L^\infty_x})t}\lw(1+\ln^+\|\w^{\e}_0\|_{H^2_x}\rw)\lw(1+\|\w^{\e}_0\|_{L^\infty_x}\rw)\bigg) \big(1+\|\nabla_x\mathfrak{s}^{\e}_0\|_{L^\infty_x}\big).
\enda
\end{align}

\item If $\sup_{\e>0}\|\w^{\e}_0\|_{L^\infty_x \cap H^2_x}<\infty$ and $\sup_{\e>0}\|\mathfrak{s}^{\e}_0\|_{H^3_x}<\infty$ , then \eqref{U.uwgrow} and \eqref{U.rtainf} hold, and
\begin{align}\label{U.rtaH3}
\bega
\|\mathfrak{s}^{\e}(t)\|_{H^3_x} \leq C\exp\bigg(e^{e^{C (1+\|\w^{\e}_0\|_{L^\infty_x})t}} \Big(1+\|\w^{\e}_0\|_{H^2_x}\Big)^{e^{C (1+\|\w^{\e}_0\|_{L^\infty_x})t}}\bigg) \|\mathfrak{s}^{\e}_0\|_{H^3_x}.
\enda
\end{align}
\end{enumerate}
\end{lemma}

\begin{proof}

We note that the growth function $\Psi(t)$ defined in \eqref{Psi-def} may become singular if the initial data 
$\|\nabla_x\mathfrak{s}^{\e}_0\|_{L^\infty_x}$ or 
$\|\w^{\e}_0\|_{H^2_x}$ is not uniformly bounded. 
Nevertheless, it is still possible to obtain estimates that are uniform in $\e$. 
Although the estimate \eqref{EDbdd2D} does not directly control this singularity in the limit, the key idea is to absorb the singular scale 
$\log\!\Big(\log\!\big(\log(1/\e)\big)\Big)$ in \eqref{ABC1} 
into a small fractional power of $\e$, as in \eqref{Boot1}. \\
\noindent (i) 
In the finite velocity energy case \eqref{caseEC}, arguing as in Lemma \ref{L.rutL2w} yields
\begin{align}\label{rutL2fin}
\bega
\frac{1}{2}\frac{d}{dt}\bigg(\|\rho^{\e}\|_{L^2_x}^2+\frac{1}{k_B}\|u^{\e}\|_{L^2_x}^2+\frac{3}{2}\|\ta^{\e}\|_{L^2_x}^2\bigg) &\leq \frac{1}{2}\|\nabla_x\cdot u^{\e}\|_{L^\infty_x} \bigg(\|\rho^{\e}\|_{L^2_x}^2+\frac{1}{k_B}\|u^{\e}\|_{L^2_x}^2+\frac{3}{2}\|\ta^{\e}\|_{L^2_x}^2\bigg) \cr 
&+\|\nabla_x(\rho^{\e}+\ta^{\e})\|_{L^\infty_x} \frac{\|\mathrm{\Theta}^{\e}-1\|_{L^2_x}}{\e}  \|u^{\e}\|_{L^2_x} +\mathfrak{S}_W^{0}(t), 
\enda
\end{align}
where $\mathfrak{S}_W^{0}(t)$ is defined in \eqref{Amadef0}.  
For $\mathfrak{S}_W^{0}(t)$, Hölder’s inequality and $\eqref{ABGscale}_1$ yield $\mathfrak{S}_W^{0}(t) \leq \kappa^{\frac{1}{2}}\mathcal{D}_G^{\frac{1}{2}}\mathcal{E}_M^{\frac{1}{2}}$.
Using \eqref{divbouclaim} and $\frac{1}{\e}\|\mathrm{\Theta}^{\e}(s)-1\|_{L^2_x}\leq \mathcal{E}_M^{\frac{1}{2}}$ from \eqref{pTaleq}, the second line of \eqref{rutL2fin} vanishes in the limit $\e\to0$. Thus, we obtain the result \eqref{U.L2}.  \\
In the infinite velocity energy case \eqref{caseECX}, applying Grönwall’s inequality to Lemma \ref{L.rutL2w} yields the result \eqref{U.L2inf}. \\
\noindent (ii) We multiply $k_B \mathrm{\Theta}^{\e}$ to $\eqref{locconNew}_2$ and then apply the Leray projector $\mathbb{P}$ to obtain
\begin{align*}
&\p_t \mathbb{P}u^{\e} +\mathbb{P}\Big(u^{\e}\cdot \nabla_x u^{\e}\Big) 
=\mathbb{P}\bigg[\frac{k_B(1-\mathrm{\Theta}^{\e})}{\e}\nabla_x(\rho^{\e}+\ta^{\e})- \frac{1}{\e^2}\frac{1}{\mathrm{P}^{\e}}\sum_{j} \p_{x_j} \mathbf{r}_{ij}^{\e}\bigg],
\end{align*}
where we have used $\mathbb{P}\nabla_x(\rho^{\e}+\ta^{\e}) = 0$.  
Taking the $L^2$–energy estimate for $\mathbb{P}u^{\e}$ yields
\begin{align}\label{Puenergy}
\bega
\frac{1}{2}\frac{d}{dt}\|\mathbb{P}u^{\e}\|_{L^2_x}^2 &\leq  \bigg|\int_{\R^2} \mathbb{P}\Big(u^{\e}\cdot \nabla_x u^{\e}\Big) \mathbb{P}u^{\e} dx\bigg| \cr 
&+ \bigg(\bigg\|\frac{k_B(1-\mathrm{\Theta}^{\e})}{\e}\nabla_x(\rho^{\e}+\ta^{\e})\bigg\|_{L^2_x} + \bigg\|\frac{1}{\e^2} \frac{1}{\mathrm{P}^{\e}}\sum_{j} \p_{x_j}\mathbf{r}_{ij}^{\e}\bigg\|_{L^2_x}\bigg)\|\mathbb{P}u^{\e}\|_{L^2_x},
\enda
\end{align}
where we used $\|\mathbb{P}(\cdot)\|_{L^2_x} \leq \|(\cdot)\|_{L^2_x}$.  
We focus on estimating the first term in \eqref{Puenergy}.  
Since the Leray projector is self-adjoint, we have
\begin{align*}
\bega
\int_{\R^2} \mathbb{P}\Big(u^{\e}\cdot \nabla_x u^{\e}\Big) \mathbb{P}u^{\e}  = \int_{\R^2} \Big(u^{\e}\cdot \nabla_x u^{\e}\Big) \mathbb{P} u^{\e} dx.
\enda
\end{align*}
Using the decomposition $u^{\e}= \mathbb{P}u^{\e} + \mathbb{P}^{\perp}u^{\e}$ twice, we obtain
\begin{align}\label{uuPu}
\bega
\int_{\R^2} &\mathbb{P}\Big(u^{\e}\cdot \nabla_x u^{\e}\Big) \mathbb{P}u^{\e} \leq 
\frac{1}{2}\|\nabla_x\cdot u^{\e}\|_{L^\infty_x}\|\mathbb{P}u^{\e}\|_{L^2_x}^2 + \int_{\R^2} \Big(u^{\e}\cdot \nabla_x \mathbb{P}^{\perp}u^{\e}\Big) \mathbb{P}u^{\e} dx \cr 
&\leq \frac{1}{2} \|\nabla_x\cdot u^{\e}\|_{L^\infty_x}\|\mathbb{P}u^{\e}\|_{L^2_x}^2 + \|\nabla_x \mathbb{P}^{\perp}u^{\e}\|_{L^\infty_x}\|\mathbb{P}u^{\e}\|_{L^2_x}^2 + \|\mathbb{P}^{\perp}u^{\e}\|_{L^4_x} \|\nabla_x \mathbb{P}^{\perp}u^{\e}\|_{L^4_x}\|\mathbb{P}u^{\e}\|_{L^2_x}.
\enda
\end{align}
Substituting \eqref{uuPu} into \eqref{Puenergy}, dividing by $\|\mathbb{P}u^{\e}\|_{L^2_x}$, and applying Grönwall’s inequality yields
\begin{align*}
\bega
\bigg|\|\mathbb{P}u^{\e}(t)\|_{L^2_x}-e^{C\int_0^t\|(\nabla_x \mathbb{P}^{\perp}u^{\e})(s)\|_{L^\infty_x}ds} \|\mathbb{P}u^{\e}_0\|_{L^2_x}\bigg| \leq C\int_0^t e^{C\int_s^t\|(\nabla_x \mathbb{P}^{\perp}u^{\e})(\tau)\|_{L^\infty_x}d\tau} \mathfrak{R}(s) ds ,
\enda
\end{align*}
where
\begin{align*}
\bega
\mathfrak{R}(t):= \|\mathbb{P}^{\perp}u^{\e}\|_{L^4_x} \|\nabla_x \mathbb{P}^{\perp}u^{\e}\|_{L^4_x}+ \frac{\|\mathrm{\Theta}^{\e}-1\|_{L^2_x}}{\e} \|\nabla_x(\rho^{\e}+\ta^{\e})\|_{L^\infty_x} +\bigg\|\frac{1}{\e^2} \frac{1}{\mathrm{P}^{\e}}\sum_{j} \p_{x_j}\mathbf{r}_{ij}^{\e}\bigg\|_{L^2_x}.
\enda
\end{align*}
From \eqref{incomp.2D} in Theorem \ref{T.2D.global}, choosing $r=8$ gives
\begin{align*}
\bega
\sum_{0\leq|\al_x|\leq1}\|\p^{\al_x}\mathbb{P}^{\perp}u^{\e}\|_{L^8_TL^4_x} \les \sum_{0\leq|\al_x|\leq1}\|\p^{\al_x}\mathbb{P}^{\perp}u^{\e}\|_{L^8_T\dot{B}_{4,1}^0} \les \e^{\frac{1}{8}-} ,
\enda
\end{align*}
where we used \eqref{Besov-Lp}.  
Since $\frac{1}{\e}\|\mathrm{\Theta}^{\e}-1\|_{L^2_x}\leq \mathcal{E}_M^{1/2}$ from \eqref{pTaleq}, and both $\int_0^t\|(\nabla_x \mathbb{P}^{\perp}u^{\e})(s)\|_{L^\infty_x}ds\to \e^{\frac{1}{4}-}$ and $\int_0^t\|\nabla_x (\rho^{\e}+\ta^{\e})(s)\|_{L^\infty_x}ds\to \e^{\frac{1}{4}-}$ by \eqref{incomp.2D} and Proposition \ref{P.div.ineq}, we finally obtain $\int_0^t \mathfrak{R}(s) ds \leq \e^{\frac{1}{4}-}+\kappa^{\frac{1}{2}-}$. This gives the result \eqref{U.PuL2}. 
\\
(iii) Applying Lemma \ref{L.rho-ta} directly yields \eqref{U.rtaL2}.  \\
\hide (iii) 
For $2\leq p<\infty$, Lemma \ref{L.w-Lp} gives
\begin{align*}
\|\w^{\e}(t)\|_{L^p_x} \leq e^{C\int_0^t\|(\nabla_x\cdot u^{\e})(s)\|_{L^\infty_x}ds}\|\w^{\e}_0\|_{L^p_x} + C\int_0^t e^{C\int_s^t\|(\nabla_x\cdot u^{\e})(\tau)\|_{L^\infty_x}d\tau}\|\varPi_{\w}^{\e}(s)\|_{L^p_x}ds ,
\end{align*}
where $\varPi_{\w}^{\e}$ is defined in \eqref{gwrta}. We observe that 
\[\int_0^t\|\varPi_{\w}^{\e}(s)\|_{L^p_x(\R^2)}ds \leq C\int_0^t\|\varPi_{\w}^{\e}(s)\|_{H^1_x(\R^2)}ds \to 0, \quad \mbox{as} \quad \e\to0,\]
where we used the Gagliardo–Nirenberg interpolation \eqref{Ga-Ni} for $p\geq2$. Then, using \eqref{divbouclaim}, we obtain the result \eqref{U.wLp>2} for $2\leq p<\infty$.  
For $p=\infty$, applying \eqref{wmax} in Proposition \ref{P.max} and noting that $\int_0^t\|\varPi_{\w}^{\e}(s)\|_{L^\infty_x}\,ds \to 0$ as $\e \to0$, we obtain the result \eqref{U.wLp>2} for $p=\infty$. \\ \unhide 
(iv) For $2 \leq p \leq \infty$, using Lemma~\ref{L.w-Lp} and~\eqref{C*def}, we obtain 
\begin{align*}
\bega
\|\w^{\e}_B(t)\|_{L^p_x} \leq \mathcal{C}_* \Big(\e^{\frac{1}{4}-} +\kappa^{\frac{1}{2}-}\Big) \to 0 , \quad \mbox{as} \quad \e\to0, \quad &\mbox{for} \quad 2\leq p \leq \infty,
\enda
\end{align*}
where $\mathcal{C}_*$ is defined in~\eqref{C*def}.  
\\
(v) Let $\w_0^{\e} \in L^\infty_x \cap H^2_x$.  
Using \eqref{nablau-T}, we obtain the uniform bound for $\|\nabla_x \mathbb{P}u^{\e}(t)\|_{L^\infty_x}$ in $\eqref{U.uwgrow}_1$. Similarly, applying the $L^\infty_x$ estimate of $\w^{\e}$ in \eqref{U.wLp<2} with $p=\infty$ to \eqref{wH2e}, and then using \eqref{GH2est-T} to bound $\int_0^t \mathfrak{O}(s)\,ds$, yields $\eqref{U.uwgrow}_2$.  
\\ (vi) Let $\w_0^{\e} \in L^\infty_x \cap H^2_x$ and assume $\|\nabla_x\mathfrak{s}^{\e}_0\|_{L^\infty_x}<\infty$. Applying \eqref{U.uwgrow} to \eqref{drta-T} proves \eqref{U.rtainf}.  \\
(vii) Let $\w_0^{\e} \in L^\infty_x \cap H^2_x$ and $\mathfrak{s}^{\e}_0 \in H^3_x$. Writing the $H^3_x$ energy estimate for \eqref{rtaeqnnew} gives
\begin{align}\label{rtaHs}
\bega
\frac{d}{dt}\|\mathfrak{s}^{\e}\|_{H^3_x}^2 &\les \Big(\|\nabla_x \mathbb{P}u^{\e}\|_{L^\infty_x}\|\nabla_x\mathfrak{s}^{\e}\|_{H^2_x} +\|\nabla_x\mathfrak{s}^{\e}\|_{L^\infty_x}\|\nabla_x\mathbb{P}u^{\e}\|_{H^2_x} \Big)\|\mathfrak{s}^{\e}\|_{H^3_x} \cr 
& + \big\|\varPi_{\mathfrak{s}}^{\e}\big\|_{H^3_x}\|\mathfrak{s}^{\e}\|_{H^3_x} + \sum_{1\leq|\al_x|\leq3}\|\p^{\al_x}\mathbb{P}^{\perp}u\|_{L^\infty_x}\|\mathfrak{s}^{\e}\|_{H^3_x}^2,
\enda
\end{align}
where $\varPi_{\mathfrak{s}}^{\e}$ is defined in \eqref{gwrta}.
Using $\eqref{ABGscale}_1$ for $1\leq |\al|\leq 4$, we obtain
\begin{align*}
\bega
\big\|\varPi_{\mathfrak{s}}^{\e}\big\|_{H^3_x} &\leq \frac{C}{\e^2}\sum_{1\leq|\al|\leq4}\bigg(\|\p^{\al}\mathbf{r}_{ij}^{\e}(t)\|_{L^2_x}+\|\p^{\al}\mathfrak{q}_j^{\e}(t)\|_{L^2_x} + \sum_{0\leq \beta<\al} \|\p^{\al-\beta}\mathrm{U}^{\e}_j \p^{\beta}\mathbf{r}_{ij}^{\e}(t)\|_{L^2_x}\bigg) \leq C\kappa^{\frac{1}{2}}\mathcal{D}_G^{\frac{1}{2}}.
\enda
\end{align*}
Dividing \eqref{rtaHs} by $\|\mathfrak{s}^{\e}\|_{H^3_x}$ gives
\begin{align*}
\bega
\frac{d}{dt}\|\mathfrak{s}^{\e}\|_{H^3_x} &\les \Big(\|\nabla_x \mathbb{P}u^{\e}\|_{L^\infty_x} +\|\w^{\e}\|_{H^2_x} +\sum_{1\leq|\al_x|\leq3}\|\p^{\al_x}\mathbb{P}^{\perp}u\|_{L^\infty_x} \Big)\|\mathfrak{s}^{\e}\|_{H^3_x} +\big\|\varPi_{\mathfrak{s}}^{\e}\big\|_{H^3_x}.
\enda
\end{align*}
We use Agmon’s inequality \eqref{Agmon}, namely $\|\nabla_x\mathfrak{s}^{\e}\|_{L^\infty_x}\leq \|\mathfrak{s}^{\e}\|_{H^3_x}$.
For the growth factors $\|\nabla_x\mathbb{P}u^{\e}\|_{L^\infty_x}$ and $\|\w^{\e}\|_{H^2_x}$, we use $\|\nabla_x\mathbb{P}u^{\e}\|_{H^2_x}  \leq C\|\w^{\e}\|_{H^2_x}$ from \eqref{pmw}, and then apply \eqref{U.uwgrow} to conclude the result. 
\end{proof}

\hide
\begin{remark}
Let us take a solution from Theorem \ref{T.2D.global} when we use the space-timde derivatives \eqref{caseB} with $\mathcal{E}_{tot}\leq \mathfrak{m}$.
Recall the scaled time derivative $\p_{\tilde{t}} =\e^{\mathfrak{n}}\p_t$. With the relation between $\e$ and $\kappa$ in Theorem \ref{T.2D.global}, even though we have \eqref{EDbdd2D}, we get
\begin{align*}
&\|\p_t(\rho^{\e},u^{\e},\ta^{\e})(t)\|_{H^s_x} \les \frac{\kappa^{-\frac{(|\al|-\mathrm{N})_+}{2}}}{\e^{\mathfrak{n}}}\mathcal{E}_{tot}^{\frac{1}{2}}(t) \rightarrow \infty, \quad \mbox{as} \quad \e \rightarrow 0,
\end{align*}
for $0\leq s\leq \mathrm{N}$.
\end{remark}

Even though the upper bound of the $H^s_x$ norm of $\p_t(\rho^{\e},u^{\e},\ta^{\e})$ blows up, we can obtain a refined scale for the quantities $\p_t \w^{\e}$ and $\p_t (\rho^{\e} - \frac{3}{2}\ta^{\e})$ from the solution constructed in Theorem \ref{T.2D.global}.

\begin{lemma}\label{L.t-refine}
For the solution constructed in solution from Theorem \ref{T.2D.global} when we use the space-timde derivatives \eqref{caseB} with $\mathcal{E}_{tot}\leq \mathfrak{m}$, we have 
\begin{align*}
\bega
&\|\p_t\w^{\e}\|_{L^2_TH^2_x} \leq C, \qquad \Big\|\p_t\Big(\rho^{\e}-\frac{3}{2}\ta^{\e}\Big)\Big\|_{L^2_TH^3_x} \leq C_T,
\enda
\end{align*}
for some positive constant $C_T>0$ depending on the right-side of \eqref{EDbdd2D}. 
\end{lemma}

\begin{proof}
Recall the vorticity equation \eqref{weqnnew}:
\begin{align}
\bega
&\p_t\w^{\e} + u^{\e}\cdot \nabla_x \w^{\e} + (\nabla_x\cdot u^{\e})\w^{\e} =\varPi_{\w}^{\e}, 
\enda
\end{align}
\begin{align}
\bega
\varPi_{\w}^{\e}(t,x) &:= -k_B\mathrm{\Theta}^{\e}\nabla_x^{\perp}\ta^{\e}\cdot\nabla_x(\rho^{\e}+\ta^{\e})
-\frac{1}{\e^2}\nabla_x^{\perp}\cdot\bigg(\frac{1}{\mathrm{P}^{\e} }\sum_{j} \p_{x_j} \mathbf{r}_{ij}^{\e}\bigg), 
\enda
\end{align}
We take the $H^1_x$ norm of the vorticity equation \eqref{weqn-g}. For the convection term, we have  
\begin{align*}
\bega
\|u^{\e} \cdot \nabla_x \w^{\e}\|_{H^1_x} \les (\|u^{\e}\|_{L^\infty_x} + \|\nabla_x u^{\e}\|_{L^\infty_x}) \|\nabla_x \w^{\e}\|_{H^1_x}  \leq C_{conv}.
\enda
\end{align*}
We also have $\|\varPi_{\w}\|_{L^2_tH^2_x} \rightarrow 0$. Thus, we conclude that $\|\p_t\w^{\e}\|_{L^2_tH^2_x} \leq C$. Similarly, we can derive the estimate for $(\rho^{\e} - \frac{3}{2}\ta^{\e})$ using the equation \eqref{rtaeqnnew}.
\end{proof}
\unhide


\section{Uniform local estimates in 3D}\label{Sec.3D}

This section is devoted to establishing a local-in-time uniform estimate in the three-dimensional case. 
The theorem proved here holds provided that $(\rho^{\e}_0,u^{\e}_0,\ta^{\e}_0)$, together with the microscopic initial data, are uniformly bounded in $\e$ at least in $H^3_x$.
The result below is stated for the purely spatial derivative case~\eqref{caseA}. 
The slight differences arising in the space--time derivative case~\eqref{caseB} with arbitrary $\mathfrak{n}\geq0$ will be discussed in a remark.


\begin{theorem}\label{T.3D.unif}
Let $\Omega=\R^3$. For some $\mathrm{N}\geq3$, 
assume that the initial data are uniformly bounded in $\e$ as follows:
\begin{align}\label{ini.unif.3D}
\bega
\sup_{\e>0}\bigg(&\sum_{0\leq|\al|\leq \mathrm{N}+1} \kappa^{(|\al|-\mathrm{N})_+}\|\p^{\al}(\rho^{\e}_0,u^{\e}_0,\ta^{\e}_0)\|_{L^2_x}^2 \bigg) < \mathfrak{m}, \cr 
\sup_{\e>0}\bigg(&\sum_{0\leq|\al|\leq\mathrm{N}+1}\frac{\kappa^{(|\al|-\mathrm{N})_+}}{\e^2} \int_{\Omega\times\R^3}\frac{|\p^{\al}\AC{\P}F^{\e}|^2}{M^{\e}}\bigg|_{t=0}dvdx\bigg) < \mathfrak{m},
\enda
\end{align}
for a free parameter $\mathfrak{m}>0$. 
Suppose that the $\mathrm{N}$-Bootstrap Assumptions~\eqref{Boot1} and~\eqref{Boot2} hold on $t \in [0,T_*]$ for $T_*\les \frac{1}{\mathfrak{m}}$. 
Then $F^{\e}$ satisfies the following uniform-in-$\e$ estimates for all $0 \leq t \leq T_*$:
\begin{itemize}
\item 
The total energy and dissipation defined in~\eqref{EDtotdef} satisfy
\begin{align}\label{EDbdd3D}
\bega
\mathcal{E}_{tot}^{\mathrm{N}}(F^\e (t))+\int_0^t\mathcal{D}_{tot}^{\mathrm{N}}(F^\e (s))ds&< \infty.
\enda
\end{align}

\item  
The acoustic/compressible components are uniformly controlled uniformly-in-$\e$ for any $s\in [0,\mathrm{N}-1)$:
\begin{align}\label{incomp.T.loc}
\bega
\Big\|\Big((\rho^{\e}+\ta^{\e}),\mathbb{P}^{\perp}u^{\e}\Big)&\Big\|_{L^r_T\dot{B}_{p,1}^{s+3(\frac{1}{p}-\frac{1}{2})+\frac{1}{r}}} \leq C \e^{\frac{1}{r}}\Big\|\Big((\rho^{\e}_0+\ta^{\e}_0),\mathbb{P}^{\perp}u^{\e}_0\Big)\Big\|_{\dot{B}_{2,1}^s} \cr 
&+ C\e^{\frac{1}{r}}\int_0^{T}\Big(\mathcal{E}_{tot}^{\mathrm{N}}(F^\e (s)) +\kappa^{\frac{1}{2}} \sqrt{\mathcal{D}_{tot}^{\mathrm{N}}(F^\e (s))}\Big) ds.
\enda
\end{align}
Here $2\leq p\leq \infty$, $\frac{1}{r} \leq \frac{1}{2}-\frac{1}{p}$ and $(r,p,d) \neq (2,\infty,3)$ while the Leray projector $\mathbb{P}$ and the Besov norms are defined in~\eqref{LerayPdef} and~\eqref{Besovdef}, respectively.

\item The microscopic energy at lower regularity is controlled by
\begin{align}\label{Grefine.T.loc}
\bega
\mathcal{E}^{\mathrm{N}-2}_G( F^\e (t) ) +\frac{\sigma_L}{C}\int_0^t\mathcal{D}^{\mathrm{N}-2}_G(F^\e (s) ) ds &\leq \mathcal{E}^{\mathrm{N}-2}_G(F^\e (0)) \cr 
&+ \kappa^{\frac{1}{2}}\mathcal P (\mathcal{E}_{tot}^{\mathrm{N}}(F^\e (t)),\mathcal{D}_{tot}^{\mathrm{N}}(F^\e (t)) ),
\enda
\end{align}
where $\mathcal P (a,b)$ is a polynomial in $a$ and $b$.

\end{itemize}

\end{theorem}

\begin{corollary}
(1) If the following estimate holds for some $T_M>0$,
\begin{align*}
\sup_{0\leq t\leq T_M}\|\nabla_x(\rho^{\e},u^{\e},\ta^{\e})(t)\|_{L^\infty_x} \leq C,
\end{align*}
then the existence time $T_*$ in Theorem~\ref{T.3D.unif} can be extended up to $T_M$. \\
(2) If the initial energy $\mathcal{E}_{tot}(0)$ is sufficiently small, then Theorem~\ref{T.3D.unif} holds for any time $T_*>0$. 
Indeed, in Theorem~\ref{T.3D.unif} we may choose $T_*< \frac{1}{C}\frac{1}{C\mathcal{E}_{tot}(0)+\e}$.
\end{corollary}

\hide
\begin{remark}\label{Rmk.divLinf}
If we use space-time derivatives \eqref{caseB}, an additional contribution 
$\|\nabla_x\nabla_x\cdot u^{\e}(t)\|_{L^\infty_x}\mathcal{E}_{M}(t)$ appears. 
To bound this term by the energy, we may use Agmon's inequality \eqref{Agmon}, $\|\nabla_x\nabla_x\cdot u^{\e}(t)\|_{L^\infty_x}\leq \|\nabla_x^2u^{\e}\|_{H^2_x} \les  \mathcal{E}_{M}^{\frac{1}{2}}(t)$.
This requires $\mathrm{N} > \frac{d}{2}+2$. Alternatively, we can estimate it using the Strichartz estimate in Proposition \ref{P.div.ineq}. 
By applying \eqref{Besov-Lp}, we have
\begin{align*}
\bega
\|\nabla_x\nabla_x\cdot u^{\e}\|_{L^\infty} &\leq C\|\nabla_x\nabla_x\cdot u^{\e}\|_{\dot{B}_{2+\frac{1}{\delta},1}^{\frac{d\delta}{1+2\delta}}}, \quad d=2,3.
\enda
\end{align*}
To control the right-hand side with Proposition \ref{P.div.ineq}, we need to choose $s= 2+\frac{d\delta}{1+2\delta} + \frac{d+1}{2}\frac{1}{2+4\delta}< \mathrm{N}-1$. Thus, for $\delta \to 0^+$, one requires $\mathrm{N} > 3^+ + \frac{d+1}{4}$ when $d=2,3$. This restriction arises because of the forcing term (see Remark \ref{Rmk.divusimil}).
\end{remark}
\unhide

\begin{remark}
To obtain the same result as in Theorem~\ref{T.3D.unif} when using space--time derivatives, we must restrict $\kappa=\e^{q}$ to $q<\frac{2}{9}$. 
This restriction originates from assumption \eqref{hassume}, which is required to control the high-order moment estimates.
In addition, the choice of $\mathrm{N}$ depends on the range of $\mathfrak{n}$. 
The only nontrivial contribution arising from the time derivative is
$\e^{\mathfrak{n}-1}\|\nabla_x\nabla_x\cdot u^{\e}(t)\|_{L^\infty_x}\mathcal{E}_{M}(t)$,
coming from $\mathcal{Z}_W^{time}(t)$ defined in \eqref{Bmadef}.
If $0<\mathfrak{n}<1$, the conservation laws \eqref{locconNew} provide the better scaling for $\|\nabla_x\nabla_x\cdot u^{\e}(t)\|_{L^\infty_x}$, yielding the same result for $\mathrm{N}\geq3$. 
If $\mathfrak{n}>1$, the factor $\e^{\mathfrak{n}-1}$ absorbs the remaining terms via the bootstrap assumption \eqref{Boot1}, again giving the same result for $\mathrm{N}\geq3$. 
However, when $\mathfrak{n}=1$, the term $\|\nabla_x\nabla_x\cdot u^{\e}(t)\|_{L^\infty_x}$ must be estimated using Strichartz estimates to extract a small power of $\e$, which requires $\mathrm{N}\geq4$.

\hide
\begin{itemize}
\item Space-time derivatives \eqref{caseB} with $\mathfrak{n} = 1$, for $\mathrm{N} > d/2+2$.
\item Space-time derivatives \eqref{caseB} with $\mathfrak{n} > 1$, for $\mathrm{N} > d/2+1$, and $(\e,\kappa)$ satisfies $\e^{\mathfrak{n}-1}\kappa^{-\frac{1}{2}} \leq 1$.
\item Space-time derivatives \eqref{caseB} with $\mathfrak{n} < 1$, for $\mathrm{N} > d/2+1$.
\end{itemize}
\unhide
\end{remark}


\begin{proof}[\textbf{Proof of Theorem \ref{T.3D.unif}}]
The proof is similar to that of Theorem~\ref{T.2D.global}, so we only sketch the main steps. 
Under the bootstrap assumptions \eqref{Boot1} and \eqref{Boot2}, we obtain the same estimate as in \eqref{EDbddG}. 
Applying the bound $\|\nabla_x(\rho^{\e},u^{\e},\ta^{\e})(t)\|_{L^\infty_x}^2 \leq \mathcal{E}_M(t)$ from \eqref{rutinf} to \eqref{EDbddG} yields
\begin{align*}
\bega
\mathcal{E}_{tot}(t)+\int_0^t\mathcal{D}_{tot}(s)ds \leq  \mathcal{E}_{tot}(0)+  C\int_0^t\Big((1+\mathcal{E}_{tot}(s))\mathcal{E}_{tot}(s)+\e\Big)ds.
\enda
\end{align*}
Consequently, for $0\leq t\leq T_*$,
\begin{align}\label{Etotgrowth}
\mathcal{E}_{tot}(t)+\int_0^t\mathcal{D}_{tot}(s)ds  \leq \frac{1}{\Big(C\mathcal{E}_{tot}(0)+\e\Big)^{-1}-Ct}, \qquad T_*< \frac{1}{C}\frac{1}{C\mathcal{E}_{tot}(0)+\e}.
\end{align}
This implies \eqref{EDbdd3D}. 
The remaining estimates \eqref{incomp.T.loc} and \eqref{Grefine.T.loc} follow by the same arguments as in the proof of Theorem~\ref{T.2D.global}.

\end{proof}

\begin{lemma}\label{L.ini3d}
Assume that a family of initial data $\{F_0^{\e}\}_{\e>0}$ is uniformly
bounded in $\e$ as in \eqref{ini.unif.3D}. 
In addition, suppose that $\{F_0^{\e}\}_{\e>0}$ satisfies
\begin{align*}
\bega
&\sum_{0 \leq |\alpha| \leq \mathrm{N}+1}\!
\left\|
\p^\alpha \left( \frac{F_0^\e - \mu}{\sqrt{\tilde{\mu}}} \right)
 \right\|_{L^\infty_{x,v}} \ll 1.
 \enda
\end{align*}
Then the family of strong Boltzmann solutions $\{F^{\e}\}_{\e>0}$ to \eqref{BE}, constructed in Theorem~\ref{T.3D.unif} on the time interval $t\in[0,T_*]$, satisfies the bootstrap assumptions \eqref{Boot1} and \eqref{Boot2} on $t\in[0,T_*]$.
\end{lemma}
\begin{proof}
Since the estimate \eqref{Etotgrowth} is not singular with respect to the scaling, for any $\mathfrak{y}>0$ the factor $\e^{\mathfrak{y}}$ absorbs the right-hand side of \eqref{Etotgrowth}. This proves the first bootstrap assumption \eqref{Boot1}. The second bootstrap assumption \eqref{Boot2} follows in the same way as in the proof 
of Lemma~\ref{L.ini}.
\end{proof}

\hide
\begin{lemma}
Assume that the family of initial data $\{F_0^{\e}\}_{\e>0}$ satisfies \eqref{ini.unif.3D} and \eqref{ini.v.3D}.
Then the family of strong Boltzmann solutions $\{F^{\e}\}_{\e>0}$ to \eqref{BE}, constructed in Theorem~\ref{T.3D.unif} on the time interval $t\in[0,T_*]$, satisfies the following estimates uniformly in $\e$:
\begin{align}\label{U.rtaH3}
\bega
\|\w^{\e}(t)\|_{H^2_x} &\leq C(\|\w^{\e}_0\|_{H^2_x},T_*,\e^{\frac{1}{2}-},\kappa^{\frac{1}{2}}-), \cr 
\Big\|\Big(\rho^{\e}-\frac{3}{2}\ta^{\e}\Big)(t)\Big\|_{H^3_x} &\leq C\bigg(\|\w^{\e}_0\|_{H^2_x\cap L^\infty_x},\Big\|\Big(\rho^{\e}_0-\frac{3}{2}\ta^{\e}_0\Big)\Big\|_{H^3_x},T_*,\e^{\frac{1}{2}-},\kappa^{\frac{1}{2}}-\bigg).
\enda
\end{align}
Here $C(a,b,c,d)$ denotes a positive constant depending on $a,b,c,d$.
\end{lemma}
\begin{proof}
Taking $\nabla_x^{\perp}$ of $\eqref{locconNew}_2$ yields
\begin{align*}
\bega
&\p_t\w^{\e} + \mathbb{P}u^{\e}\cdot \nabla_x \w^{\e} 
= \w^{\e}\cdot\nabla_x u^{\e} + \bar{\varPi}_{\w}^{\e},
\enda
\end{align*}
where $\bar{\varPi}_{\w}^{\e}$ is defined in \eqref{barvarPidef}. 
Proceeding as in the proof of $\eqref{U.uwgrow}_2$ and \eqref{U.rtaH3} 
yields the desired result.
\end{proof}
\unhide

\section{Convergence}\label{Sec.conv}

In this section, we present convergence theorems for different topologies of the initial vorticity:
\[
\w^{\e}_0 \in H^2(\R^d), \quad 
\w^{\e}_0 \in (L^\infty  \cap L^1)(\R^2), \quad
\w^{\e}_0 \in (L^p  \cap L^1) (\R^2), \quad
\w^{\e}_0 \in \mathcal{M}(\R^2) \ \text{with distinguished sign}.
\]

For later convenience, we rewrite the equations for $\w^{\e}$ and $\mathfrak{s}^{\e}=\frac{3}{2}\ta^{\e}-\rho^{\e}$ derived from the Boltzmann equation.

\begin{lemma}
For $\Omega=\R^2$, the equation for $\w^{\e}$ in \eqref{weqnnew} and the equation for $\mathfrak{s}^{\e}=\frac{3}{2}\ta^{\e}-\rho^{\e}$ in \eqref{rtaeqnnew} can be rewritten as follows:
\begin{align} 
 \p_t\w^{\e} + \mathbb{P}u^{\e}\cdot \nabla_x \w^{\e} = \bar{\varPi}_{\w}^{\e},\label{weqnP} \\ 
\p_t \mathfrak{s}^{\e}  + \mathbb{P}u^{\e} \cdot \nabla_x \mathfrak{s}^{\e} =  \bar{\varPi}_{\mathfrak{s}}^{\e},\label{rtaeqnP} 
\end{align}
where 
\begin{align}\label{barvarPidef}
\bega
\bar{\varPi}_{\w}^{\e}(t,x)&:= -\nabla_x\cdot \big( \mathbb{P}^{\perp}u^{\e} \w^{\e}\big) -k_B\mathrm{\Theta}^{\e}\nabla_x^{\perp}\ta^{\e}\cdot\nabla_x(\rho^{\e}+\ta^{\e})
-\frac{1}{\e^2}\nabla_x^{\perp}\cdot\bigg(\frac{1}{\mathrm{P}^{\e} }\nabla_x \cdot \mathbf{r}^{\e}\bigg), \cr 
\bar{\varPi}_{\mathfrak{s}}^{\e}(t,x)&:= -\mathbb{P}^{\perp}u^{\e} \cdot \nabla_x \mathfrak{s}^{\e} - \frac{1}{\e^2}\frac{1}{k_B\mathrm{P}^{\e}\mathrm{\Theta}^{\e} } \big(\nabla_x\cdot \mathfrak{q}^{\e} + \nabla_x \mathrm{U}^{\e} : \mathbf{r}^{\e} \big).
\enda
\end{align}
\end{lemma}

\subsection{$H^2$-vorticity (classical solution)}

In this section, we prove that the vorticity $\w^{\e}$, associated with the Boltzmann solutions constructed in Theorem~\ref{T.3D.unif} for $d=3$ and Theorem~\ref{T.2D.global} for $d=2$, converges to that of the incompressible Euler equations in $H^3(\R^d)$ for $d=2,3$, as long as the Boltzmann solution exists.
To avoid the singular terms $\frac{1}{\e}\nabla_x \cdot u^{\e}$ and $\frac{1}{\e}\nabla_x(\rho^{\e}+\ta^{\e})$ arising in the local conservation laws \eqref{locconNew}, we establish convergence at the level of the differences $(\w^{\e} - \w^E)$ and $(\tfrac{3}{2}\ta^{\e}-\rho^{\e}) - (\tfrac{3}{2}\ta^E-\rho^E)$.
Moreover, we obtain a convergence rate determined by the microscopic dissipation $\kappa^{\frac{1}{2}}$ and the acoustic dispersion $\e^{\frac{d-1}{4}}$.

From the perspective of the zero-viscosity limit in incompressible flows, convergence has been extensively studied in the literature. In two dimensions with bounded vorticity, convergence of the velocity was established in \cite{Chemin}, while convergence rates for smooth vorticity in dimensions $d=2,3$ were obtained in \cite{Masmoudi}. For results on the low Mach number limit, we refer to \cite{Danchin,Danchin1,Me-Sch}. 

Recall the incompressible Euler equations \eqref{incompE} for $(\rho^{E},u^{E},\ta^{E})$. 
\hide 
\begin{align}\label{incompE}
\bega
&\p_t \rho^E + u^E \cdot \nabla_x \rho^E = 0, \cr 
&\p_t u^E + u^E \cdot \nabla_x u^E + \nabla_x p^E = 0, \cr 
&\p_t \ta^E + u^E \cdot \nabla_x \ta^E = 0, \cr 
&\nabla_x \cdot u^E = 0, \qquad \nabla_x(\rho^E + \ta^E) = 0.
\enda
\end{align} \unhide
In the two–dimensional case, the vorticity satisfies a pure transport equation, whereas in the three–dimensional case it additionally includes the vortex–stretching term:
\begin{align}\label{weqnE}
\bega
\begin{cases}\p_t \w^E + u^E \cdot \nabla_x \w^E = 0,  &\mbox{when} \quad d=2, \\ 
\p_t \w^E + u^E \cdot \nabla_x \w^E = \w^E\cdot\nabla_x u^E,  &\mbox{when} \quad d=3 , \end{cases}
\enda
\end{align}
where $\w^E(t,x) = \nabla_x^\perp \cdot u^E(t,x) = -\p_2 u^E_1(t,x) + \p_1 u^E_2(t,x)$ for $d=2$, and $\w^E(t,x) = \nabla_x \times u^E(t,x)$ for $d=3$.
Since the case $d=3$ subsumes $d=2$, we consider only the three-dimensional setting in this subsection.

\begin{definition}\label{D.macrod}
We define the differences between the macroscopic fields of the Boltzmann equation and the incompressible Euler equations as follows:
\begin{align*}
\bega
(\rho^\mathsf{d}, u^\mathsf{d}, \ta^\mathsf{d},\w^\mathsf{d})(t,x) &:= (\rho^{\e}, \mathbb{P}u^{\e}, \ta^{\e},\w^{\e})(t,x) - (\rho^E, u^E, \ta^E,\w^E)(t,x).
\enda
\end{align*}
\end{definition}

We note that, in the three–dimensional case, the vorticity equation contains the vortex stretching term $\w\cdot\nabla_x u$. Since the case $d=3$ subsumes the estimates for $d=2$, we consider only the three–dimensional setting for simplicity.

\begin{lemma}\label{L.deqn}
The differences $u^{\mathsf{d}}$, $\w^\mathsf{d}$ and $\rho^\mathsf{d} - \frac{3}{2}\ta^\mathsf{d}$ defined in Definition \ref{D.macrod} satisfy
\begin{align}\label{barueqn}
\bega
\p_tu^{\mathsf{d}} &+\mathbb{P}\big(\mathbb{P}u^{\e}\cdot \nabla_x \mathbb{P}u^{\e}\big) - u^E\cdot\nabla_x u^E 
=\mathbb{P}\bigg[\frac{k_B(1-\mathrm{\Theta}^{\e})}{\e}\nabla_x(\rho^{\e}+\ta^{\e})- \frac{1}{\e^2}\frac{1}{\mathrm{P}^{\e}}\sum_{j} \p_{x_j} \mathbf{r}_{ij}^{\e}\bigg] \cr 
& + \nabla_xp^E - \mathbb{P}\Big(\mathbb{P}^{\perp}u^{\e}\cdot \nabla_x \mathbb{P}u^{\e} +\mathbb{P}u^{\e}\cdot \nabla_x \mathbb{P}^{\perp}u^{\e} + \mathbb{P}^{\perp}u^{\e}\cdot \nabla_x \mathbb{P}^{\perp}u^{\e}\Big),
\enda
\end{align}
and
\begin{align}\label{barweqn}
\bega
\p_t \w^\mathsf{d} &+ \mathbb{P}u^{\e} \cdot \nabla_x \w^\mathsf{d} + u^\mathsf{d} \cdot \nabla_x \w^E = \bar{\varPi}_{\w}^{\e}, \quad &\mbox{for} \quad d=2, \cr
\p_t \w^\mathsf{d} &+ \mathbb{P}u^{\e} \cdot \nabla_x \w^\mathsf{d} + u^\mathsf{d} \cdot \nabla_x \w^E =   \w^\mathsf{d}\cdot\nabla_x u^{\e} + \w^E\cdot \nabla_x(u^{\e}-u^E) + \bar{\varPi}_{\w}^{\e}, \quad &\mbox{for} \quad d=3,
\enda
\end{align}
and
\begin{align}\label{barrtaeqn}
\bega
\p_t \Big( \frac{3}{2}\ta^\mathsf{d}-\rho^\mathsf{d}\Big) + \mathbb{P}u^{\e} \cdot \nabla_x \Big( \frac{3}{2}\ta^\mathsf{d}-\rho^\mathsf{d}\Big) + u^\mathsf{d} \cdot \nabla_x \Big(\frac{3}{2}\ta^E-\rho^E\Big) &=  \bar{\varPi}_{\mathfrak{s}}^{\e},
\enda
\end{align}
where $\bar{\varPi}_{\w}^{\e}$ and $\bar{\varPi}_{\mathfrak{s}}^{\e}$ are defined in \eqref{barvarPidef}.
\end{lemma}
\begin{remark}
In the infinite velocity energy case \eqref{caseECX}, the radial eddy $\bar{u}$ cancels out in the difference $u^\mathsf{d}=\mathbb{P}u^{\e}-u^E$.
\end{remark}

\begin{proof}[Proof of Lemma \ref{L.deqn}]
For \eqref{barueqn}, we apply the Leray projector $\mathbb{P}$ to the equation for $u^{\e}$ in $\eqref{locconNew}_2$ and subtract the Euler equation for $u^{E}$ in \eqref{incompE} to have
\begin{align*}
\p_t(\mathbb{P}u^{\e}-u^{E}) &+\mathbb{P}\big(\mathbb{P}u^{\e}\cdot \nabla_x \mathbb{P}u^{\e}\big) - u^E\cdot\nabla_x u^E 
=\mathbb{P}\bigg[\frac{k_B(1-\mathrm{\Theta}^{\e})}{\e}\nabla_x(\rho^{\e}+\ta^{\e})- \frac{1}{\e^2}\frac{1}{\mathrm{P}^{\e}}\sum_{j} \p_{x_j} \mathbf{r}_{ij}^{\e}\bigg] \cr 
& + \nabla_xp^E - \mathbb{P}\Big(\mathbb{P}^{\perp}u^{\e}\cdot \nabla_x \mathbb{P}u^{\e} +\mathbb{P}u^{\e}\cdot \nabla_x \mathbb{P}^{\perp}u^{\e} + \mathbb{P}^{\perp}u^{\e}\cdot \nabla_x \mathbb{P}^{\perp}u^{\e}\Big),
\end{align*}
where we decomposed $u^{\e} = \mathbb{P}u^{\e} + \mathbb{P}^{\perp}u^{\e}$ for both factors in $u^{\e}\cdot\nabla_x u^{\e}$. 

Subtracting \eqref{weqnE} from \eqref{weqnP} yields \eqref{barweqn}. Adding $\eqref{incompE}_2$ and \eqref{rtaeqnP} gives \eqref{barrtaeqn}.

\end{proof}

\begin{definition}
We define the local Maxwellian associated with the macroscopic variables of the incompressible Euler equations as
\begin{align}\label{MEdef}
(\mathrm{P}^E, \mathrm{U}^E, \mathrm{\Theta}^E) := (e^{\e\rho^E}, \e u^E, e^{\e\ta^E}), \qquad M^E(t,x,v) 
:= \frac{\mathrm{P}^E}{(2\pi k_B \mathrm{\Theta}^E)^{\frac{3}{2}}} e^{-\frac{|v - \mathrm{U}^E|^2}{2k_B \mathrm{\Theta}^E}}.
\end{align}
\end{definition}
Recall that we decomposed $F^{\e}$ into a macroscopic part and a microscopic part as $F^{\e}=M^{\e}+\AC{\P}F^{\e}$. From Lemma \ref{L.MFME}, we have
\begin{align}\label{F-ME}
\bega
\frac{1}{\e^2}\frac{|\p^{\al}(F^{\e}(t)-M^E(t))|^2}{M^{\e}(t)} &\leq \frac{1}{\e^2}\frac{|\p^{\al}\AC{\P}F^{\e}(t)|^2}{M^{\e}(t)} +|\p^{\al}(\rho^{\e}-\rho^E,u^{\e}-u^E,\ta^{\e}-\ta^E)|^2 + O(\e^2).
\enda
\end{align}
For the macroscopic variable $u^{\e}$, we decompose it using the Hodge decomposition \ref{Hodge}. For $(\rho^{\e},\ta^{\e})$, we decompose them into the primitive variables and the thermodynamic variable via \eqref{thermoV}, respectively:
\begin{align*}
\bega
&u^{\e}-u^E = (\mathbb{P}u^{\e}-u^{E}) + \mathbb{P}^{\perp}u^{\e}, \cr
&|\rho^{\e}-\rho^E|^2+|\ta^{\e}-\ta^E|^2 \leq C\bigg(\Big(\frac{3}{2}\ta^{\e}-\rho^{\e}\Big)-\Big(\frac{3}{2}\ta^E-\rho^E\Big)\bigg)^2 + C\Big((\rho^{\e}+\ta^{\e})-(\rho^{E}+\ta^{E})\Big)^2.
\enda
\end{align*}
We present the convergence of each component in different topologies in the following theorem.  
Recall the radial energy decomposition in Definition \ref{D.Ra-E} and that $\bar{u}\cdot\nabla_x \bar{\w}=0$.

\hide
{\color{red}
수정: 

1) 일단 $F^\e=M^\e+\AC{\P}F^{\e}$ 를 다시 이쿠에이션 넘버로 리마인드해주고  
2) $\AC{\P}F^{\e} $ 가 어떤 센스로 $0$으로 가는지 정확하게 써주세요.  
3) 그리고 primitive variables 과 thermodynamic variables \eqref{thermoV} 을 리마인드해주세요. 이어서 써놓은 $u, P^\perp u^\e$-컨벌젼스를 써주세요
4) 그리고 나머지 macroscopic variables 의 컨벌젼스는 primitive variables 에 대해 아래와 같이 보였다: (힌트: 간결하게 
\[
\frac{3}{2} \theta^\e - \rho^\e \to \frac{3}{2} \theta^E - \rho^E  \ \ \text{in $L^\infty(0,T; H^3(\Omega))$}
\]
\[
\rho^\e + \theta^\e \to \rho^E + \theta^E
\ \ \text{in $L^r(0,T; B.....)$}
\]
5) Theorem 이 너무 복잡하지 않게, 컨벌젼스나 바운드가 퀀티파이 되면 "with quantitative rate" 이라고 명시해주고, 그게 아니면 그냥 아무말 안하면 됩니다. "with quantitative rate" 인 경우는 리마크에 정확히 레잇을 써주세요
) 
}
\unhide

\hide
{\color{red} The following theorem would be perfect if we mention/correct a few things:
\begin{itemize}
    \begin{itemize}
        \item $H^2$ convergence of $u_0$ is enough or do we really need $H^3$ convergence of $u_0$ (so that $H^2$ convergence of $\omega^\e_0$)?

        \begin{itemize}
            \item detailed question1: For the convergence what do we need? $H^2$ convergence of $u_0$? or $H^3$ convergence of $u_0$? 

            \item For "By classical well-posedness for \eqref{incompE} with these initial data,
this limit is the unique classical solution.", what do we need? $H^2$ convergence of $u_0$? or $H^3$ convergence of $u_0$? 
        \end{itemize}

        \item We don't need convergence assumption for $(3/2 \theta^\e_0 - \rho^\e_0, \rho^\e + \theta^\e)$, right? If it is the case we should comment this (sort of no wellpreparedness) in the remark starts with "The initial data are said to be well-prepared"

        \item in the same remark starts with "The initial data are said to be well-prepared", we should say that we don't assume time derivative control at the initial data, which is another wellpreparedness.  
    \end{itemize}


\end{itemize}

}
\unhide

\hide
We define 
\begin{align}\label{mupm}
\mu_-(v) &:= M_{[1,0,1-c_0]} =  \frac{1}{(2\pi k_B (1-c_0))^{\frac{3}{2}}}\exp\left(-\frac{|v|^2}{2k_B(1-c_0)}\right), \cr 
\mu_+(v) &:= M_{[1,0,1+c_0]} =  \frac{1}{(2\pi k_B (1+c_0))^{\frac{3}{2}}}\exp\left(-\frac{|v|^2}{2k_B(1+c_0)}\right)
\end{align}
so that 
\begin{align*}
C_1\mu_-(v) \leq M^{\e} \leq C_2 \mu_+(v) 
\end{align*}
\unhide

\begin{theorem}\label{T.C.Hk}
Let $\Omega=\R^d$ with $d=2$ or $3$.
Assume a family of initial data $\{F_0^{\e}\}_{\e>0}$ is uniformly
bounded in $\e$ componentwisely, in the following sense:
\begin{align*}
\sup_{\e>0}
\bigg(
\|\w_0^{\e}\|_{H^2(\R^d)}
+ \|u_0^{\e}-\bar{u}\|_{L^2(\R^d)}
+ \|(\rho_0^{\e},\ta_0^{\e})\|_{H^3(\R^d)}
+ \sum_{0\leq \ell \leq3}\left\|\frac{1}{\e}\frac{ \nabla_x^\ell \AC{\P}F^{\e}_0}{\sqrt{M^{\e}_0}}\right\|_{L^2(\R^d_x  \times \R^3_v)} 
\bigg)
<+\infty.
\end{align*}
Here, 
$\bar{u}$
denotes the radial eddy defined in Definition~\ref{D.Ra-E} when $d=2$,
while we set $\bar{u}=0$ for $d=3$. Additionally, the family $\{F_0^{\e}\}_{\e>0}$ satisfies the 
$4$-Admissible Blow-up Condition for $d=2$ or the $5$-Admissible Blow-up Condition for $d=3$ in~\eqref{ABC1} of Definition~\ref{ABC}.


We assume that, as $\e \to 0$,
\begin{equation*}
 \big(\mathbb{P}u^\e_0 -\bar{u} ,
 \ \tfrac{3}{2}\theta^\e_0-\rho^\e_0
 \big)
 \to 
  \big(\mathbb{P}u_0 -\bar{u},
  \ \tfrac{3}{2}\theta_0-\rho_0
  \big) \ \   \text{in $H^3(\R^d_x) 
  \times 
  H^3(\R^d_x). 
  $} 
\end{equation*}

Then the following statements hold for a family of Boltzmann solutions $\{F^{\e}\}_{\e>0}$ to~\eqref{BE} with $\kappa=\e^q$ for some $0<q<2$, on
the time interval $t\in[0,T_\e]$ with $T_\e\to\infty$ as $\e\to0$ when $d=2$, and for $t\in[0,T]$ with $T\ll1$ when $d=3$ (see our constructions in Theorem~\ref{T.2D.global} for $d=2$ and Theorem~\ref{T.3D.unif} for $d=3$).

\begin{enumerate}

\item[(1)]
The family of macroscopic fields associated with the Boltzmann solution $\big\{\big(\mathbb{P}u^\e,\ \mathbb{P}^{\perp}u^\e,\ \tfrac{3}{2}\theta^\e-\rho^\e,\ \rho^\e+\theta^\e\big)\big\}_{\e>0}$ converges, as $\e\to 0$, to $\big(u^E,\ 0,\ \tfrac{3}{2}\theta^E-\rho^E,\ 0\big)$ in the respective topologies specified below.
This convergence yields (and hence provides a construction of) a limiting triple
$(\rho^E,u^E,\ta^E)$ which is a unique classical solution of the incompressible Euler system \eqref{incompE} 
with initial data $(\rho^E,u^E,\ta^E)|_{t=0}=(\rho_0,\mathbb{P}u_0,\ta_0)$.

\item[(2)] The topological spaces in which convergence holds are as follows:
\begin{itemize}
\item[(i)] The solenoidal part of the velocity field converges as follows:
\begin{align}
&\mathbb{P}u^{\e} -u^E \to 0  \ \ \text{in  $L^\infty (0, T; H^2(\R^d))$ with rate, }
\label{TC.wH2-0}\\
&\mathbb{P}u^{\e} -u^E \to 0  \ \ \text{in  $L^\infty (0, T; \dot{H}^3(\R^d))$. }
\label{TC.wH3-0}
\end{align}

\item[(ii)] The entropic fluctuation converges as follows:
\begin{equation}
     \frac{3}{2} \theta^\e - \rho^\e  \to \frac{3}{2} \theta^E - \rho^E  \ \  \text{in  $L^\infty (0, T; {H}^{k}(\R^d))$,}  \label{conv_entr}
\end{equation}
for $k=2,3$ (in particular, with a quantitative rate for $k=2$). 

\item[(iii)] The irrotational part and pressure fluctuation vanish as follows:
\begin{align}
 \mathbb{P}^{\perp}u^{\e}, \ \rho^{\e}+\ta^{\e}    \to  0  \ \ &\text{in $ L^r (0,T; \dot{B}_{p,1}^{s+d(\frac{1}{p}-\frac{1}{2})+\frac{1}{r}}(\R^d))$ with rate $\e^{\frac{(d-1)-}{4}}$} ,\label{conv_irrorhota} 
\end{align}
for $2<p<\infty$ and $\frac{1}{r}\leq d(\frac{1}{2}-\frac{1}{p})$ and for any $s\in(0,3)$ in \eqref{conv_irrorhota}.

\end{itemize}
\item[(3)] The microscopic part converges to zero:
\begin{equation}\label{TC.H2.G}
  \frac{1}{\e^2}  \sqrt{\nu^{\e}} \frac{\nabla_x^\ell \AC{\P}F^{\e} }{\sqrt{M^{\e}}}   \to   0   
  \ \ \text{in $L^2\big(0,T;\, L^2( \R^d_x \times  \R^3_v)\big)
$ for all $0 \leq \ell \leq 3$ with rate $\kappa^{\frac{1}{2}}$}.
\end{equation}

\end{enumerate}
\end{theorem}

\begin{remark}
We summarize the quantitative convergence rates in the Big-O sense as follows:
\begin{align*}
 \|\mathbb{P}u^{\e}_0-u^E_0\|_{H^2_x}   
+  \Big(\e^{\frac{(d-1)-}{4}-} +\kappa^{\frac{1}{2}-}\Big)  
\ \ &\text{in \eqref{TC.wH2-0}},\\ 
 \bigg\|\Big(\frac{3}{2}\ta^{\e}_0-\rho^{\e}_0\Big)
-\Big(\frac{3}{2}\ta^E_0-\rho^E_0\Big)\bigg\|_{H^2_x}  
+   \Big(\e^{\frac{(d-1)-}{4}-} +\kappa^{\frac{1}{2}-}\Big) 
\ \ &\text{in \eqref{conv_entr} for $H^2(\R^d)$}.
\end{align*}
\end{remark}

\begin{remark}
Our approach does not assume the a priori existence of solutions to the incompressible Euler equations. This result provides a existence theorem that constructs fluid solutions directly from kinetic solutions.
\end{remark}

\begin{remark}
The initial data are said to be \emph{well-prepared} if the microscopic part converges to zero faster than the Mach number in any topology, namely, if $\frac{1}{\e}\AC{\P}F_0 \to 0$. 
In this sense, our analysis does not impose such well-preparedness assumption on the initial data. 
We do not assume that the acoustic variables $\mathbb{P}^{\perp} u^{\e}_0$ and $\rho^{\e}_0+\theta^{\e}_0$ vanish at the initial time; hence, the initial data are not well-prepared in this sense. Moreover, we do not impose any control on the time derivatives at the initial data, which is another form of well-preparedness commonly assumed in the literature.
\end{remark}

\begin{corollary} For $d=2$ and $d=3$, consider the following initial data satisfying
\begin{align*}
\bega
&\w_0\in H^2(\R^d), \quad  (u_0-\bar{u}) \in L^2(\R^d), \quad \rho_0, \ta_0 \in H^3(\R^d), \quad \frac{1}{\e}\frac{\AC{\P}F_0}{\sqrt{\tilde{\mu}}}\in H^{3}_x(\R^d ; L^2_v(\R^3)), \cr 
&\bigg\|\frac{1}{\e}\frac{(F_0-\mu)}{\sqrt{\tilde{\mu}}}\bigg\|_{L^\infty_{x,v}}\les \frac{1}{\e^{1-}},
\enda
\end{align*}
where $\tilde{\mu}=M_{[1,0,1-c_0]}$ for some $0<c_0\ll1$.
The mollified sequence $\{F^{\e}_0\}_{\e>0}$ constructed from Lemma \ref{L.molli} satisfies the assumption of Theorem \ref{T.C.Hk}, and hence all the results of Theorem \ref{T.C.Hk} hold. 
\end{corollary}

\hide
{\color{blue}[찬우: 윗부분 스테이트먼트가 의미를 정확히 전달하지 않습니다.  우리가 보인것은 “볼츠만의 $u^\e, \rho^\e + \theta^\e$ 등등이 유체 해의 존재성과 관련없이 일단 수렴을 하고, 그 수렴한 $(u, \rho, \theta)$ 가 유체의 해다” 라는 것을 보인거죠. 2D 에서는 임의의 시간ㅇ에 대해, 그리고 3D 에서는 짧은 시간에서. 아래 정리에는 써놓은 방식은 마치 유체해의 존재를 가정해야 컨벌젼스가 있는 것 처럼 써놓아서 미스리딩하네요. Golse 가 우리랑 비슷한 개념의 수렴성을 많이  공부했으니, 그분이 어떻게 스테이트했는지 찾아보고 스테이트먼트 써보세요.  ]}

{\color{red}이런씩으로 써보죠. 

1) $F_0$ 의 컴포넌트와이즈 조건을 써주고 ($G|_{t=0}$ 의 조건 포함) 각각의 컴포넌트의 조건을 써주기
2) 해당되는 클래시칼 볼츠만 솔루션은 존재하여 (여기서는 존재-시간에 대해 써주지 맙시다)
3) 다음을 $ t\in [0,T]$ 에서 만족한다 where Here, $T$ is arbitrary large for $d=2$, while $T \ll1 $ for $d=3.$ 
4) $\AC{\P}F^{\e}\to 0$ 이고 (여기에 (1)을 써주면 되겠습니다)
5) $(P u^\e, P^\perp u^\e, \rho^\e+ \theta^\e, \frac{3}{2} \theta^\e - \rho^\e) \to (P u^E, P^\perp u^E, \rho^E+ \theta^E, \frac{3}{2} \theta^E - \rho^E)$ 로 수렴하고 이들은 다음의 incompressible Euler equations 를 만족한다:  방정식 써주기

6) Moreover, (2) 와 (3) 을 써주면 되겠습니다. 

띠어럼 아래에 첫 리마크: 우리의 접근법은 인컴프레서블 오일러 해의 존재성을 어 프라리오리 가정하지 않는다. 따라서 이 결과는 유체 해를 키네틱 해로부터 컨스트럭트하는 새로운 존재성 띠어럼이다.

두번째 리마크: 웰프리페어드니스를 얼마나 피했는지 논의}
\unhide

\hide
\begin{theorem}
Then the following statements hold: Fix $T>0$. As $\e \to 0$,
\begin{enumerate}

\item[(1)] The microscopic part converges to zero:
\begin{equation}\label{TC.H2.G}
  \frac{1}{\e^2} \frac{\p^\alpha \AC{\P}F^{\e} }{\sqrt{M^\e}}  \to   0   
  \ \ \text{in $L^2 (0,T; L^2 (\Omega \times \R^d))$ for all $|\alpha| \leq 3$ with rate}.
\end{equation}

\item[(2)] The solenoidal part of the velocity field converges as follows:
\begin{align}
&\mathbb{P}u^{\e} \to u^E  \ \ \text{in  $L^\infty (0, T; H^2(\R^d))$ with rate, }
\label{TC.wH2}\\
&\mathbb{P}u^{\e} \to u^E  \ \ \text{in  $L^\infty (0, T; \dot{H}^3(\R^d))$. }
\label{TC.wH3}
\end{align}
\item[(3)] The irrotational part and thermodynamic variables (see \eqref{thermoV}) converges as follows:
\begin{align}
 \mathbb{P}^{\perp}u^{\e}, \ \rho^{\e}+\ta^{\e}    \to  0  \ \ &\text{in $ L^r (0,T; \dot{B}_{p,1}^{s+2(\frac{1}{p}-\frac{1}{2})+\frac{1}{r}}(\R^d))$ with rate(?)} ,\label{conv_irrorhota}\\
 \frac{3}{2} \theta^\e - \rho^\e  \to \frac{3}{2} \theta^E - \rho^E  \ \ &\text{in  $L^\infty (0, T; {H}^{s}(\R^d))$ with rate(?),}\label{conv_entr}
\end{align}
for $s=2,3$, with a quantitative rate for $s=2$, and for $2<p<\infty$ and $\frac{1}{r}=d(\frac{1}{2}-\frac{1}{p})$. {\color{blue}여기까지}

\begin{align*}
\bega
\sum_{0\leq|\al|\leq 3}\frac{1}{\e^2}\int_0^T\int_{\Omega \times \R^3} \frac{|\p^{\al}\AC{\P}F^{\e}(t)|^2}{M^{\e}(t)} dvdxdt \leq C_T \e^2\kappa^{\frac{1}{2}}, \quad \text{as} \quad \e \rightarrow 0,
\enda
\end{align*}
where the constant \(C_T\) depends on $\int_0^T \mathcal{D}_{G}(t)dt$.

\item[(2)] The velocity converges as follows:
\begin{itemize}
\item[(i)] $H^2_x$ convergence. For all $t\in[0,T]$,
\begin{align}\label{TC.wH2}
\bega
\sup_{0 \leq t \leq T}&\|(\mathbb{P}u^{\e}-u^E)(t)\|_{H^2_x} \leq C_T\|(\mathbb{P}u^{\e}-u^E)(0)\|_{H^2_x}  + \mathcal{C}_* \Big(\e^{\frac{(d-1)-}{4}-} +\kappa^{\frac{1}{2}-}\Big).
\enda
\end{align}
Here, the constant $C_T$ depends only on $\|\mathbb{P}u_0-\bar{u}\|_{L^2_x}$ and $\|\w_0\|_{H^2_x}$,
and the constant $\mathcal{C}_*$ is defined in \eqref{C*def}.
\item[(ii)] $H^3_x$ convergence. For all $t\in[0,T]$,
\begin{align}\label{TC.wH3}
\bega
\sup_{0 \leq t \leq T}&\|(\mathbb{P}u^{\e}-u^E)(t)\|_{H^3_x} \to 0, \quad \mbox{as} \quad \e\to 0.
\enda
\end{align}
\item[(iii)] The compressible part $\mathbb{P}^{\perp}u^{\e}$ converges to zero:
\begin{align*}
\bega
\mathbb{P}^{\perp}u^{\e} \quad \to \quad 0 \quad \mbox{in} \quad L^r_T\dot{B}_{p,1}^{s+2(\frac{1}{p}-\frac{1}{2})+\frac{1}{r}}(\R^d)
\enda
\end{align*}
for $2<p<\infty$ and $\frac{1}{r}=d(\frac{1}{2}-\frac{1}{p})$.
\end{itemize}

\item[(3)] For the remaining macroscopic variables, the primitive variables and thermodynamic variables defined in \eqref{thermoV} satisfy

\begin{align*}
\bega
&\sup_{0 \leq t \leq T}\bigg\|\Big(\rho^{\e}-\frac{3}{2}\ta^{\e}\Big)(t)-\Big(\rho^E-\frac{3}{2}\ta^E\Big)(t)\bigg\|_{H^s_x} \to 0, \quad \mbox{as} \quad \e\to 0, \cr 
&\rho^{\e}+\ta^{\e} \quad \to \quad 0 \quad \mbox{in} \quad L^r_T\dot{B}_{p,1}^{s+2(\frac{1}{p}-\frac{1}{2})+\frac{1}{r}}(\R^d),
\enda
\end{align*}
for $s=2,3$, with a quantitative rate for $s=2$, and for $2<p<\infty$ and $\frac{1}{r}=d(\frac{1}{2}-\frac{1}{p})$.

\end{enumerate}

If we further assume $\rho_0, \ta_0 \in H^3(\R^d)$, let $(\rho^{E},u^E,\ta^E)$ be the unique solution of the incompressible Euler system \eqref{incompE} with initial data $(\rho_0,u_0,\ta_0)$, and let $F^{\e}$ be the Boltzmann solution corresponding to the mollified initial data $F^{\e}_0$ in the same sense as above. Then we have the following convergences:
\begin{enumerate}
\item[(2)] The incompressible part converges in $H^2(\R^d)$:
\begin{align*}
\bega
\sup_{0 \leq t \leq T}&\|(\mathbb{P}u^{\e}-u^E)(t)\|_{H^s_x} +\sup_{0 \leq t \leq T}\bigg\|\Big(\rho^{\e}-\frac{3}{2}\ta^{\e}\Big)(t)-\Big(\rho^E-\frac{3}{2}\ta^E\Big)(t)\bigg\|_{H^s_x} \to 0, \quad \mbox{as} \quad \e\to 0.
\enda
\end{align*}
for $s=2,3$ and with quantitative rate for $s=2$.
\item[(3)] The compressible part converges to 0.
\begin{align*}
\bega
(\rho^{\e}+\ta^{\e}),\mathbb{P}^{\perp}u^{\e} \quad \to 0 \quad \mbox{in} \quad L^r_T\dot{B}_{p,1}^{s+2(\frac{1}{p}-\frac{1}{2})+\frac{1}{r}}(\R^d)
\enda
\end{align*}
for $2<p<\infty$ and $\frac{1}{r}=d(\frac{1}{2}-\frac{1}{p})$.
\end{enumerate}
\end{theorem}
\unhide 

\hide
We can conclude the claim because the microscopic part converges to $0$ by \eqref{Gto0}, and the macroscopic part satisfies
\begin{align*}
\bega
\sum_{0\leq|\al_x|\leq s}\frac{1}{\e^2}\int_{\Omega \times \R^3} \frac{|\p^{\al_x}(M^{\e}(t)-M^E(t))|^2}{M^{\e}(t)} dvdx \les (1+C\e)\|(\rho^\mathsf{d},u^\mathsf{d},\ta^\mathsf{d})(t)\|_{H^s_x}^2,
\enda
\end{align*}
for \(s = 2,3\).
\unhide 


\begin{remark}
To obtain $H^3_x$ convergence, it is necessary to control $\mathbb{P}^{\perp}u^{\e}$ and $\rho^{\e}+\ta^{\e}$ up to third-order derivatives.  
The number of required derivatives depends on the spatial dimension: the $4$-admissible blow-up condition in two dimensions and the $5$-admissible blow-up condition in three dimensions, due to the Strichartz estimates in Proposition~\ref{P.div.ineq}.
\end{remark}

In Theorem~\ref{T.2D.global}, the solution satisfying the Admissible Blow-up 
Condition~\eqref{ABC1} also satisfies \eqref{incomp.2D} and 
\eqref{Gto0}. These correspond to \eqref{TC.H2.G} and 
\eqref{conv_irrorhota}, respectively.
Hence, it suffices to focus on the convergence of the vorticity and the thermodynamic variable 
$\mathfrak{s}^{\e}=\frac{3}{2}\ta^{\e}-\rho^{\e}$. 
Once we establish the following proposition, we obtain \eqref{TC.wH2-0}, \eqref{TC.wH3-0}, and \eqref{conv_entr} in Theorem \ref{T.C.Hk}.




\begin{proposition}\label{P.C.Hk}
Suppose the same assumptions as in Theorem \ref{T.C.Hk}. 
Then the following convergences hold:
\begin{enumerate}
\item[(1)] $H^2_x$ convergence. For all $t\in[0,T]$,
\begin{align}\label{TC.wH2}
\bega
&\sup_{0 \leq t \leq T}\Big(\|(\mathbb{P}u^{\e}-u^E)(t)\|_{L^2_x}+\|(\w^{\e}-\w^E)(t)\|_{H^1_x}\Big) \cr 
&\leq C_T\Big(\|(\mathbb{P}u^{\e}-u^E)(0)\|_{L^2_x}+\|(\w^{\e}-\w^E)(0)\|_{H^1_x}\Big)   + \mathcal{C}_* \Big(\e^{\frac{(d-1)-}{4}-} +\kappa^{\frac{1}{2}-}\Big).
\enda
\end{align}
Here, $C_T$ depends only on the right-hand sides of \eqref{BEC}, and $\mathcal{C}_*$ is defined in \eqref{C*def}.

\item[(2)] $H^3_x$ convergence. For all $t\in[0,T]$,
\begin{align}\label{TC.wH3}
\bega
\sup_{0 \leq t \leq T}\big(&\|(\mathbb{P}u^{\e}-u^E)(t)\|_{L^2_x} +\|(\w^{\e}-\w^E)(t)\|_{H^2_x}\big) \to 0, \quad \mbox{as} \quad \e\to 0.
\enda
\end{align}
\end{enumerate}
For the thermodynamic variable, the following convergence results hold:
\begin{enumerate}
\item[(1)] $H^2_x$ convergence. For all $t\in[0,T]$,
\begin{align*}
\bega
\sup_{0 \leq t \leq T}&\bigg\|\Big(\frac{3}{2}\ta^{\e}-\rho^{\e}\Big)(t)-\Big(\frac{3}{2}\ta^E-\rho^E\Big)(t)\bigg\|_{H^2_x} \cr 
&\leq C_T \bigg\|\Big(\frac{3}{2}\ta^{\e}_0-\rho^{\e}_0\Big)-\Big(\frac{3}{2}\ta^E_0-\rho^E_0\Big)\bigg\|_{H^2_x}  + \mathcal{C}_* \Big(\e^{\frac{(d-1)-}{4}-} +\kappa^{\frac{1}{2}-}\Big).
\enda
\end{align*}
Here, $C_T$ depends only on the right-hand sides of \eqref{BEC} and \eqref{C}, and $\mathcal{C}_*$ is defined in \eqref{C*def}.

\item[(2)] $H^3_x$ convergence. For all $t\in[0,T]$,
\begin{align*}
\bega
\sup_{0 \leq t \leq T}\bigg\|\Big(\frac{3}{2}\ta^{\e}-\rho^{\e}\Big)(t)-\Big(\frac{3}{2}\ta^E-\rho^E\Big)(t)\bigg\|_{H^3_x} \to 0, \quad \mbox{as} \quad \e\to 0.
\enda
\end{align*}
\end{enumerate}
\end{proposition}

\hide
\textcolor{blue}{ (Convergence in $L^\infty_t$ )
\begin{remark}\label{R.M-M}
If we have $\|\nabla_x(\rho+\ta)\|_{L^\infty_tH^s_x}\les \kappa^{\delta-\frac{s+1}{2}}$, then we can obtain the following $L^\infty_t$ convergence for the macroscopic parts: 
\begin{enumerate}
\item The macroscopic part of the Boltzmann equation converges to $M^E$ in $H^3_x$ as follows: 
\begin{align}\label{M-M,CH3}
\bega
\sum_{0\leq |\al_x| \leq 3}&\frac{1}{\e^2}\bigg\| \int_{\Omega \times \R^3} \frac{|\p^{\al_x}(M^{\e}-M^E)|^2}{M^{\e}} dvdx\bigg\|_{L^\infty_t} \rightarrow 0, \quad \mbox{as} \quad \kappa \rightarrow 0.
\enda
\end{align}
\item The solution of the Boltzmann equation converges to $M^E$ in $H^2_x$ with the following rate
\begin{align}\label{M-M,CH2}
\bega
\sum_{0\leq |\al_x| \leq 2}&\frac{1}{\e^2}\bigg\| \int_{\Omega \times \R^3} \frac{|\p^{\al_x}(M^{\e}-M^E)|^2}{M^{\e}} dvdx\bigg\|_{L^\infty_t} \les e^{Ct}\bigg(\|\w^\mathsf{d}(0)\|_{H^1_x}^2 +\Big\|\Big(\rho^\mathsf{d}-\frac{3}{2}\ta^\mathsf{d}\Big)(0)\Big\|_{H^2_x}^2 \bigg) \cr 
&+e^{Ct}(\kappa^{2\delta-1}+\kappa) +\kappa^{2\delta-2} + \int_0^t e^{C(t-s)} BE~part(s)ds.
\enda
\end{align}
\end{enumerate}
Remark: $L^\infty_t$ convergence is $\kappa^{1/2}$ worse than $L^2_t$ estimate. 
\end{remark} }
\unhide

\subsubsection{Convergence in \(H^2\)}

Before proving Proposition \ref{P.C.Hk}, we first note that, by \eqref{U.L2}, \eqref{U.L2inf}, \eqref{U.PuL2}, \eqref{U.uwgrow}, and \eqref{U.rtaH3} in Lemma \ref{L.unif}, together with the assumptions $\w_0^{\e}\in H^2(\R^d)$ and $\w_0^E \in H^2(\R^d)$, we obtain   
\begin{align}\label{BEC}
\bega
&\sup_{0 \leq t \leq T} \|\w^{\e}(t)\|_{H^2_x} + \sup_{0 \leq t \leq T} \|\w^E(t)\|_{H^2_x} \leq C_T(\|\mathbb{P}u^{\e}_0-\bar{u}\|_{L^2_x},\|\w^{\e}_0\|_{H^2_x},\|u^{E}_0-\bar{u}\|_{L^2_x},\|\w^{E}_0\|_{H^2_x}).
\enda
\end{align}
If we further assume $\mathfrak{s}^{\e}_0= \frac{3}{2}\ta^{\e}_0-\rho^{\e}_0 \in H^3(\R^d)$ and $\mathfrak{s}^E_0= \frac{3}{2}\ta^E_0-\rho^E_0 \in H^3(\R^d)$, then we also have 
\begin{align}\label{C}
\bega
&\sup_{0 \leq t \leq T} \|\mathfrak{s}^{\e}(t)\|_{H^3_x} + \sup_{0 \leq t \leq T} \|\mathfrak{s}^E(t)\|_{H^3_x}  \leq C_T(\|(u^{\e}_0-\bar{u},\mathfrak{s}^{\e}_0)\|_{H^3_x},\|(u^{E}_0-\bar{u},\mathfrak{s}^E_0)\|_{H^3_x}).
\enda
\end{align}

We begin the proof of Proposition \ref{P.C.Hk} by deriving energy estimates for 
$\|u^{\mathsf{d}}(t)\|_{L^2_x}$, $\|\w^\mathsf{d}(t)\|_{H^{s-1}_x}$ and 
$\|(\tfrac{3}{2}\ta^\mathsf{d}-\rho^\mathsf{d})(t)\|_{H^s_x}$ 
for $s=2,3$.

\begin{lemma}\label{L.Pu-uconv}
Let $\Omega=\R^d$ with $d=2$ or $3$, and let $u^\mathsf{d}$ be the solution to \eqref{barueqn}. Then the following inequality holds:
\begin{align}\label{u-uL2C}
\bega
\frac{d}{dt}\|u^{\mathsf{d}}(t)\|_{L^2_x} \leq \|\w^{\e}(t)\|_{H^2_x}\|u^{\mathsf{d}}(t)\|_{L^2_x} + \mathfrak{Q}(t),
\enda
\end{align}
where
\begin{align}\label{udforce}
\bega
\mathfrak{Q}(t)  &:= \bigg\| \frac{k_B(1-\mathrm{\Theta}^{\e}(t))}{\e}\bigg\|_{L^2_x} \|\nabla_x(\rho^{\e}+\ta^{\e})(t)\|_{L^\infty_x}+ \frac{1}{\e^2}\bigg\|\frac{1}{\mathrm{P}^{\e}}\sum_{j} \p_{x_j} \mathbf{r}_{ij}^{\e}(t)\bigg\|_{L^2_x} \cr
&+ \Big\|\Big(\mathbb{P}^{\perp}u^{\e}\cdot \nabla_x \mathbb{P}u^{\e} +\mathbb{P}u^{\e}\cdot \nabla_x \mathbb{P}^{\perp}u^{\e} + \mathbb{P}^{\perp}u^{\e}\cdot \nabla_x \mathbb{P}^{\perp}u^{\e}\Big)(t) \Big\|_{L^2_x} .
\enda
\end{align}
\end{lemma}
\begin{proof}
Writing the energy estimate for $u^{\mathsf{d}}$ from \eqref{barueqn} gives 
\begin{align}\label{udE}
\bega
\frac{1}{2}\frac{d}{dt}\|u^{\mathsf{d}}\|_{L^2_x}^2 &\leq  \int_{\Omega} \bigg(\mathbb{P}\big(\mathbb{P}u^{\e}\cdot \nabla_x \mathbb{P}u^{\e}\big) - u^E\cdot\nabla_x u^E\bigg) u^{\mathsf{d}} dx + \mathfrak{Q}(t) \|u^{\mathsf{d}}\|_{L^2_x}.
\enda
\end{align}
For the first term on the right-hand side of \eqref{udE}, since $\mathbb{P}$ is self-adjoint and $\mathbb{P}(\mathbb{P}u^{\e}-u^E)=\mathbb{P}u^{\e}-u^E$, we have
\begin{align*}
\bega
\int_{\Omega} \bigg(\mathbb{P}\big(\mathbb{P}u^{\e}\cdot \nabla_x \mathbb{P}u^{\e}\big) - u^E\cdot\nabla_x u^E\bigg) u^{\mathsf{d}} dx 
&= \int_{\Omega} \big( u^{\mathsf{d}}\cdot \nabla_x \mathbb{P}u^{\e}\big) u^{\mathsf{d}} dx + \int_{\Omega} \big(u^E\cdot \nabla_x u^{\mathsf{d}} \big)u^{\mathsf{d}} dx \cr 
&\leq \int_{\Omega} |\nabla_x \mathbb{P}u^{\e}| |u^{\mathsf{d}}|^2 dx.
\enda
\end{align*}
Here, the last term on the first line vanishes because $u^E$ is divergence-free.
Then, taking $L^\infty_x$ to $\nabla_x\mathbb{P}u^{\e}$ then applying from Agmon's inequality \eqref{Agmon}, we get 
\begin{align*}
\bega
\int_{\Omega} |\nabla_x \mathbb{P}u^{\e}| |u^{\mathsf{d}}|^2 dx &\leq C\|\nabla_x\mathbb{P}u^{\e}\|_{H^2_x} \|u^{\mathsf{d}}\|_{L^2_x}^2 \leq C\|\w^{\e}\|_{H^2_x} \|u^{\mathsf{d}}\|_{L^2_x}^2,
\enda
\end{align*}
where we used $\|\nabla_x\mathbb{P}u^{\e}\|_{H^2_x}\leq C\|\w^{\e}\|_{H^2_x}$ from \eqref{PotenLp}. Then, dividing each side of \eqref{udE} by $\|u^{\mathsf{d}}\|_{L^2_x}$ gives the result. 
\end{proof}

\begin{lemma}\label{L.wrtaHsC}
Let $\Omega=\R^d$ with $d=2$ or $3$.
Let $\w^\mathsf{d}$ and $\rho^\mathsf{d}-\tfrac{3}{2}\ta^\mathsf{d}$ denote the solutions to \eqref{barweqn} and \eqref{barrtaeqn}, respectively.  
If these solutions satisfy \eqref{BEC} and \eqref{C}, then the following inequalities hold:
\begin{align}\label{wH2C}
\bega
\frac{d}{dt}&\|\w^\mathsf{d}\|_{H^{s}_x} \les  \|\nabla_x\mathbb{P}u^{\e}\|_{L^\infty_x} \|\w^\mathsf{d}\|_{H^s_x} + \|\w^{\e}\|_{L^\infty_x} \|\nabla_xu^\mathsf{d}\|_{H^s_x} + \Big(\|u^\mathsf{d}\|_{H^{s+1}_x}\|\w^E\|_{L^\infty_x}+\|u^\mathsf{d}\|_{L^\infty_x}\|\w^E\|_{H^{s+1}_x}\Big)   \cr 
&+ \delta_{s=1}\Big(\|\w^\mathsf{d}\|_{H^1_x}\|\nabla_x\mathbb{P}u^{\e}\|_{H^2_x} + \|\nabla_x\w^{\e}\|_{H^1_x} \|\nabla_xu^\mathsf{d}\|_{H^1_x} \Big) \cr 
&+\delta_{s=2}\Big(\|\nabla_x\mathbb{P}u^{\e}\|_{H^s_x}\|\w^\mathsf{d}\|_{L^\infty_x} +\|\w^\mathsf{d}\|_{L^\infty_x}\|\nabla_x\mathbb{P}u^{\e}\|_{H^s_x} + \|\w^{\e}\|_{H^s_x} \|\nabla_xu^\mathsf{d}\|_{L^\infty_x}\Big)\cr 
&+ \sum_{0\leq|\al_x|\leq s}\|\p^{\al_x}\nabla_x\mathbb{P}^{\perp}u^{\e}\|_{L^\infty_x}\big(\|\w^{\e}\|_{H^s_x}+\|\w^\mathsf{d}\|_{H^s_x}\big) + \|\bar{\varPi}_{\w}^{\e}\|_{H^s_x},
\enda
\end{align}
for \(s=1,2\), and 
\begin{align}\label{rtaH3C} 
\frac{d}{dt}&\Big\|\frac{3}{2}\ta^\mathsf{d}-\rho^\mathsf{d}\Big\|_{H^s_x}  \les \|\nabla_x \mathbb{P}u^{\e}\|_{L^\infty_x}\Big\|\Big(\frac{3}{2}\ta^\mathsf{d}-\rho^\mathsf{d}\Big)\Big\|_{H^s_x} \cr 
&+\delta_{s=2}\|\nabla_x \mathbb{P}u^{\e}\|_{H^2_x}\Big\|\nabla_x\Big(\frac{3}{2}\ta^\mathsf{d}-\rho^\mathsf{d}\Big)\Big\|_{H^2_x} +\delta_{s=3}\|\nabla_x\mathbb{P}u^{\e}\|_{H^{s-1}_x}\Big\|\nabla_x\Big(\frac{3}{2}\ta^\mathsf{d}-\rho^\mathsf{d}\Big)\Big\|_{L^\infty_x}  \cr 
&+ \bigg(\|u^\mathsf{d}\|_{L^\infty_x}\Big\|\nabla_x\Big(\frac{3}{2}\ta^E-\rho^E\Big)\Big\|_{H^{s}_x} +\|u^\mathsf{d}\|_{H^{s}_x}\Big\|\nabla_x\Big(\frac{3}{2}\ta^E-\rho^E\Big)\Big\|_{L^\infty_x} \bigg) +\big\|\bar{\varPi}_{\mathfrak{s}}^{\e} \big\|_{H^s_x}, 
\end{align}
for \(s=2,3\), where $\bar{\varPi}_{\w}^{\e}$ and $\bar{\varPi}_{\mathfrak{s}}^{\e}$ are defined in \eqref{barvarPidef}.
\end{lemma}
\begin{proof}
(1) Writing the $H^s_x$ energy estimate for $\w^\mathsf{d}$ from \eqref{barweqn} yields
\begin{align}\label{wdHs}
\bega
\frac{d}{dt}\frac{1}{2}\|\w^\mathsf{d}\|_{H^{s}_x}^2  &\leq I_{\w^\mathsf{d}}+II_{\w^\mathsf{d}}+III_{\w^\mathsf{d}}+IV_{\w^\mathsf{d}}+ \|\w^\mathsf{d}\|_{H^2_x} \|\bar{\varPi}_{\w}^{\e}\|_{H^2_x},
\enda
\end{align}
where 
\begin{align*}
\bega
&I_{\w^\mathsf{d}}:= \sum_{0\leq|\al_x|\leq s}\int_{\Omega} \p^{\al_x}(\mathbb{P}u^{\e}\cdot\nabla_x \w^\mathsf{d})\p^{\al_x}\w^\mathsf{d}dx, \quad II_{\w^\mathsf{d}}:= \sum_{0\leq|\al_x|\leq s}\int_{\Omega} \p^{\al_x}(u^\mathsf{d}\cdot\nabla_x \w^E)\p^{\al_x}\w^\mathsf{d}dx, \cr
&III_{\w^\mathsf{d}}:= \sum_{0\leq|\al_x|\leq s}\int_{\Omega} \p^{\al_x}(\w^\mathsf{d}\cdot\nabla_x u^{\e})\p^{\al_x}\w^\mathsf{d}dx, \quad IV_{\w^\mathsf{d}}:= \sum_{0\leq|\al_x|\leq s}\int_{\Omega} \p^{\al_x}(\w^{\e}\cdot\nabla_x (u^{\e}-u^E))\p^{\al_x}\w^\mathsf{d}dx.
\enda
\end{align*}
For $I_{\w^\mathsf{d}}$, since $\nabla_x\cdot \mathbb{P}u^{\e}=0$, the contribution where all derivatives fall on $\nabla_x\w^{\mathsf{d}}$ cancels. Using \eqref{uvHk}, we get
\begin{align*}
\bega
|I_{\w^\mathsf{d}}| &\les \begin{cases}
\|\nabla_x\mathbb{P}u^{\e}\|_{L^\infty_x}\|\w^\mathsf{d}\|_{H^s_x}\|\w^\mathsf{d}\|_{H^s_x} , \quad &s=1, \\
\Big(\|\nabla_x\mathbb{P}u^{\e}\|_{H^s_x}\|\w^\mathsf{d}\|_{L^\infty_x} +\|\nabla_x\mathbb{P}u^{\e}\|_{L^\infty_x}\|\w^\mathsf{d}\|_{H^s_x}\Big)\|\w^\mathsf{d}\|_{H^s_x} , \quad &s=2,
\end{cases} \cr 
&\les \|\nabla_x\mathbb{P}u^{\e}\|_{L^\infty_x} \|\w^\mathsf{d}\|_{H^s_x}^2 + \delta_{s=2}\|\nabla_x\mathbb{P}u^{\e}\|_{H^s_x}\|\w^\mathsf{d}\|_{L^\infty_x}\|\w^\mathsf{d}\|_{H^s_x}.
\enda
\end{align*}
For $II_{\w^\mathsf{d}}$, applying \eqref{uvHk} to $\|u^\mathsf{d}\w^E\|_{H^{s+1}_x}$ yields
\begin{align}\label{uwE}
\bega
|II_{\w^\mathsf{d}}| &\les (\|u^\mathsf{d}\|_{H^{s+1}_x}\|\w^E\|_{L^\infty_x}+\|u^\mathsf{d}\|_{L^\infty_x}\|\w^E\|_{H^{s+1}_x})\|\w^\mathsf{d}\|_{H^{s}_x}. 
\enda
\end{align}
For $III_{\w^\mathsf{d}}$, decompose $u^{\e}= \mathbb{P}u^{\e}+ \mathbb{P}^{\perp}u^{\e}$ and use \eqref{uvHk}:
\begin{align*}
\bega
|III_{\w^\mathsf{d}}| &\leq
\begin{cases}
\big(\|\w^\mathsf{d}\|_{H^1_x}\|\nabla_x\mathbb{P}u^{\e}\|_{H^2_x}+\|\w^\mathsf{d}\|_{H^s_x}\|\nabla_x\mathbb{P}u^{\e}\|_{L^\infty_x}\big)\|\w^\mathsf{d}\|_{H^s_x} , \quad s=1,  \\
\big(\|\w^\mathsf{d}\|_{L^\infty_x}\|\nabla_x\mathbb{P}u^{\e}\|_{H^s_x}+\|\w^\mathsf{d}\|_{H^s_x}\|\nabla_x\mathbb{P}u^{\e}\|_{L^\infty_x}\big)\|\w^\mathsf{d}\|_{H^s_x} , \quad s=2, \end{cases} \cr 
&+ \sum_{0\leq|\al_x|\leq s}\|\p^{\al_x}\nabla_x\mathbb{P}^{\perp}u^{\e}\|_{L^\infty_x}\|\w^\mathsf{d}\|_{H^s_x}^2,
\enda
\end{align*}
where for $s=1$ we used Ladyzhenskaya’s inequality \eqref{Lady} $\|\w^\mathsf{d}\nabla_x^2\mathbb{P}u^{\e}\|_{L^2_x} \leq \|\w^\mathsf{d}\|_{L^4_x} \|\nabla_x^2\mathbb{P}u^{\e}\|_{L^4_x} \leq \|\w^\mathsf{d}\|_{H^1_x}\|\nabla_x\mathbb{P}u^{\e}\|_{H^2_x}$.
For $IV_{\w^\mathsf{d}}$, decompose $u^{\e}-u^E= u^\mathsf{d} + \mathbb{P}^{\perp}u^{\e}$ and apply \eqref{uvHk}:
\begin{align*}
\bega
|IV_{\w^\mathsf{d}}| &\leq \begin{cases} 
\big(\|\w^{\e}\|_{L^\infty_x} \|\nabla_xu^\mathsf{d}\|_{H^s_x}+\|\nabla_x\w^{\e}\|_{H^1_x} \|\nabla_xu^\mathsf{d}\|_{H^1_x}\big) \|\w^\mathsf{d}\|_{H^s_x}, \quad &s=1, \\
\big(\|\w^{\e}\|_{L^\infty_x} \|\nabla_xu^\mathsf{d}\|_{H^s_x}+\|\w^{\e}\|_{H^s_x} \|\nabla_xu^\mathsf{d}\|_{L^\infty_x}\big) \|\w^\mathsf{d}\|_{H^s_x}, \quad &s=2, 
\end{cases} \cr 
&+ \sum_{0\leq|\al_x|\leq s}\|\p^{\al_x}\nabla_x\mathbb{P}^{\perp}u^{\e}\|_{L^\infty_x}\|\w^{\e}\|_{H^2_x}\|\w^\mathsf{d}\|_{H^s_x},
\enda
\end{align*}
where, for $s=1$, we again used Ladyzhenskaya’s inequality \eqref{Lady}  $\|\nabla_x\w^{\e}\nabla_xu^\mathsf{d}\|_{L^2_x} \leq \|\nabla_x\w^{\e}\|_{L^4_x} \|\nabla_xu^\mathsf{d}\|_{L^4_x} \leq \|\nabla_x\w^{\e}\|_{H^1_x} \|\nabla_xu^\mathsf{d}\|_{H^1_x}$. 
Combining the bounds for $I_{\w^\mathsf{d}},\dots, IV_{\w^\mathsf{d}}$ in \eqref{wdHs} and dividing by $\|\w^\mathsf{d}\|_{H^{s}_x}$ gives \eqref{wH2C}.
\\(2) Next, we write the $H^s_x$ energy estimate for $\frac{3}{2}\ta^\mathsf{d}-\rho^\mathsf{d}$ from \eqref{barrtaeqn}:
\begin{align*}
\bega
\frac{d}{dt}\Big\|\frac{3}{2}\ta^\mathsf{d}-\rho^\mathsf{d}\Big\|_{H^s_x}^2 & \les  \sum_{0\leq|\al_x|\leq s} \int_{\Omega}\p^{\al_x}\Big(\mathbb{P}u^{\e} \cdot \nabla_x \Big(\frac{3}{2}\ta^\mathsf{d}-\rho^\mathsf{d}\Big)\Big) \p^{\al_x}\Big(\frac{3}{2}\ta^\mathsf{d}-\rho^\mathsf{d}\Big) dx \cr 
&+ \sum_{0\leq|\al_x|\leq s} \int_{\Omega} \p^{\al_x}\Big(u^\mathsf{d} \cdot \nabla_x \Big(\frac{3}{2}\ta^E-\rho^E\Big)\Big)\p^{\al_x}\Big(\frac{3}{2}\ta^\mathsf{d}-\rho^\mathsf{d}\Big) dx \cr 
& +\big\|\bar{\varPi}_{\mathfrak{s}}^{\e} \big\|_{H^s_x} \Big\|\frac{3}{2}\ta^\mathsf{d}-\rho^\mathsf{d}\Big\|_{H^s_x}.
\enda
\end{align*}
Using the same argument as in \eqref{uwE}, we obtain the estimate \eqref{rtaH3C}.
\end{proof}

\begin{lemma}\label{L.forcingC}
For the Boltzmann solution constructed in Theorem~\ref{T.2D.global} (for $d=2$ and $\mathrm{N}\geq 4$) or Theorem~\ref{T.3D.unif} (for $d=3$ and $\mathrm{N}\geq 5$) on $t\in[0,T]$ satisfying \eqref{BEC}, the following estimates hold:
\begin{align*}
\bega
\int_0^t \mathfrak{O}(s) ds \leq \mathcal{C}_* \Big(\e^{\frac{(d-1)-}{4}-} +\kappa^{\frac{1}{2}-}\Big),
\enda
\end{align*}
and
\begin{align*}
\bega
\int_0^t\Big(\|\bar{\varPi}_{\w}^{\e}(s)\|_{H^{s-1}_x} + \big\|\bar{\varPi}_{\mathfrak{s}}^{\e}(s) \big\|_{H^s_x}\Big)ds &\leq \mathcal{C}_* \Big(\e^{\frac{(d-1)-}{4}-} +\kappa^{\frac{1}{2}-}\Big),
\enda
\end{align*}
where $\mathfrak{O}(s)$ and $\bar{\varPi}_{\w}^{\e},\bar{\varPi}_{\mathfrak{s}}^{\e}$ are defined in \eqref{udforce} and \eqref{barvarPidef}, respectively, and the constant $\mathcal{C}_*$ is defined in~\eqref{C*def}.
\end{lemma}
\begin{proof}
(1) Considering the fact that $\mathbb{P}u^{\e} \in L^{2+ }$, we apply the H\"{o}lder inequality to the second line of \eqref{udforce} as follows: for any $0<\delta \ll1$, 
\begin{align*}
\bega
&\Big\|\Big(\mathbb{P}^{\perp}u^{\e}\cdot \nabla_x \mathbb{P}u^{\e} +\mathbb{P}u^{\e}\cdot \nabla_x \mathbb{P}^{\perp}u^{\e} + \mathbb{P}^{\perp}u^{\e}\cdot \nabla_x \mathbb{P}^{\perp}u^{\e}\Big)  \Big\|_{L^2_x} \cr 
&\leq \|\mathbb{P}^{\perp}u^{\e}\|_{L^\infty_x}\big(\|\nabla_x\mathbb{P}u^{\e}\|_{L^2_x}+ \|\nabla_x\mathbb{P}^{\perp}u^{\e}\|_{L^2_x}\big) + \|\mathbb{P}u^{\e}\|_{L^4_x} \|\nabla_x \mathbb{P}^{\perp}u^{\e}\|_{L^4_x}.
\enda
\end{align*}
Then, by the same argument as in the estimate of \eqref{GH2est-T} and \eqref{C*def}, we obtain
\begin{align*}
\bega
\int_0^t \mathfrak{O}(s) ds \leq \mathcal{C}_* \Big(\e^{\frac{(d-1)-}{4}-} +\kappa^{\frac{1}{2}-}\Big).
\enda
\end{align*}
(2) Using the same argument as in the proof of \eqref{GH2esti}, we obtain 
\begin{align*}
\bega
\|\bar{\varPi}_{\w}^{\e}(t)\|_{H^{s-1}_x} &\les \sum_{0\leq|\al_x|\leq s}\|\p^{\al_x}\mathbb{P}^{\perp}u^{\e}(t)\|_{L^\infty_x} + \sum_{0\leq|\al_x|\leq s}\|\p^{\al_x}(\rho^{\e}+\ta^{\e})(t)\|_{L^\infty_x} +  \kappa^{\frac{1}{2}}\mathcal{D}_G^{\frac{1}{2}}(t), \cr 
\big\|\bar{\varPi}_{\mathfrak{s}}^{\e}(t) \big\|_{H^s_x}&\les \sum_{0\leq|\al_x|\leq s}\|\p^{\al_x}\mathbb{P}^{\perp}u^{\e}(t)\|_{L^\infty_x} + \kappa^{\frac{1}{2}}\mathcal{D}_G^{\frac{1}{2}}(t).
\enda
\end{align*}
Then, applying~\eqref{EDbdd2D} and~\eqref{incomp.2D} from Theorem~\ref{T.2D.global}, and using~\eqref{C*def} together with~\eqref{BEC}, we obtain the desired estimates.
\end{proof}

\begin{proof}[{\bf Proof of \eqref{TC.wH2} in Proposition \ref{P.C.Hk}}]
We first note that, since $\nabla_x\cdot u^\mathsf{d}=0$, applying \eqref{puw} to $u^\mathsf{d} := \mathbb{P}u^{\e}-u^E$ yields
\begin{align}\label{dpmw}
\bega
\|\nabla_xu^\mathsf{d}\|_{H^s_x} \les \|\w^\mathsf{d}\|_{H^s_x}, \qquad \|u^\mathsf{d}\|_{H^s_x} \les \|u^\mathsf{d}\|_{L^2_x} + \|\w^\mathsf{d}\|_{H^{s-1}_x}, \quad \mbox{for} \quad s\geq1.
\enda
\end{align}
We now combine the estimate \eqref{u-uL2C} with \eqref{wH2C} for $s=1$ from Lemma~\ref{L.Pu-uconv} and Lemma~\ref{L.wrtaHsC}, respectively. Using \eqref{BEC} and \eqref{C}, we obtain
\begin{align*}
\bega
\frac{d}{dt}\big(\|u^{\mathsf{d}}\|_{L^2_x}+\|\w^\mathsf{d}\|_{H^1_x}\big) &\les  \Big(1+\sum_{0\leq|\al_x|\leq 1}\|\p^{\al_x}\nabla_x\mathbb{P}^{\perp}u^{\e}\|_{L^\infty_x}\Big)\big(\|u^{\mathsf{d}}\|_{L^2_x}+\|\w^\mathsf{d}\|_{H^1_x}\big) \cr 
&+ \sum_{0\leq|\al_x|\leq 1}\|\p^{\al_x}\nabla_x\mathbb{P}^{\perp}u^{\e}\|_{L^\infty_x} + \mathfrak{Q} + \|\bar{\varPi}_{\w}^{\e}\|_{H^1_x},
\enda
\end{align*}
where we have used $\|u^\mathsf{d}\|_{L^\infty_x} \les \|u^\mathsf{d}\|_{H^2_x}$ from Agmon's inequality \eqref{Agmon}, together with \eqref{dpmw}.
Applying Gr\"{o}nwall's inequality and Lemma~\ref{L.forcingC}, we get
\begin{align*}
\bega
\|u^{\mathsf{d}}(t)\|_{L^2_x} +\|\w^\mathsf{d}(t)\|_{H^1_x} &\les e^{Ct}\big(\|u^{\mathsf{d}}_0\|_{L^2_x}+\|\w^\mathsf{d}_0\|_{H^1_x}\big) \cr 
&+e^{Ct}\int_0^t \bigg(\sum_{0\leq|\al_x|\leq 1}\|\p^{\al_x}\nabla_x\mathbb{P}^{\perp}u^{\e}\|_{L^\infty_x} + \mathfrak{Q}(s)+\|\bar{\varPi}_{\w}^{\e}(s)\|_{H^1_x}\bigg) ds.
\enda
\end{align*}
For the terms $\mathfrak{Q}(s)$ and $\|\bar{\varPi}_{\w}^{\e}\|_{H^1_x}$, using Lemma~\ref{L.forcingC}, we conclude the proof of \eqref{TC.wH2}.
By exactly the same argument, we obtain the corresponding convergence for $\|(\tfrac{3}{2}\ta^\mathsf{d}-\rho^\mathsf{d})(t)\|_{H^2_x}^2$.
\end{proof}

\subsubsection{Convergence in $H^3$}

In this part, we prove the $H^3_x$ convergence \eqref{TC.wH3} by using the $H^2_x$ estimate \eqref{TC.wH2}. The proof is inspired by the inviscid limit of the Navier--Stokes equations presented in \cite{Masmoudi}. 

To establish the convergence, we introduce a mollified solution sequence and apply the triangle inequality.

\begin{definition}\label{D.macrod2}
Let \((\rho^{E,\lambda}(t), u^{E,\lambda}(t), \ta^{E,\lambda}(t))\) denote the solution of the incompressible Euler equations \eqref{incompE} corresponding to the mollified initial data
\begin{align}\label{ud0}
\bega
(\rho_0^{E,\lambda}(x),u_0^{E,\lambda}(x),\ta_0^{E,\lambda}(x)) := \mathcal{F}^{-1}(1_{|\xi|\leq 1/\lambda}\mathcal{F}(\rho^E_0,u^E_0,\ta^E_0)).
\enda
\end{align}
We then split the difference into two parts:
\begin{align*}
\bega
(\rho^{\mathsf{d}_1}, u^{\mathsf{d}_1}, \ta^{\mathsf{d}_1},\w^{\mathsf{d}_1})(t,x) &:= (\rho^{\e}, \mathbb{P}u^{\e}, \ta^{\e},\w^{\e})(t,x) - (\rho^{E,\lambda}, u^{E,\lambda}, \ta^{E,\lambda},\w^{E,\lambda})(t,x), \cr 
(\rho^{\mathsf{d}_2}, u^{\mathsf{d}_2}, \ta^{\mathsf{d}_2},\w^{\mathsf{d}_2})(t,x) &:= (\rho^{E,\lambda}, u^{E,\lambda}, \ta^{E,\lambda},\w^{E,\lambda})(t,x) - (\rho^E, u^E, \ta^E,\w^E)(t,x),
\enda
\end{align*}
where $\w^{E,\lambda}:=\nabla_x\times u^{E,\lambda}$.
\end{definition}

We next present some properties of \((\rho^{\mathsf{d}_2}, u^{\mathsf{d}_2}, \ta^{\mathsf{d}_2}, \w^{\mathsf{d}_2})\), which will be used in the proof.

\begin{lemma}\label{L14.incomp} 
Assume \eqref{C}. Then, the following estimates hold.
\begin{enumerate}
\item The solution of the incompressible Euler equations with mollified initial data satisfies
\begin{align}\label{H4,1/lam}
\bega
\sup_{0\leq t \leq T}\|(\rho^{E,\lambda},u^{E,\lambda},\ta^{E,\lambda})(t)\|_{H^4_x} \les \frac{1}{\lambda},
\enda
\end{align}
where \(\lambda\) is the mollification parameter.

\item The difference between the mollified and original Euler solutions satisfies
\begin{align}\label{tilu}
\bega
\sup_{0\leq t \leq T}\|u^{\mathsf{d}_2}(t)\|_{H^{s'}_x} \leq C\lambda^{3-s'}, \quad \mbox{for} \quad 0\leq s' \leq 3.
\enda
\end{align}

\item The difference is bounded in \(H^3_x\):
\begin{align}\label{tilrut}
\bega
&\|\w^{\mathsf{d}_2}(t)\|_{H^2_x} \les  e^{CT}\big(\|\w^{\mathsf{d}_2}(0)\|_{H^2_x} +\lambda^{2-s'} \big),
\enda
\end{align}
for $d/2 < s' < 2$.
\end{enumerate}
\end{lemma}

\begin{proof}
(i) By the definition \eqref{ud0}, we have
\begin{align}\label{E-lam-0bdd}
\bega
&\|(\rho^{E,\lambda}_0,u^{E,\lambda}_0,\ta^{E,\lambda}_0)\|_{H^{s'}_x} \leq \frac{C}{\lambda^{3-s'}}, \quad \mbox{for} \quad s'=3,4, \cr
&\|(\rho^E_0-\rho^{E,\lambda}_0,u^E_0-u^{E,\lambda}_0,\ta^E_0-\ta^{E,\lambda}_0)\|_{H^{s'}_x}\leq C\lambda^{3-s'}, \quad \mbox{for} \quad s'\leq 3.
\enda
\end{align}
Applying the commutator estimate \eqref{commutator} to \eqref{incompE}, we derive
\begin{align*}
\bega
\frac{d}{dt}\|(\rho^{E,\lambda},u^{E,\lambda},\ta^{E,\lambda})\|_{H^{\tilde{s}}_x}^2 &\leq C \|\nabla_x(\rho^{E,\lambda},u^{E,\lambda},\ta^{E,\lambda})\|_{L^\infty_x}\|(\rho^{E,\lambda},u^{E,\lambda},\ta^{E,\lambda})\|_{H^{\tilde{s}}_x}^2,
\enda
\end{align*}
for $\tilde{s}=3,4$. Using \(\|\nabla_x(\rho^{E,\lambda},u^{E,\lambda},\ta^{E,\lambda})\|_{L^\infty_x} \leq C\) from \eqref{C} and applying Gr\"{o}nwall’s inequality, we obtain
\begin{align*}
\bega
\sup_{0 \leq t \leq T}\|(\rho^{E,\lambda},u^{E,\lambda},\ta^{E,\lambda})(t)\|_{H^{\tilde{s}}_x} \leq \frac{C}{\lambda^{\tilde{s}-3}}, \quad \mbox{for} \quad \tilde{s} = 3, 4.
\enda
\end{align*}
(ii) Subtracting the equations \eqref{incompE} for $(\rho^E,u^E,\ta^E)$ from those for $(\rho^{E,\lambda},u^{E,\lambda},\ta^{E,\lambda})$, we get
\begin{align*}
\bega
&\p_t \rho^{\mathsf{d}_2} +u^{\mathsf{d}_2}\cdot\nabla_x \rho^{E,\lambda} + u^{E}\cdot\nabla_x \rho^{\mathsf{d}_2} = 0, \cr
&\p_t u^{\mathsf{d}_2} +u^{\mathsf{d}_2}\cdot\nabla_x u^{E,\lambda} + u^{E}\cdot\nabla_x u^{\mathsf{d}_2} + \nabla_x(p^{E,\lambda}-p^E)= 0, \cr 
&\p_t \ta^{\mathsf{d}_2} +u^{\mathsf{d}_2}\cdot\nabla_x \ta^{E,\lambda} + u^{E}\cdot\nabla_x \ta^{\mathsf{d}_2} = 0.
\enda
\end{align*}
From the energy estimate we get
\begin{align*}
\bega
\frac{d}{dt} \|u^{\mathsf{d}_2}\|_{H^{s'}_x}^2 &\leq (\|u^{E}\|_{H^3_x}+\|u^{E,\lambda}\|_{H^3_x})\|u^{\mathsf{d}_2}\|_{H^{s'}_x}^2, \quad \mbox{for} \quad d/2 < s' \leq 2.
\enda
\end{align*}
Since \(\|u^{E}\|_{H^3_x} \leq C\), \(\|u^{E,\lambda}\|_{H^3_x} \leq C\), and \(\eqref{E-lam-0bdd}_2\), inequality \eqref{tilu} follows.
 \\
(iii) From the $H^3_x$ boundedness in \eqref{C} and \eqref{E-lam-0bdd}, and by the same method as in Lemma \ref{L.wrtaHsC}, we obtain the energy inequality
\begin{align*}
\bega
\frac{d}{dt}\|\w^{\mathsf{d}_2}\|_{H^2_x} &\les  \|\w^{\mathsf{d}_2}\|_{H^2_x} + \|u^{\mathsf{d}_2}\|_{L^2_x} +\|u^{\mathsf{d}_2}\|_{L^\infty_x}\|\w^{E,\lambda}\|_{H^3_x}.
\enda
\end{align*}
With Sobolev embedding \eqref{embed}, $\|u^{\mathsf{d}_2}\|_{L^\infty_x}\leq\|u^{\mathsf{d}_2}\|_{H^{s'}_x}$ for $d/2<s'<2$. Combining this with \eqref{H4,1/lam} and \eqref{tilu}, we deduce
\begin{align*}
\bega
\frac{d}{dt}\|\w^{\mathsf{d}_2}\|_{H^2_x} &\les  \|\w^{\mathsf{d}_2}\|_{H^2_x} + \lambda^{3-s'-1}.
\enda
\end{align*}
Applying Gr\"{o}nwall’s inequality yields the result \eqref{tilrut}.
\end{proof}

\begin{proof}[{\bf Proof of \eqref{TC.wH3} in Proposition \ref{P.C.Hk}}] 
Combining \eqref{u-uL2C} in Lemma~\ref{L.Pu-uconv} and \eqref{wH2C} in Lemma~\ref{L.wrtaHsC} for \(s=2\), and invoking \eqref{BEC}--\eqref{C}, we obtain
\begin{align*}
\bega
\frac{d}{dt}\big(\|u^{\mathsf{d}_1}\|_{L^2_x}+\|\w^{\mathsf{d}_1}\|_{H^2_x}\big) &\les \Big(1+\sum_{0\leq|\al_x|\leq 2}\|\p^{\al_x}\nabla_x\mathbb{P}^{\perp}u^{\e}\|_{L^\infty_x}\Big) \big(\|u^{\mathsf{d}_1}\|_{L^2_x}+\|\w^{\mathsf{d}_1}\|_{H^2_x}\big)  \cr 
& +\|u^{\mathsf{d}_1}\|_{L^\infty_x}\|\w^{E,\lambda}\|_{H^3_x} + \sum_{0\leq|\al_x|\leq 2}\|\p^{\al_x}\nabla_x\mathbb{P}^{\perp}u^{\e}\|_{L^\infty_x} + \mathfrak{Q} + \|\bar{\varPi}_{\w}^{\e}\|_{H^2_x}.
\enda
\end{align*}
Since \(\|\nabla_x \w^{E,\lambda}\|_{H^2_x} \les \|u^{E,\lambda}\|_{H^4_x}\) by \eqref{puw}, and 
\(\|u^{E,\lambda}\|_{H^4_x} \les 1/\lambda\) by \eqref{H4,1/lam}, applying Gr\"{o}nwall’s inequality, we deduce
\begin{align}\label{wrtH3C-a}
\bega
\|u^{\mathsf{d}_1}(t)\|_{L^2_x}&+\|\w^{\mathsf{d}_1}(t)\|_{H^2_x} \les e^{Ct}\big(\|u^{\mathsf{d}_1}(0)\|_{L^2_x}+\|\w^{\mathsf{d}_1}(0)\|_{H^2_x} \big) +e^{Ct}\frac{\int_0^t\|u^{\mathsf{d}_1}(s)\|_{L^\infty_x}ds}{\lambda} \cr 
& + e^{Ct}\int_0^t \bigg( \sum_{0\leq|\al_x|\leq 2}\|\p^{\al_x}\nabla_x\mathbb{P}^{\perp}u^{\e}(s)\|_{L^\infty_x} + \mathfrak{Q}(s) + \|\bar{\varPi}_{\w}^{\e}(s)\|_{H^2_x} \bigg) ds .
\enda
\end{align}
Combining \eqref{wrtH3C-a} with the bounds for $u^{\mathsf{d}_2}$ and $\w^{\mathsf{d}_2}$ in \eqref{tilu} and \eqref{tilrut}, respectively, we find
\begin{align}\label{rutaH3C}
\bega
\|u^\mathsf{d}(t)\|_{L^2_x} &+\|\w^\mathsf{d}(t)\|_{H^2_x} \les e^{Ct}\big(\|u^{\mathsf{d}_1}(0)\|_{L^2_x}+\|\w^{\mathsf{d}_1}(0)\|_{H^2_x} \big) \cr 
&+e^{Ct}\frac{\int_0^t\|u^{\mathsf{d}_1}(s)\|_{L^\infty_x}ds}{\lambda} +e^{Ct}\bigg(\|\w^{\mathsf{d}_2}(0)\|_{H^2_x} + \lambda^{2-s'} + C\big(\e^{\frac{d-1}{4}-}\kappa^{-\frac{1}{2}}+\kappa^{\frac{1}{2}-}\big)\bigg) ,
\enda
\end{align}
for $d/2<s'<2$, where Lemma~\ref{L.forcingC} was applied to handle $\|\bar{\varPi}_{\w}^{\e}\|_{H^2_x}$.
It remains to estimate $\int_0^t\|u^{\mathsf{d}_1}(s)\|_{L^\infty_x}\,ds$. Decomposing $u^{\mathsf{d}_1} = \mathbb{P}u^{\e}-u^{E,\lambda}=(\mathbb{P}u^{\e}-u^E)+(u^E-u^{E,\lambda})$ and applying the Sobolev embedding \eqref{embed}, we obtain, for any $s'$ with $d/2 < s' < 2$,
\begin{align*}
\bega
\int_0^t\|u^{\mathsf{d}_1}(s)\|_{L^\infty_x}ds &\leq \int_0^t\|(\mathbb{P}u^{\e}-u^E)(s)\|_{H^2_x}ds + \int_0^t\|(u^E-u^{E,\lambda})(s)\|_{H^{s'}_x}ds \cr 
&\les \int_0^t\big(\|u^{\mathsf{d}}(s)\|_{L^2_x}+\|\w^{\mathsf{d}}(s)\|_{H^1_x}\big)ds+ \int_0^t\|(u^E-u^{E,\lambda})(s)\|_{H^{s'}_x}ds ,
\enda
\end{align*}
where we used $\|\nabla_x(\mathbb{P}u^{\e}-u^E)\|_{H^1_x}\leq \|\w^{\mathsf{d}}\|_{H^1_x}$ from \eqref{dpmw}. 
For the term $\|u^{\mathsf{d}}\|_{L^2_x}+\|\w^{\mathsf{d}}\|_{H^1_x}$, we apply the convergence result \eqref{TC.wH2}, while $\|(u^E-u^{E,\lambda})(s)\|_{H^{s'}_x}$ is controlled by \eqref{tilu}.
Altogether, this yields
\begin{align}\label{baruconv}
\bega
\int_0^t\|u^{\mathsf{d}_1}(s)\|_{L^\infty_x}ds &\les_T \big(\|u^{\mathsf{d}}_0\|_{L^2_x}+\|\w^\mathsf{d}_0\|_{H^1_x}\big) +\e^{\frac{d-1}{4}-} + (\kappa T)^{\frac{1}{2}} +\lambda^{3-s'}t.
\enda
\end{align}
Substituting \eqref{baruconv} into \eqref{rutaH3C} yields
\begin{align*}
\bega
\sup_{0 \leq t \leq T}\big(&\|u^\mathsf{d}(t)\|_{L^2_x} +\|\w^\mathsf{d}(t)\|_{H^2_x}\big) \les C_T \bigg(\big(\|u^{\mathsf{d}_1}(0)\|_{L^2_x}+\|\w^{\mathsf{d}_1}(0)\|_{H^2_x} \big) +\|\w^{\mathsf{d}_2}(0)\|_{H^2_x} \cr 
& + \e^{\frac{d-1}{4}-}\kappa^{-\frac{1}{2}}+\kappa^{\frac{1}{2}-} +\frac{1}{\lambda}\Big(\big(\|u^{\mathsf{d}}_0\|_{L^2_x}+\|\w^\mathsf{d}_0\|_{H^1_x}\big) +\e^{\frac{d-1}{4}-} + (\kappa T)^{\frac{1}{2}} +\lambda^{3-s'}T\Big)\bigg).
\enda
\end{align*}
By choosing proper $\lambda$, we obtain the desired estimate \eqref{TC.wH3}.
By the same argument, using the estimates \eqref{rtaH3C} and the bound on $\int_0^t\big\|\bar{\varPi}_{\mathfrak{s}}^{\e}(s) \big\|_{H^3_x}\,ds$ in Lemma~\ref{L.forcingC}, we obtain the corresponding convergence for $\|(\tfrac{3}{2}\ta^\mathsf{d}-\rho^\mathsf{d})(t)\|_{H^3_x}^2$.

\end{proof}

\hide
\begin{proof}[{\bf (Proof of Remark \ref{R.M-M})}]
(Proof of \eqref{M-M,CH2}) Using $\|\nabla_x(\rho+\ta)\|_{L^\infty_tH^s_x}\les \kappa^{\delta-\frac{s+1}{2}}$, we can have 
\begin{align}\label{divbou+2}
\bega
\|\nabla_x \cdot \bar{m}\|_{L^\infty(0,T;H^s_x)} =\|\nabla_x \cdot u\|_{L^\infty(0,T;H^s_x)} &\les \kappa^{\delta-\frac{s}{2}}, \qquad \|(\rho^\mathsf{d}+\ta^\mathsf{d})\|_{L^\infty(0,T;H^{s+1}_x)} \les \kappa^{\delta-\frac{s+1}{2}},
\enda
\end{align}
Then we can use \eqref{divbou+2} to the inequality \eqref{wrtH2C-c} as follows: 
\begin{align}\label{CH2-a}
\bega
\|&(\rho^\mathsf{d},u^\mathsf{d},\ta^\mathsf{d})(t)\|_{H^2_x}^2 \les e^{Ct}\bigg(\|\w^\mathsf{d}(0)\|_{H^1_x}^2 +\Big\|\Big(\rho^\mathsf{d}-\frac{3}{2}\ta^\mathsf{d}\Big)(0)\Big\|_{H^2_x}^2 \bigg) \cr 
&+e^{Ct}(\kappa^{2\delta-1}+\kappa) +\kappa^{2\delta-2} + \int_0^t e^{C(t-s)} BE~part(s)ds.
\enda
\end{align}
Applying Lemma \ref{MFMEdif}, we have \eqref{M-M,CH2}. \\ 
(Proof of \eqref{M-M,CH3}) In the inequality \eqref{rutaH3C}, we can change the quantity $\|u^{\mathsf{d}_1}\|_{L^2_tL^\infty_x}$ to $\|u^{\mathsf{d}_1}\|_{L^\infty_tL^\infty_x}$ as follows:
\begin{align}\label{CH3-a}
\bega
&\|(\rho,u,\ta)(t)-(\rho^E,u^E,\ta^E)(t)\|_{H^3_x}^2 \les e^{Ct}\bigg(\|\w^{\mathsf{d}_1}(0)\|_{H^2_x}^2 +\Big\|\Big(\rho^{\mathsf{d}_1}-\frac{3}{2}\ta^{\mathsf{d}_1}\Big)(0)\Big\|_{H^3_x}^2 \bigg) \cr 
&+e^{Ct}\frac{\|u^{\mathsf{d}_1}\|_{L^\infty_tL^\infty_x}^2}{\lambda^2} +e^{Ct}\Big(\frac{\kappa}{\lambda^2} +\kappa^{2\delta-2}+\e^2\kappa^{-1}\Big) + \int_0^t e^{C(t-s)} BE~part(s)ds \cr 
&+\|\nabla_x \cdot u(t)\|_{H^2_x}^2 +\|(\rho^{\mathsf{d}_1}+\ta^{\mathsf{d}_1})(t)\|_{H^3_x}^2 + e^{Ct}(\|(\rho^{\mathsf{d}_2}_0,u^{\mathsf{d}_2}_0,\ta^{\mathsf{d}_2}_0)\|_{H^3_x}^2+ \lambda^{4-2s'}),
\enda
\end{align}
for $d/2<s'<2$. In addition, we can modify the inequality \eqref{baruconv} to $L^\infty_t$ as follows:
\begin{align}\label{CH3baru}
\bega
\|u^{\mathsf{d}_1}\|_{L^\infty_tL^\infty_x}^2 &\leq \|u-u^E\|_{L^\infty_tH^2_x}^2 + \|u^E-u^{E,\lambda}\|_{L^\infty_tH^{s'}_x}^2 \cr 
&\les \bigg[e^{Ct}\bigg(\|\w^\mathsf{d}(0)\|_{H^1_x}^2 +\Big\|\Big(\rho^\mathsf{d}-\frac{3}{2}\ta^\mathsf{d}\Big)(0)\Big\|_{H^2_x}^2 \bigg) \cr 
&+e^{Ct}(\kappa^{2\delta-1}+\kappa) +\kappa^{2\delta-2} + \int_0^t e^{C(t-s)} BE~part(s)ds \bigg] + \lambda^{6-2s'}t,
\enda
\end{align}
where we applied \eqref{CH2-a} and \eqref{tilu} for $\|u-u^E\|_{L^\infty_tH^2_x}^2$ and $\|u^E-u^{E,\lambda}\|_{L^\infty_tH^{s'}_x}^2$, respectively. We apply \eqref{divbou+2} and \eqref{CH3baru} to \eqref{CH3-a}. 
\begin{align}\label{CH3-b}
\bega
&\|(\rho,u,\ta)(t)-(\rho^E,u^E,\ta^E)(t)\|_{H^3_x}^2 \les e^{Ct}\bigg(\|\w^{\mathsf{d}_1}(0)\|_{H^2_x}^2 +\Big\|\Big(\rho^{\mathsf{d}_1}-\frac{3}{2}\ta^{\mathsf{d}_1}\Big)(0)\Big\|_{H^3_x}^2 \bigg) \cr 
&+e^{Ct}\frac{1}{\lambda^2}\bigg[\bigg[e^{Ct}\bigg(\|\w^\mathsf{d}(0)\|_{H^1_x}^2 +\Big\|\Big(\rho^\mathsf{d}-\frac{3}{2}\ta^\mathsf{d}\Big)(0)\Big\|_{H^2_x}^2 \bigg) \cr 
&+e^{Ct}(\kappa^{2\delta-1}+\kappa) +\kappa^{2\delta-2} + \int_0^t e^{C(t-s)} BE~part(s)ds \bigg] + \lambda^{6-2s'}t\bigg] \cr 
&+\kappa^{2\delta-3} + e^{Ct}(\|(\rho^{\mathsf{d}_2}_0,u^{\mathsf{d}_2}_0,\ta^{\mathsf{d}_2}_0)\|_{H^3_x}^2+ \lambda^{4-2s'}),
\enda
\end{align}
If we choose $\delta>3/2$ and $\lambda:= \kappa^a$ for $a<1/2$, then when $\kappa\rightarrow0$, the right-side goes to $0$. Using Lemma \ref{MFMEdif}, we have \eqref{M-M,CH3}.
\end{proof}
\unhide

\subsection{$L^\infty$-vorticity (Yudovich solution) and Besov-regular vorticity}

In the case of bounded vorticity (the Yudovich class \cite{Yudovich}), the inviscid limit of the incompressible Navier--Stokes equations has been extensively studied. In two dimensions, Chemin \cite{Chemin} established convergence of the velocity by exploiting an Osgood-type argument. Subsequently, Masmoudi \cite{Masmoudi} obtained convergence results for the vorticity in both $d=2$ and $d=3$.
More recently, convergence of the vorticity in $L^p$ for all $p \in [1,\infty)$ on $\T^2$ was established in \cite{2020Constan}, where convergence rates were also obtained under additional regularity assumptions.
For further developments, see \cite{2020Constan} and the references therein.

In this section, we prove the convergence of the vorticity $\w^{\e}$ when the initial vorticity belongs to \(L^\infty \cap L^1(\R^2)\). The limit becomes the unique weak solution of the incompressible Euler equation. 
Under a Besov-type regularity condition, we also obtain a convergence rate. 

Recall the radial energy decomposion in Definition \ref{D.Ra-E} and $\bar{u}\cdot\nabla_x \bar{\w}=0$. 
Recall also that we decomposed $F^{\e}=M^{\e}+\AC{\P}F^{\e}$ as in \eqref{F-ME}.
We now present the main theorem of this section.

\begin{theorem}\label{T.C.Linf}
Let $\Omega=\R^2$.
Consider a family of initial data $\{F_0^{\e}\}_{\e>0}$ satisfies the $4$-ABC~\eqref{ABC1}, a uniform-in-$\e$ modulated entropy bound \eqref{L2unif}, and 
\begin{align*}
\sup_{\e>0} 
\|\w_0^{\e}\|_{L^\infty\cap L^1(\R^2)}
<+\infty.
\end{align*}\hide 
\begin{equation}\label{L2unif}
\sup_{\e>0}
\bigg(  \|( \rho_0^{\e}, u_0^{\e}-\bar{u}, \ta_0^{\e}
)\|_{L^2(\R^2)}
+ \left\|\frac{1}{\e}\frac{\AC{\P}F^{\e}_0}{\sqrt{\tilde{\mu}}}\right\|_{L^2_x(\R^2;L_v^2(\R^3))}
\bigg)
<+\infty.
\end{equation}
Here $\tilde{\mu}=M_{[1,0,1-c_0]}$ for some $0<c_0\ll1$, and $\bar{u}$
denotes the radial eddy defined in Definition~\ref{D.Ra-E}. \unhide 


We assume that, as $\e \to 0$,
\begin{align} 
\w^{\e}_0 \quad &\to \quad \w_0  \qquad \text{strongly in $L^2(\R^2)$} , \label{conv_initial_Yudo}  \\
\frac{3}{2} \theta^\e_0 - \rho^\e_0  \quad &\to \quad \frac{3}{2} \theta_0 - \rho_0  \qquad \text{strongly in $L^2(\R^2)$}. \label{conv_s_initial} 
\end{align}

Then the following statements hold for a family of Boltzmann solutions $\{F^{\e}\}_{\e>0}$ to~\eqref{BE} with $\kappa=\e^q$ for some $0<q<2$, on
the time interval $t\in[0,T_\e]$ with $T_\e\to\infty$ as $\e\to0$, which is constructed in Theorem~\ref{T.2D.global}.

\begin{enumerate}

\item[(1)] The family of macroscopic fields of the Boltzmann solution $\big\{(\mathbb{P} u^\e, \mathbb{P}^{\perp}u^\e, \frac{3}{2} \theta^\e - \rho^\e, \rho^\e+ \theta^\e) \big\}_{\e>0}$ converges to $(u^E, 0, \frac{3}{2} \theta^E - \rho^E, 0)$ in the following different topologies, respectively, and the limit vorticity $\w^E:=\nabla_x^{\perp}\cdot u^E$ is the unique weak solution to the vorticity--stream formulation of the Euler equation~\eqref{incompE} with initial data $\w^E|_{t=0}=\w_0$ and 
\begin{align}\label{wEYudospace}
\w^{E}  \in L^\infty(0,T;L^\infty\cap L^1 (\R^2)). 
\end{align}

\begin{itemize}
\item[(i)] The vorticities satisfy, for any $2\leq p<\infty$,
\begin{align}
&\hspace{3cm} \w^{\e} \quad \rightarrow \quad \w^{E} \quad \mbox{in} \quad  L^\infty(0,T;L^p(\R^2)), \label{wstrongYudo} 
\\ 
&\| (\mathbb{P}u^{\e}-u^E)(t) \|_{L^2 (\R^2)}
\lesssim e^{2-2e^{-Ct}} \bigg(\|\mathbb{P}u^{\e}_0-u^{E}_0\|_{L^2_x}+\mathcal{C}_* \Big(\e^{\frac{1}{4}-} +\kappa^{\frac{1}{2}-}\Big) \bigg)^{e^{-Ct}}.
\notag 
\end{align}

\item[(ii)] The entropic fluctuation $\mathfrak{s}^{\e}:=3/2\ta^{\e}-\rho^{\e}$ converges as follows:
\begin{align*}
\frac{3}{2} \theta^\e - \rho^\e  \to \frac{3}{2} \theta^E - \rho^E  \ \ &\text{in  $L^\infty (0, T; L^2(\R^2))$ }. 
\end{align*}

\item[(iii)] The irrotational part and pressure fluctuation vanish as follows:
\begin{align}
\mathbb{P}^{\perp}u^{\e}, \ \rho^{\e}+\ta^{\e}    \to  0  \ \ &\text{in $ L^r (0,T; \dot{B}_{p,1}^{s+2(\frac{1}{p}-\frac{1}{2})+\frac{1}{r}}(\R^2))$ with rate $\e^{\frac{1}{4}-}$} ,\label{TC.Yudo.irrorhota}
\end{align}
for $2\leq p\leq\infty$, $\frac{1}{r}\leq 2(\frac{1}{2}-\frac{1}{p})$, and any $s\in[0,3)$. 
\end{itemize}

\item[(2)] The microscopic part converges to zero:
\begin{equation}\label{TC.Yudo.G}
\frac{1}{\e^2} \sqrt{\nu^{\e}}\frac{\AC{\P}F^{\e} }{\sqrt{M^\e}} =  O(\kappa^{\frac{1}{2}})   \ \ \text{in $L^2 (0,T; L^2 (\R^2_x \times \R^3_v))$}.  
\end{equation}
 
\end{enumerate}
\end{theorem}

\hide
\begin{remark}
The Boltzmann solution obtained from Theorem~\ref{T.2D.global}, satisfying the 
$\mathrm{N}$-Admissible Blow-up Condition~\eqref{ABC1}, obeys the following Strichartz estimate:
\begin{align}
\mathbb{P}^{\perp}u^{\e}, \ \rho^{\e}+\ta^{\e}    \to  0  \ \ &\text{in $ L^r (0,T; \dot{B}_{p,1}^{s+2(\frac{1}{p}-\frac{1}{2})+\frac{1}{r}}(\R^2))$ with rate $\e^{\frac{1}{4}-}$} ,\label{TC.Yudo.irrorhota}
\end{align}
for $2\leq p\leq\infty$, $\frac{1}{r}\leq 2(\frac{1}{2}-\frac{1}{p})$, and any $s\in[0,\mathrm{N}-1)$.  
Here $\mathrm{N}$ can be chosen arbitrarily within the range $\mathrm{N}>1+\frac{d}{2}$.  
Such an estimate is possible because the initial data satisfy the Admissible Blow-up Condition~\eqref{ABC1}; 
a small power of $\e$ absorbs the singular scale in \eqref{ABC1} and yields a factor of order $\e^{\frac{1}{4}-}$.
\end{remark}
\unhide

\hide

\begin{remark}
We summarize the quantitative convergence rates in the Big-O sense as follows:
\begin{align*}
e^{2-2e^{-Ct}} \bigg(\|\mathbb{P}u^{\e}_0-u^{E}_0\|_{L^2_x}+\mathcal{C}_* \Big(\e^{\frac{1}{4}-} +\kappa^{\frac{1}{2}-}\Big) \bigg)^{e^{-Ct}} \ \ &\text{in \eqref{Pu-uYudo}}, \\
\bigg( \bigg(\|\mathbb{P}u^{\e}_0-u^{E}_0\|_{L^2_x}+\mathcal{C}_* \Big(\e^{\frac{1}{4}-} +\kappa^{\frac{1}{2}-}\Big) \bigg)^{e^{-Ct}}\bigg)^{\frac{2}{p}\frac{se^{-CT}}{1+se^{-CT}}-}  
\ \ &\text{in \eqref{wrate}}.
\end{align*}
Here, \(\mathcal{C}_*\) is defined in \eqref{C*def}. 
By additionally assuming that the first spatial derivative of the microscopic part is bounded, namely $\left\|\frac{1}{\e}\frac{\nabla_x\AC{\P}F^{\e}_0}{\sqrt{\tilde{\mu}}}\right\|_{L^2_x(\R^2;L_v^2(\R^3))}\leq \infty$, we can improve the scaling $\mathcal{C}_*\kappa^{\frac{1}{2}-}$ to $C_T\kappa^{\frac{1}{2}}$.
See Remark~\ref{Rmk.Grefine}.
\end{remark}

\unhide

\begin{corollary}
Consider the following initial data satisfying
\begin{align*}
&\w_0\in L^\infty\cap L^1(\R^2), \quad  (u_0-\bar{u}) \in L^2(\R^2), \quad \rho_0, \ta_0 \in L^2(\R^2), \quad \frac{1}{\e}\frac{\AC{\P}F_0}{\sqrt{\tilde{\mu}}}\in L^2_x(\R^2 ; L^2_v(\R^3)),  \cr 
&\bigg\|\frac{1}{\e}\frac{(F_0-\mu)}{\sqrt{\tilde{\mu}}}\bigg\|_{L^\infty_{x,v}}\les \frac{1}{\e^{1-}},
\end{align*}
where $\tilde{\mu}=M_{[1,0,1-c_0]}$ for some $0<c_0\ll1$.
The mollified sequence $\{F^{\e}_0\}_{\e>0}$ constructed from Lemma \ref{L.molli} satisfies the assumption of Theorem \ref{T.C.Linf}, and hence all the results of Theorem \ref{T.C.Linf} hold. 
\end{corollary}

\begin{definition}\label{D.Fsdef}
We define the norm of the Triebel--Lizorkin-type space for $0<s<1$ as follows:
\begin{align*}
\bega
[\w]_{\mathcal{F}_2^{s}(\R^2)} := \|\w\|_{L^2(\R^2)} + \inf_{\mathfrak{V}\in L^2(\R^2)} \left\{\|\mathfrak{V}\|_{L^2(\R^2)}  :  |\w(x)-\w(y)| \leq |x-y|^s(\mathfrak{V}(x)+\mathfrak{V}(y)) \right\}.
\enda
\end{align*}
\end{definition}
\begin{theorem}\label{T.C.TL}Assume all the assumptions in Theorem \ref{T.C.Linf}. If we further assume that $\sup_{\e>0}[\w^{\e}_0]_{\mathcal{F}_2^{s}(\R^2)} <+\infty$ for some $0<s<1$, where the Triebel--Lizorkin-type space $\mathcal{F}^s_2$ is defined in Definition~\ref{D.Fsdef}, then for any $2\leq p<\infty$, we obtain the following convergence rate:
\begin{align}\label{wrate}
\bega
\| (\w^{\e}-   \w^{E} ) (t)\|_{L^p(\R^2)}
 \lesssim \bigg(\|\mathbb{P}u^{\e}_0-u^{E}_0\|_{L^2_x}+\mathcal{C}_* \Big(\e^{\frac{1}{4}-} +\kappa^{\frac{1}{2}-}\Big)  \bigg)^{\frac{2}{p} e^{-Ct}\frac{se^{-CT}}{1+se^{-CT}}-}  .
\enda
\end{align}
Here, \(\mathcal{C}_*\) is defined in \eqref{C*def}. 
\end{theorem}

\hide
 {\color{red}[will be deleted]
이건 너무 마이너한 어려움이라 써놓는 것이 오히려 손해입니다. 
$\bullet$ Mechanism is different from NS. See Remark \ref{Rmk.NS,BE}. So that we should divide each side by $\|u^\mathsf{d}\|_{L^2}$  이건 중요한 포인트인데, 이런식으로 묘사를 하게되면 굉장히 기초적으로 들리고, 또한 본질을 이야기하는 것도 아니죠. 제가 인트로에 잘 써놓을게요. }
\unhide

In Theorem~\ref{T.2D.global}, the solution satisfying the Admissible Blow-up Condition~\eqref{ABC1} also satisfies \eqref{incomp.2D} and \eqref{Gto0}. These correspond to \eqref{TC.Yudo.G} and \eqref{TC.Yudo.irrorhota}, respectively.
If we establish the convergence of the vorticity, then the thermodynamic variable
$\mathfrak{s}^{\e} := \frac{3}{2}\ta^{\e} - \rho^{\e}$,
which satisfies the continuity equation~\eqref{incompE}, also converges.
Hence, it suffices to focus on the proof of the convergence of the vorticity.

To be a weak solution, the vorticity must belong to \(L^\infty\cap L^1(\R^2)\). 
However, since the total vorticity $\w^{\e}$ contains a microscopic component, it is not clear whether it lies in \(L^1(\R^2)\). 
Therefore, as defined in Definition~\ref{D.wAwB}, we decompose the vorticity into the initial data part $\w^{\e}_A$ and the forcing part $\w^{\e}_B$. We then prove the convergence of each part separately in order to obtain the result \eqref{wstrongYudo}.

Recall that in Definition~\ref{D.wAwB}, we defined $\w^{\e}_A$ and $\w^{\e}_B$ as the solutions of the following equations:
\begin{align*}
&\p_t\w^{\e}_A + \mathbb{P}u^{\e}\cdot \nabla_x \w^{\e}_A = 0 , \qquad \hspace{3mm} \w^{\e}_A(t,x)|_{t=0} = \w^{\e}_0(x), \cr 
&\p_t\w^{\e}_B + \mathbb{P}u^{\e}\cdot \nabla_x \w^{\e}_B = \bar{\varPi}_{\w}^{\e} , \qquad \w^{\e}_B(t,x)|_{t=0} = 0 . 
\end{align*}
We also defined the associated divergence-free velocity fields via the Biot--Savart law:
\begin{align*}
\bega
u^{\e}_A(t,x) := {\bf K} \ast \w^{\e}_A(t,x), \qquad 
u^{\e}_B(t,x) := {\bf K} \ast \w^{\e}_B(t,x), \qquad {\bf K}(x):= \frac{1}{2\pi}\frac{x^{\perp}}{|x|^2}.
\enda
\end{align*}

\begin{proposition} 
Suppose that the same assumptions as in Theorem~\ref{T.C.Linf} hold.  
Then each component of the vorticity, $\w^{\e}_A$ and $\w^{\e}_B$, defined in Definition~\ref{D.wAwB}, satisfies 
\begin{align}
&\w^{\e}_A \quad \rightarrow \quad \w^{E} \quad \mbox{strongly in} \quad  L^\infty(0,T;L^p(\R^2)), \quad \mbox{for} \quad  2\leq p<\infty, \label{wAstrong}
\\
&\w^{\e}_B \quad \rightarrow \quad 0 \quad  \hspace{3mm} \mbox{strongly in} \quad L^\infty(0,T;L^p(\R^2)), \quad \mbox{for} \quad  2\leq p<\infty, \label{wBto0Yudo}
\end{align}
where $w^{E}$ is the unique weak solution to the vorticity--stream formulation of the Euler equation~\eqref{incompE} with initial data $\w^E|_{t=0}=\w_0$ satisfying \eqref{wEYudospace}.
\end{proposition}

If we prove this proposition, we directly obtain the result~\eqref{wstrongYudo}.
Before we prove \eqref{wAstrong}, we first observe that for the solution to the Boltzmann equation constructed in Theorem~\ref{T.2D.global} under the assumption $\|\w_0\|_{L^\infty\cap L^1(\R^2)}<\infty$, the vorticity $\w^{\e}_A$, defined in Definition~\ref{D.wAwB}, is uniformly bounded in both $L^1(\R^2)$ and $L^\infty(\R^2)$ (see \eqref{U.wLp<2} in Lemma \ref{L.unif}):
\begin{align*}
\sup_{0\leq t\leq T}\|\w^{\e}_A(t)\|_{L^1\cap L^\infty(\R^2)} &\leq  \|\w^{\e}_0\|_{L^1\cap L^\infty(\R^2)} \leq \|\w_0\|_{L^1\cap L^\infty(\R^2)} .
\end{align*}
On the other hand, considering the full vorticity $\w^{\e}=\w^{\e}_A+\w^{\e}_B$, we obtain a uniform bound in $L^2(\R^2)\cap L^\infty(\R^2)$:
\begin{align}\label{wL2infbdd}
\bega
\sup_{0\leq t\leq T}\|\w^{\e}(t)\|_{L^2\cap L^\infty(\R^2)} &\leq \sup_{0\leq t\leq T}\|\w^{\e}_A(t)\|_{L^2\cap L^\infty(\R^2)} + \sup_{0\leq t\leq T}\|\w^{\e}_B(t)\|_{L^2\cap L^\infty(\R^2)} \cr 
&\leq C(\|\w_0\|_{L^2(\R^2)},\|\w_0\|_{L^\infty(\R^2)},T),
\enda
\end{align}
where we used \eqref{U.wLp<2} in Lemma \ref{L.unif}.

To prove \eqref{wAstrong}, we apply the interpolation inequality, to have 
\begin{align}\label{wdinterp}
\bega
\|(\w^{\e}_A-\w^E)(t)\|_{L^p(\R^2)} &\leq \|(\w^{\e}_A-\w^E)(t)\|_{L^2(\R^2)}^{\frac{2}{p}} \|(\w^{\e}_A-\w^E)(t)\|_{L^\infty(\R^2)}^{1-\frac{2}{p}} \cr 
&\leq \big(2\|\w_0\|_{L^\infty(\R^2)}\big)^{1-\frac{2}{p}} \|(\w^{\e}_A-\w^E)(t)\|_{L^2(\R^2)}^{\frac{2}{p}}.
\enda
\end{align}
Thus, it suffices to prove convergence in $L^2(\R^2)$. 
To establish the $L^2(\R^2)$ convergence of $(\w^{\e}_A-\w^E)$, let $\w^{\e,\lambda}_A$ and $\w^{E,\lambda}$ be the unique solutions of the following problems with mollified initial data, respectively:
\begin{align}\label{weqneE}
\bega
&\p_t \w^{\e,\lambda}_A + \mathbb{P}u^{\e}\cdot\nabla_x \w^{\e,\lambda}_A = 0, \qquad \hspace{3mm}
\w^{\e,\lambda}_A(0,x) = \w^{\e}_0 \ast \varphi^\lambda(x), \cr 
&\p_t \w^{E,\lambda} + u^{E} \cdot \w^{E,\lambda} = 0, \qquad \hspace{8mm} \w^{E,\lambda}(0,x) = \w^E_0 \ast \varphi^\lambda(x),
\enda
\end{align}
where $\varphi^\lambda$ is mollifier with rate $\lambda$ defined in \eqref{molli}.
Using $\w^{\e,\lambda}_A$ and $\w^{E,\lambda}$ as intermediate approximations, we obtain
\begin{align}\label{wthree}
\|(\w^{\e}_A-\w^E)(t)\|_{L^2(\R^2)} \leq \|(\w^{\e}_A-\w^{\e,\lambda}_A)(t)\|_{L^2(\R^2)} + \|(\w^{\e,\lambda}_A-\w^{E,\lambda})(t)\|_{L^2(\R^2)} + \|(\w^{E,\lambda}-\w^{E})(t)\|_{L^2(\R^2)}.
\end{align}
Since the first and the last terms converge to zero as $\lambda \to 0$ for each fixed $\e$, namely,
\begin{align*}
\|(\w^{\e}_A-\w^{\e,\lambda}_A)(t)\|_{L^2(\R^2)} &\leq \|(\w^{\e}_0-\w^{\e}_0 \ast \varphi^\lambda)\|_{L^2(\R^2)}, \cr
\|(\w^{E,\lambda}-\w^{E})(t)\|_{L^2(\R^2)} &\leq \|(\w^{E}_0-\w^E_0 \ast \varphi^\lambda)\|_{L^2(\R^2)},
\end{align*}
it remains to prove the convergence of $\|(\w^{\e,\lambda}_A-\w^{E,\lambda})(t)\|_{L^2(\R^2)}$.

Before proceeding to the proofs of Theorem~\ref{T.C.Linf}, we first establish a lemma that describes the fundamental properties of 
$(u^{\e}_A, u^{\e}_B)$ and $(\w^{\e}_A, \w^{\e}_B)$. 
Although \eqref{wBto0Yudo} was already proved in \eqref{U.wLp<2} of Lemma \ref{L.unif}, we also show that $u^{\e}_B$ converges to $0$.

\begin{lemma}\label{L.uABprop}
For $(u^{\e}_A, u^{\e}_B)$ and $(\w^{\e}_A, \w^{\e}_B)$ defined in Definitions~\ref{D.wAwB}, the following properties hold:
\begin{enumerate}

\item Both $u^{\e}_A$ and $u^{\e}_B$ are incompressible vector fields.

\item The Leray projection of $u^{\e}$ decomposes as $\mathbb{P}u^{\e}(t,x) = u^{\e}_A(t,x) + u^{\e}_B(t,x).$

\item The velocity fields $u^{\e}_A$ and $u^{\e}_B$ satisfy
\begin{align}\label{uAeqn}
\bega
\begin{cases}
\p_tu^{\e}_A + u^{\e}_A \cdot \nabla_x u^{\e}_A + \nabla_x p^{\e}_A = -\nabla_x^{\perp}(-\Delta_x)^{-1}\nabla_x \cdot \big(u^{\e}_B\w^{\e}_A\big), \\
u^{\e}_A(t,x)|_{t=0} = {\bf K} \ast \w^{\e}_0 , \end{cases}
\enda
\end{align}
and
\begin{align}\label{uBeqn}
\bega
\begin{cases} 
\p_tu^{\e}_B + u^{\e}_B \cdot \nabla_x u^{\e}_B + \nabla_x p^{\e}_B = -\nabla_x^{\perp}(-\Delta_x)^{-1}\nabla_x \cdot(u^{\e}_A\w^{\e}_B)+\nabla_x^{\perp}(-\Delta_x)^{-1}\bar{\varPi}_{\w}^{\e}, \\
u^{\e}_B(t,x)|_{t=0} = 0 , \end{cases}
\enda
\end{align}
where $\Delta_x p^{\e}_A = \sum_{i,j} \p_i \p_j (u^{\e}_{A,i} u^{\e}_{A,j})$ and 
$\Delta_x p^{\e}_B = \sum_{i,j} \p_i \p_j (u^{\e}_{B,i} u^{\e}_{B,j})$.

\item $\w^{\e}_A$ preserves the total circulation, whereas $\w^{\e}_B$ has zero circulation:
\begin{align*}
\int_{\R^2} \w^{\e}_A dx = \int_{\R^2} \w^{\e}_0 dx, \qquad \int_{\R^2} \w^{\e}_B dx = 0.
\end{align*}

\item The forcing component $u^{\e}_B$ vanishes in the limit:
\begin{align}\label{uBto0}
\sup_{0\leq t\leq T_*}\|u^{\e}_B(t)\|_{L^p_x} \leq \mathcal{C}_* \Big(\e^{\frac{1}{4}-} +\kappa^{\frac{1}{2}-}\Big) \to 0, \quad \mbox{as} \quad \e\to 0, \quad 2\leq p < \infty,
\end{align}
where $\mathcal{C}_*$ is defined in \eqref{C*def}.

\end{enumerate}
\end{lemma}

\hide
\begin{remark}
We can rewrite equation \eqref{weqnA} in the following form:
\begin{align}\label{weqnAB}
\bega
\p_t\w^{\e}_A + u^{\e}_A\cdot \nabla_x \w^{\e}_A = -u^{\e}_B\cdot \nabla_x \w^{\e}_A , \qquad \w^{\e}_A(t,x)|_{t=0} = \w^{\e}_0(x).
\enda
\end{align}
\end{remark}
\unhide

\begin{proof}
Properties~(i)–(iii) follow directly from the definitions.  
For~(iv), by Lemma~\ref{L.util} we have $\int_{\R^2} \w^{\e} dx = \int_{\R^2} \w^{\e}_0 dx $. 
From the equation for $\w^{\e}_A$ in \eqref{weqnA}, it follows that 
$\int_{\R^2} \w^{\e}_A\, dx = \int_{\R^2} \w^{\e}_0\, dx$. 
Since $\w^{\e} = \w^{\e}_A + \w^{\e}_B$, we immediately obtain property~(iv).  
 \\
 For the proof of~(v), we apply \eqref{puw} and \eqref{U.wLp<2} from Lemma \ref{L.unif} to deduce
\begin{align}\label{duBto0}
\bega
\sup_{t\in[0,T]}\|\nabla_xu^{\e}_B(t)\|_{L^p_x} \leq C\sup_{t\in[0,T]}\|\w^{\e}_B(t)\|_{L^p_x} \leq \mathcal{C}_* \Big(\e^{\frac{1}{4}-} +\kappa^{\frac{1}{2}-}\Big) \to 0 , \quad \mbox{for} \quad 2\leq p < \infty.
\enda
\end{align}
To handle the low-frequency part of $u^{\e}_B$, starting from \eqref{uBeqn} we derive the energy estimate
\begin{align}\label{uBL2R}
\bega
\frac{d}{dt}\|u^{\e}_B\|_{L^2_x} \leq \|\nabla_x^{\perp}(-\Delta_x)^{-1}\nabla_x \cdot(u^{\e}_A\w^{\e}_B)\|_{L^2_x} + \| \nabla_x^{\perp}(-\Delta_x)^{-1}\bar{\varPi}_{\w}^{\e}\|_{L^2_x},
\enda
\end{align}
where we have used $\nabla_x \cdot u^{\e}_B = 0$. 
For the first term in \eqref{uBL2R}, we apply Plancherel’s theorem.  
From~(iv), we know that $\w^{\e}_B$ has zero circulation, hence the radial eddy $\bar{u}$ must belong to $u^{\e}_A$. 
Decomposing $u^{\e}_A = (u^{\e}_A - \bar{u}) + \bar{u}$, we obtain
\begin{align*}
\bega
\|\nabla_x^{\perp}(-\Delta_x)^{-1}\nabla_x \cdot(u^{\e}_A\w^{\e}_B)\|_{L^2_x} \leq \|u^{\e}_A\w^{\e}_B\|_{L^2_x} &\leq \|(u^{\e}_A-\bar{u})\w^{\e}_B\|_{L^2_x} + \|\bar{u}\w^{\e}_B\|_{L^2_x} \cr 
&\leq \|(u^{\e}_A-\bar{u})\|_{L^2_x}\|\w^{\e}_B\|_{L^\infty_x} + \|\bar{u}\|_{L^3_x}\|\w^{\e}_B\|_{L^6_x} \cr 
&\les \mathcal{C}_* \Big(\e^{\frac{1}{4}-} +\kappa^{\frac{1}{2}-}\Big) \to 0 , \quad \mbox{as} \quad \e\to0,
\enda
\end{align*}
where we have used \eqref{U.wLp<2} and Lemma \ref{L.barubdd}.
For the second term in \eqref{uBL2R}, note that $\bar{\varPi}_{\w}^{\e}$ is a gradient form by its definition in \eqref{barvarPidef}. 
Hence, applying Plancherel’s theorem yields
\begin{align*}
\bega
\big\|\nabla_x^{\perp}(\Delta_x)^{-1}\bar{\varPi}_{\w}^{\e}\big\|_{L^2_x} &\leq \bigg\| (\mathbb{P}^{\perp}u^{\e})\w^{\e} + \bigg(\frac{k_B}{\e}(\mathrm{\Theta}^{\e}-1)\nabla_x(\rho^{\e}+\ta^{\e})
+\frac{1}{\e^2}\frac{1}{\mathrm{P}^{\e} }\sum_{j} \p_{x_j} \mathbf{r}_{ij}^{\e}\bigg) \bigg\|_{L^2_x} \cr 
&\les \|\mathbb{P}^{\perp}u^{\e}\|_{L^\infty_x}\|\w^{\e}\|_{L^2_x} + \frac{1}{\e}\|(\mathrm{\Theta}^{\e}-1)\|_{L^2_x}\|\nabla_x(\rho^{\e}+\ta^{\e})\|_{L^\infty_x}
+\frac{1}{\e^2}\sum_{j} \|\p_{x_j} \mathbf{r}_{ij}^{\e}\|_{L^2_x} \cr 
&\les \e^{\frac{1}{4}-}\mathcal{E}_M^{\frac{1}{2}} + \kappa^{\frac{1}{2}}\mathcal{D}_G^{\frac{1}{2}}.
\enda
\end{align*}
Thus, in the same way as in~\eqref{C*def}, we obtain
\begin{align}\label{uBL2to0}
\bega
\sup_{0\leq t\leq T}\|u^{\e}_B(t)\|_{L^2_x}  &\leq \|u^{\e}_B(0)\|_{L^2_x} + T \kappa^{\frac{1}{2}-} + T\e^{\frac{1}{4}-}\sup_{t\in[0,T]}\mathcal{E}_M^{\frac{1}{2}}(t) + \kappa^{\frac{1}{2}}\sqrt{T}\Big(\int_0^T\mathcal{D}_G(s)ds\Big)^{\frac{1}{2}} \cr
&\leq \mathcal{C}_* \Big(\e^{\frac{1}{4}-} +\kappa^{\frac{1}{2}-}\Big) \to 0 , \quad \mbox{as} \quad \e\to0,
\enda
\end{align}
since $u^{\e}_B(0) = 0$.  
We now claim that 
\begin{align}\label{uBLinfto0}
\bega
\sup_{0\leq t\leq T}\|u^{\e}_B(t)\|_{L^\infty_x} \leq \mathcal{C}_* \Big(\e^{\frac{1}{4}-} +\kappa^{\frac{1}{2}-}\Big) \to 0 , \quad \mbox{as} \quad \e\to0.
\enda
\end{align}
If this claim holds, then combining \eqref{uBL2to0} and \eqref{uBLinfto0} yields the desired result~\eqref{uBto0}.  
To prove the claim~\eqref{uBLinfto0}, we decompose
\begin{align*}
u^{\e}_B= {\bf K}\ast \w^{\e}_B = ({\bf K}\mathbf{1}_{|x|\leq1})\ast \w^{\e}_B + ({\bf K}\mathbf{1}_{|x|>1})\ast \w^{\e}_B,
\end{align*}
and apply Young’s convolution inequality:
\begin{align*}
\bega
\|u^{\e}_B\|_{L^\infty_x} &= \|{\bf K}_1 \ast \w^{\e}_B\|_{L^\infty_x} + \|{\bf K}_2 \ast \w^{\e}_B\|_{L^\infty_x}\leq \|{\bf K}_1 \|_{L^1_x} \|\w^{\e}_B\|_{L^\infty_x} + \|{\bf K}_2\|_{L^3_x}  \|\w^{\e}_B\|_{L^{\frac{3}{2}}_x}.
\enda
\end{align*}
The first term is bounded by $\kappa^{\frac{1}{2}-}$, since ${\bf K}_1$ is bounded in $L^1$ and \eqref{duBto0} holds.  
For the second term, by applying \eqref{puw} together with the Gagliardo–Nirenberg interpolation inequality, we obtain
\begin{align*}
\bega
\|\w^{\e}_B\|_{L^{\frac{3}{2}}_x}\leq C\|\nabla_xu^{\e}_B\|_{L^\frac{3}{2}_x} &\leq C \|\nabla_xu^{\e}_B\|_{L^2_x}^{\frac{2}{3}} \|u^{\e}_B\|_{L^2_x}^{\frac{1}{3}} \les \mathcal{C}_* \Big(\e^{\frac{1}{4}-} +\kappa^{\frac{1}{2}-}\Big),
\enda
\end{align*}
where we used \eqref{duBto0} and \eqref{uBL2to0}.  
This completes the proof of the claim~\eqref{uBLinfto0}, and hence of the lemma.
\end{proof}

We note that the inviscid limit for Yudovich solutions on the domain $\T^2$ was established in \cite{2020Constan}, where an explicit convergence rate was obtained under additional regularity assumptions. 
To prove the convergence of $\w^{\e}_A$, we require the following three lemmas. In particular, we extend the result of \cite{2020Constan}, originally established on $\T^2$, to $\R^2$ for the equation of $\w^{\e}_A$ in the presence of microscopic effects in the velocity field $\mathbb{P}u^{\e}$.

\begin{lemma}\label{L.expdu} 
For the divergence-free velocity field $u = {\bf K} \ast \w$ with ${\bf K}(x):= \frac{1}{2\pi}\frac{x^{\perp}}{|x|^2}$ associated with a vorticity $\w \in (L^2 \cap L^\infty)(\R^2)$ with $\|\w\|_{L^\infty(\R^2)}\neq 0$, the following bound holds:
\begin{align}\label{eduest}
\int_{\R^2}\Big(e^{\mathbf{b}|\nabla_x u|}-1\Big)dx &\leq C(\mathbf{b},\|\w\|_{L^2\cap L^\infty(\R^2)}), \quad \mbox{for} \quad \mathbf{b} \in \bigg(0,\frac{1}{C_*e\|\w\|_{L^\infty(\R^2)}}\bigg).
\end{align}
\end{lemma}
\begin{proof}
Using Taylor expansion and the estimate $\|\nabla_x u\|_{L^p} \leq C_* p \|\w\|_{L^p}$ in \eqref{PotenLp}, we have
\begin{align*}
\int_{\R^2}(e^{\mathbf{b}|\nabla_x u|}-1)dx &= \sum_{p=1}^{\infty}\int_{\R^2} \frac{\mathbf{b}^p|\nabla_x u|^p}{p!}dx = \sum_{p=1}^{\infty} \frac{\mathbf{b}^p\|\nabla_x u\|_{L^p(\R^2)}^p}{p!} \leq \sum_{p=1}^{\infty} \frac{\Big(C_*p\mathbf{b}\|\w\|_{L^p(\R^2)}\Big)^p}{p!} \cr 
&\leq \sum_{p=1}^{\infty} \frac{C_*^pp^p\mathbf{b}^p}{p!} \bigg(\|\w\|_{L^2(\R^2)}^\frac{2}{p}\|\w\|_{L^\infty(\R^2)}^{1-\frac{2}{p}}\bigg)^{p}  \leq \frac{\|\w\|_{L^2(\R^2)}^2}{\|\w\|_{L^\infty(\R^2)}^2}\sum_{p=1}^{\infty} \frac{\Big(C_*p\mathbf{b}\|\w\|_{L^\infty(\R^2)}\Big)^p}{p!} ,
\end{align*}
where we used the interpolation inequality. By Stirling’s bound $p! \geq \sqrt{2\pi}p^{p+\frac{1}{2}}e^{-p}$, the following series is convergent (for $C_*e \mathbf{b}\|\w\|_{L^\infty(\R^2)}<1$):
\begin{align*}
\sum_{p=1}^{\infty} \frac{\Big(C_*p\mathbf{b}\|\w\|_{L^\infty}\Big)^p}{p!} 
\leq \sum_{p=1}^{\infty} \frac{1}{\sqrt{2\pi}} p^{-\frac{1}{2}}\Big(C_*e \mathbf{b}\|\w\|_{L^\infty(\R^2)}\Big)^p \leq \frac{1}{\sqrt{2\pi}} \frac{C_*e\mathbf{b}\|\w\|_{L^\infty(\R^2)}}{1-C_*e\mathbf{b}\|\w\|_{L^\infty(\R^2)}}.
\end{align*}
Finally, we conclude that, for $\|\w\|_{L^\infty}\neq 0$,
\begin{align*}
\int_{\R^2}\Big(e^{\mathbf{b}|\nabla_x u|}-1\Big)dx &\leq \frac{\|\w\|_{L^2(\R^2)}^2}{\|\w\|_{L^\infty(\R^2)}^2} \frac{1}{1-C_*e\mathbf{b}\|\w\|_{L^\infty}}, \quad \mbox{for} \quad \mathbf{b} \in \bigg(0,\frac{1}{C_*e\|\w\|_{L^\infty(\R^2)}}\bigg).
\end{align*} 
\end{proof}

\begin{lemma}\label{L.dwesti}
Let $\{\w^{\e,\lambda}_A\}_{\lambda}$ be the family of solutions to $\eqref{weqneE}_1$ satisfying $\w^{\e,\lambda}_A \in (L^1\cap L^\infty)(\R^2)$. Then, we have
\begin{align*}
\sup_{0\leq t \leq \bar{T}}\|\nabla_x\w^{\e,\lambda}_A(t)\|_{L^{2p(t)}}^{2p(t)} \leq C _1(\mathbf{b},p_0,\bar{T})\lambda^{-4p_0+2}  + C_2(\mathbf{b},p_0,\bar{T},\|\w^{\e}\|_{L^2\cap L^\infty(\R^2)}),
\end{align*}
where $\lambda$ denotes the mollification scale of the initial data, and the monotone decreasing function $p(t)$ and the time $\bar{T}$ are defined by 
\begin{align}\label{p(t)def}
p(t) := \frac{\mathbf{b} p_0}{\mathbf{b}+2p_0 t}, \qquad p(0)=p_0>1, \qquad \mbox{for} \quad t\in [0, \bar{T}], \quad \mbox{and} \quad \bar{T}:= \frac{\mathbf{b}(p_0-1)}{2p_0}.
\end{align}
Here $\mathbf{b}$ is the constant in \eqref{eduest}. 
\end{lemma}
\begin{proof}
We first apply $\nabla_x$ to $\eqref{weqneE}_1$ to obtain
\begin{align}\label{dweqnP}
\bega
&\p_t \nabla_x\w^{\e,\lambda}_A + \mathbb{P}u^{\e}\cdot \nabla_x (\nabla_x\w^{\e,\lambda}_A) = (\nabla_x\mathbb{P}u^{\e})\cdot \nabla_x \w^{\e,\lambda}_A .
\enda
\end{align}
For notational simplicity, we write $\mathrm{Z}:= \nabla_x\w^{\e,\lambda}_A$.  
For $p(t)$ defined in \eqref{p(t)def}, starting from \eqref{dweqnP} we have 
\begin{align}\label{Z2pE}
\bega
\frac{1}{2}\frac{d}{dt}&\int_{\R^2}|\mathrm{Z}|^{2p(t)}dx  =  p'(t) \int_{\R^2}(\ln |\mathrm{Z}|) |\mathrm{Z}|^{2p(t)}dx  + p(t) \int_{\R^2}(\p_t \mathrm{Z}) \cdot  \mathrm{Z}  |\mathrm{Z}|^{2p(t)-2}  dx 
 \cr 
&= p'(t) \int_{\R^2}(\ln |\mathrm{Z}|) |\mathrm{Z}|^{2p(t)} dx + p(t) \int_{\R^2}\Big( - \mathbb{P}u^{\e}\cdot \nabla_x\mathrm{Z} +  \nabla_x \mathbb{P}u^{\e} \cdot  \mathrm{Z} \Big) \cdot  \mathrm{Z}  |\mathrm{Z}|^{2p(t)-2}dx.
\enda
\end{align}
Since $\mathbb{P}u^{\e}$ is divergence-free, the following term vanishes:
\begin{align}\label{Z2pE-2}
\int_{\R^2}(\mathbb{P}u^{\e}\cdot \nabla_x\mathrm{Z})  \cdot \mathrm{Z}  |\mathrm{Z}|^{2p(t)-2}dx = \frac{1}{2}\int_{\R^2}(\mathbb{P}u^{\e}\cdot \nabla_x|\mathrm{Z}|^{2p(t)})dx = 0.
\end{align}
For the term $\int |\nabla_x \mathbb{P}u^{\e}|\, |\mathrm{Z}|^{2p(t)}ds$ in \eqref{Z2pE}, we use the inequality
\begin{align*}
 XY \leq \begin{cases} (e^{X}-1)  + (Y\ln Y-Y+1) , \quad &Y\geq 1, \\ 
(e^{X}-1), \quad &0\leq Y< 1,
\end{cases} 
\end{align*}
and $Y\ln Y-Y+1 \leq Y\ln Y$ for $Y\geq 1$.  
Then, for $X= \mathbf{b}|\nabla_x \mathbb{P}u^{\e}|$ and $Y=\frac{1}{\mathbf{b}}|\mathrm{Z}|^{2p(t)}$, with $\mathbf{b}$ satisfying the condition \eqref{eduest}, we obtain
\begin{align}\label{PuZ}
\bega
\int_{\R^2}&|\nabla_x \mathbb{P}u^{\e}||\mathrm{Z}|^{2p(t)}dx \leq 2\int_{\R^2}(e^{\mathbf{b}|\nabla_x \mathbb{P}u^{\e}|}-1)dx + \int_{\R^2}\frac{1}{\mathbf{b}}|\mathrm{Z}|^{2p(t)}\ln \Big(\frac{1}{\mathbf{b}}|\mathrm{Z}|^{2p(t)}\Big) dx \cr 
&\leq 2 C(\mathbf{b},\|\w^{\e}\|_{L^2\cap L^\infty(\R^2)}) + \frac{2p(t)}{\mathbf{b}} \int_{\R^2}(\ln|\mathrm{Z}|) |\mathrm{Z}|^{2p(t)} dx + \mathbf{b}^{-1} \ln (\mathbf{b}^{-1}) \int_{\R^2}|\mathrm{Z}|^{2p(t)} dx,
\enda
\end{align}
where we have used Lemma \ref{L.expdu}. 
Combining \eqref{Z2pE-2} and \eqref{PuZ} with \eqref{Z2pE}, we obtain
\begin{align*}
\frac{1}{2}\frac{d}{dt}\|\mathrm{Z}\|_{L^{2p(t)}}^{2p(t)} &\leq \bigg(p'(t)+ \frac{2p^2(t)}{\mathbf{b}}\bigg) \int_{\R^2}(\ln |\mathrm{Z}|) |\mathrm{Z}|^{2p(t)}dx + p(t)\mathbf{b}^{-1} \ln (\mathbf{b}^{-1})\int_{\R^2}|\mathrm{Z}|^{2p(t)} dx \cr 
&+ p(t) \left( C(\mathbf{b},\|\w^{\e}\|_{L^2\cap L^\infty(\R^2)})\right) \cr 
&\leq p(t)\mathbf{b}^{-1} \ln (\mathbf{b}^{-1})\int_{\R^2}|\mathrm{Z}|^{2p(t)} dx + p(t) \left( C(\mathbf{b},\|\w^{\e}\|_{L^2\cap L^\infty(\R^2)})\right).
\end{align*}
Here we used that $p(t)$ is defined in~\eqref{p(t)def} so that 
$p'(t)=-2\mathbf{b}^{-1}p(t)^2$. 
Hence, the first term on the right-hand side of the first line vanishes.
Applying Gr\"{o}nwall's inequality for $t\in[0,\bar{T}]$, and using the identities $\int_0^t p(s) ds = \ln(1+\frac{2p_0t}{\mathbf{b}})^{2\mathbf{b}}$, we get
\begin{align*}
\|\mathrm{Z}(t)\|_{L^{2p(t)}}^{2p(t)} &\leq \bigg(1+\frac{2p_0t}{\mathbf{b}}\bigg)^{4\ln (\mathbf{b}^{-1})} \bigg(\|\mathrm{Z}_0\|_{L^{2p_0}}^{2p_0} + C(\mathbf{b},\|\w^{\e}\|_{L^2\cap L^\infty(\R^2)}) \bigg).
\end{align*} 
Recovering the notation $\mathrm{Z}:= \nabla_x\w^{\e,\lambda}_A$ gives 
\begin{align*}
\sup_{0\leq t \leq \bar{T}}\|\nabla_x\w^{\e,\lambda}_A(t)\|_{L^{2p(t)}}^{2p(t)} \leq C_1\|\nabla_x\w^{\e,\lambda}_A(0)\|_{L^{2p_0}}^{2p_0} + C_2.
\end{align*}
Since we mollified the initial data with respect to the rate $\lambda$, using \eqref{Ga-Ni}, we get 
\begin{align*}
\|\nabla_x\w^{\e,\lambda}_A(0)\|_{L^{2p_0}}^{2p_0} \leq C \Big(\|\nabla_x^2\w^{\e,\lambda}_A(0)\|_{L^2}^{1-\frac{1}{p_0}} \|\nabla_x\w^{\e,\lambda}_A(0)\|_{L^2}^{\frac{1}{p_0}} \Big)^{2p_0} \leq  C\Big( \lambda^{-2(1-\frac{1}{p_0})} \lambda^{-\frac{1}{p_0}} \Big)^{2p_0} \leq C \lambda^{-4p_0+2}. 
\end{align*}
This completes the proof. 
\hide
Let $A(t):= \frac{1}{2}\|\mathrm{Z}(t)\|_{L^{2p(t)}}^{2p(t)}$ to have 
\begin{align}
A'(t) &\leq p(t) \Big(C_1 A(t) + C_2 + \|\nabla_x\bar{\varPi}_{\w}^{\e}\|_{L^{2p(t)}}A^{\frac{2p(t)-1}{2p(t)}} \Big) \cr 
&\leq  p(t) \Big(\big(C_1+\|\nabla_x\bar{\varPi}_{\w}^{\e}\|_{L^{2p(t)}}\big) A(t) + C_2+\|\nabla_x\bar{\varPi}_{\w}^{\e}\|_{L^{2p(t)}} \Big).
\end{align}
(Note. To estimate $\bar{\varPi}_{\w}^{\e}$ we need $L^2_T$.)
where we used $A^{\frac{2p(t)-1}{2p(t)}} \leq 1+A(t)$ since $\frac{1}{2}\leq\frac{2p(t)-1}{2p(t)}< 1$. 
\begin{align}
A(t) \leq e^{C_1\int_0^tp(s)ds}A_0 + C_2 \int_0^t e^{C_1\int_s^tp(\tau)d\tau}  p(s) ds
\end{align}
\begin{align}
\frac{1}{2}\|\mathrm{Z}(t)\|_{L^{2p(t)}}^{2p(t)} &\leq e^{C_1\int_0^tp(s)ds}\frac{1}{2}\|\mathrm{Z}\|_{L^{2p_0}}^{2p_0} + C_2 \int_0^t e^{C_1\int_s^tp(\tau)d\tau}  p(s) ds \cr 
&= \bigg(1+\frac{2p_0t}{\mathbf{b}}\bigg)^{\frac{C_1\mathbf{b}}{2}}\frac{1}{2}\|\mathrm{Z}\|_{L^{2p_0}}^{2p_0} + C_2 \int_0^t \bigg(\frac{1+\frac{2p_0t}{\mathbf{b}}}{1+\frac{2p_0s}{\mathbf{b}}}\bigg)^{\frac{C_1\mathbf{b}}{2}}  p(s) ds
\end{align}
\begin{align}
\frac{1}{2}\|\nabla_x\w^{\e}(t)\|_{L^{2p(t)}}^{2p(t)} &\leq  \bigg(1+\frac{2p_0t}{\mathbf{b}}\bigg)^{\frac{C_1\mathbf{b}}{2}}\frac{1}{2}\|\nabla_x\w^{\e}_0\|_{L^{2p_0}}^{2p_0} + C_2 \int_0^t \bigg(\frac{1+\frac{2p_0t}{\mathbf{b}}}{1+\frac{2p_0s}{\mathbf{b}}}\bigg)^{\frac{C_1\mathbf{b}}{2}}  p(s) ds
\end{align}
where we used $\int_0^t p(s) ds = \ln(1+\frac{2p_0t}{\mathbf{b}})^{\frac{2}{\mathbf{b}}}. $
We take power $\frac{1}{2p(t)}$. 
\unhide
\end{proof}

\begin{lemma}\label{L.Puconv}
Suppose the same assumptions as in Theorem \ref{T.C.Linf}. Then the following estimate holds:
\begin{align*}
\bega
\|(\mathbb{P}u^{\e}-u^{E})(t)\|_{L^2_x} \les e^{2-2e^{-Ct}} \bigg(\|\mathbb{P}u^{\e}_0-u^{E}_0\|_{L^2_x}+\mathcal{C}_* \Big(\e^{\frac{1}{4}-} +\kappa^{\frac{1}{2}-}\Big) \bigg)^{e^{-Ct}},
\enda
\end{align*}
where $C>0$ depends only on
$\|\w^{\e}_0\|_{L^2\cap L^\infty(\R^2)}$ and $\|\w^E_0\|_{L^2\cap L^\infty(\R^2)}$, and $\mathcal{C}_*$ is defined in~\eqref{C*def}.
\end{lemma}

\hide
\begin{remark}\label{Rmk.NS,BE}
For the inviscid limit case, we will have 
\begin{align*}
\bega
\frac{d}{dt}\|u^{\mathsf{d}}(t)\|_{L^2_x}^2 \kappa \|\nabla_xu^{\mathsf{d}}(t)\|_{L^2_x}^2 \leq \int_{\R^2} |\nabla_x u^{\nu}| |u^{\mathsf{d}}|^2 dx .
\enda
\end{align*}
On the other hand, for BE case, we will have
\begin{align*}
\bega
\frac{d}{dt}\|u^{\mathsf{d}}(t)\|_{L^2_x}^2 \leq \int_{\R^2} |\nabla_x \mathbb{P}u^{\e}| |u^{\mathsf{d}}|^2 dx + \frac{1}{\e^2}\|\p_{x_j} \mathbf{r}_{ij}^{\e}\|_{L^2_x} \|u^{\mathsf{d}}(t)\|_{L^2_x}.
\enda
\end{align*}
Here the scale of the term $\|\p_{x_j} \mathbf{r}_{ij}^{\e}\|_{L^2_TL^2_x} \sim \kappa^{\frac{1}{2}}$. Hence, to get the correct convergence rate, we should divide each side by $\|u^{\mathsf{d}}(t)\|_{L^2_x}$ for BE case. {\color{blue} 무슨 말을 하고 싶은지는 알겠는데, 너무 엘러멘터리하게 코멘트가 되어있네요. 제가 나중에 정리해서 인트로에 섹시해보이게 쓸게요. 그러면 이 리마크는 지워야하겠습니다. }
\end{remark}
\unhide

\begin{proof}
We follow the argument in \cite{Chemin}.
By the same way to Lemma \ref{L.Pu-uconv}, writing the energy estimate for $u^{\mathsf{d}}$ from \eqref{barueqn} gives 
\begin{align}\label{udE2}
\bega
\frac{1}{2}\frac{d}{dt}\|u^{\mathsf{d}}\|_{L^2_x}^2 &\leq  \int_{\R^2} \bigg(\mathbb{P}\big(\mathbb{P}u^{\e}\cdot \nabla_x \mathbb{P}u^{\e}\big) - u^E\cdot\nabla_x u^E\bigg) u^{\mathsf{d}} dx + \mathfrak{Q}(t) \|u^{\mathsf{d}}\|_{L^2_x}.
\enda
\end{align}
We consider the first term of \eqref{udE2}. 
Using that $\mathbb{P}$ is self-adjoint and that $\mathbb{P}(\mathbb{P}u^{\e}-u^E)=\mathbb{P}u^{\e}-u^E$, we have
\begin{align*}
\bega
\int_{\R^2} \bigg(\mathbb{P}\big(\mathbb{P}u^{\e}\cdot \nabla_x \mathbb{P}u^{\e}\big) - u^E\cdot\nabla_x u^E\bigg) u^{\mathsf{d}} dx 
&= \int_{\R^2} \big( u^{\mathsf{d}}\cdot \nabla_x \mathbb{P}u^{\e}\big) u^{\mathsf{d}} dx + \int_{\R^2} \big(u^E\cdot \nabla_x u^{\mathsf{d}} \big)u^{\mathsf{d}} dx \cr 
&\leq \int_{\R^2} |\nabla_x \mathbb{P}u^{\e}| |u^{\mathsf{d}}|^2 dx.
\enda
\end{align*}
Here, the last term on the first line vanishes since $u^E$ is divergence-free.
Applying H\"{o}lder's inequality for any $2< \gamma < \infty$, we obtain 
\begin{align}\label{Itineq}
\bega
\int_{\R^2} |\nabla_x \mathbb{P}u^{\e}| |u^{\mathsf{d}}|^2 dx &\les \|\nabla_x\mathbb{P}u^{\e}\|_{L^\gamma_x} \|u^{\mathsf{d}}\|_{L^{\frac{2\gamma}{\gamma-2}}_x} \|u^{\mathsf{d}}\|_{L^2_x}  
\leq C\gamma \|\w^{\e}_0\|_{L^2_x\cap L^\infty_x} \|u^{\mathsf{d}}\|_{L^{\frac{2\gamma}{\gamma-2}}_x} \|u^{\mathsf{d}}\|_{L^2_x},
\enda
\end{align}
where we used $\|\nabla_x\mathbb{P}u^{\e}\|_{L^\gamma_x}\leq C\gamma \|\w^{\e}\|_{L^\gamma_x}$ from \eqref{PotenLp} and the interpolation inequality.
Next, by the Gagliardo--Nirenberg interpolation inequality in \eqref{Ga-Ni}, we have  
\begin{align}\label{u-uLp}
\bega
\|u^{\mathsf{d}}\|_{L^{\frac{2\gamma}{\gamma-2}}_x}&\leq C\|\nabla_x u^{\mathsf{d}}\|_{L^2_x}^{\frac{2}{\gamma}}\|u^{\mathsf{d}}\|_{L^2_x}^{\frac{\gamma-2}{\gamma}} \leq \big(C\|\w^{\mathsf{d}}\|_{L^2_x}\big)^{\frac{2}{\gamma}}\|u^{\mathsf{d}}\|_{L^2_x}^{\frac{\gamma-2}{\gamma}} ,
\enda
\end{align}
where we also used $\|\nabla_x u^{\mathsf{d}}\|_{L^2_x}\leq C\|\w^{\mathsf{d}}\|_{L^2_x}$ from \eqref{PotenLp}.
We note that $\|\w^{\e}\|_{L^2_x}$ is uniformly bounded; see~\eqref{wL2infbdd}.
For brevity, we define an upper bound for $\|\w^{\mathsf{d}}\|_{L^2_x}$ by $C_{\infty}$ as follows:
\begin{align*}
\sup_{0\leq t\leq T}\big(\|\w^{\e}(t)\|_{L^2_x}+\|\w^E(t)\|_{L^2_x}\big) &\leq \sup_{0\leq t\leq T}\|\w^{\e}_A(t)\|_{L^2(\R^2)} + \sup_{0\leq t\leq T}\|\w^{\e}_B(t)\|_{L^2(\R^2)} + \sup_{0\leq t\leq T}\|\w^E(t)\|_{L^2(\R^2)} \cr 
&\leq C(\|\w_0\|_{L^2(\R^2)},T) + \mathcal{C}_* \Big(\e^{\frac{1}{4}-} +\kappa^{\frac{1}{2}-}\Big) := C_{\infty},
\end{align*}
where $\mathcal{C}_*$ is defined in~\eqref{C*def}.
Substituting \eqref{u-uLp} into \eqref{Itineq} yields
\begin{align}\label{udmain}
\bega
\int_{\R^2} |\nabla_x \mathbb{P}u^{\e}| |u^{\mathsf{d}}|^2 dx &\les C C_{\infty}^{1+\frac{2}{\gamma}} \gamma  \Big(\|u^{\mathsf{d}}\|_{L^2_x}^{1+\frac{\gamma-2}{\gamma}}\Big).
\enda
\end{align}
Combining \eqref{udE2} and \eqref{udmain}, and then dividing both sides by $\|u^{\mathsf{d}}\|_{L^2_x}$, we obtain 
\begin{align}\label{udE3}
\bega
\frac{d}{dt}\|u^{\mathsf{d}}\|_{L^2_x} &\leq  C C_{\infty}^{1+\frac{2}{\gamma}} \gamma \|u^{\mathsf{d}}\|_{L^2_x}^{\frac{\gamma-2}{\gamma}} + \mathfrak{Q},
\enda
\end{align}
where $\mathfrak{Q}$ is defined in \eqref{udforce}. 
In fact, \eqref{udE3} holds for any $\gamma>2$. In particular, we may choose $\gamma(t) = 2- \log \frac{\|u^{\mathsf{d}} (t)\|_{L^2_x} }{C_2}>2 $ as long as $0<\|u^{\mathsf{d}}(t)\|_{L^2_x}<C_2$. Then we have 
\begin{align*}
\bega
\frac{d}{dt} \|u^{\mathsf{d}}\|_{L^2_x} \leq CC_{\infty}^2 \gamma(t) \|u^{\mathsf{d}}\|_{L^2_x}^{1-\frac{2}{\gamma(t)}} + \mathfrak{Q} &\leq CC_{\infty}^2 \left(2-\log \frac{\|u^{\mathsf{d}}\|_{L^2_x}}{C_2}\right) \|u^{\mathsf{d}}\|_{L^2_x} \|u^{\mathsf{d}}\|_{L^2_x}^{-\frac{2}{\gamma}} + \mathfrak{Q} \cr 
&\leq C (eC_{\infty})^2   \left(2-\log \frac{\|u^{\mathsf{d}}\|_{L^2_x}}{C_2}\right)\|u^{\mathsf{d}}\|_{L^2_x} + \mathfrak{Q},
\enda
\end{align*}
where we have used $\left(\frac{\|u^{\mathsf{d}}\|_{L^2_x}}{C_2}\right)^{-\frac{2}{\gamma}}=X^{-\frac{1}{2-\log X}} \leq X^{\frac{1}{\log X}} = e$ for $X(t)=\frac{\|u^{\mathsf{d}}(t)\|_{L^2_x}}{C_2}<1 $. Therefore, by Osgood's lemma, we obtain
\begin{equation}\label{Osgood}
   \|u^{\mathsf d}(t)\|_{L^2_x}
\le
C_2\,e^{\,2 - 2 e^{-(eC_\infty)^2 C t}}\,
\Bigg(
\frac{\|u^{\mathsf d}_0\|_{L^2_x}}{C_2}
+ \int_0^t \frac{\mathfrak{Q}(s)}{C_2}\,ds
\Bigg)^{e^{-(eC_\infty)^2 C t}}. 
\end{equation}
Applying the estimate of $\int_0^t \mathfrak{O}(s)\,ds$ from Lemma~\ref{L.forcingC} to \eqref{Osgood} completes the proof.
\end{proof}

\begin{lemma}\label{L.wAwB}
Let $\w^{\e}_A$ and $\w^{\e}_B$ be the solutions of \eqref{weqnA} and \eqref{weqnB}, respectively. Then $\w^{\e} = \w^{\e}_A + \w^{\e}_B$ is a solution of \eqref{weqnP}, and it can be represented as
\begin{align*}
\bega
\w^{\e}_A(t,x):= \w^{\e}_0((X^{\e})^{-1}(t,x)), \qquad  \w^{\e}_B(t,x):=\int_0^t \bar{\varPi}_{\w}^{\e}(s,X^{\e}(s,\mathbf{a})) ds.
\enda
\end{align*}
Here, $X^{\e}$ denotes the flow map associated with the divergence-free vector field $\mathbb{P}u^{\e}$, defined by
\begin{align}\label{traju}
\bega
\frac{d}{dt}X^{\e}(t,\mathbf{a}) &= \mathbb{P}u^{\e}(t,X^{\e}(t,\mathbf{a})), \qquad  X^{\e}(0,\mathbf{a}) = \mathbf{a}.
\enda
\end{align}
\end{lemma}

\begin{proof}[{\bf Proof of Theorem \ref{T.C.Linf}}]
Since we have already proved the convergence of $\w^{\e}_B$ in~\eqref{uBto0}, it suffices to establish the convergence of
$\|(\w^{\e,\lambda}_A-\w^{E,\lambda})(t)\|_{L^2(\R^2)}$ in \eqref{wthree}. For simplicity, set
\begin{align*}
\w^\mathsf{d,\lambda}(t,x):=\w^{\e,\lambda}_A(t,x)-\w^{E,\lambda}(t,x).
\end{align*}
Subtracting $\eqref{weqneE}_2$ from $\eqref{weqneE}_1$ yields
\begin{align*}
\bega
\p_t \w^\mathsf{d,\lambda} &+ u^{E,\lambda} \cdot \nabla_x \w^\mathsf{d,\lambda} + (\mathbb{P}u^{\e}-u^E) \cdot \nabla_x \w^{\e,\lambda}_A = 0.
\enda
\end{align*}
The corresponding $L^2$ energy estimate gives
\begin{align*}
\bega
\frac{1}{2}\frac{d}{dt}\|\w^\mathsf{d,\lambda}\|_{L^2(\R^2)}^2 &\leq 
-\int_{\R^2}((\mathbb{P}u^{\e}-u^E) \cdot \nabla_x \w^{\e,\lambda}_A) \w^{\mathsf{d},\lambda}dx
 \cr 
&\leq \|\mathbb{P}u^{\e}-u^E\|_{L^2(\R^2)} \|\nabla_x \w^{\e,\lambda}_A\|_{L^{2p(t)}(\R^2)} \|\w^{\mathsf{d},\lambda}\|_{L^{\frac{2p(t)}{p(t)-1}}(\R^2)} ,
\enda
\end{align*}
where we used H\"older's inequality with $p(t)\in(1,p(0)]$ defined in \eqref{p(t)def}. Integrating in time, we obtain
\begin{align}\label{wdL2}
\bega
\|\w^{\mathsf{d},\lambda}(t)\|_{L^2 }^2 &\leq \|\w^{\mathsf{d},\lambda}_0\|_{L^2 }^2 + \int_0^t \underbrace{\|(\mathbb{P}u^{\e}-u^E)(s) \|_{L^2 } \|\nabla_x \w^{\e,\lambda}_A(s) \|_{L^{2p(s)} } \|\w^{\mathsf{d},\lambda}(s)\|_{L^{\frac{2p(s)}{p(s)-1}}} }_{\eqref{wdL2}_*}  ds .
\enda
\end{align}
The first term on the right-hand side of \eqref{wdL2} converges to $0$ because of the convergence of the initial data. 

For the second term of \eqref{wdL2}, we apply Lemma \ref{L.dwesti} to the factor $\|\nabla_x \w^{\e,\lambda}_A\|_{L^{2p(t)}(\R^2)}$ and Lemma \ref{L.Puconv} to the factor $\|\mathbb{P}u^{\e}-u^E\|_{L^2(\R^2)}$. For the last factor in $\eqref{wdL2}_*$, using \eqref{wdinterp}, we obtain $\|\w^{\mathsf{d},\lambda}(s)\|_{L^{\frac{2p(s)}{p(s)-1}}} \leq 2\|\w_0\|_{L^2\cap L^\infty(\R^2)}$, since $\frac{2p(s)}{p(s)-1}>2$. 
Then, on the time interval $[0,\bar{T}]$ with $\bar{T}$ defined in \eqref{p(t)def}, we can bound the integrand as
\begin{align*}
\bega
\eqref{wdL2}_*  \les e^{2-2e^{-Ct}} \bigg(\|\mathbb{P}u^{\e}_0-u^{E}_0\|_{L^2_x}+\mathcal{C}_* \Big(\e^{\frac{1}{4}-} +\kappa^{\frac{1}{2}-}\Big) \bigg)^{e^{-Ct}} \lambda^{-4p_0+2}.
\enda
\end{align*}
Hence, we conclude that the right-hand side of \eqref{wdL2} converges to $0$ for $t\in[0,\bar{T}]$ as $\e,\lambda \to0$ by choosing $\lambda$ appropriately. Since $\bar{T}$ in \eqref{p(t)def} is independent of the initial data and of $\e,\lambda$, we may iterate the above argument on any finite time interval $[0,T]= [0,\bar{T}] \cup [\bar{T},2\bar{T}] \cup [2\bar{T},3\bar{T}] \cdots$. This completes the proof of \eqref{wAstrong}. 
\end{proof}

\begin{proof}[{\bf Proof of Theorem \ref{T.C.TL}}]
Next, we prove the convergence rate in \eqref{wrate} when the initial vorticity belongs to $(L^1 \cap L^\infty \cap \mathcal{F}^s_2)(\R^2)$. Since $\w^{\e}$ is uniformly bounded in $(L^2 \cap L^\infty)(\R^2)$, 
the velocity field $\mathbb{P}u^{\e} = {\bf K} \ast \w^{\e}$ is locally log-Lipschitz: 
\begin{align}\label{log-Lip}
\bega
|\mathbb{P}u^{\e}(t,x)-\mathbb{P}u^{\e}(t,y)| \leq C(\|\w^{\e}\|_{L^2\cap L^\infty}) (1+|x-y|)\ln\frac{1}{|x-y|}, \quad \mbox{for} \quad |x-y|<1, \quad x,y \in \R^2.
\enda
\end{align}
Taking the time derivative of the flow map in \eqref{traju}, we have 
\begin{align*}
\bega
|X^{\e}(t,x)&-X^{\e}(t,y)| = \bigg|\int_0^t\mathbb{P}u^{\e}(s,X^{\e}(s,x))-\mathbb{P}u^{\e}(s,X^{\e}(s,y))ds \bigg| \cr 
&\leq |x-y| + C(\|\w^{\e}\|_{L^2\cap L^\infty}) \int_0^t \Big(1+ |X^{\e}(s,x)-X^{\e}(s,y)|\Big) \ln\frac{1}{X^{\e}(s,x)-X^{\e}(s,y)} ds ,
\enda
\end{align*}
where we have used \eqref{log-Lip}. Hence, we obtain
\begin{align*}
\bega
|X^{\e}(t,x)-X^{\e}(t,y)| \leq C^{1+e^{-Ct}} |x-y|^{e^{-Ct}}.
\enda
\end{align*}
Writing the backward flow map as $A^{\e}(t,x):= (X^{\e}(t,\cdot))^{-1}(x)$, we also have the same estimate:
\begin{align*}
\bega
|A^{\e}(t,x)-A^{\e}(t,y)| \leq C^{1+e^{-Ct}} |x-y|^{e^{-Ct}}.
\enda
\end{align*}
Let us define a decreasing function $s(t)$ by
\begin{align*}
\bega
s(t):= s_0 \mathfrak{a}(t), \quad \mathfrak{a}(t) := e^{-Ct}, \quad s_0< 1. 
\enda
\end{align*}
We estimate 
\begin{align*}
\bega
\frac{|\w^{\e}_A(t,x)-\w^{\e}_A(t,y)|}{d(x,y)^{s(t)}} 
&\leq \frac{|\w^{\e}_0(A^{\e}(t,x))-\w^{\e}_0(A^{\e}(t,y))|}{d(A^{\e}(t,x),A^{\e}(t,y))^{s_0}} \frac{d(A^{\e}(t,x),A^{\e}(t,y))^{s_0}}{d(x,y)^{s(t)}} \cr 
&\leq \|A^{\e}(t)\|_{C^{0,\mathfrak{a}(t)}}^{s_0}  \big(\mathfrak{V}(A^{\e}(t,x)) + \mathfrak{V}(A^{\e}(t,y))\big),
\enda
\end{align*}
for some function $\mathfrak{V}\in L^2(\R^2)$. 
Then, taking the infimum over $\mathfrak{V} \in L^2(\R^2)$ gives
\begin{align}\label{wFsbdd}
\bega
[\w^{\e}(t)]_{\mathcal{F}_2^{s(t)}} \leq \|A(t)\|_{C^{0,\mathfrak{a}(t)}}^{s_0} [\w^{\e}_0]_{\mathcal{F}_2^{s_0}}.
\enda
\end{align}

We now return to the estimate of the vorticity difference in $L^2$. By H\"{o}lder's inequality, we have 
\begin{align*}
\bega
\|(\w^{\e}_A-\w^E)(t)\|_{L^2(\R^2)} \leq \|(\w^{\e}_A-\w^E)(t)\|_{H^{-1}(\R^2)}^{\frac{s'}{1+s'}}\|(\w^{\e}_A-\w^E)(t)\|_{H^{s'}(\R^2)}^{\frac{1}{1+s'}}, \quad \mbox{for any} \quad s'>0 .
\enda
\end{align*}
Then, using \eqref{wFsbdd}, we obtain a uniform bound for the $H^{s(t)-}$ norm as follows:
\begin{align*}
\bega
\|(\w^{\e}_A-\w^E)\|_{H^{s(t)-}} \leq C\|(\w^{\e}_A-\w^E)\|_{B^{s(t)-}_{2,2}} \leq C\|(\w^{\e}_A-\w^E)\|_{B^{s(t)}_{2,\infty}} \leq C\|(\w^{\e}_A-\w^E)\|_{F^{s(t)}_2} \leq C, 
\enda
\end{align*}
where we used \eqref{BesovFp}, and \eqref{Besov-p,q} in Lemma \ref{L.Besov}.
Hence, applying Lemma \ref{L.Puconv} yields
\begin{align}\label{w-w1}
\bega
\|(\w^{\e}_A-\w^E)(t)\|_{L^2(\R^2)} &\leq C \|(\w^{\e}_A-\w^E)(t)\|_{H^{-1}(\R^2)}^{\frac{s(t)}{1+s(t)}-} \leq C \|(u^{\e}_A-u^E)(t)\|_{L^2(\R^2)}^{\frac{s(t)}{1+s(t)}-}. 
\enda
\end{align}
Since $\mathbb{P}u^{\e} = u^{\e}_A+u^{\e}_B$ by Lemma \ref{L.uABprop}, we have from Lemma \ref{L.Puconv} and \eqref{uBto0} that 
\begin{align}\label{w-w2}
\bega
\|(u^{\e}_A-u^E)&(t)\|_{L^2(\R^2)} \leq \|(\mathbb{P}u^{\e}-u^E)(t)\|_{L^2(\R^2)} + \|u^{\e}_B(t)\|_{L^2(\R^2)} \cr 
&\les e^{2-2e^{-Ct}} \bigg(\|\mathbb{P}u^{\e}_0-u^{E}_0\|_{L^2_x}+\mathcal{C}_* \Big(\e^{\frac{1}{4}-} +\kappa^{\frac{1}{2}-}\Big) \bigg)^{e^{-Ct}} + \mathcal{C}_* \Big(\e^{\frac{1}{4}-} +\kappa^{\frac{1}{2}-}\Big).
\enda
\end{align}
Combining \eqref{w-w1} and \eqref{w-w2} with \eqref{wdinterp} proves the estimate \eqref{wrate}. 
\end{proof}

\subsection{$L^p$-vorticity (DiPerna--Lions--Majda solutions)}\label{sec:DLM}

\hide
 By the uniqueness of the weak limit, 완전 틀림 완전 틀림 this allows us to upgrade the convergence~\eqref{phiweak} to완전 틀림 완전 틀림 
 \todo{The following argument is faulty. \eqref{phistrong} cannot work in $C([0,T])$, You need to use the equation and Aubin-Lions.}
\begin{align}\label{phistrong}
\bega
\int_{\R^2}\phi^{\e}(t,x)\psi(x) dx \rightarrow \int_{\R^2}\phi^{\#}(t,x)\psi(x) dx \quad \mbox{in} \quad C([0,T]).
\enda
\end{align}

이렇게 할수 있을것 같으면, Aubin-Lions 는 사람들이 왜 쓰겠어요. 시간 오실레이션을 잡아줘야겠죠

\unhide

The notion of renormalized solutions for transport equations was introduced by DiPerna and Lions in their seminal work \cite{DiLi}, where they established well-posedness for transport equations associated with Sobolev vector fields and developed the renormalization framework that has become fundamental in fluid dynamics. 
Building on this theory, the renormalized formulation has been further investigated in the context of the vorticity equation for the incompressible Euler system in \cite{Filho2006,B-Crippa,Crippa2021}.

In two dimensions, the vorticity equation enjoys the special structure of a continuity (transport) equation, which allows for various convergence results under low regularity assumptions. In \cite{Filho2006}, it was shown that every weak solution $\w \in L^{\infty}(0,T;L^p(\R^2))$  with $p \geq 2$ to the incompressible Euler equations is a renormalized solution.
In the context of the inviscid limit, it was proved in \cite{Crippa2021} that the viscous solutions $(u^{\nu},\w^{\nu})$ of the incompressible Navier–Stokes equations converge to a limit that satisfies the Euler equations in the renormalized sense for $1 \leq p \leq 2$, together with strong convergence of the vorticity $\w^{\nu}$.
Moreover, energy conservation in the inviscid limit was established in \cite{Econ}.

In this section, we establish the strong convergence of $u^{\e}$ and $\w^{\e}$ constructed in Theorem~\ref{T.2D.global}. 
The limit is shown to be a global renormalized solution in the sense of DiPerna--Lions--Majda, as in Definition~\ref{D.soluw}, provided that the initial vorticity sequence is uniformly bounded in $(L^p\cap L^1)(\R^2)$ for $1 \le p < \infty$. 
The macroscopic fields in the Boltzmann equation interact with the microscopic component. The reason we are able to work at the level of the vorticity equation is, as emphasized earlier, that the macro–micro cancellation yields strong dissipation. This enhanced dissipation ensures that, except at the top order, the Burnett functionals vanish at the rate $\kappa^{1/2}$. 
Nevertheless, an additional difficulty is that  the $L^p$-norm of the microscopic term $\nabla_x^{\perp}\cdot(\nabla_x \cdot \mathbf{r}^{\e})$ in the vorticity equation \eqref{weqnnew} cannot be controlled by the energy and dissipation for $1 \leq p < 2$.
To overcome this difficulty, we decompose the vorticity into the initial data part and the forcing part, as introduced in Definition~\eqref{D.wAwB}.

For the incompressible vector fields $u^{\e}_A={\bf K}\ast \w^{\e}_A$ and $u^{\e}_B={\bf K}\ast \w^{\e}_B$, with ${\bf K}(x):= \frac{1}{2\pi}\frac{x^{\perp}}{|x|^2}$, we prove the strong convergence of $\w^{\e}_A$ in $L^p(\R^2)$, and that the limit $(u^{E},\w^{E})$ of $(u^{\e}_A,\w^{\e}_A)$ satisfies the incompressible Euler equations in the renormalized sense. Moreover, the forcing part $u^{\e}_B$ converges to $0$ in $L^\infty(0,T;L^\infty(\R^2))$.


The main theorem of this section shows that even without passing through the Navier–Stokes equations, the limit of the solutions $(u^{\e},\w^{\e})$ obtained from the Boltzmann equation satisfies the Euler equations in the renormalized sense.


\hide
Since the $L^p$ norm of the microscopic part of the equation \eqref{weqnP} $\nabla_x^{\perp}\cdot \p_{x_j} \mathbf{r}_{ij}^{\e}$ is not able to control for $1\leq p<2$, recall that in Definition \ref{D.wAwB} , we decomposed the solution of the vorticity equation \eqref{weqnP} $\w^{\e}$ by the initial data part and forcing term part in Definition \ref{D.wAwB}:
\begin{align*}
\bega
\w^{\e}(t,x)&= \w^{\e}_A(t,x) + \w^{\e}_B(t,x),
\enda
\end{align*}
where
\begin{align*}
\bega
\w^{\e}_A(t,x)&:= \w^{\e}_0((X^{\e})^{-1}(t,x)), \qquad \w^{\e}_B(t,x):=\int_0^t \bar{\varPi}_{\w}^{\e}(s,X^{\e}(s,\mathbf{a})) ds. 
\enda
\end{align*}
(Note that $\nabla_x^{\perp}\cdot\mathbb{P}u^{\e} \neq \w^{\e}_A$.)
Since $\mathbb{P}u^{\e}$ is divergence-free, the $L^p_x$–norm of $\w^{\e}_A$ is conserved:
\begin{align*}
&\frac{d}{dt}\|\w^{\e}_A\|_{L^p_x}^p =0, \quad 1\leq p<\infty.
\end{align*}
For the vorticity $\w^{\e}_A$, we define corresponding vector field $u^{\e}_A$ satisfying 
\begin{align*}
\nabla_x\cdot u^{\e}_A = 0, \qquad \nabla_x^{\perp}\cdot u^{\e}_A = \w^{\e}_A.
\end{align*}
\unhide


\hide
\begin{lemma}\label{L.Lpunif}
We take the solution of the Boltzmann equation $F^{\e}$ constructed in Theorem \ref{T.2D.global}. In addition, by further assuming $\|\w^{\e}_0\|_{L^p_x} \leq C$ in Theorem \ref{T.2D.unif} (iv) for $2\leq p<\infty$ and (iv) for $1\leq p<2$, we have the following estimates: 
\begin{itemize}
\item If $\|\w^{\e}_0\|_{L^p_x}$ is uniformly bounded, then $\|\w^{\e}(t)\|_{L^p_x}$ is uniformly bounded for each $2 \leq p \leq \infty$:
\begin{align}
\|\w^{\e}(t)\|_{L^p_x} \leq C\big(\|\w^{\e}_0\|_{L^p_x}+1\big) , \quad \mbox{for each} \quad 2 \leq p\leq \infty.
\end{align}
\begin{align}
\sup_{t\in[0,T]}\|\mathbb{P}u^{\e}\|_{W^{1,p}_{loc}} \leq C
\end{align}
\item If $\|\w^{\e}_0\|_{L^p_x}$ is uniformly bounded for $1 \leq p < 2$ then we have
\begin{align}
\bega
&\big\|\w^{\e}_0((X^{\e})^{-1}(t,x))\big\|_{L^p_x} \leq C\|\w^{\e}_0\|_{L^p_x}, \qquad 1\leq p < 2, \cr 
&\bigg\|\int_0^t \bar{\varPi}_{\w}^{\e}(s,X^{\e}(s,\mathbf{a})) ds\bigg\|_{L^2_x} \to 0 , \quad \mbox{as} \quad \e\to0, \cr 
&\sup_{t\in[0,T]}\|u^{\e}_A\|_{W^{1,p}_{loc}} \leq C
\enda
\end{align}
\end{itemize}
\end{lemma}
\begin{lemma} For $\w^{\e}$ satisfying (iv) in Theorem \ref{T.2D.unif} satisfies
\begin{align}
\bega
\sup_{0\leq t \leq T}\|\w^{\e}(t)\|_{L^p_x} \leq C, \quad \sup_{0\leq t\leq T}\|\mathbb{P}u^{\e}(t)\|_{W^{1,p}_{loc}} \leq C
\enda
\end{align}
\end{lemma}
\unhide

 \begin{definition}[DiPerna--Lions--Majda \cite{DiLi,DiMa}]\label{D.soluw}
A pair $(u^E,\omega^E)$ is a \emph{renormalized solution} of \eqref{weqnE2D} on $[0,T]$ with initial vorticity
$\omega^E|_{t=0}=\omega^E_0$ if:
\begin{itemize}
\item $u^E={\bf K}\ast \omega^E$ and $\nabla_x\cdot u^E = 0$ in the sense of distributions, where
${\bf K}(x):= \frac{1}{2\pi}\frac{x^{\perp}}{|x|^2}$;
\item for every $\psi \in C_c^{\infty}([0,T)\times \mathbb{R}^2)$,
\[
\int_0^T \int_{\mathbb{R}^2} \boldsymbol{\beta}(\omega^E)\,(\partial_t\psi + u^E \cdot \nabla_x \psi)\,dxdt
+ \int_{\mathbb{R}^2}\boldsymbol{\beta}(\omega^E_0)\,\psi(0,x)\,dx = 0,
\]
for all $\boldsymbol{\beta} \in C^1(\mathbb{R}) \cap L^\infty(\mathbb{R})$ that vanish in a neighborhood of $0$.
\end{itemize}
\end{definition}

\begin{definition}\cite{DiLi, B-Crippa} \label{D.sol}
We consider the continuity equation with an incompressible vector field $\mathfrak{b}$:
\begin{align}\label{rhob}
\p_t \varrho + \mathfrak{b} \cdot \nabla_x \varrho = 0, \qquad  \varrho(0,x) = \varrho_0(x), \qquad \div(\mathfrak{b})=0.
\end{align}
\begin{itemize}

\item 
A function $\varrho$ is called a \emph{renormalized solution} of \eqref{rhob} if, for every test function $\psi \in C_c^{\infty}([0,T)\times \R^2)$, the following identity holds:
\begin{align*}
\int_0^T \int_{\R^2} \boldsymbol{\beta}(\varrho)(\p_t\psi + \mathfrak{b} \cdot \nabla_x \psi) dxdt + \int_{\R^2} \boldsymbol{\beta}(\varrho_0)\psi(0,x) dx = 0,
\end{align*}
for all $\boldsymbol{\beta} \in C^1(\R) \cap L^\infty(\R)$ that vanish in a neighborhood of zero.

\item 
A function $\varrho$ is called a \emph{distributional solution} of \eqref{rhob} if $\mathfrak{b} \in L^1(0,T;L^p_{loc}(\R^2))$, 
$\varrho_0 \in L^q_{loc}(\R^2)$, and $\frac{1}{p} + \frac{1}{q} \leq 1$, 
and if $\varrho \in L^{\infty}(0,T;L^q_{loc}(\R^2))$ satisfies 
\begin{align*}
\bega
\int_0^T \int_{\R^2} \varrho (\p_t \psi + \mathfrak{b} \cdot \nabla_x \psi) \, dxdt  + \int_{\R^2} \varrho_0 \psi|_{t=0} dx =0,
\enda
\end{align*}
for any $\psi \in C_c^{\infty}([0,T)\times \R^2)$. Here, $u \in L^p_{\mathrm{loc}}(\R^2)$ means that $u \in L^p(K)$ for every compact set $K \subset \R^2$.

\end{itemize} 
\end{definition}

\begin{definition}
Let $X(t,x)$ denote the \emph{regular Lagrangian flow} satisfying 
$\frac{d}{ds}X(s,x) = \mathfrak{b}(s,X(s,x))$ for $s \in [0,T]$ 
in the renormalized sense, with $X(0,x) = x$ 
(its precise definition is given in Definition~5.2 of \cite{B-Crippa}). 
Then a function $\varrho(t,x) = \varrho_0(X(0,x))$ is called a 
\emph{Lagrangian solution} of \eqref{rhob}. 
\end{definition}

\hide
\begin{definition}
(1) We define some function spaces\\
(i) $L^0(\R^2)$: the space of real-valued measurable functions on $\R^2$ with respect to Lebesgue measure $\mathfrak{L}^2$. \\
(ii) $\mathcal{B}(X,Y)$ : the space of bounded functions between the sets $X$ and $Y$.\\
(iii) $\log L(\R^2)$ space of measurable functions $\varrho~:~\R^2\to\R$ such that $\int_{\R^2} \log(1+|\varrho(x)|)dx$ is finite \\
(2) We define the regular Lagrangian flow $X(t,x)$
\begin{align}
\bega
X \in C(0,T;L^0_{loc}(\R^2) \cap \mathcal{B} (0,T;\log L_{loc}(\R^2))
\enda
\end{align}
if it satisfies (i)
\begin{align}
\bega
\int_0^t\int_{\R^2} \p_t \beta(X(t,x)) \psi(t,x) dxdt = \int_0^t\int_{\R^2} \beta' (X(t,x)) \mathfrak{b}(t,X(t,x))) \psi(t,x) dxdt , 
\enda
\end{align}
for any $\beta \in C^1(\R\times \R^2)$ such that $|\beta(z)| \leq C(1+\log(1+|z|))$ and $|\beta'(z)|\leq \frac{C}{1+|z|}$ and $\psi\in \mathbb{D}'((0,T)\times\R^2)$.
(ii) $X(t,x)=x$ for $\mathfrak{L}^2$- a.e. $x\in\R^2$.
(iii) There exists constant $L\geq 0$ such that for $t\in[0,T]$, 
\begin{align*}
\bega
\mathfrak{L}^2 \Big(\{x\in\R^2~:~ X(t,x) \in B \}\Big) \leq L \mathfrak{L}^2(B), \quad \mbox{for every Borel set} \quad B\subset \R^2,
\enda
\end{align*}
where $\mathfrak{L}^2$ is 2-dimensional Lebesgue measure.\\
(3) We call the following solution is Lagrangian solution 
\begin{align}\label{Lagsol}
\bega
\varrho(t,x) = \varrho_0(X^{-1}(t,x)).
\enda
\end{align}
\end{definition}
\unhide


We present two theorems from \cite{DiLi,B-Crippa}. 
For simplicity, we assume that the velocity field is divergence-free, that is, $\div(\mathfrak{b}) = 0$.
We also impose the general condition required in all three theorems:
\begin{align}\label{bassume}
\bega
\frac{|\mathfrak{b}(t,x)|}{1 + |x|} 
\in L^1(0,T;L^1(\R^2)) + L^1(0,T;L^\infty(\R^2)).
\enda
\end{align}

\begin{theorem}[Theorem II.3 of \cite{DiLi}]\label{T.Di-Li}
Let $\mathfrak{b}(t,x) \in L^1(0,T;W^{1,q}_{loc}(\R^2))$ satisfy \eqref{bassume}. Then the following statements hold:
\begin{enumerate}
\item If $\varrho_0 \in L^p(\R^2)$ with $1 \leq p < \infty$, then there exists a unique renormalized solution $\varrho \in C(0,T; L^p(\R^2))$ of \eqref{rhob}.
\item Every distributional solution $\varrho \in L^\infty(0,T;L^p(\R^2))$ of \eqref{rhob}, with $\frac{1}{p}+\frac{1}{q}=1$, is also a renormalized solution.
\end{enumerate}
\end{theorem}

\begin{theorem}[\cite{B-Crippa}]\label{T.B-Crippa}
(1) Let $\mathfrak{b}$ satisfy 
\begin{align}\label{bassume2}
\bega
&\mathfrak{b}(t,x) = {\bf K} \ast \mathfrak{g}(t,x), 
\quad \text{with} \quad \mathfrak{g} \in L^1((0,T)\times \R^2), \cr
&\mathfrak{b}(t,x) \in L^p_{loc}([0,T]\times \R^2) 
\quad \text{for some} \quad p>1.
\enda
\end{align}
Then there exists a unique Lagrangian solution of \eqref{rhob}, which is also a renormalized solution.  \\
(2) Let $\varrho$ and $\varrho_n$ be the Lagrangian solutions corresponding to the divergence-free vector fields $\mathfrak{b}$ and $\mathfrak{b}_n$, respectively, both satisfying \eqref{bassume2}. 
Assume that $\mathfrak{b}_n \to \mathfrak{b}$ strongly in $L^1_{loc}([0,T]\times \R^2)$, and that $\mathfrak{b}_{n1}$ and $\mathfrak{b}_{n2}$ are equibounded in $L^1(0,T;L^1(\R^2))$ and $L^1(0,T;L^\infty(\R^2))$, respectively, where $\mathfrak{b}_n = \mathfrak{b}_{n1} + \mathfrak{b}_{n2}$. 
Then, the strong convergence of the initial data $\varrho_{n0} \to \varrho_0$ in $L^1(\R^2)$ implies the strong convergence of the Lagrangian solutions $\varrho_n \to \varrho$ in $C(0,T;L^1(\R^2))$.
\end{theorem}
\begin{proof}
(1) For the existence and uniqueness of the regular Lagrangian flow, see Theorems~6.1 and~6.4 in \cite{B-Crippa}. 
For the statement that $\varrho$ is a renormalized solution, see Proposition~7.2 therein. \\
(2) For the stability result, see Proposition~7.6 in \cite{B-Crippa}.
\end{proof}

We define the weak $L^p$-norm (the quasi-norm of the Lorentz space $L^{p, \infty}$) as
\begin{align*}
\|u\|_{L^{p,\infty}} := \sup_{\lambda \geq0} \lambda |m\{ x: |u(x)|>\lambda \}|^{\frac{1}{p}}.
\end{align*}

\hide
\begin{theorem}[Theorem 4.1 of \cite{Crippa2017}]\label{T.Crippa}
Let $\mathfrak{b}$ satisfy 
\begin{align}\label{bassume.dist}
\bega
&\mathfrak{b}(t,x) = {\bf K} \ast \mathfrak{g}(t,x), 
\quad \text{with} \quad \mathfrak{g} \in L^1((0,T)\times \R^2), \cr 
&\mathfrak{b}(t,x) \in L^{p,\infty}((0,T)\times \R^2),
\quad \text{for some} \quad p>1.
\enda
\end{align}
If $\varrho_0 \in L^\infty \cap L^1(\R^2)$, then there exists a unique distributional solution $\varrho \in L^\infty(0,T;L^\infty\cap L^1(\R^2))$ of \eqref{rhob}.
\end{theorem}
\unhide

\hide
\begin{definition}\label{D.soluw}
A pair $(u,\w)$ is called a \emph{renormalized solution} of the vorticity formulation of the Euler equations corresponding to the initial data $\w|_{t=0}=\w_0$ if the following conditions hold:
\begin{itemize}
\item $u={\bf K}\ast \w$ and $\div(u) = 0$ in the sense of distributions for ${\bf K}(x):= \frac{1}{2\pi}\frac{x^{\perp}}{|x|^2}$.
\item For every test function $\psi \in C_c^{\infty}([0,T)\times \R^2)$, the following identity holds:
\begin{align*}
\int_0^T \int_{\R^2} \boldsymbol{\beta}(\w)(\p_t\psi + u \cdot \nabla_x \psi) dxdt + \int_{\R^2} \boldsymbol{\beta}(\w_0)\psi(0,x) dx = 0,
\end{align*}
for all $\boldsymbol{\beta} \in C^1(\R) \cap L^\infty(\R)$ that vanish in a neighborhood of zero.
\end{itemize}
\end{definition}
\unhide

Recall that in Definition~\ref{D.wAwB}, we defined $\w^{\e}_A$ and $\w^{\e}_B$ as the solutions of the following equations:
\begin{align*}
&\p_t\w^{\e}_A + \mathbb{P}u^{\e}\cdot \nabla_x \w^{\e}_A = 0 , \qquad \w^{\e}_A(t,x)|_{t=0} = \w^{\e}_0(x), \cr 
&\p_t\w^{\e}_B + \mathbb{P}u^{\e}\cdot \nabla_x \w^{\e}_B = \bar{\varPi}_{\w}^{\e} , \qquad \w^{\e}_B(t,x)|_{t=0} = 0 . 
\end{align*}
We also defined the associated divergence-free velocity fields via the Biot--Savart law:
\begin{align*}
\bega
u^{\e}_A(t,x) := {\bf K} \ast \w^{\e}_A(t,x), \qquad 
u^{\e}_B(t,x) := {\bf K} \ast \w^{\e}_B(t,x), \qquad {\bf K}(x):= \frac{1}{2\pi}\frac{x^{\perp}}{|x|^2}.
\enda
\end{align*}
Recall the radial energy decomposion in Definition \ref{D.Ra-E} and $\bar{u}\cdot\nabla_x \bar{\w}=0$. 
Recall also that we decomposed $F^{\e}=M^{\e}+\AC{\P}F^{\e}$ as in \eqref{F-ME}.
We now present the main theorem of this section.

\begin{theorem}\label{T.C.Lp}
Let $\Omega=\R^2$. For $p \in [1, \infty)$, consider a family of initial data $\{F_0^{\e}\}_{\e>0}$ satisfies \eqref{L2unif}, the $4$-ABC~\eqref{ABC1}, and 
\begin{align*}
\sup_{\e>0} 
\|\w_0^{\e}\|_{L^p \cap L^1 \cap H^{-1}_{loc}(\R^2)}
<+\infty.
\end{align*}
For $p=1$, we further assume that $\{\w^{\e}_0\}_{\e>0}$ is supported in a fixed ball. 

We assume that, as $\e \to 0$, \eqref{conv_s_initial} and 
\begin{align*}
\bega
\w^{\e}_0 \quad &\to \quad \w_0  \qquad  \text{strongly in $L^p\cap L^1(\R^2)$}.
\enda
\end{align*}


Then the following statements hold for a family of Boltzmann solutions $\{F^{\e}\}_{\e>0}$ to~\eqref{BE} with $\kappa=\e^q$ for some $0<q<2$, on
the time interval $t\in[0,T_\e]$ with $T_\e\to\infty$ as $\e\to0$, which is constructed in Theorem~\ref{T.2D.global}.
\begin{enumerate}
\item[(1)] The family of macroscopic fields of the Boltzmann solution $\big\{(\mathbb{P} u^\e, \mathbb{P}^{\perp}u^\e, \frac{3}{2} \theta^\e - \rho^\e, \rho^\e+ \theta^\e) \big\}_{\e>0}$ converges to $(u^E, 0, \frac{3}{2} \theta^E - \rho^E, 0)$ along a subsequence in the following different topologies up to subsequence. The pair $(u^{E},\w^{E})$ satisfies the incompressible Euler equations in the renormalized sense as in Definition \ref{D.soluw} and
\begin{align*}
\w^{E}  \in L^\infty(0,T;L^p(\R^2)), \qquad \begin{cases} u^{E} \in L^\infty(0,T;W^{1,p}_{loc}(\R^2)), &\mbox{for} \quad 1<p<\infty, \\ 
u^{E} \in L^\infty(0,T;L^2_{loc}\cap L^{2,\infty}(\R^2)) , &\mbox{for} \quad p=1.
\end{cases}
\end{align*}

\item[(2)] The velocity variables converge as follows:
\begin{align}
&u^{\e}_A   \rightarrow   u^{E} \quad \mbox{strongly in} \quad  \begin{cases}L^2(0,T;L^2_{loc}(\R^2)), \quad &\mbox{for} \quad 1<p<\infty, \\ 
L^r(0,T;L^r_{loc}(\R^2)), \quad \mbox{for} \quad 1\leq r<2 , \quad &\mbox{for} \quad p=1,
\end{cases} \label{mconv;0} \\ 
&\w^{\e}_A   \rightarrow   \w^{E} \quad \mbox{strongly in} \quad C(0,T;L^p(\R^2)), \quad \mbox{for} \quad 1\leq p<\infty, \label{wconv;0}  \\
&u^{\e}_B, \ \w^{\e}_B   \rightarrow   0 \ \ \text{strongly in} \quad L^\infty(0,T; L^p(\R^2)), \quad \mbox{for} \quad  2\leq p<\infty \quad \mbox{with rate $\e^{\frac{1}{4}-} +\kappa^{\frac{1}{2}-}$.}  \notag 
\end{align}

\item[(3)] The entropic fluctuation $\mathfrak{s}^{\e}:=3/2\ta^{\e}-\rho^{\e}$ converges as follows:
\begin{align}\label{TC.Lp.ent}
\frac{3}{2} \theta^\e - \rho^\e  \to \frac{3}{2} \theta^E - \rho^E  \ \ &\text{in  $L^\infty (0, T; L^2(\R^2))$ }. 
\end{align}

\item[(4)] The irrotational part $\mathbb{P}^{\perp}u^{\e}$ and pressure fluctuation $\rho^{\e}+\ta^{\e}$ vanish as \eqref{TC.Yudo.irrorhota}.
\hide\begin{align}
\mathbb{P}^{\perp}u^{\e}, \ \rho^{\e}+\ta^{\e}    \to  0  \ \ &\text{in $ L^r (0,T; \dot{B}_{p,1}^{s+2(\frac{1}{p}-\frac{1}{2})+\frac{1}{r}}(\R^2))$ with rate $\e^{\frac{1}{4}-}$} ,\label{TC.Yudo.irrorhota}
\end{align}
for $2\leq p\leq\infty$, $\frac{1}{r}\leq 2(\frac{1}{2}-\frac{1}{p})$, and any $s\in[0,3)$. \unhide

\item[(5)] The microscopic part converges to zero, as quantified in \eqref{TC.Yudo.G}.
\end{enumerate}
\end{theorem}

\hide
\begin{remark}
The compressible components converge as follows:
\begin{align}\label{TC.Lp.irrorhota}
\mathbb{P}^{\perp}u^{\e}, \ \rho^{\e}+\ta^{\e}    \to  0  \ \ &\text{in $ L^r (0,T; \dot{B}_{p,1}^{s+2(\frac{1}{p}-\frac{1}{2})+\frac{1}{r}}(\R^2))$ with rate}.
1\end{align}
This holds for $2\leq p\leq\infty$, $\frac{1}{r}\leq 2(\frac{1}{2}-\frac{1}{p})$, and any $s\in[0,\mathrm{N}-1)$.
\end{remark}
\unhide

\hide
\begin{remark}
We summarize the quantitative convergence rates in the Big-O sense as follows:
\begin{align*}
\e^{\frac{1}{4}-} +\kappa^{\frac{1}{2}-}
\ \ &\text{in \eqref{wBto0Lp}},\\
\e^{\frac{1}{4}-}
\ \ &\text{in \eqref{TC.Lp.irrorhota}}.
\end{align*}
\end{remark}
\unhide

\begin{corollary}
Assume
\begin{align}
&\w_0\in L^p\cap L^1 \cap H^{-1}_{loc}(\R^2), \notag 
\\
& (\rho_0, u_0-\bar{u},\ta_0) \in L^2(\R^2) ,
\quad \frac{1}{\e}\frac{\AC{\P}F_0}{\sqrt{\tilde{\mu}}}\in L^2_x(\R^2 ; L^2_v(\R^3)), \quad 
\bigg\|\frac{1}{\e}\frac{(F_0-\mu)}{\sqrt{\tilde{\mu}}}\bigg\|_{L^\infty_{x,v}}\les \frac{1}{\e^{1-}},\label{initialL2}
\end{align}
where $\tilde{\mu}=M_{[1,0,1-c_0]}$ for some $0<c_0\ll1$.
Then the mollified sequence $\{F^{\e}_0\}_{\e>0}$ constructed from Lemma \ref{L.molli} satisfies the assumption of Theorem \ref{T.C.Lp}, and hence all the results of Theorem \ref{T.C.Lp} hold. 
\end{corollary}

\begin{remark}
Once the limit velocity field $u^{E}$ is obtained, the continuity equation with velocity field $\mathbf{b}=u^{E}$ admits a unique renormalized (or Lagrangian) solution by Theorem \ref{T.Di-Li} for $1<p<\infty$ and Theorem \ref{T.B-Crippa} for $p=1$.
\end{remark}

\hide
\begin{remark}
To prove Theorem~\ref{T.C.Lp}, it suffices to assume that 
$\sup_{\e>0}\|\w^{\e}_0\|_{L^p\cap L^1(\R^2)}<+\infty$ for $1<p<\infty$, 
or that $\sup_{\e>0}\|\w^{\e}_0\|_{L^1 \cap H^{-1}_{loc}(\R^2)} <+\infty$ with compact support when $p=1$. 
Under these assumptions, we can show that $\sup_{0\leq t\leq T}\|u^{\e}_A(t)\|_{L^2_{\mathrm{loc}}(\R^2)}$ is uniformly bounded. 
For the sake of simplicity, we assume $(\mathbb{P}u^{\e}_0-\bar{u}) \in L^2(\R^2)$.
\end{remark}
\unhide

\begin{proposition}\label{P.Econs}
Let $\Omega = \R^2$. 
Suppose the initial data $u_0^{\e}$ satisfies the finite velocity energy condition \eqref{caseEC}. 
Assume that all hypotheses of Theorem~\ref{T.2D.global} hold, and consider the corresponding solution constructed therein. 
In addition, suppose that $\mathbb{P}u^{\e}_0 \in L^2(\R^2)$ and 
\begin{align}\label{Pu0assume}
\mathbb{P}u^{\e}_0 \rightarrow u_0 
\quad \text{strongly in} \quad L^2(\R^2).
\end{align}
Then, the following statements hold:
\begin{enumerate}
\item 
The divergence-free component $\mathbb{P}u^{\e}$ is conserved in the limit sense:
\begin{align}\label{E-conservlim}
\bega
\lim_{\e \to 0} \|\mathbb{P}u^{\e}(t)\|_{L^2_x} 
= \|u_0\|_{L^2_x}.
\enda
\end{align}

\item 
For $p>1$, $u^{\e}_A$ converges to $u^{E}$ on the whole space:
\begin{align}\label{ustrong;0}
u^{\e}_A \;\rightarrow\; u^{E} 
\quad \text{strongly in} \quad C(0,T;L^2(\R^2)),
\end{align}
and the limiting velocity $u^{E}$ conserves the energy:
\begin{align}\label{E-conserv0}
\bega
\|u^{E}(t)\|_{L^2_x} = \|u_0\|_{L^2_x}.
\enda
\end{align}
\end{enumerate}
\end{proposition}

\begin{remark}
For the inviscid (zero-viscosity) limit solutions, the energy conservation \eqref{E-conserv0} was established in \cite{Econ} by showing that 
$\nu \int_0^t \|\nabla_x u^{\nu}(s)\|_{L^2_x}^2 ds \to 0$ as $\nu \to 0$. 
In contrast, for the hydrodynamic limit solution constructed in Theorem~\ref{T.2D.global}, the energy conservation \eqref{E-conservlim} follows from the vanishing microscopic dissipation $\frac{1}{\e^2}\|\nabla_x\cdot\mathbf{r}^{\e}\|_{L^2_x} \les \kappa^{\frac{1}{2}-}$ and acoustic dispersion $\e^{\frac{d-1}{4}}$.
(See the proof of \eqref{U.PuL2}.)
\end{remark}



\hide
\begin{lemma} For $\beta \in C^1(\R) \cap L^\infty(\R)$ vanishing in a neighbourhood of zero, If $\rho$ satisfies the first equation, then we have the second, third equation:
\Be\bega \label{beta;E}
&\p_t\rho + \mathfrak{b} \cdot\nabla_x \rho =0 \cr 
&\p_t \beta(\rho) + \mathfrak{b} \cdot\nabla_x \beta(\rho) =0 \cr 
&\p_t |\rho|^p + \mathfrak{b} \cdot\nabla_x |\rho|^p =0 , \quad \mbox{for} \quad p > 1.
\enda
\Ee
If $\beta \in C^2(\R) \cap L^\infty(\R)$, then for NS: 
\Be\bega \label{beta;NS}
&\p_t\rho^{\nu} + b^{\nu}\cdot\nabla_x \rho^{\nu} -\nu \Delta_x \rho^{\nu} = g^{\e} \cr 
&\p_t\beta(\rho^{\nu}) + b^{\nu}\cdot\nabla_x \beta(\rho^{\nu}) -\nu \Delta_x \beta(\rho^{\nu}) = - \nu \beta''(\rho^{\nu})|\nabla_x \rho|^2 + \beta'(\rho^{\nu})g^{\e} \cr 
&\p_t |\rho^{\nu}|^p + b^{\nu}\cdot\nabla_x |\rho^{\nu}|^p -\nu \Delta_x |\rho^{\nu}|^p = - \nu p(p-1)|\rho^{\nu}|^{p-2}|\nabla_x \rho^{\nu}|^2 + (p\rho^{\nu}|\rho^{\nu}|^{p-2})g^{\e}, \quad \mbox{for} \quad p\geq 2.
\enda
\Ee
$\bullet$ $L^p$ conservation of Euler eqn
\begin{align}
&\|\rho(t)\|_{L^p_x}=\|\rho_0\|_{L^p_x}, \quad \mbox{for} \quad p\geq1 .
\end{align}
$\bullet$ $L^2$ estimate of N-S eqn (dissipation)
\begin{align}\label{NS-L2}
\|\rho^{\nu}(t)\|_{L^2}^2 + \nu\int_0^t \|\nabla_x\rho^{\nu}(s)\|_{L^2}^2ds &\leq \|\rho^{\nu}(0)\|_{L^2}^2 + \int_0^t \|\nabla_x\cdot b^{\nu}(s)\|_{L^\infty_x} \|\rho^{\nu}(s)\|_{L^2}^2ds \cr 
&+2 \int_0^t \|g^{\e}_{\w}(s)\|_{L^2}\|\rho^{\nu}(s)\|_{L^2} ds.
\end{align}
$\bullet$ $L^p$ estimate of N-S eqn for $p\geq1$. 
\begin{align}\label{NS-Lp}
&\|\rho^{\nu}(t)\|_{L^p}^p \leq \|\rho^{\nu}(0)\|_{L^p}^p + \int_0^t \|\nabla_x\cdot b^{\nu}(s)\|_{L^\infty_x} \|\rho^{\nu}(s)\|_{L^p}^pds +p \int_0^t \|g^{\e}(s)\|_{L^p}\|\rho^{\nu}(s)\|_{L^p}^{p-1} ds.
\end{align}
\end{lemma}
\begin{proof}
\textcolor{blue}{Derivation of the third equation of \eqref{beta;NS}: 
We multiply $p\rho^{\nu}|\rho^{\nu}|^{p-2}$ to the first equation of \eqref{beta;NS} for $p\geq2$:
\begin{align*}
&\p_t |\rho^{\nu}|^p +b^{\nu}\cdot\nabla_x |\rho^{\nu}|^p -\nu  (\Delta_x \rho^{\nu}) (p\rho^{\nu}|\rho^{\nu}|^{p-2}) = g^{\e}_{\w}p\rho^{\nu}|\rho^{\nu}|^{p-2}.
\end{align*}
by using the following equality, we obtain the third equation of \eqref{beta;NS}.
\begin{align*}
\Delta_x |\rho|^p &= \nabla_x(p|\rho|^{-1}\nabla_x\rho \times sgn(\rho)) = \nabla_x (p |\rho|^{p-2}\nabla_x\rho \times \rho) \cr 
&= p(p-2)|\rho|^{p-2}\nabla_x\rho \times sgn(\rho)\nabla_x\rho \times \rho + p|\rho|^{p-2}\Delta_x\rho \times \rho + p|\rho|^{p-2}|\nabla_x\rho|^2  \cr 
&= p\rho|\rho|^{p-2}\Delta_x\rho +  p(p-1)|\rho|^{p-2}|\nabla_x \rho|^2 
\end{align*}}
\textcolor{blue}{ Proof of \eqref{NS-Lp}: 
Taking $\int_{\Omega}dx$ to the third equation of \eqref{beta;NS} gives 
\begin{align}\label{dw-Lp}
&\frac{d}{dt}\int_{\Omega}|\rho^{\nu}|^pdx -\int_{\Omega}(\nabla_x\cdot b^{\nu}) |\rho^{\nu}|^p dx  +\nu p(p-1)\int_{\Omega} |\nabla_x\rho^{\nu}|^2 |\rho^{\nu}|^{p-2}dx = p\int_{\Omega}g^{\e}_{\w}\rho^{\nu}|\rho^{\nu}|^{p-2} dx.
\end{align}
\begin{align*}
&\frac{d}{dt}\int_{\Omega}|\rho^{\nu}|^pdx +\nu p(p-1)\int_{\Omega} |\nabla_x\rho^{\nu}|^2 |\rho^{\nu}|^{p-2}dx \leq \int_{\Omega}(\nabla_x\cdot b^{\nu}) |\rho^{\nu}|^p dx +p\|g^{\e}_{\w}\|_{L^\infty_x}\|\rho^{\nu}\|_{L^p}^p.
\end{align*}
where we used $\int_{\Omega}|\rho^{\nu}|^{p-1} dx\leq \|\rho^{\nu}\|_{L^p}^p$ by the H\"{o}lder inequality.
to the last integral of \eqref{dw-Lp} gives \eqref{NS-Lp}.}
\end{proof}
\unhide

In Theorem~\ref{T.2D.global}, the solution satisfying the Admissible Blow-up 
Condition~\eqref{ABC1} also satisfies \eqref{incomp.2D} and 
\eqref{Gto0}. These imply the convergence of the microscopic part 
and the compressible part, respectively.
Since the equation for the thermodynamic variable 
$\frac{3}{2}\ta^{\e}-\rho^{\e}$ in \eqref{rtaeqnP} has a 
structure similar to that of the vorticity equation 
$\w^{\e}$~\eqref{weqnP}, it suffices to focus on the convergence 
of the vorticity.

We prove Theorem~\ref{T.C.Lp} in the following three subsections.


\subsubsection{{\bf Proof of Theorem \ref{T.C.Lp} for $1<p<\infty$}}

\hide
Statement: For $1<p<\infty$, if $\w^{\e}_0 \in L^p_c$, then $u^{\e}_A \in L^2_{loc}(\R^2)$. (2) For $p=1$, if $\w^{\e}_0 \in L^1\cap H^{-1}_{loc}(\R^2)$, then we have $u^{\e}_A\in L^2_{loc}$. (Note $\w^{\e}_A(0)=\w^{\e}_0$)
For $1<p<2$ By the Hardy-Littlewood-Sobolev inequality, we have 
\begin{align*}
\bega
\|u^{\e}_A\|_{L^q_x}\leq \|\w^{\e}_A\|_{L^p_x}, \quad q=\frac{2p}{2-p}>2. 
\enda
\end{align*}
Hence, 
\begin{align*}
\bega
\|u^{\e}_A\|_{L^2_{loc}}\leq C\|u^{\e}_A\|_{L^q_x}. 
\enda
\end{align*}
For $p=2$ Using the stream function form $u^{\e}_A = \nabla_x^{\perp}\psi^{\e}_A$, $\w^{\e} = \Delta_x \psi^{\e}_A$, we have 
\begin{align*}
\bega
\|\Delta_x\psi^{\e}_A\|_{L^2_x}= \|\w^{\e}_A\|_{L^2_x}, \qquad \|\psi^{\e}_A\|_{H^2(B_R)} \leq C(\|\w^{\e}_A\|_{L^2(B_R)} + \|\psi^{\e}_A\|_{L^2(B_R)} )
\enda
\end{align*}
For $p>2$
\begin{align*}
\bega
\|u^{\e}_A\|_{L^2_{loc}}\leq C\|u^{\e}_A\|_{L^\infty_{loc}} \leq C\|u^{\e}_A\|_{W^{1,p}_{loc}} \leq C\|\nabla_xu^{\e}_A\|_{L^p}
\enda
\end{align*}
For $p=1$, If $\w^{\e}_0 \in L^1\cap H^{-1}_{loc}(\R^2)$, then we have $u^{\e}_A\in L^2_{loc}$. Let us decompose ($\chi=1$ for $|x|\leq 2R$ and 0 for $|x|>4R$)
\begin{align*}
\bega
u^{\e}_A = {\bf K} \ast \w^{\e}_A = {\bf K} \ast (\chi \w^{\e}_A) + {\bf K} \ast ((1-\chi) \w^{\e}_A).
\enda
\end{align*}
First part: 
\begin{align*}
\bega
\|u^{\e}_{A0}\|_{L^2(\R^2)} = \| \chi \w^{\e}_A\|_{H^{-1}(\R^2)} \leq \|\w^{\e}_A\|_{H^{-1}_{loc}(\R^2)}
\enda
\end{align*}
Second part: $(1-\chi)\w^{\e}_A$ is supported in $|y|\geq 2R$. For $|x|\leq R$, we get $|x-y|\geq |y|/2$. Hence,  
\begin{align*}
\bega
|u^{\e}_{A\infty}(x)| = \int_{|y|\geq 2R} |{\bf K}(x-y)| |\w_{far}(y)| dy \leq \int_{|y|\geq 2R} \frac{|\w_{far}(y)|}{|y|} dy \leq \frac{C}{R}\|\w_{\far}(\R^2)\|_{L^1} \leq \frac{C}{R}\|\w^{\e}_A\|_{L^1(\R^2)}
\enda
\end{align*}
so that, 
\begin{align*}
\bega
\|u^{\e}_{A\infty}\|_{L^2(B_R)} \leq |B_R| \|u^{\e}_{A\infty}\|_{L^\infty(B_R)} \leq |B_R| \frac{C}{R}\|\w^{\e}_A\|_{L^1(\R^2)}
\enda
\end{align*}
Combining the estimate of $u^{\e}_{A0}$ and $u^{\e}_{A\infty}$, we have 
\begin{align*}
\bega
\|u^{\e}_A\|_{L^2_{loc}} \leq \|u^{\e}_{A0}\|_{L^2_{loc}} +\|u^{\e}_{A\infty}\|_{L^2_{loc}} \leq C
\enda
\end{align*}
\unhide
For simplicity, we assume that $\w^{\e}_0$ and $\w_0$ have compact support and belong to $L^p(\R^2)$. 
This assumption can be easily relaxed by assuming $\w^{\e}_0, \w_0 \in L^p \cap L^1(\R^2)$.
We divide the proof into several steps.  \\
{\bf (Step 1: Strong convergence of $u^{\e}_A$ in \eqref{mconv;0})}
Recall that we obtained the uniform boundedness of $\w^{\e}_A(t)$ in \eqref{U.wLp<2} of Lemma \ref{L.unif}. Then, using $u^{\e}_A={\bf K}\ast \w^{\e}_A$, we have
\begin{align}\label{uwAbdd}
\bega
&\sup_{0\leq t\leq T}\|\w^{\e}_A(t)\|_{L^p(\R^2)} = \sup_{0\leq t\leq T}\|\w^{\e}_A(0)\|_{L^p(\R^2)} =\sup_{0\leq t\leq T}\|\w^{\e}_0\|_{L^p(\R^2)}, \cr 
&\sup_{0\leq t\leq T}\|u^{\e}_A(t)\|_{L^p_{loc}(\R^2)} + \sup_{0\leq t\leq T}\|\nabla_xu^{\e}_A(t)\|_{L^p(\R^2)} \leq C,
\enda 
\end{align}
where we used \eqref{PotenLp} and the assumption $\w^{\e}_A(0)=\w^{\e}_0\in L^p_c(\R^2)$.  
In addition, since $\mathbb{P}u^{\e}\in L^\infty(0,T;L^2_{loc}(\R^2))$ by \eqref{U.L2inf} in Lemma \ref{L.unif}, and since $u^{\e}_B\to0$ in $L^\infty(0,T;L^p(\R^2))$ for $2\leq p\leq\infty$ as in \eqref{uBto0}, we see that $u^{\e}_A$ also belongs to $L^\infty(0,T;L^2_{loc}(\R^2))$:
\begin{align}\label{uAbddL2}
\bega
\sup_{0\leq t\leq T}\|u^{\e}_A(t)\|_{L^2_{loc}(\R^2)} &= \sup_{0\leq t\leq T}\|(\mathbb{P}u^{\e}-u^{\e}_B)(t)\|_{L^2_{loc}(\R^2)} \cr 
&\leq \sup_{0\leq t\leq T}\|\mathbb{P}u^{\e}(t)\|_{L^2_{loc}(\R^2)} + \sup_{0\leq t\leq T}\|u^{\e}_B(t)\|_{L^2_{loc}(\R^2)} \leq C.
\enda 
\end{align}
Hence, there exists $(u^{E},\w^{E})$ such that
\begin{align}\label{wconv;w}
\bega
u^{\e}_A &\quad \overset{\ast}{\rightharpoonup} \quad u^{E} \quad \mbox{weakly* in} \quad L^\infty(0,T;W^{1,p}_{loc}(\R^2)), \cr
\w^{\e}_A &\quad \overset{\ast}{\rightharpoonup} \quad \w^{E} \quad \mbox{weakly* in} \quad L^\infty(0,T;L^p(\R^2)).
\enda
\end{align}
To obtain strong convergence of $u^{\e}_A$, we claim that
\begin{align}\label{AubinLionp}
\bega
\p_t u^{\e}_A \in L^\infty(0,T;H^{-s}_{loc}(\R^2)), \quad \mbox{for some} \quad s>0.
\enda 
\end{align}
Using \eqref{uAeqn}, we have
\begin{align}\label{psiu-uesti}
\bega
\big\|\p_t u^{\e}_A(t)\big\|_{H^{-s}_x} &\leq \sup_{\|\psi\|_{H^s}\leq 1} \int \psi \Big(u^{\e}_A \cdot \nabla_x u^{\e}_A + \nabla_x p^{\e}_A +\nabla_x^{\perp}(-\Delta_x)^{-1}\nabla_x \cdot \big(u^{\e}_B\w^{\e}_A\big) \Big) dx.
\enda 
\end{align}

Using $u^{\e}_A\in L^\infty(0,T;L^2_{loc}(\R^2))$, we can bound the first two terms on the right-hand side of~\eqref{psiu-uesti} in a straightforward manner. We therefore focus on the term
\begin{align}\label{IAdef}
\bega
I_A(t):= \int_{\R^2} \psi \big(\nabla_x^{\perp}(-\Delta_x)^{-1}\nabla_x \cdot \big(u^{\e}_B\w^{\e}_A\big)\big)  dx .
\enda 
\end{align}
Integrating by parts, we obtain
\begin{align*}
\bega
I_A(t) \leq \sup_{\|\psi\|_{H^2(\R^2)}\leq 1} \|u^{\e}_B(t)\|_{L^\infty_x}\|\w^{\e}_A(t)\|_{L^1_{loc}} \| \nabla_x (-\Delta_x)^{-1}\nabla_x^{\perp} \psi\|_{L^\infty_x}.
\enda 
\end{align*}
Applying Agmon’s inequality~\eqref{Agmon} to control
$\big\|\nabla_x (-\Delta_x)^{-1}\nabla_x^{\perp} \psi\big\|_{L^\infty_x}$,
and then using~\eqref{uBto0}, we deduce
\begin{align*}
\bega
\sup_{0\leq t\leq T}I_A(t)\leq \sup_{0\leq t\leq T} \|u^{\e}_B(t)\|_{L^\infty_x} \|\w^{\e}_0\|_{L^p_{loc}} \|\psi\|_{H^2_x} \to 0, \quad \mbox{as} \quad \e\to0.
\enda 
\end{align*}
This proves~\eqref{AubinLionp}. 
By the Sobolev embedding in two dimensions, \eqref{uwAbdd} implies that
$u^{\e}_A \in L^\infty(0,T;L^{p_*}_{\mathrm{loc}}(\R^2))$ with $p_*=2p/(2-p)>2$ for $1<p<2$. Hence, combining \eqref{AubinLionp} and \eqref{uwAbdd}, we obtain the strong convergence of $u^{\e}_A$ in~\eqref{mconv;0}.
  \\
To prove that $u^{E}= {\bf K}\ast\w^{E}$, we consider a test function $\psi \in C_c^{\infty}((0,T)\times \R^2)$. Using \eqref{mconv;0} and \eqref{wconv;w}, we compute
\begin{align*}
\bega
\int_0^T \int_{\R^2} &u^{E} \psi dxdt = \lim_{\e\rightarrow 0}\int_0^T \int_{\R^2} u^{\e}_A \psi dxdt = \lim_{\e\rightarrow 0}\int_0^T \int_{\R^2} ({\bf K}\ast \w^{\e}_A) \psi dxdt \cr 
&= -\lim_{\e\rightarrow 0}\int_0^T \int_{\R^2} \w^{\e}_A ({\bf K}\ast\psi) dxdt = -\int_0^T \int_{\R^2} \w^{E} ({\bf K}\ast\psi) dxdt = \int_0^T \int_{\R^2} ({\bf K}\ast\w^{E}) \psi dxdt,
\enda
\end{align*}
where we used that ${\bf K} \ast \psi \in L^\infty((0,T) \times \R^2)$.  
Since the above equality holds for every test function $\psi$, we conclude that 
$u^{E} = {\bf K} \ast \w^{E}$, a.e in $(0,T) \times \R^2$. 
\\

{\bf (Step 2: $(u^{E},\w^{E})$ solves the Euler equations in the sense of renormalized solutions)}
We first consider the dual problem associated with equation~\eqref{weqnA} and its limit.  
We rewrite equation~\eqref{weqnA} in divergence form as
\begin{align}\label{weqndiv}
\p_t \w^{\e}_A + \nabla_x \cdot (\mathbb{P}u^{\e}\w^{\e}_A)  = 0.
\end{align}
For $\chi \in C_c^{\infty}((0,T)\times \R^2)$, we consider the corresponding dual problem:
\begin{align}\label{dual;BE}
\bega
-\p_t \phi^{\e} - \nabla_x \cdot (\mathbb{P}u^{\e} \phi^{\e}) = \chi, \qquad \phi^{\e}(T,x) = 0.
\enda
\end{align}
Multiplying~\eqref{weqndiv} by $\phi^{\e}$ and using~\eqref{dual;BE}, we obtain
\begin{align}\label{dual;form}
\bega
\int_0^T \int_{\R^2} \w^{\e}_A(t,x) \chi(t,x) \, dxdt &= \int_{\R^2} \w^{\e}_0(x) \phi^{\e}(0,x) \, dx .
\enda
\end{align}
From~\eqref{dual;BE}, since $\mathbb{P}u^{\e}$ is divergence-free, we have
\begin{align}\label{phiebdd}
\bega
\{\phi^{\e}\}_{\e} \subset L^{\infty}(0,T;L^p(\R^2)), \quad \mbox{for all} \quad 1\leq p\leq \infty. 
\enda
\end{align}

Hence, there exists $\phi^{\#}$ (up to a subsequence) such that 
\begin{align}\label{phiweak}
\bega
\phi^{\e} &\quad \overset{\ast}{\rightharpoonup} \quad \phi^{\#} \quad \mbox{weakly* in} \quad L^\infty(0,T;L^\infty \cap L^1(\R^2)).
\enda
\end{align}
We now consider the dual problem associated with the Euler equation:
\begin{align}\label{dual;E}
\bega
-\p_t \phi - \nabla_x \cdot (u^{E} \phi) = \chi, \qquad \phi(T,x) = 0.
\enda
\end{align}
To compare $\phi^{\e}$ and $\phi$, define for any $\psi \in C_c^{\infty}(\R^2)$
\begin{align*}
\bega
A_{\psi}^{\e}(t) := \int_{\R^2} \phi^{\e}(t,x) \psi(x) \, dx, \qquad A_{\psi}(t) := \int_{\R^2} \phi(t,x) \psi(x) \, dx.
\enda
\end{align*}

By the fundamental theorem of calculus, for $0 \le t \le T$ we have
\begin{align*}
\bega
A_{\psi}^{\e}(T) - A_{\psi}(T) =  A_{\psi}^{\e}(t) - A_{\psi}(t) + \int_t^T \frac{d}{dt} \big(A_{\psi}^{\e} - A_{\psi}\big)(s) \, ds, \quad \mbox{for } \quad 0 \leq t \leq T.
\enda
\end{align*}
Since
\begin{align}\label{PuL2conv}
\sup_{t\in[0,T]}\|\mathbb{P}u^{\e}-u^{E}\|_{L^2_{loc}} \leq \sup_{t\in[0,T]}\|u^{\e}_A-u^{E}\|_{L^2_{loc}} + \sup_{t\in[0,T]}\|u^{\e}_B\|_{L^2} \to 0, \quad \mbox{as} \quad \e\to0, 
\end{align}
by~\eqref{mconv;0} and~\eqref{uBto0}, we obtain
\begin{align*}
\bega
\int_t^T &\frac{d}{dt} (A_{\psi}^{\e}-A_{\psi})(s) ds = \int_t^T \int_{\R^2} (\phi^{\e}-\phi)u^{E} \cdot \nabla_x\psi dxdt +\int_t^T \int_{\R^2} \phi^{\e}(\mathbb{P}u^{\e}-u^{E}) \cdot \nabla_x\psi dxdt \cr 
&\leq  \int_t^T \int_{\R^2} (\phi^{\e}-\phi)u^{E} \cdot \nabla_x\psi dxdt + C\sup_{0\leq t\leq T}\|\phi^{\e}(t)\|_{L^2_x}\sup_{0\leq t\leq T}\|\mathbb{P}u^{\e}(t)-u^{E}(t)\|_{L^2_x}\|\nabla_x\psi\|_{L^\infty_x},
\enda
\end{align*}
where we used~\eqref{phiebdd},~\eqref{phiweak}, and the fact that 
$u^{E}\in L^\infty(0,T;L^p_{\mathrm{loc}}(\R^2))$.
In addition, we have $A_{\psi}^{\e}(T) - A_{\psi}(T) = 0$, since 
$\phi^{\e}(T,x) = 0$ and $\phi^{\#}(T,x) = 0$.  
Because $C_c^{\infty}(\R^2)$ is dense in $L^p(\R^2)$ for any $1<p<\infty$, it follows that
\begin{align*}
\bega
\int_{\R^2}\phi^{\e}(t,x)\psi(x) dx \rightarrow \int_{\R^2}\phi(t,x)\psi(x) dx \quad \mbox{in} \quad C([0,T]),
\enda
\end{align*}
for every $\psi \in L^q(\R^2)$ with $1<q<\infty$.  
By the uniqueness of the weak limit, this allows us to upgrade the convergence~\eqref{phiweak} to
\begin{align}\label{phistrong}
\bega
\int_{\R^2}\phi^{\e}(t,x)\psi(x) dx \rightarrow \int_{\R^2}\phi^{\#}(t,x)\psi(x) dx \quad \mbox{in} \quad C([0,T]).
\enda
\end{align}

Next, for any $\psi \in C_c^{\infty}([0,T)\times\R^2)$, 
multiplying~\eqref{dual;BE} by $\psi$ and integrating gives
\begin{align*}
\int_0^T\int_{\R^2} \phi^{\e}(\p_t\psi + \mathbb{P}u^{\e} \cdot\nabla_x\psi) dxdt + \int_{\R^2} \phi^{\e}_0 \psi_0 dx = \int_0^T\int_{\R^2} \chi \psi dxdt .
\end{align*}
Passing to the limit as $\e\to0$ and using~\eqref{PuL2conv} together with~\eqref{phistrong}, we obtain
\begin{align*}
\int_0^T\int_{\R^2} \phi^{\#}(\p_t\psi + u^{E} \cdot\nabla_x\psi) dxdt + \int_{\R^2} \phi^{\#}_0 \psi_0 dx = \int_0^T\int_{\R^2} \chi \psi dxdt .
\end{align*}

Applying the convergence of the initial data, \eqref{wconv;w}, and \eqref{phistrong}, we can pass to the limit in~\eqref{dual;form} and obtain
\begin{align}\label{wchi}
\bega
\int_0^T \int_{\R^2} \w^{E}(t,x) \chi(t,x) \, dxdt = \int_{\R^2} \w_0(x) \phi^{\#}(0,x) \, dx.
\enda
\end{align}
Following the approach developed in \cite{DiLi,Crippa2021}, we first observe that $u^{E}$ satisfies~\eqref{bassume} by decomposing
\begin{align}\label{bhold}
u^{E}= {\bf K}\ast \w^{E} = ({\bf K}\mathbf{1}_{|x|\leq1})\ast \w^{E} + ({\bf K}\mathbf{1}_{|x|>1})\ast \w^{E},
\end{align}
since $({\bf K}\mathbf{1}_{|x|\leq1}) \in L^1$, 
$({\bf K}\mathbf{1}_{|x|>1}) \in L^\infty$, 
and $\w^{E} \in L^1$ (as $\w_0 \in L^p_c(\R^2)$).  
By~\eqref{mconv;0}, we also have 
$u^{E} \in L^\infty(0,T;W^{1,p}_{loc}(\R^2))$.  
Thus, by Theorem~\ref{T.Di-Li}, there exists a unique renormalized solution 
$\w^{*} \in L^\infty(0,T;L^q(\R^2))$ of the continuity equation
\begin{align*}
\p_t \w^{*} + u^{E}\cdot \nabla_x \w^{*} = 0, \qquad \w^{*}(0,x) = \w_0(x).
\end{align*}
Moreover, since $\phi^{\#}$ is uniformly bounded in 
$L^\infty(0,T;L^p(\R^2))$ for any $1 \leq p \leq \infty$, 
$\phi^{\#}$ is also the unique renormalized solution of the dual problem~\eqref{dual;E} by Theorem~\ref{T.Di-Li}. 
Using the same dual formulation as in~\eqref{dual;E}, 
the function $\w^{*}$ satisfies
\begin{align}\label{bwchi}
\bega
\int_0^T \int_{\R^2} \w^{*}(t,x) \chi(t,x) \, dxdt = \int_{\R^2} \w_0(x) \phi^{\#}(0,x) \, dx.
\enda
\end{align}
Subtracting~\eqref{wchi} from~\eqref{bwchi} and using the uniqueness of $\phi^{\#}$ as the renormalized solution, we obtain
\begin{align*}
\bega
\int_0^T \int_{\R^2} (\w^{*} - \w^{E})(t,x) \chi(t,x) \, dxdt = 0, \quad \mbox{for any} \quad \chi \in C_c^{\infty}((0,T)\times \R^2).
\enda
\end{align*}
This implies that
\begin{align*}
\bega
\w^{E}(t,x) = \w^{*}(t,x), \quad \mbox{almost everywhere in} \quad (0,T) \times \R^2.
\enda
\end{align*}

\subsubsection{{\bf Proof of Theorem \ref{T.C.Lp} for $p=1$}}
We prove the case $p=1$ by following the same steps as in the case $1<p<\infty$.\\
{\bf (Step 1: Strong convergence of $u^{\e}_A$ in \eqref{mconv;0})}
For $p=1$, we first observe that, in the same way as in \eqref{uwAbdd}, since $\mathbb{P}u^{\e}$ is divergence-free and thus measure-preserving, it follows from \eqref{weqnA} that
\begin{align*}
\bega
&\sup_{0\leq t\leq T}\|\w^{\e}_A(t)\|_{L^1(\R^2)} = \sup_{0\leq t\leq T}\|\w^{\e}_A(0)\|_{L^1(\R^2)} =\sup_{0\leq t\leq T}\|\w^{\e}_0\|_{L^1(\R^2)}, \cr 
&\sup_{0\leq t\leq T}\|u^{\e}_A(t)\|_{L^2_{loc}(\R^2)} \leq C,
\enda
\end{align*}
where we used \eqref{uAbddL2}.  
For $p=1$, to prove the weak convergence of $\w^{\e}_A$ in
\begin{align}\label{wconv;wL1}
\bega
\w^{\e}_A &\quad \rightharpoonup \quad \w^{E} \quad \mbox{weakly in} \quad L^\infty(0,T;L^1(\R^2)),
\enda
\end{align}
it suffices to establish the equi-integrability of the family
$\{\w^{\e}_A\}_{\e>0}$ in $L^1((0,T)\times \R^2)$.  
Following the argument in~\cite{Crippa2021}, we use the fact that
$\w^{\e}_A(0) \to \w_0$ strongly in $L^1_x$ as $\e \to 0$.
Let $\sup_{0\leq t\leq T}\|\w^{\e}_A(t)\|_{L^1(\R^2)}:= C_{\w 0}$, and decompose
\begin{align*}
\bega
\w^{\e}_A(0) = \w^{\e,\eta}_{A,1}(0) +\w^{\e,\eta}_{A,\infty}(0), \quad \w^{\e,\eta}_{A,1}(0) := \w^{\e}_A(0) \mathbf{1}_{\{\w^{\e}_A(0)|> \eta/C_{\w 0}\}}, \quad \w^{\e,\eta}_{A,\infty}(0) := \w^{\e}_A(0) \mathbf{1}_{\{|\w^{\e}_A(0)|\leq \eta/C_{\w 0}\}},
\enda
\end{align*}
so that
\begin{align*}
\bega
\|\w^{\e,\eta}_{A,1}(0)\|_{L^1_x} \leq \eta, \quad \mbox{and} \quad   \|\w^{\e,\eta}_{A,\infty}(0)\|_{L^\infty_x} \leq C_{\eta}.
\enda
\end{align*}
Then both $\w^{\e,\eta}_{A,1}(0)$ and $\w^{\e,\eta}_{A,\infty}(0)$ belong to
$L^2(\R^2)\cap L^\infty(\R^2)$.  
Let $\w^{\e,\eta}_{A,1}, \w^{\e,\eta}_{A,\infty}\in
L^\infty(0,T;L^1(\R^2)\cap L^\infty(\R^2))$ be the unique weak solutions of
\begin{align*}
\p_t \w^{\e,\eta}_{A,1} &+ \mathbb{P}u^{\e}\cdot\nabla_x \w^{\e,\eta}_{A,1} = 0, \qquad 
\w^{\e,\eta}_{A,1}(0,x) = \w^{\e,\eta}_{A,1}(0), \cr
\p_t \w^{\e,\eta}_{A,\infty} &+ \mathbb{P}u^{\e}\cdot\nabla_x \w^{\e,\eta}_{A,\infty} = 0, \qquad \w^{\e,\eta}_{A,\infty}(0,x) = \w^{\e,\eta}_{A,\infty}(0).
\end{align*}
Since $\mathbb{P}u^{\e}$ is divergence-free, we have
\begin{align*}
\|\w^{\e,\eta}_{A,1}(t)\|_{L^1_x} =\|\w^{\e,\eta}_{A,1}(0)\|_{L^1_x} \leq \eta, \qquad \|\w^{\e,\eta}_{A,\infty}(t)\|_{L^\infty_x} =\|\w^{\e,\eta}_{A,\infty}(0)\|_{L^\infty_x} \leq C_{\eta}.
\end{align*}
We approximate the absolute value function $|\cdot|$ by
$\phi_{\delta}(s):=\sqrt{s^2+\delta^2}$ and multiply the equation by
$\phi'_{\delta}(\w)=\frac{\w}{\sqrt{\w^2+\delta^2}}$ to obtain
\begin{align*}
\p_t \phi_{\delta}(\w^{\e,\eta}_{A,\infty}) + \mathbb{P}u^{\e}\cdot\nabla_x \phi_{\delta}(\w^{\e,\eta}_{A,\infty}) = 0.
\end{align*}
Next, we multiply by a smooth cutoff function
$\psi_r^R\in C_c^{\infty}(\R^2)$ satisfying
$\psi_r^R=0$ for $0<|x|<r$ and $|x|>2R$, and $\psi_r^R=1$ for $2r<|x|<R$.
Letting $\delta\to 0$, we deduce
\begin{align*}
\int_{|x|\geq r} |\w^{\e,\eta}_{A,\infty}(t)|dx \leq \frac{C}{r}.
\end{align*}
By uniqueness of the linear transport equation, $\w^{\e}_A(t)=\w^{\e,\eta}_{A,1}(t)+\w^{\e,\eta}_{A,\infty}(t)$ solves \eqref{weqnA}. Moreover,
\begin{align*}
\int_{E} |\w^{\e}_{A}(t)|dx \leq \int_{E} |\w^{\e,\eta}_{A,1}(t)|dx +\int_{E \cap B_r(0)} |\w^{\e,\eta}_{A,\infty}(t)|dx + \int_{E \cap B_r^c(0)} |\w^{\e,\eta}_{A,\infty}(t)|dx \leq \eta + |E|C_{\eta} + \frac{C}{r}.
\end{align*}
This shows that $\{\w^{\e}_A\}_{\e}$ is equi-integrable in
$L^1((0,T)\times \R^2)$. Therefore, by the Dunford--Pettis theorem,
we obtain \eqref{wconv;wL1} for $p=1$.

\hide
For the function $\chi \in C_c^{\infty}(\R^2)$ with $\chi=1$ for $|x|\leq 2R$ and 0 for $|x|>4R$, we decompose
\begin{align*}
\bega
u^{\e}_A = {\bf K} \ast \w^{\e}_A = {\bf K} \ast (\chi \w^{\e}_A) + {\bf K} \ast ((1-\chi) \w^{\e}_A) := u^{\e}_{A0} + u^{\e}_{A\infty}.
\enda
\end{align*}
For the part $u^{\e}_{A0}$, since $\w^{\e}_A$ is in $H^{-1}_{loc}(\R^2)$, we have 
\begin{align*}
\bega
\sup_{0\leq t\leq T}\|u^{\e}_{A0}(t)\|_{L^2(\R^2)} = \| \chi \w^{\e}_A\|_{H^{-1}(\R^2)} \leq \|\w^{\e}_A\|_{H^{-1}_{loc}(\R^2)}
\enda
\end{align*}
Second part: $(1-\chi)\w^{\e}_A$ is supported in $|y|\geq 2R$. For $|x|\leq R$, we get $|x-y|\geq |y|/2$. Hence,  
\begin{align*}
\bega
|u^{\e}_{A\infty}(x)| = \int_{|y|\geq 2R} |{\bf K}(x-y)| |\w_{far}(y)| dy \leq \int_{|y|\geq 2R} \frac{|\w_{far}(y)|}{|y|} dy \leq \frac{C}{R}\|\w_{\far}(\R^2)\|_{L^1} \leq \frac{C}{R}\|\w^{\e}_A\|_{L^1(\R^2)}.
\enda
\end{align*}
so that, 
\begin{align*}
\bega
\|u^{\e}_{A\infty}\|_{L^2(B_R)} \leq |B_R| \|u^{\e}_{A\infty}\|_{L^\infty(B_R)} \leq |B_R| \frac{C}{R}\|\w^{\e}_A\|_{L^1(\R^2)}.
\enda
\end{align*}
Combining the estimate of $u^{\e}_{A0}$ and $u^{\e}_{A\infty}$, we have 
\begin{align*}
\bega
\|u^{\e}_A\|_{L^2_{loc}} \leq \|u^{\e}_{A0}\|_{L^2_{loc}} +\|u^{\e}_{A\infty}\|_{L^2_{loc}} \leq C.
\enda
\end{align*}
\unhide

Next, we consider the strong convergence of $u^{\e}_A$ in \eqref{mconv;0}.
For $p=1$, the convergence in \eqref{mconv;0} is already standard for the Euler equation by Theorem~1.1 in \cite{DiMa}. We therefore only sketch the proof, focusing on the microscopic part.
Since $u^{\e}_A$ and $u^{E}$ belong to $L^2_{\mathrm{loc}}([0,T]\times \R^2)$, by interpolation for $1<r<2$ we get
\begin{align*}
\bigg(\int_0^T \int_{\Omega} |u^{\e}_A-u^{E}|^r dxdt\bigg)^{\frac{1}{r}} \leq \|u^{\e}_A-u^{E}\|_{L^1_{loc}([0,T]\times \R^2)}^{\frac{2}{r}-1} \|u^{\e}_A-u^{E}\|_{L^2_{loc}([0,T]\times \R^2)}^{2-\frac{2}{r}}, 
\end{align*}
for any compact set $\Omega \subset \R^2$. Hence, it suffices to prove convergence in $L^1$, namely,
\begin{align}\label{claimuconvL1}
\int_0^T\int_{\Omega}|u^{\e}_A(t) - u^{E}(t)| dx \to 0, \quad \mbox{as} \quad \e \to 0.
\end{align}
For $\psi \in C_c^{\infty}(\R^2)$, we first claim that
\begin{align}\label{AubinLion}
\bega
\big\|\psi u^{\e}_A(t_2) -\psi u^{\e}_A(t_1)\big\|_{H^{-s-1}_x}  \leq C |t_2-t_1|, \quad \mbox{for} \quad \psi \in C_c^{\infty}(\R^2), \quad s>1, \quad 0\leq t_1,t_2\leq T.
\enda 
\end{align}
To prove \eqref{AubinLion}, we use the equation for $u^{\e}_A$ in \eqref{uAeqn}.
By the fundamental theorem of calculus in time, we obtain
\begin{align}\label{psiu-uesti}
\bega
\big\|\psi u^{\e}_A(t_2) -&\psi u^{\e}_A(t_2)\big\|_{H^{-s}_x}  \leq \bigg\| \psi \int_{t_1}^{t_2} \p_t u^{\e} dt \bigg\|_{H^{-s}_x} \leq |t_2-t_1| \max_{0\leq t\leq T} \big\| \psi \p_t u^{\e} \big\|_{H^{-s}_x} \cr 
&\leq \sup_{0\leq t\leq T} \bigg\| \psi \Big(u^{\e}_A \cdot \nabla_x u^{\e}_A + \nabla_x p^{\e}_A +\nabla_x^{\perp}(-\Delta_x)^{-1}\nabla_x \cdot \big(u^{\e}_B\w^{\e}_A\big) \Big) \bigg\|_{H^{-s}_x},
\enda 
\end{align}
for some $s>0$. By the same argument as in the estimate of $I_A(t)$ defined in \eqref{IAdef}, and using the uniform boundedness of $u^{\e}_A$ in $L^2_{loc}$, we obtain \eqref{AubinLion}.
From \eqref{AubinLion}, it follows that
\begin{align}\label{wtimecpt}
\bega
\big\|\psi \w^{\e}_A(t_2) -\psi \w^{\e}_A(t_1)\big\|_{H^{-s-1}_x}  \leq C |t_2-t_1|, \quad \mbox{for} \quad \psi \in C_c^{\infty}(\R^2), \quad 0\leq t_1,t_2\leq T.
\enda 
\end{align}
Moreover, for $k>1$, we have
\begin{align}\label{wH-sbdd}
\bega
\sup_{0\leq t\leq T}\big\|\psi \w^{\e}_A(t)\big\|_{H^{-k}(\R^2)}  \leq \sup_{0\leq t\leq T}\big\|\psi \w^{\e}_A(t)\big\|_{L^1(\R^2)} \leq \|\psi\|_{L^\infty(\R^2)} \sup_{0\leq t\leq T}\big\|\w^{\e}_A(t)\big\|_{L^1(\R^2)} \leq C.
\enda 
\end{align}
Combining \eqref{wtimecpt} and \eqref{wH-sbdd} with the Aubin--Lions lemma, we obtain
\begin{align}\label{wH-sconv}
\sup_{0\leq t\leq T}\|\psi\w^{\e}_A(t) - \psi\w^{E}(t)\|_{H^{-k}(\R^2)} \to 0, \quad \mbox{as} \quad  \e \to 0, \quad \mbox{for} \quad k>1.
\end{align}
Finally, to prove \eqref{claimuconvL1}, we show that $\{u^{\e}_A\}$ is a Cauchy sequence.
Let $\chi\in C_c^{\infty}(\R^2)$ be a smooth cutoff function such that
$\chi(x)=1$ for $|x|<1$ and $\chi(x)=0$ for $|x|>2$, and define $\chi_r(|x|):=\chi(|x|/r)$.
For $0<r<2r<R$, we decompose
\begin{align*}
\bega
(u^{\e_1}_A -u^{\e_2}_A)(t,x) &= \int_{\R^2} \bigg(\chi_r{\bf K}+ \big(\chi_R -\chi_r\big){\bf K} + \big(1-\chi_R\big){\bf K} \bigg)(x-y) (\w^{\e_1}_A-\w^{\e_2}_A)(t,y) dy \cr 
&= I_{uu}(t,x)+II_{uu}(t,x)+III_{uu}(t,x).
\enda
\end{align*}
\hide
where
\begin{align*}
\bega
I_{uu}(t,x) &=\int_{\R^2} \chi\Big(\frac{|x-y|}{r}\Big){\bf K}(x-y) (\w^{\e_1}_A-\w^{\e_2}_A)dy, \cr 
II_{uu}(t,x) &=\int_{\R^2} \bigg(\chi\Big(\frac{|x-y|}{R}\Big) -\chi\Big(\frac{|x-y|}{r}\Big) \bigg) {\bf K}(x-y) (\w^{\e_1}_A-\w^{\e_2}_A) dy, \cr 
III_{uu}(t,x) &=\int_{\R^2} \bigg(1-\chi\Big(\frac{|x-y|}{R}\Big)\bigg){\bf K}(x-y) (\w^{\e_1}_A-\w^{\e_2}_A) dy.
\enda
\end{align*}
\unhide
For $I_{uu}$ and $III_{uu}$, using $\sup_{0\leq t\leq T}\|\w^{\e}(t)\|_{L^1(\R^2)} \leq C$, we have
\begin{align*}
\bega
\|I_{uu}(t,x)\|_{L^1([0,T]\times\R^2)} &\les \big\|\chi_r{\bf K}\big\|_{L^1(\R^2)} \|(\w^{\e_1}_A-\w^{\e_2}_A)\|_{L^1([0,T]\times\R^2)} \leq Cr , \cr 
\|III_{uu}(t,x)\|_{L^1_{loc}([0,T]\times\R^2)} &\leq C\|III_{uu}(t,x)\|_{L^\infty([0,T]\times\R^2)} \les \frac{C}{R} \|(\w^{\e_1}_A-\w^{\e_2}_A)\|_{L^1([0,T]\times\R^2)} \leq \frac{C}{R}.
\enda
\end{align*}
For $II_{uu}$, applying H\"older's inequality and using \eqref{wH-sconv}, we obtain
\begin{align*}
\bega
|II_{uu}(t,x)| \leq \big\|(\chi_R -\chi_r)\big\|_{H^k(\R^2)} \|\chi_{2R}(\w^{\e_1}_A-\w^{\e_2}_A)\|_{H^{-k}(\R^2)} \to 0, \quad \mbox{as} \quad \e\to 0.
\enda
\end{align*}
This proves \eqref{claimuconvL1} and completes the convergence \eqref{mconv;0} for $p=1$.
\\ {\bf (Step 2: $(u^{E},\w^{E})$ solves the Euler equations in the sense of renormalized solutions)}
For $p=1$, we follow the argument of~\cite{DiLi,B-Crippa}.  
Since the estimate $\|\nabla_x u\|_{L^p_x} \leq C_* p \|\w\|_{L^p_x}$ in~\eqref{puw} does not hold for $p=1$, 
we instead consider the continuity equation with velocity field $u^{E}$:
\begin{align}\label{rhou3}
\p_t \varrho + u^{E} \cdot\nabla_x \varrho = 0, \qquad 
\varrho|_{t=0} = \w_0.
\end{align}
We note that the decomposition~\eqref{bhold} also holds for $p=1$.  
Since we already proved that $u^{E} = {\bf K} \ast \w^{E}$ and 
$u^{E} \in L^r(0,T;L^r_{loc}(\R^2))$ for $1 \leq r < 2$, 
the assumptions~\eqref{bassume} and~\eqref{bassume2} are satisfied.  
Therefore, by Theorem~\ref{T.B-Crippa}, the function $\varrho$ is the unique Lagrangian solution of~\eqref{rhou3}, and it is also a renormalized solution.  
Next, we mollify both the velocity field $u^{E}$ and the initial data 
using the standard mollifier~\eqref{molli} with parameter $\lambda$:
\begin{align*}
u^{E,\lambda}(t,x):= (u^{E}\ast \varphi^{\lambda})(t,x), \qquad \varrho_0^{\lambda}(x) : = (\w_0 \ast \varphi^{\lambda})(x).
\end{align*}
By the standard properties of mollification and the uniform boundedness of 
$u^{E}$ in $L^\infty(0,T;L^2_{loc}(\R^2))$, we have
\begin{align}\label{u3lu3} 
u^{E,\lambda} \quad \rightarrow \quad u^{E} \quad \mbox{strongly in} \quad L^1_{loc}([0,T]\times \R^2).
\end{align}
For each $\lambda > 0$, consider the problem
\begin{align}\label{rhou3lam}
\p_t \varrho^{\lambda} + u^{E,\lambda} \cdot\nabla_x \varrho^{\lambda} = 0, \qquad 
\varrho^{\lambda}|_{t=0} = \w_0 \ast \varphi^{\lambda}(x).
\end{align}
By the Cauchy–Lipschitz theorem, there exists a sequence of smooth solutions 
$\{\varrho^{\lambda}\}_{\lambda}$ satisfying
\begin{align*}
\|\varrho^{\lambda}(t)\|_{L^1_x} = \|\varrho^{\lambda}(0)\|_{L^1_x}= \|\w_0 \ast \varphi^{\lambda}\|_{L^1_x}= \|\w_0\|_{L^1_x},
\end{align*}
where we used $\div(u^{E,\lambda}) = 0$.  
Hence, $\{\varrho^{\lambda}\}_{\lambda}$ is uniformly bounded in 
$L^\infty(0,T;L^1(\R^2))$.  
By the same reasoning as for $u^{E}$, the field $u^{E,\lambda}$ also satisfies 
assumptions~\eqref{bassume} and~\eqref{bassume2}.  
Decomposing $u^{E,\lambda}$ as in~\eqref{bhold},  
$u^{E,\lambda} = u^{E,\lambda}_1 + u^{E,\lambda}_2$, 
we find that $u^{E,\lambda}_1$ and $u^{E,\lambda}_2$ are equibounded in 
$L^1(0,T;L^1(\R^2))$ and $L^1(0,T;L^\infty(\R^2))$, respectively.  
Thus, applying (2) of Theorem~\ref{T.B-Crippa} gives
\begin{align}\label{varrhostrong} 
\varrho^{\lambda} \quad \rightarrow \quad \varrho \quad \mbox{strongly in} \quad C(0,T;L^1(\R^2)).
\end{align}
Now, consider the dual problem of~\eqref{rhou3lam}:
\begin{align}\label{philamu}
\bega
-\p_t \phi^{\lambda} - \nabla_x \cdot (u^{E,\lambda} \phi^{\lambda}) = \chi, \qquad \phi^{\lambda}(T,x) = 0,
\enda
\end{align}
where $\chi \in C_c^{\infty}((0,T)\times \R^2)$.  
Multiplying~\eqref{rhou3lam} by $\phi^{\lambda}$ yields
\begin{align*}
\bega
\int_0^T \int_{\R^2} \varrho^{\lambda} \chi dxdt = \int_{\R^2} \varrho_0^{\lambda} \phi_0^{\lambda} dx.
\enda
\end{align*}
Since $u^{E,\lambda}$ is divergence-free, we have
\begin{align*}
\bega
\{\phi^{\lambda}\}_{\lambda} \subset L^{\infty}(0,T;L^p(\R^2)), \quad \mbox{for all} \quad 1\leq p\leq \infty. 
\enda
\end{align*}
Hence, there exists $\bar{\phi}$ (up to a subsequence) such that
\begin{align}\label{philamweak}
\bega
\phi^{\lambda} &\quad \overset{\ast}{\rightharpoonup} \quad \bar{\phi} \quad \mbox{weakly* in} \quad L^\infty(0,T;L^\infty \cap L^1(\R^2)).
\enda
\end{align}
For any $\psi \in C_c^{\infty}([0,T)\times \R^2)$, multiplying $\psi$ by~\eqref{philamu} gives
\begin{align*}
\int_0^T\int_{\R^2} \phi^{\lambda}(\p_t\psi + u^{E,\lambda} \cdot\nabla_x\psi) dxdt + \int_{\R^2} \phi^{\lambda}_0 \psi_0 dx = \int_0^T\int_{\R^2} \chi \psi dxdt .
\end{align*}
Passing to the limit and using~\eqref{u3lu3} and~\eqref{philamweak}, we obtain
\begin{align*}
\int_0^T\int_{\R^2} \bar{\phi}(\p_t\psi + u^{E} \cdot\nabla_x\psi) dxdt + \int_{\R^2} \bar{\phi}_0 \psi_0 dx = \int_0^T\int_{\R^2} \chi \psi dxdt .
\end{align*}
Hence, $\bar{\phi}$ is a distributional solution in 
$L^\infty(0,T;L^\infty \cap L^1(\R^2))$ of
\begin{align}\label{dual;Lag,lim}
\bega
-\p_t \bar{\phi} - \nabla_x \cdot (u^{E} \bar{\phi}) = \chi, \qquad \bar{\phi}(T,x) = 0.
\enda
\end{align}
Since $\phi^{\lambda}_0 \in L^\infty(\R^2)$ and $\varrho_0 \in L^1(\R^2)$, 
using~\eqref{varrhostrong} and~\eqref{philamweak}, we deduce
\begin{align}\label{rhochi}
\bega
\int_0^T \int_{\R^2} \varrho \chi dxdt = \int_{\R^2} \w_0 \bar{\phi}_0 dx,
\enda
\end{align}
where we used $\varrho^{\lambda}_0 = \w_0 \ast \varphi^{\lambda} \to \w_0$ 
strongly in $L^1(\R^2)$.  
\hide
By using the convergence of $u^{E,\lambda}$ in \eqref{u3lu3}, we can upgrade the weak star convergence to 
\begin{align*}
\bega
\phi^{\lambda} &\quad \to \quad \bar{\phi} \quad \mbox{in} \quad C(0,T;L^p_{weak}(\R^2)), \quad \mbox{for all} \quad  1<p<\infty. 
\enda
\end{align*}
\unhide
Subtracting~\eqref{dual;E} from~\eqref{dual;Lag,lim}, we find that 
$\phi - \bar{\phi}$ is a distributional solution in 
$L^\infty(0,T;L^\infty \cap L^1(\R^2))$ satisfying
\begin{align*}
\bega
\p_t (\phi-\bar{\phi}) + \nabla_x \cdot (u^{E}(\phi-\bar{\phi})) = 0, \qquad (\phi-\bar{\phi})(T,x) = 0.
\enda
\end{align*}
We note that $u^{E}$ satisfies~\eqref{bassume2} because 
$u^{E} \in L^2_{loc}([0,T]\times \R^2)$ and 
$u^{E} = {\bf K} \ast \w^{E}$.  
Thus, by Theorem~\ref{T.B-Crippa}, such a solution is unique, 
and hence $\phi^{\#} = \bar{\phi}$ almost everywhere in $(0,T)\times \R^2$.  
We note that, by the same argument as in the case $1<p<\infty$, we have \eqref{wchi}. Hence, subtracting~\eqref{wchi} from~\eqref{rhochi}, we then obtain
\begin{align*}
\bega
\int_0^T \int_{\R^2} (\w^{E}-\varrho)(t,x) \chi(t,x) \, dxdt = \int_{\R^2} \w_0(x) (\phi^{\#}-\bar{\phi})(0,x) \, dx = 0, 
\enda
\end{align*}
for all $\chi \in C_c^{\infty}((0,T)\times \R^2)$.  
Thus, we conclude that $\w^{E} = \varrho$ almost everywhere in $(0,T)\times \R^2$.

\subsubsection{ {\bf Strong convergence of $\w^{\e}_A$ in \eqref{wconv;0} for $1\leq p<\infty$ } }

Let $\varphi_n$ be a mollifier in $\R^2$. We define $\w^{\e}_n$ as the solution of the following problem, which has the same convection term as \eqref{weqnP}:
\begin{align}\label{weqn-BE-n}
\begin{cases}
\p_t \w^{\e}_{A,n} + \mathbb{P}u^{\e}\cdot\nabla_x \w^{\e}_{A,n} = 0, \\ 
\w^{\e}_{A,n}(0,x) = \w^{\e}_0 \ast \varphi_n(x).
\end{cases}
\end{align}
We also define $\w^{E}_n$ as the solution of the corresponding problem with convection velocity $u^{E}$:
\begin{align}\label{weqn-E-n}
\begin{cases}
\p_t \w^{E}_n + u^{E} \cdot \nabla_x \w^{E}_n = 0, \\ 
\w^{E}_n(0,x) = \w_0 \ast \varphi_n(x).
\end{cases}
\end{align}
By splitting $\|\w^{\e}_A - \w^{E}\|_{L^p_x}$ according to the equations \eqref{weqn-BE-n} and \eqref{weqn-E-n}, we have
\begin{align}\label{3-decomp}
\bega
\sup_{t\in(0,T)}\|\w^{\e}_A(t)- \w^{E}(t)\|_{L^p_x} &\leq I+II+III,
\enda
\end{align}
where
\begin{align*}
\bega
I = \sup_{t\in(0,T)}\|\w^{\e}_A(t)- \w^{\e}_{A,n}(t)\|_{L^p_x}, \hspace{2mm} II= \sup_{t\in(0,T)}\|\w^{\e}_{A,n}(t)- \w^{E}_n(t)\|_{L^p_x}, \hspace{2mm} III= \sup_{t\in(0,T)}\|\w^{E}_n(t)- \w^{E}(t)\|_{L^p_x}.
\enda
\end{align*}
{\bf (Estimate of $I$)} 
Since $\mathbb{P}u^{\e}$ is divergence-free, subtracting \eqref{weqn-BE-n} from \eqref{weqnP} yields
\begin{align*}
\|\w^{\e}_A(t)-\w^{\e}_{A,n}(t)\|_{L^p_x} = \|\w^{\e}_A(0)-\w^{\e}_{A,n}(0)\|_{L^p_x}, \quad \mbox{for} \quad 1\leq p <\infty. 
\end{align*}
{\bf (Estimate of $III$)}
Since both $(u^{E},\w^{E})$ and $(u^{E}_n,\w^{E}_n)$ are unique renormalized (hence Lagrangian) solutions, we have
\begin{align*}
\|\w^{E}(t)-\w^{E}_n(t)\|_{L^p_x} = \|\w^{E}(0)-\w^{E}_n(0)\|_{L^p_x}, \quad \mbox{for} \quad 1\leq p <\infty. 
\end{align*}
{\bf (Estimate of $II$)}  
By Young's convolution inequality, we obtain
\begin{align}\label{n-ini-conv}
\bega
\int_{\R^2} |\w^{\e}_{A,n}(0) - \w^{E}_n(0)| dx &= \int_{\R^2} |(\w^{\e}(0) - \w(0)) \ast \varphi_n| dx \cr
&\les \|\w^{\e}(0) - \w(0)\|_{L^1(\R^2)} \|\varphi_n\|_{L^1(\R^2)} \rightarrow 0 \quad \mbox{as} \quad \e \rightarrow 0.
\enda
\end{align}
From \eqref{weqn-BE-n} and \eqref{weqn-E-n} we note
\begin{align*}
\|\w^{\e}_{A,n}(t)\|_{L^p_x} = \|\w^{\e}_{A,n}(0)\|_{L^p_x}, \qquad \|\w^{E}_n(t)\|_{L^p_x} = \|\w^{E}_n(0)\|_{L^p_x}.
\end{align*}
Thus, there exists a subsequence such that
\begin{align*}
\w^{\e}_{A,n} \rightharpoonup \tilde{\w}_{A,n} \quad \mbox{weakly in} \quad L^r((0,T) \times \R^2), \quad 1<r<\infty. 
\end{align*}
For each fixed $n$, by the strong convergence of $u^{\e}_A$ in \eqref{mconv;0} and the decay $\sup_{0\leq t\leq T}\|u^{\e}_B(t)\|_{L^\infty_x}\to0$ in \eqref{uBto0}, the limit $\tilde{\w}_{A,n}$ satisfies \eqref{weqn-E-n}. By uniqueness, $\tilde{\w}_{A,n}=\w^{E}_n$. Moreover,
\begin{align*}
\|\w^{E}_n\|_{L^r_TL^r_x} &= \liminf_{\e \to 0} \|\w^{\e}_{A,n}\|_{L^r_TL^r_x} 
\leq \limsup_{\e \to 0} \|\w^{\e}_{A,n}\|_{L^r_TL^r_x} = T^{\frac{1}{p}} \lim_{\e \to 0} \|\w^{\e}_{A,n}(0)\|_{L^r_x} \\
&= T^{\frac{1}{p}} \|\w^{E}_n(0)\|_{L^r_x} = \|\w^{E}_n\|_{L^r_TL^r_x}.
\end{align*}
Hence, combining weak convergence with norm convergence yields
\begin{align}\label{wLptx}
\w^{\e}_{A,n} \rightarrow \w^{E}_n \quad \mbox{strongly in} \quad L^r((0,T) \times \R^2), \quad \mbox{for} \quad 1 < r < \infty.
\end{align}
Next, we claim
\begin{align}\label{wLpCw}
\w^{\e}_{A,n} \rightarrow \w^{E}_n \quad \mbox{in} \quad C(0,T;L^p_{\text{weak}}(B_R)), \quad \mbox{for} \quad 1 < p < \infty.
\end{align}
Indeed, for $\psi \in C_c^\infty(\R^2)$, multiplying $\psi(x)$ to \eqref{weqn-BE-n} and \eqref{weqn-E-n} gives
\begin{align*}
\frac{d}{dt}\int_{\R^2} (\w^{\e}_{A,n}-\w^{E}_n) \psi(x) dx &= \int_{\R^2} ( \mathbb{P}u^{\e} \w^{\e}_{A,n}  - u^{E} \w^{E}_n) \cdot\nabla_x\psi dx \cr 
&\to 0 \quad \mbox{as} \quad  \e\to0, \quad \mbox{in} \quad L^1(0,T),
\end{align*}
where we used the convergence of $u^{\e}_A$ and $u^{\e}_B$ from \eqref{mconv;0} and \eqref{uBto0}, respectively. Since $C_c^\infty(\R^2)$ is dense in $L^q(\R^2)$, this proves \eqref{wLpCw}.
Next, we claim 
\begin{align}\label{wn-normconv}
\bega
\int_{B_R} |\w^{\e}_{A,n}(t)|^r dx \rightarrow \int_{B_R} |\w^{E}_n(t)|^r dx, \quad \mbox{strongly in } \quad C[0,T], \quad \mbox{for} \quad 1 < r < \infty.
\enda
\end{align}
For $\beta(s) = |s|^{r}$ with $r > 1$, we define 
\begin{align*}
B_{\psi}^{\e}(t) := \int_{\R^2} \beta(\w^{\e}_{A,n})(t,x) \psi(x) dx, \qquad B_{\psi}(t) := \int_{\R^2} \beta(\w^{E}_n)(t,x) \psi(x) dx,
\end{align*}
where $\psi \in C_c^{\infty}(\R^2)$. 
Multiplying $\beta'(\w^{\e}_{A,n})\psi$ and $\beta'(\w^{E}_n)\psi$ to \eqref{weqn-BE-n} and \eqref{weqn-E-n}, respectively, we have
\begin{align*}
\frac{d}{dt}B_{\psi}(t) &= \int_{\R^2} \beta(\w^{E}_n(t,x)) u^{E}(t,x) \cdot \nabla_x \psi(x) dx, \\
\frac{d}{dt}B_{\psi}^{\e}(t) &= \int_{\R^2} \beta(\w^{\e}_{A,n}(t,x)) (\mathbb{P}u^{\e}(t,x)) \cdot \nabla_x \psi(x) dx.
\end{align*}
From \eqref{wLptx}, we know $\w^{\e}_{A,n} \rightarrow \w^{E}_n$ almost everywhere, up to a subsequence. By the dominated convergence theorem, since $|\beta(\cdot)| \leq C$ for $1 \leq r < \infty$, we conclude
\begin{align*}
\bega
\beta(\w^{\e}_{A,n}) \rightarrow \beta(\w^{E}_n) \quad \mbox{strongly in} \quad L^r_{loc}((0,T)\times\R^2), \quad \mbox{for any} \quad 1\leq r<\infty.
\enda
\end{align*}
Together with the convergence of $u^{\e}_A$ in \eqref{mconv;0} and the fact that $\sup_{0\leq t\leq T}\|u^{\e}_B(t)\|_{L^\infty_x} \to 0$ as $\e\to 0$ by \eqref{uBto0}, we deduce that $\frac{d}{dt}(B_{\psi}^{\e} - B_{\psi}) \to 0$ in $L^1(0,T)$.  
Moreover, by \eqref{n-ini-conv}, we have $B_{\psi}^{\e}(0) - B_{\psi}(0) \to 0$ as $\e \to 0$. Hence, by the fundamental theorem of calculus,
\begin{align*}
B_{\psi}^{\e}(t) - B_{\psi}(t) = B_{\psi}^{\e}(0) - B_{\psi}(0) + \int_0^t \frac{d}{dt} \big(B_{\psi}^{\e} - B_{\psi}\big)(s) ds,
\end{align*}
which implies the uniform convergence of $B_{\psi}^{\e}(t)$ toward $B_{\psi}(t)$:
\begin{align*}
\bega
\int_{\R^2} \beta(\w^{\e}_{A,n}(t,x)) \psi(x) dx \rightarrow \int_{\R^2} \beta(\w^{E}_n(t,x)) \psi(x) dx, \quad \mbox{uniformly in} \quad [0,T].
\enda
\end{align*}
For $\beta(s) = |s|^r$, $r>1$, and by approximating the test function $\psi_k \in C_c^{\infty}(\R^2)$ to $\psi(x) = \mathbf{1}_{[-\frac{1}{2},\frac{1}{2}]^2}$, we establish the claim \eqref{wn-normconv}. 
Combining the weak convergence \eqref{wLpCw} with the norm convergence \eqref{wn-normconv}, we conclude
\begin{align*}
\bega
\w^{\e}_{A,n} \rightarrow \w^{E}_n, \quad \mbox{strongly in } \quad C(0,T;L^r_{loc}(\R^2)), \quad \mbox{for} \quad 1 < r < \infty. 
\enda
\end{align*}
Since this convergence holds in $L^r_{loc}$, we also obtain strong convergence in $C(0,T;L^1_{loc}(\R^2))$.
Finally, to extend this convergence globally, we show that for any $\eta>0$ there exists $R>0$ such that
\begin{align}\label{claimBRc}
\bega
\sup_{t\in[0,T]}\int_{B_R^c} |\w^{\e}_{A,n}(t,x)|^r dx + \sup_{t\in[0,T]}\int_{B_R^c} |\w^{E}_n(t,x)|^r dx < \eta.
\enda
\end{align}
Let $\psi_r^R \in C_c^{\infty}(\R^2)$ be a cutoff function with $\psi_r^R=0$ for $0<|x|<r$ and $|x|>2R$, and $\psi_r^R=1$ for $2r<|x|<R$. Multiplying $r|\w^{\e}_{A,n}|^{r-1}\psi_r^R$ to \eqref{weqn-BE-n}, we get
\begin{align*}
\bega
\int_{\R^2} |\w^{\e}_{A,n}(t,x)|^r \psi_r^R(x) dx &\leq \int_{\R^2} |\w^{\e}_{A,n}(0,x)|^r \psi_r^R(x) dx + \int_{\R^2} |\w^{\e}_{A,n}(t,x)|^r |\mathbb{P}u^{\e}| |\nabla_x \psi_r^R(x)| dx.
\enda
\end{align*}
As in \eqref{bhold}, we decompose $u^{\e}_A$ into $L^1_x$ and $L^\infty_x$ parts:
\begin{align*}
\mathbb{P}u^{\e} = u^{\e}_A+u^{\e}_B = ({\bf K}\mathbf{1}_{|x|\leq1})\ast \w^{\e}_A + ({\bf K}\mathbf{1}_{|x|>1})\ast \w^{\e}_A +u^{\e}_B.
\end{align*}
Letting $R\to\infty$ and renaming $r$ by $R$, we deduce
\begin{align*}
\bega
\int_{B_R^c} |\w^{\e}_{A,n}(t,x)|^r dx &\leq \int_{B_R^c} |\w^{\e}_{A,n}(0,x)|^r  dx + \frac{C}{R}\sup_{t\in[0,T]}\|\w^{\e}_{A,n}\|_{L^\infty_x}\int_0^T\int_{\R^2} ({\bf K}\mathbf{1}_{|x|\leq1})\ast \w^{\e}_A dxdt \cr 
&+ \frac{C}{R}\big\|({\bf K}\mathbf{1}_{|x|>1})\ast \w^{\e}_A +u^{\e}_B\big\|_{L^\infty_x} \int_0^T\int_{\R^2} |\w^{\e}_{A,n}(t,x)|^r dxdt,
\enda
\end{align*}
where we used $|\nabla_x \psi_r^R(x)|\leq\frac{1}{r}$.
Using $\w^{\e}_A \in L^\infty(0,T;L^1(\R^2))$ and \eqref{uBto0}, we conclude \eqref{claimBRc}. 
Finally, we combine the estimates of $I$, $II$, and $III$ from \eqref{3-decomp}.  
Given any $\eta_0>0$, we first choose $n$ sufficiently small so that $I \leq \eta_0/4$ and $III \leq \eta_0/4$.  
Then, for this fixed $n$, there exists $\e_0>0$ such that for all $0<\e<\e_0$ one has $II \leq \eta_0/4$.  
Therefore,
\[
\sup_{t\in(0,T)}\|\w^{\e}_A(t)-\w^{E}(t)\|_{L^p_x} \leq \eta_0,
\]
which proves the strong convergence \eqref{wconv;0}.

\subsubsection{Conservation of the Energy}

In this part, we show that the solution constructed in Theorem~\ref{T.2D.global} 
conserves the $L^2$–norm of $\mathbb{P}u^{\e}$.

\hide
To establish this energy conservation, we first prove the strong convergence of $\rho^{\e}-\frac{3}{2}\ta^{\e}$ in $L^2(\R^2)$. 
\begin{lemma}\label{L.rtastrong}
Suppose the assumptions of Theorem \ref{T.C.Lp} hold.
If the initial data satisfy \eqref{rta0assume} then there exist $\rho^{E}$ and $\ta^{E}$ such that, up to a subsequence, the following convergence holds:
\begin{align}\label{rtastrong}
\bega
\rho^{\e} - \frac{3}{2}\ta^{\e} \rightarrow \rho^{E} - \frac{3}{2}\ta^{E}, \quad \mbox{strongly in } \quad L^\infty(0,T;L^2(\R^2)).
\enda
\end{align}
\end{lemma}
\begin{proof}
From the uniform boundedness $\sup_{0 \leq t \leq T}\|(\rho^{\e},\ta^{\e})(t)\|_{L^2_x}^2 \leq C$ given in \eqref{U.L2} in Lemma \ref{L.unif}, there exist $\rho^{E},\ta^{E}$ such that, up to a subsequence,
\begin{align*}
\bega
&\rho^{\e} \overset{\ast}{\rightharpoonup} \rho^{E}  \quad \mbox{weakly* in} \quad L^\infty(0,T;L^2(\R^2)), \cr 
&\ta^{\e} \overset{\ast}{\rightharpoonup} \ta^{E}  \quad \mbox{weakly* in} \quad L^\infty(0,T;L^2(\R^2)).
\enda
\end{align*}
Analogously to the decomposition of $\w^{\e}$ in Definition \ref{D.wAwB}, we decompose the Lagrangian solution of $\rho^{\e}-\frac{3}{2}\ta^{\e}$ along the trajectories generated by the divergence-free vector field $\mathbb{P}u^{\e}$:
\begin{align*}
\bega
\Big(\rho^{\e} - \frac{3}{2}\ta^{\e}\Big)(t,x)&= \mathfrak{T}^{\e}_A(t,x) + \mathfrak{T}^{\e}_B(t,x),
\enda
\end{align*}
where
\begin{align*}
\bega
\mathfrak{T}^{\e}_A(t,x):= \Big(\rho^{\e}_0 - \frac{3}{2}\ta^{\e}_0\Big)((X^{\e})^{-1}(t,x)), \qquad  \mathfrak{T}^{\e}_B(t,x):=\int_0^t \bar{\varPi}_{\mathfrak{s}}^{\e}(s,X^{\e}(s,\mathbf{a})) ds,
\enda
\end{align*}
and $X^{\e}$ denotes the flow map associated with $\mathbb{P}u^{\e}$ defined in \eqref{traju}. The term $\mathfrak{T}^{\e}_A$ satisfies the transport equation
\begin{align*}
\bega
\p_t\mathfrak{T}^{\e}_A + \mathbb{P}u^{\e}\cdot \nabla_x \mathfrak{T}^{\e}_A = 0 , \qquad \mathfrak{T}^{\e}_A(t,x)|_{t=0} = \Big(\rho^{\e}_0-\frac{3}{2}\ta^{\e}_0\Big)(x).
\enda
\end{align*}
By \eqref{U.rtaL2} in Lemma \ref{L.unif}, $\rho^{\e}-\tfrac{3}{2}\ta^{\e}$ is uniformly bounded in $L^\infty(0,T;L^2(\R^2))$.  
Moreover, the forcing part satisfies $\sup_{0\leq t\leq T}\|\mathfrak{T}^{\e}_B(t)\|_{L^\infty_x} \to 0$  as $\e \to 0$. 
Hence, by repeating the argument used in proving the convergence of $\w^{\e}_A$ (cf. \eqref{wconv;w}–\eqref{wconv;0}) with the special case $p=2$, we obtain
\begin{align*}
\bega
\mathfrak{T}^{\e}_A \rightarrow \mathfrak{T}^{\#}, \quad \mbox{strongly in } \quad C(0,T;L^2(\R^2)).
\enda
\end{align*}
In addition, $(u^{E},\mathfrak{T}^{\#})$ satisfies the transport equation in the sense of renormalized solutions.
\begin{align*}
\bega
\sup_{0\leq t\leq T}\|\mathfrak{T}^{\e}_B(t)\|_{L^\infty_x} \to 0, \quad \mbox{as} \quad \e \to0 .
\enda
\end{align*}
Finally, combining the convergence of $\mathfrak{T}^{\e}_A$ and the vanishing of $\mathfrak{T}^{\e}_B$, we conclude the strong convergence \eqref{rtastrong}.
\hide
In addition, the Boussinesq relation \eqref{incomp.2D} in Theorem \ref{T.2D.global} implies
\begin{align*}
(\rho^{\e} + \ta^{\e}) \rightarrow 0, \quad \mbox{strongly in } \quad L^{2+4\delta}(0,T;\dot{B}_{2+\frac{1}{\delta},1}^{s-\frac{1}{1+2\delta}}) \quad \mbox{as} \quad \e \to 0, \quad \mbox{for} \quad 0\leq|\al_x|\leq2 ,
\end{align*}
for $s \in (-1, \mathrm{N}-1)$ and any $0 < \delta\ll 1$.
\unhide
\end{proof}
\unhide

\hide
\begin{lemma}\label{L.E-conse} We have 
\begin{align}\label{E-decrease}
\bega
\lim_{\kappa\rightarrow 0} \bigg(k_B\|\rho^{\e}(t)\|_{L^2_x}^2+\|u^{\e}(t)\|_{L^2_x}^2+\frac{3}{2}k_B\|\ta^{\e}(t)\|_{L^2_x}^2\bigg) \leq k_B\|\rho_0\|_{L^2_x}^2+\|u_0\|_{L^2_x}^2+\frac{3}{2}k_B\|\ta_0\|_{L^2_x}^2. 
\enda
\end{align}
If we futher assume $\nabla_x(\rho_0-\frac{3}{2}\ta_0)=0$, then we have
\begin{align}\label{E-conserv}
\bega
\|u(t)\|_{L^2_x}^2=\|u_0\|_{L^2_x}^2.
\enda
\end{align}
\end{lemma}
\begin{proof}
Standard energy estimate to \eqref{locconD} gives 
\begin{align}\label{rhoL2}
\bega
\frac{1}{2}\frac{d}{dt}\|\rho\|_{L^2_x}^2 = \frac{1}{2}\int_{\Omega} (\nabla_x\cdot u) |\rho|^2 dx -\frac{1}{\e} \int_{\Omega} (\nabla_x\cdot u) \rho dx -\int_{\Omega} (\nabla_x\cdot u) |\rho|^2 dx
\enda
\end{align}
\begin{align}\label{uL2}
\bega
&\frac{1}{2}\frac{d}{dt}\|u\|_{L^2_x}^2 + \eta_0 k_B^{\frac{1}{2}} \kappa (\|\nabla_xu\|_{L^2_x}^2 + \|\nabla_x\cdot u\|_{L^2_x}^2) = \frac{1}{2}\int_{\Omega} (\nabla_x\cdot u) |u|^2 dx +\frac{k_B}{\e} \int_{\Omega} (\nabla_x\cdot u) (\rho+\ta) dx -\int_{\Omega}  \frac{k_B}{\mathrm{P}^{\e}}(\ta-\rho) \nabla_x \rho \cdot u dx + O(\e) \cr 
&= \frac{1}{2}\int_{\Omega} (\nabla_x\cdot u) |u|^2 dx +\frac{k_B}{\e} \int_{\Omega} (\nabla_x\cdot u) (\rho+\ta) dx -k_B\int_{\Omega} \rho\nabla_x(\rho+\ta)\cdot u dx +k_B\int_{\Omega} \rho(\rho+\ta) (\nabla_x\cdot u) dx -k_B\int_{\Omega} |\rho|^2 (\nabla_x\cdot u) dx +O(\e) 
\enda
\end{align}
where we used the following integration by parts:
\begin{align*}
\bega
\int_{\Omega} (\ta-\rho) \nabla_x \rho \cdot u dx &= \int_{\Omega} (\rho+\ta-2\rho) \nabla_x \rho \cdot u dx \cr 
&= -\int_{\Omega} \rho\nabla_x(\rho+\ta)\cdot u dx -\int_{\Omega} \rho(\rho+\ta) (\nabla_x\cdot u) dx +\int_{\Omega} |\rho|^2 (\nabla_x\cdot u) dx.
\enda
\end{align*}
\begin{align}\label{taL2}
\bega
\frac{1}{2}\frac{d}{dt}\|\ta\|_{L^2_x}^2 + \frac{5}{3}\eta_1 k_B^{\frac{1}{2}} \kappa \|\nabla_x\ta\|_{L^2_x}^2 = \frac{1}{2}\int_{\Omega} (\nabla_x\cdot u) |\ta|^2 dx -\frac{2}{3\e} \int_{\Omega} (\nabla_x\cdot u) \ta dx -\frac{2}{3}\int_{\Omega} (\nabla_x\cdot u) |\ta|^2 dx + O(\e)
\enda
\end{align}
Combining $k_B\eqref{rhoL2}+\eqref{uL2}+\frac{3}{2}k_B\eqref{taL2}$ gives 
\begin{align}\label{rutaL2}
\bega
&\frac{1}{2}\frac{d}{dt}\bigg(k_B\|\rho\|_{L^2_x}^2+\|u\|_{L^2_x}^2+\frac{3}{2}k_B\|\ta\|_{L^2_x}^2\bigg) +\eta_0 k_B^{\frac{1}{2}} \kappa (\|\nabla_xu\|_{L^2_x}^2 + \|\nabla_x\cdot u\|_{L^2_x}^2) + \frac{5}{2}\eta_1 k_B^{\frac{3}{2}} \kappa \|\nabla_x\ta\|_{L^2_x}^2 \cr 
&\les  \frac{1}{2}\|\nabla_x\cdot u\|_{L^\infty_x}\bigg(k_B\|\rho\|_{L^2_x}^2+\|u\|_{L^2_x}^2+\frac{3}{2}k_B\|\ta\|_{L^2_x}^2\bigg)  \cr 
&-k_B\int_{\Omega} (\nabla_x\cdot u) (|\rho|^2-\rho\ta+|\ta|^2) dx -k_B\int_{\Omega} \rho\nabla_x(\rho+\ta)\cdot u dx+ O(\e).
\enda
\end{align}
Taking time integral, we can have \eqref{E-decrease}. Now, we prove \eqref{E-conserv}. Let 
\begin{align*}
\bega
A^{\e}(t):=\bigg(k_B\|\rho^{\e}(t)\|_{L^2_x}^2+\|u^{\e}(t)\|_{L^2_x}^2+\frac{3}{2}k_B\|\ta^{\e}(t)\|_{L^2_x}^2\bigg)
\enda
\end{align*}
Using the assumption $\nabla_x(\rho_0-\frac{3}{2}\ta_0)=0$, we can have 
\begin{align*}
\bega
\lim_{\kappa\rightarrow 0}\|\nabla_x(\rho-\frac{3}{2}\ta)(t)\|_{L^\infty_x} =0 , \qquad \lim_{\kappa\rightarrow 0}\|\nabla_x(\rho+\ta)(t)\|_{L^\infty_x} =0 
\enda
\end{align*}
This implies $\nabla_x\rho=0$ and $\nabla_x\ta=0$ almost everywhere in $([0,T]\times\R^2)$. By the conservation laws of the macroscopic fields Boltzmann \eqref{conservPUTA}, we have 
\begin{align*}
\bega
\int_{\R^2} \rho^{\e}(t) dx = \int_{\R^2} \rho^{\e}_0 dx, \qquad \lim_{\e\rightarrow 0}\int_{\R^2} \rho^{\e}(t)+\ta^{\e}(t) dx = \int_{\R^2} \rho^{\e}_0+\ta^{\e}_0 dx
\enda
\end{align*}
Thus, $\rho(t,x)=\rho_0$ and $\ta(t,x)=\ta_0$ almost everywhere in $([0,T]\times\R^2)$.
We want to prove that 
\begin{align}\label{claim;A-A}
\bega
A(t)-A(0) \geq 0
\enda
\end{align}
From \eqref{rutaL2}, we have 
\begin{align}
\bega
&A^{\e}(t)-A^{\e}(0) \geq -\eta_0 k_B^{\frac{1}{2}} \kappa \int_0^t\|\nabla_xu(s)\|_{L^2_x}^2 + \|(\nabla_x\cdot u)(s)\|_{L^2_x}^2ds - \frac{5}{2}\eta_1 k_B^{\frac{3}{2}} \kappa \int_0^t \|\nabla_x\ta(s)\|_{L^2_x}^2ds +O(\log(\e))
\enda
\end{align}
Since $\nabla_x\ta=0$ a.e, we have $\lim_{\kappa\rightarrow 0}\|\nabla_x\ta(t)\|_{L^2_x}^2=0$. And $\|\nabla_xu\|_{L^2_x}^2 = \|\nabla_xu\|_{L^2_x}^2 + O(\e) = \|\w\|_{L^2_x}^2 + O(\e)$. Thus we have
\begin{align}\label{A-Alower}
\bega
&A^{\e}(t)-A^{\e}(0) \geq -\eta_0 k_B^{\frac{1}{2}} \kappa \int_0^t\|\w(s)\|_{L^2_x}^2 ds +O(\log(\e))
\enda
\end{align}
From the vorticity equation, we also have 
\begin{align}
\bega
\frac{d}{dt}\|\w\|_{L^2}^2 = -2\kappa \|\nabla_x \w\|_{L^2}^2 + O(\log(\e))
\enda
\end{align}
By the interpolation inequality \eqref{Ga-Ni}, we have 
\begin{align}
\bega
\|\w\|_{L^2} \leq \|\nabla_x\w\|_{L^2}^{\frac{2-p}{2}} \|\w\|_{L^p}^{\frac{p}{2}}, \quad \mbox{for} \quad 1<p<2
\enda
\end{align}
Combining above two things, we have 
\begin{align}
\bega
\frac{d}{dt}\|\w\|_{L^2}^2 = -2\kappa \|\w\|_{L^2}^{\frac{4}{2-p}}\|\w\|_{L^p}^{\frac{2p}{2-p}} + O(\log(\e)), \quad \mbox{for} \quad 1<p<2
\enda
\end{align}
Thus,
\begin{align}
\bega
\|\w(t)\|_{L^2}^{-\frac{2p}{2-p}} - \|\w(0)\|_{L^2}^{-\frac{2p}{2-p}} \geq  2\kappa t\frac{p}{2-p} \|\w\|_{L^p}^{-\frac{2p}{2-p}} + O(\log(\e)), \quad \mbox{for} \quad 1<p<2
\enda
\end{align}
Since $\|\w_0\|_{L^2}\rightarrow \infty$ as $\kappa\rightarrow 0$, we have 
\begin{align}\label{wL2-blow}
\bega
\|\w(t)\|_{L^2}^2 &\leq \bigg(2\kappa t\frac{p}{2-p} \bigg)^{-\frac{2-p}{p}}\|\w_0\|_{L^p}^2 + O(\log(\e))^{\frac{2-p}{p}}, \quad \mbox{for} \quad 1<p<2
\enda
\end{align}
where we used $\|\w(t)\|_{L^p}\leq \|\w_0\|_{L^p}$.
Combining \eqref{A-Alower} and \eqref{wL2-blow}, we obtain 
\begin{align}\label{A-Alower2}
\bega
A(t)-A(0) &\geq -\eta_0 k_B^{\frac{1}{2}} \kappa \int_0^t \bigg(2\kappa t\frac{p}{2-p} \bigg)^{-\frac{2-p}{p}}\|\w_0\|_{L^p}^2 + O(\log(\e))^{\frac{2-p}{p}} ds +O(\log(\e)) \cr 
&\geq - \kappa^{\frac{2p-2}{p}} \frac{2-p}{2p-2}\Big(\frac{2p}{2-p}t\Big)^{\frac{2p-2}{p}} \cr 
&\rightarrow 0 \quad \mbox{as} \quad \kappa \rightarrow 0, \quad \mbox{for} \quad 1<p<2.
\enda
\end{align}
This proves the claim \eqref{claim;A-A}. Moreover, since $\rho(t,x)=\rho_0$ and $\ta(t,x)=\ta_0$ almost everywhere in $([0,T]\times\R^2)$, we have.
\begin{align*}
\bega
\lim_{\kappa\rightarrow 0}(\|u^{\e}(t)\|_{L^2_x}^2-\|u^{\e}_0\|_{L^2_x}^2) = 0.
\enda
\end{align*}
Then the strong convergence  of $u^{\e}$ in \eqref{mwconv;0} gives the result. 
\end{proof}
\unhide

\hide
\begin{lemma}\label{L.E-conse}
Let $d=2$. In the finite velocity energy case \eqref{caseEC}, and for $1<p<\infty$, if we further assume $(\rho_0,u_0,\ta_0)\in L^2(\R^2)$ and \eqref{rta0assume}, then we have
\begin{align}\label{E-conserv0}
\bega
\lim_{\e\to0}\bigg(\|\rho^{\e}\|_{L^2_x}^2+\frac{1}{k_B}\|u^{\e}\|_{L^2_x}^2+\frac{3}{2}\|\ta^{\e}\|_{L^2_x}^2\bigg) = \lim_{\e\to0}\bigg(\|\rho^{\e}_0\|_{L^2_x}^2+\frac{1}{k_B}\|u^{\e}_0\|_{L^2_x}^2+\frac{3}{2}\|\ta^{\e}_0\|_{L^2_x}^2\bigg).
\enda
\end{align}
\end{lemma}

\begin{proof}
In the finite velocity energy case, recall the proof of the $L^2_x$ conservation in Theorem \ref{T.2D.unif} (i).  
By integrating \eqref{rutL2fin} over $[0,t]$, we obtain
\begin{align*}
\bega
&\bigg|\mathfrak{E}(t)-e^{\frac{1}{2}\int_0^t\|(\nabla_x\cdot u^{\e})(s)\|_{L^\infty_x}ds} \mathfrak{E}(0)\bigg| \cr 
&\leq \int_0^t e^{\frac{1}{2}\int_s^t\|(\nabla_x\cdot u^{\e})(\tau)\|_{L^\infty_x}d\tau} \bigg(\|\nabla_x(\rho^{\e}+\ta^{\e})(s)\|_{L^\infty_x} \frac{\|\mathrm{\Theta}^{\e}(s)-1\|_{L^2_x}}{\e} +C_2\kappa^{\frac{1}{2}}\mathcal{D}_G^{\frac{1}{2}}(s) \bigg)ds.
\enda
\end{align*}
where
\begin{align*}
\bega
\mathfrak{E}(t) := \bigg(\|\rho^{\e}(t)\|_{L^2_x}^2+\frac{1}{k_B}\|u^{\e}(t)\|_{L^2_x}^2+\frac{3}{2}\|\ta^{\e}(t)\|_{L^2_x}^2\bigg)^{\frac{1}{2}}.
\enda
\end{align*}
Since $\int_0^t\|(\nabla_x\cdot u^{\e})(s)\|_{L^\infty_x}\,ds \to 0$ as $\e\to0$, and the right-hand side also vanishes in the limit, we obtain the desired result \eqref{E-conserv0}.
\hide \\
(2) By the same way, from Lemma \ref{L.rho-ta}, we can have 
\begin{align*}
\bega
&\bigg|\Big\|\Big(\rho^{\e}-\frac{3}{2}\ta^{\e}\Big)(t)\Big\|_{L^2_x} - e^{\frac{1}{2}\int_0^t \|(\nabla_x\cdot u^{\e})(s)\|_{L^\infty_x}ds}\Big\|\Big(\rho^{\e}_0-\frac{3}{2}\ta^{\e}_0\Big)\Big\|_{L^2_x}\bigg| \les \int_0^t e^{\frac{1}{2}\int_s^t \|(\nabla_x\cdot u^{\e})(\tau)\|_{L^\infty_x}d\tau} \kappa^{\frac{1}{2}} \mathcal{D}_{G}^{\frac{1}{2}}(s) ds.
\enda
\end{align*}
Thus we get
\begin{align*}
\lim_{\e \to 0}\Big\|\Big(\rho^{\e}-\frac{3}{2}\ta^{\e}\Big)(t)\Big\|_{L^2_x}  = \lim_{\e \to 0}\Big\|\Big(\rho^{\e}_0-\frac{3}{2}\ta^{\e}_0\Big)\Big\|_{L^2_x} = 0.
\end{align*}
In addition, the Boussinesq relation \eqref{incomp.2D} in Theorem \ref{T.2D.global} implies
\begin{align*}
\p^{\al_x}(\rho^{\e} + \ta^{\e}) \rightarrow 0, \quad \mbox{almost everywhere in } \quad ((0,T)\times\R^2) \quad \mbox{as} \quad \e \to 0, \quad \mbox{for} \quad 0\leq|\al_x|\leq2 
\end{align*}
Since, $\rho^{\e}$ and $\ta^{\e}$ are linear combination of $\rho^{\e}-\frac{3}{2}\ta^{\e}$ and $\rho^{\e} + \ta^{\e}$, we have 
\begin{align}\label{rtaL2=0}
\lim_{\e \to 0}\|\rho^{\e}(t)\|_{L^2_x}= 0, \qquad \lim_{\e \to 0}\|\ta^{\e}(t)\|_{L^2_x}= 0.
\end{align}
Combining \eqref{rtaL2=0} and the initial data $\lim_{\e\to0}\|(\rho_0^{\e},\ta_0^{\e})\|_{L^2_x}^2=0$ with \eqref{E-conserv0} gives the result \eqref{E-conserv}.
\unhide
\end{proof}
\unhide

\begin{proof}[{\bf Proof of Proposition \ref{P.Econs}}]
(1) The equality~\eqref{E-conservlim} follows directly from the conservation of the $L^2$–norm of $\mathbb{P}u^{\e}$ established in~\eqref{U.PuL2} of 
Lemma \ref{L.unif}:
\begin{align*}
\bega
\lim_{\e\to0} \|\mathbb{P}u^{\e}\|_{L^2_x}^2  =  \lim_{\e\to0} \|\mathbb{P}u^{\e}_0\|_{L^2_x}^2 = \|u_0\|_{L^2_x}^2.
\enda
\end{align*}
(2) To prove~\eqref{E-conserv0}, we first observe that 
\begin{align}\label{Pu-u=u}
\bega
\Big| \lim_{\e\to0} \Big(\|\mathbb{P}u^{\e}\|_{L^2_x}^2-\|u^{\e}_A\|_{L^2_x}^2 \Big) \Big|  \leq  \lim_{\e\to0} \|u^{\e}_B\|_{L^2_x}(2\|\mathbb{P}u^{\e}\|_{L^2_x}+\|u^{\e}_B\|_{L^2_x}) = 0, 
\enda
\end{align}
since $\mathbb{P}u^{\e}=u^{\e}_A+u^{\e}_B$ by Lemma~\ref{L.uABprop} (i), 
and $\sup_{0\leq t\leq T_*}\|u^{\e}_B\|_{L^2_x} \to 0$ by~\eqref{uBto0}.  
Therefore, it suffices to prove the strong convergence of $u^{\e}_A$ in \eqref{ustrong;0}. 
Following the argument of \cite{Econ,Crippa2021}, We claim that $\{u^{\e}_A\}_{\e>0}$ is a Cauchy sequence in $L^2(\R^2)$, 
that is, for any $\eta>0$, there exist $\e_1,\e_2 >0$ such that
\begin{align}\label{uAcauchy0}
\bega
\|u^{\e_1}_A-u^{\e_2}_A\|_{L^2(\R^2)} \leq \eta.
\enda
\end{align}
We decompose the kernel ${\bf K} = {\bf K}_1 + {\bf K}_2$ smoothly 
using a cutoff function $\chi \in C^{\infty}([0,\infty))$ satisfying 
$\chi(|x|)=1$ for $|x|<1$, $\chi(|x|)\in(0,1)$ for $1\leq |x|<2$, 
and $\chi(|x|)=0$ for $|x|\geq 2$:
\begin{align*}
{\bf K}_1(x) :=\chi(\delta|x|)\frac{x^{\perp}}{|x|^2}, \qquad {\bf K}_2(x) := (1-\chi(\delta|x|))\frac{x^{\perp}}{|x|^2}.
\end{align*}
Then, using the equation for $\w^{\e}_A$ in~\eqref{weqnA}, we write
\begin{align*}
\bega
\p_t u^{\e}_A &= {\bf K}\ast(\p_t \w^{\e}_A)= {\bf K}_1\ast(\p_t \w^{\e}_A) + {\bf K}_2\ast \big( - \mathbb{P}u^{\e}\cdot\nabla_x\w^{\e}_A\big ).
\enda
\end{align*}
Integrating in time for both sequences $\p_tu^{\e_1}_A$ and $\p_tu^{\e_2}_A$, 
we estimate 
\begin{align*}
\bega
\|u^{\e_1}_A(t)-u^{\e_2}_A(t)\|_{L^2} &\leq I+II+III+IV,
\enda
\end{align*}
where
\begin{align*}
\bega
I&:= \|u^{\e_1}_A(0)-u^{\e_2}_A(0)\|_{L^2}, \cr 
II&:= \big\|{\bf K}_1\ast\big(\w^{\e_1}_A(t)-\w^{\e_2}_A(t)\big)\big\|_{L^2} + \big\|{\bf K}_1\ast\big(\w^{\e_1}_A(0)-\w^{\e_2}_A(0)\big)\big\|_{L^2}, \cr 
III&:= \bigg\| \int_0^t {\bf K}_2\ast \big(\nabla_x^{\perp}\nabla_x\cdot(u^{\e_1}_A\otimes u^{\e_1}_A-u^{\e_2}_A\otimes u^{\e_2}_A)\big) ds \bigg\|_{L^2}, \cr 
IV&:= \bigg\| \int_0^t {\bf K}_2\ast \big(\nabla_x\cdot(u^{\e_1}_B \w^{\e_1}_A-u^{\e_2}_B \w^{\e_2}_A)\big) ds \bigg\|_{L^2}.
\enda
\end{align*}
Here, we used $\mathbb{P}u^{\e}= u^{\e}_A+u^{\e}_B$ 
(from Lemma~\ref{L.uABprop}) and the identity 
$u^{\e}_A\cdot\nabla_x \w^{\e}_A 
= \nabla_x^{\perp}\div(u^{\e}_A\otimes u^{\e}_A)$. 
For $I$, since the initial data converge strongly in $L^2(\R^2)$ by assumption~\eqref{Pu0assume}, they form a Cauchy sequence; hence $I \leq \eta$. 
For $II$, applying Young’s convolution inequality gives
\begin{align*}
\bega
\big\|{\bf K}_1\ast\big(\w^{\e_1}_A(t)-\w^{\e_2}_A(t)\big)\big\|_{L^2} \leq 
\begin{cases}
\|{\bf K}_1\|_{L^q} \|\w^{\e_1}_A(t)-\w^{\e_2}_A(t)\|_{L^p}, \quad \mbox{for} \quad 1\leq p<2 , \\ 
\|{\bf K}_1\|_{L^1} \|\w^{\e_1}_A(t)-\w^{\e_2}_A(t)\|_{L^2}, \quad \mbox{for} \quad  p \geq 2 ,
\end{cases}
\enda
\end{align*}
where $\frac{1}{p} + \frac{1}{q} = 1 + \frac{1}{2}$.  
Since ${\bf K}_1 \in L^q(\R^2)$ for $1 \leq q < 2$, 
and $\w^{\e}_A$ converges strongly by~\eqref{wconv;0}, we have $II \leq \eta$.
For $III$, we integrate by parts and then apply Minkowski’s inequality to obtain
\begin{align*}
\bega
III&\leq 
\int_0^t \|\nabla_x \nabla_x^{\perp} {\bf K}_2\|_{L^2} \big\| \big(u^{\e_1}_A\otimes u^{\e_1}_A-u^{\e_2}_A\otimes u^{\e_2}_A\big)(s) \big\|_{L^1} ds \cr 
&\leq C \int_0^t \big(\|u^{\e_1}_A(s)\|_{L^2}+\|u^{\e_2}_A(s)\|_{L^2}\big) \big\| (u^{\e_1}_A-u^{\e_2}_A)(s) \big\|_{L^2} ds,
\enda
\end{align*}
where we used $\nabla_x \nabla_x^{\perp} {\bf K}_2 \in L^2(\R^2)$ 
and inserted and subtracted $u^{\e_1}_A \otimes u^{\e_2}_A$.  
Since $u^{\e}_A$ is uniformly bounded in $L^2(\R^2)$ by~\eqref{Pu-u=u}, we deduce
\begin{align*}
\bega
III &\leq C \int_0^t \big\| (u^{\e_1}_A-u^{\e_2}_A)(s) \big\|_{L^2} ds.
\enda
\end{align*}
For $IV$, integrating by parts and applying Minkowski’s inequality gives
\begin{align*}
\bega
IV\leq 
\int_0^t \|\nabla_x{\bf K}_2\|_{L^2} \big(\|u^{\e_1}_B\|_{L^\infty} \|\w^{\e_1}_A\|_{L^1} + \|u^{\e_2}_B\|_{L^\infty} \|\w^{\e_2}_A\|_{L^1}\big) ds \to 0,
\enda
\end{align*}
where we used~\eqref{uBto0}, the fact that $\nabla_x {\bf K}_2 \in L^2(\R^2)$, 
and the uniform $L^1$–boundedness of $\w^{\e}_A$.  
Combining the above estimates for $I$–$IV$, we obtain
\begin{align*}
\bega
\|u^{\e_1}_A(t)-u^{\e_2}_A(t)\|_{L^2} &\leq C\eta + C\int_0^t \|u^{\e_1}_A(s)-u^{\e_2}_A(s)\|_{L^2} ds.
\enda
\end{align*}
Applying Gr\"{o}nwall inequality, we get 
\begin{align*}
\bega
\|u^{\e_1}_A(t)-u^{\e_2}_A(t)\|_{L^2} &\leq C\eta.
\enda
\end{align*}
This proves the claim \eqref{uAcauchy0}. 
\end{proof}

\subsection{Radon measure vorticity with distinguished sign (Delort solutions)}\label{sec:RM}

The existence of the solution for vortex-sheet initial data in the incompressible Euler equations has long been a central problem. While many works have addressed this issue, DiPerna and Majda \cite{DiMa} showed that concentration may occur in approximate-solution sequences. Subsequently, when the initial vorticity has a distinguished sign with $\w_0 \in \mathcal{M}(\R^2) \cap H^{-1}_{\mathrm{loc}}(\R^2)$, Delort in \cite{Delort} established the existence of weak solutions to the Euler equations.

In this section, we show that if the initial vorticity belongs to the space of Radon measures with a distinguished sign, then the momentum obtained from the Boltzmann solution in Theorem~\ref{T.2D.global} converges strongly, and its limit is a weak solution of the incompressible Euler equation. It is worth noting that, unlike for the incompressible Euler equation, the vorticity generated by the Boltzmann equation does not preserve its sign. Moreover, the microscopic part of the vorticity equation in~\eqref{weqnnew} interferes with the estimate of $\int_{\R^2} |x|^2 \w^{\e},dx$ because of the factor $\frac{1}{\mathrm{P}^{\e}}$ infront of microscopic part. 
To circumvent this difficulty, we redefine the vorticity in this section: instead of using $\w^{\e}:=\nabla_x^{\perp} u^{\e}$, we take the curl of the momentum as the vorticity.

\begin{definition}\label{D.mdef}
We define the momentum $m^{\e}$ and its corresponding m-vorticity as follows:
\begin{equation*}
m^{\e}(t,x):= \frac{1}{\e}\mathrm{P}^{\e}\mathrm{U}^{\e}, \qquad \bw^{\e}(t,x) := \nabla_x^\perp \cdot m^{\e}(t,x) = -\p_2 m^{\e}_1(t,x) + \p_1 m^{\e}_2(t,x).
\end{equation*}
\end{definition}
Interestingly, even though the equation contains microscopic terms and $m^{\e}$ is not divergence-free, the vorticity equation for $\bw^{\e}$ still preserves the quantities $\int_{\R^2} \bw^{\e}\,dx$ and $\int_{\R^2} |x|^2 \bw^{\e}\,dx$ (see Lemma~\ref{L.intw}). 
Afterwards, if we write the Lagrangian form of the vorticity equation and decompose it into the initial data part $\bw^{\e}_A$ and the forcing part $\bw^{\e}_B$ in the same manner as in~\eqref{D.wAwB}, we can observe that $\bw^{\e}_A$ preserves its sign. Moreover, the forcing part arising from the maximal vorticity function converges to zero.

\hide
Then the mollified vorticity is uniformly bounded in $L^1_x$, $\sup_{\e>0}\|\bw^{\e}_0\|_{L^1_x}\leq C$, so that from Theorem ...(\textcolor{red}{Need to make the theorem. construction}), we can have the solution satisfying $\|\bw^{\e}(t)\|_{L^1_x}\leq C$.
Then, we define $u_0^{\e}$ from the Biot-Savart law plus divergence part: 
\begin{align}\label{u0e}
\bega
u_0^{\e} := K \ast \bw_0^{\e} + \psi_{\e}\ast(\mathbb{P}^{\perp}u_0). 
\enda
\end{align}
(See page 434 of the book) Mollification of the vorticity satisfies 
\begin{align}
\bega
\bw_0^{\e} \geq 0, \qquad \mbox{supp}~ \bw_0 \subset \{x~|~|x|<CR\}, \cr 
\int_{|x|\leq R}|\mathbb{P}u_0^{\e}|^2 dx + \int_{\R^2}\bw_0^{\e}dx \leq C(R^2).
\enda
\end{align}
Note that $m^{\e}$ is close to $u^{\e}$: $|m^{\e}-u^{\e}| \leq |\mathrm{P}^{\e}-1||u^{\e}|$, because $|\mathrm{P}^{\e}-1|$ gives one epsilon. Moreover, $m^{\e}$ and $\bw^{\e}$ satisfies the following equations: (by the same way to the proof of Lemma \ref{meqnlem})
Now, since it is just constants multiplied in front of the term $\nabla_x^{\perp}\p_{x_j} \mathbf{r}_{ij}^{\e}$, by the integration by parts, we can prove almost conservation of the following quantities: 
\begin{align}\label{pseudodef}
\bega
\int_{\R^2}\bw^{\e}(t,x)dx, \qquad \int_{\R^2}|x|^2\bw^{\e}(t,x)dx, \qquad H(\bw^{\e}(t))= \frac{1}{2\pi} \int_{\R^2}\int_{\R^2} \log|x-y|\bw^{\e}(x)\bw^{\e}(y) dydx.
\enda
\end{align}
By using such conservation, we can prove that 
\begin{align*}
\bega
\int_{\R^2}\int_{\R^2} \big|\log^-(|x-y|)\big| \bw^{\e}(y)\bw^{\e}(x) dydx \leq C
\enda
\end{align*}
Howevere, we cannot have the decay estimate of the vorticity maximal function since $\bw^{\e}$ may not preserve it's sign because of the forcing term in \eqref{mGeqn}. So that, along the trajectory, we decompose the solution of the vorticity equation to initial data part and forcing term part. 
\unhide

Before stating the main theorem of this section, we introduce several necessary definitions and equations.

\begin{lemma} The momentum $m^{\e}$ and m-vorticity $\bw^{\e}$ defined in~\eqref{D.mdef} satisfy the following equations:
\begin{align}\label{mGeqn}
&\p_tm^{\e} +m^{\e}\cdot\nabla_x m^{\e} +\frac{k_B}{\e^2}\nabla_x(\mathrm{P}^{\e} \mathrm{\Theta}^{\e})
+\frac{1}{\e^2}\nabla_x \cdot \mathbf{r}^{\e} = -u^{\e} \nabla_x \cdot m^{\e} + m^{\e}\cdot\nabla_x ((\mathrm{P}^{\e}-1) u^{\e}), 
\end{align}
\begin{align}\label{weqn-m}
\bega
&\p_t\bw^{\e} + \mathbb{P}m^{\e}\cdot \nabla_x \bw^{\e} = \varPi_{\bw}^{\e},
\enda
\end{align}
where
\begin{align*}
\bega
\varPi_{\bw}^{\e}(t,x) :=& - \mathbb{P}^{\perp}m^{\e}\cdot \nabla_x \bw^{\e} - (\nabla_x\cdot m^{\e})\bw^{\e} -\frac{1}{\e^2}\nabla_x^{\perp}\cdot (\nabla_x\cdot \mathbf{r}^{\e}) \cr 
 &-\nabla_x^{\perp}\cdot (u^{\e} \nabla_x \cdot m^{\e}) + \nabla_x^{\perp}\cdot(m^{\e}\cdot\nabla_x ((\mathrm{P}^{\e}-1) u^{\e})).
\enda
\end{align*}
\end{lemma}
\begin{proof}
Dividing the equation~$\eqref{loccon}_2$ by $\e^2$, we get
\begin{align}\label{mweqn}
\bega
& \p_t m^{\e} + m^{\e}\cdot\nabla_x u^{\e} + u^{\e}(\nabla_x\cdot m^{\e}) +\frac{k_B}{\e^2}\nabla_x(\mathrm{P}^{\e} \mathrm{\Theta}^{\e}) + \frac{1}{\e^2}\nabla_x \cdot \mathbf{r}^{\e} =0 , \cr 
& \p_t \bw^{\e} + \nabla_x^{\perp}\cdot\big(m^{\e}\cdot\nabla_x u^{\e} + u^{\e}(\nabla_x\cdot m^{\e})\big) + \frac{1}{\e^2} \nabla_x^{\perp}\cdot (\nabla_x\cdot \mathbf{r}^{\e}) =0 .
\enda
\end{align} A straightforward algebraic computation yields the desired result.
\end{proof}

\begin{remark}
Unlike the microscopic term $\nabla_x^{\perp}\cdot\big(\frac{1}{\mathrm{P}^{\e}}\p_{x_j} \mathbf{r}_{ij}^{\e}\big)$ in the equation for $\w^{\e}$ in~\eqref{weqnnew}, the microscopic term in the equation for $\bw^{\e}$ in~\eqref{weqn-m} is $\nabla_x^{\perp}\cdot\p_{x_j} \mathbf{r}_{ij}^{\e}$. 
The presence of these two derivatives, without the factor $1/\mathrm{P}^{\e}$, makes it possible to conserve the quantity $\int_{\R^2}|x|^2 \bw^{\e}\,dx$.
\end{remark}

As in Definition~\ref{D.wAwB}, we decompose $\bw^{\e}$ into its initial data part and its forcing part.

\begin{definition}\label{D.wAwBm}
For each $\e>0$, we define $\bw^{\e}_A$ and $\bw^{\e}_B$ as the solutions of the following equations:
\begin{align}
&\p_t\bw^{\e}_A + \mathbb{P}m^{\e}\cdot \nabla_x \bw^{\e}_A = 0 , \qquad \hspace{3mm} \bw^{\e}_A(t,x)|_{t=0} = \bw^{\e}_0(x), \label{weqnmA} \\ 
&\p_t\bw^{\e}_B + \mathbb{P}m^{\e}\cdot \nabla_x \bw^{\e}_B = \bar{\varPi}_{\bw}^{\e} , \qquad \bw^{\e}_B(t,x)|_{t=0} = 0 . \label{weqnmB}
\end{align}
From \(\bw^{\e}_A(t,x)\) and \(\bw^{\e}_B(t,x)\), we define \(m^{\e}_A\) and \(m^{\e}_B\) using the Biot–Savart law:
\begin{align*}
\bega
m^{\e}_A(t,x) := {\bf K} \ast \bw^{\e}_A(t,x), \qquad m^{\e}_B(t,x) := {\bf K} \ast \bw^{\e}_B(t,x),
\enda
\end{align*}
where ${\bf K}(x):= \frac{1}{2\pi}\frac{x^{\perp}}{|x|^2}$. \\
We define the pseudo–energy by
\begin{align}\label{pseudodef}
\bega
H(\bw(t))= \frac{1}{2\pi} \int_{\R^2}\int_{\R^2} \log|x-y|\bw(x)\bw(y) dydx.
\enda
\end{align}
\end{definition}

\begin{lemma}
Let $\bw^{\e}_A$ and $\bw^{\e}_B$ be the solutions of \eqref{weqnmA} and \eqref{weqnmB}, respectively. Then $\bw^{\e} = \bw^{\e}_A + \bw^{\e}_B$ is a solution of \eqref{weqn-m}, and it can be represented as
\begin{align*}
\bega
\bw^{\e}_A(t,x):= \bw^{\e}_0((X^{\e})^{-1}(t,x)), \qquad  \bw^{\e}_B(t,x):=\int_0^t \varPi_{\bw}^{\e}(s,X^{\e}(s,\mathbf{a})) ds.
\enda
\end{align*}
Here, \(X^{\e}\) denotes the flow map associated with the divergence-free vector field \(\mathbb{P}m^{\e}\):
\begin{align*}
\bega
\frac{d}{dt}X^{\e}(t,\mathbf{a}) &= \mathbb{P}m^{\e}(t,X^{\e}(t,\mathbf{a})), \qquad  X^{\e}(0,\mathbf{a}) = \mathbf{a}.
\enda
\end{align*}
\end{lemma}

\begin{definition}\label{D.weaksol} 
A function $u(t,x)$ is called a \emph{weak solution} of the  Euler equation in primitive variable form if the following conditions hold:
\begin{itemize}
\item $u\in L^1([0,T]\times B_R)$ and $u\otimes u \in L^1([0,T]\times B_R)$. 
\item $u\in Lip([0,T];H^{-s}_{loc}(\R^2))$ and $u_0 \in H^{-s}_{loc}(\R^2)$ for some $s>0$.
\item $\div(u) = 0$ in the sense of distributions.
\item For every test function $\Lambda_1 \in C_c^{\infty}(\R^+\times \R^2;\R^2)$ with $\nabla_x \cdot \Lambda_1 = 0$, 
\begin{align*}
\int_0^T \int_{\R^2} (\p_t\Lambda_1 \cdot u + \nabla_x \Lambda_1 : u\otimes u ) dxdt = 0,
\end{align*}
where $X:Y$ denotes $\sum_{i,j}X_{ij}Y_{ij}$. 
\end{itemize}
A function $\varrho(t,x)$ is called a \emph{weak solution} of the 
continuity equation with velocity field $u$ satisfying the above 
conditions if, for every test function $\Lambda_2 \in C_c^{\infty}(\R^+\times \R^2)$
\begin{align*}
\int_0^T \int_{\R^2} 
\big( \partial_t \Lambda_2 + u\cdot\nabla_x \Lambda_2 \big)\,\varrho
\,dxdt = 0.
\end{align*}
\end{definition}

\begin{definition}
We use the following notation in this section:
\begin{itemize}
\item $\mathcal{M}(\R^2)$ denotes the space of Radon measures on $\R^2$.
\item Let $\{\bw_n\}_{n\ge1}$ and $\bw$ be Radon measures on $\R^2$. 
We say that $\bw_n$ converges to $\bw$ in the weak-$*$ sense, 
denoted by $\bw_n \overset{*}{\rightharpoonup} \bw$ weakly* in $\mathcal{M}(\R^2)$, if
\[
\int_{\R^2} \varphi \, d\bw_n
\;\longrightarrow\;
\int_{\R^2} \varphi \, d\bw
\qquad\text{for all } \varphi \in C_c(\R^2).
\]
\end{itemize}
\end{definition}

\begin{theorem}\label{T.Radon}
Let $\Omega=\R^2$.
Consider a family of initial data $\{F_0^{\e}\}_{\e>0}$ satisfies \eqref{L2unif}, the $4$-ABC~\eqref{ABC1}, and 
\begin{align*}
\sup_{\e>0}
\|\bw_0^{\e}\|_{L^1\cap H^{-1}_{loc}(\R^2)}
<+\infty.
\end{align*}
\hide and which additionally satisfy the $4$-Admissible Blow-up Condition~\eqref{ABC1}.
Here $\tilde{\mu}=M_{[1,0,1-c_0]}$ for some $0<c_0\ll1$, and $\bar{u}$
denotes the radial eddy defined in Definition~\ref{D.Ra-E}.
Then there exists a family of strong Boltzmann solutions $\{F^{\e}\}_{\e>0}$ to~\eqref{BE} with $\kappa=\e^q$ for some $q<2$, and $F^{\e}\big|_{t=0}=F^{\e}_0$, and defined on
the time interval $t\in[0,T]$, where $T>0$ can be chosen arbitrarily large (see our construction in Theorem~\ref{T.2D.global}).\unhide 
Furthermore, we assume that, as $\e \to 0$, \eqref{conv_s_initial} and 
\begin{align*}
\bw^{\e}_0 \quad &\overset{\ast}{\rightharpoonup} \quad \bw_0 \in \mathcal{M}(\R^2)\cap H^{-1}_{loc}(\R^2) \qquad \text{weakly* in $\mathcal{M}(\R^2)$}, 
\end{align*}
where 
\begin{align*}
\bega
&
\quad \bw_0,\bw^{\e}_0 \geq 0  , \quad \mbox{supp}~ \bw_0, \mbox{supp}~ \bw^{\e}_0 \subset \{x~|~|x|<R\}, \quad u_0 \in L^2_{\mathrm{loc}}(\R^2), \cr 
&
\int_{\R^2} |x|^2\bw^{\e}_0 dx \leq C\int_{\R^2} |x|^2\bw_0 dx< \infty, \quad  H(\bw^{\e}_0) \leq CH(\bw_0) < \infty.
\enda
\end{align*}
Here, $H(\bw)$ is defined in \eqref{pseudodef}.

Then the following statements hold for a family of Boltzmann solutions $\{F^{\e}\}_{\e>0}$ to~\eqref{BE} with $\kappa=\e^q$ for some $0<q<2$, on
the time interval $t\in[0,T_\e]$ with $T_\e\to\infty$ as $\e\to0$, which is constructed in Theorem~\ref{T.2D.global}.


\begin{enumerate}
\item[(1)] The momentum and its vorticity satisfy the following convergence up to a subsequence (not relabeled):
\begin{align*}
\bega
m^{\e}_A \quad &\rightarrow \quad m^{\#} \qquad \text{strongly in $L^r([0,T]\times B_R(0))$}, \quad \mbox{for} \quad 1\leq r<2, ~ \text{ for } R>0, \cr 
\bw_A^{\e} \quad &\overset{\ast}{\rightharpoonup} \quad \bw^{\#} \qquad \text{weakly* in $L^\infty(0,T;\mathcal{M}(\R^2))$}.
\enda
\end{align*}
Here $m^{\#}$ is a weak solution of the Euler equations in primitive variable form in the sense of Definition~\ref{D.weaksol}. Moreover,
$m^{\#} \in L^\infty(0,T;L^2_{\mathrm{loc}}(\R^2))$ and 
$\bw^{\#} \in L^\infty(0,T;\mathcal{M}(\R^2))$.

\item[(2)] The forcing components $m^{\e}_B$ and $\bw^{\e}_B$ vanish in the limit:
\begin{align}\label{wBto0Ra}
&m^{\e}_B, \ \bw^{\e}_B   \rightarrow   0 \ \ \text{strongly in $L^\infty(0,T; L^p(\R^2))$}, \quad \mbox{for} \quad  2\leq p \leq \infty,
\end{align}
with rate $\e^{\frac{1}{4}-} +\kappa^{\frac{1}{2}-}$.

\item[(3)] The entropic fluctuation $\mathfrak{s}^{\e}:=3/2\ta^{\e}-\rho^{\e}$ converges as \eqref{TC.Lp.ent}. \hide follows:
\begin{align}\label{TC.Ra.ent}
\frac{3}{2} \theta^\e - \rho^\e  \to \frac{3}{2} \theta^E - \rho^E  \ \ &\text{in  $L^\infty (0, T; L^2(\R^2))$ }. 
\end{align}
Here, $\frac{3}{2}\theta^E - \rho^E$ is a weak solution to the
continuity equation with the divergence--free velocity field $u^E$.\unhide
\item[(4)] The irrotational part $\mathbb{P}^{\perp}u^{\e}$ and pressure fluctuation $\rho^{\e}+\ta^{\e}$ vanish as \eqref{TC.Yudo.irrorhota}.

\item[(5)] The microscopic part converges to zero, as quantified in \eqref{TC.Yudo.G}.
\end{enumerate}
\end{theorem}

\hide
\begin{remark}
The compressible components converge as follows:
\begin{align}\label{TC.Ra.irrorhota}
\mathbb{P}^{\perp}u^{\e}, \ \rho^{\e}+\ta^{\e}    \to  0  \ \ &\text{in $ L^r (0,T; \dot{B}_{p,1}^{s+2(\frac{1}{p}-\frac{1}{2})+\frac{1}{r}}(\R^2))$ with rate}.
\end{align}
This holds for $2\leq p\leq\infty$, $\frac{1}{r}\leq 2(\frac{1}{2}-\frac{1}{p})$, and any $s\in[0,\mathrm{N}-1)$.
\end{remark}
\unhide

\hide
\begin{remark}
We summarize the quantitative convergence rates in the Big-O sense as follows:
\begin{align*}
\e^{\frac{1}{4}-} +\kappa^{\frac{1}{2}-}
\ \ &\text{in \eqref{wBto0Ra}},\\
\e^{\frac{1}{4}-}
\ \ &\text{in \eqref{TC.Ra.irrorhota}}.
\end{align*}
\end{remark}
\unhide

\begin{corollary}
Under the assumption \eqref{initialL2}, we replace the condition $\w_0\in L^p\cap L^1 \cap H^{-1}_{loc}(\R^2)$ by
$\bw_0 \in \mathcal{M}(\R^2)\cap H^{-1}_{loc}(\R^2)$ with $\bw_0\geq 0$, $\mbox{supp}~ \bw_0 \subset \{x~|~|x|<R\}$, $\int_{\R^2} |x|^2\bw_0 dx< \infty$ and $H(\bw_0) < \infty$.
Then the mollified sequence $\{F^{\e}_0\}_{\e>0}$ constructed from Lemma \ref{L.molli} satisfies the assumption of Theorem \ref{T.Radon}, and hence all the results of Theorem \ref{T.Radon} hold.
\end{corollary}


In Theorem~\ref{T.2D.global}, the solution satisfying the Admissible Blow-up Condition~\eqref{ABC1} also satisfies \eqref{incomp.2D} and 
\eqref{Gto0}, which imply the convergence of the microscopic part 
and the compressible part, respectively.
Since the thermodynamic variable 
$\frac{3}{2}\ta^{\e}-\rho^{\e}$~\eqref{rtaeqnP} satisfies 
a continuity equation with velocity field $\mathbb{P}u^{\e}$, 
it suffices to focus on the convergence of the momentum.

Before proving Theorem~\ref{T.Radon}, we first examine the properties 
of $(m^{\e}_A, m^{\e}_B)$ and $(\bw^{\e}_A, \bw^{\e}_B)$, as well as 
the rotational invariance of the Boltzmann equation and its momentum equation.

\begin{lemma}\label{L.mABprop}
For $(m^{\e}_A, m^{\e}_B)$ and $(\bw^{\e}_A, \bw^{\e}_B)$ defined in Definitions~\ref{D.wAwBm}, the following properties hold:
\begin{enumerate}

\item Both $m^{\e}_A$ and $m^{\e}_B$ are incompressible vector fields.

\item The Leray projection of $m^{\e}$ decomposes as
$\mathbb{P}m^{\e}(t,x) = m^{\e}_A(t,x) + m^{\e}_B(t,x).$

\item The corresponding velocity fields $m^{\e}_A$ and $m^{\e}_B$ satisfy
\begin{align}\label{mAeqn}
\bega
\begin{cases}
\p_tm^{\e}_A + m^{\e}_A \cdot \nabla_x m^{\e}_A + \nabla_x p^{\e}_A = -\nabla_x^{\perp}(-\Delta_x)^{-1}\nabla_x \cdot \big(m^{\e}_B\bw^{\e}_A\big), \\
m^{\e}_A(t,x)|_{t=0} = {\bf K} \ast \bw^{\e}_0 , \end{cases}
\enda
\end{align}
and
\begin{align*}
\bega
\begin{cases} 
\p_tm^{\e}_B + m^{\e}_B \cdot \nabla_x m^{\e}_B + \nabla_x p^{\e}_B = -\nabla_x^{\perp}(-\Delta_x)^{-1}\nabla_x \cdot(m^{\e}_A\bw^{\e}_B)+\nabla_x^{\perp}(-\Delta_x)^{-1}\bar{\varPi}_{\bw}^{\e}, \\
m^{\e}_B(t,x)|_{t=0} = 0 , \end{cases}
\enda
\end{align*}
where $\Delta_x p^{\e}_A = \sum_{i,j} \p_i \p_j (m^{\e}_{A,i} m^{\e}_{A,j})$ and 
$\Delta_x p^{\e}_B = \sum_{i,j} \p_i \p_j (m^{\e}_{B,i} m^{\e}_{B,j})$.

\item $\bw^{\e}_A$ preserves the total circulation, whereas $\bw^{\e}_B$ has zero circulation:
\begin{align*}
\int_{\R^2} \bw^{\e}_A dx = \int_{\R^2} \bw^{\e}_0 dx, \qquad \int_{\R^2} \bw^{\e}_B dx = 0.
\end{align*}

\item If $\bw_0^{\e} \geq 0$ (resp., $\bw_0^{\e} \leq 0$), then $\bw^{\e}_A(t,x) \geq 0$ (resp., $\bw^{\e}_A(t,x) \leq 0$). 

\item The total vorticity of $\bw^{\e}_A$ is conserved:
\begin{align}\label{wAL1}
\bega
\frac{d}{dt}\int_{\R^2}\bw^{\e}_A(t,x)dx &= 0, \qquad \|\bw^{\e}_A(t)\|_{L^1_x} \leq C. 
\enda
\end{align}

\item The $L^2_x$ norm of $\bw_A^{\e}$ is singular at the level of the Admissible blow-up rate in \eqref{ABC1}:
\begin{align}\label{wAL2}
\bega
\frac{d}{dt}\int_{\R^2}|\bw^{\e}_A(t,x)|^2dx =0, \qquad \|\bw^{\e}_A(t)\|_{L^2_x} \leq  C\log\!\Big(\log\!\big(\log(1/\e)\big)\Big). 
\enda
\end{align}

\item The forcing component $m^{\e}_B$ vanishes in the limit:
\begin{align}\label{mBto0}
\sup_{0\leq t\leq T_*}\|m^{\e}_B(t)\|_{L^p_x} \les \kappa^{\frac{1}{2}-} \to 0, \quad \mbox{as} \quad \e\to 0, \quad 2\leq p \leq \infty.
\end{align}

\end{enumerate}
\end{lemma}

\begin{proof}
Most of the proof is identical to that of Lemma~\ref{L.uABprop}; we therefore only present the parts that differ from Lemma~\ref{L.uABprop}. \\
(v) This follows directly from the definition of $\bw^{\e}_A(t,x)$ in Definition~\ref{D.wAwBm}. \\
(vi) From~\eqref{weqnmA}, integration by parts yields
\begin{align*}
\bega
\frac{d}{dt}\int_{\R^2}\bw^{\e}_A(t,x)dx &= -\int_{\R^2} \mathbb{P}m^{\e}\cdot \nabla_x \bw^{\e}_A dx = \int_{\R^2} (\nabla_x \cdot \mathbb{P}m^{\e}) \bw^{\e}_A dx =0 .
\enda
\end{align*}
Since $\bw^{\e}_A(0)=\bw^{\e}_0$ by~\eqref{weqnmA} and $\bw^{\e}_0$ is uniformly bounded in $L^1_x$, we have
\begin{align*}
\bega
\|\bw^{\e}_A(t)\|_{L^1_x} = \|\bw^{\e}_A(0)\|_{L^1_x} = \|\bw^{\e}_0\|_{L^1_x} \leq C.
\enda
\end{align*}
(vii) Multiplying~\eqref{weqnmA} by $\bw^{\e}_A$ gives
\begin{align*}
\bega
\frac{d}{dt}\int_{\R^2}|\bw^{\e}_A(t,x)|^2dx &= -\int_{\R^2} (\mathbb{P}m^{\e}\cdot \nabla_x \bw^{\e}_A)\bw^{\e}_A dx = \frac{1}{2}\int_{\R^2} (\nabla_x \cdot \mathbb{P}m^{\e}) |\bw^{\e}_A|^2 dx =0 .
\enda
\end{align*}
By the admissible blow-up condition~\eqref{ABC1}, we have
\begin{align*}
\bega
\|\bw^{\e}_A(t)\|_{L^2_x} = \|\bw^{\e}_A(0)\|_{L^2_x} = \|\bw^{\e}_0\|_{L^2_x} \leq C\log\!\Big(\log\!\big(\log(1/\e)\big)\Big). 
\enda
\end{align*}
\hide
(ix) By using the Minkowski inequality, we have 
\begin{align*}
\bega
\|\bw^{\e}_B(t)\|_{L^2_x} &= \bigg\|\int_0^t \varPi_{\bw}^{\e}(s,X^{\e}(s,\mathbf{a})) ds \bigg\|_{L^2_x} \leq  \int_0^t \big\| \varPi_{\bw}^{\e}(s,X^{\e}(s,\mathbf{a}))  \big\|_{L^2_x}ds,
\enda
\end{align*}
Since $\mathbb{P}m^{\e}$ is divergence free, the particle trajectory \eqref{traj} is volume preserving and the jacobian satisfies  $J(\mathbf{a},t)=\det(\nabla_\mathbf{a}X(t,\mathbf{a}))=1$. Thus we have 
\begin{align*}
\bega
\|\bw^{\e}_B(t)\|_{L^2_x} &\leq \int_0^t \big\| \varPi_{\bw}^{\e}(s)  \big\|_{L^2_x}ds \leq \sqrt{T} \big\| \varPi_{\bw}^{\e}\big\|_{L^2_TL^2_x},
\enda
\end{align*}
for some $T>0$. 
\unhide
\end{proof}

\begin{lemma}[Rotation invariance]\label{L.rotation}
For any rotation matrix $\mathcal{O}\in SO(2)$, the Boltzmann equation and its momentum equation are rotation invariant in the following sense:
\begin{itemize}
\item Let $F^{\e}(t,x,v)$, for $(t,x,v)\in \R^+\times\R^2\times\R^3$, be a solution of the Boltzmann equation~\eqref{BE}. Then
$F^{\e}(t,\mathcal{O}(x_1,x_2),\mathcal{O}(v_1,v_2),v_3)$ is also a solution of the Boltzmann equation.
\item Let $m^{\e}_A(t,x)$ be a solution of~\eqref{mAeqn}. Then $\mathcal{O}^{-1}m^{\e}_A(t,\mathcal{O}x)$ is also a solution of~\eqref{mAeqn}.
\end{itemize}
\end{lemma}
\begin{proof}
(1) For a rotation matrix $\mathcal{O}\in SO(2)$, define $\bar{F}^{\e}(t,x,v) = F^{\e}(t,\mathcal{O}x,\mathcal{O}(v_1,v_2),v_3)$. 
A direct computation shows that $\bar{F}^{\e}$ satisfies the same Boltzmann equation, hence the equation is rotation invariant. \\
(2) We define 
\begin{align*}
\bar{m}^{\e}(t,x) := \frac{1}{\e}\int_{\R^3} v\bar{F}^{\e}(t,x,v) dv, \qquad \bar{u}^{\e}(t,x) := \frac{1}{\e}\frac{\int_{\R^3} v\bar{F}^{\e}(t,x,v) dv}{\int_{\R^3} \bar{F}^{\e}(t,x,v) dv}, \qquad \bar{\bw}^{\e}(t,x):= \nabla_x^{\perp}\cdot \bar{m}^{\e}.
\end{align*}
Then, we have
\begin{align*}
\bar{m}^{\e}(t,x) = \mathcal{O}^{-1} m^{\e}(t,\mathcal{O}x), \qquad \bar{u}^{\e}(t,x) = \mathcal{O}^{-1} u^{\e}(t,\mathcal{O}x), \qquad \bar{\bw}^{\e}(t,x) = \bw^{\e}(t,\mathcal{O}x).
\end{align*}
It follows that $\bar{m}^{\e}$ and $\bar{\bw}^{\e}$ satisfy~\eqref{mGeqn} and~\eqref{weqn-m}, respectively. 
Defining $(\bar{\bw}^{\e}_A,\bar{\bw}^{\e}_B)$ and $(\bar{m}^{\e}_A,\bar{m}^{\e}_B)$ in the same way as in Definition~\ref{D.wAwBm}, we conclude that $\bar{m}^{\e}_A$ is also a solution of~\eqref{mAeqn}.

\hide
For the rotation matrix $\mathcal{O}\in SO(2)$, We define 
\begin{align*}
\bar{F}^{\e}(t,x,v) = F^{\e}(t,\mathcal{O}x,\mathcal{O}(v_1,v_2),v_3).
\end{align*}
By the change of variable $y=\mathcal{O}x$ and $z=(\mathcal{O}(v_1,v_2),v_3)$, we have 
\begin{align*}
\p_t\bar{F}^{\e}(t,x,v) &= \p_tF^{\e}(t,\mathcal{O}x,\mathcal{O}(v_1,v_2),v_3) = \p_tF^{\e}(t,y,z), \cr 
\nabla_x\bar{F}^{\e}(t,x,v) &= (\p_{x_1}F^{\e}, \p_{x_2}F^{\e})(t,y,z) \cr 
&= \lw[\begin{matrix}
\frac{\p y_1}{\p x_1}\frac{\p F^{\e}}{\p y_1} + \frac{\p y_2}{\p x_1}\frac{\p F^{\e}}{\p y_2} \\ 
\frac{\p y_1}{\p x_2}\frac{\p F^{\e}}{\p y_1} + \frac{\p y_2}{\p x_2}\frac{\p F^{\e}}{\p y_2} 
\end{matrix} \rw]
= \lw[\begin{matrix}
\frac{\p y_1}{\p x_1} & \frac{\p y_2}{\p x_1} \\ 
\frac{\p y_1}{\p x_2} & \frac{\p y_2}{\p x_2}
\end{matrix}\rw]  
\lw[\begin{matrix}
\frac{\p F^{\e}}{\p y_1} \\ 
\frac{\p F^{\e}}{\p y_2} 
\end{matrix}\rw]
\end{align*}
So that 
\begin{align*}
v\cdot \nabla_x\bar{F}^{\e}(t,x,v) &= v\cdot \mathcal{O}^{T} \nabla_y F^{\e}(t,y,z) = (\mathcal{O}(v_1,v_2)) \cdot \nabla_y F^{\e}(t,y,z)
\end{align*}
Hence, 
\begin{align*}
\big(\p_t F^{\e}+ z \cdot \nabla_y F^{\e}\big)(t,y,z) = \mathcal{N}(F^{\e},F^{\e})(t,y,z).
\end{align*}
Since the collision operator is invariant under rotation $SO(3)$, we can see that $F_\mathcal{O}^{\e}$ is also solution of the Boltzmann equation. \\
(2) Compare 
\begin{align*}
\bar{\mathrm{P}}^{\e}\bar{\mathrm{U}}^{\e}(t,x) &= \int_{\R^3}v \bar{F}^{\e}(t,x,v) dv = \int_{\R^3}v F^{\e}(t,\mathcal{O}x,\mathcal{O}(v_1,v_2),v_3) dv \cr 
&= \int_{\R^3} \mathcal{O}^{-1}z F^{\e}(t,\mathcal{O}x,z) dz = \mathcal{O}^{-1}(\mathrm{P}^{\e}\mathrm{U}^{\e})(t,\mathcal{O}x).
\end{align*}
This implies $\bar{m}^{\e}(t,x) = \mathcal{O}^{-1}m^{\e}(t,\mathcal{O}x)$. 
Since $\mathbb{P}m^{\e}_{\mathcal{O}}(t,x)=\mathcal{O}^{-1}(\mathbb{P}m^{\e})(t,\mathcal{O}x)$, and $\bw^{\e}_R(t,x)=\bw(t,\mathcal{O}x)$, $(\bw^{\e}_A(t,\mathcal{O}x), \mathcal{O}^{-1}(\mathbb{P}m^{\e})(t,\mathcal{O}x))$ solves the equation \eqref{weqnmA}.
\unhide
\end{proof}

We note that the quantities $(m^{\e},\bw^{\e})$ can be controlled by the solution constructed in Theorem~\ref{T.2D.global}.  

\begin{lemma}\label{L.mwest} 
For the solution constructed in Theorem~\ref{T.2D.global}, the pair $(m^{\e},\bw^{\e})$ defined in Definition~\ref{D.mdef} satisfies the following estimates:
\begin{align*}
\bega
\|(m^{\e}-u^{\e})(t)\|_{L^\infty_x} \les \e\mathcal{E}_M^{\frac{1}{2}}(t)\big(\mathcal{E}_{M}^{\frac{1}{2}}(t) +1\big), \quad \mbox{and} \quad \|\bw^{\e}(t)\|_{L^2_x} \les \mathcal{E}_{M}^{\frac{1}{2}}(t).
\enda
\end{align*}
\end{lemma}
\begin{proof}
Using \eqref{uLinfEC}, we get 
\begin{align*}
\bega
\|m^{\e}-u^{\e}\|_{L^\infty_x} = \|\mathrm{P}^{\e}u^{\e}-u^{\e}\|_{L^\infty_x} \leq \|\mathrm{P}^{\e}-1\|_{L^\infty_x}\|u^{\e}\|_{L^\infty_x} \leq \e\mathcal{E}_M^{\frac{1}{2}}\big(\mathcal{E}_{M}^{\frac{1}{2}}(t) +1\big).
\enda
\end{align*}
Similarly, we have
\begin{align*}
\bega
\|\bw^{\e}\|_{L^2_x} &= \|\nabla_x^{\perp}\cdot(\mathrm{P}^{\e}u^{\e})\|_{L^2_x} \leq \|\nabla_x^{\perp}\mathrm{P}^{\e}\|_{L^2_x}\|u^{\e}\|_{L^\infty_x} + \|\mathrm{P}^{\e}\|_{L^\infty_x}\|\w^{\e}\|_{L^2_x} \leq \e\mathcal{E}_M^{\frac{1}{2}}\big(\mathcal{E}_{M}^{\frac{1}{2}}(t) +1\big) + \mathcal{E}_M^{\frac{1}{2}}.
\enda
\end{align*}
\end{proof}

\begin{lemma}\label{L.intw} 
For the solution $\bw^{\e}$ of~\eqref{weqn-m}, we have
\begin{align}\label{dwint}
&\frac{d}{dt}\int_{\R^2}\bw^{\e}(t,x)dx = 0,  \qquad
\frac{d}{dt}\int_{\R^2}|x|^2\bw^{\e}(t,x)dx =0, 
\end{align}
and
\begin{align}\label{Hwlim}
\bega
\lim_{\e\to0}H(\bw^{\e}(t))=\lim_{\e\to0}H(\bw^{\e}_0),
\enda
\end{align}
where the pseudo--energy is defined in~\eqref{pseudodef}.
\end{lemma}

\begin{remark}\label{Rmk.x2wA}
From the equation for $\bw^{\e}_A$ in~\eqref{weqnmA}, the function $\bw^{\e}_A$ does not, in general, conserve the second moment:
\begin{align*}
\bega
\frac{d}{dt}\int_{\R^2}|x|^2\bw^{\e}_A(t,x)dx = -2\int_{\R^2} (x\cdot \mathbb{P}m^{\e}) \bw^{\e}_B dx.
\enda
\end{align*}
\end{remark}

\begin{proof}[Proof of Lemma \ref{L.intw}]
(Proof of $\eqref{dwint}_1$) 
Starting from $\eqref{mweqn}_2$, integration by parts immediately yields the claim, since
\begin{align*}
\bega
\frac{d}{dt}\int_{\R^2}\bw^{\e}(t,x)dx &= - \int_{\R^2} \nabla_x^{\perp}\cdot\big(m^{\e}\cdot\nabla_x u^{\e} + u^{\e}(\nabla_x\cdot m^{\e})\big) dx = 0 .
\enda
\end{align*}
(Proof of $\eqref{dwint}_2$) 
Multiplying $\eqref{mweqn}_2$ by $|x|^2$ and integrating over $\R^2$, we obtain
\begin{align}\label{x2w=}
\bega
\frac{d}{dt}\int_{\R^2}|x|^2&\bw^{\e}(t,x)dx = I_{\bw}(t) + II_{\bw}(t),
\enda
\end{align}
where
\begin{align*}
\bega
I_{\bw}(t):= - \int_{\R^2} |x|^2 \nabla_x^{\perp}\cdot\big(m^{\e}\cdot\nabla_x u^{\e} + u^{\e}(\nabla_x\cdot m^{\e})\big) dx, \quad II_{\bw}(t):=-\frac{1}{\e^2}\sum_{j} \int_{\R^2} |x|^2 \nabla_x^{\perp}\cdot\p_{x_j} \mathbf{r}_{ij}^{\e}  dx .
\enda
\end{align*}
For $II_{\bw}(t)$, applying integration by parts to each term gives
\begin{align*}
\bega
II_{\bw}(t) &= \frac{2}{\e^2}\sum_{j} \int_{\R^2} \bigg( x_2\p_{x_j} \mathbf{r}_{1j}^{\e} -x_1 \p_{x_j} \mathbf{r}_{2j}^{\e}\bigg) dx =\frac{2}{\e^2}\int_{\R^2} \big( \mathbf{r}_{12}^{\e} -\mathbf{r}_{21}^{\e}\big) dx = 0.
\enda
\end{align*}
For $I_{\bw}(t)$, integrating by parts twice yields
\begin{align*}
\bega
I_{\bw}(t)&= 2\int_{\R^2} x_1\big(m^{\e}\cdot\nabla_x u^{\e}_2 + u^{\e}_2(\nabla_x\cdot m^{\e})\big) - x_2\big(m^{\e}\cdot\nabla_x u^{\e}_1 + u^{\e}_1(\nabla_x\cdot m^{\e})\big) dx \cr 
&= \int_{\R^2} \bigg(-2m^{\e}_1u^{\e}_2 + 2\sum_j\big(-x_1 (\p_jm^{\e}_j) u^{\e}_2 + x_1u^{\e}_2(\p_jm^{\e}_j)\big) \bigg)dx  \cr 
&\quad+ \int_{\R^2} \bigg(2m^{\e}_2u^{\e}_1 + 2\sum_j\big(x_2 (\p_jm^{\e}_j) u^{\e}_1 - x_2 u^{\e}_1(\p_jm^{\e}_j) \big)  \bigg)dx \cr 
&= 0,
\enda
\end{align*}
where we used $m^{\e}_1u^{\e}_2-m^{\e}_2u^{\e}_1=0$ since $m^{\e}_iu^{\e}_j= \mathrm{P}^{\e}u^{\e}_iu^{\e}_j$. 
Combining the identities for $I_{\bw}(t)$ and $II_{\bw}(t)$ in~\eqref{x2w=} completes the proof. \\

\noindent (Proof of \eqref{dHw})
Since $\mathbb{P}m^{\e}$ is divergence-free, there exists a stream function $\psi^{\e}$ such that 
\begin{align}\label{m-stream}
\bega
\mathbb{P}m^{\e} = \nabla_x^{\perp}\psi^{\e} = (-\p_2\psi^{\e},\p_1\psi^{\e}), \qquad \bw^{\e} = \Delta_x \psi^{\e}.
\enda
\end{align}
We split the proof into two steps. \\
(Step 1) We first claim that 
\begin{align}\label{dHw}
\bega
\frac{d}{dt}H(\bw^{\e})& \les \|\bw^{\e}\|_{L^2_x} \kappa^{\frac{1}{2}}\mathcal{D}_{G}^{\frac{1}{2}} + \|\mathbb{P}^{\perp}m^{\e}\|_{L^{2+\frac{1}{\delta}}_x} \|\nabla_x\psi^{\e}\|_{L^{2+4\delta}_x}\|\w^{\e}\|_{L^2_x} + \|\nabla_x \cdot m^{\e}\|_{L^{2+\frac{1}{\delta}}_x} \|\mathbb{P}m^{\e}\|_{L^{\frac{2+4\delta}{1+\delta}}_x}^2 \cr 
&+ \|\nabla_x\psi^{\e}\|_{L^4_x}\|m^{\e}-u^{\e}\|_{L^4_x}\|\nabla_x m^{\e}\|_{L^2_x}.
\enda
\end{align}
Differentiating $H(\bw^{\e})$ in time and using the stream formulation~\eqref{m-stream}, we obtain
\begin{align*}
\bega
\frac{d}{dt}H(\bw^{\e})&= \frac{1}{2\pi} \int_{\R^2}\int_{\R^2} \log|x-y|\Big(\p_t\bw^{\e}(x)\bw^{\e}(y)+\bw^{\e}(x)\p_t\bw^{\e}(y)\Big) dydx 
=2\int_{\R^2} \psi^{\e}(x)\p_t\bw^{\e}(x) dx.
\enda
\end{align*}
Using the equation for $\bw^{\e}$ in~$\eqref{mweqn}_2$, we can write
\begin{align*}
\bega
\frac{d}{dt}H(\bw^{\e})&= I_{H(\bw)}+ II_{H(\bw)}, 
\enda
\end{align*}
where
\begin{align*}
\bega
I_{H(\bw)}&=-\frac{2}{\e^2}\sum_{j} \int_{\R^2} \psi^{\e} \Big(\nabla_x^{\perp}\cdot\p_{x_j} \mathbf{r}_{ij}^{\e}\Big) dx , \cr 
II_{H(\bw)}&= -2\int_{\R^2} \psi^{\e}\bigg(\nabla_x^{\perp}\cdot\big(m^{\e}\cdot\nabla_x u^{\e} + u^{\e}(\nabla_x\cdot m^{\e})\big)\bigg)dx.
\enda
\end{align*}
For the microscopic term $I_{H(\bw)}$, integration by parts gives
\begin{align}\label{IIHest}
\bega
I_{H(\bw)} &= -2\sum_{j} \int_{\R^2} \bigg(\p_{x_j}\nabla_x^{\perp}\psi^{\e}\cdot \frac{1}{\e^2} \mathbf{r}_{ij}^{\e} \bigg) dx \les \|\Delta_x\psi^{\e}\|_{L^2_x} \frac{1}{\e^2}\|\mathbf{r}_{ij}^{\e}\|_{L^2_x} \les \|\bw^{\e}\|_{L^2_x} \kappa^{\frac{1}{2}}\mathcal{D}_{G}^{\frac{1}{2}} ,
\enda
\end{align}
where we used~\eqref{ABGscale}.
For the main term $II_{H(\bw)}$, in order to obtain full cancellation of the contribution $u^{\e}(\nabla_x\cdot m^{\e})$, we apply integration by parts to the term $m^{\e}\cdot\nabla_x u^{\e}$:
\begin{align*}
\bega
II_{H(\bw)}&= 2\int_{\R^2} \sum_j \p_j\psi^{\e}\bigg(\nabla_x^{\perp}\cdot\big(u^{\e}m^{\e}_j\big)\bigg)dx =  II_{H(\bw)}^1 +II_{H(\bw)}^2 ,
\enda
\end{align*}
where
\begin{align*}
\bega
II_{H(\bw)}^1&= 2\int_{\R^2} \Big(m^{\e}\cdot\nabla_x \psi^{\e}\Big) \w^{\e}dx, \qquad II_{H(\bw)}^2= 2\sum_j\int_{\R^2} \p_j\psi^{\e} \big(u^{\e}_1(-\p_2m^{\e}_j)+u^{\e}_2(\p_1m^{\e}_j)\big)dx .
\enda
\end{align*}
For $II_{H(\bw)}^1$, using $\mathbb{P}m^{\e}\cdot\nabla_x \psi^{\e}
=\nabla_x^{\perp}\psi^{\e}\cdot\nabla_x\psi^{\e}=0$ from~\eqref{m-stream} and applying H\"older's inequality, we have
\begin{align}\label{I-1Hest}
\bega
II_{H(\bw)}^1 &= 2\int_{\R^2} \Big(\mathbb{P}^{\perp}m^{\e}\cdot\nabla_x \psi^{\e}\Big) \w^{\e}dx \les \|\mathbb{P}^{\perp}m^{\e}\|_{L^{2+\frac{1}{\delta}}_x} \|\nabla_x\psi^{\e}\|_{L^{2+4\delta}_x}\|\w^{\e}\|_{L^2_x}.
\enda
\end{align}
For $II_{H(\bw)}^2$, we again use the stream form
$\mathbb{P}m^{\e}=(-\p_2\psi^{\e},\p_1\psi^{\e})$ and decompose
$u^{\e}= m^{\e}-(m^{\e}-u^{\e})$ to obtain
\begin{align*}
\bega
II_{H(\bw)}^2 
&= 2\int_{\R^2} \Big(\mathbb{P}m^{\e}_2m^{\e}_1\p_2m^{\e}_1 - \mathbb{P}m^{\e}_2m^{\e}_2\p_1m^{\e}_1 -\mathbb{P}m^{\e}_1m^{\e}_1\p_2m^{\e}_2 + \mathbb{P}m^{\e}_1m^{\e}_2\p_1m^{\e}_2 \Big)dx \cr 
&- 2\sum_j\int_{\R^2} \p_j\psi^{\e} \big((m^{\e}_1-u^{\e}_1)(-\p_2m^{\e}_j)+(m^{\e}_2-u^{\e}_2)(\p_1m^{\e}_j)\big)dx.
\enda
\end{align*}
For the first line of $II_{H(\bw)}^2$, integration by parts yields
\begin{align*}
\bega
2\int_{\R^2} (\nabla_x\cdot m^{\e})|\mathbb{P}m^{\e}|^2dx -  \int_{\R^2} (\nabla_x\cdot \mathbb{P}m^{\e})|\mathbb{P}m^{\e}|^2dx = 2\int_{\R^2} (\nabla_x\cdot m^{\e})|\mathbb{P}m^{\e}|^2dx.
\enda
\end{align*}
\hide
\begin{align*}
\bega
&= 2\int_{\R^2} (\nabla_x\cdot m^{\e}) \big(\mathbb{P}m^{\e}_1m^{\e}_1+\mathbb{P}m^{\e}_2m^{\e}_2\big)dx +  2\int_{\R^2} \Big(\mathbb{P}m^{\e}_2m^{\e}_1\p_2m^{\e}_1 +\mathbb{P}m^{\e}_2m^{\e}_2\p_2m^{\e}_2 +\mathbb{P}m^{\e}_1m^{\e}_1\p_1m^{\e}_1 + \mathbb{P}m^{\e}_1m^{\e}_2\p_1m^{\e}_2 \Big)dx \cr 
&= 2\int_{\R^2} (\nabla_x\cdot m^{\e})|\mathbb{P}m^{\e}|^2dx -  \int_{\R^2} (\nabla_x\cdot \mathbb{P}m^{\e})|\mathbb{P}m^{\e}|^2dx \cr 
&= 2\int_{\R^2} (\nabla_x\cdot m^{\e})|\mathbb{P}m^{\e}|^2dx
\enda
\end{align*}
\unhide
Since $\mathbb{P}m^{\e}\notin L^2(\R^2)$, applying H\"older's inequality to both lines of $II_{H(\bw)}^2$ gives 
\begin{align}\label{I-2Hest}
\bega
II_{H(\bw)}^2 \les \|\nabla_x \cdot m^{\e}\|_{L^{2+\frac{1}{\delta}}_x} \|\mathbb{P}m^{\e}\|_{L^{\frac{2+4\delta}{1+\delta}}_x}^2 + \|\nabla_x\psi^{\e}\|_{L^4_x}\|m^{\e}-u^{\e}\|_{L^4_x}\|\nabla_x m^{\e}\|_{L^2_x},
\enda
\end{align}
for any $\delta>0$. 
Combining~\eqref{IIHest}, \eqref{I-1Hest}, and~\eqref{I-2Hest} proves the claim~\eqref{dHw}. \\
(Step 2) 
To prove~\eqref{Hwlim}, we claim that the time integral of the right-hand side of~\eqref{dHw} vanishes:
\begin{align}\label{claimdw2}
\bega
&\lim_{\e\to0}\int_0^T \bigg[\|\bw^{\e}(t)\|_{L^2_x} \kappa^{\frac{1}{2}}\mathcal{D}_{G}^{\frac{1}{2}}(t) + \|\mathbb{P}^{\perp}m^{\e}(t)\|_{L^{2+\frac{1}{\delta}}_x} \|\nabla_x\psi^{\e}(t)\|_{L^{2+4\delta}_x}\|\w^{\e}(t)\|_{L^2_x} \cr 
& + \|\nabla_x \cdot m^{\e}(t)\|_{L^{2+\frac{1}{\delta}}_x} \|\mathbb{P}m^{\e}(t)\|_{L^{\frac{2+4\delta}{1+\delta}}_x}^2 + \|\nabla_x\psi^{\e}(t)\|_{L^4_x}\|(m^{\e}-u^{\e})(t)\|_{L^4_x}\|\nabla_x m^{\e}(t)\|_{L^2_x}\bigg] dt =0 .
\enda
\end{align}
For the first term in~\eqref{claimdw2}, we decompose $\bw^{\e}= \bw^{\e}_A+\bw^{\e}_B$ and use~\eqref{wAL2} and~\eqref{mBto0} to obtain
\begin{align*}
\bega
\int_0^T \|\bw^{\e}(s)\|_{L^2_x} \kappa^{\frac{1}{2}}\mathcal{D}_{G}^{\frac{1}{2}}(s)ds &\les \kappa^{\frac{1}{2}} \sqrt{T} \sup_{t\in[0,T]}\big(\|\bw^{\e}_A(t)\|_{L^2_x}+\|\bw^{\e}_B(t)\|_{L^2_x}\big) \Big(\int_0^T \mathcal{D}_{G}(t)dt \Big)^{\frac{1}{2}} \cr 
&\les \kappa^{\frac{1}{2}} \sqrt{T} \bigg(\log\!\Big(\log\!\big(\log(1/\e)\big)\Big)
+\kappa^{\frac{1}{2}-}\bigg) \Big(\int_0^T \mathcal{D}_{G}(t)dt \Big)^{\frac{1}{2}}.
\enda
\end{align*}
Since the factor $\kappa^{\frac{1}{2}}$ absorbs both 
$\log\!\big(\log\!\big(\log(1/\e)\big)\big)$ and 
$\int_0^T\mathcal{D}_G(t)\,dt$ for any $\kappa=\e^q$, 
this term converges to $0$ as $\e\to0$. \\
For the remaining three terms in~\eqref{claimdw2}, we estimate $\|\mathbb{P}m^{\e}\|_{L^p_x}$ for $p>2$
(because $\mathbb{P}m^{\e}\notin L^2(\R^2)$) and
$\|\mathbb{P}^{\perp}m^{\e}\|_{L^{2+\frac{1}{\delta}}_x}$.
The stream formulation~\eqref{m-stream} implies that
$\|\nabla_x\psi^{\e}\|_{L^p_x} \leq \|\mathbb{P}m^{\e}\|_{L^p_x}$ for $p>2$.
\\ To bound $\mathbb{P}m^{\e}$, we use $\|\mathbb{P}(\cdot)\|_{L^p_x}\leq\|(\cdot)\|_{L^p_x}$ for $p>1$,
decompose $u^{\e}=(u^{\e}-\bar{u})+\bar{u}$, and invoke the estimate for $\|\bar{u}\|_{L^p_x}$ in Lemma~\ref{L.barubdd}:
\begin{align}\label{dpsiest}
\bega
\|\mathbb{P}m^{\e}(t)\|_{L^p_x} &\leq C\|\mathrm{P}^{\e}(t)\|_{L^\infty_x}\|u^{\e}(t)\|_{L^p_x} \leq C\big(\|(u^{\e}-\bar{u})(t)\|_{L^p_x} + \|\bar{u}\|_{L^p_x} \big) \cr 
&\leq C\bigg(\|\nabla_x(u^{\e}-\bar{u})(t)\|_{L^2_x}^{1-\frac{2}{p}}\|(u^{\e}-\bar{u})(t)\|_{L^2_x}^{\frac{2}{p}} + \Big(\frac{1}{p-2}\Big)^{\frac{1}{p}} \bigg) \cr 
&\les 1+\mathcal{E}_M(t) + \Big(\frac{1}{p-2}\Big)^{\frac{1}{p}}, \quad \mbox{for} \quad p>2,
\enda
\end{align}
where we used the Gagliardo--Nirenberg interpolation inequality~\eqref{Ga-Ni},
$\|(u^{\e}-\bar{u})\|_{L^2_x} \leq \mathcal{E}_M^{\frac{1}{2}}$, and
$\|\nabla_x(u^{\e}-\bar{u})\|_{L^2_x} \leq \mathcal{E}_M^{\frac{1}{2}}+C$. \\
To control $\mathbb{P}^{\perp}m^{\e}$, we decompose $m^{\e}= u^{\e}-(u^{\e}-m^{\e})$ and apply
$\|\mathbb{P}^{\perp}u^{\e}\|_{L^p_x} \leq C\|\mathbb{P}^{\perp}u^{\e}\|_{\dot{B}_{p,1}^{0}}$
via the Besov embedding~\eqref{Besov-Lp}:
\begin{align}\label{incomp.m}
\bega
\|\mathbb{P}^{\perp}m^{\e}(t)\|_{L^{2+\frac{1}{\delta}}_x} &\leq \|\mathbb{P}^{\perp}u^{\e}(t)\|_{L^{2+\frac{1}{\delta}}_x} + \|\mathbb{P}^{\perp}(m^{\e}-u^{\e})(t)\|_{L^{2+\frac{1}{\delta}}_x} \cr 
&\leq \|\mathbb{P}^{\perp}u^{\e}(t)\|_{\dot{B}_{p,1}^0} + \|\mathrm{P}^{\e}(t)-1\|_{L^\infty_x}\|u^{\e}(t)\|_{L^{2+\frac{1}{\delta}}_x} \cr 
&\leq \|\mathbb{P}^{\perp}u^{\e}(t)\|_{\dot{B}_{p,1}^0} + \e (\mathcal{E}_M(t)+1).
\enda
\end{align}
Since the second and third terms in~\eqref{claimdw2} have the same structure, we estimate only the second one.
Using~\eqref{dpsiest}, \eqref{incomp.m}, and applying~\eqref{incomp.2D} with
$s=\frac{3}{4+8\delta}$, we obtain
\begin{align*}
\bega
&\int_0^T \|\mathbb{P}^{\perp}m^{\e}(t)\|_{L^{2+\frac{1}{\delta}}_x} \|\nabla_x\psi^{\e}(t)\|_{L^{2+4\delta}_x}\|\w^{\e}(t)\|_{L^2_x}dt \cr 
&\leq (1+T) \Big(\|\mathbb{P}^{\perp}u^{\e}\|_{L^{4+8\delta}_T\dot{B}_{p,1}^0} + \e \Big(\sup_{0\leq t\leq T}\mathcal{E}_M(t)+1\Big)\Big)\Big(1+ \sup_{0\leq t\leq T}\mathcal{E}_M(t) +\Big(\frac{1}{4\delta}\Big)^{\frac{1}{2+4\delta}} \Big)\sup_{0\leq t\leq T}\mathcal{E}_M^{\frac{1}{2}}(t).
\enda
\end{align*}
For some $\delta< \frac{1}{4}$, the right-hand side tends to $0$ as $\e\to0$, since
$\e^{\frac{1}{4}-}\kappa^{-\frac{1}{2}}(1+T)\frac{1}{\sqrt{\delta}}\mathcal{E}_M \to 0$.
Indeed, a small power of $\e$ absorbs the growth of $\mathcal{E}_M(t)$, since the solution constructed in Theorem~\ref{T.2D.global} satisfies the inequality \eqref{Boot1}.
Finally, for the last term in~\eqref{claimdw2}, the factor $m^{\e}-u^{\e}$ provides one power of $\e$ as in Lemma~\ref{L.mwest}. 
Using $\|(m^{\e}-u^{\e})(t)\|_{L^4_x}\leq \e\mathcal{E}_M$, 
$\|\nabla_x m^{\e}(t)\|_{L^2_x} \les \mathcal{E}_M^{\frac{1}{2}}$, and~\eqref{dpsiest}, we conclude that this term also tends to $0$. 
This completes the proof of the claim~\eqref{claimdw2}.
\end{proof}

\begin{lemma}\label{L.wmaximal} 
For any $T>0$ and $0<R<\frac{1}{2}$, we have 
\begin{align*}
\bega
\max_{0\leq t\leq \frac{T}{x_0}\in\R^2}|\bw^{\e}_A|(B_R(x_0),t) &\leq \bigg(\log\frac{1}{2R}\bigg)^{-\frac{1}{2}} \bigg[|H(\bw^{\e}(t))| + C\int_{\R^2}\bw^{\e}(t,x)dx\int_{\R^2}|x|^2\bw^{\e}(t,x)dx \cr 
&+ \Big(\|\bw^{\e}_A\|_{L^2_x}+ \|\bw^{\e}_B\|_{L^2_x}\Big)\|\bw^{\e}_B\|_{L^2_x} \bigg]^{\frac{1}{2}},
\enda
\end{align*}
where the vorticity maximal function is defined by
\begin{align*}
\bega
|\bw| (B_R(x)) := \int_{B_R(x)}|\bw(y)|dy. 
\enda
\end{align*}
\end{lemma}
\begin{remark}
We note that the left-hand side is the vorticity maximal function associated with $\bw^{\e}_A$, whereas the right-hand side involves $\bw^{\e}$. As explained in Remark~\ref{Rmk.x2wA}, since $\int_{\R^2} |x|^2 \bw^{\e}_A\,dx$ is not conserved, we estimate the right-hand side in terms of $\bw^{\e}$ rather than $\bw^{\e}_A$.
\end{remark}
\begin{proof}[Proof of Lemma \ref{L.wmaximal}]
We use the notation $\log^+(x)=\log(x)$ for $x\geq1$ and $0$ for $0<x<1$, and set $\log^-(x):= \log(x)-\log^+(x)$. 
Since $|\log^-(x)|= -\log(x)+ \log^+(x)$, we have 
\begin{align}\label{log-ww}
\bega
\int_{\R^2}\int_{\R^2} \big|\log^-(|x-y|)\big| \bw^{\e}(y)&\bw^{\e}(x) dydx = \int_{\R^2}\int_{\R^2} -\log(|x-y|)\big| \bw^{\e}(y)\bw^{\e}(x) dydx \cr 
&+ \int_{\R^2}\int_{\R^2} \log^+(|x-y|)\big| \bw^{\e}(y)\bw^{\e}(x) dydx \cr 
&\leq - H(\bw^{\e}) + \int_{\R^2}\int_{\R^2} C(|x|^2+|y|^2) \bw^{\e}(y)\bw^{\e}(x) dydx \cr 
&\leq |H(\bw^{\e}(t))| + C\bigg(\int_{\R^2}\bw^{\e}(t,x)dx\bigg)\bigg(\int_{\R^2}|x|^2\bw^{\e}(t,x)dx\bigg),
\enda
\end{align}
where we used $\log^+(|x-y|) \leq C(|x|^2+|y|^2)$. 
Next, we decompose the left-hand side of~\eqref{log-ww} using $\bw^{\e}= \bw^{\e}_A+\bw^{\e}_B$:
\begin{align}\label{log-wwdecomp}
\bega
&\int_{\R^2}\int_{\R^2} \big|\log^-(|x-y|)\big| \bw^{\e}(y)\bw^{\e}(x) dydx \cr 
&= \int_{\R^2}\int_{\R^2} \big|\log^-(|x-y|)\big| \Big(\bw^{\e}_A(y)\bw^{\e}_A(x)+\bw^{\e}_A(y)\bw^{\e}_B(x)+\bw^{\e}_B(y)\bw^{\e}_A(x)+\bw^{\e}_B(y)\bw^{\e}_B(x) \Big) dydx.
\enda
\end{align}
For the $\bw^{\e}_A(y)\bw^{\e}_A(x)$ interaction in~\eqref{log-wwdecomp}, since $\bw^{\e}_A$ has a distinguished sign by Lemma~\ref{L.mABprop}, we obtain
\begin{align}\label{wAwA}
\bega
\int_{\R^2}\int_{\R^2} \big|\log^-(|x-y|)\big| \bw^{\e}_A(y)\bw^{\e}_A(x) dydx 
&\geq \int_{B_R(x_0)}\int_{B_R(x_0)} (-\log(2R)) \bw^{\e}_A(y)\bw^{\e}_A(x) dydx \cr 
&\geq \log\frac{1}{2R} \Big(|\bw^{\e}_A|(B_R(x_0),t)\Big)^2,
\enda
\end{align}
for $0<R<\frac{1}{2}$. 
For the $\bw^{\e}_B\bw^{\e}_B$ interaction, applying H\"older's inequality and Young's convolution inequality gives
\begin{align}\label{wAwB}
\bega
\int_{\R^2}\int_{\R^2} \big|\log^-(|x-y|)\big| \bw^{\e}_B(y)\bw^{\e}_B(x) dydx 
&\les \|\bw^{\e}_B\|_{L^2_x}\|\bw^{\e}_B\|_{L^2_x} \|\log^-(|x|)\|_{L^1_x} \les \|\bw^{\e}_B\|_{L^2_x}\|\bw^{\e}_B\|_{L^2_x}.
\enda
\end{align}
The $\bw^{\e}_A\bw^{\e}_B$ interactions can be estimated in the same way.
Combining~\eqref{log-wwdecomp}, \eqref{wAwA}, and~\eqref{wAwB} with~\eqref{log-ww} yields the desired estimate.
\end{proof}

\begin{lemma}\cite{Delort}\label{L.Hconti} 
The function
\begin{align}\label{Hpsidef}
\bega
H_{\varphi}(x,y)&:= -\frac{1}{4\pi^2} \frac{\p^2}{\p_{x_1}\p_{x_2}}\int_{\R^2} \log|x-z|\log|y-z| \varphi(z) dz,
\enda
\end{align}
is bounded on $\R^2\times \R^2$ and continuous up to the diagonal $x=y$. Moreover,
\begin{align*}
\bega
H_{\varphi}(x,y)&= \frac{1}{2}(\varphi(x)+\varphi(y))\frac{(x_1-y_1)(x_2-y_2)}{4\pi|x-y|^2} + r(x,y),
\enda
\end{align*}
where $r(x,y)$ is a bounded continuous function.
\end{lemma}
\begin{proof}
For the proof, see \cite{Delort} and \cite{MaBe}. 
\end{proof}

\begin{proof}[\textbf{Proof of Theorem \ref{T.Radon}}]

We first observe that the solution constructed in Theorem~\ref{T.2D.global} satisfies
$u^{\e} \in L^\infty([0,T],L^2_{loc}(\R^2))$ by~\eqref{U.L2inf} in Lemma \ref{L.unif}. 
Using the decomposition $\mathbb{P}m^{\e} = m^{\e}_A + m^{\e}_B$, and the fact that $\mathbb{P}m^{\e}_B \to 0$ in $L^\infty([0,T],L^2_{loc}(\R^2))$ by~\eqref{mBto0}, we obtain
\begin{align*}
\bega
\sup_{0\leq t\leq T}\|\mathbb{P}m^{\e}_A(t)\|_{L^2_{loc}} \leq \sup_{0\leq t\leq T}\|\mathbb{P}(\mathrm{P}^{\e}u^{\e})(t)\|_{L^2_{loc}} + \sup_{0\leq t\leq T}\|\mathbb{P}m^{\e}_B(t)\|_{L^2_{loc}} \leq C.
\enda
\end{align*}
We also note that $\bw^{\e}_A$ is uniformly bounded in $L^\infty(0,T;L^1(\R^2))$ by~\eqref{wAL1} in Lemma~\ref{L.mABprop}. 
Hence, up to a subsequence, there exist 
$m^{\#} \in L^\infty([0,T],L^2_{loc}(\R^2))$ and
$\bw^{\#} \in L^\infty([0,T],\mathcal{M}(\R^2))$ such that
\begin{align}\label{bwconv}
\bega
m^{\e}_A(t,x) &\rightharpoonup m^{\#}(t,x), \qquad \mbox{weakly in} \qquad L^2_{loc}([0,T]\times(\R^2)), \cr 
\bw_A^{\e}(t,x) &\overset{\ast}{\rightharpoonup} d\bw^{\#}(t,x), \qquad \mbox{weakly* in} \qquad L^\infty(0,T;\mathcal{M}(\R^2)).
\enda
\end{align}
Next, arguing as in~\eqref{AubinLion} and using the equation for $m^{\e}_A$ in~\eqref{mAeqn}, we can prove
\begin{align*}
\bega
\big\|\psi m^{\e}_A(t_2) -\psi m^{\e}_A(t_2)\big\|_{H^{-s-1}_x}  \leq C |t_2-t_1|, \quad \mbox{for} \quad \psi \in C_c^{\infty}(\R^2), \quad s>1, \quad 0\leq t_1,t_2\leq T.
\enda 
\end{align*}
Then, by the standard argument in Theorem~1.1 of~\cite{DiMa}, we obtain strong convergence in $L^r([0,T]\times B_R(0))$ for $1 \leq r < 2$.
\hide
Note that proving convergence of $u^{\e}_1u^{\e}_2$ is equivalent to convergence of $m^{\e}_1m^{\e}_2$. Recall
$m^{\e} = \mathrm{P}^{\e}u^{\e}$. 
\begin{align}
\bega
\int_0^T\int_{\R^2} \varphi m^{\e}_1 m^{\e}_2 dxdt - \int_0^T\int_{\R^2} \varphi u^{\e}_1 u^{\e}_2 dxdt &= \int_0^T\int_{\R^2}\big(|\mathrm{P}^{\e}|^2-1\big) \varphi u^{\e}_1 u^{\e}_2 dxdt \cr 
&\leq \e \|\mathrm{P}^{\e}-1\|_{L^\infty_x}\|\mathrm{P}^{\e}+1\|_{L^\infty_x} \|\varphi\|_{L^2_x}\|u^{\e}\|_{L^4_x}^2 \cr 
&\to 0
\enda
\end{align}
(We should avoid $\|u\|_{L^2_x}$.)
\unhide
To show that $m^{\#}$ is a weak solution of the incompressible Euler equation, 
let $\Lambda\in C_c^{\infty}([0,T]\times\R^2)$ be divergence-free.
Then there exists $\eta \in C_c^{\infty}([0,T]\times\R^2)$ such that
$\Lambda=\nabla_x^{\perp}\eta=(-\p_2\eta,\p_1\eta)$. Hence,
\begin{align*}
\bega
\int_0^T\int_{\R^2} \Big(m^{\e}_A\p_t \Lambda &+ \nabla_x \Lambda : m^{\e}_A \otimes m^{\e}_A \Big) dxdt =  \int_0^T\int_{\R^2} (-m^{\e}_{A,1}\p_t\p_2\eta+m^{\e}_{A,2}\p_t\p_1\eta) dxdt \cr 
&+ \int_0^T\int_{\R^2}\Big(\p_1\p_2\eta(|m^{\e}_{A,2}|^2-|m^{\e}_{A,1}|^2)+(\p_1^2\eta-\p_2\eta^2)m^{\e}_{A,1}m^{\e}_{A,2}\Big) dxdt.
\enda
\end{align*}
By Lemma~\ref{L.rotation}, the equation for $m^{\e}_A$ in~\eqref{mAeqn} is rotation invariant.
Moreover, rotating the term $|m^{\e}_{A,2}|^2-|m^{\e}_{A,1}|^2$ by $\pi/4$ transforms it into $m^{\e}_{A,1}m^{\e}_{A,2}$.
Therefore, to prove that the limit $m^{\#}$ is a weak solution, i.e.,
\begin{align*}
\int_0^T \int_{\R^2} (m^{\#} \p_t\Lambda  + \nabla_x \Lambda : m^{\#}\otimes m^{\#} ) dxdt = 0,
\end{align*}
it suffices to establish the following claim:
\begin{align}\label{claim.mm}
\bega
\lim_{\e\to 0} \int_0^T\int_{\R^2} \Lambda m^{\e}_{A,1} m^{\e}_{A,2} dxdt = \int_0^T\int_{\R^2} \Lambda m^{\#}_1 m^{\#}_2 dxdt.
\enda
\end{align}
To prove~\eqref{claim.mm}, without loss of generality we take a separable test function
$\Lambda(t,x)=\psi(t)\varphi(x)$ with $\psi(t)\in C_c^{\infty}(\R_+)$ and $\varphi(x)\in C_c^{\infty}(\R^2)$.
By the Biot--Savart law $m^{\e}_A= {\bf K} \ast \bw^{\e}_A$ (Definition~\ref{D.wAwBm}), we have
\begin{align*}
\bega
\int_{\R^2} \Lambda(t,x) m^{\e}_{A,1} m^{\e}_{A,2} dxdt &= \int_{\R^2} \psi(t)\varphi(x) ({\bf K}\ast \bw^{\e}_A)_1 ({\bf K}\ast \bw^{\e}_A)_2 dx \cr 
&= \int_{\R^2} \psi(t)\varphi(x) \int_{\R^2}-\p_2 \frac{1}{2\pi}\log|z-x|\bw^{\e}_A(z)dz \int_{\R^2}\p_1 \frac{1}{2\pi}\log|x-y|\bw^{\e}_A(y)dy dx \cr 
&= \int_{\R^2}\int_{\R^2} \psi(t) H_{\varphi}(x,y) \bw^{\e}_A(x)\bw^{\e}_A(y) dxdy,
\enda
\end{align*}
where $H_{\varphi}(x,y)$ is defined in~\eqref{Hpsidef}.
Let $\chi\in C_c^{\infty}(\R^2)$ be a smooth cutoff such that $\chi(x)=1$ for $|x|<1$ and $\chi(x)=0$ for $|x|>2$.
For $\mathfrak{k}<\frac{1}{2}$, we decompose
\begin{align*}
\bega
\int_0^T\int_{\R^2}\int_{\R^2} \psi(t) H_{\varphi}(x,y) \bw^{\e}_A(x)\bw^{\e}_A(y) dxdydt&= I+II,
\enda
\end{align*}
where
\begin{align*}
\bega
I&:= \int_0^T\int_{\R^2}\int_{\R^2} \chi\Big(\frac{|x-y|}{\mathfrak{k}}\Big)H_{\varphi}(x,y) \bw^{\e}_A(t,x)\bw^{\e}_A(t,y) dxdydt, \cr
II&:= \int_0^T\int_{\R^2}\int_{\R^2} \bigg(1-\chi\Big(\frac{|x-y|}{\mathfrak{k}}\Big)\bigg)H_{\varphi}(x,y) \bw^{\e}_A(t,x)\bw^{\e}_A(t,y) dxdydt .
\enda
\end{align*}
For $II$, Lemma~\ref{L.Hconti} implies that $H_{\varphi}(x,y)$ is continuous on $\{|x-y|>\mathfrak{k}\}$. 
Hence, using the weak* convergence of $\bw^{\e}_A$ in~\eqref{bwconv}, we obtain
\begin{align*}
\bega
\bw_A^{\e}(t,x)\otimes \bw_A^{\e}(t,y) \rightharpoonup d\bw^{\#}(t,x) \otimes d\bw^{\#}(t,y), \qquad \mbox{weakly in} \qquad L^\infty([0,T],\mathcal{M}(\R^2\times\R^2)).
\enda
\end{align*}
Consequently,
\begin{align*}
\bega
\lim_{\e\to0}II &= 
\int_0^T\int_{\R^2}\int_{\R^2} \bigg(1-\chi\Big(\frac{|x-y|}{\mathfrak{k}}\Big)\bigg)H_{\varphi}(x,y) d\bw^{\#}(t,x)d\bw^{\#}(t,y)dt.
\enda
\end{align*}
For $I$, since $H_{\varphi}(x,y)$ is bounded by Lemma~\ref{L.Hconti}, 
the estimate of the vorticity maximal function in Lemma~\ref{L.wmaximal} yields
\begin{align}\label{HIest0}
\bega
I \leq \int_0^T\bigg(\log\frac{1}{2\mathfrak{k}}\bigg)^{-\frac{1}{2}} &\bigg[|H(\bw^{\e}(t))| + C\int_{\R^2}\bw^{\e}(t,x)dx\int_{\R^2}|x|^2\bw^{\e}(t,x)dx \cr 
&+ \Big(\|\bw^{\e}_A\|_{L^2_x}+ \|\bw^{\e}_B\|_{L^2_x}\Big)\|\bw^{\e}_B\|_{L^2_x} \bigg]^{\frac{1}{2}} \int_{\R^2}|\bw_A^{\e}(t,x)|dx dt.
\enda
\end{align}
In the first line of~\eqref{HIest0} we use~\eqref{dwint} and~\eqref{Hwlim} in Lemma~\ref{L.intw}.
For the last factor we use $\|\bw^{\e}_A\|_{L^1_x} \leq \|\bw_0\|_{\mathcal{M}}$ from~\eqref{wAL1}.
For the terms $\|\bw^{\e}_A\|_{L^2_x}$ and $\|\bw^{\e}_B\|_{L^2_x}$, applying~\eqref{wAL2} and~\eqref{mBto0} gives
\begin{align*}
\bega
\sup_{0\leq t\leq T}\Big(\|\bw^{\e}_A(t)\|_{L^2_x}+ \|\bw^{\e}_B(t)\|_{L^2_x}\Big)\|\bw^{\e}_B(t)\|_{L^2_x} \les \Big(\log\!\Big(\log\!\big(\log(1/\e)\big)\Big) +\kappa^{\frac{1}{2}-}\Big)\kappa^{\frac{1}{2}-}   \to 0, \quad \mbox{as} \quad \e\to0, 
\enda
\end{align*}
since a small power of $\kappa$ absorbs the singular admissible blow-up rate for any $\kappa=\e^q$.
Therefore,
\begin{align}\label{HIest}
\bega
\lim_{\e\to0} I &\leq C_T  \bigg(\log\frac{1}{2\mathfrak{k}}\bigg)^{-\frac{1}{2}} \lim_{\e\to0} \bigg[|H(\bw^{\e}_0)| + \int_{\R^2}\bw^{\e}_0dx\int_{\R^2}|x|^2\bw^{\e}_0dx \bigg]^{\frac{1}{2}} \|\bw_0\|_{\mathcal{M}}.
\enda
\end{align}
By the weak* convergence in~\eqref{bwconv}, we also have
\begin{align*}
\bega
\int_{B_R(x_0)} d|\bw^{\#}(t,x)| \leq \liminf_{\e\to0} \int_{B_R(x_0)} \bw^{\e}_A(t,x) dx.
\enda
\end{align*}
Applying Lemma~\ref{L.wmaximal} again, we obtain the same bound as in~\eqref{HIest} for $\bw^{\#}$:
\begin{align}\label{HIestbw}
\bega
\int_0^T\int_{\R^2}\int_{\R^2} \chi\Big(\frac{|x-y|}{\mathfrak{k}}\Big)H_{\varphi}(x,y) d\bw^{\#}(t,x)d\bw^{\#}(t,y) dt \leq C_T  \bigg(\log\frac{1}{2\mathfrak{k}}\bigg)^{-\frac{1}{2}}. 
\enda
\end{align}
Since $\mathfrak{k}>0$ can be chosen arbitrarily small, the right-hand sides of~\eqref{HIest} and~\eqref{HIestbw} can be made arbitrarily small. 
This proves the claim~\eqref{claim.mm}, and therefore completes the proof of Theorem~\ref{T.Radon}.\end{proof}

\StartNoTOC

\section*{Acknowledgment}
Part of this result was presented in the seminar series \href{https://www.siam.org/publications/siam-news/articles/seminar-in-the-analysis-and-methods-of-pde}{\textit{Analysis and Methods of PDE}} hosted by SIAM. 
GCB is supported by BK21 SNU Mathematical Sciences Division and acknowledges support from the National Research Foundation of Korea (NRF-2021R1C1C2094843). CK is partially supported by NSF-CAREER 2047681, the Simons Fellowship in Mathematics, and the Brain Pool Fellowship funded by the Korean Ministry of Science and ICT through the National Research Foundation of Korea (NRF-2021H1D3A2A01039047). This material is partly based upon work supported by the National Science Foundation under Grant No. DMS-2424139, while one of the  authors (C.K.) was in residence at the Simons Laufer Mathematical Sciences Institute in Berkeley, California, during the Fall 2025 semester. GCB also gratefully acknowledges the hospitality of the Department of Mathematics at the University of Wisconsin–Madison during his several long-term visits. The authors thank Hongxu Chen for carefully checking the manuscript and for valuable suggestions that improved the argument.

\StopNoTOC


\hide
\section{Discussion (NS)}

GLOBAL HILBERT EXPANSION FOR SOME NON-RELATIVISTIC KINETIC EQUATIONS, YUANJIE LEI, SHUANGQIAN LIU, QINGHUA XIAO, AND HUIJIANG ZHAO

Discussion 1:
\textcolor{red}{Note that we used $\p_{\tilde{t}} = \e^{\mathfrak{n}}\p_t$. We need to explain how uniform estimates can be obtained for $0<\mathfrak{n}<1$. When $0<\mathfrak{n}<1$, the estimates for certain terms are different, and the associated exponential growth behavior also changes.}

Essential: 
\begin{enumerate}
\item Scale of $div(u)$ and $(\rho+\ta)$ becomes $\e^{1-\mathfrak{n}}$ (Benefit of well-preparedness) (There can be some issues if we can only estimate $\div(u)$ instead of $\mathbb{P}^{\perp}u$).
\item Some terms need to change the estimates
\end{enumerate}

Discussion 2:
\textcolor{red}{Need to discuss the convergence rate when the time derivative is included and use NS}

Essential: 
\begin{enumerate}
\item G- expansion and Equation for local conservation law with dissipation
\item Vorticity equation and $\rho-3/2\ta$ equation with dissipation 
\item (Technical) Computation for $\mathcal{L}^{-1}$
\item (Technical) Embedding for time 
\item Convergence rate in the convergence section. (Some part $\kappa^{\frac{1}{2}-}$ will be changed to $\kappa^{\frac{1}{2}}$)
\end{enumerate}

\begin{lemma}
For $t \in [0,T]$, where $T > 0$ is defined in \eqref{condition}, and for both types of derivatives $\p^{\al}$ defined in \eqref{caseA} and \eqref{caseB}, and for $\mathrm{N} > \frac{d}{2} + 1$, the following estimates hold:
\begin{align}
\bega
&(3) ~ \sup_{t \in [0,T]} \TbT(t) 
\leq \frac{C}{\sqrt{T}}\e^2\kappa^{\frac{1}{2}}\Big(\int_0^T\mathcal{D}_G(s)ds\Big)^{\frac{1}{2}} + \sqrt{T}  \frac{\e^2}{\e^{\mathfrak{n}}}\bigg[\Big(\int_0^T\mathcal{D}_{tot}(s)ds\Big)^{\frac{1}{2}} + \Big(\int_0^T \mathcal{E}_{M}^2(s)ds\Big)^{\frac{1}{2}}\bigg],
\enda
\end{align}
for a positive constant $C_2 > 0$, where $\mathbf{r}_{ij}^{\e}$, $\mathfrak{q}_j^{\e}$, and $\TbT(t)$ are defined in \eqref{albe-def} and \eqref{RSdef}, respectively.
\end{lemma}
\begin{proof}
(3) 
Taking the $L^\infty_t$ norm of $\|\p_{x_j} \mathbf{r}_{ij}^{\e}\|_{L^\infty_x}$ and applying the temporal embedding \eqref{tembedding} along with $\eqref{ABGscale}_2$, we obtain
\begin{align*}
\bega
\|\p_{x_j} \mathbf{r}_{ij}^{\e}\|_{L^\infty_TL^\infty_x} &\leq \frac{1}{\sqrt{T}}\|\p_{x_j} \mathbf{r}_{ij}^{\e}\|_{L^2_TL^\infty_x} + \sqrt{T}\|\p_t\p_{x_j} \mathbf{r}_{ij}^{\e}\|_{L^2_TL^\infty_x} \cr 
&\les \frac{1}{\sqrt{T}}\e^2\kappa^{\frac{1}{2}}\Big(\int_0^T\mathcal{D}_G(s)ds\Big)^{\frac{1}{2}} + \frac{\sqrt{T}}{\e^{\mathfrak{n}}} \|(\e^{\mathfrak{n}}\p_t)\p_{x_j}\mathbf{r}_{ij}^{\e}\|_{L^2_TL^\infty_x} \cr 
&\les \frac{1}{\sqrt{T}}\e^2\kappa^{\frac{1}{2}}\Big(\int_0^T\mathcal{D}_G(s)ds\Big)^{\frac{1}{2}} + \frac{\sqrt{T}}{\e^{\mathfrak{n}}}  \e^2\bigg[\Big(\int_0^T\mathcal{D}_{tot}(s)ds\Big)^{\frac{1}{2}} + \Big(\int_0^T \mathcal{E}_{M}^2(s)ds\Big)^{\frac{1}{2}}\bigg].
\enda
\end{align*}
This completes the proof of $\eqref{ABGscale}_3$.
\end{proof}

\begin{lemma}
For $t \in [0,T]$ where $T>0$ defined in \eqref{condition}, for both types of derivatives $\p^{\al}$ defined in \eqref{caseA} and \eqref{caseB}, and for $\mathrm{N} > \frac{d}{2} + 1$, we have the following estimates: 
$\bullet$ $L^\infty_tL^2_x$ estimate: For $0\leq|\al_x|\leq\mathrm{N}$, we have 
\begin{align}\label{Gembedt}
\Big\|& \|\la v \ra^{\frac{1}{2}} \p^{\al_x}\AC{\P}F^{\e}|M^{\e}|^{-1/2} \|_{L^2_{x,v}} \Big\|_{L^\infty_t} \cr 
&\les \begin{cases} \frac{1}{\sqrt{T}}\e^2\kappa^{\frac{1}{2}}\|\mathcal{D}_G^{\frac{1}{2}}\|_{L^2_T}+\frac{\sqrt{T}}{\e^{\mathfrak{n}}}\Big(\e^2\kappa^{\frac{1}{2}}\|\mathcal{D}_G^{\frac{1}{2}}\|_{L^2_T}+\e^2\|\mathcal{E}_M^{\frac{1}{2}}\mathcal{V}_{5}^{\frac{1}{2}}\|\|_{L^2_T}\Big), \quad \mbox{when} \quad 0\leq |\al|\leq \mathrm{N}-1, \\ 
\frac{1}{\sqrt{T}}\e^2\kappa^{\frac{1}{2}}\|\mathcal{D}_G^{\frac{1}{2}}\|_{L^2_T}+\frac{\sqrt{T}}{\e^{\mathfrak{n}}}\Big(\e^2\|\mathcal{D}_{top}^{\frac{1}{2}}\|_{L^2_T} + \e^2\|\mathcal{E}_M\|_{L^2_T} +\e^2\kappa^{-\frac{1}{4}}\|\mathcal{E}_M^{\frac{1}{2}}\mathcal{V}_{5}^{\frac{1}{2}}\|_{L^2_T}\Big)
, \quad \mbox{when} \quad |\al|=\mathrm{N}.
\end{cases}
\end{align}
\end{lemma}
\begin{proof}
(Proof of \eqref{Gembedt}) We proceed similarly to the proof of \eqref{Gembed}. Applying Minkowski's inequality, we have
\begin{align*}
\Big\| \|\la v \ra^{\frac{1}{2}} \p^{\al_x}\AC{\P}F^{\e}|M^{\e}|^{-1/2} \|_{L^2_{x,v}} \Big\|_{L^\infty_t}
&\les \Big\| \| \la v \ra^{\frac{1}{2}} \p^{\al_x}\AC{\P}F^{\e}|M^{\e}|^{-1/2}\|_{L^\infty_t} \Big\|_{L^2_{x,v}}.
\end{align*}
For the space-time derivative defined in \eqref{caseB}, we note that both the energy and the dissipation involve the scaled time derivative $\p_{\tilde{t}} = \e^{\mathfrak{n}} \p_t$, as specified in \eqref{N-EDdef} and \eqref{N-EDdef2}. Therefore, after applying the time embedding \eqref{tembedding}, we multiply and divide by $\e^{\mathfrak{n}}$:
\begin{align*}
\Big\| \| \la v \ra^{\frac{1}{2}} \p^{\al_x}\AC{\P}F^{\e}|M^{\e}|^{-1/2}\|_{L^\infty_t} \Big\|_{L^2_{x,v}} &\les \frac{1}{\sqrt{T}}\Big\| \la v \ra^{\frac{1}{2}} \p^{\al_x}\AC{\P}F^{\e}|M^{\e}|^{-1/2} \Big\|_{L^2_TL^2_{x,v}} \cr 
&+ \sqrt{T}\frac{1}{\e^{\mathfrak{n}}} \Big\| \la v \ra^{\frac{1}{2}} (\e^{\mathfrak{n}}\p_t)\Big(\p^{\al_x}\AC{\P}F^{\e}|M^{\e}|^{-1/2}\Big) \Big\|_{L^2_TL^2_{x,v}}.
\end{align*}
To estimate the second term, we distribute the time derivative $\p_t$ to both $\p^{\al_x}\AC{\P}F^{\e}$ and $|M^{\e}|^{-1/2}$, following the approach used in \eqref{P12def}. Then we have
\begin{align*}
\bega
\Big\|& \|\la v \ra^{\frac{1}{2}} \p^{\al_x}\AC{\P}F^{\e}|M^{\e}|^{-1/2} \|_{L^2_{x,v}} \Big\|_{L^\infty_t} \cr 
&\les \begin{cases} \frac{1}{\sqrt{T}}\e^2\kappa^{\frac{1}{2}}\|\mathcal{D}_G^{\frac{1}{2}}\|_{L^2_T}+\frac{\sqrt{T}}{\e^{\mathfrak{n}}}\Big\|\e^2\kappa^{\frac{1}{2}}\mathcal{D}_G^{\frac{1}{2}}+\e^2\mathcal{E}_M^{\frac{1}{2}}\mathcal{V}_{5}^{\frac{1}{2}}\Big\|_{L^2_T},  &\mbox{when} \quad 0\leq |\al|\leq \mathrm{N}-1, \\ 
\frac{1}{\sqrt{T}}\e^2\kappa^{\frac{1}{2}}\|\mathcal{D}_G^{\frac{1}{2}}\|_{L^2_T}+\frac{\sqrt{T}}{\e^{\mathfrak{n}}}\Big\|\e^2\mathcal{D}_{top}^{\frac{1}{2}} + \e^2\mathcal{E}_M +\e^2\kappa^{-\frac{1}{4}}\mathcal{E}_M^{\frac{1}{2}}\mathcal{V}_{5}^{\frac{1}{2}}\Big\|_{L^2_T}
,  &\mbox{when} \quad |\al|=\mathrm{N}.
\end{cases}
\enda
\end{align*}
This completes the proof of \eqref{Gembedt}.
\end{proof}

\begin{lemma}\label{L.locconD} We have the local conservation laws in Laplacian form:
\begin{align}\label{locconD}
\bega
&\p_t \rho^{\e} +u^{\e}\cdot\nabla_x\rho^{\e} + \frac{1}{\e}\nabla_x\cdot u^{\e} =0 ,\cr
&\frac{1}{k_B\mathrm{\Theta}^{\e}}\bigg[\p_tu^{\e}_i + u^{\e}\cdot \nabla_x u^{\e}_i \bigg] - \frac{1}{\mathrm{P}^{\e}|\mathrm{\Theta}^{\e}|^{\frac{1}{2}}} \eta_0 k_B^{-\frac{1}{2}} \kappa \lw(\Delta_x u^{\e}_i + \frac{1}{3}\p_i\nabla_x\cdot u^{\e}\rw) +\frac{1}{\e}\p_i(\rho^{\e}+\ta^{\e})=g_{u^{\e}_i} , \cr 
&\frac{3}{2}\bigg[\p_t \ta^{\e} +u^{\e}\cdot\nabla_x\ta^{\e} \bigg] -\frac{5}{2}\frac{|\mathrm{\Theta}^{\e}|^{\frac 1 2}}{\mathrm{P}^{\e}}\eta_1 k_B^{\frac 1 2} \kappa \Delta_x\ta^{\e} + \frac{1}{\e}\nabla_x\cdot u^{\e}  =g_{\ta^{\e}},
\enda
\end{align}
where the terms $g_{u^{\e}}=(g_{u^{\e}_1},g_{u^{\e}_2})$ and $g_{\ta^{\e}}$ are defined as follows:
\begin{align}\label{gutadef}
\bega
g_{u^{\e}_i}&:= \frac{1}{2}\eta_0 k_B^{-\frac{1}{2}} \kappa \e \frac{1}{\mathrm{P}^{\e}|\mathrm{\Theta}^{\e}|^{\frac{1}{2}}} \sum_j \p_j\ta^{\e}\lw(\p_i u^{\e}_j+\p_j u^{\e}_i-\frac 2 3 \delta_{ij}(\nabla_x\cdot u^{\e})\rw)  +\frac{1}{\e^2}\frac{1}{k_B\mathrm{P}^{\e}\mathrm{\Theta}^{\e}}\sum_j \p_j\varXi_{ij}^A , \cr 
g_{\ta^{\e}} &:= \frac{15}{4}\eta_1 k_B^{\frac 1 2} \kappa \e \frac{|\mathrm{\Theta}^{\e}|^{\frac 1 2}}{\mathrm{P}^{\e}}|\nabla_x\ta^{\e}|^2 + \frac{1}{\e^2 k_B \mathrm{P}^{\e}\mathrm{\Theta}^{\e}} \bigg(\sum_j \p_j \varXi_{j}^B-\sum_{i,j} \p_{x_i}\mathrm{U}^{\e}_j \mathbf{r}_{ij}^{\e}\bigg) .
\enda
\end{align}
For brevity, we use the following notations:
\begin{align}\label{microAB}
\bega
\varXi_{ij}^A(t,x)&:= \lw\la \mathcal{L}^{-1}\lw\{\kappa\eps \lw(\eps\p_t \AC{\P}F^{\e}+\AC{\P}(v\cdot\nabla_x \AC{\P}F^{\e})\rw)-\mathcal{N}(\AC{\P}F^{\e},\AC{\P}F^{\e})\rw\},\mathfrak{R}^{\e}_{ij}\rw\ra_{L^2_v} , \cr 
\varXi_{j}^B(t,x)&:= \lw\la \mathcal{L}^{-1}\lw\{\kappa\eps \lw(\eps\p_t \AC{\P}F^{\e}+\AC{\P}(v\cdot\nabla_x \AC{\P}F^{\e})\rw)-\mathcal{N}(\AC{\P}F^{\e},\AC{\P}F^{\e})\rw\},\mathcal{Q}^{\e}_j\rw\ra_{L^2_v}.
\enda
\end{align}
\end{lemma}
\begin{lemma}\label{lem:Fluid_dissipative}
Let $F$ be a solution to \eqref{BE}. For $F = M^{\e} + \AC{\P}F^{\e}$, we have 
\begin{align}\label{Gform}
\bega
\AC{\P}F^{\e}&=-\kappa\eps \bigg\{ \sum_{i,j}\frac{\p_i \mathrm{U}^{\e}_j}{k_B\mathrm{\Theta}^{\e}} \mathcal{L}^{-1} (\mathfrak{R}^{\e}_{ij}M^{\e})+\sum_j \frac{\p_j\mathrm{\Theta}^{\e}}{k_B|\mathrm{\Theta}^{\e}|^2}\mathcal{L}^{-1} (\mathcal{Q}^{\e}_j M^{\e}) \bigg\}\cr 
&\quad -\mathcal{L}^{-1}\lw\{\kappa\eps \lw(\eps\p_t \AC{\P}F^{\e}+\AC{\P}(v\cdot\nabla_x \AC{\P}F^{\e})\rw)-\mathcal{N}(\AC{\P}F^{\e},\AC{\P}F^{\e})\rw\} .
\enda
\end{align}
\end{lemma}
\begin{proof}
While the detailed proof can be found in \cite{LYY}, we provide a concise outline in Appendix \ref{A.Coll} for the reader's convenience.
\end{proof}
\begin{proof}[Proof of Lemma \ref{L.locconD}]
(Proof of $\eqref{locconD}_2$)  
From the definition of $\mathbf{r}_{ij}^{\e}$, note that $\langle \mathfrak{R}^{\e}_{ij},\P F \rangle = 0$. Substituting \eqref{Gform} for $\AC{\P}F^{\e}$ and using Lemma \ref{ABcomp}, we obtain
\begin{align}\label{GAform1}
\bega
\mathbf{r}_{ij}^{\e} = \langle \mathfrak{R}^{\e}_{ij},\AC{\P} F \rangle 
&= -\eta_0 k_B^{\frac 1 2} \kappa\eps^2 |\mathrm{\Theta}^{\e}|^{\frac 1 2} \lw(\p_i u^{\e}_j + \p_j u^{\e}_i - \frac 2 3 \delta_{ij}(\nabla_x \cdot u^{\e})\rw) - \varXi_{ij}^A.
\enda
\end{align}
Applying $\sum_j \p_j$ to $\mathbf{r}_{ij}^{\e}$, we get
\begin{align}\label{GA1comp}
\bega
\sum_j \p_j \mathbf{r}_{ij}^{\e} &= -\eta_0 k_B^{\frac 1 2} \kappa\eps^2 |\mathrm{\Theta}^{\e}|^{\frac 1 2} \lw(\Delta_x u^{\e}_i + \frac{1}{3}\p_i \nabla_x \cdot u^{\e} \rw) \cr 
&\quad - \eta_0 \kappa\eps^2 \sum_j \p_j \big(k_B^{\frac 1 2} |\mathrm{\Theta}^{\e}|^{\frac 1 2}\big) \lw(\p_i u^{\e}_j + \p_j u^{\e}_i - \frac 2 3 \delta_{ij} (\nabla_x \cdot u^{\e})\rw) - \sum_j \p_j \varXi_{ij}^A.
\enda
\end{align}
Substituting \eqref{GA1comp} for the term $\sum_j \p_j \mathbf{r}_{ij}^{\e}$ in $\eqref{loccon}_2$ and using $\p|\mathrm{\Theta}^{\e}|^{\frac 1 2} = \p e^{\frac{1}{2}\e\ta^{\e}} =  \frac{\e}{2}\p\ta^{\e} e^{\frac{1}{2}\e\ta^{\e}} = \frac{\e}{2} \p \ta^{\e}|\mathrm{\Theta}^{\e}|^{\frac 1 2}$ proves $\eqref{locconD}_2$. \\
(Proof of $\eqref{locconD}_3$)  
Similarly to \eqref{GAform1}, applying \eqref{Gform} and Lemma \ref{ABcomp}, we have
\begin{align*}
\bega
\mathfrak{q}_j^{\e} = \langle \mathcal{Q}^{\e}_j,\AC{\P} F \rangle 
&= -\frac{5}{2}\eta_1 \kappa\eps^2 k_B^{\frac{3}{2}} |\mathrm{\Theta}^{\e}|^{\frac{3}{2}} \p_j \ta^{\e} - \varXi_{j}^B.
\enda
\end{align*}
Applying $\sum_j \p_j$ to both sides gives
\Be\bega \label{GBcomp}
\sum_j \p_j \mathfrak{q}_j^{\e} 
&= -\frac{5}{2}\eta_1 \kappa\eps^2 k_B^{\frac{3}{2}}|\mathrm{\Theta}^{\e}|^{\frac{3}{2}} \Delta_x \ta^{\e}  -\frac{15}{4}\eta_1 \kappa\eps^3 k_B^{\frac{3}{2}}|\mathrm{\Theta}^{\e}|^{\frac{3}{2}} |\nabla_x\ta^{\e}|^2 - \sum_j \p_j \varXi_{j}^B,
\enda\Ee
where we used 
$\sum_j \p_j(|\mathrm{\Theta}^{\e}|^{\frac{3}{2}} \p_j \ta^{\e}) = |\mathrm{\Theta}^{\e}|^{\frac{3}{2}} \Delta_x \ta^{\e} + \frac{3}{2}\e|\mathrm{\Theta}^{\e}|^{\frac{3}{2}}|\nabla_x\ta^{\e}|^2$. Substituting \eqref{GBcomp} for $\sum_j \p_j \mathfrak{q}_j^{\e}$ in $\eqref{locconS}_3$ proves $\eqref{locconD}_3$.
\end{proof}

\begin{proposition}\label{meqnlem}
The vorticity $\w^{\e}$, defined in \eqref{wdef}, and $\rho^{\e} - \frac{3}{2} \ta^{\e}$ satisfy the following equations:
\begin{align}
&\p_t\w^{\e} + u^{\e}\cdot \nabla_x \w^{\e} -\frac{|\mathrm{\Theta}^{\e}|^{\frac 1 2}}{\mathrm{P}^{\e}} \eta_0 k_B^{\frac{1}{2}} \kappa \Delta_x \w^{\e} = \varPi_{\w}^{\e}, \label{weqn-g} \\
&\p_t \Big(\rho^{\e}-\frac{3}{2}\ta^{\e}\Big) + u^{\e}\cdot \nabla_x \Big(\rho^{\e}-\frac{3}{2}\ta^{\e}\Big)-\frac{|\mathrm{\Theta}^{\e}|^{\frac 1 2}}{\mathrm{P}^{\e}}\eta_1 k_B^{\frac 1 2} \kappa \Delta_x\Big(\rho^{\e}-\frac{3}{2}\ta^{\e}\Big) = \varPi_{\rho^{\e}-\frac{3}{2}\ta^{\e}}, \label{rtaeqn-g} 
\end{align}
where 
\begin{align*}
\bega
\varPi_{\w}^{\e} &:=-(\nabla_x\cdot u^{\e})\w^{\e} + k_B\mathrm{\Theta}^{\e}\nabla_x^{\perp}\frac{1}{\mathrm{P}^{\e}|\mathrm{\Theta}^{\e}|^{\frac{1}{2}}}\eta_0 k_B^{-\frac{1}{2}} \kappa \lw(\Delta_x u^{\e}_i + \frac{1}{3}\p_i\nabla_x\cdot u^{\e}\rw) \cr 
&- k_B\mathrm{\Theta}^{\e}\nabla_x^{\perp}\ta^{\e} \cdot \bigg[\nabla_x(\rho^{\e}+\ta^{\e})
+\frac{1}{\e}\frac{1}{k_B\mathrm{P}^{\e} \mathrm{\Theta}^{\e}}\sum_{j} \p_{x_j} \mathbf{r}_{ij}^{\e}\bigg] +k_B\mathrm{\Theta}^{\e}\nabla_x^{\perp}\cdot g_{u^{\e}} , \cr
\varPi_{\rho^{\e}-\frac{3}{2}\ta^{\e}} &:= -\frac{|\mathrm{\Theta}^{\e}|^{\frac 1 2}}{\mathrm{P}^{\e}}\eta_1 k_B^{\frac 1 2} \kappa \Delta_x (\rho^{\e}+\ta^{\e}) - g_{\ta^{\e}}.
\enda
\end{align*}
Here, $g_{u^{\e}}$ and $g_{\ta^{\e}}$ are defined in \eqref{gutadef}. 
\end{proposition}

\begin{proof}
{(Proof of \eqref{weqn-g})} 
\begin{align}
\bega
&\frac{1}{k_B\mathrm{\Theta}^{\e}}\bigg[\p_tu^{\e}_i + u^{\e}\cdot \nabla_x u^{\e}_i \bigg] - \frac{1}{\mathrm{P}^{\e}|\mathrm{\Theta}^{\e}|^{\frac{1}{2}}} \eta_0 k_B^{-\frac{1}{2}} \kappa \lw(\Delta_x u^{\e}_i + \frac{1}{3}\p_i\nabla_x\cdot u^{\e}\rw) +\frac{1}{\e}\p_i(\rho^{\e}+\ta^{\e})=g_{u^{\e}_i} ,
\enda
\end{align}

We take $\nabla_x^{\perp}\cdot$ to the second equation of \eqref{locconD}, then use $
\nabla_x^{\perp}\cdot (u^{\e}\cdot\nabla_x u^{\e}) 
 = u^{\e}\cdot \nabla_x \w^{\e} + (\nabla_x \cdot u^{\e})\w^{\e}$ to have
\begin{align*}
\bega
&\frac{1}{k_B\mathrm{\Theta}^{\e}}\bigg[\p_t\w^{\e} + u^{\e}\cdot \nabla_x \w^{\e} + (\nabla_x\cdot u^{\e})\w^{\e}\bigg] + \nabla_x^{\perp}\frac{1}{k_B\mathrm{\Theta}^{\e}}\cdot \bigg[\p_tu^{\e} + u^{\e}\cdot \nabla_x u^{\e}  \bigg] \cr 
&-\frac{1}{\mathrm{P}^{\e}|\mathrm{\Theta}^{\e}|^{\frac{1}{2}}}\eta_0 k_B^{-\frac{1}{2}} \kappa \Delta_x \w^{\e} -\nabla_x^{\perp}\frac{1}{\mathrm{P}^{\e}|\mathrm{\Theta}^{\e}|^{\frac{1}{2}}}\eta_0 k_B^{-\frac{1}{2}} \kappa \lw(\Delta_x u^{\e}_i + \frac{1}{3}\p_i\nabla_x\cdot u^{\e}\rw)
= \nabla_x^{\perp}\cdot g_{u^{\e}} .
\enda
\end{align*}
To second term in the first line, we use $\p\frac{1}{\mathrm{\Theta}^{\e}} = \p e^{-\e\ta^{\e}} =  -\e\p\ta^{\e} e^{-\e\ta^{\e}} = -\e \p \ta^{\e} \frac{1}{\mathrm{\Theta}^{\e}}$ and apply the second equation of \eqref{locconNew} to get 
\begin{align*}
\bega
\nabla_x^{\perp}\frac{1}{k_B\mathrm{\Theta}^{\e}}\cdot \bigg[\p_tu^{\e} + u^{\e}\cdot \nabla_x u^{\e}  \bigg] &= -\frac{\e \nabla_x^{\perp}\ta^{\e}}{k_B\mathrm{\Theta}^{\e}} \cdot \bigg[ \p_tu^{\e} + u^{\e}\cdot \nabla_x u^{\e}\bigg] \cr 
&= \nabla_x^{\perp}\ta^{\e} \cdot \bigg[\nabla_x(\rho^{\e}+\ta^{\e})
+\frac{1}{\e}\frac{1}{k_B\mathrm{P}^{\e} \mathrm{\Theta}^{\e}}\sum_{j} \p_{x_j} \mathbf{r}_{ij}^{\e}\bigg].
\enda
\end{align*}
Thus, we have 
\begin{align*}
\bega
&\frac{1}{k_B\mathrm{\Theta}^{\e}}\bigg[\p_t\w^{\e} + u^{\e}\cdot \nabla_x \w^{\e} + (\nabla_x\cdot u^{\e})\w^{\e}\bigg] -\frac{1}{\mathrm{P}^{\e}|\mathrm{\Theta}^{\e}|^{\frac{1}{2}}}\eta_0 k_B^{-\frac{1}{2}} \kappa \Delta_x \w^{\e} =\nabla_x^{\perp}\cdot g_{u^{\e}}
\cr 
&+ \nabla_x^{\perp}\frac{1}{\mathrm{P}^{\e}|\mathrm{\Theta}^{\e}|^{\frac{1}{2}}}\eta_0 k_B^{-\frac{1}{2}} \kappa \lw(\Delta_x u^{\e}_i + \frac{1}{3}\p_i\nabla_x\cdot u^{\e}\rw) - \nabla_x^{\perp}\ta^{\e} \cdot \bigg[\nabla_x(\rho^{\e}+\ta^{\e})
+\frac{1}{\e}\frac{1}{k_B\mathrm{P}^{\e} \mathrm{\Theta}^{\e}}\sum_{j} \p_{x_j} \mathbf{r}_{ij}^{\e}\bigg].
\enda
\end{align*}
Multiplying $k_B\mathrm{\Theta}^{\e}$, gives desired result \eqref{weqn-g}. \\
{(Proof of \eqref{rtaeqn-g})} We subtract  $\eqref{locconD}_3$ from $\eqref{locconD}_1$ to have 
\begin{align*}
\bega
\p_t \Big(\rho^{\e}-\frac{3}{2}\ta^{\e}\Big) + u^{\e}\cdot \nabla_x \Big(\rho^{\e}-\frac{3}{2}\ta^{\e}\Big)+\frac{5}{2}\eta_1 k_B^{\frac 1 2} \kappa \Delta_x\ta^{\e} = -g_{\ta^{\e}}.
\enda
\end{align*}
For the term \( \frac{5}{2} \Delta_x \ta^{\e} \), decomposing \( \frac{5}{2} \ta^{\e} = -(\rho^{\e} - \frac{3}{2} \ta^{\e}) + (\rho^{\e} + \ta^{\e}) \) gives \eqref{rtaeqn-g}. 
\end{proof}

\begin{lemma}\label{L.rutaL2} For the solution $(\rho^{\e},u^{\e},\ta^{\e})$ of equations \eqref{locconD}, the following equality holds:
\begin{align*}
\bega
&\frac{1}{2}\frac{d}{dt}\bigg(\|\rho^{\e}\|_{L^2_x}^2+\bigg\|\frac{u^{\e}}{\sqrt{k_B\mathrm{\Theta}^{\e}}}\bigg\|_{L^2_x}^2+\frac{3}{2}\|\ta^{\e}\|_{L^2_x}^2\bigg) +\eta_0 k_B^{-\frac{1}{2}} \kappa \Big(\|\nabla_x u^{\e}\|_{L^2_x}^2 + \frac{1}{3}\|\nabla_x \cdot u^{\e}\|_{L^2_x}^2\Big) + \frac{5}{2}\eta_1 k_B^{\frac{1}{2}} \kappa \|\nabla_x \ta^{\e}\|_{L^2_x}^2 \cr 
&=  \frac{1}{2}\int_{\Omega}(\nabla_x\cdot u^{\e})\bigg(|\rho^{\e}|^2+\frac{|u^{\e}|^2}{k_B\mathrm{\Theta}^{\e}}+\frac{3}{2}|\ta^{\e}|^2\bigg) dx +\int_{\Omega} u^{\e} g^{\e}_{u} dx + \int_{\Omega} \ta^{\e} g^{\e}_{\ta} dx \cr 
& + \frac{1}{2}\int_{\Omega} \frac{2}{3}\frac{1}{k_B\mathrm{\Theta}^{\e}}\bigg(\nabla_x\cdot u^{\e}  +\frac{1}{\e}\frac{1}{k_B\mathrm{P}^{\e}\mathrm{\Theta}^{\e} }\sum_{j} \p_{x_j} \mathfrak{q}_j^{\e} + \frac{1}{\e}\frac{1}{k_B\mathrm{P}^{\e}\mathrm{\Theta}^{\e} }\sum_{i,j} \p_{x_i}\mathrm{U}^{\e}_j \mathbf{r}_{ij}^{\e}\bigg)|u^{\e}|^2 dx .
\enda
\end{align*}
Here, $g_{u}$ and $g_{\ta^{\e}}$ are defined in \eqref{gutadef}.
\end{lemma}
\begin{proof}
A standard energy estimate applied to \eqref{locconD} yields 
\begin{align}\label{rutaL2-a}
\bega
&\frac{1}{2}\frac{d}{dt}\|\rho^{\e}\|_{L^2_x}^2 = \frac{1}{2}\int_{\Omega} (\nabla_x \cdot u^{\e}) |\rho^{\e}|^2 dx -\frac{1}{\e} \int_{\Omega} (\nabla_x \cdot u^{\e}) \rho^{\e} dx , \cr
&\frac{1}{2}\frac{d}{dt}\frac{3}{2}\|\ta^{\e}\|_{L^2_x}^2 + \frac{5}{2}\eta_1 k_B^{\frac{1}{2}} \kappa \|\nabla_x \ta^{\e}\|_{L^2_x}^2 = \frac{1}{2}\int_{\Omega} (\nabla_x \cdot u^{\e}) \frac{3}{2}|\ta^{\e}|^2 dx -\frac{1}{\e} \int_{\Omega} (\nabla_x \cdot u^{\e}) \ta^{\e} dx + \int_{\Omega} \ta^{\e} g^{\e}_{\ta} dx,
\enda
\end{align}
and
\begin{align*}
\bega
&\sum_i \int_{\Omega} \frac{1}{k_B\mathrm{\Theta}^{\e}}\bigg[\p_tu^{\e}_i + u^{\e}\cdot \nabla_x u^{\e}_i \bigg]u^{\e}_i dx + \eta_0 k_B^{-\frac{1}{2}} \kappa (\|\nabla_x u^{\e}\|_{L^2_x}^2 + \frac{1}{3}\|\nabla_x \cdot u^{\e}\|_{L^2_x}^2) \cr 
&= \frac{1}{\e} \int_{\Omega} (\nabla_x \cdot u^{\e}) (\rho^{\e}+\ta^{\e}) dx + \int_{\Omega} u^{\e} g^{\e}_{u} dx . 
\enda
\end{align*}
For the first term of the energy estimate of $u_i$, we use the following integration by parts
\begin{align*}
\bega
&\sum_i \int_{\Omega} \frac{1}{k_B\mathrm{\Theta}^{\e}}(\p_tu^{\e}_i)u^{\e}_i dx 
= \frac{d}{dt}\frac{1}{2}\int_{\Omega} \frac{1}{k_B\mathrm{\Theta}^{\e}}|u^{\e}|^2 dx -\frac{1}{2}\int_{\Omega} \p_t\bigg(\frac{1}{k_B\mathrm{\Theta}^{\e}}\bigg)|u^{\e}|^2 dx, \cr 
&\sum_i \int_{\Omega} \frac{1}{k_B\mathrm{\Theta}^{\e}}(u^{\e}\cdot \nabla_x u^{\e}_i)u^{\e}_i dx 
= -\frac{1}{2}\int_{\Omega} \frac{1}{k_B\mathrm{\Theta}^{\e}}(\nabla_x\cdot u^{\e})  |u^{\e}|^2 dx -\frac{1}{2}\int_{\Omega} u^{\e}\cdot \nabla_x\bigg(\frac{1}{k_B\mathrm{\Theta}^{\e}}\bigg)  |u^{\e}|^2 dx,
\enda
\end{align*}
to have
\begin{align}\label{rutaL2-b}
\bega
&\frac{1}{2}\bigg\|\frac{u^{\e}}{\sqrt{k_B\mathrm{\Theta}^{\e}}}\bigg\|_{L^2_x}^2 + \eta_0 k_B^{-\frac{1}{2}} \kappa \Big(\|\nabla_x u^{\e}\|_{L^2_x}^2 + \frac{1}{3}\|\nabla_x \cdot u^{\e}\|_{L^2_x}^2\Big)= \frac{1}{2} \int_{\Omega} (\nabla_x\cdot u^{\e})\frac{|u^{\e}|^2}{k_B\mathrm{\Theta}^{\e}} dx \cr 
&\quad + \frac{1}{\e} \int_{\Omega} (\nabla_x \cdot u^{\e}) (\rho^{\e}+\ta^{\e}) dx +\frac{1}{2}\int_{\Omega} (\p_t+u^{\e}\cdot \nabla_x)\bigg(\frac{1}{k_B\mathrm{\Theta}^{\e}}\bigg)|u^{\e}|^2 dx+ \int_{\Omega} u^{\e} g^{\e}_{u} dx . 
\enda
\end{align}
For the term $(\p_t+u^{\e}\cdot \nabla_x)\frac{1}{k_B\mathrm{\Theta}^{\e}}$, using $\eqref{locconNew}_3$, we have 
\begin{align}\label{rutaL2-c}
\bega
(\p_t+u^{\e}\cdot \nabla_x)\bigg(\frac{1}{k_B\mathrm{\Theta}^{\e}}\bigg) &= -\frac{1}{k_B\mathrm{\Theta}^{\e}} (\e\p_t\ta^{\e}+\e u^{\e}\cdot \nabla_x\ta^{\e}) \cr 
&= \frac{2}{3}\frac{1}{k_B\mathrm{\Theta}^{\e}}\Big(\nabla_x\cdot u^{\e}  +\frac{1}{\e}\frac{1}{k_B\mathrm{P}^{\e}\mathrm{\Theta}^{\e} }\sum_{j} \p_{x_j} \mathfrak{q}_j^{\e} + \frac{1}{\e}\frac{1}{k_B\mathrm{P}^{\e}\mathrm{\Theta}^{\e} }\sum_{i,j} \p_{x_i}\mathrm{U}^{\e}_j \mathbf{r}_{ij}^{\e}\Big).
\enda
\end{align}
Combining \eqref{rutaL2-a} and \eqref{rutaL2-b} then using \eqref{rutaL2-c}, we obtain the desired result. 
\end{proof}

\subsection{Macroscopic Forcing Terms}\label{Sec.remain}

In this subsection, we estimate \(\varXi_{ij}^A\) and \(\varXi_{j}^B\). Recall the definitions of \(\varXi_{ij}^A\) and \(\varXi_{j}^B\) in \eqref{microAB}. Note that \(\mathcal{L}^{-1}\) is a symmetric operator in the following sense $\int \mathcal{L}^{-1}(H) G |M^{\e}|^{-1} dv = \int H \mathcal{L}^{-1}(G) |M^{\e}|^{-1} dv$ as shown in Lemma \ref{Llem1} (i) and (iii). 
By transferring the \(\mathcal{L}^{-1}\) operator to \(\mathfrak{R}^{\e}_{ij}M^{\e}\) and \(\mathcal{Q}^{\e}_{j}M^{\e}\), we can express \(\varXi_{ij}^A\) and \(\varXi_{j}^B\) as:
\begin{align}\label{varXidef}
\bega
\varXi_{ij}^A(t,x) &= \int_{\R^3} \bigg(\kappa\eps \big[\eps \p_t \AC{\P} F + \AC{\P} (v \cdot \nabla_x \AC{\P} F)\big] - \mathcal{N}(\AC{\P} F, \AC{\P} F)\bigg) \mathcal{L}^{-1}(\mathfrak{R}^{\e}_{ij}M^{\e}) |M^{\e}|^{-1} \, dv, \cr
\varXi_{j}^B(t,x) &= \int_{\R^3} \bigg(\kappa\eps \big[\eps \p_t \AC{\P} F + \AC{\P} (v \cdot \nabla_x \AC{\P} F)\big] - \mathcal{N}(\AC{\P} F, \AC{\P} F)\bigg) \mathcal{L}^{-1}(\mathcal{Q}^{\e}_{j}M^{\e}) |M^{\e}|^{-1} \, dv.
\enda
\end{align}

The goal of this subsection is to estimate \(\varPi_{\w}\) and \(\varPi_{\rho - \frac{3}{2}\ta}\), as defined in \eqref{gwrta}, using the following proposition.


\begin{proposition}\label{P.Xi}
Assume \eqref{condition} holds. For \(1 \leq |\al_x| \leq \mathrm{N}\), the quantities \(\varXi_{ij}^A\) and \(\varXi_{j}^B\), as defined in \eqref{microAB}, satisfy
\begin{align*}
\bega
\|\p^{\al_x}&\varXi_{ij}^A\|_{L^2_tL^2_x},~ \|\p^{\al_x}\varXi_{j}^B\|_{L^2_tL^2_x} \les  \e^3\kappa\varpi(\e,\kappa)\bigg(\Big(\int_0^t\mathcal{D}_{tot}(s)ds\Big)^{\frac{1}{2}}+\Big(\int_0^t\mathcal{E}_M^2(s)ds\Big)^{\frac{1}{2}}\bigg) \cr 
&+ \e^3\varpi(\e,\kappa)\bigg(\frac{1}{\sqrt{T}}+\sqrt{T}\bigg)\bigg(\int_0^t\mathcal{D}_{tot}(s)ds + \int_0^t\mathcal{E}_M^2(s)ds\bigg) \cr 
&+ \e^3\varpi(\e,\kappa)\kappa^{-\frac{1}{2}}\bigg(\frac{1}{\sqrt{T}}+\sqrt{T}\bigg) \bigg[\int_0^t \bigg(\mathcal{E}_{tot}^2(t) +\e^{\frac{3}{2}}\kappa \mathcal{E}_{tot}\Big(\mathcal{D}_{tot}(t)+\mathcal{E}_M^2(t)\Big)\bigg)ds \cr 
&+\sup_{0\leq s\leq t}\mathcal{E}_{tot}(s) e^{-c_1\e^{-{\frac{1}{\ell}}}}\bigg(\e^2\kappa\sum_{0\leq|\al|\leq\mathrm{N}+1}\kappa^{(|\al|-\mathrm{N})_+}\|\p^{\al}h_0\|_{L^\infty_{x,v}}^2 + \frac{1}{(\e\kappa)^{d}}\int_0^t \mathcal{E}_{tot}(s) ds\bigg)\bigg],
\enda
\end{align*}
for any \(t \in [0, T]\).
\end{proposition}

To prove Proposition \ref{P.Xi}, we divide the proof into the following three lemmas.

\begin{lemma}\label{L.Xi1}
Assume \eqref{condition} holds. For \(0 \leq |\al_x| \leq \mathrm{N}\), the quantities \(\varXi_{ij}^A\) and \(\varXi_{j}^B\) satisfy:
\begin{align}\label{XiL-1}
\bega
\|\p^{\al_x} \varXi_{ij}^A\|_{L^2_x},~ &\|\p^{\al_x} \varXi_{j}^B\|_{L^2_x} 
\les \sum_{0 \leq \beta_x \leq \al_x} \bigg[\e^3\kappa\varpi(\e,\kappa)\Big(\mathcal{D}_{tot}^{\frac{1}{2}}(t)+\mathbf{1}_{|\al|\geq2}\mathcal{E}_M(t)\Big) \cr 
&\quad + \sum_{0\leq|\gamma_x|\leq\lfloor|\beta_x|/2\rfloor} \big\|\la v \ra^{\frac{1}{2}} \p^{\gamma_x}\AC{\P}F^{\e}|M^{\e}|^{-1/2} \big\|_{L^\infty_xL^2_v}  \big\|\la v \ra^{\frac{1}{2}} \p^{\beta_x-\gamma_x}\AC{\P}F^{\e}|M^{\e}|^{-1/2} \big\|_{L^2_{x,v}} \bigg] \cr
& \times \bigg(\mathbf{1}_{\al_x = \beta_x} + \mathbf{1}_{\al_x > \beta_x} \bigg(\int_{\Omega \times \R^3} \la v \ra \Big|\p^{\al_x - \beta_x}\Big(\mathcal{L}^{-1}(\hat{X}M^{\e}) |M^{\e}|^{-1}\Big)\Big|^2 M^{\e} \, dv \, dx \bigg)^{\frac{1}{2}} \bigg),
\enda
\end{align}
where \(\hat{X} = \mathfrak{R}^{\e}_{ij}\) or \(\mathcal{Q}^{\e}_j\), depending on whether the estimate is for \(\varXi_{ij}^A\) or \(\varXi_{j}^B\), respectively.
\end{lemma}

\begin{proof}
The proofs for \(\varXi_{ij}^A\) and \(\varXi_{j}^B\) are identical. For simplicity, we consider only the \(\varXi_{ij}^A\) case. 
We distribute the derivative \(\p^{\al_x}\) to \(\varXi_{ij}^A\) defined in \eqref{varXidef}. Using \eqref{nonlin} for the \(\mathcal{N}(\AC{\P} F, \AC{\P} F)\) term and applying Hölder's inequality, we obtain
\[
\|\p^{\al_x} \varXi_{ij}^A\|_{L^2_x} \les \sum_{0 \leq \beta_x \leq \al_x} \big(A_1^{\beta_x} + A_2^{\beta_x} + A_3^{\beta_x}\big) A_4^{\al_x - \beta_x},
\]
where the terms are defined as follows:
\begin{align*}
&A_1^{\beta_x} =  \e^2\kappa \lw(\int_{\Omega\times\R^3} |\p^{\beta_x}\p_t \AC{\P}F^{\e} |^2 |M^{\e}|^{-1} dvdx \rw)^{\frac{1}{2}} , \cr
&A_2^{\beta_x} =  \e\kappa \lw(\int_{\Omega\times\R^3} \frac{1}{\la v\ra}|\p^{\beta_x}\AC{\P}(v\cdot\nabla_x \AC{\P}F^{\e})|^2 |M^{\e}|^{-1} dvdx \rw)^{\frac{1}{2}}, \cr
&A_3^{\beta_x} = \sum_{0\leq\gamma_x\leq\beta_x} \lw(\int_{\Omega}\Big(\int_{\R^3} \la v\ra |\p^{\gamma_x}\AC{\P}F^{\e} |^2 |M^{\e}|^{-1} dv \Big) \Big(\int_{\R^3} \la v\ra |\p^{\beta_x-\gamma_x}\AC{\P}F^{\e} |^2 |M^{\e}|^{-1} dv\Big) dx\rw)^{\frac{1}{2}}, \cr
&A_4^{\al_x-\beta_x} = \lw(\int_{\Omega\times\R^3} \la v\ra \Big|\p^{\al_x-\beta_x}\Big(\mathcal{L}^{-1}(\mathfrak{R}^{\e}_{ij}M^{\e}) |M^{\e}|^{-1}\Big)\Big|^2 M^{\e} dvdx \rw)^{\frac{1}{2}} . 
\end{align*}
Recall that since we used different dissipation depending on the number of derivatives $\quad 0\leq|\al|\leq\mathrm{N}$ and $|\al|=\mathrm{N}+1$, we obtained the following estimates in \eqref{albe.claim}: 
\begin{align}\label{Xi-a}
\bega
\|\p^{\al}\AC{\P}F^{\e}|M^{\e}|^{-\frac{1}{2}}\|_{L^2_{x,v}} \leq \begin{cases} \e^2\kappa^{\frac{1}{2}}\mathcal{D}_G^{\frac{1}{2}}, \quad &\mbox{when} \quad 0\leq|\al|\leq\mathrm{N}, \\ 
\e^2(\mathcal{D}_{top}^{\frac{1}{2}}(t)+\mathbf{1}_{|\al|\geq2}\mathcal{E}_M(t)), \quad &\mbox{when} \quad |\al|=\mathrm{N}+1.
\end{cases}
\enda
\end{align}
For the term $A_1^{\beta_x}$, we multiply and divide by $\e \varpi^{-1}(\e, \kappa)$ in front of $\p_t$. Using \eqref{Xi-a}, we derive
\begin{align}\label{Xi-1r}
\bega
A_1^{\beta_x} &\les  \frac{\e^2\kappa}{\e\varpi^{-1}(\e,\kappa)} \lw(\int_{\Omega\times\R^3} |\p^{\beta_x}(\e\varpi^{-1}(\e,\kappa)\p_t) \AC{\P}F^{\e} |^2 |M^{\e}|^{-1} dvdx \rw)^{\frac{1}{2}} \cr 
&\les \begin{cases}
\e^3\kappa^{\frac{3}{2}}\varpi(\e,\kappa)\mathcal{D}_G^{\frac{1}{2}}, \quad &\mbox{when} \quad 0\leq|\beta_x|\leq\mathrm{N}-1, \\ 
\e^3\kappa\varpi(\e,\kappa)(\mathcal{D}_{top}^{\frac{1}{2}}(t)+\mathbf{1}_{|\al|\geq2}\mathcal{E}_M(t)), \quad &\mbox{when} \quad |\beta_x|=\mathrm{N}. \end{cases}
\enda
\end{align}
For the term $A_3^{\beta_x}$, we take the $L^\infty_x$-norm for the lower-regularity part. Without loss of generality, we assume $|\gamma_x| \leq \lfloor|\beta_x|/2\rfloor$. Then, we can write the upper bound of $A_3^{\beta_x}$ as
\begin{align}\label{Xi-3r}
\bega
A_3^{\beta_x} &\les \sum_{0\leq|\gamma_x|\leq\lfloor|\beta_x|/2\rfloor} \big\|\la v \ra^{\frac{1}{2}} \p^{\gamma_x}\AC{\P}F^{\e}|M^{\e}|^{-1/2} \big\|_{L^\infty_xL^2_v}  \big\|\la v \ra^{\frac{1}{2}} \p^{\beta_x-\gamma_x}\AC{\P}F^{\e}|M^{\e}|^{-1/2} \big\|_{L^2_{x,v}}. 
\enda
\end{align}
For the term $A_2^{\beta_x}$, we claim that $A_2^{\beta_x}$ satisfies the following estimates:
\begin{align}\label{A2claim}
\bega
A_2^{\beta_x}&\les \begin{cases}
\e^3\kappa^{\frac{3}{2}}\mathcal{D}_G^{\frac{1}{2}}, \quad &\mbox{when} \quad 0\leq|\beta_x|\leq\mathrm{N}-1, \\ 
\e^3\kappa(\mathcal{D}_{top}^{\frac{1}{2}}(t)+\mathbf{1}_{|\al|\geq2}\mathcal{E}_M(t)), \quad &\mbox{when} \quad |\beta_x|=\mathrm{N}, \end{cases}
\enda
\end{align}
We first decompose $A_2^{\beta_x}=A_{2,2}^{\beta_x}+A_{2,3}^{\beta_x}+A_{2,1}^{\beta_x}$ where
\begin{align*}
A_{2,1}^{\beta_x} &= \e \kappa \lw(\int_{\Omega \times \R^3} \frac{1}{\la v \ra} |\p^{\beta_x} (v \cdot \nabla_x \AC{\P} F)|^2 |M^{\e}|^{-1} dv dx \rw)^{\frac{1}{2}}, \cr
A_{2,2}^{\beta_x} &= \e \kappa \lw(\int_{\Omega \times \R^3} \frac{1}{\la v \ra} |\P (v \cdot \nabla_x \p^{\beta_x} \AC{\P} F)|^2 |M^{\e}|^{-1} dv dx \rw)^{\frac{1}{2}}, \cr
A_{2,3}^{\beta_x} &= \e \kappa \lw(\int_{\Omega \times \R^3} \frac{1}{\la v \ra} \Big|\P (v \cdot \nabla_x \p^{\beta_x} \AC{\P} F) - \p^{\beta_x} \P (v \cdot \nabla_x \AC{\P} F)\Big|^2 |M^{\e}|^{-1} dv dx \rw)^{\frac{1}{2}}.
\end{align*}
Here, we used the following computation by the definition $\AC{\P} = \II - \P$:
\begin{align*}
\bega
\p^{\beta_x}\AC{\P}(v\cdot\nabla_x \AC{\P}F^{\e}) = \p^{\beta_x}(v\cdot\nabla_x \AC{\P}F^{\e}) -\P(v\cdot\nabla_x \p^{\beta_x}\AC{\P}F^{\e}) +\Big[\P(v\cdot\nabla_x \p^{\beta_x}\AC{\P}F^{\e}) -\p^{\beta_x}\P(v\cdot\nabla_x \AC{\P}F^{\e})\Big],
\enda
\end{align*}
Using \eqref{Xi-a}, we can easily see that the term $A_{2,1}^{\beta_x}$ can be bounded by the right-side of \eqref{A2claim}. 
For the term $A_{2,2}^{\beta_x}$, applying the H\"{o}lder inequality to the velocity integral of $\P$ gives
\begin{align*}
\bega
A_{2,2}^{\beta_x} &\leq C \e\kappa\lw(\int_{\Omega\times\R^3} |\nabla_x \p^{\beta_x}\AC{\P}F^{\e}|^2 |M^{\e}|^{-1} dvdx\rw)^{\frac{1}{2}}.
\enda
\end{align*}
Again using \eqref{Xi-a}, the term $A_{2,2}^{\beta_x}$ also satisfies the inequality \eqref{A2claim}. 
For the third term $A_{2,3}^{\beta_x}$, we compute the commutator between $\p^{\beta_x}$ and $\P$ based on the definition of $\P$ in \eqref{Pdef}. Setting $H= v\cdot\nabla_x \AC{\P}F^{\e}$, we write
\begin{align*}
\bega
\P\p^{\al_x}H&-\p^{\al_x}\P H= - \sum_{i=0}^4 \sum_{0<\beta_x\leq\al_x}\sum_{0\leq\gamma_x\leq\beta_x} \binom{\al_x}{\beta_x}\binom{\beta_x}{\gamma_x}\lw\la \p^{\al_x-\beta_x}H, \p^{\beta_x-\gamma_x}\bigg(\frac{e^{\e}_i}{M^{\e}}\bigg) \rw\ra \p^{\gamma_x}e^{\e}_i \cr 
&\les \sum_{i=0}^4 \sum_{\substack{0<\beta_x\leq\al_x \\ 0\leq\gamma_x\leq\beta_x}} \lw(\int_{\R^3} \frac{1}{\la v\ra}|\p^{\al_x-\beta_x}H|^2 |M^{\e}|^{-1} dv\rw)^{\frac{1}{2}} \lw(\int_{\R^3} \la v\ra\bigg|\p^{\beta_x-\gamma_x}\bigg(\frac{e^{\e}_i}{M^{\e}}\bigg)\bigg|^2 M^{\e} dv\rw)^{\frac{1}{2}}|\p^{\gamma_x}e^{\e}_i|,
\enda
\end{align*}
where we applied the H\"{o}lder inequality.
Next, we note that either $\p^{\beta_x}$ or $\p^{\beta_x-\gamma_x}$ must include at least one derivative, as $0<\beta_x$. Based on the definition of $e^{\e}_i$ in \eqref{5mom}, $e^{\e}_i/M^{\e}$ is a velocity polynomial, \textcolor{blue}{and $e^{\e}_i$ has a form as polynomial multiplied by $M^{\e}$. We estimated $\p^{\al}M^{\e}$ in Lemma \ref{Mal}. By the same way, we estimated $\p^{\al}$ of polynomial multiplied by $M^{\e}$ in \eqref{al-vM}. }
Following a similar argument as in \eqref{Phi2}, we obtain
\begin{align*}
\bega
\lw(\int_{\R^3} \la v\ra\bigg|\p^{\beta_x-\gamma_x}\bigg(\frac{e^{\e}_i}{M^{\e}}\bigg)\bigg|^2 M^{\e} dv\rw)^{\frac{1}{2}}|\p^{\gamma_x}e^{\e}_i| \les \e \mathcal{E}_M^{\frac{1}{2}} (1+|v|^{2+2|\gamma_x|})M^{\e} \les (1+|v|^{2+2|\gamma_x|})M^{\e},
\enda
\end{align*}
where we used $\e \mathcal{E}_M^{\frac{1}{2}} \leq C$ by \eqref{condition}. Reverting the notation $H= v\cdot\nabla_x \AC{\P}F^{\e}$ gives 
\begin{align*}
A_{2,3}^{\beta_x} 
&\les \e \kappa \lw(\int_{\Omega} \sum_{0<\gamma_x\leq\beta_x} \int_{\R^3} \frac{1}{\la v\ra}|\p^{\beta_x-\gamma_x}H|^2 |M^{\e}|^{-1} dv \int_{\R^3}(1+|v|^{2+2|\gamma_x|})^2|M^{\e}|^2 dvdx \rw)^{\frac{1}{2}}
\end{align*}
Since $\int_{\R^3}(1+|v|^{2+2|\gamma_x|})^2|M^{\e}|^2 dv\leq C$, by using \eqref{Xi-a}, $A_{2,3}^{\beta_x}$ satisfies the upper bound \eqref{A2claim}. 
For the last term $A_4^{\al_x-\beta_x}$, when $\al_x=\beta_x$, applying \eqref{L-1L2} gives
\begin{align}\label{Xi-4r}
\bega
\int_{\Omega\times\R^3} \la v\ra \Big|\Big(\mathcal{L}^{-1}(\mathfrak{R}^{\e}_{ij}M^{\e}) |M^{\e}|^{-1}\Big)\Big|^2 M^{\e} dvdx \leq \int_{\Omega\times\R^3}  |\mathfrak{R}^{\e}_{ij}|^2 M^{\e} dvdx \leq C.
\enda
\end{align}
Combining \eqref{Xi-1r}, \eqref{Xi-3r}, \eqref{A2claim}, and \eqref{Xi-4r}, we conclude the proof.
\end{proof}

We note that in the last line of \eqref{XiL-1}, when $\al_x > \beta_x$, the term $\p^{\al_x-\beta_x}\mathcal{L}^{-1}(\mathfrak{R}^{\e}_{ij}M^{\e}) |M^{\e}|^{-1}$ generates at least one factor of $\e$. However, it is essential to carefully compute the commutator between $\p^{\al_x}$ and $\mathcal{L}^{-1}$. In the following lemma, we iteratively estimate this commutator. Here, \(H\)  can denote either \(\mathfrak{R}^{\e}_{ij}M^{\e}\) or \(\mathcal{Q}^{\e}_jM^{\e}\).


\begin{lemma}\label{L-1G2}
Assume \eqref{condition} holds. Then for any \(n \geq 0\) and any function \(H \in (\mathrm{Ker}(\mathcal{L}))^\perp\) that ensures the right-hand side is bounded, we have
\begin{align}\label{L-1}
\bega
\int_{\Omega\times\R^3} \la v \ra^{n+2}|\p^{\al_x} \mathcal{L}^{-1}H|^2& |M^{\e}|^{-1} dvdx \les \int_{\Omega\times\R^3} \la v \ra^n|\p^{\al_x}H|^2 |M^{\e}|^{-1} dvdx \cr
&+ \sum_{0<|\beta_x|\leq\mathrm{N}-1}\Big(\e^2\kappa^{-1}\mathcal{E}_{top}\Big) \int_{\Omega\times\R^3}\la v \ra^{n}|\p^{\al_x-\beta_x}H(v)|^2|M^{\e}|^{-1}dvdx \cr 
&+ \mathbf{1}_{|\al_x|=\mathrm{N}}\Big(\e^2\mathcal{E}_M\Big)\Big\|\int_{\R^3}\la v \ra^{n}|H(v)|^2|M^{\e}|^{-1}dv\Big\|_{L^\infty_x},
\enda
\end{align}
for \(1 \leq |\al_x| \leq \mathrm{N}\).
\end{lemma}

\begin{proof}
We consider a function \( H \in (\mathrm{Ker}(\mathcal{L}))^\perp \) satisfying \( \mathcal{L}^{-1}H = G \) for \( G \in (\mathrm{Ker}(\mathcal{L}))^\perp \). Taking \(\p^{\al_x}\) of \(\mathcal{L}G\), we have
\begin{align*}
\p^{\al_x}\mathcal{L}G = \mathcal{L} \p^{\al_x} G -2\sum_{0<\beta_x\leq\al_x}\binom{\al_x}{\beta_x}\mathcal{N}(\p^{\beta_x}M,\p^{\al_x-\beta_x}G).
\end{align*}
Since \(G \in (\mathrm{Ker}(\mathcal{L}))^\perp\) implies \(\p^{\al}G \in (\mathrm{Ker}(\mathcal{L}))^\perp\), applying \(\mathcal{L}^{-1}\AC{\P}\) to both sides gives
\begin{align*}
\mathcal{L}^{-1}\AC{\P}\p^{\al_x}\mathcal{L}G = \p^{\al_x} G -2\sum_{0<\beta_x\leq\al_x}\binom{\al_x}{\beta_x}\mathcal{L}^{-1}\mathcal{N}(\p^{\beta_x}M,\p^{\al_x-\beta_x}G).
\end{align*}
Substituting \(\mathcal{L}^{-1}H\) for \(G\), we obtain
\begin{align*}
\mathcal{L}^{-1}\AC{\P}\p^{\al_x}H = \p^{\al_x} \mathcal{L}^{-1}H  +2\sum_{0<\beta_x\leq\al_x}\binom{\al_x}{\beta_x}\mathcal{L}^{-1}\mathcal{N}(\p^{\beta_x}M^{\e},\p^{\al_x-\beta_x}\mathcal{L}^{-1}H).
\end{align*}
Since \(H \in (\mathrm{Ker}(\mathcal{L}))^\perp\), we also have \(\p^{\al_x}H \in (\mathrm{Ker}(\mathcal{L}))^\perp\) by Lemma \ref{HH}, which implies \(\AC{\P}\p^{\al_x}H = \p^{\al_x}H\). Squaring and multiplying both sides by \(\la v \ra^{n+2}|M^{\e}|^{-1}\), we get
\begin{align*}
\int_{\R^3} &\la v \ra^{n+2}|\p^{\al_x} \mathcal{L}^{-1}H|^2 |M^{\e}|^{-1} dv \les \int_{\R^3} \la v \ra^{n+2}|\mathcal{L}^{-1}\p^{\al_x}H|^2 |M^{\e}|^{-1} dv \\
&+ \sum_{0<\beta_x\leq\al_x}\int_{\R^3} \la v \ra^{n+2}|\mathcal{L}^{-1}\mathcal{N}(\p^{\beta_x}M^{\e},\p^{\al_x-\beta_x}\mathcal{L}^{-1}H)|^2 |M^{\e}|^{-1} dv.
\end{align*}
Using \eqref{L-1L2}, we have
\begin{align*}
\int_{\Omega\times\R^3} \la v \ra^{n+2}&|\p^{\al_x} \mathcal{L}^{-1}H|^2 |M^{\e}|^{-1} dvdx \les B_1^{\al_x} + B_2^{\al_x},
\end{align*}
where
\begin{align*}
B_1^{\al_x} &= \int_{\Omega\times\R^3} \la v \ra^n|\p^{\al_x}H|^2 |M^{\e}|^{-1} dvdx, \cr 
B_2^{\al_x} &= \sum_{0<\beta_x\leq\al_x}\int_{\Omega\times\R^3} \la v \ra^n|\mathcal{N}(\p^{\beta_x}M^{\e},\p^{\al_x-\beta_x}\mathcal{L}^{-1}H)|^2 |M^{\e}|^{-1} dvdx.
\end{align*}
We note that \(B_1^{\al_x}\) is the leading-order term, while \(B_2^{\al_x}\) generates additional factors of \(\e\) due to \(\p^{\beta_x}M^{\e}\). Applying \eqref{nonlin2} and \eqref{L-1L2}, we obtain
\begin{align}\label{L-13a}
\bega
B_2^{\al_x} &\les \sum_{0<\beta_x\leq\al_x} \int_{\Omega}\int_{\R^3}\la v \ra^{n+2}|\p^{\beta_x}M^{\e}|^2|M^{\e}|^{-1}dv\int_{\R^3}\la v \ra^{n+2}|\p^{\al_x-\beta_x}\mathcal{L}^{-1}H(v)|^2|M^{\e}|^{-1}dvdx .
\enda
\end{align}
When \(|\beta_x| = \mathrm{N}\), we cannot take the \(L^\infty_x\) norm of \(\p^{\beta_x}M^{\e}\) because we lack an upper bound for \(\|\nabla_x^{\mathrm{N}}(\rho,u,\ta)\|_{L^\infty_x}\). For \(\p^{\beta_x}M^{\e}\), we take \(L^\infty_x\) for \(0 \leq |\beta_x| \leq \mathrm{N}-1\) and \(L^2_x\) for \(|\beta_x| = \mathrm{N}\), applying \eqref{Phiinfscale} and \eqref{Phiscale}, respectively:
\begin{align}\label{L-13a}
\bega
B_2^{\al_x} &\les \begin{cases} 
\e^2\kappa^{-1}\mathcal{E}_{tot} \int_{\Omega\times\R^3}\la v \ra^{n+2}|\p^{\al_x-\beta_x}\mathcal{L}^{-1}H(v)|^2|M^{\e}|^{-1}dvdx, & \text{if } 0 \leq |\beta_x| \leq \mathrm{N}-1, \\ 
\e^2\mathcal{E}_M\Big\|\int_{\R^3}\la v \ra^{n+2}|\mathcal{L}^{-1}H(v)|^2|M^{\e}|^{-1}dv\Big\|_{L^\infty_x}, & \text{if } |\beta_x| = \mathrm{N}.
\end{cases}
\enda
\end{align}
Combining the estimate of $B_1^{\al_x}$ and $B_2^{\al_x}$ gives 
\begin{align}\label{L-13}
\bega
\int_{\Omega\times\R^3} \la v \ra^{n+2}&|\p^{\al_x} \mathcal{L}^{-1}H|^2 |M^{\e}|^{-1} dvdx \les \int_{\Omega\times\R^3} \la v \ra^n|\p^{\al_x}H|^2 |M^{\e}|^{-1} dvdx \cr
&+ \sum_{0<|\beta_x|\leq\mathrm{N}-1}\e^2\kappa^{-1}\mathcal{E}_{tot} \int_{\Omega\times\R^3}\la v \ra^{n+2}|\p^{\al_x-\beta_x}\mathcal{L}^{-1}H(v)|^2|M^{\e}|^{-1}dvdx \cr 
&+ 1_{|\al_x|=\mathrm{N}}\e^2\mathcal{E}_M\Big\|\int_{\R^3}\la v \ra^{n+2}|\mathcal{L}^{-1}H(v)|^2|M^{\e}|^{-1}dv\Big\|_{L^\infty_x}.
\enda
\end{align}
The derivatives in the second line of \eqref{L-13} are fewer than \(|\al_x|\) due to \(0 < |\beta_x|\). Applying \eqref{L-13} recursively to \(\p^{\al_x-\beta_x}\mathcal{L}^{-1}\), and using \(\e^2\kappa^{-1}\mathcal{E}_{top} \leq C\) from \eqref{condition}, we establish \eqref{L-1}.
\end{proof}


\begin{lemma}\label{L.Xi2} 
Assume \eqref{condition}. Then, we have
\begin{align*}
\bega
\int_{\Omega\times\R^3} \la v \ra \Big|\p^{\al_x}\Big(\mathcal{L}^{-1}(\hat{X}M^{\e}) |M^{\e}|^{-1}\Big)\Big|^2 M^{\e} dv &\leq \e^2\mathcal{E}_M, \quad \text{for} \quad 1 \leq |\al_x| \leq \mathrm{N},
\enda
\end{align*}
for any \(n \geq 0\), where \(\hat{X} = \mathfrak{R}^{\e}_{ij}\) or \(\mathcal{Q}^{\e}_j\).
\end{lemma}

\begin{proof}
The proof for \(\hat{X} = \mathfrak{R}^{\e}_{ij}\) and \(\hat{X} = \mathcal{Q}^{\e}_j\) is identical. We will consider the case where \(\hat{X} = \mathfrak{R}^{\e}_{ij}\). We distribute the derivative \(\p^{\al_x}\) as follows:
\begin{align*}
\int_{\Omega\times\R^3} \la v \ra \Big|\p^{\al_x}\Big(\mathcal{L}^{-1}(\mathfrak{R}^{\e}_{ij}M^{\e}) &|M^{\e}|^{-1}\Big)\Big|^2 M^{\e} dv = I_1^{\al_x} + I_2^{\al_x},
\end{align*}
where
\begin{align*}
I_1^{\al_x} &= \int_{\Omega\times\R^3} \la v \ra \Big|\p^{\al_x}\Big(\mathcal{L}^{-1}(\mathfrak{R}^{\e}_{ij}M^{\e})\Big)\Big|^2 |M^{\e}|^{-1} dv, \cr
I_2^{\al_x} &= \sum_{0<\beta_x\leq\al_x}\binom{\al_x}{\beta_x} \int_{\Omega\times\R^3} \la v \ra \Big|\p^{\al_x-\beta_x}\Big(\mathcal{L}^{-1}(\mathfrak{R}^{\e}_{ij}M^{\e})\Big) \p^{\beta_x}|M^{\e}|^{-1}\Big|^2 M^{\e} dv.
\end{align*}
For \(I_1^{\al_x}\), applying Lemma \ref{L-1G2} for \(H = \mathfrak{R}^{\e}_{ij}M^{\e}\), we have 
\begin{align*}
\bega
I_1^{\al_x} &\les \int_{\Omega\times\R^3} \la v \ra|\p^{\al_x}(\mathfrak{R}^{\e}_{ij}M^{\e})|^2 |M^{\e}|^{-1} dvdx \cr
&+ \sum_{0<|\beta_x|\leq\mathrm{N}-1}\Big(\e^2\kappa^{-1}\mathcal{E}_{tot}\Big) \int_{\Omega\times\R^3}\la v \ra^{n-2}|\p^{\al_x-\beta_x}(\mathfrak{R}^{\e}_{ij}M^{\e})|^2|M^{\e}|^{-1}dvdx \cr 
&+ 1_{|\al_x|=\mathrm{N}}\Big(\e^2\mathcal{E}_M\Big)\Big\|\int_{\R^3}\la v \ra^{n}|\mathfrak{R}^{\e}_{ij}M^{\e}|^2|M^{\e}|^{-1}dv\Big\|_{L^\infty_x}.
\enda
\end{align*}
Applying \eqref{al-vM} and \eqref{Phiscale}, we have
\begin{align*}
\int_{\Omega\times\R^3} \la v \ra|\p^{\al_x}(\mathfrak{R}^{\e}_{ij}M^{\e})|^2 |M^{\e}|^{-1} dvdx \les \Big(\e\mathcal{E}_M^{\frac{1}{2}}\Big)^2 , \quad \text{for} \quad 1\leq|\al_x|\leq \mathrm{N}.
\end{align*}
Thus, we obtain
\begin{align}\label{L-1H1r}
\bega
I_1^{\al_x} &\les \e^2\mathcal{E}_M + \sum_{0<|\beta_x|\leq\mathrm{N}-1}\Big(\e^2\kappa^{-1}\mathcal{E}_{top}\Big) \Big(\e^2\mathcal{E}_M\Big) + 1_{|\al_x|=\mathrm{N}}\Big(\e^2\mathcal{E}_M\Big) \les \e^2\mathcal{E}_M,
\enda
\end{align}
where we used \(\e^2\kappa^{-1}\mathcal{E}_{top} \leq C\) by \eqref{condition}. 
For \(I_2^{\al_x}\), using \eqref{al-1/M}, we get
\begin{align}\label{L-1H2-a}
\bega
I_2^{\al_x} &\les \sum_{1\leq i\leq |\beta_x|}\eps^{2i} \int_{\R^3} \la v \ra^{n+4|\beta_x|} \Big|\p^{\al_x-\beta_x}\Big(\mathcal{L}^{-1}(\mathfrak{R}^{\e}_{ij}M^{\e})\Big)\Big|^2 |M^{\e}|^{-1} dv \cr 
& \times  \sum_{\substack{\beta_{x_1}+\cdots+\beta_{x_i}=\beta_x \\ \beta_{x_i}>0}}|\p^{\beta_{x_1}}(\rho,u,\ta)|^2\times \cdots \times |\p^{\beta_{x_i}}(\rho,u,\ta)|^2.
\enda
\end{align}
We integrate over \(\Omega\) and then take the \(L^\infty_x\) norm for the \((\rho,u,\ta)\) terms. The first line of \eqref{L-1H2-a} can be estimated similarly to \eqref{L-1H1r}. The second line of \eqref{L-1H2-a} can be estimated similarly to the proof of \eqref{Phi2}:
\begin{align}\label{L-1H2r}
\bega
I_2^{\al_x} 
&\les \Big(\e^2\mathcal{E}_M\Big) \Big(\e^2\mathcal{E}_M\Big) \les \e^4\mathcal{E}_M^2.
\enda
\end{align}
Combining \eqref{L-1H1r} and \eqref{L-1H2r}, and using \(\e\mathcal{E}_M^{\frac{1}{2}}\leq C\) by \eqref{condition}, we obtain the desired result.
\end{proof}

\begin{proof}[\textbf{Proof of Proposition \ref{P.Xi}}]
We apply Lemma \ref{L.Xi2} to the last line of \eqref{XiL-1} in Lemma \ref{L.Xi1} and use $\e\mathcal{E}_M^{\frac{1}{2}}\leq C$ by \eqref{condition}, obtaining
\begin{align}\label{Xi-proof-1}
\bega
\|\p^{\al_x}\varXi_{ij}^A\|_{L^2_x}&,~ \|\p^{\al_x}\varXi_{j}^B\|_{L^2_x} \les  \e^3\kappa\varpi(\e,\kappa)\Big(\mathcal{D}_{tot}^{\frac{1}{2}}(t)+\mathbf{1}_{|\al|\geq2}\mathcal{E}_M(t)\Big) \cr 
&+ \sum_{0\leq|\beta_x|\leq\lfloor|\al_x|/2\rfloor} \big\|\la v \ra^{\frac{1}{2}} \p^{\beta_x}\AC{\P}F^{\e}|M^{\e}|^{-1/2} \big\|_{L^\infty_xL^2_v}  \big\|\la v \ra^{\frac{1}{2}} \p^{\al_x-\beta_x}\AC{\P}F^{\e}|M^{\e}|^{-1/2} \big\|_{L^2_{x,v}}.
\enda
\end{align}
Next, we take the \(L^2_t\)-norm over \(t \in [0,T]\). Then, the last line of \eqref{Xi-proof-1} can be estimed as 
\begin{align*}
\bega
\sum_{0\leq|\beta_x|\leq\lfloor\mathrm{N}/2\rfloor} \Big\|\|\la v \ra^{\frac{1}{2}} \p^{\beta_x}\AC{\P}F^{\e}|M^{\e}|^{-1/2} \|_{L^\infty_xL^2_v} \Big\|_{L^2_t} \Big\| \|\la v \ra^{\frac{1}{2}} \p^{\al_x-\beta_x}\AC{\P}F^{\e}|M^{\e}|^{-1/2} \|_{L^2_{x,v}} \Big\|_{L^\infty_T}
\enda
\end{align*}
where we took $L^\infty_T$ to the last norm. Then applying \eqref{Gembed} and \eqref{Gembedt} for $L^2_t$ norm and $L^\infty_T$ norm, respectively, we obtain
\begin{align}\label{Xi-proof-2}
\bega
&\|\p^{\al_x}\varXi_{ij}^A\|_{L^2_tL^2_x},~ \|\p^{\al_x}\varXi_{j}^B\|_{L^2_tL^2_x} \les  \e^3\kappa\varpi(\e,\kappa)\bigg(\Big(\int_0^t\mathcal{D}_{tot}(s)ds\Big)^{\frac{1}{2}}+\mathbf{1}_{|\al|\geq2}\Big(\int_0^t\mathcal{E}_M^2(s)ds\Big)^{\frac{1}{2}}\bigg) \cr 
&+ \e^3\varpi(\e,\kappa)\bigg(\frac{1}{\sqrt{T}}+\sqrt{T}\bigg)\bigg(\int_0^t\mathcal{D}_{tot}(s)ds + \int_0^t\mathcal{E}_M^2(s)ds +\sum_{i=1,2}\kappa^{-\frac{1}{2}} \int_0^t \mathcal{E}_M^i(s)\mathcal{V}_{4i+1}(s)ds\bigg).
\enda
\end{align}
For the term $\int_0^t \mathcal{E}_M^i(s)\mathcal{V}_{4i+1}(s)ds$, we use $\e\mathcal{E}_M^{\frac{1}{2}}\leq C$ by \eqref{condition} and apply \eqref{L.Vdecomp} and \eqref{hL2} to get
\begin{align}\label{Xi-proof-3}
\bega
\int_0^t &\mathcal{E}_M(s)\mathcal{V}_{4i+1}(s)ds 
\leq C \int_0^t \bigg[\mathcal{E}_{tot}^2(t) +\e^{\frac{3}{2}}\kappa \mathcal{E}_{tot}\Big(\mathcal{D}_{tot}(t)+\mathbf{1}_{|\al|\geq2}\mathcal{E}_M^2(t)\Big)\bigg]ds \cr 
&+\sup_{0\leq s\leq t}\mathcal{E}_{tot}(s) e^{-c_1\e^{-{\frac{1}{\ell}}}}\bigg[\e^2\kappa\sum_{0\leq|\al|\leq\mathrm{N}+1}\kappa^{(|\al|-\mathrm{N})_+}\|\p^{\al}h_0\|_{L^\infty_{x,v}}^2 + \frac{1}{(\e\kappa)^{d}}\int_0^t \mathcal{E}_{tot}(s) ds\bigg].
\enda
\end{align}
Combining \eqref{Xi-proof-2} and \eqref{Xi-proof-3}, we obtain the desired result. 
\end{proof}

We estimate the forcing terms derived from Proposition \ref{P.max}, Proposition \ref{P.macro.u}, and Proposition \ref{P.macro.rta}.

\begin{proposition}\label{P.forcing}
\textcolor{red}{We need to estimate $\|g^{\e}_{\w}\|_{L^p_x}$ for $1<p<2$ by using
\begin{align}
\|AB\|_{L^p_x} &\leq \|A\|_{L^{\frac{2p}{2-p}}_x} \|B\|_{L^2_x}.
\end{align}
to the nonlinear terms
}
Assume \eqref{condition} holds. The forcing terms \(\varPi_{\w}\) and \(\varPi_{\mathfrak{s}}\), as defined in \eqref{gwrta}, satisfy
\begin{align*}
\bega
\sum_{0\leq|\al_x|\leq \mathrm{N}-1}&\int_{\Omega}\p^{\al_x}\w \p^{\al_x}\varPi_{\w} dx \leq C \bigg[ \Big(\|\nabla_x \cdot u(t)\|_{H^{\mathrm{N}-1}_x} + \|\nabla_x(\rho+\ta)(t)\|_{H^{\mathrm{N}-1}_x}\Big) \mathcal{E}_M(t)  \cr
&+ \e\mathcal{E}_{tot}^{\frac{3}{2}}(t) +\e\mathcal{E}_{tot}^2(t) +\e\kappa^{\frac{1}{2}}\mathcal{D}_G^{\frac{1}{2}}(t)\mathcal{E}_M(t) + \frac{1}{\e^2\kappa^{\frac{1}{2}}}\|\nabla_x\varXi_{ij}^A(t)\|_{H^{\mathrm{N}-1}_x}\mathcal{E}_{top}^{\frac{1}{2}}(t) \bigg],
\enda
\end{align*}
\begin{align*}
\bega
\sum_{0\leq|\al_x|\leq \mathrm{N}}&\int_{\Omega}\p^{\al_x}\Big(\rho-\frac{3}{2}\ta\Big) \p^{\al_x}(\varPi_{\mathfrak{s}}) dx \leq C \bigg[  \kappa^{\frac{1}{2}} \|\nabla_x (\rho+\ta)(t)\|_{H^\mathrm{N}_x}  \mathcal{E}_{top}^{\frac{1}{2}}(t)   \cr
&+ \e\mathcal{E}_{tot}^{\frac{3}{2}}(t) +\e\mathcal{E}_{tot}^2(t) +\e\kappa^{\frac{1}{2}}\mathcal{D}_G^{\frac{1}{2}}(t)\mathcal{E}_M(t) + \frac{1}{\e^2\kappa^{\frac{1}{2}}}\|\nabla_x\varXi_{j}^B(t)\|_{H^{\mathrm{N}-1}_x}\mathcal{E}_{top}^{\frac{1}{2}}(t)\bigg] ,
\enda
\end{align*}
and
\begin{align*}
\bega
&\Big\|(\varPi_{\w}+(\nabla_x \cdot u)\w)(t)\Big\|_{L^\infty_x}  \leq C  \bigg[ \|\nabla_x(\rho+\ta)(t)\|_{L^\infty_x} \mathcal{E}_M^{\frac{1}{2}}(t)  +\e\kappa^{\frac{1}{2}}\mathcal{E}_{tot}(t) \cr 
&\hspace{4cm}+ \e\kappa^{\frac{1}{2}}\mathcal{D}_G^{\frac{1}{2}}(t)\mathcal{E}_M^{\frac{1}{2}}(t) + \frac{1}{\e^2}\|\nabla_x\varXi_{ij}^A(t)\|_{L^\infty_x}\bigg], \cr 
&\Big\|\nabla_x\varPi_{\mathfrak{s}}(t)\Big\|_{L^\infty_x} \leq  C \bigg[ \kappa\|\Delta_x (\rho+\ta)(t)\|_{L^\infty_x} +\e\kappa^{\frac{1}{2}}\mathcal{E}_{tot}(t) + \e\kappa^{\frac{1}{2}}\mathcal{D}_G^{\frac{1}{2}}(t)\mathcal{E}_M^{\frac{1}{2}}(t) +  \frac{1}{\e^2}\|\nabla_x\varXi_{j}^B(t)\|_{L^\infty_x} \bigg]
\enda
\end{align*}
for a positive constant $C>0$.
\end{proposition}
\begin{proof}
By the definition of $\mathcal{E}_M$ and $\mathcal{E}_{top}$ in \eqref{N-EDdef} and \eqref{N-EDdef2}, we have  
\begin{align}\label{wleq}
\bega
\|(\rho,u,\ta)\|_{H^{\mathrm{N}}_x} + \|\w\|_{H^{\mathrm{N}-1}_x} &\leq \mathcal{E}_M^{\frac{1}{2}} , \cr 
\|(\rho,u,\ta)\|_{\dot{H}^{\mathrm{N}+1}_x} + \|\w\|_{\dot{H}^{\mathrm{N}}_x} &\leq  \kappa^{-\frac{1}{2}}\mathcal{E}_{top}^{\frac{1}{2}} ,
\enda
\end{align}
where we used $\|\w\|_{H^{s}_x} \leq \|\nabla_x u\|_{H^{s}_x}$ by \eqref{pmw}. \\
{\bf (Estimate of $\varPi_{\w}$)} Recall the definition of \(\varPi_{\w}\) defined in \eqref{gwrta}.
We estimate the $\varPi_{\w}$ term component by component. For the first term of $\varPi_{\w}$, we apply integration by parts to $\Delta_x \w$,
%
using \eqref{uvHk}, to obtain
\begin{align}\label{Vw--1}
\bega
\kappa&\sum_{0\leq|\al_x|\leq \mathrm{N}-1}\bigg|\int_{\Omega}\p^{\al_x}\Big((\mathrm{\Theta}^{\e}-1) \Delta_x \w\Big) \p^{\al_x}\w dx\bigg| \cr 
&\les \kappa\|(\mathrm{\Theta}^{\e}-1)\|_{L^\infty_x}\|\nabla_x \w\|_{H^{\mathrm{N}-1}_x}^2+ \kappa\|\nabla_x(\mathrm{\Theta}^{\e}-1)\|_{H^{\mathrm{N}-1}_x}\|\nabla_x \w\|_{H^{\mathrm{N}-1}_x}\|\w\|_{H^{\mathrm{N}-1}_x} \cr 
&\les \kappa (\e\mathcal{E}_M^{\frac{1}{2}})\kappa^{-1}\mathcal{E}_{top} + \kappa (\e\mathcal{E}_M^{\frac{1}{2}})\kappa^{-\frac{1}{2}}\mathcal{E}_{top}^{\frac{1}{2}}\mathcal{E}_M^{\frac{1}{2}} \les \e\mathcal{E}_{tot}^{\frac{3}{2}},
\enda
\end{align}
where we used \eqref{pTaleq} and \eqref{wleq}. 
For the second term of $\varPi_{\w}$, using \eqref{uvHk}, we obtain
\begin{align}\label{Vw-0}
\bega
\|(\nabla_x \cdot u)\w\|_{H^{\mathrm{N}-1}_x} \|\w\|_{H^{\mathrm{N}-1}_x} &\les \|\nabla_x \cdot u\|_{L^\infty_x}\|\w\|_{H^{\mathrm{N}-1}_x}^2 + \|\nabla_x \cdot u\|_{H^{\mathrm{N}-1}_x}\|\w\|_{L^\infty_x}\|\w\|_{H^{\mathrm{N}-1}_x} \cr 
&\les \|\nabla_x \cdot u\|_{H^{\mathrm{N}-1}_x} \mathcal{E}_M,
\enda
\end{align}
where we used $\|\w\|_{L^\infty_x} \leq \|\w\|_{H^2_x}$ by \eqref{Agmon}. 
For the second line of $\varPi_{\w}$, since the Sobolev space $H^2_x$ is algebra, for $\mathrm{N}\geq3$, we find
\begin{align}\label{Vw-1}
\bega
&\bigg\|k_B\mathrm{\Theta}^{\e}\nabla_x^{\perp}\ta \cdot \bigg[\nabla_x(\rho+\ta)
-\frac{1}{\e}\frac{1}{k_B\mathrm{P}^{\e} \mathrm{\Theta}^{\e}}\sum_{j} \p_{x_j} \mathbf{r}_{ij}^{\e}\bigg] \bigg\|_{H^{\mathrm{N}-1}_x}\|\w\|_{H^{\mathrm{N}-1}_x} \cr  
&\les \|k_B\mathrm{\Theta}^{\e}\|_{H^{\mathrm{N}-1}_x} \|\nabla_x^{\perp}\ta \|_{H^{\mathrm{N}-1}_x} \bigg[\|\nabla_x(\rho+\ta)\|_{H^{\mathrm{N}-1}_x}
+\frac{1}{\e}\bigg\|\frac{1}{k_B\mathrm{P}^{\e} \mathrm{\Theta}^{\e}}\bigg\|_{H^{\mathrm{N}-1}_x}\sum_{j} \|\p_{x_j} \mathbf{r}_{ij}^{\e}\|_{H^{\mathrm{N}-1}_x}\bigg]\|\w\|_{H^{\mathrm{N}-1}_x}  \cr 
&\les  (1+\e\mathcal{E}_M^{\frac{1}{2}})\mathcal{E}_M^{\frac{1}{2}}\bigg[\|\nabla_x(\rho+\ta)\|_{H^{\mathrm{N}-1}_x} + \frac{1}{\e}\e^2\kappa^{\frac{1}{2}}\mathcal{D}_G^{\frac{1}{2}}\bigg]\mathcal{E}_M^{\frac{1}{2}} \les \|\nabla_x(\rho+\ta)\|_{H^{\mathrm{N}-1}_x}\mathcal{E}_M + \e\kappa^{\frac{1}{2}}\mathcal{D}_G^{\frac{1}{2}}\mathcal{E}_M ,
\enda
\end{align}
where we used $\|\mathrm{\Theta}^{\e}\|_{H^{\mathrm{N}-1}_x} \les 1+\e\mathcal{E}_M^{\frac{1}{2}}$, \eqref{ABGscale} and $\|\frac{1}{k\mathrm{P}^{\e} \mathrm{\Theta}^{\e}}\|_{H^{\mathrm{N}-1}_x} \les 1+\e\mathcal{E}_M^{\frac{1}{2}}$ by \eqref{pTaleq}  and \eqref{condition}.
Next, we estimate the $k_B\mathrm{\Theta}^{\e} \nabla_x^{\perp}\cdot g_{u^{\e}}$ term. Recall the definition of $g_{u^{\e}_i}$ in \eqref{gutadef}:
\begin{align*}
\bega
g_{u^{\e}_i}&:= \eta_0 k_B^{-\frac{1}{2}} \kappa \Big(\frac{1}{\mathrm{P}^{\e}|\mathrm{\Theta}^{\e}|^{\frac{1}{2}}}-1\Big)  \lw(\Delta_x u_i + \frac{1}{3}\p_i\nabla_x\cdot u\rw) \cr 
&\quad + \frac{1}{2}\eta_0 k_B^{-\frac{1}{2}} \kappa \e \frac{1}{\mathrm{P}^{\e}|\mathrm{\Theta}^{\e}|^{\frac{1}{2}}} \sum_j \p_j\ta\lw(\p_i u_j+\p_j u_i-\frac 2 3 \delta_{ij}(\nabla_x\cdot u)\rw)  +\frac{1}{\e^2}\frac{1}{k_B\mathrm{P}^{\e}\mathrm{\Theta}^{\e}}\sum_j \p_j\varXi_{ij}^A.
\enda
\end{align*}
For the first line of $\nabla_x^{\perp}\cdot g_{u^{\e}}$, we use integration by parts. By the same way to \eqref{Vw--1}, we obtain
\begin{align}\label{Vw-2}
\bega
\kappa\sum_{0\leq|\al_x|\leq \mathrm{N}-1} \bigg|\int_{\Omega}\nabla_x^\perp \cdot\p^{\al_x}\bigg[\Big(\frac{1}{\mathrm{P}^{\e}|\mathrm{\Theta}^{\e}|^{\frac{1}{2}}}-1\Big)  \lw(\Delta_x u_i + \frac{1}{3}\p_i\nabla_x\cdot u\rw)\bigg] \p^{\al_x}\w dx\bigg| \les \e\mathcal{E}_{tot}^{\frac{3}{2}}.
\enda
\end{align}
For the second term of $k_B\mathrm{\Theta}^{\e} \nabla_x^{\perp}\cdot g_{u^{\e}}$, by the same way to \eqref{Vw-1}, we can easily see that 
\begin{align}\label{Vw-2.5}
\bega
\kappa \e &\bigg\| k_B\mathrm{\Theta}^{\e} \nabla_x^{\perp}\cdot \bigg[\frac{1}{\mathrm{P}^{\e}|\mathrm{\Theta}^{\e}|^{\frac{1}{2}}} \sum_j \p_j\ta\lw(\p_i u_j+\p_j u_i-\frac 2 3 \delta_{ij}(\nabla_x\cdot u)\rw)\bigg] \bigg\|_{H^{\mathrm{N}-1}_x}\|\w\|_{H^{\mathrm{N}-1}_x} \cr 
&\les \kappa \e \bigg[\|k_B\mathrm{\Theta}^{\e}\|_{H^{\mathrm{N}-1}_x}\bigg\|\frac{1}{\mathrm{P}^{\e}|\mathrm{\Theta}^{\e}|^{\frac{1}{2}}}\bigg\|_{H^{\mathrm{N}}_x} \sum_j \|\p_j\ta\|_{H^{\mathrm{N}}_x}\Big\|\p_i u_j+\p_j u_i-\frac 2 3 \delta_{ij}(\nabla_x\cdot u)\Big\|_{H^{\mathrm{N}}_x}\bigg]\|\w\|_{H^{\mathrm{N}-1}_x} \cr 
&\les \kappa \e (1+\e\mathcal{E}_M^{\frac{1}{2}})\mathcal{E}_M^{\frac{1}{2}}(\kappa^{-\frac{1}{2}}\mathcal{E}_{top}^{\frac{1}{2}})(\kappa^{-\frac{1}{2}}\mathcal{E}_{top}^{\frac{1}{2}})\mathcal{E}_M^{\frac{1}{2}} \les \e\mathcal{E}_{tot}^2.
\enda
\end{align}
For the last term of $k_B\mathrm{\Theta}^{\e} \nabla_x^{\perp}\cdot g_{u^{\e}}$, we apply integration by parts:
\begin{align}\label{Vw-3}
\bega
\int_{\Omega}\p^{\al_x}&\nabla_x^{\perp}\cdot\bigg(\frac{1}{\e^2}\frac{1}{k_B\mathrm{P}^{\e}\mathrm{\Theta}^{\e}}\sum_j \p_j\varXi_{ij}^A\bigg) \p^{\al_x}\w dx \leq \frac{1}{\e^2}\bigg\|\frac{1}{k_B\mathrm{P}^{\e}\mathrm{\Theta}^{\e}}\sum_j \p_j\varXi_{ij}^A\bigg\|_{H^{\mathrm{N}-1}_x} \|\nabla_x\w\|_{H^{\mathrm{N}-1}_x} \cr 
&\les \frac{1}{\e^2}(1+\e\mathcal{E}_M^{\frac{1}{2}}) \|\nabla_x\varXi_{ij}^A\|_{H^{\mathrm{N}-1}_x}\|\nabla_x\w\|_{H^{\mathrm{N}-1}_x} \leq \frac{1}{\e^2}\|\nabla_x\varXi_{ij}^A\|_{H^{\mathrm{N}-1}_x} \Big(\kappa^{-\frac{1}{2}}\mathcal{E}_{top}^{\frac{1}{2}}\Big).
\enda
\end{align}
Combining \eqref{Vw--1}, \eqref{Vw-0}, \eqref{Vw-1}, \eqref{Vw-2}, \eqref{Vw-2.5}, and \eqref{Vw-3} gives 
\begin{align*}
\sum_{0\leq|\al_x|\leq \mathrm{N}-1}\int_{\Omega}\p^{\al_x}\w \p^{\al_x}\varPi_{\w} dx &\les \Big(\|\nabla_x \cdot u\|_{H^{\mathrm{N}-1}_x} + \|\nabla_x(\rho+\ta)\|_{H^{\mathrm{N}-1}_x}\Big) \mathcal{E}_M  + \e\mathcal{E}_{tot}^{\frac{3}{2}} +\e\mathcal{E}_{tot}^2 \cr 
&+ \e\kappa^{\frac{1}{2}}\mathcal{D}_G^{\frac{1}{2}}\mathcal{E}_M + \frac{1}{\e^2}\|\nabla_x\varXi_{ij}^A\|_{H^{\mathrm{N}-1}_x} \Big(\kappa^{-\frac{1}{2}}\mathcal{E}_{top}^{\frac{1}{2}}\Big)
\end{align*}
This gives the desired result. 
By a similar approach, we can estimate $\|(\varPi_{\w}+(\nabla_x \cdot u)\w)\|_{L^\infty_x}$. 
\\
{\bf (Estimate of $\varPi_{\mathfrak{s}}$)} 
Recall the definition of \(\varPi_{\mathfrak{s}}\) defined in \eqref{gwrta}.
For the first terms of $\varPi_{\mathfrak{s}}$, we use integration by parts:
\begin{align}\label{Vrta-0}
\bega
\kappa \sum_{0\leq|\al_x|\leq \mathrm{N}}\int_{\Omega}\p^{\al_x} \Delta_x (\rho+\ta) \p^{\al_x}\Big(\rho-\frac{3}{2}\ta\Big) dx &\les \kappa \sum_{0\leq|\al_x|\leq \mathrm{N}}\bigg|\int_{\Omega}\p^{\al_x} \nabla_x (\rho+\ta) \p^{\al_x}\nabla_x\Big(\rho-\frac{3}{2}\ta\Big) dx\bigg| \cr 
&\les \kappa \|\nabla_x (\rho+\ta)\|_{H^\mathrm{N}_x}  (\kappa^{-\frac{1}{2}}\mathcal{E}_{top}^{\frac{1}{2}}).
\enda
\end{align}
Recall the definition of $g_{\ta}$ is \eqref{gutadef}:
\begin{align*}
\bega
g_{\ta} &:= \frac{5}{2}\eta_1 k_B^{\frac 1 2} \kappa \Big(\frac{|\mathrm{\Theta}^{\e}|^{\frac 1 2}}{\mathrm{P}^{\e}}-1\Big) \Delta_x\ta + \frac{15}{4}\eta_1 k_B^{\frac 1 2} \kappa \e \frac{|\mathrm{\Theta}^{\e}|^{\frac 1 2}}{\mathrm{P}^{\e}}|\nabla_x\ta|^2 \cr 
&\quad+ \frac{1}{\e^2 k_B \mathrm{P}^{\e}\mathrm{\Theta}^{\e}} \bigg(\sum_j \p_j \varXi_{j}^B-\sum_{i,j} \p_{x_i}\mathrm{U}^{\e}_j \mathbf{r}_{ij}^{\e}\bigg) .
\enda
\end{align*}
For the first and second terms of $g_{\ta}$, we use integration by parts. By the same way to \eqref{Vw--1}, we obtain
\begin{align}\label{Vrta-1}
\bega
&\kappa \sum_{0\leq|\al_x|\leq \mathrm{N}}\bigg|\int_{\Omega} \p^{\al_x}\bigg[\Big(\frac{|\mathrm{\Theta}^{\e}|^{\frac 1 2}}{\mathrm{P}^{\e}}-1\Big)  \Delta_x\ta\bigg] \p^{\al_x}\Big(\rho-\frac{3}{2}\ta\Big) dx\bigg| 
\les \e\mathcal{E}_{tot}^{\frac{3}{2}} , \cr 
&\kappa \e \sum_{0\leq|\al_x|\leq \mathrm{N}}\bigg|\int_{\Omega} \p^{\al_x}\bigg[ \frac{|\mathrm{\Theta}^{\e}|^{\frac 1 2}}{\mathrm{P}^{\e}} |\nabla_x\ta|^2\bigg] \p^{\al_x}\Big(\rho-\frac{3}{2}\ta\Big) dx\bigg| 
\les \e\mathcal{E}_{tot}^2.
\enda
\end{align}
For the last line of $g_{\ta}$, by the same way to \eqref{Vw-3}, we can have the following estimate for the term $\p_j \varXi_{j}^B$:
\begin{align}\label{Vrta-4}
\bega
\int_{\Omega}\p^{\al_x} \bigg( \frac{1}{\e^2 k_B \mathrm{P}^{\e}\mathrm{\Theta}^{\e}}\p_j \varXi_{j}^B \bigg) \p^{\al_x}\Big(\rho-\frac{3}{2}\ta\Big) dx 
\les \frac{1}{\e^2}\|\nabla_x\varXi_{j}^B\|_{H^{\mathrm{N}-1}_x}\Big(\kappa^{-\frac{1}{2}}\mathcal{E}_{top}^{\frac{1}{2}}\Big).
\enda
\end{align}
For the term $\p_{x_i}\mathrm{U}^{\e}_j \mathbf{r}_{ij}^{\e}$ applying \eqref{ABGscale} and \eqref{uvHk} gives
\begin{align}\label{Vrta-3.5}
\bega
\Big\|\frac{1}{\e^2 k_B \mathrm{P}^{\e}\mathrm{\Theta}^{\e}}\sum_{i,j} \p_{x_i}\mathrm{U}^{\e}_j \mathbf{r}_{ij}^{\e}\Big\|_{H^{\mathrm{N}}_x} 
&\les \frac{1}{\e}(1+\e\mathcal{E}_M^{\frac{1}{2}})\Big(\|\nabla_xu\|_{L^\infty_x}\|\mathbf{r}_{ij}^{\e}\|_{H^{\mathrm{N}}_x} + \|\nabla_xu\|_{H^{\mathrm{N}}_x}\|\mathbf{r}_{ij}^{\e}\|_{L^\infty_x}\Big) \cr 
&\les \frac{1}{\e}(1+\e\mathcal{E}_M^{\frac{1}{2}}) \mathcal{E}_M^{\frac{1}{2}}(\e^2\kappa^{\frac{1}{2}}\mathcal{D}_G^{\frac{1}{2}}) \les \e\kappa^{\frac{1}{2}}\mathcal{D}_G^{\frac{1}{2}}\mathcal{E}_M^{\frac{1}{2}}.
\enda
\end{align}
Combining \eqref{Vrta-1}, \eqref{Vrta-0}, \eqref{Vrta-4}, and \eqref{Vrta-3.5}, we have 
\begin{align*}
\sum_{0\leq|\al_x|\leq \mathrm{N}}\int_{\Omega}\p^{\al_x}\Big(\rho-\frac{3}{2}\ta\Big)\p^{\al_x}(\varPi_{\mathfrak{s}}) dx &\les \e\mathcal{E}_{tot}^{\frac{3}{2}} +\e\mathcal{E}_{tot}^2 + \kappa^{\frac{1}{2}} \|\nabla_x (\rho+\ta)\|_{H^\mathrm{N}_x}  \mathcal{E}_{top}^{\frac{1}{2}} \cr 
&+ \e\kappa^{\frac{1}{2}}\mathcal{D}_G^{\frac{1}{2}}\mathcal{E}_M + \frac{1}{\e^2}\|\nabla_x\varXi_{j}^B\|_{H^{\mathrm{N}-1}_x} \Big(\kappa^{-\frac{1}{2}}\mathcal{E}_{top}^{\frac{1}{2}}\Big).
\end{align*}
Thus, we have the result. Similarly, we can estimate $\|\nabla_x\varPi_{\mathfrak{s}}(t)\|_{L^\infty_x}$.
\end{proof}

\unhide

\appendix

\StartNoTOC

\section{Basic properties of operators and functions}\label{A.Coll}


In this part, we present the fundamental properties of the projection operator $\P$, the collision operator $\mathcal{N}(\cdot,\cdot)$, and the linear operator $\mathcal{L}$, as defined in \eqref{Pdef}, \eqref{Qdef}, and \eqref{Ldef}, respectively. Additionally, we provide proofs for Lemma \ref{Llem1} and Lemma \ref{L.coer}.

\begin{lemma}\label{Pprop} 
For any given functions $F$, $G$, and $H$, such that the right-hand side is bounded, the operators $\P$ and $\AC{\P}$, defined in \eqref{Pdef}, satisfy the following properties:
\begin{align}\label{Pproj}
\int_{\R^3} (\P H) H |M^{\e}|^{-1} \, dv &= \int_{\R^3} |\P H|^2 |M^{\e}|^{-1} \, dv, \qquad 
\int_{\R^3} (\P G) (\AC{\P} H) |M^{\e}|^{-1} \, dv = 0,
\end{align}
and
\begin{align*}
\int_{\R^3} |H|^2 |M^{\e}|^{-1} \, dv 
= \int_{\R^3} |\P H|^2 |M^{\e}|^{-1} \, dv 
+ \int_{\R^3} |\AC{\P} H|^2 |M^{\e}|^{-1} \, dv,
\end{align*}
with the following inequalities:
\begin{align*}
\int_{\R^3}|\P H|^2 |M^{\e}|^{-1}dv \leq \int_{\R^3}|H|^2 |M^{\e}|^{-1}dv , \qquad \int_{\R^3}|\AC{\P} H|^2 |M^{\e}|^{-1}dv \leq \int_{\R^3}|H|^2 |M^{\e}|^{-1}dv.
\end{align*}
If we further assume that $
\sup_{t \in [0,T]} \big(|\mathrm{P}^{\e}(t,x) -1|, |\mathrm{U}^{\e}(t,x)|, |\mathrm{\Theta}^{\e}(t,x) -1| \big) < 1,$ then the following inequality holds:
\begin{align}\label{ACPX}
\bega
\int_{\R^3} \la v \ra^n |\AC{\P} H(t,x,v)|^2 |M^{\e}|^{-1} dv \les \int_{\R^3} \la v \ra^n |H(t,x,v)|^2 |M^{\e}|^{-1} dv,
\enda
\end{align}
for any positive constant $n \geq 0$.
\end{lemma}
\hide
\begin{proof}
By definition of $\AC{\P}$, we have 
\begin{align*}
\bega
\int_{\R^3} \la v \ra^n |\AC{\P}H|^2 |M^{\e}|^{-1} dv \les \int_{\R^3} \la v \ra^n |H|^2 M^{\e} dv + \int_{\R^3} \la v \ra^n |\PH|^2 |M^{\e}|^{-1} dv.
\enda
\end{align*}
Applying the H\"{o}lder inequality to the inner product of $\P$ in \eqref{Pdef}, we have 
\begin{align*}
\bega
|\PH| \les \sum_{i=0}^4 \lw(\int_{\R^3}|H|^2 |M^{\e}|^{-1} dv\rw)^{\frac{1}{2}} \lw(\int_{\R^3}\Big|\frac{e_i}{M^{\e}}\Big|^2 M^{\e} dv\rw)^{\frac{1}{2}} |e_i|.
\enda
\end{align*}
Since $e_i$ has exponential decay $M^{\e}$, we have the result. 
\end{proof}
\unhide


\hide
\begin{lemma}
For $Lf=-\frac{2}{\sqrt{\mu}}\mathcal{N}(\mu,\sqrt{\mu}f)$, we have 
\begin{align}\label{coerLmu}
\la Lf, f \ra_{L^2_v} \geq \sigma_L \|\sqrt{\nu}(\mathbf{I}-\P_{\mu})f\|_{L^2_v}
\end{align}
where 
\begin{align*}
\P_{\mu}:= \int_{\R^3}f\sqrt{\mu}dv\sqrt{\mu}+\sum_i\int_{\R^3}\frac{v_i}{\sqrt{k_B}}f\sqrt{\mu}dv\frac{v_i}{\sqrt{k_B}}\sqrt{\mu} + \int_{\R^3}\frac{|v|^2-3k_B}{\sqrt{6}k_B}f\sqrt{\mu}dv\frac{|v|^2-3k_B}{\sqrt{6}k_B}\sqrt{\mu}
\end{align*}
\end{lemma}

We can write 
\begin{align*}
M^{\e}(t,x,v)=\frac{\mathrm{P}^{\e}}{(2\pi k_B \mathrm{\Theta}^{\e})^{\frac{3}{2}}}e^{-\frac{|v-\mathrm{U}^{\e}|^2}{2k_B\mathrm{\Theta}^{\e}}} = \frac{\mathrm{P}^{\e}}{\sqrt{\mathrm{\Theta}^{\e}}^3}\frac{1}{(2\pi k_B )^{\frac{3}{2}}}e^{-\frac{|\tilde{v}|^2}{2k_B}} = \frac{\mathrm{P}^{\e}}{\sqrt{\mathrm{\Theta}^{\e}}^3}\mu(\tilde{v}), \quad \mbox{where} \quad \frac{v-\mathrm{U}^{\e}}{\sqrt{\mathrm{\Theta}^{\e}}}=\tilde{v}.
\end{align*}
If we write $Lf=-\frac{2}{\sqrt{\mu}}\mathcal{N}(\mu,\sqrt{\mu}f)$ then, similarity: 
\begin{align*}
\la \mathcal{L}(H), H|M^{\e}|^{-1} \ra_{L^2_v} \approx \bigg\la L\bigg(\frac{H}{\sqrt{M^{\e}}}\bigg), \frac{H}{\sqrt{M^{\e}}} \bigg\ra_{L^2_v}
\end{align*}
where $L(\cdot)$ is the linear operator perturbed from the global Maxwellian $L(\cdot)=-\frac{2}{\sqrt{\mu}}\mathcal{N}(\mu,\sqrt{\mu}(\cdot))$.
\unhide

\begin{proof}[\textbf{Proof of Lemma \ref{Llem1}}]
(1) By the definition of $\mathcal{L}$ in \eqref{Ldef}, and using the property $M^{\e}(v)M^{\e}(v_*)=M^{\e}(v')M^{\e}(v_*')$, we have
\begin{align}\label{FGsym}
\bega
\int_{\R^3}\mathcal{L}(F) &G|M^{\e}|^{-1} dv = \frac{1}{4}\int_{\R^3} \int_{\mathbb{R}^3}\int_{\mathbb{S}^2_+}|(v-v_*)\cdot w|M^{\e}(v)M^{\e}(v_*)\bigg[ \frac{F(v)}{M^{\e}(v)}+\frac{F(v_*)}{M^{\e}(v_*)} \cr
&-\frac{F(v')}{M^{\e}(v')}-\frac{F(v_*')}{M^{\e}(v_*')}\bigg] \times \bigg[ \frac{G(v)}{M^{\e}(v)}+\frac{G(v_*)}{M^{\e}(v_*)} -\frac{G(v')}{M^{\e}(v')}-\frac{G(v_*')}{M^{\e}(v_*')}\bigg] dwdv_* dv.
\enda
\end{align}
(2) From \eqref{FGsym}, it follows that $\mathcal{L}(F)=0$ if and only if $F \in \operatorname{span} \left\{ 1,v,|v|^2 \right\} M^{\e}$. \\
(3) A similar proof for the linear operator $L$ near the global Maxwellian can be found in Lemma 3.2 in \cite{Guo-NS}, page 638. From \eqref{FGsym}, the operator $\mathcal{L}$ satisfies $\la (1,v,|v|^2),\mathcal{L}F \ra = \int_{\R^3} \mathcal{L}((1,v,|v|^2)M^{\e}) F |M^{\e}|^{-1} dv = 0$ for any $F$. This implies that $\mathcal{L}: (\mathrm{Ker}(\mathcal{L}))^\perp \rightarrow (\mathrm{Ker}(\mathcal{L}))^\perp$. 
Next, consider the bilinear mapping $B(F,G):=\int_{\R^3}(\mathcal{L}F) G |M^{\e}|^{-1} dv$. Since $\int_{\R^3}F G |M^{\e}|^{-1} dv \leq C_F \|G|M^{\e}|^{-\frac{1}{2}}\|_{L^2_v}$, the term $\int_{\R^3} (\cdot) G |M^{\e}|^{-1} dv$ can be regarded as a bounded linear functional for a given $F$. Thus, by the Lax-Milgram theorem, the coercivity property \eqref{coer1} implies the existence of a unique function $F \in (\mathrm{Ker}(\mathcal{L}))^\perp$ such that $\la \mathcal{L}F, G|M^{\e}|^{-1} \ra_{L^2_v} = \la J,G|M^{\e}|^{-1} \ra_{L^2_v}$ for any $J \in (\mathrm{Ker}(\mathcal{L}))^\perp$. \\
(4) Recall the definition $\mathcal{L}(F) = -2\mathcal{N}(M^{\e},F)$. Applying \eqref{nonlin2}, we obtain
\begin{align*}
\int_{\R^3}\la v \ra^n |\mathcal{L}(F)|^2 |M^{\e}|^{-1} dv &\leq \int_{\R^3}(1+|v|^{n+2})|M^{\e}|^2|M^{\e}|^{-1} dv \int_{\R^3}(1+|v|^{n+2})|F|^2|M^{\e}|^{-1} dv.
\end{align*}
Since $|\mathrm{P}^{\e}(t,x) -1|, |\mathrm{U}^{\e}(t,x)|, |\mathrm{\Theta}^{\e}(t,x) -1| < 1$, we have $\int_{\R^3}(1+|v|^{n+2})M^{\e} dv \leq C$. This establishes the result $\eqref{L-1L2}_1$. By the open mapping theorem, we also deduce $\eqref{L-1L2}_2$.
\end{proof}

\begin{proof}[\textbf{Proof of Lemma \ref{L.coer}}]

We only prove \eqref{coer1} since the proofs of \eqref{nonlin} and \eqref{nonlin2} are similar. By the definition \(\mathcal{L}(G) = -2\mathcal{N}(M^{\e}, G)\) and using \eqref{FGsym}, we have
\begin{align}\label{LGG}
\bega
&\int_{\R^3}\mathcal{L}(G) G|M^{\e}|^{-1} dv \cr 
&= \frac{1}{4}\int_{\R^3} \int_{\mathbb{R}^3}\int_{\mathbb{S}^2_+}|(v-v_*)\cdot w|M^{\e}(v)M^{\e}(v_*)\bigg[ \frac{G(v)}{M^{\e}(v)}+\frac{G(v_*)}{M^{\e}(v_*)} -\frac{G(v')}{M^{\e}(v')}-\frac{G(v_*')}{M^{\e}(v_*')}\bigg]^2 dwdv_* dv \cr 
&= \frac{1}{4}\int_{\R^3} \int_{\mathbb{R}^3}\int_{\mathbb{S}^2_+}|(v-v_*)\cdot w|M^{\e}(v)M^{\e}(v_*)\bigg[ \frac{\AC{\P}G(v)}{M^{\e}(v)}+\frac{\AC{\P}G(v_*)}{M^{\e}(v_*)} -\frac{\AC{\P}G(v')}{M^{\e}(v')}-\frac{\AC{\P}G(v_*')}{M^{\e}(v_*')}\bigg]^2 dwdv_* dv,
\enda
\end{align}
where we used the definition of $\P$ in \eqref{Pdef}.
Because of the coercivity, there exists a constant \(\sigma_L > 0\) such that
\begin{align*}
\la Lf, f \ra_{L^2_v} \geq \sigma_L \|\sqrt{\nu}\AC{\P}_{\mu}f\|_{L^2_v}^2, \quad \mbox{where} \quad L(f) = -\frac{2}{\sqrt{\mu}}\mathcal{N}(\mu, \sqrt{\mu}f),
\end{align*}
and $\AC{\P}_{\mu}= (\mathbf{I}-\P_{\mu})$, and $\P_{\mu}$ is defined as 
\begin{align*}
\bega
\P_{\mu}(f)&= \lw\{\int_{\R^3} \sqrt{\mu}f dv\rw\} \sqrt{\mu} + \lw\{\int_{\R^3} \frac{v}{\sqrt{k_B}} \sqrt{\mu}f dv \rw\} \cdot \frac{v}{\sqrt{k_B}}\sqrt{\mu} \cr
&\quad + \lw\{ \int_{\R^3} \frac{1}{\sqrt{6 }}\lw(\frac{|v|^2-3k_B}{k_B}\rw)\sqrt{\mu}f dv \rw\} \frac{1}{\sqrt{6 }}\lw(\frac{|v|^2-3k_B}{k_B}\rw)\sqrt{\mu}.
\enda
\end{align*}
Let us define a function \(g\) as \(g(\tilde{v}) = \frac{G(\sqrt{\mathrm{\Theta}^{\e}}\tilde{v}+\mathrm{U}^{\e})}{\sqrt{\mu(\tilde{v})}}\). 
Once we write \(\frac{v-\mathrm{U}^{\e}}{\sqrt{\mathrm{\Theta}^{\e}}} = \tilde{v}\), then we have \(M^{\e}(v) = \frac{\mathrm{P}^{\e}}{\sqrt{\mathrm{\Theta}^{\e}}^3}\mu(\tilde{v})\) and
\begin{align}\label{Gtog}
\bega
\frac{\AC{\P}G(v)}{M^{\e}(v)} = \frac{G(v)-\P G(v)}{M^{\e}(v)} = \frac{\sqrt{\mathrm{\Theta}^{\e}}^3}{\mathrm{P}^{\e}}\AC{\P}_{\mu}\bigg(\frac{G(\sqrt{\mathrm{\Theta}^{\e}}\tilde{v}+\mathrm{U}^{\e})}{\sqrt{\mu(\tilde{v})}}\bigg) \frac{1}{\sqrt{\mu(\tilde{v})}} = \frac{\sqrt{\mathrm{\Theta}^{\e}}^3}{\mathrm{P}^{\e}} \frac{\AC{\P}_{\mu}g(\tilde{v})}{\sqrt{\mu(\tilde{v})}}.
\enda
\end{align}
By applying the change of variables \(\frac{v-\mathrm{U}^{\e}}{\sqrt{\mathrm{\Theta}^{\e}}} = \tilde{v}\) and \(\frac{v_*-\mathrm{U}^{\e}}{\sqrt{\mathrm{\Theta}^{\e}}} = \tilde{v}_*\) to \eqref{LGG}, we get
\begin{align*}
\bega
&\int_{\R^3}\mathcal{L}(G) G|M^{\e}|^{-1} dv \cr 
&= \frac{1}{4}\sqrt{\mathrm{\Theta}^{\e}}^7 \int_{\R^3} \int_{\mathbb{R}^3}\int_{\mathbb{S}^2_+}|(\tilde{v}-\tilde{v}_*)\cdot w| \mu(\tilde{v})\mu(\tilde{v}_*)\bigg[ \frac{\AC{\P}_{\mu}g(\tilde{v})}{\sqrt{\mu(\tilde{v})}} +\frac{\AC{\P}_{\mu}g(\tilde{v}_*)}{\sqrt{\mu(\tilde{v}_*)}} -\frac{\AC{\P}_{\mu}g(\tilde{v}')}{\sqrt{\mu(\tilde{v}')}} -\frac{\AC{\P}_{\mu}g(\tilde{v}_*')}{\sqrt{\mu(\tilde{v}_*')}} \bigg]^2 dwd\tilde{v}_* d\tilde{v} \cr 
&= \sqrt{\mathrm{\Theta}^{\e}}^7 \bigg\la L\bigg(\frac{G(\sqrt{\mathrm{\Theta}^{\e}}\tilde{v}+\mathrm{U}^{\e})}{\sqrt{\mu(\tilde{v})}}\bigg), \bigg(\frac{G(\sqrt{\mathrm{\Theta}^{\e}}\tilde{v}+\mathrm{U}^{\e})}{\sqrt{\mu(\tilde{v})}}\bigg) \bigg\ra_{L^2_{\tilde{v}}}.
\enda
\end{align*}
Thus, we have
\begin{align*}
\la \mathcal{L}(G), G|M^{\e}|^{-1} \ra_{L^2_v} &\geq \sigma_L \sqrt{\mathrm{\Theta}^{\e}}^7 \int_{\R^3}\nu(\tilde{v}) |\AC{\P}_{\mu}g(\tilde{v})|^2 d\tilde{v}  
= \sigma_L \mathrm{P}^{\e}\sqrt{\mathrm{\Theta}^{\e}} \int_{\R^3}\nu\lw(\frac{v-\mathrm{U}^{\e}}{\sqrt{\mathrm{\Theta}^{\e}}}\rw) \bigg|\frac{\AC{\P}G(v)}{\sqrt{M^{\e}(v)}}\bigg|^2dv,
\end{align*}
where we used the change of variables \(\tilde{v} = \frac{v-\mathrm{U}^{\e}}{\sqrt{\mathrm{\Theta}^{\e}}}\) with \(d\tilde{v} = \frac{1}{\sqrt{\mathrm{\Theta}^{\e}}^3} dv\) and $M^{\e}(v) = \frac{\mathrm{P}^{\e}}{\sqrt{\mathrm{\Theta}^{\e}}^3}\mu(\tilde{v})$ and \eqref{Gtog}.

\hide
Define gain term and loss term, to denote $\mathcal{N}(F,G)=\mathcal{N}^+(F,G)-\mathcal{N}^-(F,G)$: 
\begin{align}\label{Qdefgl}
\bega
\mathcal{N}^+(F,G)&:=   \frac{1}{2}\int_{\mathbb{R}^3}\int_{\mathbb{S}^2_+}B(v-v_*,w)\big[ F(v_*')G(v')+F(v')G(v_*')\big]dwdv_* \cr 
\mathcal{N}^-(F,G)&:=   \frac{1}{2}\int_{\mathbb{R}^3}\int_{\mathbb{S}^2_+}B(v-v_*,w)\big[F(v_*)G(v)+F(v)G(v_*)\big]dwdv_*
\enda
\end{align}
{\bf (Proof of \eqref{nonlin}) }
We only see the loss term in \eqref{Qdefgl}. 
Applying the same symmetry in \eqref{FGsym}, we have 
\begin{align*}
\int_{\R^3}&\mathcal{N}^-(F,G)H |M^{\e}|^{-1} dv = \frac{1}{2}\int_{\R^3} \int_{\mathbb{R}^3}\int_{\mathbb{S}^2_+}|(v-v_*)\cdot w|\cr
& \times \bigg[\frac{F(v)}{\sqrt{M^{\e}(v)}}\frac{G(v_*)}{\sqrt{M^{\e}(v_*)}}+ \frac{F(v_*)}{\sqrt{M^{\e}(v_*)}}\frac{G(v)}{\sqrt{M^{\e}(v)}} \bigg]\frac{H(v)}{\sqrt{M^{\e}(v)}} \sqrt{M^{\e}(v_*)} dwdv_*dv
\end{align*}
We only consider the first term. Applying H\"{o}lder inequality for $dv_*$-integral:
\begin{align*}
\bigg|\int_{\R^3}&\mathcal{N}^-(F,G)H |M^{\e}|^{-1} dv\bigg| \les \int_{\R^3}  \bigg( \int_{\R^3} \frac{|G(v_*)|^2}{M^{\e}(v_*)}\bigg)^{\frac{1}{2}} \bigg( \int_{\R^3} |v-v_*|^2 M^{\e}(v_*)\bigg)^{\frac{1}{2}}
\bigg[ \frac{F(v)}{\sqrt{M^{\e}(v)}}\frac{H(v)}{\sqrt{M^{\e}(v)}} \bigg]dv
\end{align*}
Explicit computation gives  
\begin{align*}
\bigg( \int_{\R^3} |v-v_*|^2 M^{\e}(v_*)\bigg)^{\frac{1}{2}} &= \bigg( \int_{\R^3} \big(|v-\mathrm{U}^{\e}|^2+|v_*-\mathrm{U}^{\e}|^2\big) M^{\e}(v_*)\bigg)^{\frac{1}{2}} \cr 
&= \bigg(|v-\mathrm{U}^{\e}|^2\mathrm{P}^{\e} + 3k_B \mathrm{P}^{\e}\mathrm{\Theta}^{\e}\bigg)^{\frac{1}{2}} \cr 
&\leq \sqrt{3k_B\mathrm{P}^{\e}\mathrm{\Theta}^{\e}} \bigg( 1+\frac{|v-\mathrm{U}^{\e}|}{\sqrt{3k_B\mathrm{\Theta}^{\e}}} \bigg)
\end{align*}
Next, applying H\"{o}lder inequality for $dv$-integral gives 
\begin{align*}
\bigg|\int_{\R^3}\mathcal{N}^-(F,G)H |M^{\e}|^{-1} dv\bigg| &\les  \sqrt{k_B\mathrm{P}^{\e}\mathrm{\Theta}^{\e}} \bigg( \int_{\R^3} \frac{|G(v_*)|^2}{M^{\e}(v_*)}\bigg)^{\frac{1}{2}} \bigg(\int_{\R^3} \bigg( 1+\frac{|v-\mathrm{U}^{\e}|}{\sqrt{\mathrm{\Theta}^{\e}}} \bigg)\frac{|F(v)|^2}{M^{\e}(v)}dv \bigg)^{\frac{1}{2}} \cr 
&\quad \times 
\bigg(\int_{\R^2}\bigg( 1+\frac{|v-\mathrm{U}^{\e}|}{\sqrt{\mathrm{\Theta}^{\e}}} \bigg)\frac{|H(v)|^2}{ {M^{\e}(v)}} dv \bigg)^{\frac{1}{2}}
\end{align*}
{\bf (Proof of \eqref{nonlin2}) }
We only consider the gain term defined in \eqref{Qdefgl}. We use $|v-v_*|=|v'-v_*'|$ and apply change of variable $(v,v_*)\rightarrow(v',v_*')$ to have 
\begin{align*}
\int_{\R^3}&\la v \ra^n|\mathcal{N}^+(F,G)|^2|M^{\e}|^{-1} dv \les  \int_{\R^3} \int_{\R^3}(1+|v'|^n+|v_*'|^n)|v'-v_*'|^2 \bigg[\frac{|F(v')|^2}{M^{\e}(v')}\frac{|G(v_*')|^2}{M^{\e}(v_*')}+ \frac{|F(v_*')|^2}{M^{\e}(v_*')}\frac{|G(v')|^2}{M^{\e}(v')}\bigg]M^{\e}(v_*) dv_*'dv',
\end{align*}
where we used $(1+|v|)^n \les 1+|v'|^n+|v_*'|^n$. 
Using $(1+|v'|^n+|v_*'|^n)|v'-v_*'|^2 \les (1+|v'|^{n+2})+(1+|v_*'|^{n+2})$, we get
\begin{align*}
\int_{\R^3}&\la v \ra^n|\mathcal{N}^+(F,G)|^2|M^{\e}|^{-1} dv \cr 
&\les  \int_{\R^3}(1+|v'|^{n+2})\frac{|F(v')|^2}{M^{\e}(v')}dv'\int_{\R^3}\frac{|G(v_*')|^2}{M^{\e}(v_*')}dv_*' + \int_{\R^3}\frac{|F(v')|^2}{M^{\e}(v')}dv'\int_{\R^3}(1+|v_*'|^{n+2})\frac{|G(v_*')|^2}{M^{\e}(v_*')}dv_*' . 
\end{align*}
The loss term $\mathcal{N}^-(F,G)$ can be estimated similarly.
\unhide
\end{proof}

\begin{lemma} 
Let $w(v) = e^{c_1 |v|^2}$ for a positive constant $c_1 > 0$. Then, the kernel $
\mathbf{k}_w(v, v_*) = \mathbf{k}(v, v_*) \frac{w(v)}{w(v_*)}$, where $\mathbf{k}(v, v_*)$ is defined in \eqref{kdef}, satisfies
\begin{align}\label{kw}
\mathbf{k}_w(v, v_*) \leq \frac{1}{|v - v_*|} e^{-C |v - v_*|^2}, \qquad \int_{\R^3} \mathbf{k}_w(v, v_*) dv_* \les \frac{1}{1 + |v|}.
\end{align}
Furthermore, $\Gamma(f,g)= \frac{1}{\sqrt{\mu}}\mathcal{N}(\sqrt{\mu}f,\sqrt{\mu}g)$, satisfies
\begin{align}\label{wgammav}
\bigg| w \Gamma\bigg(\frac{\p^{\beta} h}{w}, \frac{\p^{\al - \beta} h}{w}\bigg) \bigg| \leq C \nu(v) \|\p^{\beta} h\|_{L^\infty_v} \|\p^{\al - \beta} h\|_{L^\infty_v},
\end{align}
where $C \geq 1$ is a positive constant depending on $c_1$ in $w(v) = e^{c_1 |v|^2}$. \\
For $\bar{\mu}(v)$, $\bar{\mathbf{k}}(v, v_*)$ and $\bar{\mathbf{k}}_{w}(v, v_*)$ defined in \eqref{muMdef} and \eqref{LkMdef} we have 
\begin{align}\label{kw2}
\bega
&|\bar{\mathbf{k}}(v, v_*)|\leq C\Big(|v-v_*|+\frac{1}{|v-v_*|}\Big)\exp\bigg(-C(|v-v_*|^2)-C\frac{||v|^2-|v_*|^2|^2}{|v-v_*|^2}\bigg), \cr
&\sup_{v\in\R^3}\int_{\R^3} \bar{\mathbf{k}}_{w}(v, v_*) dv_* \leq C, \qquad
\int_{\R^3} \bar{\mathbf{k}}_{w}(v, v_*)e^{c|v-v_*|^2} dv_* \leq C, \quad \mbox{for small constant } \quad c>0,
\enda
\end{align}
and
\begin{align}\label{wgammav2}
\bigg| \frac{w}{\sqrt{\bar{\mu}}} \mathcal{N}\bigg(\frac{\p^{\beta} h\sqrt{\bar{\mu}}}{w}, \frac{\p^{\al - \beta} h\sqrt{\bar{\mu}}}{w}\bigg) \bigg| \leq C \nu(v) \|\p^{\beta} h\|_{L^\infty_v} \|\p^{\al - \beta} h\|_{L^\infty_v}.
\end{align} 
\end{lemma}

\begin{proof}
For \eqref{kw} and \eqref{wgammav}, see (3.22), (4.139), (4.137), and Lemma 2 in \cite{JangKim}.
For \eqref{kw2} and \eqref{wgammav2}, see \cite{GuoJJCPAM}.
\end{proof}

\hide
\begin{lemma}
The local Maxwellian \eqref{M-def} satisfies the following moment integrals:
\begin{align}\label{M-moment}
\int_{\R^3}M\bmx 
1\\ 
|v-\mathrm{U}|^2\\ 
|v-\mathrm{U}|^4\\ 
|v-\mathrm{U}|^6\\ 
(v_i-\mathrm{U}_i)(v_j-\mathrm{U}_j)\\
(v_j-\mathrm{U}_j)(v_i-\mathrm{U}_i)|v-\mathrm{U}|^2 \\ 
(v_i-\mathrm{U}_i)^2(v_j-\mathrm{U}_j)^2 \\ 
(v_i-\mathrm{U}_i)^4
\emx dv
=
\bmx 
\mathrm{P}\\
3k_B \mathrm{P}\mathrm{\Theta}\\
15k_B^2 \mathrm{P}|\mathrm{\Theta}|^2\\
105 k_B^3 \mathrm{P} |\mathrm{\Theta}|^3\\
\delta_{ij}k_B\mathrm{P} \mathrm{\Theta}\\
5\delta_{ij}k_B^2\mathrm{P} |\mathrm{\Theta}|^2\\
k_B^2\mathrm{P} |\mathrm{\Theta}|^2 \\ 
3k_B^2\mathrm{P} |\mathrm{\Theta}|^2
\emx.
\end{align}
\end{lemma}
\unhide

\hide
\begin{lemma}\label{ABcomp}
There exist universal constants $\eta_0 > 0$ and $\eta_1 > 0$ such that
\begin{align}
&(1)~ \la \mathcal{L}^{-1}(\mathfrak{R}^{\e}_{ij} M^{\e}), \mathfrak{R}^{\e}_{kl} \ra_{L^2_v} = \eta_0 |k_B \mathrm{\Theta}^{\e}|^{\frac{3}{2}} 
\lw(\delta_{ik} \delta_{jl} + \delta_{il} \delta_{kj} - \frac{2}{3} \delta_{ij} \delta_{kl}\rw), \label{<AA>} \\ 
&(2)~ \la \mathcal{L}^{-1}(\mathcal{Q}^{\e}_i M^{\e}), \mathcal{Q}^{\e}_k \ra_{L^2_v} = \frac{5}{2} \eta_1 |k_B \mathrm{\Theta}^{\e}|^{\frac{5}{2}} \delta_{ik}, \label{<BB>} \\ 
&(3)~ \la \mathcal{L}^{-1}(\mathcal{Q}^{\e}_k M^{\e}), \mathfrak{R}^{\e}_{ij} \ra_{L^2_v} = 0, \label{<AB>}
\end{align}
where $\mathfrak{R}^{\e}_{ij}$ and $\mathcal{Q}^{\e}_i$ are defined in \eqref{ABdef}.
\end{lemma}

\begin{proof}
In the regime near a global Maxwellian $\mu_*(v):=\frac{1}{\sqrt{2\pi}}e^{-\frac{|v|^2}{2}}$, we define
\beq \bega \label{ABs-def}
&\hat{A}_{ij}^{\mu} = \lw(v_i v_j - \frac{|v|^2}{3}\delta_{ij}\rw)\sqrt{\mu}, \qquad A_{ij}^{\mu} = L^{-1}\hat{A}_{ij}^{\mu}, \qquad \hat{B}_i^{\mu} = \lw(v_i \frac{|v|^2 - 5}{2}\rw)\sqrt{\mu}, \qquad B_i^{\mu} = L^{-1}\hat{B}_i^{\mu},
\enda \eeq
then, from \cite{Guo-NS}, we have
\begin{align*}
\la A_{ij}^{\mu} ,\hat{A}_{kl}^{\mu} \ra_{L^2_v}&=\eta_0 \lw(\delta_{ik}\delta_{jl}+\delta_{il}\delta_{kj}-\frac 2 3 \delta_{ij}\delta_{kl}
\rw), \qquad
\la B_i^{\mu} ,\hat{B}_k^{\mu} \ra_{L^2_v}
= \frac{5}{2}\eta_1\delta_{ik}, \qquad 
\la B_l^{\mu} ,\hat{A}_{ij}^{\mu}\ra_{L^2_v}=0,
\end{align*}
for any \(1 \leq i,j,k,l \leq 3\), where \(\eta_0\) and \(\eta_1\) correspond to viscosity and heat conductivity, respectively.
Writing \(\tilde{v} = \frac{v - \mathrm{U}^{\e}}{\sqrt{k_B\mathrm{\Theta}^{\e}}}\), we obtain the following relations between $(\mathfrak{R}^{\e}_{ij},\mathcal{Q}^{\e}_i)$ and $(\hat{A}_{ij}^{\mu},\hat{B}_i^{\mu})$:
\begin{align*}
\bega
&\mathfrak{R}^{\e}_{ij}(v) 
= \frac{k_B\mathrm{\Theta}^{\e}}{\sqrt{\mu(\tilde{v})}} \hat{A}_{ij}^{\mu}(\tilde{v}), \qquad 
\mathcal{L}^{-1}(\mathfrak{R}^{\e}_{ij} M^{\e})(v) = \frac{1}{k_B\mathrm{\Theta}^{\e}} A_{ij}^{\mu}(\tilde{v}) \sqrt{\mu(\tilde{v})}, \\ 
&\mathcal{Q}^{\e}_i(v) 
= \frac{(k_B\mathrm{\Theta}^{\e})^{3/2}}{\sqrt{\mu(\tilde{v})}} \hat{B}_i^{\mu}(\tilde{v}), \quad 
\mathcal{L}^{-1}(\mathcal{Q}^{\e}_i M^{\e})(v) = \frac{1}{\sqrt{k_B\mathrm{\Theta}^{\e}}} B_i^{\mu}(\tilde{v}) \sqrt{\mu(\tilde{v})}.
\enda
\end{align*}
Using these relations, we derive the desired result.
\end{proof}
\unhide

\section{Basic Inequalities and Computations}\label{A.B}

\begin{lemma} The following inequalities hold:
\begin{itemize}
\item Sobolev Embedding: 
\begin{align}\label{embed}
\bega
&\|u\|_{L^\infty(\R^d)} \leq \|u\|_{H^k(\R^d)}, \quad \text{for} \quad k > d/2, \cr
&\|u\|_{L^{p^*}(\R^d)} \leq \|u\|_{W^{1,p}(\R^d)}, \quad \text{where} \quad \frac{1}{p^*} = \frac{1}{p} - \frac{1}{d}.
\enda
\end{align}

\item Gagliardo-Nirenberg interpolation inequality:
\begin{align}\label{Ga-Ni}
\|D^j u\|_{L^p(\R^d)} \leq C \|D^m u \|_{L^r(\R^d)}^{\theta} \|u\|_{L^q(\R^d)}^{1-\theta},
\end{align}
where \(1 \leq q, r \leq \infty\), \(j < m\), \(1 \leq p <\infty \), \(\theta \in [0,1]\), and \(\sigma\) satisfies:
\begin{align*}
\frac{1}{p} = \frac{j}{d} + \theta\lw(\frac{1}{r} - \frac{m}{d} \rw) + \frac{1-\theta}{q}.
\end{align*}
If \(r > 1\) and \(m - j - n/r\) is a non-negative integer, then an additional assumption \(j/m \leq \theta < 1\) is required.

\item Agmon's Inequality for $x\in\R^d$:
\begin{align}\label{Agmon}
\bega
&\|u\|_{L^\infty(\R^2)} \leq C\|u\|_{L^2(\R^2)}^{\frac{1}{2}} \|u\|_{H^2(\R^2)}^{\frac{1}{2}}, \quad 
&\|u\|_{L^\infty(\R^3)} \leq C\|u\|_{H^1(\R^3)}^{\frac{1}{2}} \|u\|_{H^2(\R^3)}^{\frac{1}{2}}.
\enda
\end{align}

\item Ladyzhenskaya's Inequality:
\begin{align}\label{Lady}
\|u\|_{L^4(\R^2)} \leq C\|u\|_{L^2(\R^2)}^{\frac{1}{2}} \|\nabla_xu\|_{L^2(\R^2)}^{\frac{1}{2}}, 
\quad
\|u\|_{L^4(\R^3)} \leq C\|u\|_{L^2(\R^3)}^{\frac{1}{4}} \|\nabla_xu\|_{L^2(\R^3)}^{\frac{3}{4}}. 
\end{align}




\end{itemize}
\end{lemma}

\hide \\
(Proof of \eqref{Agmon}:) Applying H\"{o}lder inequality, we have for $s_1 < \frac{d}{2} < s_2$, 
\begin{align*}
u(x) &= \int_{\R^3} e^{2\pi i x \cdot \xi} \hat{u}(\xi) d\xi \cr 
&= \int_{|\xi|\leq L} e^{2\pi i x \cdot \xi} \hat{u}(\xi) \frac{(1+|\xi|)^{s_1}}{(1+|\xi|)^{s_1}}d\xi + \int_{|\xi|\leq L} e^{2\pi i x \cdot \xi} \hat{u}(\xi) \frac{|\xi|^{s_2}}{|\xi|^{s_2}}d\xi \cr 
&\leq C\|u\|_{H^{s_1}_x} \bigg(\int_{|\xi|\leq L}\frac{1}{(1+|\xi|)^{2s_1}}d\xi\bigg)^{\frac{1}{2}} + C\|u\|_{H^{s_2}_x} \bigg(\int_{|\xi|\geq L}\frac{1}{|\xi|^{2s_2}}d\xi\bigg)^{\frac{1}{2}} \cr 
&\leq C\|u\|_{H^{s_1}_x} \bigg(4\pi\int_{r\leq L}\frac{r^{d-1}}{(1+r)^{2s_1}}dr\bigg)^{\frac{1}{2}} + C\|u\|_{H^{s_2}_x} \bigg(4\pi\int_{r\geq L}\frac{r^{d-1}}{r^{2s_2}}dr\bigg)^{\frac{1}{2}} \cr 
&\leq C\|u\|_{H^{s_1}_x} \bigg(L^{d-2s_1}\bigg)^{\frac{1}{2}} + C\|u\|_{H^{s_2}_x}\bigg(\frac{1}{2s_2-d}\frac{1}{L^{2s_2-d}}\bigg)^{\frac{1}{2}} \cr 
&\leq C\|u\|_{H^{s_1}_x} L^{\frac{d-2s_1}{2}} + C\|u\|_{H^{s_2}_x} L^{-\frac{2s_2-d}{2}}
\end{align*}
Choosing $L^{\frac{d-2s_1}{2}+\frac{2s_2-d}{2}}:= \|u\|_{H^1_x}^{-1}\|u\|_{H^2_x}$ gives 
\begin{align*}
u(x) &\leq C\|u\|_{H^{s_1}_x}^{\frac{s_2-\frac{d}{2}}{s_2-s_1}}\|u\|_{H^{s_2}_x}^{\frac{\frac{d}{2}-s_1}{s_2-s_1}}
\end{align*}
Let $\vartheta =\frac{\frac{d}{2}-s_1}{s_2-s_1}$, we have 
\begin{align*}
u(x) &\leq C\|u\|_{H^{s_1}_x}^{1-\vartheta}\|u\|_{H^{s_2}_x}^{\vartheta}, \quad \frac{d}{2}  = (1-\vartheta)s_1 + \vartheta s_2, \quad s_1<\frac{d}{2}<s_2.
\end{align*}
Then the inequality \eqref{Agmon} is the special case of it.
\textcolor{blue}{To get the Local existence for general dimension. Note for the estimate \eqref{pMinf}. Remark: $\kappa$ singularity comes from the following: $\p^2M$ produces $\|\nabla_x^2 u\|_{L^\infty_x}$. In the worst case, $\|\nabla_x^2 u\|_{L^\infty_x} \leq \|u\|_{H^4_x} \leq \kappa^{-\frac{1}{4}}\mathcal{E}_{tot}^{\frac{1}{2}} $. We should check that, the situation is the same for higher dimension, $d=4,5,6,\cdots$.
What we have in mind: For any $d\in \mathbb{N}$, choose $\mathrm{N} = \lfloor\frac{d}{2}+1 \rfloor$. Then 
Agmon's inequality \eqref{Agmon} gives 
\begin{align}
\bega
\|u\|_{L^\infty_x} \les \begin{cases}
\|u\|_{H^{\frac{d}{2}-1}}^{\frac{1}{2}} \|u\|_{H^{\frac{d}{2}+1}}^{\frac{1}{2}}, \quad \mbox{when} \quad d: \mbox{even} \\ 
\|u\|_{H^{\frac{d}{2}-\frac{1}{2}}}^{\frac{1}{2}} \|u\|_{H^{\frac{d}{2}+\frac{1}{2}}}^{\frac{1}{2}}, \quad \mbox{when} \quad d: \mbox{odd}
\end{cases} 
\enda
\end{align}
Then we have 
\begin{align}
\bega
\|\nabla_x^2u\|_{L^\infty_x} \les \begin{cases} \kappa^{-\frac{1}{4}}\mathcal{E}_M^{\frac{1}{4}} \mathcal{E}_{top}^{\frac{1}{4}},\quad \mbox{if we choose}\quad \mathrm{N} = \lfloor\frac{d}{2}+1 \rfloor \\ 
\mathcal{E}_M^{\frac{1}{2}},\quad \mbox{if we choose}\quad \mathrm{N} > \lfloor\frac{d}{2}+1 \rfloor
\end{cases}
\enda
\end{align}
}
\unhide

\hide
\begin{lemma}[Brezis–Gallou\"et inequality] For a bounded domain $\Omega \subset \mathbb{R}^2$, and for a function $u\in H^2(\Omega)$, we have 
\begin{align}\label{Bre-Gal}
\|u\|_{L^\infty_x} \leq C\|u\|_{H^1_x}\bigg(1+ \Big(\log \big(1+\frac{\|u\|_{H^2_x}}{\|u\|_{H^1_x}}\big)\Big)^{\frac{1}{2}}\bigg).
\end{align}
\end{lemma}
\unhide

\hide
\begin{lemma}[Time Embedding] 
For any function \(f \in H^1_t\), the following inequality holds:
\begin{align}\label{tembedding}
\sup_{0 \leq t \leq T} |f(t)| \les \frac{1}{\sqrt{T}} \|f\|_{L^2_T} + \sqrt{T} \|\p_t f\|_{L^2_T}.
\end{align}
\end{lemma}
\begin{proof}
The proof can be found in Appendix A of \cite{JangKim}.
\hide
We claim the following inequality
\Be\label{SE1}
|f(t)| \leq \frac{2}{T}\int_0^t |f(s)|ds + \int_0^t |f'(s)| ds.
\Ee
If \eqref{SE1} holds, then we have a temporal embedding \eqref{tembedding} by H\"{o}lder inequality. Now, we prove \eqref{SE1}.
When $t\in [0,T/2]$, we have
\begin{align*}
f(t) = \frac{1}{T/2} \int_t^{t+\frac{T}{2}} \left(f(s) - \int_t^s f'(\tau)d\tau\right) ds.
\end{align*}
Hence, for $t \leq s \leq t+T/2$
\begin{align*}
|f(t)| &\leq \frac{1}{T/2} \int_t^{t+\frac{T}{2}}|f(s)|ds + \frac{1}{T/2} \int_t^{t+\frac{T}{2}}\int_t^s |f'(\tau)|d\tau ds \cr
&\leq \frac{1}{T/2} \int_t^{t+\frac{T}{2}}|f(s)|ds + \frac{1}{T/2} \int_t^{t+\frac{T}{2}} \int_t^{t+\frac{T}{2}} |f'(\tau)|d\tau ds \cr
&\leq \frac{1}{T/2} \int_t^{t+\frac{T}{2}}|f(s)|ds + \int_t^{t+\frac{T}{2}} |f'(\tau)|d\tau \cr
&\leq \frac{1}{T/2} \int_0^T|f(s)|ds + \int_0^T |f'(s)|ds.
\end{align*}
On the other hand, when $t\in (T/2,T]$, we have
\begin{align*}
f(t) = \frac{1}{T/2} \int_{t-\frac{T}{2}}^t \left(f(s) + \int_s^t f'(\tau)d\tau\right) ds.
\end{align*}
Therefore, for $t-T/2 \leq s \leq T$,
\begin{align*}
|f(t)| &\leq \frac{1}{T/2} \int_{t-\frac{T}{2}}^t|f(s)|ds + \frac{1}{T/2} \int_{t-\frac{T}{2}}^t\int_s^t |f'(\tau)|d\tau ds \cr
&\leq \frac{1}{T/2} \int_{t-\frac{T}{2}}^t|f(s)|ds + \frac{1}{T/2} \int_{t-\frac{T}{2}}^t \int_{t-\frac{T}{2}}^t |f'(\tau)|d\tau ds \cr
&\leq \frac{1}{T/2} \int_{t-\frac{T}{2}}^t|f(s)|ds + \int_{t-\frac{T}{2}}^t |f'(\tau)|d\tau \cr
&\leq \frac{1}{T/2} \int_0^T|f(s)|ds + \int_0^T |f'(s)|ds.
\end{align*}
\unhide
\end{proof}
\unhide

\begin{lemma}\cite{MaBe} 
Let \(u_*, v_* \in L^\infty_x \cap H^N_x(\mathbb{R}^2)\) and \(N \in \mathbb{N} \cup \{0\}\). For any multi-indices \(\al_x\) and \(\beta_x\) such that \(0 \leq \beta_x \leq \al_x\) and \(|\al_x| \leq N\), the following hold:
\begin{align}\label{uvHk}
\|\p^{\beta_x} u_* \p^{\al_x - \beta_x} v_*\|_{L^2_x} &\leq C \Big(\|u_*\|_{L^\infty_x} \|v_*\|_{H^N_x} + \|v_*\|_{L^\infty_x} \|u_*\|_{H^N_x} \Big),
\end{align}
and
\begin{align}\label{commutator}
\sum_{0 \leq |\al_x| \leq N} \|\p^{\al_x}(u_* \cdot \nabla_x v_*) - u_* \cdot \p^{\al_x} \nabla_x v_*\|_{L^2_x} &\leq C \Big( \|\nabla_x u_*\|_{L^\infty_x} \|v_*\|_{H^N_x} + \|\nabla_x v_*\|_{L^\infty_x} \|u_*\|_{H^N_x} \Big).
\end{align}
If \(u_* \in H^3_x(\mathbb{R}^2)\) is divergence-free and satisfies \(\nabla_x \times u_* = \w_*\), then
\begin{align}\label{puw}
\bega
\|\w_*\|_{H^k_x} \les \|\nabla_x u_* \|_{H^k_x} \les \|\w_*\|_{H^k_x},
\enda
\end{align}
\begin{align}\label{PotenLp}
\|\nabla_x u_*\|_{L^p_x} &\leq C_p \|\w_*\|_{L^p_x}, \quad \text{for} \quad 1 < p < \infty.
\end{align}
Here $C_p \leq Cp$ for $2\leq p<\infty$.
And
\begin{align}\label{Be-Ma}
\|\nabla_x u_*\|_{L^\infty_x} &\les \lw(1 + \ln^+ \|\w_*\|_{H^2_x} \rw) \lw(1 + \|\w_*\|_{L^\infty_x} \rw),
\end{align}
where \(\ln^+(x)\) equals \(\ln(x)\) for \(x > 1\), and \(0\) otherwise.
\end{lemma}

\begin{proof}
We refer the reader to \cite{MaBe} for a detailed proof. In particular, for the estimate \eqref{Be-Ma}, see the proof of Proposition 3.8 in \cite{MaBe}.
\hide
(Proof of \eqref{uvHk}) Let $|\al_x|=N$. Then the H\"{o}lder inequality gives 
\begin{align}\label{uvHol}
\bega
\|\p^{\beta_x}u\p^{\al_x-\beta_x}v\|_{L^2_x}&= \lw(\int_{\Omega}|\p^{\beta_x}u|^2|\p^{\al_x-\beta_x}v|^2dx\rw)^{\frac{1}{2}} \cr 
&\leq \lw(\int_{\Omega}|\p^{\beta_x}u|^{\frac{2N}{|\beta_x|}}dx\rw)^{\frac{|\beta_x|}{2N}} \lw(\int_{\Omega} |\p^{\al_x-\beta_x}v|^{\frac{2N}{|\al_x-\beta_x|}}dx\rw)^{\frac{|\al_x-\beta_x|}{2N}} \cr 
&=  \|\p^{\beta_x}u\|_{L^{\frac{2N}{|\beta_x|}}_x} \|\p^{\al_x-\beta_x}v\|_{L^{\frac{2N}{|\al_x-\beta_x|}}_x}.
\enda
\end{align}
We use the Gagliardo–Nirenberg interpolation inequality: 
\begin{align}\label{Ga-Ni}
\|\p^j u \|_{L^p_x} \leq C \|u\|_{L^q_x}^{1-\ta} \|\p^m u \|_{L^r_x}^{\ta}
\end{align}
for $1\leq q \leq \infty$ and $1\leq r \leq \infty$, $j\leq m$, $p\geq1$, $\theta\in[0,1]$ satisfying 
\begin{align*}
\frac{1}{p}=\frac{j}{n}+\theta\lw(\frac{1}{r}-\frac{m}{n} \rw) +\frac{1-\ta}{q}.
\end{align*}
Once we choose $\theta = \frac{j}{m}$, $q=\infty$, $r=2$ with dimension $n=2$, then we have $p= \frac{2m}{j}$ as follows:
\begin{align*}
\frac{1}{p}=\frac{j}{2}+\frac{j}{m}\lw(\frac{1}{2}-\frac{m}{2} \rw), \qquad \mbox{so that} \qquad p= \frac{2m}{j}.
\end{align*}
In this case, the inequality \eqref{Ga-Ni} becomes
\begin{align}\label{Ga-Ni2}
\|D^j u \|_{L^{\frac{2m}{j}}_x} \leq C \|u\|_{L^\infty_x}^{1-\frac{j}{m}} \|D^m u \|_{L^2_x}^{\frac{j}{m}}
\end{align}
Choosing $(j,m)=(|\beta_x|,|\al_x|)$ and $(j,m)=(|\al_x-\beta_x|,|\al_x|)$, we have 
\begin{align}\label{Ga-Ni3}
\|\p^{\beta_x} u \|_{L^{\frac{2|\al_x|}{|\beta_x|}}_x} \leq C \|u\|_{L^\infty_x}^{1-\frac{|\beta_x|}{|\al_x|}} \|\p^{|\al_x|}u\|_{L^2_x}^{\frac{|\beta_x|}{|\al_x|}}, \qquad \|\p^{\al_x-\beta_x} u \|_{L^{\frac{2|\al_x|}{|\al_x-\beta_x|}}_x} \leq C \|u\|_{L^\infty_x}^{1-\frac{|\al_x-\beta_x|}{|\al_x|}} \|\p^{|\al_x|}u\|_{L^2_x}^{\frac{|\al_x-\beta_x|}{|\al_x|}}.
\end{align}
where we wrote $\|\p^{|\al_x|}u\|_{L^p_x} = \sum_{\al_x=N}\|\p^{\al_x}u\|_{L^p_x}$ for $|\al_x|=N$.

We apply \eqref{Ga-Ni3} to \eqref{uvHol} to get
\begin{align*}
\|\p^{\beta_x}u\p^{\al_x-\beta_x}v\|_{L^2_x} &\leq C \|u\|_{L^\infty_x}^{1-\frac{|\beta_x|}{N}} \|\p^{|\al_x|} u \|_{L^2_x}^{\frac{|\beta_x|}{N}} \times \|v\|_{L^\infty_x}^{1-\frac{|\al_x-\beta_x|}{N}} \|\p^{|\al_x|} v \|_{L^2_x}^{\frac{|\al_x-\beta_x|}{N}} \cr 
&\leq C \bigg(\|u\|_{L^\infty_x}\|\p^{|\al_x|} v \|_{L^2_x} \bigg)^{\frac{|\al_x-\beta_x|}{N}} \bigg( \|v\|_{L^\infty_x} \|\p^{|\al_x|} u \|_{L^2_x}\bigg)^{\frac{|\beta_x|}{N}} \cr 
&\leq C \|u\|_{L^\infty_x}\|\p^{|\al_x|} v \|_{L^2_x}  + C \|v\|_{L^\infty_x} \|\p^{|\al_x|} u \|_{L^2_x},
\end{align*}
where we used the Young's inequality
\begin{align*}
AB \leq \frac{1}{p}A^p + \frac{1}{q}B^q, \qquad \mbox{for} \qquad p=\frac{N}{|\al_x-\beta_x|}, \quad q=\frac{N}{|\beta_x|}.
\end{align*}
Taking summation for $|\al_x|\leq N$, we prove \eqref{uvHk}. \\
(Proof of \eqref{commutator}) Once we compute the following term $\p^{\al_x}(u_*\cdot\nabla_xv_*)-u_*\cdot\p^{\al_x}\nabla_xv_*$, then both of $u_*$ and $v_*$ have at least one derivative. Thus, we can write 
\begin{align*}
\sum_{0\leq |\al_x|\leq N}\|\p^{\al_x}(u_*\cdot\nabla_xv_*)-u_*\cdot\p^{\al_x}\nabla_xv_*\|_{L^2_x}&\leq \sum_{0\leq |\al'|\leq N-1}\sum_{0\leq \beta' \leq \al'}\|\p^{\beta'}(\nabla_xu)\p^{\al'-\beta'}(\nabla_xv)\|_{L^2_x} 
\end{align*}
Applying \eqref{uvHk}, we have the result \eqref{commutator}. 
\unhide
\end{proof}

\begin{lemma}[Faà di Bruno's formula]
For a multi-index \(\al\) with \(|\al| \geq 1\), the derivative of a composite function \(h(x_1, \cdots, x_d) = f(g_1(x_1, \cdots, x_d), \cdots, g_m(x_1, \cdots, x_d))\) satisfies
\begin{align}\label{Faadi}
\p^{\al}h = \sum_{1\leq|{\bf r}|\leq |\al|} \frac{\p^{|{\bf r}|}}{\p y_1^{r_1}\cdots \p y_m^{r_m}}f \sum_{s=1}^{|\al|} \sum_{p_s(\al,{\bf r})} (\al!) \prod_{j=1}^{s} \frac{[(\p^{\bl_j}g_1,\cdots,\p^{\bl_j}g_m)]^{\bk_j}}{\bk_j! (\bl_j!)^{|\bk_j|}},
\end{align}
where
\begin{align*}
\al! = \prod_{i=0}^2 \al_i!, \qquad 
\mathbf{v}^{\bk} = \prod_{i=1}^d v_i^{k_i}, \qquad \binom{\al}{{\bf r}} = \prod_{i=1}^d \binom{\al_i}{r_i} = \frac{\al!}{{\bf r}!(\al-{\bf r})!},
\end{align*}
for a \(d\)-dimensional vector \(\mathbf{v} = (v_1, \cdots, v_d)\) and \(\kappa = (k_1, \cdots, k_d)\).
The notation \(\mathbf{v} \prec \mathbf{u}\) means one of the following holds:
\begin{enumerate}
\item \(|\mathbf{v}| < |\mathbf{u}|\).
\item \(|\mathbf{v}| = |\mathbf{u}|\) and \(v_1 < u_1\).
\item \(|\mathbf{v}| = |\mathbf{u}|\), \(v_1 = u_1, \cdots, v_k = u_k\), and \(v_{k+1} < u_{k+1}\) for some \(1 \leq k < d\).
\end{enumerate}
The set \(p_s(\al, {\bf r})\) is defined as
\begin{align}\label{ps}
p_s(\al,{\bf r}) = \left\{\bk_1,\cdots,\bk_s;\bl_1,\cdots,\bl_s : |\bk_i|>0,~ 0 \prec \bl_1 \prec \cdots \prec \bl_s, ~\sum_{i=1}^{s}\bk_i = {\bf r}, ~\sum_{i=1}^{s}|\bk_i|\bl_i=\al \right\}.
\end{align}
Here, the vectors \(\bk_j\) and \(\bl_j\) are \(m\)-dimensional and \(d\)-dimensional, respectively, where each element is a non-negative integer.
\end{lemma}

\begin{proof}[\textbf{Proof of Lemma \ref{Mal}}]
Note that $M^{\e}$ is a function depending on the macroscopic fields \((\rho^{\e},u^{\e},\ta^{\e})\). Faà di Bruno's formula \eqref{Faadi} explicitly provides the coefficients arising from the chain rule due to the multi-index derivative \(\p^{\al}\). Applying Faà di Bruno's formula \eqref{Faadi}, we have 
\begin{align*}
\bega
\p^{\al} M^{\e} &= \sum_{1\leq i\leq |\al|} \sum_{|{\bf r}|=i}
\p_{\rho^{\e}}^{r_1} \p_{u^{\e}_1}^{r_2} \p_{u^{\e}_2}^{r_3} \p_{u^{\e}_3}^{r_4} \p_{\ta^{\e}}^{r_5} M^{\e}
\sum_{s=1}^{|\al|} \sum_{p_s(\al,{\bf r})} (\al!) \prod_{j=1}^{s} \frac{[\p^{\bl_j}(\rho^{\e},u^{\e}_1,u^{\e}_2,u^{\e}_3,\ta^{\e})]^{\bk_j}}{\bk_j! (\bl_j!)^{|\bk_j|}}.
\enda
\end{align*}
Here, ${\bf r}=(r_1,r_2, r_3,r_4,r_5)$ is a 5-vector with non-negative integer components \(r_i\). Similarly, $\bk_j$ is a 5-dimensional vector and $\bl_j$ is a 4-dimensional vector, both with non-negative integer components. The derivative $\p^{\bl_j}$ represents $\p^{\bl_j}=\p_{\tilde{t}}^{l_j^0}\p_{x_1}^{l_j^1}\p_{x_2}^{l_j^2}\p_{x_3}^{l_j^3}$.
It is worth noting that the term $|{\bf r}|$ corresponds to the power of \(\eps\), since each $(\rho^{\e},u^{\e},\ta^{\e})$ derivative of $M^{\e}$ gives $\e$, for example
\begin{align*}
\bega
\frac{\p M^{\e}}{\p \rho^{\e}} = \frac{\p \mathrm{P}^{\e}}{\p \rho^{\e}} \frac{\p M^{\e}}{\p \mathrm{P}^{\e}} = \e\mathrm{P}^{\e} \frac{\p M^{\e}}{\p \mathrm{P}^{\e}}.
\enda
\end{align*}
Thus, \(\p^{\al}M^{\e}\) can be decomposed based on the order of \(\eps^i\) as \(\p^{\al} M^{\e} = \sum_{1 \leq i \leq |\al|}\eps^i \Phi_{\al}^i M^{\e}\), where 
\begin{align}\label{Phiidef}
\bega
\e^i\Phi_{\al}^i M^{\e} := \sum_{|{\bf r}|=i}
\p_{\rho^{\e}}^{r_1} \p_{u^{\e}_1}^{r_2} \p_{u^{\e}_2}^{r_3} \p_{u^{\e}_3}^{r_4} \p_{\ta^{\e}}^{r_5} M^{\e}
\sum_{s=1}^{|\al|} \sum_{p_s(\al,{\bf r})} (\al!) \prod_{j=1}^{s} \frac{[\p^{\bl_j}(\rho^{\e},u^{\e}_1,u^{\e}_2,u^{\e}_3,\ta^{\e})]^{\bk_j}}{\bk_j! (\bl_j!)^{|\bk_j|}}.
\enda
\end{align}
This provides the result \eqref{RpM-def}.
Next, we verify \eqref{Phi1def}. When $|{\bf r}|=1$, we have
\begin{align*}
\e\Phi_{\al}^1 M^{\e} &= \sum_{|{\bf r}|=1}\frac{\p^1 M^{\e}}{\p( \rho^{\e})^{r_1}\p(u^{\e}_1)^{r_2}\p(u^{\e}_2)^{r_3}\p(u^{\e}_3)^{r_4} \p(\ta^{\e})^{r_5}}  \sum_{p_1(\al,{\bf r})} (\al!)  \frac{[\p^{\al}(\rho^{\e},u^{\e}_1,u^{\e}_2,u^{\e}_3,\ta^{\e})]^{{\bf r}}}{{\bf r}! \al!} \cr 
&= \lw(\frac{\p M^{\e}}{\p \rho^{\e}},\frac{\p M^{\e}}{\p u^{\e}_1},\frac{\p M^{\e}}{\p u^{\e}_2},\frac{\p M^{\e}}{\p u^{\e}_3},\frac{\p M^{\e}}{\p \ta^{\e}}\rw)\cdot \p^{\al}(\rho^{\e},u^{\e}_1,u^{\e}_2,u^{\e}_3,\ta^{\e}) \cr 
&=\lw(\e\mathrm{P}^{\e}\frac{\p M^{\e}}{\p \mathrm{P}^{\e}},\e\frac{\p M^{\e}}{\p u^{\e}_1},\e\frac{\p M^{\e}}{\p u^{\e}_2},\e\frac{\p M^{\e}}{\p u^{\e}_3},\e\mathrm{\Theta}^{\e}\frac{\p M^{\e}}{\p \mathrm{\Theta}^{\e}}\rw)\cdot \p^{\al}(\rho^{\e},u^{\e}_1,u^{\e}_2,u^{\e}_3,\ta^{\e}) .
\end{align*}
Then, dividing both sides by \(\e M^{\e}\), and using  
\begin{align}\label{MTAderi}
\bega
\frac{\p M^{\e}}{\p (\mathrm{P}^{\e}, \mathrm{U}^{\e}, \mathrm{\Theta}^{\e})} = \bigg( \frac{1}{\mathrm{P}^{\e}}, \frac{(v - \mathrm{U}^{\e})}{k_B \mathrm{\Theta}^{\e}}, \frac{1}{\mathrm{\Theta}^{\e}} \lw( \frac{|v - \mathrm{U}^{\e}|^2}{2k_B \mathrm{\Theta}^{\e}} - \frac{3}{2} \rw) \bigg) M^{\e},
\enda
\end{align}
we obtain \eqref{Phi1def}. By an explicit computation, for $|\beta|=1$ and $|\gamma|=1$, we also have 
\begin{align*}
\bega
&\p^{\gamma}\p^{\beta}M^{\e}= \e\Phi_{\beta+\gamma}^1M^{\e} + \e^2 \Phi_{\beta}^1\Phi_{\gamma}^1M^{\e} \cr 
&-\e^2\bigg(\frac{\p^{\gamma}u^{\e}\cdot\p^{\beta}u^{\e}}{k_B \mathrm{\Theta}^{\e}} +\p^{\gamma}\ta^{\e}\frac{\p^{\beta}u^{\e}\cdot(v - \mathrm{U}^{\e})}{k_B \mathrm{\Theta}^{\e}} + \p^{\gamma}\ta^{\e}\p^{\beta}\ta^{\e}\frac{|v - \mathrm{U}^{\e}|^2}{2k_B \mathrm{\Theta}^{\e}} + \p^{\beta}\ta^{\e}\frac{\p^{\gamma}u^{\e}\cdot(v - \mathrm{U}^{\e})}{k_B \mathrm{\Theta}^{\e}} \bigg)  M^{\e}.
\enda
\end{align*}
By iteration, we have the result \eqref{Phi2def}. \\
Next, we prove \eqref{Phiinf0}. We split \eqref{Phiidef} for \( |\al| \geq 1 \) into the following two parts:
\begin{align}\label{PhiAB}
\bega
\e^i\Phi_{\al}^i M^{\e} = \underbrace{\sum_{|{\bf r}|=i}
\p_{\rho^{\e}}^{r_1} \p_{u^{\e}_1}^{r_2} \p_{u^{\e}_2}^{r_3} \p_{\ta^{\e}}^{r_4} M^{\e}}_{\mathcal{I}_{\Phi}}
\underbrace{\sum_{s=1}^{|\al|} \sum_{p_s(\al,{\bf r})} (\al!) \prod_{j=1}^{s} \frac{[\p^{\bl_j}(\rho^{\e},u^{\e}_1,u^{\e}_2,\ta^{\e})]^{\bk_j}}{\bk_j! (\bl_j!)^{|\bk_j|}}}_{\mathcal{II}_{\Phi}}.
\enda
\end{align}
We note that the derivative of \( \ta^{\e} \) acting on \( M^{\e} \) produces a second-order velocity polynomial \( |v - \mathrm{U}^{\e}|^2 \) with $\e$ since $\frac{\p \mathrm{\Theta}^{\e}}{\p \ta^{\e}} = \e\mathrm{\Theta}^{\e} $.
When every derivative acts on \( \mathrm{\Theta}^{\e} \), i.e., \( r_4 = i \) in \eqref{PhiAB}, we obtain
\begin{align}\label{MTAderi2}
\bega
\frac{\p^i M^{\e}}{\p (\mathrm{\Theta}^{\e})^i} = \sum_{0 \leq k \leq i} \binom{i}{k} \p^k \frac{\mathrm{P}^{\e}}{(2 \pi k_B \mathrm{\Theta}^{\e})^{\frac{3}{2}}} \p^{i-k} e^{-\frac{|v - \mathrm{U}^{\e}|^2}{2 k_B \mathrm{\Theta}^{\e}}}.
\enda
\end{align}
Because of \eqref{MTAderi} and \eqref{MTAderi2}, the maximum velocity growth coming from the \(i\)-th derivative of the \(\mathcal{I}_{\Phi}\)-part is as follows:
\begin{align}\label{PhiA}
\mathcal{I}_{\Phi} \les \e^i(1 + |v - \mathrm{U}^{\e}|^{2i}) M^{\e} \les \e^i(1 + |v|^{2i}) M^{\e},
\end{align}
where we used \eqref{condition}.
For the \(\mathcal{II}_{\Phi}\) part, by definition in \eqref{ps}, every vector $\bk_j$ and $\bl_j$ must be non-zero, and we have $\sum_j |\bk_j| = |{\bf r}| = i$. Thus, we get
\begin{align*}
\bega
\mathcal{II}_{\Phi} &\les  [\p^{\bl_1}(\rho,u,\ta)]^{\bk_1}\times \cdots \times [\p^{\bl_s}(\rho,u,\ta)]^{\bk_s}.
\enda
\end{align*}
From the condition $\sum_{i=1}^{s}|\bk_i|\bl_i=\al$ in the definition of $p_s(\al,{\bf r})$ in \eqref{ps}, we have
\begin{align}\label{PhiB}
\bega
\mathcal{II}_{\Phi} &\les \sum_{\substack{\al_1 + \cdots + \al_i = \al \\ \al_i > 0}} |\p^{\al_1}(\rho, u, \ta)| \times \cdots \times |\p^{\al_i}(\rho, u, \ta)|.
\enda
\end{align}
Combining \eqref{PhiA} and \eqref{PhiB}, we obtain \eqref{Phiinf0}. Then, \eqref{al-M} is a direct consequence of \eqref{Phiinf0}. \\
For the proof of \eqref{Phiscale} to \eqref{Rinfscale}, we take \(\|\langle v \rangle^n \cdot |M^{\e}|^{-n_0}\|_{L^2_x L^p_v}\) for \(n_0 < 1\) to \(\partial^{\al} M^{\e}\).
Since the velocity part has exponential decay, we have
\begin{align*}
\eps^i \|\la v \ra^n |\Phi_{\al}^i| |M^{\e}|^{1-n_0}\|_{L^2_x L^p_v} \les \eps^i \sum_{\substack{\al_1 + \cdots + \al_i = \al \\ \al_i > 0}} \big\| |\p^{\al_1}(\rho, u, \ta)| \times \cdots \times |\p^{\al_i}(\rho, u, \ta)| \big\|_{L^2_x}.
\end{align*}
Assuming, without loss of generality, that \( \al_1 \geq \al_k \) for all \( 1 \leq k \leq i \), we take \(L^2_x\) on the \(\p^{\al_1}\)-part and \(L^\infty_x\) on the remaining terms to get
\begin{align}\label{Phi2}
\bega
\e^i\|\la v\ra^{n}|\Phi_{\al}^i|&|M^{\e}|^{1-n_0}\|_{L^2_xL^p_v} \les  \e^i\sum_{\substack{\al_1+\cdots+\al_i=\al \\ \al_i>0}} \|\p^{\al_1}(\rho,u,\ta)\|_{L^2_x} \|\p^{\al_2}(\rho,u,\ta)\|_{L^\infty_x}\times \cdots \times \|\p^{\al_i}(\rho,u,\ta)\|_{L^\infty_x}.
\enda
\end{align}
We note that the top order derivative only appear in $\Phi_{\al}^1$. Thus, by using 
\begin{align*}
\bega
\sum_{0\leq|\al|\leq\mathrm{N}}\|\p^{\al}(\rho,u,\ta)\|_{L^2_x} \leq \mathcal{E}_M^{\frac{1}{2}}, \qquad \sum_{|\al|=\mathrm{N}+1}\|\p^{\al}(\rho,u,\ta)\|_{L^2_x} \leq \kappa^{-\frac{1}{2}}\mathcal{E}_{top}^{\frac{1}{2}}, \cr 
\sum_{0\leq|\al|\leq\mathrm{N}-2}\|\p^{\al}(\rho,u,\ta)\|_{L^\infty_x} \leq \mathcal{E}_M^{\frac{1}{2}}, \qquad \sum_{|\al|=\mathrm{N}-1}\|\p^{\al}(\rho,u,\ta)\|_{L^\infty_x} \leq \kappa^{-\frac{1}{2}}\mathcal{E}_{top}^{\frac{1}{2}},
\enda
\end{align*}
by Agmon's inequality \eqref{Agmon}, we have 
\begin{align*}
\bega
\|\la v\ra^n |\p^{\al}M^{\e}||M^{\e}|^{-n_0}\|_{L^2_xL^p_v} &\les \sum_{1\leq i\leq |\al|}\e^i\mathcal{E}_{M}^{\frac{i}{2}}, \quad &\mbox{for} \quad 1\leq|\al|\leq\mathrm{N}, \cr
\|\la v\ra^n |\p^{\al}M^{\e}||M^{\e}|^{-n_0}\|_{L^2_xL^p_v} &\les \e\kappa^{-\frac{1}{2}}\mathcal{E}_{top}^{\frac{1}{2}} + \sum_{2\leq i\leq |\al|}\e^i\mathcal{E}_{M}^{\frac{i}{2}}, \quad &\mbox{for} \quad |\al|=\mathrm{N}+1.
\enda
\end{align*}
By using \( \eps\mathcal{E}_{M}^{\frac{1}{2}} \leq 1/4 \) from \eqref{condition}, we have the result \eqref{Phiscale}. By the same way, taking $L^\infty_x$ instead of $L^2_x$ in \eqref{Phi2}, we can have \eqref{Phiinfscale}. Similarly, taking the sums \( \sum_{2 \leq i \leq |\al|} \) to \eqref{Phi2}, we obtain the results \eqref{Rscale}.
\end{proof}

\begin{lemma} 
Assume \eqref{condition}. Let \(n\) be the degree of the velocity polynomial of \(\mathbf{A}(x,v)\), where
\begin{align*}
\mathbf{A}(x,v) = \frac{1}{\mathrm{P}^{\e}},~\text{or}~\frac{v_i-\mathrm{U}^{\e}_i}{\sqrt{k_B\mathrm{P}^{\e} \mathrm{\Theta}^{\e}}},~\text{or}~\frac{1}{\sqrt{6 \mathrm{P}^{\e}}}\lw(\frac{|v-\mathrm{U}^{\e}|^2}{k_B\mathrm{\Theta}^{\e}}-3\rw),~\text{or}~\mathfrak{R}^{\e}_{ij},~\text{or}~\mathcal{Q}^{\e}_j.
\end{align*}
Then, for \( |\al| \geq 1 \), the following inequalities hold:
\begin{align}
&|\p^{\al}(\mathbf{A}(x,v)M^{\e})| \les \sum_{1\leq i\leq |\al|}\eps^i (1+|v|^{2i+n})\sum_{\substack{\al_1+\cdots+\al_i=\al \\ \al_i>0}}|\p^{\al_1}(\rho,u,\ta)|\times \cdots \times |\p^{\al_i}(\rho,u,\ta)| M^{\e} , \label{al-vM}  \\[10pt]
&|\p^{\al}|M^{\e}|^{-1/2}| \les \sum_{1\leq i\leq |\al|}\eps^i (1+|v|^{2i})\sum_{\substack{\al_1+\cdots+\al_i=\al \\ \al_i>0}}|\p^{\al_1}(\rho,u,\ta)|\times \cdots \times |\p^{\al_i}(\rho,u,\ta)| |M^{\e}|^{-1/2}. \label{al-1/M}
\end{align}
\end{lemma}

\begin{proof}
Using the same approach as in the proof of \eqref{al-M} in Lemma \ref{Mal}, we can derive \eqref{al-vM} and \eqref{al-1/M}. This can be achieved by applying the relation 
\(\partial |M^{\e}|^{-1/2} = -\frac{1}{2} |M^{\e}|^{-3/2} \partial M^{\e}\),
inductively.
\end{proof}

\hide
\begin{proof}[\textbf{Proof of Lemma \ref{lem:Fluid_dissipative}}]
We substitute $M^{\e}+\AC{\P}F^{\e}$ for $F$ in \eqref{BE} to obtain
\beq \notag
\eps\p_t (M^{\e}+\AC{\P}F^{\e})+v\cdot\nabla_x (M^{\e}+\AC{\P}F^{\e})=\frac{2}{\kappa\eps}\mathcal{N}(\AC{\P}F^{\e},M^{\e})+\frac 1 {\kappa\eps}\mathcal{N}(\AC{\P}F^{\e},\AC{\P}F^{\e}).
\eeq
Note that $\p_t M^{\e} = \e \Phi_{(1,0,0)}^1 M^{\e}$ is spanned by $(1, v, |v|^2)$, as shown in \eqref{Phi1def}, while $\mathcal{N}(\cdot, \cdot)$ is microscopic. Thus, taking $\AC{\P}$ on both sides gives
\[
\eps \p_t \AC{\P} F 
+ \AC{\P} \lw\{v \cdot \nabla_x \AC{\P} F + v \cdot \nabla_x M^{\e} \rw\}
= \frac{2}{\kappa \eps} \mathcal{N}(\AC{\P} F, M^{\e}) + \frac{1}{\kappa \eps} \mathcal{N}(\AC{\P} F, \AC{\P} F).
\]
By the definition of $\mathcal{L}$ in \eqref{Ldef}, we have $-2\mathcal{N}(\AC{\P} F, M^{\e}) = \mathcal{L}(\AC{\P} F)$. Applying $\kappa \e \mathcal{L}^{-1}$ to both sides, we get  
\begin{align*}
\AC{\P} F = -\kappa \eps \mathcal{L}^{-1} \lw\{\AC{\P} \big(v \cdot \nabla_x M^{\e}\big)\rw\} 
-\mathcal{L}^{-1} \lw\{\kappa \eps \big(\eps \p_t \AC{\P} F + \AC{\P} \big(v \cdot \nabla_x \AC{\P} F \big)\big) 
- \mathcal{N}(\AC{\P} F, \AC{\P} F)\rw\}.
\end{align*}
For the term $\AC{\P} \big(v \cdot \nabla_x M^{\e}\big)$, we use \eqref{pM} to compute
\[\bega
\sum_{i}\AC{\P}(v_i \p_i M^{\e})
&=\sum_{i,j}\frac{\p_i \mathrm{U}^{\e}_j}{k_B\mathrm{\Theta}^{\e}}\AC{\P}\lw\{ \lw((v_i-\mathrm{U}^{\e}_i)(v_j-\mathrm{U}^{\e}_j)-\frac 1 3 \delta_{ij}|v-\mathrm{U}^{\e}|^2\rw)M^{\e}\rw\} \\
&+\sum_i \frac{\p_i\mathrm{\Theta}^{\e}}{k_B|\mathrm{\Theta}^{\e}|^2}\AC{\P}\lw\{(v_i-\mathrm{U}^{\e}_i)\lw(\frac{|v-\mathrm{U}^{\e}|^2-5k_B\mathrm{\Theta}^{\e}}{2}\rw) M^{\e}\rw\}.
\enda
\]
With this computation, we obtain the result in \eqref{Gform}.
\end{proof}
\unhide

\begin{lemma}\label{L.1/P}
Assume \eqref{condition}. For \((z_1,z_2) = (1,\mathrm{P}^{\e})\), \((\mathrm{\Theta}^{\e},1)\), or \((\mathrm{\Theta}^{\e},\mathrm{P}^{\e})\), the following estimates hold:
\begin{align*}
&(1)~ \bigg\|\p^{\al}\bigg(\frac{\sqrt{z_1}}{z_2}\bigg)\bigg\|_{L^2_x} \les \e \kappa^{-\frac{(|\al|-\mathrm{N})_+}{2}}\mathcal{E}_{top}^{\frac{1}{2}}, \quad \text{for} \quad 1 \leq |\al| \leq \mathrm{N}+1, \cr 
&(2)~ \bigg\|\p^{\al}\bigg(\frac{\sqrt{z_1}}{z_2}\bigg)\bigg\|_{L^\infty_x} \les \e \kappa^{-\frac{(|\al|-(\mathrm{N}-2))_+}{2}}\mathcal{E}_{top}^{\frac{1}{2}}, \quad \text{for} \quad 1 \leq |\al| \leq \mathrm{N}-1. 
\end{align*}
\end{lemma}

\begin{proof}
Using Faà di Bruno's formula \eqref{Faadi} and following the same procedure as in \eqref{Phi2}, we can derive the desired results.
\hide
(1) By the Faà di Bruno's formula \eqref{Faadi}, we have 
\begin{align}\label{z/z}
\p^{\al}\bigg(\frac{\sqrt{z_1}}{z_2}\bigg) = \sum_{1\leq|{\bf r}|\leq |\al|} \frac{\p^{|{\bf r}|}}{\p z_1^{r_1}\p z_2^{r_2}}\bigg(\frac{\sqrt{z_1}}{z_2}\bigg) \sum_{s=1}^{|\al|} \sum_{p_s(\al,{\bf r})} (\al!) \prod_{j=1}^{s} \frac{[\p^{\bl_j}(z_1,z_2)]^{\bk_j}}{\bk_j! (\bl_j!)^{|\bk_j|}}.
\end{align}
Using \eqref{condition}, we have 
\begin{align}\label{z/z1}
\bigg|\frac{\p^{|{\bf r}|}}{\p z_1^{r_1}\p z_2^{r_2}}\bigg(\frac{z_1}{z_2}\bigg)\bigg| \leq C. 
\end{align}
Since $\p^{\bl_j}\mathrm{P}^{\e} = \e \p^{\bl_j}\rho$, and $\p^{\bl_j}\mathrm{\Theta}^{\e} = \e \p^{\bl_j}\ta$, we have 
\begin{align}\label{z/z2}
\bigg|\sum_{s=1}^{|\al|} \sum_{p_s(\al,{\bf r})} (\al!) \prod_{j=1}^{s} \frac{[\p^{\bl_j}(z_1,z_2)]^{\bk_j}}{\bk_j! (\bl_j!)^{|\bk_j|}}\bigg| \les \e^i \sum_{\substack{\al_1+\cdots+\al_i=\al \\ \al_i>0}} |\p^{\al_1}(\rho,\ta)|\times \cdots \times |\p^{\al_i}(\rho,\ta)| .
\end{align}
We take $L^2_x$ on \eqref{z/z} and combine the estimate \eqref{z/z1} and \eqref{z/z2}. Then, by the same way of \eqref{Phi2}, we obtain
\begin{align*} \bigg\|\p^{\al}\bigg(\frac{\sqrt{z_1}}{z_2}\bigg)\bigg\|_{L^2_x} \les \sum_{i=1}^{|\al|}\e^i\kappa^{-\frac{|\al|-1}{2}}\mathcal{E}_{top}^{\frac{i}{2}} \les \e \kappa^{-\frac{|\al|-1}{2}}\mathcal{E}_{top}^{\frac{1}{2}},
\end{align*}
for $1\leq|\al|\leq 5$. Agmon's inequality \eqref{Agmon} gives the result \eqref{1/Pinf}. 
\unhide
\end{proof}

\begin{lemma}\label{L.pta} For the derivatives $\p^{\al}$ defined in \eqref{mindex}, and for $\mathrm{Z}_1 \in \left\{ \mathrm{\Theta}^{\e}, 1/\mathrm{\Theta}^{\e}\right\}$, $z_1 = \ta^{\e}$ and $\mathrm{Z}_2 \in \left\{\mathrm{P}^{\e}, 1/\mathrm{P}^{\e} \right\}$, $z_2=\rho^{\e}$, and for $k=1,2$, we have the following estimates:
\begin{align}\label{pTa}
&|\p^{\al}\mathrm{Z}_k| \leq \sum_{1\leq i\leq |\al|} C\e^i \sum_{\substack{\al_1+\cdots+\al_i=\al \\ \al_i>0}}|\p^{\al_1}z_k|\times \cdots \times |\p^{\al_i}z_k| ,
\end{align}
and
\begin{align}\label{pTaleq}
\bega
\|\p^{\al}\mathrm{Z}\|_{L^2_x} &\leq C \e \mathcal{E}_{M}^{\frac{1}{2}}, \quad \mbox{for} \quad 1 \leq |\al| \leq \mathrm{N}, \cr \|\p^{\al}\mathrm{Z}\|_{L^\infty_x} &\leq C \e \mathcal{E}_{M}^{\frac{1}{2}}, \quad \mbox{for} \quad 1 \leq |\al| \leq \mathrm{N}-2,
\enda
\end{align}
for some positive constants $C>0$.
\end{lemma}
\begin{proof}
Applying Faà di Bruno's formula and $|\mathrm{\Theta}^{\e}|\leq 3/2$, we have \eqref{pTa}. 
\hide 
(Below is explicity computation:)
\begin{align*}
&\p_i \mathrm{\Theta}^{\e} = \p_i e^{\e\ta} = \e\p_i\ta e^{\e\ta} = \e\p_i\ta \mathrm{\Theta}^{\e} \cr 
&\p_j\p_i\mathrm{\Theta}^{\e} = \e\p_j\p_i\ta \mathrm{\Theta}^{\e} + \e^2\p_i\ta \p_j \ta \mathrm{\Theta}^{\e} \cr 
&\p_k\p_j\p_i\mathrm{\Theta}^{\e} = \e\p_k\p_j\p_i\ta \mathrm{\Theta}^{\e} + \e^2(\p_k \ta \p_j\p_i\ta +\p_k\p_i\ta \p_j \ta+\p_k\p_j\ta \p_i \ta) \mathrm{\Theta}^{\e} +\e^3 (\p_k\ta)(\p_j\ta)(\p_i\ta)\mathrm{\Theta}^{\e}
\end{align*}
\unhide
Without loss of generality, let us assume that $|\al_1|\geq |\al_2|\geq \cdots \geq |\al_i|$. Taking $L^2_x$ or $L^\infty_x$, we get 
\begin{align*}
\|\p^{\al}\mathrm{\Theta}^{\e}\|_{L^2_x} &\leq \sum_{1\leq i\leq |\al|} C\e^i \sum_{\substack{\al_1+\cdots+\al_i=\al \\ \al_i>0}}\|\p^{\al_1}\ta\|_{L^2_x} \|\p^{\al_2}\ta\|_{L^\infty_x} \times \cdots \times \|\p^{\al_i}\ta\|_{L^\infty_x}  \les \sum_{1\leq i\leq |\al|} C\e^i \mathcal{E}_{M}^{\frac{i}{2}} \leq C \e \mathcal{E}_{M}^{\frac{1}{2}},
\end{align*}
and
\begin{align*}
\|\p^{\al}\mathrm{\Theta}^{\e}\|_{L^\infty_x} &\leq \sum_{1\leq i\leq |\al|} C\e^i \sum_{\substack{\al_1+\cdots+\al_i=\al \\ \al_i>0}}\|\p^{\al_1}\ta\|_{L^\infty_x}\times \cdots \times \|\p^{\al_i}\ta\|_{L^\infty_x}  \les \sum_{1\leq i\leq |\al|} C\e^i \mathcal{E}_{M}^{\frac{i}{2}} \leq C \e \mathcal{E}_{M}^{\frac{1}{2}}.
\end{align*}

\end{proof}

\begin{lemma}\label{L.MFME}
Taylor's theorem provides the following expansions for \(M^{\e}\) and \(M^E\): \\
(1) Taylor expansion of \(M^E\), as defined in \eqref{MEdef}, around the local Maxwellian \(M^{\e}\):
\begin{align}\label{MFME}
\bega
M^E &= M^{\e} + \e(\bar{\rho}, \bar{u}, \bar{\ta}) \cdot \bigg(1,\frac{(v-\mathrm{U}^{\e})}{k_B\mathrm{\Theta}^{\e}},\frac{|v-\mathrm{U}^{\e}|^2-3k_B\mathrm{\Theta}^{\e}}{2k_B\mathrm{\Theta}^{\e}}\bigg) M^{\e} \\
&\quad + \e^2 (\mathrm{P}^{\vartheta}\bar{\rho}, \bar{u}, \mathrm{\Theta}^{\vartheta}\bar{\ta})^T \int_0^1 \left\{D^2_{(\mathrm{P}^{\vartheta}, \mathrm{U}^{\vartheta}, \mathrm{\Theta}^{\vartheta})}M(\vartheta)\right\} (1 - \vartheta) d\vartheta (\mathrm{P}^{\vartheta}\bar{\rho}, \bar{u}, \mathrm{\Theta}^{\vartheta}\bar{\ta}),
\enda
\end{align}
where
\begin{align*}
\mathrm{P}^{\vartheta} := e^{\e(\vartheta\rho^E+(1-\vartheta)\rho^{\e})}, \quad  
\mathrm{U}^{\vartheta}:= \vartheta\mathrm{U}^E +(1-\vartheta)\mathrm{U}^{\e}, \quad 
\mathrm{\Theta}^{\vartheta} := e^{\e(\vartheta\ta^E+(1-\vartheta)\ta^{\e})},
\end{align*}
and \(M(\vartheta) := M_{[\mathrm{P}^{\vartheta}, \mathrm{U}^{\vartheta}, \mathrm{\Theta}^{\vartheta}]}\). Here, \((\bar{\rho}, \bar{u}, \bar{\ta}) = (\rho^E-\rho^{\e}, u^E-u^{\e}, \ta^E-\ta^{\e})\). \\
(2) Taylor expansion of \(M^{\e}\) around the global Maxwellian \(\mu\): 
\begin{align}\label{MFtaylor}
\bega
M^{\e} &= \mu + \e\left(\rho^{\e} + u^{\e} \cdot \frac{v}{k_B} + \ta^{\e} \frac{|v|^2 - 3k_B}{2k_B}\right) \mu \\
&\quad + \e^2 (\mathrm{P}^{\vartheta}\rho^{\e}, u^{\e}, \mathrm{\Theta}^{\vartheta}\ta^{\e})^T \int_0^1 \left\{D^2_{(\mathrm{P}^{\vartheta}, \mathrm{U}^{\vartheta}, \mathrm{\Theta}^{\vartheta})}M(\vartheta)\right\} (1 - \vartheta) d\vartheta (\mathrm{P}^{\vartheta}\rho^{\e}, u^{\e}, \mathrm{\Theta}^{\vartheta}\ta^{\e})^T,
\enda
\end{align}
where
\begin{align*}
\mathrm{P}^{\vartheta} := e^{\e(\vartheta\rho^{\e})}, \quad  
\mathrm{U}^{\vartheta}:= \vartheta\mathrm{U}^{\e}, \quad 
\mathrm{\Theta}^{\vartheta} := e^{\e(\vartheta\ta^{\e})},
\end{align*}
and \(M(\vartheta) = M_{[\mathrm{P}^{\vartheta}, \mathrm{U}^{\vartheta}, \mathrm{\Theta}^{\vartheta}]}\).
 \\
In both cases, the following notation is used for \(D^2_{(\mathrm{P}^{\vartheta}, \mathrm{U}^{\vartheta}, \mathrm{\Theta}^{\vartheta})}M(\vartheta)\):
\begin{align*}
\bega
D^2_{(\mathrm{P}^{\vartheta}, \mathrm{U}^{\vartheta}, \mathrm{\Theta}^{\vartheta})}M(\vartheta) := \lw[ \begin{array}{ccc}
\frac{\p^2 M(\vartheta)}{\p \mathrm{P}^{\vartheta} \p \mathrm{P}^{\vartheta}} & \frac{\p^2 M(\vartheta)}{\p \mathrm{U}_i^{\vartheta} \p \mathrm{P}^{\vartheta}} & \frac{\p^2 M(\vartheta)}{\p \mathrm{\Theta}^{\vartheta} \p \mathrm{P}^{\vartheta}} \\[5pt]
\frac{\p^2 M(\vartheta)}{\p \mathrm{P}^{\vartheta} \p \mathrm{U}_i^{\vartheta}} & \frac{\p^2 M(\vartheta)}{\p \mathrm{U}_j^{\vartheta} \p \mathrm{U}_i^{\vartheta}} & \frac{\p^2 M(\vartheta)}{\p \mathrm{\Theta}^{\vartheta} \p \mathrm{U}_i^{\vartheta}} \\[5pt]
\frac{\p^2 M(\vartheta)}{\p \mathrm{P}^{\vartheta} \p \mathrm{\Theta}^{\vartheta}} & \frac{\p^2 M(\vartheta)}{\p \mathrm{U}_i^{\vartheta} \p \mathrm{\Theta}^{\vartheta}} & \frac{\p^2 M(\vartheta)}{\p \mathrm{\Theta}^{\vartheta} \p \mathrm{\Theta}^{\vartheta}}
\end{array} \rw].
\enda
\end{align*}
\end{lemma}

\begin{proof}
We only prove \eqref{MFME}. By Taylor's theorem, we have
\begin{align}\label{Taylor}
M(1)= M(0) + \frac{d M(\vartheta)}{d\vartheta}\bigg|_{\vartheta=0} + \int_0^1 \frac{d^2 M(\vartheta)}{d\vartheta^2}(1-\vartheta)d\vartheta.
\end{align}
Note that $M(1)=M^E$ and $M(0)=M^{\e}$. By the chain rule, we obtain
\begin{align}\label{M'}
\bega
\frac{dM(\vartheta)}{d\vartheta}\bigg|_{\vartheta=0}&= \nabla_{(\mathrm{P}^{\vartheta},\mathrm{U}^{\vartheta},\mathrm{\Theta}^{\vartheta})}M(\vartheta)\bigg|_{\vartheta=0}  \e(\mathrm{P}^{\e}(\rho^E-\rho^{\e}),u^E-u^{\e},\mathrm{\Theta}^{\e}(\ta^E-\ta^{\e}))^{T} \cr 
&=\e(\rho^E-\rho^{\e},u^E-u^{\e},\ta^E-\ta^{\e}) \cdot \lw(1,\frac{(v-\mathrm{U}^{\e})}{k_B\mathrm{\Theta}^{\e}}, \frac{|v-\mathrm{U}^{\e}|^2-3k_B\mathrm{\Theta}^{\e}}{2k_B\mathrm{\Theta}^{\e}}\rw)M^{\e}, 
\enda
\end{align}
where we used
\begin{align*}
\nabla_{(\mathrm{P},\mathrm{U},\mathrm{\Theta})} M&=\lw(\frac{1}{\mathrm{P}},\frac{(v-\mathrm{U})}{k_B\mathrm{\Theta}}, \frac{|v-\mathrm{U}|^2-3k_B\mathrm{\Theta}}{2k_B|\mathrm{\Theta}|^2}\rw)M.
\end{align*}
For the last term in \eqref{Taylor}, applying the chain rule twice yields
\begin{align}\label{M''}
\bega
\frac{d^2M(\vartheta)}{d\vartheta^2}
&=\e^2\Big(\mathrm{P}^{\vartheta}(\rho^E-\rho^{\e}),u^E-u^{\e},\mathrm{\Theta}^{\vartheta}(\ta^E-\ta^{\e})\Big)^T\left\{D^2_{(\mathrm{P}^{\vartheta},\mathrm{U}^{\vartheta},\mathrm{\Theta}^{\vartheta})}M(\vartheta)\right\} \cr 
&\quad \times \Big(\mathrm{P}^{\vartheta}(\rho^E-\rho^{\e}),u^E-u^{\e},\mathrm{\Theta}^{\vartheta}(\ta^E-\ta^{\e})\Big).
\enda
\end{align}
Combining \eqref{M'} and \eqref{M''} with \eqref{Taylor} yields the result.
\end{proof}

\hide
\begin{proof}[{\bf Proof of Lemma \ref{MFMEdif}}]
Using \eqref{MFME}, we can decompose as follows: 
\begin{align}\label{MMs}
\bega
\frac{1}{\e^2}&\int_{\Omega \times \R^3} \frac{|\p^{\al}(M^{\e}-M^E)|^2}{M^{\e}} dvdx = \frac{1}{\e^2}\int_{\Omega \times \R^3}|X_1^{\al}+X_2^{\al}+X_3^{\al}|^2|M^{\e}|^{-1} dvdx,
\enda
\end{align}
where
\begin{align}\label{X123}
\bega
X_1^{\al}&= 
\e\p^{\al}(\bar{\rho},\bar{u},\bar{\ta})\cdot \bigg(\frac{1}{\mathrm{P}^{\e}},\frac{(v-\mathrm{U}^{\e})}{k_B\mathrm{\Theta}^{\e}},\frac{|v-\mathrm{U}^{\e}|^2-3k_B\mathrm{\Theta}^{\e}}{2k_B|\mathrm{\Theta}^{\e}|^2}\bigg)M^{\e}, \cr  
X_2^{\al}&= \e\sum_{1\leq\beta\leq\al}\sum_{0\leq\gamma\leq\beta}\binom{\al}{\beta}\binom{\beta}{\gamma}\p^{\al-\beta}(\bar{\rho},\bar{u},\bar{\ta})\cdot\p^{\beta-\gamma}\bigg(\frac{1}{\mathrm{P}^{\e}},\frac{(v-\mathrm{U}^{\e})}{k_B\mathrm{\Theta}^{\e}},\frac{|v-\mathrm{U}^{\e}|^2-3k_B\mathrm{\Theta}^{\e}}{2k_B|\mathrm{\Theta}^{\e}|^2}\bigg)\p^{\gamma}M^{\e}, \cr 
X_3^{\al}&= \e^2 \p^{\al}\bigg( (\bar{\rho},\bar{u},\bar{\ta})^T \int_0^1 \left\{D^2_{(\mathrm{P}^{\vartheta},\mathrm{U}^{\vartheta},\mathrm{\Theta}^{\vartheta})}M(\vartheta)\right\} (1-\vartheta) d\vartheta (\bar{\rho},\bar{u},\bar{\ta}) \bigg) .
\enda
\end{align}
By using 
We note that $X_1$ is the main order term. The terms $X_2$ and $X_3$ is order $\e^2$ because $\p^{\beta-\gamma}$ or $\p^{\gamma}$ in $X_2$ makes one $\e$. For the main order term $X_1$, we have
\begin{align}\label{X1}
\bega
\frac{1}{\e^2}\int_{\Omega \times \R^3}X_1^2|M^{\e}|^{-1}dvdx = \int_{\Omega} \frac{1}{\mathrm{P}^{\e}}|\p^{\al}\bar{\rho}|^2 dx+\int_{\Omega} \frac{\mathrm{P}^{\e}}{k_B\mathrm{\Theta}^{\e}}|\p^{\al}\bar{u}|^2 dx + \int_{\Omega} \frac{3\mathrm{P}^{\e}}{2|\mathrm{\Theta}^{\e}|^2}|\p^{\al}\bar{\ta}|^2 dx.
\enda
\end{align}
For $s=2,3$, by using the $H^3_x$ boundedness in \eqref{BEC} and \eqref{C}, we can have the following estimates.
\begin{align}\label{Xesti}
\bega
\sum_{0\leq|\al|\leq s}\|X_1^{\al}\|_{L^2_x}&\les \e\|(\bar{\rho},\bar{u},\bar{\ta})\|_{H^s_x}(1+|v|^2)M^{\e}, \cr 
\sum_{0\leq|\al|\leq s}\|X_2^{\al}\|_{L^2_x} &\les \e^2\|(\bar{\rho},\bar{u},\bar{\ta})\|_{H^s_x}(1+|v|^{2+2|\al|})M^{\e}, \cr 
\sum_{0\leq|\al|\leq s}\|X_3^{\al}\|_{L^2_x} &\les \e^2 \|(\bar{\rho},\bar{u},\bar{\ta})\|_{H^s_x}^2 (1+|v|^{4+2|\al|}) e^{C\e|v|^2+C\e^2}M^{\e} .
\enda
\end{align}
Applying \eqref{Xesti} to \eqref{MMs} we have 
\begin{align*}
\bega
&\bigg|\frac{1}{\e^2}\int_{\Omega \times \R^3} \frac{|\p^{\al}(M^{\e}-M^E)|^2}{M^{\e}} dvdx - \int_{\Omega} \frac{1}{\mathrm{P}^{\e}}|\p^{\al}\bar{\rho}|^2 dx+\int_{\Omega} \frac{\mathrm{P}^{\e}}{k_B\mathrm{\Theta}^{\e}}|\p^{\al}\bar{u}|^2 dx + \int_{\Omega} \frac{3\mathrm{P}^{\e}}{2|\mathrm{\Theta}^{\e}|^2}|\p^{\al}\bar{\ta}|^2 dx \bigg| \cr 
&\les \e \|(\bar{\rho},\bar{u},\bar{\ta})\|_{H^s_x}^2 +\e^4 \|(\bar{\rho},\bar{u},\bar{\ta})\|_{H^s_x}^4 .
\enda
\end{align*}
Using $\|(\bar{\rho},\bar{u},\bar{\ta})\|_{H^s_x} \leq C$ for $s=2,3$ because of \eqref{BEC} and \eqref{C}, this guarantees Lemma \eqref{MFMEdif}. Now, we should prove \eqref{Xesti}. For brevity, we only consider $X_3^{\al}$ case. We explicitly distribute the derivative $\p^{\al}$ in $X_3^{\al}$.
\begin{align}\label{X3=}
\bega
X_3^{\al} &= \e^2 \sum_{0\leq\beta\leq\al}\sum_{0\leq\gamma\leq\beta}\binom{\al}{\beta}\binom{\beta}{\gamma} \p^{\al-\beta}(\bar{\rho},\bar{u},\bar{\ta})^T \bigg(\int_0^1 \p^{\gamma}\left\{D^2_{(\mathrm{P}^{\vartheta},\mathrm{U}^{\vartheta},\mathrm{\Theta}^{\vartheta})}M(\vartheta)\right\} (1-\vartheta) d\vartheta \bigg) \p^{\beta-\gamma} (\bar{\rho},\bar{u},\bar{\ta}).
\enda
\end{align}
Applying Faà di Bruno's formula \eqref{Faadi} to $\p^{\al}D^2_{(\mathrm{P}^{\vartheta},\mathrm{U}^{\vartheta},\mathrm{\Theta}^{\vartheta})}M(\vartheta)$, gives
\begin{align*}
\p^{\al}(D^2_{(\mathrm{P}^{\vartheta},\mathrm{U}^{\vartheta},\mathrm{\Theta}^{\vartheta})}M(\vartheta)) = \sum_{1\leq|\beta|\leq |\al|} \frac{\p^{|\beta|}D^2_{(\mathrm{P}^{\vartheta},\mathrm{U}^{\vartheta},\mathrm{\Theta}^{\vartheta})}M(\vartheta)}{\p( \mathrm{P}^{\vartheta})^{\beta_1}\p(\mathrm{U}_{1\vartheta})^{\beta_2}\p(\mathrm{U}_{2\vartheta})^{\beta_3} \p(\mathrm{\Theta}^{\vartheta})^{\beta_4}} \sum_{s=1}^{|\al|} \sum_{p_s(\al,\beta)} (\al!) \prod_{j=1}^{s} \frac{[\p^{\bl_j}(\mathrm{P}^{\vartheta},\mathrm{U}_{1\vartheta},\mathrm{U}_{2\vartheta},\mathrm{\Theta}^{\vartheta})]^{\bk_j}}{\bk_j! (\bl_j!)^{|\bk_j|}}.
\end{align*}
Since 
$(|\mathrm{P}^{\vartheta}-1|,|\mathrm{U}^{\vartheta}|,|\mathrm{\Theta}^{\vartheta}-1|)<\delta$, by the same way to \eqref{Phiinf0}, we can have 
\begin{align}\label{X3-a}
\bega
&|(D^2_{(\mathrm{P}^{\vartheta},\mathrm{U}^{\vartheta},\mathrm{\Theta}^{\vartheta})}M(\vartheta))| \les (1+|v|^4) M(\vartheta), \cr 
&|\p^{\al}(D^2_{(\mathrm{P}^{\vartheta},\mathrm{U}^{\vartheta},\mathrm{\Theta}^{\vartheta})}M(\vartheta))| \les \sum_{i=1}^{|\al|}\e^i  (1+|v|^{4+2i}) M(\vartheta) \cr 
&\times \sum_{\substack{\al_1+\cdots+\al_i=\al \\ \al_i>0}}\vartheta^{i} |\p^{\al_1}\big(\vartheta (\rho^E,u^E,\ta^E)+ (1-\vartheta)(\rho,u,\ta) \big)|\times \cdots \times |\p^{\al_i}\big(\vartheta (\rho^E,u^E,\ta^E)+ (1-\vartheta)(\rho,u,\ta) \big)| ,
\enda
\end{align}
where we used
\begin{align*}
\p^{\al}(\mathrm{P}^{\vartheta},\mathrm{U}^{\vartheta},\mathrm{\Theta}^{\vartheta}) =  \e \p^{\al}\big(\vartheta (\rho^E,u^E,\ta^E)+ (1-\vartheta)(\rho,u,\ta) \big), \quad \mbox{for} \quad |\al|\geq1.
\end{align*}
Comparing $M(\vartheta)$ and $M^{\e}$, we have 
\begin{align}\label{X3-b}
\bega
\frac{M(\vartheta)}{M^{\e}} &= \frac{\mathrm{P}^{\vartheta}}{(2\pi k_B\mathrm{\Theta}^{\vartheta})^{\frac{3}{2}}} \frac{(2\pi k_B\mathrm{\Theta}^{\e} )^{\frac{3}{2}}}{\mathrm{P}^{\e}} e^{-\frac{|v-\mathrm{U}^{\vartheta}|^2}{2k_B\mathrm{\Theta}^{\vartheta}}} e^{\frac{|v-\mathrm{U}^{\e} |^2}{2k_B\mathrm{\Theta}^{\e}}} 
\leq C e^{C\e|v|^2+C\e^2},
\enda
\end{align}
where we used $\mathrm{\Theta}^{\vartheta}-\mathrm{\Theta}^{\e} = \e\vartheta(\ta^E-\ta)$ and $\|(\rho,u,\ta)\|_{L^\infty_x} , \|(\rho^E,u^E,\ta^E)\|_{L^\infty_x} \leq C$. Applying \eqref{X3-a} and \eqref{X3-b} to \eqref{X3=} gives 
\begin{align*}
\bega
X_3^{\al} &\les \e^2 \sum_{0\leq\beta\leq\al}\sum_{0\leq\gamma\leq\beta}|\p^{\al-\beta}(\bar{\rho},\bar{u},\bar{\ta})||\p^{\beta-\gamma}(\bar{\rho},\bar{u},\bar{\ta})|  \sum_{i=0}^{|\gamma|}\e^i  (1+|v|^{4+2i}) e^{C\e|v|^2+C\e^2} M^{\e} \cr 
&\times \sum_{\substack{\gamma_1+\cdots+\gamma_i=\gamma \\ \gamma_i>0}} |\p^{\gamma_1}\big(\vartheta (\rho^E,u^E,\ta^E)+ (1-\vartheta)(\rho,u,\ta) \big)|\times \cdots \times |\p^{\gamma_i}\big(\vartheta (\rho^E,u^E,\ta^E)+ (1-\vartheta)(\rho,u,\ta) \big)|.
\enda
\end{align*}
Then using \eqref{BEC} and \eqref{C}, we can obtain the estimate for $X_3^{\al}$ in \eqref{Xesti}. 
\end{proof}
\unhide

\hide
\begin{lemma}\label{L.MFME}
For $M^{\e}$ and $M^E$ defined in \eqref{M-def} and \eqref{MEdef}, respectively, the following expansion holds:
\begin{align}\label{MFME}
\bega
M^E
&= M^{\e}+\e(\rho^E-\rho^{\e},u^E-u^{\e},\ta^E-\ta^{\e})\cdot\bigg(1,\frac{(v-\mathrm{U}^{\e})}{k_B\mathrm{\Theta}^{\e}},\frac{|v-\mathrm{U}^{\e}|^2-3k_B\mathrm{\Theta}^{\e}}{2k_B\mathrm{\Theta}^{\e}}\bigg) M^{\e} \cr 
&+\e^2\Big(\mathrm{P}^{\vartheta}(\rho^E-\rho^{\e}),u^E-u^{\e},\mathrm{\Theta}^{\vartheta}(\ta^E-\ta^{\e})\Big)^T\left\{D^2_{(\mathrm{P}^{\vartheta},\mathrm{U}^{\vartheta},\mathrm{\Theta}^{\vartheta})}M(\vartheta)\right\} \Big(\mathrm{P}^{\vartheta}(\rho^E-\rho^{\e}),u^E-u^{\e},\mathrm{\Theta}^{\vartheta}(\ta^E-\ta^{\e})\Big),
\enda
\end{align}
where

\end{lemma}
\unhide

\hide
\section{Relative entropy}\label{A.re;entp}

We define the relative entropy:
\begin{align}\label{relHdef}
\mathcal{H}(F|M)= \int_{\Omega\times \R^3} \bigg( F \ln \frac{F}{M} - F + M \bigg) dvdx
\end{align}

\begin{lemma}\label{L.H.FME} For relative entropy defined in \eqref{relHdef} and $M^{E}$ defined in \eqref{MEdef}, we have 
\begin{align}\label{F|ME}
\bega
\frac{1}{\e^2}\frac{d}{dt}& \mathcal{H}(F^{\e}|M^E) +\frac{1}{\e^4\kappa}\mathcal{D}(\mathcal{H})(F^{\e}) =\frac{1}{\e^2}\int_{\Omega} \bigg[\frac{\mathrm{P}^{\e}}{k_B\mathrm{\Theta}^{E}}(u^{\e}-u^{E})\cdot \nabla_x p^E \bigg]dx \cr 
&-\frac{1}{\e^2} \int_{\Omega} \bigg[\frac{1}{k_B\mathrm{\Theta}^{E}}\sum_{i,j}\p_iu^{E}_j\la\mathfrak{R}_{ij}^E,F^{\e}\ra + \frac{1}{k_B\mathrm{\Theta}^{E}}\sum_i\p_i\ta^{E}\la\mathcal{Q}_i^E,F^{\e}\ra\bigg] dx.
\enda
\end{align}
where $\hat A_{ij}^E$ and $\mathcal{Q}_i^E$ are defined by 
\beq\label{AB-defE}
\hat A_{ij}^E=(v_i-U_i^E)(v_j-U_j^E)-\frac{|v-U^E |^2}{3}\delta_{ij},\qquad \mathcal{Q}_i^E = (v_i-U_i^E)\frac{(|v-U^E|^2-5k_B\mathrm{\Theta}^E)}{2}.
\eeq
Here, $\mathcal{D}(\mathcal{H})(F^{\e})$ denotes  \begin{align}\label{DHdef}
\bega
\mathcal{D}(\mathcal{H})(F^{\e}) &:=\frac{1}{4}\int_{\Omega \times \R^6}\int_{\mathbb{S}^2}B(v-v_*,w)\big[ F^{\e}(v)F^{\e}(v_*)-F^{\e}(v')F^{\e}(v_*')\big]\ln \frac{F^{\e}(v)F^{\e}(v_*)}{F^{\e}(v')F^{\e}(v_*')} dwdv_*dvdx .
\enda
\end{align}
\end{lemma}
\begin{proof}
By the definition of relative entropy \eqref{relHdef}, we have 
\begin{align} \label{F|M}
\bega
\frac{d}{dt} \mathcal{H}(F^{\e}|M^{E}) -\frac{1}{\e^2\kappa}\int_{\Omega\times\R^3} \mathcal{N}(F^{\e},F^{\e})\ln \frac{F^{\e}}{M^{E}} dvdx &= \int_{\Omega\times\R^3} \bigg(-\frac{v}{\e}\cdot\nabla_x F^{\e}\bigg)  \ln \frac{F^{\e}}{M^{E}} dvdx \cr 
&+ \int_{\Omega\times\R^3}\p_tM^{E}\bigg(-\frac{F^{\e}}{M^{E}}+1\bigg) dvdx.
\enda
\end{align}
For the first term of the right-side of \eqref{F|M}, we use 
\begin{align*} 
\nabla_x F^{\e} \ln \frac{F^{\e}}{M^{E}} = \nabla_x(F^{\e}\ln F^{\e}- F^{\e}) - \nabla_xF^{\e} \ln M^{E},
\end{align*}
to have 
\begin{align*} 
\int_{\Omega\times\R^3} \bigg(-\frac{v}{\e}\cdot\nabla_x F^{\e}\bigg)  \ln \frac{F^{\e}}{M^{E}} dvdx &= -\int_{\Omega\times\R^3} \frac{v}{\e}\cdot \nabla_x(F^{\e}\ln F^{\e} -F^{\e}) dvdx + \frac{v}{\e}\cdot \nabla_x F^{\e} \ln M^{E} dvdx \cr 
&= - \int_{\Omega\times\R^3} F^{\e} \frac{v}{\e}\cdot \frac{\nabla_x M^{E}}{M^{E}} dvdx .
\end{align*}
In addition, using $M^E(v)M^E(v_*)=M^E(v')M^E(v_*')$ symmetry of the velocity variables gives
\begin{align*} 
\int_{\Omega\times\R^3} \mathcal{N}(F^{\e},F^{\e})\ln \frac{F^{\e}}{M^E}dvdx = -\mathcal{D}(\mathcal{H})(F^{\e}).
\end{align*}
Then, the equality \eqref{F|M} becomes
\begin{align}\label{F|M3}
\frac{d}{dt} \mathcal{H}(F^{\e}|M^{E}) +\frac{1}{\e^2\kappa}\mathcal{D}(\mathcal{H})(F^{\e}) = \int_{\Omega\times\R^3} \bigg[ \p_t M^{E} - \frac{F^{\e}}{M^{E}}\bigg(\p_tM^{E} + \frac{v}{\e}\cdot \nabla_x M^{E}\bigg) \bigg] dvdx.
\end{align}
By the same way to the proof of \eqref{Phitx}, for the $(\rho^E,u^E,\ta^E)$ solves the incompressible Euler equation \eqref{incompE}, we have 
\begin{align}\label{PhitxE2}
\bega
\left(\p_tM^{E} + \frac{v}{\e}\cdot\nabla_x M^{E}\right) &=(v-\mathrm{U}^{E}) \bigg(\frac{ - \e \nabla_x p^E}{k_B\mathrm{\Theta}^{E}}\bigg)M^E \cr 
&+ \frac{1}{k_B\mathrm{\Theta}^{E}}\sum_{i,j}\p_iu^{E}_j\mathfrak{R}_{ij}^EM^E + \frac{1}{k_B\mathrm{\Theta}^{E}}\sum_i\p_i\ta^{E}\mathcal{Q}_i^EM^E.
\enda
\end{align}
Applying \eqref{PhitxE2} to the last term of \eqref{F|M3} gives 
\begin{align}\label{rhs.3}
\bega
\int_{\Omega\times\R^3} \frac{F^{\e}}{M^E}\bigg(\p_tM^E &+ \frac{v}{\e}\cdot \nabla_x M^E\bigg) dvdx = \int_{\Omega} \bigg[(\mathrm{P}^{\e}\mathrm{U}^{\e}-\mathrm{P}^{\e}\mathrm{U}^{E}) \frac{ -\e\nabla_x p^E}{k_B\mathrm{\Theta}^{E}} \bigg]dx \cr 
&+ \int_{\Omega} \bigg[\frac{1}{k_B\mathrm{\Theta}^{E}}\sum_{i,j}\p_iu^{E}_j\la\mathfrak{R}_{ij}^E,F^{\e}\ra + \frac{1}{k_B\mathrm{\Theta}^{E}}\sum_i\p_i\ta^{E}\la\mathcal{Q}_i^E,F^{\e}\ra\bigg] dx.
\enda
\end{align}
We also note that the $\mathrm{P}^E$ satisfies the continuity equation from \eqref{incompE}:
\begin{align*}
\bega
\p_t \mathrm{P}^E + u^E \cdot \nabla_x \mathrm{P}^E = \mathrm{P}^E\Big(\p_t \rho^E + u^E \cdot \nabla_x \rho^E\Big) = 0.
\enda
\end{align*}
Thus, we get
\begin{align}\label{rhs.4}
\frac{d}{dt}\int_{\Omega\times\R^3}  M^E  dvdx = \frac{d}{dt}\int_{\Omega}  \mathrm{P}^E dx =0. 
\end{align}
Applying \eqref{rhs.3} and \eqref{rhs.4} to $\frac{1}{\e^2}\times\eqref{F|M3}$ gives the result. 
\end{proof}

\begin{lemma}\label{L.M|M} For the relative entropy defined in \eqref{relHdef} and for $M^{\e}$ and $M^{E}$ defined in \eqref{M-def} and \eqref{MEdef}, respectively, we have  
\begin{align*}
\begin{split}
\frac{1}{\e^2}\mathcal{H}(M^{\e}|M^{E}) 
&\geq \int_{\Omega} \frac{1}{2}\mathrm{P}^{\e}|\rho^{\e}-\rho^{E}|^2 dx + \int_{\Omega}\frac{\mathrm{P}^{\e}}{2k_B\mathrm{\Theta}^{E}} |u^{\e}-u^{E}|^2  dx + \int_{\Omega}\frac{3}{4}\mathrm{P}^{\e}|\ta^{\e}-\ta^E|^2 dx \cr 
&- \frac{\e}{6}\|\mathrm{P}^{\e}\|_{L^\infty_x} \bigg[\|\rho^E-\rho^{\tau_1}\|_{L^2_x}^2\|\rho^E-\rho^{\tau_1}\|_{L^\infty_x}\|e^{\tau_1(\rho^E-\rho^{\tau_1})}\|_{L^\infty_x} \bigg] \cr 
&- \frac{\e}{4}\|\mathrm{P}^{\e}\|_{L^\infty_x} \bigg[\|\ta^E-\ta^{\tau_2}\|_{L^2_x}^2\|\ta^E-\ta^{\tau_2}\|_{L^\infty_x}\|e^{\tau_2(\ta^E-\ta^{\tau_2})}\|_{L^\infty_x} \bigg],
\end{split}
\end{align*}
for some $\tau_1\in(0,\e)$ and $\tau_2\in(0,\e)$.
\end{lemma}
\begin{proof}
By an explicit computation, we can easily have 
\begin{align}\label{M|M=}
\begin{split}
\mathcal{H}(M^{\e}|M^{E})&= \int_{\Omega} \mathrm{P}^{\e}\bigg[ \e(\rho^{\e}-\rho^E) -\frac{3}{2}\e(\ta^{\e}-\ta^E) \bigg] dx + \int_{\Omega} (-\mathrm{P}^{\e}+\mathrm{P}^{E})dx \cr 
&\quad + \int_{\Omega} \bigg[\frac{3\mathrm{P}^{\e}}{2}\bigg(\frac{\mathrm{\Theta}^{\e}}{\mathrm{\Theta}^{E}}-1\bigg) + \frac{\mathrm{P}^{\e}|\mathrm{U}^{\e}-\mathrm{U}^{E}|^2}{2k_B\mathrm{\Theta}^{E}} \bigg] dx.
\end{split}
\end{align}
We claim that 
\begin{align}\label{rho-rho}
\begin{split}
&\bigg|\frac{1}{\e^2}\int_{\Omega} \mathrm{P}^{\e}\Big( \e(\rho^{\e}-\rho^E) + \frac{\mathrm{P}^{E}}{\mathrm{P}^{\e}}-1\Big) dx - \int_{\Omega} \mathrm{P}^{\e} \frac{|\rho^{\e}-\rho^{E}|^2}{2}dx \bigg|\cr 
&\leq \frac{\e}{6}\|\mathrm{P}^{\e}\|_{L^\infty_x} \bigg[\|\rho^E-\rho^{\tau_1}\|_{L^2_x}^2\|\rho^E-\rho^{\tau_1}\|_{L^\infty_x}\|e^{\tau_1(\rho^E-\rho^{\tau_1})}\|_{L^\infty_x} \bigg] ,
\end{split}
\end{align}
for some $\tau_1\in(0,\e)$. We define 
\begin{align*}
q(\e) := \big(\e(\rho^{\e}-\rho^E)\big)+ e^{\e(\rho^E-\rho^{\e})}-1.
\end{align*}
By the Taylor's theorem, we have 
\begin{align*}
\frac{1}{\e^2}q(\e) - \frac{1}{2}|\rho^{\e}-\rho^{E}|^2 = \frac{1}{\e^2} \frac{\e^3}{3!}(\rho^E-\rho^{\tau_1})^3e^{\tau_1(\rho^E-\rho^{\tau_1})},
\end{align*}
for some $\tau_1\in(0,\e)$. 
Taking $\int_{\Omega} \mathrm{P}^{\e}(\cdot)dx$ to both sides, we can prove the claim \eqref{rho-rho}. By the same way, we can have 
\begin{align}\label{ta-ta}
\begin{split}
&\bigg|\frac{1}{\e^2}\int_{\Omega} \frac{3}{2}\mathrm{P}^{\e}\bigg(-\e(\ta^{\e}-\ta^E) +\bigg(\frac{\mathrm{\Theta}^{\e}}{\mathrm{\Theta}^{E}}-1\bigg) \bigg)  dx - \int_{\Omega}\frac{3}{4}\mathrm{P}^{\e}|\ta^{\e}-\ta^E|^2 dx\bigg| \cr 
&\leq \frac{\e}{4}\|\mathrm{P}^{\e}\|_{L^\infty_x} \bigg[\|\ta^E-\ta^{\tau_2}\|_{L^2_x}^2\|\ta^E-\ta^{\tau_2}\|_{L^\infty_x}\|e^{\tau_2(\ta^E-\ta^{\tau_2})}\|_{L^\infty_x} \bigg],
\end{split}
\end{align}
for some $\tau_2\in(0,\e)$.
Combining \eqref{rho-rho} and \eqref{ta-ta} to \eqref{M|M=} gives the result. 
\end{proof}

For notational simplicity, we define
\begin{align*}
\begin{split}
\mathbf{d}_{(\e,E)}(t) &:= \|\rho^{\e}(t)-\rho^{E}(t)\|_{L^2_x} +\|u^{\e}(t)-u^{E}(t)\|_{L^2_x} +\|\ta^{\e}(t)-\ta^{E}(t)\|_{L^2_x}  \cr 
&+ \|\rho^{\e}(t)-\rho^{E}(t)\|_{L^\infty_x} +\|u^{\e}(t)-u^{E}(t)\|_{L^\infty_x} +\|\ta^{\e}(t)-\ta^{E}(t)\|_{L^\infty_x}.
\end{split}
\end{align*}

\begin{lemma}\label{L.relH} We have 
\begin{align*}
\begin{split}
&\int_{\Omega} \frac{\mathrm{P}^{\e}}{2}\bigg[|\rho^{\e}-\rho^{E}|^2 +\frac{|u^{\e}-u^{E}|^2}{k_B\mathrm{\Theta}^{E}} +\frac{3}{2}|\ta^{\e}-\ta^E|^2\bigg]dx  +\frac{1}{\e^4\kappa}\int_0^t\mathcal{D}(\mathcal{H})(F^{\e})(s)ds \cr 
&\leq \frac{1}{\e^2}\mathcal{H}(F^{\e}|M^E)(0) \cr 
&+ C\int_0^t \bigg\|\frac{\mathrm{P}^{\e}(s)}{k_B\mathrm{\Theta}^{E}(s)}\bigg\|_{L^\infty_x}\Big(\|\nabla_xu^E(s)\|_{L^\infty_x}+\|\nabla_x\ta^E(s)\|_{L^\infty_x} \Big) \bigg(\|u^E(s)-u^{\e}(s)\|_{L^2_x}^2 + \bigg\|\frac{\mathrm{\Theta}^{\e}(s)-\mathrm{\Theta}^{E}(s)}{\e}\bigg\|_{L^2_x}^2 \bigg) ds \cr
& + \int_0^t 
\bigg[\bigg\|\frac{1}{k_B\mathrm{\Theta}^{E}(s)}\bigg\|_{L^\infty_x} \kappa^{\frac{1}{2}}\mathcal{D}_G^{\frac{1}{2}}(s) \bigg(\|\nabla_xu^E(s)\|_{L^2_x} + \|\nabla_x\ta^E(s)\|_{L^2_x} + \e\mathbf{d}_{(\e,E)}(s)\|\nabla_x\ta^E(s)\|_{L^\infty_x}\bigg) \bigg] ds \cr
&+ \int_0^t \bigg[\bigg\|\frac{\mathrm{P}^{\e}(s)}{k_B\mathrm{\Theta}^{E}(s)}\bigg\|_{L^\infty_x} \big\|\div(u^{\e})\big\|_{L^{2+\frac{1}{\delta}}_x}\|u^E(s)\|_{H^1_x}^2 + \bigg\|\nabla_x\frac{\mathrm{P}^{\e}(s)}{k_B\mathrm{\Theta}^{E}(s)}\bigg\|_{L^\infty_x}\mathbf{d}_{(\e,E)}(s)\|u^E(s)\|_{H^1_x}^2 \bigg] ds \cr
&+ \e\|\mathrm{P}^{\e}(t)\|_{L^\infty_x} \bigg[\frac{1}{6}\mathbf{d}_{(\tau_1,E)}^3e^{\tau_1\mathbf{d}_{(\tau_1,E)}} + \frac{1}{4}\mathbf{d}_{(\tau_2,E)}^3e^{\tau_2\mathbf{d}_{(\tau_2,E)}} \bigg] + \e \int_0^t \bigg\|\frac{\mathrm{P}^{\e}(s)}{2k_B\mathrm{\Theta}^{E}(s)}\bigg\|_{L^\infty_x}\|\nabla_x\ta^E(s)\|_{L^\infty_x} \mathbf{d}_{(\e,E)}^3(s)  ds.
\end{split}
\end{align*}
for some $\tau_1\in(0,\e)$ and $\tau_2\in(0,\e)$ and positive constant $C>0$. Here, the relative entropy $H(\cdot|\cdot)$ and relative entropy dissipation $\mathcal{D}(\mathcal{H})(F^{\e})$ defined in  \eqref{relHdef} and \eqref{DHdef}.
\end{lemma}
\begin{proof}
We take time integral to \eqref{F|ME} in Lemma \ref{L.H.FME} to have
\begin{align}\label{intH}
\bega
\frac{1}{\e^2}\mathcal{H}(F^{\e}|M^E)(t) +\frac{1}{\e^4\kappa}\int_0^t\mathcal{D}(\mathcal{H})(F^{\e})(s)ds = \frac{1}{\e^2}\mathcal{H}(F^{\e}|M^E)(0) +  \int_0^t \Big(H_1 + H_2 + H_3\Big) ds,
\enda
\end{align}
where
\begin{align*}
\bega
H_1 &:= \int_{\Omega} \bigg[\frac{\mathrm{P}^{\e}}{k_B\mathrm{\Theta}^{E}}(u^{\e}-u^{E})\cdot \nabla_x p^E \bigg]dx, \cr 
H_2&:=-\frac{1}{\e^2} \int_{\Omega} \frac{1}{k_B\mathrm{\Theta}^{E}}\sum_{i,j}\p_iu^{E}_j\la\mathfrak{R}_{ij}^E,F^{\e}\ra dx, \qquad H_3:=-\frac{1}{\e^2} \int_{\Omega}  \frac{1}{k_B\mathrm{\Theta}^{E}}\sum_i\p_i\ta^{E}\la\mathcal{Q}_i^E,F^{\e}\ra dx.
\enda
\end{align*}
We first estimate $H_2$ and $H_3$. For the estimate of $H_2$, by using 
\begin{align*}
\bega
\mathfrak{R}_{ij}^{\e}-\mathfrak{R}_{ij}^E &= (v_i-\e u^{\e}_i)(\e u_j^E-\e u^{\e}_j)+(v_j-\e u^{\e}_j)(\e u_i^E-\e u^{\e}_i) - \e^2(u_i^E-u^{\e}_i)(u_j^E-u^{\e}_j) \cr 
& -\frac{2(v-\e u^{\e})\cdot(\e u^E-\e u^{\e})-\e^2|u^E-u^{\e}|^2}{3}\delta_{ij},
\enda
\end{align*}
we have
\begin{align*}
\bega
&H_2= \frac{1}{\e^2} \int_{\Omega} \bigg[\frac{1}{k_B\mathrm{\Theta}^{E}}\sum_{i,j}\p_iu^{E}_j\la(\mathfrak{R}_{ij}^{\e}-\mathfrak{R}_{ij}^E),F^{\e}\ra\bigg] dx - \frac{1}{\e^2} \int_{\Omega} \bigg[\frac{1}{k_B\mathrm{\Theta}^{E}}\sum_{i,j}\p_iu^{E}_j\la\mathfrak{R}_{ij}^{\e},F^{\e}\ra\bigg] dx \cr 
&= -\int_{\Omega} \bigg[\frac{\mathrm{P}^{\e}}{k_B\mathrm{\Theta}^{E}}\sum_{i,j}\p_iu^{E}_j\Big((u_i^E-u^{\e}_i)(u_j^E-u^{\e}_j)-\frac{|u^E-u^{\e}|^2}{3}\delta_{ij}\Big)\bigg] dx - \frac{1}{\e^2} \int_{\Omega} \bigg[\frac{1}{k_B\mathrm{\Theta}^{E}}\sum_{i,j}\p_iu^{E}_j\mathbf{r}_{ij}^{\e}\bigg] dx.
\enda
\end{align*}
Applying $\eqref{ABGscale}_1$, we have
\begin{align}\label{H2est}
\bega
H_2&\leq \bigg\|\frac{\mathrm{P}^{\e}}{k_B\mathrm{\Theta}^{E}}\bigg\|_{L^\infty_x}\|\nabla_xu^E\|_{L^\infty_x} \|u^E-u^{\e}\|_{L^2_x}^2 + \frac{1}{\e^2}\bigg\|\frac{1}{k_B\mathrm{\Theta}^{E}}\bigg\|_{L^\infty_x}\|\nabla_xu^E\|_{L^2_x}\|\mathbf{r}_{ij}^{\e}\|_{L^2_x} \cr 
&\leq \bigg\|\frac{\mathrm{P}^{\e}}{k_B\mathrm{\Theta}^{E}}\bigg\|_{L^\infty_x}\|\nabla_xu^E\|_{L^\infty_x} \|u^E-u^{\e}\|_{L^2_x}^2 + \kappa^{\frac{1}{2}}\bigg\|\frac{1}{k_B\mathrm{\Theta}^{E}}\bigg\|_{L^\infty_x}\|\nabla_xu^E\|_{L^2_x}\mathcal{D}_G^{\frac{1}{2}}. 
\enda
\end{align}
Similarly, we use  
\begin{align*}
\bega
\mathcal{Q}_i^E -\mathcal{Q}_i^{\e} &= \frac{1}{2}(v_i-\e u^{\e}_i)\Big[2(v-\e u^{\e})\cdot(\e u^{\e}-\e u^E)+|\e u^{\e}-\e u^E|^2+5k_B(\mathrm{\Theta}^{\e}-\mathrm{\Theta}^{E})\Big] \cr 
&+ \frac{1}{2}(\e u^{\e}_i-\e u_i^E)\Big[|v-\e u^{\e}|^2+2(v-\e u^{\e})\cdot(\e u^{\e}-\e u^E)+|\e u^{\e}-\e u^E|^2\Big] \cr 
&+ \frac{1}{2}(\e u^{\e}_i-\e u_i^E)\Big[-5k_B\mathrm{\Theta}^{\e}+5k_B(\mathrm{\Theta}^{\e}-\mathrm{\Theta}^{E})\Big],
\enda
\end{align*}
to estimate $H_3$ as folllows:
\begin{align*}
\bega
H_3&=\frac{1}{\e^2} \int_{\Omega} \bigg[\frac{1}{k_B\mathrm{\Theta}^{E}}\sum_i\p_i\ta^{E}\la(\mathcal{Q}_i^{\e}-\mathcal{Q}_i^E),F^{\e}\ra \bigg] dx - \frac{1}{\e^2} \int_{\Omega} \bigg[\frac{1}{k_B\mathrm{\Theta}^{E}}\sum_i\p_i\ta^{E}\la \mathcal{Q}_i^{\e},F^{\e}\ra \bigg] dx \cr 
&=\frac{1}{\e^2} \int_{\Omega} \bigg[\frac{\mathrm{P}^{\e}}{k_B\mathrm{\Theta}^{E}}\sum_i\p_i\ta^{E}\bigg(\e(u^{\e}_i-u^{E}_i)\frac{|\e u^{\e}-\e u^E|^2+5k_B(\mathrm{\Theta}^{\e}-\mathrm{\Theta}^{E})}{2}\bigg) \bigg] dx \cr 
&\quad - \frac{1}{\e^2} \int_{\Omega} \bigg[\frac{1}{k_B\mathrm{\Theta}^{E}}\sum_i\p_i\ta^{E} \e(u^{\e}_j-u^{E}_j)\mathbf{r}_{ij}^{\e} \bigg] dx - \frac{1}{\e^2} \int_{\Omega} \bigg[\frac{1}{k_B\mathrm{\Theta}^{E}}\sum_i\p_i\ta^{E}\la \mathcal{Q}_i^{\e},F^{\e}\ra \bigg] dx.
\enda
\end{align*}
Applying $\eqref{ABGscale}_1$, we have
\begin{align}\label{H3est}
\bega
H_3&\leq \bigg\|\frac{\mathrm{P}^{\e}}{2k_B\mathrm{\Theta}^{E}}\bigg\|_{L^\infty_x}\|\nabla_x\ta^E\|_{L^\infty_x}\bigg[ \|u^E-u^{\e}\|_{L^2_x}\frac{5k_B}{2}\bigg\|\frac{\mathrm{\Theta}^{\e}-\mathrm{\Theta}^{E}}{\e}\bigg\|_{L^2_x} + \e\|u^{\e}-u^{E}\|_{L^\infty_x}\|u^{\e}-u^{E}\|_{L^2_x}^2 \bigg] \cr 
& + \bigg\|\frac{1}{k_B\mathrm{\Theta}^{E}}\bigg\|_{L^\infty_x} \bigg(\|\nabla_x\ta^E\|_{L^\infty_x} \e\|u^{\e}-u^{E}\|_{L^2_x}\kappa^{\frac{1}{2}}\mathcal{D}_G^{\frac{1}{2}} + \|\nabla_x\ta^E\|_{L^2_x}\kappa^{\frac{1}{2}}\mathcal{D}_G^{\frac{1}{2}}\bigg).  
\enda
\end{align}
For the estimate of $H_1$, applying integration by parts gives 
\begin{align*}
\bega
H_1 &= -\int_{\Omega} \bigg[\frac{\mathrm{P}^{\e}}{k_B\mathrm{\Theta}^{E}}\Big(\nabla_x \cdot (u^{\e}-u^{E})\Big) p^E \bigg]dx -\int_{\Omega} \bigg[(u^{\e}-u^{E})\cdot\nabla_x\frac{\mathrm{P}^{\e}}{k_B\mathrm{\Theta}^{E}} p^E \bigg]dx \cr 
&\leq \bigg\|\frac{\mathrm{P}^{\e}}{k_B\mathrm{\Theta}^{E}}\bigg\|_{L^\infty_x} \big\|(\div(u^{\e})-\div(u^{E}))\big\|_{L^{2+\frac{1}{\delta}}_x}\|p^E\|_{L^{\frac{2\delta+1}{\delta+1}}_x} + \bigg\|\nabla_x\frac{\mathrm{P}^{\e}}{k_B\mathrm{\Theta}^{E}}\bigg\|_{L^\infty_x}\|u^{\e}-u^{E}\|_{L^2_x}\|p^E\|_{L^2_x},
\enda
\end{align*}
where we used H\"{o}lder inequality for any $0<\delta<1$.
For the estimate of $p^E$, taking fourier transform to $\Delta_x p^E$ gives 
\begin{align*}
\bega
-\Delta_x p^E = \sum_{i,j}\p_i\p_j(u^E_iu^E_j), \qquad \widehat{p^E} = -\frac{\xi_i\xi_j}{|\xi|^2} \widehat{u^E_iu^E_j}.
\enda
\end{align*}
So that, we have $\|p^E\|_{L^p_x} \leq C\|u^E_iu^E_j\|_{L^p_x}$ for $1<p<\infty$, hence,
\begin{align*}
\bega
\|p^E\|_{L^{\frac{2\delta+1}{\delta+1}}_x} \leq C\|u^E_iu^E_j\|_{L^{\frac{2\delta+1}{\delta+1}}_x} \leq C \|u^E\|_{H^1_x},
\enda
\end{align*}
where we used Gagliardo-Nirenberg interpolation inequality in \eqref{Ga-Ni}. This gives
\begin{align}\label{H1est}
\bega
H_1 &\leq \bigg\|\frac{\mathrm{P}^{\e}}{k_B\mathrm{\Theta}^{E}}\bigg\|_{L^\infty_x} \big\|\div(u^{\e})\big\|_{L^{2+\frac{1}{\delta}}_x}\|u^E\|_{H^1_x}^2 + \bigg\|\nabla_x\frac{\mathrm{P}^{\e}}{k_B\mathrm{\Theta}^{E}}\bigg\|_{L^\infty_x}\|u^{\e}-u^{E}\|_{L^2_x}\|u^E\|_{H^1_x}^2.
\enda
\end{align}
where we used $\div(u^E)=0$. 
For the left-side of \eqref{intH}, we use
\begin{align*}
\bega
\frac{1}{\e^2}\mathcal{H}(F^{\e}|M^E) = \frac{1}{\e^2}\mathcal{H}(F^{\e}|M^{\e}) + \frac{1}{\e^2}\mathcal{H}(M^{\e}|M^E) \geq \frac{1}{\e^2}\mathcal{H}(M^{\e}|M^E) ,
\enda
\end{align*}
and Lemma \ref{L.M|M}. And for the right-side of \eqref{intH}, combining the estimates \eqref{H2est} \eqref{H3est} and \eqref{H1est} yields the result. 
\end{proof}
\unhide

\hide
\section{Proof of the Strichartz estimate for 3D}\label{A.Stri}

\subsubsection{Estimate of \cite{JMR}}
The following limit case of the Strichartz estimate is not true: 
\[
\| u \|_{L^2_t L^\infty_x} \lesssim _t \| \Box u \|_{L^1_t L^2_x}
\]
This is only true if the Fourier transform of $u $ is compactly supported in $\xi$. \cite{JMR} provides the sharp growth estimate \eqref{est:JMR} when this support $|\xi| \leq \lambda$ grows as $\lambda \rightarrow \infty$:

\begin{proposition}[Proposition 6.3 of \cite{JMR}]
Consider $f \in C^\infty (\R; \mathcal S(\R^3)), v_0 \in \mathcal S(\R^3)$ and $v_1 \in \mathcal S (\R^3)$. Let $v \in L^\infty_{loc} (\R_+; L^2(\R^3))$ solves 
\begin{equation}
\Box v = f, \ \ v|_{t=0}= v_0 , \ \ \p_t v |_{t=0} = v_1 .
\end{equation}
Define 
\begin{equation}
    S^\lambda = \varphi (\lambda^{-1} D_x),
\end{equation}
where $\varphi \in C^\infty_0 (\R^3)$ real, $\varphi \equiv 1$ when $|\xi| \leq 1$ and $\text{spt} (\varphi) \in \{|\xi| \leq  2\}$. $S^\lambda$ is a pseudo-differential operator such that $S^\lambda v = f$ means that 
\begin{equation}
    \varphi (\frac{\xi}{\lambda}) \hat v (\xi) = \hat{f}(\xi),
\end{equation}
or
\begin{equation}
    S^\lambda v(x) = \frac{1}{(2 \pi)^n} \int_{\R^n} e^{i x \cdot \xi}
    \varphi (\frac{\xi}{\lambda}) \hat v (\xi) d \xi. 
\end{equation}

Then there exists $C(\varphi)$ does not depend of $\tau $ and $v$:
    \begin{equation}\label{est:JMR}
        \sqrt{\int^\tau_0
        \| S^\lambda (v) (t) \|_{L^\infty_x}^2 dt 
        } \leq C(\varphi) \sqrt{\ln (1+ 2 \lambda \tau )}
        \left(\| v_0 \|_{\dot{H}^1_x} + \| v_1 \|_{L^2_x} + \int^\tau_0 \| f(t) \|_{L^2_x}dt\right)
    \end{equation}
\end{proposition}

\begin{lemma}[Lemma 8.1 of \cite{JMR}]

  Consider $v \in C^\infty (\R ; S (\R^3))$ such that 
\[
\Box v =g , \ \  v|_{t=0} = 0 , \ \  \p_t v |_{t=0 }=0.
\]

Suppose the support of the spatial Fourier transform $\hat g (t, \xi)$ is contained in $\{|\xi| \leq \lambda \}$. Then 
\begin{align}
    \| v(t) \|_{L^2_x}^2 \lesssim  \ln (1+ \lambda t) \int^t_0 \| g(s) \|_{L^1_x}^2 ds \label{est:v}\\
     \| \p_t v(t) \|_{L^2_x}^2 \lesssim  \ln (1+ \lambda t) \int^t_0 \| g(s) \|_{L^1_x}^2 ds\label{est:v_t}
\end{align}
\end{lemma}
\begin{proof}

Step 1. Fourier transform and the inverse transform:
\begin{equation}
\begin{split}
    \hat g (\xi) &= \int_{\R^d} e^{- i x \cdot \xi} g( x) dx\\
    g(x)  & = \frac{1}{(2 \pi)^d} \int_{\R^d} e^{i x\cdot \xi} \hat g (\xi) d \xi 
\end{split}
\end{equation}

We first prove the lemma for $\lambda =1$. Assume $\hat g(t, \xi) $ supported in $|\xi| \leq 1$.

Solution of the linear wave has the formula: 
\begin{equation}
    \begin{split}
        \hat v(t, \xi) = \int^t_0 \sin ( (t-s) |\xi|) \hat g (s, \xi) \frac{ds}{|\xi|}
        = \frac{1}{(2 \pi )^d}\int^t_0  \int_{\R^d} e^{- i x \cdot \xi}\sin ( (t-s) |\xi|)   g(s,  x) dx   \frac{ds}{|\xi|}
        ,\\
        \frac{\p_t \hat v (t, \xi)}{|\xi|} = \int^t_0 \cos ((t-s) |\xi|) \hat g (s, \xi) \frac{ds}{|\xi|}.
    \end{split}
\end{equation}

Thus 
\begin{align*}
    &\| v(t) \|_{L^2_x}^2 \\
    & = \| \hat v (t) \|_{L^2_\xi}^2 \\
    & = \int_{\R^d} \hat{v} (t, \xi) \overline{\hat{v}} (t, \xi)d \xi\\ 
    & 
        =\int_{\R^d}   \int^t_0  \int_{\R^d} e^{- i x \cdot \xi}\sin ( (t-s) |\xi|)   g(s,  x) dx \frac{ds}{|\xi|}
         \int^t_0  \int_{\R^d} e^{ i y \cdot \xi }\sin ( (t-\tau ) |\xi|)   \overline{g(\tau ,  y) }dy\frac{d\tau }{|\xi|}
    d \xi \\
    & = \int^t_0 \int^t_0  \iint_{\R^d \times \R^d}
    K_- (\tau, s, x- y )
    g(s,  x)\overline{g(\tau ,  y) }
    d x d y ds d \tau \\ 
    & \leq   \int_s     \int_x 
   |g(s,  x)|  \int_\tau \sqrt{ \sup_z | K_- (\tau, s,z )|}  \sqrt{ \sup_z | K_- (\tau, s,z )|} \int_y | \overline{g(\tau ,  y) }|
     \\ 
   & \leq  \sqrt{  \int_s    \left(  \int_x 
   |g(s,  x)|\right)^2}
   \sqrt{  \int_s    \left(
   \int_\tau \sqrt{ \sup_z | K_- (\tau, s,z )|}  \sqrt{ \sup_z | K_- (\tau, s,z )|} \int_y | \overline{g(\tau ,  y) }|
   \right)^2}
\end{align*}
Similarly,
\begin{equation}
    \| \frac{\p_t v (t)}{|\xi|}\|_{L^2}^2
    = \int^t_0 \int^t_0 \iint_{\R^d \times \R^d} K_+ (\tau, s, x-y)
    g(s,x) \overline{g(\tau, x)} dx dy ds d \tau 
\end{equation}
where 
\begin{align}
    K_- (\tau, s, x-y) : = C \int_{\xi \in \R^d} \mathbf{1}_{|\xi| \leq 1 }
    \sin ((t-s) |\xi|) \sin ((t-\tau) |\xi|) e^{- i (x-y) \cdot \xi} \frac{1}{|\xi|^2}d \xi  \\
      K_+ (\tau, s, x-y) : = C \int_{\xi \in \R^d} \mathbf{1}_{|\xi| \leq 1 } 
    \cos ((t-s) |\xi|) \cos ((t-\tau) |\xi|) e^{- i (x-y) \cdot \xi} \frac{1}{|\xi|^2}d \xi 
\end{align}

Now using the Cauchy-Schwarz inequality to the term in the last square root: 
\begin{align*}
    & \int_s    \left(
   \int_\tau \sqrt{ \sup_z | K_- (\tau, s,z )|}  \sqrt{ \sup_z | K_- (\tau, s,z )|} \int_y | \overline{g(\tau ,  y) }|
   \right)^2\\
   & \leq  \int_s   
  \int_\tau   \sup_z | K_- (\tau, s,z )| 
  \int_\tau \sup_z | K_- (\tau, s,z )|
 \left( \int_y | \overline{g(\tau ,  y) }|\right)^2 \\
  & \leq    
 \sup_s \int_\tau  \sup_z | K_- (\tau, s,z )| 
 \int_\tau   \left( \int_y | \overline{g(\tau ,  y) }|\right)^2   \int_s \sup_z | K_- (\tau, s,z )|\\
   & \leq    
 \sup_s \int_\tau  \sup_z | K_- (\tau, s,z )| 
 \int_\tau   \left( \int_y | \overline{g(\tau ,  y) }|\right)^2   \int_s \sup_z | K_- (\tau, s,z )|\\
 & \leq  \left( \sup_s \int_\tau  \sup_z | K_- (\tau, s,z )|  \right)^2 \int^t_0 \| g(s) \|_{L^1}^2 ds 
\end{align*}
Therefore we conclude that 
\begin{equation}\label{est:v(t)}
\begin{split}
    \| v(t) \|_{L^2_x}^2 \leq  \left( \sup_s \int_\tau  \sup_z | K_- (\tau, s,z )|  \right)\int^t_0 \| g(s) \|_{L^1}^2 ds \\
        \| \p_t  v(t) \|_{\dot{H}^{-1}}^2 \leq  \left( \sup_s \int_\tau  \sup_z | K_+ (\tau, s,z )|  \right)\int^t_0 \| g(s) \|_{L^1}^2 ds 
        \end{split}
\end{equation}

Step 2. Introduce that 
\begin{equation}
    M(\lambda, z ) = \int_{|\xi| \leq 1} \cos (\lambda |\xi|) e^{- i z \cdot \xi} \frac{d \xi}{|\xi|^2}
\end{equation}
Then 
\begin{equation}
    K_\pm (\tau, s, z) = \frac{C}{2} \{
    M(\tau - s, z) \pm M (2t -s - \tau, z)
    \}
\end{equation}

\hide
We will apply the Schur's lemma: 
\begin{lemma}
For $Tf (x) = \int K(x,y) f(y) dy$ we have that 
\begin{equation}
    2 \| T \|_{L^2 \to L^2}  \leq \sup_y \int |K(x,y)| d x  + \sup_x \int |K(x,y)| d y  
\end{equation}
\end{lemma}
\unhide

Suppose \eqref{est:M} given:
\begin{equation}\label{est:M}
    \sup_z |M(\lambda, z) | \leq \frac{20 \pi}{1 + |\lambda|}
\end{equation}

 Then we have that 
\begin{equation}
    \sup_{z} |K_\pm (\tau,s, z)| \lesssim \frac{1}{1+ |\tau -s|} + \frac{1}{1 + |2t - \tau - s|}
\end{equation}
and hence 
\begin{equation}
    \sup_\tau \int^t_0  \sup_{z} |K_\pm (\tau,s, z)|  \leq  \sup_\tau \left(\ln (1+\tau) + \ln (1+ 2t - \tau)\right) \lesssim \ln (1+ t)
\end{equation}
Then by the Step 1 we are done.

Step 3. We claim \eqref{est:M}.

Clearly $|M(\lambda, z)| \lesssim 1 $.

In polar coordinate 
\begin{equation}
\begin{split}
    M(\lambda, z) &= 2 \pi \int^1_0 d r \int^1_{-1} d \omega \cos (\lambda r) e^{i |z| r \omega}\\
    & =  \pi \int^1_{-1} \frac{\sin (|z|\omega + \lambda)}{ |z| \omega + \lambda}
 + \frac{\sin (|z|\omega -\lambda)}{ |z| \omega - \lambda}  d \omega 
\end{split}
\end{equation}
But 
\begin{equation}
    \int^1_{-1} \frac{\sin (|z|\omega \pm \lambda)}{ |z| \omega \pm \lambda} d \omega
    = \frac{1}{|z|} \int^{|z| \pm \lambda}_{-|z| \pm \lambda} \frac{\sin (a)}{a} da \leq \frac{8 \pi}{|z|}
\end{equation}
Hence
\begin{equation}
   \sup_{|z| \geq \frac{\lambda}{2}} |M(\lambda, z)| \lesssim \frac{1}{|\lambda|}
\end{equation}

Combining with this and the boundedness at the beginning of this step, we prove the claim \eqref{est:M}.

Step 4. For $\lambda \geq 0$, and $\hat g(t, \xi)$ supported in $\{|\xi| \leq \lambda\}$. Introduce 
\begin{equation}
    g_\lambda (t,x) = \frac{1}{\lambda^2} g (\frac{t}{\lambda}, \frac{x}{\lambda}), \ \   v_\lambda (t,x) =   v (\frac{t}{\lambda}, \frac{x}{\lambda})
\end{equation}
which solves 
\begin{equation}
    \Box v_\lambda = g_\lambda , \ \ v_\lambda |_{t=0} =0 ,  \ \  \p_t v_\lambda |_{t=0} = 0. 
\end{equation}
where $\hat{g}_\lambda (t, \xi) = \lambda \hat{g} (\frac{t}{\lambda},  \lambda \xi)$ supported in $\{ |\xi| \leq 1 \}$. Applying the result of previous steps we proof the lemma.

\end{proof}

\begin{proposition}[Proposition 6.3 of \cite{JMR}]
Consider $f \in C^\infty (\R; \mathcal S(\R^3)), v_0 \in \mathcal S(\R^3)$ and $v_1 \in \mathcal S (\R^3)$. Let $v \in L^\infty_{loc} (\R_+; L^2(\R^3))$ solves 
\begin{equation}
\Box v = f, \ \ v|_{t=0}= v_0 , \ \ \p_t v |_{t=0} = v_1 .
\end{equation}
Define 
\begin{equation}
    S^\lambda = \varphi (\lambda^{-1} D_x),
\end{equation}
where $\varphi \in C^\infty_0 (\R^3)$ real, $\varphi \equiv 1$ when $|\xi| \leq 1$ and $\text{spt} (\varphi) \in \{|\xi| \leq  2\}$. $S^\lambda$ is a pseudo-differential operator such that $S^\lambda v = f$ means that 
\begin{equation}
    \varphi (\frac{\xi}{\lambda}) \hat v (\xi) = \hat{f}(\xi),
\end{equation}
or
\begin{equation}
    S^\lambda v(x) = \frac{1}{(2 \pi)^n} \int_{\R^n} e^{i x \cdot \xi}
    \varphi (\frac{\xi}{\lambda}) \hat v (\xi) d \xi. 
\end{equation}

Then there exists $C(\varphi)$ does not depend of $\tau $ and $v$:
    \begin{equation}
        \sqrt{\int^\tau_0
        \| S^\lambda (v) (t) \|_{L^\infty_x}^2 dt 
        } \leq C(\varphi) \sqrt{\ln (1+ 2 \lambda \tau )}
        \left(\| v_0 \|_{\dot{H}^1_x} + \| v_1 \|_{L^2_x} + \int^\tau_0 \| f(t) \|_{L^2_x}dt\right)
    \end{equation}
\end{proposition}

\begin{corollary}

Apply the result with $\tau = T/\e$ and $v(t, x) = u^\e (\e t, x)$. If $\e \p_t u^\e, \p_x u^\e =O(1) $ and if 
\begin{equation}\label{wave_e}
    \e^2 \p_t^2 u^\e - \Delta u^\e = \e O(1)  \ \ \ L^1(0,T; L^2 (\R^3))
\end{equation} 
then 
\begin{equation}\label{est:JMR_e}
  \sqrt{\int^T_0
        \| S^\lambda (u) (t) \|_{L^\infty_x}^2 dt 
        }  
        \leq C  (\varphi)\sqrt \e  \sqrt{\ln (1+ 2 \lambda  \frac{T}{\e} )}
        \left(\|  u^\e_0  \|_{\dot{H}^1_x} + \| \e \p_t u^\e_0 \|_{L^2_x} + O(1)\right) 
\end{equation}
\end{corollary}

\begin{proof}[Proof of Proposition]
Suppose Lemma and \eqref{est:v}, \eqref{est:v_t} are given. 

Note that 
\begin{equation}\label{norm:SL}
    \| S^\lambda (u) (t) \|_{L^2 (0,T; L^\infty)} = \sup_{\| g \|_{L^2 (0,T; L^1)} \leq 1} \int_{[0,T] \times \R^3} S^\lambda (u) g dx dt
\end{equation}
For each $g$, consider a linear wave solution 
\begin{equation}
    \Box v = g , \ \ v|_{t=T} = 0, \ \ \p_t v |_{t=T} = 0.
\end{equation}

We also note that 
\begin{equation}
    \int_{[0,T] \times \R^3} S^\lambda (u) g dx d t = \int_{[0,T] \times \R^3}u S^\lambda ( g) dx d t
\end{equation}

Note that $S^\lambda$ and $\Box$ commute so that 
\begin{equation}
    \Box S^\lambda v = S^\lambda g,  \ \ S^\lambda v|_{t=T} =0 , \ \ \p_t S^\lambda v|_{t=T} =0. 
\end{equation}

Using the integration by parts, we identity the RHS of \eqref{norm:SL} to 
\begin{align*}
   & \int_{[0,T] \times \R^3}
    S^\lambda (u)g dx dt  \\
    & = \int_{[0,T] \times \R^3}
    u S^\lambda (g) dx dt \\
    & = \int_{[0,T] \times \R^3}
    u \Box S^\lambda (v) dx dt\\
     & = -\int_{[0,T] \times \R^3}
  \p_t   u  \p_t  S^\lambda (v) dx dt + \int_{\R^3} u \p_t S^\lambda (v)  dx  |_{t=0}^{t=T} - \int_{[0,T] \times \R^3}
 \Delta  u    S^\lambda (v) dx dt \\
  & =  \int_{[0,T] \times \R^3}
  \p_t\p_t    u    S^\lambda (v) dx dt 
  - \int_{\R^3} \p_t   u S^\lambda (v) dx \Big|_{t=0}^{t=T}
  + \int_{\R^3} u \p_t S^\lambda (v)  dx  \Big|_{t=0}^{t=T} - \int_{[0,T] \times \R^3}
 \Delta  u    S^\lambda (v) dx dt \\
 & =  \int_{[0,T] \times \R^3}
  \Box    u    S^\lambda (v) dx dt  + \int_{\R^3}  \p_t  u (0) S^\lambda (v)(0) dx  
  -  \int_{\R^3} u (0) \p_t S^\lambda (v) (0)  dx  
\end{align*}
Therefore 
\begin{equation}
  \left|  \int_{[0,T] \times \R^3}
    S^\lambda (u)g dx dt\right| \leq \| f \|_{L^1 (0,T; L^2_x)} \| S^\lambda (v) \|_{L^\infty (0,T; L^2_x)} + \| u_0 \|_{\dot{H}^1} \| \p_t S^\lambda (v) (0) \|_{\dot{H}^{-1}}
 + \| u_1 \|_{L^2_x} \| \p_t S^\lambda (v)(0) \|_{L^2_x}
\end{equation}

Now we use \eqref{est:v} and \eqref{est:v_t}, and the fact that $\| S^\lambda w\|_{L^2_x} \leq \|  w \|_{L^2_x}$. Therefore, we derive that 
\begin{align}
&\| S^\lambda (u) \|_{L^2 (0,T; L^\infty_x)}\\
& \leq \max\{\| v \|_{L^\infty(0,T; L^2_x)},\| \p_t v \|_{L^\infty(0,T; L^2_x)}\}
\left(
\| f \|_{L^1 (0,T; L^2_x)} +  \| u_0 \|_{\dot{H}^1}  + \| u_1 \|_{L^2_x}
\right)\\
& \leq \sqrt{C\ln (1+ \lambda T)}  \left(
\| f \|_{L^1 (0,T; L^2_x)} +  \| u_0 \|_{\dot{H}^1}  + \| u_1 \|_{L^2_x}
\right),
\end{align}

\end{proof}

The following lemma provides the bound for the high frequency part $|\xi| \geq \lambda $:
\begin{lemma}
$u$ solving \eqref{wave_e} should satisfies that 
\begin{align}
    \| u  \|_{L^2_TL_x^\infty}  &\leq  \sqrt{\e } \left(\sqrt{ \ln (1+ 2   T/\e ) }    \left(
\| O(1)\|_{L^1 (0,T; L^2_x)} +  \| u_0 \|_{\dot{H}^1}  + \| u_1 \|_{L^2_x}
\right)  + 
    \| \p^2_x u \|_{L^2_TL^2_x} \right)
\end{align}
\end{lemma}
\begin{proof}
Since $u^\e = S^\lambda u^\e + (I - S^\lambda )u^\e$ and 
\[
\| (I - S^\lambda ) u \|_{L^\infty_x} \lesssim
\|  \widehat{(I - S^\lambda ) u} \|_{L^1_\xi} \lesssim  
2^{-k/2} \| \p^2 u \|_{L^2}
\]
for $\lambda = 2^k$. (See the Besov space part)

Therefore we derive that 
\begin{align}
    \| u  \|_{L^2_TL_x^\infty}  &\leq \sqrt{ \ln (1+ 2 \lambda T ) }    \left(
\| f \|_{L^1 (0,T; L^2_x)} +  \| u_0 \|_{\dot{H}^1}  + \| u_1 \|_{L^2_x}
\right)  + 
    \frac{1}{\sqrt \lambda} \| \p^2_x u \|_{L^2_TL^2_x}
\end{align}

For our case, we choose $\lambda= 1/\e$ and hence $u$ solving \eqref{wave_e} should satisfies that 
\begin{align}
    \| u  \|_{L^2_TL_x^\infty}  &\leq  \sqrt{\e } \left(\sqrt{ \ln (1+ 2   T/\e ) }    \left(
\| O(1)\|_{L^1 (0,T; L^2_x)} +  \| u_0 \|_{\dot{H}^1}  + \| u_1 \|_{L^2_x}
\right)  + 
    \| \p^2_x u \|_{L^2_TL^2_x} \right)
\end{align}

\end{proof}
\unhide

\section{Example of the Initial Data via Mollification}\label{A.molli}

In this section, we give an example of a family of initial data $\{F_0^{\e}\}_{\e>0}$ such that, after mollification, it satisfies the Admissible Blow-up Condition~\eqref{ABC1}.

\begin{definition}
We denote by $\varphi^{\mathfrak{d}}$ a standard mollifier satisfying
\begin{align}\label{molli}
\bega
\int_{\R^d} \varphi(x) dx=1, \quad \varphi(x) \geq0, \quad \varphi(x) \in C_c^{\infty}(\R^d), \quad \varphi^{\mathfrak{d}}(x) := \frac{1}{\mathfrak{d}^{d}}\varphi(x/\mathfrak{d}).
\enda
\end{align}
\end{definition}

\hide
\begin{lemma}\label{L.molli}
For a given function $u(x)$, we define
\begin{align}
\bega
u^{\e} := \varphi^{\mathfrak{d}} \ast u
\enda
\end{align}
where $\int_{\R^d} \varphi dx=1$, $\varphi(x)\geq0$, $\varphi \in C_c^{\infty}(\R^d)$, $\varphi^{\mathfrak{d}}(x) := \frac{1}{\mathfrak{d}^d}\varphi(x/\mathfrak{d})$. Then the followings hold:
\begin{itemize}
\item If $u \in H^k(\R^d)$, then $\|u^{\e}\|_{H^{k+m}(\R^d)} \leq \frac{C}{\mathfrak{d}^{m}}\|u\|_{H^k(\R^d)}$, for $m,k \in \mathbb{Z}^+\cup \{0\}$.
\item If $u \in L^p$, then $\|u^{\e}\|_{L^p}\leq C\|u\|_{L^p}$ for $1\leq p\leq \infty$.
\item If $u\in \mathcal{M}(\R^d)$, then $\|u^{\e}\|_{L^1}\leq C\|u\|_{\mathcal{M}}$, where $\mathcal{M}(\R^d)$ is space of the Radon measures on $\R^d$.
\end{itemize}
\end{lemma}
\begin{proof}
(1) Applying Young's convolution inequality gives 
\begin{align*}
\bega
\|u^{\e}\|_{H^{k+m}} &\leq \sum_{|\beta_x|=m}\sum_{|\gamma_x|=k}\|\p^{\beta_x}\varphi^{\mathfrak{d}}\|_{L^1} \|\p^{\gamma_x}u\|_{L^2} \leq  \frac{C}{\mathfrak{d}^m}\|u\|_{H^k}.
\enda
\end{align*} 
(2) Applying Young's convolution inequality gives 
\begin{align*}
\bega
\|u^{\e}\|_{L^p}  = \|\varphi^{\mathfrak{d}} \ast u\|_{L^p} \leq \|\varphi^{\mathfrak{d}}\|_{L^1} \|u\|_{L^p} \leq \|u\|_{L^p}.
\enda
\end{align*}
(3) By definition 
\begin{align*}
\bega
\|u^{\e}\|_{L^1}  = \int_{\R^d}\bigg|\int_{\R^d}\varphi^{\mathfrak{d}}(x-y) du(y)\bigg|dx \leq \int_{\R^d}\int_{\R^d}\varphi^{\mathfrak{d}}(x-y)dx d|u|(y) \leq \|u\|_{\mathcal{M}}.
\enda
\end{align*}

Example: 
\begin{align}
\bega
\|\bar{u}^{\e}\|_{L^\infty_x} \leq \e^{-\frac{4}{3}} \|\bar{u}\|_{L^3_x} \leq C 
\enda
\end{align}
Note that we chose radial vorticity to decay near $r=0$. So that, 
\begin{align}
\bega
|\bar{u}^{\e}(t,x)| \leq \frac{1}{|x|} \int_0^{|x|} r\bar{\w}(r) dr \leq C|x|^n, \quad \mbox{near} \quad |x|=0
\enda
\end{align}
Thus, we have 
\begin{align}
\bega
\|\bar{u}^{\e}\|_{L^p_x}^p &\leq \int_{|x|\leq 1} |\bar{u}^{\e}|^p dx + \int_{|x|>1} |\bar{u}^{\e}|^p dx \cr 
&\leq C\int_{|x|\leq 1} |x|^{np} dx + C\int_{|x|>1} \frac{1}{|x|^p} dx \bigg(\int_{\R^2}\w(x)dx\bigg)^p \cr 
&\leq C,
\enda
\end{align}
for $p>2$.
\end{proof}
\unhide

\begin{lemma}\label{L.molli}
Let $\Omega = \R^2$. We choose initial data 
$F_0= M_{[e^{\e\rho_0},\e u_0, e^{\e\ta_0}]}+\AC{\P}F_0$ satisfying the following conditions:
\begin{itemize}

\item The initial macroscopic and microscopic parts belong to the following spaces:
\begin{align*}
(\rho_0,u_0-\bar{u},\ta_0)\in L^2_x(\R^2), \qquad \frac{1}{\e}\frac{(F_0-\mu)}{\sqrt{\tilde{\mu}}}\in L^2_x(\R^2 ; L^2_v(\R^3)),
\end{align*}
where $\tilde{\mu}=M_{[1,0,1-c_0]}$ for any $0<c_0\ll 1$.

\item The initial velocity tail is controlled in a weak $L^\infty_{x,v}$ sense:
\begin{align}\label{ini.tail}
\bega
&\sup_{0<\e\ll1}\bigg\|\frac{\e^{1-}}{\e}\frac{(F|_{t=0}-\mu)}{\sqrt{\tilde{\mu}}}\bigg\|_{L^\infty_{x,v}}\les 1, \quad \mbox{for the purely spatial derivative case $\eqref{caseA}$}, \cr 
&\sup_{0<\e\ll1}\bigg\|\frac{\e^{1-}}{\e}\frac{(\e^{\mathfrak{n}}\p_tF|_{t=0}-\mu)}{\sqrt{\tilde{\mu}}}\bigg\|_{L^\infty_{x,v}}\les 1, \quad \mbox{for the space-time derivative case $\eqref{caseB}$}.
\enda
\end{align}
\end{itemize}

For the standard mollifier \eqref{molli} with the following mollification rate,
\begin{align*}
\bega
\mathfrak{d}(\e):= \log\!\Big(\log\!\big(\log(1/\e)\big)\Big)^{-\frac{1}{10}},
 \enda
\end{align*}
we define the mollified macroscopic and microscopic fields by
\begin{align*}
\bega
(\rho^{\e}_0,u^{\e}_0,\ta^{\e}_0)
&:= \big(\rho_0 \ast \varphi^{\mathfrak{d}},\,
\widetilde{u}_0 \ast \varphi^{\mathfrak{d}} + \bar{u},\,
\ta_0 \ast \varphi^{\mathfrak{d}}\big), \qquad \AC{\P}F^{\e}_0 &:= \AC{\P}F_0 \ast \varphi^{\mathfrak{d}}.
\enda
\end{align*}
Here, $\bar{u}$ denotes the radial eddy defined in Definition~\ref{D.Ra-E} for $d=2$.
Then, the following mollified sequence of initial data

\begin{align}\label{molliseq}
\bega
F^{\e}_0(x,v) := M_{[e^{\e\rho^{\e}_0},\,\e u^{\e}_0,\,e^{\e\ta^{\e}_0}]} + \AC{\P}F^{\e}_0,
\enda
\end{align}
satisfies the Admissible Blow-up Condition~\eqref{ABC1}.
\hide
\begin{align}
\bega
F^{\e}_0(x,v) := M_{[e^{\e(\rho_0\ast \varphi^{\mathfrak{d}})},\,\e ((u_0-\bar{u})\ast \varphi^{\mathfrak{d}}+\bar{u}),\,e^{\e(\ta_0\ast \varphi^{\mathfrak{d}})}]} + \AC{\P}F_0\ast \varphi^{\mathfrak{d}},
\enda
\end{align}
\unhide
\hide
(1) For arbitrary time $T_*>1$, once we choose $\mathfrak{d}(\e)$ as 
\begin{align}\label{ek-rel:ex}
\mathfrak{d}(\e) := 
\left[\log\!\left(\log\!\Big(\log\!\big((1+T_*)^5\e^{-\mathfrak{s}}\big)\Big)\right)\right]^{-\frac{1}{\mathrm{N}+\mathfrak{j}}}, \quad \textcolor{red}{\mbox{Need to revise later $T_*^5$}}
\end{align}
for any $0<\mathfrak{s}\ll1$, then, $F_0^{\e}$ satisfies the two bootstrap assumption \eqref{Boot1} and \eqref{Boot2}. 
\unhide
\end{lemma}
\begin{proof}
By elementary properties of mollification, we readily obtain
\begin{align*}
\bega
&\sum_{1 \leq |\alpha| \leq 5}
\| \p^\alpha (\rho^\e_0, u^\e_0, \theta^\e_0) \|_{L^2_x}^2
\leq C \mathfrak{d}^{-10}\|(\rho_0,u_0,\theta_0)\|_{L^2_x}^2 \les \log\!\Big(\log\!\big(\log(1/\e)\big)\Big).
\enda
\end{align*}
For the microscopic part, the same argument yields
\begin{align*}
\bega
\sum_{1 \leq |\alpha| \leq 5}\frac{1}{\e^2}\int_{\Omega\times\R^3} \frac{|\p^{\al}\AC{\P}F^{\e}_0|^2}{\tilde{\mu}} dvdx &\leq \sum_{1 \leq |\alpha| \leq 5} \frac{C}{\e^2}\int_{\Omega\times\R^3} \frac{|\p^{\al}(\AC{\P}F_0 \ast \varphi^{\mathfrak{d}})|^2}{\tilde{\mu}} dvdx \leq \mathfrak{d}^{-10} \frac{1}{\e^2} \left\| \frac{\AC{\P}F_0}{\sqrt{\tilde{\mu}}} \right\|_{L^2_{x,v}}^2.
\enda
\end{align*}
We decompose $\AC{\P}F_0= F_0-M_0$, where $M_0=M_{[e^{\e\rho_0},\e u_0, e^{\e\ta_0}]}$, to obtain
\begin{align*}
\bega
\mathfrak{d}^{-10} \frac{1}{\e^2} \left\| \frac{\AC{\P}F_0}{\sqrt{\tilde{\mu}}} \right\|_{L^2_{x,v}}^2 &\leq \mathfrak{d}^{-10} \bigg(\frac{1}{\e^2} \left\| \frac{F_0-\mu}{\sqrt{\tilde{\mu}}} \right\|_{L^2_{x,v}}^2 + \frac{1}{\e^2} \left\| \frac{M_0-\mu}{\sqrt{\tilde{\mu}}} \right\|_{L^2_{x,v}}^2 \bigg) \les \log\!\Big(\log\!\big(\log(1/\e)\big)\Big),
\enda
\end{align*}
where we used the Taylor expansion for $M_0-\mu$ in \eqref{MFtaylor}.
For the $L^\infty_{x,v}$ part of \eqref{ABC1}, applying Young's convolution inequality gives
\begin{align*}
\bega
\left\| \p^{\al} \left(\frac{F^\e_0 - \mu}{ \sqrt{\tilde{\mu}} }
\right) \right\|_{L^\infty_{x,v}} 
&= \left\| \left(\frac{F_0 - \mu}{ \sqrt{\tilde{\mu}} }\right)\ast \p^{\al}\varphi^{\mathfrak{d}} \right\|_{L^\infty_{x,v}} \cr 
&\leq \left\| \left(\frac{F_0 - \mu}{ \sqrt{\tilde{\mu}} }\right)\right\|_{L^\infty_{x,v}}\|\p^{\al}\varphi^{\mathfrak{d}} \|_{L^1_x} \leq \frac{C}{\mathfrak{d}^{5}} \leq \log\!\Big(\log\!\big(\log(1/\e)\big)\Big).
\enda
\end{align*} 
\end{proof}

\hide
(If we mollified $F_0^{\e} = F_0\ast \varphi^{\mathfrak{d}}$)
Since $(\mathrm{P}^{\e}, \mathrm{U}^{\e}, \mathrm{\Theta}^{\e}) := (e^{\e\rho^{\e}}, \e u^{\e}, e^{\e\ta^{\e}})$, we write 
\begin{align*}
\rho^{\e}_0= \frac{1}{\e}\ln\bigg(\int_{\R^3}F_0^{\e}dv\bigg), \quad
u^{\e}_0 = \frac{1}{\e}\bigg(\frac{\int_{\R^3}vF_0^{\e}dv}{\int_{\R^3}F_0^{\e}dv}\bigg), \quad 
\ta^{\e}_0=\frac{1}{\e}\ln\bigg(\frac{\int_{\R^3}\frac{|v|^2}{2}F_0^{\e}dv-\frac{1}{2}\frac{(\int_{\R^3}vF_0^{\e}dv)^2}{\int_{\R^3}F_0^{\e}dv} }{\frac{3k_B}{2}\int_{\R^3}F_0^{\e}dv}\bigg).
\end{align*}
By Leibniz’s rule, we have 
\begin{align}\label{x-ini}
&\|\p^{\al_x}(\rho^{\e}_0, u^{\e}_0, \ta^{\e}_0)\|_{L^2_x}^2 \leq C\mathfrak{d}^{-2(|\al_x|+m)} \|(\rho_0, u_0, \ta_0)\|_{H^{-\mathfrak{j}}_x}^2 . 
\end{align}
\unhide

\hide
Let $\rho_0 \in H^{-m}$ and $u_0 \in H^3$. Can we have $\|u^{\e}_0-u_0\|_{H^3} \to 0$?
\begin{align*}
\rho^{\e}_0= \frac{1}{\e}\ln\bigg(\int_{\R^3}F_0^{\e}dv\bigg), \quad
u^{\e}_0 = \frac{1}{\e}\bigg(\frac{\int_{\R^3}vF_0^{\e}dv}{\int_{\R^3}F_0^{\e}dv}\bigg), \quad 
\ta^{\e}_0=\frac{1}{\e}\ln\bigg(\frac{\int_{\R^3}\frac{|v|^2}{2}F_0^{\e}dv-\frac{1}{2}\frac{(\int_{\R^3}vF_0^{\e}dv)^2}{\int_{\R^3}F_0^{\e}dv} }{\frac{3k_B}{2}\int_{\R^3}F_0^{\e}dv}\bigg).
\end{align*}
\begin{align*}
\|u^{\e}_0-u_0\|_{H^3_x} &= \frac{1}{\e} \bigg\| \bigg(\frac{\int_{\R^3}vF_0dv \ast \varphi^{\mathfrak{d}} }{\int_{\R^3}F_0dv \ast \varphi^{\mathfrak{d}} } - \frac{\int_{\R^3}vF_0dv}{\int_{\R^3}F_0dv}\bigg) \bigg\|_{H^3_x} \cr 
&\leq \frac{1}{\e} \bigg\| \frac{1}{\int_{\R^3}F_0dv \ast \varphi^{\mathfrak{d}}} \bigg(\int_{\R^3}vF_0dv \ast \varphi^{\mathfrak{d}}- \int_{\R^3}vF_0dv\bigg) \bigg\|_{H^3_x} + \frac{1}{\e} \bigg\| \textcolor{red}{\int_{\R^3}vF_0dv}\bigg(\frac{1}{\int_{\R^3}F_0dv \ast \varphi^{\mathfrak{d}}} - \frac{1}{\textcolor{red}{\int_{\R^3}F_0dv}}\bigg) \bigg\|_{H^3_x}
\end{align*}
\unhide

\begin{lemma}\label{L.molli3d}
Let $\Omega = \R^3$. We choose initial data
$F_0= M_{[e^{\e\rho_0},\e u_0, e^{\e\ta_0}]}+\AC{\P}F_0$ satisfying the following conditions:
\begin{itemize}
\item For some $\mathrm{N}\geq3$, the initial macroscopic and microscopic parts satisfy
\begin{align*}
(\rho_0,u_0,\ta_0)\in H^{\mathrm{N}}_x(\R^2), \qquad \frac{1}{\e}\frac{(F_0-\mu)}{\sqrt{\tilde{\mu}}}\in H^{\mathrm{N}}_x(\R^2 ; L^2_v(\R^3)).
\end{align*}
where $\tilde{\mu}=M_{[1,0,1-c_0]}$ for any $0<c_0\ll 1$.

\item The initial velocity tail is controlled in a weak $L^\infty_{x,v}$ sense as in \eqref{ini.tail}.

\end{itemize}
Then the mollified initial data defined as in \eqref{molliseq} satisfies the assumptions of Theorem~\ref{T.3D.unif}.
\end{lemma}
\begin{proof}
The conclusion follows immediately since the mollified initial data satisfy the assumptions of Lemma~\ref{L.ini3d}.
\end{proof}

\section{Strichartz estimate}\label{A.Stri}

In this part, we introduce several notations needed for the Strichartz estimate and then state the Strichartz estimate used in Section \ref{Sec.macro.Stri}.

\begin{definition}
We define the Fourier transform and its inverse by
\begin{align*}
\mathcal{F}(g)(\xi) = \hat{g}(\xi) := \int_{\R^d} e^{-ix\cdot \xi} g(x) dx, \qquad \mathcal{F}^{-1}(\hat{g})(x) := \frac{1}{(2\pi)^d}\int_{\R^d} e^{ix\cdot \xi} \hat{g}(\xi) d\xi.
\end{align*}
For $s\in\R$, we define the Fourier multiplier $|D|^s$ by
\begin{align}\label{fourierDdef}
|D|^s f := \mathcal{F}^{-1}\big(|\xi|^s\hat{f}(\xi)\big).
\end{align}
Let $\chi,\varphi \in C_c^\infty(\R^d)$ be radial functions such that
\[
\mathrm{supp}\,\chi \subset \{|\xi|\le \tfrac43\}, \qquad
\mathrm{supp}\,\varphi \subset 
\{\tfrac34 \le |\xi| \le \tfrac83\},
\]
and
\[
\chi(\xi)+\sum_{j\ge0}\varphi(2^{-j}\xi)=1
\quad \text{for all } \xi\in\R^d, \qquad  
\sum_{j\in\Z}\varphi(2^{-j}\xi)=1
\quad \text{for } \xi\neq0.
\]
The inhomogeneous and homogeneous dyadic blocks are defined by
\begin{align*}
&\Delta_jf=0, \quad \mbox{if} \quad j\leq -2, \quad \Delta_{-1}f:=\mathcal{F}^{-1}\!\big(\chi(\xi)\hat f(\xi)\big),
\quad
\Delta_j f:=\mathcal{F}^{-1}\!\big(\varphi(2^{-j}\xi)\hat f(\xi)\big),
\quad j\ge0, \cr 
&\dot{\Delta}_j u
:= \mathcal{F}^{-1}\!\big(\varphi(2^{-j}\xi)\hat{u}(\xi)\big),
\qquad j\in\Z.
\end{align*}
For $s\in\R$ and $p,q \in [1,\infty]$, the homogeneous and inhomogeneous Besov norms are defined by
\begin{align}\label{Besovdef}
\bega 
\|u\|_{\dot{B}_{p,q}^{s}} := \bigg(\sum_{j\in \Z} 2^{qjs}\|\dot{\Delta}_j u\|_{L^p}^q \bigg)^{\frac{1}{q}}, \qquad
\|u\|_{B_{p,q}^{s}} := \bigg(\sum_{j\in \Z} 2^{qjs}\|\Delta_j u\|_{L^p}^q \bigg)^{\frac{1}{q}}.
\enda
\end{align}
For $T>0$, $s\in \R$, and $p,q,\varrho \in [1,\infty]$, we define the time--space Besov norm by
\begin{align}\label{timeBesovdef}
\bega 
\|u\|_{\widetilde{L}^{\varrho}_T(\dot{B}_{p,q}^{s})} := \bigg(\sum_{j\in \Z} 2^{qjs}\|\dot{\Delta}_j u\|_{L^{\varrho}_T(L^p)}^q \bigg)^{\frac{1}{q}}.
\enda
\end{align}
We denote by $\mathcal{S}_h'(\R^d)$ the space of tempered distributions satisfying
\begin{align*}
\bega
\lim_{\lambda \rightarrow \infty}\|\tau (\lambda D)u\|_{L^\infty} =0,
\enda
\end{align*}
for any $\tau \in \mathcal{D}(\R^d)$, where $\mathcal{D}(\R^d)$ denotes the space of smooth compactly supported functions on $\R^d$. 
\end{definition}

\begin{lemma}\label{L.Besov}\cite{Danchin,Triebel}
The following properties hold for Besov spaces.
\begin{itemize}
\item Let $1\leq p_1 \leq p_2 \leq \infty$ and $1\leq q_1 \leq q_2 \leq \infty$. 
For any $s\in\mathbb{R}$, the space $B_{p_1,q_1}^s$ is continuously embedded in $B_{p_2,q_2}^{s-d(\frac{1}{p_1}-\frac{1}{p_2})}$. In particular,
\begin{align}\label{Besov-H}
\bega
\|u\|_{\dot{B}_{2,2}^{s}} \leq C\|u\|_{H^{s}}.
\enda 
\end{align}

\item  For any $p,q\in[1,\infty]$ with $p\leq q$, the space $\dot{B}_{p,1}^{\frac{d}{p}-\frac{d}{q}}$ is continuously embedded into $L^q$.
\begin{align}\label{Besov-Lp}
\bega
\|u\|_{L^q} \leq C\|u\|_{\dot{B}_{p,1}^{\frac{d}{p}-\frac{d}{q}}},
\enda
\end{align}
for some constant $C>0$.

\item Let $s_1,s_2\in \mathbb{R}$ with $s_1<s_2$, $\tau\in(0,1)$, and $p,q\in[1,\infty]$. 
For any $u\in \mathcal{S}_h'$, the interpolation inequality holds for some constant $C>0$. 
\begin{align}\label{Besov-inter}
\bega
\|u\|_{\dot{B}_{p,1}^{\tau s_1+(1-\tau)s_2}} &\leq \frac{C}{s_2-s_1}\bigg(\frac{1}{\tau}+\frac{1}{(1-\tau)}\bigg) \|u\|_{\dot{B}_{p,\infty}^{s_1}}^{\tau} \|u\|_{\dot{B}_{p,\infty}^{s_2}}^{1-\tau} \cr 
&\leq \frac{C}{s_2-s_1}\bigg(\frac{1}{\tau}+\frac{1}{(1-\tau)}\bigg) \|u\|_{\dot{B}_{p,2}^{s_1}}^{\tau} \|u\|_{\dot{B}_{p,2}^{s_2}}^{1-\tau}.
\enda
\end{align}

\item The time-space Besov norms satisfy 
\begin{align}\label{Besovineq}
\bega 
\|u\|_{\widetilde{L}^{\varrho}_T(\dot{B}_{p,q}^{s})} \leq \|u\|_{L^{\varrho}_T(\dot{B}_{p,q}^{s})}, \quad q \geq \varrho, \qquad \|u\|_{\widetilde{L}^{\varrho}_T(\dot{B}_{p,q}^{s})} \geq \|u\|_{L^{\varrho}_T(\dot{B}_{p,q}^{s})}, \quad \varrho \geq q.
\enda
\end{align}

\item For the space $\mathcal{F}^s_p$ defined in Definition \ref{D.Fsdef}, if $p\in(1,\infty]$ and $s\in(0,1)$, then
\begin{align}\label{BesovFp}
\bega
\|u\|_{B^s_{p,\infty}(\R^2)} \leq \|u\|_{\mathcal{F}^s_p(\R^2)} \leq \|u\|_{B^s_{p,1}(\R^2)}.
\enda
\end{align}

\item Let $0<p\leq\infty$, $0<q_1,q_2<\infty$, $s\in\mathbb{R}$ and $\e>0$. Then
\begin{align}\label{Besov-p,q}
\bega
\|u\|_{B^{s}_{p,q_2}(\R^d)} \leq C\|u\|_{B^{s+\e}_{p,q_1}(\R^d)}.
\enda
\end{align}

\end{itemize}
\end{lemma}

\begin{proof}
The embeddings in the first statement and \eqref{Besov-Lp} follow from Propositions 2.71 and 2.39 in \cite{Danchin}. 
The interpolation inequality \eqref{Besov-inter} is given in Proposition 2.22 of \cite{Danchin}. 
The time–Besov inequalities \eqref{Besovineq} can be found on page 98 of \cite{Danchin}. 
Estimate \eqref{BesovFp} follows from Definition 3.30 in \cite{Danchin}. 
Finally, \eqref{Besov-p,q} is stated in Proposition 2 of Section 2.3.2 in \cite{Triebel}.
\end{proof}

\begin{proposition}\cite{Danchin}\label{P.Book1}
Let $(\varsigma_1,\varsigma_2)$ be the solution of \eqref{varsigeqn}. Then, for any $s \in \mathbb{R}$, we have 
\begin{align*}
\bega
\|(\varsigma_1,\varsigma_2)\|_{\widetilde{L}^r_T\dot{B}_{p,1}^{s+d(\frac{1}{p}-\frac{1}{2})+\frac{1}{r}}} \leq C \e^{\frac{1}{r}}\|(\varsigma_1(0),\varsigma_2(0))\|_{\dot{B}_{2,1}^s} + C\e^{1+\frac{1}{r}-\frac{1}{\bar{r}'}}\|(\bPhi_{\varsigma_1},\bPhi_{\varsigma_2})\|_{\widetilde{L}^{\bar{r}'}_T\dot{B}_{\bar{p}',1}^{s+d(\frac{1}{\bar{p}'}-\frac{1}{2})+\frac{1}{\bar{r}'}-1}},
\enda
\end{align*}
\begin{align*}
\bega
p \geq 2, \qquad \frac{2}{r} \leq \min \left\{ 1, (d-1)\bigg(\frac{1}{2}-\frac{1}{p}\bigg)\right\}, \qquad (r,p,d) \neq (2,\infty,3), \cr 
\bar{p} \geq 2, \qquad \frac{2}{\bar{r}} \leq \min \left\{ 1, (d-1)\bigg(\frac{1}{2}-\frac{1}{\bar{p}}\bigg)\right\}, \qquad (\bar{r},\bar{p},d) \neq (2,\infty,3), \cr 
\enda
\end{align*}
where $\frac{1}{\bar{p}} + \frac{1}{\bar{p}'} = 1$ and $\frac{1}{\bar{r}} + \frac{1}{\bar{r}'} = 1$. Here, the Besov norm and the time-space Besov norm are defined in \eqref{Besovdef} and \eqref{timeBesovdef}, respectively.
\end{proposition}

\begin{proof}[Proof of Proposition \ref{P.Book}]
We choose $\bar{r}'=1$ and $\bar{p}'=2$ in Proposition \ref{P.Book1}. Then, applying the inequality \eqref{Besovineq} to the space-time Besov norm
\begin{align*}
\bega
\|\cdot\|_{\widetilde{L}^{1}_T\dot{B}_{2,1}^{s}} \leq \|\cdot\|_{L^{1}_T\dot{B}_{2,1}^{s}} , \qquad \|\cdot\|_{\widetilde{L}^{r}_T\dot{B}_{p,1}^{s}} \geq \|\cdot\|_{L^{r}_T\dot{B}_{p,1}^{s}},
\enda
\end{align*}
gives the result. 
\end{proof}

\begin{proof}[Proof of Proposition \ref{P.Book1}]
Although the proof can be found in Proposition 10.30 of \cite{Danchin}, we briefly sketch the idea.
Let $\varsigma := (\varsigma_1,\varsigma_2)^T$ and 
$\bPhi := (\bPhi_{\varsigma_1},\bPhi_{\varsigma_2})^T$. 
Then we can rewrite~\eqref{varsigeqn} in the following form:
\begin{align*}
\bega
\p_t \varsigma = i \e^{-1} |D| J \varsigma + \bPhi, \qquad \varsigma(t,x)|_{t=0} = \varsigma_0(x) \qquad J:= \left( \begin{array}{cc} 0 & -1 \\ 1 & 0  \end{array} \right).
\enda
\end{align*}
We introduce the rescaled function
\[
\bar{\varsigma} (t,x)
:= \varsigma\Big(\frac{\varepsilon}{N}t, \frac{1}{N}x\Big), \qquad  \bar{\bPhi}(t,x):= \bPhi\Big(\frac{\varepsilon}{N}t, \frac{1}{N}x\Big).
\]
Then $\bar{\varsigma}$ satisfies the unit–frequency equation
\begin{align*}
\bega
\p_t \bar{\varsigma} = i |D| J \bar{\varsigma} +  \frac{\e}{N}\bar{\bPhi}, \qquad \bar{\varsigma}(t,x)|_{t=0} = \varsigma_0(x/N).
\enda
\end{align*}
Let $\mathbf{U}(t) := e^{\,it|D|J}$ denote the half–wave propagator on $\R^2$.  
Choose a smooth bump function $\varphi \in C_c^{\infty}((0,\infty))$ supported in $\{\,1/2 \le |\xi| \le 2\,\}$, and define the Littlewood–Paley projection
\[
P_N \varsigma=\mathcal{F}^{-1}\big(\varphi(|\xi|/N)  \hat{\varsigma}(\xi)\big),
\quad N\in2^{\mathbb N}.
\]
Then, by the Keel–Tao theorem for dispersive kernels with $(1+|t|)^{-1/2}$ decay (in the non-endpoint case), we have
\begin{align*}
\|\mathbf{U}(\cdot)P_1 \bar{\varsigma}_0\|_{L_t^{r}(L_x^{p})}
\les \|P_1\bar{\varsigma}_0\|_{L_x^2}, \qquad \Big\|\int_{0}^{t} \mathbf{U}(t-s)P_1\bar{\bPhi}(s)\,ds\Big\|_{L_t^{r}(L_x^{p})}
&\les
\|P_1\bar{\bPhi}\|_{L_t^{\bar{r}'}(L_x^{\bar{p}'})}.
\end{align*}
Consequently,
\[
\|P_1\bar{\varsigma}\|_{L_t^{r}L_x^{p}}
\ \les
\|P_1\bar{\varsigma}_0\|_{L_x^2}+\frac{\e}{N}\|P_1\bar{\bPhi}\|_{L_t^{\bar{r}'}L_x^{\bar{p}'}}.
\]
Applying the change of variables $y = x/N$ and $\tau = (\e/N)t$, we obtain the scaling relations
\begin{align*}
\|P_1\bar{\varsigma}\|_{L_t^{r}(L_x^{p})}
= N^{\frac{d}{p}}\Big(\frac{N}{\e}\Big)^{\!\frac1q}\|P_N\varsigma\|_{L_\tau^{r}L_y^{p}},
\qquad
\|P_1\bar{\varsigma}_0\|_{L_x^2}=N^{\frac d2}\|P_N\varsigma_0\|_{L_y^2}.
\end{align*}
Hence,
\begin{align*}
\bega
\Big(\frac{N}{\e}\Big)^{\frac{1}{r}} N^{\frac{d}{p}}\|P_N \varsigma \|_{L^r_tL^p_x} \leq C N^{\frac{d}{2}}\|P_N\varsigma_0\|_{L^2_x} + \Big(\frac{N}{\e}\Big)^{\frac{1}{\bar{r}'}} N^{\frac{d}{\bar{p}'}}\frac{\e}{N}\|P_N\bPhi\|_{L^{\bar{r}'}_tL^{\bar{p}'}_x}.
\enda
\end{align*}
Dividing both sides by $N^{\frac{d}{2}}$ and multiplying by $\e^{\frac{1}{r}}N^{s}$ yields
\begin{align*}
\bega
N^{s+d(\frac{1}{p}-\frac{1}{2})+\frac{1}{r}}\|P_N \varsigma \|_{L^r_tL^p_x} \leq C \e^{\frac{1}{r}}N^{s}\|P_N\varsigma_0\|_{L^2_x} + \e^{1+\frac{1}{r}-\frac{1}{\bar{r}'}}N^{s+d(\frac{1}{\bar{p}'}-\frac{1}{2})+\frac{1}{\bar{r}'}-1}\|P_N\bPhi\|_{L^{\bar{r}'}_tL^{\bar{p}'}_x}.
\enda
\end{align*}
Finally, summing over dyadic frequencies $N=2^j$ with $j\in\mathbb{Z}$ gives the desired estimate.
\hide
Taking $l^q(\mathbb{Z})$ gives 
\begin{align*}
\bega
\|\varsigma\|_{\widetilde{L}^r_T\dot{B}_{p,q}^{s+d(\frac{1}{p}-\frac{1}{2})+\frac{1}{r}}} \leq C \e^{\frac{1}{r}}\|\varsigma(0)\|_{\dot{B}_{2,q}^s} + C\e^{1+\frac{1}{r}-\frac{1}{\bar{r}'}}\|\bPhi\|_{\widetilde{L}^{\bar{r}'}_T\dot{B}_{\bar{p}',q}^{s+d(\frac{1}{\bar{p}'}-\frac{1}{2})+\frac{1}{\bar{r}'}-1}}.
\enda
\end{align*}
\unhide
\end{proof}

\hide
\begin{proposition} Let $d\geq 2$. We define $\mathcal{C}:=\{\xi\in\R^d~|~ r\leq|\xi|\leq R\}$ for positive constants such that $r<R$.
For the solution of the wave equation 
\begin{align}
\bega
\begin{cases} \p_t^2 u - \Delta_x u =0 , \cr 
(u,\p_tu)|_{t=0} = u_0, u_1, \end{cases}
\enda
\end{align}
if $\hat{u}_0$ and $\hat{u}_1$ are supported in $\mathcal{C}$, then we have 
\begin{align*}
\bega
\|u(t)\|_{L^\infty_x} \leq \frac{C}{t^{\frac{d-1}{2}}}(\|u_0\|_{L^1_x}+\|u_1\|_{L^1_x}), \qquad t>0. 
\enda
\end{align*}
\end{proposition}
\begin{proof}
We write that solution of the wave equation 
\begin{align*}
\bega
\hat{u}(t,\xi) = e^{it|\xi|}\frac{1}{2}\bigg(\hat{u}_0(\xi)+\frac{1}{i|\xi|}\hat{u}_1(\xi)\bigg) + e^{-it|\xi|}\frac{1}{2}\bigg(\hat{u}_0(\xi)-\frac{1}{i|\xi|}\hat{u}_1(\xi)\bigg)
\enda
\end{align*}
Let $\varphi(\xi)$ be a smooth compactly supported function in $\R^2/\{0\}$ with value 1 near $\mathcal{C}$. Multiplying $\varphi^2(\xi)$ to $\hat{u}$, then taking inverse fourier transform gives 
\begin{align}
\bega
u(t,x) &= \mathcal{F}^{-1}(e^{it|\xi|}\varphi) \ast \mathcal{F}^{-1}\bigg(\frac{1}{2}\bigg(\hat{u}_0(\xi)+\frac{1}{i|\xi|}\hat{u}_1(\xi)\bigg)\varphi \bigg) \cr 
&+ \mathcal{F}^{-1}(e^{-it|\xi|}\varphi) \ast \mathcal{F}^{-1}\bigg(\frac{1}{2}\bigg(\hat{u}_0(\xi)-\frac{1}{i|\xi|}\hat{u}_1(\xi)\bigg)\varphi \bigg)
\enda
\end{align}
Applying Young's convolution inequality, we have 
\begin{align}\label{Wuinf}
\bega
\|u(t)\|_{L^\infty_x} &\leq  \Big\|\mathcal{F}^{-1}(e^{it|\xi|}\varphi)\Big\|_{L^\infty_x} \bigg\| \mathcal{F}^{-1}\bigg(\frac{1}{2}\bigg(\hat{u}_0(\xi)+\frac{1}{i|\xi|}\hat{u}_1(\xi)\bigg)\varphi \bigg) \bigg\|_{L^1_x} \cr 
&+ \Big\|\mathcal{F}^{-1}(e^{-it|\xi|}\varphi)\Big\|_{L^\infty_x} \bigg\|\mathcal{F}^{-1}\bigg(\frac{1}{2}\bigg(\hat{u}_0(\xi)-\frac{1}{i|\xi|}\hat{u}_1(\xi)\bigg)\varphi \bigg)\bigg\|_{L^1_x}
\enda
\end{align}
For the latter term in both lines, we have 
\begin{align*}
\bega
\bigg\|\mathcal{F}^{-1}\bigg(\frac{1}{2}\bigg(\hat{u}_0(\xi)\pm \frac{1}{i|\xi|}\hat{u}_1(\xi)\bigg)\varphi \bigg)\bigg\|_{L^1_x} &\leq \frac{1}{2}\|u_0\|_{L^1_x} + \frac{1}{2} \bigg\|\mathcal{F}^{-1}\bigg(\frac{\varphi}{|\xi|}\bigg)\bigg\|_{L^\infty_x} \|u_1\|_{L^1_x} \cr 
&\leq C(\|u_0\|_{L^1_x}+\|u_1\|_{L^1_x}),
\enda
\end{align*}
where we used $\|\mathcal{F}^{-1}(\varphi/|\xi|)\|_{L^\infty_x} \leq \|(\varphi/|\xi|)\|_{L^1_{\xi}} \leq C$ by the Haudorff-Young inequality. 
Let 
\begin{align*}
\bega
K^{\pm}(t,x):= \int_{\R^d}e^{ix\cdot\xi}e^{\pm it|\xi|}\varphi d\xi.
\enda
\end{align*}
For the $L^\infty_x$ norm of \eqref{Wuinf}, we claim that 
\begin{align}\label{Kdecay}
\bega
\Big\|\int_{\R^d}e^{ix\cdot\xi}(e^{\pm it|\xi|}\varphi)d\xi\Big\|_{L^\infty_x} \leq \frac{C}{t^{\frac{d-1}{2}}}
\enda
\end{align}
Since $L^\infty_x$ does not depend on the dilation, we consider $K^{\pm}(t,tx)$: 
\begin{align*}
\bega
K^{\pm}(t,tx)= \int_{\R^d}e^{it (x\cdot\xi\pm|\xi|)}\varphi d\xi := \int_{\R^d}e^{it \Phi(\xi)}\varphi d\xi.
\enda
\end{align*}
We denote $\Phi(\xi) := x\cdot\xi\pm|\xi|$. $\nabla_{\xi}\Phi(\xi) = x\pm \frac{\xi}{|\xi|}$. We prove that 
\begin{align}
|K^{\pm}(t,tx)| &\leq \frac{C_N}{(c_0t)^N}, \qquad \mbox{if} \qquad |\nabla_{\xi}\Phi(\xi)|\geq c_0 \label{Kbdd1} \\
|K^{\pm}(t,tx)| &\leq C_{N'}\int_{K}\frac{1}{(1+c_0t|\nabla_{\xi}\Phi(\xi)|^2)^{N'}}d\xi, \qquad \mbox{if} \qquad |\nabla_{\xi}\Phi(\xi)|\leq c_0, \quad \xi\in K \label{Kbdd2}
\end{align}
for $c_0 \in (0,1]$ and compact set $K\subset \R^d$.
(Proof of \eqref{Kbdd1}) We define an operator 
\begin{align*}
La = -i \sum_{j=1}^d \frac{\p_j\Phi}{|\nabla_{\xi}\Phi|^2} \p_j a
\end{align*}
Then we have 
\begin{align*}
Le^{it\Phi} = -i \sum_{j=1}^d \frac{\p_j\Phi}{|\nabla_{\xi}\Phi|^2} it \p_j \Phi e^{it\Phi} = te^{it\Phi}
\end{align*}
By the integration by parts, we have 
\begin{align*}
K^{\pm}(t,tx) &= \frac{1}{t^N}\int_{\R^d}L^Ne^{it \Phi(\xi)}\varphi d\xi = (-1)^N\frac{1}{t^N}\int_{\R^d}e^{it \Phi(\xi)} (L^T)^N\varphi d\xi 
\end{align*}
Then we can prove that $|\int_{\R^d}e^{it \Phi(\xi)} (L^T)^N\varphi d\xi | \leq C$ when $|\nabla_{\xi}\Phi(\xi)|\geq c_0$.  \\
(Proof of \eqref{Kbdd2}) Similarly, we define an operator 
\begin{align*}
L_t = \frac{1}{1+t|\nabla_{\xi}\Phi(\xi)|^2}(Id - i \sum_{j=1}^d \p_j \Phi \p_j)
\end{align*}
Then we have 
\begin{align*}
L_te^{it\Phi} = e^{it\Phi}
\end{align*}
By the same way, we have 
\begin{align*}
K^{\pm}(t,tx) &= \int_{\R^d}(L_t)^Ne^{it \Phi(\xi)}\varphi d\xi = (-1)^N\int_{\R^d}e^{it \Phi(\xi)} (L_t^T)^N\varphi d\xi 
\end{align*}
We can have that 
\begin{align*}
|(L_t^T)^N\varphi(\xi)| \leq \frac{1}{(1+t|\nabla_{\xi}\Phi(\xi)|^2)^N}.
\end{align*}
when $|\nabla_{\xi}\Phi(\xi)|\leq c_0$. \\
For $\Phi(\xi)= x\cdot \xi \pm |\xi|$, combining \eqref{Kbdd1} for $N=\frac{d-1}{2}$ and \eqref{Kbdd2} $N'=d$, we obtain 
\begin{align*}
\bega
|K^{\pm}(t,tx)| \leq \frac{C}{t^{\frac{d-1}{2}}} + \int_{\mathcal{C}_x} \frac{C}{\big(1+t\big|x\pm\frac{\xi}{|\xi|}\big|^2\big)^d} d\xi
\enda
\end{align*}
where $\mathcal{C}_x:= \{ \xi \in \R^d ~|~ r\leq |\xi| \leq R,~ |x\pm\frac{\xi}{|\xi|}| \leq 1\}$. We then decompose $\xi$ into components parallel and perpendicular to $x$:
\begin{align*}
\bega
\xi = \xi_{||} + \xi_{\perp}, \qquad \xi_{||}:= \lw\la \xi, \frac{x}{|x|} \rw\ra \frac{x}{|x|} , \qquad \xi_{\perp}:= \xi- \lw\la \xi, \frac{x}{|x|} \rw\ra \frac{x}{|x|}.
\enda
\end{align*}
Note that $\xi_{||}$ and $\xi_{\perp}$ are 1-dimensional and $(d-1)$ dimensional vector, respectively. Since 
\begin{align*}
\bega
\bigg|x\pm\frac{\xi}{|\xi|}\bigg| = \bigg|x\pm\frac{\xi_{||}}{|\xi|} \pm\frac{\xi_{\perp}}{|\xi|} \bigg| \geq \bigg|\frac{\xi_{\perp}}{|\xi|}\bigg|,
\enda
\end{align*}
we have 
\begin{align*}
\bega
|K^{\pm}(t,tx)| \leq \frac{C}{t^{\frac{d-1}{2}}} + \int_{\mathcal{C}_x} \frac{C}{\big(1+t\big|\frac{\xi_{\perp}}{|\xi|}\big|^2\big)^d} d\xi_{||} d\xi_{\perp}.
\enda
\end{align*}
Since $r\leq|\xi|\leq R$, applying change of the variable $t^{\frac{1}{2}}\xi_{\perp} = \xi_{\perp}'$ with $t^{\frac{d-1}{2}}d\xi_{\perp} = d\xi_{\perp}'$ gives 
\begin{align*}
\bega
|K^{\pm}(t,tx)| \leq \frac{C}{t^{\frac{d-1}{2}}} + Ct^{-\frac{d-1}{2}} \leq \frac{C}{t^{\frac{d-1}{2}}}.
\enda
\end{align*}
\textcolor{red}{The constant $C$ depends on $r$ and $R$.} Thus we proved the claim \eqref{Kdecay}.

\end{proof}
\unhide

\StopNoTOC

\bibliographystyle{abbrv}

\def\cprime{$'$}

\end{document}